\documentclass{article}

\usepackage{rotating}
\usepackage{tikz-cd}
\usepackage[all]{xy}
\usepackage{amssymb}
\usepackage{amsmath}
\usepackage{amsthm}
\newtheorem{thm}{Theorem}
\newtheorem{prop}{Proposition}
\newtheorem{lem}{Lemma}
\newtheorem{rem}{Remark}
\newtheorem{cor}{Corollary}
\newtheorem{defi}{Definition}
\newtheorem{defiprop}{Definition-Proposition}

\def\Set{\mathop{\rm Set}\nolimits}
\def\Top{\mathop{\rm Top}\nolimits}
\def\RTop{\mathop{\rm RTop}\nolimits}
\def\RvTop{\mathop{\rm RvTop}\nolimits}
\def\Sch{\mathop{\rm Sch}\nolimits}
\def\AdSp{\mathop{\rm AdSp}\nolimits}

\def\CH{\mathop{\rm CH}\nolimits}

\def\alg{\mathop{\rm alg}\nolimits}
\def\hom{\mathop{\rm hom}\nolimits}

\def\codim{\mathop{\rm codim}\nolimits}
\def\dim{\mathop{\rm dim}\nolimits}

\def\Tr{\mathop{\rm Tr}\nolimits}
\def\Cor{\mathop{\rm Cor}\nolimits}

\def\SmVar{\mathop{\rm SmVar}\nolimits}
\def\Gr{\mathop{\rm Gr}\nolimits}

\def\PSmVar{\mathop{\rm PSmVar}\nolimits}
\def\Var{\mathop{\rm Var}\nolimits}
\def\Hom{\mathop{\rm Hom}\nolimits}
\def\Spec{\mathop{\rm Spec}\nolimits}

\def\supp{\mathop{\rm supp}\nolimits}
\def\QPVar{\mathop{\rm QPVar}\nolimits}
\def\AnSp{\mathop{\rm AnSp}\nolimits}
\def\CW{\mathop{\rm CW}\nolimits}
\def\Cw{\mathop{\rm Cw}\nolimits}
\def\PVar{\mathop{\rm PVar}\nolimits}

\def\Tot{\mathop{\rm Tot}\nolimits}
\def\sing{\mathop{\rm sing}\nolimits}

\def\Im{\mathop{\rm Im}\nolimits}

\def\Cone{\mathop{\rm Cone}\nolimits}
\def\ad{\mathop{\rm ad}\nolimits}

\def\nul{\mathop{\rm nul}\nolimits}
\def\log{\mathop{\rm log}\nolimits}
\def\Diff{\mathop{\rm Diff}\nolimits}
\def\An{\mathop{\rm An}\nolimits}
\def\Ho{\mathop{\rm Ho}\nolimits}
\def\PSh{\mathop{\rm PSh}\nolimits}
\def\AnSm{\mathop{\rm AnSm}\nolimits}
\def\Tr{\mathop{\rm Tr}\nolimits}

\def\pt{\mathop{\rm pt}\nolimits}
\def\DA{\mathop{\rm DA}\nolimits}

\def\Fun{\mathop{\rm Fun}\nolimits}

\def\Ab{\mathop{\rm Ab}\nolimits}
\def\Bti{\mathop{\rm Bti}\nolimits}
\def\TM{\mathop{\rm TM}\nolimits}
\def\Ouv{\mathop{\rm Ouv}\nolimits}
\def\coker{\mathop{\rm coker}\nolimits}
\def\Ext{\mathop{\rm Ext}\nolimits}
\def\CS{\mathop{\rm CS}\nolimits}

\def\Cat{\mathop{\rm Cat}\nolimits}
\def\RCat{\mathop{\rm RCat}\nolimits}
\def\TriCat{\mathop{\rm TriCat}\nolimits}
\def\an{\mathop{\rm an}\nolimits}
\def\holim{\mathop{\rm holim}\nolimits}

\def\Ring{\mathop{\rm Ring}\nolimits}
\def\cRing{\mathop{\rm cRing}\nolimits}
\def\TRing{\mathop{\rm TRing}\nolimits}
\def\cTRing{\mathop{\rm cTRing}\nolimits}
\def\Mod{\mathop{\rm Mod}\nolimits}

\def\Shv{\mathop{\rm Shv}\nolimits}
\def\Proj{\mathop{\rm Proj}\nolimits}
\def\Inj{\mathop{\rm Inj}\nolimits}

\def\nul{\mathop{\rm nul}\nolimits}

\def\Vect{\mathop{\rm Vect}\nolimits}
\def\Spv{\mathop{\rm Spv}\nolimits}
\def\Spa{\mathop{\rm Spa}\nolimits}
\def\Gal{\mathop{\rm Gal}\nolimits}
\def\AbCat{\mathop{\rm AbCat}\nolimits}

\def\dar[#1]{\ar@<2pt>[#1]\ar@<-2pt>[#1]}

\title{The De Rham, complex Hodge and $p$-adic Hodge realization functors 
on the derived category of relative motives over a field of characteristic zero}

\author{Johann Bouali}

\begin{document}

\maketitle

\begin{abstract}
We introduce the categories of geometric complex mixed Hodge modules on algebraic varieties 
over a subfield $k\subset\mathbb C$,
and for a prime number $p$, the categories of $p$-adic mixed Hodge modules on algebraic varieties
over a subfield $k\subset\mathbb C_p$.
We then give a complex Hodge realization functor on the derived category of relative motives over $k\subset\mathbb C$
and a $p$-adic Hodge realization functor on the derived category of relative motives over $k\subset\mathbb C_p$.
\end{abstract}

\tableofcontents

\section{Introduction}

\hskip 1cm Let $k$ be a field of characteristic zero and $S$ be a scheme of finite type over $k$. 
Let $DA_c(S)$ be the derived category of constructible motives over $S$. 
In a previous work (\cite{B4}) when $k={\mathbb C}$ is the field of complex numbers, we built a Hodge realization functor :
$$\mathcal F^{Hdg}_S :\DA_c(S)\rightarrow D(MHM(S))$$
where $D(MHM(S))$ is the derived category of mixed Hodge modules introduced by Morihiko Saito. 
This functor commute to the sixth operations formalism and define a 2-functor morphism on the category of schemes over $\mathbb C$.\\

\hskip 1cm In this work, we extend this realization functor 
to the general case of any field $k$ of characteristic zero embedded in $\mathbb C$ 
and we develop a $p$-adic analog of this realization. \\

\hskip 1cm The first step of the construction involves a rational version (i.e. over $k$) of the category of mixed Hodge modules. 
The key point of the construction is 
the existence of the Kashiwara-Malgrange $V$-filtration over $k$ for regular holonomic $\mathcal D$-modules 
(theorem \ref{HSk}) which is proved by induction on dimension using the complex case, 
and Saito's theorem on the strictness and stability of the direct image for proper morphisms in the complex case. 
The mixed hodge modules over $k$ are then the mixed Hodge modules whose regular holonomic sheaf is defined over $k$.
To develop a $p$-adic analog, we also introduce the full subcategory of geometric mixed Hodge modules.

\hskip 1cm For simplicity, we assume $S$ smooth and let $(M,F^\bullet M)$ be a filtered regular holonomic $\mathcal D_S$-module. 
We say $(M,F^\bullet M)$ is pure de Rham if it belongs to the full abelian category
generated by the successive higher direct images of the structural $\mathcal D_X$-module $(O_X,F_b)$ 
of $S$-schemes $X$ proper over $S$ and smooth over $k$, $F_b$ being the trivial filtration. \\
Let $(M, F^\bullet M, W^\bullet M)$ be a bi-filtered regular holonomic $\mathcal D_S$-module 
($F^\bullet M$ is the Hodge filtration and $W^\bullet M$ is the weight filtration). 
We say $(M,F,W)$ is de Rham if the associated graded module $Gr_W(M,F)$ is a pure de Rham module and if the weight filtration
is finite and satisfy an admissibility condition with respect to the Cartier divisors of $S$.
The de Rham modules over $S\in\Var(k)$ are introduced in definition \ref{DRMdef}.

\hskip 1cm A constructible sheaf over $S$ is given by a constructible sheaf $K$ over the analytic complex space $S_{\mathbb C}^{an}$ 
such that there exist a stratification $(S_i)$ of $S$ over $k$ 
such that $K_{|S_{i, {\mathbb C}}^{an}}$ is a $\mathbb Q$-local system. \\ 
We denote by $D_{c,k}(S_{\mathbb C}^{an})\subset D(S_{\mathbb C}^{an})$ the full subcategory
of the derived category of sheaves (or presheaves) on $S_{\mathbb C}^{an}$ whose cohomology sheaves are constructible over $S$
and by $P_k(S_{\mathbb C}^{an}):=P(S_{\mathbb C}^{an})\cap D_{c,k}(S_{\mathbb C}^{an})$ 
the full subcategory of perverse sheaves with constructible cohomology sheaves over $S$.

\hskip 1cm A geometric mixed Hodge module over $S$, assumed to be smooth for simplicity, 
is a triple $((M,F,W),(K,W),\alpha)$ where $(M,F,W)$ is a de Rham module over $S$, 
$(K,W)\in P_k(S_{\mathbb C}^{an})$ is a filtered perverse sheaf over $S$ 
and $\alpha$ is an isomorphism $\alpha:(K,W)\otimes {\mathbb C}\simeq DR(S)((M,W)^{an})$ 
compatible with the de Rham comparison theorem and where $DR(S)(M^{an})$ is the de Rham complex associated to $M^{an}$.
The geometric mixed Hodge module over $S\in\Var(k)$ are introduced in definition \ref{MHMgmkdef}. 

\hskip 1cm Let $D(MHM_{gm,k,\mathbb C}(S))$ be the derived category of 
the category of complexes of geometric mixed Hodge modules over $S$. 
This category can be defined for any scheme $S$ of finite type over $k$. We prove the following theorem :

\begin{thm}\label{thm1} 
Let $\Var(k)$ be the category of schemes of finite type over a subfield $k\subset\mathbb C$. Then :
\vskip 2mm
\begin{itemize}
\item The categories $D(MHM_{gm,k,\mathbb C}(S))$, for $S\in\Var(k)$, are endowed with the formalism of the sixth operation.
\item There exist a Hodge realization functor :
\vskip -2mm
$$\mathcal F^{Hdg}_S : \DA_c(S)\rightarrow D(MHM_{gm,k,\mathbb C}(S))$$
\vskip -1mm
compatible with the sixth operations formalism.
\end{itemize}
\end{thm}

\hskip 1cm Let $p$ be a prime number, ${\mathbb C}_p$ be the completion of an algebraic closure of $\mathbb Q_p$ 
and let $k\subset K\subset\mathbb C_p$ be a subfield of a $p$-adic field $K$. 
Let $S$ be a smooth scheme over $k$, 
$\mathbb B_{dR,S}$ be the sheaf of relative de Rham $p$-adic periods over the pro-étale site $S^{an}_{proet}$
of the $p$-adic analytic Huber space $S^{an}$ (introduced by Fontaine, Faltings and Scholze) 
and $C_{B_{dR}}(S)$  be the category of complexes of ${\mathbb B}_{dR,S}$-modules.

\hskip 1cm For any lisse ${\mathbb Q}_p$-sheaf $L$ over $S_{et}$, 
the sheaf $L\otimes {\mathbb B}_{dR,S}$ over $S^{an}_{proet}$ has a Poincaré resolution 
by the de Rham complex ${\cal O} {\mathbb B}_{dR,S}\otimes\Omega^{\bullet}_{S^{an}}$. 
We obtain a functor from the category of lisse ${\mathbb Q}_p$-sheaves over $S_{et}$ 
to the category of complexes of ${\mathbb B}_{dR,S}$-modules.
Using, a Beilinson's devissage by nearby and vanishing cycles functors, we extend this functor to perverse sheaves :
$$\mathbb B_{dr,S} : P_k(S_{et}) \rightarrow C_{B_{dR}}(S).$$

\hskip 1cm A geometric $p$-adic mixed Hodge module over $S$, assumed to be smooth for simplicity, is a triple $((M,F,W),(K,W),\alpha)$ 
where $(M,F,W)$ is a de Rham module over $S$, $(K,W)$ is a filtered perverse sheaf over $S_{et}$ and
$$\alpha :\mathbb B_{dr,S}(K,W)\simeq F^0DR(S)((O\mathbb B_{dR,S},F)\otimes_{O_S}(M,F,W)^{an})$$
is an isomorphism of complexes of $W$-filtered $({\mathbb B}_{dR,S},G)$-modules over $S^{an}_{proet}$, 
compatible with the $p$-adic de Rham comparison theorem (Faltings and Scholze)
and where $DR(S)$ is the de Rham complex associated to an analytic $D_S$-module and $G:=\Gal(\bar K/K)$ is the Galois group of $K$.
The geometric $p$-adic mixed Hodge module over $S\in\Var(k)$ are introduced in definition \ref{MHMgmkpdef}. \\

\hskip 1cm Let $D(MHM_{gm,k,\mathbb C_p}(S))$ 
be the derived category of the category of complexes of $p$-adic geometric Hodge modules over $S$. 
In fact $D(MHM_{gm,k,\mathbb C_p}(S))$  can be defined for any scheme $S$ of finite type over $k$. 
We prove the following $p$-adic version of theorem 1 :

\begin{thm}\label{thm2} 
Let $k\hookrightarrow K\hookrightarrow\mathbb C_p$ be a subfield of a $p$-adic field. Then :
\vskip 2mm

\begin{itemize}
\item The categories $D(MHM_{gm,k,\mathbb C_p}(S))$, for $S\in\Var(k)$, are endowed with the formalism of the sixth operation.
\vskip 2mm
\item There exist a $p$-adic Hodge realization functor :
\vskip -2mm
$$\mathcal F^{Hdg}_S : DA_c(S)\rightarrow D(MHM_{gm,k,\mathbb C_p}(S))$$
\vskip -1mm
compatible with the sixth operations formalism
\end{itemize}
\end{thm}

The proof of the first part of the theorems \ref{thm1} and \ref{thm2} is similar to the proof of the classical complex case, 
the crucial point is to show in the $p$-adic case that the isomorphism $\alpha$ is functorial. 
We are then reduced to prove that the functor $\mathbb B_{dr}$ commute to direct images in the proper case (theorem \ref{TfBdrthm}).
To do this, we use $p$-adic Hodge theory comparison theorems in the open case (\cite{Chinois}).

The proof of the second part of the theorems follows our strategy of \cite{B4} : \\ 
In the case $k\subset\mathbb C$, 
we first fully faithfully embed $D(MHM_{gm,k,\mathbb C}(S))$ into the fiber product of 
the derived category of bi-filtered regular holonomic $\mathcal D_S$-modules over $k$ 
and the derived category of filtered constructible $\mathbb Q$-sheaves with $k$-rational stratification, 
and construct the realization functor inside this big category (definition \ref{HodgeRealDAsing}), 
then we check that the image is contained in $DMHM_{gm,k,\mathbb C}(S)$ and commutes with the six operation (theorem \ref{main}). \\
In the case $k\subset K\subset\mathbb C_p$, 
we first fully faithfully embed $D(MHM_{gm,k,\mathbb C_p}(S))$ into the fiber product of 
the derived category of the category of bi-filtered regular holonomic $\mathcal D_S$-modules over $k$ 
and the derived category of filtered constructible $\mathbb Q_p$-etale sheaves with $k$-rational stratification 
and construct the realization functor inside this big category (definition \ref{HodgeRealDAsingp}), 
then we check that the image is contained in $D(MHM_{gm,k,\mathbb C_p}(S))$ and commutes with the six operation 
(theorem \ref{mainpadic}).

I am grateful to F.Mokrane for his help and support during the preparation of this work 
as well as J.Wildeshaus for the interest he has brought to this work.
I also thank A.C.Le Bras, P.Boyer and S.Morra for the help and interest they have brought to this work.

\section{Preliminaries and Notations}

\subsection{Notations}

\begin{itemize}

\item After fixing a universe, we denote by
\begin{itemize}
\item $\Set$ the category of sets, 
\item $\Top$ the category of topological spaces, 
\item $\Ring$ the category of rings and $\cRing\subset\Ring$ the full subscategory of commutative rings,
\item $\TRing$ the category of topological rings and 
$\cTRing\subset\TRing$ the full subscategory of commutative topological rings,
\item $\RTop$ the category of ringed spaces, 
\begin{itemize}
\item whose set of objects is
$\RTop:=\left\{(X,O_X), \; X\in\Top, \, O_X\in\PSh(X,\Ring)\right\}$
\item whose set of morphism is 
$\Hom((T,O_T),(S,O_S)):=\left\{((f:T\to S),(a_f:f^*O_S\to O_T))\right\}$
\end{itemize}
and by $ts:\RTop\to\Top$ the forgetfull functor.
\item $\RvTop$ the category of valued ringed spaces
\begin{itemize}
\item whose set of objects is
\begin{equation*}
\RvTop:=\left\{(X,O_X,(v_x,x\in X)), \; X\in\Top, \, O_X\in\PSh(X,\cTRing), \, v_x\in\Spv(O_{X,x})\right\}
\end{equation*}
where $O_X$ is a sheaf of complete topological commutative ring for a non archimedean semi-norm, 
$\Spv(-)$ denote the set of continous valuations of a topological commutative ring for a non archimedean semi-norm
\item whose set of morphism is 
\begin{eqnarray*}
\Hom((T,O_T,(v_x,x\in T)),(S,O_S,(v_x,x\in S))):= \\
\left\{((f:T\to S),(a_f:f^*O_S\to O_T)), \, a_f^*(v_{f(x)})=v_x \, \mbox{for all} \, x\in T\right\}.
\end{eqnarray*}
\end{itemize}

\item $\Cat$ the category of small categories which comes with the forgetful functor $o:\Cat\to\Fun(\Delta^1,\Set)$, 
where $\Fun(\Delta^1,\Set)$ is the category of simplicial sets, 
\item $\RCat$ the category of ringed topos
\begin{itemize}
\item whose set of objects is 
$\RCat:=\left\{(\mathcal X,O_X), \; \mathcal X\in\Cat, \, O_X\in\PSh(\mathcal X,\Ring)\right\}$,
\item whose set of morphism is 
$\Hom((\mathcal T,O_T),(\mathcal S,O_S)):=\left\{((f:\mathcal T\to \mathcal S),(a_f:f^*O_S\to O_T)),\right\}$
\end{itemize}
and by $tc:\RCat\to\Cat$ the forgetfull functor,
\item $\AbCat$ a category consisting of a small set of abelian categories,
\item $\TriCat$ a category constisting of a small set of triangulated categories.
\end{itemize}

\item Let $F:\mathcal C\to\mathcal C'$ be a functor with $\mathcal C,\mathcal C'\in\Cat$. 
For $X\in\mathcal C$, we denote by $F(X)\in\mathcal C'$ the image of $X$, 
and for $X,Y\in\mathcal C$,  we denote by $F^{X,Y}:\Hom(X,Y)\to\Hom(F(X),F(Y))$ the corresponding map.

\item For $\mathcal C\in\Cat$, we denote by $\mathcal C^{op}\in\Cat$ the opposite category whose set of object
is the one of $\mathcal C$ : $(\mathcal C^{op})^0=\mathcal C^0$, and whose morphisms are the morphisms of $\mathcal C$
with reversed arrows.

\item Let $\mathcal C\in\Cat$. For $S\in\mathcal C$, we denote by $\mathcal C/S$ the category 
\begin{itemize}
\item whose set of objects $(\mathcal C/S)^0=\left\{X/S=(X,h)\right\}$ consist of the morphisms $h:X\to S$ with $X\in\mathcal C$,
\item whose set of morphism $\Hom(X'/S,X/S)$ between $X'/S=(X',h'),X/S=(X,h)\in\mathcal C/S$ 
consits of the morphisms $(g:X'\to X)\in\Hom(X',X)$ such that $h\circ g=h'$. 
\end{itemize}
We have then, for $S\in\mathcal C$, the canonical forgetful functor 
\begin{eqnarray*}
r(S):\mathcal C/S\to\mathcal C, \; \; X/S\mapsto r(S)(X/S)=X, \; (g:X'/S\to X/S)\mapsto r(S)(g)=g
\end{eqnarray*}
and we denote again $r(S):\mathcal C\to\mathcal C/S$ the corresponding morphism of (pre)sites.
\begin{itemize}
\item Let $F:\mathcal C\to\mathcal C'$ be a functor with $\mathcal C,\mathcal C'\in\Cat$. 
Then for $S\in\mathcal C$, we have the canonical functor 
\begin{eqnarray*}
F_S:\mathcal C/S\to\mathcal C'/F(S), \; \; X/S\mapsto F(X/S)=F(X)/F(S), \\ (g:X'/S\to X/S)\mapsto (F(g):F(X')/F(S)\to F(X)/F(S))
\end{eqnarray*}
\item Let $\mathcal S\in\Cat$. Then, for a morphism $f:X'\to X$ with $X,X'\in\mathcal S$
we have the functor  
\begin{eqnarray*}
C(f):\mathcal S/X'\to\mathcal S/X, \; \; Y/X'=(Y,f_1)\mapsto C(f)(Y/X'):=(Y,f\circ f_1)\in\mathcal S/X, \\ 
(g:Y_1/X'\to Y_2/X')\mapsto (C(f)(g):=g:Y_1/X\to Y_2/X)
\end{eqnarray*}
\item Let $\mathcal S\in\Cat$ a category which admits fiber products. 
Then, for a morphism $f:X'\to X$ with $X,X'\in\mathcal S$, we have the pullback functor
\begin{eqnarray*}
P(f):\mathcal S/X\to\mathcal S/X', \; \; Y/X\mapsto P(f)(Y/X):=Y\times_X X'/X'\in\mathcal S/X', \\ 
(g:Y_1/X\to Y_2/X)\mapsto (P(f)(g):=(g\times I):Y_1\times_X X'\to Y_2\times_X X')
\end{eqnarray*} 
which is right adjoint to $C(f):\mathcal S/X'\to\mathcal S/X$,
and we denote again $P(f):\mathcal S/X'\to\mathcal S/X$ the corresponding morphism of (pre)sites.
\end{itemize}

\item Let $\mathcal C,\mathcal I\in\Cat$. Assume that $\mathcal C$ admits fiber products.
For $(S_{\bullet})\in\Fun(\mathcal I^{op},\mathcal C)$, 
we denote by $\mathcal C/(S_{\bullet})\in\Fun(\mathcal I,\Cat)$ the diagram of category given by
\begin{itemize}
\item for $I\in\mathcal I$, $\mathcal C/(S_{\bullet})(I):=\mathcal C/S_I$,
\item for $r_{IJ}:I\to J$, $\mathcal C/(S_{\bullet})(r_{IJ}):=P(r_{IJ}):\mathcal C/S_I\to\mathcal C/S_J$,
where we denoted again $r_{IJ}:S_J\to S_I$ the associated morphism in $\mathcal C$. 
\end{itemize}

\item Let $(F,G):\mathcal C\leftrightarrows\mathcal C'$ an adjonction between two categories.
\begin{itemize}
\item For $X\in C$ and $Y\in C'$, we consider the adjonction isomorphisms
\begin{itemize}
\item $I(F,G)(X,Y):\Hom(F(X),Y)\to\Hom(X,G(Y)), \; (u:F(X)\to Y)\mapsto (I(F,G)(X,Y)(u):X\to G(Y))$
\item $I(F,G)(X,Y):\Hom(X,G(Y))\to\Hom(F(X),Y), \; (v:X\to G(Y))\mapsto (I(F,G)(X,Y)(v):F(X)\to Y)$. 
\end{itemize}
\item For $X\in\mathcal C$, we denote by 
$\ad(F,G)(X):=I(F,G)(X,F(X))(I_{F(X)}):X\to G\circ F(X)$. 
\item For $Y\in\mathcal C'$ we denote also by 
$\ad(F,G)(Y):=I(F,G)(G(Y),Y)(I_{G(Y)}):F\circ G(Y)\to Y$.
\end{itemize}
Hence, 
\begin{itemize}
\item for $u:F(X)\to Y$ a morphism with $X\in C$ and $Y\in C'$, we have 
$I(F,G)(X,Y)(u)=G(u)\circ\ad(F,G)(X)$,
\item for $v:X\to G(Y)$ a morphism with $X\in C$ and $Y\in C'$, we have 
$I(F,G)(X,Y)(v)=\ad(F,G)(Y)\circ F(v)$.
\end{itemize}

\item Let $\mathcal C$ a category. 
\begin{itemize}
\item We denote by $(\mathcal C,F)$ the category of filtered objects :
$(X,F)\in(\mathcal C,F)$ is a sequence $(F^{\bullet}X)_{\bullet\in\mathbb Z}$ indexed by $\mathbb Z$
with value in $\mathcal C$ together with monomorphisms $a_p:F^pX\hookrightarrow F^{p-1}X\hookrightarrow X$.
\item We denote by $(\mathcal C,F,W)$ the category of bifiltered objects : $(X,F,W)\in(\mathcal C,F,W)$ 
is a sequence $(W^{\bullet}F^{\bullet}X)_{\bullet,\bullet}\in\mathbb Z^2$ indexed by $\mathbb Z^2$ with value in $\mathcal C$ 
together with monomorphisms $W^qF^pX\hookrightarrow F^{p-1}X$, $W^qF^pX\hookrightarrow W^{q-1}F^pX$.
\end{itemize}

\item Let $\mathcal A$ an additive category.
\begin{itemize}
\item We denote by $C(\mathcal A):=\Fun(\mathbb Z,\mathcal A)$ 
the category of (unbounded) complexes with value in $\mathcal A$, 
where we have denoted $\mathbb Z$ the category whose set of objects is $\mathbb Z$,
and whose set of morphism between $m,n\in\mathbb Z$ consists of one element (identity) if $n=m$,
of one elemement if $n=m+1$ and is $\emptyset$ in the other cases. 
\item We have the full subcategories $C^b(\mathcal A)$, $C^-(\mathcal A)$, $C^+(\mathcal A)$ of 
$C(\mathcal A)$ consisting of bounded, resp. bounded above, resp. bounded below complexes.
\item We denote by $K(\mathcal A):=\Ho(C(\mathcal A))$ the homotopy category of $C(\mathcal A)$
whose morphisms are equivalent homotopic classes of morphism and by
$Ho:C(\mathcal A)\to K(\mathcal A)$ the full homotopy functor. 
The category $K(\mathcal A)$ is in the standard way a triangulated category.
\end{itemize}

\item Let $\mathcal A$ an additive category.
\begin{itemize}
\item We denote by $C_{fil}(\mathcal A)\subset(C(\mathcal A),F)=C(\mathcal A,F)$ the full
additive subcategory of filtered complexes of $\mathcal A$ such that the filtration is biregular : 
for $(A^{\bullet},F)\in(C(\mathcal A),F)$, we say that $F$ is biregular if $F^{\bullet}A^r$ is finite for all $r\in\mathbb Z$.
\item We denote by $C_{2fil}(\mathcal A)\subset(C(\mathcal A),F,W)=C(\mathcal A,F,W)$ the full
subcategory of bifiltered complexes of $\mathcal A$ such that the filtration is biregular.
\item For $A^{\bullet}\in C(\mathcal A)$, we denote by $(A^{\bullet},F_b)\in(C(\mathcal A),F)$ the complex endowed
with the trivial filtration (filtration bete) : $F^pA^n=0$ if $p\geq n+1$ and $F^pA^n=A^n$ if $p\leq n$. 
Obviously, a morphism $\phi:A^{\bullet}\to B^{\bullet}$, with $A^{\bullet},B^{\bullet}\in C(\mathcal A)$
induces a morphism $\phi:(A^{\bullet},F_b)\to (B^{\bullet},F_b)$.
\item For $(A^{\bullet},F)\in C(\mathcal A,F)$, we denote by $(A^{\bullet},F(r))\in C(\mathcal A,F)$ the filtered
complex where the filtration is given by $F^p(A^{\bullet},F(r)):=F^{p+r}(A^{\bullet},F)$.
\end{itemize}

\item Let $\mathcal A$ be an abelian category. 
Then the additive category $(\mathcal A,F)$ is an exact category which admits kernel and cokernel 
(but is NOT an abelian category). 
A morphism $\phi:(M,F)\to(N,F)$ with $(M,F)\in(\mathcal A,F)$ is strict if 
the inclusion $\phi(F^nM)\subset F^nN\cap\Im(\phi)$ is an equality, i.e. if $\phi(F^nM)=F^nN\cap\Im(\phi)$.

\item Let $\mathcal A$ be an abelian category. 
\begin{itemize}
\item For $(A^{\bullet},F)\in C(\mathcal A,F)$, considering $a_p:F^pA^{\bullet}\hookrightarrow A^{\bullet}$
the structural monomorphism of of the filtration, we denote by, for $n\in\mathbb N$, 
\begin{equation*}
H^n(A^{\bullet},F)\in(\mathcal A,F), \; 
F^pH^n(A^{\bullet},F):=\Im(H^n(a_p):H^n(F^pA^{\bullet})\to H^n(A^{\bullet}))\subset H^n(A^{\bullet})
\end{equation*}
the filtration induced on the cohomology objects of the complex.
In the case $(A^{\bullet},F)\in C_{fil}(\mathcal A)$, 
the spectral sequence $E^{p,q}_r(A^{\bullet},F)$ associated to $(A^{\bullet},F)$ 
converge to $\Gr_F^pH^{p+q}(A^{\bullet},F)$, that is for all $p,q\in\mathbb Z$, there exist $r_{p+q}\in\mathbb N$,
such that $E^{p,q}_s(A^{\bullet},F)=\Gr_F^pH^{p+q}(A^{\bullet},F)$ for all $s\leq r_{p+q}$.
\item A morphism $m:(A^{\bullet},F)\to(B^{\bullet},F)$ with $(A^{\bullet},F),(B^{\bullet},F)\in C(\mathcal A,F)$
is said to be a filtered quasi-isomorphism if for all $n,p\in\mathbb Z$, 
\begin{equation*}
H^n\Gr^p_F(m):H^n(\Gr^p_FA^{\bullet})\xrightarrow{\sim}H^n(\Gr^p_FB^{\bullet})
\end{equation*}
is an isomorphism in $\mathcal A$. 
Consider a commutative diagram in $C(\mathcal A,F)$
\begin{equation*}
\xymatrix{(A^{\bullet},F)\ar[r]^{m}\ar[d]_{\phi} & (B^{\bullet},F)\ar[r]^{i_2}\ar[d]_{\psi} & 
\Cone(m)=((A^{\bullet},F)[1]\oplus(B^{\bullet},F),d,d'-m)\ar[r]^{p_1}\ar[d]^{(\phi[1],\psi)} & 
(A^{\bullet},F)[1]\ar[d]^{\phi[1]} \\
(A^{'\bullet},F)\ar[r]^{m'} & (B^{'\bullet},F)\ar[r]^{i_2} & 
\Cone(m')=((A^{'\bullet},F)[1]\oplus(B^{'\bullet},F),d,d'-m')\ar[r]^{p_1} & (A^{'\bullet},F)[1]}
\end{equation*}
If $\phi$ and $\psi$ are filtered quasi-isomorphisms, then $(\phi[1],\psi)$ is an filtered quasi-isomorphism.
That is, the filtered quasi-isomorphism satisfy the 2 of 3 property for canonical triangles.
\end{itemize}

\item Let $\mathcal A$ be an abelian category.
\begin{itemize}
\item We denote by $D(\mathcal A)$ the localization of $K(\mathcal A)$
with respect to the quasi-isomorphisms and by $D:K(\mathcal A)\to D(\mathcal A)$ the localization functor. 
The category $D(\mathcal A)$ is a triangulated category in the unique way such that $D$ a triangulated functor.
\item We denote by $D_{fil}(\mathcal A)$ the localization of $K_{fil}(\mathcal A)$
with respect to the filtered quasi-isomorphisms and by 
$D:K_{fil}(\mathcal A)\to D_{fil}(\mathcal A)$ the localization functor.
\end{itemize}

\item Let $\mathcal A$ be an abelian category.
We denote by $\Inj(A)\subset A$ the full subcategory of injective objects,
and by $\Proj(A)\subset A$ the full subcategory of projective objects.

\item For $\mathcal S\in\Cat$ a small category, we denote by 
\begin{itemize}
\item $\PSh(\mathcal S):=\PSh(\mathcal S,\Ab):=\Fun(\mathcal S,\Ab)$ the category of presheaves on $\mathcal S$,
i.e. the category of presheaves of abelian groups on $\mathcal S$,
\item $K(\mathcal S):=K(\PSh(\mathcal S))=\Ho(C(\mathcal S))$
In particular, we have the full homotopy functor $Ho:C(\mathcal S)\to K(\mathcal S)$,
\item $C_{(2)fil}(\mathcal S):=C_{(2)fil}(\PSh(\mathcal S))\subset C(\PSh(\mathcal S),F,W)$
the big abelian category of (bi)filtered complexes of presheaves on $\mathcal S$ with value in abelian groups
such that the filtration is biregular, and $\PSh_{(2)fil}(\mathcal S)=(\PSh(\mathcal S),F,W)$,
\item $K_{fil}(\mathcal S):=K_{fil}(\PSh(\mathcal S))=\Ho(C_{fil}(\mathcal S))$,
\item $K_{fil,r}(\mathcal S):=K_{fil,r}(\PSh(\mathcal S))=\Ho_r(C_{fil}(\mathcal S))$,
$K_{fil,\infty}(\mathcal S):=K_{fil,\infty}(\PSh(\mathcal S))=\Ho_{\infty}(C_{fil}(\mathcal S))$.
\end{itemize}
For $f:\mathcal T\to\mathcal S$ a morphism a presite with $\mathcal T,\mathcal S\in\Cat$,
given by the functor $P(f):\mathcal S\to\mathcal T$,
we will consider the adjonctions given by the direct and inverse image functors : 
\begin{itemize}
\item $(f^*,f_*)=(f^{-1},f_*):\PSh(\mathcal S)\leftrightarrows\PSh(\mathcal T)$,
which induces $(f^*,f_*):C(\mathcal S)\leftrightarrows C(\mathcal T)$, we denote, 
for $F\in C(\mathcal S)$ and $G\in C(\mathcal T)$ by
\begin{equation*}
\ad(f^*,f_*)(F):F\to f_*f^*F \; , \; \ad(f^*,f_*)(G):f^*f_*G\to G
\end{equation*}
the adjonction maps,
\item $(f_*,f^{\bot}):\PSh(\mathcal T)\leftrightarrows\PSh(\mathcal S)$,
which induces $(f_*,f^{\bot}):C(\mathcal T)\leftrightarrows C(\mathcal S)$, we denote
for $F\in C(\mathcal S)$ and $G\in C(\mathcal T)$ by
\begin{equation*}
\ad(f_*,f^{\bot})(F):G\to f^{\bot}f_* G \; , \; \ad(f_*,f^{\bot})(G):f_*f^{\bot}F\to F
\end{equation*}
the adjonction maps. 
\end{itemize}

\item For $(\mathcal S,O_S)\in\RCat$ a ringed topos, we denote by 
\begin{itemize}
\item $\PSh_{O_S}(\mathcal S)$ the category of presheaves of $O_S$ modules on $\mathcal S$, 
whose objects are $\PSh_{O_S}(\mathcal S)^0:=\left\{(M,m),M\in\PSh(\mathcal S),m:M\otimes O_S\to M\right\}$,
together with the forgetful functor $o:\PSh(\mathcal S)\to \PSh_{O_S}(\mathcal S)$,
\item $C_{O_S}(\mathcal S)=C(\PSh_{O_S}(\mathcal S))$ 
the big abelian category of complexes of presheaves of $O_S$ modules on $\mathcal S$,
\item $K_{O_S}(\mathcal S):=K(\PSh_{O_S}(\mathcal S))=\Ho(C_{O_S}(\mathcal S))$,
in particular, we have the full homotopy functor $Ho:C_{O_S}(\mathcal S)\to K_{O_S}(\mathcal S)$,
\item $C_{O_S(2)fil}(\mathcal S):=C_{(2)fil}(\PSh_{O_S}(\mathcal S))\subset C(\PSh_{O_S}(\mathcal S),F,W)$,
the big abelian category of (bi)filtered complexes of presheaves of $O_S$ modules on $\mathcal S$ such that the filtration is biregular
and $\PSh_{O_S(2)fil}(\mathcal S)=(\PSh_{O_S}(\mathcal S),F,W)$,
\item $K_{O_Sfil}(\mathcal S):=K_{fil}(\PSh_{O_S}(\mathcal S))=\Ho(C_{O_Sfil}(\mathcal S))$,
\item $K_{O_Sfil,r}(\mathcal S):=K_{fil,r}(\PSh_{O_S}(\mathcal S))=\Ho_r(C_{O_Sfil}(\mathcal S))$,
$K_{O_Sfil,\infty}(\mathcal S):=K_{fil,\infty}(\PSh_{O_S}(\mathcal S))=\Ho_{\infty}(C_{O_Sfil}(\mathcal S))$.
\end{itemize}

\item For $\mathcal S\in\Cat$ a small category and $n\in\mathbb N$, 
$\PSh_{\mathbb Z/n\mathbb Z}(\mathcal S)\subset\PSh(\mathcal S)$ is the full subcategory of $n$-torsion presheaves. 
The functor 
\begin{equation*}
(-)\otimes\mathbb Z/n\mathbb Z:\PSh(\mathcal S)\to\PSh_{\mathbb Z/n\mathbb Z}(\mathcal S), 
F\mapsto F\otimes\mathbb Z/n\mathbb Z
\end{equation*}
is right exact and its restriction the full subcategory $\PSh(\mathcal S)_L\subset\PSh(\mathcal S)$ 
of torsion free presheaves is exact.
For $\mathcal S\in\Cat$ and $p\in\mathbb N$ a prime number,
\begin{equation*}
\PSh_{\mathbb Z_p}(\mathcal S)\subset\PSh(\mathbb N\times\mathcal S)=\PSh(\mathcal S,\Fun(\mathbb N,\Ab)) 
\end{equation*}
is category whose objects are $(F_l)_{l\in\mathbb N}$ with $F_l\in\PSh_{\mathbb Z/p^l\mathbb Z}(\mathcal S)$
such that $F_l\to F_{l+1}/p^lF_{l+1}$ is an isomorphism.
We then have 
\begin{equation*}
C_{\mathbb Z_p}(\mathcal S):=C(\PSh_{\mathbb Z_p}(\mathcal S))\subset 
C(\mathbb N\times\mathcal S)=\PSh(\mathcal S,\Fun(\mathbb N,C(\mathbb Z))).
\end{equation*}
We get the functor 
\begin{equation*}
(-)\otimes\mathbb Z_p:C(\mathcal S)\to C_{\mathbb Z_p}(\mathcal S), 
F\mapsto (F\otimes\mathbb Z/p^l\mathbb Z)_{l\in\mathbb N}
\end{equation*}
which is right exact.
For $\mathcal S\in\Cat$ a site with topology $\tau$, we have the localization 
\begin{equation*}
D_{\mathbb Z_p}(\mathcal S):=\Ho_{\tau}C(\PSh_{\mathbb Z_p}(\mathcal S))
\end{equation*}
of $\tau$ local equivalences of $C(\PSh_{\mathbb Z_p}(\mathcal S))\subset\PSh(\mathcal S,\Fun(\mathbb N,C(\mathbb Z)))$.

\item For $\mathcal S_{\bullet}\in\Fun(\mathcal I,\Cat)$ a diagram of (pre)sites, with $\mathcal I\in\Cat$ a small category, 
we denote by $S_{\bullet}:=\Gamma\mathcal S_{\bullet}\in\Cat$ the associated diagram category
\begin{itemize}
\item whose objects are $\Gamma\mathcal S_{\bullet}^0:=\left\{(X_I,u_{IJ})_{I\in\mathcal I}\right\}$, 
with $X_I\in\mathcal S_I$, and for $r_{IJ}:I\to J$ with $I,J\in\mathcal I$,  
$u_{IJ}:X_J\to r_{IJ}(X_I)$ are morphism in $\mathcal S_J$
noting again $r_{IJ}:\mathcal S_I\to\mathcal S_J$ the associated functor,
\item whose morphism are $m=(m_I):(X_I,u_{IJ})\to(X'_I,v_{IJ})$ satisfying 
$v_{IJ}\circ m_I =r_{IJ}(m_J)\circ u_{IJ}$ in $\mathcal S_J$, 
\end{itemize}
We have then 
$\PSh(\mathcal S_{\bullet})=\PSh(\Gamma\mathcal S_{\bullet})$ the category of presheaves on $\mathcal S_{\bullet}$,
\begin{itemize}
\item whose objects are $\PSh(\mathcal S_{\bullet})^0:=\left\{(F_I,u_{IJ})_{I\in\mathcal I}\right\}$, 
with $F_I\in\PSh(\mathcal S_I)$, and for $r_{IJ}:I\to J$ with $I,J\in\mathcal I$, 
$u_{IJ}:F_I\to r_{IJ*}F_J$ are morphism in $\PSh(\mathcal S_I)$, 
noting again $r_{IJ}:\mathcal S_J\to\mathcal S_I$ the associated morphism of presite,
\item whose morphism are $m=(m_I):(F_I,u_{IJ})\to(G_I,v_{IJ})$ satisfying 
$v_{IJ}\circ m_I =r_{IJ*}m_J\circ u_{IJ}$ in $\PSh(\mathcal S_I)$, 
\end{itemize}
Let $\mathcal I,\mathcal I'\in\Cat$ be small categories.
Let $(f_{\bullet},s):\mathcal T_{\bullet}\to\mathcal S_{\bullet}$ a morphism a diagrams of (pre)site with 
$\mathcal T_{\bullet}\in\Fun(\mathcal I,\Cat),\mathcal S_{\bullet}\in\Fun(\mathcal I',\Cat)$,
which is by definition given by a functor $s:\mathcal I\to \mathcal I'$ and morphism of functor 
$P(f_{\bullet}):\mathcal S_{s(\bullet)}:=\mathcal S_{\bullet}\circ s\to\mathcal T_{\bullet}$.
Here, we denote for short, $\mathcal S_{s(\bullet)}:=\mathcal S_{\bullet}\circ s\in\Fun(\mathcal I,\Cat)$.
We have then, for $r_{IJ}:I\to J$ a morphism, with $I,J\in\mathcal I$, a commutative diagram in $\Cat$ 
\begin{equation*}
D_{fIJ}:=\xymatrix{\mathcal S_{s(J)}\ar[r]^{r^s_{IJ}} & \mathcal S_{s(I)} \\
\mathcal T_J\ar[r]^{r^t_{IJ}}\ar[u]^{f_J} & \mathcal T_I\ar[u]^{f_I}}.
\end{equation*}
In particular we get the adjonction given by the direct and inverse image functors :
\begin{eqnarray*}
((f_{\bullet},s)^*,(f_{\bullet},s)^*)=((f_{\bullet},s)^{-1},(f_{\bullet},s)_*):
\PSh(\mathcal S_{s(\bullet)})\leftrightarrows\PSh(\mathcal T_{\bullet}), \\
F=(F_I,u_{IJ})\mapsto (f_{\bullet},s)^*F:=(f_I^*F_I,T(D_{fIJ})(F_J)\circ f_I^*u_{IJ}), \\ 
G=(G_I,v_{IJ})\mapsto (f_{\bullet},s)_*G:=(f_{I*}G_I,f_{I*}v_{IJ}).
\end{eqnarray*} 

\item Let $\mathcal I\in\Cat$ a small category.
For $(\mathcal S_{\bullet},O_{S_{\bullet}})\in\Fun(\mathcal I,\RCat)$ a diagram of ringed topos, we denote by 
\begin{equation*}
(\mathcal S_{\bullet},O_{S_{\bullet}}):=(\Gamma\mathcal S_{\bullet},O_{\Gamma\mathcal S_{\bullet}})\in\RCat.
\end{equation*}
We have then
$\PSh_{O_{S_{\bullet}}}(\mathcal S_{\bullet})=\PSh_{O_{\Gamma\mathcal S_{\bullet}}}(\Gamma\mathcal S_{\bullet})$ 
the category of presheaves of modules on $(\mathcal S_{\bullet},O_{S_{\bullet}})$,
\begin{itemize}
\item whose objects are $\PSh_{O_{S_{\bullet}}}(\mathcal S_{\bullet})^0:=\left\{(F_I,u_{IJ})_{I\in\mathcal I}\right\}$, 
with $F_I\in\PSh_{O_{S_I}}(\mathcal S_I)$, and for $r_{IJ}:I\to J$ with $I,J\in\mathcal I$, 
$u_{IJ}:F_I\to r_{IJ*}F_J$ are morphism in $\PSh_{O_{S_I}}(\mathcal S_I)$, 
noting again $r_{IJ}:\mathcal S_J\to\mathcal S_I$ the associated morphism of presite,
\item whose morphism are $m=(m_I):(F_I,u_{IJ})\to(G_I,v_{IJ})$ satisfying 
$v_{IJ}\circ m_I =r_{IJ*}m_J\circ u_{IJ}$ in $\PSh_{O_{S_I}}(\mathcal S_I)$, 
\end{itemize}

\item For $A\in\Ring$, $\dim_K(A)$ denote the Krull dimension of $A$.
For $\sigma:A\to B$ a morphism with $A,B\in\cRing$, we have the extention of scalar functor
\begin{eqnarray*}
\otimes_AB:(-)\otimes_AB:\Mod(A)\to\Mod(B), \; M\mapsto M\otimes_AB \\ 
(m:M'\to M)\mapsto (m_B:=m\otimes I:M'\otimes_AB\to M\otimes_AB).
\end{eqnarray*}
which is left ajoint to the restriction of scalar 
\begin{eqnarray*}
Res_{A/B}:\Mod(B)\to\Mod(A), \; M=(M,a_M)\mapsto M=(M,a_M\circ\sigma), \; (m:M'\to M)\mapsto (m:M'\to M)
\end{eqnarray*}
The adjonction maps are 
\begin{itemize}
\item for $M\in\Mod(A)$, the canonical map in $\Mod(A)$
\begin{equation*}
n_{A/B}(M):M\to M\otimes_A B, \; n_{A/B}(M)(m):=m\otimes 1, 
\end{equation*}
\item for $M\in\Mod(B)$, 
\begin{equation*}
I\times\Delta_B:M\otimes_AB=M\otimes_BB\otimes_AB
\end{equation*}
in $\Mod(B)$, where $\Delta_B:B\otimes_AB\to B$ is given by for $x,y\in B$, $\Delta_B(x,y)=x-y$.
\end{itemize}
Let $\sigma:A\to B$ a morphism with $A,B\in\cRing$. A module $M\in\Mod(B)$ is said to be defined over $A$
if there exist a module $M_0\in\Mod(A)$ and an isomorphism $M\simeq M_0\otimes_AB$ in $\Mod(B)$.
A module $M\in\Mod(B)$ is defined over $A$ if and only if
there exist a presentation of $M$, that is an exact sequence in $\Mod(B)$,
$B^{\oplus^J}\xrightarrow{\phi}B^{\oplus^I}\to M\to 0$,
such that $\phi\circ\sigma(A^{\oplus^J})\subset A^{\oplus^I}$. 

\item For $f=(f,a_f):(\mathcal T,O_T)\to(\mathcal S,O_S)$ a morphism of ringed topos 
with $(\mathcal S,O_S),(\mathcal T,O_T)\in\RCat$, $a_f:f^*O_S\to O_T$, we have the pull-back of presheaves of modules
\begin{eqnarray*}
f^{*mod}:\PSh_{O_S}(\mathcal S)\to\PSh_{O_T}(\mathcal T), \; M\mapsto f^{*mod}M:=f^*M\otimes_{f^*O_S}O_T \\ 
(m:M'\to M)\mapsto (f^{*mod}M:=f^*m\otimes I:f^{*mod}M'\to f^{*mod}M).
\end{eqnarray*}
which is left ajoint to 
\begin{eqnarray*}
f_*:\PSh_{O_T}(\mathcal T)\to\PSh_{O_S}(\mathcal S), \; M=(M,a_M)\mapsto f_*M=(f_*M,a_M\circ a_f), \\
(m:M'\to M)\mapsto (f_*m:f_*M'\to f_*M)
\end{eqnarray*}
The adjonction maps are 
\begin{itemize}
\item for $M\in\PSh_{O_S}(\mathcal S)$, the canonical map in $\PSh_{O_S}(\mathcal S)$
\begin{eqnarray*}
\ad(f^{*mod},f_*)(M):=n_{f^*O_S/O_T}(M):M\to f_*f^{*mod}M=f_*f^*M\otimes_{f_*f^*O_S}f_*O_T, \\ 
n_{f^*O_S/O_T}(M)(m):=\ad(f_*,f^*)(M)(m)\otimes 1, 
\end{eqnarray*}
\item for $M\in\PSh_{O_T}(\mathcal S)$, the canonical map 
\begin{equation*}
\ad(f^{*mod},f_*)(M):=I\times\Delta_{O_T}:f^{*mod}f_*M=f^*f_*M\otimes_{O_T}O_T\otimes_{f^*O_S}O_T\to M,  
\end{equation*}
where $\Delta_{O_T}:O_T\otimes_{f^*O_S}O_T\to O_T$ is given by for 
$x,y\in\Gamma(T,O_T)$, $T\in\mathcal T$, $\Delta_{O_T}(x,y)=x-y$.
\end{itemize}
Let $f=(f,a_f):(\mathcal T,O_T)\to(\mathcal S,O_S)$ a morphism of ringed topos 
with $(\mathcal S,O_S),(\mathcal T,O_T)\in\RCat$, $a_f:f^*O_S\to O_T$
A presheaf $M\in\PSh_{O_T}(\mathcal T)$ is said to be defined over $(\mathcal S,O_S)$
if there exist a  $M_0\in\in\PSh_{O_S}(\mathcal S)$ such that $M\simeq f^{*mod}M_0$ in $\PSh_{O_T}(\mathcal T)$.
For $M\in\PSh_{O_T}(\mathcal T)$ quasi-coherent, $M$ is locally defined over $(\mathcal S,O_S)$
if and only if there exists locally a presentation of $M$, 
that is an exact sequence in $\PSh_{O_T}(\mathcal T')$, $\mathcal T'\subset T$,
$O_T^{\oplus^J}\xrightarrow{\phi}O_T^{\oplus^I}\to M\to 0$,
such that $\phi\circ a_f(O_S^{\oplus^J})\subset O_S^{\oplus^I}$. 

\item For $X\in\Top$, we denote by $\dim_F(X)$ its Krull dimension and $\dim_L(X)$ its Lebegue dimension.
Note that if $X$ is Hausdorf $\dim_F(X)=0$ and if $X$ is everywhere not Hausdorf $\dim_L(X)=0$.
For $X\in\Top$ and $x\in X$, we denote by $\dim_{F,x}(X):=inf_{x\in U}\dim_F(U)$ its Krull dimension at $x$ 
and $\dim_{L,x}(X):=inf_{x\in U}\dim_L(U)$ its Lebegue dimension at $x$.

\item Denote by $\Sch\subset\RTop$ the full subcategory of schemes. 
For $X\in\Sch$, $\dim(X):=\dim_F(X)$. For $X=\Spec A\in\Sch$ an affine scheme, $\dim(X)=\dim_K(A)$.
For $X\in\Sch$ and $x\in X$, $\dim_x(X):=\dim_{F,x}(X)=\dim(O_{X,x})$.
A morphism $h:U\to S$ with $U,S\in\Sch$ is said to be smooth if it is flat with smooth fibers geometric fibers.
A morphism $r:U\to X$ with $U,X\in\Sch$ is said to be etale if it is non ramified and flat.
In particular an etale morphism $r:U\to X$ with $U,X\in\Sch$ 
is smooth and quasi-finite (i.e. the fibers are either the empty set or a finite subset of $X$)
For $X\in\Sch$, we denote by 
\begin{itemize}
\item $\Sch^{ft}/X\subset\Sch/X$ the full subcategory consisting of objects
$X'/X=(X',f)\in\Sch/X$ such that $f:X'\to X$ is an morphism of finite type
\item $X^{et}\subset\Sch^{ft}/X$ the full subcategory consisting of objects
$U/X=(X,h)\in\Sch/X$ such that $h:U\to X$ is an etale morphism.
\item $X^{sm}\subset\Sch^{ft}/X$ the full subcategory consisting of objects
$U/X=(X,h)\in\Sch/X$ such that $h:U\to X$ is a smooth morphism.
\end{itemize}
For a field $k$, we consider $\Sch/k:=\Sch/\Spec k$ the category of schemes over $\Spec k$. We then denote by
\begin{itemize}
\item $\Var(k)=\Sch^{ft}/k\subset\Sch/k$ the full subcategory consisting of algebraic varieties over $k$, 
i.e. schemes of finite type over $k$,
\item $\PVar(k)\subset\QPVar(k)\subset\Var(k)$ the full subcategories consisting of quasi-projective varieties and projective varieties respectively, 
\item $\PSmVar(k)\subset\SmVar(k)\subset\Var(k)$ the full subcategories consisting of smooth varieties and smooth projective varieties respectively.
\end{itemize}
For a morphism of field $\sigma:k\hookrightarrow K$, we have the extention of scalar functor
\begin{eqnarray*}
\otimes_kK:\Sch/k\to\Sch/K, \; X\mapsto X_K:=X_{K,\sigma}:=X\otimes_kK, \; (f:X'\to X)\mapsto (f_K:=f\otimes I:X'_K\to X_K).
\end{eqnarray*}
which is left ajoint to the restriction of scalar 
\begin{eqnarray*}
Res_{k/K}:\Sch/K\to\Sch/k, \; X=(X,a_X)\mapsto X=(X,\sigma\circ a_X), \; (f:X'\to X)\mapsto (f:X'\to X)
\end{eqnarray*}
The adjonction maps are 
\begin{itemize}
\item for $X\in\Sch/k$, the projection $\pi_{k/K}(X):X_K\to X$ in $\Sch/k$,
for $X=\cup_iX_i$ an affine open cover with $X_i=\Spec(A_i)$ we have by definition $\pi_{k/K}(X_i)=n_{k/K}(A_i)$,
\item for $X\in\Sch/K$, $I\times\Delta_K:X\hookrightarrow X_K=X\times_KK\otimes_kK$ in $\Sch/K$,
where $\Delta_K:K\otimes_kK\to K$ is the diagonal which is given by for $x,y\in K$, $\Delta_K(x,y)=x-y$.
\end{itemize}
The extention of scalar functor restrict to a functor
\begin{eqnarray*}
\otimes_kK:\Var(k)\to\Var(K), \; X\mapsto X_K:=X_{K,\sigma}:=X\otimes_kK, \; (f:X'\to X)\mapsto (f_K:=f\otimes I:X'_K\to X_K).
\end{eqnarray*}
and for $X\in\Var(k)$ we have $\pi_{k/K}(X):X_K\to X$ the projection in $\Sch/k$.
An algebraic variety $X\in\Var(K)$ is said to be defined over $k$ if there exists $X_0\in\Var(k)$
such that $X\simeq X_0\otimes_kK$ in $\Var(K)$. 
By definition, 
\begin{itemize}
\item for $X=\Spec(A)\in\Var(K)$ an affine variety, $X$ is defined over $K$ if $A\in\Mod(K)$ is defined over $k$, 
that is if $A=K[x_1,\ldots,x_N]/I$ is a presentation of $A$, $I=<f_1,\cdots f_r>\subset K[x_1,\ldots,x_N]$ 
with $f_1,\ldots,f_r\in k[x_1,\cdots,x_N]$ is generated by elements over $k$.
\item for $X=\Proj(B)\in\PVar(K)$ an projective variety, $X$ is defined over $K$ if $B\in\Mod(K)$ is defined over $k$, 
that is if $B=K[x_0,\ldots,x_N]/I$ is a presentation of $B$ with $I$ generated by homogeneous elements, 
$I=<f_1,\cdots f_r>\subset K[x_0,\ldots,x_N]$ with $f_1,\ldots,f_r\in k[x_0,\cdots,x_N]$ homogeneous.
\end{itemize}
For $X=(X,a_X)\in\Var(k)$, we have $\Sch^{ft}/X=\Var(k)/X$ since for $f:X'\to X$ a morphism of schemes of finite type,
$(X',a_X\circ f)\in\Var(k)$ is the unique structure of variety over $k$ of $X'\in\Sch$ such that $f$ becomes a morphism 
of algebraic varieties over $k$, in particular we have
\begin{itemize}
\item $X^{et}\subset\Sch^{ft}/X=\Var(k)/X$,
\item $X^{sm}\subset\Sch^{ft}/X=\Var(k)/X$.
\end{itemize}

\item Denote by $\CW\subset\Top$ the full subcategory of $CW$ complexes, by $\CS\subset\CW$ the full subcategory of $\Delta$ complexes,
by $\TM(\mathbb R)\subset\CW$ the full subcategory of topological (real) manifolds
which admits a CW structure (a topological manifold admits a CW structure if it admits a differential structure)
and by $\Diff(\mathbb R)\subset\RTop$ the full subcategory of differentiable (real) manifold.  

\item Denote by $\AnSp(\mathbb C)\subset\RTop$ the full subcategory of analytic spaces over $\mathbb C$,
and by $\AnSm(\mathbb C)\subset\AnSp(\mathbb C)$ the full subcategory of smooth analytic spaces (i.e. complex analytic manifold).
For $X\in\AnSp(\mathbb C)$, we set $\dim(X):=1/2\dim_L(X)$, and for $x\in X$ $\dim_x(X):=1/2\dim_{L,x}(X)$.
For $X\in\AnSp(\mathbb C)$ and $x\in X$, there exist by Weirstrass preparation theorem a finite surjective morphism 
$r:X_x\to\mathbb D^n_0$ where $\mathbb D^n=D(0,1)^n\subset\mathbb C^n$ is the open ball
and $\dim_x(X):=1/2\dim_{L,x}(X)=\dim_K(O_{X,x})=n$. 
For $X\in\AnSm(\mathbb C)$ and $x\in X$, there exist an isomorphism $r:X_x\xrightarrow{\sim}\mathbb D^n_0$,
hence there exist a covering by open subsets $X=\cup_i X_i$ such that $r_i:X_i\xrightarrow{\sim}\mathbb D^{n_i}$.
If $X\in\AnSm(\mathbb C)$ is connected then $\dim(X):=\dim_L(X)=2n$ where $r:X_x\xrightarrow{\sim}\mathbb D^n_0$ for $x\in X$ .
For $X\in\AnSp(\mathbb C)$, $\dim(X):=\dim_L(X)=\dim_L(X_{reg})$ where $X_{reg}\subset X$ is the smooth locus of $X$. 
A morphism $h:U\to S$ with $U,S\in\AnSp(\mathbb C)$ is said to be smooth if it is flat with smooth fibers.
A morphism $r:U\to X$ with $U,X\in\AnSp(\mathbb C)$ is said to be etale if it is non ramified and flat.
For $X\in\AnSp(\mathbb C)$, we denote by 
\begin{itemize}
\item $X^{et}\subset\AnSp(\mathbb C)/X$ the full subcategory consisting of objects
$U/X=(X,h)\in\AnSp(\mathbb C)/X$ such that $h:U\to X$ is an etale morphism.
\item $X^{sm}\subset\AnSp(\mathbb C)/X$ the full subcategory consisting of objects
$U/X=(X,h)\in\AnSp(\mathbb C)/X$ such that $h:U\to X$ is a smooth morphism.
\end{itemize}
By the Weirstrass preparation theorem (or the implicit function theorem if $U$ and $X$ are smooth),
a morphism $r:U\to X$ with $U,X\in\AnSp(\mathbb C)$ is etale if and only if it is an isomorphism local.
Hence for $X\in\AnSp(\mathbb C)$, the morphism of site $\pi_X:X^{et}\to X$ is an isomorphism of site.

\item For $V\in\Var(\mathbb C)$, we denote by $V^{an}\in\AnSp(\mathbb C)$   
the complex analytic space associated to $V$ with the usual topology induced by the usual topology of $\mathbb C^N$. 
For $W\in\AnSp(\mathbb C)$, we denote by $W^{cw}\in\AnSp(\mathbb C)$ the topological space given by $W$ which is a $CW$ complex.
For simplicity, for $V\in\Var(\mathbb C)$, we denote by $V^{cw}:=(V^{an})^{cw}\in\CW$. We have then  
\begin{itemize}
\item the analytical functor $\An:\Var(\mathbb C)\to\AnSp(\mathbb C)$, $\An(V)=V^{an}$,
\item the forgetful functor $\Cw=tp:\AnSp(\mathbb C)\to\CW$, $\Cw(W)=W^{cw}$,
\item the composite of these two functors $\widetilde\Cw=\Cw\circ\An:\Var(\mathbb C)\to\CW$, $\widetilde\Cw(V)=V^{cw}$. 
\end{itemize}

\item Let $S\in\RTop$. Let $S=\cup_{i\in L}S_i$ open cover and
$i_i:S_i\hookrightarrow\tilde S_i$ closed embedding with $\tilde S_i\in\RTop$, $L\in\Set$.
We denote by $(\tilde S_I)\in\Fun(\mathcal P(L)^{op},\RTop)$ the diagram given by for $I\in L$
$\tilde S_L:=\Pi_{i\in I}\tilde S_I$ and for $I\subset J$, $p_{IJ}:\tilde S_J\to\tilde S_I$ is the projection.
We have then open embeddings $j_I:S_I:=\cap_{i\in I}S_i\hookrightarrow S$ and
closed embeddings $i_I:S_I\hookrightarrow\tilde S_I$. We consider the functor
\begin{eqnarray*}
T(S/(\tilde S_I)):C(S)\to C(S/(\tilde S_I))\hookrightarrow C((\tilde S_I)), \; 
K\mapsto T(S/(\tilde S_I))(K):=(i_{I*}j_I^*K,I).
\end{eqnarray*}

\item Let $k\subset\mathbb C$ a subfield. For $S\in\Var(k)$, let $S=\cup_iS_i$ affine open cover and
$i_i:S_i\hookrightarrow\tilde S_i$ closed embedding with $\tilde S_i\in\SmVar(k)$ connected.
We denote by $DR(S):=DR(S)^{[-]}$ the De Rham functor
\begin{eqnarray*}
DR(S):=DR(S)^{[-]}:C_{\mathcal Dfil}(S/(\tilde S_I))\to C(S_{\mathbb C}^{an}/(\tilde S_{I,\mathbb C}^{an})), \\
((M_I,F),u_{IJ})\mapsto DR(S)((M_I,F),u_{IJ}):=(DR(\tilde S_{I,\mathbb C}^{an})(M_I,F)[-d_{\tilde S_I}],DR(u_{IJ}))
\end{eqnarray*}

\item We denote by $\AdSp\subset\RvTop$ the full subcategory of adic spaces.
By definition, for $X=(X,O_X,O_X^+)\in\AdSp$ there exist an open cover $X=\cup_iX_i$ such that 
\begin{equation*}
X_i=\Spa(R_i,R_i^+):=\left\{v\in\Spv(R_i), \, v(f)\leq 1 \, \mbox{for all} \, f\in R_i^+\right\}
\end{equation*}
with $R_i\in\cTRing$ for a non archimedean semi-norm and $R_i^+\subset R_i^o\subset R_i$ a subring, 
where $R_i^o=\left\{f\in R_i s.t |f_i|\leq 1\right\}$.
We then have
\begin{equation*}
R_i^+=\left\{f\in R_i, \, v(f)\leq 1 \, \mbox{for all} \, v\in\Spa(R_i,R_i^+)\right\}.
\end{equation*}

\item For $K\subset\mathbb C_p$ a p-adic field,
denote by $\AnSp(K)\subset\RTop$ the full subcategory of analytic spaces over $K$,
By definition, for $X\in\AnSp(K)$ there exist an open cover $X=\cup_iX_i$ such that $X_i=\Spv(R_i)$
where $R_i\in\cTRing$ is a Tate algebra over $K$. Note that we have the map $m:X\to M(X)$ in $\RTop$
with $m_{|X_i}:X_i=\Spv(R_i)\to M(X_i)=\Spec(R_i)$ which sends a valuation $v$ to its support 
$p=\left\{f\in R_i, \, v(f)=0\right\}$ and where $M(X_i)$ is endowed with the standard G-topology.
By definition, we have
\begin{itemize}
\item the forgetfull functor $o_K:\AdSp/(K,K^+)\to\AnSp(K)$, such that $o_K(\Spa(R,R^+))=\Spv(R)$,
\item the canonical functor $R_K:\AnSp(K)\to\AdSp/(K,O_K)$ such that $R_K(\Spv(R))=\Spa(R,R^o)$.
\end{itemize}
We denote by $\AnSm(K)\subset\AnSp(K)$ the full subcategory of smooth analytic spaces.
For $X\in\AnSp(K)$, we set $\dim(X):=\dim_L(X)$, and for $x\in X$ $\dim_x(X):=\dim_{L,x}(X)$.
For $X\in\AnSp(K)$ and $x\in X$, there exist by Weirstrass preparation theorem a finite surjective morphism 
$r:X_x\to\mathbb D^n_0$ where $\mathbb D^n=D(0,1)^n\subset K^n$ is the open ball
and $\dim_x(X):=\dim_{L,x}(X)=\dim_K(O_{X,x})=n$. 
For $X=\Spv(A)\in\AnSp(K)$ affinoid, $\dim(X):=\dim_L(X)=\dim_K(A)$ where
for the last equality note that for $0<r<1$, 
$D(0,r)^n=\Spv(K<x_1,\cdots,x_n>_r)\subset D(0,1)^n=\Spv(K<x_1,\cdots,x_n>)$ is a rational open subset 
since the norm is ultrametric in contrast to the complex case.
A morphism $h:U\to S$ with $U,S\in\AnSp(K)$ is said to be smooth if it is flat with smooth geometric fibers.
A morphism $r:U\to X$ with $U,X\in\AnSp(K)$ is said to be etale if it is non ramified and flat.
For $X\in\AnSp(K)$, we denote by 
\begin{itemize}
\item $X^{et}\subset\AnSp(K)/X$ the full subcategory consisting of objects
$U/X=(X,h)\in\AnSp(K)/X$ such that $h:U\to X$ is an etale morphism.
\item $X^{sm}\subset\AnSp(K)/X$ the full subcategory consisting of objects
$U/X=(X,h)\in\AnSp(K)/X$ such that $h:U\to X$ is a smooth morphism.
\end{itemize}
For $X\in\AnSp(K)$, we have the morphism of site $\pi_X:X^{et}\to X$.

\item Let $K\subset\mathbb C_p$ a p-adic field.
For $V\in\Var(K)$, we denote by $V^{an}\in\AnSp(K)$. We have then  
the analytical functor $\An:\Var(K)\to\AnSp(K)$, $\An(V)=V^{an}$, $\An(f)=f^{an}$.
We will also consider the canonical functor $R_K:\AnSp(K)\to\AdSp/(K,O_K)$,
which sends by definition $X=\Spv(R)$ affinoid with $R$ a Tate algebra over $K$ to
$X=\Spa(R,R^o)$ with $R^o:=\left\{f\in R, \, |f|_p\leq 1\right\}$.

\item Denote by $\Top^2$ the category whose set of objects is 
\begin{equation*}
(\Top^2)^0:=\left\{(X,Z), \; Z\subset X \; \mbox{closed}\right\}\subset\Top\times\Top
\end{equation*}
and whose set of morphism between $(X_1,Z_1),(X_2,Z_2)\in\Top^2$ is
\begin{eqnarray*}
\Hom_{\Top^2}((X_1,Z_1),(X_2,Z_2)):= 
\left\{(f:X_1\to X_2), \; \mbox{s.t.} \; Z_1\subset f^{-1}(Z_2)\right\}\subset\Hom_{\Top}(X_1,X_2)
\end{eqnarray*}
For $S\in\Top$, $\Top^2/S:=\Top^2/(S,S)$ is then by definition the category whose set of objects is 
\begin{eqnarray*}
(\Top^2/S)^0:=\left\{((X,Z),h), h:X\to S, \; Z\subset X \; \mbox{closed} \;\right\}\subset\Top/S\times\Top
\end{eqnarray*}
and whose set of morphisms between $(X_1,Z_1)/S=((X_1,Z_1),h_1),(X_2,Z_2)/S=((X_2,Z_2),h_2)\in\Top^2/S$
is the subset
\begin{eqnarray*}
\Hom_{\Top^2/S}((X_1,Z_1)/S,(X_2,Z_2)/S):= \\
\left\{(f:X_1\to X_2), \; \mbox{s.t.} \; h_1\circ f=h_2 \; \mbox{and} \; Z_1\subset f^{-1}(Z_2)\right\}
\subset\Hom_{\RTop}(X_1,X_2)
\end{eqnarray*}
We denote by
\begin{eqnarray*}
\mu_S:\Top^{2,pr}/S:=\left\{((Y\times S,Z),p), p:Y\times S\to S, \; Z\subset Y\times S \; \mbox{closed} \;\right\}
\hookrightarrow\Top^2/S
\end{eqnarray*}
the full subcategory whose objects are those with $p:Y\times S\to S$ a projection,
and again $\mu_S:\Top^2/S\to\Top^{2,pr}/S$ the corresponding morphism of sites. 
We denote by 
\begin{eqnarray*}
\Gr_S^{12}:\Top/S\to\Top^{2,pr}/S, \; X/S\mapsto\Gr_S^{12}(X/S):=(X\times S,\bar X)/S, \\
(g:X/S\to X'/S)\mapsto\Gr_S^{12}(g):=(g\times I_S:(X\times S,\bar X)\to(X'\times S,\bar X'))
\end{eqnarray*}
the graph functor, $X\hookrightarrow X\times S$ being the graph embedding (which is a closed embedding if $X$ is separated),
and again $\Gr_S^{12}:\Top^{2,pr}/S\to\Top/S$ the corresponding morphism of sites.

\item Denote by $\RTop^2$ the category whose set of objects is 
\begin{equation*}
(\RTop^2)^0:=\left\{((X,O_X),Z), \; Z\subset X \; \mbox{closed}\right\}\subset\RTop\times\Top
\end{equation*}
and whose set of morphism between $((X_1,O_{X_1}),Z_1),((X_2,O_{X_2}),Z_2)\in\RTop^2$ is
\begin{eqnarray*}
\Hom_{\RTop^2}(((X_1,O_{X_1}),Z_1),((X_2,O_{X_2}),Z_2)):= \\
\left\{(f:(X_1,O_{X_1})\to (X_2,O_{X_2})), \; \mbox{s.t.} \; Z_1\subset f^{-1}(Z_2)\right\}
\subset\Hom_{\RTop}((X_1,O_{X_1}),(X_2,O_{X_2}))
\end{eqnarray*}
For $(S,O_S)\in\RTop$, $\RTop^2/(S,O_S):=\RTop^2/((S,O_S),S)$ is then by definition the category whose set of objects is 
\begin{eqnarray*}
(\RTop^2/(S,O_S))^0:= \\
\left\{(((X,O_X),Z),h), h:(X,O_X)\to(S,O_S), \; Z\subset X \; \mbox{closed} \;\right\}\subset\RTop/(S,O_S)\times\Top
\end{eqnarray*}
and whose set of morphisms between $(((X_1,O_{X_1}),Z_1),h_1),(((X_2,O_{X_2}),Z_2),h_2)\in\RTop^2/(S,O_S)$
is the subset
\begin{eqnarray*}
\Hom_{\RTop^2/(S,O_S)}(((X_1,O_{X_1}),Z_1)/(S,O_S),((X_2,O_{X_2}),Z_2)/(S,O_S)):= \\
\left\{(f:(X_1,O_{X_1})\to (X_2,O_{X_2})), \; \mbox{s.t.} \; h_1\circ f=h_2 \; \mbox{and} \; Z_1\subset f^{-1}(Z_2)\right\} \\ 
\subset\Hom_{\RTop}((X_1,O_{X_1}),(X_2,O_{X_2}))
\end{eqnarray*}
We denote by 
\begin{eqnarray*}
\mu_S:\RTop^{2,pr}/S:=\left\{(((Y\times S,q^*O_Y\otimes p^*O_S),Z),p), p:Y\times S\to S, \; 
Z\subset Y\times S \; \mbox{closed} \;\right\}\hookrightarrow\RTop^2/S
\end{eqnarray*}
the full subcategory whose objects are those with $p:Y\times S\to S$ is a projection,
and again $\mu_S:\RTop^2/S\to\RTop^{2,pr}/S$ the corresponding morphism of sites. 
We denote by 
\begin{eqnarray*}
\Gr_S^{12}:\RTop/S\to\RTop^{2,pr}/S, \\ 
(X,O_X)/(S,O_S)\mapsto\Gr_S^{12}((X,O_X)/(S,O_S)):=((X\times S,q^*O_X\otimes p^*O_S),\bar X)/(S,O_S), \\
(g:(X,O_X)/(S,O_S)\to (X',O_{X'})/(S,O_S))\mapsto \\
\Gr_S^{12}(g):=(g\times I_S:((X\times S,q^*O_X\otimes p^*O_S),\bar X)\to((X'\times S,q^*O_X\otimes p^*O_S),\bar X'))
\end{eqnarray*}
the graph functor, $X\hookrightarrow X\times S$ being the graph embedding (which is a closed embedding if $X$ is separated),
$p:X\times S\to S$, $q:X\times S\to X$ the projections,
and again $\Gr_S^{12}:\RTop^{2,pr}/S\to\RTop/S$ the corresponding morphism of sites.

\item We denote by $\Sch^2\subset\RTop^2$ the full subcategory such that the first factors are schemes. 
For a field $k$, we denote by $\Sch^2/k:=\Sch^2/(\Spec k,\left\{\pt\right\})$ and by
\begin{itemize}
\item $\Var(k)^2\subset\Sch^2/k$ the full subcategory such that the first factors are algebraic varieties over $k$, 
i.e. schemes of finite type over $k$,
\item $\PVar(k)^2\subset\QPVar(k)^2\subset\Var(k)^2$ the full subcategories such that the first factors are
quasi-projective varieties and projective varieties respectively, 
\item $\PSmVar(k)^2\subset\SmVar(k)^2\subset\Var(k)^2$ the full subcategories such that the first factors are
smooth varieties and smooth projective varieties respectively.
\end{itemize}
In particular we have, for $S\in\Var(k)$, the graph functor 
\begin{eqnarray*}
\Gr_S^{12}:\Var(k)/S\to\Var(k)^{2,pr}/S, \; X/S\mapsto\Gr_S^{12}(X/S):=(X\times S,X)/S, \\
(g:X/S\to X'/S)\mapsto\Gr_S^{12}(g):=(g\times I_S:(X\times S,X)\to(X'\times S,X'))
\end{eqnarray*}
the graph embedding $X\hookrightarrow X\times S$ is a closed embedding 
since $X$ is separated in the subcategory of schemes $\Sch\subset\RTop$,
and again $\Gr_S^{12}:\Var(k)^{2,pr}/S\to\Var(k)/S$ the corresponding morphism of sites.

\item We denote by $\CW^2\subset\Top^2$ the full subcategory such that the first factors are $CW$ complexes,
by $\TM(\mathbb R)^2\subset\CW^2$ the full subcategory such that the first factors are topological (real) manifolds
and by $\Diff(\mathbb R)^2\subset\RTop^2$ the full subcategory such that the first factors are differentiable (real) manifold. 

\end{itemize}

\subsection{The $p$-adic de Rham period sheaves for adic spaces over a $p$-adic field}

Let $p$ a prime number.
For $X=(X,O_X,O_X^+)\in\AdSp/(K,K^+)$ an adic space over a p-adic field $K\subset\mathbb C_p$, we consider
\begin{itemize}
\item the map $W_X:\mathbb A_{inf,X}:=W(\hat O_X^{b+})\to\hat O_X^+$ 
where $\hat O_X$ denote the completion of $O_X$ with respect to
the ideal $pO_X\subset O_X$, $b$ the tilting functor and $W$ the Witt vectors,
\item the map $W_X:\mathbb B_{inf,X}:=W(\hat O_X^{b+})[p^{-1}]\to\hat O_X:=\hat O_X^+[p^{-1}]$
\item the integral period sheaf $\mathbb B^+_{dr,X}:=\varprojlim_{n\in\mathbb N}\mathbb B_{inf,X}/(\ker W_X)^n$
with the filtration $F^k\mathbb B^+_{dr,X}:=(\ker W_X)^k\mathbb B^+_{dr,X}\subset\mathbb B^+_{dr,X}$.
\item the period sheaf $\mathbb B_{dr,X}:=\mathbb B^+_{dr,X}[t^{-1}]$ 
where $t$ is a generator of the ideal $\ker W_X\subset\mathbb B^+_{dr,X}$
with the filtration $F^k\mathbb B_{dr,X}:=\sum_j t^{-j}F^{k+j}\mathbb B^+_{dr,X}\subset\mathbb B_{dr,X}$.
\item the integral sheaf 
\begin{equation*}
O\mathbb B_{dr,X}^+:=\varprojlim_{n\in\mathbb N}
(O_X^+\hat\otimes_{W(O_K/pO_K)}\mathbb A_{inf,X}[p^{-1}])/(\ker I\otimes W_X)^n
\end{equation*}
with the filtration $F^kO\mathbb B^+_{dr,X}:=(\ker W_X)^kO\mathbb B^+_{dr,X}\subset O\mathbb B^+_{dr,X}$.
\item the period sheaf $O\mathbb B_{dr,X}:=O\mathbb B^+_{dr,X}[t^{-1}]$ 
where $t$ is a generator of the ideal $\ker W_X\subset O\mathbb B^+_{dr,X}$
with the filtration $F^kO\mathbb B_{dr,X}:=\sum_j t^{-j}F^{k+j}O\mathbb B^+_{dr,X}\subset O\mathbb B_{dr,X}$.
\end{itemize}

\subsection{The classical theorems for etale topology on schemes, for CW complexes,
and the comparaison theorem with the analytic topologies for algebraic varieties over local fields}

We first recall the smooth base change theorem

\begin{thm}
\begin{itemize}
\item[(i)]Consider a commutative diagram in $\Sch$ which is cartesian 
\begin{equation*}
\xymatrix{X_T\ar[r]^{f'}\ar[d]^{g'} & T\ar[d]^{g} \\
X\ar[r]^f & S}
\end{equation*}
such that $g$ is smooth or more generally locally acyclic. 
Let $F\in C(X^{et})$ be a torsion sheaf where we recall that $X^{et}\subset\Sch^{ft}/X$ is the small etale site.
Then the transformation map (see \cite{B4} section 2) in $D(T^{et})$
\begin{equation*}
T(f,g)(F):g^*Rf_*F\to Rf'_*g^{'*}F
\end{equation*}
is an isomorphism.
\item[(ii)] Let $k'/k$ an extention of field of characteristic zero.
Let $f:X\to S$ a morphism in $\Var(k)$. Let $F\in C(X^{et})$ be a torsion sheaf
Then the transformation map (see \cite{B4} section 2) in $D(S_{k'}^{et})$
\begin{equation*}
T(f,\pi_{k/k'})(F):\pi_{k/k'}^*Rf_*F\to Rf_{k'*}\pi_{k/k'}^*F
\end{equation*}
is an isomorphism where we recall (see section 2)
$\pi_{k/k'}=\pi_{k/k'}(X):X_{k'}\to X$ and $\pi_{k/k'}=\pi_{k/k'}(S):S_{k'}\to S$
are the projections.
\end{itemize}
\end{thm}

\begin{proof}
\noindent(i): Standard : see \cite{milnes} for example.

\noindent(ii):Follows from (i).
\end{proof}

We now recall the proper base change theorem :

\begin{thm}\label{PBCetth}
Consider a commutative diagram in $\Sch$ which is cartesian 
\begin{equation*}
\xymatrix{X_T\ar[r]^{f'}\ar[d]^{g'} & T\ar[d]^{g} \\
X\ar[r]^f & S}
\end{equation*}
such that $f$ is proper. 
Let $F\in C(X^{et})$ be a torsion sheaf where we recall that $X^{et}\subset\Sch^{ft}/X$ is the small etale site.
Then the transformation map (see \cite{B4} section 2) in $D(T^{et})$
\begin{equation*}
T(f,g)(F):g^*Rf_*F\to Rf'_*g^{'*}F
\end{equation*}
is an isomorphism.
\end{thm}

\begin{proof}
Standard : see \cite{milnes} for example.
\end{proof}

We deduce from the proper base change theorem the projection formula and the Kunneth formula for 
the cohomology of etale sheaves:

\begin{thm}\label{PFet}
\begin{itemize}
\item[(i)]Let $f:X\to S$ a proper morphism with $S,X\in\Sch$. Let $n\in\mathbb N$.
Let $F\in C_{\mathbb Z/n\mathbb Z}(X^{et})$ and $G\in C_{\mathbb Z/n\mathbb Z}(S^{et})$ be $n$-torsion sheaves.
Then the transformation map (see \cite{B4} section 2) in $D_{\mathbb Z/n\mathbb Z}(S^{et})$
\begin{equation*}
T(f,\otimes)(F,G):Rf_*F\otimes_{\mathbb Z/n\mathbb Z}^L G\to Rf_*(F\otimes_{\mathbb Z/n\mathbb Z}^Lf^*G)
\end{equation*}
is an isomorphism.
\item[(ii)]Let $f:X\to S$ a morphism with $S,X\in\Sch$. Let $n\in\mathbb N$.
Let $F\in C_{\mathbb Z/n\mathbb Z}(X^{et})$ and $G\in C_{\mathbb Z/n\mathbb Z}(S^{et})$ be $n$-torsion sheaves.
Then the transformation map (see \cite{B4} section 2) in $D_{\mathbb Z/n\mathbb Z}(S^{et})$
given by (i) and the open embedding case after taking a compactification of $f$
\begin{equation*}
T_!(f,\otimes)(F,G):Rf_!F\otimes_{\mathbb Z/n\mathbb Z}^L G\to Rf_!(F\otimes_{\mathbb Z/n\mathbb Z}^Lf^*G)
\end{equation*}
is an isomorphism.
\end{itemize}
\end{thm}

\begin{proof}
\noindent(i):Follows from theorem \ref{PBCetth}: see \cite{milnes} for example.

\noindent(ii):Follows from (i) by taking a compactification $\bar f:\bar X\to\bar S$ of $f:X\to S$.
\end{proof}

\begin{rem}\label{PFetrem}
Let $f:X\to S$ a morphism with $S,X\in\Sch$. Let $n\in\mathbb N$.
Let $F\in C_{\mathbb Z/n\mathbb Z}(X^{et})$ and $G\in C_{\mathbb Z/n\mathbb Z}(S^{et})$ be $n$-torsion sheaves.
Then, if $f$ is not proper, $H^k(Rf_*F\otimes_{\mathbb Z/n\mathbb Z}^L G)$
is NOT isomorphic in $\Shv_{\mathbb Z/n\mathbb Z}(S^{et})$ to $H^kRf_*(F\otimes_{\mathbb Z/n\mathbb Z}^LG)$
in general.
\end{rem}

Let $f_1:X_1\to S$ and $f_2:X_2\to S$ two morphisms with $X_1,X_2,S\in\Sch$.
Denote $p_1:X_1\times_SX_2\to X_1$ and $p_2:X_1\times_SX_2\to X_2$ the base change maps.
We have then 
\begin{equation*}
f_1\otimes f_2=f_1\circ p_1=f_2\circ p_2:X_1\times_SX_2\to S.
\end{equation*}
Let $F_1\in C_{\mathbb Z/n\mathbb Z}(X_1^{et})$ and $F_2\in C_{\mathbb Z/n\mathbb Z}(X_2^{et})$ be $n$-torsion sheaves.
Then the canonical map in $C_{\mathbb Z/n\mathbb Z}(S^{et})$ (see \cite{B4} section 2)
\begin{eqnarray*}
T(f_1,f_2,\otimes)(F_1,F_2):Rf_{1*}F_1\otimes_{\mathbb Z/n\mathbb Z}^LRf_{2*}F_2
\xrightarrow{\ad(p_1^*,Rp_{1*})(F_1)\otimes\ad(p_2^*,Rp_{2*})(F_2)} \\
Rf_{1*}Rp_{1*}p_1^*F_1\otimes_{\mathbb Z/n\mathbb Z}^LRf_{2*}Rp_{2*}p_2^*F_2 
\xrightarrow{=}R(f_1\otimes f_2)_*p_1^*F_1\otimes_{\mathbb Z/n\mathbb Z}^LR(f_1\otimes f_2)_*p_2^*F_2 \\
\xrightarrow{T(\otimes,E)(-,-)}R(f_1\times f_2)_*(p_1^*F\otimes_{\mathbb Z/n\mathbb Z}^Lp_2^*F_2).
\end{eqnarray*}

\begin{thm}\label{KUNetth}
\begin{itemize}
\item[(i)] Let $f_1:X_1\to S$ and $f_2:X_2\to S$ two proper morphisms with $X_1,X_2,S\in\Sch$.
Denote $p_1:X_1\times_SX_2\to X_1$ and $p_2:X_1\times_SX_2\to X_2$ the base change maps.
Let $F_1\in C_{\mathbb Z/n\mathbb Z}(X_1^{et})$ and $F_2\in C_{\mathbb Z/n\mathbb Z}(X_2^{et})$ be $n$-torsion sheaves.
Then the canonical map in $C_{\mathbb Z/n\mathbb Z}(S^{et})$ given above
\begin{equation*}
T(f_1,f_2,\otimes)(F_1,F_2):Rf_{1*}F_1\otimes_{\mathbb Z/n\mathbb Z}^LRf_{2*}F_2\to 
R(f_1\times f_2)_*(p_1^*F\otimes_{\mathbb Z/n\mathbb Z}^Lp_2^*F_2)
\end{equation*}
is an isomorphism.
\item[(ii)] Let $f_1:X_1\to S$ and $f_2:X_2\to S$ two morphisms with $X_1,X_2,S\in\Sch$.
Denote $p_1:X_1\times_SX_2\to X_1$ and $p_2:X_1\times_SX_2\to X_2$ the base change maps.
Let $F_1\in C_{\mathbb Z/n\mathbb Z}(X_1^{et})$ and $F_2\in C_{\mathbb Z/n\mathbb Z}(X_2^{et})$ be $n$-torsion sheaves.
Then the canonical map in $C_{\mathbb Z/n\mathbb Z}(S^{et})$ given after taking compactification of $f_1$ and $f_1$
by the one of (i) for the compactifications and on the other hand by the open embedding case
\begin{equation*}
T(f_1,f_2,\otimes)(F_1,F_2):Rf_{1!}F_1\otimes_{\mathbb Z/n\mathbb Z}^LRf_{2!}F_2\to 
R(f_1\times f_2)_!(p_1^*F\otimes_{\mathbb Z/n\mathbb Z}^Lp_2^*F_2)
\end{equation*}
is an isomorphism.
\end{itemize}
\end{thm}

\begin{proof}
\noindent(i): Follows from theorem \ref{PBCetth}: see \cite{milnes} for example.

\noindent(ii):Follows from (i) by taking a compactification $\bar f_1:\bar X_1\to\bar S$ of $f_1:X_1\to S$
and a compactification $\bar f_2:\bar X_2\to\bar S$ of $f_2:X_2\to S$.
\end{proof}

\begin{rem}\label{KUNetthrem}
Let $f_1:X_1\to S$ and $f_2:X_2\to S$ two morphisms with $S,X_1,X_2\in\Sch$. Let $n\in\mathbb N$.
Denote $p_1:X_1\times_SX_2\to X_1$ and $p_2:X_1\times_SX_2\to X_2$ the base change maps.
Let $F_1\in C_{\mathbb Z/n\mathbb Z}(X_1^{et})$ and $F_2\in C_{\mathbb Z/n\mathbb Z}(X_2^{et})$ be $n$-torsion sheaves.
Then, if $f_1$ or $f_2$ is not proper, $H^k(Rf_{1*}F_1\otimes_{\mathbb Z/n\mathbb Z}^LRf_{2*}F_2)$
is NOT isomorphic in $\Shv_{\mathbb Z/n\mathbb Z}(S^{et})$ to 
$H^kR(f_1\times f_2)_*(p_1^*F_1\otimes_{\mathbb Z/n\mathbb Z}^Lp_2^*F_2)$
in general.
\end{rem}

We now recall the comparaison theorems :

\begin{thm}\label{CPetusu}
Let $f:X\to S$ a morphism with $X,S\in\Var(\mathbb C)$.
Let $F\in C(X^{et})$ be a torsion sheaf where we recall that 
$X^{et}\subset\Sch^{ft}/X=\Var(\mathbb C)/X$ is the small etale site.
Then the transformation map (see \cite{B4} section 2) in $D(S^{an})$
\begin{equation*}
T(f,an)(F):\an_S^*Rf_*F\to Rf_*\an_X^*F
\end{equation*}
is an isomorphism.
\end{thm}

\begin{proof}
Standard, see \cite{milnes} for example : follows from theorem \ref{PBCetth} (i) and the open embedding case.
\end{proof}

\begin{thm}\label{CPetbet}
Let $K\subset\mathbb C_p$ be a $p$-adic field. Let $f:X\to S$ a morphism with $X,S\in\Var(K)$.
Let $F\in C_{\mathbb Z/n\mathbb Z}(X^{et})$ be an $n$-torsion sheaf where we recall that 
$X^{et}\subset\Sch^{ft}/X=\Var(\mathbb C)/X$ is the small etale site.
Then the transformation map (see \cite{B4} section 2) in $D(S^{an,et})$
\begin{equation*}
T(f,an)(F):\an_S^*Rf_*F\to Rf_*\an_X^*F
\end{equation*}
is an isomorphism.
\end{thm}

\begin{proof}
Standard, see \cite{huber} for example : follows from theorem \ref{PBCetth} (i) and the open embedding case.
\end{proof}

On the other hand for CW complexes, we have the followings :

\begin{thm}\label{PBCusu}
Consider a commutative diagram in $\CW$ which is cartesian 
\begin{equation*}
\xymatrix{X_T\ar[r]^{f'}\ar[d]^{g'} & T\ar[d]^{g} \\
X\ar[r]^f & S}
\end{equation*}
such that $f$ is proper. Let $F\in C(X)$. Then the transformation map (see \cite{B4} section 2) in $D(T)$
\begin{equation*}
T(f,g)(F):g^*Rf_*F\to Rf'_*g^{'*}F
\end{equation*}
is an isomorphism.
\end{thm}

\begin{proof}
Standard.
\end{proof}

\begin{thm}\label{PFusu}
\begin{itemize}
\item[(i)]Let $f:X\to S$ a proper morphism with $S,X\in\CW$. Let $F\in C(X)$ and $G\in C(S)$.
Then the transformation map (see \cite{B4} section 2) in $D(S)$
\begin{equation*}
T(f,\otimes)(F,G):Rf_*F\otimes^L G\to Rf_*(F\otimes^Lf^*G)
\end{equation*}
is an isomorphism.
\item[(ii)]Let $f:X\to S$ a morphism with $S,X\in\CW$. Let $F\in C(X)$ and $G\in C(S)$.
Then the transformation map (see \cite{B4} section 2) in $D(S)$
\begin{equation*}
T_!(f,\otimes)(F,G):Rf_!F\otimes^L G\to Rf_!(F\otimes^Lf^*G)
\end{equation*}
is an isomorphism.
\end{itemize}
\end{thm}

\begin{proof}
Standard.
\end{proof}

\begin{thm}\label{KUNusu}
\begin{itemize}
\item[(i)] Let $f_1:X_1\to S$ and $f_2:X_2\to S$ two proper morphisms with $X_1,X_2,S\in\CW$.
Denote $p_1:X_1\times_SX_2\to X_1$ and $p_2:X_1\times_SX_2\to X_2$ the base change maps.
Let $F_1\in C(X_1)$ and $F_2\in C(X_2)$. Then the canonical map in $C(S)$ given as above
\begin{equation*}
T(f_1,f_2,\otimes)(F_1,F_2):Rf_{1*}F_1\otimes^LRf_{2*}F_2\to R(f_1\times f_2)_*(p_1^*F\otimes^Lp_2^*F_2)
\end{equation*}
is an isomorphism.
\item[(ii)] Let $f_1:X_1\to S$ and $f_2:X_2\to S$ two morphisms with $X_1,X_2,S\in\CW$.
Denote $p_1:X_1\times_SX_2\to X_1$ and $p_2:X_1\times_SX_2\to X_2$ the base change maps.
Let $F_1\in C(X_1)$ and $F_2\in C(X_2)$.
Then the canonical map in $C(S)$ given after taking compactification of $f_1$ and $f_1$
by the one of (i) for the compactifications and on the other hand by the open embedding case
\begin{equation*}
T(f_1,f_2,\otimes)(F_1,F_2):Rf_{1!}F_1\otimes^LRf_{2!}F_2\to R(f_1\times f_2)_!(p_1^*F\otimes^Lp_2^*F_2)
\end{equation*}
is an isomorphism.
\end{itemize}
\end{thm}

\begin{proof}
Standard.
\end{proof}

\subsection{Constructible and perverse sheaves on algebraic varieties over a subfield $k\subset\mathbb C$}

Let $S\in\AnSp(\mathbb C)$. We have
\begin{itemize}
\item the classical dual functor
\begin{eqnarray*}
\mathbb D^0_S:C(S)\to C(S), \; \; K\mapsto\mathbb D_S^0K:=\mathcal Hom(LK,E_{usu}(\mathbb Z_S))
\end{eqnarray*}
which induces in the derived category
\begin{eqnarray*}
\mathbb D^0_S:D(S)\to D(S), \; \;
K\mapsto\mathbb D_S^0K:=\mathcal Hom(LK,E_{usu}(\mathbb Z_S))=R\mathcal Hom(K,\mathbb Z_S)
\end{eqnarray*}
\item the Verdier dual functor
\begin{eqnarray*}
\mathbb D^v_S:D(S)\to D(S), \; \; K\mapsto\mathbb D_S^vK:=R\mathcal Hom(K,w_S)
\end{eqnarray*}
where $w_S:=a_S^!\mathbb Z_S$ is the dualizing complex, $a_S:S\to\left\{\pt\right\}$ being the terminal map.
\end{itemize}
For $S\in\AnSm(\mathbb C)$ smooth connected of dimension $\dim(S)=d_S$ we have $\mathbb D^v_S=\mathbb D^0_S[2d_S]$.

We recall the definition of constructible sheaves on algebraic varieties over a subfield $k\subset\mathbb C$:

\begin{defi}\label{constr1}
Let $k\subset\mathbb C$ a subfield. Let $S\in\Var(k)$.
\begin{itemize}
\item[(i)]A sheaf $K\in\Shv(S_{\mathbb C}^{an})$ is called constructible 
(with respect to a Zariski stratification over $k$) if there exists a stratification 
$S=\sqcup_{\alpha} S_{\alpha}$ with $l_{\alpha}:S_{\alpha}\hookrightarrow S$ locally closed subsets
such that $l_{\alpha}^*K\in\Shv(S_{\alpha,\mathbb C}^{an})$ 
are (finite dimensional) local systems (for the usual topology) for all $\alpha$.
Note that we make the hypothesis that the strata $S_{\alpha}$ are defined over $k$. 
\item[(ii)] We denote by 
\begin{equation*}
C_{c,k}(S_{\mathbb C}^{an})\subset C(S_{\mathbb C}^{an}) \; \mbox{and} \;
D_{c,k}(S_{\mathbb C}^{an})\subset D(S_{\mathbb C}^{an})
\end{equation*}
the full subcategories consisting of $K\in C(S_{\mathbb C}^{an})$ 
such that $a_{usu}H^nK\in\Shv(S_{\mathbb C}^{an})$ are constructible 
with respect to a Zariski stratification over $k$ for all $n\in\mathbb Z$ (see (i)).
\item[(ii)'] We denote by 
\begin{equation*}
C_{fil,c,k}(S_{\mathbb C}^{an})\subset C_{fil}(S_{\mathbb C}^{an}) \; \mbox{and} \;
D_{fil,c,k}(S_{\mathbb C}^{an})\subset D_{fil}(S_{\mathbb C}^{an})
\end{equation*}
the full subcategories consisting of $(K,W)\in C_{fil}(S_{\mathbb C}^{an})$ 
such that $a_{usu}H^n\Gr_k^WK\in\Shv(S_{\mathbb C}^{an})$ are constructible 
with respect to a Zariski stratification over $k$ for all $n,k\in\mathbb Z$ (see (i)).
\item[(iii)] We denote by 
\begin{equation*}
P_k(S_{\mathbb C}^{an}):=P(S_{\mathbb C}^{an})\cap D_{c,k}(S_{\mathbb C}^{an})
\subset D_{c,k}(S_{\mathbb C}^{an})\otimes\mathbb Q\subset D_{c,k}(S_{\mathbb C}^{an}) 
\end{equation*}
the full subcategory of perverse sheaves (which are by definition torsion free) whose cohomology sheaves are
constructible with respect to a Zariski stratification (defined over $k$), see (ii).
\item[(iii)'] We denote by 
\begin{equation*}
P_{fil,k}(S_{\mathbb C}^{an}):=P_{fil}(S_{\mathbb C}^{an})\cap D_{fil,c,k}(S_{\mathbb C}^{an})
\subset D_{fil,c,k}(S_{\mathbb C}^{an})\otimes\mathbb Q\subset D_{fil,c,k}(S_{\mathbb C}^{an}) 
\end{equation*} 
the full subcategory of filtered perverse sheaves (which are by definition torsion free)
whose cohomology sheaves are constructible with respect to a Zariski stratification (defined over $k$), see (ii)'.
\end{itemize}
\end{defi}

\begin{thm}\label{constrthm1}
Let $k\subset\mathbb C$ a subfield. 
\begin{itemize}
\item Let $S\in\Var(k)$. Then for $K\in D_{c,k}(S_{\mathbb C}^{an})$, $\mathbb D^v_SK\in D_{c,k}(S_{\mathbb C}^{an})$.
For $K,K'\in D_{c,k}(S_{\mathbb C}^{an})$, $K\otimes^L K'\in D_{c,k}(S_{\mathbb C}^{an})$.
\item Let $f:X\to S$ a morphism with $S,X\in\Var(k)$.
Then for $K\in D_{c,k}(X_{\mathbb C}^{an})$, $Rf_*K\in D_{c,k}(S_{\mathbb C}^{an})$ and 
$Rf_!K=\mathbb D^v_SRf_*\mathbb D^v_XK\in D_{c,k}(S_{\mathbb C}^{an})$.
\item Let $f:X\to S$ a morphism with $S,X\in\Var(k)$.
Then for $K\in D_{c,k}(S_{\mathbb C}^{an})$, $f^*K\in D_{c,k}(X_{\mathbb C}^{an})$ and 
$f^!K:=\mathbb D^v_Xf^*\mathbb D^v_SK\in D_{c,k}(X_{\mathbb C}^{an})$.
\end{itemize}
For $S\in\Var(k)$, we have thus the full subcategory
\begin{eqnarray*}
D_{c,k}(S_{\mathbb C}^{an})^{gm}:&=&<Rf_*\mathbb Z_{X_{\mathbb C}^{an}},(f:X\to S)\in\Var(k)> \\
&=&<Rf_*\mathbb Z_{X_{\mathbb C}^{an}},(f:X\to S)\in\Var(k) \; \mbox{proper}, \; X \mbox{smooth}>
\subset D_{c,k}(S_{\mathbb C}^{an})
\end{eqnarray*}
where $<>$ means the full triangulated category generated by.
\end{thm}

\begin{proof}
Standard : follows from the fact that a morphism $f:X\to S$ with $X,S\in\Var(k)$ admits a Whitney stratification
whose strata $X=\sqcup_{\alpha,\beta} X_{\alpha,\beta}$, $S=\sqcup_{\alpha} S_{\alpha}$ 
are Zariski locally closed subset, that is locally closed subvarieties defined over $k$.
\end{proof}

\begin{defi}
Let $k\subset\mathbb C$ a subfield
We have then, using theorem \ref{constrthm1}, for $S\in\Var(k)$, the full subcategory
\begin{eqnarray*}
D_{fil,c,k}(S_{\mathbb C}^{an})^{gm}:=<(K,W), \, \mbox{s.t.,for all} \, n\in\mathbb Z, \\ 
\Gr_n^WK=Rf_{n*}\mathbb Z_{X_{n,\mathbb C}^{'an}}, \, 
X_n\in\SmVar(k), \, f_n:X_n\to S \, \mbox{proper}>\subset D_{fil,c,k}(S_{\mathbb C}^{an})
\end{eqnarray*}
where $<>$ means the full triangulated category generated by.
\end{defi}

Let $k\subset\mathbb C$ a subfield.

Let $S\in\Var(k)$ and $D=V(s)\subset S$ a Cartier divisor. 
Denote $i:D\hookrightarrow S$ the closed embedding and $j:S^o:=S\backslash D\hookrightarrow S$ the open embedding.
Let $\pi:\tilde S_{\mathbb C}^{o,an}\to S_{\mathbb C}^{o,an}$ the universal covering.
Denote by $T:\tilde S_{\mathbb C}^{o,an}\to\tilde S_{\mathbb C}^{o,an}$ the monodromy automorphism.
We then consider,
\begin{itemize}
\item for $K\in C(S_{\mathbb C}^{o,an})$, the nearby cycle functor 
\begin{equation*}
\psi_DK:=i^*R(j\circ\pi)_*\pi^*K\in D(D_{\mathbb C}^{an}), 
\end{equation*}
we write again $\psi_DK:=i_*\psi_DK\in D(S_{\mathbb C}^{an})$,
\item for $K\in C(S_{\mathbb C}^{o,an})$, the vanishing cycle functor 
\begin{equation*}
\phi_DK:=\Cone(\ad(j\circ\pi^*,j\circ\pi_*)(K):i^*K\to i^*R(j\circ\pi)_*\pi^*K=:\psi_DK\in D(D_{\mathbb C}^{an})
\end{equation*}
together with the canonical map $c(\phi_DK):\psi_DK\to\phi_DK$ in $D(D_{\mathbb C}^{o,an})$, 
we write again $\phi_DK:=i_*\phi_DK\in D(S_{\mathbb C}^{an})$,
\item for $K\in C_c(S_{\mathbb C}^{o,an})$, the canonical morphisms in $D_c(D_{\mathbb C}^{an})$
\begin{eqnarray*}
can(K):=c(\phi_DK):\psi_DK\to\phi_DK \; , \; var(K):=(0,T-I):\phi_DK\to\psi_DK.
\end{eqnarray*}
\end{itemize}
For $K\in C_{c,k}(S_{\mathbb C}^{o,an})$, 
we have $\psi_DK\in C_{c,k}(D_{\mathbb C}^{o,an})$ since it preserve local systems : 
if $(S_{\alpha})$ is a stratification for $K$ 
then $(S_{\alpha}\cap S^o,S_{\alpha}\cap D)$ is a stratification for $\psi_DK$.
This implies, for $K\in C_{c,k}(S_{\mathbb C}^{o,an})$, $\phi_DK\in C_{c,k}(D_{\mathbb C}^{o,an})$.

Let $S\in\Var(k)$ and $D=V(s)\subset S$ a Cartier divisor. 
Denote $i:D\hookrightarrow S$ the closed embedding and $j:S^o:=S\backslash D\hookrightarrow S$ the open embedding.
Let $\pi:\tilde S_{\mathbb C}^{o,an}\to S_{\mathbb C}^{o,an}$ the universal covering. 
\begin{itemize}
\item For $K,K'\in C_c(S_{\mathbb C}^{o,an})$, we have by (Verdier) duality and theorem \ref{KUNusu}
$\psi_D(K\otimes^L K')=\psi_D(K)\otimes^L\psi_D(K')$,
\item For $K,K'\in C_c(S_{\mathbb C}^{o,an})$ we have by (Verdier) duality and the preceding point 
$\phi_D(K\otimes^L K')=\phi_D(K)\otimes^L\phi_D(K')$
\item For $K,K'\in C(S_{\mathbb C}^{o,an})$, we have the tranformation map in $D(S_{\mathbb C}^{o,an})$
\begin{eqnarray*}
T(\otimes,\psi_D)(K,K'):(\psi_DK)\otimes^L K'\xrightarrow{I\otimes(\ad(i^*,i_*)(-)\circ\ad(\pi^*,\pi_*)(-))}
\psi_DK\otimes^L\psi_DK'=\psi_D(K\otimes^L K').
\end{eqnarray*}
\item For $K,K'\in C_c(S_{\mathbb C}^{o,an})$ the tranformation map in $D_c(S_{\mathbb C}^{o,an})$
\begin{eqnarray*}
T(\otimes,\phi_D)(K,K'):\phi_D(K\otimes^L K')=\mathbb D_S^v\psi_D(\mathbb D_S^vK\otimes^L\mathbb D_S^vK') \\
\xrightarrow{\mathbb D_S^vT(\otimes,\psi_D)(\mathbb D_S^vK,\mathbb D_S^vK')}
\mathbb D_S^v(\psi_D\mathbb D_S^v(K)\otimes^L\mathbb D_S^vK')=\phi_DK\otimes^L K'.
\end{eqnarray*}
\item for $K\in C_c(S_{\mathbb C}^{o,an})$ we have a canonical isomorphism in $D_c(D_{\mathbb C}^{an})$
\begin{equation*}
T(D,\psi_D)(K):\phi_D\mathbb D_S^vK\xrightarrow{\sim}\mathbb D_S^v\psi_DK[1]
\end{equation*}
\end{itemize}

\begin{defi}\label{ukC}
Let $S\in\Var(k)$ and $D=V(s)\subset S$ a Cartier divisor. 
Denote $i:D\hookrightarrow S$ the closed embedding and $j:S^o:=S\backslash D\hookrightarrow S$ the open embedding.
Let $\pi:\tilde S_{\mathbb C}^{o,an}\to S_{\mathbb C}^{o,an}$ the universal covering.
Denote by $T:\tilde S_{\mathbb C}^{o,an}\to\tilde S_{\mathbb C}^{o,an}$ the monodromy automorphism.
For $K\in P_k(S_{\mathbb C}^{o,an})$, we consider 
\begin{itemize}
\item the canonical decomposition
\begin{equation*}
\psi_DK:=\psi_D^uK\oplus\psi_D^tK\in P_k(D_{\mathbb C}^{an}), 
\end{equation*}
with $\psi_D^uK:=\ker(T-I)^N\subset\psi_DK$ with $N$ such that $\ker(T-I)^{N+1}=\ker(T-I)^N$
and $\psi_D^tK:=\ker Q(T)\subset\psi_DK$ where $P(X)=(X-1)^NQ(X)$ is the minimal polynomial of $T:\psi_DK\to\psi_DK$, 
note that we have $\psi^u_DK,\psi_D^tK\in P_{k}(D_{\mathbb C}^{o,an})$ 
since $\psi^u_D$ and $\psi^t_D$ preserve local systems : if $(S_{\alpha})$ is a stratification for $K$ 
then $(S_{\alpha}\cap S^o,S_{\alpha}\cap D)$ is a stratification for $\psi^u_DK$ and $\psi^t_DK$,
\item the canonical decomposition
\begin{equation*}
\phi_DK:=\phi_D^uK\oplus\phi_D^tK\in P_k(D_{\mathbb C}^{an}), 
\end{equation*}
with $\phi_D^uK:=\ker(T-I)^N\subset\psi_DK$ with $N$ such that $\ker(T-I)^{N+1}=\ker(T-I)^N$
and $\phi_D^tK:=\ker Q(T)\subset\phi_DK$ where $P(X)=(X-1)^NQ(X)$ is the minimal polynomial of $T:\phi_DK\to\phi_DK$, 
\item the canonical morphisms in $D_c(D_{\mathbb C}^{an})$
\begin{eqnarray*}
can(K):=c(\phi^u_DK):\psi^u_DK\to\phi^u_DK \; , \; var(K):=(0,T-I):\phi^u_DK\to\psi^u_DK.
\end{eqnarray*}
\item the maximal extension
\begin{eqnarray*}
x_{S^o/S}(K):=\Cone(p_u\circ\ad(i^*,i_*)(-)\circ\ad(\pi^*,\pi_*)(K): \\
Rj_*K\to i_*R(j\circ\pi)_*\pi^*K=:\psi_DK\to\psi_D^uK)\in D_k(S_{\mathbb C}^{an}).
\end{eqnarray*}
where $p_u:\psi_DK\to\psi_D^uK$ is the projection, note that $p_t\circ\ad(i^*,i_*)(-)\circ\ad(\pi^*,\pi_*)(K)=0$
where $p_t:\psi_DK\to\psi_D^tK$ is the projection.
\end{itemize}
\end{defi}

We will use a version of a result of Beilison on perverse sheaves.

\begin{defi}
Let $k\subset\mathbb C$ a subfield. Let $S\in\Var(k)$ and $D=V(s)\subset S$ a Cartier divisor. 
Denote $S^o:=S\backslash D$. We denote by $P_k(S_{\mathbb C}^{o,an})\times_JP_k(D_{\mathbb C}^{an})$ the category 
\begin{itemize}
\item whose objects are $(K',K'',u,v)$ where 
$K'\in P_k(S_{\mathbb C}^{o,an})$ and $K''\in P_k(D_{\mathbb C}^{an})$ are perverse sheaves
and $u:\psi^u_DK^{'an}\to K^{''an}$ and $v:K^{''an}\to\psi^u_DK^{'an}$ are morphism in $D_c(D_{\mathbb C}^{an})$ 
such that $v\circ u=T-I$, 
\item whose morphisms are $m=(m',m''):(K'_1,K''_1,u_1,v_1)\to(K'_2,K''_2,u_2,v_2)$ 
such that $u_2\circ\psi^u_Dm'=m''\circ u_1$ and $\psi^u_Dm'\circ v_1=v_2\circ m''$.
\end{itemize}
\end{defi}

We give the following version of the well known theorem for perserse sheaves

\begin{thm}\label{PSk0}
Let $k\subset\mathbb C$ a subfield. Let $S\in\Var(k)$ and $D=V(s)\subset S$ a Cartier divisor. 
Denote $i:D\hookrightarrow S$ the closed embedding and $j:S^o:=S\backslash D\hookrightarrow S$ the open embedding.
Let $\pi:\tilde S_{\mathbb C}^{o,an}\to S_{\mathbb C}^{o,an}$ the universal covering.
Denote by $T:\tilde S_{\mathbb C}^{o,an}\to\tilde S_{\mathbb C}^{o,an}$ the monodromy automorphism.
Then the functor
\begin{eqnarray*}
(j^*,\phi^u_D[-1],can,var):P_k(S_{\mathbb C}^{an})\to P_k(S_{\mathbb C}^{o,an})\times_JP_k(D_{\mathbb C}^{an})
\end{eqnarray*}
is an equivalence of category whose inverse is
\begin{eqnarray*}
P_k(S_{\mathbb C}^{o,an})\times_JP_k(D_{\mathbb C}^{an})\to P_k(S_{\mathbb C}^{an}), \\ 
(K',K'',u,v)\mapsto 
H^1(\psi^u_DK'\xrightarrow{(c(x_{S^o/S}(K')),u)}x_{S^o/S}(K')\oplus i_*K''\xrightarrow{((0,T-I),v)}\psi^u_DK').
\end{eqnarray*}
We denote, for $K\in P_k(S_{\mathbb C}^{an})$ by
\begin{eqnarray*}
Is(K):=(0,(\ad(j^*,j_*)(K),\ad(j\circ\pi^*,j\circ\pi_*)(K)),0): \\
K\xrightarrow{\sim}
(\psi^u_DK\xrightarrow{(c(x_{S^o/S}(K)),can(K))}x_{S^o/S}(K)\oplus i_*\phi^u_DK
\xrightarrow{((0,T-I),var(K):=(0,T-I))}\psi^u_DK)[-1]
\end{eqnarray*}
the canonical isomorphism in $D_{c,k}(S_{\mathbb C}^{an})$.
\end{thm}

\begin{proof}
Similar to the proof of \cite{Be2} : follows from the fact that $P_k(S_{\mathbb C}^{an})$
form an abelian category stable by the nearby and vanishing cycle functors.
\end{proof}

In the filtered case, we will consider the weight monodromy filtration 
for open embeddings and nearby and vanishing cycle functors :

\begin{defi}
Let $k\subset\mathbb C$ a subfield.
Let $S\in\Var(k)$ and $j:S^o\hookrightarrow S$ an open embedding such that 
$D:=S\backslash S^o=V(s)\subset S$ is a Cartier divisor. We then consider
\begin{equation*}
P_{fil,k}(S_{\mathbb C}^{o,an})^{ad,D}\subset P_{fil,k}(S_{\mathbb C}^{o,an})
\end{equation*}
the full subcategory consiting of $(K,W)\in P_{fil,k}(S_{\mathbb C}^{o,an})$
such that the relative weight monodromy filtration of $W$ with respect to $D\subset S$ exists,
in particular, for $(K,W)\in P_{fil,k}(S_{\mathbb C}^{o,an})^{ad,D}$ we have
$T(\psi_DW_kK)\subset\psi_DW_kK$ for all $k\in\mathbb Z$ and $T:\psi_DK\to\psi_D K$ is quasi-unipotent.
\end{defi}

\begin{defi}\label{jw}
Let $k\subset\mathbb C$ a subfield.
\begin{itemize}
\item[(i)]Let $S\in\Var(k)$ and $j:S^o\hookrightarrow S$ an open embedding such that 
$D:=S\backslash S^o=V(s)\subset S$ is a Cartier divisor.
\begin{itemize}
\item For $(K,W)\in P_{fil,k}(S_{\mathbb C}^{o,an})^{ad,D}$, we consider as in \cite{Saito}
\begin{equation*}
j_{*w}(K,W):=(Rj_*K,W)\in P_{fil,k}(S_{\mathbb C}^{an}), \; W_kRj_*K:=<Rj_*W_kK,W(N)_kK>\subset Rj_*K 
\end{equation*}
so that $j^*j_{*w}(K,W)=(K,W)$,
where $W_kRj_*K\subset Rj_*K$ is given by $W$ and the weight monodromy filtration $W(N)$ of 
the universal cover $\pi:\tilde S_{\mathbb C}^{o,an}\to S_{\mathbb C}^{o,an}$.
Note that a stratitification of $W_kRj_*K$ is given by the closure of a stratification of $W_kK$ 
and $D:=S\backslash S^o$. 
\item For $(K,W)\in P_{fil,k}(S_{\mathbb C}^{o,an})^{ad,D}$, we consider
\begin{equation*}
j_{!w}(K,W):=\mathbb D_S^vj_{*w}\mathbb D_S^v(K,W)\in P_{fil,k}(S_{\mathbb C}^{an}) 
\end{equation*}
so that $j^*j_{!w}(K,W)=(K,W)$.
\end{itemize}
For $(K',W)\in P_{fil,k}(S_{\mathbb C}^{an})^{ad,D}$, there is, by construction,
\begin{itemize}
\item a canonical map $\ad(j^*,j_{*w})(K',W)=\ad(j^*,j_*)(K'):(K',W)\to j_{*w}j^*(K',W)$ 
in $P_{fil,k}(S_{\mathbb C}^{an})$,
\item a canonical map $\ad(j_{!w},j^*)(K',W)=\ad(j_!,j^*)(K'):j_{!w}j^*(K',W)\to(K',W)$ 
in $P_{fil,k}(S_{\mathbb C}^{an})$.
\end{itemize}
\item[(ii)] Let $S\in\Var(k)$. Let $j:S^o:=S\backslash Z\hookrightarrow S$ an open embedding with
$Z=V(\mathcal I)\subset S$ an arbitrary closed subset, $\mathcal I\subset O_S$ being an ideal subsheaf. 
Taking generators $\mathcal I=(s_1,\ldots,s_r)$, we get $Z=V(s_1,\ldots,s_r)=\cap^r_{i=1}Z_i\subset S$ with 
$Z_i=V(s_i)\subset S$, $s_i\in\Gamma(S,\mathcal L_i)$ and $L_i$ a line bundle. 
Note that $Z$ is an arbitrary closed subset, $d_Z\geq d_X-r$ needing not be a complete intersection. Denote by 
$j_I:S^{o,I}:=\cap_{i\in I}(S\backslash Z_i)=S\backslash(\cup_{i\in I}Z_i)\xrightarrow{j_I^o}S^o\xrightarrow{j} S$ 
the open complementary embeddings, where $I\subset\left\{1,\cdots,r\right\}$. Denote 
\begin{equation*}
\mathcal D(Z/S):=\left\{(Z_i)_{i\in[1,\ldots r]},Z_i\subset S,\cap Z_i=Z\right\},Z'_i\subset Z_i
\end{equation*}
the flag category.
Let $P_{fil,k}(S_{\mathbb C}^{o,an})^{ad,(Z_i)}\subset P_{fil,k}(S_{\mathbb C}^{o,an})$
the full subcategory such that the relative weight monodromy filtration of $W$ with respect to the $Z_i\subset S$ exists.
For $(K,W)\in C(P_{fil,k}(S_{\mathbb C}^{o,an}))^{ad,(Z_i)}$, we define by (i) 
\begin{itemize}
\item the (bi)-filtered complex of $D_S$-modules
\begin{equation*}
j_{*w}(K,W):=\varinjlim_{\mathcal D(Z/S)}\Tot_{card I=\bullet}(j_{I*w}j_I^{o*}(K,W))\in C(P_{fil,k}(S_{\mathbb C}^{an})) 
\end{equation*}
where the horizontal differential are given by, 
if $I\subset J$, $d_{IJ}:=\ad(j^*_{IJ},j_{IJ*w})(j_I^{o*}(K,W))$, 
$j_{IJ}:S^{oJ}\hookrightarrow S^{oI}$ being the open embedding, 
and $d_{IJ}=0$ if $I\notin J$,
\item the (bi)-filtered complex of $D_S$-modules
\begin{eqnarray*}
j_{!w}(K,W):=\varprojlim_{\mathcal D(Z/S)}\Tot_{card I=-\bullet}(j_{I!w}j_I^{o*}(K,W)) 
=\mathbb D_S^vj_{*w}\mathbb D_S^v(K,W)\in C(P_{fil,k}(S_{\mathbb C}^{an})),
\end{eqnarray*}
where the horizontal differential are given by, 
if $I\subset J$, $d_{IJ}:=\ad(j_{IJ!w},j^*_{IJ})(j_I^{o*}(K,W))$, 
$j_{IJ}:S^{oJ}\hookrightarrow S^{oI}$ being the open embedding, and $d_{IJ}=0$ if $I\notin J$.
\end{itemize}
By definition, we have for $(K,W)\in C(P_{fil,k}(S_{\mathbb C}^{o,an}))^{ad,(Z_i)}$, 
$j^*j_{*w}(K,W)=(K,W)$ and $j^*j_{!w}(K,W)=(K,W)$.
For $(K',W)\in C(P_{fil,k}(S_{\mathbb C}^{an}))^{ad,(Z_i)}$, there is, by (i),
\begin{itemize}
\item a canonical map $\ad(j^*,j_{*w})(K',W):(K',W)\to j_{*w}j^*(K',W)$ in $C(P_{fil,k}(S_{\mathbb C}^{an}))$,
\item a canonical map $\ad(j_{!w},j^*)(K',W):j_{!w}j^*(K',W)\to(K',W)$ in $C(P_{fil,k}(S_{\mathbb C}^{an}))$.
\end{itemize}
\end{itemize}
\end{defi}

\begin{defi}\label{gammaw}
Let $S\in\Var(k)$. Let $Z\subset S$ a closed subset.
Denote by $j:S\backslash Z\hookrightarrow S$ the complementary open embedding. 
\begin{itemize}
\item[(i)] We define using definition \ref{jw}, the filtered Hodge support section functor
\begin{eqnarray*}
\Gamma^w_Z:C(P_{fil,k}(S_{\mathbb C}^{an})^{ad,(Z_i)})\to C(P_{fil,k}(S_{\mathbb C}^{an})), \\ 
(K,W)\mapsto\Gamma^w_Z(K,W):=\Cone(\ad(j^*,j_{*w})(K,W):(K,W)\to j_{*w}j^*(K,W))[-1],
\end{eqnarray*}
together we the canonical map $\gamma^w_Z(K,W):\Gamma^w_Z(K,W)\to (K,W)$.
\item[(i)'] Since $j_{*w}:C(P_{fil,k}(S_{\mathbb C}^{o,an})^{ad,(Z_i)})\to C(P_{fil,k}(S_{\mathbb C}^{an}))$ 
is an exact functor, $\Gamma^w_Z$ induces the functor
\begin{eqnarray*}
\Gamma^w_Z:D_{fil,c,k}(S_{\mathbb C}^{an})^{ad,(Z_i)}\to D_{fil,c,k}(S_{\mathbb C}^{an}), \; 
(K,W)\mapsto\Gamma^w_Z(K,W)
\end{eqnarray*}
\item[(ii)] We define using definition \ref{jw}, the dual filtered Hodge support section functor
\begin{eqnarray*}
\Gamma^{\vee,w}_Z:C(P_{fil,k}(S_{\mathbb C}^{an})^{ad,(Z_i)})\to C(P_{fil,k}(S_{\mathbb C}^{an})), \\ 
(K,W)\mapsto\Gamma^{\vee,w}_Z(K,W):=\Cone(\ad(j_{!w},j^*)(K,W):j_{!w},j^*(K,W)\to (K,W)),
\end{eqnarray*}
together we the canonical map $\gamma^{\vee,Hdg}_Z(K,W):(K,W)\to\Gamma_Z^{\vee,w}(K,W)$.
\item[(ii)'] Since $j_{!w}:C(P_{fil,k}(S_{\mathbb C}^{o,an})^{ad,(Z_i)})\to C(P_{fil,k}(S_{\mathbb C}^{an}))$ 
is an exact functor, $\Gamma^{\vee,w}_Z$ induces the functor
\begin{eqnarray*}
\Gamma^{\vee,w}_Z:D_{fil,c,k}(S_{\mathbb C}^{an})^{ad,(Z_i)}\to D_{fil,c,k}(S_{\mathbb C}^{an}), \; 
(K,W)\mapsto\Gamma^{\vee,w}_Z(K,W)
\end{eqnarray*}
\end{itemize}
\end{defi}

Let $S\in\Var(k)$ and $D=V(s)\subset S$ a Cartier divisor. 
Denote $i:D\hookrightarrow S$ the closed embedding and $j:S^o:=S\backslash D\hookrightarrow S$ the open embedding.
Let $\pi:\tilde S_{\mathbb C}^{o,an}\to S_{\mathbb C}^{o,an}$ the universal covering.
We then consider, for $(K,W)\in D_{fil,c}(S_{\mathbb C}^{o,an})^{ad,D}=\Ho(C(P_{fil,k}(S_{\mathbb C}^{o,an})^{ad,D})$, 
\begin{itemize}
\item the filtered nearby cycle functor 
\begin{equation*}
\psi_D(K,W):=(\psi_DK,W)\in D_{fil,c}(D_{\mathbb C}^{an}), \; W_k(\psi_D(K,W)):=<W_k\psi_DK,W(N)_k\psi_DK>\subset\psi_DK,
\end{equation*}
\item the vanishing cycle functor 
\begin{equation*}
\phi_D(K,W):=\Cone(\ad(j\circ\pi^*,j\circ\pi_*)(K,W):i^*(K,W)\to\psi_D(K,W))\in D_{fil,c,k}(D_{\mathbb C}^{an}),
\end{equation*}
where the morphism
\begin{equation*}
\ad(j\circ\pi^*,j\circ\pi_*)(K,W):i^*(K,W)\to i^*R(j\circ\pi)^*(j\circ\pi)^*(K,W)
\end{equation*}
being compatible with the weight monodromy filtration induces the morphism
\begin{equation*}
\ad(j\circ\pi^*,j\circ\pi_*)(K,W):i^*(K,W)\to\psi_D(K,W),
\end{equation*}
\item the canonical morphisms in $D_{fil,c,k}(D_{\mathbb C}^{an})$
\begin{eqnarray*}
can(K,W):=c(\phi_D(K,W)):\psi_D(K,W)\to\phi_D(K,W), \\ 
var(K,W):=(0,T-I):\phi_D(K,W)\to\psi_D(K,W).
\end{eqnarray*}
\end{itemize}

\begin{defi}\label{ukCw}
Let $S\in\Var(k)$ and $D=V(s)\subset S$ a Cartier divisor. 
Denote $i:D\hookrightarrow S$ the closed embedding and $j:S^o:=S\backslash D\hookrightarrow S$ the open embedding.
Let $\pi:\tilde S_{\mathbb C}^{o,an}\to S_{\mathbb C}^{o,an}$ the universal covering.
For $(K,W)\in P_{fil,k}(S_{\mathbb C}^{o,an})^{ad,D}$, we consider 
\begin{itemize}
\item the canonical decomposition
\begin{equation*}
\psi_D(K,W):=\psi_D^u(K,W)\oplus\psi_D^t(K,W)\in P_{fil,k}(D_{\mathbb C}^{an}), 
\end{equation*} 
\item the canonical decomposition
\begin{equation*}
\phi_D(K,W):=\phi_D^u(K,W)\oplus\phi_D^t(K,W)\in P_{fil,k}(D_{\mathbb C}^{an}), 
\end{equation*} 
\item the canonical morphisms in $D_{fil}(D_{\mathbb C}^{an})$
\begin{eqnarray*}
can(K,W):=c(\phi^u_D(K,W)):\psi^u_D(K,W)\to\phi^u_D(K,W), \\ 
var(K,W):=(0,T-I):\phi^u_D(K,W)\to\psi^u_D(K,W).
\end{eqnarray*}
\item the maximal extension
\begin{eqnarray*}
x_{S^o/S}(K,W):=\Cone(p_u\circ\ad(i^*,i_*)(-)\circ\ad(\pi^*,\pi_*)(K,W): \\
Rj_*(K,W)\to\psi_D(K,W)\to\psi_D^u(K,W))\in D_{fil,k}(S_{\mathbb C}^{an}).
\end{eqnarray*}
where the morphism
\begin{equation*}
\ad(i^*,i_*)(-)\circ\ad(\pi^*,\pi_*)(K,W):j_*(K,W)\to i_*i^*R(j\circ\pi)^*(j\circ\pi)^*(K,W)
\end{equation*}
being compatible with the weight monodromy filtration induces the morphism
\begin{equation*}
\ad(i^*,i_*)(-)\circ\ad(\pi^*,\pi_*)(K,W):j_{*w}(K,W)\to\psi_D(K,W).
\end{equation*}
\end{itemize}
\end{defi}

\begin{defi}\label{fw}
Let $k\subset\mathbb C$ a subfield.
\begin{itemize}
\item[(i)]Let $f:X\to S$ a morphism with $S,X\in\Var(k)$. 
Consider the graph factorization $f:X\xrightarrow{l}X\times S\xrightarrow{p}S$ of $f$ 
where $l$ the the graph closed embedding and $p$ is the projection.
We have, using definition \ref{gammaw},
\begin{itemize}
\item the inverse image functor
\begin{eqnarray*}
f^{*w}:D_{fil,c,k}(S_{\mathbb C}^{an})^{ad,(\Gamma_{f,i})}\to D_{fil,c,k}(X_{\mathbb C}^{an}), \;
(K,W)\mapsto f^{*w}(K,W):=l^*\Gamma_X^{\vee,w}p^*(K,W)
\end{eqnarray*}
\item the exceptional inverse image functor
\begin{eqnarray*}
f^{!w}:D_{fil,c,k}(S_{\mathbb C}^{an})^{ad,(\Gamma_{f,i})}\to D_{fil,c,k}(X_{\mathbb C}^{an}), \;
(K,W)\mapsto f^{!w}(K,W):=l^*\Gamma_X^wp^*(K,W).
\end{eqnarray*}
\end{itemize}
\item[(ii)]Let $f:X\to S$ a morphism with $S,X\in\Var(k)$.
Consider a compactification $f:X\hookrightarrow{j}\bar X\xrightarrow{\bar f}S$ of $f$
with $\bar X\in\Var(k)$, $j$ an open embedding and $\bar f$ a proper morphism.
Denote $Z=\bar X\backslash X$. We have, using definition \ref{jw},
\begin{itemize}
\item the direct image functor
\begin{eqnarray*}
Rf_{*w}:D_{fil,c,k}(X_{\mathbb C}^{an})^{ad,(Z_i)}\to D_{fil,c,k}(S_{\mathbb C}^{an}), \;
(K,W)\mapsto Rf_{*w}(K,W):=R\bar f_*j_{*w}(K,W)
\end{eqnarray*}
\item the proper direct image functor
\begin{eqnarray*}
Rf_{!w}:D_{fil,c,k}(X_{\mathbb C}^{an})^{ad,(Z_i)}\to D_{fil,c,k}(S_{\mathbb C}^{an}), \;
(K,W)\mapsto Rf_{!w}(K,W):=R\bar f_*j_{!w}(K,W).
\end{eqnarray*}
\end{itemize}
\item[(iii)] Let $S\in\Var(k)$. Denote by $\Delta_S:S\hookrightarrow S\times S$ the diagonal closed embedding and
$p_1:S\times S\to S$, $p_2:S\times S\to S$ the projections. We have by (i) the functor
\begin{eqnarray*}
\otimes^{Lw}:D_{fil,c,k}(S_{\mathbb C}^{an})^{ad,(S_i)}\times D_{fil,c,k}(S_{\mathbb C}^{an})^{ad,(S_i)}
\to D_{fil,c,k}(S_{\mathbb C}^{an}), \\ 
((K_1,W),(K_2,W))\mapsto (K_1,W)\otimes^{L,w}(K_2,W):=\Delta_S^{!w}(p_1^*(K_1,W)\otimes^L p_2^*(K_2,W)).
\end{eqnarray*}
\end{itemize}
\end{defi}

Let $S\in\Var(k)$ and $D=V(s)\subset S$ a Cartier divisor. 
Denote $i:D\hookrightarrow S$ the closed embedding and $j:S^o:=S\backslash D\hookrightarrow S$ the open embedding.
In the filtered case, we get, for $(K,W)\in P_{fil,k}(S_{\mathbb C}^{an})^{ad,D}$ 
the map in $D_{fil,c,k}(S_{\mathbb C}^{an})$
\begin{eqnarray*}
Is(K,W):=(0,(\ad(j^*,j_*)(K),\ad(j\circ\pi^*,j\circ\pi_*)(K)),0): \\
(K,W)\to(\psi^u_D(K,W)\xrightarrow{(c(x_{S^o/S}(K,W)),can(K,W))}x_{S^o/S}(K,W)\oplus i_*\phi^u_D(K,W) \\
\xrightarrow{((0,T-I),var(K,W):=(0,T-I))}\psi^u_D(K,W))[-1]
\end{eqnarray*}
which is NOT an isomorphism in general (it leads to different $W$-filtration on perverse cohomology).

\subsection{Constructible and perverse etale sheaves on algebraic varieties over a field $k$ of charactersitic $0$} 

Let $k$ a field of caracteristic zero. Let $S\in\Var(k)$. We have
\begin{itemize}
\item the classical dual functor
\begin{eqnarray*}
\mathbb D^0_S:C(S^{et})\to C(S^{et}), \; \; K\mapsto\mathbb D_S^0K:=\mathcal Hom(LK,E_{et}(\mathbb Z_S))
\end{eqnarray*}
which induces in the derived category
\begin{eqnarray*}
\mathbb D^0_S:D(S^{et})\to D(S^{et}), \; \;
K\mapsto\mathbb D_S^0K:=\mathcal Hom(LK,E_{et}(\mathbb Z_S))=R\mathcal Hom(K,\mathbb Z_S)
\end{eqnarray*}
\item the Verdier dual functor
\begin{eqnarray*}
\mathbb D^v_S:D(S^{et})\to D(S^{et}), \; \; K\mapsto\mathbb D_S^vK:=R\mathcal Hom(K,w_S)
\end{eqnarray*}
where $w_S:=a_S^!\mathbb Z_S$ is the dualizing complex, $a_S:S\to\left\{\pt\right\}$ being the terminal map.
\end{itemize}
For $S\in\SmVar(k)$ smooth connected of dimension $\dim(S)=d_S$ we have $\mathbb D^v_S=\mathbb D^0_S[2d_S]$.

We recall the definition of constructible etale sheaves on algebraic varieties over a field $k$ of caracteristic zero:

\begin{defi}\label{constr2}
Let $k$ a field of caracteristic zero. Let $S\in\Var(k)$. Let $l$ a prime number.
Recall (see section 2.1) that 
\begin{equation*}
K=(K_n)_{n\in\mathbb N}\in\Shv_{\mathbb Z_l}(S^{et})\subset\PSh(S^{et},\Fun(\mathbb N,\Ab))
\end{equation*}
is a projective system with $K_n\in\Shv_{\mathbb Z/l^n\mathbb Z}(S^{et})$ 
such that $K_n\to K_{n+1}/l^nK_{n+1}$ is an isomorphism.
\begin{itemize}
\item[(i)]A sheaf $K\in\Shv_{\mathbb Z_l}(S^{et})$ 
is called constructible if there exists a stratification $S=\sqcup_{\alpha} S_{\alpha}$ 
with $l_{\alpha}:S_{\alpha}\hookrightarrow S$ locally closed subsets such that 
$l_{\alpha}^*K\in\Shv_{\mathbb Z_l}(S_{\alpha}^{et})$ 
are (finite dimensional) local systems (for the etale topology) for all $\alpha$. 
\item[(ii)]We denote by 
\begin{equation*}
C_{\mathbb Z_l,c,k}(S^{et})\subset C_{\mathbb Z_l}(S^{et}) \; \mbox{and} \; 
D_{\mathbb Z_l,c,k}(S^{et})\subset D_{\mathbb Z_l}(S^{et})
\end{equation*}
the full subcategories consisting of $K\in C_{\mathbb Z_l}(S^{et})$ 
such that $a_{et}H^nK\in\Shv_{\mathbb Z_l}(S^{et})$ are constructible for all $n\in\mathbb Z$ (see (i)).
\item[(ii)']We denote by 
\begin{equation*}
C_{\mathbb Z_lfil,c,k}(S^{et})\subset C_{\mathbb Z_lfil}(S^{et}) \; \mbox{and} \; 
D_{\mathbb Z_lfil,c,k}(S^{et})\subset D_{\mathbb Z_lfil}(S^{et})
\end{equation*}
the full subcategories consisting of $(K,W)\in C_{\mathbb Z_lfil}(S^{et})$ 
such that $a_{et}\Gr_k^WH^nK\in\Shv_{\mathbb Z_l}(S^{et})$ are constructible for all $n,k\in\mathbb Z$ (see (i)).
\item[(iii)] We denote by 
\begin{equation*}
P_{\mathbb Z_l,k}(S^{et})\subset D_{\mathbb Z_l,c,k}(S^{et})\otimes_{\mathbb Z_l}\mathbb Q_l 
\subset D_{\mathbb Z_l,c,k}(S^{et})
\end{equation*}
the full subcategory of perverse sheaves (which are by definition torsion free).
\item[(iii)'] We denote by 
\begin{equation*}
P_{\mathbb Z_l,fil,k}(S^{et})\subset D_{\mathbb Z_lfil,c,k}(S^{et})\otimes_{\mathbb Z_l}\mathbb Q_l 
\subset D_{\mathbb Z_lfil,c,k}(S^{et})
\end{equation*}
the full subcategory of filtered perverse sheaves (which are by definition torsion free).
\end{itemize}
Let $k/k'$ a field extention. Let $S\in\Var(k)$. Let $l$ a prime number.
\begin{itemize}
\item[(i)]A sheaf $K\in\Shv_{\mathbb Z_l}(S_{k'}^{et})$ 
is called constructible (over k) if there exists a stratification $S=\sqcup_{\alpha} S_{\alpha}$ 
with $l_{\alpha}:S_{\alpha}\hookrightarrow S$ locally closed subsets such that 
$l_{\alpha}^*K\in\Shv_{\mathbb Z_l}(S_{\alpha}^{et})$ 
are (finite dimensional) local systems (for the etale topology) for all $\alpha$. 
\item[(ii)]We denote by 
\begin{equation*}
C_{\mathbb Z_l,c,k}(S_{k'}^{et})\subset C_{\mathbb Z_l}(S_{k'}^{et}) \; \mbox{and} \; 
D_{\mathbb Z_l,c,k}(S_{k'}^{et})\subset D_{\mathbb Z_l}(S_{k'}^{et})
\end{equation*}
the full subcategories consisting of $K\in C_{\mathbb Z_l}(S_{k'}^{et})$ 
such that $a_{et}H^nK\in\Shv_{\mathbb Z_l}(S_{k'}^{et})$ are constructible over $k$ for all $n\in\mathbb Z$ (see (i)).
\item[(ii)']We denote by 
\begin{equation*}
C_{\mathbb Z_lfil,c,k}(S_{k'}^{et})\subset C_{\mathbb Z_lfil}(S_{k'}^{et}) \; \mbox{and} \; 
D_{\mathbb Z_lfil,c,k}(S_{k'}^{et})\subset D_{\mathbb Z_lfil}(S_{k'}^{et})
\end{equation*}
the full subcategories consisting of $(K,W)\in C_{\mathbb Z_lfil}(S_{k'}^{et})$ 
such that $a_{et}\Gr_k^WH^nK\in\Shv_{\mathbb Z_l}(S_{k'}^{et})$ are constructible over $k$ for all $n,k\in\mathbb Z$ (see (i)).
\item[(iii)] We denote by 
\begin{equation*}
P_{\mathbb Z_l,k}(S_{k'}^{et}):=P_{\mathbb Z_l,k'}(S_{k'}^{et})\cap D_{\mathbb Z_l,c,k}(S_{k'}^{et})
\subset D_{\mathbb Z_l,c,k}(S_{k'}^{et}) 
\end{equation*}
the full subcategory of perverse sheaves whose stratification is defined over $k$.
\item[(iii)'] We denote by 
\begin{equation*}
P_{\mathbb Z_l,fil,k}(S_{k'}^{et}):=P_{\mathbb Z_l,fil,k'}(S_{k'}^{et})\cap D_{\mathbb Z_lfil,c,k}(S_{k'}^{et})
\subset D_{\mathbb Z_lfil,c,k}(S_{k'}^{et}) 
\end{equation*}
the full subcategory of filtered perverse sheaves whose stratification is defined over $k$.
\end{itemize}
\end{defi}

\begin{thm}\label{constrthm2}
Let $k$ a field of caracteristic zero. Let $l$ a prime number.
\begin{itemize}
\item Let $S\in\Var(k)$. Then for $K\in D_{\mathbb Z_l,c,k}(S^{et})$, $\mathbb D^v_SK\in D_{\mathbb Z_l,c,k}(S^{et})$.
For $K,K'\in D_{\mathbb Z_l,c,k}(S^{et})$, $K\otimes^L K'\in D_{\mathbb Z_l,c,k}(S^{et})$.
\item Let $f:X\to S$ a morphism with $S,X\in\Var(k)$.
Then for $K\in D_{\mathbb Z_l,c,k}(X^{et})$, $Rf_*K\in D_{\mathbb Z_l,c,k}(S^{et})$ 
and $Rf_!K=\mathbb D^v_SRf_*\mathbb D^v_XK\in D_{\mathbb Z_l,c,k}(S^{et})$.
\item Let $f:X\to S$ a morphism with $S,X\in\Var(k)$.
Then for $K\in D_{\mathbb Z_l,c,k}(S^{et})$, $f^*K\in D_{\mathbb Z_l,c,k}(X^{et})$ 
and $f^!K:=\mathbb D^v_Xf^*\mathbb D^v_SK\in D_{\mathbb Z_l,c,k}(X^{et})$.
\end{itemize}
We have then, for $S\in\Var(k)$, the full subcategory
\begin{eqnarray*}
D_{\mathbb Z_l,c,k,gm}(S^{et}):&=&<Rf_*\mathbb Z_{X^{et}},(f:X\to S)\in\Var(k)> \\
&=&<Rf_*\mathbb Z_{X^{et}},(f:X\to S)\in\Var(k) \; \mbox{proper}, \; X \mbox{smooth}>\subset D_{\mathbb Z_l,c,k}(S^{et})
\end{eqnarray*}
where $<>$ means the full triangulated category generated by.
\end{thm}

\begin{proof}
Standard : see \cite{milnes} for example.
\end{proof}

\begin{defi}
We have then, using theorem \ref{constrthm2}, for $S\in\Var(k)$, the full subcategory
\begin{eqnarray*}
D_{\mathbb Z_lfil,c,k,gm}(S^{et}):=<(K,W), \, \mbox{s.t., for all} \, n\in\mathbb Z
\Gr_n^WK=Rf_{n*}\mathbb Z_{X_n^{et}}, \, X_n\in\SmVar(k), \,  f_n:X\to S \, \mbox{proper}>
\subset D_{\mathbb Z_lfil,c,k}(S^{et})
\end{eqnarray*}
where $<>$ means the full triangulated category generated by.
\end{defi}

We now give, following Ayoub, the definition of the nearby and vanishing cycle functors for constructible etale sheaves :

Let $k$ a field of characteristic zero. Let $l$ a prime number.
Let $S\in\Var(k)$ and $D=V(s)\subset S$ a Cartier divisor with $s\in\Gamma(S,L)$. 
Denote $i:D\hookrightarrow S$ the closed embedding, $j:S^o:=S\backslash D\hookrightarrow S$ the open embedding,
so that $s$ induces locally for the Zariski topology on $S$ a morphism $s:S^o\to\mathbb G_m$. Let 
\begin{equation*}
\pi:\tilde S^o:=\varinjlim_{n\in\mathbb N}\tilde S^o_n:=\Spec(L[t]/(t^n-s))\to S^o 
\end{equation*}
the universal covering given by taking the $n$-th roots of $s$ 
and $\mathcal A^{\bullet}_D(S^o):=((S^o\times_{S^o\times S^o}S^o)/\mathbb G_m)\in\Fun(\Delta^{\bullet},\Var(k)^{sm}/S^o)$
the diagram of lattices (see \cite{AyoubT}). We then consider, 
\begin{itemize}
\item for $K\in C_{\mathbb Z_l}(S^{o,et})$, the nearby cycle functor 
\begin{eqnarray*}
\psi_DK:&=&i^*e(S^{et})_*R(j\circ\pi)_*\pi^*\mathcal Hom(\mathcal A^{\bullet}_D(S^o),e(S^{et})^*K) \\
&=&i^*R(j\circ\pi)_*\pi^*\mathcal Hom(\mathcal A^{\bullet,e}_D(S^o),K)\in D_{\mathbb Z_l}(D^{et}), 
\end{eqnarray*}
together with the monodromy morphism $T:\psi_DK\to\psi_DK$ in $D_{\mathbb Z_l}(D^{et})$, where 
\begin{equation*}
\mathcal A^{\bullet,e}_D(S^o):=(\cdots\to Rp_{S^o!}\mathbb Z_{l,S^o\times S^o}\to\cdots)\in D_{\mathbb Z_l}(S^{o,et}) 
\end{equation*}
and we write again $\psi_DK:=i_*\psi_DK\in D_{\mathbb Z_l,c}(S^{et})$,
note that by adjonction $(Rp_{S^o!},p_{S^o}^*)$ we have $e(S^{et})^*\psi_DK=\psi_De(S^{et})^*K$,
\item for $K\in C_{\mathbb Z_l}(S^{o,et})$, the vanishing cycle functor 
\begin{equation*}
\phi_DK:=\Cone(\ad(j\circ\pi^*,j\circ\pi_*)(-)\circ ev(K):i^*K\to\psi_D K\in D_{\mathbb Z_l}(D^{et})
\end{equation*}
together with the canonical map $c(\phi_DK):\psi_DK\to\phi_DK$ in $D_{\mathbb Z_l}(D^{et})$,
we write again $\phi_DK:=i_*\phi_DK\in D_{\mathbb Z_l}(S^{et})$,
\item for $K\in C_{\mathbb Z_l,c,k}(S^{o,et})$, the canonical morphisms in $D_{\mathbb Z_l,c,k}(D^{et})$
\begin{eqnarray*}
can(K):=c(\phi_DK):\psi_DK\to\phi_DK \; , \; var(K):=((0,T-I)):\phi_DK\to\psi_DK.
\end{eqnarray*}
\end{itemize}

Let $S\in\Var(k)$ and $D=V(s)\subset S$ a Cartier divisor with $s\in\Gamma(S,L)$. 
Denote $i:D\hookrightarrow S$ the closed embedding and $j:S^o:=S\backslash D\hookrightarrow S$ the open embedding. Let 
\begin{equation*}
\pi:\tilde S^o:=\varinjlim_{n\in\mathbb N} \tilde S^o_n:=\Spec(L[t]/(t^n-s))\to S^o 
\end{equation*}
the universal covering
and $\mathcal A^{\bullet}_D(S^o):=((S^o\times_{S^o\times S^o}S^o)/\mathbb G_m)\in\Fun(\Delta^{\bullet},\Var(k)^{sm}/S^o)$
the diagram of lattices (see \cite{AyoubT}).
\begin{itemize}
\item For $K,K'\in C_{\mathbb Z_l,c,k}(S^{o,et})$, we have by (Verdier) duality and theorem \ref{KUNetth}(ii)
 $\psi_D(K\otimes^L K')=\psi_D(K)\otimes^L\psi_D(K')$,
\item For $K,K'\in C_{\mathbb Z_l,c,k}(S^{o,et})$ we have by (Verdier) duality and the preceding point 
$\phi_D(K\otimes^L K')=\phi_D(K)\otimes^L\phi_D(K')$
\item For $K,K'\in C_{\mathbb Z_l}(S^{o,et})$, we have the transformation map in $D_{\mathbb Z_l}(S^{o,et})$
\begin{eqnarray*}
T(\otimes,\psi_D)(K,K'):(\psi_DK)\otimes^L K'\xrightarrow{I\otimes(\ad(i^*,i_*)(-)\circ\ad(\pi^*,\pi_*)(-))}
\psi_DK\otimes^L\psi_DK'=\psi_D(K\otimes^L K').
\end{eqnarray*}
\item For $K,K'\in C_{\mathbb Z_l,c,k}(S^{o,et})$ the transformation map in $D_{\mathbb Z_l,c,k}(S^{o,et})$
\begin{eqnarray*}
T(\otimes,\phi_D)(K,K'):\phi_D(K\otimes^L K')=\mathbb D_S^v\psi_D(\mathbb D_S^vK\otimes^L\mathbb D_S^vK') \\
\xrightarrow{\mathbb D_S^vT(\otimes,\psi_D)(\mathbb D_S^vK,\mathbb D_S^vK')}
\mathbb D_S^v(\psi_D\mathbb D_S^v(K)\otimes^L\mathbb D_S^vK')=\phi_DK\otimes^L K'.
\end{eqnarray*}
\item For $K\in C_{\mathbb Z_l,c,k}(S^{o,et})$ we have a canonical isomorphism in $D_{\mathbb Z_l,c,k}(D^{et})$
\begin{equation*}
T(D,\psi_D)(K):\phi_D\mathbb D_S^vK\xrightarrow{\sim}\mathbb D_S^v\psi_DK.
\end{equation*}
\end{itemize}

We have then the following :

\begin{prop}\label{phipsiperv}
Let $k$ a field of characteristic zero. Let $S\in\Var(k)$ and $D\subset S$ a Cartier divisor. 
For $K\in P_{\mathbb Z_l,k}(S^{et})$, we have $\psi_DK[-1],\phi_DK[-1]\in P_{\mathbb Z_l,k}(S^{et})$.
\end{prop}

\begin{proof}
Take an embedding $k\hookrightarrow\mathbb C$ and consider the morphism of site
\begin{equation*}
\an_S:S_{\mathbb C}^{an}=S_{\mathbb C}^{an,et}\to S_{\mathbb C}^{et}\xrightarrow{\pi_{k/\mathbb C}(S)} S^{et} 
\end{equation*}
given by the analytical functor.
By \cite{AyoubB}, we have a canonical isomorphism in $D_{c,k}(S_{\mathbb C}^{an})\otimes\mathbb Q$
\begin{equation*}
T(an,\psi):\an_S^*\psi_DK\xrightarrow{\sim}\psi_D\an_S^*K
\end{equation*}
where we recall (see section 2) that $D_{c,k}(S_{\mathbb C}^{an})\subset D(S_{\mathbb C}^{an})$
is the full subcategory consisting of classes of complexes of presheaves whose cohomology sheaves are
constructible with respect to a Zariski stratification of $S$ (in particular defined over $k$). Since 
\begin{equation*}
K=((K_n)_{n\in\mathbb N})\otimes\mathbb Q_p\in P_{\mathbb Z_l,k}(S^{et}), 
\end{equation*}
$\an_S^*K\in P_k(S_{\mathbb C}^{an})$, 
where we recall (see section 2) that $P_k(S_{\mathbb C}^{an})\subset D_{c,k}(S_{\mathbb C}^{an})\otimes\mathbb Q$
is the full subcategory consisting of presheaves sheaves whose cohomology sheaves are constructible
with respect to a Zariski stratification of $S$ (defined over $k$).
Thus by the complex case (see e.g. \cite{PS}) $\psi_D\an_S^*K[-1]\in P_k(S_{\mathbb C}^{an})$.
Hence $\an_S^*\psi_DK[-1]\in P_k(S_{\mathbb C}^{an})$ and thus $\psi_DK[-1],\phi_DK[-1]\in P_{\mathbb Z_l,k}(S^{et})$.
\end{proof}

\begin{defi}\label{uket}
Let $S\in\Var(k)$ and $D=V(s)\subset S$ a Cartier divisor. 
Denote $i:D\hookrightarrow S$ the closed embedding and $j:S^o:=S\backslash D\hookrightarrow S$ the open embedding.
Let $\pi:\tilde S^{o}\to S^{o}$ the universal covering. For $K\in P_{\mathbb Z_l,k}(S^{o,et})$, we consider 
\begin{itemize}
\item the canonical decomposition
\begin{equation*}
\psi_DK:=\psi_D^uK\oplus\psi_D^tK\in P_{\mathbb Z_l,k}(D^{et}), 
\end{equation*}
with $\psi_D^uK:=\ker(T-I)^N\subset\psi_DK$ with $N$ such that $\ker(T-I)^{N+1}=\ker(T-I)^N$
and $\psi_D^tK:=\ker Q(T)\subset\psi_DK$ where $P(X)=(X-1)^NQ(X)$ is the minimal polynomial of $T:\psi_DK\to\psi_DK$, 
note that we have $\psi^u_DK,\psi_D^tK\in P_{\mathbb Z_l,k}(D^{o,et})$ 
since $\psi^u_D$ and $\psi^t_D$ preserve local systems : if $(S_{\alpha})$ is a stratification for $K$ 
then $(S_{\alpha}\cap S^o,S_{\alpha}\cap D)$ is a stratification for $\psi^u_DK$ and $\psi^t_DK$,
\item the canonical decomposition
\begin{equation*}
\phi_DK:=\phi_D^uK\oplus\phi_D^tK\in P_{\mathbb Z_l,k}(D^{et}), 
\end{equation*}
with $\phi_D^uK:=\ker(T-I)^N\subset\psi_DK$ with $N$ such that $\ker(T-I)^{N+1}=\ker(T-I)^N$
and $\phi_D^tK:=\ker Q(T)\subset\phi_DK$ where $P(X)=(X-1)^NQ(X)$ is the minimal polynomial of $T:\phi_DK\to\phi_DK$, 
\item the canonical morphisms in $D_c(D^{et})$
\begin{eqnarray*}
can(K):=c(\phi^u_DK):\psi^u_DK\to\phi^u_DK \; , \; var(K):=(0,T-I):\phi^u_DK\to\psi^u_DK.
\end{eqnarray*}
\item the maximal extension
\begin{eqnarray*}
x_{S^o/S}(K):=\Cone(p_u\circ\ad(i^*,i_*)(-)\circ\ad(\pi^*,\pi_*)(-)\circ ev(K): \\
Rj_*K\to i_*R(j\circ\pi)_*\pi^*\mathcal Hom(\mathcal A^{\bullet}_D(S^o),K)=:\psi_DK\to\psi_D^uK)
\in P_{\mathbb Z_l,k}(S^{et}).
\end{eqnarray*}
where $p_u:\psi_DK\to\psi_D^uK$ is the projection, 
note that $p_t\circ\ad(i^*,i_*)(-)\circ\ad(\pi^*,\pi_*)(K)\circ ev(K)=0$ 
where $p_t:\psi_DK\to\psi_D^tK$ is the projection.
\end{itemize}
\end{defi}

We will use the etale version of a result of Beilinson on perverse sheaves.

\begin{defi}
Let $k$ a field of characteristic zero. Let $S\in\Var(k)$ and $D=V(s)\subset S$ a Cartier divisor. 
Denote $S^o:=S\backslash D$. We denote by $P_{\mathbb Z_l}(S^{o,et})\times_JP_{\mathbb Z_l}(D^{et})$ the category 
\begin{itemize}
\item whose objects are $(K',K'',u,v)$ 
where $K'\in P_{\mathbb Z_l}(S^{o,et})$ and $K''\in P_{\mathbb Z_l}(D^{et})$ are perverse sheaves
and $u:\psi^u_DK^{'an}\to K^{''an}$ and $v:K^{''an}\to\psi^u_DK^{'an}$ are morphism in $D_c(D^{et})$ such that
$v\circ u=T-I$, 
\item whose morphisms are $m=(m',m''):(K'_1,K''_1,u_1,v_1)\to(K'_2,K''_2,u_2,v_2)$ 
such that $u_2\circ\psi^u_Dm'=m''\circ u_1$ and $\psi^u_Dm'\circ v_1=v_2\circ m''$.
\end{itemize}
\end{defi}

We give the version for etale construcible sheaves of the well known theorem \ref{PSk0} for perserse sheaves

\begin{thm}\label{PSk}
Let $k$ a field of characteristic zero. Let $S\in\Var(k)$ and $D=V(s)\subset S$ a Cartier divisor. 
Denote $i:D\hookrightarrow S$ the closed embedding and $j:S^o:=S\backslash D\hookrightarrow S$ the open embedding.
Then the functor
\begin{eqnarray*}
(j^*,\phi^u_D[-1],can,var):P_{\mathbb Z_l,k}(S^{et})\to P_{\mathbb Z_l,k}(S^{o,et})\times_JP_{\mathbb Z_l,k}(D^{et})
\end{eqnarray*}
is an equivalence of category whose inverse is
\begin{eqnarray*}
P_{\mathbb Z_l,k}(S^{o,et})\times_JP_{\mathbb Z_l,k}(D^{et})\to P_{\mathbb Z_l,k}(S^{et}), \\ 
(K',K'',u,v)\mapsto H^1(\psi^u_DK'\xrightarrow{c(x_{S^o/S}(K')),u)}x_{S^o/S}(K')\oplus i_*K''
\xrightarrow{((0,T-I),v)}\psi^u_DK').
\end{eqnarray*}
We denote, for $K\in P_{\mathbb Z_l}(S^{et})$ by
\begin{eqnarray*}
Is(K):=(0,(\ad(j^*,j_*)(K),\ad((j\circ\pi)^*,(j\circ\pi)_*)(K)),0): \\
K\xrightarrow{\sim}
(\psi^u_DK\xrightarrow{(c(x_{S^o/S}(K)),can(K))}x_{S^o/S}(K)\oplus i_*\phi^u_DK
\xrightarrow{((0,T-I),var(K):=(0,T-I))}\psi^u_DK)[-1]
\end{eqnarray*}
the canonical isomorphism in $D_{\mathbb Z_l,c,k}(S^{et})$.
\end{thm}

\begin{proof}
Similar to the proof of theorem \ref{PSk0}. 
Follows from the fact that $P_{\mathbb Z_l,k}(S^{et})\subset D_{\mathbb Z_l,c,k}(S^{et})$ is a full abelian subcategory
(as the heart of the perverse t-structure)
which is stable by the functors $\phi_D$ and $\phi_D$ for $D\subset S$ a Cartier divisor by proposition \ref{phipsiperv}.
\end{proof}

In the filtered case, we will consider the weight monodromy filtration for open embeddings 
and nearby and vanishing cycle functors:

\begin{defi}
Let $k$ a field of characteristic zero. Let $l$ a prime number.
Let $S\in\Var(k)$ and $j:S^o\hookrightarrow S$ an open embedding such that 
$D:=S\backslash S^o=V(s)\subset S$ is a Cartier divisor. We consider
\begin{equation*}
P_{\mathbb Z_lfil,k}(S^{o,et})^{ad,D}\subset P_{\mathbb Z_lfil,k}(S^{o,et})
\end{equation*}
the full subcategory consisting of $(K,W)\in P_{\mathbb Z_lfil,k}(S^{o,et})$
the relative weight monodromy filtration of $W$ with respect to $D\subset S$ exists,
in particular, for $(K,W)\in P_{\mathbb Z_lfil,k}(S^{o,et})^{ad,D}$ we have
$T(\psi_DW_kK)\subset\psi_DW_kK$ for all $k\in\mathbb Z$ and $T:\psi_DK\to\psi_D K$ is quasi-unipotent.
\end{defi}

\begin{defi}\label{jwet}
Let $k$ a field of characteristic zero. Let $l$ a prime number.
\begin{itemize}
\item[(i)] Let $S\in\Var(k)$ and $j:S^o\hookrightarrow S$ an open embedding such that 
$D:=S\backslash S^o=V(s)\subset S$ is a Cartier divisor.
\begin{itemize}
\item For $(K,W)\in P_{\mathbb Z_lfil,k}(S^{o,et})^{ad,D}$, we consider as in \cite{Saito}
\begin{equation*}
j_{*w}(K,W):=(Rj_*K,W)\in P_{\mathbb Z_lfil,k}(S^{et}), \; W_kRj_*K:=<Rj_*W_kK,W(N)_kK>\subset Rj_*K 
\end{equation*}
so that $j^*j_{*w}(K,W)=(K,W)$,
where $W_kRj_*K\subset Rj_*K$ is given by $W$ and the weight monodromy filtration $W(N)$ of 
the universal cover $\pi:\tilde S^o\to S^o$, see \cite{AyoubT}.
Note that a stratitification of $W_kRj_*K$ is given by the closure of a stratification of $W_kK$ 
and $D:=S\backslash S^o$. 
\item For $(K,W)\in P_{\mathbb Z_l,fil,k}(S^{o,et})^{ad,D}$, we consider
\begin{equation*}
j_{!w}(K,W):=\mathbb D_S^vj_{*w}\mathbb D_S^v(K,W)\in P_{\mathbb Z_l,fil}(S^{et}) 
\end{equation*}
so that $j^*j_{!w}(K,W)=(K,W)$.
\end{itemize}
For $(K',W)\in P_{\mathbb Z_lfil,k}(S^{et})^{ad,D}$, there is, by construction,
\begin{itemize}
\item a canonical map $\ad(j^*,j_{*w})(K',W)=\ad(j^*,j_*)(K'):(K',W)\to j_{*w}j^*(K',W)$ in $P_{\mathbb Z_l,fil,k}(S^{et})$,
\item a canonical map $\ad(j_{!w},j^*)(K',W)=\ad(j_!,j^*)(K'):j_{!w}j^*(K',W)\to(K',W)$ in $P_{\mathbb Z_l,fil,k}(S^{et})$.
\end{itemize}
\item[(ii)] Let $S\in\Var(k)$. Let $j:S^o:=S\backslash Z\hookrightarrow S$ an open embedding with
$Z=V(\mathcal I)\subset S$ an arbitrary closed subset, $\mathcal I\subset O_S$ being an ideal subsheaf. 
Taking generators $\mathcal I=(s_1,\ldots,s_r)$, we get $Z=V(s_1,\ldots,s_r)=\cap^r_{i=1}Z_i\subset S$ with 
$Z_i=V(s_i)\subset S$, $s_i\in\Gamma(S,\mathcal L_i)$ and $L_i$ a line bundle. 
Note that $Z$ is an arbitrary closed subset, $d_Z\geq d_X-r$ needing not be a complete intersection. Denote by 
$j_I:S^{o,I}:=\cap_{i\in I}(S\backslash Z_i)=S\backslash(\cup_{i\in I}Z_i)\xrightarrow{j_I^o}S^o\xrightarrow{j} S$ 
the open complementary embeddings, where $I\subset\left\{1,\cdots,r\right\}$. Denote 
\begin{equation*}
\mathcal D(Z/S):=\left\{(Z_i)_{i\in[1,\ldots r]},Z_i\subset S,\cap Z_i=Z\right\},Z'_i\subset Z_i
\end{equation*}
the flag category. 
Let $P_{\mathbb Z_lfil,k}(S^{o,et})^{ad,(Z_i)}\subset P_{\mathbb Z_lfil,k}(S^{o,et})$
the full subcategory such that the relative weight monodromy filtration of $W$ with respect to the $Z_i\subset S$ exists.
For $(K,W)\in C(P_{\mathbb Z_l,fil,k}(S^{o,et})^{ad,(Z_i)})$, we define by (i) 
\begin{itemize}
\item the (bi)-filtered complex of $D_S$-modules
\begin{equation*}
j_{*w}(K,W):=\varinjlim_{\mathcal D(Z/S)}\Tot_{card I=\bullet}(j_{I*w}j_I^{o*}(K,W))\in C(P_{\mathbb Z_l,fil}(S^{et})) 
\end{equation*}
where the horizontal differential are given by, 
if $I\subset J$, $d_{IJ}:=\ad(j^*_{IJ},j_{IJ*w})(j_I^{o*}(K,W))$, 
$j_{IJ}:S^{oJ}\hookrightarrow S^{oI}$ being the open embedding, 
and $d_{IJ}=0$ if $I\notin J$,
\item the (bi)-filtered complex of $D_S$-modules
\begin{eqnarray*}
j_{!w}(K,W):=\varprojlim_{\mathcal D(Z/S)}\Tot_{card I=-\bullet}(j_{I!w}j_I^{o*}(K,W)) 
=\mathbb D_S^vj_{*w}\mathbb D_S^v(K,W)\in C(P_{\mathbb Z_l,fil}(S^{et})),
\end{eqnarray*}
where the horizontal differential are given by, 
if $I\subset J$, $d_{IJ}:=\ad(j_{IJ!w},j^*_{IJ})(j_I^{o*}(K,W))$, 
$j_{IJ}:S^{oJ}\hookrightarrow S^{oI}$ being the open embedding, and $d_{IJ}=0$ if $I\notin J$.
\end{itemize}
By definition, we have for $(K,W)\in C(P_{\mathbb Z_l,fil,k}(S^{o,et})^{ad,(Z_i)})$, 
$j^*j_{*w}(K,W)=(K,W)$ and $j^*j_{!w}(K,W)=(K,W)$.
For $(K',W)\in C(P_{\mathbb Z_l,fil}(S^{et})^{ad,(Z_i)})$, there is, by (i),
\begin{itemize}
\item a canonical map $\ad(j^*,j_{*w})(K',W):(K',W)\to j_{*w}j^*(K',W)$ in $C(P_{\mathbb Z_l,fil,k}(S^{et}))$,
\item a canonical map $\ad(j_{!w},j^*)(K',W):j_{!w}j^*(K',W)\to(K',W)$ in $C(P_{\mathbb Z_l,fil,k}(S^{et}))$.
\end{itemize}
\item[(iii)] Let $S\in\Var(k)$. Let $j:S^o:=S\backslash Z\hookrightarrow S$ an open embedding with
$Z=V(\mathcal I)\subset S$ an arbitrary closed subset (over $k$), $\mathcal I\subset O_S$ being an ideal subsheaf.
Let $k/k'$ a field extension. 
For $(K,W)\in C(P_{\mathbb Z_l,fil,k}(S_{k'}^{o,et})^{ad,(Z_i)})$, (ii) gives
\begin{itemize}
\item the (bi)-filtered complex of $D_S$-modules
\begin{equation*}
j_{*w}(K,W):=\varinjlim_{\mathcal D(Z/S)}\Tot_{card I=\bullet}(j_{I*w}j_I^{o*}(K,W))\in C(P_{\mathbb Z_l,fil,k}(S_{k'}^{et})) 
\end{equation*}
\item the (bi)-filtered complex of $D_S$-modules
\begin{eqnarray*}
j_{!w}(K,W):=\varprojlim_{\mathcal D(Z/S)}\Tot_{card I=-\bullet}(j_{I!w}j_I^{o*}(K,W)) 
=\mathbb D_S^vj_{*w}\mathbb D_S^v(K,W)\in C(P_{\mathbb Z_l,fil}(S_{k'}^{et})),
\end{eqnarray*}
\end{itemize}
By definition, we have for $(K,W)\in C(P_{\mathbb Z_l,fil,k}(S_{k'}^{o,et})^{ad,(Z_i)})$, 
$j^*j_{*w}(K,W)=(K,W)$ and $j^*j_{!w}(K,W)=(K,W)$.
For $(K',W)\in C(P_{\mathbb Z_l,fil,k}(S_{k'}^{et})^{ad,(Z_i)})$, we have, see (ii),
\begin{itemize}
\item the canonical map $\ad(j^*,j_{*w})(K',W):(K',W)\to j_{*w}j^*(K',W)$ in $C(P_{\mathbb Z_l,fil,k}(S_{k'}^{et}))$,
\item the canonical map $\ad(j_{!w},j^*)(K',W):j_{!w}j^*(K',W)\to(K',W)$ in $C(P_{\mathbb Z_l,fil,k}(S_{k'}^{et}))$.
\end{itemize}
\end{itemize}
\end{defi}

\begin{defi}\label{gammawet}
Let $S\in\Var(k)$. Let $Z\subset S$ a closed subset. 
Denote by $j:S\backslash Z\hookrightarrow S$ the complementary open embedding. 
Let $l$ a prime number.
\begin{itemize}
\item[(i)] We define using definition \ref{jwet}, the filtered Hodge support section functor
\begin{eqnarray*}
\Gamma^w_Z:C(P_{\mathbb Z_lfil,k}(S^{et})^{ad,(Z_i)})\to C(P_{\mathbb Z_lfil,k}(S^{et})), \\ 
(K,W)\mapsto\Gamma^w_Z(K,W):=\Cone(\ad(j^*,j_{*w})(K,W):(K,W)\to j_{*w}j^*(K,W))[-1],
\end{eqnarray*}
together we the canonical map $\gamma^w_Z(K,W):\Gamma^w_Z(K,W)\to (K,W)$. Since 
\begin{equation*}
j_{*w}:C(P_{\mathbb Z_lfil,k}(S^{o,et})^{ad,(Z_i)})\to C(P_{\mathbb Z_lfil,k}(S^{et}))
\end{equation*}
is an exact functor, $\Gamma^w_Z$ induces the functor
\begin{eqnarray*}
\Gamma^w_Z:D_{\mathbb Z_lfil,c,k}(S^{et})^{ad,(Z_i)}\to D_{\mathbb Z_lfil,c,k}(S^{et}), \; (K,W)\mapsto\Gamma^w_Z(K,W)
\end{eqnarray*}
\item[(ii)] We define using definition \ref{jw}, the dual filtered Hodge support section functor
\begin{eqnarray*}
\Gamma^{\vee,w}_Z:C(P_{\mathbb Z_l,fil}(S^{et})^{ad,(Z_i)})\to C(P_{\mathbb Z_l,fil,k}(S^{et})), \\ 
(K,W)\mapsto\Gamma^{\vee,w}_Z(K,W):=\Cone(\ad(j_{!w},j^*)(K,W):j_{!w},j^*(K,W)\to (K,W)),
\end{eqnarray*}
together we the canonical map $\gamma^{\vee,Hdg}_Z(K,W):(K,W)\to\Gamma_Z^{\vee,w}(K,W)$. Since 
\begin{equation*}
j_{!w}:C(P_{\mathbb Z_l,fil,k}(S^{o,et})^{ad,(Z_i)})\to C(P_{\mathbb Z_l,fil,k}(S^{et}))
\end{equation*}
is an exact functor, $\Gamma^{\vee,w}_Z$ induces the functor
\begin{eqnarray*}
\Gamma^{\vee,w}_Z:D_{\mathbb Z_lfil,c,k}(S^{et})^{ad,(Z_i)}\to D_{\mathbb Z_lfil,c,k}(S^{et}), \; 
(K,W)\mapsto\Gamma^{\vee,w}_Z(K,W)
\end{eqnarray*}
\end{itemize}
Let $k/k'$ a field extension. Let $S\in\Var(k)$. Let $Z\subset S$ a closed subset (over $k$). 
Denote by $j:S\backslash Z\hookrightarrow S$ the complementary open embedding. Let $l$ a prime number.
\begin{itemize}
\item[(i)'] Then (i) gives the filtered Hodge support section functor
\begin{eqnarray*}
\Gamma^w_Z:C(P_{\mathbb Z_l,fil,k}(S_{k'}^{et})^{ad,(Z_i)})\to C(P_{\mathbb Z_l,fil,k}(S_{k'}^{et})), \\ 
(K,W)\mapsto\Gamma^w_Z(K,W):=\Cone(\ad(j^*,j_{*w})(K,W):(K,W)\to j_{*w}j^*(K,W))[-1],
\end{eqnarray*}
together we the canonical map $\gamma^w_Z(K,W):\Gamma^w_Z(K,W)\to (K,W)$, which induces
\begin{eqnarray*}
\Gamma^w_Z:D_{\mathbb Z_lfil,c,k}(S_{k'}^{et})^{ad,(Z_i)}\to D_{\mathbb Z_lfil,c,k}(S_{k'}^{et}), \; 
(K,W)\mapsto\Gamma^w_Z(K,W)
\end{eqnarray*}
\item[(ii)'] Then, (ii) gives the dual filtered Hodge support section functor
\begin{eqnarray*}
\Gamma^{\vee,w}_Z:C(P_{\mathbb Z_l,fil,k}(S_{k'}^{et})^{ad,(Z_i)})\to C(P_{\mathbb Z_l,fil,k}(S_{k'}^{et})), \\ 
(K,W)\mapsto\Gamma^{\vee,w}_Z(K,W):=\Cone(\ad(j_{!w},j^*)(K,W):j_{!w},j^*(K,W)\to (K,W)),
\end{eqnarray*}
together we the canonical map $\gamma^{\vee,Hdg}_Z(K,W):(K,W)\to\Gamma_Z^{\vee,w}(K,W)$, which induces
\begin{eqnarray*}
\Gamma^{\vee,w}_Z:D_{\mathbb Z_lfil,c,k}(S_{k'}^{et})^{ad,(Z_i)}\to D_{\mathbb Z_lfil,c,k}(S_{k'}^{et}), \; 
(K,W)\mapsto\Gamma^{\vee,w}_Z(K,W)
\end{eqnarray*}
\end{itemize}
\end{defi}

Let $S\in\Var(k)$ and $D=V(s)\subset S$ a Cartier divisor. 
Denote $i:D\hookrightarrow S$ the closed embedding and $j:S^o:=S\backslash D\hookrightarrow S$ the open embedding.
Let $\pi:\tilde S^{o}\to S^{o}$ the universal covering. Let $l$ a prime number.
We then consider, for $(K,W)\in D_{\mathbb Z_lfil,c,k}(S^{o,et})^{ad,D}=\Ho(C(P_{\mathbb Z_lfil,k}(S^{o,et})^{ad,D})$, 
\begin{itemize}
\item the filtered nearby cycle functor 
\begin{equation*}
\psi_D(K,W):=(\psi_DK,W)\in D_{\mathbb Z_lfil,c,k}(D^{et}), \; W_k(\psi_D(K,W)):=<W_k\psi_DK,W(N)_k\psi_DK>\subset\psi_DK,
\end{equation*}
\item the vanishing cycle functor 
\begin{equation*}
\phi_D(K,W):=\Cone(\ad(j\circ\pi^*,j\circ\pi_*)(K,W)\circ ev(K,W):i^*(K,W)\to\psi_D(K,W))
\in D_{\mathbb Z_lfil,c}(D^{et}),
\end{equation*}
where the morphism 
\begin{equation*}
\ad(j\circ\pi^*,j\circ\pi_*)(K,W)\circ ev(K,W):i^*(K,W)\to 
i^*R(j\circ\pi)_*\pi^*\mathcal Hom(\mathcal A^{\bullet}_D(S^o),(K,W)) 
\end{equation*}
being by definition compatible with the weight monodromy filtration induces the morphism 
\begin{equation*}
\ad(j\circ\pi^*,j\circ\pi_*)(K,W)\circ ev(K,W):i^*(K,W)\to\psi_D(K,W) 
\end{equation*}
\item the canonical morphisms in $D_{\mathbb Z_lfil,c,k}(D^{et})$
\begin{eqnarray*}
can(K,W):=c(\phi_D(K,W)):\psi_D(K,W)\to\phi_D(K,W), \\
var(K,W):=(0,T-I):\phi_D(K,W)\to\psi_D(K,W),
\end{eqnarray*}
\end{itemize}
Let $k/k'$ a field extension.
For $(K,W)\in D_{\mathbb Z_lfil,c,k}(S_{k'}^{o,et})^{ad,D}=\Ho(C(P_{\mathbb Z_lfil,k}(S_{k'}^{o,et})^{ad,D})$, we get
\begin{itemize}
\item the filtered nearby cycle functor 
\begin{equation*}
\psi_D(K,W):=(\psi_DK,W)\in D_{\mathbb Z_lfil,c,k}(D_{k'}^{et}), \; 
W_k(\psi_D(K,W)):=<W_k\psi_DK,W(N)_k\psi_DK>\subset\psi_DK,
\end{equation*}
\item the vanishing cycle functor 
\begin{equation*}
\phi_D(K,W):=\Cone(\ad(j\circ\pi^*,j\circ\pi_*)(K,W)\circ ev(K,W):
i^*(K,W)\to\psi_D(K,W))\in D_{\mathbb Z_lfil,c,k}(D_{k'}^{et}),
\end{equation*}
\item the canonical morphisms in $D_{\mathbb Z_lfil,c,k}(D_{k'}^{et})$
\begin{eqnarray*}
can(K,W):=c(\phi_D(K,W)):\psi_D(K,W)\to\phi_D(K,W), \\
var(K,W):=(0,T-I):\phi_D(K,W)\to\psi_D(K,W).
\end{eqnarray*}
\end{itemize}

\begin{defi}\label{uketw}
Let $S\in\Var(k)$ and $D=V(s)\subset S$ a Cartier divisor. 
Denote $i:D\hookrightarrow S$ the closed embedding and $j:S^o:=S\backslash D\hookrightarrow S$ the open embedding.
Let $\pi:\tilde S^{o}\to S^{o}$ the universal covering.
For $(K,W)\in P_{fil,k}(S_{\mathbb C}^{o,an})^{ad,D}$, we consider 
\begin{itemize}
\item the canonical decomposition
\begin{equation*}
\psi_D(K,W):=\psi_D^u(K,W)\oplus\psi_D^t(K,W)\in P_{\mathbb Z_lfil,k}(D^{et}), 
\end{equation*} 
\item the canonical decomposition
\begin{equation*}
\phi_D(K,W):=\phi_D^u(K,W)\oplus\phi_D^t(K,W)\in P_{\mathbb Z_lfil,k}(D^{et}), 
\end{equation*} 
\item the canonical morphisms in $D_{\mathbb Z_lfil}(D^{et})$
\begin{eqnarray*}
can(K,W):=c(\phi^u_D(K,W)):\psi^u_D(K,W)\to\phi^u_D(K,W) \; , \; var(K,W):=(0,T-I):\phi^u_D(K,W)\to\psi^u_D(K,W).
\end{eqnarray*}
\item the maximal extension
\begin{eqnarray*}
x_{S^o/S}(K,W):=\Cone(p_u\circ\ad(i^*,i_*)(-)\circ\ad(\pi^*,\pi_*)(K,W)\circ ev(K,W): \\
j_{*w}(K,W)\to\psi^u_D(K,W))\in D_{\mathbb Z_lfil,c,k}(S^{et}),
\end{eqnarray*}
where the morphism 
\begin{equation*}
\ad(\pi^*,\pi_*)(K,W)\circ ev(K,W):Rj_*(K,W)\to R(j\circ\pi)_*\pi^*\mathcal Hom(\mathcal A^{\bullet}_D(S^o),(K,W)) 
\end{equation*}
being by definition compatible with the weight monodromy filtration induces the morphism 
\begin{equation*}
\ad(i^*,i_*)(-)\circ\ad(\pi^*,\pi_*)(K,W)\circ ev(K,W):j_{*w}(K,W)\to\psi_D(K,W).
\end{equation*}
\end{itemize}
\end{defi}

\begin{defi}\label{fwet}
Let $k$ a field of characteristic zero. Let $l$ a prime number.
\begin{itemize}
\item[(i)]Let $f:X\to S$ a morphism with $S,X\in\Var(k)$. 
Consider the graph factorization $f:X\xrightarrow{l}X\times S\xrightarrow{p}S$ of $f$ 
where $l$ the the graph closed embedding and $p$ is the projection.
We have, using definition \ref{gammawet},
\begin{itemize}
\item the inverse image functor
\begin{eqnarray*}
f^{*w}:D_{\mathbb Z_lfil,c,k}(S^{et})^{ad,(\Gamma_{f,i})}\to D_{\mathbb Z_lfil,c,k}(X^{et}), \; 
(K,W)\mapsto f^{*w}(K,W):=l^*\Gamma_X^{\vee,w}p^*(K,W)
\end{eqnarray*}
\item the exceptional inverse image functor
\begin{eqnarray*}
f^{!w}:D_{\mathbb Z_lfil,c,k}(S^{et})^{ad,(\Gamma_{f,i})}\to D_{\mathbb Z_lfil,c,k}(X^{et}), \; 
(K,W)\mapsto f^{!w}(K,W):=l^*\Gamma_X^wp^*(K,W).
\end{eqnarray*}
\end{itemize}
\item[(i)']Let $f:X\to S$ a morphism with $S,X\in\Var(k)$. 
Consider the graph factorization $f:X\xrightarrow{l}X\times S\xrightarrow{p}S$ of $f$ 
where $l$ the the graph closed embedding and $p$ is the projection. Let $k/k'$ a field extension. We have, by (i),
\begin{itemize}
\item the inverse image functor
\begin{eqnarray*}
f^{*w}:D_{\mathbb Z_lfil,c,k}(S_{k'}^{et})^{ad,(\Gamma_{f,i})}\to D_{\mathbb Z_lfil,c,k}(X_{k'}^{et}), \; 
(K,W)\mapsto f^{*w}(K,W):=l^*\Gamma_X^{\vee,w}p^*(K,W)
\end{eqnarray*}
\item the exceptional inverse image functor
\begin{eqnarray*}
f^{!w}:D_{\mathbb Z_lfil,c,k}(S_{k'}^{et})^{ad,(\Gamma_{f,i})}\to D_{\mathbb Z_lfil,c,k}(X_{k'}^{et}), \; 
(K,W)\mapsto f^{!w}(K,W):=l^*\Gamma_X^wp^*(K,W).
\end{eqnarray*}
\end{itemize}
\item[(ii)]Let $f:X\to S$ a morphism with $S,X\in\Var(k)$.
Consider a compactification $f:X\hookrightarrow{j}\bar X\xrightarrow{\bar f}S$ of $f$
with $\bar X\in\Var(k)$, $j$ an open embedding and $\bar f$ a proper morphism.
Denote $Z:=\bar X\backslash X$. We have, using definition \ref{jw},
\begin{itemize}
\item the direct image functor
\begin{eqnarray*}
Rf_{*w}:D_{\mathbb Z_lfil,c,k}(X^{et})^{ad,(Z_i)}\to D_{\mathbb Z_lfil,c,k}(S^{et}), \;
(K,W)\mapsto Rf_{*w}(K,W):=R\bar f_*j_{*w}(K,W)
\end{eqnarray*}
\item the proper direct image functor
\begin{eqnarray*}
Rf_{!w}:D_{\mathbb Z_lfil,c,k}(X^{et})^{ad,(Z_i)}\to D_{\mathbb Z_lfil,c,k}(S^{et}), \; 
(K,W)\mapsto Rf_{!w}(K,W):=R\bar f_*j_{!w}(K,W).
\end{eqnarray*}
\end{itemize}
\item[(ii)']Let $f:X\to S$ a morphism with $S,X\in\Var(k)$.
Consider a compactification $f:X\hookrightarrow{j}\bar X\xrightarrow{\bar f}S$ of $f$
with $\bar X\in\Var(k)$, $j$ an open embedding and $\bar f$ a proper morphism. Let $k/k'$ a field extension. 
We have, by (ii),
\begin{itemize}
\item the direct image functor
\begin{eqnarray*}
Rf_{*w}:D_{\mathbb Z_lfil,c,k}(X_{k'}^{et})^{ad,(Z_i)}\to D_{\mathbb Z_lfil,c,k}(S_{k'}^{et}), \; 
(K,W)\mapsto Rf_{*w}(K,W):=R\bar f_*j_{*w}(K,W)
\end{eqnarray*}
\item the proper direct image functor
\begin{eqnarray*}
Rf_{!w}:D_{\mathbb Z_lfil,c,k}(X_{k'}^{et})^{ad,(Z_i)}\to D_{\mathbb Z_lfil,c,k}(S_{k'}^{et}), \; 
(K,W)\mapsto Rf_{!w}(K,W):=R\bar f_*j_{!w}(K,W).
\end{eqnarray*}
\end{itemize}
\item[(iii)] Let $S\in\Var(k)$. Denote by $\Delta_S:S\hookrightarrow S\times S$ the diagonal closed embedding and
$p_1:S\times S\to S$, $p_2:S\times S\to S$ the projections.
We have by (i) the functor
\begin{eqnarray*}
\otimes^{Lw}:D_{\mathbb Z_lfil,c,k}(S^{et})^{ad,(S_i)}\times D_{\mathbb Z_lfil,c,k}(S^{et})^{ad,(S_i)}
\to D_{\mathbb Z_lfil,c,k}(S^{et}), \\ 
((K_1,W),(K_2,W))\mapsto (K_1,W)\otimes^{L,w}(K_2,W):=\Delta_S^{!w}(p_1^*(K_1,W)\otimes^L p_2^*(K_2,W)).
\end{eqnarray*}
\end{itemize}
\end{defi}

Let $S\in\Var(k)$ and $D=V(s)\subset S$ a Cartier divisor. 
Denote $i:D\hookrightarrow S$ the closed embedding and $j:S^o:=S\backslash D\hookrightarrow S$ the open embedding.
In the filtered case, we get, for $(K,W)\in P_{\mathbb Z_lfil,k}(S^{et})^{ad,D}$ the map in $D_{\mathbb Z_lfil,c,k}(S^{et})$
\begin{eqnarray*}
Is(K,W):=(0,(\ad(j^*,j_*)(K),\ad(j\circ\pi^*,j\circ\pi_*)(K)),0): \\
(K,W)\to(\psi^u_D(K,W)\xrightarrow{(c(x_{S^o/S}(K,W)),can(K,W))}x_{S^o/S}(K,W)\oplus i_*\phi^u_D(K,W) \\
\xrightarrow{((0,T-I),var(K,W):=(0,T-I))}\psi^u_DK)[-1]
\end{eqnarray*}
which is NOT an isomorphism in general (it leads to different $W$-filtration on perverse cohomology).

We recall the definition of constructible pro etale sheaves on algebraic varieties over a subfield $k\subset K$
of a $p$ adic field:

\begin{defi}\label{constr3}
Let $k\subset K\subset\mathbb C_p$ a subfield of a p adic field. Let $l$ a prime number. 
Let $S\in\Var(k)$.
\begin{itemize}
\item[(i)]A sheaf $K\in\Shv_{\mathbb Z_l}(S_K^{an,pet})$ is called constructible 
if there exists a stratification $S=\sqcup_{\alpha} S_{\alpha}$ 
with $l_{\alpha}:S_{\alpha}\hookrightarrow S$ locally closed subsets (defined over $k$)
such that $l_{\alpha}^*K\in\Shv(S_{\alpha}^{an,pet})$ 
are (finite dimensional) local systems (for the etale topology) for all $\alpha$. 
\item[(ii)]We denote by 
\begin{equation*}
C_{\mathbb Z_l,c,k}(S^{an,pet})\subset C_{\mathbb Z_l}(S^{an,pet}) \; \mbox{and} \; 
D_{\mathbb Z_l,c,k}(S^{an,pet})\subset D_{\mathbb Z_l}(S^{an,pet})
\end{equation*}
the full subcategories consisting of $K\in C(S_K^{an,pet})$ 
such that $a_{et}H^nK\in\Shv(S_K^{an,pet})$ are constructible for all $n\in\mathbb Z$.
\item[(ii)']We denote by 
\begin{equation*}
C_{\mathbb Z_lfil,c,k}(S^{an,pet})\subset C_{\mathbb Z_lfil}(S^{an,pet}) \; \mbox{and} \; 
D_{\mathbb Z_lfil,c,k}(S^{an,pet})\subset D_{\mathbb Z_lfil}(S^{an,pet})
\end{equation*}
the full subcategories consisting of $(K,W)\in C_{\mathbb Z_lfil}(S_K^{an,pet})$ 
such that $a_{et}H^n\Gr_k^WK\in\Shv_{\mathbb Z_l}(S_K^{an,pet})$ are constructible for all $n,k\in\mathbb Z$.
\item[(iii)] We denote by $P_{\mathbb Z_l,k}(S^{an,pet})\subset D_{\mathbb Z_l,c,k}(S^{an,pet})$ 
the full subcategory of perverse sheaves.
\item[(iii)'] We denote by $P_{\mathbb Z_l,fil,k}(S^{an,pet})\subset D_{\mathbb Z_lfil,c,k}(S^{an,pet})$ 
the full subcategory of filtered perverse sheaves.
\end{itemize}
\end{defi}

Let $K\subset\mathbb C_p$ a p adic field.
Let $S\in\Var(K)$ and $D=V(s)\subset S$ a Cartier divisor. 
Denote $i:D\hookrightarrow S$ the closed embedding and $j:S^o:=S\backslash D\hookrightarrow S$ the open embedding.
Let $\pi:\tilde S^{o,an}\to S^{o,an}$ the perfectoid universal covering (see \cite{Scholze}.
We then consider, 
\begin{itemize}
\item for $K\in C_{\mathbb Z_l}(S^{o,an,pet})$, the nearby cycle functor 
\begin{equation*}
\psi_DK:=i^*R(j\circ\pi)_*\pi^*K=i^*(j\circ\pi)_*\pi^*K\in D_{\mathbb Z_l}(D^{an,pet}), 
\end{equation*}
we write again $\psi_DK:=i_*\psi_DK\in D_{\mathbb Z_l}(S^{an,pet})$,
\item for $K\in C_{\mathbb Z_l}(S^{o,an,pet})$, the vanishing cycle functor 
\begin{equation*}
\phi_DK:=\Cone(\ad(j\circ\pi^*,j\circ\pi_*)(K):i^*K\to\psi_D K\in D_{\mathbb Z_l}(D^{an,pet})
\end{equation*}
together with the canonical map $c(\phi_DK):\psi_DK\to\phi_DK$ in $D_{\mathbb Z_l}(S^{an,pet})$
we write again $\phi_DK:=i_*\phi_DK\in D_{\mathbb Z_l}(S^{an,pet})$,
\item for $K\in C_{\mathbb Z_l,c,k}(S^{o,an,pet})$, the canonical morphisms in $D_{\mathbb Z_l,c,k}(D^{an,pet})$
\begin{eqnarray*}
can(K):=c(\phi_DK):\psi_DK\to\phi_DK \; , \; var(K):=(0,T-I):\phi_DK\to\psi_DK,
\end{eqnarray*}
\item for $K\in P_{\mathbb Z_l}(S^{o,an,pet})$, the maximal extension
\begin{eqnarray*}
x_{S^o/S}(K):=\Cone(p_u\circ\ad(i^*,i_*)(-)\circ\ad(\pi^*,\pi_*)(K): \\
Rj_*K\to i_*R(j\circ\pi)_*\pi^*K=:\psi_DK\to\psi_D^uK)\in P_{\mathbb Z_l}(S^{an,pet}).
\end{eqnarray*}
\end{itemize}
 
Let $k\subset K\subset\mathbb C_p$ a subfield of a p adic field.
Let $S\in\Var(k)$ and $D=V(s)\subset S$ a Cartier divisor. Denote by $S^o:=S\backslash D$.
By definition, we have for $K\in C_{\mathbb Z_l}(S^{o,et})$
\begin{eqnarray*}
\an_S^*(\psi_DK)=\psi_D(\an_S^*K)\in D_{\mathbb Z_l}(S_K^{an,pet}) \; \mbox{and} \;
\an_S^*(\phi_DK)=\phi_D(\an_S^*K)\in D_{\mathbb Z_l}(S_K^{an,pet}).
\end{eqnarray*}
where $\an_S:S_K^{an,pet}\xrightarrow{\an_S}S_K^{et}\xrightarrow{\otimes_kK}S^{et}$
is the morphism of site induced by the analytical functor.

We then deduce from the algebraic case the following :

\begin{cor}\label{PSkp}
Let $k\subset K\subset\mathbb C_p$ a subfield of a p adic field.
Let $S\in\Var(k)$ and $D=V(s)\subset S$ a Cartier divisor. 
Denote $i:D\hookrightarrow S$ the closed embedding and $j:S^o:=S\backslash D\hookrightarrow S$ the open embedding.
\begin{itemize} 
\item[(i)]For $K\in P_{\mathbb Z_l,k}(S^{et})$, we have $\psi_D\an_S^*K,\phi_D\an_S^*K\in P_{\mathbb Z_l,k}(S_K^{an,pet})$
\item[(ii)]We have, for $K\in P_{\mathbb Z_l,k}(S^{et})$, denoting again $K=\an_S^*K\in P_{\mathbb Z_l,k}(S_K^{an,pet})$,
\begin{eqnarray*}
Is(K):=(0,(\ad(j^*,j_*)(K),\ad((j\circ\pi)^*,(j\circ\pi)_*)(K)),0): \\
K\xrightarrow{\sim}
(\psi^u_DK\xrightarrow{(c(x_{S^o/S}(K)),can(K))}x_{S^o/S}(K')\oplus i_*\phi^u_DK\xrightarrow{((0,T-I),var(K))}\psi^u_DK)[-1]
\end{eqnarray*}
the canonical isomorphism in $D_{\mathbb Z_l,c,k}(S_K^{an,pet})$.
\end{itemize}
\end{cor}

\begin{proof}
\noindent(i):Follows from proposition \ref{phipsiperv}.

\noindent(ii):Follows immediately from theorem \ref{PSk}.
\end{proof}

\subsection{Presheaves on the big Zariski site or on the big etale site}

Let $k$ a field of caracteristic $0$.
For $S\in\Var(k)$, we denote by $\rho_S:\Var(k)^{sm}/S\hookrightarrow\Var(k)/S$ be the full subcategory 
consisting of the objects $U/S=(U,h)\in\Var(k)/S$ such that the morphism $h:U\to S$ is smooth. 
That is, $\Var(k)^{sm}/S$ is the category  
\begin{itemize}
\item whose objects are smooth morphisms 
$U/S=(U,h)$, $h:U\to S$ with $U\in\Var(k)$, 
\item whose morphisms $g:U/S=(U,h_1)\to V/S=(V,h_2)$ 
is a morphism $g:U\to V$ of complex algebraic varieties such that $h_2\circ g=h_1$. 
\end{itemize}
We denote again $\rho_S:\Var(k)/S\to\Var(k)^{sm}/S$ the associated morphism of site.
We will consider 
\begin{equation*}
r^s(S):\Var(k)\xrightarrow{r(S)}\Var(k)/S\xrightarrow{\rho_S}\Var(k)^{sm}/S
\end{equation*}
the composite morphism of site.
For $S\in\Var(k)$, we denote by $\mathbb Z_S:=\mathbb Z(S/S)\in \PSh(\Var(k)^{sm}/S)$ the constant presheaf
By Yoneda lemma, we have for $F\in C(\Var(k)^{sm}/S)$, $\mathcal Hom(\mathbb Z_S,F)=F$.
For $f:T\to S$ a morphism, with $T,S\in\Var(k)$, we have the following commutative diagram of sites
\begin{equation}\label{pf0}
\xymatrix{\Var(k)/T\ar[d]^{P(f)}\ar[r]^{\rho_T} & \Var(k)^{sm}/T\ar[d]^{P(f)} \\ 
\Var(k)/S\ar[r]^{\rho_S} & \Var(k)^{sm}/S} 
\end{equation}
We denote, for $S\in\Var(k)$, the obvious morphism of sites 
\begin{equation*}
\tilde e(S):\Var(k)/S\xrightarrow{\rho_S}\Var(k)^{sm}/S\xrightarrow{e(S)}\Ouv(S)  
\end{equation*}
where $\Ouv(S)$ is the set of the Zariski open subsets of $S$, given by the inclusion functors
$\tilde e(S):\Ouv(S)\hookrightarrow\Var(k)^{sm}/S\hookrightarrow\Var(k)/S$.
By definition, for $f:T\to S$ a morphism with $S,T\in\Var(k)$, the commutative diagram of sites (\ref{pf0})
extend a commutative diagram of sites :
\begin{equation}\label{empf0}
\xymatrix{
\tilde e(T):\Var(k)/T\ar[d]^{P(f)}\ar[rr]^{\rho_T} & \, & \Var(k)^{sm}/T\ar[d]^{P(f)}\ar[rr]^{e(T)} & \, & 
\Ouv(T)\ar[d]^{P(f)} \\
\tilde e(S):\Var(k)/S\ar[rr]^{\rho_S} & \, & \Var(k)^{sm}/S\ar[rr]^{e(S)} & \, & \Ouv(S)}
\end{equation}

\begin{itemize}
\item As usual, we denote by
\begin{equation*}
(f^*,f_*):=(P(f)^*,P(f)_*):C(\Var(k)^{sm}/S)\to C(\Var(k)^{sm}/T)
\end{equation*}
the adjonction induced by $P(f):\Var(k)^{sm}/T\to\Var(k)^{sm}/S$.
Since the colimits involved in the definition of $f^*=P(f)^*$ are filtered, $f^*$ also preserve monomorphism. 
Hence, we get an adjonction
\begin{equation*}
(f^*,f_*):C_{fil}(\Var(k)^{sm}/S)\leftrightarrows C_{fil}(\Var(k)^{sm}/T), \; f^*(G,F):=(f^*G,f^*F)
\end{equation*}
\item As usual, we denote by
\begin{equation*}
(f^*,f_*):=(P(f)^*,P(f)_*):C(\Var(k)/S)\to C(\Var(k)/T)
\end{equation*}
the adjonction induced by $P(f):\Var(k)/T\to \Var(k)/S$.
Since the colimits involved in the definition of $f^*=P(f)^*$ are filtered, $f^*$ also preserve monomorphism. 
Hence, we get an adjonction
\begin{equation*}
(f^*,f_*):C_{fil}(\Var(k)/S)\leftrightarrows C_{fil}(\Var(k)/T), \; f^*(G,F):=(f^*G,f^*F)
\end{equation*}
\end{itemize}

For $h:U\to S$ a smooth morphism with $U,S\in\Var(k)$, 
the pullback functor $P(h):\Var(k)^{sm}/S\to \Var(k)^{sm}/U$ 
admits a left adjoint $C(h)(X\to U)=(X\to U\to S)$.
Hence, $h^*:C(\Var(k)^{sm}/S)\to C(\Var(k)^{sm}/U)$ admits a left adjoint
\begin{equation*}
h_{\sharp}:C(\Var(k)^{sm}/U)\to C(\Var(k)^{sm}/S), \; 
F\mapsto((V,h_0)\mapsto\lim_{(V',h\circ h')\to(V,h_0)}F(V',h'))
\end{equation*}
Note that we have for $V/U=(V,h')$ with $h':V\to U$ a smooth morphism 
we have $h_{\sharp}(\mathbb Z(V/U))=\mathbb Z(V'/S)$ with $V'/S=(V',h\circ h')$.
Hence, since projective presheaves are the direct summands of the representable presheaves,
$h_{\sharp}$ sends projective presheaves to projective presheaves.

We have the support section functors of a closed embedding $i:Z\hookrightarrow S$ for presheaves on the big Zariski site.
\begin{defi}\label{gamma}
Let $i:Z\hookrightarrow S$ be a closed embedding with $S,Z\in\Var(k)$ 
and $j:S\backslash Z\hookrightarrow S$ be the open complementary subset.
\begin{itemize}
\item[(i)] We define the functor
\begin{equation*}
\Gamma_Z:C(\Var(k)^{sm}/S)\to C(\Var(k)^{sm}/S), \;
G^{\bullet}\mapsto\Gamma_Z G^{\bullet}:=\Cone(\ad(j^*,j_*)(G^{\bullet}):G^{\bullet}\to j_*j^*G^{\bullet})[-1],
\end{equation*}
so that there is then a canonical map $\gamma_Z(G^{\bullet}):\Gamma_ZG^{\bullet}\to G^{\bullet}$.
\item[(ii)] We have the dual functor of (i) :
\begin{equation*}
\Gamma^{\vee}_Z:C(\Var(k)^{sm}/S)\to C(\Var(k)^{sm}/S), \; 
F\mapsto\Gamma^{\vee}_Z(F^{\bullet}):=\Cone(\ad(j_{\sharp},j^*)(G^{\bullet}):j_{\sharp}j^*G^{\bullet}\to G^{\bullet}), 
\end{equation*}
together with the canonical map $\gamma^{\vee}_Z(G):F\to\Gamma^{\vee}_Z(G)$.
\item[(iii)] For $F,G\in C(\Var(k)^{sm}/S)$, we denote by 
\begin{equation*}
I(\gamma,hom)(F,G):=(I,I(j_{\sharp},j^*)(F,G)^{-1}):\Gamma_Z\mathcal Hom(F,G)\xrightarrow{\sim}\mathcal Hom(\Gamma^{\vee}_ZF,G)
\end{equation*}
the canonical isomorphism given by adjonction.
\end{itemize}
\end{defi}

Let $S_{\bullet}\in\Fun(\mathcal I,\Var(k))$ with $\mathcal I\in\Cat$, a diagram of algebraic varieties.
It gives the diagram of sites $\Var(k)^2/S_{\bullet}\in\Fun(\mathcal I,\Cat)$.  
\begin{itemize}
\item Then $C_{fil}(\Var(k)/S_{\bullet})$ is the category  
\begin{itemize}
\item whose objects $(G,F)=((G_I,F)_{I\in\mathcal I},u_{IJ})$,
with $(G_I,F)\in C_{fil}(\Var(k)/S_I)$,
and $u_{IJ}:(G_I,F)\to r_{IJ*}(G_J,F)$ for $r_{IJ}:I\to J$, denoting again $r_{IJ}:S_I\to S_J$, are morphisms
satisfying for $I\to J\to K$, $r_{IJ*}u_{JK}\circ u_{IJ}=u_{IK}$ in $C_{fil}(\Var(k)/S_I)$,
\item the morphisms $m:((G,F),u_{IJ})\to((H,F),v_{IJ})$ being (see section 2.1) a family of morphisms of complexes,  
\begin{equation*}
m=(m_I:(G_I,F)\to (H_I,F))_{I\in\mathcal I}
\end{equation*}
such that $v_{IJ}\circ m_I=p_{IJ*}m_J\circ u_{IJ}$ in $C_{fil}(\Var(k)/S_I)$.
\end{itemize}
\item Then $C_{fil}(\Var(k)^{sm}/S_{\bullet})$ is the category  
\begin{itemize}
\item whose objects $(G,F)=((G_I,F)_{I\in\mathcal I},u_{IJ})$,
with $(G_I,F)\in C_{fil}(\Var(k)^{sm}/S_I)$,
and $u_{IJ}:(G_I,F)\to r_{IJ*}(G_J,F)$ for $r_{IJ}:I\to J$, denoting again $r_{IJ}:S_I\to S_J$, are morphisms
satisfying for $I\to J\to K$, $r_{IJ*}u_{JK}\circ u_{IJ}=u_{IK}$ in $C_{fil}(\Var(k)^{sm}/S_I)$,
\item the morphisms $m:((G,F),u_{IJ})\to((H,F),v_{IJ})$ being (see section 2.1) a family of morphisms of complexes,  
\begin{equation*}
m=(m_I:(G_I,F)\to (H_I,F))_{I\in\mathcal I}
\end{equation*}
such that $v_{IJ}\circ m_I=p_{IJ*}m_J\circ u_{IJ}$ in $C_{fil}(\Var(k)^{sm}/S_I)$.
\end{itemize}
\end{itemize}
As usual, we denote by
\begin{eqnarray*}
(f_{\bullet}^*,f_{\bullet*}):=(P(f_{\bullet})^*,P(f_{\bullet})_*): 
C(\Var(k)^{(sm)}/S_{\bullet})\to C(\Var(k)^{(sm)}/T_{\bullet})
\end{eqnarray*}
the adjonction induced by 
$P(f_{\bullet}):\Var(k)^{(sm)}/T_{\bullet}\to \Var(k)^{(sm)}/S_{\bullet}$.
Since the colimits involved in the definition of $f_{\bullet}^*=P(f_{\bullet})^*$ are filtered, 
$f_{\bullet}^*$ also preserve monomorphism. Hence, we get an adjonction
\begin{eqnarray*}
(f_{\bullet}^*,f_{\bullet*}):
C_{fil}(\Var(k)^{(sm)}/S_{\bullet})\leftrightarrows C_{fil}(\Var(k)^{(sm)}/T_{\bullet}), \\
f_{\bullet}^*((G_I,F),u_{IJ}):=((f_I^*G_I,f_I^*F),T(f_I,r_{IJ})(-)\circ f_I^*u_{IJ}).
\end{eqnarray*}

Let $S\in\Var(k)$. Let $S=\cup_{i=1}^l S_i$ an open affine cover and denote by $S_I=\cap_{i\in I} S_i$.
Let $i_i:S_i\hookrightarrow\tilde S_i$ closed embeddings, with $\tilde S_i\in\Var(k)$. 
For $I\subset\left[1,\cdots l\right]$, denote by $\tilde S_I=\Pi_{i\in I}\tilde S_i$.
We then have closed embeddings $i_I:S_I\hookrightarrow\tilde S_I$ and for $J\subset I$ the following commutative diagram
\begin{equation*}
D_{IJ}=\xymatrix{ S_I\ar[r]^{i_I} & \tilde S_I \\
S_J\ar[u]^{j_{IJ}}\ar[r]^{i_J} & \tilde S_J\ar[u]^{p_{IJ}}}  
\end{equation*}
where $p_{IJ}:\tilde S_J\to\tilde S_I$ is the projection
and $j_{IJ}:S_J\hookrightarrow S_I$ is the open embedding so that $j_I\circ j_{IJ}=j_J$.
This gives the diagram of algebraic varieties $(\tilde S_I)\in\Fun(\mathcal P(\mathbb N),\Var(k))$ which
the diagram of sites $\Var(k)^{sm}/(\tilde S_I)\in\Fun(\mathcal P(\mathbb N),\Cat)$. 
Denote by $m:\tilde S_I\backslash(S_I\backslash S_J)\hookrightarrow\tilde S_I$ the open embedding.
Then $C_{fil}(\Var(k)^{sm}/(\tilde S_I))$ is the category  
\begin{itemize}
\item whose objects $(G,F)=((G_I,F),u_{IJ})$ with $(G_I,F)\in C_{fil}(\Var(k)^{sm}/\tilde S_I)$,
and $u_{IJ}:(G_I,F)\to p_{IJ*}(G_J,F)$ are morphisms
satisfying for $I\subset J\subset K$, $p_{IJ*}u_{JK}\circ u_{IJ}=u_{IK}$ in $C_{fil}(\Var(k)^{sm}/\tilde S_I)$,
\item the morphisms $m:((G,F),u_{IJ})\to((H,F),v_{IJ})$ being a family of morphisms of complexes,  
\begin{equation*}
m=(m_I:(G_I,F)\to (H_I,F))_{I\in\mathcal I}
\end{equation*}
such that $v_{IJ}\circ m_I=p_{IJ*}m_J\circ u_{IJ}$ in $C_{fil}(\Var(k)^{sm}/\tilde S_I)$.
\end{itemize}
Similarly, $C_{fil}(\Var(k)^{sm}/(\tilde S_I)^{op})$ 
is the category  
\begin{itemize}
\item whose objects $(G,F)=((G_I,F)_{I\subset\left[1,\cdots l\right]},u_{IJ})$,
with $(G_I,F)\in C_{fil}(\Var(k)^{(sm)}/\tilde S_I)$,
and $u_{IJ}:(G_J,F)\to p_{IJ}^*(G_I,F)$ for $I\subset J$, are morphisms
satisfying for $I\subset J\subset K$, $p_{JK}^*u_{IJ}\circ u_{JK}=u_{IK}$ in $C_{fil}(\Var(k)^{(sm)}/\tilde S_K)$,
\item the morphisms $m:((G,F),u_{IJ})\to((H,F),v_{IJ})$ being (see section 2.1) a family of morphisms of complexes,  
\begin{equation*}
m=(m_I:(G_I,F)\to (H_I,F))_{I\subset\left[1,\cdots l\right]}
\end{equation*}
such that $v_{IJ}\circ m_J=p_{IJ}^*m_I\circ u_{IJ}$ in $C_{fil}(\Var(k)^{(sm)}/\tilde S_J)$.
\end{itemize}

\begin{defi}
Let $S\in\Var(k)$. Let $S=\cup_{i=1}^l S_i$ an open cover and denote by $S_I=\cap_{i\in I} S_i$.
Let $i_i:S_i\hookrightarrow\tilde S_i$ closed embeddings, with $\tilde S_i\in\Var(k)$. 
We will denote by 
$C_{fil}(\Var(k)^{sm}/(S/(\tilde S_I)))\subset C_{fil}(\Var(k)^{sm}/(\tilde S_I))$ 
the full subcategory  
whose objects $(G,F)=((G_I,F)_{I\subset\left[1,\cdots l\right]},u_{IJ})$,
with $(G_I,F)\in C_{fil,S_I}(\Var(k)^{sm}/\tilde S_I)$,
and $u_{IJ}:m^*(G_I,F)\to m^*p_{IJ*}(G_J,F)$ for $I\subset J$, are $\infty$-filtered Zariski local equivalence,
\end{defi}

We now give the definition of the $\mathbb A^1$ local property :

Denote by
\begin{eqnarray*}
p_a:\Var(k)^{(sm)}/S\to\Var(k)^{(sm)}/S, \;  
X/S=(X,h)\mapsto (X\times\mathbb A^1)/S=(X\times\mathbb A^1,h\circ p_X), \\ 
(g:X/S\to X'/S)\mapsto ((g\times I_{\mathbb A^1}):X\times\mathbb A^1/S\to X'\times\mathbb A^1/S)
\end{eqnarray*}
the projection functor and again by $p_a:\Var(k)^{(sm)}/S\to\Var(k)^{(sm)}/S$
the corresponding morphism of site.

\begin{defi}\label{a1locdef}
Let $S\in\Var(k)$. Denote for short $\Var(k)^{(sm)}/S$ 
either the category $\Var(k)/S$ or the category $\Var(k)^{sm}/S$.
\begin{itemize}
\item[(i0)] A complex $F\in C(\Var(k)^{(sm)}/S)$ is said to be $\mathbb A^1$ homotopic if
$\ad(p_a^*,p_{a*})(F):F\to p_{a*}p_a^*F$ is an homotopy equivalence.
\item[(i)] A complex $F\in C(\Var(k)^{(sm)}/S)$ is said to be $\mathbb A^1$ invariant 
if for all $U/S\in\Var(k)^{(sm)}/S$,
\begin{equation*}
F(p_U):F(U/S)\to F(U\times\mathbb A^1/S) 
\end{equation*}
is a quasi-isomorphism, where $p_U:U\times\mathbb A^1\to U$ is the projection.
Obviously, if a complex $F\in C(\Var(k)^{(sm)}/S)$ is $\mathbb A^1$ homotopic then it is $\mathbb A^1$ invariant.
\item[(ii)] Let $\tau$ a topology on $\Var(k)$. 
A complex $F\in C(\Var(k)^{(sm)}/S)$ is said to be $\mathbb A^1$ local for
the topology $\tau$, if for a (hence every) $\tau$ local equivalence $k:F\to G$ with $k$ injective and
$G\in C(\Var(k)^{(sm)}/S)$ $\tau$ fibrant, e.g. $k:F\to E_{\tau}(F)$, 
$G$ is $\mathbb A^1$ invariant for all $n\in\mathbb Z$.
\item[(iii)] A morphism $m:F\to G$ with $F,G\in C(\Var(k)^{(sm)}/S)$ is said to an $(\mathbb A^1,et)$ local equivalence 
if for all $H\in C(\Var(k)^{(sm)}/S)$ which is $\mathbb A^1$ local for the etale topology
\begin{equation*}
\Hom(L(m),E_{et}(H)):\Hom(L(G),E_{et}(H))\to\Hom(L(F),E_{et}(H)) 
\end{equation*}
is a quasi-isomorphism.
\end{itemize}
\end{defi}
Denote $\square^*:=\mathbb P^*\backslash\left\{1\right\}$
\begin{itemize}
\item Let $S\in\Var(k)$. For $U/S=(U,h)\in\Var(k)^{sm}/S$, we consider 
\begin{equation*}
\square^*\times U/S=(\square^*\times U,h\circ p)\in\Fun(\Delta,\Var(k)^{sm}/S).
\end{equation*}
For $F\in C^-(\Var(k)^{sm}/S)$, it gives the complex
\begin{equation*}
C_*F\in C^-(\Var(k)^{sm}/S), U/S=(U,h)\mapsto C_*F(U/S):=\Tot F(\square^*\times U/S)
\end{equation*}
together with the canonical map $c_F:=(0,I_F):F\to C_*F$.
For $F\in C(\Var(k)^{sm}/S)$, we get
\begin{equation*}
C_*F:=\holim_n C_*F^{\leq n}\in C(\Var(k)^{sm}/S),
\end{equation*}
together with the canonical map $c_F:=(0,I_F):F\to C_*F$.
For $m:F\to G$ a morphism, with $F,G\in C(\Var(k)^{sm}/S)$,
we get by functoriality the morphism $C_*m:C_*F\to C_*G$.
\item Let $S\in\Var(k)$. For $U/S=(U,h)\in\Var(k)/S$, we consider 
\begin{equation*}
\square^*\times U/S=(\mathbb A^*\times U,h\circ p)\in\Fun(\Delta,\Var(k)/S).
\end{equation*}
For $F\in C^-(\Var(k)/S)$, it gives the complex
\begin{equation*}
C_*F\in C^-(\Var(k)/S), U/S=(U,h)\mapsto C_*F(U/S):=\Tot F(\square^*\times U/S)
\end{equation*}
together with the canonical map $c=c(F):=(0,I_F):F\to C_*F$.
For $F\in C(\Var(k)/S)$, we get
\begin{equation*}
C_*F:=\holim_n C_*F^{\leq n}\in C(\Var(k)/S),
\end{equation*}
together with the canonical map $c_F:=(0,I_F):F\to C_*F$.
For $m:F\to G$ a morphism, with $F,G\in C(\Var(k)/S)$,
we get by functoriality the morphism $C_*m:C_*F\to C_*G$.
\end{itemize}

\begin{prop}\label{ca1Var}
\begin{itemize}
\item[(i)]Let $S\in\Var(k)$.
Then for $F\in C(\Var(k)^{sm}/S)$, $C_*F$ is $\mathbb A^1$ local for the etale topology
and  $c(F):F\to C_*F$ is an equivalence $(\mathbb A^1,et)$ local.
\item[(ii)]A morphism $m:F\to G$ with $F,G\in C(\Var(k)^{(sm)}/S)$ is an $(\mathbb A^1,et)$ local equivalence
if and only if there exists 
\begin{equation*}
\left\{X_{1,\alpha}/S,\alpha\in\Lambda_1\right\},\ldots,\left\{X_{r,\alpha}/S,\alpha\in\Lambda_r\right\}
\subset\Var(k)^{(sm)}/S 
\end{equation*}
such that we have in $\Ho_{et}(C(\Var(k)^{(sm)}/S))$
\begin{eqnarray*}
\Cone(m)\xrightarrow{\sim}
\Cone(\oplus_{\alpha\in\Lambda_1}\Cone(\mathbb Z(X_{1,\alpha}\times\mathbb A^1/S)\to\mathbb Z(X_{1,\alpha}/S)) \\
\to\cdots\to\oplus_{\alpha\in\Lambda_r}\Cone(\mathbb Z(X_{r,\alpha}\times\mathbb A^1/S)\to\mathbb Z(X_{r,\alpha}/S)))
\end{eqnarray*}
\end{itemize}
\end{prop}

\begin{proof}
Standard : see Ayoub's thesis for example.
\end{proof}

\begin{defiprop}\label{projmodstr}
Let $S\in\Var(k)$.
\begin{itemize}
\item[(i)]With the weak equivalence the $(\mathbb A^1,et)$ local equivalence and 
the fibration the epimorphism with $\mathbb A^1_S$ local and etale fibrant kernels gives
a model structure on  $C(\Var(k)^{sm}/S)$ : the left bousfield localization
of the projective model structure of $C(\Var(k)^{sm}/S)$. 
We call it the projective $(\mathbb A^1,et)$ model structure.
\item[(ii)]With the weak equivalence the $(\mathbb A^1,et)$ local equivalence and 
the fibration the epimorphism with $\mathbb A^1_S$ local and etale fibrant kernels gives
a model structure on  $C(\Var(k)/S)$ : the left bousfield localization
of the projective model structure of $C(\Var(k)/S)$. 
We call it the projective $(\mathbb A^1,et)$ model structure.
\end{itemize}
\end{defiprop}

\begin{proof}
See \cite{C.D}.
\end{proof}

\begin{prop}\label{g1}
Let $g:T\to S$ a morphism with $T,S\in\Var(k)$.
\begin{itemize}
\item[(i)] The adjonction $(g^*,g_*):C(\Var(k)^{sm}/S)\leftrightarrows C(\Var(k)^{sm}/T)$
is a Quillen adjonction for the $(\mathbb A^1,et)$ projective model structure (see definition-proposition \ref{projmodstr}).
\item[(i)'] Let $h:U\to S$ a smooth morphism with $U,S\in\Var(k)$.
The adjonction $(h_{\sharp},h^*):C(\Var(k)^{sm}/U)\leftrightarrows C(\Var(k)^{sm}/S)$
is a Quillen adjonction for the $(\mathbb A^1,et)$ projective model structure.
\item[(i)''] The functor $g^*:C(\Var(k)^{sm}/S)\to C(\Var(k)^{sm}/T)$
sends quasi-isomorphism to quasi-isomorphism,
sends equivalence Zariski local to equivalence Zariski local, and equivalence etale local to equivalence etale local,
sends $(\mathbb A^1,et)$ local equivalence to $(\mathbb A^1,et)$ local equivalence.
\item[(ii)] The adjonction $(g^*,g_*):C(\Var(k)/S)\leftrightarrows C(\Var(k)/T)$
is a Quillen adjonction for the $(\mathbb A^1,et)$ projective model structure (see definition-proposition \ref{projmodstr}).
\item[(ii)'] The adjonction $(g_{\sharp},g^*):C(\Var(k)/T)\leftrightarrows C(\Var(k)/S)$
is a Quillen adjonction for the $(\mathbb A^1,et)$ projective model structure.
\item[(ii)''] The functor $g^*:C(\Var(k)/S)\to C(\Var(k)/T)$
sends quasi-isomorphism to quasi-isomorphism,
sends equivalence Zariski local to equivalence Zariski local, and equivalence etale local to equivalence etale local,
sends $(\mathbb A^1,et)$ local equivalence to $(\mathbb A^1,et)$ local equivalence.
\end{itemize}
\end{prop}

\begin{proof}
Standard : see \cite{C.D} for example.
\end{proof}

\begin{prop}\label{rho1}
Let $S\in\Var(k)$. 
\begin{itemize}
\item[(i)] The adjonction $(\rho_S^*,\rho_{S*}):C(\Var(k)^{sm}/S)\leftrightarrows C(\Var(k)/S)$
is a Quillen adjonction for the $(\mathbb A^1,et)$ projective model structure.
\item[(ii)]The functor $\rho_{S*}:C(\Var(k)/S)\to C(\Var(k)^{sm}/S)$
sends quasi-isomorphism to quasi-isomorphism,
sends equivalence Zariski local to equivalence Zariski local, and equivalence etale local to equivalence etale local,
sends $(\mathbb A^1,et)$ local equivalence to $(\mathbb A^1,et)$ local equivalence.
\end{itemize}
\end{prop}

\begin{proof}
Standard : see \cite{C.D} for example.
\end{proof}

Let $S\in\Var(k)$. Let $S=\cup_{i=1}^l S_i$ an open affine cover and denote by $S_I=\cap_{i\in I} S_i$.
Let $i_i:S_i\hookrightarrow\tilde S_i$ closed embeddings, with $\tilde S_i\in\Var(k)$.
\begin{itemize}
\item For $(G_I,K_{IJ})\in C(\Var(k)^{(sm)}/(\tilde S_I)^{op})$ and 
$(H_I,T_{IJ})\in C(\Var(k)^{(sm)}/(\tilde S_I))$, we denote
\begin{eqnarray*}
\mathcal Hom((G_I,K_{IJ}),(H_I,T_{IJ})):=(\mathcal Hom(G_I,H_I),u_{IJ}((G_I,K_{IJ}),(H_I,T_{IJ})))
\in C(\Var(k)^{(sm)}/(\tilde S_I))
\end{eqnarray*}
with
\begin{eqnarray*}
u_{IJ}((G_I,K_{IJ})(H_I,T_{IJ})):\mathcal Hom(G_I,H_I) \\
\xrightarrow{\ad(p_{IJ}^*,p_{IJ*})(-)}p_{IJ*}p_{IJ}^*\mathcal Hom(G_I,H_I)
\xrightarrow{T(p_{IJ},hom)(-,-)}p_{IJ*}\mathcal Hom(p_{IJ}^*G_I,p_{IJ}^*H_I) \\
\xrightarrow{\mathcal Hom(p_{IJ}^*G_I,T_{IJ})}p_{IJ*}\mathcal Hom(p_{IJ}^*G_I,H_J)
\xrightarrow{\mathcal Hom(K_{IJ},H_J)}p_{IJ*}\mathcal Hom(G_J,H_J).
\end{eqnarray*}
This gives in particular the functor
\begin{eqnarray*}
\mathbb D^0_{(\tilde S_I)}:C(\Var(k)^{(sm)}/(\tilde S_I)^{op})\to C(\Var(k)^{(sm)}/(\tilde S_I)), \\
(H_I,T_{IJ})\mapsto(\mathbb D^0_{\tilde S_I}LH_I,T^d_{IJ}):=\mathcal Hom((LH_I,T_{IJ}),(E_{et}(\mathbb Z_{\tilde S_I}),I_{IJ})).
\end{eqnarray*}
\item For $(G_I,K_{IJ})\in C(\Var(k)^{(sm)}/(\tilde S_I))$ and 
$(H_I,T_{IJ})\in C(\Var(k)^{(sm)}/(\tilde S_I)^{op})$, we denote
\begin{eqnarray*}
\mathcal Hom((G_I,K_{IJ}),(H_I,T_{IJ})):=(\mathcal Hom(G_I,H_I),u_{IJ}((G_I,K_{IJ}),(H_I,T_{IJ})))
\in C(\Var(k)^{(sm)}/(\tilde S_I)^{op})
\end{eqnarray*}
with
\begin{eqnarray*}
u_{IJ}((G_I,K_{IJ})(H_I,T_{IJ})):\mathcal Hom(G_J,H_J) \\
\xrightarrow{\mathcal Hom(\ad(p_{IJ}^*,p_{IJ*})(G_J),H_J)}\mathcal Hom(p_{IJ}^*p_{IJ*}G_J,H_J) 
\xrightarrow{\mathcal Hom(p_{IJ}^*K_{IJ},H_J)}\mathcal Hom(p_{IJ}^*G_I,H_J) \\
\xrightarrow{\mathcal Hom(p_{IJ}^*G_I,T_{IJ})}\mathcal Hom(p_{IJ}^*G_I,p_{IJ}^*H_I) 
\xrightarrow{T(p_{IJ},\hom)(-,-)^{-1}}p_{IJ}^*\mathcal Hom(G_I,H_I).
\end{eqnarray*}
This gives in particular the functor
\begin{eqnarray*}
\mathbb D^0_{(\tilde S_I)}:C(\Var(k)^{(sm)}/(\tilde S_I))\to C(\Var(k)^{(sm)}/(\tilde S_I)^{op}), \\
(H_I,T_{IJ})\mapsto(\mathbb D^0_{\tilde S_I}LH_I,T^d_{IJ}):=(\mathcal Hom((LH_I,T_{IJ}),(E_{et}(\mathbb Z_{\tilde S_I}),I_{IJ}))).
\end{eqnarray*}
\end{itemize}

\begin{defi}\label{a1locdefIJ}
Let $S\in\Var(k)$. Let $S=\cup_{i=1}^l S_i$ an open affine cover and denote by $S_I=\cap_{i\in I} S_i$.
Let $i_i:S_i\hookrightarrow\tilde S_i$ closed embeddings, with $\tilde S_i\in\Var(k)$.
\begin{itemize}
\item[(i0)]A complex $(F_I,u_{IJ})\in C(\Var(k)^{(sm)}/(\tilde S_I))$ is said to be $\mathbb A^1$ homotopic if 
$\ad(p_a^*,p_{a*})((F_I,u_{IJ})):(F_I,u_{IJ})\to p_{a*}p_a^*(F_I,u_{IJ})$ is an homotopy equivalence.
\item[(i)] A complex  $(F_I,u_{IJ})\in C(\Var(k)^{(sm)}/(\tilde S_I))$ is said to be $\mathbb A^1$ invariant 
if for all $(X_I/\tilde S_I,s_{IJ})\in\Var(k)^{(sm)}/(\tilde S_I)$ 
\begin{equation*}
(F_I(p_{X_I})):(F_I(X_I/\tilde S_I),F_J(s_{IJ})\circ u_{IJ}(-))\to 
(F_I(X_I\times\mathbb A^1/\tilde S_I),F_J(s_{IJ}\times I)\circ u_{IJ}(-)) 
\end{equation*}
is a quasi-isomorphism, where $p_{X_I}:X_I\times\mathbb A^1\to X_I$ are the projection,
and $s_{IJ}:X_I\times\tilde S_{J\backslash I}/\tilde S_J\to X_J/\tilde S_J$.
Obviously a complex $(F_I,u_{IJ})\in C(\Var(k)^{(sm)}/(\tilde S_I))$ is $\mathbb A^1$ invariant
if and only if all the $F_I$ are $\mathbb A^1$ invariant.
\item[(ii)]Let $\tau$ a topology on $\Var(k)$. 
A complex $F=(F_I,u_{IJ})\in C(\Var(k)^{(sm)}/(\tilde S_I))$ is said to be $\mathbb A^1$ local 
for the $\tau$ topology induced on $\Var(k)/(\tilde S_I)$, 
if for an (hence every) $\tau$ local equivalence $k:F\to G$ with $k$ injective and 
$G=(G_I,v_{IJ})\in C(\Var(k)^{(sm)}/(\tilde S_I))$ $\tau$ fibrant,
e.g. $k:(F_I,u_{IJ})\to (E_{\tau}(F_I),E(u_{IJ}))$, $G$ is $\mathbb A^1$ invariant.
\item[(iii)] A morphism $m=(m_I):(F_I,u_{IJ})\to (G_I,v_{IJ})$ with 
$(F_I,u_{IJ}),(G_I,v_{IJ})\in C(\Var(k)^{(sm)}/(\tilde S_I))$ 
is said to be an $(\mathbb A^1,et)$ local equivalence 
if for all $(H_I,w_{IJ})\in C(\Var(k)^{(sm)}/(\tilde S_I))$ which is $\mathbb A^1$ local for the etale topology
\begin{eqnarray*}
(\Hom(L(m_I),E_{et}(H_I))):\Hom(L(G_I,v_{IJ}),E_{et}(H_I,w_{IJ}))\to\Hom(L(F_I,u_{IJ}),E_{et}(H_I,w_{IJ})) 
\end{eqnarray*}
is a quasi-isomorphism (of complexes of abelian groups).
Obviously, if a morphism $m=(m_I):(F_I,u_{IJ})\to (G_I,v_{IJ})$ with 
$(F_I,u_{IJ}),(G_I,u_{IJ})\in C(\Var(k)^{(sm)}/(\tilde S_I))$ is an $(\mathbb A^1,et)$ local equivalence, 
then all the $m_I:F_I\to G_I$ are $(\mathbb A^1,et)$ local equivalence.
\item[(iv)] A morphism $m=(m_I):(F_I,u_{IJ})\to (G_I,v_{IJ})$ with 
$(F_I,u_{IJ}),(G_I,v_{IJ})\in C(\Var(k)^{(sm)}/(\tilde S_I)^{op})$ is said to be an $(\mathbb A^1,et)$ local equivalence 
if for all $(H_I,w_{IJ})\in C(\Var(k)^{(sm)}/(\tilde S_I)^{op})$ which is $\mathbb A^1$ local for the etale topology
\begin{eqnarray*}
(\Hom(L(m_I),E_{et}(H_I))):\Hom(L(G_I,v_{IJ}),E_{et}(H_I,w_{IJ}))\to\Hom(L(F_I,u_{IJ}),E_{et}(H_I,w_{IJ})) 
\end{eqnarray*}
is a quasi-isomorphism (of complexes of abelian groups).
Obviously, if a morphism $m=(m_I):(F_I,u_{IJ})\to (G_I,v_{IJ})$ with 
$(F_I,u_{IJ}),(G_I,u_{IJ})\in C(\Var(k)^{(sm)}/(\tilde S_I)^{op})$ is an $(\mathbb A^1,et)$ local equivalence, 
then all the $m_I:F_I\to G_I$ are $(\mathbb A^1,et)$ local equivalence and for all 
$(H_I,w_{IJ})\in C(\Var(k)^{(sm)}/(\tilde S_I))$,
\begin{eqnarray*}
(\Hom(L(m_I),E_{et}(H_I))):\Hom(L(G_I,v_{IJ}),E_{et}(H_I,w_{IJ}))\to\Hom(L(F_I,u_{IJ}),E_{et}(H_I,w_{IJ})) 
\end{eqnarray*}
is a quasi-isomorphism (of diagrams of complexes of abelian groups).
\end{itemize}
\end{defi}

\begin{prop}\label{ca1VarIJ}
Let $S\in\Var(k)$. Let $S=\cup_{i=1}^l S_i$ an open affine cover and denote by $S_I=\cap_{i\in I} S_i$.
Let $i_i:S_i\hookrightarrow\tilde S_i$ closed embeddings, with $\tilde S_i\in\Var(k)$.
\begin{itemize}
\item[(i)]A morphism $m:F\to G$ with $F,G\in C(\Var(k)^{(sm)}/(\tilde S_I))$ 
is an $(\mathbb A^1,et)$ local equivalence if and only if there exists 
\begin{eqnarray*}
\left\{(X_{1,\alpha,I}/\tilde S_I,u^1_{IJ}),\alpha\in\Lambda_1\right\},\ldots,
\left\{(X_{r,\alpha,I}/\tilde S_I,u^r_{IJ}),\alpha\in\Lambda_r\right\}
\subset\Var(k)^{(sm)}/(\tilde S_I)
\end{eqnarray*}
with $u^l_{IJ}:X_{l,\alpha,I}\times\tilde S_{J\backslash I}/\tilde S_J\to X_{l,\alpha,J}/\tilde S_J$,
such that we have in $\Ho_{et}(C(\Var(k)^{(sm)}/(\tilde S_I)))$
\begin{eqnarray*}
\Cone(m)\xrightarrow{\sim}\Cone( \oplus_{\alpha\in\Lambda_1} 
\Cone((\mathbb Z(X_{1,\alpha,I}\times\mathbb A^1/\tilde S_I),\mathbb Z(u_{IJ}^1\times I)) 
\to(\mathbb Z(X_{1,\alpha,I}/\tilde S_I),\mathbb Z(u_{IJ}^1))) \\
\to\cdots\to\oplus_{\alpha\in\Lambda_r}
\Cone((\mathbb Z(X_{r,\alpha,I}\times\mathbb A^1/\tilde S_I),\mathbb Z(u_{IJ}^r\times I)) 
\to(\mathbb Z(X_{r,\alpha,I}/\tilde S_I),\mathbb Z(u^r_{IJ}))))
\end{eqnarray*}
\item[(ii)]A morphism $m:F\to G$ with $F,G\in C(\Var(k)^{(sm)}/(\tilde S_I)^{op})$ 
is an $(\mathbb A^1,et)$ local equivalence if and only if there exists 
\begin{eqnarray*}
\left\{(X_{1,\alpha,I}/\tilde S_I,u^1_{IJ}),\alpha\in\Lambda_1\right\},\ldots,
\left\{(X_{r,\alpha,I}/\tilde S_I,u^r_{IJ}),\alpha\in\Lambda_r\right\}
\subset\Var(k)^{(sm)}/(\tilde S_I)^{op}
\end{eqnarray*}
with $u^l_{IJ}:X_{l,\alpha,J}/\tilde S_J\to X_{l,\alpha,I}\times\tilde S_{J\backslash I}/\tilde S_J$,
such that we have in $\Ho_{et}(C(\Var(k)^{(sm)}/(\tilde S_I)^{op}))$
\begin{eqnarray*}
\Cone(m)\xrightarrow{\sim}\Cone( \oplus_{\alpha\in\Lambda_1} 
\Cone((\mathbb Z(X_{1,\alpha,I}\times\mathbb A^1/\tilde S_I),\mathbb Z(u_{IJ}^1\times I)) 
\to(\mathbb Z(X_{1,\alpha,I}/\tilde S_I),\mathbb Z(u_{IJ}^1))) \\
\to\cdots\to\oplus_{\alpha\in\Lambda_r}
\Cone((\mathbb Z(X_{r,\alpha,I}\times\mathbb A^1/\tilde S_I),\mathbb Z(u_{IJ}^r\times I)) 
\to(\mathbb Z(X_{r,\alpha,I}/\tilde S_I),\mathbb Z(u^r_{IJ}))))
\end{eqnarray*}
\end{itemize}
\end{prop}

\begin{proof}
Standard. See Ayoub's thesis for example.
\end{proof}

\begin{itemize}
\item For $f:X\to S$ a morphism with $X,S\in\Var(k)$, we denote as usual (see \cite{C.D} for example),
$\mathbb Z^{tr}(X/S)\in\PSh(\Var(k)/S)$ the presheaf given by
\begin{itemize}
\item for $X'/S\in\Var(k)/S$, with $X'$ irreducible, 
$\mathbb Z^{tr}(X/S)(X'/S):=\mathcal Z^{fs/X}(X'\times_S X)\subset\mathcal Z_{d_{X'}}(X'\times_S X)$
which consist of algebraic cycles $\alpha=\sum_in_i\alpha_i\in\mathcal Z_{d_{X'}}(X'\times_S X)$ such that,
denoting $\supp(\alpha)=\cup_i\alpha_i\subset X'\times_S X$ its support and $f':X'\times_S X\to X'$ the projection,
$f'_{|\supp(\alpha)}:\supp(\alpha)\to X'$ is finite surjective, 
\item for $g:X_2/S\to X_1/S$ a morphism, with $X_1/S,X_2/S\in\Var(k)/S$,
\begin{equation*}
\mathbb Z^{tr}(X/S)(g):\mathbb Z^{tr}(X/S)(X_1/S)\to\mathbb Z^{tr}(X/S)(X_2/S), \;
\alpha\mapsto (g\times I)^{-1}(\alpha)
\end{equation*}
with $g\times I:X_2\times_S X\to X_1\times_S X$, noting that, by base change,
$f_{2|\supp((g\times I)^{-1}(\alpha))}:\supp((g\times I)^{-1}(\alpha))\to X_2$ is finite surjective,
$f_2:X_2\times_S X\to X_2$ being the projection.
\end{itemize}
\item For $f:X\to S$ a morphism with $X,S\in\Var(k)$ and $r\in\mathbb N$, 
we denote as usual (see \cite{C.D} for example),
$\mathbb Z^{equir}(X/S)\in\PSh(\Var(k)/S)$ the presheaf given by
\begin{itemize}
\item for $X'/S\in\Var(k)/S$, with $X'$ irreducible, 
$\mathbb Z^{equir}(X/S)(X'/S):=\mathcal Z^{equir/X}(X'\times_S X)\subset\mathcal Z_{d_{X'}}(X'\times_S X)$
which consist of algebraic cycles $\alpha=\sum_in_i\alpha_i\in\mathcal Z_{d_{X'}}(X'\times_S X)$ such that,
denoting $\supp(\alpha)=\cup_i\alpha_i\subset X'\times_S X$ its support and $f':X'\times_S X\to X'$ the projection,
$f'_{|\supp(\alpha)}:\supp(\alpha)\to X'$ is dominant, with fibers either empty or of dimension $r$, 
\item for $g:X_2/S\to X_1/S$ a morphism, with $X_1/S,X_2/S\in\Var(k)/S$,
\begin{equation*}
\mathbb Z^{equir}(X/S)(g):\mathbb Z^{equir}(X/S)(X_1/S)\to\mathbb Z^{equir}(X/S)(X_2/S), \;
\alpha\mapsto (g\times I)^{-1}(\alpha)
\end{equation*}
with $g\times I:X_2\times_S X\to X_1\times_S X$, noting that, by base change,
$f_{2|\supp((g\times I)^{-1}(\alpha))}:\supp((g\times I)^{-1}(\alpha))\to X_2$ is obviously dominant,
with fibers either empty or of dimension $r$, $f_2:X_2\times_S X\to X_2$ being the projection.
\end{itemize}
\item Let $S\in\Var(k)$. We denote by 
$\mathbb Z_S(d):=\mathbb Z^{equi0}(S\times\mathbb A^d/S)[-2d]$
the Tate twist. For $F\in C(\Var(k)/S)$, we denote by $F(d):=F\otimes\mathbb Z_S(d)$.
\end{itemize}

For $S\in\Var(k)$, let $\Cor(\Var(k)^{sm}/S)$ be the category 
\begin{itemize}
\item whose objects are smooth morphisms 
$U/S=(U,h)$, $h:U\to S$ with $U\in\Var(k)$, 
\item whose morphisms $\alpha:U/S=(U,h_1)\to V/S=(V,h_2)$ 
is finite correspondence that is $\alpha\in\oplus_i\mathcal Z^{fs}(U_i\times_S V)$, 
where $U=\sqcup_i U_i$, with $U_i$ connected (hence irreducible by smoothness), and $\mathcal Z^{fs}(U_i\times_S V)$ 
is the abelian group of cycle finite and surjective over $U_i$. 
\end{itemize}
We denote by 
$\Tr(S):\Cor(\Var(k)^{sm}/S)\to\Var(k)^{sm}/S$ 
the morphism of site
given by the inclusion functor
$\Tr(S):\Var(k)^{sm}/S\hookrightarrow\Cor(\Var(k)^{sm}/S)$
It induces an adjonction
\begin{equation*}
(\Tr(S)^*\Tr(S)_*):C(\Var(k)^{sm}/S)\leftrightarrows C(\Cor(\Var(k)^{sm}/S))
\end{equation*}
A complex of preheaves $G\in C(\Var(k)^{sm}/S)$ is said to admit transferts
if it is in the image of the embedding
\begin{equation*}
\Tr(S)_*:C(\Cor(\Var(k)^{sm}/S)\hookrightarrow C(\Var(k)^{sm}/S),
\end{equation*}
that is $G=\Tr(S)_*\Tr(S)^*G$.

We will use to compute the algebraic Gauss-Manin realization functor the following

\begin{thm}\label{DDADM}
Let $\phi:F^{\bullet}\to G^{\bullet}$ an etale local equivalence with $F^{\bullet},G^{\bullet}\in C(\Var(k)^{sm}/S)$.
If $F^{\bullet}$ and $G^{\bullet}$ are $\mathbb A^1$ local and admit tranferts 
then $\phi:F^{\bullet}\to G^{\bullet}$ is a Zariski local equivalence.
Hence if $F\in C(\Var(k)^{sm}/S)$ is $\mathbb A^1$ local and admits transfert 
\begin{equation*}
k:E_{zar}(F)\to E_{et}(E_{zar}(F))=E_{et}(F) 
\end{equation*}
is a Zariski local equivalence.
\end{thm}

\begin{proof}
See \cite{C.D}.
\end{proof}

\subsection{Presheaves on the big Zariski site or the big etale site of pairs}

We recall the definition given in subsection 5.1 :
For $S\in\Var(k)$, $\Var(k)^2/S:=\Var(k)^2/(S,S)$ is by definition (see subsection 2.1)
the category whose set of objects is 
\begin{eqnarray*}
(\Var(k)^2/S)^0:=
\left\{((X,Z),h), h:X\to S, \; Z\subset X \; \mbox{closed} \;\right\}\subset\Var(k)/S\times\Top
\end{eqnarray*}
and whose set of morphisms between $(X_1,Z_1)/S=((X_1,Z_1),h_1),(X_1,Z_1)/S=((X_2,Z_2),h_2)\in\Var(k)^2/S$
is the subset
\begin{eqnarray*}
\Hom_{\Var(k)^2/S}((X_1,Z_1)/S,(X_2,Z_2)/S):= \\
\left\{(f:X_2\to X_2), \; \mbox{s.t.} \; h_1\circ f=h_2 \; \mbox{and} \; Z_1\subset f^{-1}(Z_2)\right\} 
\subset\Hom_{\Var(k)}(X_1,X_2)
\end{eqnarray*}
The category $\Var(k)^2$ admits fiber products : $(X_1,Z_1)\times_{(S,Z)}(X_2,Z_2)=(X_1\times_S X_2,Z_1\times_Z Z_2)$.
In particular, for $f:T\to S$ a morphism with $S,T\in\Var(k)$, we have the pullback functor
\begin{equation*}
P(f):\Var(k)^2/S\to\Var(k)^2/T, P(f)((X,Z)/S):=(X_T,Z_T)/T, P(f)(g):=(g\times_S f)
\end{equation*}
and we note again $P(f):\Var(k)^2/T\to\Var(k)^2/S$ the corresponding morphism of sites.

We will consider in the construction of the filtered De Rham realization functor the
full subcategory $\Var(k)^{2,sm}/S\subset\Var(k)^2/S$ such that the first factor is a smooth morphism :
We will also consider, in order to obtain a complex of D modules in the construction of the filtered De Rham realization functor,
the restriction to the full subcategory $\Var(k)^{2,pr}/S\subset\Var(k)^2/S$ 
such that the first factor is a projection :

\begin{defi}\label{PVar12S}
\begin{itemize}
\item[(i)]Let $S\in\Var(k)$. We denote by
\begin{equation*}
\rho_S:\Var(k)^{2,sm}/S\hookrightarrow\Var(k)^2/S 
\end{equation*}
the full subcategory
consisting of the objects $(U,Z)/S=((U,Z),h)\in\Var(k)^2/S$ such that the morphism $h:U\to S$ is smooth.
That is, $\Var(k)^{2,sm}/S$ is the category  
\begin{itemize}
\item whose objects are $(U,Z)/S=((U,Z),h)$, with $U\in\Var(k)$, $Z\subset U$ a closed subset, 
and $h:U\to S$ a smooth morphism,
\item whose morphisms $g:(U,Z)/S=((U,Z),h_1)\to (U',Z')/S=((U',Z'),h_2)$ 
is a morphism $g:U\to U'$ of complex algebraic varieties such that $Z\subset g^{-1}(Z')$ and  $h_2\circ g=h_1$. 
\end{itemize}
We denote again $\rho_S:\Var(k)^2/S\to\Var(k)^{2,sm}/S$ the associated morphism of site.We have 
\begin{equation*}
r^s(S):\Var(k)^2\xrightarrow{r(S):=r(S,S)}\Var(k)^2/S\xrightarrow{\rho_S}\Var(k)^{2,sm}/S
\end{equation*}
the composite morphism of site.
\item[(ii)]Let $S\in\Var(k)$. We will consider the full subcategory 
\begin{equation*}
\mu_S:\Var(k)^{2,pr}/S\hookrightarrow\Var(k)^2/S
\end{equation*}
whose subset of object consist of those whose morphism is a projection to $S$ : 
\begin{eqnarray*}
(\Var(k)^{2,pr}/S)^0:=\left\{((Y\times S,X),p), \; Y\in\Var(k), \;
 p:Y\times S\to S \; \mbox{the projection}\right\}\subset(\Var(k)^2/S)^0.
\end{eqnarray*}
\item[(iii)]We will consider the full subcategory 
\begin{equation*}
\mu_S:(\Var(k)^{2,smpr}/S)\hookrightarrow\Var(k)^{2,sm}/S
\end{equation*}
whose subset of object consist of those whose morphism is a smooth projection to $S$ : 
\begin{eqnarray*}
(\Var(k)^{2,smpr}/S)^0:=\left\{((Y\times S,X),p), \; Y\in\SmVar(k), \;
 p:Y\times S\to S \; \mbox{the projection}\right\}\subset(\Var(k)^2/S)^0
\end{eqnarray*}
\end{itemize}
\end{defi}
For $f:T\to S$ a morphism with $T,S\in\Var(k)$, we have by definition, the following commutative diagram of sites
\begin{equation}\label{muf}
\xymatrix{\Var(k)^2/T\ar[rr]^{\mu_T}\ar[dd]_{P(f)}\ar[rd]^{\rho_T} & \, & 
\Var(k)^{2,pr}/T\ar[dd]^{P(f)}\ar[rd]^{\rho_T} & \, \\
\, & \Var(k)^{2,sm}/T\ar[rr]^{\mu_T}\ar[dd]_{P(f)} & \, & \Var(k)^{2,smpr}/T\ar[dd]^{P(f)} \\
\Var(k)^2/S\ar[rr]^{\mu_S}\ar[rd]^{\rho_S} & \, & \Var(k)^{2,pr}/S\ar[rd]^{\rho_S} & \, \\
\, & \Var(k)^{2,sm}/S\ar[rr]^{\mu_S} & \, & \Var(k)^{2,smpr}/S}.
\end{equation}

Recall we have (see subsection 2.1), for $S\in\Var(k)$, the graph functor 
\begin{eqnarray*}
\Gr_S^{12}:\Var(k)/S\to\Var(k)^{2,pr}/S, \; X/S\mapsto\Gr_S^{12}(X/S):=(X\times S,X)/S, \\
(g:X/S\to X'/S)\mapsto\Gr_S^{12}(g):=(g\times I_S:(X\times S,X)\to(X'\times S,X'))
\end{eqnarray*}
Note that $\Gr_S^{12}$ is fully faithfull.
For $f:T\to S$ a morphism with $T,S\in\Var(k)$, we have by definition, the following commutative diagram of sites
\begin{equation}\label{Grf}
\xymatrix{\Var(k)^{2,pr}/T\ar[rr]^{\Gr_T^{12}}\ar[dd]_{P(f)}\ar[rd]^{\rho_T} & \, & 
\Var(k)/T\ar[dd]^{P(f)}\ar[rd]^{\rho_T} & \, \\
\, & \Var(k)^{2,smpr}/T\ar[rr]^{\Gr_T^{12}}\ar[dd]_{P(f)} & \, & \Var(k)^{sm}/T\ar[dd]^{P(f)} \\
\Var(k)^{2,pr}/S\ar[rr]^{\Gr_S^{12}}\ar[rd]^{\rho_S} & \, & \Var(k)/S\ar[rd]^{\rho_S} & \, \\
\, & \Var(k)^{2,sm}/S\ar[rr]^{\Gr_S^{12}} & \, & \Var(k)^{sm}/S}.
\end{equation}
where we recall that $P(f)((X,Z)/S):=((X_T,Z_T)/T)$, since smooth morphisms are preserved by base change.

\begin{itemize}
\item As usual, we denote by
\begin{equation*}
(f^*,f_*):=(P(f)^*,P(f)_*):C(\Var(k)^{2,(sm)}/S)\to C(\Var(k)^{2,(sm)}/T)
\end{equation*}
the adjonction induced by $P(f):\Var(k)^{2,(sm)}/T\to \Var(k)^{2,(sm)}/S$.
Since the colimits involved in the definition of $f^*=P(f)^*$ are filtered, $f^*$ also preserve monomorphism. 
Hence, we get an adjonction
\begin{equation*}
(f^*,f_*):C_{fil}(\Var(k)^{2,(sm)}/S)\leftrightarrows C_{fil}(\Var(k)^{2,(sm)}/T), \; f^*(G,F):=(f^*G,f^*F)
\end{equation*}
For $S\in\Var(k)$, 
we denote by $\mathbb Z_S:=\mathbb Z((S,S)/(S,S))\in\PSh(\Var(k)^{2,(sm)}/S)$ the constant presheaf.
By Yoneda lemma, we have for $F\in C(\Var(k)^{2,(sm)}/S)$, $\mathcal Hom(\mathbb Z_S,F)=F$.
\item As usual, we denote by
\begin{equation*}
(f^*,f_*):=(P(f)^*,P(f)_*):C(\Var(k)^{2,(sm)pr}/S)\to C(\Var(k)^{2,(sm)pr}/T)
\end{equation*}
the adjonction induced by $P(f):\Var(k)^{2,(sm)pr}/T\to \Var(k)^{2,(sm)pr}/S$.
Since the colimits involved in the definition of $f^*=P(f)^*$ are filtered, $f^*$ also preserve monomorphism. 
Hence, we get an adjonction
\begin{equation*}
(f^*,f_*):C_{fil}(\Var(k)^{2,(sm)pr}/S)\leftrightarrows C_{fil}(\Var(k)^{2,(sm)pr}/T), \; f^*(G,F):=(f^*G,f^*F)
\end{equation*}
For $S\in\Var(k)$, 
we denote by $\mathbb Z_S:=\mathbb Z((S,S)/(S,S))\in\PSh(\Var(k)^{2,sm}/S)$ the constant presheaf.
By Yoneda lemma, we have for $F\in C(\Var(k)^{2,sm}/S)$, $\mathcal Hom(\mathbb Z_S,F)=F$.
\end{itemize}

\begin{itemize}
\item For $h:U\to S$ a smooth morphism with $U,S\in\Var(k)$, 
$P(h):\Var(k)^{2,sm}/S\to\Var(k)^{2,sm}/U$ admits a left adjoint
\begin{equation*}
C(h):\Var(k)^{2,sm}/U\to\Var(k)^{2,sm}/S, \; C(h)((U',Z'),h')=((U',Z'),h\circ h').
\end{equation*}
Hence $h^*:C(\Var(k)^{2,sm}/S)\to C(\Var(k)^{2,sm}/U)$ admits a left adjoint
\begin{eqnarray*}
h_{\sharp}:C(\Var(k)^{2,sm}/U)\to C(\Var(k)^{2,sm}/S), \\
F\mapsto (h_{\sharp}F:((U,Z),h_0)\mapsto\lim_{((U',Z'),h\circ h')\to ((U,Z),h_0)} F((U',Z')/U))
\end{eqnarray*}
\item For $h:X\to S$ a morphism with $X,S\in\Var(k)$, 
$P(h):\Var(k)^2/S\to\Var(k)^2/X$ admits a left adjoint
\begin{equation*}
C(h):\Var(k)^2/X\to\Var(k)^2/S, \; C(h)((X',Z'),h')=((X',Z'),h\circ h').
\end{equation*}
Hence $h^*:C(\Var(k)^2/S)\to C(\Var(k)^2/X)$ admits a left adjoint
\begin{eqnarray*}
h_{\sharp}:C(\Var(k)^2/X)\to C(\Var(k)^{2,sm}/S), \\
F\mapsto (h_{\sharp}F:((X,Z),h_0)\mapsto\lim_{((X',Z'),h\circ h')\to ((X,Z),h_0)} F((X',Z')/X))
\end{eqnarray*}
\item For $p:Y\times S\to S$ a projection with $Y,S\in\Var(k)$ with $Y$ smooth, 
$P(p):\Var(k)^{2,smpr}/S\to\Var(k)^{2,smpr}/Y\times S$ admits a left adjoint
\begin{eqnarray*}
C(p):\Var(k)^{2,smpr}/Y\times S\to\Var(k)^{2,smpr}/S, \\ 
C(p)((Y'\times S,Z'),p')=((Y'\times S,Z'),p\circ p').
\end{eqnarray*}
Hence $p^*:C(\Var(k)^{2,smpr}/S)\to C(\Var(k)^{2,smpr}/Y\times S)$ admits a left adjoint
\begin{eqnarray*}
p_{\sharp}:C(\Var(k)^{2,smpr}/Y\times S)\to C(\Var(k)^{2,smpr}/S), \\
F\mapsto (p_{\sharp}F:((Y_0\times S,Z),p_0)\mapsto
\lim_{((Y'\times Y\times S,Z'),p\circ p')\to ((Y_0\times S,Z),p_0)} F((Y'\times Y\times S,Z')/Y\times S))
\end{eqnarray*}
\item For $p:Y\times S\to S$ a projection with $Y,S\in\Var(k)$, 
$P(p):\Var(k)^{2,pr}/S\to\Var(k)^{2,pr}/Y\times S$ admits a left adjoint
\begin{equation*}
C(p):\Var(k)^{2,pr}/Y\times S\to\Var(k)^{2,pr}/S, \; C(p)((Y'\times S,Z'),p')=((Y'\times S,Z'),p\circ p').
\end{equation*}
Hence $p^*:C(\Var(k)^{2,pr}/S)\to C(\Var(k)^{2,pr}/Y\times S)$ admits a left adjoint
\begin{eqnarray*}
p_{\sharp}:C(\Var(k)^{2,pr}/Y\times S)\to C(\Var(k)^{2,pr}/S), \\
F\mapsto (p_{\sharp}F:((Y_0\times S,Z),p_0)\mapsto
\lim_{((Y'\times Y\times S,Z'),p\circ p')\to ((Y_0\times S,Z),p_0)} F((Y'\times Y\times S,Z')/Y\times S))
\end{eqnarray*}
\end{itemize}

Let $S_{\bullet}\in\Fun(\mathcal I,\Var(k))$ with $\mathcal I\in\Cat$, a diagram of algebraic varieties.
It gives the diagram of sites $\Var(k)^2/S_{\bullet}\in\Fun(\mathcal I,\Cat)$.  
\begin{itemize}
\item Then $C_{fil}(\Var(k)^{2,(sm)}/S_{\bullet})$ is the category  
\begin{itemize}
\item whose objects $(G,F)=((G_I,F)_{I\in\mathcal I},u_{IJ})$,
with $(G_I,F)\in C_{fil}(\Var(k)^{2,(sm)}/S_I)$,
and $u_{IJ}:(G_I,F)\to r_{IJ*}(G_J,F)$ for $r_{IJ}:I\to J$, denoting again $r_{IJ}:S_I\to S_J$, are morphisms
satisfying for $I\to J\to K$, $r_{IJ*}u_{JK}\circ u_{IJ}=u_{IK}$ in $C_{fil}(\Var(k)^{2,(sm)}/S_I)$,
\item the morphisms $m:((G,F),u_{IJ})\to((H,F),v_{IJ})$ being (see section 2.1) a family of morphisms of complexes,  
\begin{equation*}
m=(m_I:(G_I,F)\to (H_I,F))_{I\in\mathcal I}
\end{equation*}
such that $v_{IJ}\circ m_I=p_{IJ*}m_J\circ u_{IJ}$ in $C_{fil}(\Var(k)^{2,(sm)}/S_I)$.
\end{itemize}
\item Then $C_{fil}(\Var(k)^{2,(sm)pr}/S_{\bullet})$ is the category  
\begin{itemize}
\item whose objects $(G,F)=((G_I,F)_{I\in\mathcal I},u_{IJ})$,
with $(G_I,F)\in C_{fil}(\Var(k)^{2,(sm)pr}/S_I)$,
and $u_{IJ}:(G_I,F)\to r_{IJ*}(G_J,F)$ for $r_{IJ}:I\to J$, denoting again $r_{IJ}:S_I\to S_J$, are morphisms
satisfying for $I\to J\to K$, $r_{IJ*}u_{JK}\circ u_{IJ}=u_{IK}$ in $C_{fil}(\Var(k)^{2,(sm)}/S_I)$,
\item the morphisms $m:((G,F),u_{IJ})\to((H,F),v_{IJ})$ being (see section 2.1) a family of morphisms of complexes,  
\begin{equation*}
m=(m_I:(G_I,F)\to (H_I,F))_{I\in\mathcal I}
\end{equation*}
such that $v_{IJ}\circ m_I=p_{IJ*}m_J\circ u_{IJ}$ in $C_{fil}(\Var(k)^{2,(sm)pr}/S_I)$.
\end{itemize}
\end{itemize}
For $s:\mathcal I\to\mathcal J$ a functor, with $\mathcal I,\mathcal J\in\Cat$, and
$f_{\bullet}:T_{\bullet}\to S_{s(\bullet)}$ a morphism with 
$T_{\bullet}\in\Fun(\mathcal J,\Var(k))$ and $S_{\bullet}\in\Fun(\mathcal I,\Var(k))$, 
we have by definition, the following commutative diagrams of sites
\begin{equation}\label{mufIJ}
\xymatrix{\Var(k)^2/T_{\bullet}\ar[rr]^{\mu_{T_{\bullet}}}\ar[dd]_{P(f_{\bullet})}\ar[rd]^{\rho_{T_{\bullet}}} & \, & 
\Var(k)^{2,pr}/T_{\bullet}\ar[dd]^{P(f_{\bullet})}\ar[rd]^{\rho_{T_{\bullet}}} & \, \\
\, & \Var(k)^{2,sm}/T_{\bullet}\ar[rr]^{\mu_{T_{\bullet}}}\ar[dd]_{P(f_{\bullet})} & \, & 
\Var(k)^{2,smpr}/T_{\bullet}\ar[dd]^{P(f_{\bullet})} \\
\Var(k)^2/S_{\bullet}\ar[rr]^{\mu_{S_{\bullet}}}\ar[rd]^{\rho_{S_{\bullet}}} & \, & 
\Var(k)^{2,pr}/S_{\bullet}\ar[rd]^{\rho_{S_{\bullet}}} & \, \\
\, & \Var(k)^{2,sm}/S_{\bullet}\ar[rr]^{\mu_{S_{\bullet}}} & \, & \Var(k)^{2,smpr}/S_{\bullet}}.
\end{equation}
and
\begin{equation}\label{GrfIJ}
\xymatrix{\Var(k)^{2,pr}/T_{\bullet}
\ar[rr]^{\Gr_{T_{\bullet}}^{12}}\ar[dd]_{P(f_{\bullet})}\ar[rd]^{\rho_{T_{\bullet}}} & \, & 
\Var(k)/T\ar[dd]^{P(f_{\bullet})}\ar[rd]^{\rho_{T_{\bullet}}} & \, \\
\, & \Var(k)^{2,smpr}/T_{\bullet}\ar[rr]^{\Gr_{T}^{12}}\ar[dd]_{P(f_{\bullet})} & \, & 
\Var(k)^{sm}/T_{\bullet}\ar[dd]^{P(f_{\bullet})} \\
\Var(k)^{2,pr}/S_{\bullet}\ar[rr]^{\Gr_{S_{\bullet}}^{12}}\ar[rd]^{\rho_{S_{\bullet}}} & \, & 
\Var(k)/S_{\bullet}\ar[rd]^{\rho_{S_{\bullet}}} & \, \\
\, & \Var(k)^{2,sm}/S_{\bullet}\ar[rr]^{\Gr_{S_{\bullet}}^{12}} & \, & 
\Var(k)^{sm}/S_{\bullet}}.
\end{equation}
Let $s:\mathcal I\to\mathcal J$ a functor, with $\mathcal I,\mathcal J\in\Cat$, and
$f_{\bullet}:T_{\bullet}\to S_{s(\bullet)}$ a morphism with 
$T_{\bullet}\in\Fun(\mathcal J,\Var(k))$ and $S_{\bullet}\in\Fun(\mathcal I,\Var(k))$.
\begin{itemize}
\item As usual, we denote by
\begin{eqnarray*}
(f_{\bullet}^*,f_{\bullet*}):=(P(f_{\bullet})^*,P(f_{\bullet})_*): 
C(\Var(k)^{2,(sm)}/S_{\bullet})\to C(\Var(k)^{2,(sm)}/T_{\bullet})
\end{eqnarray*}
the adjonction induced by 
$P(f_{\bullet}):\Var(k)^{2,(sm)}/T_{\bullet}\to \Var(k)^{2,(sm)}/S_{\bullet}$.
Since the colimits involved in the definition of $f_{\bullet}^*=P(f_{\bullet})^*$ are filtered, 
$f_{\bullet}^*$ also preserve monomorphism. Hence, we get an adjonction
\begin{eqnarray*}
(f_{\bullet}^*,f_{\bullet*}):
C_{fil}(\Var(k)^{2,(sm)}/S_{\bullet})\leftrightarrows C_{fil}(\Var(k)^{2,(sm)}/T_{\bullet}), \\
f_{\bullet}^*((G_I,F),u_{IJ}):=((f_I^*G_I,f_I^*F),T(f_I,r_{IJ})(-)\circ f_I^*u_{IJ})
\end{eqnarray*}
\item As usual, we denote by
\begin{eqnarray*}
(f_{\bullet}^*,f_{\bullet*}):=(P(f_{\bullet})^*,P(f_{\bullet})_*):
C(\Var(k)^{2,(sm)pr}/S_{\bullet})\to C(\Var(k)^{2,(sm)pr}/T_{\bullet})
\end{eqnarray*}
the adjonction induced by 
$P(f_{\bullet}):\Var(k)^{2,(sm)pr}/T_{\bullet}\to \Var(k)^{2,(sm)pr}/S_{\bullet}$.
Since the colimits involved in the definition of $f_{\bullet}^*=P(f_{\bullet})^*$ are filtered, 
$f_{\bullet}^*$ also preserve monomorphism. Hence, we get an adjonction
\begin{eqnarray*}
(f_{\bullet}^*,f_{\bullet*}):
C_{fil}(\Var(k)^{2,(sm)pr}/S_{\bullet})\leftrightarrows C_{fil}(\Var(k)^{2,(sm)pr}/T_{\bullet}), \\
f_{\bullet}^*((G_I,F),u_{IJ}):=((f_I^*G_I,f_I^*F),T(f_I,r_{IJ})(-)\circ f_I^*u_{IJ})
\end{eqnarray*}
\end{itemize}

Let $S\in\Var(k)$. Let $S=\cup_{i=1}^l S_i$ an open affine cover and denote by $S_I=\cap_{i\in I} S_i$.
Let $i_i:S_i\hookrightarrow\tilde S_i$ closed embeddings, with $\tilde S_i\in\Var(k)$. 
For $I\subset\left[1,\cdots l\right]$, denote by $\tilde S_I=\Pi_{i\in I}\tilde S_i$.
We then have closed embeddings $i_I:S_I\hookrightarrow\tilde S_I$ and for $J\subset I$ the following commutative diagram
\begin{equation*}
D_{IJ}=\xymatrix{ S_I\ar[r]^{i_I} & \tilde S_I \\
S_J\ar[u]^{j_{IJ}}\ar[r]^{i_J} & \tilde S_J\ar[u]^{p_{IJ}}}  
\end{equation*}
where $p_{IJ}:\tilde S_J\to\tilde S_I$ is the projection
and $j_{IJ}:S_J\hookrightarrow S_I$ is the open embedding so that $j_I\circ j_{IJ}=j_J$.
This gives the diagram of algebraic varieties $(\tilde S_I)\in\Fun(\mathcal P(\mathbb N),\Var(k))$ 
which gives the diagram of sites $\Var(k)^2/(\tilde S_I)\in\Fun(\mathcal P(\mathbb N),\Cat)$.  
This gives also the diagram of algebraic varieties $(\tilde S_I)^{op}\in\Fun(\mathcal P(\mathbb N)^{op},\Var(k))$
which gives the diagram of sites $\Var(k)^2/(\tilde S_I)^{op}\in\Fun(\mathcal P(\mathbb N)^{op},\Cat)$.  
\begin{itemize}
\item Then $C_{fil}(\Var(k)^{2,(sm)}/(\tilde S_I))$ is the category  
\begin{itemize}
\item whose objects $(G,F)=((G_I,F)_{I\subset\left[1,\cdots l\right]},u_{IJ})$,
with $(G_I,F)\in C_{fil}(\Var(k)^{2,(sm)}/\tilde S_I)$,
and $u_{IJ}:(G_I,F)\to p_{IJ*}(G_J,F)$ for $I\subset J$, are morphisms
satisfying for $I\subset J\subset K$, $p_{IJ*}u_{JK}\circ u_{IJ}=u_{IK}$ in $C_{fil}(\Var(k)^{2,(sm)}/\tilde S_I)$,
\item the morphisms $m:((G,F),u_{IJ})\to((H,F),v_{IJ})$ being (see section 2.1) a family of morphisms of complexes,  
\begin{equation*}
m=(m_I:(G_I,F)\to (H_I,F))_{I\subset\left[1,\cdots l\right]}
\end{equation*}
such that $v_{IJ}\circ m_I=p_{IJ*}m_J\circ u_{IJ}$ in $C_{fil}(\Var(k)^{2,(sm)}/\tilde S_I)$.
\end{itemize}
\item Then $C_{fil}(\Var(k)^{2,(sm)pr}/(\tilde S_I))$ is the category  
\begin{itemize}
\item whose objects $(G,F)=((G_I,F)_{I\subset\left[1,\cdots l\right]},u_{IJ})$,
with $(G_I,F)\in C_{fil}(\Var(k)^{2,(sm)pr}/\tilde S_I)$,
and $u_{IJ}:(G_I,F)\to p_{IJ*}(G_J,F)$ for $I\subset J$, are morphisms
satisfying for $I\subset J\subset K$, $p_{IJ*}u_{JK}\circ u_{IJ}=u_{IK}$ in $C_{fil}(\Var(k)^{2,(sm)pr}/\tilde S_I)$,
\item the morphisms $m:((G,F),u_{IJ})\to((H,F),v_{IJ})$ being (see section 2.1) a family of morphisms of complexes,  
\begin{equation*}
m=(m_I:(G_I,F)\to (H_I,F))_{I\subset\left[1,\cdots l\right]}
\end{equation*}
such that $v_{IJ}\circ m_I=p_{IJ*}m_J\circ u_{IJ}$ in $C_{fil}(\Var(k)^{2,(sm)pr}/\tilde S_I)$.
\end{itemize}
\item Then $C_{fil}(\Var(k)^{2,(sm)}/(\tilde S_I)^{op})$ 
is the category  
\begin{itemize}
\item whose objects $(G,F)=((G_I,F)_{I\subset\left[1,\cdots l\right]},u_{IJ})$,
with $(G_I,F)\in C_{fil}(\Var(k)^{2,(sm)}/\tilde S_I)$,
and $u_{IJ}:(G_J,F)\to p_{IJ}^*(G_I,F)$ for $I\subset J$, are morphisms
satisfying for $I\subset J\subset K$, $p_{JK}^*u_{IJ}\circ u_{JK}=u_{IK}$ in $C_{fil}(\Var(k)^{2,(sm)}/\tilde S_K)$,
\item the morphisms $m:((G,F),u_{IJ})\to((H,F),v_{IJ})$ being (see section 2.1) a family of morphisms of complexes,  
\begin{equation*}
m=(m_I:(G_I,F)\to (H_I,F))_{I\subset\left[1,\cdots l\right]}
\end{equation*}
such that $v_{IJ}\circ m_J=p_{IJ}^*m_I\circ u_{IJ}$ in $C_{fil}(\Var(k)^{2,(sm)}/\tilde S_J)$.
\end{itemize}
\item Then $C_{fil}(\Var(k)^{2,(sm)pr}/(\tilde S_I)^{op})$ 
is the category  
\begin{itemize}
\item whose objects $(G,F)=((G_I,F)_{I\subset\left[1,\cdots l\right]},u_{IJ})$,
with $(G_I,F)\in C_{fil}(\Var(k)^{2,(sm)pr}/\tilde S_I)$,
and $u_{IJ}:(G_J,F)\to p_{IJ}^*(G_I,F)$ for $I\subset J$, are morphisms
satisfying for $I\subset J\subset K$, $p_{JK}^*u_{IJ}\circ u_{JK}=u_{IK}$ in $C_{fil}(\Var(k)^{2,(sm)pr}/\tilde S_K)$,
\item the morphisms $m:((G,F),u_{IJ})\to((H,F),v_{IJ})$ being (see section 2.1) a family of morphisms of complexes,  
\begin{equation*}
m=(m_I:(G_I,F)\to (H_I,F))_{I\subset\left[1,\cdots l\right]}
\end{equation*}
such that $v_{IJ}\circ m_J=p_{IJ}^*m_I\circ u_{IJ}$ in $C_{fil}(\Var(k)^{2,(sm)pr}/\tilde S_J)$.
\end{itemize}
\end{itemize}

We now define the Zariski and the etale topology on $\Var(k)^2/S$.

\begin{defi}\label{tau12}
Let $S\in\Var(k)$. 
\begin{itemize}
\item[(i)]Denote by $\tau$ a topology on $\Var(k)$, e.g. the Zariski or the etale topology. 
The $\tau$ covers in $\Var(k)^2/S$ of $(X,Z)/S$ are the families of morphisms 
\begin{eqnarray*}
\left\{(c_i:(U_i,Z\times_X U_i)/S\to(X,Z)/S)_{i\in I}, \; 
\mbox{with} \; (c_i:U_i\to X)_{i\in I} \, \tau \, \mbox{cover of} \, X \, \mbox{in} \, \Var(k)\right\}
\end{eqnarray*}
\item[(ii)]Denote by $\tau$ the Zariski or the etale topology on $\Var(k)$. 
The $\tau$ covers in $\Var(k)^{2,sm}/S$ of $(U,Z)/S$ are the families of morphisms 
\begin{eqnarray*}
\left\{(c_i:(U_i,Z\times_U U_i)/S\to(U,Z)/S)_{i\in I}, \; 
\mbox{with} \; (c_i:U_i\to U)_{i\in I} \, \tau \, \mbox{cover of} \, U \, \mbox{in} \, \Var(k)\right\}
\end{eqnarray*}
\item[(iii)]Denote by $\tau$ the Zariski or the etale topology on $\Var(k)$. 
The $\tau$ covers in $\Var(k)^{2,(sm)pr}/S$ of $(Y\times S,Z)/S$ are the families of morphisms 
\begin{eqnarray*}
\left\{(c_i\times I_S:(U_i\times S,Z\times_{Y\times S} U_i\times S)/S\to(Y\times S,Z)/S)_{i\in I}, \; 
\mbox{with} \; (c_i:U_i\to Y)_{i\in I} \, \tau \, \mbox{cover of} \, Y \, \mbox{in} \, \Var(k)\right\}
\end{eqnarray*}
\end{itemize}
\end{defi}

Will now define the $\mathbb A^1$ local property on $\Var(k)^2/S$.

Denote $\square^*:=\mathbb P_{\mathbb C}^*\backslash\left\{1\right\}$
\begin{itemize}
\item Let $S\in\Var(k)$. For $(X,Z)/S=((X,Z),h)\in\Var(k)^{2,(sm)}/S$, we consider 
\begin{equation*}
(\square^*\times X,\square^*\times Z)/S=((\square^*\times X,\square^*\times Z,h\circ p)\in\Fun(\Delta,\Var(k)^{2,(sm)}/S).
\end{equation*}
For $F\in C^-(\Var(k)^{2,(sm)}/S)$, it gives the complex
\begin{equation*}
C_*F\in C^-(\Var(k)^{2,(sm)}/S), (X,Z)/S=((X,Z),h)\mapsto C_*F((X,Z)/S):=\Tot F((\square^*\times X,\square^*\times Z/S)
\end{equation*}
together with the canonical map $c_F:=(0,I_F):F\to C_*F$.
For $F\in C(\Var(k)^{2,(sm)}/S)$, we get
\begin{equation*}
C_*F:=\holim_n C_*F^{\leq n}\in C(\Var(k)^{2,(sm)}/S),
\end{equation*}
together with the canonical map $c_F:=(0,I_F):F\to C_*F$.
For $m:F\to G$ a morphism, with $F,G\in C(\Var(k)^{2,(sm)}/S)$,
we get by functoriality the morphism $C_*m:C_*F\to C_*G$.
\item Let $S\in\Var(k)$. For $(Y\times S,Z)/S=((Y\times S,Z),h)\in\Var(k)^{2,(sm)pr}/S$, we consider 
\begin{equation*}
(\square^*\times Y\times S,\square^*\times Z)/S=(\square^*\times Y\times S,\square^*\times Z,h\circ p)
\in\Fun(\Delta,\Var(k)/S).
\end{equation*}
For $F\in C^-(\Var(k)^{2,(sm)pr}/S)$, it gives the complex
\begin{eqnarray*}
C_*F\in C^-(\Var(k)^{2,(sm)pr}/S), \\ 
(Y\times S,Z)/S=((Y\times S,Z),h)\mapsto C_*F((Y\times S,Z)/S):=\Tot F(\square^*\times Y\times S,\square^*\times Z)/S)
\end{eqnarray*}
together with the canonical map $c=c(F):=(0,I_F):F\to C_*F$.
For $F\in C(\Var(k)^{2,(sm)pr}/S)$, we get
\begin{equation*}
C_*F:=\holim_n C_*F^{\leq n}\in C(\Var(k)^{2,(sm)pr}/S),
\end{equation*}
together with the canonical map $c=c(F):=(0,I_F):F\to C_*F$.
For $m:F\to G$ a morphism, with $F,G\in C(\Var(k)^{2,(sm)pr}/S)$,
we get by functoriality the morphism $C_*m:C_*F\to C_*G$.
\item Let $S\in\Var(k)$. Let $S=\cup_{i=1}^l S_i$ an open affine cover and denote by $S_I=\cap_{i\in I} S_i$.
Let $i_i:S_i\hookrightarrow\tilde S_i$ closed embeddings, with $\tilde S_i\in\Var(k)$.
For $F=(F_I,u_{IJ})\in C(\Var(k)^{2,(sm)}/(\tilde S_I))$, it gives the complex
\begin{equation*}
C_*F=(C_*F_I,C_*u_{IJ})\in C(\Var(k)^{2,(sm)}/(\tilde S_I)), 
\end{equation*}
together with the canonical map $c_F:=(0,I_F):F\to C_*F$.
\item Let $S\in\Var(k)$. Let $S=\cup_{i=1}^l S_i$ an open affine cover and denote by $S_I=\cap_{i\in I} S_i$.
Let $i_i:S_i\hookrightarrow\tilde S_i$ closed embeddings, with $\tilde S_i\in\Var(k)$.
For $F=(F_I,u_{IJ})\in C(\Var(k)^{2,(sm)}/(\tilde S_I))$, it gives the complex
\begin{equation*}
C_*F=(C_*F_I,C_*u_{IJ})\in C(\Var(k)^{2,(sm)}/(\tilde S_I)^{op}), 
\end{equation*}
together with the canonical map $c_F:=(0,I_F):F\to C_*F$.
\item Let $S\in\Var(k)$. Let $S=\cup_{i=1}^l S_i$ an open affine cover and denote by $S_I=\cap_{i\in I} S_i$.
Let $i_i:S_i\hookrightarrow\tilde S_i$ closed embeddings, with $\tilde S_i\in\Var(k)$.
For $F=(F_I,u_{IJ})\in C(\Var(k)^{2,(sm)pr}/(\tilde S_I))$, it gives the complex
\begin{equation*}
C_*F=(C_*F_I,C_*u_{IJ})\in C(\Var(k)^{2,(sm)pr}/(\tilde S_I)), 
\end{equation*}
together with the canonical map $c_F:=(0,I_F):F\to C_*F$.
\item Let $S\in\Var(k)$. Let $S=\cup_{i=1}^l S_i$ an open affine cover and denote by $S_I=\cap_{i\in I} S_i$.
Let $i_i:S_i\hookrightarrow\tilde S_i$ closed embeddings, with $\tilde S_i\in\Var(k)$.
For $F=(F_I,u_{IJ})\in C(\Var(k)^{2,(sm)pr}/(\tilde S_I))$, it gives the complex
\begin{equation*}
C_*F=(C_*F_I,C_*u_{IJ})\in C(\Var(k)^{2,(sm)pr}/(\tilde S_I)^{op}), 
\end{equation*}
together with the canonical map $c_F:=(0,I_F):F\to C_*F$.
\end{itemize}

Let $S\in\Var(k)$. Denote for short $\Var(k)^{2,(sm)}/S$ 
either the category $\Var(k)^2/S$ or the category $\Var(k)^{2,sm}/S$. Denote by
\begin{eqnarray*}
p_a:\Var(k)^{2,(sm)}/S\to\Var(k)^{2,(sm)}/S, \\ 
(X,Z)/S=((X,Z),h)\mapsto (X\times\mathbb A^1,Z\times\mathbb A^1)/S=((X\times\mathbb A^1,Z\times\mathbb A^1,h\circ p_X), \\ 
(g:(X,Z)/S\to (X',Z')/S)\mapsto 
((g\times I_{\mathbb A^1}):(X\times\mathbb A^1,Z\times\mathbb A^1)/S\to (X'\times\mathbb A^1,Z'\times\mathbb A^1)/S)
\end{eqnarray*}
the projection functor and again by $p_a:\Var(k)^{2,(sm)}/S\to\Var(k)^{2,(sm)}/S$
the corresponding morphism of site.
Let $S\in\Var(k)$.Denote for short $\Var(k)^{2,(sm)}/S$ 
either the category $\Var(k)^2/S$ or the category $\Var(k)^{2,sm}/S$. 
Denote for short $\Var(k)^{2,(sm)pr}/S$ 
either the category $\Var(k)^{2,pr}/S$ or the category $\Var(k)^{2,smpr}/S$. Denote by
\begin{eqnarray*}
p_a:\Var(k)^{2,(sm)pr}/S\to\Var(k)^{2,(sm)pr}/S, \\ 
(Y\times S,Z)/S=((Y\times S,Z),p_S)\mapsto 
(Y\times S\times\mathbb A^1,Z\times\mathbb A^1)/S=((Y\times S\times\mathbb A^1,Z\times\mathbb A^1,p_S\circ p_{Y\times S}), \\ 
(g:(Y\times S,Z)/S\to (Y'\times S,Z')/S)\mapsto 
((g\times I_{\mathbb A^1}):(Y\times S\times\mathbb A^1,Z\times\mathbb A^1)/S\to 
(Y'\times S\times\mathbb A^1,Z'\times\mathbb A^1)/S)
\end{eqnarray*}
the projection functor and again by $p_a:\Var(k)^{2,(sm)pr}/S\to\Var(k)^{2,(sm)pr}/S$
the corresponding morphism of site.

\begin{defi}\label{a1loc12def}
\begin{itemize}
\item[(i0)]A complex $F\in C(\Var(k)^{2,(sm)}/S)$ is said to be $\mathbb A^1$ homotopic if 
$\ad(p_a^*,p_{a*})(F):F\to p_{a*}p_a^*F$ is an homotopy equivalence.
\item[(i0)']A complex $F\in C(\Var(k)^{2,(sm)pr}/S)$ is said to be $\mathbb A^1$ homotopic if 
$\ad(p_a^*,p_{a*})(F):F\to p_{a*}p_a^*F$ is an homotopy equivalence.
\item[(i)] A complex  $F\in C(\Var(k)^{2,(sm)}/S)$, is said to be $\mathbb A^1$ invariant 
if for all $(X,Z)/S\in\Var(k)^{2,(sm)}/S$ 
\begin{equation*}
F(p_X):F((X,Z)/S)\to F((X\times\mathbb A^1,(Z\times\mathbb A^1))/S) 
\end{equation*}
is a quasi-isomorphism, where $p_X:(X\times\mathbb A^1,(Z\times\mathbb A^1))\to (X,Z)$ is the projection.
Obviously, if a complex  $F\in C(\Var(k)^{2,(sm)}/S)$ is $\mathbb A^1$ homotopic, then it is $\mathbb A^1$ invariant.
\item[(i)'] A complex $G\in C(\Var(k)^{2,(sm)pr}/S)$, 
is said to be $\mathbb A^1$ invariant if for all $(Y\times S,Z)/S\in\Var(k)^{2,(sm)pr}/S$ 
\begin{equation*}
G(p_{Y\times S}):G((Y\times S,Z)/S)\to G((Y\times\mathbb A^1\times S,(Z\times\mathbb A^1))/S) 
\end{equation*}
is a quasi-isomorphism of abelian group.
Obviously, if a complex  $F\in C(\Var(k)^{2,(sm)pr}/S)$ is $\mathbb A^1$ homotopic, then it is $\mathbb A^1$ invariant.
\item[(ii)]Let $\tau$ a topology on $\Var(k)$. 
A complex $F\in C(\Var(k)^{2,(sm)}/S)$ is said to be $\mathbb A^1$ local 
for the $\tau$ topology induced on $\Var(k)^2/S$, 
if for an (hence every) $\tau$ local equivalence $k:F\to G$ with $k$ injective and $G\in C(\Var(k)^{2,(sm)}/S)$ $\tau$ fibrant,
e.g. $k:F\to E_{\tau}(F)$, $G$ is $\mathbb A^1$ invariant.
\item[(ii)']Let $\tau$ a topology on $\Var(k)$. 
A complex $F\in C(\Var(k)^{2,(sm)pr}/S)$ is said to be $\mathbb A^1$ local 
for the $\tau$ topology induced on $\Var(k)^{2,pr}/S$, 
if for an (hence every) $\tau$ local equivalence $k:F\to G$ with $k$ injective and $G\in C(\Var(k)^{2,(sm)pr}/S)$ $\tau$ fibrant,
e.g. $k:F\to E_{\tau}(F)$, $G$ is $\mathbb A^1$ invariant.
\item[(iii)] A morphism $m:F\to G$ with $F,G\in C(\Var(k)^{2,(sm)}/S)$ is said to an $(\mathbb A^1,et)$ local equivalence 
if for all $H\in C(\Var(k)^{2,(sm)}/S)$ which is $\mathbb A^1$ local for the etale topology
\begin{equation*}
\Hom(L(m),E_{et}(H)):\Hom(L(G),E_{et}(H))\to\Hom(L(F),E_{et}(H)) 
\end{equation*}
is a quasi-isomorphism.
\item[(iii)'] A morphism $m:F\to G$ with $F,G\in C(\Var(k)^{2,(sm)pr}/S)$ is said to an $(\mathbb A^1,et)$ local equivalence 
if for all $H\in C(\Var(k)^{2,(sm)pr}/S)$ which is $\mathbb A^1$ local for the etale topology
\begin{equation*}
\Hom(L(m),E_{et}(H)):\Hom(L(G),E_{et}(H))\to\Hom(L(F),E_{et}(H)) 
\end{equation*}
is a quasi-isomorphism.
\end{itemize}
\end{defi}

\begin{prop}\label{ca1Var12}
\begin{itemize}
\item[(i)]Let $S\in\Var(k)$.
Then for $F\in C(\Var(k)^{2,(sm)}/S)$, $C_*F$ is $\mathbb A^1$ local for the etale topology
and  $c(F):F\to C_*F$ is an equivalence $(\mathbb A^1,et)$ local.
\item[(i)']Let $S\in\Var(k)$.
Then for $F\in C(\Var(k)^{2,(sm)pr}/S)$, $C_*F$ is $\mathbb A^1$ local for the etale topology
and  $c(F):F\to C_*F$ is an equivalence $(\mathbb A^1,et)$ local.
\item[(ii)]A morphism $m:F\to G$ with $F,G\in C(\Var(k)^{2,(sm)}/S)$ is an $(\mathbb A^1,et)$ local equivalence
if and only if $a_{et}H^nC_*\Cone(m)=0$ for all $n\in\mathbb Z$.
\item[(ii)']A morphism $m:F\to G$ with $F,G\in C(\Var(k)^{2,(sm)pr}/S)$ is an $(\mathbb A^1,et)$ local equivalence
if and only if $a_{et}H^nC_*\Cone(m)=0$ for all $n\in\mathbb Z$.
\item[(iii)]A morphism $m:F\to G$ with $F,G\in C(\Var(k)^{2,(sm)}/S)$ is an $(\mathbb A^1,et)$ local equivalence
if and only if there exists 
\begin{eqnarray*}
\left\{(X_{1,\alpha},Z_{1,\alpha})/S,\alpha\in\Lambda_1\right\},\ldots,
\left\{(X_{r,\alpha},Z_{r,\alpha})/S,\alpha\in\Lambda_r\right\}\subset\Var(k)^{2,(sm)}/S
\end{eqnarray*}
such that we have in $\Ho_{et}(C(\Var(k)^{2,(sm)}/S))$
\begin{eqnarray*}
\Cone(m)\xrightarrow{\sim}\Cone(\oplus_{\alpha\in\Lambda_1}
\Cone(\mathbb Z((X_{1,\alpha}\times\mathbb A^1,Z_{1,\alpha}\times\mathbb A^1)/S)\to\mathbb Z((X_{1,\alpha},Z_{1,\alpha})/S)) \\
\to\cdots\to\oplus_{\alpha\in\Lambda_r}
\Cone(\mathbb Z((X_{r,\alpha}\times\mathbb A^1,Z_{r,\alpha}\times\mathbb A^1)/S)\to\mathbb Z((X_{r,\alpha},Z_{r,\alpha})/S)))
\end{eqnarray*}
\item[(iii)']A morphism $m:F\to G$ with $F,G\in C(\Var(k)^{2,(sm)pr}/S)$ is an $(\mathbb A^1,et)$ local equivalence
if and only if there exists 
\begin{eqnarray*}
\left\{(Y_{1,\alpha}\times S,Z_{1,\alpha})/S,\alpha\in\Lambda_1\right\},\ldots,
\left\{(Y_{r,\alpha}\times S,Z_{r,\alpha})/S,\alpha\in\Lambda_r\right\}\subset\Var(k)^{2,(sm)pr}/S
\end{eqnarray*}
such that we have in $\Ho_{et}(C(\Var(k)^{2,(sm)}/S))$
\begin{eqnarray*}
\Cone(m)\xrightarrow{\sim}\Cone(\oplus_{\alpha\in\Lambda_1}
\Cone(\mathbb Z((Y_{1,\alpha}\times\mathbb A^1\times S,Z_{1,\alpha}\times\mathbb A^1)/S)\to
\mathbb Z((Y_{1,\alpha}\times S,Z_{1,\alpha})/S)) \\
\to\cdots\to\oplus_{\alpha\in\Lambda_r}
\Cone(\mathbb Z((Y_{r,\alpha}\times\mathbb A^1\times S,Z_{r,\alpha}\times\mathbb A^1)/S)\to
\mathbb Z((Y_{r,\alpha}\times S,Z_{r,\alpha})/S)))
\end{eqnarray*}
\end{itemize}
\end{prop}

\begin{proof}
Standard : see Ayoub's thesis section 4 for example. 
Indeed, for (iii), by definition, if $\Cone(m)$ is of the given form, 
then it is an equivalence $(\mathbb A^1,et)$ local, on the other hand if $m$ is an equivalence $(\mathbb A^1,et)$ local,
we consider the commutative diagram
\begin{equation*}
\xymatrix{F\ar[r]^{c(F)}\ar[d]^{m} & C_*F\ar[d]^{c_*m} \\ G\ar[r]^{c(G)} & C_*G}
\end{equation*}
to deduce that $\Cone(m)$ is of the given form.
\end{proof}

\begin{defiprop}\label{projmodstr12}
Let $S\in\Var(k)$.
\begin{itemize}
\item[(i)]With the weak equivalence the $(\mathbb A^1,et)$ local equivalence and 
the fibration the epimorphism with $\mathbb A^1_S$ local and etale fibrant kernels gives
a model structure on  $C(\Var(k)^{2,(sm)}/S)$ : the left bousfield localization
of the projective model structure of $C(\Var(k)^{2,(sm)}/S)$. 
We call it the projective $(\mathbb A^1,et)$ model structure.
\item[(ii)]With the weak equivalence the $(\mathbb A^1,et)$ local equivalence and 
the fibration the epimorphism with $\mathbb A^1_S$ local and etale fibrant kernels gives
a model structure on  $C(\Var(k)^{2,(sm)pr}/S)$ : the left bousfield localization
of the projective model structure of $C(\Var(k)^{2,(sm)pr}/S)$. 
We call it the projective $(\mathbb A^1,et)$ model structure.
\end{itemize}
\end{defiprop}

\begin{proof}
Similar to the proof of proposition \ref{projmodstr}.
\end{proof}

We have, similarly to the case of single varieties the following :

\begin{prop}\label{g12}
Let $g:T\to S$ a morphism with $T,S\in\Var(k)$.
\begin{itemize}
\item[(i)] The adjonction $(g^*,g_*):C(\Var(k)^{2,(sm)}/S)\leftrightarrows C(\Var(k)^{2,(sm)}/T)$
is a Quillen adjonction for the projective $(\mathbb A^1,et)$ model structure (see definition-proposition \ref{projmodstr12})
\item[(i)'] The functor $g^*:C(\Var(k)^{2,(sm)}/S)\to C(\Var(k)^{2,(sm)}/T)$
sends quasi-isomorphism to quasi-isomorphism,
sends equivalence Zariski local to equivalence Zariski local, and equivalence etale local to equivalence etale local,
sends $(\mathbb A^1,et)$ local equivalence to $(\mathbb A^1,et)$ local equivalence.
\item[(ii)] The adjonction $(g^*,g_*):C(\Var(k)^{2,(sm)pr}/S)\leftrightarrows C(\Var(k)^{2,(sm)pr}/T)$
is a Quillen adjonction for the projective $(\mathbb A^1,et)$ model structure (see definition-proposition \ref{projmodstr12})
\item[(ii)'] The functor $g^*:C(\Var(k)^{2,(sm)pr}/S)\to C(\Var(k)^{2,(sm)pr}/T)$
sends quasi-isomorphism to quasi-isomorphism,
sends equivalence Zariski local to equivalence Zariski local, and equivalence etale local to equivalence etale local,
sends $(\mathbb A^1,et)$ local equivalence to $(\mathbb A^1,et)$ local equivalence.
\end{itemize}
\end{prop}

\begin{proof}
\noindent(i):Follows immediately from definition.
\noindent(i)': Since the functor $g^*$ preserve epimorphism and also monomorphism (the colimits involved being filetered),
$g^*$ sends quasi-isomorphism to quasi-isomorphism. Hence it preserve Zariski and etale local equivalence.
The fact that it preserve $(\mathbb A^1,et)$ local equivalence then follows similarly to the single case by the fact
that $g_*$ preserve by definition $\mathbb A^1$ equivariant presheaves.

\noindent(ii) and (ii)': Similar to (i) and (i)'.
\end{proof}

\begin{prop}\label{rho12}
Let $S\in\Var(k)$. 
\begin{itemize}
\item[(i)] The adjonction $(\rho_S^*,\rho_{S*}):C(\Var(k)^{2,sm}/S)\leftrightarrows C(\Var(k)^2/S)$
is a Quillen adjonction for the $(\mathbb A^1,et)$ projective model structure.
\item[(i)']The functor $\rho_{S*}:C(\Var(k)^2/S)\to C(\Var(k)^{2,sm}/S)$
sends quasi-isomorphism to quasi-isomorphism,
sends equivalence Zariski local to equivalence Zariski local, and equivalence etale local to equivalence etale local,
sends $(\mathbb A^1,et)$ local equivalence to $(\mathbb A^1,et)$ local equivalence.
\item[(ii)] The adjonction $(\rho_S^*,\rho_{S*}):C(\Var(k)^{2,smpr}/S)\leftrightarrows C(\Var(k)^{2,pr}/S)$
is a Quillen adjonction for the $(\mathbb A^1,et)$ projective model structure.
\item[(ii)']The functor $\rho_{S*}:C(\Var(k)^{2,pr}/S)\to C(\Var(k)^{2,smpr}/S)$
sends quasi-isomorphism to quasi-isomorphism,
sends equivalence Zariski local to equivalence Zariski local, and equivalence etale local to equivalence etale local,
sends $(\mathbb A^1,et)$ local equivalence to $(\mathbb A^1,et)$ local equivalence.
\end{itemize}
\end{prop}

\begin{proof}
Similar to the proof of proposition \ref{rho1}.
\end{proof}

\begin{prop}\label{mu12}
Let $S\in\Var(k)$. 
\begin{itemize}
\item[(i)] The adjonction $(\mu_S^*,\mu_{S*}):C(\Var(k)^{2,pr}/S)\leftrightarrows C(\Var(k)^2/S)$
is a Quillen adjonction for the $(\mathbb A^1,et)$ projective model structure.
\item[(i)']The functor $\mu_{S*}:C(\Var(k)^2/S)\to C(\Var(k)^{2,pr}/S)$
sends quasi-isomorphism to quasi-isomorphism,
sends equivalence Zariski local to equivalence Zariski local, and equivalence etale local to equivalence etale local,
sends $(\mathbb A^1,et)$ local equivalence to $(\mathbb A^1,et)$ local equivalence.
\item[(ii)] The adjonction $(\mu_S^*,\mu_{S*}):C(\Var(k)^{2,smpr}/S)\leftrightarrows C(\Var(k)^{2,pr}/S)$
is a Quillen adjonction for the $(\mathbb A^1,et)$ projective model structure.
\item[(ii)']The functor $\mu_{S*}:C(\Var(k)^{2,sm}/S)\to C(\Var(k)^{2,smpr}/S)$
sends quasi-isomorphism to quasi-isomorphism,
sends equivalence Zariski local to equivalence Zariski local, and equivalence etale local to equivalence etale local,
sends $(\mathbb A^1,et)$ local equivalence to $(\mathbb A^1,et)$ local equivalence.
\end{itemize}
\end{prop}

\begin{proof}
Similar to the proof of proposition \ref{rho1}. 
Indeed, for (i)' or (ii)', if $m:F\to G$ with $F,G\in C(\Var(k)^{2,(sm)})$ is an equivalence $(\mathbb A^1,et)$
local then (see proposition \ref{ca1Var12}), there exists 
\begin{equation*}
\left\{(X_{1,\alpha},Z_{1,\alpha})/S,\alpha\in\Lambda_1\right\},\ldots,
\left\{(X_{r,\alpha},Z_{r,\alpha})/S,\alpha\in\Lambda_r\right\}
\subset\Var(k)^{2,(sm)}/S 
\end{equation*}
such that we have in $\Ho_{et}(C(\Var(k)^{2,(sm)}/S))$
\begin{eqnarray*}
\Cone(m)\xrightarrow{\sim}\Cone(\oplus_{\alpha\in\Lambda_1}
\Cone(\mathbb Z((X_{1,\alpha}\times\mathbb A^1,Z_{1,\alpha}\times\mathbb A^1)/S)
\to\mathbb Z((X_{1,\alpha},Z_{1,\alpha})/S)) \\
\to\cdots\to\oplus_{\alpha\in\Lambda_r}
\Cone(\mathbb Z((X_{r,\alpha}\times\mathbb A^1,Z_{r,\alpha}\times\mathbb A^1)/S)\to
\mathbb Z((X_{r,\alpha},Z_{r,\alpha})/S))) \\
\xrightarrow{\sim}\Cone(
\Cone(\oplus_{\alpha\in\Lambda_1}\mathbb Z((X_{1,\alpha},Z_{1,\alpha})/S)\otimes\mathbb Z(\mathbb A^1,\mathbb A^1)/S\to
\oplus_{\alpha\in\Lambda_1}\mathbb Z((X_{1,\alpha},Z_{1,\alpha})/S)) \\
\to\cdots\to
\Cone(\oplus_{\alpha\in\Lambda_r}\mathbb Z((X_{r,\alpha},Z_{r,\alpha})/S)\otimes\mathbb Z((\mathbb A^1,\mathbb A^1)/S)\to
\oplus_{\alpha\in\Lambda_r}\mathbb Z((X_{r,\alpha},Z_{r,\alpha})/S)),
\end{eqnarray*}
this gives in $\Ho_{et}(C(\Var(k)^{2,(sm)pr}/S))$
\begin{eqnarray*}
\Cone(\mu_{S*}m)\xrightarrow{\sim}\Cone( \\
\Cone((L\mu_{S*}\oplus_{\alpha\in\Lambda_1}\mathbb Z((X_{1,\alpha},Z_{1,\alpha})/S))\otimes\mathbb Z((\mathbb A^1,\mathbb A^1)/S)
\to(L\mu_{S*}\oplus_{\alpha\in\Lambda_1}\mathbb Z((X_{1,\alpha},Z_{1,\alpha})/S)) \\
\to\cdots\to
\Cone((L\mu_{S*}\oplus_{\alpha\in\Lambda_r}\mathbb Z((X_{r,\alpha},Z_{r,\alpha})/S))\otimes\mathbb Z((\mathbb A^1,\mathbb A^1)/S)
\to(L\mu_{S*}\oplus_{\alpha\in\Lambda_r}\mathbb Z((X_{1,\alpha},Z_{1,\alpha})/S))))
\end{eqnarray*}
hence $\mu_{S*}m:\mu_{S*}F\to\mu_{S*}G$ is an equivalence $(\mathbb A^1,et)$ local.
\end{proof}

We also have

\begin{prop}\label{Gra1}
Let $S\in\Var(k)$. 
\begin{itemize}
\item[(i)] The adjonction $(\Gr_S^{12*},\Gr_{S*}^{12}):C(\Var(k)/S)\leftrightarrows C(\Var(k)^{2,pr}/S)$
is a Quillen adjonction for the $(\mathbb A^1,et)$ projective model structure.
\item[(ii)] The adjonction $(\Gr_S^{12*}\Gr_{S*}^{12}:C(\Var(k)^{sm}/S)\leftrightarrows C(\Var(k)^{2,smpr}/S)$
is a Quillen adjonction for the $(\mathbb A^1,et)$ projective model structure.
\end{itemize}
\end{prop}

\begin{proof}
Immediate from definition.
\end{proof}

\begin{itemize}
\item For $f:X\to S$ a morphism with $X,S\in\Var(k)$ and $Z\subset X$ a closed subset, we denote
$\mathbb Z^{tr}((X,Z)/S)\in\PSh(\Var(k)^2/S)$ the presheaf given by
\begin{itemize}
\item for $(X',Z')/S\in\Var(k)^2/S$, with $X'$ irreducible, 
\begin{equation*}
\mathbb Z^{tr}((X,Z)/S)((X',Z')/S):=\left\{\alpha\in\mathcal Z^{fs/X}(X'\times_S X),s.t. p_X(p_{X'}^{-1}(Z'))\subset Z\right\}
\subset\mathcal Z_{d_{X'}}(X'\times_S X) 
\end{equation*}
\item for $g:(X_2,Z_2)/S\to (X_1,Z_1)/S$ a morphism, with $(X_1,Z_1)/S,(X_2,Z_2)/S\in\Var(k)^2/S$,
\begin{equation*}
\mathbb Z^{tr}((X,Z)/S)(g):\mathbb Z^{tr}((X,Z)/S)((X_1,Z_1)/S)\to\mathbb Z^{tr}((X,Z)/S)((X_2,Z_2)/S), \;
\alpha\mapsto (g\times I)^{-1}(\alpha)
\end{equation*}
with $g\times I:X_2\times_S X\to X_1\times_S X$.
\end{itemize}
\item For $f:X\to S$ a morphism with $X,S\in\Var(k)$, $Z\subset X$ a closed subset and $r\in\mathbb N$, we denote
$\mathbb Z^{equir}((X,Z)/S)\in\PSh(\Var(k)^2/S)$ the presheaf given by
\begin{itemize}
\item for $(X',Z')/S\in\Var(k)^2/S$, with $X'$ irreducible,
\begin{equation*}
\mathbb Z^{equir}((X,Z)/S)((X',Z')/S):=\left\{\alpha\in\mathcal Z^{equir/X}(X'\times_S X),s.t. p_X(p_{X'}^{-1}(Z'))\right\}
\subset\mathcal Z_{d_{X'}}(X'\times_S X) 
\end{equation*}
\item for $g:(X_2,Z_2)/S\to (X_1,Z_1)/S$ a morphism, with $(X_1,Z_1)/S,(X_2,Z_2)/S\in\Var(k)^2/S$,
\begin{equation*}
\mathbb Z^{equir}((X,Z)/S)(g):\mathbb Z^{equir}((X,Z)/S)((X_1,Z_1)/S)\to\mathbb Z^{equir}((X,Z)/S)((X_2,Z_2)/S), \;
\alpha\mapsto (g\times I)^{-1}(\alpha)
\end{equation*}
with $g\times I:X_2\times_S X\to X_1\times_S X$.
\end{itemize}
\item Let $S\in\Var(k)$. We denote by 
$\mathbb Z_S(d):=\mathbb Z^{equi0}((S\times\mathbb A^d,S\times\mathbb A^d)/S)[-2d]$
the Tate twist. For $F\in C(\Var(k)^2/S)$, we denote by $F(d):=F\otimes\mathbb Z_S(d)$.
\end{itemize}

For $S\in\Var(k)$, let $\Cor(\Var(k)^{2,(sm)}/S)$ be the category 
\begin{itemize}
\item whose objects are those of $\Var(k)^{2,(sm)}/S$, i.e. 
$(X,Z)/S=((X,Z),h)$, $h:X\to S$ with $X\in\Var(k)$, $Z\subset X$ a closed subset, 
\item whose morphisms $\alpha:(X',Z)/S=((X',Z),h_1)\to (X,Z)/S=((X,Z),h_2)$ 
is finite correspondence that is $\alpha\in\oplus_i\mathbb Z^{tr}((X_i,Z)/S)((X',Z')/S)$, 
where $X'=\sqcup_i X'_i$, with $X'_i$ connected, 
the composition being defined in the same way as the morphism $\Cor(\Var(k)^{(sm)}/S)$. 
\end{itemize}
We denote by 
$\Tr(S):\Cor(\Var(k)^{2,(sm)}/S)\to\Var(k)^{2,(sm)}/S$ 
the morphism of site
given by the inclusion functor
$\Tr(S):\Var(k)^{2,(sm)}/S\hookrightarrow\Cor(\Var(k)^{2,(sm)}/S)$
It induces an adjonction
\begin{equation*}
(\Tr(S)^*\Tr(S)_*):C(\Var(k)^{2,(sm)}/S)\leftrightarrows C(\Cor(\Var(k)^{2,(sm)}/S))
\end{equation*}
A complex of preheaves $G\in C(\Var(k)^{2,(sm)}/S)$ is said to admit transferts
if it is in the image of the embedding
\begin{equation*}
\Tr(S)_*:C(\Cor(\Var(k)^{2,(sm)}/S)\hookrightarrow C(\Var(k)^{2,(sm)}/S),
\end{equation*}
that is $G=\Tr(S)_*\Tr(S)^*G$.
We then have the full subcategory $\Cor(\Var(k)^{2,(sm)pr}/S)\subset\Cor(\Var(k)^{2,(sm)}/S)$
consisting of the objects of $\Var(k)^{2,(sm)pr}/S)$.
We have the adjonction
\begin{equation*}
(\Tr(S)^*\Tr(S)_*):C(\Var(k)^{2,(sm)pr}/S)\leftrightarrows C(\Cor(\Var(k)^{2,(sm)pr}/S))
\end{equation*}
A complex of preheaves $G\in C(\Var(k)^{2,(sm)pr}/S)$ is said to admit transferts
if it is in the image of the embedding
\begin{equation*}
\Tr(S)_*:C(\Cor(\Var(k)^{2,(sm)pr}/S)\hookrightarrow C(\Var(k)^{2,(sm)pr}/S),
\end{equation*}
that is $G=\Tr(S)_*\Tr(S)^*G$.

Let $S\in\Var(k)$. Let $S=\cup_{i=1}^l S_i$ an open affine cover and denote by $S_I=\cap_{i\in I} S_i$.
Let $i_i:S_i\hookrightarrow\tilde S_i$ closed embeddings, with $\tilde S_i\in\Var(k)$.
\begin{itemize}
\item For $(G_I,K_{IJ})\in C(\Var(k)^{2,(sm)}/(\tilde S_I)^{op})$ and 
$(H_I,T_{IJ})\in C(\Var(k)^{2,(sm)}/(\tilde S_I))$, we denote
\begin{eqnarray*}
\mathcal Hom((G_I,K_{IJ}),(H_I,T_{IJ})):=(\mathcal Hom(G_I,H_I),u_{IJ}((G_I,K_{IJ}),(H_I,T_{IJ})))
\in C(\Var(k)^{2,(sm)}/(\tilde S_I))
\end{eqnarray*}
with
\begin{eqnarray*}
u_{IJ}((G_I,K_{IJ})(H_I,T_{IJ})):\mathcal Hom(G_I,H_I) \\
\xrightarrow{\ad(p_{IJ}^*,p_{IJ*})(-)}p_{IJ*}p_{IJ}^*\mathcal Hom(G_I,H_I)
\xrightarrow{T(p_{IJ},hom)(-,-)}p_{IJ*}\mathcal Hom(p_{IJ}^*G_I,p_{IJ}^*H_I) \\
\xrightarrow{\mathcal Hom(p_{IJ}^*G_I,T_{IJ})}p_{IJ*}\mathcal Hom(p_{IJ}^*G_I,H_J)
\xrightarrow{\mathcal Hom(K_{IJ},H_J)}p_{IJ*}\mathcal Hom(G_J,H_J).
\end{eqnarray*}
This gives in particular the functor
\begin{eqnarray*}
\mathbb D_{(\tilde S_I)}^{12}:C(\Var(k)^{2,(sm)}/(\tilde S_I)^{op})\to C(\Var(k)^{2,(sm)}/(\tilde S_I)), \\
(H_I,T_{IJ})\mapsto\mathbb D_{(\tilde S_I)}^{12}L(H_I,T_{IJ}):=
\mathcal Hom((LH_I,T^q_{IJ}),(E_{et}\mathbb Z_{\tilde S_I},I_{IJ}))=(\mathbb D_{\tilde S_I}^{12}LH_I,T_{IJ}^d)
\end{eqnarray*}
\item For $(G_I,K_{IJ})\in C(\Var(k)^{2,(sm)pr}/(\tilde S_I)^{op})$ and 
$(H_I,T_{IJ})\in C(\Var(k)^{2,(sm)pr}/(\tilde S_I))$, we denote
\begin{eqnarray*}
\mathcal Hom((G_I,K_{IJ}),(H_I,T_{IJ})):=(\mathcal Hom(G_I,H_I),u_{IJ}((G_I,K_{IJ}),(H_I,T_{IJ})))
\in C(\Var(k)^{2,(sm)pr}/(\tilde S_I))
\end{eqnarray*}
with
\begin{eqnarray*}
u_{IJ}((G_I,K_{IJ})(H_I,T_{IJ})):\mathcal Hom(G_I,H_I) \\
\xrightarrow{\ad(p_{IJ}^*,p_{IJ*})(-)}p_{IJ*}p_{IJ}^*\mathcal Hom(G_I,H_I)
\xrightarrow{T(p_{IJ},hom)(-,-)}p_{IJ*}\mathcal Hom(p_{IJ}^*G_I,p_{IJ}^*H_I) \\
\xrightarrow{\mathcal Hom(p_{IJ}^*G_I,T_{IJ})}p_{IJ*}\mathcal Hom(p_{IJ}^*G_I,H_J)
\xrightarrow{\mathcal Hom(K_{IJ},H_J)}p_{IJ*}\mathcal Hom(G_J,H_J).
\end{eqnarray*}
This gives in particular the functor
\begin{eqnarray*}
\mathbb D_{(\tilde S_I)}^{12}:C(\Var(k)^{2,(sm)pr}/(\tilde S_I)^{op})\to C(\Var(k)^{2,(sm)pr}/(\tilde S_I)), 
(H_I,T_{IJ})\mapsto\mathbb D_{(\tilde S_I)}^{12}L(H_I,T_{IJ})
\end{eqnarray*}
\end{itemize}

The functors $p_a$ naturally extend to functors
\begin{eqnarray*}
p_a:\Var(k)^{2,(sm)}/(\tilde S_I)\to\Var(k)^{2,(sm)}/(\tilde S_I), \\ 
((X,Z)/\tilde S_I,u_{IJ})\mapsto ((X\times\mathbb A^1,Z\times\mathbb A^1)/\tilde S_I,u_{IJ}\times I), \\
(g:((X,Z)/\tilde S_I,u_{IJ})\to ((X',Z')/\tilde S_I,u_{IJ}))\mapsto \\
((g\times I_{\mathbb A^1}):((X\times\mathbb A^1,Z\times\mathbb A^1)/\tilde S_I,u_{IJ}\times I)\to 
((X'\times\mathbb A^1,Z'\times\mathbb A^1)/\tilde S_I,u_{IJ}\times I))
\end{eqnarray*}
the projection functor and again by $p_a:\Var(k)^{2,(sm)}/(\tilde S_I)\to\Var(k)^{2,(sm)}/(\tilde S_I)$
the corresponding morphism of site, and
\begin{eqnarray*}
p_a:\Var(k)^{2,(sm)pr}/(\tilde S_I)\to\Var(k)^{2,(sm)pr}/(\tilde S_I), \\ 
((Y\times\tilde S_I,Z)/\tilde S_I,u_{IJ})\mapsto 
((Y\times\tilde S_I\times\mathbb A^1,Z\times\mathbb A^1)/\tilde S_I,u_{IJ}\times I), \\ 
(g:((Y\times\tilde S_I,Z)/\tilde S_I,u_{IJ})\to ((Y'\times\tilde S_I,Z')/\tilde S_I,u_{IJ}))\mapsto \\
((g\times I_{\mathbb A^1}):((Y\times\tilde S_I\times\mathbb A^1,Z\times\mathbb A^1)/\tilde S_I,u_{IJ}\times I),  
((Y'\times\tilde S_I\times\mathbb A^1,Z'\times\mathbb A^1)/\tilde S_I,u_{IJ}\times I)), 
\end{eqnarray*}
the projection functor and again by $p_a:\Var(k)^{2,(sm)pr}/(\tilde S_I)\to\Var(k)^{2,(sm)pr}/(\tilde S_I)$
the corresponding morphism of site.
These functors also gives the morphisms of sites
$p_a:\Var(k)^{2,(sm)}/(\tilde S_I)^{op}\to\Var(k)^{2,(sm)}/(\tilde S_I)^{op}$ and
$p_a:\Var(k)^{2,(sm)pr}/(\tilde S_I)^{op}\to\Var(k)^{2,(sm)pr}/(\tilde S_I)^{op}$.

\begin{defi}\label{a1loc12defIJ}
Let $S\in\Var(k)$. Let $S=\cup_{i=1}^l S_i$ an open affine cover and denote by $S_I=\cap_{i\in I} S_i$.
Let $i_i:S_i\hookrightarrow\tilde S_i$ closed embeddings, with $\tilde S_i\in\Var(k)$.
\begin{itemize}
\item[(i0)]A complex $(F_I,u_{IJ})\in C(\Var(k)^{2,(sm)}/(\tilde S_I))$ is said to be $\mathbb A^1$ homotopic if 
$\ad(p_a^*,p_{a*})((F_I,u_{IJ})):(F_I,u_{IJ})\to p_{a*}p_a^*(F_I,u_{IJ})$ is an homotopy equivalence.
\item[(i0)']A complex $(F_I,u_{IJ})\in C(\Var(k)^{2,(sm)pr}/(\tilde S_I))$ is said to be $\mathbb A^1$ homotopic if 
$\ad(p_a^*,p_{a*})((F_I,u_{IJ})):(F_I,u_{IJ})\to p_{a*}p_a^*(F_I,u_{IJ})$ is an homotopy equivalence.
\item[(i)] A complex  $(F_I,u_{IJ})\in C(\Var(k)^{2,(sm)}/(\tilde S_I))$ is said to be $\mathbb A^1$ invariant 
if for all $((X_I,Z_I)/\tilde S_I,s_{IJ})\in\Var(k)^{2,(sm)}/(\tilde S_I)$ 
\begin{equation*}
(F_I(p_{X_I})):(F_I((X_I,Z_I)/\tilde S_I),F_J(s_{IJ})\circ u_{IJ}(-)\to 
(F_I((X_I\times\mathbb A^1,(Z_I\times\mathbb A^1))/\tilde S_I),F_J(s_{IJ}\times I)\circ u_{IJ}(-)) 
\end{equation*}
is a quasi-isomorphism, where $p_{X_I}:(X_I\times\mathbb A^1,(Z_I\times\mathbb A^1))\to (X_I,Z_I)$ are the projection,
and $s_{IJ}:(X_I\times\tilde S_{J\backslash I},Z_I)/\tilde S_J\to(X_J,Z_J)/\tilde S_J$.
Obviously a complex $(F_I,u_{IJ})\in C(\Var(k)^{2,(sm)}/(\tilde S_I))$ is $\mathbb A^1$ invariant
if and only if all the $F_I$ are $\mathbb A^1$ invariant.
\item[(i)'] A complex  $(G_I,u_{IJ})\in C(\Var(k)^{2,(sm)pr}/(\tilde S_I))$ is said to be $\mathbb A^1$ invariant 
if for all $((Y\times\tilde S_I,Z_I)/\tilde S_I,s_{IJ})\in\Var(k)^{2,(sm)pr}/(\tilde S_I)$ 
\begin{eqnarray*}
(G_I(p_{Y\times\tilde S_I})):
(G_I((Y\times\tilde S_I,Z_I)/\tilde S_I),G_J(s_{IJ})\circ u_{IJ}(-))\to \\
(G_I((Y\times\tilde S_I\times\mathbb A^1,(Z_I\times\mathbb A^1))/\tilde S_I),G_J(s_{IJ}\times I)\circ u_{IJ}(-)) 
\end{eqnarray*}
is a quasi-isomorphism. 
Obviously a complex  $(G_I,u_{IJ})\in C(\Var(k)^{2,(sm)pr}/(\tilde S_I))$ is $\mathbb A^1$ invariant
if and only if all the $G_I$ are $\mathbb A^1$ invariant.
\item[(ii)]Let $\tau$ a topology on $\Var(k)$. 
A complex $F=(F_I,u_{IJ})\in C(\Var(k)^{2,(sm)}/(\tilde S_I))$ is said to be $\mathbb A^1$ local 
for the $\tau$ topology induced on $\Var(k)^2/(\tilde S_I)$, 
if for an (hence every) $\tau$ local equivalence $k:F\to G$ with $k$ injective and 
$G=(G_I,v_{IJ})\in C(\Var(k)^{2,(sm)}/(\tilde S_I))$ $\tau$ fibrant,
e.g. $k:(F_I,u_{IJ})\to (E_{\tau}(F_I),E(u_{IJ}))$, $G$ is $\mathbb A^1$ invariant.
\item[(ii)']Let $\tau$ a topology on $\Var(k)$. 
A complex $F=(F_I,u_{IJ})\in C(\Var(k)^{2,(sm)pr}/(\tilde S_I))$ is said to be $\mathbb A^1$ local 
for the $\tau$ topology induced on $\Var(k)^2/(\tilde S_I)$, 
if for an (hence every) $\tau$ local equivalence $k:F\to G$ with $k$ injective and 
$G=(G_I,u_{IJ})\in C(\Var(k)^{2,(sm)pr}/(\tilde S_I))$ $\tau$ fibrant,
e.g. $k:(F_I,u_{IJ})\to (E_{\tau}(F_I),E(u_{IJ}))$, $G$ is $\mathbb A^1$ invariant.
\item[(iii)] A morphism $m=(m_I):(F_I,u_{IJ})\to (G_I,v_{IJ})$ with 
$(F_I,u_{IJ}),(G_I,v_{IJ})\in C(\Var(k)^{2,(sm)}/(\tilde S_I))$ is said to be an $(\mathbb A^1,et)$ local equivalence 
if for all $H=(H_I,w_{IJ})\in C(\Var(k)^{2,(sm)}/(\tilde S_I))$ which is $\mathbb A^1$ local for the etale topology 
\begin{eqnarray*}
(\Hom(L(m_I),E_{et}(H_I))):\Hom(L(G_I,v_{IJ}),E_{et}(H_I,w_{IJ}))\to\Hom(L(F_I,u_{IJ}),E_{et}(H_I,w_{IJ})) 
\end{eqnarray*}
is a quasi-isomorphism (of complexes of abelian groups).
Obviously, if a morphism $m=(m_I):(F_I,u_{IJ})\to (G_I,v_{IJ})$ with 
$(F_I,u_{IJ}),(G_I,u_{IJ})\in C(\Var(k)^{2,(sm)}/(\tilde S_I))$ is an $(\mathbb A^1,et)$ local equivalence, 
then all the $m_I:F_I\to G_I$ are $(\mathbb A^1,et)$ local equivalence.
\item[(iii)'] A morphism $m=(m_I):(F_I,u_{IJ})\to (G_I,v_{IJ})$ with 
$(F_I,u_{IJ}),(G_I,v_{IJ})\in C(\Var(k)^{2,(sm)pr}/(\tilde S_I))$ is said to be an $(\mathbb A^1,et)$ local equivalence 
if for all $(H_I,w_{IJ})\in C(\Var(k)^{2,(sm)pr}/(\tilde S_I))$ which is $\mathbb A^1$ local for the etale topology
\begin{eqnarray*}
(\Hom(L(m_I),E_{et}(H_I))):\Hom(L(G_I,v_{IJ}),E_{et}(H_I,w_{IJ}))\to\Hom(L(F_I,u_{IJ}),E_{et}(H_I,w_{IJ})) 
\end{eqnarray*}
is a quasi-isomorphism (of complexes of abelian groups).
Obviously, if a morphism $m=(m_I):(F_I,u_{IJ})\to (G_I,v_{IJ})$ with 
$(F_I,u_{IJ}),(G_I,u_{IJ})\in C(\Var(k)^{2,(sm)pr}/(\tilde S_I))$ is an $(\mathbb A^1,et)$ local equivalence, 
then all the $m_I:F_I\to G_I$ are $(\mathbb A^1,et)$ local equivalence.
\item[(iv)] A morphism $m=(m_I):(F_I,u_{IJ})\to (G_I,v_{IJ})$ with 
$(F_I,u_{IJ}),(G_I,v_{IJ})\in C(\Var(k)^{2,(sm)}/(\tilde S_I)^{op})$ is said to be an $(\mathbb A^1,et)$ local equivalence 
if for all $H=(H_I,w_{IJ})\in C(\Var(k)^{2,(sm)}/(\tilde S_I)^{op})$ which is $\mathbb A^1$ local for the etale topology 
\begin{eqnarray*}
(\Hom(L(m_I),E_{et}(H_I))):\Hom(L(G_I,v_{IJ}),E_{et}(H_I,w_{IJ}))\to\Hom(L(F_I,u_{IJ}),E_{et}(H_I,w_{IJ})) 
\end{eqnarray*}
is a quasi-isomorphism (of complexes of abelian groups).
Obviously, if a morphism $m=(m_I):(F_I,u_{IJ})\to (G_I,v_{IJ})$ with 
$(F_I,u_{IJ}),(G_I,u_{IJ})\in C(\Var(k)^{2,(sm)}/(\tilde S_I)^{op})$ is an $(\mathbb A^1,et)$ local equivalence, 
then all the $m_I:F_I\to G_I$ are $(\mathbb A^1,et)$ local equivalence and
for all $H=(H_I,w_{IJ})\in C(\Var(k)^{2,(sm)}/(\tilde S_I))$ which is $\mathbb A^1$ local for the etale topology 
\begin{eqnarray*}
(\Hom(L(m_I),E_{et}(H_I))):\Hom(L(G_I,v_{IJ}),E_{et}(H_I,w_{IJ}))\to\Hom(L(F_I,u_{IJ}),E_{et}(H_I,w_{IJ})) 
\end{eqnarray*}
is a quasi-isomorphism (of diagrams of complexes of abelian groups)
\item[(iv)'] A morphism $m=(m_I):(F_I,u_{IJ})\to (G_I,v_{IJ})$ with 
$(F_I,u_{IJ}),(G_I,v_{IJ})\in C(\Var(k)^{2,(sm)pr}/(\tilde S_I)^{op})$ is said to be an $(\mathbb A^1,et)$ local equivalence 
if for all $(H_I,w_{IJ})\in C(\Var(k)^{2,(sm)pr}/(\tilde S_I)^{op})$ which is $\mathbb A^1$ local for the etale topology
\begin{eqnarray*}
(\Hom(L(m_I),E_{et}(H_I))):\Hom(L(G_I,v_{IJ}),E_{et}(H_I,w_{IJ}))\to\Hom(L(F_I,u_{IJ}),E_{et}(H_I,w_{IJ})) 
\end{eqnarray*}
is a quasi-isomorphism (of complexes of abelian groups).
Obviously, if a morphism $m=(m_I):(F_I,u_{IJ})\to (G_I,v_{IJ})$ with 
$(F_I,u_{IJ}),(G_I,u_{IJ})\in C(\Var(k)^{2,(sm)pr}/(\tilde S_I)^{op})$ is an $(\mathbb A^1,et)$ local equivalence, 
then all the $m_I:F_I\to G_I$ are $(\mathbb A^1,et)$ local equivalence and
for all $(H_I,w_{IJ})\in C(\Var(k)^{2,(sm)pr}/(\tilde S_I))$ which is $\mathbb A^1$ local for the etale topology
\begin{eqnarray*}
(\Hom(L(m_I),E_{et}(H_I))):\Hom(L(G_I,v_{IJ}),E_{et}(H_I,w_{IJ}))\to\Hom(L(F_I,u_{IJ}),E_{et}(H_I,w_{IJ})) 
\end{eqnarray*}
is a quasi-isomorphism (of diagrams of complexes of abelian groups).
\end{itemize}
\end{defi}

\begin{prop}\label{ca1Var12IJ}
Let $S\in\Var(k)$. Let $S=\cup_{i=1}^l S_i$ an open affine cover and denote by $S_I=\cap_{i\in I} S_i$.
Let $i_i:S_i\hookrightarrow\tilde S_i$ closed embeddings, with $\tilde S_i\in\Var(k)$.
\begin{itemize}
\item[(i)]Then for $F\in C(\Var(k)^{2,(sm)}/(\tilde S_I)^{op})$, 
$C_*F$ is $\mathbb A^1$ local for the etale topology
and  $c(F):F\to C_*F$ is an equivalence $(\mathbb A^1,et)$ local.
\item[(i)']Then for $F\in C(\Var(k)^{2,(sm)pr}/(\tilde S_I)^{op})$, 
$C_*F$ is $\mathbb A^1$ local for the etale topology
and  $c(F):F\to C_*F$ is an equivalence $(\mathbb A^1,et)$ local.
\item[(ii)]A morphism $m:F\to G$ with $F,G\in C(\Var(k)^{2,(sm)}/(\tilde S_I)^{op})$ 
is an $(\mathbb A^1,et)$ local equivalence
if and only if $a_{et}H^nC_*\Cone(m)=0$ for all $n\in\mathbb Z$.
\item[(ii)']A morphism $m:F\to G$ with $F,G\in C(\Var(k)^{2,(sm)pr}/(\tilde S_I)^{op})$ 
is an $(\mathbb A^1,et)$ local equivalence
if and only if $a_{et}H^nC_*\Cone(m)=0$ for all $n\in\mathbb Z$.
\item[(iii)]A morphism $m:F\to G$ with $F,G\in C(\Var(k)^{2,(sm)}/(\tilde S_I)^{op})$ 
is an $(\mathbb A^1,et)$ local equivalence if and only if there exists 
\begin{eqnarray*}
\left\{((X_{1,\alpha,I},Z_{1,\alpha,I})/\tilde S_I,u^1_{IJ}),\alpha\in\Lambda_1\right\},\ldots,
\left\{((X_{r,\alpha,I},Z_{r,\alpha,I})/\tilde S_I,u^r_{IJ}),\alpha\in\Lambda_r\right\}
\subset\Var(k)^{2,(sm)}/(\tilde S_I)^{op}
\end{eqnarray*}
with 
\begin{equation*}
u^l_{IJ}:(X_{l,\alpha,J},Z_{l,\alpha,J})/\tilde S_J\to 
(X_{l,\alpha,I}\times\tilde S_{J\backslash I},Z_{l,\alpha,I}\times\tilde S_{J\backslash I})/\tilde S_J
\end{equation*}
such that we have in $\Ho_{et}(C(\Var(k)^{2,(sm)}/(\tilde S_I)^{op}))$
\begin{eqnarray*}
\Cone(m)\xrightarrow{\sim}\Cone( \\ \oplus_{\alpha\in\Lambda_1} 
\Cone((\mathbb Z((X_{1,\alpha,I}\times\mathbb A^1,Z_{1,\alpha,I}\times\mathbb A^1)/\tilde S_I),\mathbb Z(u_{IJ}^1\times I)) 
\to(\mathbb Z((X_{1,\alpha,I},Z_{1,\alpha,I})/\tilde S_I),\mathbb Z(u_{IJ}^1))) \\
\to\cdots\to \\ \oplus_{\alpha\in\Lambda_r}
\Cone((\mathbb Z((X_{r,\alpha,I}\times\mathbb A^1,Z_{r,\alpha,I}\times\mathbb A^1)/\tilde S_I),\mathbb Z(u_{IJ}^r\times I)) 
\to(\mathbb Z((X_{r,\alpha,I},Z_{r,\alpha,I})/\tilde S_I),\mathbb Z(u^r_{IJ}))))
\end{eqnarray*}
\item[(iii)']A morphism $m:F\to G$ with $F,G\in C(\Var(k)^{2,(sm)pr}/(\tilde S_I)^{op})$ 
is an $(\mathbb A^1,et)$ local equivalence if and only if there exists 
\begin{eqnarray*}
\left\{((Y_{1,\alpha,I}\times\tilde S_I,Z_{1,\alpha,I})/\tilde S_I,u^1_{IJ}),\alpha\in\Lambda_1\right\},\ldots,
\left\{((Y_{r,\alpha,I}\times\tilde S_I,Z_{r,\alpha,I})/\tilde S_I,u^r_{IJ}),\alpha\in\Lambda_r\right\} \\
\subset\Var(k)^{2,(sm)pr}/(\tilde S_I)
\end{eqnarray*}
with 
\begin{equation*}
u^l_{IJ}:(Y_{l,\alpha,J}\times\tilde S_J,Z_{l,\alpha,J})/\tilde S_J\to 
(Y_{l,\alpha,I}\times\tilde S_J,Z_{l,\alpha,I}\times\tilde S_{J\backslash I})/\tilde S_J
\end{equation*}
such that we have in $\Ho_{et}(C(\Var(k)^{2,(sm)}/(\tilde S_I)^{op}))$
\begin{eqnarray*}
\Cone(m)\xrightarrow{\sim}\Cone(\oplus_{\alpha\in\Lambda_1} \\
\Cone((\mathbb Z((Y_{1,\alpha,I}\times\mathbb A^1\times\tilde S_I,Z_{1,\alpha,I}\times\mathbb A^1)/\tilde S_I),
\mathbb Z(u_{IJ}^1\times I)) 
\to(\mathbb Z((Y_{1,\alpha,I}\times S,Z_{1,\alpha,I})/\tilde S_I),\mathbb Z(u_{IJ}))) \\
\to\cdots\to\oplus_{\alpha\in\Lambda_r} \\
\Cone((\mathbb Z((Y_{r,\alpha,I}\times\mathbb A^1\times\tilde S_I,Z_{r,\alpha,I}\times\mathbb A^1)/\tilde S_I),
\mathbb Z(u_{IJ}^r\times I)) 
\to(\mathbb Z((Y_{r,\alpha,I}\times\tilde S_I,Z_{r,\alpha})/\tilde S_I),\mathbb Z(u_{IJ}^r)))
\end{eqnarray*}
\item[(iv)] A similar statement then (iii) holds for equivalence $(\mathbb A^1,et)$ local
$m:F\to G$ with $F,G\in C(\Var(k)^{2,(sm)}/(\tilde S_I))$
\item[(iv)'] A similar statement then (iii) holds for equivalence $(\mathbb A^1,et)$ local
$m:F\to G$ with $F,G\in C(\Var(k)^{2,(sm)pr}/(\tilde S_I))$
\end{itemize}
\end{prop}

\begin{proof}
Similar to the proof of proposition \ref{ca1Var12}. See Ayoub's thesis for example.
\end{proof}

In the filtered case we also consider :

\begin{defi}
Let $S\in\Var(k)$. Let $S=\cup_{i=1}^l S_i$ an open affine cover and denote by $S_I=\cap_{i\in I} S_i$.
Let $i_i:S_i\hookrightarrow\tilde S_i$ closed embeddings, with $\tilde S_i\in\SmVar(k)$. 
\begin{itemize}
\item[(i)]A filtered complex $(G,F)\in C_{fil}(\Var(k)^{2,(sm)}/S)$ 
is said to be $r$-filtered $\mathbb A^1$ homotopic if 
$\ad(p_a^*,p_{a*})(G,F):(G,F)\to p_{a*}p_a^*(G,F)$ is an $r$-filtered homotopy equivalence.
\item[(i)']A filtered complex $(G,F)\in C_{fil}(\Var(k)^{2,(sm)}/(\tilde S_I))$ 
is said to be $r$-filtered $\mathbb A^1$ homotopic if 
$\ad(p_a^*,p_{a*})(G,F):(G,F)\to p_{a*}p_a^*(G,F)$ is an $r$-filtered homotopy equivalence.
\item[(ii)]A filtered complex $(G,F)\in C_{fil}(\Var(k)^{2,(sm)pr}/S)$ 
is said to be $r$-filtered $\mathbb A^1$ homotopic if 
$\ad(p_a^*,p_{a*})(G,F):(G,F)\to p_{a*}p_a^*(G,F)$ is an $r$-filtered homotopy equivalence.
\item[(ii)']A filtered complex $(G,F)\in C_{fil}(\Var(k)^{2,(sm)pr}/(\tilde S_I))$ 
is said to be $r$-filtered $\mathbb A^1$ homotopic if 
$\ad(p_a^*,p_{a*})(G,F):(G,F)\to p_{a*}p_a^*(G,F)$ is an $r$-filtered homotopy equivalence.
\end{itemize}
\end{defi}

We will use to compute the algebraic De Rahm realization functor the followings

\begin{thm}\label{DDADM12}
Let $S\in\Var(k)$. 
\begin{itemize}
\item[(i)]Let $\phi:F^{\bullet}\to G^{\bullet}$ an etale local equivalence 
with $F^{\bullet},G^{\bullet}\in C(\Var(k)^{2,sm}/S)$.
If $F^{\bullet}$ and $G^{\bullet}$ are $\mathbb A^1$ local and admit tranferts 
then $\phi:F^{\bullet}\to G^{\bullet}$ is a Zariski local equivalence.
Hence if $F\in C(\Var(k)^{2,sm}/S)$ is $\mathbb A^1$ local and admits transfert 
\begin{equation*}
k:E_{zar}(F)\to E_{et}(E_{zar}(F))=E_{et}(F) 
\end{equation*}
is a Zariski local equivalence.
\item[(ii)]Let $\phi:F^{\bullet}\to G^{\bullet}$ an etale local equivalence
with $F^{\bullet},G^{\bullet}\in C(\Var(k)^{2,smpr}/S)$.
If $F^{\bullet}$ and $G^{\bullet}$ are $\mathbb A^1$ local and admit tranferts 
then $\phi:F^{\bullet}\to G^{\bullet}$ is a Zariski local equivalence.
Hence if $F\in C(\Var(k)^{2,smpr}/S)$ is $\mathbb A^1$ local and admits transfert 
\begin{equation*}
k:E_{zar}(F)\to E_{et}(E_{zar}(F))=E_{et}(F) 
\end{equation*}
is a Zariski local equivalence.
\end{itemize}
\end{thm}

\begin{proof}
Similar to the proof of theorem \ref{DDADM}.
\end{proof}

\begin{thm}\label{DDADM12fil}
Let $S\in\Var(k)$. Let $S=\cup_{i=1}^l S_i$ an open affine cover and denote by $S_I=\cap_{i\in I} S_i$.
Let $i_i:S_i\hookrightarrow\tilde S_i$ closed embeddings, with $\tilde S_i\in\Var(k)$. 
\begin{itemize}
\item[(i)]Let $\phi:(F^{\bullet},F)\to (G^{\bullet},F)$ a filtered etale local equivalence 
with $(F^{\bullet},F),(G^{\bullet},F)\in C_{fil}(\Var(k)^{2,sm}/S)$.
If $(F^{\bullet},F)$ and $(G^{\bullet},F)$ are $r$-filtered $\mathbb A^1$ homotopic and admit tranferts 
then $\phi:(F^{\bullet},F)\to (G^{\bullet},F)$ is an $r$-filtered Zariski local equivalence.
Hence if $(G,F)\in C(\Var(k)^{2,sm}/S)$ is $r$-filtered $\mathbb A^1$ homotopic and admits transfert 
\begin{equation*}
k:E_{zar}(G,F)\to E_{et}(E_{zar}(G,F))=E_{et}(G,F) 
\end{equation*}
is an $r$-filtered Zariski local equivalence.
\item[(i)']Let $\phi:(F^{\bullet},F)\to (G^{\bullet},F)$ a filtered etale local equivalence 
with $((F_I^{\bullet},F),u_{IJ}),((G_I^{\bullet},F),v_{IJ})\in C_{fil}(\Var(k)^{2,sm}/(\tilde S_I))$.
If $((F^{\bullet},F),u_{IJ})$ and $((G^{\bullet},F),v_{IJ})$ are $r$-filtered $\mathbb A^1$ homotopic and admit tranferts 
then $\phi:((F^{\bullet},F),u_{IJ})\to ((G^{\bullet},F),v_{IJ})$ is an $r$-filtered Zariski local equivalence.
Hence if $((G_I,F),u_{IJ})\in C(\Var(k)^{2,sm}/S)$ is $r$-filtered $\mathbb A^1$ homotopic and admits transfert 
\begin{equation*}
k:(E_{zar}(G_I,F),u_{IJ})\to (E_{et}(E_{zar}(G_I,F)),u_{IJ})=(E_{et}(G,F),u_{IJ}) 
\end{equation*}
is an $r$-filtered Zariski local equivalence.
\item[(ii)]Let $\phi:(F^{\bullet},F)\to (G^{\bullet},F)$ a filtered etale local equivalence
with $(F^{\bullet},F),(G^{\bullet},F)\in C_{fil}(\Var(k)^{2,smpr}/S)$.
If $F^{\bullet}$ and $G^{\bullet}$ are $r$-filtered $\mathbb A^1$ homotopic and admit tranferts 
then $\phi:(F^{\bullet},F)\to (G^{\bullet},F)$ is an $r$-filtered Zariski local equivalence.
Hence if $(G,F)\in C(\Var(k)^{2,smpr}/S)$ is $r$-filtered $\mathbb A^1$ homotopic and admits transfert 
\begin{equation*}
k:E_{zar}(F)\to E_{et}(E_{zar}(F))=E_{et}(F) 
\end{equation*}
is an $r$-filtered Zariski local equivalence.
\item[(ii)']Let $\phi:(F^{\bullet},F)\to (G^{\bullet},F)$ a filtered etale local equivalence 
with $((F_I^{\bullet},F),u_{IJ}),((G_I^{\bullet},F),v_{IJ})\in C_{fil}(\Var(k)^{2,smpr}/(\tilde S_I))$.
If $((F^{\bullet},F),u_{IJ})$ and $((G^{\bullet},F),v_{IJ})$ are $r$-filtered $\mathbb A^1$ homotopic and admit tranferts 
then $\phi:((F^{\bullet},F),u_{IJ})\to ((G^{\bullet},F),v_{IJ})$ is an $r$-filtered Zariski local equivalence.
Hence if $((G_I,F),u_{IJ})\in C(\Var(k)^{2,smpr}/S)$ is $r$-filtered $\mathbb A^1$ homotopic and admits transfert 
\begin{equation*}
k:(E_{zar}(G_I,F),u_{IJ})\to (E_{et}(E_{zar}(G_I,F)),u_{IJ})=(E_{et}(G,F),u_{IJ}) 
\end{equation*}
is an $r$-filtered Zariski local equivalence.
\end{itemize}
\end{thm}

\begin{proof}
Similar to the proof of theorem \ref{DDADM12}.
\end{proof}

\subsection{The Borel-Moore Corti-Hanamura resolution functors 
$R^{CH}$, $\hat R^{CH}$, $R^{0CH}$, and $\hat R^{0CH}$}

Let $k$ a field of caracteristic zero.

\begin{defi}\label{desVardef}
\begin{itemize}
\item[(i)]Let $X_0\in\Var(k)$ and $Z\subset X_0$ a closed subset.
A desingularization of $(X_0,Z)$ is a pair of complex varieties $(X,D)\in\Var^2(k))$,
together with a morphism of pair of varieties $\epsilon:(X,D)\to(X_0,\Delta)$ with $Z\subset\Delta$ such that 
\begin{itemize}
\item $X\in\SmVar(k)$ and 
$D:=\epsilon^{-1}(\Delta)=\epsilon^{-1}(Z)\cup(\cup_i E_i)\subset X$ is a normal crossing divisor
\item $\epsilon:X\to X_0$ is a proper modification with discriminant $\Delta$, 
that is $\epsilon:X\to X_0$ is proper and $\epsilon:X\backslash D\xrightarrow{\sim}X\backslash\Delta$ is an isomorphism.
\end{itemize}
\item[(ii)]Let $X_0\in\Var(k)$ and $Z\subset X_0$ a closed subset such that $X_0\backslash Z$ is smooth.
A strict desingularization of $(X_0,Z)$ is a pair of complex varieties $(X,D)\in\Var^2(\mathbb C))$,
together with a morphism of pair of varieties $\epsilon:(X,D)\to(X_0,Z)$ such that 
\begin{itemize}
\item $X\in\SmVar(k)$ and $D:=\epsilon^{-1}(Z)\subset X$ is a normal crossing divisor
\item $\epsilon:X\to X_0$ is a proper modification with discriminant $Z$, 
that is $\epsilon:X\to X_0$ is proper and $\epsilon:X\backslash D\xrightarrow{\sim}X\backslash Z$ is an isomorphism.
\end{itemize}
\end{itemize}
\end{defi}

We have the following well known resolution of singularities of complex algebraic varieties
and their functorialities :

\begin{thm}\label{desVar}
\begin{itemize}
\item[(i)] Let $X_0\in\Var(k)$ and $Z\subset X_0$ a closed subset.
There exists a desingularization of $(X_0,Z)$, that is a pair of complex varieties $(X,D)\in\Var^2(k))$,
together with a morphism of pair of varieties $\epsilon:(X,D)\to(X_0,\Delta)$ with $Z\subset\Delta$ such that 
\begin{itemize}
\item $X\in\SmVar(k)$ and 
$D:=\epsilon^{-1}(\Delta)=\epsilon^{-1}(Z)\cup(\cup_iE_i)\subset X$ is a normal crossing divisor
\item $\epsilon:X\to X_0$ is a proper modification with discriminant $\Delta$, 
that is $\epsilon:X\to X_0$ is proper and $\epsilon:X\backslash D\xrightarrow{\sim}X\backslash\Delta$ is an isomorphism.
\end{itemize}
\item[(ii)] Let $X_0\in\PVar(k)$ and $Z\subset X_0$ a closed subset such that $X_0\backslash Z$ is smooth.
There exists a strict desingularization of $(X_0,Z)$, that is a pair of complex varieties $(X,D)\in\PVar^2(k))$,
together with a morphism of pair of varieties $\epsilon:(X,D)\to(X_0,Z)$ such that 
\begin{itemize}
\item $X\in\PSmVar(k)$ and $D:=\epsilon^{-1}(Z)\subset X$ is a normal crossing divisor
\item $\epsilon:X\to X_0$ is a proper modification with discriminant $Z$, 
that is $\epsilon:X\to X_0$ is proper and $\epsilon:X\backslash D\xrightarrow{\sim}X\backslash Z$ is an isomorphism.
\end{itemize}
\end{itemize}
\end{thm}

\begin{proof}
\noindent(i):Standard. See \cite{PS} for example.

\noindent(ii):Follows immediately from (i).
\end{proof}

We use this theorem to construct a resolution of a morphism by Corti-Hanamura morphisms,
we will need these resolution in the definition of the filtered De Rham realization functor :

\begin{defiprop}\label{RCHdef0}
\begin{itemize}
\item[(i)]Let $h:V\to S$ a morphism, with $V,S\in\Var(k)$. 
Let $\bar{S}\in\PVar(k)$ be a compactification of $S$.
\begin{itemize}
\item There exist a compactification $\bar X_0\in\PVar(\mathbb C)$ of $V$ such that $h:V\to S$
extend to a morphism $\bar f_0=\bar{h}_0:\bar X_0\to\bar{S}$. Denote by $\bar Z=\bar X_0\backslash V$.
We denote by $j:V\hookrightarrow\bar X_0$ the open embedding and by 
$i_0:\bar Z\hookrightarrow\bar X_0$ the complementary closed embedding.
We then consider $X_0:=\bar f_0^{-1}(S)\subset\bar X_0$ the open subset, $f_0:=\bar f_{0|X_0}:X_0\to S$, $Z=\bar Z\cap X_0$,
and we denote again $j:V\hookrightarrow X_0$ the open embedding and by $i_0:Z\hookrightarrow X_0$ the complementary closed embedding.
\item In the case $V$ is smooth, we take, using theorem \ref{desVar}(ii), a strict desingularization 
$\bar\epsilon:(\bar X,\bar D)\to(\bar X_0,\bar Z)$ of the pair $(\bar X_0,\bar Z)$,
with $\bar X\in\PSmVar(\mathbb C)$ and $\bar D=\cup_{i=1}^s\bar D_i\subset\bar X$ a normal crossing divisor.  
We denote by $i_{\bullet}:\bar D_{\bullet}\hookrightarrow\bar X=\bar X_{c(\bullet)}$ the morphism of simplicial varieties
given by the closed embeddings $i_I:\bar D_I=\cap_{i\in I}\bar D_i\hookrightarrow\bar X$.
Then the morphisms $\bar f:=\bar f_0\circ\bar\epsilon:\bar X\to\bar S$ and 
$\bar f_{D_{\bullet}}:=\bar f\circ i_{\bullet}:\bar D_{\bullet}\to\bar S$ are projective since
$\bar X$ and $\bar D_I$ are projective varieties.  
We then consider $(X,D):=\bar\epsilon^{-1}(X_0,Z)$, $\epsilon:=\bar\epsilon_{|X}:(X,D)\to(X_0,Z)$
We denote again by $i_{\bullet}:D_{\bullet}\hookrightarrow X=X_{c(\bullet)}$ the morphism of simplicial varieties
given by the closed embeddings $i_I:D_I=\cap_{i\in I}D_i\hookrightarrow X$.
Then the morphisms $f:=f_0\circ\epsilon:X\to S$ and 
$f_{D_{\bullet}}:=f\circ i_{\bullet}:D_{\bullet}\to S$ are projective since $\bar f:\bar X_0\to\bar S$ is projective.
\end{itemize}
\item[(ii)]Let $g:V'/S\to V/S$ a morphism, with $V'/S=(V',h'),V/S=(V,h)\in\Var(k)/S$
\begin{itemize}
\item Take (see (i)) a compactification $\bar X_0\in\PVar(\mathbb C)$ of $V$ such that $h:V\to S$
extend to a morphism $\bar f_0=\bar{h}_0:\bar X_0\to\bar{S}$. Denote by $\bar Z=\bar X_0\backslash V$.
Then, there exist a compactification $\bar X'_0\in\PVar(\mathbb C)$ of $V'$ such that $h':V'\to S$
extend to a morphism $\bar f'_0=\bar{h}'_0:\bar X'_0\to\bar{S}$, $g:V'\to V$ extend to a morphism 
$\bar{g}_0:\bar X'_0\to\bar X_0$ and $\bar f_0\circ\bar{g}_0=\bar f'_0$ 
that is $\bar{g}_0$ is gives a morphism $\bar{g}_0:\bar X'_0/\bar{S}\to\bar X_0/\bar{S}$. 
Denote by $\bar Z'=\bar X'_0\backslash V'$. We then have the following commutative diagram
\begin{equation*}
\xymatrix{V\ar[r]^j & \bar X_0 & \, & \bar Z\ar[ll]_i \\
V'\ar[r]^{j'}\ar[u]^{\bar g} & \bar X'_0\ar[u]^{\bar{g}_0} & \bar Z'\ar[l]_{i'} & 
\bar{g_0}^{-1}(\bar Z)\ar[l]_{i''_{g,0}}\ar[u]^{\bar{g}'}:i'_{g,0}}
\end{equation*}
It gives the following commutative diagram
\begin{equation*}
\xymatrix{V\ar[r]^j & X_0:=\bar f_0^{-1}(S) & \, & Z\ar[ll]_i \\
V'\ar[r]^{j'}\ar[u]^g & X'_0:=\bar f_0^{'-1}(S)\ar[u]^{\bar{g}_0} & Z'\ar[l]_{i'} & 
\bar{g_0}^{-1}(Z)\ar[l]_{i''_{g,0}}\ar[u]^{\bar{g}'_0}:i'_{g,0}}
\end{equation*}
\item In the case $V$ and $V'$ are smooth, we take using theorem \ref{desVar} a strict desingularization 
$\bar\epsilon:(\bar X,\bar D)\to(\bar X_0,\bar Z)$ of $(\bar X_0,\bar Z)$.  
Then there exist a strict desingularization
$\bar\epsilon'_{\bullet}:(\bar X',\bar D')\to(\bar X'_0,\bar Z')$ of $(\bar X'_0,\bar Z')$
and a morphism $\bar g:\bar X'\to\bar X$ such that the following diagram commutes
\begin{equation*}
\xymatrix{\bar X'_0\ar[r]^{\bar{g}_0} & \bar X_0 \\
\bar X'\ar[u]^{\bar\epsilon'}\ar[r]^{\bar g} & \bar X\ar[u]^{\bar\epsilon}}.
\end{equation*}  
We then have the following commutative diagram in $\Fun(\Delta,\Var(k))$
\begin{equation*}
\xymatrix{V=V_{c(\bullet)}\ar[r]^j & \bar X=\bar X_{c(\bullet)} & \, & \bar D_{s_g(\bullet)}\ar[ll]_{i_{\bullet}} \\
V'=V'_{c(\bullet)}\ar[r]^{j'}\ar[u]^g & \bar X'=\bar X'_{c(\bullet)}\ar[u]^{\bar{g}} & \bar D'_{\bullet}\ar[l]_{i'_{\bullet}} & 
\bar{g}^{-1}(\bar D_{s_g(\bullet)})\ar[l]_{i''_{g\bullet}}\ar[u]^{\bar{g}'_{\bullet}}:i'_{g\bullet}}
\end{equation*}
where $i_{\bullet}:\bar D_{\bullet}\hookrightarrow\bar X_{\bullet}$ the morphism of simplicial varieties
given by the closed embeddings $i_n:\bar D_n\hookrightarrow\bar X_n$,
and $i'_{\bullet}:\bar D'_{\bullet}\hookrightarrow\bar X'_{\bullet}$ the morphism of simplicial varieties
given by the closed embeddings $i'_n:\bar D'_n\hookrightarrow\bar X'_n$. 
It gives the commutative diagram in $\Fun(\Delta,\Var(k))$
\begin{equation*}
\xymatrix{V=V_{c(\bullet)}\ar[r]^j & X:=\bar\epsilon^{-1}(X_0)=X_{c(\bullet)} & \, & D_{s_g(\bullet)}\ar[ll]_{i_{\bullet}} \\
V'=V'_{c(\bullet)}\ar[r]^{j'}\ar[u]^g & X':=\bar\epsilon^{',-1}(X'_0)=X'_{c(\bullet)}\ar[u]^{\bar{g}} & 
D'_{\bullet}\ar[l]_{i'_{\bullet}} & \bar{g}^{-1}(D_{s_g(\bullet)})\ar[l]_{i''_{g\bullet}}\ar[u]^{\bar{g}'_{\bullet}}:i'_{g\bullet}}
\end{equation*}
\end{itemize}
\end{itemize}
\end{defiprop}

\begin{proof}
\noindent(i): Let $\bar X_{00}\in\PVar(\mathbb C)$ be a compactification of $V$. 
Let $l_0:\bar X_0=\bar{\Gamma_h}\hookrightarrow\bar X_{00}\times\bar{S}$ be the closure of the graph of $h$ and
$\bar f_0:=p_{\bar{S}}\circ l_0:\bar X_0\hookrightarrow\bar X_{00}\times\bar{S}\to\bar{S}$, 
$\epsilon_{\bar X_0}:=p_{\bar X_{00}}\circ l_0:\bar X_0\hookrightarrow\bar X_{00}\times\bar{S}\to\bar X_{00}$
be the restriction to $\bar X_0$ of the projections.
Then, $\bar X\in\PVar(\mathbb C)$, $\epsilon_{\bar X_0}:\bar X_0\to\bar X_{00}$ is a proper modification which does not affect
the open subset $V\subset\bar X_0$, and $\bar f_0=\bar{h}_0:\bar X_0\to\bar{S}$ is a compactification of $h$.

\noindent(ii): There is two things to prove: 
\begin{itemize}
\item Let $\bar f_0:\bar X_0\to\bar S$ a compactification of $h:V\to S$ and 
$\bar f'_{00}:\bar X'_{00}\to\bar S$ a compactification of $h':V'\to S$ (see (i)). 
Let $l_0:\bar X'_0\hookrightarrow\bar{\Gamma_g}\subset\bar X'_{00}\times_{\bar S}\bar X_0$ be the closure of the graph of $g$,
$\bar f'_0:=(\bar f'_{00},\bar f_0)\circ l_0:\bar X'_0\hookrightarrow\bar X'_{00}\times_S\bar X_0\to\bar S$ and
$\bar{g}_0:=p_{\bar X_0}\circ l_0:\bar X'_0\hookrightarrow\bar X'_{00}\times_{\bar S}\bar X_0\to\bar X_0$, 
$\epsilon_{\bar X'_{00}}:=p_{\bar X'_0}\circ i:\bar X'_0\hookrightarrow\bar X'_{00}\times_{\bar S}\bar X_0\to\bar X'_{00}$
be the restriction to $X$ of the projections.
Then $\epsilon_{\bar X'_{00}}:\bar X'_0\to\bar X'_{00}$ is a proper modification which does not affect
the open subset $V'\subset\bar X'_0$, $\bar f'_0:\bar X'_0\to\bar S$ is an other compactification of $h':V'\to S$
and $\bar{g}_0:\bar X'_0\to\bar X_0$ is a compactification of $g$.
\item In the case $V$ and $V'$ are smooth, we take, using theorem \ref{desVar}, a strict desingularization 
$\bar\epsilon:(\bar X,\bar D)\to(\bar X_0,\bar Z)$ of the pair $(\bar X_0,\bar Z)$.
Take then, using theorem \ref{desVar}, a strict desingularization 
$\bar\epsilon'_1:(\bar X',\bar D')\to(\bar X\times_{\bar X_0}\bar X'_0,\bar X\times_{\bar X_0}\bar Z')$ 
of the pair $(\bar X\times_{\bar X_0}\bar X'_0,\bar X\times_{\bar X_0}\bar Z')$.
We consider then following commutative diagram whose square is cartesian :
\begin{equation*}
\xymatrix{ \, & \bar X'_0\ar[r]^{\bar g_0} &  X_0 \\
\, & \bar X\times_{\bar X_0}\bar X'_0\ar[u]^{\bar\epsilon'_0}\ar[r]^{\bar{g}'_0} & \bar X\ar[u]^{\epsilon} \\
\bar X'\ar[ruu]^{\epsilon'}\ar[rru]^{\bar{g}}\ar[ru]^{\bar\epsilon'_1} & \, & \,}
\end{equation*}  
and $\bar\epsilon':=\bar\epsilon'_0\circ\bar\epsilon'_1:(\bar X',\bar D')\to(\bar X'_0,\bar Z')$
is a strict desingularization of the pair $(\bar X\times_{\bar X_0}\bar X'_0, \bar X\times_{\bar X_0}\bar Z')$.
\end{itemize}
\end{proof}

Let $S\in\Var(k)$. Recall we have the dual functor
\begin{eqnarray*}
\mathbb D^0_S:C(\Var(k)/S)\to C(\Var(k)/S), \; F\mapsto\mathbb D^0_S(F):=\mathcal Hom(F,E_{et}(\mathbb Z(S/S)))
\end{eqnarray*}
which induces the functor
\begin{eqnarray*}
L\mathbb D_S:D_{\tau}(\Var(k)/S)\to D_{\tau}(\Var(k)/S), \; 
F\mapsto L\mathbb D_S(F):=\mathbb D^0_S(LF):=\mathcal Hom(LF,E_{et}(\mathbb Z(S/S)))
\end{eqnarray*}
with $\tau$ a topology on $\Var(k)$.

We will use the following resolutions of representable presheaves by Corti-Hanamura presheaves and
their the functorialities.

\begin{defi}\label{RCHdef} 
\begin{itemize}
\item[(i)]Let $h:U\to S$ a morphism, with $U,S\in\Var(k)$ and $U$ smooth.
Take, see definition-proposition \ref{RCHdef0},
$\bar f_0=\bar{h}_0:\bar X_0\to\bar S$ a compactification of $h:U\to S$ and denote by $\bar Z=\bar X_0\backslash U$.
Take, using theorem \ref{desVar}(ii), a strict desingularization 
$\bar\epsilon:(\bar X,\bar D)\to(\bar X_0,\bar Z)$ of the pair $(\bar X_0,\bar Z)$, with $\bar X\in\PSmVar(k)$ and 
$\bar D:=\epsilon^{-1}(\bar Z)=\cup_{i=1}^s\bar D_i\subset\bar X$ a normal crossing divisor.  
We denote by $i_{\bullet}:\bar D_{\bullet}\hookrightarrow\bar X=\bar X_{c(\bullet)}$ the morphism of simplicial varieties
given by the closed embeddings $i_I:\bar D_I=\cap_{i\in I}\bar D_i\hookrightarrow\bar X$
We denote by $j:U\hookrightarrow\bar X$ the open embedding and by $p_S:\bar X\times S\to S$ 
and $p_S:U\times S\to S$ the projections.
Considering the graph factorization $\bar f:\bar X\xrightarrow{\bar l}\bar X\times\bar S\xrightarrow{p_{\bar S}}\bar S$
of $\bar f:\bar X\to\bar S$, where $\bar l$ is the graph embedding and $p_{\bar S}$ the projection,
we get closed embeddings $l:=\bar l\times_{\bar S}S:X\hookrightarrow\bar X\times S$ and 
$l_{D_I}:=\bar D_I\times_{\bar X} l:D_I\hookrightarrow\bar D_I\times S$.
We then consider the following map in $C(\Var(k)^2/S)$
\begin{eqnarray*}
r_{(\bar X,\bar D)/S}(\mathbb Z(U/S)):R_{(\bar X,\bar D)/S}(\mathbb Z(U/S)) \\
\xrightarrow{:=}p_{S*}E_{et}(\Cone(\mathbb Z(i_{\bullet}\times I):
(\mathbb Z((\bar D_{\bullet}\times S,D_{\bullet})/\bar X\times S),u_{IJ})\to
\mathbb Z((\bar X\times S,X)/\bar X\times S))) \\
\xrightarrow{p_{S*}E_{et}(0,k\circ\ad((j\times I)^*,(j\times I)_*)(\mathbb Z((\bar X\times S,X)/\bar X\times S)))}
p_{S*}E_{et}(\mathbb Z((U\times S,U)/U\times S))=:\mathbb D_S^{12}(\mathbb Z(U/S)).
\end{eqnarray*}
Note that $\mathbb Z((\bar D_I\times S,D_I)/\bar X\times S)$ and $\mathbb Z((\bar X\times S,X)/\bar X\times S))$
are obviously $\mathbb A^1$ invariant.
Note that $r_{(X,D)/S}$ is NOT an equivalence $(\mathbb A^1,et)$ local by proposition \ref{rho12} since
$\rho_{\bar X\times S*}\mathbb Z((\bar D_{\bullet}\times S,D_{\bullet})/\bar X\times S)=0$ and 
$\rho_{\bar X\times S*}\ad((j\times I)^*,(j\times I)_*)(\mathbb Z((\bar X\times S,X)/\bar X\times S)))$ 
is not an equivalence $(\mathbb A^1,et)$ local.
\item[(ii)]Let $g:U'/S\to U/S$ a morphism, with $U'/S=(U',h'),U/S=(U,h)\in\Var(k)/S$, with $U$ and $U'$ smooth.
Take, see definition-proposition \ref{RCHdef0}(ii),a compactification $\bar f_0=\bar h:\bar X_0\to\bar S$ of $h:U\to S$ 
and a compactification $\bar f'_0=\bar{h}':\bar X'_0\to\bar S$ of $h':U'\to S$ such that 
$g:U'/S\to U/S$ extend to a morphism $\bar g_0:\bar X'_0/\bar S\to\bar X_0/\bar S$. 
Denote $\bar Z=\bar X_0\backslash U$ and $\bar Z'=\bar X'_0\backslash U'$.
Take, see definition-proposition \ref{RCHdef0}(ii), a strict desingularization 
$\bar\epsilon:(\bar X,\bar D)\to(\bar X_0,\bar Z)$ of $(\bar X_0,\bar Z)$,
a strict desingularization $\bar\epsilon'_{\bullet}:(\bar X',\bar D')\to(\bar X'_0,\bar Z')$ of $(\bar X'_0,\bar Z')$
and a morphism $\bar g:\bar X'\to\bar X$ such that the following diagram commutes
\begin{equation*}
\xymatrix{\bar X'_0\ar[r]^{\bar{g}_0} & \bar X_0 \\
\bar X'\ar[u]^{\bar\epsilon'}\ar[r]^{\bar g} & \bar X\ar[u]^{\bar\epsilon}}.
\end{equation*} 
We then have, see definition-proposition \ref{RCHdef0}(ii),
the following commutative diagram in $\Fun(\Delta,\Var(k))$
\begin{equation}\label{RCHdia}
\xymatrix{U=U_{c(\bullet)}\ar[r]^j & \bar X=\bar X_{c(\bullet)} & \, & \bar D_{s_g(\bullet)}\ar[ll]_{i_{\bullet}} \\
U'=U'_{c(\bullet)}\ar[r]^{j'}\ar[u]^g & \bar X'=\bar X'_{c(\bullet)}\ar[u]^{\bar{g}} & \bar D'_{\bullet}\ar[l]_{i'_{\bullet}} & 
\bar{g}^{-1}(\bar D_{s_g(\bullet)})\ar[l]_{i''_{g\bullet}}\ar[u]^{\bar{g}'_{\bullet}}:i'_{g\bullet}}
\end{equation}
Denote by $p_S:\bar X\times S\to S$ and $p'_S:\bar X'\times S\to S$ the projections
We then consider the following map in $C(\Var(k)^2/S)$
\begin{eqnarray*}
R_S^{CH}(g):R_{(\bar X,\bar D)/S}(\mathbb Z(U/S))\xrightarrow{:=} \\
p_{S*}E_{et}(\Cone(\mathbb Z(i_{\bullet}\times I):
(\mathbb Z((\bar D_{s_g(\bullet)}\times S,D_{s_g(\bullet)})/\bar X\times S),u_{IJ})\to
\mathbb Z((\bar X\times S,X)/\bar X\times S))) \\
\xrightarrow{T((\bar g\times I),E)(-)\circ p_{S*}\ad((\bar{g}\times I)^*,(\bar{g}\times I)_*)(-)} \\
p'_{S*}E_{et}(\Cone(\mathbb Z(i'_{g\bullet}\times I): \\
(\mathbb Z((\bar{g}^{-1}(\bar D_{s_g(\bullet)})\times S,\bar g^{-1}(D_{s_g(\bullet)})/\bar X'\times S),u_{IJ})
\to\mathbb Z((\bar X'\times S,X')/\bar X'\times S)))) \\
\xrightarrow{p'_{S*}E_{et}(\mathbb Z(i''_{g\bullet}\times I),I)} \\
p'_{S*}E_{et}(\Cone(\mathbb Z(i'_{\bullet}\times I):
((\mathbb Z((\bar D'_{\bullet}\times S,D'_{\bullet)})/\bar X'\times S),u_{IJ})
\to\mathbb Z((\bar X'\times S,X')/\bar X'\times S))) \\
\xrightarrow{=:}R_{(\bar X',\bar D')/S}(\mathbb Z(U'/S))
\end{eqnarray*}
Then by the diagram (\ref{RCHdia}) and adjonction, the following diagram in $C(\Var(k)^2/S)$ obviously commutes
\begin{equation*}
\xymatrix{R_{(\bar X,\bar D)/S}(\mathbb Z(U/S))\ar[rr]^{r_{(\bar X,\bar D)/S}(\mathbb Z(U/S))}\ar[d]_{R_S^{CH}(g)} & \, & 
p_{S*}E_{et}(\mathbb Z((U\times S,U)/U\times S))=:\mathbb D_S^{12}(\mathbb Z(U/S))
\ar[d]^{D^{12}_S(g):=T(g\times I,E)(-)\circ\ad((g\times I)^*,(g\times I)_*)(E_{et}(\mathbb Z((U\times S,U)/U\times S)))} \\
R_{(\bar X',\bar D')/S}(\mathbb Z(U'/S))\ar[rr]^{r_{(\bar X',\bar D')/S}(\mathbb Z(U'/S))} & \, & 
p'_{S*}E_{et}(\mathbb Z((U'\times S,U')/U'\times S))=:\mathbb D^{12}_S(\mathbb Z(U'/S))}.
\end{equation*}
\item[(iii)]For $g_1:U''/S\to U'/S$, $g_2:U'/S\to U/S$ two morphisms
with $U''/S=(U',h''),U'/S=(U',h'),U/S=(U,h)\in\Var(k)/S$, with $U$, $U'$ and $U''$ smooth.
We get from (i) and (ii) 
a compactification $\bar f=\bar{h}:\bar X\to\bar S$ of $h:U\to S$,
a compactification $\bar f'=\bar{h}':\bar X'\to\bar S$ of $h':U'\to S$,
and a compactification $\bar f''=\bar{h}'':\bar X''\to\bar S$ of $h'':U''\to S$,
with $\bar X,\bar X',\bar X''\in\PSmVar(k)$, 
$\bar D:=\bar X\backslash U\subset\bar X$ $\bar D':=\bar X'\backslash U'\subset\bar X'$, 
and $\bar D'':=\bar X''\backslash U''\subset\bar X''$ normal crossing divisors, 
such that $g_1:U''/S\to U'/S$ extend to $\bar g_1:\bar X''/\bar S\to\bar X'/\bar S$,
$g_2:U'/S\to U/S$ extend to $\bar g_2:\bar X'/\bar S\to\bar X/\bar S$, and
\begin{eqnarray*}
R_S^{CH}(g_2\circ g_1)=R_S^{CH}(g_1)\circ R_S^{CH}(g_2):R_{(\bar X,\bar D)/S}\to R_{(\bar X'',\bar D'')/S} 
\end{eqnarray*}
\item[(iv)] For  
\begin{eqnarray*}
Q^*:=(\cdots\to\oplus_{\alpha\in\Lambda^n}\mathbb Z(U^n_{\alpha}/S)
\xrightarrow{(\mathbb Z(g^n_{\alpha,\beta}))}\oplus_{\beta\in\Lambda^{n-1}}\mathbb Z(U^{n-1}_{\beta}/S)\to\cdots)
\in C(\Var(k)/S)
\end{eqnarray*}
a complex of (maybe infinite) direct sum of representable presheaves with $U^*_{\alpha}$ smooth,
we get from (i), (ii) and (iii) the map in $C(\Var(k)^2/S)$
\begin{eqnarray*}
r_S^{CH}(Q^*):R^{CH}(Q^*):=
(\cdots\to\oplus_{\beta\in\Lambda^{n-1}}\varinjlim_{(\bar X^{n-1}_{\beta},\bar D^{n-1}_{\beta})/S}
R_{(\bar X^{n-1}_{\beta},\bar D^{n-1}_{\beta})/S}(\mathbb Z(U^{n-1}_{\beta}/S)) \\
\xrightarrow{(R_S^{CH}(g^n_{\alpha,\beta}))}\oplus_{\alpha\in\Lambda^n}\varinjlim_{(\bar X^n_{\alpha},\bar D^n_{\alpha})/S}
R_{(\bar X^n_{\alpha},\bar D^n_{\alpha})/S}(\mathbb Z(U^n_{\alpha}/S))\to\cdots)
\to\mathbb D^{12}_S(Q^*),
\end{eqnarray*}
where for $(U^n_{\alpha},h^n_{\alpha})\in\Var(k)/S$, the inductive limit run over all the compactifications 
$\bar f_{\alpha}:\bar X_{\alpha}\to\bar S$ of $h_{\alpha}:U_{\alpha}\to S$ with $\bar X_{\alpha}\in\PSmVar(k)$
and $\bar D_{\alpha}:=\bar X_{\alpha}\backslash U_{\alpha}$ a normal crossing divisor.
For $m=(m^*):Q_1^*\to Q_2^*$ a morphism with 
\begin{eqnarray*}
Q_1^*:=(\cdots\to\oplus_{\alpha\in\Lambda^n}\mathbb Z(U^n_{1,\alpha}/S)
\xrightarrow{(\mathbb Z(g^n_{\alpha,\beta}))}\oplus_{\beta\in\Lambda^{n-1}}\mathbb Z(U^{n-1}_{1,\beta}/S)\to\cdots), \\
Q_2^*:=(\cdots\to\oplus_{\alpha\in\Lambda^n}\mathbb Z(U^n_{2,\alpha}/S)
\xrightarrow{(\mathbb Z(g^n_{\alpha,\beta}))}\oplus_{\beta\in\Lambda^{n-1}}\mathbb Z(U^{n-1}_{2,\beta}/S)\to\cdots)
\in C(\Var(k)/S)
\end{eqnarray*}
complexes of (maybe infinite) direct sum of representable presheaves with $U^*_{1,\alpha}$ and $U^*_{2,\alpha}$ smooth,
we get again from (i), (ii) and (iii) a commutative diagram in $C(\Var(k)^2/S)$
\begin{equation*}
\xymatrix{R^{CH}(Q_2^*)\ar[rr]^{r_S^{CH}(Q_2^*)}\ar[d]_{R_S^{CH}(m):=(R_S^{CH}(m^*))} & \, & 
\mathbb D^{12}_S(Q_2^*)\ar[d]^{\mathbb D^{12}_S(m):=(\mathbb D^{12}_S(m^*))} \\
R^{CH}(Q_1^*)\ar[rr]^{r_S^{CH}(Q_1^*)} & \, & \mathbb D^{12}_S(Q_1^*)}.
\end{equation*}
\end{itemize}
\end{defi}

\begin{itemize}
\item Let $S\in\Var(k)$
For $(h,m,m')=(h^*,m^*,m^{'*}):Q_1^*[1]\to Q_2^*$ an homotopy with $Q_1^*,Q_2^*\in C(\Var(k)/S)$
complexes of (maybe infinite) direct sum of representable presheaves with $U^*_{1,\alpha}$ and $U^*_{2,\alpha}$ smooth,
\begin{equation*}
(R_S^{CH}(h),R_S^{CH}(m),R_S^{CH}(m'))=(R_S^{CH}(h^*),R_S^{CH}(m^*),R_S^{CH}(m^{'*})):
R^{CH}(Q_2^*)[1]\to R^{CH}(Q_1^*)
\end{equation*}
is an homotopy in $C(\Var(k)^2/S)$ using definition \ref{RCHdef} (iii). 
In particular if $m:Q_1^*\to Q_2^*$ with $Q_1^*,Q_2^*\in C(\Var(k)/S)$
complexes of (maybe infinite) direct sum of representable presheaves with $U^*_{1,\alpha}$ and $U^*_{2,\alpha}$ smooth
is an homotopy equivalence, then $R_S^{CH}(m):R^{CH}(Q_2^*)\to R^{CH}(Q_1^*)$ is an homotopy equivalence.
\item Let $S\in\SmVar(k)$. Let $F\in\PSh(\Var(k)^{sm}/S)$. Consider 
\begin{eqnarray*}
q:LF:=(\cdots\to\oplus_{(U_{\alpha},h_{\alpha})\in\Var(k)^{sm}/S}\mathbb Z(U_{\alpha}/S)
\xrightarrow{(\mathbb Z(g^n_{\alpha,\beta}))}
\oplus_{(U_{\alpha},h_{\alpha})\in\Var(k)^{sm}/S}\mathbb Z(U_{\alpha}/S)\to\cdots)\to F
\end{eqnarray*}
the canonical projective resolution given in subsection 2.3.3.
Note that the $U_{\alpha}$ are smooth since $S$ is smooth and $h_{\alpha}$ are smooth morphism.
Definition \ref{RCHdef}(iv) gives in this particular case the map in $C(\Var(k)^2/S)$
\begin{eqnarray*}
r_S^{CH}(\rho_S^*LF):R^{CH}(\rho_S^*LF):=
(\cdots\to\oplus_{(U_{\alpha},h_{\alpha})\in\Var(k)^{sm}/S}\varinjlim_{(\bar X_{\alpha},\bar D_{\alpha})/S}
R_{(\bar X_{\alpha},\bar D_{\alpha})/S}(\mathbb Z(U_{\alpha}/S)) \\
\xrightarrow{(R_S^{CH}(g^n_{\alpha,\beta}))}
\oplus_{(U_{\alpha},h_{\alpha})\in\Var(k)^{sm}/S}\varinjlim_{(\bar X_{\alpha},\bar D_{\alpha})/S}
R_{(\bar X_{\alpha},\bar D_{\alpha})/S}(\mathbb Z(U_{\alpha}/S))\to\cdots)
\to\mathbb D^{12}_S(\rho_S^*LF),
\end{eqnarray*}
where for $(U_{\alpha},h_{\alpha})\in\Var(k)^{sm}/S$, the inductive limit run over all the compactifications 
$\bar f_{\alpha}:\bar X_{\alpha}\to\bar S$ of $h_{\alpha}:U_{\alpha}\to S$ with $\bar X_{\alpha}\in\PSmVar(k)$
and $\bar D_{\alpha}:=\bar X_{\alpha}\backslash U_{\alpha}$ a normal crossing divisor.
Definition \ref{RCHdef}(iv) gives then by functoriality in particular, for $F=F^{\bullet}\in C(\Var(k)^{sm}/S)$, 
the map in $C(\Var(k)^2/S)$
\begin{eqnarray*}
r_S^{CH}(\rho_S^*LF)=(r_S^{CH}(\rho_S^*LF^*)):R^{CH}(\rho_S^*LF)\to\mathbb D^{12}_S(\rho_S^*LF).
\end{eqnarray*}

\item Let $g:T\to S$ a morphism with $T,S\in\SmVar(k)$. Let $h:U\to S$ a smooth morphism with $U\in\Var(k)$.
Consider the cartesian square
\begin{equation*}
\xymatrix{U_T\ar[r]^{h'}\ar[d]^{g'} & T\ar[d]^g \\
U\ar[r]^h & S}
\end{equation*}
Note that $U$ is smooth since $S$ and $h$ are smooth, and $U_T$ is smooth since $T$ and $h'$ are smooth. 
Take, see definition-proposition \ref{RCHdef0}(ii),a compactification $\bar f_0=\bar h:\bar X_0\to\bar S$ of $h:U\to S$ 
and a compactification $\bar f'_0=\bar{g\circ h'}:\bar X'_0\to\bar S$ of $g\circ h':U'\to S$ such that 
$g':U_T/S\to U/S$ extend to a morphism $\bar{g}'_0:\bar X'_0/\bar S\to \bar X_0/\bar S$. 
Denote $\bar Z=\bar X_0\backslash U$ and $\bar Z'=\bar X'_0\backslash U_T$.
Take, see definition-proposition \ref{RCHdef0}(ii), a strict desingularization 
$\bar\epsilon:(\bar X,\bar D)\to(\bar X_0,\bar Z)$ of $(\bar X_0,\bar Z)$,
a desingularization $\bar\epsilon'_{\bullet}:(\bar X',\bar D')\to(\bar X'_0,\bar Z')$ of $(\bar X'_0,\bar Z')$
and a morphism $\bar g':\bar X'\to\bar X$ such that the following diagram commutes
\begin{equation*}
\xymatrix{\bar X'_0\ar[r]^{\bar{g}'_0} & \bar X_0 \\
\bar X'\ar[u]^{\bar\epsilon'}\ar[r]^{\bar{g}'} & \bar X\ar[u]^{\bar\epsilon}}.
\end{equation*} 
We then have, see definition-proposition \ref{RCHdef0}(ii),
the following commutative diagram in $\Fun(\Delta,\Var(k))$
\begin{equation*}
\xymatrix{U=U_{c(\bullet)}\ar[r]^j & \bar X=\bar X_{c(\bullet)} & \, & \bar D_{s_{g'}(\bullet)}\ar[ll]_{i_{\bullet}} \\
U_T=U_{T,c(\bullet)}\ar[r]^{j'}\ar[u]^{g'} & \bar X'=X'_{c(\bullet)}\ar[u]^{\bar{g}'} & \bar D'_{\bullet}\ar[l]_{i'_{\bullet}} & 
\bar{g}^{'-1}(\bar D_{s_{g'}(\bullet)})\ar[l]_{i''_{g'\bullet}}\ar[u]^{(\bar{g}')'_{\bullet}}:i'_{g\bullet}}
\end{equation*}
We then consider the following map in $C(\Var(k)^2/T)$, see definition \ref{RCHdef}(ii)
\begin{eqnarray*}
T(g,R^{CH})(\mathbb Z(U/S)):g^*R_{(\bar X,\bar D)/S}(\mathbb Z(U/S)) \\ 
\xrightarrow{g^*R^{CH}_S(g')}g^*R_{(\bar X',\bar D')/S}(\mathbb Z(U_T/S))=g^*g_*R_{(\bar X',\bar D')/T}(\mathbb Z(U_T/T)) \\
\xrightarrow{\ad(g^*,g_*)(R_{(\bar X',\bar D')/T}(\mathbb Z(U_T/T)))}R_{(\bar X',\bar D')/T}(\mathbb Z(U_T/T))
\end{eqnarray*}
For  
\begin{eqnarray*}
Q^*:=(\cdots\to\oplus_{\alpha\in\Lambda^n}\mathbb Z(U^n_{\alpha}/S)
\xrightarrow{(\mathbb Z(g^n_{\alpha,\beta}))}\oplus_{\beta\in\Lambda^{n-1}}\mathbb Z(U^{n-1}_{\beta}/S)\to\cdots)
\in C(\Var(k)/S)
\end{eqnarray*}
a complex of (maybe infinite) direct sum of representable presheaves with $h^n_{\alpha}:U^n_{\alpha}\to S$ smooth,
we get the map in $C(\Var(k)^2/T)$
\begin{eqnarray*}
T(g,R^{CH})(Q^*):g^*R^{CH}(Q^*)=
(\cdots\to\oplus_{\alpha\in\Lambda^n}\varinjlim_{(\bar X^n_{\alpha},\bar D^n_{\alpha})/S}
g^*R_{(\bar X^n_{\alpha},\bar D^n_{\alpha})/S}(\mathbb Z(U^n_{\alpha}/S))\to\cdots) \\
\xrightarrow{(T(g,R^{CH})(\mathbb Z(U^n_{\alpha}/S)))} 
(\cdots\to\oplus_{\alpha\in\Lambda^n}\varinjlim_{(\bar X^{n'}_{\alpha},\bar D^{n'}_{\alpha})/T}
R_{(\bar X^{n'}_{\alpha},\bar D^{n'}_{\alpha})/T}(\mathbb Z(U^n_{\alpha,T}/S))\to\cdots)=:R^{CH}(g^*Q^*).
\end{eqnarray*}
Let $F\in\PSh(\Var(k)^{sm}/S)$. Consider 
\begin{eqnarray*}
q:LF:=(\cdots\to\oplus_{(U_{\alpha},h_{\alpha})\in\Var(k)^{sm}/S}\mathbb Z(U_{\alpha}/S)\to\cdots)\to F
\end{eqnarray*}
the canonical projective resolution given in subsection 2.3.3.
We then get in particular the map in $C(\Var(k)^2/T)$
\begin{eqnarray*}
T(g,R^{CH})(\rho_S^*LF):g^*R^{CH}(\rho_S^*LF)= \\
(\cdots\to\oplus_{(U_{\alpha},h_{\alpha})\in\Var(k)^{sm}/S}\varinjlim_{(\bar X_{\alpha},\bar D_{\alpha})/S}
g^*R_{(\bar X_{\alpha},\bar D_{\alpha})/S}(\mathbb Z(U_{\alpha}/S))\to\cdots) 
\xrightarrow{(T(g,R^{CH})(\mathbb Z(U_{\alpha}/S)))} \\ 
(\cdots\to\oplus_{(U_{\alpha},h_{\alpha})\in\Var(k)^{sm}/S}\varinjlim_{(\bar X'_{\alpha},\bar D'_{\alpha})/T}
R_{(\bar X'_{\alpha},\bar D'_{\alpha})/T}(\mathbb Z(U_{\alpha,T}/S))\to\cdots)=:R^{CH}(\rho_T^*g^*LF). 
\end{eqnarray*}
By functoriality, we get in particular for $F=F^{\bullet}\in C(\Var(k)^{sm}/S)$, the map in $C(\Var(k)^2/T)$
\begin{eqnarray*}
T(g,R^{CH})(\rho_S^*LF):g^*R^{CH}(\rho_S^*LF)\to R^{CH}(\rho_T^*g^*LF).
\end{eqnarray*}

\item Let $S_1,S_2\in\SmVar(k)$ and $p:S_1\times S_2\to S_1$ the projection. 
Let $h:U\to S_1$ a smooth morphism with $U\in\Var(k)$. Consider the cartesian square
\begin{equation*}
\xymatrix{U\times S_2\ar[r]^{h\times I}\ar[d]^{p'} & S_1\times S_2\ar[d]^p \\
U\ar[r]^h & S_1}
\end{equation*}
Take, see definition-proposition \ref{RCHdef0}(i),a compactification $\bar f_0=\bar h:\bar X_0\to\bar S_1$ of $h:U\to S_1$.
Then $\bar f_0\times I:\bar X_0\times S_2\to\bar S_1\times S_2$ is a compactification of $h\times I:U\times S_2\to S_1\times S_2$ 
and $p':U\times S_2\to U$ extend to $\bar{p}'_0:=p_{X_0}:\bar X_0\times S_2\to\bar X_0$. Denote $Z=X_0\backslash U$. 
Take see theorem \ref{desVar}(i), a strict desingularization 
$\bar\epsilon:(\bar X,\bar D)\to(\bar X_0,\bar Z)$ of the pair $(\bar X_0,\bar Z)$.  
We then have the following commutative diagram in $\Fun(\Delta,\Var(k))$ whose squares are cartesian
\begin{equation}\label{RCHdiap}
\xymatrix{U=U_{c(\bullet)}\ar[r]^j & \bar X & \bar D_{\bullet}\ar[l]_{i_{\bullet}} \\
U\times S_2=(U\times S_2)_{c(\bullet)}\ar[r]^{j\times I}\ar[u]^g & \bar X\times S_2\ar[u]^{\bar{p}':=p_{\bar X}} & 
\bar D_{\bullet}\times S_2\ar[l]_{i'_{\bullet}}\ar[u]^{\bar{p'}'_{\bullet}}}
\end{equation}
Then the map in $C(\Var(k)^2/S_1\times S_2)$
\begin{eqnarray*}
T(p,R^{CH})(\mathbb Z(U/S_1)):p^*R_{(\bar X,\bar D)/S_1}(\mathbb Z(U/S_1))\xrightarrow{\sim} 
R_{(\bar X\times S_2,\bar D_{\bullet}\times S_2)/S_1\times S_2}(\mathbb Z(U\times S_2/S_1\times S_2))
\end{eqnarray*}
is an isomorphism.
Hence, for $Q^*\in C(\Var(k)/S_1)$ a complex of (maybe infinite) direct sum of representable presheaves of smooth morphism,
the map in $C(\Var(k)^2/S_1\times S_2)$
\begin{eqnarray*}
T(p,R^{CH})(Q^*):p^*R^{CH}(Q^*)\xrightarrow{\sim}R^{CH}(p^*Q^*)
\end{eqnarray*}
is an isomorphism.
In particular, for $F\in C(\Var(k)^{sm}/S_1)$ the map in $C(\Var(k)^2/S_1\times S_2)$
\begin{eqnarray*}
T(p,R^{CH})(\rho_{S_1}^*LF):p^*R^{CH}(\rho_{S_1}^*LF)\xrightarrow{\sim}R^{CH}(\rho_{S_1\times S_2}^*p^*LF)
\end{eqnarray*}
is an isomorphism.

\item Let $h_1:U_1\to S$, $h_2:U_2\to S$ two morphisms with $U_1,U_2,S\in\Var(k)$, $U_1,U_2$ smooth.
Denote by $p_1:U_1\times_SU_2\to U_1$ and $p_2:U_1\times_S U_2\to U_2$ the projections.
Take, see definition-proposition \ref{RCHdef0}(i)),
a compactification $\bar f_{10}=\bar{h}_1:\bar X_{10}\to\bar S$ of $h_1:U_1\to S$
and a compactification $\bar f_{20}=\bar{h}_2:\bar X_{20}\to\bar S$ of $h_2:U_2\to S$. Then, 
\begin{itemize}
\item $\bar f_{10}\times\bar f_{20}:\bar X_{10}\times_{\bar S}\bar X_{20}\to S$ 
is a compactification of $h_1\times h_2:U_1\times_SU_2\to S$.
\item $\bar p_{10}:=p_{X_{10}}:\bar X_{10}\times_{\bar S} \bar X_{20}\to\bar X_{10}$ 
is a compactification of $p_1:U_1\times_SU_2\to U_1$.
\item $\bar p_{20}:=p_{X_{20}}:\bar X_{10}\times_{\bar S}\bar X_{20}\to\bar X_{20}$ 
is a compactification of $p_2:U_1\times_SU_2\to U_2$.
\end{itemize}
Denote $\bar Z_1=\bar X_{10}\backslash U_1$ and $\bar Z_2=\bar X_{20}\backslash U_2$.
Take, see theorem \ref{desVar}(i), a strict desingularization 
$\bar\epsilon_1:(\bar X_1,\bar D)\to(\bar X_{10},Z_1)$ of the pair $(\bar X_{10},\bar Z_1)$
and a  strictdesingularization 
$\bar\epsilon_2:(\bar X_2,\bar E)\to(\bar X_{20},Z_2)$ of the pair $(\bar X_{20},\bar Z_2)$. 
Take then a strict desingularization
\begin{equation*}
\bar\epsilon_{12}:((\bar X_1\times_{\bar S}\bar X_2)^N,\bar F)\to
(\bar X_1\times_{\bar S}\bar X_2,(D\times_{\bar S}\bar X_2)\cup(\bar X_1\times_{\bar S}\bar E)) 
\end{equation*}
of the pair $(\bar X_1\times_{\bar S}\bar X_2,(\bar D\times_{\bar S}\bar X_2)\cup(\bar X_1\times_{\bar S}\bar E))$. 
We have then the following commutative diagram 
\begin{equation*}
\xymatrix{\, & \bar X_1\ar[r]^{\bar f_1} & \bar S \\
\, & \bar X_1\times_{\bar S}\bar X_2\ar[r]^{\bar p_1}\ar[u]^{\bar p_2} &\bar X_2\ar[u]^{\bar f_2} \\
(\bar X_1\times_{\bar S}\bar X_2)^N\ar[ru]^{\bar\epsilon_{12}}\ar[rru]^{(\bar p_1)^N}\ar[ruu]^{(\bar p_2)^N} & \, & \,}
\end{equation*}
and
\begin{itemize}
\item $\bar f_1\times\bar f_2:\bar X_1\times_{\bar S}\bar X_2\to\bar S$ 
is a compactification of $h_1\times h_2:U_1\times_SU_2\to S$.
\item $(\bar p_1)^N:=\bar p_1\circ\epsilon_{12}:(\bar X_1\times_{\bar S}\bar X_2)^N\to\bar X_1$ is a compactification of 
$p_1:U_1\times_SU_2\to U_1$.
\item $(\bar p_2)^N:=\bar p_2\circ\epsilon_{12}:(\bar X_1\times_{\bar S}\bar X_2)^N\to\bar X_2$ is a compactification of 
$p_2:U_1\times_SU_2\to U_2$.
\end{itemize}
We have then the morphism in $C(\Var(k)^2/S)$
\begin{eqnarray*}
T(\otimes,R_S^{CH})(\mathbb Z(U_1/S),\mathbb Z(U_2/S)):=R_S^{CH}(p_1)\otimes R_S^{CH}(p_2): \\
R_{(\bar X_1,\bar D)/S}(\mathbb Z(U_1/S))\otimes R_{(X_2,E))/S}(\mathbb Z(U_2/S)) 
\xrightarrow{\sim} R_{(\bar X_1\times_{\bar S}\bar X_2)^N,\bar F)/S}(\mathbb Z(U_1\times_S U_2/S))
\end{eqnarray*}
For   
\begin{eqnarray*}
Q_1^*:=(\cdots\to\oplus_{\alpha\in\Lambda^n}\mathbb Z(U^n_{1,\alpha}/S)
\xrightarrow{(\mathbb Z(g^n_{\alpha,\beta}))}\oplus_{\beta\in\Lambda^{n-1}}\mathbb Z(U^{n-1}_{1,\beta}/S)\to\cdots), \\
Q_2^*:=(\cdots\to\oplus_{\alpha\in\Lambda^n}\mathbb Z(U^n_{2,\alpha}/S)
\xrightarrow{(\mathbb Z(g^n_{\alpha,\beta}))}\oplus_{\beta\in\Lambda^{n-1}}\mathbb Z(U^{n-1}_{2,\beta}/S)\to\cdots)
\in C(\Var(k)/S)
\end{eqnarray*}
complexes of (maybe infinite) direct sum of representable presheaves with $U^*_{\alpha}$ smooth, 
we get the morphism in $C(\Var(k)^2/S)$
\begin{eqnarray*}
T(\otimes,R_S^{CH})(Q^*_1,Q^*_2): R^{CH}(Q^*_1)\otimes R^{CH}(Q^*_2) 
\xrightarrow{(T(\otimes,R_S^{CH})(\mathbb Z(U_{1,\alpha}^m),\mathbb Z(U_{2,\beta}^n))} R^{CH}(Q^*_1\otimes Q^*_2)).
\end{eqnarray*}
For $F_1,F_2\in C(\Var(k)^{sm}/S)$, we get in particular the morphism in $C(\Var(k)^2/S)$
\begin{eqnarray*}
T(\otimes,R_S^{CH})(\rho_S^*LF_1,\rho_S^*LF_2):R^{CH}(\rho_S^*LF_1)\otimes R^{CH}(\rho_S^*LF_2) 
\to R^{CH}(\rho_S^*(LF_1\otimes LF_2)).
\end{eqnarray*}

\end{itemize}

\begin{defi}\label{sharpstar}
Let $h:U\to S$ a morphism, with $U,S\in\Var(k)$, $U$ irreducible.
Take, see definition-proposition \ref{RCHdef0},
$\bar f_0=\bar{h}_0:\bar X_0\to\bar S$ a compactification of $h:U\to S$ and denote by $\bar Z=\bar X_0\backslash U$.
Take, using theorem \ref{desVar}, a desingularization 
$\bar\epsilon:(\bar X,\bar D)\to(\bar X_0,\Delta)$ of the pair $(\bar X_0,\Delta)$, $\bar Z\subset\Delta$,
with $\bar X\in\PSmVar(k)$ and 
$\bar D:=\bar\epsilon^{-1}(\Delta)=\cup_{i=1}^s\bar D_i\subset\bar X$ a normal crossing divisor.
Denote $d_X:=\dim(\bar X)=\dim(U)$.
\begin{itemize}
\item[(i)] The cycle $(\Delta_{\bar D_{\bullet}}\times S)\subset\bar D_{\bullet}\times\bar D_{\bullet}\times S$
induces by the diagonal $\Delta_{\bar D_{\bullet}}\subset\bar D_{\bullet}\times\bar D_{\bullet}$
gives the morphism in $C(\Var(k)^2/S)$
\begin{eqnarray*}
[\Delta_{\bar D_{\bullet}}]\in\Hom(\mathbb Z^{tr}((\bar D_{\bullet}\times S,D_{\bullet})/S),
p_{S*}E_{et}(\mathbb Z((\bar D_{\bullet}\times S,D_{\bullet})/\bar X\times S)(d_X)[2d_X]))
\xrightarrow{\sim} \\
\Hom(\mathbb Z((\bar D_{\bullet}\times S\times\bar X,D_{\bullet})/\bar X\times S), \\
\mathbb Z^{tr}((\bar D_{\bullet}\times S\times\mathbb P^{d_X},D_{\bullet}\times\mathbb P^{d_X})/\bar X\times S)/
\mathbb Z^{tr}((-)\times\mathbb P^{d_X-1},(-)\times\mathbb P^{d_X-1})) \\
\subset H^0(\mathcal Z_{d_{D_{\bullet}}+d_S}(\square^*\times\bar D_{\bullet}\times\bar D_{\bullet}\times S),
\mbox{s.t.}\alpha_*(\times D_{\bullet})=D_{\bullet})
\end{eqnarray*}
\item[(ii)] The cycle $(\Delta_{\bar X}\times S)\subset\bar X\times\bar X\times S$
induces by the diagonal $\Delta_{\bar X}\subset\bar X\times\bar X$
gives the morphism in $C(\Var(k)^2/S)$
\begin{eqnarray*}
[\Delta_{\bar X}]\in\Hom(\mathbb Z^{tr}((\bar X\times S,X)/S),
p_{S*}E_{et}(\mathbb Z((\bar X\times S,X)/\bar X\times S)(d_X)[2d_X]))
\xrightarrow{\sim} \\
\Hom(\mathbb Z((\bar X\times S\times\bar X,X)/\bar X\times S), \\
\mathbb Z^{tr}((\bar X\times S\times\mathbb P^{d_X},X\times\mathbb P^{d_X})/\bar X\times S)/
\mathbb Z^{tr}((-)\times\mathbb P^{d_X-1},(-)\times\mathbb P^{d_X-1})) \\
\subset H^0(\mathcal Z_{d_X+d_S}(\square^*\times\bar X\times\bar X\times S), \mbox{s.t.}\alpha_*(\times X)=X)
\end{eqnarray*}
\end{itemize}
Let $h:U\to S$ a morphism, with $U,S\in\Var(k)$, $U$ smooth connected (hence irreducible by smoothness).
Take, see definition-proposition \ref{RCHdef0},
$\bar f_0=\bar{h}_0:\bar X_0\to\bar S$ a compactification of $h:U\to S$ and denote by $\bar Z=\bar X_0\backslash U$.
Take, using theorem \ref{desVar}(ii), a strict desingularization 
$\bar\epsilon:(\bar X,\bar D)\to(\bar X_0,\bar Z)$ of the pair $(\bar X_0,\bar Z)$
with $\bar X\in\PSmVar(k)$ and 
$\bar D:=\bar\epsilon^{-1}(\bar Z)=\cup_{i=1}^s\bar D_i\subset\bar X$ a normal crossing divisor.
Denote $d_X:=\dim(\bar X)=\dim(U)$.
\begin{itemize}
\item[(iii)]We get from (i) and (ii) the morphism in $C(\Var(k)^2/S)$
\begin{eqnarray*}
T(p_{S\sharp},p_{S*})(\mathbb Z((\bar D_{\bullet}\times S,D_{\bullet})/\bar X\times S),
\mathbb Z((\bar X\times S,X)/\bar X\times S)):=([\Delta_{\bar D_{\bullet}}],[\Delta_{\bar X}]): \\
\Cone(\mathbb Z(i_{\bullet}\times I):
(\mathbb Z^{tr}((\bar D_{\bullet}\times S,D_{\bullet})/S),u_{IJ})\to\mathbb Z^{tr}((\bar X\times S,X)/S))\to \\
p_{S*}E_{et}(\Cone(\mathbb Z(i_{\bullet}\times I):
(\mathbb Z((\bar D_{\bullet}\times S,D_{\bullet})/\bar X\times S),u_{IJ})\to \\
\mathbb Z((\bar X\times S,X)/\bar X\times S)))(d_X)[2d_X]=:R_{(\bar X,\bar D)/S}(\mathbb Z(U/S))(d_X)[2d_X]
\end{eqnarray*}
\item[(iii)']which gives the map in $C(\Var(k)^{2,smpr}/S)$
\begin{eqnarray*}
T^{\mu,q}(p_{S\sharp},p_{S*})(\mathbb Z((\bar D_{\bullet}\times S,D_{\bullet})/\bar X\times S),
\mathbb Z((\bar X\times S,X)/\bar X\times S)): \\
\Cone(\mathbb Z(i_{\bullet}\times I):
(\mathbb Z^{tr}((\bar D_{\bullet}\times S,D_{\bullet})/S),u_{IJ})\to\mathbb Z^{tr}((\bar X\times S,X)/S))= \\
L\rho_{S*}\mu_{S*}\Cone(\mathbb Z(i_{\bullet}\times I):
(\mathbb Z^{tr}((\bar D_{\bullet}\times S,D_{\bullet})/S),u_{IJ})\to\mathbb Z^{tr}((\bar X\times S,X)/S)) \\
\xrightarrow{L\rho_{S*}\mu_{S*}T(p_{S\sharp},p_{S*})(\mathbb Z((\bar D_{\bullet}\times S,D_{\bullet})/\bar X\times S),
\mathbb Z((\bar X\times S,X)/\bar X\times S))}
L\rho_{S*}\mu_{S*}R_{(\bar X,\bar D)/S}(\mathbb Z(U/S))(d_X)[2d_X]
\end{eqnarray*}
\end{itemize}
\end{defi}

\begin{prop}\label{sharpstarprop}
Let $h:U\to S$ a morphism, with $U,S\in\Var(k)$, $U$ irreducible.
Take, see definition-proposition \ref{RCHdef0},
$\bar f_0=\bar{h}_0:\bar X_0\to\bar S$ a compactification of $h:U\to S$ and denote by $\bar Z=\bar X_0\backslash U$.
Take, using theorem \ref{desVar}(ii), a desingularization 
$\bar\epsilon:(\bar X,\bar D)\to(\bar X_0,\Delta)$ of the pair $(\bar X_0,\Delta)$, $\bar Z\subset\Delta$
with $\bar X\in\PSmVar(k)$ and 
$\bar D:=\bar\epsilon^{-1}(\Delta)=\cup_{i=1}^s\bar D_i\subset\bar X$ a normal crossing divisor.
Denote $d_X:=\dim(\bar X)=\dim(U)$.
\begin{itemize}
\item[(i)] The morphism
\begin{eqnarray*}
[\Delta_{\bar D_{\bullet}}]:\mathbb Z^{tr}((\bar D_{\bullet}\times S,D_{\bullet})/S),
\to p_{S*}E_{et}(\mathbb Z((\bar D_{\bullet}\times S,D_{\bullet})/\bar X\times S)(d_X)[2d_X])
\end{eqnarray*}
given in definition \ref{sharpstar}(i) is an equivalence $(\mathbb A^1,et)$ local.
\item[(ii)] The morphism
\begin{eqnarray*}
[\Delta_{\bar X}]:\mathbb Z^{tr}((\bar X\times S,X)/S),
\to p_{S*}E_{et}(\mathbb Z((\bar X\times S,X)/\bar X\times S)(d_X)[2d_X])
\end{eqnarray*}
given in definition \ref{sharpstar}(ii) is an equivalence $(\mathbb A^1,et)$ local.
\end{itemize}
Let $h:U\to S$ a morphism, with $U,S\in\Var(k)$, $U$ smooth connected (hence irreducible by smoothness).
Take, see definition-proposition \ref{RCHdef0},
$\bar f_0=\bar{h}_0:\bar X_0\to\bar S$ a compactification of $h:U\to S$ and denote by $\bar Z=\bar X_0\backslash U$.
Take, using theorem \ref{desVar}(ii), a strict desingularization 
$\bar\epsilon:(\bar X,\bar D)\to(\bar X_0,\bar Z)$ of the pair $(\bar X_0,\bar Z)$,
with $\bar X\in\PSmVar(k)$ and 
$\bar D:=\bar\epsilon^{-1}(Z)=\cup_{i=1}^s\bar D_i\subset\bar X$ a normal crossing divisor.
\begin{itemize}
\item[(iii)] The morphism
\begin{eqnarray*}
T(p_{S\sharp},p_{S*})(\mathbb Z((\bar D_{\bullet}\times S,D_{\bullet})/\bar X\times S),
\mathbb Z((\bar X\times S,X)/\bar X\times S)):=([\Delta_{\bar D_{\bullet}}],[\Delta_{\bar X}]): \\
\Cone(\mathbb Z(i_{\bullet}\times I):
(\mathbb Z^{tr}((\bar D_{\bullet}\times S,D_{\bullet})/S),u_{IJ})\to\mathbb Z^{tr}((\bar X\times S,X)/S))\to \\
p_{S*}E_{et}(\Cone(\mathbb Z(i_{\bullet}\times I):
(\mathbb Z((\bar D_{\bullet}\times S,D_{\bullet})/\bar X\times S),u_{IJ})\to \\
\mathbb Z((\bar X\times S,X)/\bar X\times S)))(d_X)[2d_X]=:R_{(\bar X,\bar D)/S}(\mathbb Z(U/S))(d_X)[2d_X]
\end{eqnarray*}
given in definition \ref{sharpstar}(iii)' is an equivalence $(\mathbb A^1,et)$ local.
\item[(iii)']The morphism
\begin{eqnarray*}
T^{\mu,q}(p_{S\sharp},p_{S*})(\mathbb Z((\bar D_{\bullet}\times S,D_{\bullet})/\bar X\times S),
\mathbb Z((\bar X\times S,X)/\bar X\times S)): \\
\Cone(\mathbb Z(i_{\bullet}\times I):
(\mathbb Z^{tr}((\bar D_{\bullet}\times S,D_{\bullet})/S),u_{IJ})\to\mathbb Z^{tr}((\bar X\times S,X)/S)) \\
\to L\rho_{S*}\mu_{S*}R_{(\bar X,\bar D)/S}(\mathbb Z(U/S))(d_X)[2d_X]
\end{eqnarray*}
given in definition \ref{sharpstar}(iii)' is an equivalence $(\mathbb A^1,et)$ local.
\end{itemize}
\end{prop}

\begin{proof}
\noindent(i):
By Yoneda lemma, it is equivalent to show that for every morphism $g:T\to S$
with $T\in\Var(k)$ and every closed subset $E\subset T$, the composition morphism
\begin{eqnarray*}
[\Delta_{\bar D_{\bullet}}]:
\Hom^{\bullet}(\mathbb Z((T,E)/S),C_*\mathbb Z^{tr}((\bar D_{\bullet}\times S,D_{\bullet})/S))
\xrightarrow{\Hom^{\bullet}(\mathbb Z((T,E)/S),C_*\Delta_{\bar D_{\bullet}})} \\
\Hom^{\bullet}(\mathbb Z((T,E)/S),p_{S*}E_{et}(\mathbb Z((\bar D_{\bullet}\times S,D_{\bullet})/\bar X\times S)(d_X)[2d_X]))
\end{eqnarray*}
is a quasi-isomorphism of abelian groups.
But this map is the composite
\begin{eqnarray*}
\Hom^{\bullet}(\mathbb Z((T,E)/S),\mathbb Z^{tr}((\bar D_{\bullet}\times S,D_{\bullet})/S))
\xrightarrow{[\Delta_{\bar D_{\bullet}}]} \\ 
\Hom^{\bullet}(\mathbb Z((T,E)/S),p_{S*}E_{et}(\mathbb Z((\bar D_{\bullet}\times S,D_{\bullet})/\bar X\times S)(d_X)[2d_X]))
\xrightarrow{\sim} \\
\Hom^{\bullet}(\mathbb Z((T\times\bar X,E)/S\times\bar X), \\
C_*\mathbb Z^{tr}((\bar D_{\bullet}\times S\times\mathbb P^{d_X},D_{\bullet}\times\mathbb P^{d_X})/\bar X\times S)/
\mathbb Z^{tr}((-)\times\mathbb P^{d_X-1},(-)\times\mathbb P^{d_X-1}))
\end{eqnarray*}
which is clearly a quasi-isomorphism.

\noindent(ii): Similar to (i).

\noindent(iii):Follows from (i) and (ii).

\noindent(iii)':Follows from (iii) and the fact that 
$\mu_{S*}$ preserve $(\mathbb A^1,et)$ local equivalence (see proposition \ref{mu12}) and the fact that
$\rho_{S*}$ preserve $(\mathbb A^1,et)$ local equivalence (see proposition \ref{rho12}).
\end{proof}

\begin{defi}\label{tus}
\begin{itemize}
\item[(i)]Let $h:U\to S$ a morphism, with $U,S\in\Var(k)$, $U$ smooth.
Take, see definition-proposition \ref{RCHdef0},
$\bar f_0=\bar{h}_0:\bar X_0\to\bar S$ a compactification of $h:U\to S$ and denote by $\bar Z=\bar X_0\backslash U$.
Take, using theorem \ref{desVar}(ii), a strict desingularization 
$\bar\epsilon:(\bar X,\bar D)\to(\bar X_0,\bar Z)$ of the pair $(\bar X_0,\bar Z)$,
with $\bar X\in\PSmVar(k)$ and 
$\bar D:=\bar\epsilon^{-1}(\bar Z)=\cup_{i=1}^s\bar D_i\subset\bar X$ a normal crossing divisor.
We will consider the following canonical map in $C(\Var(k)^{sm}/S)$
\begin{eqnarray*}
T_{(\bar X,\bar D)/S}(U/S):\Gr_{S*}^{12}L\rho_{S*}\mu_{S*}R_{(\bar X,\bar D)/S}(\mathbb Z(U/S)) 
\xrightarrow{q}
\Gr_{S*}^{12}\rho_{S*}\mu_{S*}R_{(\bar X,\bar D)/S}(\mathbb Z(U/S)) \\
\xrightarrow{r_{(\bar X,\bar D)/S}(\mathbb Z(U/S))} 
\Gr_{S*}^{12}\rho_{S*}\mu_{S*}p_{S*}E_{et}(\mathbb Z((U\times S,U)/U\times S)) 
\xrightarrow{l(U/S)}h_*E_{et}(\mathbb Z(U/U))=:\mathbb D_S^0(\mathbb Z(U/S))
\end{eqnarray*}
where, for $h':V\to S$ a smooth morphism with $V\in\Var(k)$,
\begin{eqnarray*}
l^{00}(U/S)(V/S):\mathbb Z((U\times S,U)/U\times S)(V\times U\times S,V\times_SU/U\times S)
\to\mathbb Z(U/U)(V\times_S U), \alpha\mapsto\alpha_{|V\times_SU}  
\end{eqnarray*}
which gives
\begin{eqnarray*}
l^0(U/S)(V/S):E^0_{et}(\mathbb Z((U\times S,U)/U\times S))(V\times U\times S,V\times_SU/U\times S)
\to E^0_{et}(\mathbb Z(U/U))(V\times_S U), 
\end{eqnarray*}
and by induction 
\begin{equation*}
\tau^{\leq i}l(U/S):\Gr_{S*}^{12}\rho_{S*}\mu_{S*}p_{S*}E^{\leq i}_{et}(\mathbb Z((U\times S,U)/U\times S)) 
\to h_*E^{\leq i}_{et}(\mathbb Z(U/U))
\end{equation*}
where $\tau^{\leq i}$ is the cohomological truncation.
\item[(ii)]Let $g:U'/S\to U/S$ a morphism, with $U'/S=(U',h'),U/S=(U,h)\in\Var(k)/S$, $U$,$U'$ smooth.
Take, see definition-proposition \ref{RCHdef0}(ii),a compactification $\bar f_0=\bar h:\bar X_0\to\bar S$ of $h:U\to S$ 
and a compactification $\bar f'_0=\bar{h}':\bar X'_0\to S$ of $h':U'\to S$ such that 
$g:U'/S\to U/S$ extend to a morphism $\bar g_0:\bar X'_0/\bar S\to\bar X_0/\bar S$. 
Denote $\bar Z=\bar X_0\backslash U$ and $\bar Z'=\bar X'_0\backslash U'$.
Take, see definition-proposition \ref{RCHdef0}(ii), a strict desingularization 
$\bar\epsilon:(\bar X,\bar D)\to(\bar X_0,\bar Z)$ of $(\bar X_0,\bar Z)$,
a strict desingularization $\bar\epsilon'_{\bullet}:(\bar X',\bar D')\to(\bar X'_0,\bar Z')$ of $(\bar X'_0,\bar Z')$
and a morphism $\bar g:\bar X'\to\bar X$ such that the following diagram commutes
\begin{equation*}
\xymatrix{\bar X'_0\ar[r]^{\bar{g}_0} & \bar X_0 \\
\bar X'\ar[u]^{\bar\epsilon'}\ar[r]^{\bar g} & \bar X\ar[u]^{\bar\epsilon}}.
\end{equation*}
Then by the diagram given in definition \ref{RCHdef}(ii), 
the following diagram in $C(\Var(k)^{sm}/S)$ obviously commutes
\begin{equation*}
\xymatrix{\Gr_{S*}^{12}L\rho_{S*}\mu_{S*}R_{(\bar X,\bar D)/S}(\mathbb Z(U/S))
\ar[rr]^{T_{(\bar X,\bar D)/S}(U/S)}\ar[d]_{R_S^{CH}(g)} & \, & 
h_*E_{et}(\mathbb Z(U/U)):=\mathbb D_S^0(\mathbb Z(U/S))
\ar[d]^{T(g,E)(-)\circ\ad(g^*,g_*)(E_{et}(\mathbb Z(U/U))):=\mathbb D_S^0(g)} \\
\Gr_{S*}^{12}L\rho_{S*}\mu_{S*}R_{(\bar X',\bar D')/S}(\mathbb Z(U'/S))
\ar[rr]^{T_{(\bar X',\bar D')/S}(U'/S)} & \, & 
h'_*E_{et}(\mathbb Z(U'/U')):=\mathbb D_S^0(\mathbb Z(U'/S))}
\end{equation*}
where $l(U/S)$ are $l(U'/S)$ are the maps given in (i).
\item[(iii)]Let $S\in\SmVar(k)$.
Let $F\in C(\Var(k)^{sm}/S)$. We get from (i) and (ii) morphisms in $C(\Var(k)^{sm}/S)$
\begin{eqnarray*}
T_S^{CH}(LF):\Gr_{S*}^{12}L\rho_{S*}\mu_{S*}R_{(\bar X^*,\bar D^*)/S}(\rho_S^*LF) \\
\xrightarrow{r_S^{CH}(LF)}\Gr_{S*}^{12}L\rho_{S*}\mu_{S*}\mathbb D^{12}_S(\rho_S^*LF)
\xrightarrow{l(L(F)}\mathbb D_S^0(L(F))
\end{eqnarray*}
\end{itemize}
\end{defi}

\begin{lem}\label{grRCH}
\begin{itemize}
\item[(i)]Let $h:U\to S$ a morphism, with $U,S\in\Var(k)$, $U$ smooth.
Take, see definition-proposition \ref{RCHdef0},
$\bar f_0=\bar{h}_0:\bar X_0\to\bar S$ a compactification of $h:U\to S$ and denote by $\bar Z=\bar X_0\backslash U$.
Take, using theorem \ref{desVar}(ii), a  strict desingularization 
$\bar\epsilon:(\bar X,\bar D)\to(\bar X_0,\bar Z)$ of the pair $(\bar X_0,\bar Z)$,
with $\bar X\in\PSmVar(k)$ and 
$\bar D:=\bar\epsilon^{-1}(\bar Z)=\cup_{i=1}^s\bar D_i\subset\bar X$ a normal crossing divisor.
Then the map in $C(\Var(k)^{2,smpr}/S)$
\begin{eqnarray*}
\ad(\Gr_S^{12*},\Gr_{S*}^{12})(L\rho_{S*}\mu_{S*}R_{(\bar X,\bar D)/S}(\mathbb Z(U/S)))\circ q: \\
\Gr_S^{12*}L\Gr_{S*}^{12}L\rho_{S*}\mu_{S*}R_{(\bar X,\bar D)/S}(\mathbb Z(U/S))\to 
L\rho_{S*}\mu_{S*}R_{(\bar X,\bar D)/S}(\mathbb Z(U/S))
\end{eqnarray*}
is an equivalence $(\mathbb A^1,et)$ local.
\item[(ii)]Let $S\in\SmVar(k)$.
Let $F\in C(\Var(k)^{sm}/S)$. Then the map in $C(\Var(k)^{2,smpr}/S)$
\begin{eqnarray*}
\ad(\Gr_S^{12*},\Gr_{S*}^{12})(L\rho_{S*}\mu_{S*}R_{(\bar X^*,\bar D^*)/S}(\rho_S^*LF))\circ q: \\
\Gr_S^{12*}L\Gr_{S*}^{12}L\rho_{S*}\mu_{S*}R_{(\bar X^*,\bar D^*)/S}(\rho_S^*LF)\to 
L\rho_{S*}\mu_{S*}R_{(\bar X^*,\bar D^*)/S}(\rho_S^*LF)
\end{eqnarray*}
is an equivalence $(\mathbb A^1,et)$ local.
\end{itemize}
\end{lem}

\begin{proof}
\noindent(i): Follows from proposition \ref{sharpstarprop}.

\noindent(ii): Follows from (i).
\end{proof}

\begin{defi}\label{RCHhatdef} 
\begin{itemize}
\item[(i)]Let $h:U\to S$ a morphism, with $U,S\in\Var(k)$ and $U$ smooth.
Take, see definition-proposition \ref{RCHdef0},
$\bar f_0=\bar{h}_0:\bar X_0\to\bar S$ a compactification of $h:U\to S$ and denote by $\bar Z=\bar X_0\backslash U$.
Take, using theorem \ref{desVar}(ii), a strict desingularization 
$\bar\epsilon:(\bar X,\bar D)\to(\bar X_0,\bar Z)$ of the pair $(\bar X_0,\bar Z)$, with $\bar X\in\PSmVar(k)$ and 
$\bar D:=\epsilon^{-1}(\bar Z)=\cup_{i=1}^s\bar D_i\subset\bar X$ a normal crossing divisor.  
We denote by $i_{\bullet}:\bar D_{\bullet}\hookrightarrow\bar X=\bar X_{c(\bullet)}$ the morphism of simplicial varieties
given by the closed embeddings $i_I:\bar D_I=\cap_{i\in I}\bar D_i\hookrightarrow\bar X$
We denote by $j:U\hookrightarrow\bar X$ the open embedding and by $p_S:\bar X\times S\to S$ 
and $p_S:U\times S\to S$ the projections.
Considering the graph factorization $\bar f:\bar X\xrightarrow{\bar l}\bar X\times\bar S\xrightarrow{p_{\bar S}}\bar S$
of $\bar f:\bar X\to\bar S$, where $\bar l$ is the graph embedding and $p_{\bar S}$ the projection,
we get closed embeddings $l:=\bar l\times_{\bar S}S:X\hookrightarrow\bar X\times S$ and 
$l_{D_I}:=\bar D_I\times_{\bar X} l:D_I\hookrightarrow\bar D_I\times S$.
We then consider the map in $C(\Var(k)^{2,smpr}/S)$
\begin{eqnarray*}
T(\hat R^{CH},R^{CH})(\mathbb Z(U/S)):\hat R_{(\bar X,\bar D)/S}(\mathbb Z(U/S)) \\
\xrightarrow{:=}\Cone(\mathbb Z(i_{\bullet}\times I):
(\mathbb Z^{tr}((\bar D_{\bullet}\times S,D_{\bullet})/S),u_{IJ})\to\mathbb Z^{tr}((\bar X\times S,X)/S))(-d_X)[-2d_X] \\
\xrightarrow{T^{\mu,q}(p_{S\sharp},p_{S*})(\mathbb Z((\bar D_{\bullet}\times S,D_{\bullet})/\bar X\times S),
\mathbb Z((\bar X\times S,X)/\bar X\times S))(-d_X)[-2d_X]} \\
L\rho_{S*}\mu_{S*}R_{(\bar X,\bar D)/S}(\mathbb Z(U/S)).
\end{eqnarray*}
given in definition \ref{sharpstar}(iii).
\item[(ii)]Let $g:U'/S\to U/S$ a morphism, with $U'/S=(U',h'),U/S=(U,h)\in\Var(k)/S$, with $U$ and $U'$ smooth.
Take, see definition-proposition \ref{RCHdef0}(ii),a compactification $\bar f_0=\bar h:\bar X_0\to\bar S$ of $h:U\to S$ 
and a compactification $\bar f'_0=\bar{h}':\bar X'_0\to\bar S$ of $h':U'\to S$ such that 
$g:U'/S\to U/S$ extend to a morphism $\bar g_0:\bar X'_0/\bar S\to\bar X_0/\bar S$. 
Denote $\bar Z=\bar X_0\backslash U$ and $\bar Z'=\bar X'_0\backslash U'$.
Take, see definition-proposition \ref{RCHdef0}(ii), a strict desingularization 
$\bar\epsilon:(\bar X,\bar D)\to(\bar X_0,\bar Z)$ of $(\bar X_0,\bar Z)$,
a strict desingularization $\bar\epsilon'_{\bullet}:(\bar X',\bar D')\to(\bar X'_0,\bar Z')$ of $(\bar X'_0,\bar Z')$
and a morphism $\bar g:\bar X'\to\bar X$ such that the following diagram commutes
\begin{equation*}
\xymatrix{\bar X'_0\ar[r]^{\bar{g}_0} & \bar X_0 \\
\bar X'\ar[u]^{\bar\epsilon'}\ar[r]^{\bar g} & \bar X\ar[u]^{\bar\epsilon}}.
\end{equation*} 
We then have, see definition-proposition \ref{RCHdef0}(ii),
the diagram (\ref{RCHdia}) in $\Fun(\Delta,\Var(k))$
\begin{equation*}
\xymatrix{U=U_{c(\bullet)}\ar[r]^j & \bar X=\bar X_{c(\bullet)} & \, & \bar D_{s_g(\bullet)}\ar[ll]_{i_{\bullet}} \\
U'=U'_{c(\bullet)}\ar[r]^{j'}\ar[u]^g & \bar X'=\bar X'_{c(\bullet)}\ar[u]^{\bar{g}} & \bar D'_{\bullet}\ar[l]_{i'_{\bullet}} & 
\bar{g}^{-1}(\bar D_{s_g(\bullet)})\ar[l]_{i''_{g\bullet}}\ar[u]^{\bar{g}'_{\bullet}}:i'_{g\bullet}}
\end{equation*}
Consider 
\begin{eqnarray*}
[\Gamma_{\bar g}]^t\in\Hom(\mathbb Z^{tr}((\bar X\times S,X)/S)(-d_X)[-2d_X],
\mathbb Z^{tr}((\bar X'\times S,X')/S)(-d_{X'})[-2d_{X'}]) \\ 
\xrightarrow{\sim}\Hom(\mathbb Z^{tr}((\bar X\times\mathbb A^{d_{X'}}\times S,X\times\mathbb A^{d_{X'}})/S), \\
\mathbb Z_{tr}((\bar X'\times\mathbb P^{d_X}\times S,X'\times\mathbb P^{d_X})/S)/
\mathbb Z_{tr}((-)\times\mathbb P^{d_X-1},(-)\times\mathbb P^{d_X-1}))
\end{eqnarray*}
the morphism given by the transpose of the graph $\Gamma_g\subset X'\times_S X$ of $\bar g:\bar X'\to\bar X$. Then,  
since $i_{\bullet}\circ\bar g'_{\bullet}=\bar g\circ i''_{g\bullet}=\bar g\circ i'\circ\circ i'_{g\bullet}$, 
we have the factorization 
\begin{eqnarray*}
[\Gamma_g]^t\circ\mathbb Z(i_{\bullet}\times I): 
(\mathbb Z^{tr}((\bar D_{s_g(\bullet)}\times S,D_{s_g(\bullet)})/S),u_{IJ})(-d_X)[-2d_X] \\
\xrightarrow{[\Gamma_{\bar g'_{\bullet}}]^t}
(\mathbb Z^{tr}((\bar{g}^{-1}(\bar D_{s_g(\bullet)})\times S,\bar{g}^{-1}(D_{s_g(\bullet)}))/S),u_{IJ})(-d_{X'})[-2d_{X'}]) \\
\xrightarrow{\mathbb Z(i'_{g\bullet}\times I)} 
\mathbb Z^{tr}((\bar X'\times S,X')/S)(-d_{X'})[-2d_{X'}].
\end{eqnarray*}
with
\begin{eqnarray*}
[\Gamma_{\bar g'_{\bullet}}]^t\in
\Hom((\mathbb Z^{tr}((\bar D_{s_g(\bullet)}\times\mathbb A^{d_{X'}}\times S,
D_{s_g(\bullet)}\times\mathbb A^{d_{X'}})/S),u_{IJ}), \\
(\mathbb Z_{tr}((\bar{g}^{-1}(\bar D_{s_g(\bullet)})\times\mathbb P^{d_X}\times S,
\bar{g}^{-1}(D_{s_g(\bullet)})\times\mathbb P^{d_X})/S),u_{IJ})/
\mathbb Z^{tr}((-)\times\mathbb P^{d_X-1},(-)\times\mathbb P^{d_X-1})).
\end{eqnarray*}
We then consider the following map in $C(\Var(k)^{2,pr}/S)$
\begin{eqnarray*}
\hat R_S^{CH}(g):\hat R_{(\bar X,\bar D)/S}(\mathbb Z(U/S))\xrightarrow{:=} \\
\Cone(\mathbb Z(i_{\bullet}\times I):
(\mathbb Z^{tr}((\bar D_{s_g(\bullet)}\times S,D_{s_g(\bullet)})/S),u_{IJ})\to
\mathbb Z^{tr}((\bar X\times S,X)/S))(-d_X)[-2d_X] \\
\xrightarrow{([\Gamma_{\bar g'_{\bullet}}]^t,[\Gamma_{\bar g}]^t)} \\
\Cone(\mathbb Z(i'_{g\bullet}\times I): \\
(\mathbb Z^{tr}((\bar{g}^{-1}(\bar D_{s_g(\bullet)})\times S,\bar g^{-1}(D_{s_g(\bullet)})/S),u_{IJ})
\to\mathbb Z^{tr}((\bar X'\times S,X')/S)))(-d_{X'})[-2d_{X'}] \\
\xrightarrow{(\mathbb Z(i''_{g\bullet}\times I),I)(-d_{X'})[-2d_{X'}]} \\
\Cone(\mathbb Z(i'_{\bullet}\times I):
((\mathbb Z^{tr}((\bar D'_{\bullet}\times S,D'_{\bullet)})/S),u_{IJ})\to
\mathbb Z^{tr}((\bar X'\times S,X')/S)))(-d_{X'})[-2d_{X'}] \\
\xrightarrow{=:}\hat R_{(\bar X',\bar D')/S}(\mathbb Z(U'/S))
\end{eqnarray*}
Then the following diagram in $C(\Var(k)^{2,smpr}/S)$ commutes by definition
\begin{equation*}
\xymatrix{\hat R_{(\bar X,\bar D)/S}(\mathbb Z(U/S))
\ar[rr]^{T(\hat R^{CH},R^{CH})(\mathbb Z(U/S))}\ar[d]_{\hat R_S^{CH}(g)} & \, & 
L\rho_{S*}\mu_{S*}R_{(\bar X,\bar D)/S}(\mathbb Z(U/S))\ar[d]^{L\rho_{S*}\mu_{S*}R_S^{CH}(g)} \\
\hat R_{(\bar X',\bar D')/S}(\mathbb Z(U'/S))\ar[rr]^{T(\hat R^{CH},R^{CH})(\mathbb Z(U'/S))} & \, & 
L\rho_{S*}\mu_{S*}R_{(\bar X',\bar D')/S}(\mathbb Z(U'/S))}.
\end{equation*}
\item[(iii)] For $g_1:U''/S\to U'/S$, $g_2:U'/S\to U/S$ two morphisms
with $U''/S=(U',h''),U'/S=(U',h'),U/S=(U,h)\in\Var(k)/S$, with $U$, $U'$ and $U''$ smooth.
We get from (i) and (ii) 
a compactification $\bar f=\bar{h}:\bar X\to\bar S$ of $h:U\to S$,
a compactification $\bar f'=\bar{h}':\bar X'\to\bar S$ of $h':U'\to S$,
and a compactification $\bar f''=\bar{h}'':\bar X''\to\bar S$ of $h'':U''\to S$,
with $\bar X,\bar X',\bar X''\in\PSmVar(k)$, 
$\bar D:=\bar X\backslash U\subset\bar X$ $\bar D':=\bar X'\backslash U'\subset\bar X'$, 
and $\bar D'':=\bar X''\backslash U''\subset\bar X''$ normal crossing divisors, 
such that $g_1:U''/S\to U'/S$ extend to $\bar g_1:\bar X''/\bar S\to\bar X'/\bar S$,
$g_2:U'/S\to U/S$ extend to $\bar g_2:\bar X'/\bar S\to\bar X/\bar S$, and
\begin{eqnarray*}
\hat R_S^{CH}(g_2\circ g_1)=\hat R_S^{CH}(g_1)\circ \hat R_S^{CH}(g_2):
\hat R_{(\bar X,\bar D)/S}\to \hat R_{(\bar X'',\bar D'')/S} 
\end{eqnarray*}
\item[(iv)] For  
\begin{eqnarray*}
Q^*:=(\cdots\to\oplus_{\alpha\in\Lambda^n}\mathbb Z(U^n_{\alpha}/S)
\xrightarrow{(\mathbb Z(g^n_{\alpha,\beta}))}\oplus_{\beta\in\Lambda^{n-1}}\mathbb Z(U^{n-1}_{\beta}/S)\to\cdots)
\in C(\Var(k)/S)
\end{eqnarray*}
a complex of (maybe infinite) direct sum of representable presheaves with $U^*_{\alpha}$ smooth,
we get from (i),(ii) and (iii) the map in $C(\Var(k)^{2,smpr}/S)$
\begin{eqnarray*}
T(\hat R^{CH},R^{CH})(Q^*):\hat R^{CH}(Q^*):=
(\cdots\to\oplus_{\beta\in\Lambda^{n-1}}\varinjlim_{(\bar X^{n-1}_{\beta},\bar D^{n-1}_{\beta})/S}
\hat R_{(\bar X^{n-1}_{\beta},\bar D^{n-1}_{\beta})/S}(\mathbb Z(U^{n-1}_{\beta}/S)) \\
\xrightarrow{(\hat R_S^{CH}(g^n_{\alpha,\beta}))}\oplus_{\alpha\in\Lambda^n}\varinjlim_{(\bar X^n_{\alpha},\bar D^n_{\alpha})/S}
\hat R_{(\bar X^n_{\alpha},\bar D^n_{\alpha})/S}(\mathbb Z(U^n_{\alpha}/S))\to\cdots) 
\to L\rho_{S*}\mu_{S*}R^{CH}(Q^*),
\end{eqnarray*}
where for $(U^n_{\alpha},h^n_{\alpha})\in\Var(k)/S$, the inductive limit run over all the compactifications 
$\bar f_{\alpha}:\bar X_{\alpha}\to\bar S$ of $h_{\alpha}:U_{\alpha}\to S$ with $\bar X_{\alpha}\in\PSmVar(k)$
and $\bar D_{\alpha}:=\bar X_{\alpha}\backslash U_{\alpha}$ a normal crossing divisor.
For $m=(m^*):Q_1^*\to Q_2^*$ a morphism with 
\begin{eqnarray*}
Q_1^*:=(\cdots\to\oplus_{\alpha\in\Lambda^n}\mathbb Z(U^n_{1,\alpha}/S)
\xrightarrow{(\mathbb Z(g^n_{\alpha,\beta}))}\oplus_{\beta\in\Lambda^{n-1}}\mathbb Z(U^{n-1}_{1,\beta}/S)\to\cdots), \\
Q_2^*:=(\cdots\to\oplus_{\alpha\in\Lambda^n}\mathbb Z(U^n_{2,\alpha}/S)
\xrightarrow{(\mathbb Z(g^n_{\alpha,\beta}))}\oplus_{\beta\in\Lambda^{n-1}}\mathbb Z(U^{n-1}_{2,\beta}/S)\to\cdots)
\in C(\Var(k)/S)
\end{eqnarray*}
complexes of (maybe infinite) direct sum of representable presheaves with $U^*_{1,\alpha}$ and $U^*_{2,\alpha}$ smooth,
we get again from (i),(ii) and (iii) a commutative diagram in $C(\Var(k)^{2,smpr}/S)$
\begin{equation*}
\xymatrix{\hat R^{CH}(Q_2^*)\ar[rr]^{T(\hat R_S^{CH},R_S^{CH})(Q_2^*)}\ar[d]_{\hat R_S^{CH}(m):=(\hat R_S^{CH}(m^*))} & \, & 
L\rho_{S*}\mu_{S*}R^{CH}(Q_2^*)\ar[d]^{L\rho_{S*}\mu_{S*}R_S^{CH}(m):=L\rho_{S*}\mu_{S*}(R_S^{CH}(m^*))} \\
\hat R^{CH}(Q_1^*)\ar[rr]^{T(\hat R_S^{CH},R_S^{CH})(Q_1^*)} & \, & 
L\rho_{S*}\mu_{S*}R^{CH}(Q_1^*)}.
\end{equation*}
\end{itemize}
\end{defi}

\begin{itemize}
\item Let $S\in\Var(k)$
For $(h,m,m')=(h^*,m^*,m^{'*}):Q_1^*[1]\to Q_2^*$ an homotopy with $Q_1^*,Q_2^*\in C(\Var(k)/S)$
complexes of (maybe infinite) direct sum of representable presheaves with $U^*_{1,\alpha}$ and $U^*_{2,\alpha}$ smooth,
\begin{equation*}
(\hat R_S^{CH}(h),\hat R_S^{CH}(m),\hat R_S^{CH}(m'))=(\hat R_S^{CH}(h^*),\hat R_S^{CH}(m^*),\hat R_S^{CH}(m^{'*})):
R^{CH}(Q_2^*)[1]\to R^{CH}(Q_1^*)
\end{equation*}
is an homotopy in $C(\Var(k)^{2,smpr}/S)$ using definition \ref{RCHhatdef} (iii). 
In particular if $m:Q_1^*\to Q_2^*$ with $Q_1^*,Q_2^*\in C(\Var(k)/S)$
complexes of (maybe infinite) direct sum of representable presheaves with $U^*_{1,\alpha}$ and $U^*_{2,\alpha}$ smooth
is an homotopy equivalence, then $\hat R_S^{CH}(m):\hat R^{CH}(Q_2^*)\to\hat R^{CH}(Q_1^*)$ is an homotopy equivalence.
\item Let $S\in\SmVar(k)$. Let $F\in\PSh(\Var(k)^{sm}/S)$. Consider 
\begin{eqnarray*}
q:LF:=(\cdots\to\oplus_{(U_{\alpha},h_{\alpha})\in\Var(k)^{sm}/S}\mathbb Z(U_{\alpha}/S)
\xrightarrow{(\mathbb Z(g^n_{\alpha,\beta}))}
\oplus_{(U_{\alpha},h_{\alpha})\in\Var(k)^{sm}/S}\mathbb Z(U_{\alpha}/S)\to\cdots)\to F
\end{eqnarray*}
the canonical projective resolution given in subsection 2.3.3.
Note that the $U_{\alpha}$ are smooth since $S$ is smooth and $h_{\alpha}$ are smooth morphism.
Definition \ref{RCHhatdef}(iv) gives in this particular case the map in $C(\Var(k)^2/S)$
\begin{eqnarray*}
T(\hat R_S^{CH},R_S^{CH})(\rho_S^*LF):\hat R^{CH}(\rho_S^*LF):=
(\cdots\to\oplus_{(U_{\alpha},h_{\alpha})\in\Var(k)^{sm}/S}\varinjlim_{(\bar X_{\alpha},\bar D_{\alpha})/S}
\hat R_{(\bar X_{\alpha},\bar D_{\alpha})/S}(\mathbb Z(U_{\alpha}/S)) \\
\xrightarrow{(\hat R_S^{CH}(g^n_{\alpha,\beta}))}
\oplus_{(U_{\alpha},h_{\alpha})\in\Var(k)^{sm}/S}\varinjlim_{(\bar X_{\alpha},\bar D_{\alpha})/S}
\hat R_{(\bar X_{\alpha},\bar D_{\alpha})/S}(\mathbb Z(U_{\alpha}/S))\to\cdots)
\to L\rho_{S*}\mu_{S*}R^{CH}(\rho_S^*LF),
\end{eqnarray*}
where for $(U_{\alpha},h_{\alpha})\in\Var(k)^{sm}/S$, the inductive limit run over all the compactifications 
$\bar f_{\alpha}:\bar X_{\alpha}\to\bar S$ of $h_{\alpha}:U_{\alpha}\to S$ with $\bar X_{\alpha}\in\PSmVar(k)$
and $\bar D_{\alpha}:=\bar X_{\alpha}\backslash U_{\alpha}$ a normal crossing divisor.
Definition \ref{RCHhatdef}(iv) gives then by functoriality in particular, for $F=F^{\bullet}\in C(\Var(k)^{sm}/S)$, 
the map in $C(\Var(k)^{2,smpr}/S)$
\begin{eqnarray*}
T(\hat R_S^{CH},R_S^{CH})(\rho_S^*LF):\hat R^{CH}(\rho_S^*LF)\to L\rho_{S*}\mu_{S*}R^{CH}(\rho_S^*LF).
\end{eqnarray*}

\item Let $g:T\to S$ a morphism with $T,S\in\SmVar(k)$. Let $h:U\to S$ a smooth morphism with $U\in\Var(k)$.
Consider the cartesian square
\begin{equation*}
\xymatrix{U_T\ar[r]^{h'}\ar[d]^{g'} & T\ar[d]^g \\
U\ar[r]^h & S}
\end{equation*}
Note that $U$ is smooth since $S$ and $h$ are smooth, and $U_T$ is smooth since $T$ and $h'$ are smooth. 
Take, see definition-proposition \ref{RCHdef0}(ii),a compactification $\bar f_0=\bar h:\bar X_0\to\bar S$ of $h:U\to S$ 
and a compactification $\bar f'_0=\bar{g\circ h'}:\bar X'_0\to\bar S$ of $g\circ h':U'\to S$ such that 
$g':U_T/S\to U/S$ extend to a morphism $\bar{g}'_0:\bar X'_0/\bar S\to \bar X_0/\bar S$. 
Denote $\bar Z=\bar X_0\backslash U$ and $\bar Z'=\bar X'_0\backslash U_T$.
Take, see definition-proposition \ref{RCHdef0}(ii), a strict desingularization 
$\bar\epsilon:(\bar X,\bar D)\to(\bar X_0,\bar Z)$ of $(\bar X_0,\bar Z)$,
a desingularization $\bar\epsilon'_{\bullet}:(\bar X',\bar D')\to(\bar X'_0,\bar Z')$ of $(\bar X'_0,\bar Z')$
and a morphism $\bar g':\bar X'\to\bar X$ such that the following diagram commutes
\begin{equation*}
\xymatrix{\bar X'_0\ar[r]^{\bar{g}'_0} & \bar X_0 \\
\bar X'\ar[u]^{\bar\epsilon'}\ar[r]^{\bar{g}'} & \bar X\ar[u]^{\bar\epsilon}}.
\end{equation*} 
We then have, see definition-proposition \ref{RCHdef0}(ii),
the following commutative diagram in $\Fun(\Delta,\Var(k))$
\begin{equation*}
\xymatrix{U=U_{c(\bullet)}\ar[r]^j & \bar X=\bar X_{c(\bullet)} & \, & \bar D_{s_{g'}(\bullet)}\ar[ll]_{i_{\bullet}} \\
U_T=U_{T,c(\bullet)}\ar[r]^{j'}\ar[u]^{g'} & \bar X'=X'_{c(\bullet)}\ar[u]^{\bar{g}'} & \bar D'_{\bullet}\ar[l]_{i'_{\bullet}} & 
\bar{g}^{'-1}(\bar D_{s_{g'}(\bullet)})\ar[l]_{i''_{g'\bullet}}\ar[u]^{(\bar{g}')'_{\bullet}}:i'_{g\bullet}}
\end{equation*}
We then consider the following map in $C(\Var(k)^{2,pr}/T)$, 
\begin{eqnarray*}
T(g,\hat R^{CH})(\mathbb Z(U/S)):g^*\hat R_{(\bar X,\bar D)/S}(\mathbb Z(U/S)) \\
\xrightarrow{:=}g^*\Cone(\mathbb Z(i_{\bullet}\times I):
(\mathbb Z^{tr}((\bar D_{\bullet}\times S,D_{\bullet})/S),u_{IJ})\to\mathbb Z^{tr}((\bar X\times S,X)/S))(-d_X)[-2d_X] \\ 
\xrightarrow{T(g,L)(-)\circ T(g,c)(-)} \\
\Cone(\mathbb Z(i'_{g\bullet}\times I): 
(\mathbb Z^{tr}((\bar D_{\bullet}\times T,\bar g^{-1}(D_{s_g(\bullet)})/T),u_{IJ})
\to\mathbb Z^{tr}((\bar X\times T,X')/T)))(-d_X)[-2d_X] \\
\xrightarrow{([\Gamma_{\bar g'_{\bullet}}]^t,[\Gamma_{\bar g}]^t)} \\
\Cone(\mathbb Z(i'_{g\bullet}\times I): \\
(\mathbb Z^{tr}((\bar{g}^{-1}(\bar D_{s_g(\bullet)})\times T,\bar g^{-1}(D_{s_g(\bullet)})/T),u_{IJ})
\to\mathbb Z^{tr}((\bar X'\times T,X')/T)))(-d_{X'})[-2d_{X'}] \\
\xrightarrow{(\mathbb Z(i''_{g\bullet}\times I),I)(-d_{X'})[-2d_{X'}]} \\
\Cone(\mathbb Z(i'_{\bullet}\times I):
((\mathbb Z^{tr}((\bar D'_{\bullet}\times T,D'_{\bullet)})/T),u_{IJ})\to
\mathbb Z^{tr}((\bar X'\times S,X')/T)))(-d_{X'})[-2d_{X'}] \\
\xrightarrow{=:}\hat R_{(\bar X',\bar D')/T}(\mathbb Z(U_T/T))
\end{eqnarray*}
For  
\begin{eqnarray*}
Q^*:=(\cdots\to\oplus_{\alpha\in\Lambda^n}\mathbb Z(U^n_{\alpha}/S)
\xrightarrow{(\mathbb Z(g^n_{\alpha,\beta}))}\oplus_{\beta\in\Lambda^{n-1}}\mathbb Z(U^{n-1}_{\beta}/S)\to\cdots)
\in C(\Var(k)/S)
\end{eqnarray*}
a complex of (maybe infinite) direct sum of representable presheaves with $h^n_{\alpha}:U^n_{\alpha}\to S$ smooth,
we get the map in $C(\Var(k)^{2,smpr}/T)$
\begin{eqnarray*}
T(g,\hat R^{CH})(Q^*):g^*\hat R^{CH}(Q^*)=
(\cdots\to\oplus_{\alpha\in\Lambda^n}\varinjlim_{(\bar X^n_{\alpha},\bar D^n_{\alpha})/S}
g^*\hat R_{(\bar X^n_{\alpha},\bar D^n_{\alpha})/S}(\mathbb Z(U^n_{\alpha}/S))\to\cdots) \\
\xrightarrow{(T(g,\hat R^{CH})(\mathbb Z(U^n_{\alpha}/S)))} 
(\cdots\to\oplus_{\alpha\in\Lambda^n}\varinjlim_{(\bar X^{n'}_{\alpha},\bar D^{n'}_{\alpha})/T}
\hat R_{(\bar X^{n'}_{\alpha},\bar D^{n'}_{\alpha})/T}(\mathbb Z(U^n_{\alpha,T}/S))\to\cdots)=:\hat R^{CH}(g^*Q^*)
\end{eqnarray*}
together with the commutative diagram in $C(\Var(k)^{2,smpr}/T)$
\begin{equation*}
\xymatrix{g^*\hat R^{CH}(Q^*)\ar[rrr]^{T(g,\hat R^{CH})(Q^*)}\ar[d]_{g^*T(\hat R_S^{CH},R_S^{CH})(Q^*)} 
& \, & \, & \hat R^{CH}(g^*Q^*)\ar[d]^{T(\hat R_T^{CH},R_T^{CH})(g^*Q)} \\
g^*L\rho_{S*}\mu_{S*}R^{CH}(Q^*)\ar[rrr]^{T(g,R^{CH})(Q^*)\circ T(g,\mu)(-)\circ T(g,\rho)(-)\circ T(g,L)(-)} 
& \, & \, & L\rho_{T*}\mu_{T*}R^{CH}(g^*Q^*)}.
\end{equation*}
Let $F\in\PSh(\Var(k)^{sm}/S)$. Consider 
\begin{eqnarray*}
q:LF:=(\cdots\to\oplus_{(U_{\alpha},h_{\alpha})\in\Var(k)^{sm}/S}\mathbb Z(U_{\alpha}/S)\to\cdots)\to F
\end{eqnarray*}
the canonical projective resolution given in subsection 2.3.3.
We then get in particular the map in $C(\Var(k)^{2,smpr}/T)$
\begin{eqnarray*}
T(g,\hat R^{CH})(\rho_S^*LF):g^*\hat R^{CH}(\rho_S^*LF)= \\
(\cdots\to\oplus_{(U_{\alpha},h_{\alpha})\in\Var(k)^{sm}/S}\varinjlim_{(\bar X_{\alpha},\bar D_{\alpha})/S}
g^*\hat R_{(\bar X_{\alpha},\bar D_{\alpha})/S}(\mathbb Z(U_{\alpha}/S))\to\cdots) 
\xrightarrow{(T(g,\hat R^{CH})(\mathbb Z(U_{\alpha}/S)))} \\ 
(\cdots\to\oplus_{(U_{\alpha},h_{\alpha})\in\Var(k)^{sm}/S}\varinjlim_{(\bar X'_{\alpha},\bar D'_{\alpha})/T}
\hat R_{(\bar X'_{\alpha},\bar D'_{\alpha})/T}(\mathbb Z(U_{\alpha,T}/S))\to\cdots)=:\hat R^{CH}(\rho_T^*g^*LF), 
\end{eqnarray*}
and by functoriality, we get in particular for $F=F^{\bullet}\in C(\Var(k)^{sm}/S)$, 
the map in $C(\Var(k)^{2,smpr}/T)$
\begin{eqnarray*}
T(g,\hat R^{CH})(\rho_S^*LF):g^*\hat R^{CH}(\rho_S^*LF)\to\hat R^{CH}(\rho_T^*g^*LF)
\end{eqnarray*}
together with the commutative diagram in $C(\Var(k)^{2,smpr}/T)$
\begin{equation*}
\xymatrix{g^*\hat R^{CH}(\rho_S^*LF)\ar[rrr]^{T(g,\hat R^{CH})(\rho_S^*LF)}\ar[d]_{g^*T(\hat R_S^{CH},R_S^{CH})(\rho_S^*LF)} 
& \, & \, & \hat R^{CH}(\rho_T^*g^*LF)\ar[d]^{T(\hat R_T^{CH},R_T^{CH})(\rho_T^*g^*LF)} \\
g^*L\rho_{S*}\mu_{S*}R^{CH}(\rho_S^*LF)
\ar[rrr]^{L\rho_{T*}\mu_{T*}T(g,R^{CH})(\rho_S^*LF)\circ T(g,\mu)(-)\circ T(g,\rho)(-)\circ T(g,L)(-)} 
& \, & \, & L\rho_{T*}\mu_{T*}R^{CH}(\rho_T^*g^*LF)}.
\end{equation*}

\item Let $S_1,S_2\in\SmVar(k)$ and $p:S_1\times S_2\to S_1$ the projection. 
Let $h:U\to S_1$ a smooth morphism with $U\in\Var(k)$. Consider the cartesian square
\begin{equation*}
\xymatrix{U\times S_2\ar[r]^{h\times I}\ar[d]^{p'} & S_1\times S_2\ar[d]^p \\
U\ar[r]^h & S_1}
\end{equation*}
Take, see definition-proposition \ref{RCHdef0}(i),a compactification $\bar f_0=\bar h:\bar X_0\to\bar S_1$ of $h:U\to S_1$.
Then $\bar f_0\times I:\bar X_0\times S_2\to\bar S_1\times S_2$ is a compactification of $h\times I:U\times S_2\to S_1\times S_2$ 
and $p':U\times S_2\to U$ extend to $\bar{p}'_0:=p_{X_0}:\bar X_0\times S_2\to\bar X_0$. Denote $Z=X_0\backslash U$. 
Take see theorem \ref{desVar}(i), a strict desingularization 
$\bar\epsilon:(\bar X,\bar D)\to(\bar X_0,\bar Z)$ of the pair $(\bar X_0,\bar Z)$.  
We then have the commutative diagram (\ref{RCHdiap}) in $\Fun(\Delta,\Var(k))$ whose squares are cartesian
\begin{equation*}
\xymatrix{U=U_{c(\bullet)}\ar[r]^j & \bar X & \bar D_{\bullet}\ar[l]_{i_{\bullet}} \\
U\times S_2=(U\times S_2)_{c(\bullet)}\ar[r]^{j\times I}\ar[u]^g & \bar X\times S_2\ar[u]^{\bar{p}':=p_{\bar X}} & 
\bar D_{\bullet}\times S_2\ar[l]_{i'_{\bullet}}\ar[u]^{\bar{p'}'_{\bullet}}}
\end{equation*}
Then the map in $C(\Var(k)^{2,smpr}/S_1\times S_2)$
\begin{eqnarray*}
T(p,\hat R^{CH})(\mathbb Z(U/S_1)):p^*\hat R_{(\bar X,\bar D)/S_1}(\mathbb Z(U/S_1))\xrightarrow{\sim} 
\hat R_{(\bar X\times S_2,\bar D_{\bullet}\times S_2)/S_1\times S_2}(\mathbb Z(U\times S_2/S_1\times S_2))
\end{eqnarray*}
is an isomorphism.
Hence, for $Q^*\in C(\Var(k)/S_1)$ a complex of (maybe infinite) direct sum of representable presheaves of smooth morphism,
the map in $C(\Var(k)^{2,smpr}/S_1\times S_2)$
\begin{eqnarray*}
T(p,\hat R^{CH})(Q^*):p^*\hat R^{CH}(Q^*)\xrightarrow{\sim}\hat R^{CH}(p^*Q^*)
\end{eqnarray*}
is an isomorphism.
In particular, for $F\in C(\Var(k)^{sm}/S_1)$ the map in $C(\Var(k)^{2,smpr}/S_1\times S_2)$
\begin{eqnarray*}
T(p,\hat R^{CH})(\rho_{S_1}^*LF):p^*\hat R^{CH}(\rho_{S_1}^*LF)\xrightarrow{\sim}\hat R^{CH}(\rho_{S_1\times S_2}^*p^*LF)
\end{eqnarray*}
is an isomorphism.

\item Let $h_1:U_1\to S$, $h_2:U_2\to S$ two morphisms with $U_1,U_2,S\in\Var(k)$, $U_1,U_2$ smooth.
Denote by $p_1:U_1\times_SU_2\to U_1$ and $p_2:U_1\times_S U_2\to U_2$ the projections.
Take, see definition-proposition \ref{RCHdef0}(i)),
a compactification $\bar f_{10}=\bar{h}_1:\bar X_{10}\to\bar S$ of $h_1:U_1\to S$
and a compactification $\bar f_{20}=\bar{h}_2:\bar X_{20}\to\bar S$ of $h_2:U_2\to S$. Then, 
\begin{itemize}
\item $\bar f_{10}\times\bar f_{20}:\bar X_{10}\times_{\bar S}\bar X_{20}\to S$ 
is a compactification of $h_1\times h_2:U_1\times_SU_2\to S$.
\item $\bar p_{10}:=p_{X_{10}}:\bar X_{10}\times_{\bar S} \bar X_{20}\to\bar X_{10}$ 
is a compactification of $p_1:U_1\times_SU_2\to U_1$.
\item $\bar p_{20}:=p_{X_{20}}:\bar X_{10}\times_{\bar S}\bar X_{20}\to\bar X_{20}$ 
is a compactification of $p_2:U_1\times_SU_2\to U_2$.
\end{itemize}
Denote $\bar Z_1=\bar X_{10}\backslash U_1$ and $\bar Z_2=\bar X_{20}\backslash U_2$.
Take, see theorem \ref{desVar}(i), a strict desingularization 
$\bar\epsilon_1:(\bar X_1,\bar D)\to(\bar X_{10},Z_1)$ of the pair $(\bar X_{10},\bar Z_1)$
and a  strictdesingularization 
$\bar\epsilon_2:(\bar X_2,\bar E)\to(\bar X_{20},Z_2)$ of the pair $(\bar X_{20},\bar Z_2)$. 
Take then a strict desingularization
\begin{equation*}
\bar\epsilon_{12}:((\bar X_1\times_{\bar S}\bar X_2)^N,\bar F)\to
(\bar X_1\times_{\bar S}\bar X_2,(D\times_{\bar S}\bar X_2)\cup(\bar X_1\times_{\bar S}\bar E)) 
\end{equation*}
of the pair $(\bar X_1\times_{\bar S}\bar X_2,(\bar D\times_{\bar S}\bar X_2)\cup(\bar X_1\times_{\bar S}\bar E))$. 
We have then the following commutative diagram 
\begin{equation*}
\xymatrix{\, & \bar X_1\ar[r]^{\bar f_1} & \bar S \\
\, & \bar X_1\times_{\bar S}\bar X_2\ar[r]^{\bar p_1}\ar[u]^{\bar p_2} &\bar X_2\ar[u]^{\bar f_2} \\
(\bar X_1\times_{\bar S}\bar X_2)^N\ar[ru]^{\bar\epsilon_{12}}\ar[rru]^{(\bar p_1)^N}\ar[ruu]^{(\bar p_2)^N} & \, & \,}
\end{equation*}
and
\begin{itemize}
\item $\bar f_1\times\bar f_2:\bar X_1\times_{\bar S}\bar X_2\to\bar S$ 
is a compactification of $h_1\times h_2:U_1\times_SU_2\to S$.
\item $(\bar p_1)^N:=\bar p_1\circ\epsilon_{12}:(\bar X_1\times_{\bar S}\bar X_2)^N\to\bar X_1$ is a compactification of 
$p_1:U_1\times_SU_2\to U_1$.
\item $(\bar p_2)^N:=\bar p_2\circ\epsilon_{12}:(\bar X_1\times_{\bar S}\bar X_2)^N\to\bar X_2$ is a compactification of 
$p_2:U_1\times_SU_2\to U_2$.
\end{itemize}
We have then the morphism in $C(\Var(k)^{2,smpr}/S)$
\begin{eqnarray*}
T(\otimes,\hat R_S^{CH})(\mathbb Z(U_1/S),\mathbb Z(U_2/S)):=\hat R_S^{CH}(p_1)\otimes\hat R_S^{CH}(p_2): \\
\hat R_{(\bar X_1,\bar D)/S}(\mathbb Z(U_1/S))\otimes\hat R_{(X_2,E))/S}(\mathbb Z(U_2/S)) 
\xrightarrow{\sim}\hat R_{(\bar X_1\times_{\bar S}\bar X_2)^N,\bar F)/S}(\mathbb Z(U_1\times_S U_2/S))
\end{eqnarray*}
For   
\begin{eqnarray*}
Q_1^*:=(\cdots\to\oplus_{\alpha\in\Lambda^n}\mathbb Z(U^n_{1,\alpha}/S)
\xrightarrow{(\mathbb Z(g^n_{\alpha,\beta}))}\oplus_{\beta\in\Lambda^{n-1}}\mathbb Z(U^{n-1}_{1,\beta}/S)\to\cdots), \\
Q_2^*:=(\cdots\to\oplus_{\alpha\in\Lambda^n}\mathbb Z(U^n_{2,\alpha}/S)
\xrightarrow{(\mathbb Z(g^n_{\alpha,\beta}))}\oplus_{\beta\in\Lambda^{n-1}}\mathbb Z(U^{n-1}_{2,\beta}/S)\to\cdots)
\in C(\Var(k)/S)
\end{eqnarray*}
complexes of (maybe infinite) direct sum of representable presheaves with $U^*_{\alpha}$ smooth, 
we get the morphism in $C(\Var(k)^{2,smpr}/S)$
\begin{eqnarray*}
T(\otimes,\hat R_S^{CH})(Q^*_1,Q^*_2):\hat R^{CH}(Q^*_1)\otimes R^{CH}(Q^*_2) 
\xrightarrow{(T(\otimes,\hat R_S^{CH})(\mathbb Z(U_{1,\alpha}^m),\mathbb Z(U_{2,\beta}^n))}
\hat R^{CH}(Q^*_1\otimes Q^*_2))
\end{eqnarray*},
together with the commutative diagram in $C(\Var(k)^{2,smpr}/S)$
\begin{equation*}
\xymatrix{\hat R^{CH}(Q^*_1)\otimes R^{CH}(Q^*_2)
\ar[rrr]^{T(\otimes,\hat R_S^{CH})(Q_1^*,Q_2)}
\ar[d]_{T(\hat R_S^{CH},R_S^{CH})(Q_1^*)\otimes T(\hat R_S^{CH},R_S^{CH})(Q_2^*)} 
& \, & \, & \hat R^{CH}(Q^*_1\times Q_2^*)\ar[d]^{T(\hat R_S^{CH},R_S^{CH})(Q^*_1\otimes Q^*_2)} \\
L\rho_{S*}\mu_{S*}(R^{CH}(Q_1^*)\otimes R^{CH}(Q_2^*))
\ar[rrr]^{L\rho_{S*}\mu_{S*}T(\otimes,R_S^{CH})(Q_1^*,Q_2)} 
& \, & \, & L\rho_{S*}\mu_{S*}R^{CH}(Q_1^*\otimes Q_2)}.
\end{equation*}
For $F_1,F_2\in C(\Var(k)^{sm}/S)$, we get in particular the morphism in $C(\Var(k)^2/S)$
\begin{eqnarray*}
T(\otimes,R_S^{CH})(\rho_S^*LF_1,\rho_S^*LF_2):R^{CH}(\rho_S^*LF_1)\otimes R^{CH}(\rho_S^*LF_2) 
\to R^{CH}(\rho_S^*(LF_1\otimes LF_2))
\end{eqnarray*}
together with the commutative diagram in $C(\Var(k)^{2,smpr}/S)$
\begin{equation*}
\xymatrix{\hat R^{CH}(\rho_S^*LF_1)\otimes R^{CH}(\rho_S^*LF_2)
\ar[rrr]^{T(\otimes,\hat R_S^{CH})(\rho_S^*LF_1,\rho_S^*LF_2)}
\ar[d]_{T(\hat R_S^{CH},R_S^{CH})(\rho_S^*LF_1)\otimes T(\hat R_S^{CH},R_S^{CH})(\rho_S^*LF_2)} 
& \, & \, & \hat R^{CH}(\rho_S^*LF_1\times\rho_S^*LF_2)
\ar[d]^{T(\hat R_S^{CH},R_S^{CH})(\rho_S^*LF_1\otimes\rho_S^*LF_2)} \\
L\rho_{S*}\mu_{S*}(R^{CH}(\rho_S^*LF_1)\otimes R^{CH}(\rho_S^*LF_2))
\ar[rrr]^{L\rho_{S*}\mu_{S*}T(\otimes,R_S^{CH})(\rho_S^*LF_1,\rho_S^*LF_2)} 
& \, & \, & L\rho_{S*}\mu_{S*}R^{CH}(\rho_S^*LF_1\times\rho_S^*LF_2)}.
\end{equation*}

\end{itemize}


For $S\in\Var(k)$, we will use rather the functors $R^{0CH}_S$ and $\hat R^{0CH}_S$
since we are working in the image of the graph functor $\Gr_S^{12}:\Var(k)/S\to\Var(k)^2/S$.
We have the full subcategory $\SmVar(k)/S\subset\Var(k)/S$ 
whose objects are morphisms $f:X\to S$ with $X\in\SmVar(k)$. Then $\Gr_S^{12}(\SmVar(k)/S)\subset\Var(k)^{2,smpr}/S$.
If $S\in\SmVar(k)$, we have the factorization of morphism of site 
\begin{equation*}
\Gr_S^{12}:\Var(k)^{2,smpr}/S\xrightarrow{\Gr_S^{12}}\SmVar(k)/S\xrightarrow{\rho_S}\Var(k)^{sm}/S.
\end{equation*}

\begin{defi}\label{R0CHdef} 
\begin{itemize}
\item[(i)]Let $h:U\to S$ a morphism, with $U,S\in\Var(k)$ and $U$ smooth.
Take, see definition-proposition \ref{RCHdef0},
$\bar f_0=\bar{h}_0:\bar X_0\to\bar S$ a compactification of $h:U\to S$ and denote by $\bar Z=\bar X_0\backslash U$.
Take, using theorem \ref{desVar}(ii), a strict desingularization 
$\bar\epsilon:(\bar X,\bar D)\to(\bar X_0,\bar Z)$ of the pair $(\bar X_0,\bar Z)$, with $\bar X\in\PSmVar(k)$ and 
$\bar D:=\epsilon^{-1}(\bar Z)=\cup_{i=1}^s\bar D_i\subset\bar X$ a normal crossing divisor.  
We denote by $i_{\bullet}:\bar D_{\bullet}\hookrightarrow\bar X=\bar X_{c(\bullet)}$ the morphism of simplicial varieties
given by the closed embeddings $i_I:\bar D_I=\cap_{i\in I}\bar D_i\hookrightarrow\bar X$
We denote by $j:U\hookrightarrow\bar X$ the open embedding. We then consider the following map in $C(\Var(k)/S)$
\begin{eqnarray*}
r^0_{(\bar X,\bar D)/S}(\mathbb Z(U/S)):R^0_{(\bar X,\bar D)/S}(\mathbb Z(U/S)) \\
\xrightarrow{:=}
\bar f_*E_{et}(\Cone(\mathbb Z(i_{\bullet}):(\mathbb Z((\bar D_{\bullet})/\bar X),u_{IJ})\to\mathbb Z((\bar X/\bar X)))) \\
\xrightarrow{\bar f_*E_{et}(0,k\circ\ad(j^*,j_*)(\mathbb Z(\bar X/\bar X)))}
h_*E_{et}(\mathbb Z(U/U))=:\mathbb D^0_S(\mathbb Z(U/S)).
\end{eqnarray*}
Note that $\mathbb Z(\bar D_I/\bar X)$ and $\mathbb Z(\bar X/\bar X)$ are obviously $\mathbb A^1$ invariant.
Note that $r_{(X,D)/S}$ is NOT an equivalence $(\mathbb A^1,et)$ local by proposition \ref{rho1} since
$\rho_{\bar X*}\mathbb Z(\bar D_{\bullet}/\bar X)=0$, 
and $\rho_{\bar X*}\ad(j^*,j_*)(\mathbb Z(\bar X/\bar X))$ is not an equivalence $(\mathbb A^1,et)$ local.
\item[(ii)]Let $g:U'/S\to U/S$ a morphism, with $U'/S=(U',h'),U/S=(U,h)\in\Var(k)/S$, with $U$ and $U'$ smooth.
Take, see definition-proposition \ref{RCHdef0}(ii),a compactification $\bar f_0=\bar h:\bar X_0\to\bar S$ of $h:U\to S$ 
and a compactification $\bar f'_0=\bar{h}':\bar X'_0\to\bar S$ of $h':U'\to S$ such that 
$g:U'/S\to U/S$ extend to a morphism $\bar g_0:\bar X'_0/\bar S\to\bar X_0/\bar S$. 
Denote $\bar Z=\bar X_0\backslash U$ and $\bar Z'=\bar X'_0\backslash U'$.
Take, see definition-proposition \ref{RCHdef0}(ii), a strict desingularization 
$\bar\epsilon:(\bar X,\bar D)\to(\bar X_0,\bar Z)$ of $(\bar X_0,\bar Z)$,
a strict desingularization $\bar\epsilon'_{\bullet}:(\bar X',\bar D')\to(\bar X'_0,\bar Z')$ of $(\bar X'_0,\bar Z')$
and a morphism $\bar g:\bar X'\to\bar X$ such that the following diagram commutes
\begin{equation*}
\xymatrix{\bar X'_0\ar[r]^{\bar{g}_0} & \bar X_0 \\
\bar X'\ar[u]^{\bar\epsilon'}\ar[r]^{\bar g} & \bar X\ar[u]^{\bar\epsilon}}.
\end{equation*} 
We then have, see definition-proposition \ref{RCHdef0}(ii),
the commutative diagram (\ref{RCHdia}) in $\Fun(\Delta,\Var(k))$
\begin{equation*}
\xymatrix{U=U_{c(\bullet)}\ar[r]^j & \bar X=\bar X_{c(\bullet)} & \, & \bar D_{s_g(\bullet)}\ar[ll]_{i_{\bullet}} \\
U'=U'_{c(\bullet)}\ar[r]^{j'}\ar[u]^g & \bar X'=\bar X'_{c(\bullet)}\ar[u]^{\bar{g}} & \bar D'_{\bullet}\ar[l]_{i'_{\bullet}} & 
\bar{g}^{-1}(\bar D_{s_g(\bullet)})\ar[l]_{i''_{g\bullet}}\ar[u]^{\bar{g}'_{\bullet}}:i'_{g\bullet}}
\end{equation*}
We then consider the following map in $C(\Var(k)/S)$
\begin{eqnarray*}
R_S^{0CH}(g):R^0_{(\bar X,\bar D)/S}(\mathbb Z(U/S))\xrightarrow{:=} \\
\bar f_*E_{et}(\Cone(\mathbb Z(i_{\bullet}):(\mathbb Z((\bar D_{s_g(\bullet)})/\bar X),u_{IJ})\to\mathbb Z(\bar X/\bar X))) \\
\xrightarrow{T(\bar g,E)(-)\circ p_{S*}\ad(\bar{g}^*,\bar{g}_*)(-)} \\
\bar f'_*E_{et}(\Cone(\mathbb Z(i'_{g\bullet}): 
(\mathbb Z((\bar{g}^{-1}(\bar D_{s_g(\bullet)})/\bar X'),u_{IJ})\to\mathbb Z((\bar X'/\bar X'))) \\
\xrightarrow{\bar f'_*E_{et}(\mathbb Z(i''_{g\bullet}),I)} \\
\bar f'_*E_{et}(\Cone(\mathbb Z(i'_{\bullet}):(\mathbb Z(\bar D'_{\bullet}/\bar X'),u_{IJ})\to\mathbb Z(\bar X'/\bar X'))) \\
\xrightarrow{=:}R^0_{(\bar X',\bar D')/S}(\mathbb Z(U'/S))
\end{eqnarray*}
Then by the diagram (\ref{RCHdia}) and adjonction, the following diagram in $C(\Var(k)/S)$ obviously commutes
\begin{equation*}
\xymatrix{R^0_{(\bar X,\bar D)/S}(\mathbb Z(U/S))\ar[rr]^{r_{(\bar X,\bar D)/S}(\mathbb Z(U/S))}\ar[d]_{R_S^{0CH}(g)} & \, & 
h_*E_{et}(\mathbb Z(U/U)=:\mathbb D^0_S(\mathbb Z(U/S))\ar[d]^{D_S(g):=T(g,E)(-)\circ\ad(g^*,g_*)(E_{et}(\mathbb Z(U/U)))} \\
R^0_{(\bar X',\bar D')/S}(\mathbb Z(U'/S))\ar[rr]^{r^0_{(\bar X',\bar D')/S}(\mathbb Z(U'/S))} & \, & 
h'_*E_{et}(\mathbb Z(U'/U'))=:\mathbb D^0_S(\mathbb Z(U'/S))}.
\end{equation*}
\item[(iii)]For $g_1:U''/S\to U'/S$, $g_2:U'/S\to U/S$ two morphisms
with $U''/S=(U',h''),U'/S=(U',h'),U/S=(U,h)\in\Var(k)/S$, with $U$, $U'$ and $U''$ smooth.
We get from (i) and (ii) 
a compactification $\bar f=\bar{h}:\bar X\to\bar S$ of $h:U\to S$,
a compactification $\bar f'=\bar{h}':\bar X'\to\bar S$ of $h':U'\to S$,
and a compactification $\bar f''=\bar{h}'':\bar X''\to\bar S$ of $h'':U''\to S$,
with $\bar X,\bar X',\bar X''\in\PSmVar(k)$, 
$\bar D:=\bar X\backslash U\subset\bar X$ $\bar D':=\bar X'\backslash U'\subset\bar X'$, 
and $\bar D'':=\bar X''\backslash U''\subset\bar X''$ normal crossing divisors, 
such that $g_1:U''/S\to U'/S$ extend to $\bar g_1:\bar X''/\bar S\to\bar X'/\bar S$,
$g_2:U'/S\to U/S$ extend to $\bar g_2:\bar X'/\bar S\to\bar X/\bar S$, and
\begin{eqnarray*}
R_S^{0CH}(g_2\circ g_1)=R_S^{0CH}(g_1)\circ R_S^{0CH}(g_2):R^0_{(\bar X,\bar D)/S}\to R^0_{(\bar X'',\bar D'')/S} 
\end{eqnarray*}
\item[(iv)] For  
\begin{eqnarray*}
Q^*:=(\cdots\to\oplus_{\alpha\in\Lambda^n}\mathbb Z(U^n_{\alpha}/S)
\xrightarrow{(\mathbb Z(g^n_{\alpha,\beta}))}\oplus_{\beta\in\Lambda^{n-1}}\mathbb Z(U^{n-1}_{\beta}/S)\to\cdots)
\in C(\Var(k)/S)
\end{eqnarray*}
a complex of (maybe infinite) direct sum of representable presheaves with $U^*_{\alpha}$ smooth,
we get from (i), (ii) and (iii) the map in $C(\Var(k)/S)$
\begin{eqnarray*}
r_S^{0CH}(Q^*):R^{0CH}(Q^*):=
(\cdots\to\oplus_{\beta\in\Lambda^{n-1}}\varinjlim_{(\bar X^{n-1}_{\beta},\bar D^{n-1}_{\beta})/S}
R^0_{(\bar X^{n-1}_{\beta},\bar D^{n-1}_{\beta})/S}(\mathbb Z(U^{n-1}_{\beta}/S)) \\
\xrightarrow{(R_S^{CH}(g^n_{\alpha,\beta}))}\oplus_{\alpha\in\Lambda^n}\varinjlim_{(\bar X^n_{\alpha},\bar D^n_{\alpha})/S}
R^0_{(\bar X^n_{\alpha},\bar D^n_{\alpha})/S}(\mathbb Z(U^n_{\alpha}/S))\to\cdots)
\to\mathbb D_S(Q^*),
\end{eqnarray*}
where for $(U^n_{\alpha},h^n_{\alpha})\in\Var(k)/S$, the inductive limit run over all the compactifications 
$\bar f_{\alpha}:\bar X_{\alpha}\to\bar S$ of $h_{\alpha}:U_{\alpha}\to S$ with $\bar X_{\alpha}\in\PSmVar(k)$
and $\bar D_{\alpha}:=\bar X_{\alpha}\backslash U_{\alpha}$ a normal crossing divisor.
For $m=(m^*):Q_1^*\to Q_2^*$ a morphism with 
\begin{eqnarray*}
Q_1^*:=(\cdots\to\oplus_{\alpha\in\Lambda^n}\mathbb Z(U^n_{1,\alpha}/S)
\xrightarrow{(\mathbb Z(g^n_{\alpha,\beta}))}\oplus_{\beta\in\Lambda^{n-1}}\mathbb Z(U^{n-1}_{1,\beta}/S)\to\cdots), \\
Q_2^*:=(\cdots\to\oplus_{\alpha\in\Lambda^n}\mathbb Z(U^n_{2,\alpha}/S)
\xrightarrow{(\mathbb Z(g^n_{\alpha,\beta}))}\oplus_{\beta\in\Lambda^{n-1}}\mathbb Z(U^{n-1}_{2,\beta}/S)\to\cdots)
\in C(\Var(k)/S)
\end{eqnarray*}
complexes of (maybe infinite) direct sum of representable presheaves with $U^*_{1,\alpha}$ and $U^*_{2,\alpha}$ smooth,
we get again from (i), (ii) and (iii) a commutative diagram in $C(\Var(k)/S)$
\begin{equation*}
\xymatrix{R^{0CH}(Q_2^*)\ar[rr]^{r_S^{0CH}(Q_2^*)}\ar[d]_{R_S^{0CH}(m):=(R_S^{0CH}(m^*))} & \, & 
\mathbb D^0_S(Q_2^*)\ar[d]^{\mathbb D_S(m):=(\mathbb D^0_S(m^*))} \\
R^{0CH}(Q_1^*)\ar[rr]^{r_S^{0CH}(Q_1^*)} & \, & \mathbb D^0_S(Q_1^*)}.
\end{equation*}
\item[(v)] Let  
\begin{eqnarray*}
Q^*:=(\cdots\to\oplus_{\alpha\in\Lambda^n}\mathbb Z(U^n_{\alpha}/S)
\xrightarrow{(\mathbb Z(g^n_{\alpha,\beta}))}\oplus_{\beta\in\Lambda^{n-1}}\mathbb Z(U^{n-1}_{\beta}/S)\to\cdots)
\in C(\Var(k)/S)
\end{eqnarray*}
a complex of (maybe infinite) direct sum of representable presheaves with $U^*_{\alpha}$ smooth,
we have by definition 
\begin{equation*}
\Gr_S^{12*}R^{0CH}(Q^*)=R^{CH}(Q^*)\in C(\Var(k)^2/S).
\end{equation*}
\end{itemize}
\end{defi}

\begin{itemize}
\item Let $S\in\Var(k)$
For $(h,m,m')=(h^*,m^*,m^{'*}):Q_1^*[1]\to Q_2^*$ an homotopy with $Q_1^*,Q_2^*\in C(\Var(k)/S)$
complexes of (maybe infinite) direct sum of representable presheaves with $U^*_{1,\alpha}$ and $U^*_{2,\alpha}$ smooth,
\begin{equation*}
(R_S^{0CH}(h),R_S^{0CH}(m),R_S^{0CH}(m'))=(R_S^{0CH}(h^*),R_S^{0CH}(m^*),R_S^{0CH}(m^{'*})):
R^{0CH}(Q_2^*)[1]\to R^{0CH}(Q_1^*)
\end{equation*}
is an homotopy in $C(\Var(k)/S)$ using definition \ref{R0CHdef} (iii). 
In particular if $m:Q_1^*\to Q_2^*$ with $Q_1^*,Q_2^*\in C(\Var(k)/S)$
complexes of (maybe infinite) direct sum of representable presheaves with $U^*_{1,\alpha}$ and $U^*_{2,\alpha}$ smooth
is an homotopy equivalence, then $R_S^{0CH}(m):R^{0CH}(Q_2^*)\to R^{0CH}(Q_1^*)$ is an homotopy equivalence.
\item Let $S\in\SmVar(k)$. Let $F\in\PSh(\Var(k)^{sm}/S)$. Consider 
\begin{eqnarray*}
q:LF:=(\cdots\to\oplus_{(U_{\alpha},h_{\alpha})\in\Var(k)^{sm}/S}\mathbb Z(U_{\alpha}/S)
\xrightarrow{(\mathbb Z(g^n_{\alpha,\beta}))}
\oplus_{(U_{\alpha},h_{\alpha})\in\Var(k)^{sm}/S}\mathbb Z(U_{\alpha}/S)\to\cdots)\to F
\end{eqnarray*}
the canonical projective resolution given in subsection 2.3.3.
Note that the $U_{\alpha}$ are smooth since $S$ is smooth and $h_{\alpha}$ are smooth morphism.
Definition \ref{R0CHdef}(iv) gives in this particular case the map in $C(\Var(k)/S)$
\begin{eqnarray*}
r_S^{0CH}(\rho_S^*LF):R^{0CH}(\rho_S^*LF):=
(\cdots\to\oplus_{(U_{\alpha},h_{\alpha})\in\Var(k)^{sm}/S}\varinjlim_{(\bar X_{\alpha},\bar D_{\alpha})/S}
R^0_{(\bar X_{\alpha},\bar D_{\alpha})/S}(\mathbb Z(U_{\alpha}/S)) \\
\xrightarrow{(R_S^{0CH}(g^n_{\alpha,\beta}))}
\oplus_{(U_{\alpha},h_{\alpha})\in\Var(k)^{sm}/S}\varinjlim_{(\bar X_{\alpha},\bar D_{\alpha})/S}
R^0_{(\bar X_{\alpha},\bar D_{\alpha})/S}(\mathbb Z(U_{\alpha}/S))\to\cdots)
\to\mathbb D^0_S(\rho_S^*LF),
\end{eqnarray*}
where for $(U_{\alpha},h_{\alpha})\in\Var(k)^{sm}/S$, the inductive limit run over all the compactifications 
$\bar f_{\alpha}:\bar X_{\alpha}\to\bar S$ of $h_{\alpha}:U_{\alpha}\to S$ with $\bar X_{\alpha}\in\PSmVar(\mathbb C)$
and $\bar D_{\alpha}:=\bar X_{\alpha}\backslash U_{\alpha}$ a normal crossing divisor.
Definition \ref{R0CHdef}(iv) gives then by functoriality in particular, for $F=F^{\bullet}\in C(\Var(k)^{sm}/S)$, 
the map in $C(\Var(k)/S)$
\begin{eqnarray*}
r_S^{0CH}(\rho_S^*LF)=(r_S^{0CH}(\rho_S^*LF^*)):R^{0CH}(\rho_S^*LF)\to\mathbb D^0_S(\rho_S^*LF).
\end{eqnarray*}

\item Let $g:T\to S$ a morphism with $T,S\in\SmVar(k)$. Let $h:U\to S$ a smooth morphism with $U\in\Var(k)$.
Consider the cartesian square
\begin{equation*}
\xymatrix{U_T\ar[r]^{h'}\ar[d]^{g'} & T\ar[d]^g \\
U\ar[r]^h & S}
\end{equation*}
Note that $U$ is smooth since $S$ and $h$ are smooth, and $U_T$ is smooth since $T$ and $h'$ are smooth. 
Take, see definition-proposition \ref{RCHdef0}(ii),a compactification $\bar f_0=\bar h:\bar X_0\to\bar S$ of $h:U\to S$. 
Take, see definition-proposition \ref{RCHdef0}(ii), a strict desingularization 
$\bar\epsilon:(\bar X,\bar D)\to(\bar X_0,\bar Z)$ of $(\bar X_0,\bar Z)$.
Then $\bar f'_0=\bar{g\circ h'}:\bar X_T\to\bar T$ is a compactification of $g\circ h':U_T\to S$ such that 
$g':U_T/S\to U/S$ extend to a morphism $\bar{g}'_0:\bar X_T/\bar S\to \bar X/\bar S$. 
Denote $\bar Z=\bar X_0\backslash U$ and $\bar Z'=\bar X_T\backslash U_T$.
Take, see definition-proposition \ref{RCHdef0}(ii), a strict desingularization 
$\epsilon'_{\bullet}:(\bar X',\bar D')\to(\bar X_T,\bar Z')$ of $(\bar X_T,\bar Z')$.
Denote $\bar g'=\bar g'_0\circ\epsilon'_{\bullet}:\bar X'\to\bar X$. 
We then have, see definition-proposition \ref{RCHdef0}(ii),
the following commutative diagram in $\Fun(\Delta,\Var(k))$
\begin{equation*}
\xymatrix{U=U_{c(\bullet)}\ar[r]^j & \bar X=\bar X_{c(\bullet)} & \, & \bar D_{s_{g'}(\bullet)}\ar[ll]_{i_{\bullet}} \\
U_T=U_{T,c(\bullet)}\ar[r]^{j'}\ar[u]^{g'} & \bar X_T=X_{T,c(\bullet)}\ar[u]^{\bar{g}'_0} & \, & 
\bar{g}^{'-1}(\bar D_{s_{g'}(\bullet)})\ar[ll]_{i_{g\bullet}}\ar[u]^{(\bar{g}')'_{\bullet}} \\
U_T=U_{T,c(\bullet)}\ar[r]^{j'}\ar[u]^{g'} & \bar X'=X'_{c(\bullet)}\ar[u]^{\bar{g}'} & \bar D'_{\bullet}\ar[l]_{i'_{\bullet}} & 
\bar{g}^{'-1}(\bar D_{s_{g'}(\bullet)})\ar[l]_{i''_{g'\bullet}}\ar[u]^{\epsilon'_{\bullet}}:i'_{g\bullet}}
\end{equation*}
We then consider the following map in $C(\Var(k)/T)$, see definition \ref{R0CHdef}(ii)
\begin{eqnarray*}
T(g,R^{0CH})(\mathbb Z(U/S)):g^*R^0_{(\bar X,\bar D)/S}(\mathbb Z(U/S)) \\ 
\xrightarrow{g^*R^{0CH}_S(g')}g^*R^0_{(\bar X',\bar D')/S}(\mathbb Z(U_T/S))=g^*g_*R^0_{(\bar X',\bar D')/T}(\mathbb Z(U_T/T)) \\
\xrightarrow{\ad(g^*,g_*)(R^0_{(\bar X',\bar D')/T}(\mathbb Z(U_T/T)))}R^0_{(\bar X',\bar D')/T}(\mathbb Z(U_T/T))
\end{eqnarray*}
For  
\begin{eqnarray*}
Q^*:=(\cdots\to\oplus_{\alpha\in\Lambda^n}\mathbb Z(U^n_{\alpha}/S)
\xrightarrow{(\mathbb Z(g^n_{\alpha,\beta}))}\oplus_{\beta\in\Lambda^{n-1}}\mathbb Z(U^{n-1}_{\beta}/S)\to\cdots)
\in C(\Var(k)/S)
\end{eqnarray*}
a complex of (maybe infinite) direct sum of representable presheaves with $h^n_{\alpha}:U^n_{\alpha}\to S$ smooth,
we get the map in $C(\Var(k)/T)$
\begin{eqnarray*}
T(g,R^{0CH})(Q^*):g^*R^{0CH}(Q^*)=
(\cdots\to\oplus_{\alpha\in\Lambda^n}\varinjlim_{(\bar X^n_{\alpha},\bar D^n_{\alpha})/S}
g^*R^0_{(\bar X^n_{\alpha},\bar D^n_{\alpha})/S}(\mathbb Z(U^n_{\alpha}/S))\to\cdots) \\
\xrightarrow{(T(g,R^{0CH})(\mathbb Z(U^n_{\alpha}/S)))} 
(\cdots\to\oplus_{\alpha\in\Lambda^n}\varinjlim_{(\bar X^{n'}_{\alpha},\bar D^{n'}_{\alpha})/T}
R^0_{(\bar X^{n'}_{\alpha},\bar D^{n'}_{\alpha})/T}(\mathbb Z(U^n_{\alpha,T}/S))\to\cdots)=:R^{CH}(g^*Q^*).
\end{eqnarray*}
Let $F\in\PSh(\Var(k)^{sm}/S)$. Consider 
\begin{eqnarray*}
q:LF:=(\cdots\to\oplus_{(U_{\alpha},h_{\alpha})\in\Var(k)^{sm}/S}\mathbb Z(U_{\alpha}/S)\to\cdots)\to F
\end{eqnarray*}
the canonical projective resolution given in subsection 2.3.3.
We then get in particular the map in $C(\Var(k)/T)$
\begin{eqnarray*}
T(g,R^{0CH})(\rho_S^*LF):g^*R^{0CH}(\rho_S^*LF)= \\
(\cdots\to\oplus_{(U_{\alpha},h_{\alpha})\in\Var(k)^{sm}/S}\varinjlim_{(\bar X_{\alpha},\bar D_{\alpha})/S}
g^*R^0_{(\bar X_{\alpha},\bar D_{\alpha})/S}(\mathbb Z(U_{\alpha}/S))\to\cdots) 
\xrightarrow{(T(g,R^{0CH})(\mathbb Z(U_{\alpha}/S)))} \\ 
(\cdots\to\oplus_{(U_{\alpha},h_{\alpha})\in\Var(k)^{sm}/S}\varinjlim_{(\bar X'_{\alpha},\bar D'_{\alpha})/T}
R^0_{(\bar X'_{\alpha},\bar D'_{\alpha})/T}(\mathbb Z(U_{\alpha,T}/S))\to\cdots)=:R^{CH}(\rho_T^*g^*LF). 
\end{eqnarray*}
By functoriality, we get in particular for $F=F^{\bullet}\in C(\Var(k)^{sm}/S)$, the map in $C(\Var(k)/T)$
\begin{eqnarray*}
T(g,R^{0CH})(\rho_S^*LF):g^*R^{0CH}(\rho_S^*LF)\to R^{0CH}(\rho_T^*g^*LF).
\end{eqnarray*}

\item Let $S_1,S_2\in\SmVar(k)$ and $p:S_1\times S_2\to S_1$ the projection. 
Let $h:U\to S_1$ a smooth morphism with $U\in\Var(k)$. Consider the cartesian square
\begin{equation*}
\xymatrix{U\times S_2\ar[r]^{h\times I}\ar[d]^{p'} & S_1\times S_2\ar[d]^p \\
U\ar[r]^h & S_1}
\end{equation*}
Take, see definition-proposition \ref{RCHdef0}(i),a compactification $\bar f_0=\bar h:\bar X_0\to\bar S_1$ of $h:U\to S_1$.
Then $\bar f_0\times I:\bar X_0\times S_2\to\bar S_1\times S_2$ is a compactification of $h\times I:U\times S_2\to S_1\times S_2$ 
and $p':U\times S_2\to U$ extend to $\bar{p}'_0:=p_{X_0}:\bar X_0\times S_2\to\bar X_0$. Denote $Z=X_0\backslash U$. 
Take see theorem \ref{desVar}(i), a strict desingularization 
$\bar\epsilon:(\bar X,\bar D)\to(\bar X_0,\bar Z)$ of the pair $(\bar X_0,\bar Z)$.  
We then have the commutative diagram (\ref{RCHdiap}) in $\Fun(\Delta,\Var(k))$ whose squares are cartesian
\begin{equation*}
\xymatrix{U=U_{c(\bullet)}\ar[r]^j & \bar X & \bar D_{\bullet}\ar[l]_{i_{\bullet}} \\
U\times S_2=(U\times S_2)_{c(\bullet)}\ar[r]^{j\times I}\ar[u]^g & \bar X\times S_2\ar[u]^{\bar{p}':=p_{\bar X}} & 
\bar D_{\bullet}\times S_2\ar[l]_{i'_{\bullet}}\ar[u]^{\bar{p'}'_{\bullet}}}
\end{equation*}
Then the map in $C(\Var(k)/S_1\times S_2)$
\begin{eqnarray*}
T(p,R^{0CH})(\mathbb Z(U/S_1)):p^*R^0_{(\bar X,\bar D)/S_1}(\mathbb Z(U/S_1))\xrightarrow{\sim} 
R^0_{(\bar X\times S_2,\bar D_{\bullet}\times S_2)/S_1\times S_2}(\mathbb Z(U\times S_2/S_1\times S_2))
\end{eqnarray*}
is an isomorphism.
Hence, for $Q^*\in C(\Var(k)/S_1)$ a complex of (maybe infinite) direct sum of representable presheaves of smooth morphism,
the map in $C(\Var(k)/S_1\times S_2)$
\begin{eqnarray*}
T(p,R^{0CH})(Q^*):p^*R^{0CH}(Q^*)\xrightarrow{\sim}R^{0CH}(p^*Q^*)
\end{eqnarray*}
is an isomorphism.
In particular, for $F\in C(\Var(k)^{sm}/S_1)$ the map in $C(\Var(k)/S_1\times S_2)$
\begin{eqnarray*}
T(p,R^{0CH})(\rho_{S_1}^*LF):p^*R^{0CH}(\rho_{S_1}^*LF)\xrightarrow{\sim}R^{0CH}(\rho_{S_1\times S_2}^*p^*LF)
\end{eqnarray*}
is an isomorphism.

\item Let $h_1:U_1\to S$, $h_2:U_2\to S$ two morphisms with $U_1,U_2,S\in\Var(k)$, $U_1,U_2$ smooth.
Denote by $p_1:U_1\times_SU_2\to U_1$ and $p_2:U_1\times_S U_2\to U_2$ the projections.
Take, see definition-proposition \ref{RCHdef0}(i)),
a compactification $\bar f_{10}=\bar{h}_1:\bar X_{10}\to\bar S$ of $h_1:U_1\to S$
and a compactification $\bar f_{20}=\bar{h}_2:\bar X_{20}\to\bar S$ of $h_2:U_2\to S$. Then, 
\begin{itemize}
\item $\bar f_{10}\times\bar f_{20}:\bar X_{10}\times_{\bar S}\bar X_{20}\to S$ 
is a compactification of $h_1\times h_2:U_1\times_SU_2\to S$.
\item $\bar p_{10}:=p_{X_{10}}:\bar X_{10}\times_{\bar S} \bar X_{20}\to\bar X_{10}$ 
is a compactification of $p_1:U_1\times_SU_2\to U_1$.
\item $\bar p_{20}:=p_{X_{20}}:\bar X_{10}\times_{\bar S}\bar X_{20}\to\bar X_{20}$ 
is a compactification of $p_2:U_1\times_SU_2\to U_2$.
\end{itemize}
Denote $\bar Z_1=\bar X_{10}\backslash U_1$ and $\bar Z_2=\bar X_{20}\backslash U_2$.
Take, see theorem \ref{desVar}(i), a strict desingularization 
$\bar\epsilon_1:(\bar X_1,\bar D)\to(\bar X_{10},Z_1)$ of the pair $(\bar X_{10},\bar Z_1)$
and a  strictdesingularization 
$\bar\epsilon_2:(\bar X_2,\bar E)\to(\bar X_{20},Z_2)$ of the pair $(\bar X_{20},\bar Z_2)$. 
Take then a strict desingularization
\begin{equation*}
\bar\epsilon_{12}:((\bar X_1\times_{\bar S}\bar X_2)^N,\bar F)\to
(\bar X_1\times_{\bar S}\bar X_2,(D\times_{\bar S}\bar X_2)\cup(\bar X_1\times_{\bar S}\bar E)) 
\end{equation*}
of the pair $(\bar X_1\times_{\bar S}\bar X_2,(\bar D\times_{\bar S}\bar X_2)\cup(\bar X_1\times_{\bar S}\bar E))$. 
We have then the following commutative diagram 
\begin{equation*}
\xymatrix{\, & \bar X_1\ar[r]^{\bar f_1} & \bar S \\
\, & \bar X_1\times_{\bar S}\bar X_2\ar[r]^{\bar p_1}\ar[u]^{\bar p_2} &\bar X_2\ar[u]^{\bar f_2} \\
(\bar X_1\times_{\bar S}\bar X_2)^N\ar[ru]^{\bar\epsilon_{12}}\ar[rru]^{(\bar p_1)^N}\ar[ruu]^{(\bar p_2)^N} & \, & \,}
\end{equation*}
and
\begin{itemize}
\item $\bar f_1\times\bar f_2:\bar X_1\times_{\bar S}\bar X_2\to\bar S$ 
is a compactification of $h_1\times h_2:U_1\times_SU_2\to S$.
\item $(\bar p_1)^N:=\bar p_1\circ\epsilon_{12}:(\bar X_1\times_{\bar S}\bar X_2)^N\to\bar X_1$ is a compactification of 
$p_1:U_1\times_SU_2\to U_1$.
\item $(\bar p_2)^N:=\bar p_2\circ\epsilon_{12}:(\bar X_1\times_{\bar S}\bar X_2)^N\to\bar X_2$ is a compactification of 
$p_2:U_1\times_SU_2\to U_2$.
\end{itemize}
We have then the morphism in $C(\Var(k)/S)$
\begin{eqnarray*}
T(\otimes,R_S^{0CH})(\mathbb Z(U_1/S),\mathbb Z(U_2/S)):=R_S^{0CH}(p_1)\otimes R_S^{CH}(p_2): \\
R^0_{(\bar X_1,\bar D)/S}(\mathbb Z(U_1/S))\otimes R^0_{(X_2,E))/S}(\mathbb Z(U_2/S)) 
\xrightarrow{\sim} R^0_{(\bar X_1\times_{\bar S}\bar X_2)^N,\bar F)/S}(\mathbb Z(U_1\times_S U_2/S))
\end{eqnarray*}
For   
\begin{eqnarray*}
Q_1^*:=(\cdots\to\oplus_{\alpha\in\Lambda^n}\mathbb Z(U^n_{1,\alpha}/S)
\xrightarrow{(\mathbb Z(g^n_{\alpha,\beta}))}\oplus_{\beta\in\Lambda^{n-1}}\mathbb Z(U^{n-1}_{1,\beta}/S)\to\cdots), \\
Q_2^*:=(\cdots\to\oplus_{\alpha\in\Lambda^n}\mathbb Z(U^n_{2,\alpha}/S)
\xrightarrow{(\mathbb Z(g^n_{\alpha,\beta}))}\oplus_{\beta\in\Lambda^{n-1}}\mathbb Z(U^{n-1}_{2,\beta}/S)\to\cdots)
\in C(\Var(k)/S)
\end{eqnarray*}
complexes of (maybe infinite) direct sum of representable presheaves with $U^*_{\alpha}$ smooth, 
we get the morphism in $C(\Var(k)/S)$
\begin{eqnarray*}
T(\otimes,R_S^{0CH})(Q^*_1,Q^*_2): R^{0CH}(Q^*_1)\otimes R^{0CH}(Q^*_2) 
\xrightarrow{(T(\otimes,R_S^{0CH})(\mathbb Z(U_{1,\alpha}^m),\mathbb Z(U_{2,\beta}^n))} R^{0CH}(Q^*_1\otimes Q^*_2)).
\end{eqnarray*}
For $F_1,F_2\in C(\Var(k)^{sm}/S)$, we get in particular the morphism in $C(\Var(k)/S)$
\begin{eqnarray*}
T(\otimes,R_S^{0CH})(\rho_S^*LF_1,\rho_S^*LF_2):R^{0CH}(\rho_S^*LF_1)\otimes R^{0CH}(\rho_S^*LF_2) 
\to R^{0CH}(\rho_S^*(LF_1\otimes LF_2)).
\end{eqnarray*}

\end{itemize}

\begin{defi}\label{sharpstar0}
Let $h:U\to S$ a morphism, with $U,S\in\Var(k)$, $U$ irreducible.
Take, see definition-proposition \ref{RCHdef0},
$\bar f_0=\bar{h}_0:\bar X_0\to\bar S$ a compactification of $h:U\to S$ and denote by $\bar Z=\bar X_0\backslash U$.
Take, using theorem \ref{desVar}, a desingularization 
$\bar\epsilon:(\bar X,\bar D)\to(\bar X_0,\Delta)$ of the pair $(\bar X_0,\Delta)$, $\bar Z\subset\Delta$,
with $\bar X\in\PSmVar(k)$ and 
$\bar D:=\bar\epsilon^{-1}(\Delta)=\cup_{i=1}^s\bar D_i\subset\bar X$ a normal crossing divisor.
Denote $d_X:=\dim(\bar X)=\dim(U)$.
\begin{itemize}
\item[(i)] The diagonal $\Delta_{\bar D_{\bullet}}\subset\bar D_{\bullet}\times\bar D_{\bullet}$
induces the morphism in $C(\Var(k)/S)$
\begin{eqnarray*}
[\Delta_{\bar D_{\bullet}}]\in\Hom(\mathbb Z^{tr}(\bar D_{\bullet}/S),
\bar f_*E_{et}(\mathbb Z(\bar D_{\bullet}/\bar X)(d_X)[2d_X]))
\xrightarrow{\sim} \\
\Hom(\mathbb Z(\bar D_{\bullet}\times_S\bar X/\bar X),
\mathbb Z^{tr}(\bar D_{\bullet}\times\mathbb P^{d_X}/\bar X)/
\mathbb Z^{tr}(\bar D_{\bullet}\times\mathbb P^{d_X-1}/\bar X)) \\
\subset H^0(\mathcal Z_{d_{D_{\bullet}}}(\square^*\times\bar D_{\bullet}\times_S\bar D_{\bullet}))
\end{eqnarray*}
\item[(ii)] The cycle $\Delta_{\bar X}\subset\bar X\times_S\bar X$ induces by the morphism in $C(\Var(k)/S)$
\begin{eqnarray*}
[\Delta_{\bar X}]\in\Hom(\mathbb Z^{tr}(\bar X/S),
\bar f_*E_{et}(\mathbb Z(\bar X/\bar X)(d_X)[2d_X]))
\xrightarrow{\sim} \\
\Hom(\mathbb Z(\bar X\times_S\bar X/\bar X),
\mathbb Z^{tr}(\bar X\times\mathbb P^{d_X}/\bar X)/\mathbb Z^{tr}(\bar X\times\mathbb P^{d_X-1}/\bar X)) \\
\subset H^0(\mathcal Z_{d_X}(\square^*\times\bar X\times_S\bar X))
\end{eqnarray*}
\end{itemize}
Let $h:U\to S$ a morphism, with $U,S\in\Var(k)$, $U$ smooth connected (hence irreducible by smoothness).
Take, see definition-proposition \ref{RCHdef0},
$\bar f_0=\bar{h}_0:\bar X_0\to\bar S$ a compactification of $h:U\to S$ and denote by $\bar Z=\bar X_0\backslash U$.
Take, using theorem \ref{desVar}(ii), a strict desingularization 
$\bar\epsilon:(\bar X,\bar D)\to(\bar X_0,\bar Z)$ of the pair $(\bar X_0,\bar Z)$
with $\bar X\in\PSmVar(k)$ and 
$\bar D:=\bar\epsilon^{-1}(\bar Z)=\cup_{i=1}^s\bar D_i\subset\bar X$ a normal crossing divisor.
Denote $d_X:=\dim(\bar X)=\dim(U)$.
We get from (i) and (ii) the morphism in $C(\Var(k)/S)$
\begin{eqnarray*}
T(\bar f_{\sharp},\bar f_*)(\mathbb Z(D_{\bullet}/\bar X),\mathbb Z(\bar X/\bar X))
:=([\Delta_{\bar D_{\bullet}}],[\Delta_{\bar X}]): \\
\Cone(\mathbb Z(i_{\bullet}):(\mathbb Z^{tr}(\bar D_{\bullet}/S),u_{IJ})\to\mathbb Z^{tr}(\bar X/S))\to \\
\bar f_*E_{et}(\Cone(\mathbb Z(i_{\bullet}):
(\mathbb Z(\bar D_{\bullet}/\bar X),u_{IJ})\to\mathbb Z(\bar X/\bar X)))(d_X)[2d_X] \\
=:R^0_{(\bar X,\bar D)/S}(\mathbb Z(U/S))(d_X)[2d_X].
\end{eqnarray*}
\end{defi}

\begin{defi}\label{R0CHhatdef} 
\begin{itemize}
\item[(i)]Let $h:U\to S$ a morphism, with $U,S\in\Var(k)$ and $U$ smooth.
Take, see definition-proposition \ref{RCHdef0},
$\bar f_0=\bar{h}_0:\bar X_0\to\bar S$ a compactification of $h:U\to S$ and denote by $\bar Z=\bar X_0\backslash U$.
Take, using theorem \ref{desVar}(ii), a strict desingularization 
$\bar\epsilon:(\bar X,\bar D)\to(\bar X_0,\bar Z)$ of the pair $(\bar X_0,\bar Z)$, with $\bar X\in\PSmVar(k)$ and 
$\bar D:=\epsilon^{-1}(\bar Z)=\cup_{i=1}^s\bar D_i\subset\bar X$ a normal crossing divisor.  
We denote by $i_{\bullet}:\bar D_{\bullet}\hookrightarrow\bar X=\bar X_{c(\bullet)}$ the morphism of simplicial varieties
given by the closed embeddings $i_I:\bar D_I=\cap_{i\in I}\bar D_i\hookrightarrow\bar X$
We denote by $j:U\hookrightarrow\bar X$ the open embedding. We then consider the map in $C(\Var(k)/S)$
\begin{eqnarray*}
T(\hat R^{0CH},R^{0CH})(\mathbb Z(U/S)):\hat R^0_{(\bar X,\bar D)/S}(\mathbb Z(U/S)) \\
\xrightarrow{:=}\Cone(\mathbb Z(i_{\bullet}):
(\mathbb Z^{tr}(D_{\bullet}/S),u_{IJ})\to\mathbb Z^{tr}(X/S))(-d_X)[-2d_X] \\
\xrightarrow{T(\bar f_{\sharp},\bar f_*)(\mathbb Z(\bar D_{\bullet}/\bar X),\mathbb Z(\bar X/\bar X))(-d_X)[-2d_X]} \\
R^0_{(\bar X,\bar D)/S}(\mathbb Z(U/S)).
\end{eqnarray*}
given in definition \ref{sharpstar}(iii).
\item[(ii)]Let $g:U'/S\to U/S$ a morphism, with $U'/S=(U',h'),U/S=(U,h)\in\Var(k)/S$, with $U$ and $U'$ smooth.
Take, see definition-proposition \ref{RCHdef0}(ii),a compactification $\bar f_0=\bar h:\bar X_0\to\bar S$ of $h:U\to S$ 
and a compactification $\bar f'_0=\bar{h}':\bar X'_0\to\bar S$ of $h':U'\to S$ such that 
$g:U'/S\to U/S$ extend to a morphism $\bar g_0:\bar X'_0/\bar S\to\bar X_0/\bar S$. 
Denote $\bar Z=\bar X_0\backslash U$ and $\bar Z'=\bar X'_0\backslash U'$.
Take, see definition-proposition \ref{RCHdef0}(ii), a strict desingularization 
$\bar\epsilon:(\bar X,\bar D)\to(\bar X_0,\bar Z)$ of $(\bar X_0,\bar Z)$,
a strict desingularization $\bar\epsilon'_{\bullet}:(\bar X',\bar D')\to(\bar X'_0,\bar Z')$ of $(\bar X'_0,\bar Z')$
and a morphism $\bar g:\bar X'\to\bar X$ such that the following diagram commutes
\begin{equation*}
\xymatrix{\bar X'_0\ar[r]^{\bar{g}_0} & \bar X_0 \\
\bar X'\ar[u]^{\bar\epsilon'}\ar[r]^{\bar g} & \bar X\ar[u]^{\bar\epsilon}}.
\end{equation*} 
We then have, see definition-proposition \ref{RCHdef0}(ii),
the diagram (\ref{RCHdia}) in $\Fun(\Delta,\Var(k))$
\begin{equation*}
\xymatrix{U=U_{c(\bullet)}\ar[r]^j & \bar X=\bar X_{c(\bullet)} & \, & \bar D_{s_g(\bullet)}\ar[ll]_{i_{\bullet}} \\
U'=U'_{c(\bullet)}\ar[r]^{j'}\ar[u]^g & \bar X'=\bar X'_{c(\bullet)}\ar[u]^{\bar{g}} & \bar D'_{\bullet}\ar[l]_{i'_{\bullet}} & 
\bar{g}^{-1}(\bar D_{s_g(\bullet)})\ar[l]_{i''_{g\bullet}}\ar[u]^{\bar{g}'_{\bullet}}:i'_{g\bullet}}
\end{equation*}
Consider 
\begin{eqnarray*}
[\Gamma_{\bar g}]^t\in\Hom(\mathbb Z^{tr}(\bar X/S)(-d_X)[-2d_X],
\mathbb Z^{tr}(\bar X'/S)(-d_{X'})[-2d_{X'}]) \\ 
\xrightarrow{\sim}\Hom(\mathbb Z^{tr}(\bar X\times\mathbb P^{d_X}/S)/\mathbb Z_{tr}(\bar X\times\mathbb P^{d_X-1}/S), \\
\mathbb Z_{tr}(\bar X'\times\mathbb P^{d_{X'}}/S)/\mathbb Z_{tr}(\bar X'\times\mathbb P^{d_{X'}-1}/S)
\end{eqnarray*}
the morphism given by the transpose of the graph $\Gamma_g\subset X'\times_S X$ of $\bar g:\bar X'\to\bar X$. Then,  
since $i_{\bullet}\circ\bar g'_{\bullet}=\bar g\circ i''_{g\bullet}=\bar g\circ i'\circ\circ i'_{g\bullet}$, 
we have the factorization 
\begin{eqnarray*}
[\Gamma_g]^t\circ\mathbb Z(i_{\bullet}): 
(\mathbb Z^{tr}(\bar D_{s_g(\bullet)}/S),u_{IJ})(-d_X)[-2d_X] \\
\xrightarrow{[\Gamma_{\bar g'_{\bullet}}]^t}
(\mathbb Z^{tr}(\bar{g}^{-1}(\bar D_{s_g(\bullet)})/S),u_{IJ})(-d_{X'})[-2d_{X'}] \\
\xrightarrow{\mathbb Z(i'_{g\bullet})} 
\mathbb Z^{tr}(\bar X'/S)(-d_{X'})[-2d_{X'}].
\end{eqnarray*}
with
\begin{eqnarray*}
[\Gamma_{\bar g'_{\bullet}}]^t\in
\Hom((\mathbb Z^{tr}(\bar D_{s_g(\bullet)}\times\mathbb P^{d_X}/S),u_{IJ})/
(\mathbb Z^{tr}(\bar D_{s_g(\bullet)}\times\mathbb P^{d_{X-1}}/S),u_{IJ}), \\
(\mathbb Z_{tr}(\bar{g}^{-1}(\bar D_{s_g(\bullet)})\times\mathbb P^{d_{X'}}/S),u_{IJ})/
(\mathbb Z_{tr}(\bar{g}^{-1}(\bar D_{s_g(\bullet)})\times\mathbb P^{d_{X'-1}}/S),u_{IJ})).
\end{eqnarray*}
We then consider the following map in $C(\Var(k)/S)$
\begin{eqnarray*}
\hat R_S^{0CH}(g):\hat R^0_{(\bar X,\bar D)/S}(\mathbb Z(U/S))\xrightarrow{:=} \\
\Cone(\mathbb Z(i_{\bullet}):
(\mathbb Z^{tr}(\bar D_{s_g(\bullet)}/S),u_{IJ})\to\mathbb Z^{tr}(\bar X/S)(-d_X)[-2d_X] \\
\xrightarrow{([\Gamma_{\bar g'_{\bullet}}]^t,[\Gamma_{\bar g}]^t)} \\
\Cone(\mathbb Z(i'_{g\bullet}): 
(\mathbb Z^{tr}(\bar{g}^{-1}(\bar D_{s_g(\bullet)})/S),u_{IJ})\to\mathbb Z^{tr}(\bar X'/S))(-d_{X'})[-2d_{X'}] \\
\xrightarrow{(\mathbb Z(i''_{g\bullet}),I)(-d_{X'})[-2d_{X'}]} \\
\Cone(\mathbb Z(i'_{\bullet}):
((\mathbb Z^{tr}(\bar D'_{\bullet}/S),u_{IJ})\to\mathbb Z^{tr}(\bar X'/S))(-d_{X'})[-2d_{X'}] \\
\xrightarrow{=:}\hat R^0_{(\bar X',\bar D')/S}(\mathbb Z(U'/S))
\end{eqnarray*}
Then the following diagram in $C(\Var(k)/S)$ commutes by definition
\begin{equation*}
\xymatrix{\hat R^0_{(\bar X,\bar D)/S}(\mathbb Z(U/S))
\ar[rr]^{T(\hat R^{0CH},R^{0CH})(\mathbb Z(U/S))}\ar[d]_{\hat R_S^{0CH}(g)} & \, & 
R^0_{(\bar X,\bar D)/S}(\mathbb Z(U/S))\ar[d]^{R_S^{0CH}(g)} \\
\hat R^0_{(\bar X',\bar D')/S}(\mathbb Z(U'/S))\ar[rr]^{T(\hat R^{0CH},R^{0CH})(\mathbb Z(U'/S))} & \, & 
R^0_{(\bar X',\bar D')/S}(\mathbb Z(U'/S))}.
\end{equation*}
\item[(iii)] For $g_1:U''/S\to U'/S$, $g_2:U'/S\to U/S$ two morphisms
with $U''/S=(U',h''),U'/S=(U',h'),U/S=(U,h)\in\Var(k)/S$, with $U$, $U'$ and $U''$ smooth.
We get from (i) and (ii) 
a compactification $\bar f=\bar{h}:\bar X\to\bar S$ of $h:U\to S$,
a compactification $\bar f'=\bar{h}':\bar X'\to\bar S$ of $h':U'\to S$,
and a compactification $\bar f''=\bar{h}'':\bar X''\to\bar S$ of $h'':U''\to S$,
with $\bar X,\bar X',\bar X''\in\PSmVar(k)$, 
$\bar D:=\bar X\backslash U\subset\bar X$ $\bar D':=\bar X'\backslash U'\subset\bar X'$, 
and $\bar D'':=\bar X''\backslash U''\subset\bar X''$ normal crossing divisors, 
such that $g_1:U''/S\to U'/S$ extend to $\bar g_1:\bar X''/\bar S\to\bar X'/\bar S$,
$g_2:U'/S\to U/S$ extend to $\bar g_2:\bar X'/\bar S\to\bar X/\bar S$, and
\begin{eqnarray*}
\hat R_S^{0CH}(g_2\circ g_1)=\hat R_S^{0CH}(g_1)\circ \hat R_S^{0CH}(g_2):
\hat R^0_{(\bar X,\bar D)/S}\to \hat R^0_{(\bar X'',\bar D'')/S} 
\end{eqnarray*}
\item[(iv)] For  
\begin{eqnarray*}
Q^*:=(\cdots\to\oplus_{\alpha\in\Lambda^n}\mathbb Z(U^n_{\alpha}/S)
\xrightarrow{(\mathbb Z(g^n_{\alpha,\beta}))}\oplus_{\beta\in\Lambda^{n-1}}\mathbb Z(U^{n-1}_{\beta}/S)\to\cdots)
\in C(\Var(k)/S)
\end{eqnarray*}
a complex of (maybe infinite) direct sum of representable presheaves with $U^*_{\alpha}$ smooth,
we get from (i),(ii) and (iii) the map in $C(\Var(k)/S)$
\begin{eqnarray*}
T(\hat R^{0CH},R^{0CH})(Q^*):\hat R^{0CH}(Q^*):=
(\cdots\to\oplus_{\beta\in\Lambda^{n-1}}\varinjlim_{(\bar X^{n-1}_{\beta},\bar D^{n-1}_{\beta})/S}
\hat R^0_{(\bar X^{n-1}_{\beta},\bar D^{n-1}_{\beta})/S}(\mathbb Z(U^{n-1}_{\beta}/S)) \\
\xrightarrow{(\hat R_S^{0CH}(g^n_{\alpha,\beta}))}\oplus_{\alpha\in\Lambda^n}\varinjlim_{(\bar X^n_{\alpha},\bar D^n_{\alpha})/S}
\hat R^0_{(\bar X^n_{\alpha},\bar D^n_{\alpha})/S}(\mathbb Z(U^n_{\alpha}/S))\to\cdots) 
\to R^{0CH}(Q^*),
\end{eqnarray*}
where for $(U^n_{\alpha},h^n_{\alpha})\in\Var(k)/S$, the inductive limit run over all the compactifications 
$\bar f_{\alpha}:\bar X_{\alpha}\to\bar S$ of $h_{\alpha}:U_{\alpha}\to S$ with $\bar X_{\alpha}\in\PSmVar(k)$
and $\bar D_{\alpha}:=\bar X_{\alpha}\backslash U_{\alpha}$ a normal crossing divisor.
For $m=(m^*):Q_1^*\to Q_2^*$ a morphism with 
\begin{eqnarray*}
Q_1^*:=(\cdots\to\oplus_{\alpha\in\Lambda^n}\mathbb Z(U^n_{1,\alpha}/S)
\xrightarrow{(\mathbb Z(g^n_{\alpha,\beta}))}\oplus_{\beta\in\Lambda^{n-1}}\mathbb Z(U^{n-1}_{1,\beta}/S)\to\cdots), \\
Q_2^*:=(\cdots\to\oplus_{\alpha\in\Lambda^n}\mathbb Z(U^n_{2,\alpha}/S)
\xrightarrow{(\mathbb Z(g^n_{\alpha,\beta}))}\oplus_{\beta\in\Lambda^{n-1}}\mathbb Z(U^{n-1}_{2,\beta}/S)\to\cdots)
\in C(\Var(k)/S)
\end{eqnarray*}
complexes of (maybe infinite) direct sum of representable presheaves with $U^*_{1,\alpha}$ and $U^*_{2,\alpha}$ smooth,
we get again from (i),(ii) and (iii) a commutative diagram in $C(\Var(k)/S)$
\begin{equation*}
\xymatrix{\hat R^{0CH}(Q_2^*)\ar[rr]^{T(\hat R_S^{0CH},R_S^{0CH})(Q_2^*)}\ar[d]_{\hat R_S^{0CH}(m):=(\hat R_S^{0CH}(m^*))} & \, & 
R^{0CH}(Q_2^*)\ar[d]^{(R_S^{0CH}(m^*))} \\
\hat R^{0CH}(Q_1^*)\ar[rr]^{T(\hat R_S^{0CH},R_S^{0CH})(Q_1^*)} & \, & R^{0CH}(Q_1^*)}.
\end{equation*}
\item[(v)] Let  
\begin{eqnarray*}
Q^*:=(\cdots\to\oplus_{\alpha\in\Lambda^n}\mathbb Z(U^n_{\alpha}/S)
\xrightarrow{(\mathbb Z(g^n_{\alpha,\beta}))}\oplus_{\beta\in\Lambda^{n-1}}\mathbb Z(U^{n-1}_{\beta}/S)\to\cdots)
\in C(\Var(k)/S)
\end{eqnarray*}
a complex of (maybe infinite) direct sum of representable presheaves with $U^*_{\alpha}$ smooth,
we have by definition 
\begin{equation*}
\Gr_S^{12*}\hat R^{0CH}(Q^*)=\hat R^{CH}(Q^*)\in C(\Var(k)^{2,smpr}/S).
\end{equation*}
\end{itemize}
\end{defi}

\begin{itemize}
\item Let $S\in\Var(k)$
For $(h,m,m')=(h^*,m^*,m^{'*}):Q_1^*[1]\to Q_2^*$ an homotopy with $Q_1^*,Q_2^*\in C(\Var(k)/S)$
complexes of (maybe infinite) direct sum of representable presheaves with $U^*_{1,\alpha}$ and $U^*_{2,\alpha}$ smooth,
\begin{equation*}
(\hat R_S^{0CH}(h),\hat R_S^{0CH}(m),\hat R_S^{0CH}(m'))=(\hat R_S^{0CH}(h^*),\hat R_S^{0CH}(m^*),\hat R_S^{0CH}(m^{'*})):
R^{0CH}(Q_2^*)[1]\to R^{0CH}(Q_1^*)
\end{equation*}
is an homotopy in $C(\Var(k)/S)$ using definition \ref{R0CHhatdef} (iii). 
In particular if $m:Q_1^*\to Q_2^*$ with $Q_1^*,Q_2^*\in C(\Var(k)/S)$
complexes of (maybe infinite) direct sum of representable presheaves with $U^*_{1,\alpha}$ and $U^*_{2,\alpha}$ smooth
is an homotopy equivalence, then $\hat R_S^{0CH}(m):\hat R^{0CH}(Q_2^*)\to\hat R^{0CH}(Q_1^*)$ is an homotopy equivalence.
\item Let $S\in\SmVar(k)$. Let $F\in\PSh(\Var(k)^{sm}/S)$. Consider 
\begin{eqnarray*}
q:LF:=(\cdots\to\oplus_{(U_{\alpha},h_{\alpha})\in\Var(k)^{sm}/S}\mathbb Z(U_{\alpha}/S)
\xrightarrow{(\mathbb Z(g^n_{\alpha,\beta}))}
\oplus_{(U_{\alpha},h_{\alpha})\in\Var(k)^{sm}/S}\mathbb Z(U_{\alpha}/S)\to\cdots)\to F
\end{eqnarray*}
the canonical projective resolution given in subsection 2.3.3.
Note that the $U_{\alpha}$ are smooth since $S$ is smooth and $h_{\alpha}$ are smooth morphism.
Definition \ref{R0CHhatdef}(iv) gives in this particular case the map in $C(\Var(k)/S)$
\begin{eqnarray*}
T(\hat R_S^{0CH},R_S^{0CH})(\rho_S^*LF):\hat R^{0CH}(\rho_S^*LF):=
(\cdots\to\oplus_{(U_{\alpha},h_{\alpha})\in\Var(k)^{sm}/S}\varinjlim_{(\bar X_{\alpha},\bar D_{\alpha})/S}
\hat R^0_{(\bar X_{\alpha},\bar D_{\alpha})/S}(\mathbb Z(U_{\alpha}/S)) \\
\xrightarrow{(\hat R_S^{0CH}(g^n_{\alpha,\beta}))}
\oplus_{(U_{\alpha},h_{\alpha})\in\Var(k)^{sm}/S}\varinjlim_{(\bar X_{\alpha},\bar D_{\alpha})/S}
\hat R^0_{(\bar X_{\alpha},\bar D_{\alpha})/S}(\mathbb Z(U_{\alpha}/S))\to\cdots)
\to R^{0CH}(\rho_S^*LF),
\end{eqnarray*}
where for $(U_{\alpha},h_{\alpha})\in\Var(k)^{sm}/S$, the inductive limit run over all the compactifications 
$\bar f_{\alpha}:\bar X_{\alpha}\to\bar S$ of $h_{\alpha}:U_{\alpha}\to S$ with $\bar X_{\alpha}\in\PSmVar(k)$
and $\bar D_{\alpha}:=\bar X_{\alpha}\backslash U_{\alpha}$ a normal crossing divisor.
Definition \ref{R0CHhatdef}(iv) gives then by functoriality in particular, for $F=F^{\bullet}\in C(\Var(k)^{sm}/S)$, 
the map in $C(\Var(k)/S)$
\begin{eqnarray*}
T(\hat R_S^{0CH},R_S^{0CH})(\rho_S^*LF):\hat R^{0CH}(\rho_S^*LF)\to R^{0CH}(\rho_S^*LF).
\end{eqnarray*}

\item Let $g:T\to S$ a morphism with $T,S\in\SmVar(k)$. Let $h:U\to S$ a smooth morphism with $U\in\Var(k)$.
Consider the cartesian square
\begin{equation*}
\xymatrix{U_T\ar[r]^{h'}\ar[d]^{g'} & T\ar[d]^g \\
U\ar[r]^h & S}
\end{equation*}
Note that $U$ is smooth since $S$ and $h$ are smooth, and $U_T$ is smooth since $T$ and $h'$ are smooth. 
Take, see definition-proposition \ref{RCHdef0}(ii),a compactification $\bar f_0=\bar h:\bar X_0\to\bar S$ of $h:U\to S$. 
Take, see definition-proposition \ref{RCHdef0}(ii), a strict desingularization 
$\bar\epsilon:(\bar X,\bar D)\to(\bar X_0,\bar Z)$ of $(\bar X_0,\bar Z)$.
Then $\bar f'_0=\bar{g\circ h'}:\bar X_T\to\bar T$ is a compactification of $g\circ h':U_T\to S$ such that 
$g':U_T/S\to U/S$ extend to a morphism $\bar{g}'_0:\bar X_T/\bar S\to \bar X/\bar S$. 
Denote $\bar Z=\bar X_0\backslash U$ and $\bar Z'=\bar X_T\backslash U_T$.
Take, see definition-proposition \ref{RCHdef0}(ii), a strict desingularization 
$\epsilon'_{\bullet}:(\bar X',\bar D')\to(\bar X_T,\bar Z')$ of $(\bar X_T,\bar Z')$.
Denote $\bar g'=\bar g'_0\circ\epsilon'_{\bullet}:\bar X'\to\bar X$. 
We then have, see definition-proposition \ref{RCHdef0}(ii),
the following commutative diagram in $\Fun(\Delta,\Var(k))$
\begin{equation*}
\xymatrix{U=U_{c(\bullet)}\ar[r]^j & \bar X=\bar X_{c(\bullet)} & \, & \bar D_{s_{g'}(\bullet)}\ar[ll]_{i_{\bullet}} \\
U_T=U_{T,c(\bullet)}\ar[r]^{j'}\ar[u]^{g'} & \bar X_T=X_{T,c(\bullet)}\ar[u]^{\bar{g}'_0} & \, & 
\bar{g}^{'-1}(\bar D_{s_{g'}(\bullet)})\ar[ll]_{i_{g\bullet}}\ar[u]^{(\bar{g}')'_{\bullet}} \\
U_T=U_{T,c(\bullet)}\ar[r]^{j'}\ar[u]^{g'} & \bar X'=X'_{c(\bullet)}\ar[u]^{\bar{g}'} & \bar D'_{\bullet}\ar[l]_{i'_{\bullet}} & 
\bar{g}^{'-1}(\bar D_{s_{g'}(\bullet)})\ar[l]_{i''_{g'\bullet}}\ar[u]^{\epsilon'_{\bullet}}:i'_{g\bullet}}
\end{equation*}
We then consider the following map in $C(\Var(k)/T)$, 
\begin{eqnarray*}
T(g,\hat R^{0CH})(\mathbb Z(U/S)):g^*\hat R^0_{(\bar X,\bar D)/S}(\mathbb Z(U/S)) \\
\xrightarrow{:=}g^*\Cone(\mathbb Z(i_{\bullet}):
(\mathbb Z^{tr}(\bar D_{\bullet}/S),u_{IJ})\to\mathbb Z^{tr}(\bar X/S))(-d_X)[-2d_X] \\ 
\xrightarrow{=} \\
\Cone(\mathbb Z(i_{g\bullet}):(\mathbb Z^{tr}(\bar g^{-1}(D_{s_g(\bullet)})/T),u_{IJ})
\to\mathbb Z^{tr}(\bar X_T/T))(-d_X)[-2d_X] \\
\xrightarrow{(\mathbb Z(i''_{g\bullet}),[\Gamma_{\epsilon'}]^t)} \\
\Cone(\mathbb Z(i'_{\bullet}):
((\mathbb Z^{tr}(\bar D'_{\bullet}/T),u_{IJ})\to\mathbb Z^{tr}((\bar X'/T)))(-d_{X'})[-2d_{X'}] \\
\xrightarrow{=:}\hat R^0_{(\bar X',\bar D')/T}(\mathbb Z(U_T/T))
\end{eqnarray*}
For  
\begin{eqnarray*}
Q^*:=(\cdots\to\oplus_{\alpha\in\Lambda^n}\mathbb Z(U^n_{\alpha}/S)
\xrightarrow{(\mathbb Z(g^n_{\alpha,\beta}))}\oplus_{\beta\in\Lambda^{n-1}}\mathbb Z(U^{n-1}_{\beta}/S)\to\cdots)
\in C(\Var(k)/S)
\end{eqnarray*}
a complex of (maybe infinite) direct sum of representable presheaves with $h^n_{\alpha}:U^n_{\alpha}\to S$ smooth,
we get the map in $C(\Var(k)/T)$
\begin{eqnarray*}
T(g,\hat R^{0CH})(Q^*):g^*\hat R^{0CH}(Q^*)=
(\cdots\to\oplus_{\alpha\in\Lambda^n}\varinjlim_{(\bar X^n_{\alpha},\bar D^n_{\alpha})/S}
g^*\hat R^0_{(\bar X^n_{\alpha},\bar D^n_{\alpha})/S}(\mathbb Z(U^n_{\alpha}/S))\to\cdots) \\
\xrightarrow{(T(g,\hat R^{0CH})(\mathbb Z(U^n_{\alpha}/S)))} 
(\cdots\to\oplus_{\alpha\in\Lambda^n}\varinjlim_{(\bar X^{n'}_{\alpha},\bar D^{n'}_{\alpha})/T}
\hat R_{(\bar X^{n'}_{\alpha},\bar D^{n'}_{\alpha})/T}(\mathbb Z(U^n_{\alpha,T}/S))\to\cdots)=:\hat R^{CH}(g^*Q^*)
\end{eqnarray*}
together with the commutative diagram in $C(\Var(k)/T)$
\begin{equation*}
\xymatrix{g^*\hat R^{0CH}(Q^*)\ar[rrr]^{T(g,\hat R^{0CH})(Q^*)}\ar[d]_{g^*T(\hat R_S^{0CH},R_S^{0CH})(Q^*)} 
& \, & \, & \hat R^{0CH}(g^*Q^*)\ar[d]^{T(\hat R_T^{0CH},R_T^{0CH})(g^*Q)} \\
g^*R^{0CH}(Q^*)\ar[rrr]^{T(g,R^{0CH})(Q^*)} & \, & \, & R^{0CH}(g^*Q^*)}.
\end{equation*}
Let $F\in\PSh(\Var(k)^{sm}/S)$. Consider 
\begin{eqnarray*}
q:LF:=(\cdots\to\oplus_{(U_{\alpha},h_{\alpha})\in\Var(k)^{sm}/S}\mathbb Z(U_{\alpha}/S)\to\cdots)\to F
\end{eqnarray*}
the canonical projective resolution given in subsection 2.3.3.
We then get in particular the map in $C(\Var(k)/T)$
\begin{eqnarray*}
T(g,\hat R^{0CH})(\rho_S^*LF):g^*\hat R^{0CH}(\rho_S^*LF)= \\
(\cdots\to\oplus_{(U_{\alpha},h_{\alpha})\in\Var(k)^{sm}/S}\varinjlim_{(\bar X_{\alpha},\bar D_{\alpha})/S}
g^*\hat R_{(\bar X_{\alpha},\bar D_{\alpha})/S}(\mathbb Z(U_{\alpha}/S))\to\cdots) 
\xrightarrow{(T(g,\hat R^{0CH})(\mathbb Z(U_{\alpha}/S)))} \\ 
(\cdots\to\oplus_{(U_{\alpha},h_{\alpha})\in\Var(k)^{sm}/S}\varinjlim_{(\bar X'_{\alpha},\bar D'_{\alpha})/T}
\hat R_{(\bar X'_{\alpha},\bar D'_{\alpha})/T}(\mathbb Z(U_{\alpha,T}/S))\to\cdots)=:\hat R^{0CH}(\rho_T^*g^*LF), 
\end{eqnarray*}
and by functoriality, we get in particular for $F=F^{\bullet}\in C(\Var(k)^{sm}/S)$, 
the map in $C(\Var(k)/T)$
\begin{eqnarray*}
T(g,\hat R^{0CH})(\rho_S^*LF):g^*\hat R^{0CH}(\rho_S^*LF)\to\hat R^{0CH}(\rho_T^*g^*LF)
\end{eqnarray*}
together with the commutative diagram in $C(\Var(k)/T)$
\begin{equation*}
\xymatrix{g^*\hat R^{0CH}(\rho_S^*LF)\ar[rrr]^{T(g,\hat R^{0CH})(\rho_S^*LF)}\ar[d]_{g^*T(\hat R_S^{0CH},R_S^{0CH})(\rho_S^*LF)} 
& \, & \, & \hat R^{0CH}(\rho_T^*g^*LF)\ar[d]^{T(\hat R_T^{0CH},R_T^{0CH})(\rho_T^*g^*LF)} \\
g^*L\rho_{S*}\mu_{S*}R^{CH}(\rho_S^*LF)
\ar[rrr]^{T(g,R^{0CH})(\rho_S^*LF)} & \, & \, & R^{CH}(\rho_T^*g^*LF)}.
\end{equation*}

\item Let $S_1,S_2\in\SmVar(k)$ and $p:S_1\times S_2\to S_1$ the projection. 
Let $h:U\to S_1$ a smooth morphism with $U\in\Var(k)$. Consider the cartesian square
\begin{equation*}
\xymatrix{U\times S_2\ar[r]^{h\times I}\ar[d]^{p'} & S_1\times S_2\ar[d]^p \\
U\ar[r]^h & S_1}
\end{equation*}
Take, see definition-proposition \ref{RCHdef0}(i),a compactification $\bar f_0=\bar h:\bar X_0\to\bar S_1$ of $h:U\to S_1$.
Then $\bar f_0\times I:\bar X_0\times S_2\to\bar S_1\times S_2$ is a compactification of $h\times I:U\times S_2\to S_1\times S_2$ 
and $p':U\times S_2\to U$ extend to $\bar{p}'_0:=p_{X_0}:\bar X_0\times S_2\to\bar X_0$. Denote $Z=X_0\backslash U$. 
Take see theorem \ref{desVar}(i), a strict desingularization 
$\bar\epsilon:(\bar X,\bar D)\to(\bar X_0,\bar Z)$ of the pair $(\bar X_0,\bar Z)$.  
We then have the commutative diagram (\ref{RCHdiap}) in $\Fun(\Delta,\Var(k))$ whose squares are cartesian
\begin{equation*}
\xymatrix{U=U_{c(\bullet)}\ar[r]^j & \bar X & \bar D_{\bullet}\ar[l]_{i_{\bullet}} \\
U\times S_2=(U\times S_2)_{c(\bullet)}\ar[r]^{j\times I}\ar[u]^g & \bar X\times S_2\ar[u]^{\bar{p}':=p_{\bar X}} & 
\bar D_{\bullet}\times S_2\ar[l]_{i'_{\bullet}}\ar[u]^{\bar{p'}'_{\bullet}}}
\end{equation*}
Then the map in $C(\Var(k)/S_1\times S_2)$
\begin{eqnarray*}
T(p,\hat R^{0CH})(\mathbb Z(U/S_1)):p^*\hat R^0_{(\bar X,\bar D)/S_1}(\mathbb Z(U/S_1))\xrightarrow{\sim} 
\hat R^0_{(\bar X\times S_2,\bar D_{\bullet}\times S_2)/S_1\times S_2}(\mathbb Z(U\times S_2/S_1\times S_2))
\end{eqnarray*}
is an isomorphism.
Hence, for $Q^*\in C(\Var(k)/S_1)$ a complex of (maybe infinite) direct sum of representable presheaves of smooth morphism,
the map in $C(\Var(k)/S_1\times S_2)$
\begin{eqnarray*}
T(p,\hat R^{0CH})(Q^*):p^*\hat R^{0CH}(Q^*)\xrightarrow{\sim}\hat R^{0CH}(p^*Q^*)
\end{eqnarray*}
is an isomorphism.
In particular, for $F\in C(\Var(k)^{sm}/S_1)$ the map in $C(\Var(k)/S_1\times S_2)$
\begin{eqnarray*}
T(p,\hat R^{0CH})(\rho_{S_1}^*LF):p^*\hat R^{0CH}(\rho_{S_1}^*LF)\xrightarrow{\sim}\hat R^{0CH}(\rho_{S_1\times S_2}^*p^*LF)
\end{eqnarray*}
is an isomorphism.

\item Let $h_1:U_1\to S$, $h_2:U_2\to S$ two morphisms with $U_1,U_2,S\in\Var(k)$, $U_1,U_2$ smooth.
Denote by $p_1:U_1\times_SU_2\to U_1$ and $p_2:U_1\times_S U_2\to U_2$ the projections.
Take, see definition-proposition \ref{RCHdef0}(i)),
a compactification $\bar f_{10}=\bar{h}_1:\bar X_{10}\to\bar S$ of $h_1:U_1\to S$
and a compactification $\bar f_{20}=\bar{h}_2:\bar X_{20}\to\bar S$ of $h_2:U_2\to S$. Then, 
\begin{itemize}
\item $\bar f_{10}\times\bar f_{20}:\bar X_{10}\times_{\bar S}\bar X_{20}\to S$ 
is a compactification of $h_1\times h_2:U_1\times_SU_2\to S$.
\item $\bar p_{10}:=p_{X_{10}}:\bar X_{10}\times_{\bar S} \bar X_{20}\to\bar X_{10}$ 
is a compactification of $p_1:U_1\times_SU_2\to U_1$.
\item $\bar p_{20}:=p_{X_{20}}:\bar X_{10}\times_{\bar S}\bar X_{20}\to\bar X_{20}$ 
is a compactification of $p_2:U_1\times_SU_2\to U_2$.
\end{itemize}
Denote $\bar Z_1=\bar X_{10}\backslash U_1$ and $\bar Z_2=\bar X_{20}\backslash U_2$.
Take, see theorem \ref{desVar}(i), a strict desingularization 
$\bar\epsilon_1:(\bar X_1,\bar D)\to(\bar X_{10},Z_1)$ of the pair $(\bar X_{10},\bar Z_1)$
and a  strictdesingularization 
$\bar\epsilon_2:(\bar X_2,\bar E)\to(\bar X_{20},Z_2)$ of the pair $(\bar X_{20},\bar Z_2)$. 
Take then a strict desingularization
\begin{equation*}
\bar\epsilon_{12}:((\bar X_1\times_{\bar S}\bar X_2)^N,\bar F)\to
(\bar X_1\times_{\bar S}\bar X_2,(D\times_{\bar S}\bar X_2)\cup(\bar X_1\times_{\bar S}\bar E)) 
\end{equation*}
of the pair $(\bar X_1\times_{\bar S}\bar X_2,(\bar D\times_{\bar S}\bar X_2)\cup(\bar X_1\times_{\bar S}\bar E))$. 
We have then the following commutative diagram 
\begin{equation*}
\xymatrix{\, & \bar X_1\ar[r]^{\bar f_1} & \bar S \\
\, & \bar X_1\times_{\bar S}\bar X_2\ar[r]^{\bar p_1}\ar[u]^{\bar p_2} &\bar X_2\ar[u]^{\bar f_2} \\
(\bar X_1\times_{\bar S}\bar X_2)^N\ar[ru]^{\bar\epsilon_{12}}\ar[rru]^{(\bar p_1)^N}\ar[ruu]^{(\bar p_2)^N} & \, & \,}
\end{equation*}
and
\begin{itemize}
\item $\bar f_1\times\bar f_2:\bar X_1\times_{\bar S}\bar X_2\to\bar S$ 
is a compactification of $h_1\times h_2:U_1\times_SU_2\to S$.
\item $(\bar p_1)^N:=\bar p_1\circ\epsilon_{12}:(\bar X_1\times_{\bar S}\bar X_2)^N\to\bar X_1$ is a compactification of 
$p_1:U_1\times_SU_2\to U_1$.
\item $(\bar p_2)^N:=\bar p_2\circ\epsilon_{12}:(\bar X_1\times_{\bar S}\bar X_2)^N\to\bar X_2$ is a compactification of 
$p_2:U_1\times_SU_2\to U_2$.
\end{itemize}
We have then the morphism in $C(\Var(k)/S)$
\begin{eqnarray*}
T(\otimes,\hat R_S^{0CH})(\mathbb Z(U_1/S),\mathbb Z(U_2/S)):=\hat R_S^{0CH}(p_1)\otimes\hat R_S^{0CH}(p_2): \\
\hat R^0_{(\bar X_1,\bar D)/S}(\mathbb Z(U_1/S))\otimes\hat R^0_{(X_2,E))/S}(\mathbb Z(U_2/S)) 
\xrightarrow{\sim}\hat R^0_{(\bar X_1\times_{\bar S}\bar X_2)^N,\bar F)/S}(\mathbb Z(U_1\times_S U_2/S))
\end{eqnarray*}
For   
\begin{eqnarray*}
Q_1^*:=(\cdots\to\oplus_{\alpha\in\Lambda^n}\mathbb Z(U^n_{1,\alpha}/S)
\xrightarrow{(\mathbb Z(g^n_{\alpha,\beta}))}\oplus_{\beta\in\Lambda^{n-1}}\mathbb Z(U^{n-1}_{1,\beta}/S)\to\cdots), \\
Q_2^*:=(\cdots\to\oplus_{\alpha\in\Lambda^n}\mathbb Z(U^n_{2,\alpha}/S)
\xrightarrow{(\mathbb Z(g^n_{\alpha,\beta}))}\oplus_{\beta\in\Lambda^{n-1}}\mathbb Z(U^{n-1}_{2,\beta}/S)\to\cdots)
\in C(\Var(k)/S)
\end{eqnarray*}
complexes of (maybe infinite) direct sum of representable presheaves with $U^*_{\alpha}$ smooth, 
we get the morphism in $C(\Var(k)/S)$
\begin{eqnarray*}
T(\otimes,\hat R_S^{0CH})(Q^*_1,Q^*_2):\hat R^{0CH}(Q^*_1)\otimes R^{0CH}(Q^*_2) 
\xrightarrow{(T(\otimes,\hat R_S^{CH})(\mathbb Z(U_{1,\alpha}^m),\mathbb Z(U_{2,\beta}^n))}
\hat R^{0CH}(Q^*_1\otimes Q^*_2))
\end{eqnarray*},
together with the commutative diagram in $C(\Var(k)/S)$
\begin{equation*}
\xymatrix{\hat R^{0CH}(Q^*_1)\otimes R^{0CH}(Q^*_2)
\ar[rrr]^{T(\otimes,\hat R_S^{0CH})(Q_1^*,Q_2)}
\ar[d]_{T(\hat R_S^{0CH},R_S^{CH})(Q_1^*)\otimes T(\hat R_S^{0CH},R_S^{0CH})(Q_2^*)} 
& \, & \, & \hat R^{0CH}(Q^*_1\times Q_2^*)\ar[d]^{T(\hat R_S^{0CH},R_S^{0CH})(Q^*_1\otimes Q^*_2)} \\
R^{0CH}(Q_1^*)\otimes R^{0CH}(Q_2^*)
\ar[rrr]^{T(\otimes,R_S^{0CH})(Q_1^*,Q_2)} & \, & \, & R^{0CH}(Q_1^*\otimes Q_2)}.
\end{equation*}
For $F_1,F_2\in C(\Var(k)^{sm}/S)$, we get in particular the morphism in $C(\Var(k)/S)$
\begin{eqnarray*}
T(\otimes,R_S^{0CH})(\rho_S^*LF_1,\rho_S^*LF_2):R^{0CH}(\rho_S^*LF_1)\otimes R^{0CH}(\rho_S^*LF_2) 
\to R^{0CH}(\rho_S^*(LF_1\otimes LF_2))
\end{eqnarray*}
together with the commutative diagram in $C(\Var(k)/S)$
\begin{equation*}
\xymatrix{\hat R^{0CH}(\rho_S^*LF_1)\otimes R^{0CH}(\rho_S^*LF_2)
\ar[rrr]^{T(\otimes,\hat R_S^{0CH})(\rho_S^*LF_1,\rho_S^*LF_2)}
\ar[d]_{T(\hat R_S^{0CH},R_S^{0CH})(\rho_S^*LF_1)\otimes T(\hat R_S^{0CH},R_S^{0CH})(\rho_S^*LF_2)} 
& \, & \, & \hat R^{0CH}(\rho_S^*LF_1\times\rho_S^*LF_2)
\ar[d]^{T(\hat R_S^{0CH},R_S^{0CH})(\rho_S^*LF_1\otimes\rho_S^*LF_2)} \\
R^{0CH}(\rho_S^*LF_1)\otimes R^{0CH}(\rho_S^*LF_2)
\ar[rrr]^{T(\otimes,R_S^{0CH})(\rho_S^*LF_1,\rho_S^*LF_2)} 
& \, & \, & R^{0CH}(\rho_S^*LF_1\times\rho_S^*LF_2)}.
\end{equation*}

\end{itemize}

\section{Triangulated category of motives}

\subsection{Definition and the six functor formalism}

The category of motives is obtained by inverting the $(\mathbb A^1_S,et)$ equivalence.
Hence the $\mathbb A^1_S$ local complexes of presheaves plays a key role.
\begin{defi}
The derived category of motives of complex algebraic varieties over $S$ is
the category
\begin{equation*}
\DA(S):=\Ho_{\mathbb A_S^1,et}(C(\Var(k)^{sm}/S)),
\end{equation*}
which is the localization of the category of complexes of presheaves on $\Var(k)^{sm}/S$
with respect to $(\mathbb A_S^1,et)$ local equivalence and we denote by
\begin{equation*}
D(\mathbb A_S^1,et):=D(\mathbb A_S^1)\circ D(et):C(\Var(k)^{sm}/S)\to\DA(S)
\end{equation*}
the localization functor. We have 
$\DA^-(S):=D(\mathbb A_S^1,et)(\PSh(\Var(k)^{sm}/S,C^-(\mathbb Z)))\subset\DA(S)$
the full subcategory consisting of bounded above complexes.
\end{defi}

\begin{defi}
The stable derived category of motives of complex algebraic varieties over $S$ is the category
\begin{equation*}
\DA_{st}(S):=\Ho_{\mathbb A_S^1,et}(C_{\Sigma}(\Var(k)^{sm}/S)),
\end{equation*}
which is the localization of the category of 
$\mathbb G_{mS}$-spectra ($\Sigma F^{\bullet}=F^{\bullet}\otimes\mathbb G_{mS}$) 
of complexes of presheaves on $\Var(k)^{sm}/S$
with respect to $(\mathbb A_S^1,et)$ local equivalence. The functor 
\begin{equation*}
\Sigma^{\infty}:C(\Var(k)^{sm}/S)\hookrightarrow
C_{\Sigma}(\Var(k)^{sm}/S)
\end{equation*}
induces the functor $\Sigma^{\infty}:\DA(S)\to\DA_{st}(S)$.
\end{defi}

We have all the six functor formalism by \cite{C.D}. 
We give a list of the operation we will use : 

\begin{itemize}

\item For $f:T\to S$ a morphism with $S,T\in\Var(k)$, the adjonction
\begin{equation*}
(f^*,f_*):C(\Var(k)^{sm}/S)\leftrightarrows C(\Var(k)^{sm}/T)
\end{equation*}
is a Quillen adjonction which induces in the derived categories ($f^*$ derives trivially),
$(f^*,Rf_*):\DA(S)\leftrightarrows\DA(T)$.

\item For $h:V\to S$ a smooth morphism with $V,S\in\Var(k)$, the adjonction
\begin{equation*}
(h_{\sharp},h^*):C(\Var(k)^{sm}/V)\leftrightarrows C(\Var(k)^{sm}/S)
\end{equation*}
is a Quillen adjonction which induces in the derived categories ($h^*$ derive trivially) 
$(Lh_{\sharp},h^*)=:\DA(V)\leftrightarrows\DA(S)$.

\item For $i:Z\hookrightarrow S$ a closed embedding, with $Z,S\in\Var(k)$,
\begin{equation*}
(i_*,i^!):=(i_*,i^{\bot}):C(\Var(k)^{sm}/Z)\leftrightarrows C(\Var(k)^{sm}/S)
\end{equation*}
is a Quillen adjonction, which induces in the derived categories ($i_*$ derive trivially) 
$(i_*,Ri^!):\DA(Z)\leftrightarrows\DA(S)$. The fact that $i_*$ derive trivially 
(i.e. send $(\mathbb A^1,et)$ local equivalence to $(\mathbb A^1,et)$ local equivalence is proved in \cite{AyoubT}.

\item For $S\in\Var(k)$, the adjonction given by
the tensor product of complexes of abelian groups and the internal hom of presheaves
\begin{equation*}
((\cdot\otimes\cdot),\mathcal Hom^{\bullet}(\cdot,\cdot)):
C(\Var(k)^{sm}/S)^2\to C(\Var(k)^{sm}/S),
\end{equation*}
is a Quillen adjonction, which induces in the derived category
\begin{equation*},
((\cdot\otimes^L\cdot),R\mathcal Hom^{\bullet}(\cdot,\cdot)):\DA(S)^2\to\DA(S),
\end{equation*} 

\item Let $M,N\in\DA(S)$, $Q^{\bullet}$ projectively cofibrant such that $M=D(\mathbb A^1,et)(Q^{\bullet})$, 
and $G^{\bullet}$ be $\mathbb A^1$ local for the etale topology such that $N=D(\mathbb A^1,et)(G^{\bullet})$. Then, 
\begin{equation}
R\mathcal Hom^{\bullet}(M,N)=\mathcal Hom^{\bullet}(Q^{\bullet},E(G^{\bullet}))\in\DA(S).
\end{equation}
This is well defined since if $s:Q_1\to Q_2$ is a etale local equivalence, 
\begin{equation*}
\mathcal Hom(s,E(G)):\mathcal Hom(Q_1,E(G))\to\mathcal Hom(Q_2,E(G))
\end{equation*}
is a etale local equivalence for $1\leq i\leq l$. 

\end{itemize}

We get from the first point 2 functors :
\begin{itemize}
\item The 2-functor $C(\Var(k)^{sm}/\cdot):\Var(k)\to Ab\Cat$, given by
\begin{equation*}
S\mapsto C(\Var(k)^{sm}/S) \;,\;
 (f:T\to S)\mapsto(f^*:C(\Var(k)^{sm}/S)\to C(\Var(k)^{sm}/T)). 
\end{equation*}
\item The 2-functor $\DA(\cdot):\Var(k)\to\TriCat$, given by
\begin{equation*}
S\mapsto\DA(S) \;,\; (f:T\to S)\mapsto(f^*:\DA(S)\to\DA(T)). 
\end{equation*}
\end{itemize}
The main theorem is the following :
\begin{thm}\cite{AyoubT}\cite{C.D}\label{2functDM}
The 2-functor $\DA(\cdot):\Var(k)\to\TriCat$, given by
\begin{equation*}
S\mapsto\DA(S) \;,\; (f:T\to S)\mapsto(f^*:\DA(S)\to\DA(T)) 
\end{equation*}
is a 2-homotopic functor (\cite{AyoubT})
\end{thm}

From theorem \ref{2functDM}, we get in particular
\begin{itemize}

\item For $f:T\to S$ a morphism with $T,S\in\Var(k)$, there by theorem \ref{2functDM} is also a pair of adjoint functor
\begin{equation*}
(f_!,f^!):\DA(S)\leftrightarrows\DA(T) 
\end{equation*}
\begin{itemize}
\item with $f_!=Rf_*$ if $f$ is proper,
\item with $f^!=f^*[d]$ if $f$ is smooth of relative dimension $d$. 
\end{itemize}
For $h:U\to S$ a smooth morphism with $U,S\in\Var(k)$ irreducible, have, for $M\in\DA(U)$, an isomorphism  
\begin{equation}\label{If}
Lh_{\sharp}M\to h_!M[d],
\end{equation}
in $\DA(S)$.

\item The 2-functor $S\in\Var(k)\mapsto\DA(S)$ satisfy the localization property, that is
for $i:Z\hookrightarrow X$ a closed embedding with $Z,X\in\Var(k)$, 
denote by $j:S\backslash Z\hookrightarrow S$ the open complementary subset, 
we have for $M\in\DA(S)$ a distinguish triangle in $\DA(S)$
\begin{equation*}
j_{\sharp}j^*M\xrightarrow{\ad(j_{\sharp},j^*)(M)}M\xrightarrow{\ad(i^*,i_*)(M)}i_*i^*M\to j_{\sharp}j^*M[1]
\end{equation*}
equivalently, 
\begin{itemize}
\item the functor 
\begin{equation*}
(i^*,j^*):\DA(S)\xrightarrow{\sim}\DA(Z)\times\DA(S\backslash Z)
\end{equation*}
is conservative,
\item  and for $F\in C(\Var(k)^{sm}/Z)$,
the adjonction map $\ad(i^*,i_*)(F):i^*i_*F\to F$ is an equivalence Zariski local, 
hence for $M\in\DA(S)$, the induced map in the derived category 
\begin{equation*}
\ad(i^*,i_*)(M):i^*i_*M\xrightarrow{\sim} M 
\end{equation*}
is an isomorphism.
\end{itemize}

\item For $f:X\to S$ a proper map, $g:T\to S$ a morphism, with $T,X,S\in\Var(k)$, and $M\in\DA(X)$,
\begin{equation*}
T(f,g)(M):g^*Rf_*M\to Rf'_*\tilde g^{'*}M 
\end{equation*}
is an isomorphism in $\DA(T)$ if $f$ is proper.

\end{itemize}

\begin{defi}
The derived category of extended motives of complex algebraic varieties over $S$ is
the category
\begin{equation*}
\underline{\DA}(S):=\Ho_{\mathbb A_S^1,et}(C(\Var(k)/S)),
\end{equation*}
which is the localization of the category of complexes of presheaves on $\Var(k)/S$
with respect to $(\mathbb A_S^1,et)$ local equivalence and we denote by
\begin{equation*}
D(\mathbb A_S^1,et):=D(\mathbb A_S^1)\circ D(et):C(\Var(k)/S)\to\underline{\DA}(S)
\end{equation*}
the localization functor. We have 
$\underline{\DA}^-(S):=D(\mathbb A_S^1,et)(\PSh(\Var(k)/S,C^-(\mathbb Z)))\subset\underline{\DA}(S)$
the full subcategory consisting of bounded above complexes.
\end{defi}

\subsection{Constructible motives and resolution of a motive by Corti-Hanamura motives}

We now give the definition of the motives of morphisms $f:X\to S$ which are constructible motives and
the definition of the category of Corti-Hanamura motives.

\begin{defi}
Let $S\in\Var(k)$, 
\begin{itemize}
\item the homological motive functor is 
$M(/S):\Var(k)/S\to\DA(S) \;, \; (f:X\to S)\mapsto M(X/S):=f_!f^!M(S/S)$,
\item the cohomological motive functor is 
$M^{\vee}(/S):\Var(k)/S\to\DA(S) \;, \; (f:X\to S)\mapsto M(X/S)^{\vee}:=Rf_*M(X/X):=f_*E_{et}(\mathbb Z_X)$,
\item the Borel-Moore motive functor is 
$M^{BM}(/S):\Var(k)/S\to\DA(S)\;, \; (f:X\to S)\mapsto M^{BM}(X/S):=f_!M(X/X)$,
\item the (homological) motive with compact support functor is
$M_c(/S):\Var(k)/S\to\DA(S)\; , \; (f:X\to S)\mapsto M_c(X/S):=Rf_*f^!M(S/S)$.
\end{itemize}
Let $f:X\to S$ a morphism with $X,S\in\Var(k)$. Assume that there exist a factorization
$f:X\xrightarrow{i} Y\times S\xrightarrow{p} S$, with $Y\in\SmVar(k)$, 
$i:X\hookrightarrow Y$ is a closed embedding and $p$ the projection. Then, 
\begin{equation*}
Q(X/S):=p_{\sharp}\Gamma^{\vee}_X\mathbb Z_{Y\times S}\in C(\Var(k)^{sm}/S)
\end{equation*}
is projective, admits transfert, and satisfy $D(\mathbb A^1_S,et)(Q(X/S))=M(X/S)$.
\end{defi}

\begin{defi}\label{constrDM}
\begin{itemize}
\item[(i)] Let $S\in\Var(k)$. 
We define the full subcategory 
\begin{eqnarray*}
\DA_c(S):&=&<Rf_*\mathbb Z_X,(f:X\to S)\in\Var(k)> \\
&=&<Rf_*\mathbb Z_X,(f:X\to S)\in\Var(k), \; \mbox{proper}, \; X \; \mbox{smooth}\subset\DA(S) 
\end{eqnarray*}
where $<,>$ denoted the full triangulated category generated by. 
\item[(ii)]Let $X,S\in\Var(k)$. If $f:X\to S$ is proper (but not necessary smooth) and $X$ is smooth, 
$M(X/S)$ is said to be a Corti-Hanamura motive and we have by above in this case 
$M(X/S)=M^{BM}(X/S)[c]=M(X/S)^{\vee}[c]$, with $c=\codim(X,X\times S)$ where $f:X\hookrightarrow X\times S\to S$. 
We denote by 
\begin{equation*}
\mathcal{CH}(S)=\left\{M(X/S)\right\}_{\left\{X/S=(X,f),f \mbox{pr},X \mbox{sm}\right\}}^{pa}\subset DM(S) 
\end{equation*}
the full subcategory which is the pseudo-abelian completion of the full subcategory 
whose objects are Corti-Hanamura motives.
\item[(iii)] We denote by
\begin{equation*}
\mathcal {CH}^0(S)\subset\mathcal{CH}(S) 
\end{equation*}
the full subcategory which is the pseudo-abelian completion of the full subcategory 
whose objects are Corti-Hanamura motives $M(X/S)$ such that the morphism $f:X\to S$ is projective.
\end{itemize}
\end{defi}

For bounded above motives, we have

\begin{thm}\label{weightst}
Let $S\in\Var(k)$.
\begin{itemize}
\item[(i)] There exists a unique weight structure $\omega$ on $\DA^-(S)$ such that $\DA^-(S)^{\omega=0}=\mathcal{CH}(S)$
\item[(ii)] There exist a well defined functor 
\begin{equation*}
W(S):\DA^-(S)\to K^-(\mathcal{CH}(S)) \; , \; W(S)(M)=[M^{(\bullet)}]
\end{equation*}
where $M^{(\bullet)}\in C^-(\mathcal{CH}(S))$ is a bounded above weight complex, such that if $m\in\mathbb Z$ is the highest weight,
we have a generalized distinguish triangle for all $i\leq m$
\begin{equation}\label{weightDT}
T_i:M^{(i)}[i]\to M^{(i+1)}[(i+1)]\to\cdots\to M^{(m)}[m]\to M^{w\geq i}
\end{equation}
Moreover the maps $w(M)^{(\geq i)}:M^{\geq i}\to M$ induce an isomorphism
 $w(M):holim_{i} M^{\geq i}\xrightarrow{\sim} M$
in $\DA^-(S)$.
\item[(iii)]Denote by $Chow(S)$ the category of Chow motives, which is the pseudo-abelian completion of the category 
\begin{itemize}
\item whose set of objects consist of the $X/S=(X,f)\in\Var(k)/S$ such that $f$ is proper and $X$ is smooth, 
\item whose set of morphisms between $X_1/S$ and $X_2/S$ is $\CH^{d_1}(X_1\times_S X_2)$, 
and the composition law is given in \cite{C.H}.
\end{itemize} 
We have then a canonical functor $CH_S:Chow(S)\hookrightarrow\DA(S)$, with $CH_S(X/S):=M(X/S):=Rf_*\mathbb Z(X/X)$,
which is a full embedding whose image is the category $\mathcal{CH}(S)$.
\end{itemize}
\end{thm}

\begin{proof}
\noindent(i): The category $\DA(S)$ is clearly weakly generated by $\mathcal{CH}(S)$.
 Moreover $\mathcal{CH}(S)\subset\DA(S)$ is negative. Hence, the result follows from
\cite{Bondarko} theorem 4.3.2 III.

\noindent(ii): Follows from (i) by standard fact of weight structure on triangulated categories. 
See \cite{Bondarko} theorem 3.2.2 and theorem 4.3.2 V for example.

\noindent(iii): See \cite{Fangzhou}.
\end{proof}

\begin{thm}\label{weightst2}
Let $S\in\Var(k)$.
\begin{itemize}
\item[(i)] There exists a unique weight structure $\omega$ on $\DA^-(S)$ such that $\DA^-(S)^{\omega=0}=\mathcal{CH}^0(S)$
\item[(ii)] There exist a well defined functor 
\begin{equation*}
W(S):\DA^-(S)\to K^-(\mathcal{CH}^0(S)) \; , \; W(S)(M)=[M^{(\bullet)}]
\end{equation*}
where $M^{(\bullet)}\in C^-(\mathcal{CH}^0(S))$ is a bounded above weight complex, such that if $m\in\mathbb Z$ is the highest weight,
we have a generalized distinguish triangle for all $i\leq m$
\begin{equation}\label{weightDT2}
T_i:M^{(i)}[i]\to M^{(i+1)}[(i+1)]\to\cdots\to M^{(m)}[m]\to M^{w\geq i}
\end{equation}
Moreover the maps $w(M)^{(\geq i)}:M^{\geq i}\to M$ induce an isomorphism 
$w(M):holim_{i} M^{\geq i}\xrightarrow{\sim} M$ in $\DA^-(S)$.
\end{itemize}
\end{thm}

\begin{proof}
Similar to the proof of theorem \ref{weightst}.
\end{proof}

\begin{cor}\label{weightst2Cor}
Let $S\in\Var(k)$. Let $M\in\DA(S)$. 
Then there exist $(F,W)\in C_{fil}(\Var(k)^{sm}/S)$ such that $D(\mathbb A^1,et)(F)=M$
and $D(\mathbb A^1,et)(\Gr^W_pF)\in\mathcal{CH}^0(S)$.
\end{cor}

\begin{proof}
By theorem \ref{weightst2}, there exist, by induction, for $i\in\mathbb Z$, a distinguish triangle in $\DA(S)$
\begin{equation}
T_i:M^{(i)}[i]\xrightarrow{m_i} M^{(i+1)}\xrightarrow{m_{i+1}}\cdots\xrightarrow{m_{m-1}}M^{(m)}[m]\to M^{w\geq i}
\end{equation}
with $M^{(j)}[j]\in\mathcal{CH}^0(S)$ and $w(M):holim_{i} M^{\geq i}\xrightarrow{\sim} M$ in $\DA^-(S)$.
For $i\in\mathbb Z$, take $(F_j)_{j\geq i},F_{w\geq i}\in C(\Var(k)^{sm}/S)$ 
such that $D(\mathbb A^1,et)(F_j)=M^{(j)}[j]$, $D(\mathbb A^1,et)(F_{w\geq i})=M^{w\geq i}$ 
and such that we have in $C(\Var(k)^{sm}/S)$, 
\begin{equation}
F_{w\geq i}=\Cone(F_i\xrightarrow{m_i} F_{i+1}\xrightarrow{m_{i+1}}\cdots\xrightarrow{m_{m-1}}F_m)
\end{equation}
where $m_j:F_j\to F_{j+1}$ are morphisms in $C(\Var(k)^{sm}/S)$ such that $D(\mathbb A^1,et)(m_j)=m_j$.
Now set $F=\holim_iF_{w\geq i}\in C(\Var(k)^{sm}/S)$ and $W_iF:=F_{w\geq i}\hookrightarrow F$,
so that $(F,W)\in C_{fil}(\Var(k)^{sm}/S)$ satisfy $D(\mathbb A^1,et)(\Gr^W_pF)=M^{(p)}[p]\in\mathcal{CH}^0(S)$.
\end{proof}

\section{The (filtered) D modules and the (filtered) De Rham functor 
on algebraic varieties over a field $k$ of characteristic zero}

\subsection{The D-modules on smooth algebraic varieties over a field $k$ of characteristic zero 
and their functorialities}

Let $k$ a field of characteristic zero.

For $S=(S,O_S)\in\SmVar(k)$, we denote by
\begin{itemize}
\item $D_S:=D(O_S)\subset\mathcal Hom_{\mathbb C_S}(O_S,O_S)$ the subsheaf consisting of differential operators.
By a $D_S$ module, we mean a left $D_S$ module.
\item we denote by
\begin{itemize}
\item $\PSh_{\mathcal D}(S)$ the abelian category of Zariski presheaves on $S$ with a structure of left $D_S$ module,
and by $\PSh_{\mathcal D,c}(S)\subset\PSh_{\mathcal D}(S)$ 
the full subcategory whose objects are coherent sheaves of left $D_S$ modules,
\item $\PSh_{\mathcal D^{op}}(S)$ the abelian category of Zariski presheaves on $S$ with a structure of right $D_S$ module,
and by $\PSh_{\mathcal D^{op},c}(S)\subset\PSh_{\mathcal D^{op}}(S)$ 
the full subcategories whose objects are coherent sheaves of right $D_S$ modules,
\end{itemize}
\item we denote by
\begin{itemize}
\item $C_{\mathcal D}(S)=C(\PSh_{\mathcal D}(S))$ 
the category of complexes of Zariski presheaves on $S$ with a structure of $D_S$ module,
\item $C_{\mathcal D^{op}}(S)=C(\PSh_{\mathcal D^{op}}(S))$ the category of complexes of Zariski presheaves on $S$
with a structure of right $D_S$ module,
\end{itemize}
\item in the filtered case we have
\begin{itemize}
\item $C_{\mathcal D(2)fil}(S)\subset C(\PSh_{\mathcal D}(S),F,W):=C(\PSh_{D(O_S)}(S),F,W)$ 
the category of (bi)filtered complexes of algebraic $D_S$ modules such that the filtration is biregular,
$D_{\mathcal D(2)fil}(S):=\Ho_{zar}(C_{\mathcal D(2)fil}(S))$ its localization with respect
to filtered Zariski local equivalence, and more generally
$D_{\mathcal D(2)fil,r}(S):=\Ho_{zar}(C_{\mathcal D(2)fil}(S))$ its localization with respect
to $r$-filtered Zariski local equivalence for $r=1,\cdots,\infty$,
\item $C_{\mathcal D0fil}(S)\subset C_{\mathcal Dfil}(S)$ the full subcategory such that
the filtration is a filtration by $D_S$ submodule (which is stronger then Griffitz transversality), 
$C_{\mathcal D(1,0)fil}(S)\subset C_{\mathcal D2fil}(S)$
the full subcategory such that $W$ is a filtration by $D_S$ submodules,
$D_{\mathcal D(1,0)fil}(S):=\Ho_{zar}(C_{\mathcal D(1,0)fil}(S))$ its localization with respect
to filtered Zariski local equivalence, and more generally
$D_{\mathcal D(1,0)fil,r}(S):=\Ho_{zar}(C_{\mathcal D(1,0)fil}(S))$ its localization with respect
to $r$-filtered Zariski local equivalence for $r=1,\cdots,\infty$,
\item $C_{\mathcal D^{op}(2)fil}(S)\subset C(\PSh_{\mathcal D^{op}}(S),F,W):=C(\PSh_{D(O_S)^{op}}(S),F,W)$ 
the category of (bi)filtered complexes of algebraic (resp. analytic) right $D_S$ modules such that the filtration is biregular,
as in the left case we consider the subcategories as above.
\end{itemize}
\end{itemize}

\begin{defi}\label{Daffine}
An $X\in\SmVar(k)$ is said to be D-affine if the following two condition hold:
\begin{itemize}
\item[(i)] The global section functor $\Gamma(X,\cdot):\mathcal QCoh_{\mathcal D}(X)\to Mod(\Gamma(X,D_X))$ is exact.
\item[(ii)] If $\Gamma(X,M)=0$ for $M\in\mathcal QCoh_{\mathcal D}(X)$, then $M=0$.
\end{itemize} 
\end{defi}

\begin{prop}\label{Daffineprop}
If $X\in\SmVar(k)$ is D-affine, then :
\begin{itemize}
\item[(i)] Any $M\in\mathcal QCoh_{\mathcal D}(X)$ is generated by its global sections.
\item[(ii)] The functor $\Gamma(X,\cdot):\mathcal QCoh_{\mathcal D}(X)\to Mod(\Gamma(X,D_X))$ is an equivalence
of category whose inverse is $L\in Mod(\Gamma(X,D_X))\mapsto D_X\otimes_{\Gamma(X,D_X)}L\in\mathcal QCoh_{\mathcal D}(X)$.
\item[(iii)] We have $\Gamma(X,\cdot)(\mathcal Coh_{\mathcal D}(X))=Mod(\Gamma(X,D_X))_f$, that is the global sections of a
coherent $D_X$ module is a finite module over the differential operators on $X$.
\end{itemize} 
\end{prop}

\begin{proof}
Similar to the complex case : see \cite{LvDmod}.
\end{proof}

Let $f:X\to S$ be a morphism with $X,S\in\SmVar(k)$,
Then, we recall from \cite{B4} section 4.1, the transfers modules 
\begin{itemize}
\item $(D_{X\to S},F^{ord}):=f^{*mod}(D_S,F^{ord}):=f^*(D_S,F^{ord})\otimes_{f^*O_S}(O_X,F_b)$ 
which is a left $D_X$ module and a left and right $f^*D_S$ module
\item $(D_{X\leftarrow S},F^{ord}):=(K_X,F_b)\otimes_{O_X}(D_{X\to S},F^{ord})\otimes_{f^*O_S}f^*(K_S,F_b)$.
which is a right $D_X$ module and a left and right $f^*D_S$ module.
\end{itemize}

\begin{prop}
Let $i:Z\hookrightarrow S$ be a closed embedding with $Z,S\in\SmVar(k)$.
Then, $D_{Z\to S}=i^*D_S/D_S\mathcal I_Z$ and it is a locally free (left) $D_Z$ module.
Similarly, $D_{Z\leftarrow S}=i^*D_S/\mathcal I_ZD_S$ and it is a locally free right $D_Z$ module.
\end{prop}

\begin{proof}
Similar to the complex case : see \cite{LvDmod}.
\end{proof}

\begin{itemize}

\item Let $S\in\SmVar(k)$. 
\begin{itemize}
\item For $M\in C_{\mathcal D}(S)$, we have the canonical projective resolution $q:L_D(M)\to M$ of complexes
of $D_S$ modules.
\item Let $\tau$ a topology on $S$. For $M\in C_{\mathcal D}(S)$, 
there exist a unique strucure of $D_S$ module on the flasque presheaves $E_{\tau}^i(M)$ such that
$E_{\tau}(M)\in C_{\mathcal D}(S)$ (i.e. is a complex of $D_S$ modules) 
and that the map $k:M\to E(M)$ is a morphism of complexes of $D_S$ modules.
\end{itemize}

\item Let $S\in\SmVar(k)$. 
For $M\in C_{\mathcal D^{(op)}}(S)$, $N\in C(S)$, we will consider the induced D module structure
(right $D_S$ module in the case one is a left $D_S$ module and the other one is a right one)
on the presheaf $M\otimes N:=M\otimes_{\mathbb Z_S}N$ (see section 2). We get the bifunctor
\begin{eqnarray*}
C(S)\times C_{\mathcal D}(S)\to C_{\mathcal D}(S), (M,N)\mapsto M\otimes N
\end{eqnarray*}

\item Let $S\in\SmVar(k)$. 
For $M,N\in C_{\mathcal D^{(op)}}(S)$, $M\otimes_{O_S} N$ has a canonical structure of $D_S$ modules 
(right $D_S$ module in the case one is a left $D_S$ module and the other one is a right one)
given by (in the left case) for $S^o\subset S$ an open subset, 
\begin{equation*}
m\otimes n\in\Gamma(S^o,M\otimes_{O_S}N),\gamma\in\Gamma(S^o,D_S), \;
\gamma.(m\otimes n):=(\gamma.m)\otimes n-m\otimes\gamma.n
\end{equation*}
This gives the bifunctor
\begin{eqnarray*}
C_{\mathcal D^{(op)}}(S)^2\to C_{\mathcal D^{(op)}}(S), (M,N)\mapsto M\otimes_{O_S}N
\end{eqnarray*}
More generally, let $f:X\to S$ a morphism with $X,S\in\Var(k)$. Assume $S$ smooth. 
For $M,N\in C_{f^*\mathcal D^{(op)}}(X)$, $M\otimes_{f^*O_S} N$ (see section 2),
has a canonical structure of $f^*D_S$ modules 
(right $f^*D_S$ module in the case one is a left $f^*D_S$ module and the other one is a right one)
given by (in the left case) for $X^o\subset X$ an open subset, 
\begin{equation*}
m\otimes n\in\Gamma(X^o,M\otimes_{f^*O_S} N),\gamma\in\Gamma(X^o,f^*D_S), \;
\gamma.(m\otimes n):=(\gamma.m)\otimes n-m\otimes\gamma.n
\end{equation*}
This gives the bifunctor
\begin{eqnarray*}
C_{f^*\mathcal D^{(op)}}(X)^2\to C_{f^*\mathcal D^{(op)}}(X), (M,N)\mapsto M\otimes_{f^*O_S}N
\end{eqnarray*}

\item Let $S\in\SmVar(k)$. The internal hom bifunctor 
\begin{equation*}
\mathcal Hom(\cdot,\cdot):=\mathcal Hom_{\mathbb Z_S}(\cdot,\cdot):C(S)^2\to C(S)
\end{equation*}
induces a bifunctor
\begin{equation*}
\mathcal Hom(\cdot,\cdot):=\mathcal Hom_{\mathbb Z_S}(\cdot,\cdot):C(S)\times C_{\mathcal D}(S)\to C_{\mathcal D}(S)
\end{equation*}
such that, for $F\in C(S)$ and $G\in C_{\mathcal D}(S)$, the $D_S$ structure on $\mathcal Hom^{\bullet}(F,G)$ is given by
\begin{equation*}
\gamma\in\Gamma(S^o,D_S)\longmapsto 
(\phi\in\Hom^{p}(F^{\bullet}_{|S^o},G_{|S^o})\mapsto
(\gamma\cdot\phi:\alpha\in F^{\bullet}(S^o)\mapsto\gamma\cdot\phi^p(S^o)(\alpha))
\end{equation*}
where $\phi^p(S^o)(\alpha)\in\Gamma(S^o,G)$.

\item Let $S\in\SmVar(k)$. 
For $M,N\in C_{\mathcal D}(S)$, $\mathcal Hom_{O_S}(M,N)$,
has a canonical structure of $D_S$ modules given by for $S^o\subset S$ an open subset and
$\phi\in\Gamma(S^o,\mathcal Hom(M,O_S))$, $\gamma\in\Gamma(S^o,D_S)$,
$(\gamma.\phi)(m):=\gamma.(\phi(m))-\phi(\gamma.m)$
This gives the bifunctor
\begin{eqnarray*}
\Hom^{\bullet}_{O_S}(-,-):C_{\mathcal D}(S)^2\to C_{\mathcal D}(S)^{op}, (M,N)\mapsto\mathcal Hom^{\bullet}_{O_S}(M,N)
\end{eqnarray*}

\item Let $S\in\SmVar(k)$. 
We have the bifunctors
\begin{itemize}
\item $\Hom^{\bullet}_{D_S}(-,-):C_{\mathcal D}(S)^2\to C(S)$, $(M,N)\mapsto\mathcal Hom^{\bullet}_{D_S}(M,N)$,
and if $N$ is a bimodule (i.e. has a right $D_S$ module structure whose opposite coincide with the left one),
$\mathcal Hom_{D_S}(M,N)\in C_{\mathcal D^{op}}(S)$ given by for $S^o\subset S$ an open subset and
$\phi\in\Gamma(S^o,\mathcal Hom(M,N))$, $\gamma\in\Gamma(S^o,D_S)$, $(\phi.\gamma)(m):=(\phi(m)).\gamma$
\item $\Hom_{D_S}(-,-):C_{\mathcal D^{op}}(S)^2\to C(S)$, $(M,N)\mapsto\mathcal Hom_{D_S}(M,N)$
and if $N$ is a bimodule, $\mathcal Hom_{D_S}(M,N)\in C_{\mathcal D}(S)$
\end{itemize}
For $M\in C_{\mathcal D}(S)$, we get in particular the dual with respect $\mathbb D_S$,
\begin{equation*}
\mathbb D_SM:=\mathcal Hom_{D_S}(M,D_S)\in C_{\mathcal D}(S) \; ; \;
\mathbb D^K_SM:=\mathcal Hom_{D_S}(M,D_S)\otimes_{O_S}\mathbb D^O_Sw(K_S)[d_S]\in C_{\mathcal D}(S)
\end{equation*}
and we have canonical map $d:M\to\mathbb D_S^2M$.
This functor induces in the derived category, for $M\in D_{\mathcal D}(S)$,
\begin{equation*}
L\mathbb D_SM:=R\mathcal Hom_{D_S}(L_DM,D_S)\otimes_{O_S}\mathbb D^O_Sw(K_S)[d_S]=\mathbb D^K_SL_DM\in D_{\mathcal D}(S).
\end{equation*}
where $\mathbb D^O_Sw(S):\mathbb D_S^Ow(K_S)\to\mathbb D^O_SK_S=K_S^{-1}$ 
is the dual of the Koczul resolution of the canonical bundle (proposition \ref{resw}), 
and the canonical map $d:M\to L\mathbb D_S^2M$.

\item Let $f:X\to S$ a morphism with $X,S\in\SmVar(k)$. 
For $N\in C_{\mathcal D,f^*\mathcal D}(X)$ and $M\in C_{\mathcal D}(X)$, $N\otimes_{D_X}M$ 
has the canonical $f^*D_S$ module structure given by, for $X^o\subset X$ an open subset, 
\begin{equation*}
\gamma\in\Gamma(X^o,f^*D_S), m\in\Gamma(X^o,M),n\in\Gamma(X^o,N), \; \gamma.(n\otimes m)=(\gamma.n)\otimes m.
\end{equation*}
This gives the functor
\begin{equation*}
C_{\mathcal D,f^*\mathcal D}(X)\times C_{\mathcal D}(X)\to C_{f^*\mathcal D}(X), \; (M,N)\mapsto M\otimes_{D_X}N
\end{equation*}

\item Let $f:X\to S$ be a morphism with $X,S\in\SmVar(k)$,
Then, for $M\in C_{\mathcal D}(S)$, $O_X\otimes_{f^*O_S}f^*M$ has a canonical $D_X$ module structure given by
given by, for $X^o\subset X$ an open subset,
\begin{equation*}
m\otimes n\in\Gamma(X^o,O_X\otimes_{f^*O_S}f^*M), \gamma\in\Gamma(X^o,D_X), \; 
\gamma.(m\otimes n):=(\gamma.m)\otimes n-m\otimes df(\gamma).n.
\end{equation*}
This gives the inverse image functor
\begin{equation*}
f^{*mod}:\PSh_{\mathcal D}(S)\to\PSh_{\mathcal D}(X), \; \; 
M\mapsto f^{*mod}M:=O_X\otimes_{f^*O_S}f^*M=D_{X\to S}\otimes_{f^*D_S}f^*M 
\end{equation*}
which induces in the derived category the functor
\begin{equation*}
Lf^{*mod}:D_{\mathcal D}(S)\to D_{\mathcal D}(X), \; \;
M\mapsto Lf^{*mod}M:=O_X\otimes^L_{f^*O_S}f^*M=O_X\otimes_{f^*O_S}f^*L_DM,
\end{equation*}
We will also consider the shifted inverse image functor
\begin{equation*}
Lf^{*mod[-]}:=Lf^{*mod}[d_S-d_X]:D_{\mathcal D}(S)\to D_{\mathcal D}(X).
\end{equation*}

\item Let $f:X\to S$ be a morphism with $X,S\in\SmVar(k)$.
For $M\in C_{\mathcal D}(X)$, $D_{X\leftarrow S}\otimes_{D_X}M$ has the canonical $f^*D_S$ module structure given above. 
Then, the direct image functor 
\begin{equation*}
f^0_{*mod}:\PSh_{\mathcal D}(X)\to\PSh_{\mathcal D}(S), \; \; 
M\mapsto f_{*mod} M:=f_*(D_{X\leftarrow S}\otimes_{D_X}M)
\end{equation*}
induces in the derived category the functor
\begin{equation*}
\int_{f}=Rf_{*mod}:D_{\mathcal D}(X)\to D_{\mathcal D}(S), \; \; 
M\mapsto \int_{f}M=Rf_*(D_{X\leftarrow S}\otimes^L_{D_X}M).
\end{equation*}

\end{itemize}

The functorialities given above induce :

\begin{itemize}

\item Let $S\in\SmVar(k)$. For $(M,F)\in C_{fil}(S)$ and $(N,F)\in C_{fil}(S)$, recall that 
\begin{equation*}
F^p((M,F)\otimes(N,F)):=\Im(\oplus_qF^qM\otimes F^{p-q}N\to M\otimes N)
\end{equation*}
This gives the functor
\begin{equation*}
(\cdot,\cdot):C_{fil}(S)\times C_{\mathcal Dfil}(S)\to C_{\mathcal Dfil}(S) \; , 
\;((M,F),(N,F))\mapsto (M,F)\otimes(N,F).
\end{equation*}
It induces in the derived categories by taking r-projective resolutions the bifunctors, for $r=1,\ldots,\infty$,
\begin{equation*}
(\cdot,\cdot):D_{\mathcal Dfil,r}(S)\times D_{fil,r}(S)\to D_{\mathcal Dfil,r}(S) \; , \;
((M,F),(N,F))\mapsto (M,F)\otimes^L(N,F)=L_D(M,F)\otimes (N,F).
\end{equation*}

\item Let $S\in\SmVar(k)$.
For $(M,F)\in C_{O_Sfil}(S)$ and $(N,F)\in C_{O_Sfil}(S)$, recall that 
\begin{equation*}
F^p((M,F)\otimes_{O'_S}(N,F)):=\Im(\oplus_qF^qM\otimes_{O'_S} F^{p-q}N\to M\otimes_{O'_S} N)
\end{equation*}
It induces in the derived categories by taking r-projective resolutions the bifunctors, for $r=1,\ldots,\infty$,
\begin{equation*}
(\cdot,\cdot):D_{\mathcal Dfil,r}(S)^2\to D_{\mathcal Dfil,r}(S) \; , \;((M,F),(N,F))\mapsto (M,F)\otimes_{O_S}^L(N,F).
\end{equation*}
More generally, let $f:X\to S$ a morphism with $X,S\in\Var(k)$. Assume $S$ smooth. We have the bifunctors
\begin{eqnarray*}
(\cdot,\cdot):D_{f^*\mathcal Dfil,r}(X)^2\to D_{f^*\mathcal Dfil,r}(X), \\
((M,F),(N,F))\mapsto (M,F)\otimes_{f^*O_S}^L(N,F)=(M,F)\otimes_{f^*O_S}L_{f^*D}(N,F).
\end{eqnarray*}

\item Let $S\in\SmVar(k)$. The hom functor induces the bifunctor
\begin{eqnarray*}
\Hom(-,-):C_{\mathcal Dfil}(S)\times C_{fil}(S)\to C_{\mathcal D(1,0)fil}(S), ((M,W),(N,F))\mapsto\mathcal Hom((M,W),(N,F)).
\end{eqnarray*}

\item Let $S\in\SmVar(k)$. The hom functor induces the bifunctor
\begin{eqnarray*}
\Hom_{O_S}(-,-):C_{\mathcal Dfil}(S)^2\to C_{\mathcal D2fil}(S), ((M,W),(N,F))\mapsto\mathcal Hom_{O_S}((M,W),(N,F)).
\end{eqnarray*}

\item Let $S\in\SmVar(k)$. 
The hom functor induces the bifunctors
\begin{itemize}
\item $\Hom_{D_S}(-,-):C_{\mathcal Dfil}(S)^2\to C_{2fil}(S)$, $((M,W),(N,F))\mapsto\mathcal Hom_{D_S}((M,W),(N,F))$,
\item $\Hom_{D_S}(-,-):C_{\mathcal D^{op}fil}(S)^2\to C_{2fil}(S)$, $((M,W),(N,F))\mapsto\mathcal Hom_{D_S}((M,W),(N,F))$.
\end{itemize}
We get the filtered dual 
\begin{eqnarray*}
\mathbb D^K_S(\cdot):C_{\mathcal D(2)fil}(S)\to C_{\mathcal D(2)fil}(S)^{op}, \; 
(M,F)\mapsto\mathbb D^K_S(M,F):=\mathcal Hom_{D_S}((M,F),D_S)\otimes_{O_S}\mathbb D_S^Ow(K_S)[d_S]
\end{eqnarray*}
together with the canonical map $d(M,F):(M,F)\to\mathbb D^{2,K}_S(M,F)$.
Of course $\mathbb D^K_S(\cdot)(C_{\mathcal D(1,0)fil}(S))\subset C_{\mathcal D(1,0)fil}(S)$.
It induces in the derived categories $D_{\mathcal Dfil,r}(S)$, for $r=1,\ldots,\infty$, the functors
\begin{equation*}
L\mathbb D_S(\cdot):D_{\mathcal D(2)fil,r}(S)\to D_{\mathcal D(2)fil,r}(S)^{op}, 
(M,F)\mapsto L\mathbb D_S(M,F):=\mathbb D^K_S L_D(M,F).
\end{equation*}
together with the canonical map $d(M,F):L_D(M,F)\to \mathbb D^{2,K}_SL_D(M,F)$.

\item Let $f:X\to S$ be a morphism with $X,S\in\SmVar(k)$.
Then, the inverse image functor
\begin{eqnarray*}
f^{*mod}:C_{\mathcal D(2)fil}(S)\to C_{\mathcal D(2)fil}(X), \\
(M,F)\mapsto f^{*mod}(M,F):=(O_X,F_b)\otimes_{f^*O_S}f^*(M,F)=(D_{X\to S},F^{ord})\otimes_{f^*D_S}f^*(M,F), 
\end{eqnarray*} 
induces in the derived categories the functors, for $r=1,\ldots,\infty$ (resp. $r\in(1,\ldots\infty)^2$),
\begin{eqnarray*}
Lf^{*mod}:D_{\mathcal D(2)fil,r}(S)\to D_{\mathcal D(2)fil,r}(X), \\  
(M,F)\mapsto Lf^{*mod}M:=(O_X,F_b)\otimes^L_{f^*O_S}f^*(M,F)=(O_X,F_b)\otimes_{f^*O_S}f^*L_D(M,F).
\end{eqnarray*}
Of course $f^{*mod}(C_{\mathcal D(1,0)fil}(S))\subset C_{\mathcal D(1,0)fil}(X)$. Note that
\begin{itemize}
\item If the $M$ is a complex of locally free $O_S$ modules, 
then $Lf^{*mod}(M,F)=f^{*mod}(M,F)$ in $D_{\mathcal D(2)fil,\infty}(S)$.
\item If the $\Gr^p_FM$ are complexes of locally free $O_S$ modules, 
then $Lf^{*mod}(M,F)=f^{*mod}(M,F)$ in $D_{\mathcal D(2)fil}(S)$.
\end{itemize}
We will consider also the shifted inverse image functors
\begin{equation*}
Lf^{*mod[-]}:=Lf^{*mod}[d_S-d_X]:D_{\mathcal D(2)fil,r}(S)\to D_{\mathcal D(2)fil,r}(X). 
\end{equation*}

\item Let $f:X\to S$ be a morphism with $X,S\in\SmVar(k)$.
Then,the direct image functor 
\begin{equation*}
f^{00}_{*mod}:(\PSh_{\mathcal D}(X),F)\to(\PSh_{\mathcal D}(S),F), \; \;
(M,F)\mapsto f_{*mod} (M,F):=f_*((D_{S\leftarrow X},F^{ord})\otimes_{D_X}(M,F))
\end{equation*}
induces in the derived categories by taking r-injective resolutions the functors, for $r=1,\ldots,\infty$,
\begin{eqnarray*}
\int_{f}=Rf_{*mod}:D_{\mathcal D(2)fil,r}(X)\to D_{\mathcal D(2)fil,r}(S), 
(M,F)\mapsto \int_{f}(M,F)=Rf_*((D_{S\leftarrow X},F^{ord})\otimes^L_{D_X}(M,F)).
\end{eqnarray*}
Let $f_1:X\to Y$ and $f_2:Y\to S$ two morphism with $X,Y,S\in\SmVar(\mathbb C)$ or with $X,Y,S\in\AnSm(\mathbb C)$.
We have, for $(M,F)\in C_{\mathcal Dfil}(X)$, the canonical transformation map in $D_{\mathcal D(2)fil,r}(S)$
\begin{eqnarray*}
T(\int_{f_2}\circ\int_{f_1},\int_{f_2\circ f_1})(M,F): \\
\int_{f_2}\int_{f_1}(M,F):=
Rf_{2*}((D_{Y\leftarrow S},F^{ord})\otimes^L_{D_Y}Rf_{1*}((D_{X\leftarrow Y},F^{ord})\otimes^L_{D_X}(M,F))) \\
\xrightarrow{T(f_1,\otimes)(-,-)}
Rf_{2*}Rf_{1*}(f_1^*(D_{Y\leftarrow S},F^{ord})\otimes^L_{D_Y}((D_{X\leftarrow Y},F^{ord})\otimes^L_{D_X}(M,F))) \\
\xrightarrow{\sim} 
Rf_{2*}Rf_{1*}((f_1^*(D_{Y\leftarrow S},F^{ord})\otimes^L_{D_Y}(D_{X\leftarrow Y},F^{ord}))\otimes^L_{D_X}(M,F)) \\
\xrightarrow{\sim}
Rf_{2*}Rf_{1*}((D_{X\leftarrow S},F^{ord})\otimes^L_{D_X}(M,F)):=\int_{f_2\circ f_1}(M,F)
\end{eqnarray*}

\item  Let $f:X\to S$ be a morphism with $X,S\in\SmVar(k)$.
Then the functor
\begin{eqnarray*}
f^{\hat{*}mod}:C_{\mathcal D2fil}(S)\to C_{\mathcal D2fil}(X), 
(M,F)\mapsto f^{\hat{*}mod}(M,F):=\mathbb D^K_XL_Df^{*mod}L_D\mathbb D^K_S(M,F) 
\end{eqnarray*} 
induces in the derived categories the exceptional inverse image functors, 
for $r=1,\ldots,\infty$ (resp. $r\in(1,\ldots\infty)^2$),
\begin{eqnarray*}
Lf^{\hat{*}mod}:D_{\mathcal D(2)fil,r}(S)\to D_{\mathcal D(2)fil,r}(X), \\
(M,F)\mapsto Lf^{\hat{*}mod}(M,F):=L\mathbb D_XLf^{*mod}L\mathbb D_S(M,F):=f^{\hat{*}mod}L_D(M,F).
\end{eqnarray*}
Of course $f^{\hat*mod}(C_{\mathcal D(1,0)fil}(S))\subset C_{\mathcal D(1,0)fil}(X)$.
We will also consider the shifted exceptional inverse image functors
\begin{equation*}
Lf^{\hat{*}mod[-]}:=Lf^{\hat{*}mod}[d_S-d_X]:D_{\mathcal D(2)fil,r}(S)\to D_{\mathcal D(2)fil,r}(X).
\end{equation*}

\item Let $S_1,S_2\in\SmVar(k)$. Consider $p:S_1\times S_2\to S_1$ the projection. 
Since $p$ is a projection, we have a canonical embedding $p^*D_{S_1}\hookrightarrow D_{S_1\times S_2}$.
For $(M,F)\in C_{\mathcal D(2)fil}(S_1\times S_2)$, $(M,F)$ has a canonical $p^*D_{S_1}$ module structure. 
Moreover, with this structure, for $(M_1,F)\in C_{\mathcal D(2)fil}(S_1)$ 
\begin{equation*}
\ad(p^{*mod},p)(M_1,F):(M_1,F)\to p_*p^{*mod}(M_1,F)
\end{equation*}
is a map of complexes of $D_{S_1}$ modules, and for $(M_{12},F)\in C_{\mathcal D(2)fil}(S_1\times S_2))$
\begin{equation*}
\ad(p^{*mod},p)(M_{12},F):p^{*mod}p_*(M_{12},F)\to (M_{12},F)
\end{equation*}
is a map of complexes of $D_{S_1\times S_2}$ modules.

\end{itemize}

\begin{prop}\label{compDmod}
\begin{itemize}
\item[(i)] Let $f_1:X\to Y$ and $f_2:Y\to S$ two morphism with $X,Y,S\in\SmVar(k)$. 
\begin{itemize}
\item Let $(M,F)\in C_{\mathcal D(2)fil,r}(S)$. Then $(f_2\circ f_1)^{*mod}(M,F)=f_1^{*mod}f_2^{*mod}(M,F)$.
\item Let $(M,F)\in D_{\mathcal D(2)fil,r}(S)$. Then $L(f_2\circ f_1)^{*mod}(M,F)=Lf_1^{*mod}(Lf_2^{*mod}(M,F))$.
\end{itemize}
\item[(ii)] Let $f_1:X\to Y$ and $f_2:Y\to S$ two morphism with $X,Y,S\in\SmVar(k)$. 
Let $M\in D_{\mathcal D}(X)$. Then, 
\begin{equation*}
T(\int_{f_2}\circ\int_{f_1},\int_{f_2\circ f_1})(M):\int_{f_2}\int_{f_1}(M)\xrightarrow{\sim}\int_{f_2\circ f_1}(M)
\end{equation*}
is an isomorphism in $D_{\mathcal D}(S)$ (i.e. if we forget filtration).
\item[(iii)] Let $i_0:Z_2\hookrightarrow Z_1$ and $i_1:Z_1\hookrightarrow S$ two closed embedding, 
with $Z_2,Z_1,S\in\SmVar(k)$. Let $(M,F)\in C_{\mathcal D(2)fil}(Z_2)$. 
Then, $(i_1\circ i_0)_{*mod}(M,F)=i_{1*mod}(i_{0*mod}(M,F))$ in  $C_{\mathcal D(2)fil}(S)$.
\end{itemize}
\end{prop}

\begin{proof}
Similar to the complex case : see \cite{B4}. 
\end{proof}

\begin{prop}\label{otimesed}
For $X\in\SmVar(k)$, we have for $(M,F),(N,F)\in C_{O_Xfil}(X)$ or $(M,F),(N,F)\in C_{\mathcal Dfil}(X)$,
Denote by $\Delta_X:X\hookrightarrow X\times X$ the diagonal closed embedding and $p_1:X\times X\to X$,
$p_2:X\times X\to X$ the projections. We have
\begin{equation*}
(M,F)\otimes_{O_X}(N,F)=\Delta_X^{*mod}(p_1^{*mod}(M,F)\otimes_{O_{X\times X}}p_2^{*mod}(N,F))
\end{equation*}
\end{prop}

\begin{proof}
Similar to the complex case : see \cite{LvDmod}.
\end{proof}

Let $i:Z\hookrightarrow S$ a closed embedding, with $Z,S\in\SmVar(k)$. We have the functor
\begin{eqnarray*}
i^{\sharp}:C_{\mathcal Dfil}(S)\to C_{\mathcal Dfil}(Z), 
(M,F)\mapsto i^{\sharp}(M,F):=\mathcal Hom_{i^*D_S}((D_{S\leftarrow Z},F^{ord}),i^*(M,F))
\end{eqnarray*}
where the (left) $D_Z$ module structure on $i^{\sharp}M$ 
comes from the right module structure on $D_{S\leftarrow Z}$, resp. $O_Z$. We denote by
\begin{itemize}
\item for $(M,F)\in C_{\mathcal Dfil}(S)$, the canonical map in $C_{\mathcal Dfil}(S)$
\begin{eqnarray*}
\ad(i_{*mod},i^{\sharp})(M,F): i_{*mod}i^{\sharp}(M,F):=
i_*(\mathcal Hom_{i^*D_S}((D_{S\leftarrow Z},F^{ord}),i^*(M,F))\otimes_{D_Z}(D_{S\leftarrow Z},F^{ord})) \\
\to (M,F), \phi\otimes P\mapsto\phi(P)
\end{eqnarray*}
\item for $(N,F)\in C_{\mathcal Dfil}(Z)$, the canonical map in $C_{\mathcal Dfil}(Z)$
\begin{eqnarray*}
\ad(i_{*mod},i^{\sharp})(N,F):(N,F)\to i^{\sharp}i_{*mod}(N,F):=
\mathcal Hom_{i^*D_S}(D_{S\leftarrow Z},i^*i_*((N,F)\otimes_{D_Z}(D_{S\leftarrow Z},F^{ord}))) \\
n\mapsto(P\mapsto n\otimes P)
\end{eqnarray*}
\end{itemize}
The functor $i^{\sharp}$ induces in the derived category the functor :
\begin{eqnarray*}
Ri^{\sharp}:D_{\mathcal D(2)fil,r}(S)\to D_{\mathcal D(2)fil,r}(Z), \; (M,F)\mapsto \\
Ri^{\sharp}(M,F):=R\mathcal Hom_{i^*D_S}((D_{Z\leftarrow S},F^{ord}),i^*(M,F))
=\mathcal Hom_{i^*D_S}((D_{Z\leftarrow S},F^{ord}),E(i^*(M,F))).
\end{eqnarray*}

\begin{prop}\label{isharpk}
Let $i:Z\hookrightarrow S$ a closed embedding, with $Z,S\in\SmVar(k)$. 
The functor  $i_{*mod}:C_{\mathcal D}(Z)\to C_{\mathcal D}(S)$ admit a right adjoint
which is the functor $i^{\sharp}:C_{\mathcal D}(S)\to C_{\mathcal D}(Z)$ and
\begin{equation*}
\ad(i_{*mod},i^{\sharp})(N):N\to i^{\sharp}i_{*mod}N \mbox{\; and \;} \ad(i_{*mod},i^{\sharp})(M):i_{*mod}i^{\sharp}M\to M
\end{equation*}
are the adjonction maps. 
\end{prop}

\begin{proof}
Similar to the complex case : see \cite{LvDmod}.
\end{proof}

One of the main results in D modules is Kashiwara equivalence :
\begin{thm}\label{Keqk}
Let $i:Z\hookrightarrow S$ a closed embedding with $Z,S\in\SmVar(k)$. 
\begin{itemize}
\item[(i)] The functor 
$i_{*mod}:\mathcal QCoh_{\mathcal D}(Z)\to\mathcal QCoh_{\mathcal D,Z}(S)$ is an equivalence of category whose inverse is 
given by $i^{\sharp}:=a_{\tau}i^{\sharp}:\mathcal QCoh_{\mathcal D}(S)\to\mathcal QCoh_{\mathcal D}(Z)$.
That is, for $M\in QCoh_{\mathcal D,Z}(S)$ and $N\in QCoh_{\mathcal D}(Z)$, the adjonction maps
\begin{equation*}
\ad(i_{*mod},i^{\sharp})(M):i_{*mod}i^{\sharp}M\xrightarrow{\sim} M \; , \; 
\ad(i_{*mod},i^{\sharp})(N):i^{\sharp}i_{*mod}N\xrightarrow{\sim} N
\end{equation*}
are isomorphisms.
\item[(ii)] The functor 
$\int_i=i_{*mod}:D_{\mathcal D}(Z)\to D_{\mathcal D,Z}(S)$ is an equivalence of category whose inverse is 
given by $Ri^{\sharp}:D_{\mathcal D}(S)\to D_{\mathcal D}(Z)$.
That is, for $M\in D_{\mathcal D,Z}(S)$ and $N\in D_{\mathcal D}(Z)$, the adjonction maps
\begin{equation*}
\ad(\int_i,Ri^{\sharp})(M):i_{*mod}Ri^{\sharp}M\xrightarrow{\sim} M \; , \; 
\ad(\int_i,Ri^{\sharp})(N):Ri^{\sharp}i_{*mod} N\xrightarrow{\sim} N
\end{equation*}
are isomorphisms.
\end{itemize}
\end{thm}

\begin{proof}
Similar to the complex case : see \cite{LvDmod} :(ii) follows from (i).
\end{proof}

\begin{lem}\label{imodjk}
Let $i:Z\hookrightarrow S$ a closed embedding with $Z,S\in\Var(k)$. 
Denote by $j:U:=S\backslash Z\hookrightarrow Z$ the open complementary embedding.
Then, if $i$ is a locally complete intersection embedding (e.g. if $Z,S$ are smooth),
we have for $M\in C_{O_U}(U)$ quasi-coherent, $Li^{*mod}Rj_*M=0$.
\end{lem}

\begin{proof}
Similar to the complex case : see \cite{B4}.
\end{proof}

We deduce from theorem \ref{Keqk}(i) and lemma \ref{imodjk} the localization for $D$-modules 
for a closed embedding of smooth algebraic varieties:

\begin{thm}\label{KZSk}
Let $i:Z\hookrightarrow S$ a closed embedding with $Z,S\in\SmVar(k)$. Denote by $c=\codim(Z,S)$.
Then, for $M\in C_{\mathcal D}(S)$, we have by Kashiwara equivalence the following map in $C_{\mathcal D}(S)$ :
\begin{eqnarray*}
\mathcal K_{Z/S}(M):\Gamma_ZE(M)
\xrightarrow{\ad(i_{*mod},i^{\sharp})(-)^{-1}}i_{*mod}i^{\sharp}\Gamma_ZE(M) \\
\xrightarrow{\gamma_Z(-)}i_{*mod}i^{\sharp}(E(M))
\xrightarrow{\mathcal Hom(q_K,E(i^*M))\circ\mathcal Hom(O_Z,T(i,E)(M)}
i_{*mod}K^{\vee}_{i^*O_S}(O_Z)\otimes_{i^*O_S} M 
\end{eqnarray*}
which is an equivalence Zariski local. It gives the isomorphism in $D_{\mathcal D}(S)$
\begin{equation*}
\mathcal K_{Z/S}(M):R\Gamma_ZM\to i_{*mod}K^{\vee}_{i^*O_S}(O_Z)=i_{*mod}Li^{*mod}M[c]
\end{equation*}
\end{thm}

\begin{proof}
Follows from theorem \ref{Keqk} and lemma \ref{imodjk}.
\end{proof}

Let $k$ a field of caracteristic zero. Let $S\in\SmVar(k)$.
Let $M\in\PSh_{\mathcal D,c}(S)$ a coherent $D_S$ module 
so that it admits a good filtration $(M,F)$ for the filtered ring $(D_S,F^{ord})$.
We then have the characteristic variety 
\begin{equation*}
Ch(M):=\supp(cc(\Gr^FM))\subset T_S
\end{equation*}
which is the support of the characteristic cycle $cc(\Gr^FM)\in\mathcal Z(T_S)$ of the coherent sheaf $\Gr^FM\in\Shv_c(T_S)$.
Since for two good filtration $(M,F)$ and $(M,F')$ there exist $r,s\in\mathbb Z$
satisfying $F^{'i}M\subset F^{i-r}M\subset F^{'i-s}M$ for all $i$, 
$cc(\Gr^FM)\in\mathcal Z(T_S)$ and $Ch(M)\in T_S$ does NOT depend on the choice of a good filtration $F$.

For $k\subset k'$ a subfield of characteristic zero and $S\in\SmVar(k)$, we have by definition 
\begin{equation*}
cc(\Gr^F(\pi_{k/k'}(S)^{*mod}M))=cc(\Gr^FM)\otimes_kk'\in\mathcal Z(T_{S_{k'}})
\end{equation*}
and thus 
\begin{equation*}
Ch(\pi_{k/k'}(S)^{*mod}M)=Ch(M)_{k'}\subset T_{S_{k'}} 
\end{equation*}
with $S_{k'}:=S\otimes_kk'$,
since if $(M,F)$ is a good filtration then $(\pi_{k/k'}(S)^{*mod}M,\pi_{k/k'}(S)^{*mod}F)$ is a good filtration,
$\pi_{k/k'}(S):S_{k'}:=S\otimes_kk'\to S$ being the projection (see section 2).

We have the following proposition :

\begin{prop}\label{chk}
Let $k$ a field of characteristic zero.
\begin{itemize}
\item[(i)]Let $i:Z\hookrightarrow S$ a closed embedding with $S,Z\in\SmVar(k)$.
Let $M\in\PSh_{\mathcal D,c}(Z)$ a coherent $D_Z$ module 
so that it admits a good filtration $(M,F)$ for the filtered ring $(D_S,F^{ord})$.
Then $i_{*mod}(M,F)$ is a good filtration for the filtered ring $(D_Z,F^{ord})$
and $Ch(i_{*mod}M)=di(Ch(M))\subset T_S$ where $d_i:T_Z\hookrightarrow T_{S|Z}:=p_S^{-1}(Z)\hookrightarrow T_S$ 
is the closed embedding where the first embedding is given by the differential of $i:i^*O_S\to O_Z$.
\item[(ii)]Let $S\in\SmVar(k)$.
Let $M\in\PSh_{\mathcal D,c}(S)$ a coherent $D_S$ module 
so that it admits a good filtration $(M,F)$ for the filtered ring $(D_S,F^{ord})$.
Let $Ch(M)=\cup_l C^l_M$ with $C^l_M$ the irreducible components of $Ch(M)$. 
Then $\dim(C^l_M)\leq\dim(S)$ for all $l$.
\end{itemize}
\end{prop}

\begin{proof}
\noindent(i): Similar to the proof of \cite{LvDmod}.

\noindent(ii) Let $S=\cup_iS_i$ an open affine cover. Since by definition $Ch(j_i^*M)=Ch(M)\cap p_S^{-1}(S_i)$,
it is enought to prove the result for a smooth affine variety. 
So, let $S'\in\SmVar(k)$ affine and $i:S'\hookrightarrow\mathbb A^n_k$ a closed embedding.
By (i) $Ch(i_{*mod}M)=di(Ch(M))$ so the result follows from \cite{Coutino}.
\end{proof}

\begin{defi}\label{holk}
Let $k$ a field of characteristic zero.
\begin{itemize}
\item[(i)]Let $S\in\SmVar(k)$ connected (hence irreducible since $S$ is smooth). 
A coherent $D_S$ module $M\in\PSh_{\mathcal D,c}(S)$ is called holonomic if 
all the irreducible components $C^l_M$ of $\supp(Ch(M))=\cup_l C^l_M\in T_S$ are of dimension $\dim(C^l_M)=\dim(S)$.
\item[(ii)]Let $S\in\SmVar(k)$. Then $S=\sqcup_i S_i$ with $S_i\in\SmVar(k)$ connected.
A coherent $D_S$ module $M\in\PSh_{\mathcal D,c}(S)$ is called holonomic 
if $j_i^*M\in\PSh_{\mathcal D,c}(S_i)$ is holonomic for all $i$.
\end{itemize}
Let $k\subset k'$ a subfield. Consider the projection $\pi:=\pi_{k/k'}(S):S_{k'}\to S$.
By definition $M\in\PSh_{\mathcal D,c}(S)$ is holonomic if and only if 
$\pi^{*mod}M\in\PSh_{\mathcal D,c}(S_{k'})$ is holonomic since $Ch(M)_{k'}=Ch(\pi^{*mod}M)\subset T_{S_{k'}}$.
In particular considering an embedding $\sigma:k\subset\mathbb C$, $M\in\PSh_{\mathcal D,c}(S)$ is holonomic 
if and only if $\pi_{k/\mathbb C}(S)^{*mod}M\in\PSh_{\mathcal D,c}(S_{\mathbb C})$ is holonomic.
\end{defi}

\begin{prop}\label{holksub}
Let $S\in\SmVar(k)$.
\begin{itemize}
\item[(i)] Consider an exact sequence in $\PSh_{\mathcal D,c}(S)$
\begin{equation*}
0\to M_1\to M_2\to M_3\to 0.
\end{equation*}
Then $M_2\in\PSh_{\mathcal D,h}(S)$ if and only if $M_1,M_3\in\PSh_{\mathcal D,h}(S)$.
\item[(ii)] An holonomic module $M\in\PSh_{\mathcal D,h}(S)$ has finite length.
\end{itemize}
\end{prop}

\begin{proof}
Similar to the complex case : See \cite{Coutino} or \cite{LvDmod}.
\end{proof}

Let $S\in\SmVar(k)$.
A locally free $O_S$ module with a structure of $D_S$ module is called an integrable connexion.
We denote by $\Vect_{\mathcal D}(S)\subset\PSh_{\mathcal D,c}(S)$ 
the full subcategory whose set of objects consists of integrable connexions.
By definition, an integrable connexion $M\in\Vect_{\mathcal D}(S)$ is holonomic since $Ch(M)=i_0(S)\subset T_S$
where $i_0:S\hookrightarrow T_S$, $i_0(s)=(s,0)$ is the zero section embedding.
Hence $\Vect_{\mathcal D}(S)\subset\PSh_{\mathcal D,h}(S)$.

\begin{prop}\label{holkgen}
Let $k$ a field of characteristic zero. Let $S\in\SmVar(k)$. 
\begin{itemize}
\item[(i)] A coherent $D_S$ module $M\in\PSh_{\mathcal D,c}(S)$ which is a coherent $O_S$ module is a locally free $O_S$ module.
\item[(ii)] An holonomic $D_S$ module $M\in\PSh_{\mathcal D,h}(S)$ is generically an integrable connexion,
that is there exists an open subset $j:S^o\subset S$ such that $M_{|S^o}:=j^*M\in\Vect_{\mathcal D}(S^o)$.
\end{itemize}
\end{prop}

\begin{proof}
\noindent(i): Similar to the complex case : see \cite{LvDmod}.

\noindent(ii): Similar to the complex case : 
follows from (i) since there exist an open subset $S^o\subset S$ such that $ch(M)\cap p^{-1}(S^o)=T_{S^o}S^o$
where $T_SS\subset T_S$ is the zero section.
\end{proof}

Let $k$ a field of characteristic zero. Let $S\in\SmVar(k)$
\begin{itemize}
\item we consider
\begin{itemize}
\item the full subcategories
\begin{equation*}
C_{\mathcal D,h}(S)\subset C_{\mathcal D,c}(S)\subset C_{\mathcal D}(S) \; \mbox{and} \;
D_{\mathcal D,h}(S)\subset D_{\mathcal D,c}(S)\subset D_{\mathcal D}(S) 
\end{equation*}
consisting of complexes of presheaves $M$ such that 
$a_{\tau}H^n(M)$ are coherent, resp. holonomic, sheaves of $D_S$ modules, 
$a_{\tau}$ being the sheaftification functor for the Zariski topology,
\item the full subcategories
\begin{equation*}
C_{\mathcal D^{op},h}(S)\subset C_{\mathcal D^{op},c}(S)\subset C_{\mathcal D^{op}}(S)  \; \mbox{and} \;
D_{\mathcal D^{op},h}(S)\subset D_{\mathcal D^{op},c}(S)\subset D_{\mathcal D^{op}}(S)  
\end{equation*}
the full subcategories consisting of complexes of presheaves $M$ such that 
$a_{\tau}H^n(M)$ are coherent, resp. holonomic, sheaves of right $D_S$ modules,
\end{itemize}
\item in the filtered case we have 
\begin{itemize}
\item the full subcategories
\begin{equation*}
C_{\mathcal D(2)fil,h}(S)\subset C_{\mathcal D(2)fil,c}(S)\subset C_{\mathcal D(2)fil}(S), \; \mbox{and} \;
D_{\mathcal D(2)fil,h}(S)\subset D_{\mathcal D(2)fil,c}(S)\subset D_{\mathcal D(2)fil}(S), 
\end{equation*}
consisting of filtered complexes of presheaves $(M,F)$ such that 
$a_{\tau}H^n(M,F)$ are filtered coherent, resp. filtered holonomic, sheaves of $D_S$ modules, 
that is $a_{\tau}H^n(M)$ are coherent, resp. holonomic sheaves of $D_S$ modules and 
$F$ induces a good filtration on $a_{\tau}H^n(M)$ 
(in particular $F^pa_{\tau}H^n(M)\subset a_{\tau}H^n(M)$ are coherent sub $O_S$ modules),
the full subcategories
\begin{eqnarray*}
C_{\mathcal D(1,0)fil,h}(S)=C_{\mathcal D2fil,h}(S)\cap C_{\mathcal D(1,0)fil}(S)\subset C_{\mathcal D2fil,h}(S), 
\; \mbox{and} \\
D_{\mathcal D(1,0)fil,h}(S)=D_{\mathcal D2fil,h}(S)\cap D_{\mathcal D(1,0)fil}(S)\subset D_{\mathcal D2fil,h}(S), 
\end{eqnarray*}
consisting of filtered complexes of presheaves $(M,F,W)$ 
such that $a_{\tau}H^n(M,F)$ are filtered holonomic sheaves of $D_S$ modules and such that $W^pM\subset M$ are $D_S$ submodules
(recall that the $O_S$ submodules $F^pM\subset M$ are NOT $D_S$ submodules but satisfy by definition 
$md:F^rD_S\otimes F^pM\subset F^{p+r}M$),
\item and similarly the full subcategories
\begin{equation*}
C_{\mathcal D^{op}(2)fil,h}(S)\subset C_{\mathcal D^{op}(2)fil,c}(S)\subset C_{\mathcal D^{op}(2)fil}(S), 
\end{equation*}
the full subcategories consisting of filtered complexes of presheaves $(M,F)$ 
such that $a_{\tau}H^n(M,F)$ are filtered coherent, resp. filtered holonomic, sheaves of right $D_S$ modules.
\end{itemize}
\end{itemize}
Let $S\in\Var(k)$. Let $Z\subset S$ a closed subset. Denote by $j:S\backslash Z\hookrightarrow S$ the open embedding.
We denote by $C_{\mathcal Dfil,h,Z}(S)\subset C_{\mathcal Dfil,h}(S)$ the full subcategory consisting of
$(M,F)\in C_{\mathcal Dfil,h}(S)$ such that $j^*\Gr_F^pM\in C_{O_S}(S)$ is acyclic for all $p\in\mathbb Z$.

\begin{prop}\label{compDmodh}
Let $f:X\to S$ a morphism with $X,S\in\SmVar(k)$. Then,
\begin{itemize}
\item[(i)] For $(M,F)\in C_{\mathcal D(2)fil,h}(S)$, we have $L\mathbb D_S(M,F)\in D_{\mathcal D(2)fil,h}(S)$.
\item[(ii)] For $(M,F)\in C_{\mathcal D(2)fil,h}(S)$, 
we have $Lf^{*mod}(M,F)\in D_{\mathcal D(2)fil,h}(X)$ and $Lf^{\hat*mod}(M,F)\in D_{\mathcal D(2)fil,h}(X)$.
\item[(iii)] For $M\in C_{\mathcal D,h}(X)$, we have $\int_fM\in D_{\mathcal D,h}(S)$
and $\int_{f!}M:=L\mathbb D_S\int_fL\mathbb D_X\in D_{\mathcal D,h}(S)$.
\item[(iv)] If $f$ is proper, for $(M,F)\in C_{\mathcal D(2)fil,h}(X)$, 
we have $\int_f(M,F)\in D_{\mathcal D(2)fil,h}(S)$.
\item[(v)] For $(M,F),(N,F)\in C_{\mathcal D(2)fil,h}(S)$, $(M,F)\otimes^L_{O_S}(N,F)\in D_{\mathcal D(2)fil,h}(S)$
\end{itemize}
\end{prop}

\begin{proof}
Similar to the proof of the complex case in \cite{LvDmod} or simply follows from the complex case 
since $M\in\PSh_{\mathcal D,c}(S)$ is holonomic 
if and only if $\pi_{k/\mathbb C}(S)^{*mod}M\in\PSh_{\mathcal D,c}(S_{\mathbb C})$ is holonomic.
Note that (v) follows from (ii) by proposition \ref{otimesed}.
\end{proof}

\begin{prop}\label{dmodinv}
Let $S\in\SmVar(k)$. For $M\in C_{\mathcal D,c}(S)$, 
the canonical map $d(M):M\to\mathbb D^2_SL_DM$ is an equivalence Zariski local
\end{prop}

\begin{proof}
Standard.
\end{proof}

\begin{prop}\label{compDmodhat}
Let $f_1:X\to Y$ and $f_2:Y\to S$ two morphism with $X,Y,S\in\SmVar(k)$. 
Let $M\in C_{\mathcal D,h}(S)$. 
Then, we have $L(f_2\circ f_1)^{\hat{*}mod}M=Lf_1^{\hat{*}mod}(Lf_2^{\hat{*}mod}M)$ in $D_{\mathcal D,h}(X)$.
\end{prop}

\begin{proof}
Follows from proposition \ref{compDmod} (i), proposition \ref{compDmodh} and proposition \ref{dmodinv},
or directly from the complex case.
\end{proof}

\begin{thm}\label{holkstr}
Let $k$ a field of characteristic zero. Let $S\in\SmVar(k)$. 
\begin{itemize}
\item[(i)]Let $M\in D_{\mathcal D,c}(S)$. Then $M\in D_{\mathcal D,h}(S)$ if and only if
there exist a finite sequence 
\begin{equation*}
S=S_0\supset S_1\supset\cdots\supset S_r\supset S_{r+1}=\emptyset
\end{equation*}
such that for all $l$, $S_l\backslash S_{l+1}$ is smooth 
and $H^k(i_l^{*mod}M)\in\Vect_{\mathcal D}(S_l\backslash S_{l+1})$ are integrable connexion for all $k\in\mathbb Z$, 
$i_l:S_l\backslash S_{l+1}\hookrightarrow S$ being the locally closed embedding.
\item[(ii)]Let $M\in\PSh_{\mathcal D,c}(S)$. Then $M\in\PSh_{\mathcal D,h}(S)$ if and only if
there exist a finite sequence 
\begin{equation*}
S=S_0\supset S_1\supset\cdots\supset S_r\supset S_{r+1}=\emptyset
\end{equation*}
such that for all $l$, $S_l\backslash S_{l+1}$ is smooth and 
$i_l^{*mod}M\in\Vect_{\mathcal D}(S_l\backslash S_{l+1})$ is an integrable connexion, 
$i_l:S_l\backslash S_{l+1}\hookrightarrow S$ being the locally closed embedding.
\end{itemize}
\end{thm}

\begin{proof}
\noindent(i): Similar to the complex case : follows from proposition \ref{holkgen},
theorem \ref{KZSk} and proposition \ref{compDmodh}.

\noindent(ii):It is a particular case of (i).
\end{proof}

Let $k$ a field of characteristic zero.
\begin{itemize}
\item Let $C\in\SmVar(k)$ connected (hence irreducible). 
An algebraic meromorphism connexion $M=(M,\nabla)\in\Mod(K_C,D(O_{C,s}))$ at $s\in C$
is a $K_C$ module $M$ endowed with a $k$ linear map 
$\nabla:M\to\Omega^1_{C,s}\otimes_{O_{C,s}}M$ such that $\nabla(fm)=d_f\otimes m+f\nabla(m)$ for $f\in K_C$
where $K_C:=Frac(O_{C,s})$ is the field of fraction of $C$,
\item Let $S\in\SmVar(k)$. 
An algebraic meromorphism connexion $M=(M,\nabla)\in\PSh_{O_S(*D),D_S}(S)$ along a (Cartier divisor) $D\subset S$
is a coherent $O_S(*D)$ module which has a structure of $D_S$ module. 
In particular, $M_{|S\backslash D}\in\Vect_{\mathcal D}(S\backslash D)$ is an integrable connexion
since it is a $D_{S\backslash D}$ module which is a coherent $O_{S\backslash D}$ module.
\end{itemize}

\begin{lem}
Let $S\in\SmVar(k)$. Let $D\subset S$ a (Cartier) divisor. 
Denote by $j:S^o:=S\backslash D\hookrightarrow S$ the open embedding. Then, the restriction
\begin{equation*}
j^*:\PSh_{O_S(*D),D_S}(S)\to\Vect_{\mathcal D}(S^o)
\end{equation*}
is an equivalence of category whose inverse is
\begin{equation*}
j_*:\Vect_{\mathcal D}(S^o)\to\PSh_{O_S(*D),D_S}(S).
\end{equation*}
By proposition \ref{compDmodh}, we get a full subcategory $\PSh_{O_S(*D),D_S}(S)\subset\PSh_{\mathcal D,h}(S)$.
\end{lem}

\begin{proof}
Standard fact on coherent $O_S(*D)=j_*O_{S^o}$ module.
\end{proof}

We now give the definition of the regularity of integrable connexions and holonomic $D_S$-modules on $S\in\SmVar(k)$:
We first define it for integrable connexion and holonomic $D_C$-module for $C\in\SmVar(k)$ a smooth algebraic curve over $k$.

\begin{defi}\label{regholk0}
Let $k$ a field of characteristic zero.
\begin{itemize}
\item[(i)]Let $C\in\SmVar(k)$ connected (hence irreducible). 
An algebraic meromorphism connexion $M=(M,\nabla)\in\Mod(K_C,D(O_{C,s}))$ at $s\in C$ 
is called regular if there exists a finitely generated $O_{C,s}$ module $L\subset M$ such that $M=K_{C,s}L$
and $x\nabla(L)\subset\Omega^1_{C,s}\otimes_{O_{C,s}}L$ for some local parameter $x\in O_{C,s}$.
We call such an $L\subset M$ an integral lattice.
\item[(ii)]Let $C\in\SmVar(k)$. An integrable connexion $M=(M,\nabla)\in\Vect_{\mathcal D}(C)$ is called regular
if for any smooth compactification $\bar C\in\PSmVar(k)$ of $C$ with
$j:C\hookrightarrow\bar C$ denoting the open embedding the algebraic meromorphic connexion 
\begin{equation*}
j_*M:=(j_*M,j_*\nabla)\in\PSh_{O_{\bar C}(*\bar C\backslash C),D_{\bar C}}(\bar C) 
\end{equation*}
is regular at all $s\in\bar C$, that is for all $s\in\bar C$ the algebraic meromorphic connexion 
$(j_*M)_s=((j_*M)_s,(j_*\nabla)_s)\in\Mod(K_C,D(O_{\bar C,s}))$ at $s\in\bar C$ is regular (see (i)).
\item[(iii)]Let $C\in\SmVar(k)$. Let $M\in\PSh_{\mathcal D,h}(C)$ an holonomic $D_C$ module.
Then by proposition \ref{holkgen}, there exist an open subset $l:C^o\subset C$ 
such that $M_{|C^o}:=l^*M\in\Vect_{\mathcal D}(C^o)$ is an integrable connexion. 
We say that $M$ is regular if $l^*M\in\Vect_{\mathcal D}(C^o)$ is regular (see (ii)).
\end{itemize}
\end{defi}

We have the following :

\begin{prop}\label{regkeq}
Let $C\in\SmVar(k)$ connected (hence irreducible). 
Consider an algebraic meromorphism connexion $M=(M,\nabla)\in\Mod(K_C,D(O_{C,s}))$ at $s\in C$.
The following are equivalent :
\begin{itemize}
\item[(i)] $M$ is regular
\item[(ii)] For any $m\in M$, there exists a finitly generated $O_{C,s}$ submodule $L\subset M$ such that
$x\nabla(L)\subset L$.
\item[(iii)] For any $m\in M$ there exist a polynomial $F(t)=t^m+a_1t^{m-1}+\cdots+a_m\in O_{C,s}[t]$
such that $F(x\nabla)(m)=0$.
\end{itemize}
\end{prop}

\begin{proof}
Similar to the complex case.
\end{proof}

We have then the following lemma

\begin{lem}\label{regholk0lem}
\begin{itemize}
\item[(i)]Let $k:C\to C'$ a morphism with $C,C'\in\SmVar(k)$ smooth algebraic curves.
\begin{itemize}
\item Let $M\in\PSh_{\mathcal D,h}(C)$. 
Then $M$ is regular if and only if $H^k\int_kM\in\PSh_{\mathcal D,h}(C')$ are regular for all $k$.
\item Let $N\in\PSh_{\mathcal D,h}(C')$. 
Then $M$ is regular if and only if $H^kLk^{*mod}N\in\PSh_{\mathcal D,h}(C)$ are regular for all $k$.
\end{itemize}
\item[(ii)] Let $\sigma:k\hookrightarrow\mathbb C$ an embedding. 
Let $C\in\SmVar(k)$ and $M\in\PSh_{\mathcal D,h}(C)$. 
Then $M$ is regular if and only if $\pi_{k/\mathbb C}(C)^{*mod}M\in\PSh_{\mathcal D,h}(C_{\mathbb C})$ is regular
\end{itemize}

\end{lem}

\begin{proof}
\noindent(i):Similar to the complex case : see \cite{LvDmod}.

\noindent(ii): Follows from the fact that for $l:C^o\hookrightarrow C$ an open subset such that
$l^*M\in\Vect_{\mathcal D}(C^o)$ is an integral connexion and $s\in\bar C$ with $\bar C\in\SmVar(k)$
a compactification of $C$, $j:C^o\hookrightarrow C\hookrightarrow\bar C$, 
if $L\subset (j_*l^*M)_s$ is an integral lattice then 
$\pi_{k/\mathbb C}(\bar C)^{*mod}L\subset (\pi_{k/\mathbb C}(\bar C)^{*mod}j_*l^*M)_s$
is an integral lattice, and conversely if $L'\subset(\pi_{k/\mathbb C}(\bar C)^{*mod}j_*l^*M)_s$ is an integral lattice then
$L'\cap (j_*l^*M)_s\subset (j_*l^*M)_s$ is an integral lattice the canonical map
\begin{equation*}
(j_*n_{O_{C^o}/O_{C^o_{\mathbb C}}}(l^*M))_s:(j_*l^*M)_s\hookrightarrow\pi_{k/\mathbb C}(\bar C)^{*mod}(j_*l^*M)_s, 
m\mapsto m\otimes 1
\end{equation*}
being injective since $l^*M$ is a locally free $O_{C^o}$ module.
\end{proof}

For integral connexions and holonomic $D_S$ modules on $S\in\SmVar(k)$ an algebraic variety of arbitrary dimesion over $k$,
we define it by the case of curves

\begin{defi}\label{regholk}
Let $k$ a field of characteristic zero.
\begin{itemize}
\item[(i)]Let $S\in\SmVar(k)$. 
An algebraic meromorphism connexion $M=(M,\nabla)\in\PSh_{O_S(*D),D_S}(S)$ along a (Cartier divisor) $D\subset S$
is called regular, if for all morphism $i_C:C\to X$ with $C\in\SmVar(k)$ a smooth curve and all $s=D\cap i_C(C)$, 
the meromorphic connexion $(i_C^{*mod}M,\nabla)\in\Mod(K_C,D(O_{C,s}))$ is regular (see definition \ref{regholk0}).
\item[(ii)]Let $S\in\SmVar(k)$. An integrable connexion $M=(M,\nabla)\in\Vect_{\mathcal D}(S)$ is called regular
if for any smooth compactification $\bar S\in\PSmVar(k)$ of $S$ 
with $D:=\bar S\backslash S\subset\bar S$ a (Cartier) divisor, 
the algebraic meromorphic connexion $(j_*M)=((j_*M),(j_*\nabla))\in\PSh_{O_{\bar S}(*D),D_{\bar S}}(\bar S)$ 
along $D\subset S$ is regular, where $j:S\hookrightarrow\bar S$ is the open embedding (see (i)).
\item[(iii)]Let $S\in\SmVar(k)$. An holonomic $D_S$ module $M\in\PSh_{\mathcal D,h}(S)$ is called regular
if for all morphism $i_C:C\to S$ with $C\in\SmVar(k)$, 
$i_C^{*mod}M\in\PSh_{\mathcal D,h}(C)$ is regular (see definition \ref{regholk0}).
\item[(iv)]Let $\sigma:k\hookrightarrow\mathbb C$ an embedding.
Let $S\in\SmVar(k)$ and $M\in\PSh_{\mathcal D,h}(S)$. Consider the projection $\pi_{k/\mathbb C}(S):S\to S_{\mathbb C}$.
Then by lemma \ref{regholk0lem}(ii), if $\pi_{k/\mathbb C}(S)^{*mod}M\in\PSh_{\mathcal D,h}(S_{\mathbb C})$ is regular
then $M\in\PSh_{\mathcal D,h}(S)$ is regular. 
So, let $k$ a field of characteristic zero. Let $S\in\SmVar(k)$ and $M\in\PSh_{\mathcal D,h}(S)$.
We say that $M$ is regular in the strong sense if 
\begin{equation*}
\pi_{k/\mathbb C}(S_0)^{*mod}M_0\in\PSh_{\mathcal D,h}(S_{0\mathbb C})
\end{equation*}
is regular, where $k_0\subset k$ is a subfield of finite transcandence degree over $\mathbb Q$ 
such that $S$ and $M$ are defined that is $S=S_{0k}:=S_0\otimes_{k_0}k$ with $S_0\in\SmVar(k_0)$ and 
$M=\pi_{k_0/k}(S_0)^{*mod}M_0$ with $M_0\in\PSh_{\mathcal D,h}(S_0)$,
and we take an embedding $\sigma:k_0\hookrightarrow\mathbb C$.
This definition does NOT depend on the choice of the subfield $k_0$ and the embedding $\sigma:k_0\hookrightarrow\mathbb C$.
\end{itemize}
For $S\in\SmVar(k)$, we denote $\PSh_{\mathcal D,rh}(S)\subset\PSh_{\mathcal D,h}(S)$ the full subcategory
consisting of holonomic $D_S$ modules $M\in\PSh_{\mathcal D,h}(S)$ regular in the strong sense (see (iv)). 
\end{defi}

We have then the following easy proposition :

\begin{prop}\label{holregksub}
Let $S\in\SmVar(k)$. Consider an exact sequence in $\PSh_{\mathcal D,h}(S)$
\begin{equation*}
0\to M_1\to M_2\to M_3\to 0.
\end{equation*} 
\begin{itemize}
\item[(i)]Then, $M_2$ is regular if and only if $M_1$ and $M_3$ are regular
\item[(ii)]Then, $M_2\in\PSh_{\mathcal D,rh}(S)$ if and only if $M_1,M_3\in\PSh_{\mathcal D,rh}(S)$.
\end{itemize}
\end{prop}

\begin{proof}
\noindent(i):Similar to the complex case : by definition we are reduced to the case of integrable connexions on curves.
But for $0\to M_1\to M_2\xrightarrow{q} M_3\to 0$ an exact sequence of integrable connexions on a curve $C\in\SmVar(k)$
with compactification $\bar C\in\PSmVar(k)$, $M_1$ and $M_3$ are regular at $s\in\bar C$ if and only if
$M_2$ is regular at $s$ by proposition \ref{regkeq} (use (iii)).

\noindent(ii):Follows by definition from the complex case which is a particular case of (i).
\end{proof}

Let $k$ a field of characteristic zero. Let $S\in\SmVar(k)$
\begin{itemize}
\item we consider
\begin{itemize}
\item the full subcategories
\begin{equation*}
C_{\mathcal D,rh}(S)\subset C_{\mathcal D,h}(S) \; \mbox{and} \; D_{\mathcal D,rh}(S)\subset D_{\mathcal D,h}(S) 
\end{equation*}
consisting of complexes of presheaves $M$ such that $a_{\tau}H^n(M)\in\PSh_{\mathcal D,rh}(S)$ 
(see definition \ref{regholk}), $a_{\tau}$ being the sheaftification functor for the Zariski topology,
\item the full subcategories
\begin{equation*}
C_{\mathcal D^{op},rh}(S)\subset C_{\mathcal D^{op},h}(S) \; \mbox{and} \;
D_{\mathcal D^{op},rh}(S)\subset D_{\mathcal D^{op},h}(S) 
\end{equation*}
the full subcategories consisting of complexes of presheaves $M$ such that $a_{\tau}H^n(M)^{op}\in\PSh_{\mathcal D,rh}(S)$,
\end{itemize}
\item in the filtered case we have 
\begin{itemize}
\item the full subcategories
\begin{equation*}
C_{\mathcal D(2)fil,rh}(S)\subset C_{\mathcal D(2)fil,h}(S), \; \mbox{and} \;
D_{\mathcal D(2)fil,rh}(S)\subset D_{\mathcal D(2)fil,h}(S), 
\end{equation*}
consisting of filtered complexes of presheaves $(M,F)$ such that $a_{\tau}H^n(M)\in\PSh_{\mathcal D,rh}(S)$ 
(see definition \ref{regholk}), the full subcategories
\begin{eqnarray*}
C_{\mathcal D(1,0)fil,rh}(S)=C_{\mathcal D2fil,rh}(S)\cap C_{\mathcal D(1,0)fil}(S)\subset C_{\mathcal D2fil,rh}(S), 
\; \mbox{and} \\
D_{\mathcal D(1,0)fil,rh}(S)=D_{\mathcal D2fil,rh}(S)\cap D_{\mathcal D(1,0)fil}(S)\subset D_{\mathcal D2fil,rh}(S), 
\end{eqnarray*}
consisting of filtered complexes of presheaves $(M,F,W)\in C_{\mathcal D(1,0)fil,h}(S)$ 
such that $a_{\tau}H^n(M)\in\PSh_{\mathcal D,rh}(S)$ (see definition \ref{regholk})
\item and similarly the full subcategories
\begin{equation*}
C_{\mathcal D^{op}(2)fil,rh}(S)\subset C_{\mathcal D^{op}(2)fil,h}(S), \; \mbox{and} \;
D_{\mathcal D^{op}(2)fil,rh}(S)\subset D_{\mathcal D^{op}(2)fil,h}(S)
\end{equation*}
the full subcategories consisting of filtered complexes of presheaves $(M,F)\in C_{\mathcal D^{op}(2)fil,h}(S)$ 
such that $a_{\tau}H^n(M)^{op}\in\PSh_{\mathcal D,rh}(S)$.
\end{itemize}
\end{itemize}
Let $S\in\SmVar(k)$. Let $Z\subset S$ a closed subset. Denote by $j:S\backslash Z\hookrightarrow S$ the open embedding.
We denote by $C_{\mathcal Dfil,Z,rh}(S)\subset C_{\mathcal Dfil,rh}(S)$ the full subcategory consisting of
$(M,F)\in C_{\mathcal Dfil,rh}(S)$ such that $j^*\Gr_F^pM\in C_{O_S}(S)$ is acyclic for all $p\in\mathbb Z$.

We now give an equivalent definition of regular holonomic $D_S$ modules together with a result on
stability by direct, inverse image and duality.

Let $S\in\SmVar(k)$. For each $M\in\PSh_{\mathcal D,h}(S)$ there exists by proposition \ref{holksub}(ii)
a finite sequence of holonomic submodules
\begin{equation*}
0=M_{r+1}\subset M_r\subset\cdots\subset M_1\subset M_0=M
\end{equation*}
such that $M_i/M_{i+1}\in\PSh_{\mathcal D,h}(S)$ is simple. 

\begin{defi}\label{DSholksimple}
Let $k:Z^o\hookrightarrow S$ a locally closed embedding with $S,Z^o\in\SmVar(k)$,
and assume $k$ is affine. We define for $M\in\PSh_{\mathcal D,h}(Z^o)$ the minimal extension
\begin{equation*}
L_{Z^o/S}(M):=T(k_!,k_*)(\int_{k!}M)\subset\int_kM
\end{equation*}
where $T(k_!,k_*)(M):\int_{k!}M\to\int_kM$ is given by, 
using a factorization of $k$ by open embeddings and proper morphisms
\begin{itemize}
\item the adjonction map $T(j_!,j_*)(N):=\ad(j^*,j_*)(j_!N):j_!N\to j_*N$ for open embeddings $j:X^o\hookrightarrow X$
with $X\in\SmVar(k)$,
\item the trace map on proper morphisms.
\end{itemize}
By proposition \ref{holksub}(i), $L_{Z^o/S}(M)\in\PSh_{\mathcal D,h}(S)$ is holonomic.
\end{defi}

\begin{thm}\label{DSholksimplethm}
\begin{itemize}
\item[(i)]Let $k:Z^o\hookrightarrow S$ a locally closed embedding with $S,Z^o\in\SmVar(k)$,
and assume $k$ is affine. Let $M\in\PSh_{\mathcal D,h}(Z^o)$. 
If $M$ is simple, then $L_{Z^o/S}(M)$ is also simple, and is the unique simple submodule of $\int_kM$ and
the unique quotient module of $\int_{k!}M$.
\item[(ii)]Let $S\in\SmVar(k)$. Let $M\in\PSh_{\mathcal D,h}(S)$. If $M$ is a simple $D_S$ module
then there exist $k:Z^o\hookrightarrow S$ a locally closed embedding with $Z\in\SmVar(k)$, $k$ affine,
such that $M\simeq L_{Z^o/S}(N)$ with $N\in\Vect_{\mathcal D}(Z^o)$ a simple integral connexion.
\item[(iii)]Let $k:Z^o\hookrightarrow S$, $k':Z^{'o}$ locally closed embeddings with $S,Z^o,Z^{'o}\in\SmVar(k)$,
$k,k'$ affine. Let $N\in\Vect_{\mathcal D}(Z^o)$ and $N\in\Vect_{\mathcal D}(Z^{'o})$ simple integral connexions.
Then $L_{Z^o/S}(N)\simeq L_{Z^{'o}/S}(N')$ in $\PSh_{\mathcal D}(S)$ if and only if 
$\bar Z^o=\bar Z^{'o}$ and $N_{|U}\simeq N'_{|U}$ for an open dense subset $U\subset Z^o\cap Z^{'o}$.
\end{itemize}
\end{thm}

\begin{proof}
Similar to the complex case : see \cite{LvDmod}.
\end{proof}

\begin{thm}\label{holregkthmeq}
Let $S\in\SmVar(k)$. Let $M\in\PSh_{\mathcal D,h}(S)$. Take by proposition \ref{holksub}(ii)
a finite sequence of holonomic submodules
\begin{equation*}
0=M_{r+1}\subset M_r\subset\cdots\subset M_1\subset M_0=M
\end{equation*}
such that $M_i/M_{i+1}\in\PSh_{\mathcal D,h}(S)$ is simple. By theorem \ref{DSholksimplethm}
there exist locally closed embeddings $k_i:Z_i^o\hookrightarrow S$ with $Z_i^o\in\SmVar(k)$ 
and $N_i\in\Vect_{\mathcal D}(Z_i^o)$ simple integrable connexion such that
$M_i\simeq L_{S_i^o/S}(N_i)$ in $\PSh_{\mathcal D}(S)$.
Then $M$ is regular if and only if 
the simple integral connexions $N_i\in\Vect_{\mathcal D}(Z_i^o)$ are regular (see definition \ref{regholk}).
\end{thm}

\begin{proof}
Similar to the proof of the complex case in \cite{LvDmod}.
\end{proof}

\begin{thm}\label{compDmodrhthm}
Let $f:X\to S$ a morphism with $X,S\in\SmVar(k)$. Then,
\begin{itemize}
\item[(i)] For $(M,F)\in C_{\mathcal D(2)fil,rh}(S)$, we have $L\mathbb D_S(M,F)\in D_{\mathcal D(2)fil,rh}(S)$.
\item[(ii)] For $(M,F)\in C_{\mathcal D(2)fil,rh}(S)$, 
we have $Lf^{*mod}(M,F)\in D_{\mathcal D(2)fil,rh}(X)$ and $Lf^{\hat*mod}(M,F)\in D_{\mathcal D(2)fil,rh}(X)$.
\item[(iii)] For $M\in C_{\mathcal D,rh}(X)$, we have $\int_fM\in D_{\mathcal D,rh}(S)$.
and $\int_{f!}M:=L\mathbb D_S\int_fL\mathbb D_X\in D_{\mathcal D,rh}(S)$.
\item[(iv)] If $f$ is proper, for $(M,F)\in C_{\mathcal D(2)fil,rh}(X)$, 
we have $\int_f(M,F)\in D_{\mathcal D(2)fil,rh}(S)$.
\item[(v)] For $(M,F),(N,F)\in C_{\mathcal D(2)fil,rh}(S)$, $(M,F)\otimes^L_{O_S}(N,F)\in D_{\mathcal D(2)fil,rh}(S)$
\end{itemize}
\end{thm}

\begin{proof}
Follows by definition from the complex case :
\noindent(i),(ii) and (iii): See \cite{LvDmod}.

\noindent(iv): Follows from (iii) and stability of coherent $O_X$-modules by direct image of proper morphism $f:X\to S$.

\noindent(v):Follows from (ii) by proposition \ref{otimesed}.
\end{proof}

\subsection{The D modules on singular algebraic varieties over a field $k$ of characteristic zero}

In this subsection by defining the category of complexes of filtered D-modules in the singular case and
there functorialities. 

\subsubsection{Definition}

In all this subsection, we fix the notations:
Let $k$ a field of characteristic zero. For $S\in\Var(k)$, we denote by $S=\cup_{i=1}^lS_i$ an open cover
such that there exits closed embeddings $i_iS_i\hookrightarrow\tilde S_i$ with $\tilde S_i\in\SmVar(k)$.
We have then closed embeddings $i_I:S_I:=\cap_{i\in I} S_i\hookrightarrow\tilde S_I:=\Pi_{i\in I}\tilde S_I$.
Then for $I\subset J$, we denote by $j_{IJ}:S_J\hookrightarrow S_I$ the open embedding 
and $p_{IJ}:\tilde S_J\to\tilde S_I$ the projection, so that $p_{IJ}\circ i_J=i_I\circ j_{IJ}$.
This gives the diagram of algebraic varieties $(\tilde S_I)\in\Fun(\mathcal P(\mathbb N),\Var(k))$ which
gives the diagram of sites $(\tilde S_I):=\Ouv(\tilde S_I)\in\Fun(\mathcal P(\mathbb N),\Cat)$.
It also gives the diagram of sites $(\tilde S_I)^{op}:=\Ouv(\tilde S_I)^{op}\in\Fun(\mathcal P(\mathbb N),\Cat)$.
For $I\subset J$, we denote by $m:\tilde S_I\backslash(S_I\backslash S_J)\hookrightarrow\tilde S_I$ the open embedding.

\begin{defi}
Let $S\in\Var(k)$ and let $S=\cup_iS_i$ an open cover
such that there exist closed embeddings $i_iS_i\hookrightarrow\tilde S_i$ with $\tilde S_i\in\SmVar(k)$. 
Then, $\PSh_{\mathcal D(2)fil}(S/(\tilde S_I))\subset\PSh_{\mathcal D(2)fil}((\tilde S_I))$ is the full subcategory
\begin{itemize}
\item whose objects are $(M,F)=((M_I,F)_{I\subset\left[1,\cdots l\right]},s_{IJ})$, with
\begin{itemize}
\item $(M_I,F)\in\PSh_{\mathcal D(2)fil}(\tilde S_I)$ such that $\mathcal I_{S_I}M_I=0$, 
in particular $(M_I,F)\in\PSh_{\mathcal D(2)fil,S_I}(\tilde S_I)$
\item $s_{IJ}:m^*(M_I,F)\xrightarrow{\sim} m^*p_{IJ*}(M_J,F)[d_{\tilde S_I}-d_{\tilde S_J}]$ 
for $I\subset J$, are isomorphisms, $p_{IJ}:\tilde S_J\to\tilde S_I$ being the projection,
satisfying for $I\subset J\subset K$, $p_{IJ*}s_{JK}\circ s_{IJ}=s_{IK}$ ; 
\end{itemize}
\item the morphisms $m:(M,F)\to(N,F)$ between  
$(M,F)=((M_I,F)_{I\subset\left[1,\cdots l\right]},s_{IJ})$ and $(N,F)=((N_I,F)_{I\subset\left[1,\cdots l\right]},r_{IJ})$
are by definition a family of morphisms of complexes,  
\begin{equation*}
m=(m_I:(M_I,F)\to (N_I,F))_{I\subset\left[1,\cdots l\right]}
\end{equation*}
such that $r_{IJ}\circ m_J=p_{IJ*}m_J\circ s_{IJ}$ in $C_{\mathcal D,S_J}(\tilde S_J)$.
\end{itemize}
We denote by 
\begin{equation*}
\PSh_{\mathcal D(2)fil,rh}(S/(\tilde S_I))\subset
\PSh_{\mathcal D(2)fil,h}(S/(\tilde S_I))\subset\PSh_{\mathcal D(2)fil,c}(S/(\tilde S_I))
\subset\PSh_{\mathcal D(2)fil}(S/(\tilde S_I))
\end{equation*}
the full subcategory consisting of $((M_I,F),s_{IJ})\in\PSh_{\mathcal D(2)fil}(S/(\tilde S_I))$ 
such that for all $I\subset[1,\ldots,l]$, $(M_I,F)\in\PSh_{\mathcal D(2)fil,c}(\tilde S_I)$ 
resp. $(M_I,F)\in\PSh_{\mathcal D(2)fil,h}(\tilde S_I)$,
resp. $(M_I,F)\in\PSh_{\mathcal D(2)fil,h}(\tilde S_I)$ and $M_I\in\PSh_{\mathcal D,rh}(\tilde S_I)$
(see definition \ref{regholk})
We have the full subcategories
\begin{eqnarray*} 
\PSh_{\mathcal D(1,0)fil,rh}(S/(\tilde S_I))\subset\PSh_{\mathcal D2fil,rh}(S/(\tilde S_I)), \;
\PSh_{\mathcal D(1,0)fil,h}(S/(\tilde S_I))\subset\PSh_{\mathcal D2fil,h}(S/(\tilde S_I)), \\
\PSh_{\mathcal D(1,0)fil,h}(S/(\tilde S_I))\subset\PSh_{\mathcal D2fil,h}(S/(\tilde S_I)), 
\end{eqnarray*}
consisting of $((M_I,F,W),s_{IJ})$ such that $W^pM_I$ are $D_{\tilde S_I}$ submodules.
\end{defi}
We recall from \cite{B4} the following
\begin{itemize}
\item A morphism $m=(m_I):((M_I),s_{IJ})\to((N_I),r_{IJ})$ in $C(\PSh_{\mathcal D}(S/(\tilde S_I)))$ 
is a Zariski, resp. usu, local equivalence if and only if all the $m_I$ are Zariski local equivalences.
\item A morphism $m=(m_I):((M_I,F),s_{IJ}\to((N_I,F),r_{IJ}))$ in $C(\PSh_{\mathcal D(2)fil}(S/(\tilde S_I)))$
is a filtered Zariski local equivalence 
if and only if all the $m_I$ are filtered Zariski local equivalence.
\item By definition, a morphism $m=(m_I):((M_I,F),s_{IJ})\to((N_I,F),r_{IJ})$ in $C(\PSh_{\mathcal D(2)fil}(S/(\tilde S_I)))$
is an $r$-filtered Zariski local equivalence 
if there exist $m_i:((C_{iI},F),s_{iIJ})\to((C_{(i+1)I},F),s_{(i+1)IJ})$, $0\leq i\leq s$, 
with $((C_{iI},F),s_{iIJ})\in C(\PSh_{\mathcal D(2)fil}(S/(\tilde S_I)))$, 
$((C_{0I},F),s_{iIJ})=((M_I,F),s_{IJ})$, $((C_{sI},F),s_{sIJ})=((N_I,F),r_{IJ})$ such that
\begin{equation*}
m=m_s\circ\cdots\circ m_i\circ\cdots\circ m_0:((M_I,F),s_{IJ}\to((N_I,F),r_{IJ}))
\end{equation*}
with $m_i:((C_{iI},F),s_{iIJ})\to((C_{(i+1)I},F),s_{(i+1)IJ})$ either filtered Zariski  local equivalence
or $r$-filtered homotopy equivalence.
\end{itemize}

\begin{defiprop}\label{defprops} 
Let $S\in\Var(k)$ and let $S=\cup_iS_i$ an open cover
such that there exist closed embeddings $i_iS_i\hookrightarrow\tilde S_i$ with $\tilde S_i\in\SmVar(k)$.
Then $\PSh_{\mathcal D(2)fil}(S/(\tilde S_I))$ does not depend on the open covering of $S$ and the closed embeddings and we set
\begin{equation*}
\PSh_{\mathcal D(2)fil}(S):=\PSh_{\mathcal D(2)fil}(S/(\tilde S_I))
\end{equation*}
We denote by $C^0_{\mathcal D(2)fil}(S):=C(\PSh_{\mathcal D(2)fil}(S/(\tilde S_I)))$ and by 
$D^0_{\mathcal D(2)fil,r}(S):=K^0_{\mathcal D(2)fil,r}(S)([E_1]^{-1})$
the localization of the $r$-filtered homotopy category with respect to the classes of filtered Zariski local equivalences.
\end{defiprop}

\begin{proof}
Similar to the complex case : see \cite{B4}.
\end{proof}

We now give the definition of our category :

\begin{defi}\label{Dmodsingdef}
Let $S\in\Var(k)$ and let $S=\cup_iS_i$ an open cover
such that there exist closed embeddings $i_i:S_i\hookrightarrow\tilde S_i$ with $\tilde S_i\in\SmVar(k)$. 
Then, $C_{\mathcal D(2)fil}(S/(\tilde S_I))\subset C_{\mathcal D(2)fil}((\tilde S_I))$ is the full subcategory
\begin{itemize}
\item whose objects are $(M,F)=((M_I,F)_{I\subset\left[1,\cdots l\right]},u_{IJ})$, with
\begin{itemize}
\item $(M_I,F)\in C_{\mathcal D(2)fil,S_I}(\tilde S_I)$ that is $(M_I,F)\in C_{\mathcal D(2)fil}(\tilde S_I)$
satisfy $n_I^*\Gr_F^pM_I\in C_O(\tilde S_I)$ is acyclic for all $p\in\mathbb Z$, 
where $n_I:\tilde S_I\backslash S_I\hookrightarrow\tilde S_I$ is the open embedding,
\item $u_{IJ}:m^*(M_I,F)\to m^*p_{IJ*}(M_J,F)[d_{\tilde S_I}-d_{\tilde S_J}]$ 
for $J\subset I$, are morphisms, $p_{IJ}:\tilde S_J\to\tilde S_I$ being the projection, 
satisfying for $I\subset J\subset K$, $p_{IJ}*u_{JK}\circ u_{IJ}=u_{IK}$ in $C_{\mathcal Dfil}(\tilde S_I)$ ;
\end{itemize}
\item the morphisms $m:((M_I,F),u_{IJ})\to((N_I,F),v_{IJ})$ between  
$(M,F)=((M_I,F)_{I\subset\left[1,\cdots l\right]},u_{IJ})$ and $(N,F)=((N_I,F)_{I\subset\left[1,\cdots l\right]},v_{IJ})$
being a family of morphisms of complexes,  
\begin{equation*}
m=(m_I:(M_I,F)\to (N_I,F))_{I\subset\left[1,\cdots l\right]}
\end{equation*}
such that $v_{IJ}\circ m_I=p_{IJ*}m_J\circ u_{IJ}$ in $C_{\mathcal Dfil}(\tilde S_I)$.
\end{itemize}
We denote by $C^{\sim}_{\mathcal D(2)fil}(S/(\tilde S_I))\subset C_{\mathcal D(2)fil}(S/(\tilde S_I))$ the full subcategory 
consisting of objects $((M_I,F),u_{IJ})$ such that the $u_{IJ}$ are filtered Zariski local equivalences.
\end{defi}

Let $S\in\Var(k)$ and let $S=\cup_{i=1}^lS_i$ an open cover
such that there exist closed embeddings $i_i:S_i\hookrightarrow\tilde S_i$ with $\tilde S_i\in\SmVar(k)$. Then,
We denote by
\begin{equation*}
C_{\mathcal D(2)fil,rh}(S/(\tilde S_I))\subset 
C_{\mathcal D(2)fil,h}(S/(\tilde S_I))\subset C_{\mathcal D(2)fil,c}(S/(\tilde S_I))
\subset C_{\mathcal D(2)fil}(S/(\tilde S_I))
\end{equation*}
the full subcategories consisting of those $((M_I,F),u_{IJ})\in C_{\mathcal D(2)fil}(S/(\tilde S_I))$ 
such that for all $I\subset[1,\ldots,l]$, 
$(M_I,F)\in C_{\mathcal D(2)fil,S_I,c}(\tilde S_I)$,
that is such that $a_{\tau}H^n(M_I,F)\in\PSh_{\mathcal Dfil,c}(\tilde S_I)$ 
are coherent endowed with a good filtration for all $n\in\mathbb Z$, 
resp. $(M_I,F)\in C_{\mathcal D(2)fil,S_I,h}(\tilde S_I)$, 
that is such that $a_{\tau}H^n(M_I,F)\in\PSh_{\mathcal Dfil,h}(\tilde S_I)$ are filtered holonomic for all $n\in\mathbb Z$, 
resp. such that $(M_I,F)\in C_{\mathcal D(2)fil,S_I,rh}(\tilde S_I)$, 
that is such that $a_{\tau}H^n(M_I,F)\in\PSh_{\mathcal Dfil,h}(\tilde S_I)$ are filtered holonomic for all $n\in\mathbb Z$ 
and $a_{\tau}H^nM_I\in\PSh_{\mathcal D,rh}(\tilde S_I)$(see definition \ref{regholk}).

We denote by
\begin{eqnarray*}
C_{\mathcal D(1,0)fil,h}(S/(\tilde S_I))\subset C_{\mathcal D2fil,h}(S/(\tilde S_I)), \;
C_{\mathcal D(1,0)fil,rh}(S/(\tilde S_I))\subset C_{\mathcal D2fil,rh}(S/(\tilde S_I)), \\
C_{\mathcal D(1,0)fil}(S/(\tilde S_I))\subset C_{\mathcal D2fil}(S/(\tilde S_I)), 
\end{eqnarray*}
the full subcategories consisting of those $((M_I,F,W),u_{IJ})\in C_{\mathcal D2fil}(S/(\tilde S_I))$
such that $W^pM_I$ are $D_{\tilde S_I}$ submodules.

We recall from \cite{B4} the following
\begin{itemize}
\item A morphism $m=(m_I):(M_I,u_{IJ})\to(N_I,v_{IJ})$ in $C_{\mathcal D}(S/(\tilde S_I))$ 
is a Zariski local equivalence if and only if all the $m_I$ are Zariski local equivalences.
\item A morphism $m=(m_I):((M_I,F),u_{IJ}\to((N_I,F),v_{IJ}))$ in $C_{\mathcal D(2)fil}(S/(\tilde S_I))$
is a filtered Zariski local equivalence if and only if all the $m_I$ are filtered Zariski local equivalence.
\item Let $r=1,\cdots,\infty$.
By definition, a morphism $m=(m_I):((M_I,F),u_{IJ})\to((N_I,F),v_{IJ})$ in $C_{\mathcal D(2)fil}(S/(\tilde S_I))$
is an $r$-filtered Zariski local equivalence 
if there exist $m_i:(C_{iI},F),u_{iIJ})\to(C_{(i+1)I},F),u_{(i+1)IJ})$, $0\leq i\leq s$, 
with $(C_{iI},F),u_{iIJ})\in C_{\mathcal D(2)fil}(S/(\tilde S_I))$, 
$(C_{0I},F),u_{iIJ})=(M_I,F),u_{IJ})$, $(C_{sI},F),u_{sIJ})=(N_I,F),v_{IJ})$ such that
\begin{equation*}
m=m_s\circ\cdots\circ m_i\circ\cdots\circ m_0:((M_I,F),u_{IJ}\to((N_I,F),v_{IJ}))
\end{equation*}
with $m_i:(C_{iI},F),u_{iIJ})\to(C_{(i+1)I},F),u_{(i+1)IJ})$ either filtered Zariski local equivalence
or $r$-filtered homotopy equivalence (i.e. $r$-filtered for the first filtration and filtered for the second filtration).
\end{itemize}

\begin{defi}\label{defpropK} 
Let $S\in\Var(k)$ and let $S=\cup_iS_i$ an open cover
such that there exist closed embeddings $i_i:S_i\hookrightarrow\tilde S_i$ with $\tilde S_i\in\SmVar(k)$. 
\begin{itemize}
\item[(i)]We have the derived category
\begin{equation*}
D_{\mathcal D(2)fil}(S/(\tilde S_I)):=C_{\mathcal D(2)fil}(S/(\tilde S_I))^{\sim}([E_1]^{-1})
\end{equation*}
the localization with respect to the classes of filtered Zariski local equivalences, together with
the localization functor
\begin{equation*}
D(zar):C_{\mathcal D(2)fil}(S/(\tilde S_I))^{\sim}\to K_{\mathcal D(2)fil}(S/(\tilde S_I))
\to D_{\mathcal D(2)fil}(S/(\tilde S_I)).
\end{equation*}
\item[(ii)]We have the full subcategories
\begin{equation*}
D_{\mathcal D(1,0)fil,rh}(S/(\tilde S_I))\subset D_{\mathcal D(1,0)fil,h}(S/(\tilde S_I))
\subset D_{\mathcal D2fil}(S/(\tilde S_I)) 
\end{equation*}
which are the image of $C_{\mathcal D(1,0)fil,h}(S/(\tilde S_I))$, 
resp. of $C_{\mathcal D(1,0)fil,rh}(S/(\tilde S_I))$, by the localization functor
$D(zar):C_{\mathcal D(2)fil}(S/(\tilde S_I))\to D_{\mathcal D(2)fil}(S/(\tilde S_I))$.
\item[(iii)]We have, for $r=1,\ldots,\infty$, the $r$-filtered homotopy category
\begin{equation*}
K_{\mathcal D(2)fil,r}(S/(\tilde S_I)):=\Ho_{r}(C_{\mathcal D(2)fil}(S/(\tilde S_I)))
\end{equation*}
whose objects are those of $C_{\mathcal D(2)fil}(S/(\tilde S_I))$ and
whose morphism are $r$-filtered homotopy classes of morphisms 
($r$-filtered for the first filtration and filtered for the second), and
\begin{equation*}
D_{\mathcal D(2)fil,r}(S/(\tilde S_I)):=K_{\mathcal D(2)fil,r}(S/(\tilde S_I)))([E_1]^{-1})
\end{equation*}
the localization with respect to the classes of filtered Zariski local equivalences, together with
the localization functor
\begin{equation*}
D(zar):C_{\mathcal D(2)fil}(S/(\tilde S_I))\to K_{\mathcal D(2)fil,r}(S/(\tilde S_I))
\to D_{\mathcal D(2)fil,r}(S/(\tilde S_I)).
\end{equation*}
\item[(iv)]We have 
\begin{equation*}
D_{\mathcal D(1,0)fil,\infty,h}(S/(\tilde S_I))\subset D_{\mathcal D2fil,\infty,h}(S/(\tilde S_I))
\subset D_{\mathcal D2fil,\infty}(S/(\tilde S_I)) 
\end{equation*}
the full subcategories which are the image of $C_{\mathcal D2fil,h}(S/(\tilde S_I))$, 
resp. of $C_{\mathcal D(1,0)fil,rh}(S/(\tilde S_I))$, by the localization functor
$D(zar):C_{\mathcal D(2)fil}(S/(\tilde S_I))\to D_{\mathcal D(2)fil,\infty}(S/(\tilde S_I))$.
\end{itemize}
\end{defi}

Let $S\in\Var(k)$ and let $S=\cup_iS_i$ an open cover
such that there exist closed embeddings $i_i:S_i\hookrightarrow\tilde S_i$ with $\tilde S_i\in\SmVar(k)$. 
\begin{itemize}
\item We denote by
\begin{equation*}
C_{\mathcal D(2)fil}(S/(\tilde S_I))^0\subset C_{\mathcal D(2)fil}(S/(\tilde S_I)) \; \mbox{and} \;
D_{\mathcal D(2)fil}(S/(\tilde S_I))^0\subset D_{\mathcal D(2)fil}(S/(\tilde S_I))
\end{equation*}
the full subcategories consisting of $((M_I,F),u_{IJ})\in C_{\mathcal D(2)fil}(S/(\tilde S_I))$ such that 
\begin{equation*}
H^n((M_I,F),u_{IJ})=(H^n(M_I,F),H^nu_{IJ})\in\PSh^0_{\mathcal D(2)fil}(S/(\tilde S_I))
\end{equation*}
that is such that the $H^nu_{IJ}$ are isomorphism.
\item We have the full embedding functor
\begin{eqnarray*}
\iota^0_{S/(\tilde S_I)}:C^0_{\mathcal D(2)fil}(S):=C^0_{\mathcal D(2)fil}(S/(\tilde S_I))
\hookrightarrow C_{\mathcal D(2)fil}(S/(\tilde S_I))^0
\hookrightarrow C_{\mathcal D(2)fil}(S/(\tilde S_I)), \\
((M_I,F),s_{IJ})\mapsto ((M_I,F),s_{IJ})
\end{eqnarray*}
This full embedding induces in the derived category the functors
\begin{eqnarray*}
\iota^0_{S/(\tilde S_I)}:D^0_{\mathcal D(2)fil,r}(S):=\Ho_{zar}(C^0_{\mathcal D(2)fil,\infty}(S/(\tilde S_I)))\to 
D_{\mathcal D(2)fil}(S/(\tilde S_I))^0\hookrightarrow D_{\mathcal D(2)fil,r}(S/(\tilde S_I)), \\
((M_I,F),s_{IJ})\mapsto ((M_I,F),s_{IJ}).
\end{eqnarray*}
We can show that this functor is a full embedding.
\end{itemize}

\subsubsection{Duality in the singular case}

The definition of Saito's category comes with a dual functor :

\begin{defi}\label{dualsing}
Let $S\in\Var(k)$ and let $S=\cup S_i$ an open cover such that there exist closed embeddings
$i_i:S_i\hookrightarrow\tilde S_i$ with $\tilde S_i\in\SmVar(k)$.
We have the dual functor : 
\begin{equation*}
\mathbb D_S^K:C^0_{\mathcal Dfil}(S/(\tilde S_I))\to C^0_{\mathcal Dfil}(S/(\tilde S_I)), \;
((M_I,F),s_{IJ})\mapsto(\mathbb D^K_{\tilde S_I}(M_I,F),s^d_{IJ}),
\end{equation*}
with, denoting for short $d_{IJ}:=d_{\tilde S_J}-d_{\tilde S_I}$,
\begin{eqnarray*}
u^q_{IJ}:\mathbb D^K_{\tilde S_I}(M_I,F)\xrightarrow{\mathbb D^K(s_{IJ}^{-1})}
\mathbb D^K_{\tilde S_I}p_{IJ*}(M_J,F)[d_{IJ}]\xrightarrow{T_*(p_{IJ},D)(-)}
p_{IJ*}\mathbb D^K_{\tilde S_J}(M_J,F)[d_{IJ}]
\end{eqnarray*}
It induces in the derived category the functor
\begin{equation*}
L\mathbb D_S^K:D^0_{\mathcal Dfil}(S/(\tilde S_I))\to D^0_{\mathcal Dfil}(S/(\tilde S_I)), \;
((M_I,F),s_{IJ})\mapsto\mathbb D^K_SQ((M_I,F),s_{IJ}),
\end{equation*}
with $q:Q((M_I,F),s_{IJ})\to((M_I,F),s_{IJ})$ a projective resolution.
\end{defi}

For our definition, we have

\begin{defi}\label{dualsingour}
Let $S\in\Var(\mathbb C)$ and let $S=\cup S_i$ an open cover such that there exist closed embeddings
$i_i:S_i\hookrightarrow\tilde S_i$ with $\tilde S_i\in\SmVar(\mathbb C)$.
We have the dual functors : 
\begin{equation*}
\mathbb D_S:C_{\mathcal Dfil}(S/(\tilde S_I))\to C_{\mathcal Dfil}(S/(\tilde S_I)^{op}), \;
((M_I,F),s_{IJ})\mapsto(\mathbb D_{\tilde S_I}(M_I,F),s^d_{IJ}),
\end{equation*}
with, denoting for short $d_{IJ}:=d_{\tilde S_J}-d_{\tilde S_I}$,
\begin{eqnarray*}
u^q_{IJ}:p_{IJ*}\mathbb D_{\tilde S_J}(M_J,F)[d_{IJ}]
\end{eqnarray*}
and
\begin{equation*}
\mathbb D_S:C_{\mathcal Dfil}(S/(\tilde S_I)^{op})\to C_{\mathcal Dfil}(S/(\tilde S_I)), \;
((M_I,F),s_{IJ})\mapsto(\mathbb D_{\tilde S_I}(M_I,F),s^d_{IJ}),
\end{equation*}
with, denoting for short $d_{IJ}:=d_{\tilde S_J}-d_{\tilde S_I}$,
\begin{eqnarray*}
u^q_{IJ}:\mathbb D_{\tilde S_I}(M_I,F)
\end{eqnarray*}
It induces in the derived category the functors
\begin{equation*}
L\mathbb D_S:D_{\mathcal Dfil,r}(S/(\tilde S_I))\to D_{\mathcal Dfil,r}(S/(\tilde S_I)^{op}), \;
((M_I,F),s_{IJ})\mapsto\mathbb D_SQ((M_I,F),s_{IJ}),
\end{equation*}
with $q:Q((M_I,F),s_{IJ})\to((M_I,F),s_{IJ})$ a projective resolution,
and
\begin{equation*}
L\mathbb D_S:D_{\mathcal Dfil,r}(S/(\tilde S_I)^{op})\to D_{\mathcal Dfil,r}(S/(\tilde S_I)), \;
((M_I,F),s_{IJ})\mapsto\mathbb D_SQ((M_I,F),s_{IJ}),
\end{equation*}
with $q:Q((M_I,F),s_{IJ})\to((M_I,F),s_{IJ})$ a projective resolution.
\end{defi}

\subsubsection{Inverse image in the singular case}

We give in this subsection the inverse image functors between our categories.

Let $n:S^o\hookrightarrow S$ be an open embedding with $S\in\Var(k)$ and
let $S=\cup_i S_i$ an open cover such that there exist closed embeddings
$i_i:S_i\hookrightarrow\tilde S_i$ with $\tilde S_i\in\SmVar(k)$. 
Denote $S^o_I:=n^{-1}(S_I)=S_I\cap S^o$ and $n_I:=n_{|S^o_I}:S^o_I\hookrightarrow S^o$ the open embeddings. 
Consider open embeddings $\tilde n_I:\tilde S^o_I\hookrightarrow\tilde S_I$ such that $\tilde S^o_I\cap S_I=S^o_I$,
that is which are lift of $n_I$. We have the functor
\begin{eqnarray*}
n^*:C_{\mathcal Dfil}(S/(\tilde S_I))\to C_{\mathcal Dfil}(S^o/(\tilde S^o_I)), \\ 
(M,F)=((M_I,F),u_{IJ})\mapsto n^*(M,F):=(\tilde n_I)^*(M,F):=(\tilde n_I^*(M_I,F),n^*u_{IJ})
\end{eqnarray*}
which derive trivially.

Let $f:X\to S$ be a morphism, with $X,S\in\Var(k)$, 
such that there exist a factorization $f;X\xrightarrow{l}Y\times S\xrightarrow{p_S}S$
with $Y\in\SmVar(k)$, $l$ a closed embedding and $p_S$ the projection,
and consider $S=\cup_{i=1}^lS_i$ an open cover such that there exist closed embeddings 
$i_i:S_i\hookrightarrow\tilde S_i$, with $\tilde S_i\in\SmVar(k)$ ;
The (graph) inverse image functors is 
\begin{eqnarray*}
f^{*mod[-],\Gamma}:C_{\mathcal Dfil}(S/(\tilde S_I))\to C_{\mathcal Dfil}(X/(Y\times\tilde S_I)), \\  
(M,F)=((M_I,F),u_{IJ})\mapsto  
f^{*mod[-],\Gamma}(M,F):=(\Gamma_{X_I}E(p_{\tilde S_I}^{*mod[-]}(M_I,F)),\tilde f_J^{*mod[-]}u_{IJ})
\end{eqnarray*}
with $\tilde f_J^{*mod[-]}u_{IJ}$ as in \cite{B4},
It induces in the derived categories the functor
\begin{eqnarray*}
Rf^{*mod[-],\Gamma}:D_{\mathcal D(2)fil,r}(S/(\tilde S_I))\to D_{\mathcal D(2)fil,r}(X/(Y\times\tilde S_I)), \\   
(M,F)=((M_I,F),u_{IJ})\mapsto  
f^{*mod[-],\Gamma}(M,F):=(\Gamma_{X_I}E(p_{\tilde S_I}^{*mod[-]}(M_I,F)),\tilde f_J^{*mod[-]}u_{IJ}).
\end{eqnarray*}
It gives by duality the functor
\begin{eqnarray*}
Lf^{\hat*mod[-],\Gamma}:D_{\mathcal D(2)fil,r}(S/(\tilde S_I))\to D_{\mathcal D(2)fil,r}(X/(Y\times\tilde S_I)), \\  
(M,F)=((M_I,F),u_{IJ})\mapsto 
Lf^{\hat*mod[-],\Gamma}(M,F):=L\mathbb D_SRf^{*mod[-],\Gamma}L\mathbb D_S(M,F).
\end{eqnarray*}

The following proposition is easy :

\begin{prop}\label{compDImod}
Let $f_1:X\to Y$ and $f_2:Y\to S$ two morphism with $X,Y,S\in\Var(k)$. 
Assume there exist factorizations $f_1:X\xrightarrow{l_1} Y'\times Y\xrightarrow{p_Y} Y$ and 
$f_2:Y\xrightarrow{l_2} Y''\times S\xrightarrow{p_S} S$ with $Y',Y''\in\SmVar(k)$,
$l_1,l_2$ closed embeddings and $p_S,p_Y$ the projections. We have then the factorization 
\begin{equation*}
f_2\circ f_1:X\xrightarrow{(l_2\circ I_{Y'})\circ l_1}Y'\times Y''\times S\xrightarrow{p_S} S.
\end{equation*}
We have, for $(M,F)\in C^{\sim}_{\mathcal D(2)fil}(S/(\tilde S_I))$,
$R(f_2\circ f_1)^{*mod[-],\Gamma}(M,F)=Rf_{2}^{*mod[-],\Gamma}\circ Rf_{1}^{*mod[-],\Gamma}(M,F)$.
\end{prop}

\begin{proof}
Similar to the complex case : see \cite{B4}.
\end{proof}

\subsubsection{Direct image functor in the singular case}

We define the direct image functors between our category.

Let $f:X\to S$ be a morphism with $X,S\in\Var(k)$,
and assume there exist a factorization $f:X\xrightarrow{l}Y\times S\xrightarrow{p_S}S$ 
with $Y\in\SmVar(k)$, $l$ a closed embedding and $p_S$ a the projection ; 
Let $S=\cup_{i=1}^l S_i$ an open cover such that there exist closed embeddings
$i_i:S_i\hookrightarrow\tilde S_i$ with $\tilde S_i\in\SmVar(k)$.
Then $X=\cup_{i=1}^lX_i$ with $X_i:=f^{-1}(S_i)$. 
Denote, for $I\subset\left[1,\cdots l\right]$, $S_I=\cap_{i\in I} S_i$ and $X_I=\cap_{i\in I}X_i$.
For $I\subset\left[1,\cdots l\right]$, denote by $\tilde S_I=\Pi_{i\in I}\tilde S_i$,
We define the direct image functor on our category by
\begin{eqnarray*}
f^{FDR}_{*mod}:C_{\mathcal D(2)fil}(X/(Y\times\tilde S_I))\to C_{\mathcal D(2)fil}(S/(\tilde S_I)), \\ 
((M_I,F),u_{IJ})\mapsto (\tilde f^{FDR}_{I*mod}(M_I,F),f^k(u_{IJ})):=
(p_{\tilde S_I*}E((\Omega^{\bullet}_{Y\times\tilde S_I/\tilde S_I},F_b)\otimes_{O_{Y\times\tilde S_I}}(M_I,F)[d_Y]),f^k(u_{IJ}))
\end{eqnarray*}
with $f^k(u_{IJ})$ as in \cite{B4}.
It induces in the derived categories the functor
\begin{eqnarray*}
\int^{FDR}_f:D_{\mathcal D(2)fil,r}(X/(Y\times\tilde S_I))\to D_{\mathcal D(2)fil,r}(S/(\tilde S_I)), \;  
((M_I,F),u_{IJ})\mapsto (\tilde f^{FDR}_{I*mod}(M_I,F),f^k(u_{IJ})).
\end{eqnarray*}

In the algebraic case, we have the followings:

\begin{prop}\label{compDmodiihsing}
Let $f_1:X\to Y$ and $f_2:Y\to S$ two morphism with $X,Y,S\in\QPVar(\mathbb C)$ quasi-projective. 
Then there exist factorizations $f_1:X\xrightarrow{l_1} Y'\times Y\xrightarrow{p_Y} Y$ and 
$f_2:Y\xrightarrow{l_2} Y''\times S\xrightarrow{p_S} S$ 
with $Y'=\mathbb P^{N,o}\subset\mathbb P^N$,$Y''=\mathbb P^{N',o}\subset\mathbb P^{N'}$ open subsets,
$l_1,l_2$ closed embeddings and $p_S,p_Y$ the projections.
We have then the factorization 
$f_2\circ f_1:X\xrightarrow{(l_2\circ I_{Y'})\circ l_1}Y'\times Y''\times S\xrightarrow{p_S} S$.
Let $i:S\hookrightarrow\tilde S$ a closed embedding with $\tilde S=\mathbb P^{n,o}\subset\mathbb P^n$ an open subset.
\begin{itemize}
\item[(i)]Let $M\in C_{\mathcal D}(X/(Y'\times Y''\times\tilde S))$. 
Then, we have $\int^{FDR}_{f_2\circ f_1}(M)=\int^{FDR}_{f_2}(\int^{FDR}_{f_1}(M))$ 
in $D_{\mathcal D}(S/(\tilde S_I))$.
\item[(ii)]Let $M\in C_{\mathcal D(2)fil,h}(X/(Y'\times Y''\times\tilde S))$.
Then, we have $\int^{FDR}_{(f_2\circ f_1)!}(M)=\int^{FDR}_{f_2!}(\int^{FDR}_{f_1!}(M))$ 
in $D_{\mathcal D,h}(S/(\tilde S_I))$.
\end{itemize}
\end{prop}

\begin{proof}
Similar to the complex case : see \cite{B4}.
\end{proof}

\subsubsection{Tensor product in the singular case}

Let $S\in\Var(k)$. Let $S=\cup S_i$ an open cover such that there exist closed embeddings
$i_i:S_i\hookrightarrow\tilde S_i$ with $\tilde S_i\in\SmVar(k)$.
We have, as in the complex case, the tensor product functors
\begin{eqnarray*}
(-)\otimes^{[-]}_{O_S}(-):C^2_{\mathcal Dfil}(S/(\tilde S_I))\to C_{\mathcal Dfil}(S/(\tilde S_I)), \;
(((M_I,F),u_{IJ}),((N_I,F),v_{IJ}))\mapsto \\
((M_I,F),u_{IJ})\otimes_{O_S}^{[-]}((N_I,F),v_{IJ}):=
((M_I,F)\otimes_{O_{\tilde S_I}}(N_I,F)[d_{\tilde S_I}],u_{IJ}\otimes v_{IJ}),
\end{eqnarray*}
with, denoting for short $d_{IJ}:=d_{\tilde S_J}-d_{\tilde S_I}$ and $d_I:=d_{\tilde S_I}$,
\begin{eqnarray*}
u_{IJ}\otimes v_{IJ}:(M_I,F)\otimes_{O_{\tilde S_I}}(N_I,F)[d_I]
\xrightarrow{T(p_{IJ}^{*mod},p_{IJ})(-)[d_I]}p_{IJ*}p_{IJ}^{*mod}((M_I,F)\otimes_{O_{\tilde S_I}}(N_I,F))[d_I] \\
\xrightarrow{=}p_{IJ*}(p_{IJ}^{*mod}(M_I,F)\otimes_{O_{\tilde S_J}}p_{IJ}^{*mod}(N_I,F))[d_I] \\
\xrightarrow{I(p_{IJ}^{*mod},p_{IJ})(-,-)(u_{IJ})\otimes I(p_{IJ}^{*mod},p_{IJ})(-,-)(v_{IJ})[d_I]}
p_{IJ*}((M_J,F)\otimes_{O_{\tilde S_J}}(N_J,F))[d_J+d_{IJ}].
\end{eqnarray*}
It induces in the derived category, for $1\leq r\leq\infty$, the functors
\begin{eqnarray*}
(-)\otimes^L_{O_S}(-):D^2_{\mathcal Dfil,r}(S/(\tilde S_I))\to D_{\mathcal Dfil,r}(S/(\tilde S_I)), \;
(((M_I,F),u_{IJ}),((N_I,F),v_{IJ}))\mapsto \\
((M_I,F),u_{IJ})\otimes_{O_S}^L((N_I,F),v_{IJ}):=
(L_D(M_I,F)\otimes_{O_{\tilde S_I}}L_D(N_I,F)[d_{\tilde S_I}],u^q_{IJ}\otimes v^q_{IJ}).
\end{eqnarray*}

We have the following easy proposition :

\begin{prop}\label{otimesedsing}
Let $S\in\Var(k)$. Denote $\Delta_S:S\hookrightarrow S\times S$ the diagonal embedding.
Denote $p_1:S\times S\to S$ and $p_2:S\times S\to S$ the projections.
Let $S=\cup S_i$ an open cover such that there exist
closed embeddings $i_i:S_i\hookrightarrow\tilde S_i$ closed embedding with $\tilde S_i\in\SmVar(k)$.
We have, for $(M_I,u_{IJ}),(N_I,v_{IJ})\in C_{\mathcal D}(S/(\tilde S_I))$,
\begin{equation*}
(M_I,u_{IJ})\otimes^{[-]}_{O_S}(N_I,v_{IJ})=
\Delta_S^{*mod,\Gamma}((p_{1I}^{*mod}M_I,p_{1I}^{*mod}u_{IJ})\otimes_{O_{S\times S}}
(p_{2I}^{*mod}N_I,p_{2I}^{*mod}v_{IJ}))
\end{equation*}
and
\begin{equation*}
((M_I,F),u_{IJ})\otimes^L_{O_S}((N_I,F),v_{IJ})=
R\Delta_S^{*mod,\Gamma}((p_{1I}^{*mod}M_I,p_{1I}^{*mod}u_{IJ})\otimes_{O_{S\times S}}
(p_{2I}^{*mod}N_I,p_{2I}^{*mod}v_{IJ}))
\end{equation*}
\end{prop}

\begin{proof}
Follows from proposition \ref{otimesed} and theorem \ref{KZSk}. 
\end{proof}

\subsubsection{The 2 functors of D modules on the category of algebraic varieties 
over a field $k$ of characteristic zero and the transformation maps}

\begin{defi}\label{TDmodlem0sing}
Consider a commutative diagram in $\Var(k)$ which is cartesian :
\begin{equation*}
D=\xymatrix{X_T\ar[r]^{f'}\ar[d]^{g'} & T\ar[d]^{g} \\
X\ar[r]^{f} & S}.
\end{equation*}
Assume there exist factorizations $f:X\xrightarrow{l_1}Y_1\times S\xrightarrow{p_S}S$, 
$g:T\xrightarrow{l_2}Y_2\times S\xrightarrow{p_S}S$, with $Y_1,Y_2\in\SmVar(k)$, 
$l_1,l_2$ closed embeddings and $p_S$, $p_S$ the projections.
Let $S=\cup S_i$ an open cover such that there exist
closed embeddings $i_i:S_i\hookrightarrow\tilde S_i$ closed embedding with $\tilde S_i\in\SmVar(k)$.
We then have as in the complex case, for $(M,F)=((M_I,F),u_{IJ})\in C_{\mathcal D(2)fil}(X/(Y_1\times\tilde S_I))$,
the following canonical transformation map in $D_{\mathcal D(2)fil,r}(T/(Y_2\times\tilde S_I))$,
\begin{eqnarray*}
T^{\mathcal Dmod}(f,g)(M,F): \\
Rg^{*mod,\Gamma}\int^{FDR}_f(M,F):=(\Gamma_{T_I}E(\tilde g_I^{*mod}p_{\tilde S_I*}
E((\Omega^{\bullet}_{Y_1\times\tilde S_I/\tilde S_I},F_b)\otimes_{O_{Y_1\times\tilde S_I}}(M_I,F))),
\tilde g_J^{*mod}f^k(u_{IJ})) \\
\xrightarrow{(T^O_{\omega}(p_{\tilde S_I},\tilde g_I)(M_I,F))} \\
(\Gamma_{T_I}E(p_{Y_2\times\tilde S_I*}E((\Omega^{\bullet}_{Y_1\times Y_2\times\tilde S_I/Y_2\times\tilde S_I},F_b)
\otimes_{O_{Y_1\times Y_2\times\tilde S_I}}p_{Y_1\times\tilde S_I}^{*mod}(M_I,F))),
f^{'k}(p_{Y_1\times\tilde S_J}^{*mod}(u_{IJ}))) \\
\xrightarrow{(T^O_{\omega}(\gamma,\otimes)(p_{Y_1\times\tilde S_I}^{*mod}(M_I,F)))^{-1}} \\
(p_{Y_2\times\tilde S_I*}E((\Omega^{\bullet}_{Y_1\times Y_2\times\tilde S_I/Y_2\times\tilde S_I},F_b)
\otimes_{O_{Y_1\times Y_2\times\tilde S_I}}\Gamma_{Y_1\times T_I}E(p_{Y_1\times\tilde S_I}^{*mod}((M_I,F)))),
f^{'k}(\tilde g_J^{''*mod}(u^q_{IJ}))) \\
=:\int^{FDR}_{f'}Rg^{'*mod,\Gamma}(M,F).
\end{eqnarray*}
\end{defi}

\begin{prop}\label{PDmod1sing}
Consider a commutative diagram in $\Var(k)$
\begin{equation*}
D=(f,g)=\xymatrix{X_T\ar[r]^{f'}\ar[d]^{g'} & T\ar[d]^{g} \\
X\ar[r]^{f} & S}.
\end{equation*}
which is cartesian. Assume there exist factorizations $f:X\xrightarrow{l_1}Y_1\times S\xrightarrow{p_S}S$, 
$g:T\xrightarrow{l_2}Y_2\times S\xrightarrow{p_S}S$, with $Y_1,Y_2\in\SmVar(k)$, 
$l_1,l_2$ closed embeddings and $p_S$, $p_S$ the projections.
Let $S=\cup S_i$ an open cover such that there exist
closed embeddings $i_i:S_i\hookrightarrow\tilde S_i$ closed embedding with $\tilde S_i\in\SmVar(k)$.
For $(M,F)=((M_I,F),u_{IJ})\in C_{\mathcal D(2)fil,c}(X/(Y\times\tilde S_I))$,
\begin{equation*}
T^{\mathcal Dmod}(f,g):Rg^{*mod,\Gamma}\int^{FDR}_f(M,F)\to\int^{FDR}_{f'}Rg^{'*mod,\Gamma}(M,F)
\end{equation*}
is an isomorphism in $D_{\mathcal D(2)fil,r}(T/(Y_2\times\tilde S_I))$.
\end{prop}

\begin{proof}
Similar to the complex case.
\end{proof}

\begin{thm}\label{compDmodrhthmsing}
Let $f:X\to S$ a morphism with $X,S\in\Var(k)$. 
Assume there exists a factorization $f:X\xrightarrow{l}Y\times S\xrightarrow{p}S$
with $Y\in\SmVar(k)$, $l$ a closed embedding and $p$ the projection.
Let $S=\cup S_i$ an open cover such that there exist
closed embeddings $i_i:S_i\hookrightarrow\tilde S_i$ closed embedding with $\tilde S_i\in\SmVar(k)$. Then,
\begin{itemize}
\item[(i)] For $(M,F)\in C_{\mathcal D(2)fil,rh}(S/(\tilde S_I)^{op})$, 
we have $L\mathbb D_S(M,F)\in D_{\mathcal D(2)fil,rh}(S/(\tilde S_I))$.
\item[(ii)] For $M\in C_{\mathcal D,rh}(S/(\tilde S_I))$, 
$Rf^{*mod,\Gamma}(M)\in D_{\mathcal D,rh}(X/(Y\times\tilde S_I))$
and $Lf^{\hat*mod,\Gamma}M\in D_{\mathcal D,rh}(X/(Y\times\tilde S_I))$.
\item[(iii)] For $M\in C_{\mathcal D,rh}(X/(Y\times\tilde S_I))$, 
$\int_fM\in D_{\mathcal D,rh}(S/(\tilde S_I))$ and 
$\int_{f!}M:=L\mathbb D_S\int_fL\mathbb D_X\in D_{\mathcal D,rh}(S/(\tilde S_I))$.
\item[(iv)] If $f$ is proper, for $(M,F)\in C_{\mathcal D(2)fil,rh}(X/(Y\times\tilde S_I)))$, 
we have $\int_f(M,F)\in D_{\mathcal D(2)fil,rh}(S/(\tilde S_I))$.
\item[(v)] For $(M,F),(N,F)\in C_{\mathcal D(2)fil,rh}(S/(\tilde S_I))$, 
$(M,F)\otimes^L_{O_S}(N,F)\in D_{\mathcal D(2)fil,rh}(S/(\tilde S_I))$
\end{itemize}
\end{thm}

\begin{proof}
Follows from theorem \ref{compDmodrhthm}.
\end{proof}

\subsection{The category of complexes of quasi-coherent sheaves on an algebraic variety
whose cohomology sheaves has a structure of D-modules}

\subsubsection{Definition on a smooth algebraic variety and the functorialities}

\begin{defi}
Let $S\in\SmVar(k)$. Let $Z\subset S$ a closed subset.
Denote by $j:S\backslash Z\hookrightarrow S$ the open complementary embedding.
\begin{itemize}
\item[(i)] We denote by $C_{O_S,\mathcal D,Z}(S)\subset C_{O_S,\mathcal D}(S)$
the full subcategory consisting of $M\in C_{O_S,\mathcal D}(S)$ such that
such that $j^*H^nM=0$ for all $n\in\mathbb Z$.
\item[(ii)] We denote by $C_{O_Sfil,\mathcal D,Z}(S)\subset C_{O_Sfil,\mathcal D}(S)$
the full subcategory consisting of $(M,F)\in C_{O_Sfil,\mathcal D}(S)$ such that there exist $r\in\mathbb N$
and an $r$-filtered homotopy equivalence $m:(M,F)\to(M',F)$ with $(M',F)\in C_{O_Sfil,\mathcal D}(S)$
such that $j^*H^n\Gr_F^p(M',F)=0$ for all $n,p\in\mathbb Z$.
\end{itemize}
\end{defi}

\begin{defi}
Let $S\in\SmVar(k)$.
We have then (see section 2), for $r=1,\cdots,\infty$, the homotopy category 
$K_{O_Sfil,\mathcal D,r}(S)=\Ho_r(C_{O_Sfil,\mathcal D}(S))$ whose objects are those of $C_{O_Sfil,\mathcal D}(S)$
and whose morphisms are $r$-filtered homotopy classes of morphism, and
its localization $D_{O_Sfil,\mathcal D,r}(S)=K_{O_Sfil,\mathcal D,r}(S)([E_1]^{-1})$ with respect to
filtered zariski, resp. usu local equivalence.
Note that the classes of filtered $\tau$ local equivalence constitute a right multiplicative system.
\end{defi}

\begin{itemize}

\item Let $S\in\SmVar(k)$. Let $(M,F)\in C_{O_Sfil,\mathcal D}(S)$.
Then, the canonical morphism $q:L_O(M,F)\to(M,F)$ in $C_{O_Sfil}(S)$ being a quasi-isomorphism of $O_S$ modules, 
we get in a unique way $L_O(M,F)\in C_{O_Sfil,\mathcal D}(S)$ such that $q:L_O(M,F)\to(M,F)$ is a morphism
in $C_{O_Sfil,\mathcal D}(S)$

\item Let $f:X\to S$ be a morphism with $X,S\in\SmVar(k)$.
Let $(M,F)\in C_{O_Sfil,\mathcal D}(S)$. Then, $f^{*mod}H^n(M,F):=(O_X,F_b)\otimes_{f^*O_S}f^*H^n(M,F)$
is canonical a filtered $D_X$ module (see section 4.1 or 4.2).
Consider the canonical surjective map $q(f):H^nf^{*mod}(M,F)\to f^{*mod}H^n(M,F)$.
Then, $q(f)$ is an isomorphism if $f$ is smooth.
Let $h:U\to S$ be a smooth morphism with $U,S\in\SmVar(k)$. We get the functor
\begin{eqnarray*}
h^{*mod}:C_{O_Sfil,\mathcal D}(S)\to C_{O_Ufil,\mathcal D}(U),
(M,F)\mapsto h^{*mod}(M,F), 
\end{eqnarray*} 

\item Let $S\in\SmVar(k)$, and 
let $i:Z\hookrightarrow S$ a closed embedding and denote by $j:S\backslash Z\hookrightarrow S$ the open complementary. 
For $M\in C_{O_S,\mathcal D}(S)$, the cohomology presheaves of 
\begin{equation*}
\Gamma_Z M:=\Cone(\ad(j^*,j_*)(M):M\to j_*j^*M)[-1]
\end{equation*}
has a canonical $D_S$-module structure 
(as $j^*H^nM$ is a $j^*D_S$ module, $H^nj_*j^*M=j_*j^*H^nM$ has an induced structure of $D_S$ module),
and $\gamma_Z(M):\Gamma_Z M\to M$ is a map in $C_{O_S,\mathcal D}(S)$.
For $Z_2\subset Z$ a closed subset and $M\in C_{O_S,\mathcal D}(S)$, 
$T(Z_2/Z,\gamma)(M):\Gamma_{Z_2}M\to\Gamma_Z M$ is a map in $C_{O_S,\mathcal D}(S)$.
We get the functor
\begin{eqnarray*}
\Gamma_Z:C_{O_Sfil,\mathcal D}(S)\to C_{O_Sfil,\mathcal D}(S), \\ 
(M,F)\mapsto\Gamma_Z(M,F):=\Cone(\ad(j^*,j_*)((M,F)):(M,F)\to j_*j^*(M,F))[-1],
\end{eqnarray*}
together we the canonical map $\gamma_Z(M,F):\Gamma_Z(M,F)\to (M,F)$

More generally, let $h:Y\to S$ a morphism with $Y,S\in\Var(k)$, $S$ smooth, and 
let $i:X\hookrightarrow Y$ a closed embedding and denote by $j:Y\backslash X\hookrightarrow Y$ the open complementary. 
For $M\in C_{h^*O_S,h^*\mathcal D}(Y)$, 
\begin{equation*}
\Gamma_X M:=\Cone(\ad(j^*,j_*)(M):M\to j_*j^*M)[-1]
\end{equation*}
has a canonical $h^*D_S$-module structure,
(as $j^*H^nM$ is a $j^*h^*D_S$ module, $H^nj_*j^*M=j_*j^*H^nM$ has an induced structure of $j^*h^*D_S$ module),
and $\gamma_X(M):\Gamma_X M\to M$ is a map in $C_{h^*O_S,h^*\mathcal D}(Y)$.
For $X_2\subset X$ a closed subset and $M\in C_{h^*O_S,h^*\mathcal D}(Y)$, 
$T(Z_2/Z,\gamma)(M):\Gamma_{X_2}M\to\Gamma_X M$ is a map in $C_{h^*O_S,h^*\mathcal D}(Y)$. We get the functor
\begin{eqnarray*}
\Gamma_X:C_{h^*O_Sfil,h^*\mathcal D}(Y)\to C_{h^*O_Sfil,h^*\mathcal D}(Y), \\ 
(M,F)\mapsto\Gamma_X(M,F):=\Cone(\ad(j^*,j_*)((M,F)):(M,F)\to j_*j^*(M,F))[-1],
\end{eqnarray*}
together we the canonical map $\gamma_X(M,F):\Gamma_X(M,F)\to (M,F)$

\item Let $f:X\to S$ be a morphism with $X,S\in\SmVar(k)$. 
Consider the factorization $f:X\xrightarrow{l}X\times S\xrightarrow{p}S$,
where $l$ is the graph embedding and $p$ the projection. We get from the two preceding points the functor
\begin{eqnarray*}
f^{*mod,\Gamma}:C_{O_Sfil,\mathcal D}(S)\to C_{O_Xfil,\mathcal D}(X\times S),
(M,F)\mapsto f^{*mod,\Gamma}(M,F):=\Gamma_Xp^{*mod}(M,F), 
\end{eqnarray*} 
and
\begin{eqnarray*}
f^{*mod[-],\Gamma}:C_{O_Sfil,\mathcal D}(S)\to C_{O_Xfil,\mathcal D}(X\times S), \\
(M,F)\mapsto f^{*mod[-],\Gamma}(M,F):=\Gamma_XE(p^{*mod}(M,F))[-d_X], 
\end{eqnarray*} 
which induces in the derived categories the functor
\begin{eqnarray*}
Rf^{*mod[-],\Gamma}:D_{O_Sfil,\mathcal D}(S)\to D_{O_Xfil,\mathcal D}(X\times S), \\
(M,F)\mapsto Rf^{*mod[-],\Gamma}(M,F):=\Gamma_XE(p^{*mod[-]}(M,F)). 
\end{eqnarray*} 

\end{itemize}

\subsubsection{Definition on a singular algebraic variety and the functorialities}

\begin{defi}\label{DmodsingVardef}
Let $S\in\Var(k)$ and let $S=\cup_iS_i$ an open cover
such that there exist closed embeddings $i_i:S_i\hookrightarrow\tilde S_i$ with $\tilde S_i\in\SmVar(k)$. Then,
$C_{Ofil,\mathcal D}(S/(\tilde S_I))$ is the category
\begin{itemize}
\item whose objects are $(M,F)=((M_I,F)_{I\subset\left[1,\cdots l\right]},u_{IJ})$, with
\begin{itemize}
\item $(M_I,F)\in C_{O_{\tilde S_I}fil\mathcal D,S_I}(\tilde S_I)$,
\item $u_{IJ}:m^*(M_I,F)\to m^*p_{IJ*}(M_J,F)[d_{\tilde S_J}-d_{\tilde S_I}]$ 
for $J\subset I$, are morphisms, $p_{IJ}:\tilde S_J\to\tilde S_I$ being the projection, 
satisfying for $I\subset J\subset K$, $p_{IJ*}u_{JK}\circ u_{IJ}=u_{IK}$ in $C_{O_{\tilde S_I}fil,\mathcal D}(\tilde S_I)$ ;
\end{itemize}
\item whose morphisms $m:((M_I,F),u_{IJ})\to((N_I,F),v_{IJ})$ between  
$(M,F)=((M_I,F)_{I\subset\left[1,\cdots l\right]},u_{IJ})$ and $(N,F)=((N_I,F)_{I\subset\left[1,\cdots l\right]},v_{IJ})$
are a family of morphisms of complexes,  
\begin{equation*}
m=(m_I:(M_I,F)\to (N_I,F))_{I\subset\left[1,\cdots l\right]}
\end{equation*}
such that $v_{IJ}\circ m_I=p_{IJ*}m_J\circ u_{IJ}$ in $C_{O_{\tilde S_I}fil,\mathcal D}(\tilde S_I)$.
\end{itemize}
We denote by $C^{\sim}_{Ofil,\mathcal D}(S/(\tilde S_I))\subset C_{Ofil,\mathcal D}(S/(\tilde S_I))$ the full subcategory 
consisting of objects $((M_I,F),u_{IJ})$ such that the $u_{IJ}$ are $\infty$-filtered Zariski local equivalences.
\end{defi}

\begin{defi}
Let $S\in\Var(k)$ and let $S=\cup_iS_i$ an open cover
such that there exist closed embeddings $i_i:S_i\hookrightarrow\tilde S_i$ with $\tilde S_i\in\SmVar(k)$.
We have then (see \cite{B4}), for $r=1,\cdots,\infty$, the homotopy category 
\begin{equation*}
K_{Ofil,\mathcal D,r}(S/(\tilde S_I)):=\Ho_r(C_{Ofil,\mathcal D}(S/(\tilde S_I))) 
\end{equation*}
whose objects are those of $C_{Ofil,\mathcal D}(S/(\tilde S_I))$
and whose morphisms are $r$-filtered homotopy classes of morphism, and its localization 
\begin{equation*}
D_{fil,\mathcal D,r}(S/(\tilde S_I)):=K_{Ofil,\mathcal D,r}(S/(\tilde S_I))([E_1]^{-1})
\end{equation*}
with respect to the classes of filtered zariski local equivalence.
Note that the classes of filtered $\tau$ local equivalence constitute a right multiplicative system.
\end{defi}

Let $f:X\to S$ be a morphism, with $X,S\in\Var(k)$, 
such that there exist a factorization $f;X\xrightarrow{l}Y\times S\xrightarrow{p_S}S$
with $Y\in\SmVar(k)$, $l$ a closed embedding and $p_S$ the projection,
and consider $S=\cup_{i=1}^lS_i$ an open cover such that there exist closed embeddings 
$i_i:S_i\hookrightarrow\tilde S_i$, with $\tilde S_i\in\SmVar(k)$.
Then, $X=\cup^l_{i=1} X_i$ with $X_i:=f^{-1}(S_i)$.
We then have the filtered De Rham the inverse image functor :
\begin{eqnarray*}
f^{*mod[-],\Gamma}:C_{Ofil,\mathcal D}(S/(\tilde S_I))\to C_{Ofil,\mathcal D}(X/(Y\times\tilde S_I)), \; \; 
(M,F)=((M_I,F),u_{IJ})\mapsto \\
f^{*mod[-],\Gamma}(M,F):=(\Gamma_{X_I}E(p_{\tilde S_I}^{*mod[-]}(M_I,F))),\tilde f_J^{*mod[-]}u_{IJ})
\end{eqnarray*}
with $\tilde f_J^{*mod[-]}u_{IJ}$ as in the complex case
It induces in the derived categories, the functor
\begin{eqnarray*}
Rf^{*mod[-],\Gamma}:D_{Ofil,\mathcal D,r}(S/(\tilde S_I)\to D_{Ofil,\mathcal D,r}(X/(Y\times\tilde S_I)), \\  
(M,F)=((M_I,F),u_{IJ})\mapsto \\
Rf^{*mod[-],\Gamma}:=f^{*mod[-],\Gamma}(M,F):=(\Gamma_{X_I}E(p_{\tilde S_I}^{*mod[-]}(M_I,F)),\tilde f_J^{*mod[-]}u_{IJ}).
\end{eqnarray*}

\subsection{The (filtered) De Rahm functor over a field $k$ of characteristic zero 
and Riemann Hilbert for holonomic D-modules on smooth algebraic varieties over a subfield $k\subset\mathbb C$}

Let $j:S^o\hookrightarrow S$ an open embedding with $S=(S,O_S)\in\RTop$. 
Denote by $Z:=S\backslash S^o$ the closed complementary subset.
Recall that we have, see \cite{B4}, for $(M,F)\in C_{\mathcal Dfil}(S^o)$ the canonical maps in $C_{fil}(S)$
\begin{eqnarray*}
T^w(j,\otimes)(M,F):DR(S)(j_*(M,F)):=(\Omega^{\bullet}_S,F_b)\otimes_{O_S}j_*(M,F)
\xrightarrow{\ad(j^*,j_*)(\Omega^{\bullet}_S)\otimes I} \\
j_*j^*(\Omega^{\bullet}_S,F_b)\otimes_{O_S}j_*(M,F)
\xrightarrow{=}j_*((\Omega^{\bullet}_{S^o},F_b)\otimes_{O_{S^o}}(M,F))=:j_*DR(S^o)(M,F)
\end{eqnarray*}
and
\begin{eqnarray*}
T^w(\gamma_Z,\otimes)(M,F):DR(S)(\Gamma_Z(M,F)):=(\Omega^{\bullet}_S,F_b)\otimes_{O_S}\Gamma_Z(M,F) \\
\xrightarrow{(I,T^w(j,\otimes)(M,F))}\Gamma_Z((\Omega^{\bullet}_{S^o},F_b)\otimes_{O_{S^o}}(M,F))
=:\Gamma_ZDR(S^o)(M,F).
\end{eqnarray*}

Let $k$ a field of characteristic zero.

\begin{prop}\label{DRhUS}
Let $Y,S\in\SmVar(k)$. Let $p:Y\times S\to S$ the projection.
For $(M,F)\in C_{\mathcal Dfil}(Y\times S)$, 
\begin{equation*}
DR(Y\times S/S)(M,F):=(\Omega^{\bullet}_{Y\times S/S},F_b)\otimes_{O_{Y\times S}}(M,F)\in C_{p^*O_Sfil}(Y\times S)
\end{equation*}
is a naturally a complex of filtered $p^*D_S$ modules, that is 
\begin{equation*}
DR(Y\times S/S)(M,F):=(\Omega^{\bullet}_{Y\times S/S},F_b)\otimes_{O_{Y\times S}}(M,F)\in C_{p^*\mathcal Dfil}(Y\times S),
\end{equation*}
where the $p^*D_S$ module structure on $\Omega^p_{Y\times S/S}\otimes_{O_{Y\times S}}M^n$ is given by
for $(Y\times S)^o\subset Y\times S$ an open subset,
\begin{equation*}
(\gamma\in\Gamma((Y\times S)^o,T_{Y\times S}),
\hat\omega\otimes m\in\Gamma((Y\times S)^o,\Omega^p_{Y\times S/S}\otimes_{O_{Y\times S}}M^n))\mapsto 
\gamma.(\hat\omega\otimes m):=(\hat\omega\otimes(\gamma.m).
\end{equation*}
Moreover, if $\phi:(M_1,F)\to(M_2,F)$ a morphism with $(M_1,F),(M_2,F)\in C_{\mathcal Dfil}(Y\times S)$,
\begin{equation*}
DR(Y\times S/S)(\phi):=(I\otimes\phi):(\Omega^{\bullet}_{Y\times S/S},F_b)\otimes_{O_{Y\times S}}(M_1,F)
\to(\Omega^{\bullet}_{Y\times S/S},F_b)\otimes_{O_{Y\times S}}(M_2,F)
\end{equation*}
is a morphism in $C_{p^*\mathcal Dfil}(Y\times S)$.
\end{prop}

\begin{proof}
Follows imediately by definition : see \cite{B4}.
\end{proof}

\begin{prop}\label{TDhwM}
Consider a commutative diagram in $\SmVar(k)$ :
\begin{equation*}
D=\xymatrix{Y\times S\ar[r]^{p} & S \\
Y'\times T\ar[u]^{g''=(g''_0\times g)}\ar[r]^{p'} & T\ar[u]^{g}}
\end{equation*}
with $p$ and $p'$ the projections. 
For $(M,F)\in C_{\mathcal Dfil}(Y\times S)$ the map in $C_{g^{''*}p^*O_Sfil}(Y'\times T)$
\begin{equation*}
\Omega_{(Y'\times T/Y\times S)/(T/S)}(M,F):g''^*((\Omega^{\bullet}_{Y\times S/S},F_b)\otimes_{O_{Y\times S}}(M,F))
\to(\Omega^{\bullet}_{Y'\times T/T},F_b)\otimes_{O_{Y'\times T}}g^{''*mod}(M,F)
\end{equation*}
given in \cite{B4} section 4.1 is a map in $C_{g^{''*}p^*\mathcal Dfil}(Y'\times T)$. 
Hence, for $(M,F)\in C_{\mathcal Dfil}(Y\times S)$, the map in $C_{O_Tfil}(T)$ (with $L_D$ instead of $L_O$)
\begin{eqnarray*}
T^O_{\omega}(D)(M):g^{*mod}L_D(p_*E((\Omega^{\bullet}_{Y\times S/S},F_b)\otimes_{O_{Y\times S}}(M,F)))\to
p'_*E((\Omega^{\bullet}_{Y'\times T/T},F_b)\otimes_{O_{Y'\times T}}g^{''*mod}(M,F)),
\end{eqnarray*}
is a map in $C_{\mathcal Dfil}(T)$.
\end{prop}

\begin{proof}
Follows imediately by definition. 
\end{proof}

\begin{prop}\label{resw}
Let $S\in\SmVar(k)$.
\begin{itemize}
\item We have the filtered resolutions of $K_S$ by the following complex of locally free right $D_S$ modules:
$\omega(S):\omega(K_S):=(\Omega^{\bullet}_{S},F_b)[d_S]\otimes_{O_S}(D_S,F_b)\to (K_S,F_b)$ and
$\omega(S):\omega(K_S,F^{ord}):=(\Omega^{\bullet}_{S},F_b)[d_S]\otimes_{O_S}(D_S,F^{ord})\to (K_S,F^{ord})$
\item Dually, we have the filtered resolution of $O_S$ by the following complex of locally free (left) $D_S$ modules:
$\omega^{\vee}(S):\omega(O_S):=(\wedge^{\bullet}T_S,F_b)[d_S]\otimes_{O_S}(D_S,F_b)\to (O_S,F_b)$ and
$\omega^{\vee}(S):\omega(O_S,F^{ord}):=(\wedge^{\bullet}T_S,F_b)[d_S]\otimes_{O_S}(D_S,F^{ord})\to (O_S,F^{ord})$.
\end{itemize}
Let $S_1,S_2\in\SmVar(k)$. Consider the projection $p=p_1:S_1\times S_2\to S_1$.
\begin{itemize}
\item We have the filtered resolution of $D_{S_1\times S_2\to S_1}$ by the following complexes of (left)
$(p^*D_{S_1}$ and right $D_{S_1\times S_2})$ modules :
\begin{equation*}
\omega(S_1\times S_2/S_1):
(\Omega^{\bullet}_{S_1\times S_2/S_1}[d_{S_2}],F_b)\otimes_{O_{S_1\times S_2}}(D_{S_1\times S_2},F^{ord})
\to(D_{S_1\times S_2\leftarrow S_1},F^{ord}).
\end{equation*}
\item Dually, we have the filtered resolution of $D_{S_1\times S_2\to S_1}$ by the following complexes of (left)
$(p^*D_{S_1},D_{S_1\times S_2})$ modules :
\begin{equation*}
\omega^{\vee}(S_1\times S_2/S_1):
(\wedge^{\bullet}T_{S_1\times S_2/S_1}[d_{S_2}],F_b)\otimes_{O_{S_1\times S_2}}(D_{S_1\times S_2},F^{ord})
\to(D_{S_1\times S_2\to S_1},F^{ord}),
\end{equation*}
\end{itemize}
\end{prop}

\begin{proof} 
Similar to the complex case: see \cite{LvDmod}.
\end{proof}

\begin{defi}\label{directmod}
\begin{itemize}
\item[(i)] Let $i:Z\hookrightarrow S$ be a closed embedding, with $Z,S\in\SmVar(k)$.
Then, for $(M,F)\in C_{\mathcal Dfil}(Z)$, we set 
\begin{equation*}
i_{*mod}(M,F):=i^0_{*mod}(M,F):=i_*((M,F)\otimes_{D_Z}(D_{Z\leftarrow S},F^{ord}))\in C_{\mathcal Dfil}(S)
\end{equation*}
\item[(ii)] Let $S_1,S_2\in\SmVar(k)$ and $p:S_1\times S_2\to S_1$ be the projection.
Then, for $(M,F)\in C_{\mathcal Dfil}(S_1\times S_2)$, we set
\begin{itemize}
\item $p^0_{*mod}(M,F):=p_*(DR(S_1\times S_2/S_1)(M,F)):=
p_*((\Omega^{\bullet}_{S_1\times S_2/S_1},F_b)\otimes_{O_{S_1\times S_2}}(M,F))[d_{S_2}]\in C_{\mathcal Dfil}(S_1)$,
\item $p_{*mod}(M,F):=p_*E(DR(S_1\times S_2/S_1)(M,F)):=
p_*E((\Omega^{\bullet}_{S_1\times S_2/S_1},F_b)\otimes_{O_{S_1\times S_2}}(M,F))[d_{S_2}]\in C_{\mathcal Dfil}(S_1)$.
\end{itemize}
\item[(iii)] Let $f:X\to S$ be a morphism, with $X,S\in\SmVar(k)$.
Consider the factorization $f:X\xrightarrow{i} X\times S\xrightarrow{p_S}S$, 
where $i$ is the graph embedding and $p_S:X\times S\to S$ is the projection.
Then, for $(M,F)\in C_{\mathcal Dfil}(X)$ we set
\begin{itemize}
\item $f^{FDR}_{*mod}(M,F):=p_{S*mod}i_{*mod}(M,F)\in C_{\mathcal Dfil}(S)$,
\item $\int_f^{FDR}(M,F):=f^{FDR}_{*mod}(M,F):=p_{S*mod}i_{*mod}(M,F)\in D_{\mathcal Dfil}(S)$.
\end{itemize}
By proposition \ref{Pham} below, we have $\int_f^{FDR}M=\int_fM\in D_{\mathcal D}(X)$.
\item[(iii)] Let $f:X\to S$ be a morphism, with $X,S\in\SmVar(k)$.
Consider the factorization $f:X\xrightarrow{i} X\times S\xrightarrow{p_S}S$, 
where $i$ is the graph embedding and $p_S:X\times S\to S$ is the projection.
Then, for $(M,F)\in C_{\mathcal Dfil}(X)$ we set
\begin{itemize}
\item $f^{FDR}_{!mod}(M,F):=\mathbb D_S^KL_Df^{FDR}_{*mod}\mathbb D_S^KL_D(M,F):=
\mathbb D_S^KL_Dp_{S*mod}i_{*mod}\mathbb D_{X\times S}^KL_D(M,F)\in C_{\mathcal Dfil}(S)$,
\item $\int_{f!}^{FDR}(M,F):=f^{FDR}_{!mod}(M,F):=
\mathbb D_S^KL_Dp_{S*mod}i_{*mod}\mathbb D_{X\times S}^KL_D(M,F)\in D_{\mathcal Dfil}(S)$.
\end{itemize}
\end{itemize}
\end{defi}

\begin{prop}\label{Pham}
\begin{itemize}
\item[(i)] Let $i:Z\hookrightarrow S$ a closed embedding with $S,Z\in\SmVar(k)$.
Then for $(M,F)\in C_{\mathcal Dfil}(Z)$, we have
\begin{equation*}
\int_i(M,F):=Ri_*((M,F)\otimes^L_{D_Z}(D_{Z\leftarrow S},F^{ord})=i_*((M,F)\otimes_{D_Z}(D_{Z\leftarrow S},F^{ord}))=i_{*mod}(M,F).
\end{equation*}
\item[(ii)] Let $S_1,S_2\in\SmVar(k)$ and $p:S_{12}:=S_1\times S_2\to S_1$ be the projection.
Then, for $(M,F)\in C_{\mathcal Dfil}(S_1\times S_2)$ we have
\begin{eqnarray*}
\int_p(M,F):&=&Rp_*((M,F)\otimes^L_{D_{S_1\times S_2}}(D_{S_1\times S_2\leftarrow S_1},F^{ord})) \\
&=&p_*E((\Omega^{\bullet}_{S_1\times S_2/S_1},F_b)\otimes_{O_{S_1\times S_2}}
(D_{S_1\times S_2},F^{ord})\otimes_{D_{S_1\times S_2}}(M,F))[d_{S_2}] \\
&=&p_*E((\Omega^{\bullet}_{S_1\times S_2/S_1},F_b)\otimes_{O_{S_1\times S_2}}(M,F))[d_{S_2}]
=:p_{*mod}(M,F).
\end{eqnarray*}
where the second equality follows from Griffitz transversality (the canonical isomorphism map respect by definition the filtration).
\item[(iii)] Let $f:X\to S$ be a morphism with $X,S\in\SmVar(k)$.
Then for $M\in C_{\mathcal D}(X)$, we have $\int^{FDR}_fM=\int_fM$.
\end{itemize}
\end{prop}

\begin{proof}
\noindent(i):Follows from the fact that $D_{Z\leftarrow S}$ is a locally free $D_Z$ module and that $i_*$ is an exact functor.

\noindent(ii): Since $\Omega^{\bullet}_{S_{12}/S_1}[d_{S_2}],F_b)\otimes_{O_{S_{12}}}D_{S_{12}}$
is a complex of locally free $D_{S_1\times S_2}$ modules, we have in $D_{fil}(S_1\times S_2)$, using proposition \ref{resw},
\begin{eqnarray*}
(D_{S_1\times S_2\leftarrow S_1},F^{ord})\otimes^L_{D_{S_1\times S_2}}(M,F)=
(\Omega^{\bullet}_{S_{12}/S_1}[d_{S_2}],F_b)\otimes_{O_{S_{12}}}(D_{S_{12}},F^{ord})\otimes_{D_{S_{12}}}(M,F).
\end{eqnarray*} 

\noindent(iii): Follows from (i) and (ii) by proposition \ref{compDmod}(ii).
\end{proof}

Let $k$ a field of characteristic zero. Let $S\in\SmVar(k)$ connected.
We use to shift the De Rham functor 
in order to have compatibility with perverse sheaves in the complex or p-adic ananlytic case by setting
for $(M,F)\in C_{\mathcal Dfil}(S)$, $DR(S)^{[-]}(M,F):=DR(S)(M,F)[-d_S]\in C_{fil}(S)$.

\begin{itemize}
\item Let $f:X\to S$ a morphism with $S,X\in\SmVar(k)$.
Recall that we have for $(M,F)\in C_{\mathcal Dfil}(X)$ the canonical map in $D_{fil}(S)$
\begin{eqnarray*}
T_*(f,DR)(M,F):DR(S)(\int_f(M,F)):=
(\Omega_S^{\bullet},F_b)\otimes_{O_S}Rf_*((D_{X\leftarrow S},F^{ord})\otimes_{D_X}^L(M,F)) \\
\xrightarrow{\iota(S)\otimes I}
Rf_*((D_{X\leftarrow S},F^{ord})\otimes_{D_X}^L(M,F))\otimes^L_{D_S}(K_S,F^{ord}) \\
\xrightarrow{T(f,\otimes)((D_{X\leftarrow S},F^{ord})\otimes_{D_X}^L(M,F),(K_S,F^{ord}))} \\
Rf_*(f^*(K_S,F^{ord})\otimes^L_{f^*D_S}(D_{X\leftarrow S},F^{ord})\otimes_{D_X}^L(M,F)) 
\xrightarrow{=}Rf_*((K_X,F^{ord})\otimes_{D_X}^L(M,F)) \\
\xrightarrow{\iota(X)\otimes I}
Rf_*((\Omega_X^{\bullet},F_b)\otimes_{O_X}(M,F))=:Rf_*DR(X)(M,F)
\end{eqnarray*}
which is an isomorphism by the projection formula for quasi-coherent sheaves for a morphism of ringed topos
and proposition \ref{resw}.
In particular, for $S\in\SmVar(k)$ and $(M,F)\in C_{\mathcal Dfil,c}(S^o)$,  
\begin{eqnarray*}
T^w(j,\otimes)(M,F):DR(S)(j_*(M,F)):=(\Omega^{\bullet}_S,F_b)\otimes_{O_S}j_*(M,F) \\ 
\to j_*((\Omega^{\bullet}_{S^o},F_b)\otimes_{O_{S^o}}(M,F))=:j_*DR(S^o)(M,F)
\end{eqnarray*}
is a filtered quasi-isomorphism.
\item Let $f:X\to S$ a morphism with $S,X\in\Var(k)$. 
Assume there exist a factorization $f:X\xrightarrow{l} Y\times S\xrightarrow{p}S$
with $Y\in\SmVar(k)$, $l$ a closed embedding and $p$ the projection.
Let $S=\cup_iS_i$ an open cover such that there exists closed embedding 
$i_i:S_i\hookrightarrow\tilde S_i$ with $\tilde S_i\in\SmVar(k)$.
We have for $((M_I,F),u_{IJ})\in C_{\mathcal Dfil}(X/(Y\times\tilde S_I))$ the canonical map in $D_{fil}(S/(\tilde S_I))$
\begin{eqnarray*}
T_*(f,DR)((M_I,F),u_{IJ}): \\ 
DR(S)(\int_f((M_I,F),u_{IJ})):=
((\Omega_{\tilde S_I}^{\bullet},F_b)\otimes_{O_{\tilde S_I}}
p_{\tilde S_I*}E((\Omega^{\bullet}_{Y\times\tilde S_I/\tilde S_I},F_b)\otimes_{O_{Y\times\tilde S_I}}(M_I,F)),DR(fu_{IJ})) \\
\xrightarrow{(k\circ T(p_{\tilde S_I},\otimes)(-,-))} 
(p_{\tilde S_I*}E((p_{\tilde S_I}^*\Omega_{\tilde S_I}^{\bullet},F_b)\otimes_{p_{\tilde S_I}^*O_{\tilde S_I}}
(\Omega^{\bullet}_{Y\times\tilde S_I/\tilde S_I},F_b)\otimes_{O_{Y\times\tilde S_I}}(M_I,F)),DR(fu_{IJ})) \\
\xrightarrow{w(Y\times\tilde S_I)} 
p_{\tilde S_I*}E((\Omega^{\bullet}_{Y\times\tilde S_I},F_b)\otimes_{O_{Y\times\tilde S_I}}(M_I,F)),DR(fu_{IJ}))
=:Rf_*DR(X)((M_I,F),u_{IJ})
\end{eqnarray*}
$w(Y\times\tilde S_I)$ being the wedge product, which is an isomorphism by
by the projection formula for quasi-coherent sheaves for a morphism of ringed topos.
\end{itemize}

Let $k$ a field of characteristic zero. Let $S\in\SmVar(k)$.
Recall we denote for $M,N\in C_{\mathcal D}(S)$, 
\begin{eqnarray*}
m(M,N):\mathcal Hom_{D_S}(M,D_S)\otimes_{D_S} N\to\mathcal Hom_{D_S}(M,N), \;
(\phi\otimes n)\mapsto (m\mapsto\phi(m)n)
\end{eqnarray*}
the multiplication map in $C(S)$.
It induces in the derived category for $M,N\in C_{\mathcal D}(S)$ the map in $D(S)$
\begin{eqnarray*}
m(L_DM,N):R\mathcal Hom_{D_S}(M,D_S)\otimes_{D_S}^LN=\mathcal Hom_{D_S}(L_DM,D_S)\otimes_{D_S} N \\
\to\mathcal Hom_{D_S}(L_DM,N)=R\mathcal Hom_{D_S}(M,N).
\end{eqnarray*}
We use for the proof of theorem \ref{DRKk} in the next subsection the following :

\begin{prop}\label{homDS}
Let $k$ a field of characteristic zero. Let $S\in\SmVar(k)$.
\begin{itemize}
\item[(i)] Let $M,N\in C_{\mathcal D}(S)$. If $N\in C_{\mathcal D,c}(S)$,
\begin{eqnarray*}
m(L_DM,N):R\mathcal Hom_{D_S}(M,D_S)\otimes_{D_S}^LN=\mathcal Hom_{D_S}(L_DM,D_S)\otimes_{D_S} N \\
\to\mathcal Hom_{D_S}(L_DM,N)=R\mathcal Hom_{D_S}(M,N).
\end{eqnarray*}
is an isomorphism in $D(S)$.
\item[(ii)] Let $M,N\in C_{\mathcal D}(S)$. 
If $N\in C_{\mathcal D,c}(S)$, we have using (i) a canonical isomorphism in $D(S)$
\begin{eqnarray*}
D(M,N):R\mathcal Hom_{D_S}(M,N)\xrightarrow{m(L_DM,N)^{-1}}R\mathcal Hom_{D_S}(M,D_S)\otimes_{D_S}^LN \\
\xrightarrow{=:}K_S\otimes_{O_S}^LL\mathbb D_SM[-d_S]\otimes_{D_S}^LN 
\xrightarrow{=}K_S\otimes_{D_S}^L\mathbb D_SM\otimes_{O_S}^LN[-d_S]=:DR(S)^{[-]}(L\mathbb D_SM\otimes_{O_S}^LN)
\end{eqnarray*}
\end{itemize}
\end{prop}

\begin{proof}
\noindent(i):Standard.

\noindent(ii):Follows from (i).
\end{proof}

\subsubsection{Some complements on the (filtered)De Rahm functor for D modules on smooth algebraic varieties 
over a subfield $k\subset\mathbb C$}

In this section, for $S\in\AnSm(\mathbb C)$, we write for short $DR(S):=DR(S)^{[-]}$, where
we recall for $S$ connected $DR(S)^{[-]}:=DR(S)[-d_S]$.

For $S\in\AnSp(\mathbb C)$, we denote by
\begin{eqnarray*}
\alpha(S):\mathbb C_S\hookrightarrow DR(S)(O_S)
\end{eqnarray*}
the inclusion map in $C(S)$.
In particular, we get for $S\in\Var(\mathbb C)$, the inclusion map
\begin{eqnarray*}
\alpha(S):\mathbb C_S\hookrightarrow DR(S)(O_{S^{an}})
\end{eqnarray*}
in $C(S^{an})$.

For $S\in\AnSp(\mathbb C)$, we denote by
\begin{eqnarray*}
\iota(S)::DR(S)(O_S)\to K_S, \; h\otimes w\mapsto hw
\end{eqnarray*}
the canonical map in $C(S)$.
In particular, we get for $S\in\Var(\mathbb C)$, the canonical map
\begin{eqnarray*}
\iota(S):DR(S)(O_{S^{an}})\to K_S
\end{eqnarray*}
in $C(S^{an})$.

\begin{itemize}
\item Let $f:X\to S$ a morphism with $S,X\in\AnSm(\mathbb C)$.
Recall that we have for $(M,F)\in C_{\mathcal Dfil}(X)$ the canonical map in $D_{fil}(S)$
\begin{eqnarray*}
T_*(f,DR)(M,F):DR(S)(\int_f(M,F)):=(\Omega_S^{\bullet},F_b)\otimes_{O_S}Rf_*((D_{X\leftarrow S},F^{ord})\otimes_{D_X}^L(M,F)) \\
\xrightarrow{\iota(S)\otimes I}
Rf_*((D_{X\leftarrow S},F^{ord})\otimes_{D_X}^L(M,F))\otimes^L_{D_S}(K_S,F^{ord}) \\
\xrightarrow{T(f,\otimes)((D_{X\leftarrow S},F^{ord})\otimes_{D_X}^L(M,F),(K_S,F^{ord}))} \\
Rf_*(f^*(K_S,F^{ord})\otimes^L_{f^*D_S}(D_{X\leftarrow S},F^{ord})\otimes_{D_X}^L(M,F)) 
\xrightarrow{=}Rf_*(K_X\otimes_{O_X}((D_X,F^{ord})\otimes_{D_X}^L(M,F)) \\
\xrightarrow{\iota(X)\otimes I}
Rf_*((\Omega_X^{\bullet},F_b)\otimes_{O_X}(M,F))=:Rf_*DR(X)(M,F)
\end{eqnarray*}
which is an isomorphism by the projection formula for quasi-coherent sheaves for a morphism of ringed topos.
In particular, for $S\in\AnSm(\mathbb C)$ and $(M,F)\in C_{\mathcal Dfil,c}(S^o)$,  
\begin{eqnarray*}
T^w(j,\otimes)(M,F):DR(S)(j_*(M,F)):=(\Omega^{\bullet}_S,F_b)\otimes_{O_S}j_*(M,F) \\
\to j_*((\Omega^{\bullet}_{S^o},F_b)\otimes_{O_{S^o}}(M,F))=:j_*DR(S^o)(M,F)
\end{eqnarray*}
is a filtered quasi-isomorphism.
\item Let $f:X\to S$ a morphism with $S,X\in\AnSp(\mathbb C)$. 
Assume there exist a factorization $f:X\xrightarrow{l} Y\times S\xrightarrow{p}S$
with $Y\in\AnSm(\mathbb C)$, $l$ a closed embedding and $p$ the projection.
Let $S=\cup_iS_i$ an open cover such that there exists closed embedding 
$i_i:S_i\hookrightarrow\tilde S_i$ with $\tilde S_i\in\AnSm(\mathbb C)$.
We have for $((M_I,F),u_{IJ})\in C_{\mathcal Dfil}(X/(Y\times\tilde S_I))$ the canonical map in $D_{fil}(S/(\tilde S_I))$
\begin{eqnarray*}
T_*(f,DR)((M_I,F),u_{IJ}): \\ 
DR(S)(\int_f((M_I,F),u_{IJ})):=
((\Omega_{\tilde S_I}^{\bullet},F_b)\otimes_{O_{\tilde S_I}}
p_{\tilde S_I*}E((\Omega^{\bullet}_{Y\times\tilde S_I/\tilde S_I},F_b)\otimes_{O_{Y\times\tilde S_I}}(M_I,F)),DR(fu_{IJ})) \\
\xrightarrow{(k\circ T(p_{\tilde S_I},\otimes)(-,-))} 
(p_{\tilde S_I*}E((p_{\tilde S_I}^*\Omega_{\tilde S_I}^{\bullet},F_b)\otimes_{p_{\tilde S_I}^*O_{\tilde S_I}}
(\Omega^{\bullet}_{Y\times\tilde S_I/\tilde S_I},F_b)\otimes_{O_{Y\times\tilde S_I}}(M_I,F)),DR(fu_{IJ})) \\
\xrightarrow{w(Y\times\tilde S_I)} 
p_{\tilde S_I*}E((\Omega^{\bullet}_{Y\times\tilde S_I},F_b)\otimes_{O_{Y\times\tilde S_I}}(M_I,F)),DR(fu_{IJ}))
=:Rf_*DR(X)((M_I,F),u_{IJ})
\end{eqnarray*}
$w(Y\times\tilde S_I)$ being the wedge product,
which is an isomorphism by the projection formula for quasi-coherent sheaves for a morphism of ringed topos.
\item Let $S\in\AnSm(\mathbb C)$.
Recall that we have for $(M,F)\in C_{\mathcal Dfil}(S)$ the canonical map in $D_{fil}(S)$
\begin{eqnarray*}
T(D,DR)(M,F):DR(S)(L\mathbb D_S(M,F))\to\mathbb D^v_S(DR(S)(M,F))
\end{eqnarray*}
\item Let $S\in\AnSp(\mathbb C)$.
Let $S=\cup_iS_i$ an open cover such that there exists closed embedding 
$i_i:S_i\hookrightarrow\tilde S_i$ with $\tilde S_i\in\AnSm(\mathbb C)$.
Recall that we have for $((M_I,F),u_{IJ})\in C_{\mathcal Dfil}(S/(\tilde S_I))$ the canonical map in $D_{fil}(S/(\tilde S_I))$
\begin{eqnarray*}
T(D,DR)(((M_I,F),u_{IJ})):DR(S)(L\mathbb D_S((M_I,F),u_{IJ}))\to\mathbb D^v_S(DR(S)((M_I,F),u_{IJ}))
\end{eqnarray*}
\item Let $f:X\to S$ a morphism with $S,X\in\AnSm(\mathbb C)$.
Recall that we have for $(N,F)\in C_{\mathcal Dfil}(S)$ the canonical map in $D_{fil}(X)$
\begin{eqnarray*}
T^*(f,DR)(M,F):f^*DR(S)(N,F)\xrightarrow{\iota_S}f^*\mathcal Hom_{D_S}(O_S,L_D(N,F)) \\
\xrightarrow{T(f,hom)(-,-)}\mathcal Hom_{f^*D_S}(f^*O_S,f^*L_D(N,F)) \\
\xrightarrow{Tr(-,-)}
\mathcal Hom_{D_X}(f^{*mod}O_S,f^{*mod}L_D(N,F))=\mathcal Hom_{D_X}(O_X,f^{*mod}L_D(N,F)) \\
\xrightarrow{(\iota(X)\otimes I)^{-1}}
\Omega_X^{\bullet}\otimes_{O_X}f^{*mod}L_D(N,F)=:DR(X)(Lf^{*mod}(N,F)) 
\end{eqnarray*}
\end{itemize}

Let $k\subset\mathbb C$ a subfield.
\begin{itemize}
\item Let $f:X\to S$ a morphism with $S,X\in\Var(k)$. 
Assume there exist a factorization $f:X\xrightarrow{l} Y\times S\xrightarrow{p}S$
with $Y\in\SmVar(k)$, $l$ a closed embedding and $p$ the projection.
Let $S=\cup_iS_i$ an affine open cover so that there exists closed embedding 
$i_i:S_i\hookrightarrow\tilde S_i$ with $\tilde S_i\in\SmVar(k)$.
We have for $((M_I,F),u_{IJ})\in C_{\mathcal Dfil}(S/(\tilde S_I))$ the canonical map 
in $D_{fil}(X_{\mathbb C}^{an}/(Y\times\tilde S_I)_{\mathbb C}^{an})$
\begin{eqnarray*}
f^!DR(S)((M_I,F),u_{IJ})=\Gamma_X\mathbb D^vp^!\mathbb DDR(S)(((M_I,F),u_{IJ})^{an}) \\
\xrightarrow{\mathbb DT^*(p,DR)(-)}\Gamma_XDR(Y\times S)((p_{\tilde S_I}^{*mod}((M_I,F),u_{IJ})^{an}) \\
\xrightarrow{T^w(j,\otimes)(-)^{-1}}DR(Y\times S)(\Gamma_X(p_{\tilde S_I}^{*mod}(M_I,F),u_{IJ})^{an})
\xrightarrow{DR(Y\times S)(T(\gamma,an)(-):=(I,T(j,an)(-)))} \\
DR(Y\times S)((\Gamma_Xp_{\tilde S_I}^{*mod}((M_I,F),u_{IJ}))^{an})=:DR(X)(f^{*mod,\Gamma}((M_I,F),u_{IJ})).
\end{eqnarray*}
\item Let $S\in\Var(k)$. Denote by $\Delta_S:S\hookrightarrow S\times S$ the graph embedding
and $p_1:S\times S\to S$ and $p_2:S\times S\to S$ the projections.
Let $S=\cup_iS_i$ an affine open cover so that there exists closed embedding 
$i_i:S_i\hookrightarrow\tilde S_i$ with $\tilde S_i\in\SmVar(k)$.
We have for $((M_I,F),u_{IJ}),((N_I,F),u_{IJ})\in C_{\mathcal Dfil}(S/(\tilde S_I))$ the canonical map 
in $D_{fil}(S_{\mathbb C}^{an}/(\tilde S^{an}_{I,\mathbb C}))$
\begin{eqnarray*}
T(\otimes,DR)((M,F)(M_I,F),u_{IJ}),((N_I,F),u_{IJ})):DR(S)((M,F)^{an})\otimes_{\mathbb C_S}DR(S)(((N_I,F),u_{IJ})^{an}) \\
\xrightarrow{=}R\Delta_S^!((p_1^*DR(S)((M_I,F),u_{IJ})^{an})\otimes_{\mathbb C_S}p_2^*DR(S)(((N_I,F),u_{IJ})^{an})) \\
\xrightarrow{=}R\Delta_S^!((p_1^!DR(S)((M_I,F),u_{IJ})^{an})\otimes_{\mathbb C_S}p_2^!DR(S)(((N_I,F),u_{IJ})^{an}))[2d_S] \\
\xrightarrow{T^!(p_1,DR)(-)\otimes T^!(p_2,DR)(-)}
\Delta_S^!DR(S\times S)(((p_1^{*mod}(M_I,F),u_{IJ})\otimes_{O_{S\times S}}p_2^{*mod}((N_I,F),u_{IJ}))^{an}) \\
\xrightarrow{T^!(\Delta_S,DR)(-)}
DR(S)((L\Delta_S^{*mod}((p_1^{*mod}(M_I,F),u_{IJ})\otimes_{O_{S\times S}}p_2^{*mod}((N_I,F),u_{IJ})))^{an}) \\ 
\xrightarrow{=}DR(S)((((M_I,F),u_{IJ})\otimes^L_{O_S}((N_I,F),u_{IJ}))^{an}).
\end{eqnarray*}
\end{itemize}

\begin{thm}\label{DRKkfil0}
Let $k\subset\mathbb C$ a subfield.
\begin{itemize}
\item[(i)] Let $j:S^o\hookrightarrow S$ an open embedding with $S\in\SmVar(k)$.
Then, for $M\in C_{\mathcal D,rh}(S^o)$, the map in $C(S_{\mathbb C}^{an})$
\begin{eqnarray*}
DR(S)(T(j,an)(M)):DR(S)((j_*M)^{an})\to DR(S)(j_*(M^{an}))
\end{eqnarray*}
is a quasi-isomorphism.
\item[(ii)] Let $S\in\Var(k)$. 
Let $S=\cup_iS_i$ an affine open cover so that there exists closed embedding 
$i_i:S_i\hookrightarrow\tilde S_i$ with $\tilde S_i\in\SmVar(k)$.
For $((M_I,F),u_{IJ})\in C_{\mathcal Dfil}(S/(\tilde S_I))$ the canonical map in 
$D_{fil}(S_{\mathbb C}^{an}/(\tilde S_{I\mathbb C}^{an}))$
\begin{eqnarray*}
T(D,DR)(((M_I,F),u_{IJ})):DR(S)(L\mathbb D_S((M_I,F),u_{IJ}))\to\mathbb D^v_S(DR(S)((M_I,F),u_{IJ}))
\end{eqnarray*}
is an isomorphism.
\item[(iii)]Let $f:X\to S$ a morphism with $S,X\in\SmVar(k)$. 
Then, for $(M,F)\in C_{\mathcal Dfil,rh}(X)$, the map in $D_{fil}(S_{\mathbb C}^{an})$
\begin{eqnarray*}
T_*(f,DR)(M,F):DR(S)((\int_f(M,F))^{an})\xrightarrow{DR(S)(T(f,an)(-))}DR(S)(\int_f(M,F)^{an}) \\
\xrightarrow{T_*(f,DR)((M,F)^{an})}Rf_*DR(X)((M,F)^{an})
\end{eqnarray*}
is an isomorphism if $f$ is proper and $o_{fil}T(f,DR)(M,F)=:T(f,DR)(M)$ is an isomorphism in $D(S_{\mathbb C}^{an})$.
\item[(iii)'] Let $f:X\to S$ a morphism with $S,X\in\Var(k)$. 
Assume there exist a factorization $f:X\xrightarrow{l} Y\times S\xrightarrow{p}S$
with $Y\in\SmVar(k)$, $l$ a closed embedding and $p$ the projection.
Let $S=\cup_iS_i$ an affine open cover so  that there exists closed embedding 
$i_i:S_i\hookrightarrow\tilde S_i$ with $\tilde S_i\in\SmVar(k)$.
For $((M_I,F),u_{IJ})\in C_{\mathcal Dfil}(X/(Y\times\tilde S_I))$ the canonical map in 
$D_{fil}(S_{\mathbb C}^{an}/(\tilde S^{an}_{I,\mathbb C}))$ 
\begin{eqnarray*}
T_*(f,DR)((M_I,F),u_{IJ}):DR(S)(\int_f((M_I,F),u_{IJ})^{an})\xrightarrow{(T(p_{\tilde S_I},an))} \\
DR(S)(\int_f(((M_I,F),u_{IJ})^{an}))\xrightarrow{T_*(f,DR)(((M_I,F),u_{IJ})^{an})}
Rf_*DR(X)(((M_I,F),u_{IJ})^{an})
\end{eqnarray*}
is an isomorphism if $f$ is proper, and $o_{fil}T(f,DR)((M_I,F),u_{IJ})=:T(f,DR)((M_I,u_{IJ}))$ 
is an isomorphism in $D(S_{\mathbb C}^{an}/(\tilde S^{an}_{I,\mathbb C}))$.
\item[(iv)] Let $S\in\SmVar(k)$. Then, for $M,N\in C_{\mathcal D,rh}(S)$, the map in $D(S_{\mathbb C}^{an})$
\begin{eqnarray*}
T(\otimes,DR)(M,N):DR(S)(M^{an})\otimes_{\mathbb C_{S_{\mathbb C}^{an}}}DR(S)(N^{an}) \\
\xrightarrow{T(\otimes,DR)(M^{an},N^{an})}
DR(S)(M^{an}\otimes_{O_S}N^{an})=DR(S)((M\otimes_{O_S}N)^{an})
\end{eqnarray*}
is an isomorphism.
\end{itemize}
\end{thm}

\begin{proof}
\noindent(i): Follows from the complex case : see \cite{LvDmod}.

\noindent(ii): Follows from the complex case which is standard.

\noindent(iii) and (iii)': Follows from GAGA for proper morphisms of complex algebraic varieties.

\noindent(iv): Follows from (i).
\end{proof}

\begin{thm}\label{DRKk}
Let $k\subset\mathbb C$ a subfield. Let $S\in\SmVar(k)$. 
\begin{itemize}
\item[(i)] For $M,N\in D_{\mathcal D,rh}(S)$, we have using
proposition \ref{homDS} and theorem \ref{DRKkfil0}(iv) and (ii) the
following canonical isomorphism in $D(\mathbb C)$
\begin{eqnarray*}
DR(M,N):R\Hom_{D_S}(M,N)\otimes_k\mathbb C=Ra_{S*}R\mathcal Hom_{D_S}(M,N)\otimes_k\mathbb C \\
\xrightarrow{Ra_{S*}D(M,N)}
\int_{a_S}(L\mathbb D_SM\otimes_{O_S}^LN)\otimes_k\mathbb C=
DR(pt)\int_{a_S}(L\mathbb D_SM\otimes_{O_S}^LN)\otimes_k\mathbb C \\
\xrightarrow{(T(\otimes,DR)(-,-)\circ T_*(a_S,DR)(-)^{-1}}
Ra_{S*}(\mathbb D_S^vDR(S)(M^{an})\otimes DR(S)(N^{an})) \\
\xrightarrow{Ra_{S*}m(DR(S)(M^{an}),DR(S)(N^{an}))}R\Hom_{\mathbb C_{S_{\mathbb C}^{an}}}(DR(S)(M^{an}),DR(S)(N^{an}))
\end{eqnarray*}
is an isomorphism.
\item[(i)'] For $M,N\in D_{\mathcal D,rh}(S)$, the canonical map in $D(\mathbb C)$
\begin{equation*}
DR(S)^{L_DM,N}:R\Hom_{D_S}(M,N)\otimes_k\mathbb C\xrightarrow{\sim}
R\Hom_{\mathbb C_{S_{\mathbb C}^{an}}}(DR(S)(M^{an}),DR(S)(N^{an}))
\end{equation*}
is equal to $DR(M,N)$, hence is an isomorphism.
\item[(ii1)] We have $DR(S)(D_{\mathcal D,rh}(S))\subset D_{\mathbb C_S,c}(S_{\mathbb C}^{an})$,
that is the image of the class of a complex of $D_S$ module with regular holonomic cohomology sheaves
is a complex of presheaves on $S_{\mathbb C}^{an}$ whose cohomology sheaves are constructible
for a Zariski stratification of $S$ defined over $k$.
\item[(ii2)] For $M\in\PSh_{\mathcal D,rh}(S)$, $DR(S)(M)\in P(S_{\mathbb C}^{an})$
that is is a perverse sheaf for a Zariski stratification of $S$ defined over $k$.
\end{itemize}
\end{thm}

\begin{proof}
\noindent(i): Follows from proposition \ref{homDS} and theorem \ref{DRKkfil0}(iv) and (ii).

\noindent(i)': Follows from the following commutative diagram in $D(\mathbb C)$
\begin{equation*}
\xymatrix{R\Hom_{D_S}(M,N)\otimes_k\mathbb C=Ra_{S*}R\mathcal Hom_{D_S}(M,N)\otimes_k\mathbb C
\ar[r]^{DR(S)^{L_DM,N}}\ar[d]^{Ra_{S*}D(M,N)} &
R\Hom_{\mathbb C_{S_{\mathbb C}^{an}}}(DR(S)(M^{an}),DR(S)(N^{an})) \\
\int_{a_S}(L\mathbb D_SM\otimes_{O_S}^LN)\otimes_k\mathbb C=
DR(pt)\int_{a_S}(L\mathbb D_SM\otimes_{O_S}^LN)\otimes_k\mathbb C
\ar[r]^{(T(\otimes,DR)(-,-)\circ T(a_S,DR)(-)^{-1}} &
Ra_{S*}(\mathbb D_S^vDR(S)(M^{an})\otimes DR(S)(N^{an}))\ar[u]^{m(-,-)}}
\end{equation*}

\noindent(ii):Similar to the proof of the complex case in \cite{LvDmod}: 
follows by definition from the locally free case by theorem \ref{holkstr}.
\end{proof}

\subsubsection{On the De Rahm functor for D modules on smooth algebraic varieties 
over a p-adic field $K\subset\mathbb C_p$}

For $S\in\AnSp(K)$, we denote by
\begin{eqnarray*}
\alpha(S):\mathbb B_{dr,S}\hookrightarrow DR(S)(O\mathbb B_{dr,S})
\end{eqnarray*}
the inclusion map in $C_{\mathbb B_{dr,S}}(S^{pet})$.
In particular, we get for $S\in\Var(K)$, the inclusion map
\begin{eqnarray*}
\alpha(S):\mathbb B_{dr,S}\hookrightarrow DR(S)(O\mathbb B_{dr,S})
\end{eqnarray*}
in $C_{\mathbb B_{dr,S}}(S^{an,pet})$.

For $S\in\AnSp(K)$, we denote by
\begin{eqnarray*}
\iota(S)::DR(S)(O\mathbb B_{dr,S})\to\mathbb B_{dr,S}\otimes_{O_S} K_S, \;
h\otimes k\otimes w\mapsto k\otimes (hw)
\end{eqnarray*}
the canonical map in $C_{\mathbb B_{dr,S}}(S^{pet})$.
In particular, we get for $S\in\Var(K)$, the canonical map
\begin{eqnarray*}
\iota(S):DR(S)(O\mathbb B_{dr,S})\to\mathbb B_{dr,S}\otimes_{O_S} K_S
\end{eqnarray*}
in $C_{\mathbb B_{dr,S}}(S^{an,pet})$.

We have the following theorem 

\begin{thm}\label{RHpadic}
Let $K\subset\mathbb C_p$ a p adic field.  
\begin{itemize}
\item[(i)]Let $S\in\AnSm(K)$. The inclusion map in $C_{\mathbb B_{dr,S}}(S^{an,pet})$
\begin{eqnarray*}
\alpha(S):\mathbb B_{dr,S}\hookrightarrow DR(S)(O\mathbb B_{dr,S}).
\end{eqnarray*}
is a quasi-isomorphism.
\item[(i)']Let $S\in\AnSm(K)$.the canonical map in $C_{\mathbb B_{dr,S}}(S^{an,pet})$
\begin{eqnarray*}
\iota(S):DR(S)(O\mathbb B_{dr,S})\to\mathbb B_{dr,S}\otimes K_S
\end{eqnarray*}
is a quasi-isomorphism.
\item[(ii)]Let $S\in\Var(K)$. Let $D=\cup D_i\subset S$ a normal crossing divisor 
and denote $j:S^o:=S\backslash D\hookrightarrow S$ the open embedding.
The inclusion map in $C_{\mathbb B_{dr,S}}(S^{an,pet})$
\begin{eqnarray*}
\alpha(S):\mathbb B_{dr,S}(\log D)\hookrightarrow
\Omega_S^{\bullet}(\log D)\otimes_{O_S}O\mathbb B_{dr,S}(\log D)
=F^0DR(S)(j_{*Hdg}(O_{S^o},F_b)\otimes_{O_S}(O\mathbb B_{dr,S},F)).
\end{eqnarray*}
where $j_{*Hdg}(O_S,F_b)=(j_*O_{S^o},V_D)$ with $V_D$ the $V$-filtration (see section 5) that is the
filtration by order of the pole in this case, is a quasi-isomorphism.
\item[(iii)]Let $S\in\AnSm(K)$. The functor
\begin{eqnarray*}
\Vect_{\mathcal D}(S)\to D_{\mathbb B_{dr,S}}(S^{an,pet}), \;
(M,F)\mapsto F^0DR(S)((M,F)^{an}\otimes_{O_S}(O\mathbb B_{dr,S},F))
\end{eqnarray*}
is fully faithful whose inverse on the image is given by 
\begin{eqnarray*}
N\in\Shv_{\mathbb B_{dr,S}}(S^{an,pet})\mapsto Re_*((N,F)\otimes_{\mathbb B_{dr,S}}(O\mathbb B_{dr,S},F)).
\end{eqnarray*}
where $e:S^{pet}\to S^{et}$ is the morphism of site given by the inclusion functor.
\end{itemize}
\end{thm}

\begin{proof}
\noindent(i):See \cite{Scholze}.

\noindent(i)': Follows from (i) by duality.

\noindent(i):See \cite{Chinois}.

\noindent(ii):See \cite{Chinois}.
\end{proof}

\begin{defi}\label{TfDRpadic}
Let $K\subset\mathbb C_p$ a p adic field. 
\begin{itemize}
\item[(i)]Let $f:X\to S$ a morphism with $S,X\in\AnSm(K)$.
We have for $(M,F)\in C_{\mathcal Dfil}(X)$ the canonical map in $D_{\mathbb B_{dr}fil}(S)$
\begin{eqnarray*}
T^{B_{dr}}(f,DR)(M,F):DR(S)(\int_f(M,F)\otimes_{O_S}(OB_{dr,S},F)) \\
\xrightarrow{:=} 
(\Omega_S^{\bullet},F_b)\otimes_{O_S}(OB_{dr,S},F)\otimes_{O_S}Rf_*((D_{X\leftarrow S},F^{ord})\otimes_{D_X}^L(M,F)) \\
\xrightarrow{\iota(S)\otimes I}
Rf_*((D_{X\leftarrow S},F^{ord})\otimes_{D_X}^L(M,F))\otimes_{D_S}\mathbb B_{dr,S}\otimes_{O_S}K_S \\
\xrightarrow{T(f,\otimes)((D_{X\leftarrow S},F^{ord})\otimes_{D_X}^L(M,F),\mathbb B_{dr,S}\otimes_{O_S}K_S)} \\
Rf_*(f^*K_S\otimes^L_{f^*D_S}(D_{X\leftarrow S},F^{ord})\otimes_{D_X}^L(M,F)\otimes_{f^*O_S}f^*\mathbb B_{dr,S}) \\
\xrightarrow{=}Rf_*(K_X\otimes_{O_X}\mathbb B_{dr,X}\otimes_{D_X}^L(M,F))
\xrightarrow{\iota(X)}Rf_*((\Omega_X^{\bullet},F_b)\otimes_{O_X}(M,F)\otimes_{O_X}(O\mathbb B_{dr,X},F)) \\
=:Rf_*DR(X)((M,F)\otimes_{O_X}(O\mathbb B_{dr,X},F))
\end{eqnarray*}
which is an isomorphism by the projection formula for quasi-coherent modules for morphisms of ringed topos (see \cite{B4}).
\item[(i)'] Let $f:X\to S$ a morphism with $S,X\in\AnSp(K)$. 
Assume there exist a factorization $f:X\xrightarrow{l} Y\times S\xrightarrow{p}S$
with $Y\in\AnSm(K)$, $l$ a closed embedding and $p$ the projection.
Let $S=\cup_iS_i$ an affine open cover so that there exists closed embedding 
$i_i:S_i\hookrightarrow\tilde S_i$ with $\tilde S_i\in\AnSm(K)$.
We have for $((M_I,F),u_{IJ})\in C_{\mathcal Dfil}(X/(Y\times\tilde S_I))$ 
the canonical map in $D_{\mathbb B_{dr}fil}(S/(\tilde S_I))$
\begin{eqnarray*}
T^{B_{dr}}(f,DR)((M_I,F),u_{IJ}): 
DR(S)(\int_f((M_I,F),u_{IJ})\otimes_{O_S}(OB_{dr,(\tilde S_I)},F)):= \\
((\Omega_{\tilde S_I}^{\bullet},F_b)\otimes_{O_{\tilde S_I}}(OB_{dr,(\tilde S_I)},F)\otimes_{O_{\tilde S_I}}
p_{\tilde S_I*}E((\Omega^{\bullet}_{Y\times\tilde S_I/\tilde S_I},F_b)\otimes_{O_{Y\times\tilde S_I}}(M_I,F)),DR(fu_{IJ})) \\
\xrightarrow{(k\circ T(p_{\tilde S_I},\otimes)(-,-))} \\
(p_{\tilde S_I*}E((p_{\tilde S_I}^*\Omega_{\tilde S_I}^{\bullet},F_b)\otimes_{p_{\tilde S_I}^*O_{\tilde S_I}}
p_{\tilde S_I}^*(OB_{dr,(\tilde S_I)},F)\otimes_{p_{\tilde S_I}^*O_{\tilde S_I}}
(\Omega^{\bullet}_{Y\times\tilde S_I/\tilde S_I},F_b)\otimes_{O_{Y\times\tilde S_I}}(M_I,F)),DR(fu_{IJ})) \\
\xrightarrow{=} \\
(p_{\tilde S_I*}E((p_{\tilde S_I}^*\Omega_{\tilde S_I}^{\bullet},F_b)\otimes_{p_{\tilde S_I}^*O_{\tilde S_I}}
(OB_{dr,(Y\times\tilde S_I)},F)\otimes_{O_{Y\times\tilde S_I}}
(\Omega^{\bullet}_{Y\times\tilde S_I/\tilde S_I},F_b)\otimes_{O_{Y\times\tilde S_I}}(M_I,F)),DR(fu_{IJ})) \\
\xrightarrow{w(Y\times\tilde S_I)} 
(p_{\tilde S_I*}E((\Omega^{\bullet}_{Y\times\tilde S_I},F_b)\otimes_{O_{Y\times\tilde S_I}}(OB_{dr,(Y\times\tilde S_I)},F)
\otimes_{O_{Y\times\tilde S_I}}(M_I,F)),DR(u_{IJ})) \\
=:Rf_*DR(X)(((M_I,F),u_{IJ})\otimes_{O_X}(OB_{dr,(Y\times\tilde S_I)},F))
\end{eqnarray*}
where $w(Y\times\tilde S_I)$ is the wedge product, 
which is an isomorphism by the projection formula for quasi-coherent modules for morphisms of ringed topos (see \cite{B4}).
\item[(ii)] Let $S\in\AnSm(K)$.
We have for $(M,F)\in C_{\mathcal Dfil}(S)$ the canonical map in $D_{fil}(S)$
\begin{eqnarray*}
T(D,DR)(M,F):DR(S)(\mathbb D_S(M,F))\to\mathbb D^v_S(DR(S)(M,F))
\end{eqnarray*}
\item[(ii)']Let $S\in\AnSp(K)$. 
Let $S=\cup_iS_i$ an affine open cover so that there exists closed embedding 
$i_i:S_i\hookrightarrow\tilde S_i$ with $\tilde S_i\in\AnSm(K)$.
We have for $((M_I,F),u_{IJ})\in C_{\mathcal Dfil}(S/(\tilde S_I))$ the canonical map in $D_{fil}(S/(\tilde S_I))$
\begin{eqnarray*}
T(D,DR)((M_I,F),u_{IJ}):DR(S)(L\mathbb D_S((M_I,F),u_{IJ}))\to\mathbb D^v_S(DR(S)((M_I,F),u_{IJ}))
\end{eqnarray*}
\item Let $f:X\to S$ a morphism with $S,X\in\Var(K)$. 
Assume there exist a factorization $f:X\xrightarrow{l} Y\times S\xrightarrow{p}S$
with $Y\in\SmVar(K)$, $l$ a closed embedding and $p$ the projection.
Let $S=\cup_iS_i$ an affine open cover so that there exists closed embedding 
$i_i:S_i\hookrightarrow\tilde S_i$ with $\tilde S_i\in\SmVar(K)$.
We have for $((M_I,F),u_{IJ})\in C_{\mathcal Dfil}(S/(\tilde S_I))$ the canonical map 
in $D_{fil}(X^{an}/(Y\times\tilde S_I)^{an})$
\begin{eqnarray*}
f^!DR(S)((M_I,F),u_{IJ})=\Gamma_X\mathbb D^vp^!\mathbb DDR(S)(((M_I,F),u_{IJ})^{an}) \\
\xrightarrow{\mathbb DT^*(p,DR)(-)}\Gamma_XDR(Y\times S)((p_{\tilde S_I}^{*mod}((M_I,F),u_{IJ})^{an}) \\
\xrightarrow{T^w(j,\otimes)(-)^{-1}}DR(Y\times S)(\Gamma_X(p_{\tilde S_I}^{*mod}(M_I,F),u_{IJ})^{an})
\xrightarrow{DR(Y\times S)(T(\gamma,an)(-):=(I,T(j,an)(-)))} \\
DR(Y\times S)((\Gamma_Xp_{\tilde S_I}^{*mod}((M_I,F),u_{IJ}))^{an})=:DR(X)(f^{*mod,\Gamma}((M_I,F),u_{IJ})).
\end{eqnarray*}
\item Let $S\in\Var(K)$. Denote by $\Delta_S:S\hookrightarrow S\times S$ the graph embedding
and $p_1:S\times S\to S$ and $p_2:S\times S\to S$ the projections.
Let $S=\cup_iS_i$ an affine open cover so that there exists closed embedding 
$i_i:S_i\hookrightarrow\tilde S_i$ with $\tilde S_i\in\SmVar(K)$.
We have for $((M_I,F),u_{IJ}),((N_I,F),u_{IJ})\in C_{\mathcal Dfil}(S/(\tilde S_I))$ the canonical map 
in $D_{fil}(S^{an}/(\tilde S^{an}_I))$
\begin{eqnarray*}
T(\otimes,DR)((M_I,F),u_{IJ}),((N_I,F),u_{IJ})):
DR(S)(((M_I,F),u_{IJ})^{an})\otimes_{\mathbb C_S}DR(S)(((N_I,F),u_{IJ})^{an}) \\
\xrightarrow{=}R\Delta_S^!((p_1^*DR(S)((M_I,F),u_{IJ})^{an})\otimes_{\mathbb C_S}p_2^*DR(S)(((N_I,F),u_{IJ})^{an})) \\
\xrightarrow{=}R\Delta_S^!((p_1^!DR(S)((M_I,F),u_{IJ})^{an})\otimes_{\mathbb C_S}p_2^!DR(S)(((N_I,F),u_{IJ})^{an}))[2d_S] \\
\xrightarrow{T^!(p_1,DR)(-)\otimes T^!(p_2,DR)(-)}
\Delta_S^!DR(S\times S)(((p_1^{*mod}(M_I,F),u_{IJ})\otimes_{O_{S\times S}}p_2^{*mod}((N_I,F),u_{IJ}))^{an}) \\
\xrightarrow{T^!(\Delta_S,DR)(-)}
DR(S)((L\Delta_S^{*mod}((p_1^{*mod}(M_I,F),u_{IJ})\otimes_{O_{S\times S}}p_2^{*mod}((N_I,F),u_{IJ})))^{an}) \\ 
\xrightarrow{=}DR(S)((((M_I,F),u_{IJ})\otimes^L_{O_S}((N_I,F),u_{IJ}))^{an}).
\end{eqnarray*}
\end{itemize}
\end{defi}

\begin{thm}\label{DRKkfil1}
\begin{itemize}
\item[(i)] Let $j:S^o\hookrightarrow S$ an open embedding with $S\in\SmVar(K)$.
Then, for $M\in C_{\mathcal D,rh}(S^o)$, the map in $C(S^{an})$
\begin{eqnarray*}
DR(S)(T(j,an)(M)):DR(S)((j_*M)^{an})\to DR(S)(j_*(M^{an}))
\end{eqnarray*}
is a quasi-isomorphism.
\item[(ii)] Let $S\in\Var(K)$. 
Then, for $((M_I,F),u_{IJ})\in C_{\mathcal Dfil,rh}(S)$, the map in $D_{\mathbb B_{dr}fil}(S^{an})$
\begin{eqnarray*}
T(D,DR)(M,F):DR(S)(L\mathbb D_S((M_I,F),u_{IJ})^{an})\to\mathbb D_S^v(DR(S)(((M_I,F),u_{IJ})^{an}))
\end{eqnarray*}
is an isomorphism.
\item[(iii)]Let $f:X\to S$ a morphism with $S,X\in\SmVar(K)$. 
Then, for $(M,F)\in C_{\mathcal Dfil,rh}(X)$, the map in $D_{\mathbb B_{dr}fil}(S^{an})$
\begin{eqnarray*}
T^{B_{dr}}(f,DR)(M,F):DR(S)((\int_f(M,F))^{an}\otimes_{O_S}(OB_{dr,S},F)) 
\xrightarrow{DR(S)(T(\int_f,an)(M,F))} \\
DR(S)(\int_f((M,F)^{an})\otimes_{O_S}(OB_{dr,S},F)) 
\xrightarrow{T^{B_{dr}}(f,DR)((M,F)^{an})}Rf_*DR(X)((M,F)^{an}\otimes_{O_X}(O\mathbb B_{dr,X},F))
\end{eqnarray*}
is an isomorphism if $f$ is proper.
\end{itemize}
\item[(iii)']Let $f:X\to S$ a morphism with $S,X\in\Var(K)$. 
Then, for $((M_I,F),u_{IJ})\in C_{\mathcal Dfil}(X/(Y\times\tilde S_I))$ 
the map in $D_{\mathbb B_{dr}fil}(S^{an}/(\tilde S^{an}_I))$
\begin{eqnarray*}
T^{B_{dr}}(f,DR)((M_I,F),u_{IJ}): 
DR(S)(\int_f((M_I,F),u_{IJ})^{an}\otimes_{O_S}(OB_{dr,(\tilde S_I)},F))
\xrightarrow{DR(S)((T(p_{\tilde S_I},an)(-)))} \\
DR(S)(\int_f(((M_I,F),u_{IJ})^{an})\otimes_{O_S}(OB_{dr,(\tilde S_I)},F))
\xrightarrow{T^{B_{dr}}(f,DR)(((M_I,F),u_{IJ})^{an})} \\
Rf_*DR(X)(((M_I,F),u_{IJ})^{an}\otimes_{O_S}(OB_{dr,(Y\times\tilde S_I)},F))
\end{eqnarray*}
is an isomorphism if $f$ is proper.
\end{thm}

\begin{proof}
Similar to the proof of theorem \ref{DRKkfil0}. 
For (iii) and (iii)', we use GAGA for proper morphism of algebraic varieties over a $p$-adic field.
\end{proof}

\section{The De Rham modules over a field $k$ of characteristic $0$ :
the Kashiwara Malgrange $V$-filtration and the Hodge filtration in the geometric case}

\subsection{The Kashiwara Malgrange $V$ filtration for geometric D modules 
on smooth algebraic varieties over a field of characteristic zero and the nearby and vanishing cycle functors.}

Let $k\subset\mathbb C$ a subfield. 
For $S\in\Var(k)$, consider $\pi:=\pi_{k/\mathbb C}(S):S_{\mathbb C}:=S\otimes_k\mathbb C\to S$ the projection 
so that we have the injective map in $\PSh_{O_S}(S_{\mathbb C})$
\begin{equation*}
n_{k/\mathbb C}(O_S):\pi^*O_S\hookrightarrow O_{S_{\mathbb C}}:=\pi^{*mod}O_S:=\pi^*O_S\otimes_k\mathbb C, \; \;
h\mapsto n_{k/\mathbb C}(O_S)(h):=h\otimes 1.
\end{equation*}
For $M\in\PSh_{O_S}(S)$, we have (see section 2) the injective map in $\PSh_{O_S}(S_{\mathbb C})$
\begin{eqnarray*}
n_{O_S/O_{S_{\mathbb C}}}(M):
\pi^*M\hookrightarrow\pi^{*mod}M:=\pi^*M\otimes_{\pi^*O_S}O_{S_{\mathbb C}}=\pi^*M\otimes_k\mathbb C, \; \;
m\mapsto n_{O_S/O_{S_{\mathbb C}}}(M)(m):=m\otimes 1
\end{eqnarray*}
For $S\in\SmVar(k)$ and $M\in\PSh_{\mathcal D}(S)$
\begin{eqnarray*}
n_{O_S/O_{S_{\mathbb C}}}(M):
\pi^*M\hookrightarrow\pi^{*mod}M:=\pi^*M\otimes_{\pi^*O_S}O_{S_{\mathbb C}}=\pi^*M\otimes_k\mathbb C, \; \;
m\mapsto n_{O_S/O_{S_{\mathbb C}}}(M)(m):=m\otimes 1
\end{eqnarray*}
is a morphism in $\PSh_{D_S}(S_{\mathbb C})$, that is a morphism of $\pi^*D_S$ modules.

\begin{defi}\label{VfilKM0}
Let $k$ a field of characteristic zero. Let $S\in\SmVar(k)$. 
\begin{itemize}
\item[(i)] Let $D=V(s)\subset S$ be a smooth (Cartier) divisor, 
where $s\in\Gamma(S,L)$ is a section of the line bundle $L=L_D$ associated to $D$. 
Let $M\in\PSh_{\mathcal D}(S)$. 
A $V_D$-filtration $(M,V_D)$ for $M$ (see \cite{B4}) is called a Kashiwara-Malgrange $V_D$-filtration for $M$ if
\begin{itemize}
\item $V_{D,k}M$ are coherent $V_{D,0}D_S$ modules for all $k\in\mathbb Z$, that is $(M,V_D)$ is a good filtration,
in particular $M\in\PSh_{\mathcal D,c}(S)$ is coherent
\item $sV_{D,k}M=V_{D,k-1}M$ for $k<<0$,
\item all eigenvalues of $s\partial_s:\Gr_{V_D,k}:=V_{D,k}M/V_{D,k-1}M\to\Gr_{V_D,k}M:=V_{D,k}M/V_{D,k-1}M$ 
have real part between $k-1$ and $k$.
\end{itemize}
Almost by definition, a Kashiwara-Malgrange $V_D$-filtration for $M$ if it exists is unique,
so that we denote it by $(M,V_D)\in\PSh_{O_Sfil}(S)$ and $(M,V_D)$ is strict.
In particular if $m:(M_1,F)\to (M_2,F)$ a morphism with $(M_1,F),(M_2,F)\in\PSh_{\mathcal D(2)fil}(S)$
such that $M_1$ and $M_2$ admit the Kashiwara-Malgrange filtration for $D\subset S$,
we have $m(V_{D,q}F^pM_1)\subset V_{D,q}F^pM_2$, that is we get $m:(M_1,F,V_D)\to(M_2,F,V_D)$ a filtered morphism,
and if $0\to M'\to M\to M''\to 0$ is an exact sequence, 
$0\to (M',V_D)\to (M,V_D)\to (M'',V_D)\to 0$ is an exact sequence (strictness).
\item[(ii)]More generally, let $Z=V(s_1,\ldots,s_r)=D_1\cap\cdots\cap D_r\subset S$ be a smooth Zariski closed subset,
where $s_i\in\Gamma(S,L_i)$ is a section of the line bundle $L=L_{D_i}$ associated to $D_i$. 
Let $M\in\PSh_{\mathcal D}(S)$. 
A $V_Z$-filtration $(M,V_Z)$ for $M$ (see \cite{B4}) is called a Kashiwara-Malgrange $V_Z$-filtration for $M$ if
\begin{itemize}
\item $V_{Z,k}M$ are coherent $V_{Z,0}O_S$ modules for all $k\in\mathbb Z$,
\item $\sum_{i=1}^rs_iV_{Z,k}M=V_{Z,k-1}M$ for $k<<0$,
\item all eigenvalues of 
$\sum_{i=1}^rs_i\partial_{s_i}:\Gr_{V_Z,k}M:=V_{Z,k}M/V_{Z,k-1}M\to\Gr^V_{Z,k}M:=V_{Z,k}M/V_{Z,k-1}M$ 
have real part between $k-1$ and $k$.
\end{itemize}
Almost by definition, a Kashiwara-Malgrange $V_Z$-filtration for $M$ if it exists is unique (see \cite{Sa2}),
so that we denote it by $(M,V_Z)\in\PSh_{O_Sfil}(S)$ and $(M,V_Z)$ is strict.
In particular if $m:(M_1,F)\to (M_2,F)$ a morphism with $(M_1,F),(M_2,F)\in\PSh_{\mathcal D(2)fil}(S)$
such that $M_1$ and $M_2$ admit the Kashiwara-Malgrange filtration for $D\subset S$,
we have $m(V_{Z,q}F^pM_1)\subset V_{Z,q}F^pM_2$, that is we get $m:(M_1,F,V_Z)\to(M_2,F,V_Z)$ a filtered morphism,
and if $0\to M'\to M\to M''\to 0$ is an exact sequence, 
$0\to (M',V_Z)\to (M,V_Z)\to (M'',V_Z)\to 0$ is an exact sequence (strictness).
\end{itemize}
\end{defi}

\begin{defi}\label{VfilKM}
Let $k$ a field of characteristic zero. Let $S\in\SmVar(k)$. 
Let $D=V(s)\subset S$ a (Cartier) divisor, where $s\in\Gamma(S,L)$ 
is a section of the line bundle $L=L_D$ associated to $D$. 
We then have the graph embedding $i:S\hookrightarrow L$, $i(x):=(x,s(x))$. 
Let $M\in\PSh_{\mathcal D,c}(S)$.
If the Kashiwara-Malgrange $V_S$-filtration exist on $M':=i_{*mod}M\in\PSh_{\mathcal D,c}(L)$ 
(see definition \ref{VfilKM0}) and the eigeinvalues of $s\partial_s$ are rational numbers,
we refine the $V_S$-filtration to all rational numbers as follows : 
for $\alpha=k-1+r/q\in\mathbb Q$, $k,q,r\in\mathbb Z$, $q\leq 0$, $0\leq r\leq q-1$, we set
\begin{eqnarray*}
V_{S,\alpha}M':=q_{V,k}^{-1}(\oplus_{k-1<\beta\leq\alpha}\Gr_{k,\beta}^{V_S}M\subset V_{S,k}M'
\end{eqnarray*}
with $\Gr_{k,\beta}^{V_S}M':=\ker(\partial_ss-\beta I)\subset\Gr_k^{V_S}M'$ and $q_{V,k}:V_{S,k}M'\to\Gr_k^{V_S}M'$
is the projection.
We set similarly
\begin{eqnarray*}
V_{S,<\alpha}M:=q_{V,k}^{-1}(\oplus_{k-1<\beta <\alpha}\Gr_{k,\beta}^{V_S}M')\subset V_{S,k}M'
\end{eqnarray*}
The Hodge filtration induced on $\Gr^V_{\alpha}M$ is 
\begin{eqnarray*}
F^p\Gr^{V_S}_{\alpha}M':=(F^pM\cap V_{S,\alpha}M')/(F^pM\cap V_{S,<\alpha}M')
\end{eqnarray*}
We call it the rational Kashiwara-Malgrange $V_S$-filtration of $M'=i_{*mod}M\in\PSh_{\mathcal D,c}(L)$.
Using theorem \ref{Keqk}, we set for $\alpha\in\mathbb Q$, 
\begin{equation*}
V_{D,\alpha}M:=i^{\sharp}V_{S,\alpha}i_{*mod}M\subset M=i^{\sharp}i_{*mod}M
\end{equation*}
and $(M,V_D)$ is the rational Kashiwara-Malgrange $V_D$-filtration of $M\in\PSh_{\mathcal D,c}(S)$.
\end{defi}

\begin{defi}\label{DHdgpsi0}
Let $k$ a field of characteristic zero. Let $S\in\SmVar(k)$.
Let $D=V(s)\subset S$ a divisor with $s\in\Gamma(S,L)$ and $L$ a line bundle ($S$ being smooth, $D$ is Cartier).
For $M\in\PSh_{\mathcal D,c}(S)$ such that the rational Kashiwara-Malgrange $V_D$ filtration exists,
that is the $V_S$-filtration on $i_{*mod}M$ exists, we define, using definition 
\begin{itemize}
\item the nearby cycle functor
\begin{equation*}
\psi_DM:=i^{\sharp}(\oplus_{-1\leq\alpha<0}\Gr_{V_S,\alpha}i_{*mod}M)
=\oplus_{-1\leq\alpha<0}\Gr_{V_D,\alpha}M\in\PSh_{\mathcal D,D}(S),
\end{equation*}
\item the unipotent nearby cycle functor
\begin{equation*}
\psi^u_DM:=i^{\sharp}(\Gr_{V_S,-1}i_{*mod}M)=\Gr_{V_D,-1}M\in\PSh_{\mathcal D,D}(S),
\end{equation*}
\item the vanishing cycle functor
\begin{equation*}
\phi_DM:=i^{\sharp}(\oplus_{-1<\alpha\leq 0}\Gr_{V_S,\alpha}i_{*mod}M)
=\oplus_{-1<\alpha\leq 0}\Gr_{V_D,\alpha}M\in\PSh_{\mathcal D,D}(S),
\end{equation*}
\item the unipotent nearby cycle functor
\begin{equation*}
\phi^u_DM:=i^{\sharp}(\Gr_{V_S,0}i_{*mod}M)=\Gr_{V_D,0}M\in\PSh_{\mathcal D,D}(S),
\end{equation*}
\item the canonical maps in $\PSh_{\mathcal D,D}(S)$
\begin{equation*}
can(M):=(\partial_s,I):\psi_DM\to\phi_DM \; , \; var(M):=(I,s):\phi_DM\to\psi_DM.
\end{equation*}
\item the canonical maps in $\PSh_{\mathcal D,D}(S)$
\begin{equation*}
can(M):=\partial_s:\psi^u_DM\to\phi^u_DM \; , \; var(M):=s:\phi^u_DM\to\psi^u_DM.
\end{equation*}
\end{itemize}
By the complex case (see \cite{Kashiwara}), after considering a subfield $k_0\subset k$ and an embedding
$\sigma:k_0\hookrightarrow\mathbb C$ we have
\begin{itemize}
\item for $M\in\PSh_{\mathcal D,h}(S)$, $\psi_DM,\phi_DM\in\PSh_{\mathcal D,h}(S)$.
\item for $M\in\PSh_{\mathcal D,rh}(S)$, $\psi_DM,\phi_DM\in\PSh_{\mathcal D,rh}(S)$.
\end{itemize}
\end{defi}

Let $k\subset\mathbb C$ a subfield. Let $S\in\SmVar(k)$.
Let $D=V(s)\subset S$ a divisor with $s\in\Gamma(S,L)$ and $p:L\to S$ a line bundle ($S$ being smooth, $D$ is Cartier),
so that we have the closed embedding $i:S\hookrightarrow L$, $i(x)=(x,s(x))$ 
and $D=i^{-1}(s_0)$, $s_0$ being the zero section. 
Denote by $l:L^o:=L\backslash S\hookrightarrow L$ the open embedding 
which induces the open embedding $l:=l\times_LS:S^o:=S\backslash D\hookrightarrow S$.
Denote by $\pi:\tilde L^o\to L^o$ the universal covering 
which induces the universal covering $\pi:=\pi\times_{L^o}S^o:\tilde S^o\to S^o$. 
Consider $S=\cup_{i=1}^sS_i$ an open affine cover such that 
$D\cap S_i=V(f_i)\subset S_i$ is given by $f_i\in\Gamma(S_i,O_{S_i})$,
and denote $q:L_i:=p^{-1}(S_i)\to\mathbb A^1_k$ the projection and $j_i:S_i\hookrightarrow S$ the open embeddings.
For $M\in\PSh_{\mathcal D,c}(S)$ such that the rational Kashiwara-Malgrange $V_D$-filtration on $M$ exists, we consider
\begin{itemize}
\item for $-1\leq\alpha<0$, the morphism in $D_{c,k}(S_{\mathbb C}^{an})$
\begin{eqnarray*}
A_{\alpha}(M):\pi_{\alpha}(M):=(\oplus_{i=1}^s
\Cone(\partial_s:DR(L_i/\mathbb A^1_k)((V_{D\alpha}M)^{an})\otimes_{O_S}s^{\alpha+1}O_{S^{an}_{\mathbb C}}[\log s] \\
\to DR(L_i/\mathbb A^1_k)((V_{D\alpha}M)^{an})\otimes_{O_S}s^{\alpha}O_{S^{an}_{\mathbb C}}[\log s])
\xrightarrow{(j_I^*)}\cdots) \\ 
\to(\oplus_{i=1}^sDR(S_i)((\Gr_{V_D,\alpha}M)^{an})\xrightarrow{j_I^*}\cdots)[1]
\xrightarrow{((j_i^*),0)^{-1}}DR(S)((\Gr_{V_D,\alpha}M)^{an})[1], \\  
(\sum_jm_j\otimes(\log s)^j,m')\mapsto [m_0] , 
\end{eqnarray*}
and for $\alpha=0$ the morphism in $D_{c,k}(S_{\mathbb C}^{an})$
\begin{eqnarray*}
A_0(M):\pi_0(M):= \\
\Cone((\oplus_{i=1}^s\Cone(\partial_s:DR(L_i/\mathbb A^1_k)(V_{D,-1}M)\to DR(L_i/\mathbb A^1_k)(V_{D,0}M))
\xrightarrow{(j_I^*)}\cdots)\to\pi_{-1}(M)) \\ 
\xrightarrow{:=} 
(\oplus_{i=1}^s\Cone(a(M):\Cone(\partial_s:DR(L_i/\mathbb A^1_k)((V_{D,-1}M)^{an})\to DR(L_i/\mathbb A^1_k)((V_{D,0}M)^{an})) \\
\to\Cone(\partial_s:DR(L_i/\mathbb A^1_k)((V_{D,-1}M)^{an})\otimes_{O_S}O_{S^{an}_{\mathbb C}}[\log s] \\
\to DR(L_i/\mathbb A^1_k)((V_{D,-1}M)^{an})\otimes_{O_S}s^{-1}O_{S^{an}_{\mathbb C}}[\log s]))
\xrightarrow{(j_I^*)}\cdots) \\
\xrightarrow{(q_V,A_{-1}(M))} \\
(\oplus_{i=1}^s\Cone((I,0):
\Cone(\partial_s:DR(L_i/\mathbb A^1_k)((\Gr_{V_D,-1}M)^{an})\to DR(L_i/\mathbb A^1_k)((\Gr_{V_D,0}M)^{an})) \\
\to DR(S_i)((\Gr_{V_D,-1}M)^{an})[1])\xrightarrow{j_I^*}\cdots) \\
\xrightarrow{=}
(\oplus_{i=1}^sDR(S_i)((\Gr_{V_D,0}M)^{an})[1]\xrightarrow{j_I^*}\cdots)
\xrightarrow{((j_i^*),0)^{-1}}DR(S)((\Gr_{V_D,0}M)^{an})[1],  
\end{eqnarray*}
with $a(M)=(a_1(M),a_2(M))$, $a_1(M)(m)=m\otimes1$, $a_2(M)(m)=sm\otimes s^{-1}$,
\item for $-1\leq\alpha<0$, the morphism in $D_{c,k}(S_{\mathbb C}^{an})$
\begin{eqnarray*}
B_{\alpha}(M):\pi_{\alpha}(M):=(\oplus_{i=1}^s
\Cone(\partial_s:V_{D\alpha}DR(L_i/\mathbb A^1_k)(M^{an})\otimes_{O_S}s^{\alpha+1}O_{S^{an}_{\mathbb C}}[\log s] \\
\to V_{D\alpha}DR(L_i/\mathbb A^1_k)(M^{an})\otimes_{O_S}s^{\alpha}O_{S^{an}_{\mathbb C}}[\log s])
\xrightarrow{(j_I^*)}\cdots) \\
\to (\oplus_{i=1}^sDR(p^*O_{\mathbb A_1^k})(i^*(l\circ\pi)_*(l\circ\pi)^*DR(L_i/\mathbb A^1_k)(M^{an}))
\xrightarrow{(j_I^*)}\cdots) \\
\xrightarrow{=}(\oplus_{i=1}^s\psi_DDR(p^*O_{\mathbb A_1^k})(DR(L_i/\mathbb A^1_k)(M^{an}))\xrightarrow{(j_I^*)}\cdots)
\xrightarrow{((j_i^*),0)^{-1}}\psi_DDR(S)(M^{an}), \\
(\sum_jm_j\otimes(\log s)^j,m')\mapsto\sum_j(\log s)^jm_j.
\end{eqnarray*}
and for $\alpha=0$ the morphism in $D_{c,k}(S_{\mathbb C}^{an})$
\begin{eqnarray*}
B_0(M):\pi_0(M):= \\
(\oplus_{i=1}^s\Cone(\ad(\pi^{*mod},\pi_*)(-):
\Cone(\partial_s:DR(L_i/\mathbb A^1_k)((V_{D,-1}M)^{an})\to DR(L_i/\mathbb A^1_k)((V_{D,0}M)^{an})) \\
\to\Cone(\partial_s:DR(L_i/\mathbb A^1_k)((V_{D,-1}M)^{an})\otimes_{O_S}O_{S^{an}_{\mathbb C}}[\log s] \\
\to DR(L_i/\mathbb A^1_k)((V_{D,-1}M)^{an})\otimes_{O_S}s^{-1}O_{S^{an}_{\mathbb C}}[\log s]))
\xrightarrow{(j_I^*)}\cdots) \\
\xrightarrow{(I,B_{-1}(M))}\Cone(i_*i^*DR(S)(M^{an})\to\psi_DDR(S)(M^{an}))=:\phi_DDR(S)(M^{an})
\end{eqnarray*}
\end{itemize}

\begin{thm}\label{phipsi0thm0}
Let $k\subset\mathbb C$ a subfield. Let $S\in\SmVar(k)$.
Let $D=V(s)\subset S$ a divisor with $s\in\Gamma(S,L)$ and $p:L\to S$ a line bundle ($S$ being smooth, $D$ is Cartier),
so that we have the closed embedding $i:S\hookrightarrow L$, $i(x)=(x,s(x))$ 
and $D=i^{-1}(s_0)$, $s_0$ being the zero section. 
Denote by $l:L^o:=L\backslash S\hookrightarrow L$ the open embedding 
which induces the open embedding $l:=l\times_LS:S^o:=S\backslash D\hookrightarrow S$.
Denote by $\pi:\tilde L^o\to L^o$ the universal covering 
which induces the universal covering $\pi:=\pi\times_{L^o}S^o:\tilde S^o\to S^o$. 
For $M\in\PSh_{\mathcal D,c}(S)$ such that the rational Kashiwara-Malgrange $V_D$-filtration on $M$ exists, 
\begin{itemize}
\item for each $-1\leq\alpha\leq 0$, the map in $D_{c,k}(S_{\mathbb C}^{an})$
\begin{eqnarray*}
A_{\alpha}(M):\pi_{\alpha}(M)\to DR(S)((\Gr_{V_D,\alpha}M)^{an})[1]
\end{eqnarray*}
is an isomorphism,
\item the morphisms in $D_{c,k}(S_{\mathbb C}^{an})$
\begin{eqnarray*}
B(M):=\oplus_{-1\leq\alpha<0}B_{\alpha}(M):\oplus_{-1\leq\alpha<0}\pi_{\alpha}(M)\to\psi_DDR(S)(M^{an})
\end{eqnarray*}
and
\begin{eqnarray*}
B'(M):=\oplus_{-1<\alpha<0}can(DR(S)(M^{an}))\circ B_{\alpha}\oplus B_0(M):
\oplus_{-1<\alpha\leq 0}\pi_{\alpha}(M)\to\phi_DDR(S)(M^{an})
\end{eqnarray*}
are isomorphisms.
\end{itemize}
We get, for $M\in\PSh_{\mathcal D,c}(S)$ with quasi-unipotent monodromy
such that the rational Kashiwara-Malgrange $V_D$-filtration on $M$ exists,
\begin{itemize}
\item the isomorphism in $D_{c,k}(S_{\mathbb C}^{an})$
\begin{eqnarray*}
T(\psi_D,DR)(M):=B(M)\circ(\oplus_{-1\leq\alpha<0}A_{\alpha}(M)^{-1})[-1]:
DR(S)(\psi_DM^{an})\xrightarrow{\sim}\psi_DDR(S)(M^{an})[-1] 
\end{eqnarray*}
which satisfy 
$T(\psi_D,DR)(M)\circ DR(S)(\oplus_{-1\leq\alpha<0}(s\partial_s-\alpha))=N\circ T(\psi_D,DR)(M)$ where
\begin{equation*}
N:=\log T_u\in\Hom(\psi_DDR(S)(M^{an}),\psi_DDR(S)(M^{an})), \; 
T=T_uT_s, \, T_u \, \mbox{unipotent}, \, T_s \, \mbox{nilpotent} 
\end{equation*} 
is induced by the monodromy automorphism $T:\tilde S^o\xrightarrow{\sim}\tilde S^o$ 
of the universal covering $\pi:\tilde S^o\to S^o:=S\backslash D$,
\item the isomorphism in $D_{c,k}(S_{\mathbb C}^{an})$
\begin{eqnarray*}
T(\phi_D,DR)(M):=B'(M)\circ (\oplus_{-1<\alpha\leq 0}A_{\alpha}(M)^{-1})[-1]:
DR(S)(\phi_DM^{an})\xrightarrow{\sim}\phi_DDR(S)(M^{an})[-1] 
\end{eqnarray*}
which satisfy
\begin{itemize}
\item $T(\phi_D,DR)(M)\circ DR(S)(can(M))=can(DR(S)(M^{an}))\circ T(\psi_D,DR)(M)$
\item $T(\psi_D,DR)(M)\circ DR(S)(var(M))=var(DR(S)(M^{an}))\circ T(\phi_D,DR)(M)$
\end{itemize}
\end{itemize}
As is the rest of this paper, we denote for short 
$M^{an}:=\an_S^{*mod}M=(\pi_{k/\mathbb C}(S)^{*mod}M)^{an}\in\PSh_{\mathcal D,c}(S^{an}_{\mathbb C})$
with $\an_S:S_{\mathbb C}^{an}\xrightarrow{\an_S}S_{\mathbb C}\xrightarrow{\pi_{k/\mathbb C}(S)}S$.
\end{thm}

\begin{proof}
See \cite{Sa2}.
\end{proof}

\begin{defiprop}\label{phipsi0thm0l}
Let $k$ a field of characteristic zero. Let $S\in\SmVar(k)$.
Let $D=V(s)\subset S$ a divisor with $s\in\Gamma(S,L)$ and $L$ a line bundle ($S$ being smooth, $D$ is Cartier),
so that we have the closed embedding $i:S\hookrightarrow L$, $i(x)=(x,s(x))$ 
and $D=i^{-1}(s_0)$, $s_0$ being the zero section. Denote by $j:S^o:=S\backslash D\hookrightarrow S$ the open embedding.
Then $s$ induces locally on $S$ for the Zariski toplogy a morphism $s:S^o\to\mathbb G_m$.
For $K\in D(S)$, we set (see \cite{AyoubT})
\begin{equation*}
\psi_DK:=e(S)_*i^*e(S)_*R(j\circ\pi)_*\pi^*\mathcal Hom(\mathcal A^{\bullet}_D(S^o),e(S)^*K)
\end{equation*}
together with the monodromy morphism $T:\psi_DK\to\psi_DK$ where 
\begin{equation*}
\pi:\tilde S^o:=\varinjlim_{n\in\mathbb N}\Spec(L[t]/(t^n-s))\to S^o
\end{equation*}
is the universal cover and 
$\mathcal A^{\bullet}_D(S^o):=((S^o\times_{S^o\times S^o}S^o)/\mathbb G_m)\in\Fun(\Delta^{\bullet},\Var(k)^{sm}/S^o)$
the diagram of lattices.
Denote by $l:L^o:=L\backslash S\hookrightarrow L$ the open embedding 
which induces the open embedding $l:=l\times_LS:S^o:=S\backslash D\hookrightarrow S$.
Denote again by $\pi:\tilde L^o\to L^o$ the universal covering 
which induces the universal covering $\pi:=\pi\times_{L^o}S^o:\tilde S^o\to S^o$.
Consider $S=\cup_{i=1}^sS_i$ an open affine cover such that 
$D\cap S_i=V(f_i)\subset S_i$ is given by $f_i\in\Gamma(S_i,O_{S_i})$,
denote $q:L_i:=p^{-1}(S_i)\to\mathbb A^1_k$ the projection and $j_i:S_i\hookrightarrow S$ the open embeddings
and consider $\log s\in\mathcal Hom(\mathcal A^{\bullet}_D(S^o),p^*O_{\mathbb A_k^1})$.
For $M\in\PSh_{\mathcal D,c}(S)$ such that the rational Kashiwara-Malgrange $V_D$-filtration on $M$ exists, 
we have then,
\begin{itemize}
\item for $-1\leq\alpha\leq 0$ the isomorphism in $D_c(S)$
\begin{eqnarray*}
A_{alg,\alpha}(M):\pi^{\alg}_{\alpha}(M):=(\oplus_{i=1}^s
\Cone(\partial_s:DR(L_i/\mathbb A^1_k)(V_{D\alpha}M)\otimes_{O_S}s^{\alpha+1}O_S[\log s] \\
\to DR(L_i/\mathbb A^1_k)(V_{D\alpha}M)\otimes_{O_S}s^{\alpha}O_S[\log s])\xrightarrow{j_I^*}\cdots)[-1] \\
\xrightarrow{\sim}
(\oplus_{i=1}^sDR(S_i)(\Gr_{V_D,\alpha}M)\xrightarrow{j_I^*}\cdots)\xrightarrow{((j_i^*),0)^{-1}}DR(S)(\Gr_{V_D,\alpha}M), \\ 
(\sum_jm_j\otimes(\log s)^j,m')\mapsto [m_0]
\end{eqnarray*}
and for $\alpha=0$ the isomorphism in $D_c(S)$
\begin{eqnarray*}
A_{alg,0}(M):=(q_V,A_{alg,-1}(M)):\pi^{\alg}_0(M)\xrightarrow{:=} \\
\Cone((\oplus_{i=1}^s\Cone(\partial_s:DR(L_i/\mathbb A^1_k)(V_{D,-1}M)\to DR(L_i/\mathbb A^1_k)(V_{D,0}M))
\xrightarrow{(j_I^*)}\cdots)\to\pi^{alg}_{-1}(M)) \\
\xrightarrow{\sim}DR(S)(\Gr_{V_D,0}M),
\end{eqnarray*}
given similarly as in the complex analytic case,
\item for $-1\leq\alpha<0$ the morphism in $D_c(S)$
\begin{eqnarray*}
B_{alg,\alpha}(M):\pi^{\alg}_{\alpha}(M):=(\oplus_{i=1}^s
\Cone(\partial_s:V_{D\alpha}DR(L_i/\mathbb A^1_k)(M)\otimes_{O_S}s^{\alpha+1}O_S[\log s] \\
\to V_{D\alpha}DR(L_i/\mathbb A^1_k)(M)\otimes_{O_S}s^{\alpha}O_S[\log s])\xrightarrow{(j_I^*)}\cdots)[-1] \\
\to(\oplus_{i=1}^s
DR(p^*O_{\mathbb A_k^1})(i^*l_*\pi_*\pi^{*mod}\mathcal Hom(\mathcal A_D^{\bullet}(S^o),l^*DR(L_i/\mathbb A^1_k)(M)))
\xrightarrow{(j_I^*)}\cdots)[-1] \\
\xrightarrow{=}(\oplus_{i=1}^s\psi_DDR(p^*O_{\mathbb A_k^1})(DR(L_i/\mathbb A^1_k)(M))\xrightarrow{(j_I^*)}\cdots)[-1]  
\xrightarrow{((j_i^*),0)^{-1}}\psi_DDR(S)(M)[-1], \\ 
(\sum_jm_j\otimes(\log s)^j,m')\mapsto\sum_j(\log s)^jm_j.
\end{eqnarray*}
and for $\alpha=0$ the morphism in $D_c(S)$
\begin{eqnarray*}
B_{alg,0}(M):=(I,B_{alg,-1}(M)):\pi^{\alg}_0(M)\to\phi_DDR(S)(M).
\end{eqnarray*}
\end{itemize} 
It gives, for $M\in\PSh_{\mathcal D,c}(S)$ with quasi-unipotent monodromy
such that the rational Kashiwara-Malgrange $V_D$-filtration on $M$ exists,
\begin{itemize}
\item the canonical morphism in $D(S)$
\begin{eqnarray*}
T_{alg}(\psi_D,DR)(M):=\oplus_{-1\leq\alpha<0}B_{alg,\alpha}(M)\circ A_{alg,\alpha}(M)^{-1}:
DR(S)(\psi_DM)\to\psi_DDR(S)(M)[-1] 
\end{eqnarray*}
which satisfy $T_{alg}(\psi_D,DR)(M)\circ DR(S)(s\partial_s)=N\circ T_{alg}(\psi_D,DR)(M)$ where 
\begin{equation*}
N:=\log T_u\in\Hom(\psi_DDR(S)(M),\psi_DDR(S)(M)), 
T=T_uT_s, \, T_u \, \mbox{unipotent}, \, T_s \, \mbox{nilpotent}
\end{equation*}
is induced by the the monodromy morphism $T:\psi_DDR(S)(M)\to\psi_DDR(S)(M)$.
\item the canonical morphism in $D(S)$
\begin{eqnarray*}
T_{alg}(\phi_D,DR)(M):= \\
\oplus_{-1<\alpha<0}(can(DR(S)(M))\circ B_{alg,\alpha}(M)\circ A_{alg,\alpha}(M)^{-1})
\oplus (B_{alg,0}(M)\circ A_{alg,0}(M)^{-1}): \\
DR(S)(\phi_DM)\to\phi_DDR(S)(M)[-1] 
\end{eqnarray*}
which satisfy
\begin{itemize}
\item $T_{alg}(\phi_D,DR)(M)\circ DR(S)(can(M))=can(DR(S)(M^{an}))\circ T_{alg}(\psi_D,DR)(M)$
\item $T_{alg}(\psi_D,DR)(M)\circ DR(S)(var(M))=var(DR(S)(M^{an}))\circ T_{alg}(\phi_D,DR)(M)$.
\end{itemize}
\end{itemize}
By definition we have for $k\subset\mathbb C$ a subfield, the following commutaive diagram in $D(S^{an}_{\mathbb C})$
\begin{equation*}
\xymatrix{(DR(S)(\psi_DM))^{an}\ar[rr]^{T^w(\an_S,\otimes)(M)}\ar[d]_{T_{alg}(\psi_D,DR)(M)} & \, & 
DR(S)((\psi_DM)^{an})\ar[d]^{T(\psi_D,DR)(M)} \\
(\psi_DDR(S)(M))^{an}\ar[rr]^{T(\psi_D,an)(M)} & \, & \psi_D(DR(S)(M)^{an})}
\end{equation*}
where $T(\psi_D,an)(M)$ is the isomorphism given in \cite{AyoubB}.
\end{defiprop}

\begin{proof}
The proof that $A_{alg,\alpha}(M)$ is an isomorphism is similar to the proof that 
$A_{\alpha}(M)$ is an isomorphism in theorem \ref{phipsi0thm0}.
\end{proof}

\begin{defi}\label{rhoDR0}
Let $k$ a field of caracteristic $0$. 
Let $S\in\SmVar(k)$ and $D=V(s)\subset S$ a (Cartier) divisor. 
Denote by $j:S^o:=S\backslash D\hookrightarrow S$ the open embedding.
Let $M\in\PSh_{\mathcal D,rh}(S)$ which admits the rational Kashiwara-Malgrange filtration.
There exist a subfield $k_0\subset k$ of finite transcendence degree over $\mathbb Q$, 
$S_0\in\SmVar(k_0)$ and $M_0\in\PSh_{\mathcal D,rh}(S_0)$ such that 
$S=S_0\otimes_{k_0}k$ and $M=\pi_{k_0/k}(S_0)^{*mod}M_0$.
Denote again by $j:S^o:=S\backslash D\hookrightarrow S$ the open embedding.
Consider then an embedding $\sigma:k_0\hookrightarrow\mathbb C$. 
Denote by $\pi:\tilde S_0^{an(\sigma)}\to S_0^{o,an(\sigma)}:=S^{an(\sigma)}_0\backslash D^{an(\sigma)}$ 
the universal covering and by $j:S_0\hookrightarrow S_0$ the open embedding. Denote
\begin{eqnarray*}
A_{\pi}:=\ad(i^*,i_*)(-)\circ\ad(\pi^*,\pi_*)(DR(S^o_0)(M_0^{an}))
\in\Hom_{D(S_{0\mathbb C}^{an})}(j_*DR(S^o_0)(M_0^{an}),\psi_DDR(S^o_0)(M_0^{an})) \\
=\Hom_{D(S_{0\mathbb C}^{an})}(DR(S_0)(j_*M_0^{an}),DR(S_0)(\psi_DM_0^{an}))
\end{eqnarray*}
where the equality follows from the isomorphisms
\begin{itemize}
\item $T(j,DR)(M_0):=T^w(j,\otimes)(M_0)\circ DR(S)(T(j,an)(M_0)):DR(S_0)(j_*M_0^{an})\xrightarrow{\sim}j_*DR(S^o_0)(M_0^{an})$,  
\item $T(\psi_D,DR)(M_0):\psi_DDR(S_0)(M_0^{an})\xrightarrow{\sim}\psi_DDR(S_0)(M_0^{an})$ of theorem \ref{phipsi0thm0}.
\end{itemize}
We have by definition the following commutative diagram
\begin{equation}\label{rhoDRdia}
\xymatrix{\Hom_{D_{\mathcal D,rh}(S_0)}(j_*M_0,\psi_DM_0)\otimes_k\mathbb C
\ar[rr]^{\theta_3:=DR(S_0)^{j_*M_0,\psi_DM_0}}\ar[d]^{DR(S_0)^{j_*M_0,\psi_DM_0}} & \, &
\Hom_{D(S_{0\mathbb C}^{an})}(DR(S_0)(j_*M_0^{an}),DR(S_0)(\psi_DM_0^{an})) \\
\Hom_{D(S_0)}(DR(S_0)(j_*M_0),DR(S_0)(\psi_DM_0))\otimes_k\mathbb C
\ar[rr]^{\Hom(T^w(j,\otimes)(M_0),T_{alg}(\psi_D,DR)(M_0))} & \, &
\Hom_{D(S_0)}(j_*DR(S^o_0)(M_0),\psi_DDR(S^o_0)(M_0))\otimes_k\mathbb C
\ar[u]^{\theta_2:=\mathcal Hom(T^w(an,\otimes)(-),T^w(an,\otimes)(-))\circ{\an_S^*}^{-,-}}}
\end{equation}
where $T_{alg}(\psi_D,DR)(M_0)$ is given in definition-proposition \ref{phipsi0thm0l}.
By theorem \ref{DRKk}, $\theta_3$ is an isomorphism, hence for $m=\theta_2(m_k\otimes 1)$, $\theta_3^{-1}(m)=m_k\otimes 1$
by the diagram (\ref{rhoDRdia}).
In particular for $A_{\pi}=\theta_2(A_{\pi,k}\otimes 1)$ with 
\begin{eqnarray*}
A_{\pi,k}:=\ad(i^*,i_*)(-)\circ n_{\log s}(-)\circ\ad(\pi^*,\pi_*)(DR(S^o_0)(M_0))
\in\Hom_{D(S_0)}(j_*DR(S^o_0)(M_0),\psi_DDR(S^o_0)(M_0))
\end{eqnarray*}
and $n_{\log s}(K)=K\to\mathcal Hom(\mathcal A^{\bullet}_D(S^o),K)$, we get in $\PSh_{\mathcal D,rh}(S_0)$ 
\begin{eqnarray*}
\rho_{DR,D}(M_0):=(DR(S)^{j_*M_0,\psi_DM_0})^{-1}(A_{\pi})=A_{\pi,k}:j_*M_0\to\psi_D(M_0).
\end{eqnarray*}
If $\sigma':k_0\hookrightarrow\mathbb C$ is an other embedding, 
there exists $\theta:\mathbb C\to\mathbb C$ an algebraic automorphism of $\mathbb C$ such that $\theta\circ\sigma=\sigma'$
and $\sigma'(\rho_{DR,D}(M_0))=\theta(\sigma(\rho_{DR,D}(M_0)))$ by the diagram (\ref{rhoDRdia}).
This map gives in particular the map in $\PSh_{\mathcal D,rh}(S)$
\begin{eqnarray*}
\rho_{DR,D}(M):=\rho_{DR,D}(M_0)\otimes_{O_{S_0}}O_S:j_*(M)\to\psi_D(M).
\end{eqnarray*}
For $-1\leq\alpha<0$, denote by 
\begin{equation*}
p_{\alpha}:\psi_D(M):=\oplus_{-1\leq\alpha'<0}\Gr_{V_D,\alpha'}M\to\Gr_{V_D,\alpha}M 
\end{equation*}
the projection. 
We get the map in $\PSh_{\mathcal D,rh}(S)$
\begin{eqnarray*}
\rho^u_{DR,D}(M):j_*(M)\xrightarrow{\rho_{DR,D}(M)}\psi_D(M)\xrightarrow{p_{-1}}\psi^u_D(M):=\Gr_{V_D,-1}M.
\end{eqnarray*}
Consider the decomposition 
\begin{equation*}
\psi_D(DR(S)(M^{an}))=\psi^u_D(DR(S)(M^{an}))\oplus\psi^t_D(DR(S)(M^{an}))
\end{equation*}
and its projections $p_u$ and $p_t$.
Since $p^t\circ A_{\pi}=0$ we note that we have $p_{\alpha}\circ\rho_{DR,D}(M)=0$ for $\alpha\neq -1$.
\end{defi}

We now show, using the complex case,
the existence of the rational Kashiwara-Malgrange filtration in the regular holonomic case.

\begin{lem}\label{VQjOlem}
Let $S\in\SmVar(k)$ a smooth affine variety with a closed embedding $l:S\hookrightarrow\mathbb A^N_k$.
Let $D=V(f)\subset S$ a (Cartier) Divisor which is given by a $f\in\Gamma(S,O_S)$. 
Then $D=\tilde D\cap S$ with $\tilde D=V(\tilde f)\subset\mathbb A^N_k$, 
where the polynomial $\tilde f\in\Gamma(\mathbb A^N_k,O_{\mathbb A^N_k})$ is a lift of $f$ ;
denote by $j:S\backslash D\hookrightarrow S$ and
$\tilde j:\mathbb A^N_k\backslash\tilde D\hookrightarrow\mathbb A^N_k$ the open embeddings.
We then have the graph embedding $i:S\hookrightarrow S\times\mathbb A^1_k$, $i(x)=(x,f(x))$ 
and the zero section embedding $i_0:S\hookrightarrow S\times\mathbb A^1_k$, $i_0(x)=(x,0)$. 
Denote $(x,t)\in S\times\mathbb A^1_k$ the coordinates and $s=t\partial_t$.
\begin{itemize}
\item[(i)] Let $M\in\PSh_{\mathcal D}(S)$ such that the multiplication map $m_f:M\to M$, $m_f(m)=fm$, is an isomorphism.
Denote by $\delta=1/(f-t)\in O_S(*D)$. 
Consider for $i\in\mathbb N$ the polynomials $Q_i=\pi_{j=0}^{i-1}(x+j)\in\mathbb Z[x]$ for $i>0$ and $Q_0=1$. 
We have then an isomorphism of $D_S<t,t^{-1},s>$ modules
\begin{eqnarray*}
A_f(M):M[s]f^s:=i_*(M\otimes_{O_S}O_S[s])\xrightarrow{\sim}
i_{*mod}M=i_*M\otimes_k k[\partial_t]=i_*M\otimes_{O_{S}}i_{*mod}O_S(*D), \\
ms^jf^s\mapsto m\otimes(t\partial_t)^j\delta
\end{eqnarray*}
whose inverse is
\begin{eqnarray*}
B_f(M):i_{*mod}M\xrightarrow{\sim}M[s]f^s, \; m\otimes\partial_t^j\delta\mapsto m/f^jQ_j(-s)f^s
\end{eqnarray*}
where the structure of $D_S<t,t^{-1},s>$ module on $M[s]f^s$ is given by
\begin{itemize}
\item $s.(ms^jf^s)=ms^{j+1}f^s$, $t.(ms^jf^s)=m(s+1)^jf^s$,
\item $P.(ms^jf^s)=P(f)/fms^{j+1}f^s+P(m)(s+1)^jf^s$, for $P\in\Gamma(S,D_S)$.
\end{itemize}
\item[(ii)] Consider by (i) the $D_S<t,t^{-1},s>$ submodule
\begin{eqnarray*}
N_f:=D_S[s]f^s\subset O_S(*D)[s]f^s\xrightarrow{A_f(O_S(*D))\sim}i_{*mod}O_S(*D).
\end{eqnarray*}
Then $N_f/tN_f$ an holonomic $D_S$ module.
\item[(iii)] The endomorphism
\begin{equation*}
s:N_f/tN_f\to N_f/tN_f
\end{equation*}
has a minimal polynomial which is equal to the Berstein-Sato Polynomial $b_f\in\mathbb Q[x]$ of $f$.
\end{itemize}
\end{lem}

\begin{proof}
\noindent(i): See \cite{Popa}.

\noindent(ii): By proposition \ref{compDmodh}, $O_S(*D)=j_*O_{S\backslash D}$ is an holonomic $D_S$ module.
Hence, since $N_f/tN_f\subset O_S(*D)$ is a $D_S$ submodule, $N_f/tN_f$ is an holonomic $D_S$ module
by proposition \ref{holksub}.

\noindent(iii): Follows from (ii) by \cite{Coutino} theorem 3.3.

\end{proof}

For an arbitrary field of caracteristic zero, we have the following key proposition :

\begin{prop}\label{VQjO}
Let $k$ a field of characteristic zero.
\begin{itemize}
\item[(i1)]Let $S\in\SmVar(k)$. 
Let $D=V(s)\subset S$ a (Cartier) divisor, where $s\in\Gamma(S,L)$ 
is a section of the line bundle $L=L_D$ associated to $D$. 
We then have the graph embedding $i:S\hookrightarrow L$, $i(x)=(x,s(x))$ 
and the zero section embedding $i_0:S\hookrightarrow L$, $i_0(x)=(x,0)$ and $L_0=i_0(S)$. 
Denote $j:S^o:=S\backslash D\hookrightarrow S$ and $j:L^o:=L\backslash L_0\hookrightarrow L$ 
the open complementary subsets.
Then $j_*O_{S^o}=O_S(*D)\in\PSh_{\mathcal D,rh}(S)$ admits the rational Kashiwara-Malgrange $V_D$-filtration, that is 
$j_*i_{*mod}O_{S^o}=i_{*mod}O_S(*D)\in\PSh_{\mathcal D,rh}(L)$ admits the rational Kashiwara-Malgrange $V_S$-filtration.
\item[(i2)]Let $S\in\SmVar(k)$. Let $D\subset S$ a (Cartier) divisor. 
Denote $j:S^o:=S\backslash D\hookrightarrow S$ the open complementary subset.
Then for $E\in\Vect_{\mathcal D}(S^o)$ regular in the strong sense (see definition \ref{regholk}), 
$j_*E\in\PSh_{\mathcal D,rh}(S)$ admits the rational Kashiwara-Malgrange $V_E$-filtration 
for all (Cartier) divisor $E=V(s')\subset S$.
\item[(ii)] Let $f:X\to S$ a proper morphism with $X,S\in\SmVar(k)$. Let $D\subset S$ a (Cartier) divisor. 
If $M\in\PSh_{\mathcal D,rh}(X)$ admits the rational Kashiwara-Malgrange $V_{f^{-1}(D)}$ filtration, 
then $H^n\int_fM$ admits the rational Kashiwara-Malgrange $V_D$ filtration for all $n\in\mathbb Z$.
\end{itemize}
\end{prop}

\begin{proof}

\noindent(i1): Let $S=\cup_{i=1}^rS_i$ an open affine covering such that for each $i$
$D\cap S_i=V(f_i)\subset S_i$ is given by one equation $f_i\in\Gamma(S_i,O_{S_i})$. 
For each $i$, lemma \ref{VQjOlem} (iii) applied to $S_i$ and $D\cap S_i$, 
gives the rational Kashiwara-Malgrange $V_{S_i}$-filtration for $i_{*mod}O_{S_i}(*D)$.
By unicity of the rational Kashiwara-Malgrange $V$-filtration, we get the rational Kashiwara-Malgrange $V_S$-filtration
on $i_{*mod}O_S(*D)$ since for all $k\in\mathbb Q$ 
\begin{equation*}
V_{S_i,k}i_{*mod}O_S(*D)_{|S_i\cap S_j}=V_{S_j,k}i_{*mod}O_S(*D)_{|S_i\cap S_j}.
\end{equation*}

\noindent(i2): Let $E\subset S$ a (Cartier) divisor.
Consider a subfield $k_0\subset k$ of finite transcendence degree over $\mathbb Q$ such that
$S=S_{0k}:=S\otimes_{k_0}k$ with $S_0\in\SmVar(k_0)$, $D=D_{0k}:=D\otimes_{k_0}k$ with $D_0\subset S_0$, 
$E=E_{0k}:=E\otimes_{k_0}k$ with $E_0\subset S_0$,
$E=\pi_{k_0/k}(S^o_0)^{*mod}E_0$ with $E\in\Vect_{\mathcal D}(S_0^o)$, $S_0^o:=S_0\backslash D_0$,
and an embedding $\sigma:k_0\hookrightarrow\mathbb C$. 
For simplicity, we denote again $S=S_0$, $E=E_0$, $D=D_0$,
$j:S^o:=S\backslash D\hookrightarrow S$ the open complementary subset and $E=E_0\in\Vect_{\mathcal D}(S^o)$,
and denote for short 
\begin{equation*}
\pi:=\pi_{k_0/\mathbb C}(S):S_{\mathbb C}:=S\otimes_{k_0}\mathbb C\to S
\end{equation*}
the projection. The injective map
\begin{eqnarray*}
n_{O_S/O_{S_{\mathbb C}}}(E):E=\pi^*E\to \pi^{*mod}E, \; m\mapsto n_{O_S/O_{S_{\mathbb C}}}(E)(m):=m\otimes 1
\end{eqnarray*}
induces a canonical embedding
\begin{eqnarray*}
j_*n_{O_S/O_{S_{\mathbb C}}}(E):j_*\pi^*E=\pi^*j_*E\hookrightarrow j_*(\pi^{*mod}E), \; m\mapsto m\otimes 1
\end{eqnarray*}
By the complex case, see \cite{Kashiwara}, $j_*(\pi^{*mod}E)\in\PSh_{\mathcal D,rh}(S_{\mathbb C})$ 
admits the rational Kashiwara-Malgrange $V_{E'}$-filtration for all divisor $E'\subset S_{\mathbb C}$.
We then set for $\alpha\in\mathbb Q$
\begin{eqnarray*}
V_{E,\alpha}j_*E:=\pi_*(V_{E_{\mathbb C},\alpha}j_*(\pi^{*mod}E)\cap\pi^*j_*E)
\end{eqnarray*}
so that we get a strict monomorphism 
\begin{equation*}
j_*n_{O_S/O_{S_{\mathbb C}}}(E):(j_*E,V_E)\hookrightarrow\pi_*(j_*(\pi^{*mod}E),V_{E_{\mathbb C}})
\end{equation*}
and so that this filtration satisfies the property of the Kashiwara-Malgrange $V_E$-filtration :
since $E\in\Vect_{\mathcal D}(S^o)$ is a locally free $O_{S^o}$ module, $s:j_*E\to j_*E$ is an isomorphism
with inverse $1/s:j_*E\to j_*E$.
Taking back the initial notations, this means that
we get the $V_{E_0}$ filtration on $j_*E_0\subset\PSh_{\mathcal D,h}(S_0)$. Then,
\begin{eqnarray*}
V_{E,\alpha}j_*E:=\pi_{k_0/k}(S)^{*mod}V_{E_0,\alpha}j_*E_0\subset\pi_{k_0/k}(S)^{*mod}j_*E_0=j_*E
\end{eqnarray*}
gives the $V_E$ filtration on $j_*E$.

\noindent(ii): By definition $H^n\int_fM=H^nf_*E(D_{X\leftarrow S}\otimes_{D_X}M)$.
We then see immediately that $H^nf_*E(D_{X\leftarrow S}\otimes_{D_X}(M,V_{f^{-1}(D)}))$ 
satisfy the hypothesis of the $V_D$ filtration : since $f$ is proper, the $V_{D,0}O_S$ modules 
\begin{eqnarray*}
V_{D,\alpha}H^nf_*E(((D_{X\leftarrow S},V_{f^{-1}(D)})\otimes_{D_X}(M,V_{f^{-1}(D)})):= \\
\Im(H^nf_*(I\otimes\iota_{V_{f^{-1}(D)}}(M)):
H^nf_*E(V_{f^{-1}(D),\alpha}(D_{X\leftarrow S}\otimes_{D_X}M)\to H^nf_*E(D_{X\leftarrow S}\otimes_{D_X}M)) 
\end{eqnarray*}
are coherent.
\end{proof}

\begin{thm}\label{HSk}
Let $k$ a field of characteristic zero. Let $S\in\SmVar(k)$. 
Every $M\in\PSh_{\mathcal D,rh}(S)$ holonomic, regular in the strong sense (see definition \ref{regholk}),
admits the rationnal Kashiwara-Malgrange $V_E$ filtration for each (Cartier) divisor $E\subset S$
(see definition \ref{VfilKM}).
\end{thm}

\begin{proof}
We may assume without loss of generality that $S$ is connected.
We argue by induction on $\dim\supp(M)$. If $\dim\supp(M)=0$, there is noting to prove.
Denote $i:Z:=\supp M\hookrightarrow S$ the closed embedding.
There exist by proposition \ref{holkgen} an open subset $j:S^o\hookrightarrow S$ 
with $D:=S\backslash S^o\subset S$ a (Cartier) divisor 
such that $Z^o:=Z\cap S^o$ is smooth and $i^{*mod}j^*M\in\Vect_{\mathcal D}(Z^o)$ is an integral connexion. 
Then $\dim(D\cap Z_i)=\dim(Z_i)-1$ where $Z_i\subset Z$ are the irreducible component of $Z$
since an holonomic $D_S$ module is generically an integral connexion on its support by proposition \ref{holkgen}.
Take a desingularization $\epsilon:\tilde Z\to Z$ of the pair $(Z,D\cap Z)$
and denote by $l:Z^o\hookrightarrow\tilde Z$ the open embedding.
By proposition \ref{VQjO} (i0) and (i2), $l_*i^{*mod}j^*M\in\PSh_{\mathcal D,rh}(\tilde Z)$ 
admits the rational Kashiwara Malgrange $V_{\epsilon^{-1}i^{-1}(E)}$ filtration, 
hence by proposition \ref{VQjO} (ii)
\begin{equation*}
j_*j^*M=i_{*mod}\epsilon_{*mod}l_*i^{*mod}j^*M\in\PSh_{\mathcal D,rh,Z}(S) 
\end{equation*}
admits the rational Kashiwara Malgrange $V_E$ filtration for all divisor $E\subset S$. We then consider
\begin{itemize}
\item the nearby cycle functor for $M$ 
\begin{equation*}
\psi_D(M):=\psi_D(j_*j^*M):=\oplus_{-1\leq\alpha<0}\Gr_{V_D,\alpha}j_*j^*M\in\PSh_{\mathcal D,rh,D\cap Z}(S),
\end{equation*}
and we have then, by theorem \ref{phipsi0thm0}, for any embedding $k\subset\mathbb C$, the canonical isomorphism 
\begin{equation*}
T(\psi_D,DR)(j_*j^*M):DR(S)(\psi_DM^{an})\xrightarrow{\sim}\psi_DDR(S)(M^{an})
\end{equation*}
(note that $\psi_DK=\psi_D(j_*j^*K)$ for $K\in P(S_{\mathbb C}^{an})$),
\item the unipotent nearby cycle functor for $M$ 
\begin{equation*}
\psi^u_D(M):=\psi^u_D(j_*j^*M):=\Gr_{V_D,-1}j_*j^*M\in\PSh_{\mathcal D,rh,D\cap Z}(S),
\end{equation*}
and we have then, by theorem \ref{phipsi0thm0}, for any embedding $k\subset\mathbb C$, the canonical isomorphism 
\begin{equation*}
T(\psi_D,DR)(j_*j^*M):DR(S)(\psi^u_DM^{an})\xrightarrow{\sim}\psi^u_DDR(S)(M^{an})
\end{equation*}
\item the vanishing cycle functor for $M$
\begin{equation*}
\phi^{\rho}_D(M):=H^0\Cone(\theta_{DR,D}(M):\Gamma^{\vee,h}_DM\to\psi_D(M))\in\PSh_{\mathcal D,rh,D\cap Z}(S)
\end{equation*}
which is regular holonomic by proposition \ref{holregksub},
where $\theta_{DR,D}$ is the factorization in $D_{\mathcal D,rh}(S)$  
\begin{eqnarray*}
\rho_{DR,D}(j_*j^*M)\circ\ad(j^*,j_*)(M)):M\xrightarrow{\gamma^{\vee,h}(M)} 
\Gamma^{\vee,h}_DM:=L\mathbb D_SR\Gamma_DL\mathbb D_SM\xrightarrow{\theta_{DR,D}(M)}\psi_D(M)
\end{eqnarray*}
of the map in $\PSh_{\mathcal D,rh}(S)$ of definition \ref{rhoDR0}
\begin{eqnarray*}
\rho_{DR,D}(j_*j^*M)\circ\ad(j^*,j_*)(M)):=M\xrightarrow{\ad(j^*,j_*)(M)}j_*j^*M \\
\xrightarrow{\rho_{DR,D}(j_*j^*M)}\psi_D(j_*j^*M)=:\psi_D(M),
\end{eqnarray*}
we have then by theorem \ref{DRKk}, for any embedding $k\subset\mathbb C$, the canonical isomorphism 
\begin{equation*}
(T(\gamma_D^{\vee},DR)(M),T(\psi_D,DR)(j_*j^*M)):DR(S)(\phi^{\rho}_DM^{an})\xrightarrow{\sim}\phi_DDR(S)(M^{an})
\end{equation*}
\item the unipotent vanishing cycle functor for $M$
\begin{equation*}
\phi^{\rho,u}_D(M):=H^0\Cone(p_{-1}\circ\theta_{DR,D}(M):\Gamma^{\vee,h}_DM\to\psi^u_D(M))
\in\PSh_{\mathcal D,rh,D\cap Z}(S)
\end{equation*}
we have then by theorem \ref{DRKk}, for any embedding $k\subset\mathbb C$, the canonical isomorphism 
\begin{equation*}
(T(\gamma_D^{\vee},DR)(M),T(\psi_D,DR)(j_*j^*M)):DR(S)(\phi^{\rho,u}_DM^{an})\xrightarrow{\sim}\phi^u_DDR(S)(M^{an}),
\end{equation*}
\item the canonical map $can^{\rho}(M):=H^0c(\phi^{\rho,u}_D(M)):\psi^u_D(M)\to\phi^{\rho,u}_DM$ 
in $\PSh_{\mathcal D,rh,D\cap Z}(S)$,
\item the variation map $var^{\rho}(M):=(0,exp(s\partial s+1)):\phi^{\rho,u}_DM\to\psi^u_DM$ 
in $\PSh_{\mathcal D,rh,D\cap Z}(S)$.
\end{itemize}
Consider then the following canonical map in $C_{\mathcal D,rh}(S)$
\begin{eqnarray*}
Is(M):=(0,(\ad(j^*,j_*)(M),\rho^u_{DR,D}(j_*j^*M)\circ\ad(j^*,j_*)(M)),0): \\
M\to(\psi^u_DM\xrightarrow{(c(x_{S^o/S}(M)),can^{\rho}(M))}x_{S^o/S}(M)\oplus\phi^{\rho,u}_DM
\xrightarrow{(0,exp(s\partial s+1),var^{\rho}(M):=(0,s\partial s))}\psi^u_DM).
\end{eqnarray*}
with 
\begin{equation*}
x_{S^o/S}(M):=\Cone(p_{-1}\circ\rho_{DR,D}(j_*j^*M):j_*j^*M\to\psi^u_DK)\in C_{\mathcal D,rh}(S)
\end{equation*}
which is a quasi-isomorphism by theorem \ref{PSk0} and theorem \ref{DRKk}.
By induction hypothesis, $\phi^{\rho,u}_DM\in\PSh_{\mathcal D,rh,D\cap Z}(S)$ 
admits the Kashiwara-Malgrange rational $V_E$-filtration for all divisor $E\subset S$.
Let $E\subset S$ a Cartier divisor. We then set for $\alpha\in\mathbb Q$
\begin{eqnarray*}
V_{E,\alpha}M:=Is(M)^{-1}(V_{E,\alpha}H^1((\psi^u_DM,V_E)\xrightarrow{(c(x_{S^o/S}(M)),can^{\rho}(M))} \\
\Cone(\rho^u_{DR,D}(j_*j^*M):(j_*j^*M,V_E)\to(\psi^u_DK,V_E))\oplus(\phi^{\rho,u}_DM,V_E) \\
\xrightarrow{((0,exp(s\partial s+1)),var^{\rho}(M):=(0,exp(s\partial s+1)))}(\psi^u_DM,V_E))).
\end{eqnarray*} 
\end{proof}

\begin{thm}\label{phipsi0thm}
Let $k\subset\mathbb C$ a subfield. Let $S\in\SmVar(k)$.
Let $D=V(s)\subset S$ a divisor with $s\in\Gamma(S,L)$ and $L$ a line bundle ($S$ being smooth, $D$ is Cartier).
For $M\in\PSh_{\mathcal D,rh}(S)$ with quasi-unipotent monodromy, 
so that the rational Kashiwara-Malgrange $V_D$-filtration on $M$ exists by theorem \ref{HSk},
\begin{itemize}
\item we have the canonical isomorphism in $D_{c,k}(S_{\mathbb C}^{an})$ given in theorem \ref{phipsi0thm0} :
\begin{equation*}
T(\psi_D,DR)(M):=\oplus_{-1\leq\alpha<0}B_{\alpha}(M)\circ A_{\alpha}(M)^{-1}:
DR(S)(\psi_DM^{an})\xrightarrow{\sim}\psi_DDR(S)(M^{an})[-1] 
\end{equation*}
so that $T(\psi_D,DR)(M)\circ DR(S)(s\partial_s)=N\circ T(\psi_D,DR)(M)$ where
\begin{equation*}
N:=\log T_u\in\Hom(\psi_DDR(S)(M^{an}),\psi_DDR(S)(M^{an})), \; 
T=T_uT_s, \, T_u \, \mbox{unipotent}, \, T_s \, \mbox{nilpotent} 
\end{equation*} 
is induced by the monodromy automorphism $T:\tilde S^o\xrightarrow{\sim}\tilde S^o$ 
of the universal covering $\pi:\tilde S^o\to S^o:=S\backslash D$.
\item we have the canonical isomorphism in $D_{c,k}(S_{\mathbb C}^{an})$ given in theorem \ref{phipsi0thm0} :
\begin{eqnarray*}
T(\phi_D,DR)(M):=\oplus_{-1<\alpha\leq 0}B_{\alpha}(M)\circ A_{\alpha}(M)^{-1}:
DR(S)(\phi_DM^{an})\to\phi_DDR(S)(M^{an})[-1] 
\end{eqnarray*}
so that
\begin{itemize}
\item $T(\phi_D,DR)(M)\circ DR(S)(can(M))=can(DR(S)(M))\circ T(\psi_D,DR)(M)$
\item $T(\psi_D,DR)(M)\circ DR(S)(var(M))=var(DR(S)(M))\circ T(\phi_D,DR)(M)$.
\end{itemize}
\end{itemize}
\end{thm}

\begin{proof}
Follows from theorem \ref{phipsi0thm0}.
\end{proof}

\begin{thm}\label{HSk2}
Let $S\in\SmVar(k)$.
\begin{itemize}
\item[(i)] Let $M\in\PSh_{\mathcal D,rh}(S)$. 
Let $S^o\subset S$ an open subset such that $D:=S\backslash S^o=V(s)\subset S$ is a (Cartier) divisor.
Denote $i:D\hookrightarrow S$ the closed embedding and $j:S^o\hookrightarrow S$ the open embedding.
Take an embedding $\sigma:k\hookrightarrow\mathbb C$.
We have, using definition \ref{rhoDR0} and theorem \ref{HSk}, 
the canonical quasi-isomorphism in $C_{\mathcal D,rh}(S)$ :
\begin{eqnarray*}
Is(M):=(0,(\ad(j^*,j_*)(M),\rho^u_{DR,D}(M)\circ\ad(j^*,j_*)(M)),0): \\
M\to(\psi^u_DM\xrightarrow{(c(x_{S^o/S}(M)),can(M))}x_{S^o/S}(M)\oplus\phi^u_DM
\xrightarrow{((0,exp(s\partial s+1)),var(M))}\psi^u_DM).
\end{eqnarray*}
with 
\begin{equation*}
x_{S^o/S}(M):=\Cone(\rho^u_{DR,D}(M):j_*j^*M\to\psi^u_DK)\in C_{\mathcal D,rh}(S)
\end{equation*}
\item[(ii)]Let $D=V(s)\subset S$ a Cartier divisor. 
Denote $i:D\hookrightarrow S$ the closed embedding and $j:S^o:=S\backslash D\hookrightarrow S$ the open embedding.
Then the functor
\begin{eqnarray*}
(j^*,\phi^u_D,can,var):\PSh_{\mathcal D,rh}(S)\to\PSh_{\mathcal D,rh}(S^o)\times_J\PSh_{\mathcal D,rh,D}(S)
\end{eqnarray*}
is an equivalence of category whose inverse is
\begin{eqnarray*}
\PSh_{\mathcal D,rh}(S^o)\times_J\PSh_{\mathcal D,rh}(S)\to\PSh_{\mathcal D,rh,D}(S), \\ 
(M',M'',u,v)\mapsto 
H^1((\psi^u_DM')\xrightarrow{(c(x_{S^o/S}(M')),u)}x_{S^o/S}(M')\oplus M''
\xrightarrow{((0,exp(s\partial s+1)),v)}(\psi^u_DM')).
\end{eqnarray*}
\end{itemize}
\end{thm}

\begin{proof}
\noindent(i):By theorem \ref{DRKk}, theorem \ref{phipsi0thm} and definition \ref{rhoDR} we have 
$Is(M)={DR(S)^{-,-}}^{-1}(Is(DR(S)(M)))$. The result then follows from theorem \ref{PSk0}.

\noindent(ii):Follows from (i).
\end{proof}

We now give the $p$-adic version of theorem \ref{phipsi0thm} :

\begin{thm}\label{phipsi0thmp}
Let $k\subset K\subset\mathbb C_p$ a subfield with $p$ a prime number and $K$ a $p$ adic field. 
Let $S\in\SmVar(k)$.
Let $D=V(s)\subset S$ a divisor with $s\in\Gamma(S,L)$ and $L$ a line bundle ($S$ being smooth, $D$ is Cartier).
so that we have the closed embedding $i:S\hookrightarrow L$, $i(x)=(x,s(x))$ 
and $D=i^{-1}(s_0)$, $s_0$ being the zero section.
For $M\in\PSh_{\mathcal D,rh}(S)$ with quasi-unipotent monodromy, 
so that the rational Kashiwara-Malgrange $V_D$-filtration on $M$ exists by theorem \ref{HSk},
\begin{itemize}
\item we have the canonical isomorphism in $D_{\mathbb B_{dr}}(S_K^{an,pet})$
\begin{eqnarray*}
T(\psi_D,DR)^{B_{dr}}(M):=B^{B_{dr}}(M)\circ A^{B_{dr}}(M)^{-1}: \\
DR(S)(\psi_DM^{an}\otimes_{O_{S_K}}O\mathbb B_{dr,S_K})\xrightarrow{\sim}
\psi_DDR(S)(M^{an}\otimes_{O_{S_K}}O\mathbb B_{dr,S_K})[-1] 
\end{eqnarray*}
with, for $S=\cup_{i=1}^sS_i$ an open affine cover such that 
$D\cap S_i=V(f_i)\subset S_i$ is given by $f_i\in\Gamma(S_i,O_{S_i})$,
denoting $q:L_i:=p^{-1}(S_i)\to\mathbb A^1_k$ the projection and $j_i:S_i\hookrightarrow S$ the open embeddings,
\begin{itemize}
\item the isomorphism in $D_{\mathbb B_{dr}}(S_K^{an,pet})$
\begin{eqnarray*}
A^{B_{dr}}(M):(\oplus_{i=1}^s\oplus_{-1\leq\alpha<0}
\Cone(\partial_s:DR(L_i/\mathbb A^1_k)((V_{D\alpha}M)^{an})\otimes_{O_S}s^{\alpha+1}O\mathbb B_{dr,S_K} \\
\to DR(L_i/\mathbb A^1_k)((V_{D\alpha}M)^{an})\otimes_{O_S}s^{\alpha}O\mathbb B_{dr,S_K})\xrightarrow{(j_I^*)}\cdots)[-1] \\ 
\to(\oplus_{i=1}^sDR(S_i)(\psi_DM^{an}\otimes_{O_{S_K}}O\mathbb B_{dr,S_K})\xrightarrow{j_I^*}\cdots)[-1] \\
\xrightarrow{((j_i^*),0)^{-1}}(DR(S)(\psi_DM^{an})\otimes_{O_{S_K}}O\mathbb B_{dr,S_K})[-1], \;  
(\sum_jm_j\otimes(\log s)^j,m')\mapsto [m_0] , 
\end{eqnarray*}
\item and the isomorphism in $D_{\mathbb B_{dr}}(S_K^{an,pet})$
\begin{eqnarray*}
B^{B_{dr}}(M):(\oplus_{i=1}^s\oplus_{-1\leq\alpha<0}
\Cone(\partial_s:V_{D\alpha}DR(L_i/\mathbb A^1_k)(M^{an})\otimes_{O_S}s^{\alpha+1}\otimes_{O_{S_K}}O\mathbb B_{dr,S_K} \\
\to V_{D\alpha}DR(L_i/\mathbb A^1_k)(M^{an})\otimes_{O_S}s^{\alpha}\otimes_{O_{S_K}}O\mathbb B_{dr,S_K})
\xrightarrow{(j_I^*)}\cdots)[-1] \\
\to (\oplus_{i=1}^sDR(p^*O_{\mathbb A_1^k})(i^*\pi_*\pi^{*mod}DR(L_i/\mathbb A^1_k)(M^{an})\otimes_{O_{S_K}}O\mathbb B_{dr,S_K})
\xrightarrow{(j_I^*)}\cdots)[-1] \\
\xrightarrow{=:}(\oplus_{i=1}^s\psi_DDR(p^*O_{\mathbb A_1^k})(DR(L_i/\mathbb A^1_k)(M^{an})\otimes_{O_{S_K}}O\mathbb B_{dr,S_K})
\xrightarrow{(j_I^*)}\cdots)[-1] \\
\xrightarrow{((j_i^*),0)^{-1}}\psi_DDR(S)(M^{an}\otimes_{O_{S_K}}O\mathbb B_{dr,S_K})[-1], \;
(\sum_jm_j\otimes(\log s)^j,m')\mapsto\sum_j(\log s)^jm_j,
\end{eqnarray*}
\end{itemize}
so that $T_{B_{dr}}(\psi_D,DR)(M)\circ DR(S)((s\partial_s)\otimes I)=N\circ T_{B_{dr}}(\psi_D,DR)(M)$ where 
\begin{eqnarray*}
N:=\log T_u\in\Hom(\psi_DDR(S)(M^{an}\otimes_{O_{S_K}}O\mathbb B_{dr,S_K}),
\psi_DDR(S)(M^{an}\otimes_{O_{S_K}}O\mathbb B_{dr,S_K})), \\ 
T=T_uT_s, \, T_u \, \mbox{unipotent}, \, T_s \, \mbox{nilpotent}
\end{eqnarray*}
is induced by the monodromy automorphism $T:\tilde S^o\xrightarrow{\sim}\tilde S^o$ of the
perfectoid universal covering $\pi:\tilde S^o\to S^o:=S\backslash D$ (see \cite{Scholze}).
\item there is a canonical isomorphism in $D_{\mathbb B_{dr}}(S_K^{an,pet})$
\begin{eqnarray*}
T^{B_{dr}}(\phi_D,DR)(M):DR(S)(\phi_DM^{an}\otimes_{O_{S_K}}O\mathbb B_{dr,S_K})
\xrightarrow{DR(S)((0,var(M)\otimes I))} \\
DR(S)({\phi^{\rho}_D}M^{an}\otimes_{O_{S_K}}O\mathbb B_{dr,S_K}) 
\xrightarrow{(I,T^{B_{dr}}(\psi_D,DR)(M))}\phi_DDR(S)(M^{an}\otimes_{O_{S_K}}O\mathbb B_{dr,S_K})[-1]. 
\end{eqnarray*}
where 
\begin{equation*}
\phi^{\rho}_DM:=H^0\Cone(\theta_{DR,D}(M):\Gamma_D^{\vee,h}M\to\psi_DM)\in\PSh_{\mathcal D,rh,D}(S)
\end{equation*}
with $\theta_{DR,D}$ the factorization in $D_{\mathcal D,rh}(S)$  
\begin{eqnarray*}
\rho_{DR,D}(M)\circ\ad(j^*,j_*)(M)):M\xrightarrow{\gamma_D^{\vee,h}(M)} 
\Gamma^{\vee,h}_DM:=L\mathbb D_SR\Gamma_DL\mathbb D_SM\xrightarrow{\theta_{DR,D}(M)}\psi_D(M).
\end{eqnarray*}
of the map given in definition \ref{rhoDR0}. In particular 
\begin{itemize}
\item $T^{B_{dr}}(\phi_D,DR)(M)\circ DR(S)(can(M)\otimes I)=
can(DR(S)(M^{an}\otimes_{O_{S_K}}O\mathbb B_{dr,S_K}))\circ T^{B_{dr}}(\psi_D,DR)(M)$
\item $T^{B_{dr}}(\psi_D,DR)(M)\circ DR(S)(var(M)\otimes I)=
var(DR(S)(M^{an}\otimes_{O_{S_K}}O\mathbb B_{dr,S_K}))\circ T^{B_{dr}}(\phi_D,DR)(M)$.
\end{itemize}
\end{itemize}
\end{thm}

\begin{proof}
Similar to the proof of theorem \ref{phipsi0thm}.
\end{proof}

We now look at the particular case of algebraic integral connexions.
The following follows from the work of Bhatwaderkar and Rao (\cite{BR}).

\begin{thm}\label{BRthm}
Let $R$ a local finite type algebra over a field of characteristic zero $k$.
Let $f\in R$ a non zero divisor and non invertible element. 
If $M\in\Mod(R_f)$ is a projective $R_f:=R[1/f]$ module, 
then there exist a projective (hence free) module $M'\in\Mod(R)$ such that $M'_f:=M'\otimes_R R_f=M$
\end{thm}

\begin{proof}
See \cite{BR}.
\end{proof}

\begin{cor}\label{extendEk}
Let $S\in\SmVar(k)$. 
Let $S^o\subset S$ an open subset such that $D:=S\backslash S^o=V(s)\subset S$ is a Cartier divisor.
Denote $j:S^o\hookrightarrow S$ the open embedding.
Let $E=(E,\nabla)\in\Vect_{\mathcal D}(S^o)$ an integrable connexion regular along $D$ (see definition \ref{regholk}). 
\begin{itemize}
\item[(i)]For all $s\in D$, there exist an open affine neighborhood $W_s\subset S$ of $s$ in $S$ 
and a free $O_{W_s}$-submodule $E'\subset (j_*E)_{|W_s}$ such that 
$E'_{|W_s\cap S^o}=E_{|W_s\cap S^o}$, $D_SE'=(j_*E)_{|W_s}$, $s\nabla_sE'\subset E'$
and such that all the eigeinvalue of the residue matrix $\nabla_s:E'/sE'\to E'/sE'$ 
are rational numbers $k\in\mathbb Q$, $0\leq k<1$.
\item[(i)']For $s\in D$, such an $E'\subset (j_*E)_{|W_s}$ is unique.
\item[(ii)] There exist a unique locally free $O_S$-submodule $E'\subset j_*E$ such that
$E'_{|S^o}=E$, $D_SE'=j_*E$ and $s\nabla_sE'\subset E'$
and such that all the eigeinvalue of the residue matrix $\nabla_s:E'/sE'\to E'/sE'$ 
are rational numbers $k\in\mathbb Q$, $0\leq k<1$.
\end{itemize}
\end{cor}

\begin{proof}
\noindent(i):Consider an affine neighborhood $W'_s\subset S$ of $s$.
Then, $E_{|W'_s\cap S^o}\in\Vect(W'_s\cap S^o)$ is projective since it is locally free and $W'_s\cap S^o$ is affine.
By the complex case and regularity, 
there exist an integral lattice $L'\subset j_*(E\otimes_k\mathbb C)$ such that $D_SL=E$, $s\partial_sL\subset L$
and such that all the eigeinvalue of the residue matrix $\nabla_s:L'/sL'\to L'/sL'$ 
are rational numbers $k\in\mathbb Q$, $0\leq k<1$.
Then set $L:=L'\cap\subset j_*E$. we have then $D_SL=E$, $s\partial_sL\subset L$
and all the eigeinvalue of the residue matrix $\nabla_s:L/sL\to L/sL$ 
are rational numbers $k\in\mathbb Q$, $0\leq k<1$.
In particular $L_{|S^o}=E$ and $L$ is a coherent $j_*O_S$ module. Denote $R=O_{W'_s,s}$. 
Then by theorem \ref{BRthm} the projective $R_f$ module $L_{s|S^o}\in\Mod(R_f)$  
extend to a free module $\tilde L\in\Mod(R)$, that is $\tilde L\subset j_*j^*L$ with $\tilde L\otimes_RR_f=L_{s|S^o}$. 
Then take a neighborhood $W_s\subset W'_s$ of $s$ in $W'_s$ and $E'\in\Vect(W_s)$, $E'\subset (j_*E)_{|W_s}$,
such that $E'_s=\tilde L$. Then, $E'_{|W_s\cap S^o}=E_{|W_s}$, $D_SE'=(j_*E)_{|W_s}$, $s\nabla_sE'\subset E'$
and all the eigeinvalue of the residue matrix $\nabla_s:E'/sE'\to E'/sE'$ 
are rational numbers $k\in\mathbb Q$, $0\leq k<1$.

\noindent(i)':Follows the unicity of the $V_D$ filtration on $j_*E$.

\noindent(ii):Follows from (i) and (i)'.
\end{proof}

\subsection{The De Rham modules on algebraic varieties over a field of caracteristic zero}

We recall the theorem of \cite{Saito}
\begin{thm}\label{Sa12}
Let $f:X\to S$ a projective morphism with $X,S\in\Var(\mathbb C)$, where projective means that there exist a factorization
$f:X\xrightarrow{l}\mathbb P^N\times S\xrightarrow{p_S}S$ with $l$ a closed embedding and $p_S$ the projection.
Let $S=\cup_{i=1}^s S_i$ an open cover such that there exits closed embeddings 
$i_I:S_i\hookrightarrow\tilde S_i$ with $\tilde S_i\in\SmVar(\mathbb C)$. 
For $I\subset[1,\ldots,s]$, recall that we denote $S_I:=\cap_{i\in I} S_i$ and  $X_I:=f^{-1}(S_I)$.
We have then the following commutative diagram
\begin{equation*}
\xymatrix{X_I\ar[r]^{i_I\circ l_I} & \mathbb P^N\times\tilde S_I\ar[r]^{p_{\tilde S_I}} & \tilde S_I \\
X_J\ar[u]^{j'_{IJ}}\ar[r]^{i_J\circ l_J} & \mathbb P^N\times\tilde S_J\ar[r]^{p_{\tilde S_J}}\ar[u]^{p'_{IJ}} & 
\tilde S_J\ar[u]^{p_{IJ}}}
\end{equation*}
whose right square is cartesian (see section 5). 
\begin{itemize}
\item[(i)] For 
\begin{equation*}
(((M_I,F),u_{IJ}),K,\alpha)\in HM(X), 
\end{equation*}
where $((M_I,F),u_{IJ})\in\PSh_{\mathcal Dfil}(X_I/(\mathbb P^N\times\tilde S_I))$, $K\in P(X^{an})$, we have
\begin{eqnarray*}
H^n(\int^{FDR}_f((M_I,F),u_{IJ}),Rf_*W,f_*(\alpha))\in HM(S)
\end{eqnarray*}
for all $n\in\mathbb Z$, and for all $p\in\mathbb Z$, 
the differentials of $\Gr_F^p\int^{FDR}_f((M_I,F,W),u_{IJ})$ are strict for the the Hodge filtration $F$.
\item[(ii)]Then, for 
\begin{equation*}
(((M_I,F,W),u_{IJ}),(K,W),\alpha)\in D(MHM(X)), 
\end{equation*}
where $((M_I,F,W),u_{IJ})\in C_{\mathcal D(1,0)fil}(X_I/(\mathbb P^N\times\tilde S_I))$, $(K,W)\in C_{fil}(X^{an})$, we have
\begin{eqnarray*}
H^n(\int^{FDR}_f((M_I,F,W),u_{IJ}),Rf_*(K,W),f_*(\alpha))\in MHM(S)
\end{eqnarray*}
for all $n\in\mathbb Z$, and for all $p\in\mathbb Z$, 
the differentials of $\Gr_F^p\int^{FDR}_f((M_I,F,W),u_{IJ})$ are strict for the the Hodge filtration $F$.
\end{itemize}
\end{thm}

\begin{proof}
\noindent(i):See \cite{Sa2}.

\noindent(ii):Follows from (i) using the spectral sequence of the filtered complex 
$(\int^{FDR}_f((M_I,F,W),u_{IJ}),Rf_*(K,W),f_*(\alpha))$ associated to the filtration $W$: see \cite{Saito}.
\end{proof}

\begin{defi}\label{DHdgpsi}
Let $S\in\SmVar(k)$.
Let $D=V(s)\subset S$ a divisor with $s\in\Gamma(S,L)$ and $L$ a line bundle ($S$ being smooth, $D$ is Cartier).
For $(M,F,W)\in\PSh_{\mathcal D(1,0)fil,rh}(S)$, 
consider using theorem \ref{HSk} the rationnal Kashiwara-Malgrange $V_D$-filtration on $M$ 
that is the rationnal Kashiwara-Malgrange $V_S$-filtration on $i_{*mod}M$ and for $\alpha\in\mathbb Q$, 
\begin{equation*}
V_{D,\alpha}M:=i^{\sharp}V_{S,\alpha}i_{*mod}M\subset M=i^{\sharp}i_{*mod}M. 
\end{equation*}
We then define, using definition \ref{DHdgpsi0},
\begin{itemize}
\item the nearby cycle functor
\begin{equation*}
\psi_D(M,F,W):=\oplus_{-1\leq\alpha<0}(\Gr_{V_D,\alpha}M,F,W)\in\PSh_{\mathcal D(1,0)fil,D}(S)
\end{equation*}
where $\Gr_{V_D,\alpha}M$ is endowed with the induced filtration 
$F^p\Gr_{V_D,\alpha}M:=F^pV_{D,\alpha}M/F^pV_{D,<\alpha}M$,
\item the unipotent nearby cycle functor
\begin{equation*}
\psi^u_D(M,F,W):=(\Gr_{V_D,-1}M,F,W)\in\PSh_{\mathcal D(1,0)fil,D}(S),
\end{equation*}
\item the vanishing cycle functor
\begin{equation*}
\phi_D(M,F,W):=\oplus_{-1<\alpha\leq 0}(\Gr_{V_D,\alpha}M,F,W)\in\PSh_{\mathcal D(1,0)fil,D}(S)
\end{equation*}
where $\Gr_{V_D,\alpha}M$ is endowed with the induced filtration 
$F^p\Gr_{V_D,\alpha}M:=F^pV_{D,\alpha}M/F^pV_{D,<\alpha}M$,
\item the unipotent vanishing cycle functor
\begin{equation*}
\phi^u_D(M,F,W):=(\Gr_{V_D,0}M,F,W)\in\PSh_{\mathcal D(1,0)fil,D}(S),
\end{equation*}
\item the canonical maps
\begin{eqnarray*}
can(M,F,W):=(\partial_s,I):\psi_D(M,F,W)\to\phi_D(M,F,W)(-1), \\
var(M,F,W):=(I,s):\phi_D(M,F,W)\to\psi_D(M,F,W).
\end{eqnarray*}
\item the canonical maps
\begin{eqnarray*}
can(M,F,W):=\partial_s:\psi^u_D(M,F,W)\to\phi^u_D(M,F,W)(-1), \\
var(M,F,W):=s:\phi^u_D(M,F,W)\to\psi^u_D(M,F,W).
\end{eqnarray*}

\end{itemize}
\end{defi}

\begin{prop}\label{phipsigmHdgprop}
Let $S\in\SmVar(k)$. Let $D=V(s)\subset S$ a (Cartier divisor). 
Consider a composition of proper morphisms $(f:X=X_r\xrightarrow{f_r}X_{r-1}\xrightarrow{f_1}X_0=S)\in\SmVar(k)$ and 
\begin{equation*}
(M,F)=H^{n_0}\int_{f_1}\cdots H^{n_r}\int_{f_r}(O_X,F_b)\in\PSh_{\mathcal Dfil,rh}(S).
\end{equation*} 
which admits a $V_D$-filtration (see theorem \ref{HSk}). Then,
\begin{eqnarray*}
\psi_D(M,F)=H^{n_0}\int_{f_1}\cdots H^{n_r}\int_{f_r}\psi_{f^{-1}(D)}(O_X,F_b)\in\PSh_{\mathcal Dfil,rh}(S) 
\end{eqnarray*}
and
\begin{eqnarray*}
\phi_D(M,F)=H^{n_0}\int_{f_1}\cdots H^{n_r}\int_{f_r}\phi_{f^{-1}(D)}(O_X,F_b)\in\PSh_{\mathcal Dfil,rh}(S).
\end{eqnarray*}
\end{prop}

\begin{proof}
Immediate from definition.
\end{proof}

\begin{defi}\label{DRMdef}
\begin{itemize}
\item[(i)]Let $S\in\SmVar(k)$. We define, see theorem \ref{compDmodrhthm}(iv), the full subcategories  
\begin{equation*}
PDRM^1(S)\subset PDRM^2(S)\subset PDRM(S)=\cup_{i\in\mathbb N}PDRM^i(S)\subset\PSh_{\mathcal Dfil,rh}(S)
\end{equation*}
consisting of pure De Rham modules inductively. For each $S\in\SmVar(k)$, we define
\begin{eqnarray*}
PDRM^1(S):=<H^n\int_f(O_X,F_b)(d), (f:X\to S)\in\SmVar(k) \; \mbox{proper}, n,d\in\mathbb Z> 
\subset\PSh_{\mathcal Dfil,rh}(S), 
\end{eqnarray*}
the full abelian subcategory, where $<,>$ means generated by and $(-)$ is the shift of the filtration. 
Assume we have defined $PDRM^{k-1}(S)\subset\PSh_{\mathcal Dfil,rh}(S)$ for all $S\in\SmVar(k)$. 
For each $S\in\SmVar(k)$, we define
\begin{eqnarray*}
PDRM^k(S):=<H^n\int_f(M,F)(d), (f:X\to S)\in\SmVar(k) \; \mbox{proper}, \\ 
(M,F)\in PDRM^{k-1}(X), n,d\in\mathbb Z>\subset\PSh_{\mathcal Dfil,rh}(S), 
\end{eqnarray*}
the full abelian subcategory, where $<,>$ means generated by and $(-)$ is the shift of the filtration. 
By proposition \ref{phipsigmHdgprop}, for $D=V(s)\subset S$ a (Cartier) divisor and $(M,F)\in PDRM(S)$, 
we have, using theorem \ref{HSk} or proposition \ref{VQjO}(ii), for all $k\in\mathbb Z$
\begin{equation*}
\Gr^W_k\psi_D(M,F),\Gr^W_k\psi_D(M,F)\in PDRM(S).
\end{equation*}
\item[(i)']Let $S\in\Var(k)$ non smooth.
Take an open cover $S=\cup_iS_i$ such that there are closed embedding $S_I\hookrightarrow\tilde S_I$ with $S_I\in\SmVar(k)$. 
We define as in (i), see theorem \ref{compDmodrhthmsing}(iv), the full subcategories 
\begin{equation*}
PDRM^1(S)\subset PDRM^2(S)\subset PDRM(S)=\cup_{i\in\mathbb N}PDRM^i(S)\subset\PSh^0_{\mathcal Dfil,rh}(S/(\tilde S_I))
\end{equation*}
inductively. For each $S\in\Var(k)$, we define
\begin{eqnarray*}
PDRM^1(S):=<H^n\int_{p_S}(\Gamma^{\vee,Hdg}_X(O_{Y\times\tilde S_I},F_b),x_{IJ})(d), \\
(f:X\to S)\in\Var(k)\; \mbox{proper}, \; X \, \mbox{smooth}, \;
f=p_S\circ i, n,d\in\mathbb Z>\subset\PSh^0_{\mathcal Dfil,rh}(S/(\tilde S_I))
\end{eqnarray*}
the full abelian subcategory, where $<,>$ means generated by and $(-)$ is the shift of the filtration,
$i:X\hookrightarrow Y\times S$ is a closed embedding with $Y\in\PSmVar(k)$, 
\begin{itemize} 
\item for $j:X^o\hookrightarrow X$ an open embedding we set using proposition \ref{VQjO}(i1) 
$j_{!Hdg}(O_{X^o},F_b):=(j_*O_{X^o},F)$
with $F^pj_*O_{X^o}:=\oplus_k\partial^k_sF^{p+k}V_{X\backslash X^o,k}j_*O_{X^o}$, $X\backslash X^o=V(s)\subset X$,
\item we set
\begin{eqnarray*}
(\Gamma^{\vee,Hdg}_X(O_{Y\times\tilde S_I},F_b),x_{IJ})\in C_{\mathcal Dfil,rh}(X/(Y\times\tilde S_I))
\end{eqnarray*}
where
\begin{eqnarray*}
\Gamma^{\vee,Hdg}_X(O_{Y\times\tilde S_I},F_b):=
\Cone(\ad(j_{I!Hdg},j_I^*)(-):j_{I!Hdg}j_I^*(O_{Y\times\tilde S_I\backslash X},F_b)\to (O_{Y\times\tilde S_I},F_b)),
\end{eqnarray*}
together with the maps 
\begin{eqnarray*}
x_{IJ}:\Gamma^{\vee,Hdg}_X(O_{Y\times\tilde S_I},F_b)\xrightarrow{\ad(p_{IJ}^{*mod},p_{IJ*})(-)} \\ 
p_{IJ*}\Gamma^{\vee,Hdg}_{X\times\tilde S_{J\backslash I}}(O_{Y\times\tilde S_J},F_b)\xrightarrow{\ad(j_{IJ!Hdg},j^{IJ*})(-)}
p_{IJ*}\Gamma^{\vee,Hdg}_X(O_{Y\times\tilde S_J},F_b).
\end{eqnarray*}
\end{itemize}
Assume we have defined $PDRM^{k-1}(S)\subset\PSh^0_{\mathcal Dfil,rh}(S/(\tilde S_I))$ for all $S\in\Var(k)$. 
For each $S\in\Var(k)$, we define
\begin{eqnarray*}
PDRM^k(S):=<H^n\int_{p_S}(\Gamma^{Hdg}_X((M_I,F),u_{IJ})(d), 
(f:X\to S)\in\Var(k) \; \mbox{proper}, \; X \, \mbox{smooth}, \, f=p_S\circ i, \\
((M_I,F),u_{IJ})\in PDRM^{k-1}(X), n,d\in\mathbb Z>\subset\PSh^0_{\mathcal Dfil}(S/(\tilde S_I))
\end{eqnarray*}
the full abelian subcategory, where $<,>$ means generated by and $(-)$ is the shift of the filtration,
$i:X\hookrightarrow Y\times S$ is a closed embedding with $Y\in\PSmVar(k)$, $((M_I,F),u_{IJ})\in DRM^{k-1}(X)$ and 
\begin{itemize}
\item for $j:X^o\hookrightarrow X'$ an open embedding with $X'\in\SmVar(k)$ and $(M,F)\in\PSh_{\mathcal Dfil,rh}(X^o)$
we set $j_{!Hdg}(M,F):=(j_*M,F)$ using theorem \ref{HSk} or proposition \ref{VQjO}(ii)
with $F^pj_*M:=\oplus_k\partial^k_sF^{p+k}V_{X'\backslash X^o,k}j_*(M,F)$, $X'\backslash X^o=V(s)\subset X'$,
\item we set
\begin{eqnarray*}
\Gamma^{\vee,Hdg}_X((M_I,F),u_{IJ}):=(\Gamma_X^{\vee,Hdg}(M_I,F),u^q_{IJ})\in C_{\mathcal Dfil}(X/(Y\times\tilde S_I))
\end{eqnarray*}
with $\Gamma^{\vee,Hdg}_X(M_I,F):=\ad(j_{I!Hdg}j_I^*)(-):\Cone(j_{I!Hdg}j_I^*(M_I,F)\to(M_I,F))$,
together with the maps 
\begin{eqnarray*}
u^q_{IJ}:\Gamma^{\vee,Hdg}_X(M_I,F)\xrightarrow{I(p^{*mod}_{IJ},p_{IJ*})(u_{IJ})\circ\ad(p_{IJ}^{*mod},p_{IJ*})(-)} \\ 
p_{IJ*}\Gamma^{\vee,Hdg}_{X\times\tilde S_{J\backslash I}}(M_J,F)\xrightarrow{\ad(j_{IJ!Hdg},j^{IJ*})(-)}
p_{IJ*}\Gamma^{\vee,Hdg}_X(M_I,F).
\end{eqnarray*}
Note that if $S$ is smooth then this definition of $PDRM(S)$ agree with the one given in (i).
\end{itemize}
For $k\subset\mathbb C$ and $S\in\Var(k)$, 
we have by theorem \ref{Sa12} $PDRM(S_{\mathbb C})\subset\pi_S(HM(S_{\mathbb C}))$ 
which are by definition the De Rham factor of geometric pure Hodge modules.
\item[(ii)]Let $S\in\Var(k)$.
Take an open cover $S=\cup_iS_i$ such that there are closed embedding $S_I\hookrightarrow\tilde S_I$ with $S_I\in\SmVar(k)$.
We define using the pure case (i) and (i)' the full subcategory of weak mixed De Rham modules
\begin{eqnarray*}
\widetilde{DRM}(S):=\left\{((M_I,F,W),u_{IJ}), \; \Gr^W_k((M_I,F,W),u_{IJ})\in PDRM(S)\right\}
\subset\PSh_{\mathcal D(1,0)fil,rh}(S/(\tilde S_I))
\end{eqnarray*}
whose object consists of $((M_I,F,W),u_{IJ})\in\PSh_{\mathcal D(1,0)fil,rh}(S/(\tilde S_I))$ such that 
\begin{equation*}
\Gr^W_k((M_I,F,W),u_{IJ}):=((\Gr^W_kM_I,F),\Gr_W^ku_{IJ})\in PDRM(S).
\end{equation*}
For $S\in\SmVar(k)$ and $D=V(s)\subset S$ a (Cartier) divisor, we have for $(M,F,W)\in\widetilde{DRM}(S)$, using theorem \ref{HSk}, 
\begin{equation*}
\psi_D(M,F,W),\phi_D(M,F,W)\in\widetilde{DRM}(S),
\end{equation*}
by the pure case (c.f. (i) and proposition \ref{phipsigmHdgprop}) and the strictness of the $V$-filtration.
\item[(ii)']Let $S\in\Var(k)$.
Take an open cover $S=\cup_iS_i$ such that there are closed embedding $S_I\hookrightarrow\tilde S_I$ with $S_I\in\SmVar(k)$.
The category of (mixed) de Rham modules over $S$ is the full subcategory
\begin{eqnarray*}
DRM(S)\subset\widetilde{DRM}(S)\subset\PSh_{\mathcal D(1,0)fil,rh}(S/(\tilde S_I))
\end{eqnarray*}
whose object consists of $((M_I,F,W),u_{IJ})\in\widetilde{DRM}(S)$ such that for each $I\subset[1,\ldots s]$
and every Cartier divisor $D\subset\tilde S_I$, the three filtrations $F$,$W$ and $V_D$ on $M_I$ are compatible,
the filtration $W$ is admissible for $D$ (that is the relative monodromy filtration $<W,W(N)>$ exists on $M_I$ where
$N$ is the monodromy of $\tilde S_I\backslash D$ on $M_I$), 
and so on inductively on the graded $\psi_DM_I$ until $\supp(\psi_{D_1}\cdots\psi_{D_r}M_I)$ has dimension zero.
\end{itemize}
For $k\subset\mathbb C$ and $S\in\Var(k)$, 
we have by theorem \ref{Sa12} $DRM(S_{\mathbb C})\subset\pi_S(MHM(S_{\mathbb C}))$ 
which are by definition the De Rham factor of geometric mixed Hodge modules.
In particular, for $S\in\Var(k)$, a morphism 
\begin{equation*}
m:((M_I,F,W),u_{IJ})\to((N_I,F,W),u_{IJ}), \;  ((M_I,F,W),u_{IJ}),((N_I,F,W),u_{IJ})\in DRM(S) 
\end{equation*}
is strict for the Hodge filtration $F$.
For $S\in\Var(k)$ we get $D(DRM(S)):=\Ho_{zar}(C(DRM(S)))$ after localization with Zariski local equivalence.
\end{defi}

\begin{rem}
\begin{itemize}
\item[(i)]Let $S\in\SmVar(k)$. By definition,
\begin{eqnarray*}
PDRM(S):=<H^{n_1}\int_{f_1}\cdots H^{n_r}\int_{f_r}(O_X,F_b)(d), \;
(f:X=X_r\xrightarrow{f_r}X_{r-1}\xrightarrow{f_1}X_0=S)\in\SmVar(k), \\
f_i \; \mbox{proper}, 1\leq i\leq r, n_1,\ldots n_r,d\in\mathbb Z> 
\subset\PSh_{\mathcal Dfil,rh}(S).
\end{eqnarray*} 
\item[(ii)]Let $S\in\Var(k)$ non smooth.
Take an open cover $S=\cup_iS_i$ such that there are closed embedding $S_I\hookrightarrow\tilde S_I$ with $S_I\in\SmVar(k)$.
By definition,
\begin{eqnarray*}
PDRM(S):=<H^{n_1}\int_{p_S}(\Gamma_{X_1}^{\vee,Hdg}(\cdots H^{n_r}\int_{p_{X_{r-1}}}
\Gamma_X^{\vee,Hdg}(O_{Y_r\times\tilde X_{r-1}},F_b)))(d), \\
(f:X=X_r\xrightarrow{f_r}X_{r-1}\xrightarrow{f_1}X_0=S)\in\Var(k) \; \mbox{proper},
f_i=p_{S_{i-1}}\circ i_i, 1\leq i\leq r, n_1,\ldots n_r,d\in\mathbb Z> \\
\subset\PSh_{\mathcal Dfil,rh}(S), 
\end{eqnarray*}
where $i_i:X_i\hookrightarrow Y_i\times X_{i-1}$ is a closed embedding with $Y_i\in\PSmVar(k)$. 
Note that if $S$ is smooth then this definition of $PDRM(S)$ agree with the one given in (i).
\end{itemize}
\end{rem}

\begin{itemize}
\item Let $S\in\SmVar(k)$. We consider the canonical embedding
\begin{eqnarray*}
\iota_S:C(DRM(S))\hookrightarrow C_{\mathcal D(1,0)fil}(S)
\end{eqnarray*}
which induces in the derived category the functor 
\begin{eqnarray*}
\iota_S:D(DRM(S))\to D_{\mathcal D(1,0)fil}(S)\to D_{\mathcal D(1,0)fil,\infty}(S)
\end{eqnarray*}
after localization with respect to filtered Zariski local equivalences 
and $\infty$-filtered Zariski local equivalences respectively.
Note that if $m:(M,F,W)\to(N,F,W)$ with $(M,F,W),(N,F,W)\in C(DRM(S))$ is a Zariski local equivalence,
then it is a filltered Zariski local equivalence by strictness.

\item Let $S\in\Var(k)$ non smooth. 
Take an open cover $S=\cup_iS_i$ such that there are closed embedding $S_I\hookrightarrow\tilde S_I$ with $S_I\in\SmVar(k)$. 
We consider the canonical embedding
\begin{eqnarray*}
\iota_S:C(DRM(S))\hookrightarrow C_{\mathcal D(1,0)fil}(S/(\tilde S_I)),
\end{eqnarray*}
which induces in the derived category the functor 
\begin{eqnarray*}
\iota_S:D(DRM(S))\to D_{\mathcal D(1,0)fil}(S/(\tilde S_I))\to D_{\mathcal D(1,0)fil,\infty}(S/(\tilde S_I))
\end{eqnarray*}
after localization with respect to filtered Zariski local equivalences 
and $\infty$-filtered Zariski local equivalences respectively.
Note that if $m:(M,F,W)\to(N,F,W)$ with $(M,F,W),(N,F,W)\in C(DRM(S))$ is a Zariski local equivalence,
then it is a filltered Zariski local equivalence by strictness.
\end{itemize}

\begin{defi}\label{DHdgj}
\begin{itemize}
\item[(i)] Let $S\in\SmVar(k)$.
Let $D=V(s)\subset S$ a divisor with $s\in\Gamma(S,L)$ and $L$ a line bundle ($S$ being smooth, $D$ is Cartier).
Denote by $j:S^o:=S\backslash D\hookrightarrow S$ the open complementary embedding.
Let $(M,F,W)\in DRM(S^o))$. By theorem \ref{HSk}, $M$ admits the Kashiwara-Malgrange rational $V_D$-filtration, 
that is $i_{*mod}M$ admits the Kashiwara-Malgrange rational $V_S$-filtration and $V_{D,k}M:=i^*V_{S,k}i_{*mod}M$.
We then define,
\begin{itemize}
\item the canonical extension 
\begin{eqnarray*}
j_{*Hdg}(M,F,W):=(j_*M,F,W)\in DRM(S), \; \;
F^pj_*M=\sum_{k\in\mathbb N}\partial_s^kF^{p+k}V_{D,<0}j_*M\subset j_*M, \\
W_kj_*M:=W_kj_{*w}(M,W):=<j_*W_kM,W(N)_kM>\subset j_*M 
\end{eqnarray*}
and $(j_*M,W):=j_{*w}(M,W)$ is given by monodromy weight filtration similarly as in the complex case in \cite{Saito},
so that $j^*j_{*Hdg}(M,F,W)=(M,F,W)$ and $DR(S)(j_{*Hdg}(M,F,W))=j_*DR(S^o)(M,W)$,
\item the canonical extension 
\begin{eqnarray*}
j_{!Hdg}(M,F,W):=\mathbb D_S^{Hdg}j_{*Hdg}\mathbb D_S^{Hdg}(M,F,W)\in DRM(S) 
\end{eqnarray*}
so that $j^*j_{!Hdg}(M,F,W)=(M,F,W)$ and $DR(S)(j_{!Hdg}(M,F,W))=j_!DR(S^o)(M,W)$.
\end{itemize}
Moreover for $(M',F,W)\in DRM(S)$,
\begin{itemize}
\item there is a canonical map $\ad(j^*,j_{*Hdg})(M',F,W):(M',F,W)\to j_{*Hdg}j^*(M',F,W)$ in $DRM(S)$,
\item there is a canonical map $\ad(j_{!Hdg},j^*)(M',F,W):j_{!Hdg}j^*(M',F,W)\to(M',F,W)$ in $DRM(S)$.
\end{itemize}
\item[(ii)] Let $S\in\SmVar(k)$.
Let $Z=V(\mathcal I)\subset S$ an arbitrary closed subset, $\mathcal I\subset O_S$ being an ideal subsheaf. 
Taking generators $\mathcal I=(s_1,\ldots,s_r)$, we get $Z=V(s_1,\ldots,s_r)=\cap^r_{i=1}Z_i\subset S$ with 
$Z_i=V(s_i)\subset S$, $s_i\in\Gamma(S,\mathcal L_i)$ and $L_i$ a line bundle. 
Note that $Z$ is an arbitrary closed subset, $d_Z\geq d_X-r$ needing not be a complete intersection. 
Denote by $j:S^o:=S\backslash Z\hookrightarrow S$,
$j_I:S^{o,I}:=\cap_{i\in I}(S\backslash Z_i)=S\backslash(\cup_{i\in I}Z_i)\xrightarrow{j_I^o}S^o\xrightarrow{j} S$ 
the open complementary embeddings, where $I\subset\left\{1,\cdots,r\right\}$. Denote 
\begin{equation*}
\mathcal D(Z/S):=\left\{(Z_i)_{i\in[1,\ldots r]},Z_i\subset S,\cap Z_i=Z\right\},Z'_i\subset Z_i
\end{equation*}
the flag category.
For $(M,F,W)\in DRM(S^o)$, we define by (i) 
\begin{itemize}
\item the (bi)-filtered complex of $D_S$-modules
\begin{equation*}
j_{*Hdg}(M,F,W):=\varinjlim_{\mathcal D(Z/S)}\Tot_{card I=\bullet}(j_{I*}^{Hdg}j_I^{o*}(M,F,W))\in C(DRM(S)), 
\end{equation*}
where the horizontal differential are given by, 
if $I\subset J$, $d_{IJ}:=\ad(j^*_{IJ},j^{Hdg}_{IJ*})(j_I^{o*}(M,F,W))$, 
$j_{IJ}:S^{oJ}\hookrightarrow S^{oI}$ being the open embedding, 
and $d_{IJ}=0$ if $I\notin J$,
\item the (bi)-filtered complex of $D_S$-modules
\begin{eqnarray*}
j_{!Hdg}(M,F,W):=\varprojlim_{\mathcal D(Z/S)}\Tot_{card I=-\bullet}(j_{I!}^{Hdg}j_I^{o*}(M,F,W)) 
=\mathbb D_S^{Hdg}j_{*Hdg}\mathbb D_S^{Hdg}(M,F,W)\in C(DRM(S)),
\end{eqnarray*}
where the horizontal differential are given by, 
if $I\subset J$, $d_{IJ}:=\ad(j^{Hdg}_{IJ!},j^*_{IJ})(j_I^{o*}(M,F,W))$, 
$j_{IJ}:S^{oJ}\hookrightarrow S^{oI}$ being the open embedding, and $d_{IJ}=0$ if $I\notin J$.
\end{itemize}
By definition, we have for $(M,F,W)\in C(DRM(S^o))$, $j^*j_{*Hdg}(M,F,W)=(M,F,W)$ and $j^*j_{!Hdg}(M,F,W)=(M,F,W)$.
For $(M',F,W)\in C(DRM(S))$, there is, by construction,
\begin{itemize}
\item a canonical map $\ad(j^*,j_{*Hdg})(M',F,W):(M',F,W)\to j_{*Hdg}j^*(M',F,W)$ in $C(DRM(S))$,
\item a canonical map $\ad(j_{!Hdg},j^*)(M',F,W):j_{!Hdg}j^*(M',F,W)\to(M',F,W)$ in $C(DRM(S))$.
\end{itemize}
\end{itemize}
\end{defi}

Let $j:S^o\hookrightarrow S$ an open embedding with $S\in\SmVar(k)$. For $(M,F,W)\in C(DRM(S^o))$, 
\begin{itemize}
\item we have the canonical map in $C_{\mathcal D(1,0)fil}(S)$
\begin{equation*}
T(j_{*Hdg},j_*)(M,F,W):=k\circ\ad(j^*,j_*)(j_{*Hdg}(M,F,W)):j_{*Hdg}(M,F,W)\to j_*E(M,F,W),
\end{equation*}
\item we have the canonical map in $C_{\mathcal D(1,0)fil}(S)$
\begin{eqnarray*}
T(j_!,j_{!Hdg})(M,F,W):=\mathbb D^K_SL_D(k\circ\ad(j^*,j_*)(-)): \\
j_!(M,F,W):=\mathbb D_S^KL_Dj_*E(\mathbb D_S^K(M,F,W))\to\mathbb D_S^KL_Dj_{*Hdg}\mathbb D_S^K(M,F,W)=j_{!Hdg}(M,F,W).
\end{eqnarray*}
\end{itemize}

\begin{rem}\label{remHdgkey}
Let $j:S^o\hookrightarrow S$ an open embedding, with $S\in\SmVar(k)$.
Then, for $(M,F,W)\in DRM(S^o)$, 
\begin{itemize}
\item the map  $T(j_!,j_{!Hdg})(M,W):j_{!w}(M,W)\to j_{!Hdg}(M,W)$ in $C_{\mathcal D0fil}(S)$
is a filtered quasi-isomorphism (by the acyclicity of the functor $j_*$ in the divisor case).
\item the map $T(j_{*Hdg},j_*)(M,W):j_{*Hdg}(M,W)\to j_{*w}(M,W)$ in $C_{\mathcal D0fil}(S)$
is a filtered quasi-isomorphism (by the acyclicity of the functor $j_*$ in the divisor case).
\end{itemize}
Hence, for $(M,F,W)\in DRM(S^o)$,  
\begin{itemize}
\item we get, for all $p,n\in\mathbb N$, monomorphisms
\begin{equation*}
F^pH^nT(j_!,j_{!Hdg})(M,F,W):F^pH^nj_{!w}(M,F,W)\hookrightarrow F^pH^nj_{!Hdg}(M,F,W)
\end{equation*}
in $\PSh_{O_S}(S)$, but $F^pH^nj_{!w}(M,F,W)\neq F^pH^nj_{!Hdg}(M,F,W)$ (it leads to different F-filtrations),
since $F^pH^nj_!(M,F)\subset H^nj_!M$ are sub $D_S$ module 
while the F-filtration on $H^nj_{!Hdg}(M,F)$ is given by Kashiwara-Malgrange $V$-filtrations,
hence satisfy a non trivial Griffith transversality property, 
thus $H^nj_!(M,F)$ and $H^nj_{!Hdg}(M,F)$ are isomorphic as $D_S$-modules but NOT isomorphic as filtered $D_S$-modules.
\item we get, for all $p,n\in\mathbb N$, monomorphisms
\begin{equation*}
T(j_{*Hdg},j_*)(M,F,W):F^pH^nj_{*Hdg}(M,F,W)\hookrightarrow F^pH^nj_{*w}(M,F,W)
\end{equation*}
in $\PSh_{O_S}(S)$, but $F^pH^nj_{*Hdg}(M,F,W)\neq F^pH^nj_{*w}(M,F,W)$ (it leads to different F-filtrations),
since $F^pH^nj_*E(M,F)\subset H^nj_*E(M)$ are sub $D_S$ module 
while the F-filtration on $H^nj_{*Hdg}(M,F)$ is given by Kashiwara-Malgrange $V$-filtrations,
hence satisfy a non trivial Griffith transversality property, 
thus $H^nj_*E(M,F)$ and $H^nj_{*Hdg}(M,F)$ are isomorphic as $D_S$-modules but NOT isomorphic as filtered $D_S$-modules.
\end{itemize}
\end{rem}

\begin{prop}\label{jHdgpropad}
\begin{itemize}
\item[(i)] Let $S\in\SmVar(k)$.
Let $D=V(s)\subset S$ a divisor with $s\in\Gamma(S,L)$ and $L$ a line bundle ($S$ being smooth, $D$ is Cartier).
Denote by $j:S^o:=S\backslash D\hookrightarrow S$ the open complementary embedding. Then, 
\begin{itemize}
\item $(j^*,j_{*Hdg}):DRM(S)\leftrightarrows DRM(S^o)$ is a pair of adjoint functors
\item $(j_{!Hdg},j^*):DRM(S^o)\leftrightarrows DRM(S)$ is a pair of adjoint functors.
\end{itemize}
\item[(ii)] Let $S\in\SmVar(k)$.
Let $Z=V(\mathcal I)\subset S$ an arbitrary closed subset, $\mathcal I\subset O_S$ being an ideal subsheaf. 
Denote by $j:S^o:=S\backslash Z\hookrightarrow S$. Then,
\begin{itemize}
\item $(j^*,j_{*Hdg}):D(DRM(S))\leftrightarrows D(DRM(S^o))$ is a pair of adjoint functors
\item $(j_{!Hdg},j^*):D(DRM(S^o))\leftrightarrows D(DRM(S))$ is a pair of adjoint functors.
\end{itemize}
\end{itemize}
\end{prop}

\begin{proof}
\noindent(i): Follows from the fact that for $(M,F)\in DRM(S)$, we have $F^pV_{D,<0}M=j_*F^pj^*M\cap V_{D,<0}M$,
where $V_{D,p}M:=i^{\sharp}V_{S,p}i_{*mod}M$.

\noindent(ii):Follows from (i).
\end{proof}

The map given in definition \ref{rhoDR0} induces the following:

\begin{defi}\label{rhoDR}
Let $k$ a field of caracteristic $0$. 
Let $S\in\SmVar(k)$ and $D\subset S$ a (Cartier) divisor. Let $(M,F,W)\in DRM(S)$.
The map of definition \ref{rhoDR0} given by theorem \ref{DRKk} and theorem \ref{phipsi0thm0} in $\PSh_{\mathcal D,rh}(S)$
\begin{eqnarray*}
\rho_{DR,D}(M):={DR(S)^{-,-}}^{-1}(\ad(i^*,i_*)(-)\circ\ad(\pi^*,\pi_*)(DR(S)(M))):j_*M\to\psi_D(M).
\end{eqnarray*}
induces, similarly to the complex case (\cite{Saito}) by theorem \ref{DRKk}, the unicity of the $V$-filtration,
and the definition of the monomdromy filtration,
the following maps in $\PSh_{\mathcal D(1,0)fil,rh}(S)$
\begin{eqnarray*}
\rho_{DR,D}(M,F,W):=\rho_{DR,D}(M):j_{*Hdg}(M,F,W)\to\psi_D(M,F,W)
\end{eqnarray*}
and
\begin{eqnarray*}
\rho^u_{DR,D}(M,F,W):=\rho^u_{DR,D}(M):j_{*Hdg}(M,F,W)\xrightarrow{\rho_{DR,D}(M,F,W)}\psi_D(M,F,W)
\xrightarrow{p_{-1}}\to\psi^u_D(M,F,W).
\end{eqnarray*}
\end{defi}

\begin{prop}\label{PSkDRM}
\begin{itemize}
\item[(i)] Let $(M,F,W)\in DRM(S)$. 
Let $S^o\subset S$ an open subset such that $M_{|S^o}\in\Vect_{\mathcal D}(S^o)$.
Denote $i:D:=S\backslash D\hookrightarrow S$ the closed embedding and $j:S^o\hookrightarrow S$ the open embedding.
We have the canonical quasi-isomorphism in $C_{\mathcal D,rh}(S)$ given in theorem \ref{HSk2}:
\begin{eqnarray*}
Is(M):=(0,(\ad(j^*,j_*)(M),\rho^u_{DR,D}(M)\circ\ad(j^*,j_*)(M)),0): \\
M\to(\psi^u_DM\xrightarrow{(c(x_{S^o/S}(M)),can(M))}x_{S^o/S}(M)\oplus\phi^u_DM
\xrightarrow{(0,exp(s\partial s-1)),var(M))}\psi^u_DM).
\end{eqnarray*}
gives a filtered quasi-isomorphism in $C_{\mathcal D(1,0)fil,rh}(S)$
\begin{eqnarray*}
Is(M,F,W):=(0,(\ad(j^*,j_{*Hdg})(M,F,W),\rho_{DR,D}(M,F,W)\circ\ad(j^*,j_{*Hdg})(M,F,W)),0): \\
(M,F,W)\to(\psi^u_D(M,F,W)\xrightarrow{(c(x_{S^o/S}(M,F,W)),can(M,F,W))}x_{S^o/S}(M,F,W)\oplus\phi^u_D(M,F,W) \\
\xrightarrow{(0,exp(s\partial s+1)),var(M,F,W))}\psi^u_D(M,F,W)).
\end{eqnarray*}
with, see definition \ref{rhoDR},
\begin{eqnarray*}
x_{S^o/S}(M,F,W):=\Cone(\rho^u_{DR,D}(M,F,W):j_{*Hdg}(M,F,W)\to\psi^u_D(M,F,W))\in C_{\mathcal D(1,0)fil,rh}(S)
\end{eqnarray*}
\item[(ii)]Let $S\in\SmVar(k)$. Let $D=V(s)\subset S$ a (Cartier) divisor, where $s\in\Gamma(S,L)$ 
is a section of the line bundle $L=L_D$ associated to $D$. We then have the zero section embedding
$i:S\hookrightarrow L$. We denote $L_0=i(S)$ and $j:L^o:=L\backslash L_0\hookrightarrow L$ the open complementary subset.
We denote by $DRM(S\backslash D)\times_J DRM(D)$ the category whose set of objects consists of
\begin{equation*}
\left\{(\mathcal M,\mathcal N,a,b),\mathcal M\in DRM(S\backslash D),\mathcal N\in DRM(D),
a:\psi_{D1}\mathcal M\to N,b:N\to\psi_{D1}M \right\}
\end{equation*}
The functor (see definition \ref{DHdgpsi})
\begin{eqnarray*}
(j^*,\phi^u_{D},can,var):DRM(S)\to DRM(S\backslash D)\times_J DRM(D), \\
(M,F,W)\mapsto (j^*(M,F,W),\phi^u_{D}(M,F,W),can(M,F,W),var(M,F,W))
\end{eqnarray*}
is an equivalence of category.
\end{itemize}
\end{prop}

\begin{proof}
\noindent(i):Similar to the complex case (\cite{Saito}) by theorem \ref{DRKk}, the unicity of the $V$-filtration
and the definition of the monodromy weight filtration.

\noindent(ii): follows from (i).
\end{proof}

We make the following key definition

\begin{defi}\label{gammaHdg}
Let $S\in\SmVar(k)$. Let $Z\subset S$ a closed subset.
Denote by $j:S\backslash Z\hookrightarrow S$ the complementary open embedding. 
\begin{itemize}
\item[(i)] We define using definition \ref{DHdgj}, the filtered Hodge support section functor
\begin{eqnarray*}
\Gamma^{Hdg}_Z:C(DRM(S))\to C(DRM(S)), \\ 
(M,F,W)\mapsto\Gamma^{Hdg}_Z(M,F,W):=\Cone(\ad(j^*,j_{*Hdg})(M,F,W):(M,F,W)\to j_{*Hdg}j^*(M,F,W))[-1],
\end{eqnarray*}
together we the canonical map $\gamma^{Hdg}_Z(M,F,W):\Gamma^{Hdg}_Z(M,F,W)\to (M,F,W)$.
\item[(i)'] Since $j_{*Hdg}:C(DRM(S^o))\to C(DRM(S))$ is an exact functor, 
$\Gamma^{Hdg}_Z$ induces the functor
\begin{eqnarray*}
\Gamma^{Hdg}_Z:D(DRM(S))\to D(DRM(S)), \; (M,F,W)\mapsto\Gamma^{Hdg}_Z(M,F,W)
\end{eqnarray*}
\item[(ii)] We define using definition \ref{DHdgj}, the dual filtered Hodge support section functor
\begin{eqnarray*}
\Gamma^{\vee,Hdg}_Z:C(DRM(S))\to C(DRM(S)), \\ 
(M,F,W)\mapsto\Gamma^{\vee,Hdg}_Z(M,F,W):=\Cone(\ad(j_{!Hdg},j^*)(M,F,W):j_{!Hdg},j^*(M,F,W)\to (M,F,W)),
\end{eqnarray*}
together we the canonical map $\gamma^{\vee,Hdg}_Z(M,F,W):(M,F,W)\to\Gamma_Z^{\vee,Hdg}(M,F,W)$.
\item[(ii)'] Since $j_{!Hdg}:C(DRM(S^o))\to C(DRM(S))$ is an exact functor, $\Gamma^{\vee,Hdg}_Z$ induces the functor
\begin{eqnarray*}
\Gamma^{\vee,Hdg}_Z:D(DRM(S))\to D(DRM(S)), \; (M,F,W)\mapsto\Gamma^{\vee,Hdg}_Z(M,F,W)
\end{eqnarray*}
\end{itemize}
\end{defi}

We now give the definition of the filtered Hodge inverse image functor :

\begin{defi}\label{inverseHdg}
\begin{itemize}
\item[(i)] Let $i:Z\hookrightarrow S$ be a closed embedding, with $Z,S\in\SmVar(k)$.
Then, for $(M,F,W)\in C(DRM(S))$, we set 
\begin{equation*}
i_{Hdg}^{*mod}(M,F,W):=i^*\Gr_{V_Z,0}\Gamma_Z^{Hdg}(M,F,W)\in D(DRM(Z))
\end{equation*}
and
\begin{equation*}
i_{Hdg}^{\hat*mod}(M,F,W):=i^*\Gr_{V_Z,0}\Gamma_Z^{\vee,Hdg}(M,F,W)\in D(DRM(Z)),
\end{equation*}
noting that $i_{*mod}:D(DRM(Z))\to D(DRM_Z(S))$ is an equivalence of category
whose inverse is $i^*\Gr_{V_Z,0}:D(DRM_Z(S))\to D(DRM(Z))$.
\item[(ii)] Let $f:X\to S$ be a morphism, with $X,S\in\SmVar(k)$.
Consider the factorization $f:X\xrightarrow{i} X\times S\xrightarrow{p_S}S$, 
where $i$ is the graph embedding and $p_S:X\times S\to S$ is the projection.
\begin{itemize}
\item For $(M,F,W)\in C(DRM(S))$ we set
\begin{equation*}
f_{Hdg}^{*mod}(M,F,W):=i_{Hdg}^{*mod}p_S^{*mod[-]}(M,F,W)(d_X)[2d_X]\in D(DRM(X)), 
\end{equation*}
\item For $(M,F,W)\in C(DRM(S))$ we set
\begin{equation*}
f_{Hdg}^{\hat*mod}(M,F,W):=i_{Hdg}^{\hat*mod}p_S^{*mod[-]}(M,F,W)\in D(DRM(X)), 
\end{equation*}
\end{itemize}
If $j:S^o\hookrightarrow S$ is a closed embedding, we have, for $(M,F,W)\in C(DRM(S))$, 
\begin{equation*}
j_{Hdg}^{*mod}(M,F,W)=j_{Hdg}^{\hat*mod}(M,F,W)=j^*(M,F,W)\in D(DRM(S^o)))
\end{equation*}
\item[(iii)] Let $f:X\to S$ be a morphism, with $X,S\in\SmVar(k)$.
Consider the factorization $f:X\xrightarrow{i} X\times S\xrightarrow{p_S}S$, 
where $i$ is the graph embedding and $p_S:X\times S\to S$ is the projection.
\begin{itemize}
\item For $(M,F,W)\in C(DRM(S))$ we set
\begin{equation*}
f_{Hdg}^{*mod}(M,F,W):=\Gamma_X^{Hdg}p_S^{*mod[-]}(M,F,W)(d_X)[2d_X]\in C(DRM(X\times S)), 
\end{equation*}
\item For $(M,F,W)\in C(DRM(S)))$ we set
\begin{equation*}
f_{Hdg}^{\hat*mod}(M,F,W):=\Gamma_X^{\vee,Hdg}p_S^{*mod[-]}(M,F,W)\in C(DRM(X\times S)), 
\end{equation*}
\end{itemize}
\end{itemize}
\end{defi}

\begin{defiprop}\label{TgammaHdg}
\begin{itemize}
\item[(i)] Let $g:S'\to S$ a morphism with $S',S\in\SmVar(k)$ and $i:Z\hookrightarrow S$ a closed subset. 
Then, for $(M,F,W)\in C(DRM(S)))$, there is a canonical map in $C(DRM_{S'}(S'\times S))$
\begin{equation*}
T^{Hdg}(g,\gamma)(M,F,W):g_{Hdg}^{*mod,\Gamma}\Gamma^{Hdg}_{Z}(M,F,W)\to\Gamma^{Hdg}_{Z\times_S S'}g_{Hdg}^{*mod,\Gamma}(M,F,W)
\end{equation*}
unique up to homotopy such that 
\begin{equation*}
\gamma^{Hdg}_{Z\times_S S'}(g_{Hdg}^{*mod,\Gamma}(M,F,W))\circ T^{Hdg}(g,\gamma)(M,F,W)=
g^{*mod,\Gamma}_{Hdg}\gamma^{Hdg}_{Z}(M,F,W).
\end{equation*}
\item[(i)'] Let $g:S'\to S$ a morphism with $S',S\in\SmVar(k)$ and $i:Z\hookrightarrow S$ a closed subset. 
Then, for $(M,F,W)\in C(DRM(S)))$, 
there is a canonical isomorphism in $C(DRM_{S'}(S'\times S))$
\begin{equation*}
T^{Hdg}(g,\gamma^{\vee})(M,F,W):\Gamma^{Hdg}_{Z\times_S S'}g_{Hdg}^{\hat*mod,\Gamma}(M,F,W)\xrightarrow{\sim} 
g_{Hdg}^{\hat*mod,\Gamma}\Gamma^{Hdg}_{Z}(M,F,W)
\end{equation*}
unique up to homotopy such that 
\begin{equation*}
\gamma^{\vee,Hdg}_{Z\times_S S'}(g_{Hdg}^{\hat*mod,\Gamma}(M,F,W))
\circ g^{\hat*mod,\Gamma}_{Hdg}\gamma^{\vee,Hdg}_{Z}(M,F,W)=T^{Hdg}(g,\gamma)(M,F,W).
\end{equation*}
\item[(ii)] Let $S\in\SmVar(k)$ and $i_1:Z_1\hookrightarrow S$, $i_2:Z_2\hookrightarrow Z_1$ be closed embeddings.
Then, for $(M,F,W)\in C(DRM(S)))$, 
\begin{itemize}
\item there is a canonical map 
$T(Z_2/Z_1,\gamma^{Hdg})(M,F,W):\Gamma^{Hdg}_{Z_2}(M,F,W)\to\Gamma^{Hdg}_{Z_1}(M,F,W)$ 
in $C(DRM(S)))$ unique up to homotopy such that 
\begin{equation*}
\gamma^{Hdg}_{Z_1}(M,F,W)\circ T(Z_2/Z_1,\gamma^{Hdg})(M,F,W)=\gamma^{Hdg}_{Z_2}(M,F,W) 
\end{equation*}
together with a distinguish triangle in $K(DRM(S))$
\begin{eqnarray*}
\Gamma^{Hdg}_{Z_2}(M,F,W)\xrightarrow{T(Z_2/Z_1,\gamma^{Hdg})(M,F,W)}\Gamma^{Hdg}_{Z_1}(M,F,W) \\
\xrightarrow{\ad(j_2^*,j^{Hdg}_{2*})(\Gamma^{Hdg}_{Z_1}(M,F,W))}
\Gamma^{Hdg}_{Z_1/\backslash Z_2}(G,F)\to\Gamma^{Hdg}_{Z_2}(M,F,W)[1]
\end{eqnarray*} 
\item there is a canonical map 
$T(Z_2/Z_1,\gamma^{\vee,Hdg})(M,F,W):\Gamma_{Z_1}^{\vee,Hdg}(M,F,W)\to\Gamma_{Z_2}^{\vee,Hdg}(M,F,W)$ 
in $C(DRM(S)))$ unique up to homotopy such that 
\begin{equation*}
\gamma^{\vee,Hdg}_{Z_2}(M,F,W)=T(Z_2/Z_1,\gamma^{\vee,Hdg})(M,F,W)\circ\gamma^{\vee,Hdg}_{Z_1}(M,F,W). 
\end{equation*}
together with a distinguish triangle in $K(DRM(S))$
\begin{eqnarray*}
\Gamma_{Z_1\backslash Z_2}^{\vee,Hdg}(M,F,W)\xrightarrow{\ad(j^{Hdg}_{2!},j_2^*)(M,F,W)}
\Gamma_{Z_1}^{\vee,Hdg}(M,F,W) \\ 
\xrightarrow{T(Z_2/Z_1,\gamma^{\vee,Hdg})(M,F,W))}
\Gamma^{\vee,Hdg}_{Z_2}(M,F,W)\to\Gamma_{Z_2\backslash Z_1}^{\vee,Hdg}(M,F,W)[1]
\end{eqnarray*} 
\end{itemize}
\end{itemize}
\end{defiprop}

\begin{proof}
Follows from the projection case and the closed embedding case using the adjonction maps.
\end{proof}

The definitions \ref{gammaHdg} and \ref{inverseHdg} immediately extends to the non smooth case :

\begin{defi}\label{gammaHdgsing}
Let $S\in\Var(k)$. Let $Z\subset S$ a closed subset.
Let $S=\cup_iS_i$ an open cover such that there exist closed embeddings
$i_i:S_i\hookrightarrow\tilde S_i$ with $\tilde S_i\in\SmVar(k)$. Denote $Z_I:=Z\cap S_I$. 
Denote by $j:S\backslash Z\hookrightarrow S$ and $\tilde j_I:\tilde S_I\backslash Z_I\hookrightarrow\tilde S_I$ 
the complementary open embeddings. 
\begin{itemize}
\item[(i)] We define using definition \ref{DHdgj}, the filtered Hodge support section functor
\begin{eqnarray*}
\Gamma^{Hdg}_Z:C(DRM(S)))\to C(DRM(S))), \; ((M_I,F,W),u_{IJ})\mapsto\Gamma^{Hdg}_Z((M_I,F,W),u_{IJ}):= \\
\Cone(\ad(j^*,j_{*Hdg})(M_I,F,W),u_{IJ}):=(\ad(\tilde j_I^*,\tilde j_{I*Hdg})(M_I,F,W)): \\
((M_I,F,W),u_{IJ})\to(\tilde j_{I*Hdg}\tilde j_I^*(M_I,F,W),\tilde j_{J*Hdg}u_{IJ}))[-1],
\end{eqnarray*}
together with the canonical map $\gamma^{Hdg}_Z((M_I,F,W),u_{IJ}):\Gamma^{Hdg}_Z((M_I,F,W),u_{IJ})\to((M_I,F,W),u_{IJ})$.
\item[(i)'] Since 
$\tilde j_{I*Hdg}:C(DRM(\tilde S_I\backslash S_I))\to C(DRM(\tilde S_I))$
are exact functors, $\Gamma^{Hdg}_Z$ induces the functor
\begin{eqnarray*}
\Gamma^{Hdg}_Z: D(DRM(S))\to D(DRM(S)),  ((M_I,F,W),u_{IJ})\mapsto\Gamma^{Hdg}_Z((M_I,F,W),u_{IJ})
\end{eqnarray*}
\item[(ii)] We define using definition \ref{DHdgj}, the dual filtered Hodge support section functor
\begin{eqnarray*}
\Gamma^{\vee,Hdg}_Z:C(DRM(S))\to C(DRM(S)), \\ 
((M_I,F,W),u_{IJ})\mapsto\Gamma^{\vee,Hdg}_Z((M_I,F,W),u_{IJ}):=
\mathbb D^{Hdg}_S\Gamma_Z^{Hdg}\mathbb D_S^{Hdg}((M_I,F,W),u_{IJ})= \\
\Cone(\ad(j_{!Hdg},j^*)((M_I,F,W),u_{IJ}):=(\ad(\tilde j_{I!Hdg},\tilde j_I^*)(M_I,F,W)): \\
(\tilde j_{I!Hdg}\tilde j_I^*(M_I,F,W),(\tilde j_{J*Hdg}u^d_{IJ})^d)\to ((M_I,F,W),u_{IJ})),
\end{eqnarray*}
together we the canonical map 
$\gamma^{\vee,Hdg}_Z((M_I,F,W),u_{IJ}):((M_I,F,W),u_{IJ})\to\Gamma_Z^{\vee,Hdg}((M_I,F,W),u_{IJ})$.
\item[(ii)'] Since 
$\tilde j_{I!Hdg}:C(DRM(\tilde S_I\backslash S_I))\to C(DRM(\tilde S_I))$ 
are exact functors, $\Gamma^{Hdg,\vee}_Z$ induces the functor
\begin{eqnarray*}
\Gamma^{\vee,Hdg}_Z:D(DRM(S))\to D(DRM(S)), \; ((M_I,F,W),u_{IJ})\mapsto\Gamma^{\vee,Hdg}_Z((M_I,F,W),u_{IJ})
\end{eqnarray*}
\end{itemize}
\end{defi}

\begin{defi}\label{inverseHdgsing}
Let $f:X\to S$ a morphism with $X,S\in\Var(k)$.
Assume there exist a factorization $f:X\xrightarrow{l}Y\times S\xrightarrow{p_S}S$ 
with $Y\in\SmVar(k)$, $l$ a closed embedding and $p_S$ the projection.
Let $S=\cup_{i\in I}$ an open cover such that there exist closed embeddings
$i:S_i\hookrightarrow\tilde S_i$ with $\tilde S_i\in\SmVar(k)$. 
Denote $X_I:=f^{-1}(S_I)$. We have then $X=\cup_{i\in I}X_i$ and the commutative diagrams
\begin{equation*}
\xymatrix{f:X_I\ar[r]^{l_I}\ar[rd] & Y\times S_I\ar[r]^{p_{S_I}}\ar[d]^{i_I':=(I\times i_I)} & S_I\ar[d]^{i_I} \\ 
\, & Y\times\tilde S_I\ar[r]^{p_{\tilde S_I}=:\tilde f_I} & \tilde S_I} 
\end{equation*}
\begin{itemize}
\item[(i)] For $((M_I,F,W),u_{IJ})\in C(DRM(S))$ we set (see definition \ref{gammaHdgsing} for $l$)
\begin{equation*}
f_{Hdg}^{*mod}((M_I,F,W),u_{IJ}):=\Gamma_X^{Hdg}(p_{\tilde S_I}^{*mod[-]}(M_I,F,W),u_{IJ})(d_Y)[2d_Y]\in C(DRM(X)), 
\end{equation*}
\item[(ii)] For $((M_I,F,W),u_{IJ})\in C(DRM(S)))$ we set (see definition \ref{gammaHdgsing} for $l$)
\begin{equation*}
f_{Hdg}^{\hat*mod}(M,F,W):=\Gamma_X^{\vee,Hdg}(p_{\tilde S_I}^{*mod[-]}(M_I,F,W),p_{\tilde S_I}^{*mod[-]}u_{IJ})
\in C(DRM(X)), 
\end{equation*}
\end{itemize}
Let $j:S^o\hookrightarrow S$ an open embedding with $S\in\Var(k)$.
Let $S=\cup_{i\in I}$ an open cover such that there exist closed embeddings
$i:S_i\hookrightarrow\tilde S_i$ with $\tilde S_i\in\SmVar(k)$. We have then,
for $((M_I,F,W),u_{IJ})\in C(DRM(S))$, quasi-isomorphisms in $C(DRM(S))$
\begin{eqnarray*}
I(j^*,j^{*mod}_{Hdg})(-): 
j^*((M_I,F,W),u_{IJ}):=(\tilde j_I^*(M_I,F,W),\tilde j_I^*u_{IJ})\to j_{Hdg}^{*mod}((M_I,F,W),u_{IJ})
\end{eqnarray*}
and
\begin{eqnarray*}
I(j^*,j^{\hat*mod}_{Hdg})(-):j_{Hdg}^{\hat*mod}((M_I,F,W),u_{IJ})\to 
j^*((M_I,F,W),u_{IJ}):=(\tilde j_I^*(M_I,F,W),\tilde j_I^*u_{IJ}).
\end{eqnarray*}
\end{defi}

\begin{defi}\label{otimesHdg}
Let $S\in\Var(k)$.
Let $S=\cup_{i\in I}$ an open cover such that there exist closed embeddings
$i:S_i\hookrightarrow\tilde S_i$ with $\tilde S_i\in\SmVar(k)$. We have the following functor
\begin{eqnarray*}
(-)\otimes^{Hdg}(-):D(DRM(S))^2\to D(DRM(S)), \\
(((M,F,W),u_{IJ}),((N,F,W),v_{IJ}))\mapsto ((M,F,W),u_{IJ})\otimes_{O_S}^{Hdg}((N,F,W),v_{IJ}):= \\
\Delta_{S,Hdg}^{*mod}(p_{1I}^{*mod}(M_I,F,W)\otimes_{O_{\tilde S_I\times\tilde S_I}}p_{2I}^{*mod}(N_I,F,W),
p_{1I}^{*mod}u_{IJ}\otimes p_{2I}^{*mod}v_{IJ}):= \\
\Delta_S^*\Gr_{V_{\Delta_S,0}}\Gamma_{\Delta_S}^{\vee,Hdg}
(p_{1I}^{*mod}(M_I,F,W)\otimes_{O_{\tilde S_I\times\tilde S_I}}p_{2I}^{*mod}(N_I,F,W),
p_{1I}^{*mod}u_{IJ}\otimes p_{2I}^{*mod}v_{IJ})
\end{eqnarray*}
using the definition \ref{inverseHdgsing} for the diagonal closed embedding $\Delta_S:S\hookrightarrow S\times S$.
\end{defi}

\begin{prop}\label{compDmodDRHdg}
Let $f_1:X\to Y$ and $f_2:Y\to S$ two morphism with $X,Y,S\in\QPVar(k)$. 
\begin{itemize}
\item[(i)]Let $(M,F,W)\in C(DRM(S)))$. Then, 
\begin{equation*}
(f_2\circ f_1)_{Hdg}^{*mod}(M,F,W)=f_{1Hdg}^{*mod}f_{2Hdg}^{*mod}(M,F,W)\in D(DRM(X)).
\end{equation*}
\item[(ii)]Let $(M,F,W)\in C(DRM(S)))$. Then,
\begin{equation*}
(f_2\circ f_1)_{Hdg}^{\hat*mod}(M,F,W)=f_{1Hdg}^{\hat*mod}f_{2Hdg}^{\hat*mod}(M,F,W)\in D(DRM(X))
\end{equation*}
\end{itemize}
\end{prop}

\begin{proof}
Immediate from definition.
\end{proof}

\begin{thm}\label{Bek0}
\begin{itemize}
\item[(i)]Let $S\in\Var(k)$. Let $S=\cup_{i\in I}S_i$ an open cover such that there exists
closed embedding $i_i:S_i\hookrightarrow\tilde S_i$ with $\tilde S_i\in\SmVar(k)$. Then the full embedding
\begin{eqnarray*}
\iota_S:DRM(S)\hookrightarrow\PSh^0_{\mathcal D(1,0)fil,rh}(S/(\tilde S_I))
\hookrightarrow C_{\mathcal D(1,0)fil,rh}(S/(\tilde S_I)) 
\end{eqnarray*}
induces a full embedding
\begin{equation*}
\iota_S:D(DRM(S)\hookrightarrow D_{\mathcal D(1,0)fil,rh}(S/(\tilde S_I)) 
\end{equation*}
whose image consists of $((M_I,F,W),u_{IJ})\in D_{\mathcal D(1,0)fil,rh}(S/(\tilde S_I))$ such that 
$(H^n(M_I,F,W),H^n(u_{IJ}))\in DRM(S)$ for all $n\in\mathbb Z$ and such that for all $p\in\mathbb Z$,
the differentials of $\Gr_W^p(M_I,F)$ are strict for the filtrations $F$.
\item[(i)'] Let $S\in\Var(k)$. Let $S=\cup_{i\in I}S_i$ an open cover such that there exists
closed embedding $i_i:S_i\hookrightarrow\tilde S_i$ with $\tilde S_i\in\SmVar(k)$. We have then
\begin{eqnarray*}
D(DRM(S)=<\int_{p_S}(n\times I)_{!Hdg}(\Gamma_X^{\vee,Hdg}(O_{\mathbb P^{N,o}\times\tilde S_I},F_b),x_{IJ}), \;
(f:X\xrightarrow{l}\mathbb P^{N,o}\times S\xrightarrow{p_S}S)\in\QPVar(k)> \\
=<\int_{p_S}(\Gamma_X^{\vee,Hdg}(O_{\mathbb P^{N,o}\times\tilde S_I},F_b),x_{IJ})
(f:X\xrightarrow{l}\mathbb P^N\times S\xrightarrow{p_S}S)\in\QPVar(k), \; \mbox{proper}, \; X \mbox{smooth}> \\
\subset D_{\mathcal D(1,0)fil,rh}(S/(\tilde S_I))
\end{eqnarray*}
where $n:\mathbb P^{N,o}\hookrightarrow\mathbb P^N$ are open embeddings, $l$ are closed embedding
and $<,>$ means the full triangulated category generated by.
\item[(ii)]Let $S\in\Var(k)$. Let $S=\cup_{i\in I}S_i$ an open cover such that there exists
closed embedding $i_i:S_i\hookrightarrow\tilde S_i$ with $\tilde S_i\in\SmVar(k)$. Then the full embedding
\begin{eqnarray*}
\iota_S:DRM(S)\hookrightarrow\PSh^0_{\mathcal D(1,0)fil,rh}(S/(\tilde S_I))
\hookrightarrow C_{\mathcal D(1,0)fil,rh}(S/(\tilde S_I)) 
\end{eqnarray*}
induces a full embedding
\begin{equation*}
\iota_S:D(DRM(S))\hookrightarrow D_{\mathcal D(1,0)fil,\infty,rh}(S/(\tilde S_I)) 
\end{equation*}
whose image consists of $((M_I,F,W),u_{IJ})\in D_{\mathcal D(1,0)fil,\infty,rh}(S/(\tilde S_I))$ such that 
$(H^n(M_I,F,W),H^n(u_{IJ}))\in DRM(S)$ for all $n\in\mathbb Z$ 
and such that there exist $r\in\mathbb Z$ and an $r$-filtered homotopy equivalence
$((M_I,F,W),u_{IJ})\to ((M'_I,F,W),u_{IJ})$ such that for all $p\in\mathbb Z$
the differentials of $\Gr_W^p(M'_I,F)$ are strict for the filtrations $F$.
\end{itemize}
\end{thm}

\begin{proof}
\noindent(i):We first show that $\iota_S$ is fully faithfull, that is for all
$\mathcal M=((M_I,F,W),u_{IJ}),\mathcal M'=((M'_I,F,W),u_{IJ})\in DRM(S)$ 
and all $n\in\mathbb Z$,
\begin{eqnarray*}
\iota_S:\Ext_{D(DRM(S))}^n(\mathcal M,\mathcal M'):=\Hom_{D(DRM(S))}(\mathcal M,\mathcal M'[n]) \\
\to\Ext_{\mathcal D(S)_0}^n(\mathcal M,\mathcal M')
:=\Hom_{\mathcal D(S)_0:=D_{\mathcal D(1,0)fil,rh}(S/(\tilde S_I))}(\mathcal M,\mathcal M'[n])
\end{eqnarray*}
For this it is enough to assume $S$ smooth. We then proceed by induction on $max(\dim\supp(M),\dim\supp(M'))$. 
\begin{itemize}
\item For $\supp(M)=\supp(M')=\left\{s\right\}$, it is obvious. 
If $\supp(M)=\left\{s\right\}$ and $\supp(M')=\left\{s'\right\}$ wit $s'\neq s$, then by the localization exact sequence
\begin{equation*}
\Ext_{D(MHM(S))}^n(\mathcal M,\mathcal M')=0=\Ext_{\mathcal D(S)}^n(\mathcal M,\mathcal M')
\end{equation*}
\item Denote $\supp(M)=Z\subset S$ and $\supp(M')=Z'\subset S$.
There exist an open subset $S^o\subset S$ such that $Z^o:=Z\cap S^o$ and $Z^{'o}:=Z'\cap S^o$ are smooth,
and $\mathcal M_{|Z^o}:=(i^*\Gr_{V_{Z^o},0}M_{|S^o},F,W)\in DRM(Z^o)$ and 
$\mathcal M'_{|Z^{'o}}:=(i^{'*}\Gr_{V_{Z^{'o}},0}M'_{|S^o},F,W)\in DRM(Z^{'o})$ 
are filtered vector bundles, where $j:S^o\hookrightarrow S$ is the open embedding, and
$i:Z^o\hookrightarrow S^o$, $i:Z^{'o}\hookrightarrow S^o$ the closed embeddings.
Considering the connected components of $Z^o$ and $Z^{'o}$, we way assume that $Z^o$ and $Z^{'o}$ are connected.
Shrinking $S^o$ if necessary, we may assume that either $Z^o=Z^{'o}$ or $Z^o\cap Z^{'o}=\emptyset$,
We denote $D=S\backslash S^o$. Shrinking $S^o$ if necessary, 
we may assume that $D$ is a divisor and denote by $l:S\hookrightarrow L_D$ the zero section embedding.
\begin{itemize}
\item If $Z^o=Z^{'o}$, denote $i:Z^o\hookrightarrow S^o$ the closed embedding.
We have then the following commutative diagram
\begin{equation*}
\xymatrix{\Ext_{D(DRM(S^o))}^n(\mathcal M_{|S^o},\mathcal M'_{|S^o})
\ar[rr]^{\iota_{S^o}}\ar[d]_{i^*\Gr_{V_{Z^o},0}} & \, & 
\Ext_{\mathcal D(S^o)_0}^n(\mathcal M_{|S^o},\mathcal M'_{|S^o})\ar[d]^{i^*\Gr_{V_{Z^o},0}} \\
\Ext_{D(DRM(Z^o))}^n(\mathcal M_{|Z^o},\mathcal M'_{|Z^o})\ar[rr]^{\iota_{Z^o}}\ar[u]_{i_{*mod}} & \, &
\Ext_{\mathcal D(Z^o)_0}^n(\mathcal M_{|Z^o},\mathcal M'_{|Z^o})\ar[u]^{i_{*mod}}}
\end{equation*}
Now we prove that $\iota_{Z^o}$ is an isomorphism similarly to the proof the the generic case of \cite{Be}.
On the other hand the left and right colummn are isomorphisms.
Hence $\iota_{S^o}$ is an isomorphism by the diagram.
\item If $Z^o\cap Z^{'o}=\emptyset$, we consider the following commutative diagram
\begin{equation*}
\xymatrix{\Ext_{D(DRM(S^o))}^n(\mathcal M_{|S^o},\mathcal M'_{|S^o})
\ar[rr]^{\iota_{S^o}}\ar[d]_{i^*\Gr_{V_{Z^o},0}} & \, & 
\Ext_{\mathcal D(S^o)_0}^n(\mathcal M_{|S^o},\mathcal M'_{|S^o})\ar[d]^{i^*\Gr_{V_{Z^o},0}} \\
\Ext_{D(DRM(Z^o))}^n(\mathcal M_{|Z^o},0)=0\ar[rr]^{\iota_{Z^o}}\ar[u]_{i_{*mod}} & \, &
\Ext_{\mathcal D(Z^o)_0}^n(\mathcal M_{|Z^o},0)=0\ar[u]^{i_{*mod}}}
\end{equation*}
where the left and right column are isomorphism by strictness of the $V_{Z^o}$ filtration
(use a bi-filtered injective resolution with respect to $F$ and $V_{Z^o}$ for the right column).
\end{itemize}
\item We consider now the following commutative diagram in $C(\mathbb Z)$ where we denote for short $H_0:=D(DRM(S))$
\begin{equation*}
\xymatrix{0\ar[r] & \Hom_{H_0}^{\bullet}(\Gamma^{\vee,Hdg}_D\mathcal M,\Gamma^{Hdg}_D\mathcal M')
\ar[r]^{\Hom(-,\gamma^{Hdg}_D(\mathcal M'))}\ar[d]^{\iota_S} &
\Hom_{H_0}^{\bullet}(\Gamma^{\vee,Hdg}_D\mathcal M,\mathcal M')
\ar[r]^{\Hom(-,\ad(j^*,j_{*Hdg})(\mathcal M'))}\ar[d]^{\iota_S} &
\Hom_{H_0}^{\bullet}(\Gamma^{\vee,Hdg}_D\mathcal M,j_{*Hdg}j^*\mathcal M')\ar[r]\ar[d]^{\iota_S} & 0 \\
0\ar[r] & \Hom_{\mathcal D(S)_0}^{\bullet}(\Gamma^{\vee,Hdg}_D\mathcal M,\Gamma^{Hdg}_D\mathcal M')
\ar[r]^{\Hom(-,\gamma^{Hdg}_D(\mathcal M'))} &
\Hom_{\mathcal D(S)_0}^{\bullet}(\Gamma^{\vee,Hdg}_D\mathcal M,\mathcal M')\ar[r]^{\Hom(-,\ad(j^*,j_{*Hdg})(\mathcal M'))} &
\Hom_{\mathcal D(S)_0}^{\bullet}(\Gamma^{\vee,Hdg}_D\mathcal M,j_{*Hdg}j^*\mathcal M')\ar[r] & 0}
\end{equation*}
whose lines are exact sequence. We have on the one hand,
\begin{equation*}
\Hom_{D(DRM(S))}^{\bullet}(\Gamma^{\vee,Hdg}_D\mathcal M,j_{*Hdg}j^*\mathcal M')=0=
\Hom_{\mathcal D(S)}^{\bullet}(\Gamma^{\vee,Hdg}_D\mathcal M,j_{*Hdg}j^*\mathcal M')
\end{equation*}
On the other hand by induction hypothesis
\begin{equation*}
\iota_S:\Hom_{D(DRM(S))}^{\bullet}(\Gamma^{\vee,Hdg}_D\mathcal M,\Gamma^{Hdg}_D\mathcal M')\to
\Hom_{\mathcal D(S)_0}^{\bullet}(\Gamma^{\vee,Hdg}_D\mathcal M,\Gamma^{Hdg}_D\mathcal M')
\end{equation*}
is a quasi-isomorphism. Hence, by the diagram
\begin{equation*}
\iota_S:\Hom_{D(DRM(S))}^{\bullet}(\Gamma^{\vee,Hdg}_D\mathcal M,\mathcal M')\to
\Hom_{\mathcal D(S)_0}^{\bullet}(\Gamma^{\vee,Hdg}_D\mathcal M,\mathcal M')
\end{equation*}
is a quasi-isomorphism.
\item We consider now the following commutative diagram in $C(\mathbb Z)$ where we denote for short $H_0:=D(DRM(S))$
\begin{equation*}
\xymatrix{0\ar[r] & \Hom_{H_0}^{\bullet}(\Gamma^{\vee,Hdg}_D\mathcal M,\mathcal M')
\ar[r]^{\Hom(\gamma^{\vee,Hdg}_D(\mathcal M),-)}\ar[d]^{\iota_S} &
\Hom_{H_0}^{\bullet}(\mathcal M,\mathcal M')
\ar[r]^{\Hom(\ad(j_{!Hdg},j^*)(\mathcal M'),-)}\ar[d]^{\iota_S} &
\Hom_{H_0}^{\bullet}(j_{!Hdg}j^*\mathcal M,\mathcal M')\ar[r]\ar[d]^{\iota_S} & 0 \\
0\ar[r] & \Hom_{\mathcal D(S)_0}^{\bullet}(\Gamma^{\vee,Hdg}_D\mathcal M,\mathcal M')
\ar[r]^{\Hom(\gamma^{\vee,Hdg}_D(\mathcal M),-)} &
\Hom_{\mathcal D(S)_0}^{\bullet}(\mathcal M,\mathcal M')\ar[r]^{\Hom(\ad(j_{!Hdg},j^*)(\mathcal M),-)} &
\Hom_{\mathcal D(S)_0}^{\bullet}(j_{!Hdg}j^*\mathcal M,\mathcal M')\ar[r] & 0}
\end{equation*}
whose lines are exact sequence. On the one hand, the commutative diagram
\begin{equation*}
\xymatrix{\Hom_{D(DRM(S))}^{\bullet}(j_{!Hdg}j^*\mathcal M,\mathcal M')\ar[r]^{j^*}\ar[d]^{\iota_{S}} &
\Hom_{D(DRM(S^o))}^{\bullet}(j^*\mathcal M,j^*\mathcal M')\ar[d]^{\iota_{S^o}} \\
\Hom_{\mathcal D(S)_0}^{\bullet}(j_{!Hdg}j^*\mathcal M,\mathcal M')\ar[r]^{j^*} &
\Hom_{\mathcal D(S^o)_0}^{\bullet}(j^*\mathcal M,j^*\mathcal M')}
\end{equation*}
together with the fact that the horizontal arrows $j^*$ are quasi-isomorphism 
by the functoriality given the uniqueness of the $V_S$ filtration for the embedding $l:S\hookrightarrow L_D$, 
(use a bi-filtered injective resolution with respect to $F$ and $V_S$ for the lower arrow)
and the fact that $\iota_{S^o}$ is a quasi-isomorphism by the first two point, show that
\begin{equation*}
\iota_S:\Hom_{D(DRM(S))}^{\bullet}(j_{!Hdg}j^*\mathcal M,\mathcal M')\to
\Hom_{\mathcal D(S)_0}^{\bullet}(j_{!Hdg}j^*\mathcal M,\mathcal M')
\end{equation*}
is a quasi-isomorphism. On the other hand, by the third point
\begin{equation*}
\iota_S:\Hom_{D(DRM(S))}^{\bullet}(\Gamma^{\vee,Hdg}_D\mathcal M,\mathcal M')\to
\Hom_{\mathcal D(S)_0}^{\bullet}(\Gamma^{\vee,Hdg}_D\mathcal M,\mathcal M')
\end{equation*}
is a quasi-isomorphism. Hence, by the diagram
\begin{equation*}
\iota_S:\Hom_{D(DRM(S))}^{\bullet}(\Gamma^{\vee,Hdg}_D\mathcal M,\mathcal M')\to
\Hom_{\mathcal D(S)_0}^{\bullet}(\Gamma^{\vee,Hdg}_D\mathcal M,\mathcal M')
\end{equation*}
is a quasi-isomorphism.
\end{itemize}
This shows the fully faithfulness. We now prove the essential surjectivity : let
\begin{equation*}
((M_I,F,W),u_{IJ})\in C_{\mathcal D(1,0)fil,rh}(S/(\tilde S_I)) 
\end{equation*}
such that the cohomology are mixed hodge modules and such that the differential are strict.
We proceed by induction on $card\left\{n\in\mathbb Z\right\}, \, \mbox{s.t.} H^n(M_I,F,W)\neq 0$ 
by taking the cohomological troncation 
and using the fact that the differential are strict for the filtration $F$ and the fully faithfullness.

\noindent(i)':Follows from (i).

\noindent(ii):Follows from (i). Indeed, in the composition of functor
\begin{eqnarray*}
\iota_S:D(DRM(S))\xrightarrow{\iota_S}D_{\mathcal D(1,0)fil,rh}(S/(\tilde S_I))
\to D_{\mathcal D(1,0)fil,\infty,rh}(S/(\tilde S_I)) 
\end{eqnarray*}
the second functor which is the localization functor is an isomorphism on the full subcategory
$D_{\mathcal D(1,0)fil,rh}(S/(\tilde S_I))^{st}\subset D_{\mathcal D(1,0)fil,rh}(S/(\tilde S_I))$ 
constisting of complex such that the differentials are strict for $F$, 
and the first functor $\iota_S$ is a full embedding by (i) and 
$\iota_S(D(DRM(S)))\subset D_{\mathcal D(1,0)fil,rh}(S/(\tilde S_I))^{st}$.
\end{proof}


\begin{defi}\label{DHdg}
Let $f:X\to S$ a morphism with $X,S\in\SmVar(k)$. 
Consider a compactification $f:X\xrightarrow{j}\bar X\xrightarrow{\bar f} S$ of $f$,
in particular $j$ is an open embedding and $\bar f$ is proper.
\begin{itemize}
\item[(i)]For $(M,F,W)\in C(DRM(X))$, we define, using definition \ref{DHdgj}, 
\begin{eqnarray*}
\int_f^{Hdg}(M,F,W):=\int_{\bar f}^{FDR}j_{*Hdg}(M,F,W)\in D_{\mathcal D(1,0)fil}(S)
\end{eqnarray*} 
It does not depends on the choice of the compactification as in the complex case:
for two compactification $\bar f:\bar X\to S$, $\bar f':X'\to S$, there exist
a compactification $\bar f'':\bar X''\to S$ together with morphisms 
$e:\bar X''\to\bar X$ and $e':\bar X''\to\bar X'$ such that $\bar f\circ e=\bar f\circ e'=\bar f''$. 
Let $(M,F,W)\in C(DRM(X)))$, then 
\begin{itemize}
\item by definition $H^i\int_{\bar f}^{FDR}\Gr_W^kj_{*Hdg}(M,F,W)\in PDRM(S)$ for all $i,k\in\mathbb Z$, 
hence by the spectral sequence for the filtered complex $\int_{\bar f}^{FDR}j_{*Hdg}(M,W)$
\begin{equation*}
\Gr_W^k(H^i\int_{f}^{Hdg}(M,F,W))=\Gr_W^k(H^i\int_{\bar f}^{FDR}j_{*Hdg}(M,F,W))\in PDRM(S)
\end{equation*}
since it is a sub-quotient of $H^i\int_{\bar f}^{FDR}\Gr_W^kj_{*Hdg}(M,F,W)$,
this gives by definition $H^i\int_{f}^{Hdg}(M,F,W))\in DRM(S)$ for all $i\in\mathbb Z$. 
\item $\int_f^{Hdg}(M,F,W)$ is the class of a complex such that the differential are strict for $F$
by theorem \ref{Sa12} in the complex case
\end{itemize}
We then set using theorem \ref{Bek0}
\begin{eqnarray*}  
Rf^{Hdg}_*(M,F,W):=\iota_S^{-1}\int_f^{Hdg}(M,F,W)\in D(DRM(S))
\end{eqnarray*}
\item[(ii)]For $(M,F,W)\in\in C(DRM(X))$, we define, using definition \ref{DHdgj}, 
\begin{eqnarray*}
\int_{f!}^{Hdg}(M,F,W):=\int_{\bar f}^{FDR}j_{!Hdg}(M,F,W)\in D_{\mathcal D(1,0)fil}(S)
\end{eqnarray*} 
It does not depends on the choice of the compactification 
as in the complex case: for two compactification $\bar f:\bar X\to S$, $\bar f':X'\to S$, there exist
a compactification $\bar f'':\bar X''\to S$ together with morphisms 
$e:\bar X''\to\bar X$ and $e':\bar X''\to\bar X'$ such that $\bar f\circ e=\bar f\circ e'=\bar f''$.
Let $(M,F,W)\in C(DRM(X)))$, then 
\begin{itemize}
\item by definition $H^i\int_{\bar f}^{FDR}\Gr_W^kj_{!Hdg}(M,F,W)\in PDRM(S)$ for all $i,k\in\mathbb Z$, 
hence by the spectral sequence for the filtered complex $\int_{\bar f}^{FDR}j_{!Hdg}(M,W)$
\begin{equation*}
\Gr_W^k(H^i\int_{f!}^{Hdg}(M,F,W))==\Gr_W^k(H^i\int_{\bar f}^{FDR}j_{!Hdg}(M,F,W))\in PDRM(S)
\end{equation*}
since it is a sub-quotient of $H^i\int_{\bar f}^{FDR}\Gr_W^kj_{!Hdg}(M,F,W)$,
this gives by definition $H^i\int_{f!}^{Hdg}(M,F,W))\in DRM(S)$ for all $i\in\mathbb Z$.
\item $\int_{f!}^{Hdg}(M,F,W)$ is the class of a complex such that the differential are strict for $F$
by theorem \ref{Sa12} in the complex case.
\end{itemize}
We then set using theorem \ref{Bek0}
\begin{eqnarray*}  
Rf^{Hdg}_!(M,F,W):=\iota_S^{-1}\int_{f!}^{Hdg}(M,F,W)\in D(DRM(S))
\end{eqnarray*}
\end{itemize}
\end{defi}

In the singular case, we set the following

\begin{defi}\label{DHdgsing}
Let $f:X\to S$ a morphism with $X,S\in\Var(k)$. 
Assume there exist a factorization $f:X\xrightarrow{l}Y\times S\xrightarrow{p_S}S$ with $Y\in\SmVar(k)$,
$l$ a closed embedding and $p_S$ the projection.
Let $\bar Y\in\PSmVar(k)$ a compactification of $Y$ and denote by $n:Y\hookrightarrow\bar Y$ the open embedding.
Denote again $p_S:\bar Y\times S\to S$ the projection.
Let $S=\cup_{i\in I}S_i$ an open cover such that there exists
closed embedding $i_i:S_i\hookrightarrow\tilde S_i$ with $\tilde S_i\in\SmVar(k)$.
We have then the open cover $X=\cup_iX_i$ with $X_i:=f^{-1}(S_i)$ together with
closed embeddings $i'_I:X_I:\hookrightarrow Y\times\tilde S_I$.
\begin{itemize}
\item[(i)]For $((M_I,F,W),u_{IJ})\in C(DRM(X))$, we define, 
using definition \ref{DHdgj} for $n\times I:Y\times S\hookrightarrow\bar Y\times S$, 
\begin{eqnarray*}
\int_{p_S}^{Hdg}((M_I,F,W),u_{IJ}):=\int_{p_S}^{FDR}(n\times I)_{*Hdg}((M_I,F,W),u_{IJ})
\in D_{\mathcal D(1,0)fil}(S/\tilde S_I)
\end{eqnarray*} 
with 
\begin{eqnarray*}
(n\times I)_{*Hdg}((M_I,F,W),u_{IJ}):=((n\times I)_{*Hdg}(M_I,F,W),(n\times I)_{*Hdg}u_{IJ})
\in C(\bar X/(\bar Y\times\tilde S_I))
\end{eqnarray*}
We then set using theorem \ref{Bek0} and theorem \ref{Sa12}
\begin{eqnarray*}
Rf_*^{Hdg}((M_I,F,W),u_{IJ}):=\iota_S^{-1}\int_{p_S}^{Hdg}((M_I,F,W),u_{IJ})\in D(DRM(S)).
\end{eqnarray*}
\item[(ii)]For $((M_I,F,W),u_{IJ})\in C(DRM(X))$, we define, 
using definition \ref{DHdgj} for $n\times I:Y\times S\hookrightarrow\bar Y\times S$,
\begin{eqnarray*}
\int_{p_S!}^{Hdg}((M_I,F,W),u_{IJ}):=\int_{p_S}^{FDR}(n\times I)_{!Hdg}((M_I,F,W),u_{IJ})
\in D_{\mathcal D(1,0)fil}(S/(\tilde S_I))
\end{eqnarray*} 
with 
\begin{eqnarray*}
(n\times I)_{!Hdg}((M_I,F,W),u_{IJ}):=((n\times I)_{!Hdg}(M_I,F,W),((n\times I)_{*Hdg}u^d_{IJ})^d)
\in C(\bar X/(\bar Y\times\tilde S_I))
\end{eqnarray*}
We then set using theorem \ref{Bek0} and theorem \ref{Sa12}
\begin{eqnarray*}
Rf_!^{Hdg}((M_I,F,W),u_{IJ}):=\iota_S^{-1}\int_{p_S!}^{Hdg}((M_I,F,W),u_{IJ})\in D(DRM(S)).
\end{eqnarray*}
\end{itemize}
\end{defi}

\begin{prop}\label{compDmodDRHdgD}
Let $f_1:X\to Y$ and $f_2:Y\to S$ two morphism with $X,Y,S\in\QPVar(k)$ or with $X,Y,S\in\SmVar(k)$. 
\begin{itemize}
\item[(i)]Let $(M,F,W)\in C(DRM(X))$. Then, 
\begin{equation*}
R(f_2\circ f_1)^{Hdg}_*(M,F,W)=Rf^{Hdg}_{2*}Rf^{Hdg}_{1*}(M,F,W)\in D(DRM(S)).
\end{equation*}
\item[(ii)]Let $(M,F,W)\in C(DRM(X))$. Then,
\begin{equation*}
R(f_2\circ f_1)^{Hdg}_!(M,F,W)=Rf^{Hdg}_{2!}Rf^{Hdg}_{1!}(M,F,W)\in D(DRM(S))
\end{equation*}
\end{itemize}
\end{prop}

\begin{proof}
Immediate from definition.
\end{proof}

\begin{prop}\label{Hdgpropad}
Let $f:X\to S$ with $S,X\in\SmVar(k)$ or with $S,X\in\QPVar(k)$. Then
\begin{itemize}
\item[(i)] $(f_{Hdg}^{\hat*mod},Rf^{Hdg}_*):D(DRM(S))\to D(DRM(X))$ is a pair of adjoint functors.
For $(M,F,W)\in C(DRM(S))$ we denote by 
\begin{equation*}
\ad(f_{Hdg}^{\hat*mod},Rf^{Hdg}_*)(M,F,W):(M,F,W)\to Rf^{Hdg}_*f_{Hdg}^{\hat*mod}(M,F,W) 
\end{equation*}
the adjonction map in $D(DRM(S))$. For $(N,F,W)\in C(DRM(X))$, we denote by 
\begin{equation*}
\ad(f_{Hdg}^{\hat*mod},Rf^{Hdg}_*)(N,F,W):f_{Hdg}^{\hat*mod}Rf^{Hdg}_*(N,F,W)\to (N,F,W) 
\end{equation*}
the adjonction map in $D(DRM(X))$. 
\item[(ii)]$(Rf^{Hdg}_!,f_{Hdg}^{*mod}):D(DRM(X))\to D(DRM(S))$ is a pair of adjoint functors.
For $(M,F,W)\in C(DRM(S))$ we denote by 
\begin{equation*}
\ad(Rf^{Hdg}_!,f_{Hdg}^{*mod})(M,F,W):Rf^{Hdg}_!f_{Hdg}^{*mod}(M,F,W)\to(M,F,W) 
\end{equation*}
the adjonction map in $D(DRM(S))$. For $(N,F,W)\in C(DRM(X))$, we denote by 
\begin{equation*}
\ad(Rf^{Hdg}_!,f_{Hdg}^{*mod})(N,F,W):(N,F,W)\to f_{Hdg}^{*mod}Rf^{Hdg}_!(N,F,W) 
\end{equation*}
the adjonction map in $D(DRM(X))$. 
\end{itemize}
\end{prop}

\begin{proof}
Follows from proposition \ref{jHdgpropad} after considering the graph factorization 
$f:X\hookrightarrow\bar X\times S\xrightarrow{p_S} S$ with $\bar X\in\PSmVar(k)$ a compactification of $X$.
\end{proof}

We have by proposition \ref{compDmodDRHdg} and proposition \ref{compDmodDRHdgD} the 2 functors on $\QPVar(k)$ :
\begin{itemize}
\item $D(DRM(-)):\QPVar(k)\to D(DRM(-)), \; S\mapsto D(DRM(S)), \, (f:T\to S)\mapsto Rf^{Hdg}_*$,
\item $D(DRM(-)):\QPVar(k)\to D(DRM(-)), \; S\mapsto D(DRM(S)), \, (f:T\to S)\mapsto Rf^{Hdg}_!$,
\item $D(DRM(-)):\QPVar(k)\to D(DRM(-)), \; S\mapsto D(DRM(S)), \, (f:T\to S)\mapsto f_{Hdg}^{*mod}$,
\item $D(DRM(-)):\QPVar(k)\to D(DRM(-)), \; S\mapsto D(DRM(S)), \, (f:T\to S)\mapsto f_{Hdg}^{\hat*mod}$.
\end{itemize}

\begin{prop}\label{Hdgpropadsing}
Let $f:X\to S$ with $S,X\in\SmVar(k)$ or with $S,X\in\QPVar(k)$. Then
\begin{itemize}
\item[(i)] $(f_{Hdg}^{\hat*mod},Rf^{Hdg}_*):D(DRM(S))\to D(DRM(X))$ is a pair of adjoint functors.
For $((M_I,F,W),u_{IJ})\in C(DRM(S))$ we denote by 
\begin{eqnarray*}
\ad(f_{Hdg}^{\hat*mod},Rf^{Hdg}_*)((M_I,F,W),u_{IJ}):((M_I,F,W),)\to Rf^{Hdg}_*f_{Hdg}^{\hat*mod}((M_I,F,W),u_{IJ}) 
\end{eqnarray*}
the adjonction map in $D(DRM(S))$. 
For $((N_I,F,W),u_{IJ})\in C(DRM(X))$, we denote by 
\begin{eqnarray*}
\ad(f_{Hdg}^{\hat*mod},Rf^{Hdg}_*)((N_I,F,W),u_{IJ}):f_{Hdg}^{\hat*mod}Rf^{Hdg}_*((N_I,F,W),u_{IJ})\to((N_I,F,W),u_{IJ}) 
\end{eqnarray*}
the adjonction map in $D(DRM(X))$ 
\item[(ii)]$(Rf^{Hdg}_!,f_{Hdg}^{*mod}):D(DRM(X))\to D(DRM(S))$ is a pair of adjoint functors.
For $((M_I,F,W),u_{IJ})\in C(DRM(S))$ we denote by 
\begin{eqnarray*}
\ad(Rf^{Hdg}_!,f_{Hdg}^{*mod})((M_I,F,W),u_{IJ}):Rf^{Hdg}_!f_{Hdg}^{*mod}((M_I,F,W),u_{IJ})\to((M_I,F,W),u_{IJ}) 
\end{eqnarray*}
the adjonction map in $D(DRM(S))$. 
For $((N_I,F,W),u_{IJ})\in C(DRM(X))$, we denote by 
\begin{eqnarray*}
\ad(Rf^{Hdg}_!,f_{Hdg}^{*mod})((N_I,F,W),u_{IJ}):((N_I,F,W),u_{IJ})\to f_{Hdg}^{*mod}Rf^{Hdg}_!((N_I,F,W),u_{IJ}) 
\end{eqnarray*}
the adjonction map in $D(DRM(X))$. 
\end{itemize}
\end{prop}

\begin{proof}
Follows from proposition \ref{jHdgpropad} after considering a factorization 
$f:X\hookrightarrow\bar Y\times S\xrightarrow{p_S} S$ with $\bar Y\in\PSmVar(k)$.
\end{proof}

\begin{thm}\label{sixDRMk}
Let $k$ a field of characteristic zero. 
\begin{itemize}
\item[(i)]We have the six functor formalism on $D(DRM(-)):\SmVar(k)\to\TriCat$.
\item[(ii)]We have the six functor formalism on $D(DRM(-)):\QPVar(k)\to\TriCat$.
\end{itemize}
\end{thm}

\begin{proof}
Follows from proposition \ref{Hdgpropadsing}.
\end{proof}

\begin{thm}\label{phipsithmp}
Let $k\subset K\subset\mathbb C_p$ a subfield with $p$ a prime number and $K$ a $p$ adic field. 
Let $S\in\SmVar(k)$.
Let $D=V(s)\subset S$ a divisor with $s\in\Gamma(S,L)$ and $L$ a line bundle ($S$ being smooth, $D$ is Cartier).
so that we have the closed embedding $i:S\hookrightarrow L$, $i(x)=(x,s(x))$ 
and $D=i^{-1}(s_0)$, $s_0$ being the zero section.
For $(M,F,W)\in DRM(S)$,
\begin{itemize}
\item we have the canonical isomorphism in $D_{\mathbb B_{dr}}(S_K^{an,pet})$
\begin{eqnarray*}
T^{B_{dr}}(\psi_D,DR)(M,F,W):=B^{B_{dr}}(M,F,W)\circ A^{B_{dr}}(M,F,W)^{-1}: \\
DR(S)(\psi_D(M,F,W)^{an}\otimes_{O_{S_K}}(O\mathbb B_{dr,S_K},F))\xrightarrow{\sim}
\psi_DDR(S)((M,F,W)^{an}\otimes_{O_{S_K}}(O\mathbb B_{dr,S_K},F))[-1] 
\end{eqnarray*}
with, for $S=\cup_{i=1}^sS_i$ an open affine cover such that 
$D\cap S_i=V(f_i)\subset S_i$ is given by $f_i\in\Gamma(S_i,O_{S_i})$,
denoting $q:L_i:=p^{-1}(S_i)\to\mathbb A^1_k$ the projection and $j_i:S_i\hookrightarrow S$ the open embeddings,
\begin{itemize}
\item the isomorphism in $D_{\mathbb B_{dr}2fil}(S_K^{an,pet})$
\begin{eqnarray*}
A^{B_{dr}}(M,F,W):(\oplus_{i=1}^s\oplus_{-1\leq\alpha<0}
\Cone(\partial_s:DR(L_i/\mathbb A^1_k)((V_{D\alpha}(M,F^*,W))^{an})\otimes_{O_S}s^{\alpha+1}O\mathbb B_{dr,S_K} \\
\to DR(L_i/\mathbb A^1_k)((V_{D\alpha}(M,F^{*-1},W)^{an})\otimes_{O_S}s^{\alpha}O\mathbb B_{dr,S_K}))
\xrightarrow{(j_I^*)}\cdots)[-1] \\ 
\to(\oplus_{i=1}^sDR(S_i)(\psi_D(M,F,W)^{an}\otimes_{O_{S_K}}(O\mathbb B_{dr,S_K},F))\xrightarrow{j_I^*}\cdots)[-1] \\
\xrightarrow{((j_i^*),0)^{-1}}(DR(S)(\psi_D(M,F,W)^{an})\otimes_{O_{S_K}}(O\mathbb B_{dr,S_K},F)), \;  
(\sum_jm_j\otimes(\log s)^j,m')\mapsto [m_0] , 
\end{eqnarray*}
\item and the isomorphism in $D_{\mathbb B_{dr}2fil}(S_K^{an,pet})$
\begin{eqnarray*}
B^{B_{dr}}(M,F,W):(\oplus_{i=1}^s\oplus_{-1\leq\alpha<0}
\Cone(\partial_s:V_{D\alpha}DR(L_i/\mathbb A^1_k)((M,F^*,W)^{an})\otimes_{O_{S_K}}s^{\alpha+1}O\mathbb B_{dr,S_K} \\
\to V_{D\alpha}DR(L_i/\mathbb A^1_k)((M,F^{*-1},W)^{an})\otimes_{O_{S_K}}s^{\alpha}O\mathbb B_{dr,S_K})
\xrightarrow{(j_I^*)}\cdots)[-1] \\
\to (\oplus_{i=1}^sDR(p^*O_{\mathbb A_1^k})(i^*\pi_*\pi^{*mod}DR(L_i/\mathbb A^1_k)((M,F,W)^{an})
\otimes_{O_{S_K}}(O\mathbb B_{dr,S_K},F))
\xrightarrow{(j_I^*)}\cdots)[-1] \\
\xrightarrow{=}(\oplus_{i=1}^s\psi_DDR(p^*O_{\mathbb A_1^k})(DR(L_i/\mathbb A^1_k)((M,F,W)^{an})
\otimes_{O_{S_K}}(O\mathbb B_{dr,S_K},F))
\xrightarrow{(j_I^*)}\cdots)[-1] \\
\xrightarrow{((j_i^*),0)^{-1}}\psi_DDR(S)((M,F,W)^{an}\otimes_{O_{S_K}}(O\mathbb B_{dr,S_K},F))[-1], \\
(\sum_jm_j\otimes(\log s)^j,m')\mapsto\sum_j(\log s)^jm_j,
\end{eqnarray*}
\end{itemize}
so that $T_{B_{dr}}(\psi_D,DR)(M)\circ DR(S)((s\partial_s)\otimes I)=N\circ T_{B_{dr}}(\psi_D,DR)(M)$ where 
\begin{equation*}
N:=\log T_u\in\Hom(\psi_DDR(S)(M^{an}\otimes_{O_{S_K}}O\mathbb B_{dr,S_K}),
\psi_DDR(S)(M^{an}\otimes_{O_{S_K}}O\mathbb B_{dr,S_K})), 
\end{equation*}
is induced by the monodromy automorphism $T:\tilde S^o\xrightarrow{\sim}\tilde S^o$ of the
perfectoid universal covering $\pi:\tilde S^o\to S^o:=S\backslash D$ (see \cite{Scholze}).
\item there is a canonical isomorphism in $D_{\mathbb B_{dr}2fil}(S_K^{an,pet})$
\begin{eqnarray*}
T^{B_{dr}}(\phi_D,DR)(M,F,W):DR(S)(\phi_D(M,F,W)^{an}\otimes_{O_{S_K}}(O\mathbb B_{dr,S_K},F)) \\
\xrightarrow{DR(S)(0,(var(M,F,W)\otimes I))} 
DR(S)({\phi^{\rho}_D}(M,F,W)^{an}\otimes_{O_{S_K}}(O\mathbb B_{dr,S_K},F)) \\
\xrightarrow{(I,T^{B_{dr}}(M,F,W))}\phi_DDR(S)((M,F,W)^{an}\otimes_{O_{S_K}}(O\mathbb B_{dr,S_K},F))[-1]. 
\end{eqnarray*}
where 
\begin{equation*}
\phi^{\rho}_D(M,F,W):=\Cone(\theta_{DR,D}(M,F,W):\Gamma_D^{\vee,Hdg}(M,F,W)\to\psi_D(M,F,W))
\end{equation*}
with $\theta_{DR,D}(M,F,W)$ the factorization in $C_{\mathcal D(1,0)fil,rh}(S)$  
\begin{eqnarray*}
\rho_{DR,D}(M,F,W)\circ\ad(j^*,j_{*Hdg})(M,F,W)): \\
(M,F,W)\xrightarrow{\gamma_D^{\vee,Hdg}(M,F,W)}\Gamma^{\vee,Hdg}_D(M,F,W)\xrightarrow{\theta_{DR,D}(M)}\psi_D(M,F,W).
\end{eqnarray*}
of the map given in definition \ref{rhoDR}.
\end{itemize}
\end{thm}

\begin{proof}
Follows from the proof of theorem \ref{phipsi0thmp}: $j_i^*A^{B_{dr}}(M,F,W)$ are filtered quasi-isomorphism, hence 
\begin{eqnarray*}
\Cone(\partial_s:V_{D\alpha}DR(L_i/\mathbb A^1_k)((M,F^*,W)^{an})\otimes_{O_{S_K}}s^{\alpha+1}O\mathbb B_{dr,S_K} \\
\to V_{D\alpha}DR(L_i/\mathbb A^1_k)((M,F^{*-1},W)^{an})\otimes_{O_{S_K}}s^{\alpha}O\mathbb B_{dr,S_K})
\end{eqnarray*}
is strict for the $F$-filtration (i.e. the spectral sequence for the $F$-filtration is $E_1$-degenerate,
hence the fact that $j_i^*B^{B_{dr}}(M,F,W)$ is a quasi-isomorphism implies that $j_i^*B^{B_{dr}}(M,F,W)$
is a filtered quasi-isomorphism.
\end{proof}

\section{The geometric Mixed Hodge Modules over a field $k$ of characteristic $0$}

\subsection{The complex case where $k\subset\mathbb C$}

Let $k\subset\mathbb C$ a subfield.
For $S\in\Var(k)$, we denote by 
$\an_S:S^{an}:=S_{\mathbb C}^{an}\xrightarrow{\an_S}S_{\mathbb C}\xrightarrow{\pi_{k/\mathbb C}(S)} S$
the morphism of ringed spaces given by the analytical functor.
\begin{itemize}
\item For $(M,F)\in C_{O_Sfil}(S)$, we denote by $(M,F)^{an}:=\an_S^{*mod}(M,F)\in C_{O_Sfil}(S_{\mathbb C}^{an})$.
\item For $(M,F)\in C_{\mathcal Dfil}(S)$,
we denote by $(M,F)^{an}:=\an_S^{*mod}(M,F)\in C_{\mathcal Dfil}(S_{\mathbb C}^{an})$.
\end{itemize}
We denote for short 
\begin{eqnarray*}
DR(S):=DR(S_{\mathbb C}^{an})\circ\an_S^{*mod}:C_{\mathcal Dfil}(S)\to C_{fil}(S_{\mathbb C}^{an}), \;
M\mapsto DR(S)(M^{an})
\end{eqnarray*}
the De Rham functor.

\begin{itemize}
\item Let $S\in\SmVar(k)$.
The category $C_{\mathcal D(1,0)fil,rh}(S)\times_I D_{fil,c,k}(S_{\mathbb C}^{an})$ is the category 
\begin{itemize}
\item whose set of objects is the set of triples $\left\{((M,F,W),(K,W),\alpha)\right\}$ with 
\begin{eqnarray*}
(M,F,W)\in C_{\mathcal D(1,0)fil,rh}(S), \, (K,W)\in D_{fil,c,k}(S_{\mathbb C}^{an}), \; 
\alpha:(K,W)\otimes\mathbb C_{S_{\mathbb C}^{an}}\to DR(S)^{[-]}((M,W)^{an})
\end{eqnarray*}
where $DR(S)^{[-]}:=DR(S)^{[-]}(S_{\mathbb C}^{an}):C_{\mathcal D(1,0)fil,rh}(S_{\mathbb C}^{an})\to C_{fil}(S_{\mathbb C}^{an})$ 
is the De Rahm functor 
(recall for $S'\subset S$ a connected component of $S$ of dimension $d$, $DR(S)^{[-]}_{|S'}:=DR(S)_{|S'}[d]$)
and $\alpha$ is an morphism in $D_{fil}(S_{\mathbb C}^{an})$,
\item and whose set of morphisms are 
\begin{equation*}
\phi=(\phi_D,\phi_C,[\theta]):((M_1,F,W),(K_1,W),\alpha_1)\to((M_2,F,W),(K_2,W),\alpha_2)
\end{equation*}
where $\phi_D:(M_1,F,W)\to(M_2,F,W)$ and $\phi_C:(K_1,W)\to (K_2,W)$ are morphisms, and
\begin{eqnarray*}
\theta=(\theta^{\bullet},I(DR(S)(\phi^{an}_D))\circ I(\alpha_1),I(\alpha_2)\circ I(\phi_C\otimes I)):
I(K_1,W)\otimes\mathbb C_{S^{an}}[1]\to I(DR(S)(M^{an}_2,W)) 
\end{eqnarray*}
is an homotopy, i.e. for all $i\in\mathbb Z$, 
\begin{eqnarray*}
\theta^i\circ\partial^i-\partial^{i+1}\circ\theta^i=
(I(DR(S)(\phi^{an}_D))\circ I(\alpha_1))^i-(I(\alpha_2)\circ I(\phi_C\otimes I))^i,
\end{eqnarray*}
$I:C_{fil}(S_{\mathbb C}^{an})\to K_{fil}(S_{\mathbb C}^{an})$ being the injective resolution functor : for
$(K,W)\in C_{fil}(S_{\mathbb C}^{an})$, we take an injective resolution $k:(K,W)\to I(K,W)$ 
with $I(K,W)\in C_{fil}(S_{\mathbb C}^{an})$ which is unique modulo homotopy, 
and the class $[\theta]$ of $\theta$ does NOT depend of the injective resolution ;
in particular, we have 
\begin{equation*}
DR(S)^{[-]}(\phi^{an}_D)\circ\alpha_1=\alpha_2\circ(\phi_C\otimes I) 
\end{equation*}
in $D_{fil}(S_{\mathbb C}^{an})$ ; and for
\begin{itemize}
\item $\phi=(\phi_D,\phi_C,[\theta]):((M_1,F,W),(K_1,W),\alpha_1)\to((M_2,F,W),(K_2,W),\alpha_2)$
\item $\phi'=(\phi'_D,\phi'_C,[\theta']):((M_2,F,W),(K_2,W),\alpha_2)\to((M_3,F,W),(K_3,W),\alpha_3)$
\end{itemize}
the composition law is given by 
\begin{eqnarray*}
\phi'\circ\phi:=(\phi'_D\circ\phi_D,\phi'_C\circ\phi_C,
I(DR(S)(\phi^{'an}_D))\circ[\theta]+[\theta']\circ I(\phi_C\otimes I)[1]): \\
((M_1,F,W),(K_1,W),\alpha_1)\to((M_3,F,W),(K_3,W),\alpha_3),
\end{eqnarray*}
in particular for $((M,F,W),(K,W),\alpha)\in C_{\mathcal D(1,0)fil,rh}(S)\times_I D_{fil,c,k}(S_{\mathbb C}^{an})$,
\begin{equation*}
I_{((M,F,W),(K,W),\alpha)}=(I_M,I_K,0).
\end{equation*}
\end{itemize}
We have then the full embedding
\begin{eqnarray*}
\PSh_{\mathcal D(1,0)fil,rh}(S)\times_I P_{fil,k}(S_{\mathbb C}^{an})
\hookrightarrow C_{\mathcal D(1,0)fil,rh}(S)\times_I D_{fil,c,k}(S_{\mathbb C}^{an})
\end{eqnarray*}
where $\PSh_{\mathcal D(1,0)fil,rh}(S)\times_I P_{fil,k}(S_{\mathbb C}^{an})$ is the category
\begin{itemize}
\item whose set of objects is the set of triples $\left\{((M,F,W),(K,W),\alpha)\right\}$ with 
\begin{eqnarray*}
(M,F,W)\in\PSh_{\mathcal D(1,0)fil,rh}(S), \, (K,W)\in P_{fil,k}(S_{\mathbb C}^{an}), \; 
\alpha:(K,W)\otimes\mathbb C_{S_{\mathbb C}^{an}}\to DR(S)^{[-]}((M,W)^{an})
\end{eqnarray*}
where $DR(S)^{[-]}$ is the De Rahm functor and $\alpha$ is an isomorphism in $D_{fil}(S_{\mathbb C}^{an})$,
\item and whose set of morphisms are 
\begin{equation*}
\phi=(\phi_D,\phi_C)=(\phi_D,\phi_C,0):((M_1,F,W),(K_1,W),\alpha_1)\to((M_2,F,W),(K_2,W),\alpha_2)
\end{equation*}
where $\phi_D:(M_1,F,W)\to(M_2,F,W)$ and $\phi_C:(K_1,W)\to (K_2,W)$ are morphisms (of filtered sheaves)
and $DR(S)^{[-]}(\phi^{an}_D)\circ\alpha_1=\alpha_2\circ(\phi_C\otimes I)$ in $P_{fil,k}(S_{\mathbb C}^{an})$.
\end{itemize}
\item Let $S\in\Var(k)$. Let $S=\cup_{i\in I}S_i$ an open cover such that there
exists closed embeddings $i_i:S_i\hookrightarrow\tilde S_i$ with $\tilde S_I\in\SmVar(k)$.
The category $C_{\mathcal D(1,0)fil,rh}(S/(\tilde S_I))\times_I D_{fil,c,k}(S_{\mathbb C}^{an})$ is the category 
\begin{itemize}
\item whose set of objects is the set of triples $\left\{(((M_I,F,W),u_{IJ}),(K,W),\alpha)\right\}$ with
\begin{eqnarray*} 
((M_I,F,W),u_{IJ})\in C_{\mathcal D(1,0)fil,rh}(S/(\tilde S_I)), \, (K,W)\in D_{fil,c,k}(S_{\mathbb C}^{an}), \\ 
\alpha:T(S/(\tilde S_I))((K,W)\otimes\mathbb C_{S_{\mathbb C}^{an}})\to DR(S)^{[-]}(((M_I,W),u_{IJ})^{an})
\end{eqnarray*}
where 
\begin{equation*}
DR(S)^{[-]}:=DR(S_{\mathbb C}^{an})^{[-]}:C_{\mathcal D(1,0)fil,rh}(S^{an}/(\tilde S_{I,\mathbb C}^{an}))
\to C_{fil}(S_{\mathbb C}^{an}/(\tilde S_{I,\mathbb C}^{an})) 
\end{equation*}
is the De Rahm functor and $\alpha$ is a morphism in $D_{fil}(S_{\mathbb C}^{an}/(\tilde S_{I,\mathbb C}^{an}))$,
\item and whose set of morphisms consists of 
\begin{equation*}
\phi=(\phi_D,\phi_C,[\theta]):(((M_{1I},F,W),u_{IJ}),(K_1,W),\alpha_1)\to(((M_{2I},F,W),u_{IJ}),(K_2,W),\alpha_2)
\end{equation*}
where $\phi_D:((M_1,F,W),u_{IJ})\to((M_2,F,W),u_{IJ})$ and $\phi_C:(K_1,W)\to (K_2,W)$ 
are morphisms, and
\begin{eqnarray*}
\theta=(\theta^{\bullet},I(DR(S)(\phi^{an}_D))\circ I(\alpha_1),I(\alpha_2)\circ I(\phi_C\otimes I)): \\
I(T(S/(\tilde S_I))((K_1,W)\otimes\mathbb C_{S_{\mathbb C}^{an}}))[1]\to I(DR(S)(((M_{2I},W),u_{IJ})^{an}))  
\end{eqnarray*}
is an homotopy,  
$I:C_{fil}(S_{\mathbb C}^{an}/(\tilde S^{an}_{I,\mathbb C}))\to K_{fil}(S_{\mathbb C}^{an}/(\tilde S^{an}_{I,\mathbb C}))$
being the injective resolution functor : for
$((K_I,W),t_{IJ})\in C_{fil}(S_{\mathbb C}^{an}/(\tilde S^{an}_{I,\mathbb C}))$, 
we take an injective resolution 
\begin{equation*}
k:((K_I,W),t_{IJ})\to I((K_I,W),t_{IJ}) 
\end{equation*}
with $I((K,W),t_{IJ})\in C_{fil}(S_{\mathbb C}^{an}/(\tilde S^{an}_{I,\mathbb C}))$ 
which is unique modulo homotopy, and the class $[\theta]$ of $\theta$ does NOT depend of the injective resolution ;
in particular we have 
\begin{equation*}
DR(S)^{[-]}(\phi^{an}_D)\circ\alpha_1=\alpha_2\circ(\phi_C\otimes I) 
\end{equation*}
in $D_{fil}(S_{\mathbb C}^{an}/(\tilde S^{an}_{I,\mathbb C}))$ ; and for
\begin{itemize}
\item $\phi=(\phi_D,\phi_C,[\theta]):(((M_{1I},F,W),u_{IJ}),(K_1,W),\alpha_1)\to(((M_{2I},F,W),u_{IJ}),(K_2,W),\alpha_2)$
\item $\phi'=(\phi'_D,\phi'_C,[\theta']):(((M_{2I},F,W),u_{IJ}),(K_2,W),\alpha_2)\to(((M_{3I},F,W),u_{IJ}),(K_3,W),\alpha_3)$
\end{itemize}
the composition law is given by 
\begin{eqnarray*}
\phi'\circ\phi:=(\phi'_D\circ\phi_D,\phi'_C\circ\phi_C,
I(DR(S)(\phi^{'an}_D))\circ[\theta]+[\theta']\circ I(\phi_C\otimes I)[1]): \\
(((M_{1I},F,W),u_{IJ}),(K_1,W),\alpha_1)\to(((M_{3I},F,W),u_{IJ}),(K_3,W),\alpha_3)
\end{eqnarray*}
in particular for 
$(((M_I,F,W),u_{IJ}),(K,W),\alpha)\in C_{\mathcal D(1,0)fil,rh}(S/(\tilde S_I))\times_I D_{fil,c,k}(S_{\mathbb C}^{an})$,
\begin{equation*}
I_{(((M_I,F,W),u_{IJ}),(K,W),\alpha)}=((I_{M_I}),I_K,0).
\end{equation*}
\end{itemize}
We have then full embeddings
\begin{eqnarray*}
\PSh^0_{\mathcal D(1,0)fil,rh}(S/(\tilde S_I))\times_I P_{fil,k}(S_{\mathbb C}^{an})
\hookrightarrow C^0_{\mathcal D(1,0)fil,rh}(S/(\tilde S_I))\times_I D_{fil,c,k}(S_{\mathbb C}^{an}) \\
\xrightarrow{\iota^0_{S/\tilde S_I}} C_{\mathcal D(1,0)fil,rh}(S/(\tilde S_I))^0\times_I D_{fil,c,k}(S_{\mathbb C}^{an})
\hookrightarrow C_{\mathcal D(1,0)fil,rh}(S/(\tilde S_I))\times_I D_{fil,c,k}(S_{\mathbb C}^{an})
\end{eqnarray*}
where $\PSh^0_{\mathcal D(1,0)fil,rh}(S/(\tilde S_I))\times_I P_{fil,k}(S_{\mathbb C}^{an})$ is the category 
\begin{itemize}
\item whose set of objects is the set of triples $\left\{(((M_I,F,W),u_{IJ}),(K,W),\alpha)\right\}$ with
\begin{eqnarray*} 
((M_I,F,W),u_{IJ})\in\PSh^0_{\mathcal D(1,0)fil,rh}(S/(\tilde S_I)), \, (K,W)\in P_{fil,k}(S_{\mathbb C}^{an}), \\ 
\alpha:T(S/(\tilde S_I))((K,W)\otimes\mathbb C_{S_{\mathbb C}^{an}})\to DR(S)^{[-]}(((M_I,W),u_{IJ})^{an})
\end{eqnarray*}
where $DR(S)^{[-]}$ is the De Rahm functor and $\alpha$ is an isomorphism in 
$D_{fil}(S_{\mathbb C}^{an}/(\tilde S_{I,\mathbb C}^{an}))$,
\item and whose set of morphisms are 
\begin{equation*}
\phi=(\phi_D,\phi_C)=(\phi_D,\phi_C,0):
(((M_{1I},F,W),u_{IJ}),(K_1,W),\alpha_1)\to(((M_{2I},F,W),u_{IJ}),(K_2,W),\alpha_2)
\end{equation*}
where $\phi_D:((M_1,F,W),u_{IJ})\to((M_2,F,W),u_{IJ})$ and $\phi_C:(K_1,W)\to (K_2,W)$ are morphisms (of filtered sheaves)
such that $\phi_D^{an}\circ\alpha_1=\alpha_2\circ(\phi_C\otimes I)$ in $P_{fil,k}(S_{\mathbb C}^{an})$.
\end{itemize}
\end{itemize}
Moreover,
\begin{itemize}
\item For 
$(((M_{I},F,W),u_{IJ}),(K,W),\alpha)\in C_{\mathcal D(1,0)fil,rh}(S/(\tilde S_I))\times_I D_{fil,c,k}(S_{\mathbb C}^{an})$, 
we set
\begin{equation*}
(((M_{I},F,W),u_{IJ}),(K,W),\alpha)[1]:=(((M_{I},F,W),u_{IJ})[1],(K,W)[1],\alpha[1]).
\end{equation*}
\item For 
\begin{equation*}
\phi=(\phi_D,\phi_C,[\theta]):(((M_{1I},F,W),u_{IJ}),(K_1,W),\alpha_1)\to(((M_{2I},F,W),u_{IJ}),(K_2,W),\alpha_2)
\end{equation*}
a morphism in $C_{\mathcal D(1,0)fil,rh}(S/(\tilde S_I))\times_I D_{fil,c,k}(S_{\mathbb C}^{an})$, 
we set (see \cite{CG} definition 3.12)
\begin{eqnarray*}
\Cone(\phi):=(\Cone(\phi_D),\Cone(\phi_C),((\alpha_1,\theta),(\alpha_2,0)))
\in D_{\mathcal D(1,0)fil,rh}(S/(\tilde S_I))\times_I D_{fil,c,k}(S_{\mathbb C}^{an}),
\end{eqnarray*}
$((\alpha_1,\theta),(\alpha_2,0))$ being the matrix given by the composition law, together with the canonical maps
\begin{itemize}
\item $c_1(-)=(c_1(\phi_D),c_1(\phi_C),0):(((M_{2I},F,W),u_{IJ}),(K_2,W),\alpha_2)\to\Cone(\phi)$
\item $c_2(-)=(c_2(\phi_D),c_2(\phi_C),0):\Cone(\phi)\to (((M_{1I},F,W),u_{IJ}),(K_1,W),\alpha_1)[1]$.
\end{itemize}
\end{itemize}

\begin{rem}\label{CGrem}
By \cite{CG} theorem 3.25, if 
\begin{equation*}
\phi=(\phi_D,\phi_C,[\theta]):(((M_1,F,W),u_{IJ}),(K_1,W),\alpha_1)\to(((M_2,F,W),u_{IJ}),(K_2,W),\alpha_2)
\end{equation*}
is a morphism in $C_{\mathcal D(1,0)fil,rh}(S/(\tilde S_I))\times_I D_{fil,c,k}(S_{\mathbb C}^{an})$
such that $\phi_D$ is a Zariski local equivalence and $\phi_C$ is an isomorphism then $\phi$ is an isomorphism.
\end{rem}

We get from \cite{B4} the following definition :

\begin{defi}\label{falpha}
\begin{itemize}
\item[(i)] Let $f:X\to S$ a morphism with $S,X\in\SmVar(k)$. 
Let $f:X\xrightarrow{j}\bar X\xrightarrow{\bar f}S$ a compactification of $f$ with $\bar X\in\SmVar(k)$
and $j$ the open embedding. Denote $Z:=\bar X\backslash X=\cap_iZ_i$ with $Z_i\subset\bar X$ (Cartier) divisor. Let 
\begin{equation*}
\alpha:(K,W)\otimes\mathbb C_{X_{\mathbb C}^{an}}\to DR(X)((M,W)^{an}) 
\end{equation*}
a morphism in $D_{fil}(X_{\mathbb C}^{an})$, with 
\begin{equation*}
(M,F,W)\in C(DRM(X)), \; (K,W)\in D_{fil,c,k}(X_{\mathbb C}^{an})^{ad,(Z_i)}. 
\end{equation*}
We then consider the maps in $D_{fil}(S_{\mathbb C}^{an})$
\begin{eqnarray*}
f_*\alpha:Rf_{*w}(K,W)\otimes\mathbb C_{S_{\mathbb C}^{an}}:=
R\bar f_*Rj_{*w}(K,W)\otimes\mathbb C_{S_{\mathbb C}^{an}} \\
\xrightarrow{R\bar f_*j_*\alpha}R\bar f_*Rj_{*w}DR(X)((M,W)^{an}) 
\xrightarrow{T^w(j,\otimes)(-)^{-1}}R\bar f_*DR(\bar X)(j_{*Hdg}(M,W)^{an}) \\
\xrightarrow{T(\bar f,DR)(-)^{-1}}
DR(S)((\int_{\bar f}(j_{*Hdg}(M,W))^{an})=DR(S)((\int^{Hdg}_f(M,W))^{an})
\end{eqnarray*}
and
\begin{eqnarray*}
f_!\alpha:Rf_{!w}(K,W)\otimes\mathbb C_{S_{\mathbb C}^{an}}:=
R\bar f_*Rj_{!w}(K,W)\otimes\mathbb C_{S_{\mathbb C}^{an}} \\
\xrightarrow{R\bar f_*j_!\alpha}R\bar f_*Rj_{!w}DR(X)((M,W)^{an}) 
\xrightarrow{\mathbb DT^w(j,\otimes)(-)}R\bar f_*DR(\bar X)(j_{!Hdg}(M,W)^{an}) \\
\xrightarrow{T(\bar f,DR)(-)^{-1}}
DR(S)((\int_{\bar f}(j_{!Hdg}(M,W))^{an})=DR(S)((\int^{Hdg}_{f!}(M,W))^{an}),
\end{eqnarray*}
see definition \ref{fw} and definition \ref{DHdg} .
\item[(ii)] Let $f:X\to S$ a morphism with $S,X\in\QPVar(k)$.  
Consider a factorization $f:X\xrightarrow{l}Y\times S\xrightarrow{p}S$
with $Y\in\SmVar(k)$, $l$ a closed embedding and $p_S$ the projection.
Let $\bar Y\in\PSmVar(k)$ a smooth compactification of $Y$ with $j:Y\hookrightarrow\bar Y$ the open embedding.
Then $\bar f:\bar X\xrightarrow{\bar l}\bar Y\times_S\xrightarrow{\bar p}S$ is a compactification of $f$,
with $\bar X\subset\bar Y\times S$ the closure of $X$ and $\bar l$ the closed embedding.
Denote $Z:=\bar X\backslash X=\cap_iZ_i$ with $Z_i\subset\bar X$ Cartier divisors.
Let $S=\cup_i S_i$ an open affine cover and 
$i_i:S_i\hookrightarrow\tilde S_i$ closed embedding with $\tilde S_i\in\SmVar(k)$. Let 
\begin{equation*}
\alpha:T(X/(Y\times\tilde S_I))((K,W)\otimes\mathbb C_{X_{\mathbb C}^{an}})\to DR(X)((M_I,W)^{an},u^{an}_{IJ}) 
\end{equation*}
a morphism in $D_{fil}(X_{\mathbb C}^{an}/({Y\times\tilde S_I}_{\mathbb C}^{an}))$, with
\begin{equation*}
((M_I,W),u_{IJ})\in C(DRM(X))\subset C_{\mathcal D0fil,rh}(X/(Y\times\tilde S_I)), \; 
(K,W)\in D_{fil,c,k}(X_{\mathbb C}^{an})^{ad,(Z_i)}. 
\end{equation*}
We then consider the maps in $D_{fil}(S_{\mathbb C}^{an}/(\tilde S_{I,\mathbb C}^{an}))$
\begin{eqnarray*}
f_*\alpha=f_*(\alpha):T(S/\tilde S_I)(Rf_{*w}(K,W))\otimes\mathbb C_{S_{\mathbb C}^{an}} \\
\xrightarrow{:=}T(S/\tilde S_I)(R\bar p_*(I\times j)_{*w}(K,W))\otimes\mathbb C_{S_{\mathbb C}^{an}}
\xrightarrow{=}R\bar p_*(I\times j)_{*w}T(X/(Y\times\tilde S_I))((K,W)\otimes\mathbb C_{X_{\mathbb C}^{an}}) \\
\xrightarrow{Rp_*\alpha}R\bar p_*(I\times j)_{*w}DR(X)(((M_I,W),u_{IJ})^{an}) \\ 
\xrightarrow{(T^w(I\times j,\otimes)(-))}R\bar p_*DR(X)((I\times j)_{*Hdg}((M_I,W),u_{IJ})^{an}) \\
\xrightarrow{T(\bar f,DR)(-)}
DR(S)((\int_{\bar f}(I\times j)_{*Hdg}((M_I,W),u_{IJ}))^{an})=DR(S)((\int^{Hdg}_f((M_I,W),u_{IJ}))^{an})
\end{eqnarray*}
and
\begin{eqnarray*}
f_!\alpha=f_!(\alpha):T(S/\tilde S_I)(Rf_{!w}(K,W))\otimes\mathbb C_{S_{\mathbb C}^{an}} \\
\xrightarrow{:=}T(S/\tilde S_I)(R\bar p_*(I\times j)_{*w}(K,W))\otimes\mathbb C_{S_{\mathbb C}^{an}}
\xrightarrow{=}R\bar p_*(I\times j)_{!w}T(X/(Y\times\tilde S_I))((K,W)\otimes\mathbb C_{X_{\mathbb C}^{an}}) \\
\xrightarrow{R\bar p_*\mathbb D^vR(I\times j)_*\mathbb D^v\alpha}
R\bar p_*(I\times j)_{!w}DR(X)(((M_I,W),u_{IJ})^{an}) \\
\xrightarrow{T(D,DR)(-)\circ(\mathbb DT^w(I\times j,\otimes)(-))\circ T(D,DR)(-)}
R\bar p_*DR(X)((I\times j)_{!Hdg}((M_I,W),u_{IJ})^{an}) \\
\xrightarrow{T(\bar f,DR)(-)}
DR(S)((\int_{\bar f}(I\times j)_{!Hdg}((M_I,W),u_{IJ}))^{an})=DR(S)((\int^{Hdg}_{f!}((M_I,W),u_{IJ}))^{an}), 
\end{eqnarray*}
see definition \ref{fw} and definition \ref{DHdgsing} .
\item[(iii)] Let $l:S^o\hookrightarrow S$ an open embedding with $S\in\Var(k)$ and denote $Z=S\backslash S^o=\cap_i Z_i$.
with $Z_i\subset S$ Cartier divisors.
Let $S=\cup_i S_i$ an open affine cover and 
$i_i:S_i\hookrightarrow\tilde S_i$ closed embedding with $\tilde S_i\in\SmVar(k)$. 
Let $l_I:\tilde S^o_I\hookrightarrow\tilde S_I$ open embeddings such that $\tilde S^o_I\cap S=S^o\cap S_I$. Let 
\begin{equation*}
\alpha:T(S/(\tilde S_I))((K,W)\otimes\mathbb C_{S_{\mathbb C}^{an}})\to DR(S)(((M_I,W),u_{IJ})^{an}) 
\end{equation*}
a morphism in $D_{fil}(S_{\mathbb C}^{an}/(\tilde S_{I,\mathbb C}^{an}))$, with
\begin{equation*}
((M_I,W),u_{IJ})\in C(DRM(S))\subset C_{\mathcal D0fil,rh}(S/(\tilde S_I)), \; (K,W)\in D_{fil,c,k}(S_{\mathbb C}^{an}). 
\end{equation*}
We then consider the maps in $D_{fil}(S_{\mathbb C}^{an}/(\tilde S_{I,\mathbb C}^{an}))$
\begin{eqnarray*}
\Gamma_Z(\alpha):T(S/(\tilde S_I))(\Gamma^w_Z(K,W)\otimes\mathbb C_{S_{\mathbb C}^{an}})
\xrightarrow{=}\Gamma^w_ZT(S/(\tilde S_I))((K,W)\otimes\mathbb C_{S_{\mathbb C}^{an}}) \\
\xrightarrow{R\Gamma_Z\alpha}\Gamma^w_ZDR(S)(((M_I,W),u_{IJ})^{an}) \\
\xrightarrow{(T(\gamma_Z,DR)((M_I,W),u_{IJ})^{-1}:=((I,T^w(l_I,\otimes)(M_I,W))\circ (T^w(an,\otimes)(M_I,W)))^{-1}} 
DR(S)((\Gamma^{Hdg}_Z((M_I,W),u_{IJ}))^{an})
\end{eqnarray*}
and
\begin{eqnarray*}
\Gamma^{\vee}_Z(\alpha):T(S/(\tilde S_I))(\Gamma^{\vee,w}_Z(K,W)\otimes\mathbb C_{S_{\mathbb C}^{an}})
\xrightarrow{=}\Gamma^{\vee}_ZT(S/(\tilde S_I))((K,W)\otimes\mathbb C_{S_{\mathbb C}^{an}}) \\
\xrightarrow{\Gamma^{\vee}_Z\alpha}\Gamma^{\vee,w}_ZDR(S)(((M_I,W),u_{IJ})^{an}) \\
\xrightarrow{T(\gamma^{\vee}_Z,DR)((M_I,W),u_{IJ}):=
(\mathbb D(I,T^w(l_I,\otimes)(\mathbb D(M_I,W)))\circ (\mathbb DT^w(an,\otimes)(\mathbb D(M_I,W))))} 
DR(S)((\Gamma^{\vee,w}_Z((M_I,W),u_{IJ}))^{an}),
\end{eqnarray*}
see definition \ref{gammaw} and definition \ref{gammaHdg}.
\item[(iv)] Let $f:X\to S$ a morphism with $S,X\in\QPVar(k)$. 
Consider a factorization $f:X\hookrightarrow Y\times S\xrightarrow{p}S$ with $Y\in\SmVar(k)$. 
Let $S=\cup_i S_i$ an open affine cover and 
$i_i:S_i\hookrightarrow\tilde S_i$ closed embedding with $\tilde S_i\in\SmVar(k)$. Let 
\begin{equation*}
\alpha:T(S/(\tilde S_I))((K,W)\otimes\mathbb C_{S_{\mathbb C}^{an}})\to DR(S)(((M_I,W),u_{IJ})^{an}) 
\end{equation*}
a morphism in $D_{fil}(S_{\mathbb C}^{an}/(\tilde S_{I,\mathbb C}^{an}))$, with
\begin{equation*}
((M_I,W),u_{IJ})\in C(DRM(S))\subset C_{\mathcal D0fil,rh}(S/(\tilde S_I)), \; 
(K,W)\in D_{fil,c,k}(S_{\mathbb C}^{an})^{ad,(\Gamma_{f,i})}. 
\end{equation*}
We then consider, see (iii), the maps in $D_{fil}(X_{\mathbb C}^{an}/(Y\times\tilde S_{I,\mathbb C})^{an})$
\begin{eqnarray*}
f^!\alpha=f^!(\alpha):T(X/(Y\times\tilde S_I))(f^{!w}(K,W)\otimes\mathbb C_{X_{\mathbb C}^{an}})
\xrightarrow{:=}T(X/(Y\times\tilde S_I))(\Gamma^w_Xp^*(K,W))\otimes\mathbb C_{X_{\mathbb C}^{an}}) \\
\xrightarrow{=}
(\Gamma^w_Xp_{\tilde S_I}^*T(S/\tilde S_I)((K,W)_I\otimes\mathbb C_{S_{\mathbb C}^{an}}),\Gamma^w_Xp^*T(D_{IJ})(-)) \\
\xrightarrow{R\Gamma_Xp^*\alpha}
(\Gamma^w_Xp_{\tilde S_I}^*DR(S)((M_I,W),u_{IJ})_I,\Gamma^w_Xp^*DR(u_{IJ})) \\
\xrightarrow{T^!(f,DR)(-):=T(\gamma_X,DR)(-)^{-1}\circ T^!(p,DR)(-)}DR(X)(f^{*mod}_{Hdg}((M_I,W),u_{IJ})^{an}) 
\end{eqnarray*}
and
\begin{eqnarray*}
f^*\alpha=f^*(\alpha):T(X/(Y\times\tilde S_I))(f^{*w}(K,W)\otimes\mathbb C_{X_{\mathbb C}^{an}})
\xrightarrow{:=}T(X/(Y\times\tilde S_I))(\Gamma^{\vee,w}_Xp^*(K,W)\otimes\mathbb C_{X_{\mathbb C}^{an}}) \\
\xrightarrow{=}
(\Gamma^{\vee,w}_Xp_{\tilde S_I}^*T(S/(\tilde S_I))((K,W)\otimes\mathbb C_{S_{\mathbb C}^{an}})_I,
\Gamma_X^{\vee,w}p^*T(D_{IJ})(-)) \\
\xrightarrow{\Gamma^{\vee}_Xp^*\alpha}
(\Gamma^{\vee,w}_Xp_{\tilde S_I}^*DR(S)((M_I,W),u_{IJ})_I,\Gamma_X^{\vee,w}p^*DR(u_{IJ})) \\
\xrightarrow{T^*(f,DR)(-):=T(\gamma_X^{\vee},DR)(-)\circ T^*(p,DR)(-)}DR(X)(f_{Hdg}^{\hat*mod}((M_I,W),u_{IJ})^{an}) 
\end{eqnarray*}
\item[(v)] Let $S\in\Var(k)$.Let $S=\cup_i S_i$ an open affine cover and 
$i_i:S_i\hookrightarrow\tilde S_i$ closed embedding with $\tilde S_i\in\SmVar(k)$. Let 
\begin{eqnarray*}
\alpha:T(S/(\tilde S_I))((K,W)\otimes\mathbb C_{S_{\mathbb C}^{an}})\to DR(S)(((M_I,W),u_{IJ})^{an}), \\
\alpha':T(S/(\tilde S_I))((K',W)\otimes\mathbb C_{S_{\mathbb C}^{an}})\to DR(S)(((M'_I,W),v_{IJ})^{an}) 
\end{eqnarray*}
two morphism in $D_{fil}(S_{\mathbb C}^{an}/(\tilde S_{I,\mathbb C}^{an}))$, with
\begin{equation*}
((M_I,W),u_{IJ}),((M'_I,W),u_{IJ})\in C(DRM(S))\subset C_{\mathcal D0fil,rh}(S/(\tilde S_I)), \; 
(K,W),(K',W)\in D_{fil}(S_{\mathbb C}^{an}). 
\end{equation*}
We then consider the map in $D_{fil}(S_{\mathbb C}^{an}/(\tilde S_{I,\mathbb C}^{an}))$
\begin{eqnarray*}
\alpha\otimes\alpha':T(S/(\tilde S_I))((K,W)\otimes^{L,w}(K',W)\otimes\mathbb C_{S_{\mathbb C}^{an}}) \\
\xrightarrow{=}T(S/(\tilde S_I))((K,W)\otimes\mathbb C_{S_{\mathbb C}^{an}})
\otimes T(S/(\tilde S_I))((K',W)\otimes\mathbb C_{S_{\mathbb C}^{an}}) \\
\xrightarrow{\alpha\otimes\alpha'}DR(S)(((M_I,W),u_{IJ})^{an})\otimes^{L,w} DR(S)(((M'_I,W),v_{IJ})^{an}) \\
\xrightarrow{T(\otimes^w,DR)(-,-)}DR(S)((((M_I,W),u_{IJ})\otimes^{L,w}_{O_S}((M'_I,W),v_{IJ}))^{an}) \\
\xrightarrow{=}DR(S)((((M_I,W),u_{IJ})\otimes^{Hdg}_{O_S}((M'_I,W),v_{IJ}))^{an})
\end{eqnarray*}
with $T(\otimes^w,DR)(-,-):=T(\gamma_S,DR)(-)\circ T(\otimes,DR)(-,-)\circ(T(p_1,DR)(-)\otimes T(p_2,DR)(-))$.
\end{itemize}
\end{defi}

\begin{defi}\label{DHdgpsialpha}
Let $S\in\SmVar(k)$.
Let $D=V(s)\subset S$ a divisor with $s\in\Gamma(S,L)$ and $L$ a line bundle ($S$ being smooth, $D$ is Cartier).
For $\mathcal M=((M,F,W),(K,W),\alpha)\in\PSh_{\mathcal D(1,0)fil,rh}(S)\times_IP_{fil,k}(S_{\mathbb C}^{an})$, 
we then define, using definition \ref{DHdgpsi} and theorem \ref{phipsi0thm},
\begin{itemize}
\item the nearby cycle functor
\begin{equation*}
\psi_D((M,F,W),(K,W),\alpha):=(\psi_D(M,F,W),\psi_D(K,W)[-1],\psi_D\alpha)
\in\PSh_{\mathcal D(1,0)fil,rh}(S)\times_IP_{fil,k}(S_{\mathbb C}^{an}),
\end{equation*}
with $\psi_D\alpha:=T(\psi_D,DR)(M)\circ\psi_D(\alpha)$.
\item the vanishing cycle functor
\begin{equation*}
\phi_D((M,F,W),(K,W),\alpha):=(\phi_D(M,F,W),\phi_D(K,W)[-1],\phi_D\alpha)
\in\PSh_{\mathcal D(1,0)fil,rh}(S)\times_IP_{fil,k}(S_{\mathbb C}^{an}),
\end{equation*}
with $\phi_D\alpha:=T(\phi_D,DR)(M)\circ\phi_D(\alpha)$.
\item the canonical maps in $\PSh_{\mathcal D(1,0)fil,rh}(S)\times_IP_{fil,k}(S_{\mathbb C}^{an})$
\begin{eqnarray*}
can(\mathcal M):=(can(M,F,W),can(K,W)):\psi_D((M,F,W),(K,W),\alpha)\to\phi_D((M,F,W),(K,W),\alpha)(-), \\
var(\mathcal M):=(var(M,F,W),var(K,W):\phi_D((M,F,W),(K,W),\alpha)\to\psi_D((M,F,W),(K,W),\alpha).
\end{eqnarray*}
\end{itemize}
\end{defi}

\begin{prop}\label{phipsigmHdgpropalpha}
Let $S\in\SmVar(k)$. Let $D=V(s)\subset S$ a (Cartier divisor). 
Consider a composition of proper morphisms $(f:X=X_r\xrightarrow{f_r}X_{r-1}\xrightarrow{f_1}X_0=S)\in\SmVar(k)$ and 
\begin{eqnarray*}
(M,F)=H^{n_0}\int_{f_1}\cdots H^{n_r}\int_{f_r}((O_X,F_b), 
H^{n_0}Rf_{1*}\cdots H^{n_r}Rf_{r*}\mathbb Z_X, \\
H^{n_0}f_{1*}\circ\cdots\circ H^{n_r}f_{r*}\alpha(X)) 
\in\PSh_{\mathcal Dfil,rh}(S)\times_IP_k(S_{\mathbb C}^{an}).
\end{eqnarray*} 
Then,
\begin{eqnarray*}
\psi_D(M,F)=H^{n_0}\int_{f_1}\cdots H^{n_r}\int_{f_r}(\psi_{f^{-1}(D)}(O_X,F_b)), 
H^{n_0}Rf_{1*}\cdots H^{n_r}Rf_{r*}\psi_{f^{-1}(D)}\mathbb Z_X, \\
H^{n_0}f_{1*}\circ\cdots\circ H^{n_r}f_{r*}\psi_{f^{-1}(D)}\alpha(X)) 
\in\PSh_{\mathcal Dfil,rh}(S)\times_IP_k(S_{\mathbb C}^{an})
\end{eqnarray*}
and
\begin{eqnarray*}
\phi_D(M,F)=H^{n_0}\int_{f_1}\cdots H^{n_r}\int_{f_r}(\psi_{f^{-1}(D)}(O_X,F_b)), 
H^{n_0}Rf_{1*}\cdots H^{n_r}Rf_{r*}\phi_{f^{-1}(D)}\mathbb Z_X, \\
H^{n_0}f_{1*}\circ\cdots\circ H^{n_r}f_{r*}\phi_{f^{-1}(D)}\alpha(X)) 
\in\PSh_{\mathcal Dfil,rh}(S)\times_IP_k(S_{\mathbb C}^{an}). 
\end{eqnarray*}
\end{prop}

\begin{proof}
Immediate from definition.
\end{proof}

We now come to the main definition of this section :

\begin{defi}\label{MHMgmkdef}
Let $k\subset\mathbb C$ a subfield.
\begin{itemize}
\item[(i0)]Let $S\in\Var(k)$.
Take an open cover $S=\cup_iS_i$ such that there are closed embedding $S_I\hookrightarrow\tilde S_I$ with $S_I\in\SmVar(k)$.
The category of mixed hodge modules over $k$ is the full subcategory
\begin{eqnarray*}
MHM_{k,\mathbb C}(S)\subset DRM(S)\times_IP_{fil,k}(S_{\mathbb C}^{an})
\subset\PSh_{\mathcal D(1,0)fil,rh}(S/(\tilde S_I))\times_IP_{fil,k}(S_{\mathbb C}^{an})
\end{eqnarray*}
whose object consists of $(((M_I,F,W),u_{IJ}),(K,W),\alpha)\in DRM(S)\times_IP_{fil,k}(S_{\mathbb C}^{an})$ such that 
\begin{equation*}
((\pi_{k/\mathbb C}^{*mod}(M_I,F,W),u_{IJ}),(K,W),\alpha)\in MHM(S_{\mathbb C})
\end{equation*}
where
\begin{itemize}
\item $\pi_{k/\mathbb C}:=\pi_{k/\mathbb C}(S):S_{\mathbb C}\to S$ is the projection (see section 2),
\item $DRM(S)$ is the category of de Rham modules introduced in section 5 definition \ref{DRMdef},
\item $MHM(S_{\mathbb C})$ is the category of mixed hodge modules on $S_{\mathbb C}$ introduced by Saito (\cite{Saito}).
\end{itemize} 
\item[(i)]Let $S\in\SmVar(k)$. We denote by 
\begin{eqnarray*}
HM_{gm,k,\mathbb C}(S):=<(H^{n_1}\int_{f_1}\cdots H^{n_r}\int_{f_r}(O_X,F_b)(d),
R^{n_1}f_{1*}\cdots R^{n_r}f_{r*}\mathbb Z_{X_{\mathbb C}^{an}},H^{n_1}f_{1*}\cdots H^{n_r}f_{r*}\alpha(X)), \\
(f:X=X_r\xrightarrow{f_r}X_{r-1}\to\cdots\xrightarrow{f_1}X_0=S)\in\SmVar(k), \; 
\mbox{proper}, \; n_1,\ldots n_r,d\in\mathbb Z> \\
\subset PDRM(S)\times_I P_{fil,k}(S_{\mathbb C}^{an})
\subset\PSh_{\mathcal Dfil,rh}(S)\times_I P_k(S_{\mathbb C}^{an})
\end{eqnarray*}
the full abelian subcategory, where $<,>$ means generated by and $(-)$ is the shift of the filtration, 
\begin{equation*}
\alpha(X):\mathbb C_{X_{\mathbb C}^{an}}\hookrightarrow DR(X)(O^{an}_X) 
\end{equation*}
is the inclusion quasi-isomorphism in $C(X_{\mathbb C}^{an})$, and we use definition \ref{falpha}. 
We have by proposition \ref{phipsigmHdgpropalpha} for $((M,F),K,\alpha)\in HM_{gm,k,\mathbb C}(S)$,
\begin{eqnarray*}
\Gr_k^W\psi_D((M,F),K,\alpha):=\Gr_k^W\psi_D(M,F),\Gr_k^W\psi_DK,\Gr_k^W\psi_D\alpha)\in HM_{gm,k,\mathbb C}(S) .
\end{eqnarray*}
and
\begin{eqnarray*}
\Gr_k^W\psi_D((M,F),K,\alpha):=\Gr_k^W\psi_D(M,F),\Gr_k^W\psi_DK,\Gr_k^W\psi_D\alpha)\in HM_{gm,k,\mathbb C}(S) .
\end{eqnarray*}
for all $k\in\mathbb Z$. 
We have by theorem \ref{Sa12}, for $S\in\SmVar(k)$, $HM_{gm,k,\mathbb C}(S)\subset HM(S_{\mathbb C})$
which consists of geometric Hodge module defined over $k$.
\item[(i)']Let $S\in\Var(k)$. Let $S=\cup_{i\in I}S_i$ an open cover such that there
exists closed embeddings $i_i:S_i\hookrightarrow\tilde S_i$ with $\tilde S_I\in\SmVar(k)$. We denote by 
\begin{eqnarray*}
HM_{gm,k,\mathbb C}(S):= 
<(H^{n_1}\int_{p_1}\cdots H^{n_r}\int_{p_r}(\Gamma^{\vee,Hdg}_{X_I}(O_{Y\times\tilde X_{r-1,I}},F_b),x_{IJ})(d), \\
R^{n_1}p_{1*}\cdots R^{n_r}p_{r*}T(X/(Y_r\times\tilde X_{r-1,I}))(\mathbb Z_{X_{\mathbb C}^{an}}), 
H^{n_1}p_{1*}\cdots H^{n_r}p_{r*}\alpha(X)), \\
(f:X=X_r\xrightarrow{f_r}X_{r-1}\to\cdots\xrightarrow{f_1}X_0=S)\in\Var(k), \; > \\
\subset PDRM(S)\times_I P_k(S_{\mathbb C}^{an})
\subset\PSh_{\mathcal Dfil,rh}(S/(\tilde S_I))\times_I P_k(S_{\mathbb C}^{an})
\end{eqnarray*}
the full abelian subcategory, where $<,>$ means generated by and $(-)$ is the shift of the filtration, 
$f_i:X_i\hookrightarrow Y_i\times X_{i-1}\xrightarrow{p_i}X_{i-1}$ proper, $Y_i\in\PSmVar(k)$, $X_i$ smooth, 
\begin{eqnarray*}
\alpha(X):=(\Gamma_{X_I}^{\vee}\alpha(Y_r\times\tilde X_{r-1,I})):
T(X/(Y_r\times\tilde X_{r-1,I}))(\mathbb C_{X^{an}}):=
(\Gamma^{\vee}_{X_I}\mathbb C_{(Y\times\tilde X_{r-1,I})^{an}_{\mathbb C}},t_{IJ}) \\
\xrightarrow{\sim}DR(X)(o_F(\Gamma^{\vee,Hdg}_{X_I}(O_{Y\times\tilde X_{r-1,I}},F_b),x_{IJ})^{an}).
\end{eqnarray*}
is the canonical isomorphism in $D(X_{\mathbb C}^{an}/(Y_r\times\tilde X_{r-1,I})_{\mathbb C}^{an})$, 
and we use definition \ref{falpha}. 
Note that if $S$ is smooth then this definition of $HM_{gm,k,\mathbb C}(S)$ agree with the one given in (i).
We have by theorem \ref{Sa12}, for $S\in\Var(k)$, $HM_{gm,k,\mathbb C}(S)\subset HM(S_{\mathbb C})$
which consists of geometric Hodge module defined over $k$.
\item[(ii)]Let $S\in\Var(k)$.
Take an open cover $S=\cup_iS_i$ such that there are closed embedding $S_I\hookrightarrow\tilde S_I$ with $S_I\in\SmVar(k)$.
We define using the pure case (i) and (i)' the full subcategory of geometric mixed Hodge modules defined over $k$
\begin{eqnarray*}
MHM_{gm,k,\mathbb C}(S):= \\
\left\{(((M_I,F,W),u_{IJ}),(K,W),\alpha), \; \mbox{s.t.} \; 
\Gr^W_k(((M_I,F,W),u_{IJ}),(K,W),\alpha)\in HM_{gm,k,\mathbb C}(S)\right\} \\
\subset DRM(S)\times_IP_{fil,k}(S_{\mathbb C}^{an})
\subset\PSh_{\mathcal D(1,0)fil,rh}(S/(\tilde S_I))\times_IP_{fil,k}(S_{\mathbb C}^{an})
\end{eqnarray*}
whose object consists of 
$(((M_I,F,W),u_{IJ}),(K,W),\alpha)\in DRM(S)\times_IP_{fil,k}(S_{\mathbb C}^{an})$ 
such that 
\begin{equation*}
\Gr^W_k(((M_I,F,W),u_{IJ}),(K,W),\alpha):=(\Gr^W_k((M_I,F),u_{IJ}),\Gr^W_kK,\Gr^W_k\alpha)
\in HM_{gm,k,\mathbb C}(S).
\end{equation*}
We set 
\begin{eqnarray*}
\mathbb Q_S^{Hdg}:=((\Gamma_{S_I}^{\vee,Hdg}(O_{\tilde S_I},F_b),x_{IJ}),\mathbb Q_{S^{an}}^w,\alpha(S))
\in C(MHM_{gm,k,\mathbb C}(S))
\end{eqnarray*}
where $\mathbb Q_{S^{an}}^w\in C(P_{fil,k}(S_{\mathbb C}^{an}))$ 
is such that $j_I^*\mathbb Q_{S^{an}}^w=i_I^*\Gamma^{\vee,w}_{S_I}\mathbb Q_{\tilde S_I^{an}}$ and
\begin{eqnarray*}
\alpha(S):=(\Gamma_{S_I}^{\vee}\alpha(\tilde S_I)):
T(S/(\tilde S_I))((\mathbb Q_{S^{an}}^w)\otimes\mathbb C_{S^{an}})
\xrightarrow{=}(\Gamma^{\vee,w}_{S_I}(\mathbb C_{\tilde S_I^{an}}),t_{IJ})
\xrightarrow{\sim}DR(S)(o_F(\Gamma_{S_I}^{\vee,Hdg}(O_{\tilde S_I},F_b),x_{IJ})).
\end{eqnarray*}
For $S\in\SmVar(k)$ and $D=V(s)\subset S$ a (Cartier) divisor, 
we have for $((M,F,W),(K,W),\alpha)\in MHM_{gm,k,\mathbb C}(S)$, using theorem \ref{HSk}, 
\begin{equation*}
\psi_D((M,F,W),(K,W),\alpha),\phi_D((M,F,W),(K,W),\alpha)\in MHM_{gm,k,\mathbb C}(S),
\end{equation*}
by the pure case (c.f. (i) and proposition \ref{phipsigmHdgpropalpha}) and the strictness of the $V$-filtration.
\end{itemize}
For $k\subset\mathbb C$ and $S\in\Var(k)$, we have by theorem \ref{Sa12} 
\begin{equation*}
MHM_{gm,k,\mathbb C}(S)\subset MHM_{k,\mathbb C}(S)\subset MHM(S_{\mathbb C}).
\end{equation*}
For $S\in\Var(k)$ we get $D(MHM_{gm,k,\mathbb C}(S)):=\Ho_{(zar,usu)}(C(MHM_{gm,k,\mathbb C}(S)))$ 
after localization with Zariski local equivalence and usu local equivalence.
\end{defi}

We now look at functorialities :

\begin{defi}\label{DHdgjalpha}
Let $k\subset\mathbb C$ a subfield. Let $S\in\SmVar(k)$. Let $j:S^o\hookrightarrow S$ an open embedding.
Let $Z:=S\backslash S^o=V(\mathcal I)\subset S$ an the closed complementary subset, 
$\mathcal I\subset O_S$ being an ideal subsheaf. 
Taking generators $\mathcal I=(s_1,\ldots,s_r)$, we get $Z=V(s_1,\ldots,s_r)=\cap^r_{i=1}Z_i\subset S$ with 
$Z_i=V(s_i)\subset S$, $s_i\in\Gamma(S,\mathcal L_i)$ and $L_i$ a line bundle. 
Note that $Z$ is an arbitrary closed subset, $d_Z\geq d_X-r$ needing not be a complete intersection. 
Denote by $j_I:S^{o,I}:=\cap_{i\in I}(S\backslash Z_i)=S\backslash(\cup_{i\in I}Z_i)\xrightarrow{j_I^o}S^o\xrightarrow{j} S$ 
the open embeddings.
Let $(M,F,W)\in MHM_{gm,k,\mathbb C}(S^o))$. We then define, using definition \ref{DHdgj} and definition \ref{jw}
\begin{itemize}
\item the canonical extension 
\begin{eqnarray*}
j_{*Hdg}((M,F,W),(K,W),\alpha):=(j_{*Hdg}(M,F,W),j_{*w}(K,W),j_*\alpha) \\
:=\Tot((j_{I*Hdg}j_I^*(M,F,W),j_{I*w}j_I^*(K,W),j_{I*}\alpha))\in MHM_{gm,k,\mathbb C}(S), 
\end{eqnarray*}
so that $j^*(j_{*Hdg}((M,F,W),(K,W),\alpha))=((M,F,W),(K,W),\alpha)$,
\item the canonical extension 
\begin{eqnarray*}
j_{!Hdg}((M,F,W),(K,W),\alpha):=(j_{!Hdg}(M,F,W),j_{!w}(K,W),j_!\alpha) \\
:=\Tot((j_{I!Hdg}j_I^*(M,F,W),j_{I!w}j_I^*(K,W),j_{I!}\alpha))\in MHM_{gm,k,\mathbb C}(S),  
\end{eqnarray*}
so that $j^*(j_{!Hdg}((M,F,W),(K,W),\alpha))=((M,F,W),(K,W),\alpha)$.
\end{itemize}
Moreover for $((M',F,W),(K',W),\alpha')\in MHM_{gm,k,\mathbb C}(S)$,
\begin{itemize}
\item there is a canonical map in $MHM_{gm,k,\mathbb C}(S)$
\begin{equation*}
\ad(j^*,j_{*Hdg})((M',F,W),(K',W),\alpha'):((M',F,W),(K',W),\alpha')\to j_{*Hdg}j^*((M',F,W),(K',W),\alpha'), 
\end{equation*}
\item there is a canonical map in $MHM_{gm,k,\mathbb C}(S)$
\begin{equation*}
\ad(j_{!Hdg},j^*)((M',F,W),(K',W),\alpha'):j_{!Hdg}j^*((M',F,W),(K',W),\alpha')\to((M',F,W),(K',W),\alpha').
\end{equation*}
\end{itemize}
\end{defi}

Let $j:S^o\hookrightarrow S$ an open embedding with $S\in\SmVar(k)$.
For $(M,F,W)\in C(MHM_{gm,k,\mathbb C}(S^o))$, 
\begin{itemize}
\item we have the canonical map in $C_{\mathcal D(1,0)fil}(S)\times_IC_{fil}(S_{\mathbb C}^{an})$
\begin{eqnarray*}
T(j_{*Hdg},j_*)((M,F,W),(K,W),\alpha):=(k\circ\ad(j^*,j_*)(-),k\circ\ad(j^*,j_*),0): \\
j_{*Hdg}((M,F,W),(K,W),\alpha)\to (j_*E(M,F,W),j_*E(K,W),\alpha)
\end{eqnarray*}
\item we have the canonical map in $C_{\mathcal D(1,0)fil}(S)\times_IC_{fil}(S_{\mathbb C}^{an})$
\begin{eqnarray*}
T(j_!,j_{!Hdg})((M,F,W),(K,W),\alpha):=(k\circ\ad(j_!,j^*)(-),k\circ\ad(j_!,j^*)(-),0): \\
(j_!(M,F,W),j_!(K,W),j_!\alpha)\to j_{!Hdg}((M,F,W),(K,W),\alpha)
\end{eqnarray*}
\end{itemize}
the canonical maps.

\begin{prop}\label{jHdgpropadC}
\begin{itemize}
\item[(i)] Let $S\in\SmVar(k)$.
Let $D=V(s)\subset S$ a divisor with $s\in\Gamma(S,L)$ and $L$ a line bundle ($S$ being smooth, $D$ is Cartier).
Denote by $j:S^o:=S\backslash D\hookrightarrow S$ the open complementary embedding. Then, 
\begin{itemize}
\item $(j^*,j_{*Hdg}):MHM_{gm,k,\mathbb C}(S)\leftrightarrows MHM_{gm,k,\mathbb C}(S^o)$ is a pair of adjoint functors
\item $(j_{!Hdg},j^*):MHM_{gm,k,\mathbb C}(S^o)\leftrightarrows MHM_{gm,k,\mathbb C}(S)$ is a pair of adjoint functors.
\end{itemize}
\item[(ii)] Let $S\in\SmVar(k)$.
Let $Z=V(\mathcal I)\subset S$ an arbitrary closed subset, $\mathcal I\subset O_S$ being an ideal subsheaf. 
Denote by $j:S^o:=S\backslash Z\hookrightarrow S$. Then,
\begin{itemize}
\item $(j^*,j_{*Hdg}):D(MHM_{gm,k,\mathbb C}(S))\leftrightarrows D(MHM_{gm,k,\mathbb C}(S^o))$ is a pair of adjoint functors
\item $(j_{!Hdg},j^*):D(MHM_{gm,k,\mathbb C}(S^o))\leftrightarrows D(MHM_{gm,k,\mathbb C}(S))$ is a pair of adjoint functors.
\end{itemize}
\end{itemize}
\end{prop}

\begin{proof}
\noindent(i): Follows from proposition \ref{jHdgpropad}.

\noindent(ii):Follows from (i) and the exactness of $j^*$, $j_{*Hdg}$ and $j_{!Hdg}$.
\end{proof}

\begin{defi}\label{gammaHdgalpha}
Let $S\in\SmVar(k)$. Let $Z\subset S$ a closed subset.
Denote by $j:S\backslash Z\hookrightarrow S$ the complementary open embedding. 
\begin{itemize}
\item[(i)] We define using definition \ref{gammaHdg}, definition \ref{gammaw} and definition \ref{falpha}(iii), 
the filtered Hodge support section functor
\begin{eqnarray*}
\Gamma^{Hdg}_Z:C(MHM_{gm,k,\mathbb C}(S))\to C(MHM_{gm,k,\mathbb C}(S)), \; \; ((M,F,W),(K,W),\alpha)\mapsto \\
\Gamma^{Hdg}_Z((M,F,W),(K,W),\alpha):=(\Gamma_Z^{Hdg}(M,F,W),\Gamma_Z^w(K,W),\Gamma_Z(\alpha)) \\
=\Cone(\ad(j^*,j_{*Hdg})(-):j_{*Hdg},j^*((M,F,W),(K,W),\alpha)\to((M,F,W),(K,W),\alpha)[-1]
\end{eqnarray*}
see definition \ref{DHdgjalpha} for the last equality, together we the canonical map 
\begin{eqnarray*}
\gamma^{Hdg}_Z((M,F,W),(K,W),\alpha):\Gamma^{Hdg}_Z((M,F,W),(K,W),\alpha)\to((M,F,W),(K,W),\alpha).
\end{eqnarray*}
\item[(i)'] Since $j_{*Hdg}:C(MHM_{gm,k,\mathbb C}(S^o))\to C(MHM_{gm,k,\mathbb C}(S))$ is an exact functor, 
$\Gamma^{Hdg}_Z$ induces the functor
\begin{eqnarray*}
\Gamma^{Hdg}_Z:D(MHM_{gm,k,\mathbb C}(S))\to D(MHM_{gm,k,\mathbb C}(S)), \\ 
((M,F,W),(K,W),\alpha)\mapsto\Gamma^{Hdg}_Z((M,F,W),(K,W),\alpha)
\end{eqnarray*}
\item[(ii)] We define using definition \ref{gammaHdg}, definition \ref{gammaw} and definition \ref{falpha}(iii) 
the dual filtered Hodge support section functor
\begin{eqnarray*}
\Gamma^{\vee,Hdg}_Z:C(MHM_{gm,k,\mathbb C}(S))\to C(MHM_{gm,k,\mathbb C}(S)), \; \; ((M,F,W),(K,W),\alpha)\mapsto \\
\Gamma^{\vee,Hdg}_Z((M,F,W),(K,W),\alpha):=(\Gamma_Z^{\vee,Hdg}(M,F,W),\Gamma_Z^{\vee,w}(K,W),\Gamma_Z^{\vee}(\alpha)) \\
=\Cone(\ad(j_{!Hdg},j^*)(-):j_{!Hdg},j^*((M,F,W),(K,W),\alpha) \to ((M,F,W),(K,W),\alpha))
\end{eqnarray*}
see definition \ref{DHdgjalpha} for the last equality, together we the canonical map 
\begin{eqnarray*}
\gamma^{\vee,Hdg}_Z((M,F,W),(K,W),\alpha):((M,F,W),(K,W),\alpha)\to\Gamma_Z^{\vee,Hdg}((M,F,W),(K,W),\alpha).
\end{eqnarray*}
\item[(ii)'] Since $j_{!Hdg}:C(MHM_{gm,k,\mathbb C}(S^o))\to C(MHM_{gm,k,\mathbb C}(S))$ is an exact functor, 
$\Gamma^{Hdg,\vee}_Z$ induces the functor
\begin{eqnarray*}
\Gamma^{\vee,Hdg}_Z:D(MHM_{gm,k,\mathbb C}(S))\to D(MHM_{gm,k,\mathbb C}(S)), \\ 
((M,F,W),(K,W),\alpha)\mapsto\Gamma^{\vee,Hdg}_Z((M,F,W),(K,W),\alpha)
\end{eqnarray*}
\end{itemize}
\end{defi}

In the singular case it gives :

\begin{defi}\label{gammaHdgsingalpha}
Let $S\in\Var(k)$. Let $Z\subset S$ a closed subset.
Let $S=\cup_{i=1}^sS_i$ an open cover such that there exist closed embeddings
$i_i:S_i\hookrightarrow\tilde S_i$ with $\tilde S_i\in\SmVar(k)$. Denote $Z_I:=Z\cap S_I$. 
Denote by $n:S\backslash Z\hookrightarrow S$ and $\tilde n_I:\tilde S_I\backslash Z_I\hookrightarrow\tilde S_I$ 
the complementary open embeddings. 
\begin{itemize}
\item[(i)] We define using definition \ref{gammaHdgsing}, definition \ref{gammaw} and definition \ref{falpha}(iii)
the filtered Hodge support section functor
\begin{eqnarray*}
\Gamma^{Hdg}_Z:C(MHM_{gm,k,\mathbb C}(S))\to C(MHM_{gm,k,\mathbb C}(S)), \\ 
(((M_I,F,W),u_{IJ}),(K,W),\alpha)\mapsto\Gamma^{Hdg}_Z(((M_I,F,W),u_{IJ}),(K,W),\alpha):= \\
:=(\Gamma_Z^{Hdg}((M_I,F,W),u_{IJ}),\Gamma_Z^w(K,W),\Gamma_Z(\alpha))
\end{eqnarray*}
together with the canonical map
\begin{eqnarray*}
\gamma^{Hdg}_Z(((M_I,F,W),u_{IJ}),(K,W),\alpha): \\
\Gamma^{Hdg}_Z(((M_I,F,W),u_{IJ}),(K,W),\alpha)\to(((M_I,F,W),u_{IJ}),(K,W),\alpha).
\end{eqnarray*}
\item[(i)'] By exactness of $\Gamma_Z^{Hdg}$ and $\Gamma_Z^w$ it induces the functor
\begin{eqnarray*}
\Gamma^{Hdg}_Z: D(MHM_{gm,k,\mathbb C}(S))\to D(MHM_{gm,k,\mathbb C}(S)), \\  
(((M_I,F,W),u_{IJ}),(K,W),\alpha)\mapsto\Gamma^{Hdg}_Z(((M_I,F,W),u_{IJ}),(K,W),\alpha)
\end{eqnarray*}
\item[(ii)] We define using definition \ref{gammaHdgsing}, definition \ref{gammaw} and definition \ref{falpha}(iii)
the dual filtered Hodge support section functor
\begin{eqnarray*}
\Gamma^{\vee,Hdg}_Z:C(MHM_{gm,k,\mathbb C}(S))\to C(MHM_{gm,k,\mathbb C}(S)), \; \; 
(((M_I,F,W),u_{IJ}),(K,W),\alpha)\mapsto \\
\Gamma^{\vee,Hdg}_Z(((M_I,F,W),u_{IJ}),(K,W),\alpha):=
(\Gamma_Z^{\vee,Hdg}((M_I,F,W),u_{IJ}),\Gamma_Z^{\vee,w}(K,W),\Gamma^{\vee}_Z(\alpha)),
\end{eqnarray*}
together we the canonical map 
\begin{eqnarray*}
\gamma^{\vee,Hdg}_Z(((M_I,F,W),u_{IJ}),(K,W),\alpha): \\
(((M_I,F,W),u_{IJ}),(K,W),\alpha)\to\Gamma_Z^{\vee,Hdg}(((M_I,F,W),u_{IJ}),(K,W),\alpha).
\end{eqnarray*}
\item[(ii)'] By exactness of $\Gamma_Z^{\vee,Hdg}$ and $\Gamma_Z^{\vee,w}$, it induces the functor
\begin{eqnarray*}
\Gamma^{\vee,Hdg}_Z:D(MHM_{gm,k,\mathbb C}(S))\to D(MHM_{gm,k,\mathbb C}(S)), \\
(((M_I,F,W),u_{IJ}),(K,W),\alpha)\mapsto\Gamma^{\vee,Hdg}_Z(((M_I,F,W),u_{IJ}),(K,W),\alpha) \\
:=(\Gamma_Z^{\vee,Hdg}((M_I,F,W),u_{IJ}),\Gamma_Z^{\vee,w}(K,W),\Gamma^{\vee}_Z(\alpha))
\end{eqnarray*}
\end{itemize}
\end{defi}

This gives the inverse image functor :

\begin{defi}\label{inverseHdgsingalpha}
Let $f:X\to S$ a morphism with $X,S\in\Var(k)$.
Assume there exist a factorization $f:X\xrightarrow{l}Y\times S\xrightarrow{p_S}S$ 
with $Y\in\SmVar(k)$, $l$ a closed embedding and $p_S$ the projection.
Let $S=\cup_{i\in I}$ an open cover such that there exist closed embeddings
$i:S_i\hookrightarrow\tilde S_i$ with $\tilde S_i\in\SmVar(k)$. 
Denote $X_I:=f^{-1}(S_I)$. We have then $X=\cup_{i\in I}X_i$ and the commutative diagrams
\begin{equation*}
\xymatrix{f:X_I\ar[r]^{l_I}\ar[rd] & Y\times S_I\ar[r]^{p_{S_I}}\ar[d]^{i_I':=(I\times i_I)} & S_I\ar[d]^{i_I} \\ 
\, & Y\times\tilde S_I\ar[r]^{p_{\tilde S_I}=:\tilde f_I} & \tilde S_I} 
\end{equation*}
\begin{itemize}
\item[(i)] For $(((M_I,F,W),u_{IJ}),(K,W),\alpha)\in C(MHM_{gm,k,\mathbb C}(S))$ 
we set (see definition \ref{gammaHdgsingalpha} for $l$)
\begin{eqnarray*}
f^{*Hdg}(((M_I,F,W),u_{IJ}),(K,W),\alpha):= \\
\Gamma_X^{Hdg}((p_{\tilde S_I}^{*mod[-]}(M_I,F,W),p_{\tilde S_I}^{*mod[-]}u_{IJ}),p_S^*(K,W),p_S^*\alpha)(d_Y)[2d_Y]
\in C(MHM_{gm,k,\mathbb C}(X)), 
\end{eqnarray*}
\item[(ii)] For $(((M_I,F,W),u_{IJ}),(K,W),\alpha)\in C(MHM_{gm,k,\mathbb C}(S))$ 
we set (see definition \ref{gammaHdgsingalpha} for $l$)
\begin{eqnarray*}
f^{!Hdg}(((M_I,F,W),u_{IJ}),(K,W),\alpha):= \\
\Gamma_X^{\vee,Hdg}((p_{\tilde S_I}^{*mod[-]}(M_I,F,W),p_{\tilde S_I}^{*mod[-]}u_{IJ}),p_S^*(K,W),p_S^*\alpha)
\in C(MHM_{gm,k,\mathbb C}(X)), 
\end{eqnarray*}
\end{itemize}
Let $j:S^o\hookrightarrow S$ an open embedding with $S\in\Var(k)$.
Let $S=\cup_{i\in I}$ an open cover such that there exist closed embeddings
$i:S_i\hookrightarrow\tilde S_i$ with $\tilde S_i\in\SmVar(k)$. We have then,
for $(((M_I,F,W),u_{IJ}),(K,W),\alpha)\in C(MHM_{gm,k,\mathbb C}(S))$, quasi-isomorphisms in $C(MHM_{gm,k,\mathbb C}(S))$
\begin{eqnarray*}
(I(j^*,j^{*mod}_{Hdg})(-),I):j_{*Hdg}(((M_I,F,W),u_{IJ}),(K,W),\alpha)\to \\ 
j^*(((M_I,F,W),u_{IJ}),(K,W),\alpha):=(j^*((M_I,F,W),u_{IJ}),j^*(K,W),\alpha)
\end{eqnarray*}
and
\begin{eqnarray*}
(I(j^*,j^{\hat*mod}_{Hdg})(-),I):j^{!Hdg}(((M_I,F,W),u_{IJ}),(K,W),\alpha)\to \\
j^*(((M_I,F,W),u_{IJ}),(K,W),\alpha):=(j^*((M_I,F,W),u_{IJ}),j^*(K,W),\alpha)
\end{eqnarray*}
\end{defi}

\begin{defi}\label{otimesHdgalpha}
Let $S\in\Var(k)$. Let $S=\cup_{i\in I}$ an open cover such that there exist closed embeddings
$i:S_i\hookrightarrow\tilde S_i$ with $\tilde S_i\in\SmVar(k)$. 
We have, using definition \ref{otimesHdg} the following bi-functor
\begin{eqnarray*}
(-)\otimes_{O_S}^{Hdg}(-):D(MHM_{gm,k,\mathbb C}(S))^2\to D(MHM_{gm,k,\mathbb C}(S)), \\
(((M_I,F,W),u_{IJ}),(K,W),\alpha),(((M'_I,F,W),v_{IJ}),(K',W),\alpha')\mapsto \\
(((M_I,F,W),u_{IJ}),(K,W),\alpha)\otimes_{O_S}^{Hdg}(((M'_I,F,W),v_{IJ}),(K',W),\alpha'):= \\
((M_I,F,W),u_{IJ})\otimes_{O_S}^{Hdg}((M'_I,F,W),v_{IJ}),(K,W)\otimes^{L,w}(K',W),\alpha\otimes\alpha')
\end{eqnarray*}
where the map $\alpha\otimes\alpha'$ is given in definition \ref{falpha}(v).
\end{defi}

\begin{prop}\label{compDmodDRHdgalpha}
Let $f_1:X\to Y$ and $f_2:Y\to S$ two morphism with $X,Y,S\in\QPVar(k)$. 
\begin{itemize}
\item[(i)]Let $\mathcal M\in C(MHM_{gm,k,\mathbb C}(S))$. Then, 
\begin{equation*}
(f_2\circ f_1)^{*Hdg}(\mathcal M)=f_1^{*Hdg}f_2^{*Hdg}(\mathcal M)\in D(MHM_{gm,k,\mathbb C}(X)).
\end{equation*}
\item[(ii)]Let $(M,F,W)\in C(MHM_{gm,k,\mathbb C}(S))$. Then,
\begin{equation*}
(f_2\circ f_1)^{!Hdg}(\mathcal M)=f_1^{!Hdg}f_2^{!Hdg}(\mathcal M)\in D(MHM_{gm,k,\mathbb C}(X))
\end{equation*}
\end{itemize}
\end{prop}

\begin{proof}
Immediate from definition.
\end{proof}

\begin{prop}\label{PSkMHM}
Let $S\in\SmVar(k)$. Let $D=V(s)\subset S$ a (Cartier) divisor, where $s\in\Gamma(S,L)$.
Denote $i:D:=S\backslash D\hookrightarrow S$ the closed embedding and $j:S^o\hookrightarrow S$ the open embedding.
\begin{itemize}
\item[(i)] Let $((M,F,W),(K,W),\alpha)\in MHM_{gm,k,\mathbb C}(S)$.
We have, using proposition \ref{PSkDRM}, the canonical quasi-isomorphism in $C(MHM_{gm,k,\mathbb C}(S))$ :
\begin{eqnarray*}
Is(M):=(Is(M),Is(K),0):((M,F,W),(K,W),\alpha)\to \\
(\psi^u_D((M,F,W),\psi^u_D(K,W),\psi_D\alpha)\xrightarrow{((c(x_{S^o/S}(M)),can(M)),(c(x_{S^o/S}(K)),can(K)),0)} \\
(x_{S^o/S}(M,F,W),x_{S^o/S}(K,W),x_{S^o/S}(\alpha))\oplus(\phi^u_D(M,F,W),\phi^u_D(K,W),\phi_D\alpha) \\
\xrightarrow{((0,exp(s\partial s+1)),var(M)),((0,T-I),var(K)),0)}
(\psi^u_D(M,F,W),\psi^u_D(K,W),\psi_D\alpha)).
\end{eqnarray*}
\item[(ii)]We denote by $MHM_{gm,k,\mathbb C}(S\backslash D)\times_J MHM_{gm,k,\mathbb C}(D)$ 
the category whose set of objects consists of
\begin{equation*}
\left\{(\mathcal M,\mathcal N,a,b),\mathcal M\in MHM_{gm,k,\mathbb C}(S\backslash D),\mathcal N\in MHM_{gm,k,\mathbb C}(D),
a:\psi^u_{D}\mathcal M\to N,b:N\to\psi^u_{D}M \right\}
\end{equation*}
The functor (see definition \ref{DHdgpsialpha})
\begin{eqnarray*}
(j^*,\phi^u_{D},c,v):MHM_{gm,k,\mathbb C}(S)\to MHM_{gm,k,\mathbb C}(S\backslash D)\times_J MHM_{gm,k,\mathbb C}(D), \\
((M,F,W),(K,W),\alpha)\mapsto \\
((j^*(M,F,W),j^*(K,W),j^*\alpha),(\phi^u_D(M,F,W),\phi^u_D(K,W),\phi_D\alpha), can(-),var(-))
\end{eqnarray*}
is an equivalence of category.
\end{itemize}
\end{prop}

\begin{proof}
\noindent(i):Follows from proposition \ref{PSkDRM}.

\noindent(ii):Follows from (i).
\end{proof}

Let $S\in\Var(k)$. Let $S=\cup_{i\in I}S_i$ an open cover such that there
exists closed embeddings $i_i:S_i\hookrightarrow\tilde S_i$ with $\tilde S_I\in\SmVar(k)$.
We have the category $D_{\mathcal D(1,0)fil,rh}(S/(\tilde S_I))\times_I D_{fil,c,k}(S_{\mathbb C}^{an})$  
\begin{itemize}
\item whose set of objects is the set of triples $\left\{(((M_I,F,W),u_{IJ}),(K,W),\alpha)\right\}$ with
\begin{eqnarray*} 
((M_I,F,W),u_{IJ})\in D_{\mathcal D(1,0)fil,rh}(S/(\tilde S_I)), \, (K,W)\in D_{fil,c,k}(S_{\mathbb C}^{an}), \\ 
\alpha:T(S/(\tilde S_I))(K,W)\otimes\mathbb C_{S_{\mathbb C}^{an}}\to DR(S)^{[-]}(((M_I,W),u_{IJ})^{an})
\end{eqnarray*}
where $\alpha$ is an morphism in $D_{fil}(S_{\mathbb C}^{an}/(\tilde S_{I,\mathbb C}^{an}))$,
\item and whose set of morphisms consists of 
\begin{equation*}
\phi=(\phi_D,\phi_C,[\theta]):(((M_{1I},F,W),u_{IJ}),(K_1,W),\alpha_1)\to(((M_{2I},F,W),u_{IJ}),(K_2,W),\alpha_2)
\end{equation*}
where $\phi_D:((M_1,F,W),u_{IJ})\to((M_2,F,W),u_{IJ})$ and $\phi_C:(K_1,W)\to (K_2,W)$ 
are morphisms and
\begin{eqnarray*}
\theta=(\theta^{\bullet},I(DR(S)(\phi^{an}_D))\circ I(\alpha_1),I(\alpha_2)\circ I(\phi_C\otimes I)): \\
I(T(S/(\tilde S_I))(K_1,W))\otimes\mathbb C_{S^{an}}[1]\to I(DR(S)(((M_{2I},W),u_{IJ})^{an}))  
\end{eqnarray*}
is an homotopy,  
$I:D_{fil}(S_{\mathbb C}^{an}/(\tilde S^{an}_{I,\mathbb C}))\to K_{fil}(S_{\mathbb C}^{an}/(\tilde S^{an}_{I,\mathbb C}))$
being the injective resolution functor, and for
\begin{itemize}
\item $\phi=(\phi_D,\phi_C,[\theta]):(((M_{1I},F,W),u_{IJ}),(K_1,W),\alpha_1)\to(((M_{2I},F,W),u_{IJ}),(K_2,W),\alpha_2)$
\item $\phi'=(\phi'_D,\phi'_C,[\theta']):(((M_{2I},F,W),u_{IJ}),(K_2,W),\alpha_2)\to(((M_{3I},F,W),u_{IJ}),(K_3,W),\alpha_3)$
\end{itemize}
the composition law is given by 
\begin{eqnarray*}
\phi'\circ\phi:=(\phi'_D\circ\phi_D,\phi'_C\circ\phi_C,
I(DR(S)(\phi^{'an}_D))\circ[\theta]+[\theta']\circ I(\phi_C\otimes I)[1]): \\
(((M_{1I},F,W),u_{IJ}),(K_1,W),\alpha_1)\to(((M_{3I},F,W),u_{IJ}),(K_3,W),\alpha_3),
\end{eqnarray*}
in particular for 
$(((M_I,F,W),u_{IJ}),(K,W),\alpha)\in D_{\mathcal D(1,0)fil,rh}(S/(\tilde S_I))\times_I D_{fil,c,k}(S_{\mathbb C}^{an})$,
\begin{equation*}
I_{(((M_I,F,W),u_{IJ}),(K,W),\alpha)}=((I_{M_I}),I_K,0),
\end{equation*}
\end{itemize}
and also the category 
$D_{\mathcal D(1,0)fil,rh,\infty}(S/(\tilde S_I))\times_I D_{fil,c,k}(S_{\mathbb C}^{an})$ defined in the same way,
together with the localization functor
\begin{eqnarray*}
(D(zar),I):C_{\mathcal D(1,0)fil,rh}(S/(\tilde S_I))\times_I D_{fil,c,k}(S_{\mathbb C}^{an})
\to D_{\mathcal D(1,0)fil,rh}(S/(\tilde S_I))\times_I D_{fil,c,k}(S_{\mathbb C}^{an}) \\
\to D_{\mathcal D(1,0)fil,rh,\infty}(S/(\tilde S_I))\times_I D_{fil,c,k}(S_{\mathbb C}^{an}).
\end{eqnarray*}
Note that if $\phi=(\phi_D,\phi_C,[\theta]):(((M_1,F,W),u_{IJ}),(K_1,W),\alpha_1)\to(((M_2,F,W),u_{IJ}),(K_2,W),\alpha_2)$
is a morphism in $D_{\mathcal D(1,0)fil,rh}(S/(\tilde S_I))\times_I D_{fil,c,k}(S_{\mathbb C}^{an})$
such that $\phi_D$ and $\phi_C$ are isomorphism then $\phi$ is an isomorphism (see remark \ref{CGrem}).
Moreover,
\begin{itemize}
\item For 
$(((M_{I},F,W),u_{IJ}),(K,W),\alpha)\in D_{\mathcal D(1,0)fil,rh}(S/(\tilde S_I))\times_I D_{fil,c,k}(S_{\mathbb C}^{an})$, 
we set
\begin{equation*}
(((M_{I},F,W),u_{IJ}),(K,W),\alpha)[1]:=(((M_{I},F,W),u_{IJ})[1],(K,W)[1],\alpha[1]).
\end{equation*}
\item For 
\begin{equation*}
\phi=(\phi_D,\phi_C,[\theta]):(((M_{1I},F,W),u_{IJ}),(K_1,W),\alpha_1)\to(((M_{2I},F,W),u_{IJ}),(K_2,W),\alpha_2)
\end{equation*}
a morphism in $D_{\mathcal D(1,0)fil,rh}(S/(\tilde S_I))\times_I D_{fil,c,k}(S_{\mathbb C}^{an})$, 
we set (see \cite{CG} definition 3.12)
\begin{eqnarray*}
\Cone(\phi):=(\Cone(\phi_D),\Cone(\phi_C),((\alpha_1,\theta),(\alpha_2,0)))
\in D_{\mathcal D(1,0)fil,rh}(S/(\tilde S_I))\times_I D_{fil,c,k}(S_{\mathbb C}^{an}),
\end{eqnarray*}
$((\alpha_1,\theta),(\alpha_2,0))$ being the matrix given by the composition law, together with the canonical maps
\begin{itemize}
\item $c_1(-)=(c_1(\phi_D),c_1(\phi_C),0):(((M_{2I},F,W),u_{IJ}),(K_2,W),\alpha_2)\to\Cone(\phi)$
\item $c_2(-)=(c_2(\phi_D),c_2(\phi_C),0):\Cone(\phi)\to (((M_{1I},F,W),u_{IJ}),(K_1,W),\alpha_1)[1]$.
\end{itemize}
\end{itemize}

We have then the following :

\begin{thm}\label{Bek}
\begin{itemize}
\item[(i)]Let $S\in\Var(k)$. Let $S=\cup_{i\in I}S_i$ an open cover such that there exists
closed embedding $i_i:S_i\hookrightarrow\tilde S_i$ with $\tilde S_i\in\SmVar(k)$. Then the full embedding
\begin{eqnarray*}
\iota_S:MHM_{gm,k,\mathbb C}(S)\hookrightarrow
\PSh^0_{\mathcal D(1,0)fil,rh}(S/(\tilde S_I))\times_I P_{fil,k}(S_{\mathbb C}^{an}) 
\hookrightarrow C_{\mathcal D(1,0)fil,rh}(S/(\tilde S_I))\times_I D_{fil,c,k}(S_{\mathbb C}^{an}) 
\end{eqnarray*}
induces a full embedding
\begin{equation*}
\iota_S:D(MHM_{gm,k,\mathbb C}(S))\hookrightarrow D_{\mathcal D(1,0)fil,rh}(S/(\tilde S_I))\times_I D_{fil,c,k}(S_{\mathbb C}^{an}) 
\end{equation*}
whose image consists of 
$(((M_I,F,W),u_{IJ}),(K,W),\alpha)\in D_{\mathcal D(1,0)fil,rh}(S/(\tilde S_I))\times_I D_{fil,c,k}(S_{\mathbb C}^{an})$ such that 
\begin{equation*}
((H^n(M_I,F,W),H^n(u_{IJ})),H^n(K,W),H^n\alpha)\in MHM_{gm,k,\mathbb C}(S) 
\end{equation*}
for all $n\in\mathbb Z$ and such that for all $p\in\mathbb Z$,
the differentials of $\Gr_W^p(M_I,F)$ are strict for the filtrations $F$.
\item[(i)']Let $S\in\Var(k)$. Let $S=\cup_{i\in I}S_i$ an open cover such that there exists
closed embedding $i_i:S_i\hookrightarrow\tilde S_i$ with $\tilde S_i\in\SmVar(k)$. Then,
\begin{eqnarray*}
D(MHM_{gm,k,\mathbb C}(S))=
<(\int^{FDR}_f(n\times I)_{!Hdg}(\Gamma_X^{\vee,Hdg}(O_{\mathbb P^{N,o}\times\tilde S_I},F_b),x_{IJ})(d),
Rf_*\mathbb Q^w_X,f_*\alpha(X)), \\
(f:X\xrightarrow{l}\mathbb P^{N,o}\times S\xrightarrow{p}S)\in\QPVar(k), \; d\in\mathbb Z> \\
=<(\int^{FDR}_f(\Gamma_X^{\vee,Hdg}(O_{\mathbb P^{N,o}\times\tilde S_I},F_b),x_{IJ})(d),Rf_*\mathbb Q_X,f_*\alpha(X)), \\
(f:X\xrightarrow{l}\mathbb P^{N,o}\times S\xrightarrow{p}S)\in\QPVar(k), \; 
\mbox{proper}, \; X \mbox{smooth},\; d\in\mathbb Z> \\
\subset D_{\mathcal D(1,0)fil,rh}(S/(\tilde S_I))\times_I D_{fil,c,k}(S_{\mathbb C}^{an}) 
\end{eqnarray*}
where $n:\mathbb P^{N,o}\hookrightarrow\mathbb P^N$ are open embeddings, $l$ are closed embedding
and $<,>$ means the full triangulated category generated by and $(-)$ the shift of the $F$-filtration.
\item[(ii)]Let $S\in\Var(k)$. Let $S=\cup_{i\in I}S_i$ an open cover such that there exists
closed embedding $i_i:S_i\hookrightarrow\tilde S_i$ with $\tilde S_i\in\SmVar(k)$. Then the full embedding
\begin{eqnarray*}
\iota_S:MHM_{gm,k,\mathbb C}(S)\hookrightarrow
\PSh^0_{\mathcal D(1,0)fil,rh}(S/(\tilde S_I))\times_I P_{fil,k}(S_{\mathbb C}^{an})
\hookrightarrow C_{\mathcal D(1,0)fil,rh}(S/(\tilde S_I))\times_I D_{fil,c,k}(S_{\mathbb C}^{an}) 
\end{eqnarray*}
induces a full embedding
\begin{equation*}
\iota_S:D(MHM_{gm,k,\mathbb C}(S))\hookrightarrow 
D_{\mathcal D(1,0)fil,\infty,rh}(S/(\tilde S_I))\times_I D_{fil,c,k}(S_{\mathbb C}^{an}) 
\end{equation*}
whose image consists of 
$(((M_I,F,W),u_{IJ}),(K,W),\alpha)\in D_{\mathcal D(1,0)fil,\infty,rh}(S/(\tilde S_I))\times_I D_{fil}(S_{\mathbb C}^{an})$ 
such that 
\begin{equation*}
((H^n(M_I,F,W),H^n(u_{IJ})),H^n(K,W),H^n\alpha)\in MHM_{gm,k,\mathbb C}(S) 
\end{equation*}
for all $n\in\mathbb Z$ and such that there exist $r\in\mathbb Z$ and an $r$-filtered homotopy equivalence
$((M_I,F,W),u_{IJ})\to ((M'_I,F,W),u_{IJ})$ such that for all $p\in\mathbb Z$
the differentials of $\Gr_W^p(M'_I,F)$ are strict for the filtrations $F$.
\end{itemize}
\end{thm}

\begin{proof}
\noindent(i): We first show that $\iota_S$ is fully faithfull, that is for all
$\mathcal M=(((M_I,F,W),u_{IJ}),(K,W),\alpha),\mathcal M'=(((M'_I,F,W),u_{IJ}),(K',W),\alpha')\in MHM_{gm,k,\mathbb C}(S)$ 
and all $n\in\mathbb Z$,
\begin{eqnarray*}
\iota_S:\Ext_{D(MHM_{gm,k,\mathbb C}(S))}^n(\mathcal M,\mathcal M'):=
\Hom_{D(MHM_{gm,k,\mathbb C}(S))}(\mathcal M,\mathcal M'[n]) \\
\to\Ext_{\mathcal D(S)}^n(\mathcal M,\mathcal M')
:=\Hom_{\mathcal D(S):=D_{\mathcal D(1,0)fil,rh}(S/(\tilde S_I))\times_I D_{fil}(S_{\mathbb C}^{an})}(\mathcal M,\mathcal M'[n])
\end{eqnarray*}
For this it is enough to assume $S$ smooth. We then proceed by induction on $max(\dim\supp(M),\dim\supp(M'))$. 
\begin{itemize}
\item For $\supp(M)=\supp(M')=\left\{s\right\}$, 
it is the theorem for mixed hodge complexes or absolute Hodge complexes, see \cite{CG}. 
If $\supp(M)=\left\{s\right\}$ and $\supp(M')=\left\{s'\right\}$ and $s'\neq s$, 
then by the localization exact sequence
\begin{equation*}
\Ext_{D(MHM_{gm,k,\mathbb C}(S))}^n(\mathcal M,\mathcal M')=0=\Ext_{\mathcal D(S)}^n(\mathcal M,\mathcal M')
\end{equation*}
\item Denote $\supp(M)=Z\subset S$ and $\supp(M')=Z'\subset S$.
There exist an open subset $S^o\subset S$ such that $Z^o:=Z\cap S^o$ and $Z^{'o}:=Z'\cap S^o$ are smooth,
and $\mathcal M_{|Z^o}:=((i^*\Gr_{V_{Z^o},0}M_{|S^o},F,W),i^*j^*(K,W),\alpha^*(i))\in MHM_{gm,k}(Z^o)$ and 
$\mathcal M'_{|Z^{'o}}:=((i^{'*}\Gr_{V_{Z^{'o}},0}M'_{|S^o},F,W),i^{'*}j^*K,\alpha^*(i'))\in MHM_{gm,k}(Z^{'o})$ 
are variation of geometric mixed Hodge structure over $k\subset\mathbb C$, 
where $j:S^o\hookrightarrow S$ is the open embedding, and
$i:Z^o\hookrightarrow S^o$, $i':Z^{'o}\hookrightarrow S^o$ the closed embeddings.
Considering the connected components of $Z^o$ and $Z^{'o}$, we way assume that $Z^o$ and $Z^{'o}$ are connected.
Shrinking $S^o$ if necessary, we may assume that either $Z^o=Z^{'o}$ or $Z^o\cap Z^{'o}=\emptyset$,
We denote $D=S\backslash S^o$. Shrinking $S^o$ if necessary, 
we may assume that $D$ is a divisor and denote by $l:S\hookrightarrow L_D$ the zero section embedding.
\begin{itemize}
\item If $Z^o=Z^{'o}$, denote $i:Z^o\hookrightarrow S^o$ the closed embedding.
We have then the following commutative diagram
\begin{equation*}
\xymatrix{\Ext_{D(MHM_{gm,k,\mathbb C}(S^o))}^n(\mathcal M_{|S^o},\mathcal M'_{|S^o})
\ar[rr]^{\iota_{S^o}}\ar[d]_{(i^*\Gr_{V_{Z^o},0},i^*,\alpha^*(i))} & \, & 
\Ext_{\mathcal D(S^o)}^n(\mathcal M_{|S^o},\mathcal M'_{|S^o})\ar[d]^{(i^*\Gr_{V_{Z^o},0},i^*,\alpha^*(i))} \\
\Ext_{D(MHM_{gm,k,\mathbb C}(Z^o))}^n(\mathcal M_{|Z^o},\mathcal M'_{|Z^o})
\ar[rr]^{\iota_{Z^o}}\ar[u]_{(i_{*mod},i_*,\alpha_*(i))} & \, &
\Ext_{\mathcal D(Z^o)}^n(\mathcal M_{|Z^o},\mathcal M'_{|Z^o})\ar[u]^{(i_{*mod},i_*,\alpha_*(i))}}
\end{equation*}
Now we prove that $\iota_{Z^o}$ is an isomorphism similarly to the proof the the generic case of \cite{Be}.
On the other hand the left and right colummn are isomorphisms.
Hence $\iota_{S^o}$ is an isomorphism by the diagram.
\item If $Z^o\cap Z^{'o}=\emptyset$, we consider the following commutative diagram
\begin{equation*}
\xymatrix{\Ext_{D(MHM_{gm,k,\mathbb C}(S^o))}^n(\mathcal M_{|S^o},\mathcal M'_{|S^o})
\ar[rr]^{\iota_{S^o}}\ar[d]_{(i^*\Gr_{V_{Z^o},0},i^*,\alpha^*(i))} & \, & 
\Ext_{\mathcal D(S^o)}^n(\mathcal M_{|S^o},\mathcal M'_{|S^o})\ar[d]^{(i^*\Gr_{V_{Z^o},0},i^*,\alpha^*(i))} \\
\Ext_{D(MHM_{gm,k,\mathbb C}(Z^o))}^n(\mathcal M_{|Z^o},0)=0\ar[rr]^{\iota_{Z^o}}\ar[u]_{(i_{*mod},i_*,\alpha_*(i))} & \, &
\Ext_{\mathcal D(Z^o)}^n(\mathcal M_{|Z^o},0)=0\ar[u]^{(i_{*mod},i_*,\alpha_*(i))}}
\end{equation*}
where the left and right column are isomorphism by strictness of the $V_{Z^o}$ filtration
(use a bi-filtered injective resolution with respect to $F$ and $V_{Z^o}$ for the right column).
\end{itemize}
\item We consider now the following commutative diagram in $C(\mathbb Z)$ 
where we denote for short $H:=D(MHM_{gm,k,\mathbb C}(S))$
\begin{equation*}
\xymatrix{0\ar[r] & \Hom_{H}^{\bullet}(\Gamma^{\vee,Hdg}_D\mathcal M,\Gamma^{Hdg}_D\mathcal M')
\ar[r]^{\Hom(-,\gamma^{Hdg}_D(\mathcal M'))}\ar[d]^{\iota_S} &
\Hom_{H}^{\bullet}(\Gamma^{\vee,Hdg}_D\mathcal M,\mathcal M')
\ar[r]^{\Hom(-,\ad(j^*,j_{*Hdg})(\mathcal M'))}\ar[d]^{\iota_S} &
\Hom_{H}^{\bullet}(\Gamma^{\vee,Hdg}_D\mathcal M,j_{*Hdg}j^*\mathcal M')\ar[r]\ar[d]^{\iota_S} & 0 \\
0\ar[r] & \Hom_{\mathcal D(S)}^{\bullet}(\Gamma^{\vee,Hdg}_D\mathcal M,\Gamma^{Hdg}_D\mathcal M')
\ar[r]^{\Hom(-,\gamma^{Hdg}_D(\mathcal M'))} &
\Hom_{\mathcal D(S)}^{\bullet}(\Gamma^{\vee,Hdg}_D\mathcal M,\mathcal M')\ar[r]^{\Hom(-,\ad(j^*,j_{*Hdg})(\mathcal M'))} &
\Hom_{\mathcal D(S)}^{\bullet}(\Gamma^{\vee,Hdg}_D\mathcal M,j_{*Hdg}j^*\mathcal M')\ar[r] & 0}
\end{equation*}
whose lines are exact sequence. We have on the one hand,
\begin{equation*}
\Hom_{D(MHM_{gm,k,\mathbb C}(S))}^{\bullet}(\Gamma^{\vee,Hdg}_D\mathcal M,j_{*Hdg}j^*\mathcal M')=0=
\Hom_{\mathcal D(S)}^{\bullet}(\Gamma^{\vee,Hdg}_D\mathcal M,j_{*Hdg}j^*\mathcal M')
\end{equation*}
On the other hand by induction hypothesis
\begin{equation*}
\iota_S:\Hom_{D(MHM_{gm,k,\mathbb C}(S))}^{\bullet}(\Gamma^{\vee,Hdg}_D\mathcal M,\Gamma^{Hdg}_D\mathcal M')\to
\Hom_{\mathcal D(S)}^{\bullet}(\Gamma^{\vee,Hdg}_D\mathcal M,\Gamma^{Hdg}_D\mathcal M')
\end{equation*}
is a quasi-isomorphism. Hence, by the diagram
\begin{equation*}
\iota_S:\Hom_{D(MHM_{gm,k,\mathbb C}(S))}^{\bullet}(\Gamma^{\vee,Hdg}_D\mathcal M,\mathcal M')\to
\Hom_{\mathcal D(S)}^{\bullet}(\Gamma^{\vee,Hdg}_D\mathcal M,\mathcal M')
\end{equation*}
is a quasi-isomorphism.
\item We consider now the following commutative diagram in $C(\mathbb Z)$ 
where we denote for short $H:=D(MHM_{gm,k,\mathbb C}(S))$
\begin{equation*}
\xymatrix{0\ar[r] & \Hom_{H}^{\bullet}(\Gamma^{\vee,Hdg}_D\mathcal M,\mathcal M')
\ar[r]^{\Hom(\gamma^{\vee,Hdg}_D(\mathcal M),-)}\ar[d]^{\iota_S} &
\Hom_{H}^{\bullet}(\mathcal M,\mathcal M')
\ar[r]^{\Hom(\ad(j_{!Hdg},j^*)(\mathcal M'),-)}\ar[d]^{\iota_S} &
\Hom_{H}^{\bullet}(j_{!Hdg}j^*\mathcal M,\mathcal M')\ar[r]\ar[d]^{\iota_S} & 0 \\
0\ar[r] & \Hom_{\mathcal D(S)}^{\bullet}(\Gamma^{\vee,Hdg}_D\mathcal M,\mathcal M')
\ar[r]^{\Hom(\gamma^{\vee,Hdg}_D(\mathcal M),-)} &
\Hom_{\mathcal D(S)}^{\bullet}(\mathcal M,\mathcal M')\ar[r]^{\Hom(\ad(j_{!Hdg},j^*)(\mathcal M),-)} &
\Hom_{\mathcal D(S)}^{\bullet}(j_{!Hdg}j^*\mathcal M,\mathcal M')\ar[r] & 0}
\end{equation*}
whose lines are exact sequence. On the one hand, the commutative diagram
\begin{equation*}
\xymatrix{\Hom_{D(MHM_{gm,k,\mathbb C}(S))}^{\bullet}(j_{!Hdg}j^*\mathcal M,\mathcal M')\ar[r]^{j^*}\ar[d]^{\iota_{S}} &
\Hom_{D(MHM_{gm,k,\mathbb C}(S^o))}^{\bullet}(j^*\mathcal M,j^*\mathcal M')\ar[d]^{\iota_{S^o}} \\
\Hom_{\mathcal D(S)}^{\bullet}(j_{!Hdg}j^*\mathcal M,\mathcal M')\ar[r]^{j^*} &
\Hom_{\mathcal D(S^o)}^{\bullet}(j^*\mathcal M,j^*\mathcal M')}
\end{equation*}
together with the fact that the horizontal arrows $j^*$ are quasi-isomorphism 
by the functoriality given the uniqueness of the $V_S$ filtration for the embedding $l:S\hookrightarrow L_D$, 
(use a bi-filtered injective resolution with respect to $F$ and $V_S$ for the lower arrow)
and the fact that $\iota_{S^o}$ is a quasi-isomorphism by the first two point, show that
\begin{equation*}
\iota_S:\Hom_{D(MHM_{gm,k,\mathbb C}(S))}^{\bullet}(j_{!Hdg}j^*\mathcal M,\mathcal M')\to
\Hom_{\mathcal D(S)}^{\bullet}(j_{!Hdg}j^*\mathcal M,\mathcal M')
\end{equation*}
is a quasi-isomorphism. On the other hand, by the third point
\begin{equation*}
\iota_S:\Hom_{D(MHM_{gm,k,\mathbb C}(S))}^{\bullet}(\Gamma^{\vee,Hdg}_D\mathcal M,\mathcal M')\to
\Hom_{\mathcal D(S)}^{\bullet}(\Gamma^{\vee,Hdg}_D\mathcal M,\mathcal M')
\end{equation*}
is a quasi-isomorphism. Hence, by the diagram
\begin{equation*}
\iota_S:\Hom_{D(MHM_{gm,k,\mathbb C}(S))}^{\bullet}(\Gamma^{\vee,Hdg}_D\mathcal M,\mathcal M')\to
\Hom_{\mathcal D(S)}^{\bullet}(\Gamma^{\vee,Hdg}_D\mathcal M,\mathcal M')
\end{equation*}
is a quasi-isomorphism.
\end{itemize}
This shows the fully faithfulness. We now prove the essential surjectivity : let
\begin{equation*}
(((M_I,F,W),u_{IJ}),(K,W),\alpha)\in C_{\mathcal D(1,0)fil,rh}(S/(\tilde S_I))\times_I C_{fil}(S_{\mathbb C}^{an}) 
\end{equation*}
such that the cohomology are mixed hodge modules and such that the differential are strict.
We proceed by induction on $card\left\{n\in\mathbb Z\right\}, \, \mbox{s.t.} H^n(M_I,F,W)\neq 0$ by taking for 
the cohomological troncation 
\begin{equation*}
\tau^{\leq n}(((M_I,F,W),u_{IJ}),(K,W),\alpha):=
((\tau^{\leq n}(M_I,F,W),\tau^{\leq n}u_{IJ}),\tau^{\leq n}(K,W),\tau^{\leq n}\alpha)
\end{equation*}
and using the fact that the differential are strict for the filtration $F$ and the fully faithfullness.

\noindent(i)':Follows from (i).

\noindent(ii):Follows from (i).Indeed, in the composition of functor
\begin{eqnarray*}
\iota_S:D(MHM_{gm,k,\mathbb C}(S))\xrightarrow{\iota_S}
D_{\mathcal D(1,0)fil,rh}(S/(\tilde S_I))\times_I D_{fil}(S_{\mathbb C}^{an}) \\
\to D_{\mathcal D(1,0)fil,\infty,rh}(S/(\tilde S_I))\times_I D_{fil}(S_{\mathbb C}^{an}) 
\end{eqnarray*}
the second functor which is the localization functor is an isomorphism on the full subcategory
\begin{equation*}
D_{\mathcal D(1,0)fil,rh}(S/(\tilde S_I))^{st}\times_I D_{fil}(S_{\mathbb C}^{an})
\subset D_{\mathcal D(1,0)fil,rh}(S/(\tilde S_I))\times_I D_{fil}(S_{\mathbb C}^{an}) 
\end{equation*}
constisting of complex such that the differentials are strict for $F$, 
and the first functor $\iota_S$ is a full embedding by (i) and 
$\iota_S(D(MHM_{gm,k,\mathbb C}(S)))\subset D_{\mathcal D(1,0)fil,rh}(S/(\tilde S_I))^{st}\times_I D_{fil}(S_{\mathbb C}^{an})$.
\end{proof}

\begin{defi}\label{DHdgalpha}
Let $f:X\to S$ a morphism with $X,S\in\Var(k)$. 
Assume there exist a factorization $f:X\xrightarrow{l}Y\times S\xrightarrow{p_S}S$
with $Y\in\SmVar(k)$, $l$ a closed embedding and $p_S$ the projection.
Let $\bar Y\in\PSmVar(k)$ a smooth compactification of $Y$ with $n:Y\hookrightarrow\bar Y$ the open embedding.
Then $\bar f:\bar X\xrightarrow{\bar l}\bar Y\times_S\xrightarrow{\bar p_S}S$ is a compactification of $f$,
with $\bar X\subset\bar Y\times S$ the closure of $X$ and $\bar l$ the closed embedding,
and we denote by $n':X\hookrightarrow\bar X$ the open embedding so that $f=\bar f\circ n'$.
\begin{itemize}
\item[(i)]For $(((M_I,F,W),u_{IJ}),(K,W),\alpha)\in C(MHM_{gm,k,\mathbb C}(X))$, 
we define, using definition \ref{DHdgsing} and theorem \ref{Bek}
\begin{eqnarray*}
Rf_{*Hdg}(((M_I,F,W),u_{IJ}),(K,W),\alpha):=
\iota_S^{-1}(\int_{\bar f}^{FDR}(n\times I)_{*Hdg}((M_I,F,W),u_{IJ}),Rf_{*w}(K,W),f_*(\alpha)) \\
\in D(MHM_{gm,k,\mathbb C}(S))
\end{eqnarray*}
where $f_*(\alpha)$ is given in definition \ref{falpha}, and since 
\begin{itemize} 
\item by definition 
\begin{equation*}
H^i(\int_{\bar f}^{FDR}\Gr_W^k(I\times n)_{Hdg*}((M_I,F,W),u_{IJ}),
R\bar f_*\Gr_W^kn'_{*w}(K,W),\bar f_*\Gr_W^kn'_*\alpha)\in HM_{gm,k,\mathbb C}(S) 
\end{equation*}
for all $i,k\in\mathbb Z$, hence by the spectral sequence for the filtered complexes 
$\int_{\bar f}^{FDR}(I\times n)_{Hdg*}((M_I,F,W),u_{IJ})$ and $R\bar f_*((I\times n)_{*w}(K,W))$
\begin{eqnarray*}
\Gr_W^kH^i(\int_f^{Hdg}((M_I,F,W),u_{IJ}),Rf_{*w}(K,W),f_*\alpha)):= \\
(\Gr_W^kH^i\int_{\bar f}^{FDR}(I\times n)_{Hdg*}((M_I,F,W),u_{IJ}),
\Gr_W^kH^iR\bar f_*n'_{*w}(K,W),\Gr_W^kH^if_*\alpha)\in HM_{gm,k,\mathbb C}(S) 
\end{eqnarray*}
this gives by definition 
$H^i(\int_f^{Hdg}((M_I,F,W),u_{IJ}),Rf_{*w}(K,W),f_*(\alpha))\in MHM_{gm,k,\mathbb C}(S)$ for all $i\in\mathbb Z$. 
\item $\int_{f}^{Hdg}((M_I,F,W),u_{IJ})$ is the class of a complex such that the differential are strict for $F$
by theorem \ref{Sa12} in the complex case.
\end{itemize}
\item[(ii)]For $(((M_I,F,W),u_{IJ}),(K,W),\alpha)\in C(MHM_{gm,k,\mathbb C}(X))$, 
we define, using definition \ref{DHdgsing} and theorem \ref{Bek}, 
\begin{eqnarray*}
Rf_{!Hdg}(((M_I,F,W),u_{IJ}),(K,W),\alpha):=
\iota_S^{-1}(\int_{\bar f}^{FDR}(n\times I)_{!Hdg}((M_I,F,W),u_{IJ}),Rf_{!w}(K,W),f_!(\alpha)) \\
\in D(MHM_{gm,k,\mathbb C}(S))
\end{eqnarray*}
where $f_!(\alpha)$ is given in definition \ref{falpha}, and since 
\begin{itemize} 
\item by definition 
\begin{equation*}
H^i(\int_{\bar f}^{FDR}\Gr_W^k(n\times I)_{!Hdg}((M_I,F,W),u_{IJ}),
R\bar f_*\Gr_W^kn'_{!w}(K,W),\Gr_W^kf_!\alpha)\in HM_{gm,k,\mathbb C}(S) 
\end{equation*}
for all $i,k\in\mathbb Z$, 
hence by the spectral sequence for the filtered complexes 
$\int_{\bar f}^{FDR}(n\times I)_{!Hdg}((M_I,F,W),u_{IJ})$ and $R\bar f_*(n\times I)_{!w}(K,W)$
\begin{eqnarray*}
\Gr_W^kH^i(\int_{f!}^{Hdg}((M_I,F,W),u_{IJ}),Rf_{!w}(K,W),f_!\alpha):= \\
(\Gr_W^kH^i\int_{\bar f}^{FDR}(n\times I)_{!Hdg}((M_I,F,W),u_{IJ}),
\Gr_W^kH^iR\bar f_*n'_{!w}(K,W),\Gr_W^kH^if_!\alpha)\in HM_{gm,k,\mathbb C}(S) 
\end{eqnarray*}
this gives by definition 
$H^i(\int_{f!}^{Hdg}((M_I,F,W),u_{IJ}),Rf_{!w}(K,W),f_!(\alpha))\in MHM_{gm,k,\mathbb C}(S)$ for all $i\in\mathbb Z$. 
\item $\int_{f!}^{Hdg}((M_I,F,W),u_{IJ})$ is the class of a complex such that the differential are strict for $F$
by theorem \ref{Sa12} in the complex case.
\end{itemize}
\end{itemize}
\end{defi}

\begin{prop}\label{compDmodDRHdgDalpha}
Let $f_1:X\to Y$ and $f_2:Y\to S$ two morphism with $X,Y,S\in\QPVar(k)$. 
\begin{itemize}
\item[(i)]Let $\mathcal M\in C(MHM_{gm,k,\mathbb C}(X))$. Then, 
\begin{equation*}
R(f_2\circ f_1)^{Hdg}_*\mathcal M=Rf^{Hdg}_{2*}Rf^{Hdg}_{1*}\mathcal M\in D(MHM_{gm,k,\mathbb C}(S)).
\end{equation*}
\item[(ii)]Let $(M,F,W)\in C(MHM_{gm,k,\mathbb C}(X))$. Then,
\begin{equation*}
R(f_2\circ f_1)^{Hdg}_!\mathcal M=Rf^{Hdg}_{2!}Rf^{Hdg}_{1!}\mathcal M\in D(MHM_{gm,k,\mathbb C}(S))
\end{equation*}
\end{itemize}
\end{prop}

\begin{proof}
Immediate from definition.
\end{proof}

Let $k\subset\mathbb C$ a subfield.
Definition \ref{inverseHdgsingalpha}, definition \ref{DHdgalpha} and gives by proposition \ref{compDmodDRHdgalpha} 
and proposition \ref{compDmodDRHdgDalpha} respectively, the following 2 functors :
\begin{itemize}
\item We have the following 2 functor on the category of algebraic varieties over $k\subset\mathbb C$
\begin{eqnarray*}
D(MHM_{gm,k,\mathbb C}(\cdot)):\QPVar(k)\to\TriCat, \; S\mapsto D(MHM_{gm,k,\mathbb C}(S)), \\
(f:T\to S)\longmapsto (f^{*Hdg}:(((M_I,F,W),u_{IJ}),(K,W),\alpha)\mapsto \\
f^{!Hdg}(((M_I,F,W),u_{IJ}),(K,W),\alpha):=(f^{*mod}_{Hdg}(((M_I,F,W),u_{IJ})),f^{!w}(K,W),f^!\alpha)).
\end{eqnarray*}
see definition \ref{inverseHdgsing} and definition \ref{falpha} for the equality.
\item We have the following 2 functor on the category of quasi-projective algebraic varieties over $k\subset\mathbb C$
\begin{eqnarray*}
D(MHM_{gm,k,\mathbb C}(\cdot)):\QPVar(k)\to\TriCat, \; S\mapsto D(MHM_{gm,k,\mathbb C}(S)), \\
(f:T\to S)\longmapsto (f_{*Hdg}:(((M_I,F,W),u_{IJ}),(K,W),\alpha)\mapsto Rf_{*Hdg}(((M_I,F,W),u_{IJ}),(K,W),\alpha)).
\end{eqnarray*}
\item We have the following 2 functor on the category of quasi-projective algebraic varieties over $k\subset\mathbb C$
\begin{eqnarray*}
D(MHM_{gm,k,\mathbb C}(\cdot)):\QPVar(k)\to\TriCat, \; S\mapsto D(MHM_{gm,k,\mathbb C}(S)), \\
(f:T\to S)\longmapsto (f_{!Hdg}:(((M_I,F,W)),(K,W),\alpha)\mapsto f_{!Hdg}(((M_I,F,W),u_{IJ}),(K,W),\alpha)).
\end{eqnarray*}
\item We have the following 2 functor on the category of algebraic varieties over $k\subset\mathbb C$
\begin{eqnarray*}
D(MHM_{gm,k,\mathbb C}(\cdot)):\QPVar(k)\to\TriCat, \; S\mapsto D(MHM_{gm,k,\mathbb C}(S)), \\
(f:T\to S)\longmapsto (f^{!Hdg}:(((M_I,F,W),u_{IJ}),(K,W),\alpha)\mapsto \\
f^{*Hdg}(((M_I,F,W),u_{IJ}),(K,W),\alpha):=(f^{\hat*mod}_{Hdg}(((M_I,F,W),u_{IJ})),f^{*w}(K,W),f^*\alpha)).
\end{eqnarray*}
see definition \ref{inverseHdgsing} and definition \ref{falpha} for the equality.
\end{itemize}

\begin{prop}\label{HdgpropadsingC}
Let $f:X\to S$ with $S,X\in\QPVar(k)$. Then
\begin{itemize}
\item[(i)] $(f^{*Hdg},Rf^{Hdg}_*):D(MHM_{gm,k,\mathbb C}(S))\to D(MHM_{gm,k,\mathbb C}(X))$ is a pair of adjoint functors.
\begin{itemize}
\item For $(((M_I,F,W),u_{IJ}),(K,W),\alpha)\in C(MHM_{gm,k,\mathbb C}(S))$, 
\begin{eqnarray*}
\ad(f^{*Hdg},Rf^{Hdg}_*)(((M_I,F,W),u_{IJ}),(K,W),\alpha):= \\
(\ad(f_{Hdg}^{\hat*mod},Rf^{Hdg}_*)((M_I,F,W),u_{IJ}),\ad(f^{*w},Rf_{*w})(K,W)): \\
(((M_I,F,W),u_{IJ}),(K,W),\alpha)\to Rf^{Hdg}_*f^{*Hdg}(((M_I,F,W),u_{IJ}),(K,W),\alpha) 
\end{eqnarray*}
is the adjonction map in $D(MHM_{gm,k,\mathbb C}(S))$. 
\item For $(((N_I,F,W),u_{IJ}),(P,W),\beta)\in C(MHM_{gm,k,\mathbb C}(X))$,  
\begin{eqnarray*}
\ad(f^{*Hdg},Rf^{Hdg}_*)(((N_I,F,W),u_{IJ}),(P,W),\beta):= \\
(\ad(f_{Hdg}^{\hat*mod},Rf^{Hdg}_*)((N_I,F,W),u_{IJ}),\ad(f^{*w},Rf_{*w})(P,W)): \\
f^{*Hdg}Rf^{Hdg}_*(((N_I,F,W),u_{IJ}),(P,W),\beta)\to(((N_I,F,W),u_{IJ}),(P,W),\beta) 
\end{eqnarray*}
is the adjonction map in $D(MHM_{gm,k,\mathbb C}(X))$ 
\end{itemize}
\item[(ii)]$(Rf^{Hdg}_!,f^{!Hdg}):D(MHM_{gm,k,\mathbb C}(X))\to D(MHM_{gm,k,\mathbb C}(S))$ is a pair of adjoint functors.
\begin{itemize}
\item For $(((M_I,F,W),u_{IJ}),(K,W),\alpha)\in C(MHM_{gm,k,\mathbb C}(S))$, 
\begin{eqnarray*}
\ad(Rf^{Hdg}_!,f^{!Hdg})(((M_I,F,W),u_{IJ}),(K,W),\alpha):= \\
(\ad(f_{Hdg}^{*mod},Rf^{Hdg}_!)((M_I,F,W),u_{IJ}),\ad(f^{!w},Rf_{!w})(K,W)): \\
Rf^{Hdg}_!f^{!Hdg}(((M_I,F,W),u_{IJ}),(K,W),\alpha)\to(((M_I,F,W),u_{IJ}),(K,W),\alpha) 
\end{eqnarray*}
is the adjonction map in $D(MHM_{gm,k,\mathbb C}(S))$. 
\item For $(((N_I,F,W),u_{IJ}),(P,W),\beta)\in C(MHM_{gm,k,\mathbb C}(X))$, 
\begin{eqnarray*}
\ad(Rf^{Hdg}_!,f^{!Hdg})(((N_I,F,W),u_{IJ}),(P,W),\beta):= \\
(\ad(f_{Hdg}^{*mod},Rf^{Hdg}_!)((N_I,F,W),u_{IJ}),\ad(f^{!w},Rf_{!w})(P,W)): \\
(((N_I,F,W),u_{IJ}),(P,W),\beta)\to f^{!Hdg}Rf^{Hdg}_!(((N_I,F,W),u_{IJ}),(P,W),\beta) 
\end{eqnarray*}
is the adjonction map in $D(MHM_{gm,k,\mathbb C}(X))$. 
\end{itemize}
\end{itemize}
\end{prop}

\begin{proof}
Follows from proposition \ref{jHdgpropadC} after considering a factorization 
$f:X\hookrightarrow\bar Y\times S\xrightarrow{p_S} S$ with $\bar Y\in\PSmVar(k)$.
\end{proof}

\begin{thm}\label{sixMHMkC}
Let $k\subset\mathbb C$ a subfield.
\begin{itemize}
\item[(i)]We have the six functor formalism on $D(MHM_{gm,k,\mathbb C}(-)):\SmVar(k)\to\TriCat$.
\item[(ii)]We have the six functor formalism on $D(MHM_{gm,k,\mathbb C}(-)):\QPVar(k)\to\TriCat$.
\end{itemize}
\end{thm}

\begin{proof}
Follows from proposition \ref{HdgpropadsingC}.
\end{proof}

We recall the definition of the Deligne complex of a complex manifold and 
the Deligne cohomology class of an algebraic cycle of a complex algebraic variety.
\begin{defi}\label{Delkdef}
\begin{itemize}
\item[(i)] Let $X\in\AnSm(\mathbb C)$. We have for $d\in\mathbb Z$ the Deligne complex
\begin{equation*}
\mathbb Z_{\mathcal D,X}(d):=(\mathbb Z_X(d)\hookrightarrow\tau^{\leq d}DR(X))
=(\mathbb Z_X(d)\hookrightarrow (O_X\to\cdots\to\Omega_X^{d-1})\in C(X)
\end{equation*}
Let $D\subset X$ a normal crossing divisor. We have for $d\in\mathbb Z$ the Deligne complexes
\begin{equation*}
\mathbb Z_{\mathcal D,(X,D)}(d):=(\mathbb Z_X(d)\hookrightarrow(O_X\to\cdots\to\Omega_X^{d-1}(\log D)))\in C(X)
\end{equation*}
and
\begin{equation*}
\mathbb Z_{\mathcal D,(X,D)}(d)^{\vee}:=(\mathbb Z_X(d)\hookrightarrow(O_X\to\cdots\to\Omega_X^{d-1}(\nul D)))\in C(X).
\end{equation*}
Moreover we have (see \cite{DP}) canonical products 
\begin{itemize}
\item $(-)\cdot(-):\mathbb Z_{\mathcal D,(X,D)}(d)\otimes\mathbb Z_{\mathcal D,(X,D)}(d')\to\mathbb Z_{\mathcal D,(X,D)}(d+d')$
\item $(-)\cdot(-):\mathbb Z_{\mathcal D,(X,D)}(d)^{\vee}\otimes\mathbb Z_{\mathcal D,(X,D)}(d')^{\vee}
\to\mathbb Z_{\mathcal D,(X,D)}(d+d')^{\vee}$
\end{itemize}
\item[(ii)] Let $X\in\AnSm(\mathbb C)$. We have for $d\in\mathbb Z$ the Deligne (homology) complex
\begin{eqnarray*}
C^{\bullet}_{\mathcal D}(X,\mathbb Z(d)):=
\Cone(\mathbb Z\Hom_{Diff(\mathbb R)}(\Delta^{\bullet},X)\oplus\Gamma(X,F^d\mathcal D^{\bullet}_X)
\hookrightarrow\Gamma(X,\mathcal D^{\bullet}_X))\in C(\mathbb Z)
\end{eqnarray*}
Let $D\subset X$ a normal crossing divisor. Denote $U:=X\backslash D$. 
We have for $d\in\mathbb Z$ the Deligne (homology) complexes
\begin{eqnarray*}
C^{\bullet}_{\mathcal D}((X,D),\mathbb Z(d)):=
\Cone(\mathbb Z\Hom_{Diff(\mathbb R)}(\Delta^{\bullet},U)\oplus\Gamma(X,F^d\mathcal D^{\bullet}_X(\log D))
\hookrightarrow\Gamma(X,\mathcal D^{\bullet}_X(\log D)))\in C(\mathbb Z)
\end{eqnarray*}
and
\begin{eqnarray*}
C^{\bullet}_{\mathcal D}(X,D,\mathbb Z(d)):=
\Cone(\mathbb Z\Hom_{Diff(\mathbb R)}(\Delta^{\bullet},(X,D))\oplus\Gamma(X,F^d\mathcal D^{\bullet}_X(\nul D))
\hookrightarrow\Gamma(X,\mathcal D^{\bullet}_X(\nul D)))\in C(\mathbb Z).
\end{eqnarray*}
\item[(iii)] Let $k\subset\mathbb C_p$ a subfield.
Let $X\in\PSmVar(k)$. We have, for $k\in\mathbb Z$ and $d\in\mathbb Z$, the Deligne cohomology
\begin{eqnarray*}
H_{\mathcal D}^k(X_{\mathbb C}^{an},\mathbb Z(d)):=\mathbb H^k(X_{\mathbb C}^{an},\mathbb Z_{X,\mathcal D}(d))
=H^kC^{\bullet}_{\mathcal D}(X_{\mathbb C}^{an},D,\mathbb Z(d))^{\vee}
\end{eqnarray*}
Let $U\in\SmVar(k)$. Let $X\in\PSmVar(k)$ a compactification of $U$ with $D:=X\backslash U$ a normal crossing divisor.
We have, for $k\in\mathbb Z$ and $d\in\mathbb Z$, the Deligne cohomology
\begin{eqnarray*}
H_{\mathcal D}^k(U_{\mathbb C}^{an},\mathbb Z(d)):=
\mathbb H^k(X,\mathbb Z_{(X_{\mathbb C}^{an},D_{\mathbb C}^{an}),\mathcal D}(d))
=H^kC^{\bullet}_{\mathcal D}((X_{\mathbb C}^{an},D_{\mathbb C}^{an}),\mathbb Z(d))^{\vee}
\end{eqnarray*}
and
\begin{eqnarray*}
H_{\mathcal D}^k(X,D,\mathbb Z(d)):=
\mathbb H^k(X_{\mathbb C}^{an},\mathbb Z_{(X_{\mathbb C}^{an},D_{\mathbb C}^{an}),\mathcal D}(d)^{\vee})
=H^kC^{\bullet}_{\mathcal D}(X_{\mathbb C}^{an},D_{\mathbb C}^{an},\mathbb Z(d))^{\vee}.
\end{eqnarray*}
\item[(iv)] Let $k\subset\mathbb C_p$ a subfield.
Let $U\in\SmVar(k)$. Let $X\in\PSmVar(k)$ a compactification of $U$ with $D:=X\backslash U$ a normal crossing divisor.
We define the Deligne cohomology of a (higher) cycle $Z\in\mathcal Z^d(U,n)^{\partial=0}$ by
\begin{eqnarray*}
[Z]_{\mathcal D}:=\Im(H^{2d-n}(\gamma_{\supp(Z)})([Z])), \\ 
H^k(\gamma_{\supp(Z)}):
\mathbb H^{2d-n}_{\mathcal D,\supp(Z)}(X_{\mathbb C}^{an},\mathbb Z_{X_{\mathbb C}^{an},D_{\mathbb C}^{an}}(d))
\to\mathbb H^{2d-n}_{\mathcal D}(X_{\mathbb C}^{an},\mathbb Z_{X_{\mathbb C}^{an},D_{\mathbb C}^{an}}(d)) 
\end{eqnarray*}
with $\supp(Z):=p_X(\supp(Z))\subset X$, where $\supp(Z)\subset X\times\square^n$ is the support of $Z$.
\item[(v)]Let $k\subset\mathbb C_p$ a subfield.
Let $U\in\SmVar(k)$. Let $X\in\PSmVar(k)$ a compactification of $U$ with $D:=X\backslash U$ a normal crossing divisor.
We have for $d\in\mathbb Z$ the morphism of complexes
\begin{eqnarray*}
\mathcal R^d_U:\mathcal Z^d(U,\bullet)\to C^{\bullet}_{\mathcal D}(X_{\mathbb C}^{an},D_{\mathbb C}^{an},\mathbb Z(d)), \; 
Z\mapsto\mathcal R^d_U(Z):=(T_{\bar Z},\Omega_{\bar Z},R_{\bar Z})
\end{eqnarray*}
which gives for $Z\in\mathcal Z^d(U,n)^{\partial=0}$, 
\begin{equation*}
[\mathcal R^d_U(Z)]=[Z]_{\mathcal D}\in H_{\mathcal D}^{2d-n}(U_{\mathbb C}^{an},\mathbb Z(d))
\end{equation*}
\end{itemize}
\end{defi}

Let $f:X\to S$ a morphism with $S,X\in\AnSm(\mathbb C)$.
We have for $d\in\mathbb Z$ the canonical morphism of Deligne complexes
\begin{equation*}
(\ad(f^*,f_*)(\mathbb Z_S),\Omega_{X/S}^{\leq d}):\mathbb Z_{\mathcal D,S}(d)\to f_*\mathbb Z_{\mathcal D,X}(d)
\end{equation*}
which induces after taking the resolution of the Deligne complexes by differential forms
the morphism in $C(\mathbb Z)$
\begin{eqnarray*}
f^*:=(f^*,f^*,\theta(f)^t):
\Cone(\mathbb Z\Hom_{Diff(\mathbb R)}(\Delta^{\bullet},S)^{\vee}\oplus\Gamma(S,F^d\mathcal A^{\bullet}_S)
\hookrightarrow\Gamma(S,\mathcal A^{\bullet}_S)) \\
\to\Cone(\mathbb Z\Hom_{Diff(\mathbb R)}(\Delta^{\bullet},X)^{\vee}\oplus\Gamma(X,F^d\mathcal A^{\bullet}_X)
\hookrightarrow\Gamma(X,\mathcal A^{\bullet}_X))
\end{eqnarray*}
where $\theta(f)^t$ is the homotopy in the morphism in $D_{fil}(k)\otimes_ID(\mathbb Z)$
\begin{eqnarray*}
(f^*,f^*,\theta(f)^t):
(\Gamma(S,(\Omega^{\bullet}_S,F_b)),\mathbb Z\Hom_{Diff(\mathbb R)}(\Delta^{\bullet},S)^{\vee},a_{S*}\alpha(S)) \\
\to(\Gamma(X,(\Omega^{\bullet}_X,F_b)),\mathbb Z\Hom_{Diff(\mathbb R)}(\Delta^{\bullet},X)^{\vee},a_{X*}\alpha(X)),
\end{eqnarray*}
which induces in cohomology for $n\in\mathbb Z$, the morphisms of abelian groups
\begin{equation*}
f^*:H^n_{\mathcal D}(S,\mathbb Z(d))\to H^n_{\mathcal D}(X,\mathbb Z(d)) ;
\end{equation*}
we get dually, after taking the resolution of the Deligne complexes by currents the morphism in $C(\mathbb Z)$
\begin{eqnarray*}
f_*:=(f_*,f_*,\theta(f)):C^{\bullet}_{\mathcal D}(X,\mathbb Z(d)):=
\Cone(\mathbb Z\Hom_{Diff(\mathbb R)}(\Delta^{\bullet},X)\oplus\Gamma(X,F^d\mathcal D^{\bullet}_X)
\hookrightarrow\Gamma(X,\mathcal D^{\bullet}_X)) \\
\to C^{\bullet}_{\mathcal D}(S,\mathbb Z(d)):=
\Cone(\mathbb Z\Hom_{Diff(\mathbb R)}(\Delta^{\bullet},S)\oplus\Gamma(S,F^d\mathcal D^{\bullet}_S)
\hookrightarrow\Gamma(S,\mathcal D^{\bullet}_S))
\end{eqnarray*}
where $\theta(f)$ is the homotopy in the morphism in $D_{fil}(k)\otimes_ID(\mathbb Z)$
\begin{eqnarray*}
(f_*,f_*,\theta(f)):
(\Gamma(X,(\Omega^{\bullet}_X,F_b)),\mathbb Z\Hom_{Diff(\mathbb R)}(\Delta^{\bullet},X),a_{X!}\alpha(X)) \\
\to(\Gamma(S,(\Omega^{\bullet}_S,F_b)),\mathbb Z\Hom_{Diff(\mathbb R)}(\Delta^{\bullet},S),a_{S!}\alpha(S)),
\end{eqnarray*}
which induces in homology for $n\in\mathbb Z$, the morphisms of abelian groups
\begin{equation*}
f_*:H_{n,\mathcal D}(X,\mathbb Z(d))\to H_{n,\mathcal D}(S,\mathbb Z(d)).
\end{equation*}

\begin{thm}\label{Delk}
Let $k\subset\mathbb C$ a subfield.
\begin{itemize}
\item[(i)] Let $U\in\SmVar(k)$. Denote by $a_U:U\to\pt$ the terminal map. 
Let $X\in\PSmVar(k)$ a compactification of $U$ with $D:=X\backslash U$ a normal crossing divisor.
The embedding (see theorem \ref{Bek})
\begin{equation*}
\iota:D(MHM_{gm,k,\mathbb C}(\left\{\pt\right\}))\to D_{fil}(k)\times_ID(\mathbb Z) 
\end{equation*}
induces for $k\in\mathbb Z$ and $d\in\mathbb Z$, canonical isomorphisms
\begin{eqnarray*}
\iota(a_{U!Hdg}\mathbb Z^{Hdg}_U):H^k(a_{U!Hdg}\mathbb Z^{Hdg}_U)
\xrightarrow{\sim}H^k_{\mathcal D}(X_{\mathbb C}^{an},D_{\mathbb C}^{an},\mathbb Z(d)), 
\; \mbox{and} \\
\iota(a_{U*Hdg}\mathbb Z^{Hdg}_U):H^k(a_{U*Hdg}\mathbb Z^{Hdg}_U)\xrightarrow{\sim}
H^k_{\mathcal D}(U_{\mathbb C}^{an},\mathbb Z(d)).
\end{eqnarray*}
\item[(ii)] Let $h:U\to S$ and $h':U'\to S$ two morphism with $S,U,U'\in\SmVar(k)$.
Let $X\in\PSmVar(k)$ a compactification of $U$ with $D:=X\backslash U$ a normal crossing divisor
such that $h:U\to S$ extend to $f:X\to\bar S$.
Let $X'\in\PSmVar(k)$ a compactification of $U'$ with $D':=X'\backslash U'$ a normal crossing divisor
such that $h':U'\to S$ extend to $f':X'\to\bar S$.
The embedding $\iota:D(MHM_{gm,k,\mathbb C}(\pt))\to D_{fil}(k)\times_ID(\mathbb Z)$ (see theorem \ref{Bek}) 
induces for $k\in\mathbb Z$ and $d\in\mathbb Z$ a canonical isomorphism
\begin{eqnarray*}
\iota(a_{U'\times_SU!Hdg}\mathbb Z^{Hdg}_{U'\times_SU}):
\Hom_{D(MHM_{gm,k,\mathbb C}(S))}(h_{U'!Hdg}\mathbb Z^{Hdg}_{U'},h_{U!Hdg}\mathbb Z^{Hdg}_U(d)[k]) \\
\xrightarrow{RI(-,-)}
\Hom_{D(MHM_{gm,k,\mathbb C}(\pt))}(\mathbb Z^{Hdg}_{\pt},a_{U'\times_SU!Hdg}\mathbb Z^{Hdg}_{U'\times_SU}(d)[k])
=H^k(a_{U'\times_SU!Hdg}\mathbb Z^{Hdg}_{U'\times_SU}(d)) \\
\xrightarrow{\sim}
H^k_{\mathcal D}((X'\times_SX)_{\mathbb C}^{an},((X'\times_SU)\cup(U'\times_SX))_{\mathbb C}^{an},\mathbb Z(d)).
\end{eqnarray*}
\item[(iii)]Let $U\in\SmVar(k)$.  
Let $X\in\PSmVar(k)$ a compactification of $U$ with $D:=X\backslash U$ a normal crossing divisor.
For $[Z]\in\CH^d(U,n)$ and $[Z']\in\CH^{d'}(U,n')$, we have
\begin{equation*}
([Z]\cdot[Z'])_{\mathcal D}=[Z]_{\mathcal D}\cdot[Z']_{\mathcal D}\in H^{2d+2d'-n-n'}(U_{\mathbb C}^{an},\mathbb Z(d+d'))
\end{equation*}
where the product on the left is the intersection of higher Chow cycle which is well defined modulo boundary 
(they intersect properly modulo boundary) while the right product of Deligne cohomology classes is induced by 
the product of Deligne complexes 
$(-)\cdot(-):\mathbb Z_{\mathcal D,(X,D)}(d)\otimes\mathbb Z_{\mathcal D,(X,D)}(d')\to\mathbb Z_{\mathcal D,(X,D)}(d+d')$.
\item[(iv)]Let $h:U\to S$,$h':U'\to S$, $h'':U''\to S$ three morphism with $S,U,U',U''\in\SmVar(k)$.
Let $X\in\PSmVar(k)$ a compactification of $U$ with $D:=X\backslash U$ a normal crossing divisor
such that $h:U\to S$ extend to $f:X\to\bar S$.
Let $X'\in\PSmVar(k)$ a compactification of $U'$ with $D':=X'\backslash U'$ a normal crossing divisor
such that $h':U'\to S$ extend to $f':X'\to\bar S$.
Let $X'\in\PSmVar(k)$ a compactification of $U'$ with $D':=X'\backslash U'$ a normal crossing divisor
such that $h':U'\to S$ extend to $f':X'\to\bar S$.
For $[Z]\in\CH^d(U\times_SU',n)$ and $[Z']\in\CH^{d'}(U'\times_SU'',n')$, we have
\begin{equation*}
([Z]\circ[Z'])_{\mathcal D}=[Z]_{\mathcal D}\circ[Z']_{\mathcal D}
\in H^{d''-n''}((U\times_SU'')_{\mathbb C}^{an},\mathbb Z(d''-n''))
\end{equation*}
where the composition on the left is the composition of higher correspondence modulo boundary
while the composition on the right is given by (ii).
\end{itemize}
\end{thm}

\begin{proof}
\noindent(i):Standard.

\noindent(ii):Follows on the one hand from (i) and 
on the other hand the six functor formalism on the 2-functor 
$D(MHM_{gm,k,\mathbb C}(-)):\SmVar(k)\to\TriCat$ (theorem \ref{sixMHMkC}) gives the isomorphism $RI(-,-)$.

\noindent(iii):Standard.

\noindent(iv):Follows from (iii).

\end{proof}

\subsection{The $p$-adic case where $k\subset K\subset\mathbb C_p$}

Let $p$ a prime integer. Let $k\subset K\subset\mathbb C_p$ a subfield of a $p$ adic field $K$. 
Denote by $\bar k\subset\mathbb C_p$ its algebraic closure.
\begin{itemize}
\item We denote by $G:=\Gal(\bar K/K)\subset\Gal(\bar k/k)$ the Galois group of $K$.
\item For $S\in\Var(k)$, 
\begin{itemize}
\item we will consider $\mathbb B_{dr,S_K}:=\mathbb B_{dr,R_K(S^{an}_K)}$ and
$O\mathbb B_{dr,S_K}:=O\mathbb B_{dr,R_K(S^{an}_K)}$ where $R_K:\AnSp(K)\to\AdSp/(K,O_K)$ the canonical functor
\item we will consider $\mathbb B_{dr,S_{\mathbb C_p}}:=\mathbb B_{dr,R_{\mathbb C_p}(S^{an}_{\mathbb C_p})}$ and
$O\mathbb B_{dr,S_{\mathbb C_p}}:=O\mathbb B_{dr,R_{\mathbb C_p}(S^{an}_{\mathbb C_p})}$ 
where $R_{\mathbb C_p}:\AnSp(\mathbb C_p)\to\AdSp/(\mathbb C_p,O_{\mathbb C_p})$ the canonical functor.
\end{itemize}
\item Recall (see section 2) that for a prime number $l$, 
a $\mathbb Z_l$ module $K=(K_n)_{n\in\mathbb N}\in\Fun(\mathbb N,\Ab)$ is
a projective system with $K_n$ a $l^n$ torsion group such that $K_n\to K_{n+1}/l^nK_{n+1}$ is an isomorphism.
For $S\in\Var(k)$ and $l$ a prime integer, we have (see section 2) 
\begin{itemize}
\item the full subcategory $C_{\mathbb Z_lfil}(S^{et})\subset\PSh(S^{et},\Fun(\mathbb N,C(\mathbb Z)))$ and
the full subcategory
\begin{equation*}
D_{\mathbb Z_lfil,c,k}(S^{et})\subset D_{\mathbb Z_lfil}(S^{et}):=\Ho_{et}C_{\mathbb Z_lfil}(S^{et}),
\end{equation*}
whose cohomology sheaves of the graded piece are constructible with respect to a Zariski stratification of $S$,
the full subcategory $C_{\mathbb Z_lfil}(S_K^{an,pet})\subset\PSh(S_K^{an,pet},\Fun(\mathbb N,C(\mathbb Z)))$ and
the full subcategory
\begin{equation*}
D_{\mathbb Z_lfil,c,k}(S_K^{an,pet})\subset D_{\mathbb Z_lfil}(S_K^{an,pet}):=\Ho_{pet}C_{\mathbb Z_lfil}(S_K^{an,pet})
\end{equation*} 
whose cohomology sheaves of the graded piece are constructible with respect to a Zariski stratification of $S$,
\item $P_{\mathbb Z_l,fil}(S^{et})\subset D_{\mathbb Z_lfil,c}(S^{et})$ 
and $P_{\mathbb Z_l,fil,k}(S_K^{an,pet})\subset D_{\mathbb Z_lfil,c,k}(S_K^{an,pet})$ 
the full subcategories of filtered perverse sheaves.
\end{itemize}
\item Let $S\in\Var(k)$ and $D\subset S$ a Cartier divisor. For $(K,W)\in D_{\mathbb Z_lfil}(S^{et})$, we denote for short
\begin{eqnarray*}
\psi_D(K,W):=\psi_D(K,W)[-1]\in D_{\mathbb Z_lfil}(S^{et}), \phi_D(K,W):=\phi_D(K,W)[-1]\in D_{\mathbb Z_lfil}(S^{et}), \\
x_{S\backslash D/S}(K,W):=x_{S\backslash D/S}(K,W)[-1]\in D_{\mathbb Z_lfil}(S^{et})
\end{eqnarray*}
so that it sends (filtered) perverse sheaves to (filtered) perverse sheaves.
\item For $S\in\Var(k)$, we denote by 
$\an_S:S^{an}:=S_{\mathbb C_p}^{an}\xrightarrow{\an_S} S_{\mathbb C_p}\xrightarrow{\pi_{k/\mathbb C_p}(S)} S$
the morphism of ringed spaces given by the analytical functor.
\begin{itemize}
\item For $(M,F)\in C_{O_Sfil}(S)$, we denote by $(M,F)^{an}:=\an_S^{*mod}(M,F)\in C_{O_Sfil}(S_{\mathbb C_p}^{an})$.
\item For $(M,F)\in C_{\mathcal Dfil}(S)$,
we denote by $(M,F)^{an}:=\an_S^{*mod}(M,F)\in C_{\mathcal Dfil}(S_{\mathbb C_p}^{an})$.
\end{itemize}
We denote for short 
\begin{eqnarray*}
DR(S):=DR(S_{\mathbb C_p}^{an})\circ\an_S^{*mod}:C_{\mathcal Dfil}(S)\to C_{fil}(S_{\mathbb C_p}^{an,pet}), \;
M\mapsto DR(S)(M^{an})
\end{eqnarray*}
the De Rham functor.
\item Let $S\in\Var(k)$.
\begin{itemize}
\item For $K_1,K_2\in D_{\mathbb Z_p}(S_K^{an,pet})$, we denote for short 
$K_1\otimes_{\mathbb Q_p}K_2:=K_1\otimes^L_{\mathbb Q_p}K_2\in D_{\mathbb Z_p}(S_K^{an,pet})$
the derived tensor product.
\item For $(K_1,W),(K_2,W)\in D_{\mathbb Z_pfil}(S_K^{an,pet})$, we denote for short 
$(K_1,W)\otimes_{\mathbb Q_p}(K_2,W):=(K_1,W)\otimes^L_{\mathbb Q_p}(K_2,W)\in D_{\mathbb Z_pfil}(S_K^{an,pet})$
the derived tensor product.
\item For $M,N\in D_{\mathbb B_{dr,S_K}}(S_K^{an,pet})$, we denote for short 
$M\otimes_{\mathbb B_{dr,S}}N:=M\otimes^L_{\mathbb B_{dr,S_K}}N\in D_{\mathbb B_{dr,S_K}}(S_K^{an,pet})$
the derived tensor product.
\item For $(M,W),(N,W)\in D_{\mathbb B_{dr,S_K}fil}(S_K^{an,pet})$, we denote for short 
$(M,W)\otimes_{\mathbb B_{dr,S}}(N,W):=(M,W)\otimes^L_{\mathbb B_{dr,S_K}}(N,W)\in D_{\mathbb B_{dr,S_K}}(S_K^{an,pet})$
the derived tensor product.
\end{itemize}
\item Let $S\in\SmVar(k)$ and $D\subset S$ a (Cartier) divisor. 
Denote by $j:S^o:=S\backslash D\hookrightarrow S$ the open embedding.
Then $j_*:C(S_K^{o,an,pet})\to C(S_K^{an,pet})$ is preserve pro-etale equivalence, that is $Rj_*=j_*$.
\end{itemize}

\subsubsection{The $\mathbb B_{dr}$ functor}

Motivated by theorem \ref{RHpadic} and theorem \ref{PSk}, we make the following definition:

\begin{defi}\label{VfilKMmap}
Let $S\in\SmVar(k)$. 
Let $D=V(s)\subset S$ a (Cartier) divisor, where $s\in\Gamma(S,L)$ is a section of the line bundle $L=L_D$ associated to $D$.
We denote by $j:S^o:=S\backslash D\hookrightarrow S$ the open complementary subset.
Let $\pi:\tilde S_K^{o,an}\to S_K^{o,an}$ the perfectoid universal covering.
\begin{itemize}
\item[(i)] We define, using definition \ref{DHdgj},
\begin{equation*}
\mathbb B_{dr,S^o/S,K}:=F^0DR(S)((j_{*Hdg}(O_{S^o},F_b))^{an}\otimes_{O_{S_K^{an}}}(O\mathbb B_{dr,S_K},F))
\in C_{\mathbb B_{dr,S_K}}(S_K^{an,pet})
\end{equation*}
together with the canonical map in $C_{\mathbb B_{dr,S_K}}(S_K^{an,pet})$ 
\begin{equation*}
a_S(\mathbb B_{dr,S^o/S,K}):=F^0DR(S)(\ad(j^*,j_{*Hdg})(O_S,F_b)^{an}\otimes I):\mathbb B_{dr,S_K}\to\mathbb B_{dr,S^o/S,K}
\end{equation*}
\item[(ii)] We define, using definition \ref{DHdgpsi}, the nearby cycle module
\begin{equation*}
\mathbb B_{dr,\psi_D,K}:=F^0DR(S)((\psi^u_D(O_{S^o},F_b))^{an}\otimes_{O_{S_K^{an}}}(O\mathbb B_{dr,S_K},F))
\in C_{\mathbb B_{dr,S}}(S_K^{an,pet})
\end{equation*}
together with the canonical maps in $C_{\mathbb B_{dr,S_K}}(S_K^{an,pet})$ 
\begin{equation*}
\rho_{\mathbb B_{dr},D}(O_S):=F^0DR(S)(\rho^u_{DR,D}(O_{S^o},F_b)^{an}\otimes I):\mathbb B_{dr,S^o/S,K}\to\mathbb B_{dr,\psi_D,K}
\end{equation*}
and
\begin{eqnarray*}
a_S(\mathbb B_{dr,\psi_D}):=F^0DR(S)((\rho^u_{DR,D}(O_{S^o},F_b)\circ\ad(j^*,j_{*Hdg})(O_S,F_b))^{an}\otimes I): \\
\mathbb B_{dr,S_K}\xrightarrow{a_S(\mathbb B_{dr,S^o/S,K})}\mathbb B_{dr,S^o/S,K}
\xrightarrow{\rho_{\mathbb B_{dr},D}(O_S)}\mathbb B_{dr,\psi_D,K}
\end{eqnarray*}
where $\rho^u_{DR,D}(O_{S^o},F_b):j_{*Hdg}(O_{S^o},F_b)\to\psi^u_D(O_{S^o},F_b)$ is given in definition \ref{rhoDR}.
\item[(iii)] We define, using definition \ref{DHdgpsi} and (ii), the vanishing cycle module
\begin{equation*}
\mathbb B_{dr,\phi_D,K}:=F^0DR(S)((\phi^u_D(O_{S^o},F_b))^{an}\otimes_{O_{S_K^{an}}}(O\mathbb B_{dr,S_K},F))
\in C_{\mathbb B_{dr,S_K}}(S_K^{an,pet}).
\end{equation*}
together with the canonical maps in $C_{\mathbb B_{dr,S_K}}(S_K^{an,pet})$ 
We have using definition \ref{DHdgpsi} the following maps
\begin{equation*}
can_{\mathbb B_{dr},D}(O_S):=F^0DR(S)(can(O_{S^o},F_b)^{an}\otimes I):\mathbb B_{dr,\psi_D,K}\to\mathbb B_{dr,\phi_D,K},
\end{equation*}
and
\begin{equation*}
var_{\mathbb B_{dr},D}(O_S):=F^0DR(S)(var(O_{S^o},F_b)^{an}\otimes I):\mathbb B_{dr,\phi_D,K}\to\mathbb B_{dr,\psi_D,K}
\end{equation*}
and
\begin{eqnarray*}
a_S(\mathbb B_{dr,\phi_D,K}):=
F^0DR(S)((can(O_{S^o},F_b)\circ\rho^u_{DR,D}(O_{S^o},F_b)\circ\ad(j^*,j_{*Hdg})(O_S,F_b))^{an}\otimes I): \\
\mathbb B_{dr,S_K}\xrightarrow{a_S(\mathbb B_{dr,\psi_D,K})}\mathbb B_{dr,\psi_D,K}
\xrightarrow{can_{\mathbb B_{dr},D}(O_S)}\mathbb B_{dr,\phi_D,K}
\end{eqnarray*}
\item[(iv)]Using (ii), we set
\begin{eqnarray*}
\mathbb B_{dr,x_{S^o/S},K}:=\Cone(\rho_{\mathbb B_{dr},D}(O_S):\mathbb B_{dr,S^o/S,K}\to\mathbb B_{dr,\psi_D,K})
\in C_{\mathbb B_{dr,S_K}}(S_K^{an,pet})
\end{eqnarray*}
together with the canonical map in $C_{\mathbb B_{dr,S_K}}(S_K^{an,pet})$
\begin{equation*}
a_S(\mathbb B_{dr,x_{S^o/S},K}):=(a_S(\mathbb B_{dr,S^o/S,K}),a_S(\mathbb B_{dr,\psi_D,K})):
\mathbb B_{dr,S_K}\to\mathbb B_{dr,x_{S^o/S},K}.
\end{equation*}
\item[(v)]For $L\in\Shv_{\mathbb Z_p}(S^{o,et})$ a local system, we set using theorem \ref{HSk} for $j_*(L\otimes O_{S_K^o})$
\begin{eqnarray*}
V_{D0}j_*(L\otimes_{\mathbb Q_p}\mathbb B_{dr,S^o_K}):=
V_{D0}j_*(L\otimes_{\mathbb Q_p} O_{S_K^o})\otimes_{O_{S_K}}\mathbb B_{dr,S_K}
\in C_{\mathbb B_{dr}}(S_K^{an,pet})
\end{eqnarray*}
so that we have the isomorphism in $D_{\mathbb B_{dr}}(S_K^{an,pet})$
\begin{eqnarray*}
m(L\otimes_{\mathbb Q_p}\mathbb B_{dr,S^o_K}):
V_{D0}j_*(L\otimes_{\mathbb Q_p}\mathbb B_{dr,S^o_K})\otimes_{\mathbb B_{dr,S_K}}\mathbb B_{dr,S^o/S,K} \\
\xrightarrow{:=}
V_{D0}j_*(L\otimes_{\mathbb Q_p}\mathbb B_{dr,S^o_K})\otimes_{\mathbb B_{dr,S_K}}
F^0DR(S)(j_{*Hdg}(O_{S^o},F_b)\otimes_{O_S}(O\mathbb B_{dr,S},F)) \\
\to j_*L\otimes_{\mathbb Q_p}\mathbb B_{dr,S^o/S,K}, 
\xrightarrow{:=}
j_*L\otimes_{\mathbb Q_p}F^0DR(S)(j_{*Hdg}(O_{S^o},F_b)\otimes_{O_S}(O\mathbb B_{dr,S},F)), \\
s\otimes h\otimes w\mapsto s\otimes hw.
\end{eqnarray*}
More generally, for $Z\subset S$ a closed subset and $L\in\Shv_{\mathbb Z_p}(Z^{o,et})$ a local system with $Z^o:=Z\cap S^o$, 
we set using theorem \ref{HSk} for $j_*(L\otimes O_{Z_K^o})$
\begin{eqnarray*}
V_{D0}j_*(L\otimes_{\mathbb Q_p}\mathbb B_{dr,S^o_K}):=
V_{D0}j_*(L\otimes_{\mathbb Q_p} O_{Z_K^o})\otimes_{O_{S_K}}\mathbb B_{dr,S_K}
\in C_{\mathbb B_{dr}}(S_K^{an,pet})
\end{eqnarray*}
so that we have the isomorphism in $D_{\mathbb B_{dr}}(S_K^{an,pet})$
\begin{eqnarray*}
m(L\otimes_{\mathbb Q_p}\mathbb B_{dr,S^o_K}):
V_{D0}j_*(L\otimes_{\mathbb Q_p}\mathbb B_{dr,S^o_K})\otimes_{\mathbb B_{dr,S_K}}\mathbb B_{dr,S^o/S,K}
\to j_*L\otimes_{\mathbb Q_p}\mathbb B_{dr,S^o/S,K}, \\
s\otimes h\otimes w\mapsto s\otimes hw.
\end{eqnarray*}
\item[(vi)]For $L\in\Shv_{\mathbb Z_p}(S^{o,et})$ a local system, we set using theorem \ref{HSk} for $j_*(L\otimes O_{S_K^o})$
\begin{eqnarray*}
\psi_D(L\otimes_{\mathbb Q_p}\mathbb B_{dr,S^o_K}):=
\Gr^{V_D}_{-1\leq\alpha<0}j_*(L\otimes_{\mathbb Q_p} O_{S_K^o})\otimes_{O_{S_K}}\mathbb B_{dr,S_K}
\in C_{\mathbb B_{dr}}(S_K^{an,pet})
\end{eqnarray*}
so that we have the isomorphism in $D_{\mathbb B_{dr}}(S_K^{an,pet})$
\begin{eqnarray*}
m(L\otimes_{\mathbb Q_p}\mathbb B_{dr,S^o_K}):
\psi_D(L\otimes_{\mathbb Q_p}\mathbb B_{dr,S^o_K})\otimes_{\mathbb B_{dr,S_K}}\mathbb B_{dr,\psi_D,K}
\to \psi_DL\otimes_{\mathbb Q_p}\mathbb B_{dr,\psi_D,K} \\ 
\xrightarrow{:=}
\psi_DL\otimes_{\mathbb Q_p}F^0DR(S)(\Gr^{V_D}_{-1\leq\alpha<0}(O_{S^o},F_b)\otimes_{O_S}(O\mathbb B_{dr,S},F)), \;
s\otimes h\otimes w\mapsto s\otimes hw
\end{eqnarray*}
More generally, for $Z\subset S$ a closed subset and $L\in\Shv_{\mathbb Z_p}(Z^{o,et})$ a local system with $Z^o:=Z\cap S^o$, 
we set using theorem \ref{HSk} for $j_*(L\otimes O_{Z_K^o})$
\begin{eqnarray*}
\psi_D(L\otimes_{\mathbb Q_p}\mathbb B_{dr,S^o_K}):=
\Gr^{V_D}_{-1\leq\alpha<0}j_*(L\otimes_{\mathbb Q_p}O_{Z_K^o})\otimes_{O_{S_K}}\mathbb B_{dr,S_K}
\in C_{\mathbb B_{dr}}(S_K^{an,pet})
\end{eqnarray*}
so that we have the isomorphism in $D_{\mathbb B_{dr}}(S_K^{an,pet})$
\begin{eqnarray*}
m(L\otimes_{\mathbb Q_p}\mathbb B_{dr,S^o_K}):
\psi_D(L\otimes_{\mathbb Q_p}\mathbb B_{dr,S^o_K})\otimes_{\mathbb B_{dr,S_K}}\mathbb B_{dr,\psi_D,K}
\to\psi_DL\otimes_{\mathbb Q_p}\mathbb B_{dr,\psi_D,K} \\ 
\xrightarrow{:=}
\psi_DL\otimes_{\mathbb Q_p}F^0DR(S)(\Gr^{V_D}_{-1\leq\alpha<0}(O_{S^o},F_b)\otimes_{O_S}(O\mathbb B_{dr,S},F)), \;
s\otimes h\otimes w\mapsto s\otimes hw
\end{eqnarray*}
\end{itemize}
\end{defi}

We then give using definition \ref{VfilKMmap} and the local system case
an inverse functor to the De Rham functor for De Rham modules

\begin{defi}\label{Bdr}
\begin{itemize}
\item[(i0)]Let $S\in\SmVar(k)$ irreducible. 
We then have the morphism of site $\an_S:S_K^{an,pet}\to S^{et}$ given by the analytical functor.
Let $K\in P_{\mathbb Z_p,k}(S^{et})$ a perverse sheaf, 
in particular there exist an open subset $S^o\subset S$ with $D:=S\backslash S^o$ a (Cartier) divisor
such that $K_{|S^o}:=j^*K\in C(S^{o,et})$ is a local system for the etale topology,
where we denote $j:S_0\hookrightarrow S$ the open embedding and $i:D\hookrightarrow S$ the closed embedding of the Cartier divisor.
Assume first that $K_{|D}:=i^*K$ is a local system. Then, $\psi_DK,\phi_DK\in C_{\mathbb Z_p}(D^{et})$ are local systems.
We denote again $K:=\an_S^*K\in C(S_K^{an,pet})$ and $K:=j^*K\in C(S_K^{o,an,pet})$.
Denote by $\pi:\tilde S_K^{o,an}\to S_K^{o,an}$ the perfectoid universal covering.
We then have by theorem \ref{PSk} a canonical isomorphism in $D_{\mathbb Z_p,c}(S_K^{an,pet})$
\begin{eqnarray*}
Is(K):K\xrightarrow{\sim}I_D(K):=
(\psi^u_DK\xrightarrow{(c(x_{S^o/S}(K)),can(K))}x_{S^o/S}(K)\oplus\phi^u_DK\xrightarrow{((0,T-I),var(K))}\psi^u_DK)
\end{eqnarray*}
We then set using definition \ref{VfilKMmap} 
\begin{eqnarray*}
\mathbb B_{dr,S}(x_{S^o/S}(K)):=
\Cone((\rho^u_{D,DR}(K\otimes_{\mathbb Q_p}O_{S^o})\circ i_{V0})\otimes\rho_{\mathbb B_{dr,D}}(O_S): \\
V_{D0}j_*(K\otimes_{\mathbb Q_p}O_{S^o})\otimes_{O_S}\mathbb B_{dr,S^o/S,K}\to
\psi^u_D(K\otimes_{\mathbb Q_p}O_{S^o})\otimes_{O_S}\mathbb B_{dr,\psi_D,K})
\in D_{\mathbb B_{dr}}(S_K^{an,pet})
\end{eqnarray*}
and
\begin{eqnarray*}
\mathbb B_{dr,S}(K):=\mathbb B_{dr,S}(I_D(K)):= 
(\mathbb B_{dr,S}(\psi^u_DK):=\psi^u_D(K\otimes_{\mathbb Q_p}O_{S^o})\otimes_{O_S}\mathbb B_{dr,\psi_D,K} \\ 
\xrightarrow{(c(\mathbb B_{dr,S}(x_{S^o/S}(K)))\otimes c(\mathbb B_{dr,x_{S^o/S},K}),
can(K\otimes O_{S^o})\otimes can_{\mathbb B_{dr},D}(O_S))} \\
\mathbb B_{dr,S}(x_{S^o/S}(K))\oplus\mathbb B_{dr,S}(\phi^u_DK):=
\mathbb B_{dr,S}(x_{S^o/S}(K))\oplus\phi^u_D(K\otimes_{\mathbb Q_p} O_{S^o})\otimes_{O_S}\mathbb B_{dr,\phi_D,K} \\
\xrightarrow{((0,exp(s\partial s+1)\otimes(0,exp(s\partial s+1)),
var(K\otimes O_{S^o})\otimes var_{\mathbb B_{dr},D}(O_S))} \\
\mathbb B_{dr,S}(\psi^u_DK):=\psi^u_D(K\otimes_{\mathbb Q_p}O_{S^o})\otimes_{O_S}\mathbb B_{dr,\psi_D,K})
\in D_{\mathbb B_{dr}}(S_K^{an,pet}).
\end{eqnarray*}
where we write for short $O_S:=O_{S_K}$.
\item[(i)]Let $S\in\SmVar(k)$ irreducible. 
We then have the morphism of site $\an_S:S^{an,pet}\to S^{et}$ given by the analytical functor.
We define the functor
\begin{equation*}
\mathbb B_{dr,S}:D_{\mathbb Z_pfil,c,k}(S^{et})\to D_{\mathbb B_{dr,S}fil}(S_K^{an,pet})
\end{equation*}
using the nearby and vanishing cycle functors.
For $(K,W)\in P_{pfil}(S^{et})$ a filtered perverse sheaf and 
$(D_1,\ldots,D_d)\in\mathcal S(K)$ a stratification by (Cartier) divisor $D_i\subset S$, $1\leq i\leq d$ such that 
$K_{|D(r)\backslash D(r+1)}:=l_r^*K\in D_{\mathbb Z_p,c}(D(r)\backslash D(r+1)^{et})$ are local systems for all $1\leq r\leq d$, 
where $D(r):=\cap_{1\leq i\leq r}D_i$ and $l_r:D(r)\backslash D(r+1)\hookrightarrow S$ is the locally closed embedding,
denoting $S^o:=S\backslash D_1$, we have by theorem \ref{PSk} a canonical isomorphism in $D_{\mathbb Z_p,c}(S_K^{an,pet})$
\begin{eqnarray*}
Is(K):K\xrightarrow{\sim}I_{D_d}\circ\cdots\circ I_{D_1}(K):&=&(I_{D_d}(\cdots(I_{D_1}(K)) \\
&=&(\cdots\to\xi_{D_d}\cdots\xi_{D_1}(K)
\to\cdots)_{\xi:[1,\cdots,d]\to\left\{\psi^u,V_0j_*\oplus\phi^u,\psi^u\oplus\phi^u,\psi^u\right\}},
\end{eqnarray*}
where for $D\subset S$ a (Cartier) divisor, denoting $j:S\backslash D\hookrightarrow S$ the open embedding, 
\begin{equation*}
\xi_D(K)=\psi_D^u(K) \; \; \mbox{or} \; \;  \xi_D(K)=\psi_D^u(K)\oplus\phi_D^u(K)
\; \; \mbox{or} \; \; \xi_D(K)=(V_{D,0}j_*K)\oplus\phi_D^u(K).
\end{equation*}
We then define by (i0) inductively
\begin{eqnarray*}
\mathbb B_{dr,S}(K,W):=
\lim_{(D_1,\ldots,D_d)\in\mathcal S(K)}\mathbb B_{dr,S}(I_{D_d}\circ\cdots\circ I_{D_1}(K)):= \\
(\cdots\to\xi_{D_d}\cdots\xi_{D_1}((K_1,W)\otimes_{\mathbb Q_p}O_{S^o})\otimes_{O_S} 
\mathbb B_{dr,S\backslash D_1/S,K}\otimes_{\mathbb B_{dr,S_K}}\cdots \\
\otimes_{\mathbb B_{dr,S_K}}\mathbb B_{dr,S\backslash D_{i_1}/S,K}\otimes_{\mathbb B_{dr,S_K}} 
\mathbb B_{dr,\psi_{D_{i_{r+1}}},K}\otimes_{\mathbb B_{dr,S_K}}
\cdots\otimes_{\mathbb B_{dr,S_K}}\mathbb B_{dr,\psi_{D_{i_s}},K} \\
\otimes_{\mathbb B_{dr,S_K}}\mathbb B_{dr,\phi_{D_{i_{s+1}}},K}\otimes_{\mathbb B_{dr,S_K}}\cdots
\otimes_{\mathbb B_{dr,S_K}}\mathbb B_{dr,\phi_{D_{i_d}},K}
\to\cdots)_{\xi:[1,\cdots,d]\to\left\{\psi^u,V_0j_*\oplus\phi^u,\psi^u\oplus\phi^u,\psi^u\right\}}
\end{eqnarray*}
where for $D\subset S$ a (Cartier) divisor, denoting $j:S\backslash D\hookrightarrow S$ the open embedding,
and $(L,W)\in\PSh_Z(S^{et})$ with $L$ a local system on a closed subset $Z\subset S$, 
\begin{eqnarray*}
\xi_D((L,W)\otimes_{\mathbb Q_p}O_Z)=\psi_D^u((L,W)\otimes_{\mathbb Q_p}O_Z) 
\; \; \mbox{or} \\ 
\xi_D((L,W)\otimes_{\mathbb Q_p}O_Z)=\psi_D^u((L,W)\otimes_{\mathbb Q_p}O_Z)\oplus\phi_D^u((L,W)\otimes_{\mathbb Q_p}O_Z) 
\; \; \mbox{or} \\
\xi_D((L,W)\otimes_{\mathbb Q_p}O_Z)=
(V_{D,0}j_{*w}((L,W)\otimes_{\mathbb Q_p}O_Z))\oplus\phi_D^u((L,W)\otimes_{\mathbb Q_p}O_Z).
\end{eqnarray*}
where we denote for short $O_Z:=O_{Z_K}$.
For $m:(K_1,W)\to (K_2,W)$ a morphism with $(K_1,W),(K_2,W)\in P_{pfil}(S^{et})$, considering a stratification 
$(D_1,\ldots,D_d)\in\mathcal S(K_1)\cap\mathcal S(K_2)$ by (Cartier) divisor $D_i\subset S$, $1\leq i\leq d$ 
such that $K_{1|D(r)\backslash D(r+1)},K_{2|D(r)\backslash D(r+1)}\in D_{\mathbb Z_p,c}(D(r)\backslash D(r+1)^{et})$ 
are local systems for all $1\leq r\leq d$, we get
\begin{eqnarray*}
\mathbb B_{dr,S}(m): \\
\mathbb B_{dr,S}(K_1,W):=(\cdots\to\xi_{D_d}\cdots\xi_{D_1}((K_1,W)\otimes_{\mathbb Q_p}O_{S^o}) 
\otimes_{O_S}\mathbb B_{dr,S\backslash D_{i_1}/S,K}\otimes_{\mathbb B_{dr,S_K}}\cdots \\
\otimes_{\mathbb B_{dr,S_K}}\mathbb B_{dr,S\backslash D_{i_r}/S,K}\otimes_{\mathbb B_{dr,S_K}} 
\mathbb B_{dr,\psi_{D_{i_{s+1}}},K}\otimes_{\mathbb B_{dr,S_K}} 
\cdots\otimes_{\mathbb B_{dr,S_K}}\mathbb B_{dr,\psi_{D_{i_s}},K} \\
\otimes_{\mathbb B_{dr,S_K}}\mathbb B_{dr,\phi_{D_{i_{s+1}}},K}\otimes_{\mathbb B_{dr,S_K}}\cdots
\otimes_{\mathbb B_{dr,S_K}}\mathbb B_{dr,\phi_{D_{i_d}},K}\to\cdots) \\
\xrightarrow{(\xi_{D_d}\cdots\xi_{D_1}(m)\otimes I)} \\
\mathbb B_{dr,S}(K_2,W):=(\cdots\to\xi_{D_d}\cdots\xi_{D_1}((K_2,W)\otimes_{\mathbb Q_p}O_{S^o}) 
\otimes_{O_S}\mathbb B_{dr,S\backslash D_{i_1}/S},K\otimes_{\mathbb B_{dr,S_K}}\cdots \\
\otimes\mathbb B_{dr,S\backslash D_{i_r}/S,K}\otimes_{\mathbb B_{dr,S_K}} 
\mathbb B_{dr,\psi_{D_{i_{r+1}}},K}\otimes_{\mathbb B_{dr,S_K}}\cdots 
\otimes_{\mathbb B_{dr,S_K}}\mathbb B_{dr,\psi_{D_{i_s}},K} \\
\otimes_{\mathbb B_{dr,S_K}}\mathbb B_{dr,\phi_{D_{i_{s+1}}},K}\otimes_{\mathbb B_{dr,S_K}}\cdots
\otimes_{\mathbb B_{dr,S_K}}\mathbb B_{dr,\phi_{D_{i_d}},K}\to\cdots).
\end{eqnarray*} 
Note that if $L$ is a local system on $S$ then $\mathbb B_{dr,S}(L)=L\otimes\mathbb B_{dr,S,K}$, 
that is it does NOT depend on the choice of a stratification (see remark \ref{Bdrfunctrem}). This gives the functor
\begin{equation*}
\mathbb B_{dr,S}:D_{\mathbb Z_pfil,c,k}(S^{et})=D(P_{pfil}(S^{et}))\to D_{\mathbb B_{dr,S}fil}(S_K^{an,pet}).
\end{equation*}
\item[(ii)]Let $S\in\Var(k)$. Let $S=\cup_{i\in I}S_i$ an open cover such that there
exists closed embeddings $i_i:S_i\hookrightarrow\tilde S_i$ with $\tilde S_I\in\SmVar(k)$.
We define as in (i) the functor 
\begin{eqnarray*}
\mathbb B_{dr,(\tilde S_I)}:D_{\mathbb Z_pfil,c,k}(S^{et})\xrightarrow{T(S/(\tilde S_I))}
D_{\mathbb Z_pfil,c,k}(S^{et}/(\tilde S_I^{et}))\to D_{\mathbb B_{dr}fil}(S_K^{an,pet}/(\tilde S_{I,K}^{an,pet})), \\ 
(K,W)\mapsto\mathbb B_{dr,(\tilde S_I)}(K,W):=\mathbb B_{dr,(\tilde S_I)}(i_{I*}j_I^*(K,W),t_{IJ}) 
\end{eqnarray*}
with for $(K,W)\in P_{pfil}(S^{et})$,
\begin{eqnarray*}
\mathbb B_{dr,(\tilde S_I)}(K,W):=\lim_{(D_1,\ldots,D_d)\in\mathcal S(K)} 
((\xi_{\tilde D_{d,I}}\cdots\xi_{\tilde D_{1,I}}(i_{I*}j_I^*((K,W)\otimes_{\mathbb Q_p}O_{S^o})) \\
\otimes_{O_S}
\mathbb B_{dr,\tilde S_I\backslash\tilde D_{i_1,I}/S,K}\otimes_{\mathbb B_{dr,\tilde S_{I,K}}}\cdots
\otimes_{\mathbb B_{dr,\tilde S_{I,K}}}\mathbb B_{dr,\tilde S_I\backslash D_{i_r,I}/\tilde S_I,K}
\otimes_{\mathbb B_{dr,\tilde S_{I,K}}} \\
\mathbb B_{dr,\phi_{\tilde D_{i_{r+1},I}},K}\otimes_{\mathbb B_{dr,\tilde S_{I,K}}}
\cdots\otimes_{\mathbb B_{dr,\tilde S_{I,K}}} 
\mathbb B_{dr,\phi_{\tilde D_{i_s,I}},K}\otimes_{\mathbb B_{dr,\tilde S_{I,K}}} \\
\mathbb B_{dr,\psi_{\tilde D_{i_{s+1},I}},K}\otimes_{\mathbb B_{dr,\tilde S_{I,K}}}\cdots
\otimes_{\mathbb B_{dr,\tilde S_{I,K}}}\mathbb B_{dr,\psi_{\tilde D_{i_d,I}},K}
\to\cdots)_{\xi:[1,\cdots,d]\to\left\{\psi^u,V_0j_*\oplus\phi^u,\psi^u\oplus\phi^u,\psi^u\right\}},
\mathbb B_{dr}(t_{IJ}))
\end{eqnarray*}
with $(D_1,\ldots,D_d)\in\mathcal S(K)$ stratifications by Cartier divisor $D_i\subset S$, $1\leq i\leq d$ such that 
\begin{equation*}
K_{|D(r)\backslash D(r+1)}:=l_r^*K\in D_{\mathbb Z_p,c}(D(r)\backslash D(r+1)^{et})
\end{equation*}
are local systems for all $1\leq r\leq d$, 
and $\tilde D_{s,I}\subset\tilde S_I$ are (Cartier) divisor such that $D_s\cap S_I\subset\tilde D_{s,I}\cap S$
(that is $D_s\cap S_I$ is a union of irreducible components of $\tilde D_{S,I}\cap S$ which are (Cartier) divisors),
having by theorem \ref{PSk} the canonical isomorphism in $D_{\mathbb Z_p,c,k}(S_K^{an,pet})$
\begin{eqnarray*}
Is(K):K\xrightarrow{\sim}I_{D_d}\circ\cdots\circ I_{D_1}(K):=
(\cdots\to\xi_{D_d}\cdots\xi_{D_1}(K)
\to\cdots)_{\xi:[1,\cdots,d]\to\left\{\psi^u,V_0j_*\oplus\phi^u,\psi^u\oplus\phi^u,\psi^u\right\}}.
\end{eqnarray*}
\end{itemize}
\end{defi}

\begin{rem}\label{Bdrfunctrem}
Let $S\in\SmVar(k)$. Let $K\in P_{\mathbb Z_p,k}(S^{et})$.
$(D_1,\ldots,D_d)\in\mathcal S(K)$ stratifications by Cartier divisor $D_i\subset S$, $1\leq i\leq d$ such that 
$K_{|D(r)\backslash D(r+1)}:=l_r^*K\in D_{\mathbb Z_p,c}(D(r)\backslash D(r+1)^{et})$ are local systems for all $1\leq r\leq d$.
We then have the canonical map in $D_{\mathbb B_{dr}}(S_K^{an,pet})$
\begin{eqnarray*}
a_S(K,\mathbb B_{dr}):=
(I\otimes a_S(\mathbb B_{dr,S\backslash D_i/S}), 
I\otimes a_S(\mathbb B_{dr,\phi_{D_i}}),I\otimes a_S(\mathbb B_{dr,\psi_{D_i}}))\circ(I(K)\otimes I): \\
K\otimes_{\mathbb Q_p}\mathbb B_{dr,S_K}\to 
\mathbb B_{dr,S}(K):=(\cdots\to\xi_{D_d}\cdots\xi_{D_1}((K_2,W)\otimes_{\mathbb Q_p}O_{S^o}) \\
\otimes_{O_S}\mathbb B_{dr,S\backslash D_{i_1}/S,K}\otimes_{\mathbb B_{dr,S_K}}\cdots
\otimes_{\mathbb B_{dr,S_K}}\mathbb B_{dr,S\backslash D_{i_r}/S,K}\otimes_{\mathbb B_{dr,S_K}} \\
\mathbb B_{dr,\phi_{D_{i_{r+1}}}}\otimes_{\mathbb B_{dr,S_K}}\cdots
\otimes_{\mathbb B_{dr,S_K}}\mathbb B_{dr,\phi_{D_{i_s}},K}\otimes_{\mathbb B_{dr,S_K}}
\mathbb B_{dr,\psi_{D_{i_{s+1}}},K}\otimes_{\mathbb B_{dr,S_K}}\cdots
\otimes_{\mathbb B_{dr,S_K}}\mathbb B_{dr,\psi_{D_{i_d}},K}\to\cdots).
\end{eqnarray*}
On the other hand we have by theorem \ref{RHpadic},
\begin{eqnarray*}
\alpha(S):K\otimes_{\mathbb Q_p}\mathbb B_{dr,S_K}\xrightarrow{\sim}
F^0DR(S)((K\otimes O_{S_K^{an}},F_b)\otimes_{O_{S_K^{an}}}(O\mathbb B_{dr,S_K},F))
\end{eqnarray*}
in $D_{\mathbb B_{dr}}(S_K^{an,pet})$.
\begin{itemize}
\item[(i)] If $K\in P_{\mathbb Z_p,k}(S^{et})$ is a local system then the map $a_S(K,\mathbb B_{dr})$ is an isomorphism
since by proposition \ref{PSkMHM}
\begin{eqnarray*}
I(K\otimes O_S):(K\otimes O_S)\xrightarrow{\sim}I_{D_1}\circ I_{D_d}(K\otimes O_S):=
(\cdots\to\xi_{D_d}\cdots\xi_{D_1}(O_S)\to\cdots)
\end{eqnarray*}
in $D(DRM(S))$ and since the functor
\begin{equation*}
K\otimes_{\mathbb Q_p}(-):C_{\mathbb Z_p}(S_K^{an,pet})\to C_{\mathbb Z_p}(S_K^{an,pet}), \; 
N\mapsto K\otimes_{\mathbb Q_p} N
\end{equation*}
respect etale hence pro-etale equivalences.
\item[(ii)] If $K\in P_{\mathbb Z_p,k}(S^{et})$ is NOT a local system 
then the map $a_S(K,\mathbb B_{dr})$ is NOT an isomorphism is general.
For example, for $j:S^o\hookrightarrow S$ an open embedding with $D:=S\backslash S^o$ a Cartier divisor, 
we have in $D_{\mathbb B_{dr}}(S_K^{an,pet})$, by proposition \ref{PSkMHM}, 
\begin{eqnarray*}
T(j,\mathbb B_{dr})(\mathbb Z_{p,S^o}):=m\circ(0,(I,0),0)^{-1}: 
\mathbb B_{dr,S}(j_*\mathbb Z_{p,S^o})\xrightarrow{:=} 
(\phi^u_D(O_S)\otimes_{O_S}\mathbb B_{dr,\psi_D,K}\to \\
(\Cone(V_{D,0}j_*O_{S^o}\otimes\mathbb B_{dr,S^o/S,K}\to\psi^u_D(O_S)\otimes_{O_S}\mathbb B_{dr,\psi_D,K}))
\oplus(\phi^u_D(O_S)\otimes_{O_S}\mathbb B_{dr,\phi_D,K})\to \\ 
\psi^u_D(O_S)\otimes_{O_S}\mathbb B_{dr,\psi_D,K})
\xrightarrow{\sim} 
\mathbb B_{dr,S^o/S,K}:=F^0DR(S)(j_{*Hdg}(O_{S^o},F_b)^{an}\otimes_{O_{S_K^{an}}}(O\mathbb B_{dr,S_K},F)) 
\end{eqnarray*}
which is, by theorem \ref{RHpadic} 
(in the case of $D$ a normal crossing divisor, $\mathbb B_{dr,S^o/S,K}=\mathbb B_{dr,S_K}(\log D_K)$),
NOT isomorphic in $D_{\mathbb B_{dr}}(S_K^{an,pet})$ to 
\begin{eqnarray*}
j_*\alpha(S^o):j_*\mathbb B_{dr,S^o_K}\xrightarrow{\sim}
F^0DR(S)(j_*(O_S,F_b)^{an}\otimes_{O_{S_K^{an}}}(O\mathbb B_{dr,S_K},F)),
\end{eqnarray*}
and also NOT isomorphic to
\begin{eqnarray*}
\mathbb D^v_ST_!(j,\otimes)(-,-):(j_*\mathbb Z_{p,S_K^o})\otimes\mathbb B_{dr,S_K}\xrightarrow{\sim}
\mathbb D^v_S(j_!\mathbb B_{dr,S_K})
\end{eqnarray*}
see also remark \ref{PFetrem}.
If $K\in P_{\mathbb Z_p,k}(S^{et})$ is NOT a local system, the functor
\begin{equation*}
K\otimes_{\mathbb Q_p}(-):C_{\mathbb Z_p}(S_K^{an,pet})\to C_{\mathbb Z_p}(S_K^{an,pet}), \; 
N\mapsto K\otimes_{\mathbb Q_p} N
\end{equation*}
does NOT preserve etale or pro-etale equivalence.
Note also that if $K\in P_{\mathbb Z_p,k}(S^{et})$ is NOT a local system, 
$K\otimes O_S$ is NOT a quasi-coherent $O_S$ module hence NOT an holonomic $D_S$ module.
Recall also that the filtered De Rham functor does NOT commutes with filtered tensor product in general
(it may leads to different F-filtration).
\end{itemize}
\end{rem}

Let $k\subset K\subset\mathbb C_p$ a subfield of a p-adic field.
Let $S\in\SmVar(k)$ and $D\subset S$ a Cartier divisor. Denote $S^o:=S\backslash D$.
We write for simplicity, 
\begin{itemize}
\item $\mathbb B_{dr,S}:=\mathbb B_{dr,S_K}$, $O\mathbb B_{dr,S}:=O\mathbb B_{dr,S_K}$,
\item $\mathbb B_{dr,S^o/S}:=\mathbb B_{dr,S^o/S,K}$,
$\mathbb B_{dr,\psi_D}:=\mathbb B_{dr,\psi_D,K}$ and $\mathbb B_{dr,\phi_D}:=\mathbb B_{dr,\phi_D,K}$
\item $\mathbb B_{dr,x_{S^o/S}}:=\mathbb B_{dr,x_{S^o/S},K}$.
\end{itemize}

We now look at the functorialities with respect to the proper morphisms and with respect to the open embeddings :

Let $f:X\to S$ a morphism with $X,S\in\SmVar(k)$. 
Let $D\subset S$ a (Cartier) divisor and denote $S^o:=S\backslash D$ and $X^o:=X\backslash f^{-1}(D)$.  
Denote $j:S^o\hookrightarrow S$, $j':X^o\hookrightarrow X$ the open embeddings.
\begin{itemize}
\item We have the following quasi-isomorphism in $C_{\mathbb B_{dr,S}}(X_K^{an,pet})$
\begin{eqnarray*}
m_f(\mathbb B_{dr,S^o/S}):f^*\mathbb B_{dr,S^o/S}\otimes_{f^*\mathbb B_{dr,S}}\mathbb B_{dr,X} 
\xrightarrow{I\otimes\alpha(X_K)} \\
f^*F^0DR(S)(j_{*Hdg}(O_{S^o},F_b)^{an}\otimes_{O_S}(O\mathbb B_{dr,S},F))\otimes_{f^*\mathbb B_{dr,S}} 
F^0DR(X)((O_X,F_b)^{an}\otimes_{O_X}(O\mathbb B_{dr,X}),F) \\
\xrightarrow{=} 
F^0DR(f^*O_S)(f^*j_{*Hdg}(O_{S^o},F_b)^{an}\otimes_{O_S}(O\mathbb B_{dr,S},F))\otimes_{f^*\mathbb B_{dr,S}} \\
F^0DR(X)((O_X,F_b)^{an}\otimes_{O_X}(O\mathbb B_{dr,X}),F) \\
\xrightarrow{w_X\circ\Omega_{f^*O_S/O_X}(-)}
F^0DR(X)(f^{*mod}j_{*Hdg}(O_{S^o},F_b)^{an}\otimes_{O_X}(O\mathbb B_{dr,X},F)) \\
\xrightarrow{=}
F^0DR(X)(j'_{*Hdg}(O_{X^o},F_b)^{an}\otimes_{O_X}(O\mathbb B_{dr,X}),F)=:\mathbb B_{dr,X^o/X}.
\end{eqnarray*}
\item We have the following quasi-isomorphism in $C_{\mathbb B_{dr,S}}(X_K^{an,pet})$
\begin{eqnarray*}
m_f(\mathbb B_{dr,\psi_D}):f^*\mathbb B_{dr,\psi_D}\otimes_{f^*\mathbb B_{dr,S}}\mathbb B_{dr,X} 
\xrightarrow{I\otimes\alpha(X_K)} \\
f^*F^0DR(S)(\psi_D(O_{S^o},F_b)^{an}\otimes_{O_S}(O\mathbb B_{dr,S},F))\otimes_{f^*\mathbb B_{dr,S}} 
F^0DR(X)((O_X,F_b)^{an}\otimes_{O_X}(O\mathbb B_{dr,X}),F) \\
\xrightarrow{=} 
F^0DR(f^*O_S)(f^*\psi_D(O_{S^o},F_b)^{an}\otimes_{O_S}(O\mathbb B_{dr,S},F))\otimes_{f^*\mathbb B_{dr,S}} \\
F^0DR(X)((O_X,F_b)^{an}\otimes_{O_X}(O\mathbb B_{dr,X}),F) \\
\xrightarrow{w_X\circ\Omega_{f^*O_S/O_X}(-)}
F^0DR(X)(f^{*mod}\psi_D(O_{S^o},F_b)^{an}\otimes_{O_X}(O\mathbb B_{dr,X},F)) \\
\xrightarrow{=}
F^0DR(X)(\psi_{f^{-1}(D)}(O_{X^o},F_b)^{an}\otimes_{O_X}(O\mathbb B_{dr,X},F))=:\mathbb B_{dr,\psi_{f^{-1}(D)}}.
\end{eqnarray*}
\item We have the following quasi-isomorphism in $C_{\mathbb B_{dr,S}}(X_K^{an,pet})$
\begin{eqnarray*}
m_f(\mathbb B_{dr,\psi_D}):f^*\mathbb B_{dr,\psi_D}\otimes_{f^*\mathbb B_{dr,S}}\mathbb B_{dr,X} 
\xrightarrow{I\otimes\alpha(X_K)} \\
f^*F^0DR(S)(\phi_D(O_{S^o},F_b)^{an}\otimes_{O_S}(O\mathbb B_{dr,S},F))\otimes_{f^*\mathbb B_{dr,S}} 
F^0DR(X)((O_X,F_b)^{an}\otimes_{O_X}(O\mathbb B_{dr,X}),F) \\
\xrightarrow{=} 
F^0DR(f^*O_S)(f^*\phi_D(O_{S^o},F_b)^{an}\otimes_{O_S}(O\mathbb B_{dr,S},F))\otimes_{f^*\mathbb B_{dr,S}} \\
F^0DR(X)((O_X,F_b)^{an}\otimes_{O_X}(O\mathbb B_{dr,X}),F) \\
\xrightarrow{w_X\circ\Omega_{f^*O_S/O_X}(-)}
F^0DR(X)(f^{*mod}\phi_D(O_{S^o},F_b)^{an}\otimes_{O_X}(O\mathbb B_{dr,X},F)) \\
\xrightarrow{=}
F^0DR(X)(\phi_{f^{-1}(D)}(O_{X^o},F_b)^{an}\otimes_{O_X}(O\mathbb B_{dr,X},F))=:\mathbb B_{dr,\phi_{f^{-1}(D)}}.
\end{eqnarray*}
\item We have the following quasi isomorphism in $C_{\mathbb B_{dr,S}}(S_K^{an,pet})$
\begin{eqnarray*}
m_f(\mathbb B_{dr,x_{S^o/S}}):=(m_f(\mathbb B_{dr,S^o/S}),m_f(\mathbb B_{dr,\psi_D})):
f^*\mathbb B_{dr,x_{S^o/S}}\otimes_{f^*\mathbb B_{dr,S}}\mathbb B_{dr,X}\to\mathbb B_{dr,x_{X^o/X}}
\end{eqnarray*}
\end{itemize}

\begin{defi}\label{TfjBdr1}
\begin{itemize}
\item[(i)]Let $f:X\to S$ be a proper morphism with $X,S\in\SmVar(k)$. 
Let $(K,W)\in P_{\mathbb Z_pfil,k}(X^{et})\cap D_{\mathbb Z_pfil,c,k,gm}(X^{et})$ be a filtered perverse sheaf
of geometric origin, i.e. for all $n\in\mathbb Z$, there exists a proper morphism $f'_n:X'_n\to X$ with $X'_n\in\SmVar(k)$
and $k\in\mathbb Z$ such that $\Gr^n_WK={\,}^pR^kf'_*\mathbb Z_{X'_n}\in D_{\mathbb Z_p,c,k,gm}(X^{et})$. 
We have then, by the perverse hard Lefchetz thorem, a canonical isomorphism in $D_{\mathbb Z_pfil,c,k}(S^{et})$ 
\begin{equation*}
l(K,W):(K,W)\to\oplus_{k\in\mathbb Z}{\,}^pR^kf_*(K,W).
\end{equation*}
Consider a stratification $(E_1,\ldots,E_d)\in\mathcal S(K)$ by (Cartier) divisor $E_i\subset X$, $1\leq i\leq d$, such that 
\begin{equation*}
K_{|E(r)\backslash E(r+1)}:=l_r^*K\in D_{\mathbb Z_p,c}((E(r)\backslash E(r+1))^{et}) 
\end{equation*}
are local systems for all $1\leq r\leq d$, $l_r:E(r)\hookrightarrow X$ being the locally closed embeddings.
Let $k\in\mathbb Z$. 
Take a stratification $(D_1,\ldots,D_e)\in\mathcal S(K)$ by (Cartier) divisor $D_i\subset S$, $1\leq i\leq e$, such that 
\begin{eqnarray*}
({\,}^pR^kf_*\xi_{E_d}\cdots\xi_{E_1}K)_{D(r')\backslash D(r'+1)}:=m_{r'}^*{\,}^pR^kf_*\xi_{E_d}\cdots\xi_{E_1}K
\in D_{\mathbb Z_p,c}((D(r')\backslash D(r'+1))^{et})
\end{eqnarray*}
are local systems for all $\xi:[1,\cdots,d]\to\left\{\psi^u,V_0j_*\oplus\phi^u,\psi^u\oplus\phi^u,\psi^u\right\}$ 
and all $1\leq r'\leq d$, $m_{r'}:D(r')\hookrightarrow S$ being the locally closed embeddings. This implies that 
\begin{equation*}
{\,}^pR^kf_*K_{|D(r')\backslash D(r'+1)}:=m_{r'}^*{\,}^pR^kf_*K\in D_{\mathbb Z_p,c}((D(r')\backslash D(r'+1))^{et})
\end{equation*}
are local systems for all $1\leq r'\leq e$. We then define the canonical maps in $D_{\mathbb B_{dr,S}fil}(S_K^{an,pet})$
\begin{eqnarray*}
T^k(f,\mathbb B_{dr})(K,W):\mathbb B_{dr,S}({\,}^pR^kf_*(K,W))\xrightarrow{B_{dr,S}({\,}^pR^kf_*Is(K,W))} \\ 
(\cdots\to
\mathbb B_{dr,S}({\,}^pR^kf_*\xi_{E_d}\cdots\xi_{E_1}(K,W)\to\cdots) \\
\xrightarrow{:=} 
(\cdots\to\xi_{D_e}\cdots\xi_{D_1}({\,}^pR^kf_*\xi_{E_d}\cdots\xi_{E_1}(K,W)\otimes_{\mathbb Q_p}O_{S\backslash D_1})
\otimes_{O_S} 
\mathbb B_{dr,S\backslash D_{j_1}/S}\otimes\cdots\otimes\mathbb B_{dr,S\backslash D_{j_{r'}}/S} \\ 
\otimes_{\mathbb B_{dr,S}}\mathbb B_{dr,\phi_{D_{j_{r'+1}}}} 
\otimes\cdots\otimes\mathbb B_{dr,\phi_{D_{j_{s'}}}} 
\otimes_{\mathbb B_{dr,S}}\mathbb B_{dr,\psi_{D_{j_{s'+1}}}}\otimes\cdots\otimes\mathbb B_{dr,\psi_{D_{j_e}}}\to\cdots) \\
\xrightarrow{(T(f,\otimes)(-,-)\circ l_k(-))\otimes I} \\
(\cdots\to\xi_{D_e}\cdots\xi_{D_1}({\,}
Rf_*(\xi_{E_d}\cdots\xi_{E_1}(K,W)\otimes_{\mathbb Q_p}f^*O_{S\backslash D_1}\otimes_{f^*O_S} 
f^*\mathbb B_{dr,S\backslash D_{j_1}/S}\otimes\cdots\otimes f^*\mathbb B_{dr,S\backslash D_{j_{r'}}/S} \\ 
\otimes_{f^*\mathbb B_{dr,S}}f^*\mathbb B_{dr,\phi_{D_{j_{r'+1}}}} 
\otimes\cdots\otimes f^*\mathbb B_{dr,\phi_{D_{j_{s'}}}}\otimes_{f^*\mathbb B_{dr,S}}f^*\mathbb B_{dr,\psi_{D_{j_{s'+1}}}}
\otimes\cdots\otimes f^*\mathbb B_{dr,\psi_{D_{j_e}}}))\to\cdots) \\
\xrightarrow{Rf_*(I\otimes a_X(\mathbb B_{dr,X\backslash E_{i}/X})\otimes a_X(\mathbb B_{dr,\phi_{E_i}})
\otimes a_X(\mathbb B_{dr,\psi_{E_i}}))} \\
(\cdots\to\xi_{D_e}\cdots\xi_{D_1}(
Rf_*(\xi_{E_d}\cdots\xi_{E_1}((K,W)\otimes_{\mathbb Q_p}O_{X\backslash E_1})\otimes_{O_X} \\
\mathbb B_{dr,X\backslash E_{i_1}/X}\otimes_{\mathbb B_{dr,X}}\otimes\cdots\otimes\mathbb B_{dr,X\backslash E_{i_r}/X} \\
\otimes_{\mathbb B_{dr,X}}\mathbb B_{dr,\phi_{E_{i_{r+1}}}}\otimes\cdots\otimes\mathbb B_{dr,\phi_{E_{i_s}}}
\otimes_{\mathbb B_{dr,X}}\mathbb B_{dr,\psi_{E_{i_{s+1}}}}\otimes\cdots\otimes\mathbb B_{dr,\phi_{E_{i_d}}} \\
\otimes_{f^*\mathbb B_{dr,S}}f^*\mathbb B_{dr,S\backslash D_{j_1}/S}
\otimes\cdots\otimes f^*\mathbb B_{dr,S\backslash D_{j_{r'}}/S} 
\otimes_{f^*\mathbb B_{dr,S}}f^*\mathbb B_{dr,\phi_{D_{j_{r'+1}}}}\otimes\cdots\otimes f^*\mathbb B_{dr,\phi_{D_{j_{s'}}}} \\
\otimes_{f^*\mathbb B_{dr,S}}f^*\mathbb B_{dr,\psi_{D_{j_{s'+1}}}}
\otimes\cdots\otimes f^*\mathbb B_{dr,\psi_{D_{j_e}}}))\to\cdots) \\
\xrightarrow{Rf_*(I\otimes m_f(\mathbb B_{dr,S\backslash D_j/S})
\otimes m_f(\mathbb B_{dr\phi_{D_{j}}})\otimes m_f(\mathbb B_{dr,\psi_{D_{j}}}))} \\
(\cdots\to\xi_{D_e}\cdots\xi_{D_1}(
Rf_*(\xi_{E_d}\cdots\xi_{E_1}(K,W)\otimes_{\mathbb Q_p}f^*O_{S\backslash D_1}\otimes_{f^*O_S} \\
\mathbb B_{dr,x_{X\backslash E_{i_1}/X}}\otimes_{\mathbb B_{dr,X}}
\otimes\cdots\otimes\mathbb B_{dr,x_{X\backslash E_{i_r}/X}}  
\otimes_{\mathbb B_{dr,X}}\mathbb B_{dr,\phi_{E_{i_{r+1}}}}\otimes\cdots\otimes\mathbb B_{dr,\phi_{E_{i_s}}} \\
\otimes_{\mathbb B_{dr,X}}\mathbb B_{dr,\psi_{E_{i_{s+1}}}}\otimes\cdots\otimes\mathbb B_{dr,\phi_{E_{i_d}}} 
\otimes_{\mathbb B_{dr,X}}
\mathbb B_{dr,X\backslash f^{-1}(D_{j_1})/X}\otimes\cdots\otimes\mathbb B_{dr,X\backslash f^{-1}(D_{j_{r'}})/X} \\
\otimes_{\mathbb B_{dr,X}}\mathbb B_{dr,\phi_{f^{-1}(D_{j_{r'+1}})}}\otimes\cdots\otimes\mathbb B_{dr,\phi_{f^{-1}(D_{j_{s'}})}}
\otimes_{\mathbb B_{dr,X}}\mathbb B_{dr,\psi_{f^{-1}(D_{j_{s'+1}})}}\otimes\cdots\otimes\mathbb B_{dr,\psi_{f^{-1}(D_{j_e})}}))
\to\cdots) \\
\xrightarrow{=}(\cdots\to
Rf_*(\xi_{f^{-1}(D_e)}\cdots\xi_{f^{-1}(D_1)}\xi_{E_d}\cdots\xi_{E_1}(K,W)\otimes_{\mathbb Q_p}f^*O_{S\backslash D_1}
\otimes_{f^*O_S} \\
\mathbb B_{dr,X\backslash E_{i_1}/X}\otimes_{\mathbb B_{dr,X}}
\otimes\cdots\otimes\mathbb B_{dr,X\backslash E_{i_r}/X}  
\otimes_{\mathbb B_{dr,X}}\mathbb B_{dr,\phi_{E_{i_{r+1}}}}\otimes\cdots\otimes\mathbb B_{dr,\phi_{E_{i_s}}} \\
\otimes_{\mathbb B_{dr,X}}\mathbb B_{dr,\psi_{E_{i_{s+1}}}}\otimes\cdots\otimes\mathbb B_{dr,\phi_{E_{i_d}}} 
\otimes_{\mathbb B_{dr,X}}
\mathbb B_{dr,X\backslash f^{-1}(D_{j_1})/X}\otimes\cdots\otimes\mathbb B_{dr,x_{X\backslash f^{-1}(D_{j_{r'}})/X}} \\
\otimes_{\mathbb B_{dr,X}}\mathbb B_{dr,\phi_{f^{-1}(D_{j_{r'+1}})}}\otimes\cdots\otimes\mathbb B_{dr,\phi_{f^{-1}(D_{j_{s'}})}}
\otimes_{\mathbb B_{dr,X}}\mathbb B_{dr,\psi_{f^{-1}(D_{j_{s'+1}})}}\otimes\cdots\otimes\mathbb B_{dr,\psi_{f^{-1}(D_{j_e})}})
\to\cdots) \\
\xrightarrow{=:}Rf_*\mathbb B_{dr,X}(K,W),
\end{eqnarray*}
with $l_k(K,W):{\,}^pR^kf_*(K,W)\hookrightarrow Rf_*(K,W)$, which gives the canonical map in $D_{\mathbb B_{dr,S}}(S_K^{an,pet})$
\begin{eqnarray*}
T(f,\mathbb B_{dr})(K,W):\mathbb B_{dr,S}(Rf_*(K,W))\xrightarrow{\mathbb B_{dr,S}(l(K,W))}
\oplus_{k\in\mathbb Z}\mathbb B_{dr,S}({\,}^pR^kf_*(K,W)) \\
\xrightarrow{(T^k(f,\mathbb B_{dr})(K,W))}Rf_*\mathbb B_{dr,X}(K,W).
\end{eqnarray*}
It gives, by functoriality, for $(K,W)\in D_{\mathbb Z_pfil,c,k,gm}(S^{et})$, 
the canonical map in $D_{\mathbb B_{dr,S}fil}(S_K^{an,pet})$
\begin{eqnarray*}
T(f,\mathbb B_{dr})(K,W):\mathbb B_{dr,S}(Rf_*(K,W))\to Rf_*\mathbb B_{dr,X}(K,W).
\end{eqnarray*}
\item[(ii)]Let $f:X\to S$ a morphism with $S,X\in\QPVar(k)$. 
Consider a factorization $f:X\hookrightarrow Y\times S\xrightarrow{p}S$ with $Y\in\SmVar(k)$. 
Let $S=\cup_i S_i$ an open affine cover and 
$i_i:S_i\hookrightarrow\tilde S_i$ closed embedding with $\tilde S_i\in\SmVar(k)$.
Denote by $i'_I:X_I\hookrightarrow Y\times\tilde S_I$ the closed embeddings.
For $(K,W)\in P_{\mathbb Z_pfil,k}(X^{et})\cap D_{\mathbb Z_pfil,c,k,gm}(X^{et})$, 
we have as in (i) the following map in $D_{\mathbb B_{dr}fil}(S_K^{an,pet}/(\tilde S_{I,K})^{an,pet})$
\begin{eqnarray*}
T^k(f,\mathbb B_{dr})(K,W):\mathbb B_{dr,(\tilde S_I)}({\,}^pRf_*(K,W))
\xrightarrow{=} \\
(H^k(\cdots\to\xi_{\tilde D_{e,I}}\cdots\xi_{\tilde D_{1,I}}(
p_{\tilde S_I*}E(\xi_{\tilde E_{d,I}}\cdots\xi_{\tilde E_{1,I}}(K,W))
\otimes_{\mathbb Q_p}O_{\tilde S_I\backslash D_{1,I}})\otimes_{O_{\tilde S_I}} \\
\mathbb B_{dr,\tilde S_I\backslash\tilde D_{j_1,I}/\tilde S_I}\otimes\cdots\otimes
\mathbb B_{dr,\tilde S_I\backslash\tilde D_{j_{r'},I}/\tilde S_I} 
\otimes_{\mathbb B_{dr,\tilde S_I}}\mathbb B_{dr,\phi_{\tilde D_{j_{r'+1},I}}}\otimes\cdots\otimes \\
\mathbb B_{dr,\phi_{\tilde D_{j_{s'},I}}}\otimes_{\mathbb B_{dr,\tilde S_I}}
\mathbb B_{dr,\psi_{\tilde D_{j_{s'+1},I}}}\otimes\cdots\otimes\mathbb B_{dr,\psi_{\tilde D_{j_e,I}}}
\to\cdots),\mathbb B_{dr}(t_{IJ})) \\
\xrightarrow{((I\otimes m_f(-))\circ(I\otimes a_{-}(-))\circ T(p,\otimes)(-,-)\circ l_k(-)\otimes I))} \\
((\cdots\to p_{\tilde S_I*}E(\xi_{p^{-1}(\tilde D_{e,I})}\cdots\xi_{p^{-1}(\tilde D_{1,I})}
\xi_{\tilde E_{d,I}}\cdots\xi_{\tilde E_{1,I}}((K,W)\otimes_{\mathbb Q_p}O_{X\backslash E_1})
\otimes_{O_{Y\times\tilde S_I}} \\ 
\mathbb B_{dr,(Y\times\tilde S_I)\backslash\tilde E_{i_1}/Y\times\tilde S_I}\otimes_{\mathbb B_{dr,Y\times\tilde S_I}}
\otimes\cdots\otimes\mathbb B_{dr,(Y\times\tilde S_I)\backslash\tilde E_{i_r,I}/Y\times\tilde S_I} 
\otimes_{\mathbb B_{dr,Y\times\tilde S_I}} \\
\mathbb B_{dr,\phi_{\tilde E_{i_{r+1},I}}}\otimes\cdots\otimes\mathbb B_{dr,\phi_{\tilde E_{i_s,I}}}
\otimes_{\mathbb B_{dr,Y\times\tilde S_I}}\mathbb B_{dr,\psi_{\tilde E_{i_{s+1},I}}}
\otimes\cdots\otimes\mathbb B_{dr,\phi_{\tilde E_{i_d,I}}} \\
\otimes_{\mathbb B_{dr,Y\times\tilde S_I}}
\mathbb B_{dr,(Y\times\tilde S_I)\backslash p^{-1}(\tilde D_{j_1,I})/X}\otimes\cdots\otimes
\mathbb B_{dr,(Y\times\tilde S_I)\backslash p^{-1}(\tilde D_{j_{r'}})/X} \\
\otimes_{\mathbb B_{dr,Y\times\tilde S_I}}\mathbb B_{dr,\phi_{p^{-1}(\tilde D_{j_{r'+1},I})}}\otimes\cdots
\otimes\mathbb B_{dr,\phi_{p^{-1}(\tilde D_{j_{s'},I})}} \\
\otimes_{\mathbb B_{dr,Y\times\tilde S_I}}\mathbb B_{dr,\psi_{p^{-1}(\tilde D_{j_{s'+1},I})}}\otimes\cdots\otimes
\mathbb B_{dr,\psi_{p^{-1}(\tilde D_{j_e,I})}})
\to\cdots),\mathbb B_{dr}(t_{IJ}))
\xrightarrow{=:}Rp_*\mathbb B_{dr,(Y\times\tilde S_I)}(K,W) 
\end{eqnarray*}
where $(E_1,\ldots,E_d)\in\mathcal S(K)$ is a stratification by Cartier divisor $E_i\subset X$, $1\leq i\leq d$, such that 
\begin{equation*}
K_{|E(r)\backslash E(r+1)}:=l_r^*K\in D_{\mathbb Z_p,c}((E(r)\backslash E(r+1))^{et}) 
\end{equation*}
are local systems for all $1\leq r\leq d$, $l_r:E(r)\hookrightarrow X$ being the locally closed embeddings,
and $(D_1,\ldots,D_e)\in\mathcal S(K)$ is a stratification by Cartier divisor $D_i\subset S$, $1\leq i\leq e$, such that 
\begin{eqnarray*}
({\,}^pR^kf_*\xi_{E_d}\cdots\xi_{E_1}K)_{D(r')\backslash D(r'+1)}:=m_{r'}^*{\,}^pR^kf_*\xi_{E_d}\cdots\xi_{E_1}K
\in D_{\mathbb Z_p,c}((D(r')\backslash D(r'+1))^{et})
\end{eqnarray*}
are local systems for all $\xi:[1,\cdots,d]\to\left\{\psi^u,V_0j_*\oplus\phi^u,\psi^u\oplus\phi^u,\psi^u\right\}$, 
all $1\leq r'\leq e$, $k\in\mathbb Z$, $m_{r'}:D(r')\hookrightarrow S$ being the locally closed embeddings,
$\tilde D_{s,I}\subset\tilde S_I$ (Cartier) divisor such that $D_s\cap S_I\subset\tilde D_{s,I}\cap S$,
and $\tilde E_{s,I}\subset Y\times\tilde S_I$ (Cartier) divisor such that $E_s\cap X_I\subset\tilde E_{s,I}\cap X$.
It gives, by functoriality, for $(K,W)\in D_{\mathbb Z_pfil,c,k,gm}(S^{et})$, the canonical map in 
$D_{\mathbb B_{dr}fil}(S_K^{an,pet}/(\tilde S_{I,K}^{an,pet}))$
\begin{eqnarray*}
T(f,\mathbb B_{dr})(K,W):\mathbb B_{dr,(\tilde S_I)}(Rf_*(K,W))\to Rf_*\mathbb B_{dr,(Y\times\tilde S_I)}(K,W).
\end{eqnarray*}
\end{itemize}
\end{defi}

\begin{lem}\label{TfBdrlem}
Let $f:X\to S$ a proper morphism with $X,S\in\SmVar(k)$.
Let $E\subset X$ a (Cartier) divisor. Denote $j:U:=X\backslash E\hookrightarrow X$ the open embedding. 
Let $(K,W)\in P_{\mathbb Z_pfil,k}(X^{et})^{gm}$ such that $K_{|U}$ and $K_{|E}$ are local systems 
and such that ${\,}^pR^kf_*K=R^kf_*K$ are local systems for all $k\in\mathbb Z$. Then,
\begin{itemize}
\item[(i)] The map
\begin{eqnarray*}
Rf_*(I\otimes a_X(\mathbb B_{dr,U/X,\mathbb C_p})): \\ 
Rf_*Rj_*\pi_{k/\mathbb C_p}^*(K,W)\otimes\mathbb B_{dr,S_{\mathbb C_p}}\to 
Rf_*(V_{E,0}j_*(\pi_{k/\mathbb C_p}^*(K,W)\otimes_{\mathbb Q_p}O_{U_{\mathbb C_p}})
\otimes_{O_{X_{\mathbb C_p}}}\mathbb B_{dr,U/X,\mathbb C_p}) 
\end{eqnarray*}
is an isomorphism.
\item[(ii)] The map
\begin{eqnarray*}
Rf_*(I\otimes a_X(\mathbb B_{dr,\psi_E})):Rf_*\psi_E(K,W)\otimes\mathbb B_{dr,S}\to 
Rf_*(\psi_E((K,W)\otimes_{\mathbb Q_p}O_U)\otimes_{O_X}\mathbb B_{dr,\psi_E}) 
\end{eqnarray*}
is an isomorphism.
\item[(iii)] The map
\begin{eqnarray*}
Rf_*(I\otimes a_X(\mathbb B_{dr,\phi_E})):Rf_*\phi_E(K,W)\otimes\mathbb B_{dr,S}\to 
Rf_*(\phi_E((K,W)\otimes_{\mathbb Q_p}O_U)\otimes_{O_X}\mathbb B_{dr,\phi_E}) 
\end{eqnarray*}
is an isomorphism.
\end{itemize}
\end{lem}

\begin{proof}
By considering the spectral sequence associated to $(K,W)$, 
we may assume that $K=Rf'_*\mathbb Z_{X'}$ with $f':X'\to X$ a proper morphism with $X'\in\SmVar(k)$.

\noindent(i):Consider a desingularization of the pair $(X,E)$. 
Then the $E_2$ degenerescence of the perverse Leray spectral sequence, 
and the a normal crossing divisor case (see \cite{Chinois}).
shows that $H^i(\mathbb B_{dr,U/X,\mathbb C_p})=0$ for all $i\in\mathbb Z$, $i\neq 0$. 
Hence, (i) follows from theorem \ref{KUNetth} and theorem \ref{CPetbet}.

\noindent(ii):Follows from theorem \ref{phipsithmp} and 
on the other hand theorem \ref{KUNetth} and theorem \ref{CPetbet}.

\noindent(iii):Follows from theorem \ref{phipsithmp} and 
on the other hand theorem \ref{KUNetth} and theorem \ref{CPetbet}.
\end{proof}

\begin{thm}\label{TfBdrthm}
\begin{itemize}
\item[(i)]Let $X\in\PSmVar(k)$. Let $Z\subset X$ a closed subset.
Denote by $j:U:=X\backslash Z\hookrightarrow X$ the open complementary embedding. 
Take (Cartier) divisor $D_1,\ldots,D_r\subset X$ such that $Z=\cap_{i=1}^r D_i$. 
The map in $D(\mathbb B_{dr,\mathbb C_p},G)$
\begin{eqnarray*}
R\Gamma(X_{\bar k},I\otimes a_X(\mathbb B_{dr,D_i/X})): 
R\Gamma(X_{\bar k},Rj_*\mathbb Z_{U^{et},p})\otimes\mathbb B_{dr,\mathbb C_p}=
R\Gamma(U_{\bar k},\mathbb Z_{p,U^{et}})\otimes\mathbb B_{dr,\mathbb C_p} \\
\to R\Gamma(X_{\bar k},j_*\mathbb Z_{p,U^{et}}\otimes_{\mathbb Q_p}
\mathbb B_{dr,X\backslash D_1/X}\otimes_{\mathbb B_{dr,X}}\cdots
\otimes_{\mathbb B_{dr,X}}\mathbb B_{dr,X\backslash D_1/X}) \\
\xrightarrow{=} 
R\Gamma(X_{\bar k},F^0DR(X)(j_{*Hdg}(O_U,F_b)^{an}\otimes_{O_X}(O\mathbb B_{dr,X},F)))
\end{eqnarray*}
is an isomorphism.
\item[(ii)]Let $f:X\to S$ be a proper morphism with $X,S\in\SmVar(k)$. 
For $(K,W)\in D_{\mathbb Z_pfil,c,k}(X^{et})^{gm}$, the map in $D_{\mathbb B_{dr},G,fil}(S_{\mathbb C_p}^{an,pet})$
(where the a $G$ module structure is a continuous action of the Galois group)
\begin{eqnarray*}
T(f,B_{dr})(K,W):\mathbb B_{dr,S}(Rf_*(K,W))\to Rf_*\mathbb B_{dr,X}(K,W)
\end{eqnarray*} 
given in definition \ref{TfjBdr1} is an isomorphism.
\item[(ii)']Let $f:X\to S$ a morphism with $S,X\in\QPVar(k)$. 
Consider a factorization $f:X\hookrightarrow Y\times S\xrightarrow{p}S$ with $Y\in\SmVar(k)$. 
Let $S=\cup_i S_i$ an open affine cover and 
$i_i:S_i\hookrightarrow\tilde S_i$ closed embedding with $\tilde S_i\in\SmVar(k)$.
Denote by $i'_I:X_I\hookrightarrow Y\times\tilde S_I$ the closed embeddings.
For $(K,W)\in D_{\mathbb Z_pfil,c,k}(X^{et})^{gm}$, 
the map in $D_{\mathbb B_{dr},G,fil}(S_{\mathbb C_p}^{an,pet}/(\tilde S_{I,\mathbb C_p})^{an,pet})$
\begin{eqnarray*}
T(f,\mathbb B_{dr})(K,W):\mathbb B_{dr,(\tilde S_I)}(Rf_*(K,W))\to Rp_*\mathbb B_{dr,(Y\times\tilde S_I)}(K,W) 
\end{eqnarray*}
given in definition \ref{TfjBdr1} is an isomorphism.
\end{itemize}
\end{thm}

\begin{proof}
\noindent(i):Follows from lemma \ref{TfBdrlem}(i).

\noindent(ii):Follows from lemma \ref{TfBdrlem} and 
on the other hand theorem \ref{PFet} together with theorem \ref{CPetbet}.

\noindent(ii)':Follows from lemma \ref{TfBdrlem} and 
on the other hand theorem \ref{PFet} together with theorem \ref{CPetbet} as for (ii).
\end{proof}

\begin{rem}
Let $f:X\to S$ a proper morphism with $S,X\in\Var(k)$. 
Then for $K\in C_{\mathbb Z_p}(X^{et})$, the map in $D_{\mathbb Z_p}(S^{an,pet})$
\begin{eqnarray*}
Rf_*K\otimes\mathbb B_{dr,S}\xrightarrow{\ad(Lf^{*mod},Rf_*)(\mathbb B_{dr,S})}
Rf_*K\otimes Rf_*\mathbb B_{dr,X}\xrightarrow{T(f,f,\otimes)(K,\mathbb B_{dr,X})}Rf_*(K\otimes\mathbb B_{dr,X})
\end{eqnarray*}
is an isomorphism by theorem \ref{KUNetth} and theorem \ref{CPetbet}.
In the analytic case (\cite{Scholze}), for $f:X\to S$ a smooth proper morphism with $X,S\in\AnSm(K)$
and $L\in Loc_{\mathbb Z_p}(X^{et})$ an analytic local system, the map in $D_{\mathbb Z_p}(S^{pet})$
\begin{eqnarray*}
Rf_*L\otimes\mathbb B_{dr,S}\xrightarrow{\ad(Lf^{*mod},Rf_*)(\mathbb B_{dr,S})}
Rf_*L\otimes Rf_*\mathbb B_{dr,X}\xrightarrow{T(f,f,\otimes)(L,\mathbb B_{dr,X})}Rf_*(L\otimes\mathbb B_{dr,X})
\end{eqnarray*}
is an isomorphism.
\end{rem}

\begin{defi}\label{TfjBdr2}
\begin{itemize}
\item[(i)]Let $j:S^o\hookrightarrow S$ an open embedding with $S\in\SmVar(k)$ and $D:=S\backslash S^o$ a (Cartier) divisor. 
We will consider, using definition \ref{VfilKMmap}(vi)
for $(K,W)\in P_{\mathbb Z_pfil,k}(S^{o,et})$, the canonical isomorphism in $D_{\mathbb B_{dr,S}fil}(S_K^{an,pet})$
\begin{eqnarray*}
T(j,\mathbb B_{dr})(K,W):\mathbb B_{dr,S}(j_{*w}(K,W))\xrightarrow{:=} \\
((\cdots\to\xi_{\bar E_d}\cdots\xi_{\bar E_1}\psi^u_D((K,W)\otimes_{\mathbb Q_p}O_{S^o})\otimes_{O_S}
\mathbb B_{dr,\psi_D}\otimes_{\mathbb B_{dr,S}}\mathbb B_{dr,S\backslash\bar E_{i_1}/S}
\otimes_{\mathbb B_{dr,S}}\cdots\otimes \\
\mathbb B_{dr,S\backslash \bar E_{i_r}/S}  
\otimes_{\mathbb B_{dr,S}}\mathbb B_{dr,\phi_{\bar E_{i_{r+1}}}}\otimes\cdots\otimes\mathbb B_{dr,\phi_{\bar E_{i_s}}}
\otimes_{\mathbb B_{dr,S}}\mathbb B_{dr,\psi_{\bar E_{i_{s+1}}}}\otimes\cdots\otimes\mathbb B_{dr,\psi_{\bar E_{i_d}}} 
\to\cdots)\to \\
(\cdots\to\Cone(\xi_{\bar E_d}\cdots\xi_{\bar E_1}V_{D,0}j_{*w}((K,W)\otimes_{\mathbb Q_p}O_{S^o}) \\
\to\xi_{\bar E_d}\cdots\xi_{\bar E_1}\psi^u_D((K,W)\otimes_{\mathbb Q_p}O_{S^o}))
\oplus\xi_{\bar E_d}\cdots\xi_{\bar E_1}\phi_D^u((K,W)\otimes_{\mathbb Q_p}O_{S^o}) \\
\otimes_{O_S}(\mathbb B_{dr,S^o/S}\oplus\mathbb B_{dr,\phi_D})\otimes_{\mathbb B_{dr,S}}
\mathbb B_{dr,S\backslash\bar E_{i_1}/S} 
\otimes_{\mathbb B_{dr,S}}\cdots\otimes\mathbb B_{dr,S\backslash \bar E_{i_r}/S} \\ 
\otimes_{\mathbb B_{dr,S}}\mathbb B_{dr,\phi_{\bar E_{i_{r+1}}}}\otimes\cdots\otimes\mathbb B_{dr,\phi_{\bar E_{i_s}}}
\otimes_{\mathbb B_{dr,S}}\mathbb B_{dr,\psi_{\bar E_{i_{s+1}}}}\otimes\cdots\otimes\mathbb B_{dr,\psi_{\bar E_{i_d}}} 
\to\cdots)\to \\
(\cdots\to\xi_{\bar E_d}\cdots\xi_{\bar E_1}\psi^u_D((K,W)\otimes_{\mathbb Q_p}O_{S^o})\otimes_{O_S}
\mathbb B_{dr,\psi_D}\otimes_{\mathbb B_{dr,S}}\mathbb B_{dr,S\backslash\bar E_{i_1}/S}
\otimes_{\mathbb B_{dr,S}}\cdots\otimes \\
\mathbb B_{dr,S\backslash \bar E_{i_r}/S}  
\otimes_{\mathbb B_{dr,S}}\mathbb B_{dr,\phi_{\bar E_{i_{r+1}}}}\otimes\cdots\otimes\mathbb B_{dr,\phi_{\bar E_{i_s}}}
\otimes_{\mathbb B_{dr,S}}\mathbb B_{dr,\psi_{\bar E_{i_{s+1}}}}\otimes\cdots\otimes\mathbb B_{dr,\psi_{\bar E_{i_d}}} 
\to\cdots)) \\
\xrightarrow{(0,(I,0),0)} \\
(\cdots\to\xi_{\bar E_d}\cdots\xi_{\bar E_1}V_{D,0}j_{*w}((K,W)\otimes_{\mathbb Q_p}O_{S^o}) 
\otimes_{O_S}\mathbb B_{dr,S^o/S}\otimes_{\mathbb B_{dr,S}}\mathbb B_{dr,x_{S\backslash\bar E_{i_1}/S}} 
\otimes_{\mathbb B_{dr,S}}\cdots\otimes \\
\mathbb B_{dr,S\backslash \bar E_{i_r}/S}  
\otimes_{\mathbb B_{dr,S}}\mathbb B_{dr,\phi_{\bar E_{i_{r+1}}}}\otimes\cdots\otimes\mathbb B_{dr,\phi_{\bar E_{i_s}}}
\otimes_{\mathbb B_{dr,S}}\mathbb B_{dr,\psi_{\bar E_{i_{s+1}}}}\otimes\cdots\otimes\mathbb B_{dr,\psi_{\bar E_{i_d}}} 
\to\cdots) \\
\xrightarrow{=:}V_{D0}j_{*w}\mathbb B_{dr,S^o}(K,W)\otimes_{\mathbb B_{dr,S}}\mathbb B_{dr,S^o/S}.
\end{eqnarray*}
where $(E_1,\ldots,E_d)\in\mathcal S(K)$ is a stratification by (Cartier) divisor $E_i\subset S^o$, $1\leq i\leq d$, such that 
\begin{equation*}
K_{|E(r)\backslash E(r+1)}:=l_r^*K\in D_{\mathbb Z_p,c}((E(r)\backslash E(r+1))^{et}) 
\end{equation*}
are local systems for all $1\leq r\leq d$, $l_r:E(r)\hookrightarrow S^o$ being the locally closed embeddings.
\item[(ii)]Let $l:S^o\hookrightarrow S$ an open embedding with $S\in\Var(k)$ 
such that $D=S\backslash S^o\subset S$ is a Cartier divisor.
Let $S=\cup_i S_i$ an open affine cover and 
$i_i:S_i\hookrightarrow\tilde S_i$ closed embedding with $\tilde S_i\in\SmVar(k)$. 
Let $l_I:\tilde S^o_I\hookrightarrow\tilde S_I$ open embeddings such that $\tilde S^o_I\cap S=S^o\cap S_I$
and $\tilde D_I\subset\tilde S_I$ a Cartier divisor such that $D\cap S_I\subset\tilde D_I\cap S$.
We will consider, using definition \ref{VfilKMmap}(vi), 
for $(K,W)\in P_{\mathbb Z_pfil,k}(S^{o,et})$, the canonical isomorphism 
in $D_{\mathbb B_{dr}fil}(S_K^{an,pet}/(\tilde S_{I,K})^{an,pet})$
\begin{eqnarray*}
T(l,\mathbb B_{dr})(K,W):\mathbb B_{dr,(\tilde S_I)}(l_{*w}(K,W))\xrightarrow{:=} \\
(((\cdots\to
\xi_{\tilde E_{d,I}}\circ\xi_{\tilde E_{1,I}}\psi^u_{\tilde D_I}(i_{I*}j_I^*(K,W)\otimes_{\mathbb Q_p}O_{S^o}) \\
\otimes_{O_{\tilde S_I}}\mathbb B_{dr,\psi_{\tilde D_I}}\otimes_{\mathbb B_{dr,\tilde S_I}}
\mathbb B_{dr,\tilde S_I\backslash\tilde E_{i_1,I}/S}
\otimes_{\mathbb B_{dr,\tilde S_I}}\cdots\otimes\mathbb B_{dr,\tilde S_I\backslash\tilde E_{i_r,I}/S} \\
\otimes_{\mathbb B_{dr,\tilde S_I}}\mathbb B_{dr,\phi_{\tilde E_{i_{r+1},I}}}
\otimes\cdots\otimes\mathbb B_{dr,\phi_{\tilde E_{i_s,I}}}
\otimes_{\mathbb B_{dr,\tilde S_I}}\mathbb B_{dr,\psi_{\tilde E_{i_{s+1},I}}}
\otimes\cdots\otimes\mathbb B_{dr,\psi_{\tilde E_{i_d,I}}} 
\to\cdots)\to \\
(\cdots\to
\Cone(\xi_{\tilde E_{d,I}}\circ\xi_{\tilde E_{1,I}}V_{\tilde D_I,0}l_{I*w}(i_{I*}j_I^*(K,W)\otimes_{\mathbb Q_p}O_{S^o}) \\
\to\xi_{\tilde E_{d,I}}\circ\xi_{\tilde E_{1,I}}\psi^u_{\tilde D_I}(i_{I*}j_I^*(K,W)\otimes_{\mathbb Q_p}O_{S^o}))
\oplus\xi_{\tilde E_{d,I}}\circ\xi_{\tilde E_{1,I}}\phi^u_{\tilde D_I}(i_{I*}j_I^*(K,W)\otimes_{\mathbb Q_p}O_{S^o}) \\
\otimes_{O_{\tilde S_I}}(\mathbb B_{dr,\tilde S_I^o/\tilde S_I} 
\mathbb B_{dr,\phi_D})\otimes_{\mathbb B_{dr,\tilde S_I}}
\mathbb B_{dr,\tilde S_I\backslash\tilde E_{i_1,I}/S}\otimes_{\mathbb B_{dr,\tilde S_I}}\cdots\otimes
\mathbb B_{dr,\tilde S_I\backslash\tilde E_{i_r,I}/S} \\ 
\otimes_{\mathbb B_{dr,\tilde S_I}}\mathbb B_{dr,\phi_{\tilde E_{i_{r+1},I}}}
\otimes\cdots\otimes\mathbb B_{dr,\phi_{\tilde E_{i_s,I}}}
\otimes_{\mathbb B_{dr,\tilde S_I}}\mathbb B_{dr,\psi_{\tilde E_{i_{s+1},I}}}
\otimes\cdots\otimes\mathbb B_{dr,\psi_{\tilde E_{i_d,I}}} 
\to\cdots)\to \\
(\cdots\to\xi_{\tilde E_{d,I}}\circ\xi_{\tilde E_{1,I}}\psi^u_{\tilde D_I}(i_{I*}j_I^*(K,W)\otimes_{\mathbb Q_p}O_{S^o})
\otimes_{O_{\tilde S_I}} \\
\mathbb B_{dr,\psi_D}\otimes_{\mathbb B_{dr,\tilde S_I}}\mathbb B_{dr,\tilde S_I\backslash\tilde E_{i_1,I}/S} 
\otimes_{\mathbb B_{dr,\tilde S_I}}\cdots\otimes\mathbb B_{dr,\tilde S_I\backslash\tilde E_{i_r,I}/\tilde S_I} \\
\otimes_{\mathbb B_{dr,\tilde S_I}}\mathbb B_{dr,\phi_{\tilde E_{i_{r+1},I}}}\otimes\cdots
\otimes\mathbb B_{dr,\phi_{\tilde E_{i_s,I}}}
\otimes_{\mathbb B_{dr,\tilde S_I}}\mathbb B_{dr,\psi_{\tilde E_{i_{s+1},I}}}\otimes\cdots
\otimes\mathbb B_{dr,\psi_{\tilde E_{i_d,I}}} 
\to\cdots)),\mathbb B_{dr}(t_{IJ})) \\
\xrightarrow{(0,(I,0),0)} \\
(\cdots\to
\xi_{\tilde E_{d,I}}\circ\xi_{\tilde E_{1,I}}V_{\tilde D_I,0}l_{I*w}(i_{I*}j_I^*(K,W)\otimes_{\mathbb Q_p}O_{S^o}) \\
\otimes_{O_{\tilde S_I}}\mathbb B_{dr,\tilde S_I\backslash\tilde E_{i_1,I}/S} 
\otimes_{\mathbb B_{dr,\tilde S_I}}\cdots\otimes\mathbb B_{dr,\tilde S_I\backslash\tilde E_{i_r,I}/\tilde S_I} 
\otimes_{\mathbb B_{dr,\tilde S_I}}\mathbb B_{dr,\phi_{\tilde E_{i_{r+1},I}}}\otimes\cdots 
\otimes\mathbb B_{dr,\phi_{\tilde E_{i_s,I}}} \\
\otimes_{\mathbb B_{dr,\tilde S_I}}\mathbb B_{dr,\psi_{\tilde E_{i_{s+1},I}}}\otimes\cdots
\otimes\mathbb B_{dr,\psi_{\tilde E_{i_d,I}}}\otimes_{\mathbb B_{dr,\tilde S_I}}\mathbb B_{dr,\tilde S_I^o/\tilde S_I} 
\to\cdots),\mathbb B_{dr}(t_{IJ})) \\
\xrightarrow{=:}
V_{D0}l_{*w}\mathbb B_{dr,(\tilde S_I^o)}(K,W)\otimes_{\mathbb B_{dr,S}}(\mathbb B_{dr,\tilde S_I^o/\tilde S_I},t_{IJ}).
\end{eqnarray*}
where $(E_1,\ldots,E_d)\in\mathcal S(K)$ is a stratification by Cartier divisor $E_i\subset S^o$, $1\leq i\leq d$, such that 
\begin{equation*}
K_{|E(r)\backslash E(r+1)}:=l_r^*K\in D_{\mathbb Z_p,c}((E(r)\backslash E(r+1))^{et}) 
\end{equation*}
are local systems for all $1\leq r\leq d$, $l_r:E(r)\hookrightarrow S^o$ being the locally closed embeddings,
and $\tilde E_{s,I}\subset\tilde S_I$, $\tilde D_I\subset\tilde S_I$ are (Cartier) divisor 
such that $\bar E_s\cap S_I\subset\tilde E_{s,I}\cap S_I$ and $D\cap S_I\subset\tilde D_I\cap S_I$.
It gives for $(K,W)\in D_{\mathbb Z_pfil,c,k}(S^{o,et})$, the canonical isomorphism 
in $D_{\mathbb B_{dr}fil}(S_K^{an,pet}/(\tilde S_{I,K})^{an,pet})$
\begin{eqnarray*}
T(l,\mathbb B_{dr})(K,W):\mathbb B_{dr,(\tilde S_I)}(l_{*w}(K,W))\xrightarrow{\sim}
V_{D0}l_{*w}\mathbb B_{dr,(\tilde S_I^o)}(K,W)\otimes_{\mathbb B_{dr,S}}(\mathbb B_{dr,\tilde S_I^o/\tilde S_I},t_{IJ}).
\end{eqnarray*}
\end{itemize}
\end{defi}

Let $S\in\SmVar(k)$. 
Let $j:S^o\hookrightarrow S$ an open embedding such that $D=S\backslash S^o\subset S$ is a Cartier divisor.
Denote by $\Delta_S:S\hookrightarrow S\times S$ the diagonal closed embedding and
$p_1:S\times S\to S$ and $p_2:S\times S\to S$ the projections.
\begin{itemize}
\item We have the isomorphism in $C_{\mathbb B_{dr}}(S_K^{an,pet})$
\begin{eqnarray*}
m(\mathbb B_{dr,S^o/S,K}):\mathbb B_{dr,S^o/S,K}\otimes_{\mathbb B_{dr,S_K}}\mathbb B_{dr,S^o/S,K} 
\xrightarrow{:=} \\
F^0DR(S)(j_{*Hdg}(O_{S^o},F_b)\otimes_{O_S}(O\mathbb B_{dr,S_K},F))\otimes_{\mathbb B_{dr,S_K}} 
F^0DR(S)(j_{*Hdg}(O_{S^o},F_b)\otimes_{O_S}(O\mathbb B_{dr,S_K},F)) \\
\xrightarrow{F^0w_S}
\Delta_S^{*mod}F^0DR(S\times S)((p_1^{*mod}j_{*Hdg}(O_{S^o},F_b)\otimes_{O_{S\times S}} \\
p_2^{*mod}j_{*Hdg}(O_{S^o},F_b)\otimes_{O_{S\times S}}(O\mathbb B_{dr,(S\times S)_K},F))) \\
\xrightarrow{DR(S\times S)(\ad(\Delta_{S,Hdg}^{*mod},\Delta_{S*mod})(-))} \\
\Delta_S^{*mod}F^0DR(S\times S)(\Delta_{S*mod}\Delta^{*mod}_{S,Hdg}
(p_1^{*mod}j_{*Hdg}(O_{S^o},F_b)\otimes_{O_{S\times S}}p_2^{*mod}j_{*Hdg}(O_{S^o},F_b)) \\
\otimes_{O_{S\times S}}(O\mathbb B_{dr,(S\times S)_K},F)) 
\xrightarrow{T^{B_{dr}}(\Delta_S,DR)(-)} \\
\Delta_S^{*mod}\Delta_{S*}F^0DR(S)(\Delta^{*mod}_{S,Hdg}(p_1^{*mod}j_{*Hdg}(O_{S^o},F_b)
\otimes_{O_{S\times S}}p_2^{*mod}j_{*Hdg}(O_{S^o},F_b))\otimes_{O_S}(O\mathbb B_{dr,S_K},F)) \\
\xrightarrow{=}
F^0DR(S)(j_{*Hdg}(O_{S^o},F_b)\otimes^{Hdg}_{O_S}j_{*Hdg}(O_{S^o},F_b)\otimes_{O_S}(O\mathbb B_{dr,S_K},F)) \\
\xrightarrow{F^0DR(S)(m)}
F^0DR(S)(j_{*Hdg}(O_{S^o},F_b)\otimes_{O_S}(O\mathbb B_{dr,S_K},F)):=\mathbb B_{dr,S^o/S,K}
\end{eqnarray*}
where 
\begin{equation*}
m:(j_*O_{S^o},V_D)\otimes_{O_S}(j_*O_{S^o},V_D)\xrightarrow{\sim}(j_*O_{S^o},V_D), \; m(b_1\otimes b_2)=b_1b_2 
\end{equation*}
is the multiplication map whose inverse is 
\begin{equation*}
n:(j_*O_{S^o},V_D)\xrightarrow{\sim}(j_*O_{S^o},V_D)\otimes_{O_S}(j_*O_{S^o},V_D), \; n(b)=b\otimes 1. 
\end{equation*}
\item We have the isomorphism
\begin{eqnarray*}
m(\mathbb B_{dr,\psi_D,K}):\mathbb B_{dr,\psi_D,K}\otimes_{\mathbb B_{dr,S_K}}\mathbb B_{dr,\psi_D,K} 
\xrightarrow{:=} \\
F^0DR(S)(\psi_D(O_{S^o},F_b)\otimes_{O_S}(O\mathbb B_{dr,S_K},F))\otimes_{\mathbb B_{dr,S_K}} 
F^0DR(S)(\psi_D(O_{S^o},F_b)\otimes_{O_S}(O\mathbb B_{dr,S_K},F)) \\
\xrightarrow{F^0w_S}
\Delta_S^{*mod}F^0DR(S\times S)(p_1^{*mod}\psi_D(O_{S^o},F_b)\otimes_{O_{S\times S}}p_2^{*mod}\psi_D(O_{S^o},F_b)
\otimes_{O_{S\times S}}(O\mathbb B_{dr,(S\times S)_K},F)) \\
\xrightarrow{DR(S\times S)(\ad(\Delta_{S,Hdg}^{*mod},\Delta_{S*mod})(-))} \\
\Delta_S^{*mod}F^0DR(S\times S)(\Delta_{S*mod}\Delta_{S,Hdg}^{*mod}(p_1^{*mod}\psi_D(O_{S^o},F_b)
\otimes_{O_{S\times S}}p_2^{*mod}\psi_D(O_{S^o},F_b)) \\
\otimes_{O_{S\times S}}(O\mathbb B_{dr,(S\times S)_K},F)) 
\xrightarrow{T^{B_{dr}}(\Delta_S,DR)(-)} \\
\Delta_S^{*mod}\Delta_{S*}F^0DR(S)(\Delta_{S,Hdg}^{*mod}(p_1^{*mod}\psi_D(O_{S^o},F_b)
\otimes_{O_{S\times S}}p_2^{*mod}\psi_D(O_{S^o},F_b))\otimes_{O_S}(O\mathbb B_{dr,S_K},F)) \\
\xrightarrow{=}
F^0DR(S)(\psi_D(O_{S^o},F_b)\otimes^{Hdg}_{O_S}\psi_D(O_{S^o},F_b)\otimes_{O_S}(O\mathbb B_{dr,S_K},F)) \\
\xrightarrow{F^0DR(S)(m)}
F^0DR(S)(\psi_D(O_{S^o},F_b)\otimes_{O_S}(O\mathbb B_{dr,S_K},F)):=\mathbb B_{dr,\psi_D,K}
\end{eqnarray*}
where 
\begin{equation*}
m:\psi_D(O_{S^o})\otimes_{O_S}\psi_D(O_{S^o})\xrightarrow{\sim}\psi_DO_{S^o}, \; m(b_1\otimes b_2)=b_1b_2 
\end{equation*}
is the multiplication map whose inverse is 
\begin{equation*}
n:\psi_D(O_{S^o})\xrightarrow{\sim}\psi_D(O_{S^o})\otimes_{O_S}\psi_D(O_{S^o}), \; n(b)=b\otimes 1. 
\end{equation*}
\item We have similarly to $\mathbb B_{dr,\psi_D,K}$ the isomorphism
\begin{eqnarray*}
m(\mathbb B_{dr,\phi_D,K}):\mathbb B_{dr,\phi_D,K}\otimes_{\mathbb B_{dr,S_K}}\mathbb B_{dr,\phi_D,K}
\xrightarrow{\sim}\mathbb B_{dr,\phi_D,K}
\end{eqnarray*}
with 
\begin{equation*}
m:\phi_D(O_{S^o})\otimes_{O_S}\phi_D(O_{S^o})\xrightarrow{\sim}\phi_DO_{S^o}, \; m(b_1\otimes b_2)=b_1b_2 
\end{equation*}
is the multiplication map whose inverse is 
\begin{equation*}
n:\phi_D(O_{S^o})\xrightarrow{\sim}\phi_D(O_{S^o})\otimes_{O_S}\phi_D(O_{S^o}), \; n(b)=b\otimes 1. 
\end{equation*}
\item We have similarly the isomorphism
\begin{eqnarray*}
m(\mathbb B_{dr,x_{S^o/S},K}):\mathbb B_{dr,x_{S^o/S},K}\otimes_{\mathbb B_{dr,S_K}}\mathbb B_{dr,x_{S^o/S},K}
\xrightarrow{\sim}\mathbb B_{dr,x_{S^o/S},K}
\end{eqnarray*}
with
\begin{eqnarray*}
m:=(m,0,m):\Cone(j_*O_{S^o}\to\psi_D(O_{S^o}))\otimes_{O_S}\Cone(j_*O_{S^o}\to\psi_D(O_{S^o}))
\xrightarrow{\sim}\Cone(j_*O_{S^o}\to\psi_D(O_{S^o})). 
\end{eqnarray*}
\end{itemize}

\begin{defi}\label{TOBdr}
\begin{itemize}
\item[(i)] Let $S\in\SmVar(k)$. For $(K_1,W),(K_2,W)\in P_{\mathbb Zpfil,k}(S^{et})$ filtered perverse sheaves, 
we have the isomorphism in $D_{\mathbb B_{dr}fil}(S_K^{an,pet})$
\begin{eqnarray*}
T(\otimes,\mathbb B_{dr})((K_1,W),(K_2,W)):\mathbb B_{dr,S}(K_1,W)\otimes_{B_{dr,S}}\mathbb B_{dr,S}(K_2,W) \\
\xrightarrow{:=}
(\cdots\to\xi_{D_d}\cdots\xi_{D_1}((K_1,W)\otimes_{\mathbb Q_p}O_{S^o})\otimes_{O_S} \\
\mathbb B_{dr,S\backslash D_{i_1}/S}\otimes_{\mathbb B_{dr,S}}\cdots
\otimes_{\mathbb B_{dr,S}}\mathbb B_{dr,S\backslash D_{i_r}/S}\otimes_{\mathbb B_{dr,S_K}} 
\mathbb B_{dr,\phi_{D_{i_{r+1}}}}\otimes_{\mathbb B_{dr,S}}\cdots \\
\otimes_{\mathbb B_{dr,S}}\mathbb B_{dr,\phi_{D_{i_s}}}\otimes_{\mathbb B_{dr,S}} 
\mathbb B_{dr,\psi_{D_{i_{s+1}}}}\otimes_{\mathbb B_{dr,S}}\cdots
\otimes_{\mathbb B_{dr,S}}\mathbb B_{dr,\psi_{D_{i_d}}}\to\cdots)
\otimes_{B_{dr,S}} \\
(\cdots\to 
\xi_{D_d}\cdots\xi_{D_1}((K_1,W)\otimes_{\mathbb Q_p}O_{S^o})\otimes_{O_S} \\
\mathbb B_{dr,S\backslash D_{i_1}/S}\otimes_{\mathbb B_{dr,S}}\cdots
\otimes_{\mathbb B_{dr,S_K}}\mathbb B_{dr,S\backslash D_{i_r}/S}\otimes_{\mathbb B_{dr,S}} 
\mathbb B_{dr,\phi_{D_{i_{r+1}}}}\otimes_{\mathbb B_{dr,S}}\cdots \\
\otimes_{\mathbb B_{dr,S}}\mathbb B_{dr,\phi_{D_{i_s}}}\otimes_{\mathbb B_{dr,S}} 
\mathbb B_{dr,\psi_{D_{i_{s+1}}}}\otimes_{\mathbb B_{dr,S}}\cdots
\otimes_{\mathbb B_{dr,S}}\mathbb B_{dr,\psi_{D_{i_d}}}\to\cdots) \\
\xrightarrow{=}
\xi_{D_d}\cdots\xi_{D_1}((K_1,W)\otimes_{\mathbb Q_p}(K_2,W)\otimes_{\mathbb Q_p}O_{S^o})\otimes_{O_S} \\
\mathbb B_{dr,x_{S\backslash D_{i_1}/S}}\otimes_{\mathbb B_{dr,S}}\cdots
\otimes_{\mathbb B_{dr,S}}\mathbb B_{dr,S\backslash D_{i_r}/S}\otimes_{\mathbb B_{dr,S_K}} \\
\mathbb B_{dr,\phi_{D_{i_{r+1}}}}\otimes_{\mathbb B_{dr,S_K}}\cdots 
\otimes_{\mathbb B_{dr,S_K}}\mathbb B_{dr,\phi_{D_{i_s}},K}\otimes_{\mathbb B_{dr,S}} 
\mathbb B_{dr,\psi_{D_{i_{s+1}}}}\otimes_{\mathbb B_{dr,S_K}}\cdots \\
\otimes_{\mathbb B_{dr,S}}\mathbb B_{dr,\psi_{D_{i_d}},K}\otimes_{B_{dr,S}} 
\mathbb B_{dr,S\backslash D_{i_1}/S}\otimes_{\mathbb B_{dr,S}}\cdots
\otimes_{\mathbb B_{dr,S}}\mathbb B_{dr,S\backslash D_{i_r}/S}\otimes_{\mathbb B_{dr,S}} \\
\mathbb B_{dr,\phi_{D_{i_{r+1}}}}\otimes_{\mathbb B_{dr,S}}\cdots 
\otimes_{\mathbb B_{dr,S}}\mathbb B_{dr,\phi_{D_{i_s}}}\otimes_{\mathbb B_{dr,S}} 
\mathbb B_{dr,\psi_{D_{i_{s+1}}}}\otimes_{\mathbb B_{dr,S}}\cdots
\otimes_{\mathbb B_{dr,S}}\mathbb B_{dr,\psi_{D_{i_d}}}\to\cdots) \\
\xrightarrow{(m(\mathbb B_{dr,S\backslash D_{i_1}/S}),\cdots,
m(\mathbb B_{dr,\phi_{D_{i_{r+1}}}}),\cdots,m(\mathbb B_{dr,\psi_{D_{i_d}},K}))} \\
(\cdots\to\xi_{D_d}\cdots\xi_{D_1}((K_1,W)\otimes(K_2,W)\otimes_{\mathbb Q_p}O_{S^o})\otimes_{O_S} \\
\mathbb B_{dr,S\backslash D_{i_1}/S}\otimes_{\mathbb B_{dr,S}}\cdots\otimes_{\mathbb B_{dr,S}} 
\mathbb B_{dr,S\backslash D_{i_r}/S}\otimes_{\mathbb B_{dr,S}}   
\mathbb B_{dr,\phi_{D_{i_{r+1}}},K}\otimes_{\mathbb B_{dr,S}}\cdots 
\otimes_{\mathbb B_{dr,S}}\mathbb B_{dr,\phi_{D_{i_s}}}) \\
\xrightarrow{=:}
\mathbb B_{dr,S}((K_1,W)\otimes^{L,w} (K_2,W))
\end{eqnarray*}
with $(D_1,\ldots,D_d)\in\mathcal S(K)$ a stratification by (Cartier) divisor $D_i\subset S$, $1\leq i\leq d$ such that 
\begin{equation*}
K_{1|D(r)\backslash D(r+1)},K_{2|D(r)\backslash D(r+1)}\in D_{\mathbb Z_p,c}(D(r)\backslash D(r+1)^{et}) 
\end{equation*}
are local systems for all $1\leq r\leq d$.
\item[(ii)]Let $S\in\Var(k)$. Let $S=\cup_{i\in I}S_i$ an open cover such that there
exists closed embeddings $i_i:S_i\hookrightarrow\tilde S_i$ with $\tilde S_I\in\SmVar(k)$.
For $(K_1,W),(K_2,W)\in P_{\mathbb Z_pfil,k}(S^{et})$ filtered perverse sheaves, 
we have the isomorphism in $D_{\mathbb B_{dr}fil}(S_K^{an,pet}/(\tilde S_{I,K}^{an,pet}))$ given as in (i) 
\begin{eqnarray*}
T(\otimes,\mathbb B_{dr})((K_1,W),(K_2,W)):
\mathbb B_{dr,(\tilde S_I)}(K_1,W)\otimes_{\mathbb B_{dr,S}}\mathbb B_{dr,(\tilde S_I)}(K_2,W) 
\xrightarrow{:=} \\
((\cdots\to\xi_{\tilde D_{d,I}}\cdots\xi_{\tilde D_{1,I}}(i_{I*}j_I^*(K_1,W)\otimes_{\mathbb Q_p}O_{\tilde S_I^o}) 
\otimes_{O_{\tilde S_I}} \\
\mathbb B_{dr,\tilde S_I\backslash\tilde D_{i_1,I}/S}\otimes_{\mathbb B_{dr,\tilde S_{I}}}\cdots
\otimes_{\mathbb B_{dr,\tilde S_{I}}}\mathbb B_{dr,\tilde S_I\backslash D_{i_r,I}/\tilde S_I} 
\otimes_{\mathbb B_{dr,\tilde S_{I}}} 
\mathbb B_{dr,\phi_{\tilde D_{i_{r+1},I}}}\otimes_{\mathbb B_{dr,\tilde S_{I}}} \\
\cdots\otimes_{\mathbb B_{dr,\tilde S_{I}}} 
\mathbb B_{dr,\phi_{\tilde D_{i_s,I}}}\otimes_{\mathbb B_{dr,\tilde S_{I}}} 
\mathbb B_{dr,\psi_{\tilde D_{i_{s+1},I}}}\otimes_{\mathbb B_{dr,\tilde S_{I}}}\cdots
\otimes_{\mathbb B_{dr,\tilde S_{I}}}\mathbb B_{dr,\psi_{\tilde D_{i_d,I}}}\to\cdots),\mathbb B_{dr}(t_{IJ}))
\otimes_{\mathbb B_{dr,S}} \\
((\cdots\to\xi_{\tilde D_{d,I}}\cdots\xi_{\tilde D_{1,I}}(i_{I*}j_I^*(K_2,W)\otimes_{\mathbb Q_p}O_{\tilde S_I^o})
\otimes_{O_{\tilde S_I}} \\
\mathbb B_{dr,x_{\tilde S_I\backslash\tilde D_{i_1,I}/S}}\otimes_{\mathbb B_{dr,\tilde S_{I}}}\cdots
\otimes_{\mathbb B_{dr,\tilde S_{I}}}\mathbb B_{dr,x_{\tilde S_I\backslash D_{i_r,I}/\tilde S_I}} 
\otimes_{\mathbb B_{dr,\tilde S_{I}}} 
\mathbb B_{dr,\phi_{\tilde D_{i_{r+1},I}}}\otimes_{\mathbb B_{dr,\tilde S_{I}}} \\
\cdots\otimes_{\mathbb B_{dr,\tilde S_{I}}}  
\mathbb B_{dr,\phi_{\tilde D_{i_s,I}}}\otimes_{\mathbb B_{dr,\tilde S_{I}}} 
\mathbb B_{dr,\psi_{\tilde D_{i_{s+1},I}}}\otimes_{\mathbb B_{dr,\tilde S_{I}}}\cdots
\otimes_{\mathbb B_{dr,\tilde S_{I}}}\mathbb B_{dr,\psi_{\tilde D_{i_d,I}}}\to\cdots),\mathbb B_{dr}(t_{IJ})) 
\xrightarrow{=} \\
((\cdots\to\xi_{\tilde D_{d,I}}\cdots\xi_{\tilde D_{1,I}}(i_{I*}j_I^*(K_1,W)\otimes(K_2,W)\otimes_{\mathbb Q_p}O_{\tilde S_I^o}) 
\otimes_{O_{\tilde S_I}} \\
\mathbb B^2_{dr,\tilde S_I\backslash\tilde D_{i_1,I}/S}\otimes_{\mathbb B_{dr,\tilde S_{I}}}\cdots
\otimes_{\mathbb B_{dr,\tilde S_{I}}}\mathbb B^2_{dr,\tilde S_I\backslash D_{i_r,I}/\tilde S_I}
\otimes_{\mathbb B_{dr,\tilde S_{I}}} 
\mathbb B^2_{dr,\phi_{\tilde D_{i_{r+1},I}}}\otimes_{\mathbb B_{dr,\tilde S_{I}}} \\
\cdots\otimes_{\mathbb B_{dr,\tilde S_{I}}} 
\mathbb B^2_{dr,\phi_{D_{i_s,I}}}\otimes_{\mathbb B_{dr,\tilde S_{I}}} 
\mathbb B^2_{dr,\psi_{\tilde D_{i_{s+1},I}}}\otimes_{\mathbb B_{dr,\tilde S_{I}}}\cdots 
\otimes_{\mathbb B_{dr,\tilde S_{I}}}\mathbb B^2_{dr,\psi_{\tilde D_{i_d,I}}}
\to\cdots),\mathbb B_{dr}(t_{IJ})) \\
\xrightarrow{(m(\mathbb B_{dr,\tilde S_I\backslash D_{i_1,I}/S}),\cdots,
m(\mathbb B_{dr,\phi_{\tilde D_{i_{r+1},I}}}),\cdots,m(\mathbb B_{dr,\psi_{\tilde D_{i_d,I}}}))} \\
((\cdots\to\xi_{\tilde D_{d,I}}\cdots\xi_{\tilde D_{1,I}}(i_{I*}j_I^*(K_1,W)\otimes(K_2,W)\otimes_{\mathbb Q_p}O_{\tilde S_I^o}) 
\otimes_{O_{\tilde S_I}} \\ 
\mathbb B_{dr,\tilde S_I\backslash\tilde D_{i_1,I}/S}\otimes_{\mathbb B_{dr,\tilde S_{I}}}\cdots
\otimes_{\mathbb B_{dr,\tilde S_{I}}}\mathbb B_{dr,\tilde S_I\backslash D_{i_r,I}/\tilde S_I} 
\otimes_{\mathbb B_{dr,\tilde S_{I}}}\mathbb B_{dr,\phi_{\tilde D_{i_{r+1},I}}} 
\otimes_{\mathbb B_{dr,\tilde S_{I}}}\cdots\otimes_{\mathbb B_{dr,\tilde S_{I}}} \\
\mathbb B_{dr,\phi_{\tilde D_{i_s,I}}}\otimes_{\mathbb B_{dr,\tilde S_{I}}} 
\mathbb B_{dr,\psi_{\tilde D_{i_{s+1},I}}}\otimes_{\mathbb B_{dr,\tilde S_{I}}}\cdots
\otimes_{\mathbb B_{dr,\tilde S_{I}}}\mathbb B_{dr,\psi_{\tilde D_{i_d,I}}}\to\cdots),\mathbb B_{dr}(t_{IJ})) \\
\xrightarrow{=:}
\mathbb B_{dr,(\tilde S_I)}((K_1,W)\otimes^{L,w}(K_2,W))
\end{eqnarray*}
with $(D_1,\ldots,D_d)\in\mathcal S(K)$ stratifications by Cartier divisor $D_i\subset S$, $1\leq i\leq d$ such that 
\begin{equation*}
K_{1|D(r)\backslash D(r+1)},K_{2|D(r)\backslash D(r+1)}\in D_{\mathbb Z_p,c}(D(r)\backslash D(r+1)^{et}) 
\end{equation*}
are local systems for all $1\leq r\leq d$, 
and $\tilde D_{s,I}\subset\tilde S_I$ (Cartier) divisor such that $D_s\cap S_I\subset\tilde D_{s,I}\cap S$.
This gives, for $(K_1,W),(K_2,W)\in D_{\mathbb Z_pfil,c,k}(S^{et})$, 
isomorphism in $D_{\mathbb B_{dr}fil}(S_K^{an,pet}/(\tilde S_{I,K}^{an,pet})$,  
\begin{eqnarray*}
T(\otimes,\mathbb B_{dr})((K_1,W),(K_2,W)):
\mathbb B_{dr,(\tilde S_I)}(K_1,W)\otimes_{\mathbb B_{dr,S}}\mathbb B_{dr,(\tilde S_I)}(K_2,W) \\
\xrightarrow{\sim}\mathbb B_{dr,(\tilde S_I)}((K_1,W)\otimes^{L,w}(K_2,W)).
\end{eqnarray*}
\end{itemize}
\end{defi}

\begin{defi}\label{TDBdr}
\begin{itemize}
\item[(i)] Let $S\in\SmVar(k)$. For $(K,W)\in D_{\mathbb Z_pfil,c,k}(S^{et})$, we have the isomorphism
\begin{eqnarray*}
T(D,\mathbb B_{dr})(K,W):\mathbb B_{dr,S}(\mathbb D^v_S(K,W))
\xrightarrow{B((K,W),\mathbb D_S^v(K,W))}\mathbb D_S\mathbb B_{dr,S}(K,W)
\end{eqnarray*}
given by the pairing
\begin{eqnarray*}
B((K,W),\mathbb D_S^v(K,W)):\mathbb B_{dr,S}(\mathbb D^v_SK)\otimes_{\mathbb B_{dr,S_K}}\mathbb B_{dr,S}(K,W) 
\xrightarrow{T(\otimes,\mathbb B_{dr})((K,W),\mathbb D_S^v(K,W))} \\
\mathbb B_{dr,S}(\mathbb D^v_S(K,W)\otimes(K,W)) 
\xrightarrow{\mathbb B_{dr,S}(ev_K)}\mathbb B_{dr,S}(\mathbb Z_{p,S^{et}})\xrightarrow{\alpha(S_K)^{-1}}\mathbb B_{dr,S_K}
\end{eqnarray*}
using definition \ref{TOBdr}.
\item[(ii)]Let $S\in\Var(k)$. Let $S=\cup_{i\in I}S_i$ an open cover such that there
exists closed embeddings $i_i:S_i\hookrightarrow\tilde S_i$ with $\tilde S_I\in\SmVar(k)$.
For $(K,W)\in D_{\mathbb Z_pfil,c,k}(S^{et})$, we have the isomorphism 
in $D_{\mathbb B_{dr}fil}(S_K^{an,pet}/(\tilde S_{I,K}^{an,pet})$
\begin{eqnarray*}
T(D,\mathbb B_{dr})(K,W):\mathbb B_{dr,(\tilde S_I)}(\mathbb D^v_S(K,W))
\xrightarrow{B((K,W),\mathbb D_S^v(K,W))}\mathbb D_S\mathbb B_{dr,(\tilde S_I)}(K,W)
\end{eqnarray*}
given by the pairing
\begin{eqnarray*}
B((K,W),\mathbb D_S^v(K,W)):
\mathbb B_{dr,(\tilde S_I)}(\mathbb D^v_S(K,W))\otimes_{\mathbb B_{dr,S_K}}\mathbb B_{dr,(\tilde S_I)}(K,W) \\
\xrightarrow{T(\otimes,\mathbb B_{dr})((K,W),\mathbb D_S^v(K,W))}
\mathbb B_{dr,(\tilde S_I)}(\mathbb D^v_S(K,W)\otimes(K,W)) \\
\xrightarrow{\mathbb B_{dr,(\tilde S_I)}(ev_K)}\mathbb B_{dr,(\tilde S_I)}(\mathbb Z_{p,S^{et}})
\xrightarrow{(\alpha(\tilde S_{I,K}))^{-1}}(\mathbb B_{dr,\tilde S_{I,K}},t_{IJ})
\end{eqnarray*}
using definition \ref{TOBdr}.
\end{itemize}
\end{defi}

\begin{defi}\label{TfjBdr22}
\begin{itemize}
\item[(i)]Let $j:S^o\hookrightarrow S$ an open embedding with $S\in\SmVar(k)$ and $D:=S\backslash S^o$ a (Cartier) divisor. 
We will consider, for $(K,W)\in D_{\mathbb Z_pfil,c,k}(S^{o,et})$, 
the canonical isomorphism in $D_{\mathbb B_{dr}fil}(S_K^{an,pet})$
\begin{eqnarray*}
T_!(j,\mathbb B_{dr})(K,W):\mathbb B_{dr,S}(j_{!w}(K,W)):=\mathbb B_{dr,S}(\mathbb D_S^vj_{*w}\mathbb D^v_S(K,W)) \\ 
\xrightarrow{T(D,\mathbb B_{dr})(j_{*w}\mathbb D^v_S(K,W))}\mathbb D_S\mathbb B_{dr,S}(j_{*w}\mathbb D^v_{S^o}(K,W)) \\
\xrightarrow{\mathbb D_ST(j,\mathbb B_{dr})(\mathbb D_{S^o}^v(K,W))}
\mathbb D_S(V_{D0}j_{*w}\mathbb B_{dr,S^o}(\mathbb D_{S^o}^v(K,W))\otimes\mathbb B_{dr,S^o/S}) \\
\xrightarrow{T(D,\mathbb B_{dr})(K,W)}
\mathbb D_S(V_{D0}j_{*w}\mathbb D_{S^o}\mathbb B_{dr,S^o}(K,W)\otimes\mathbb B_{dr,S^o/S}) \\
\xrightarrow{=}
V_{D0}j_{!w}\mathbb B_{dr,S^o}(K,W)\otimes_{\mathbb B_{dr,S}}\mathbb D_S\mathbb B_{dr,S^o/S},
\end{eqnarray*}
using definition \ref{TfjBdr2} and definition \ref{TDBdr}.
\item[(ii)]Let $l:S^o\hookrightarrow S$ an open embedding with $S\in\Var(k)$ such that $D=S\backslash S^o$ is a Cartier divisor.
Let $S=\cup_i S_i$ an open affine cover and 
$i_i:S_i\hookrightarrow\tilde S_i$ closed embedding with $\tilde S_i\in\SmVar(k)$. 
Let $l_I:\tilde S^o_I\hookrightarrow\tilde S_I$ open embeddings such that $\tilde S^o_I\cap S=S^o\cap S_I$
and $\tilde D_I\subset\tilde S_I$ a Cartier divisor such that $D\cap S_I\subset\tilde D_I\cap S$.
We will consider, for $(K,W)\in D_{\mathbb Z_pfil,c,k}(S^{o,et})$, the canonical isomorphism 
in $D_{\mathbb B_{dr}fil}(S_K^{an,pet}/(\tilde S_{I,K})^{an,pet})$
\begin{eqnarray*}
T_!(l,\mathbb B_{dr})(K,W):\mathbb B_{dr,(\tilde S_I)}(l_{!w}(K,W)):=
\mathbb B_{dr,(\tilde S_I)}(\mathbb D^v_Sl_{*w}\mathbb D^v_S(K,W)) \\
\xrightarrow{T(D,\mathbb B_{dr})(l_{*w}\mathbb D^v_S(K,W))}
\mathbb D_S\mathbb B_{dr,(\tilde S_I)}(l_{*w}\mathbb D_{S^o}^v(K,W)) \\
\xrightarrow{\mathbb D_ST(l,\mathbb B_{dr})(\mathbb D_{S^o}^v(K,W))}
\mathbb D_S(V_{D0}l_{*w}\mathbb B_{dr,(\tilde S_I^o)}(\mathbb D_{S^o}^v(K,W))
\otimes_{\mathbb B_{dr,S}}(\mathbb B_{dr,\tilde S_I^o/\tilde S_I},t_{IJ})) \\
\xrightarrow{T(D,\mathbb B_{dr})(K,W)}
\mathbb D_S(V_{D0}l_{*w}\mathbb D_{S^o}\mathbb B_{dr,(\tilde S_I^o)}(K,W)
\otimes_{\mathbb B_{dr,S}}(\mathbb B_{dr,\tilde S_I^o/\tilde S_I},t_{IJ})) \\
\xrightarrow{=}
V_{D0}l_{!w}\mathbb B_{dr,(\tilde S_I^o}(K,W)\otimes_{\mathbb B_{dr,S}}\mathbb D_S(\mathbb B_{dr,\tilde S_I^o/\tilde S_I})
\end{eqnarray*}
using definition \ref{TfjBdr2} and definition \ref{TDBdr}.
\end{itemize}
\end{defi}

As a consequence of this formalism we have :

\begin{thm}\label{Bdrthm}
\begin{itemize}
\item[(i)]Let $S\in\SmVar(k)$ irreducible. Let $\bar S\in\PSmVar(k)$ a compactification of $S$,
with $D:=\bar S\backslash S\subset S$ a (Cartier) divisor. Denote by $j:S^o\hookrightarrow S$ the open embedding.
We have the canonical isomorphisms, given by, using definition \ref{TfjBdr1}, definition \ref{TfjBdr2},
definition \ref{TOBdr} and definition \ref{TDBdr}, for $K,K'\in D_{\mathbb Z_p,c,k}(S_{\bar k}^{et})$, 
\begin{eqnarray*}
\mathbb B_{dr,S}(K_1,K_2):R\Hom(K_1,K_2)\otimes\mathbb B_{dr,\bar k}
\xrightarrow{=}\mathbb B_{dr,\bar k}(Ra_{\bar S*}Rj_*\mathcal Hom(K_1,K_2)) \\
\xrightarrow{T(a_{\bar S},\mathbb B_{dr})(-)}Ra_{\bar S*}\mathbb B_{dr,\bar S}(Rj_*\mathcal Hom(K_1,K_2)) \\
\xrightarrow{\mathbb B_{dr,\bar S}(T(j,\mathbb B_{dr})(-))}
Ra_{\bar S*}(V_{D,0}Rj_*\mathbb B_{dr,S}(\mathcal Hom(K_1,K_2))\otimes_{\mathbb B_{dr,S}}\mathbb B_{dr,S/\bar S}) \\
\xrightarrow{\mathbb B_{dr,S}(m(K_1,K_2)^{-1})}
Ra_{\bar S*}(V_{D,0}Rj_*\mathbb B_{dr,S}(\mathbb D_S^vK_1\otimes K_2)\otimes_{\mathbb B_{dr,S}}\mathbb B_{dr,S/\bar S}) \\
\xrightarrow{(T(D,\mathbb B_{dr})(K_1)\otimes I)\circ T(\otimes,\mathbb B_{dr})(\mathbb D_SK_1,K_2)}
Ra_{\bar S*}(V_{D,0}Rj_*(\mathbb D_S\mathbb B_{dr,S}(K_1)\otimes_{\mathbb B_{dr,S}}\mathbb B_{dr,S}(K_2))
\otimes_{\mathbb B_{dr,S}}\mathbb B_{dr,S/\bar S}) \\
\xrightarrow{m(\mathbb B_{dr,S}(K_1),\mathbb B_{dr,S}(K_2))}
Ra_{\bar S*}(V_{D,0}Rj_*\mathcal Hom(\mathbb B_{dr,S}(K_1),\mathbb B_{dr,S}(K_2))
\otimes_{\mathbb B_{dr,S}}\mathbb B_{dr,S/\bar S}) \\
\xrightarrow{=}R\Hom(\mathbb B_{dr,S}(K_1),\mathbb B_{dr,S}(K_2)).
\end{eqnarray*}
\item[(ii)]Let $S\in\Var(k)$. Let $\bar S\in\PVar(k)$ a compactification of $S$,
with $D:=\bar S\backslash S\subset S$ a Cartier divisor. Denote by $j:S^o\hookrightarrow S$ the open embedding.
Let $S=\cup_{i\in I}S_i$ an open cover such that there
exists closed embeddings $i_i:S_i\hookrightarrow\tilde S_i$ with $\tilde S_I\in\SmVar(k)$. 
We have the canonical isomorphisms, given by, using definition \ref{TfjBdr1}, definition \ref{TfjBdr2},
definition \ref{TOBdr} and definition \ref{TDBdr}, for $K,K'\in D_{\mathbb Z_p,c,k}(S_{\bar k}^{et})$, 
\begin{eqnarray*}
\mathbb B_{dr,(\tilde S_I)}(K_1,K_2):R\Hom(K_1,K_2)\otimes\mathbb B_{dr,\bar k}
\xrightarrow{=}\mathbb B_{dr,\bar k}(Ra_{\bar S*}Rj_*\mathcal Hom(K_1,K_2)) \\
\xrightarrow{T(a_{\bar S},\mathbb B_{dr})(-)}Ra_{\bar S*}\mathbb B_{dr,(\bar\tilde S_I)}(Rj_*\mathcal Hom(K_1,K_2)) \\
\xrightarrow{\mathbb B_{dr,(\bar\tilde S_I)}(T(j,\mathbb B_{dr})(-))}
Ra_{\bar S*}(V_{D,0}j_*\mathbb B_{dr,(\tilde S_I)}(\mathcal Hom(K_1,K_2))
\otimes_{\mathbb B_{dr,S}}(\mathbb B_{dr,\tilde S_I/\bar\tilde S_I},t_{IJ})) \\
\xrightarrow{\mathbb B_{dr,S}(m(K_1,K_2)^{-1})}
Ra_{\bar S*}(V_{D,0}j_*\mathbb B_{dr,(\tilde S_I)}(\mathbb D_S^vK_1\otimes K_2)
\otimes_{\mathbb B_{dr,S}}(\mathbb B_{dr,\tilde S_I/\bar\tilde S_I},t_{IJ})) \\
\xrightarrow{(T(D,\mathbb B_{dr})(K_1)\otimes I)\circ T(\otimes,\mathbb B_{dr})(\mathbb D_SK_1,K_2)} \\
Ra_{\bar S*}(V_{D,0}j_*(\mathbb D_S\mathbb B_{dr,(\tilde S_I)}(K_1)\otimes_{\mathbb B_{dr,S}}\mathbb B_{dr,(\tilde S_I)}(K_2))
\otimes_{\mathbb B_{dr,S}}(\mathbb B_{dr,\tilde S_I/\bar\tilde S_I},t_{IJ})) \\
\xrightarrow{m(\mathbb B_{dr,(\tilde S_I)}(K_1),\mathbb B_{dr,(\tilde S_I)}(K_2))}
Ra_{\bar S*}(V_{D,0}j_*\mathcal Hom(\mathbb B_{dr,(\tilde S_I)}(K_1),\mathbb B_{dr,(\tilde S_I)}(K_2))
\otimes_{\mathbb B_{dr,S}}(\mathbb B_{dr,\tilde S_I/\bar\tilde S_I},t_{IJ})) \\
\xrightarrow{=}R\Hom(\mathbb B_{dr,S}(K_1),\mathbb B_{dr,S}(K_2)).
\end{eqnarray*}
\end{itemize}
\end{thm}

\begin{proof}
Follows from theorem \ref{TfBdrthm}.
\end{proof}

Let $S\in\SmVar(k)$ and $D\subset S$ a (Cartier) divisor. 
We have by theorem \ref{phipsithmp} the following isomorphisms
\begin{eqnarray*}
F^0T^{B_{dr}}(O_S,F_b):\mathbb B_{dr,\psi_D}:=F^0DR(S)(\psi_D(O_S,F_b)^{an}\otimes_{O_S}(O\mathbb B_{dr,S},F)) \\
\xrightarrow{\sim}F^0\psi_DDR(S)(O\mathbb B_{dr,S},F)
\xrightarrow{=}\psi_DF^0DR(S)(O\mathbb B_{dr,S},F)=:\psi_D\mathbb B_{dr,S}
\end{eqnarray*}
and
\begin{eqnarray*}
F^0T^{'B_{dr}}(O_S,F_b):\mathbb B_{dr,\phi_D}:=F^0DR(S)(\phi_D(O_S,F_b)^{an}\otimes_{O_S}(O\mathbb B_{dr,S},F)) \\
\xrightarrow{\sim}F^0\phi_DDR(S)(O\mathbb B_{dr,S},F)
\xrightarrow{=}\phi_DF^0DR(S)(O\mathbb B_{dr,S},F)=:\phi_D\mathbb B_{dr,S}.
\end{eqnarray*}

\begin{defi}\label{TphipsiBdr}
\begin{itemize}
\item[(i)]Let $j:S^o\hookrightarrow S$ an open embedding with $S\in\SmVar(k)$ and $D:=S\backslash S^o$ a (Cartier) divisor. 
We will consider, for $(K,W)\in P_{\mathbb Z_pfil,k}(S^{et})$, the canonical isomorphism in $D_{\mathbb B_{dr,S}fil}(S_K^{an,pet})$
\begin{eqnarray*}
T(\psi_D,\mathbb B_{dr})(K,W):\mathbb B_{dr,S}(\psi_D(K,W))\xrightarrow{:=} \\
(\cdots\to\xi_{\bar E_d}\cdots\xi_{\bar E_1}(\psi_D((K,W)\otimes_{\mathbb Q_p} O_{S^o})) 
\otimes_{O_S}\mathbb B_{dr,\psi_D}\otimes_{\mathbb B_{dr,S}}\mathbb B_{dr,S\backslash\bar E_{i_1}/S}
\otimes_{\mathbb B_{dr,S}}\cdots\otimes\mathbb B_{dr,S\backslash \bar E_{i_r}/S} \\ 
\otimes_{\mathbb B_{dr,S}}\mathbb B_{dr,\phi_{\bar E_{i_{r+1}}}}\otimes\cdots\otimes\mathbb B_{dr,\phi_{\bar E_{i_s}}}
\otimes_{\mathbb B_{dr,S}}\mathbb B_{dr,\psi_{\bar E_{i_{s+1}}}}\otimes\cdots\otimes\mathbb B_{dr,\psi_{\bar E_{i_d}}} 
\to\cdots) \\
\xrightarrow{=}
(\cdots\to\psi_D\xi_{\bar E_d}\cdots\xi_{\bar E_1}((K,W)\otimes_{\mathbb Q_p} O_{S\backslash\bar E_1})
\otimes_{O_S}\mathbb B_{dr,x_{S\backslash\bar E_{i_1}/S}}
\otimes_{\mathbb B_{dr,S}}\cdots\otimes\mathbb B_{dr,x_{S\backslash \bar E_{i_r}/S}} \\ 
\otimes_{\mathbb B_{dr,S}}\mathbb B_{dr,\phi_{\bar E_{i_{r+1}}}}\otimes\cdots\otimes\mathbb B_{dr,\phi_{\bar E_{i_s}}}
\otimes_{\mathbb B_{dr,S}}\mathbb B_{dr,\psi_{\bar E_{i_{s+1}}}}\otimes\cdots\otimes\mathbb B_{dr,\psi_{\bar E_{i_d}}} 
\to\cdots)
\xrightarrow{=:}\psi_D\mathbb B_{dr,S}((K,W))
\end{eqnarray*}
using definition \ref{VfilKMmap}(vi).
We will also consider, for $(K,W)\in D_{\mathbb Z_pfil,c,k}(S^{o,et})$, 
the canonical isomorphism in $D_{\mathbb B_{dr}fil}(S_K^{an,pet})$
\begin{eqnarray*}
T(\phi_D,\mathbb B_{dr})(K,W):\mathbb B_{dr,S}(\phi_D(K,W))=\mathbb B_{dr,S}(\mathbb D_S^v\psi_D\mathbb D_S^v(K,W)) \\ 
\xrightarrow{T(D,\mathbb B_{dr})(\psi_D\mathbb D^v_S(K,W))}\mathbb D_S\mathbb B_{dr,S}(\psi_D\mathbb D^v_{S^o}(K,W)) 
\xrightarrow{\mathbb D_ST(\psi_D,\mathbb B_{dr})(\mathbb D_S^v(K,W))} \\
\mathbb D_S(\psi_D\mathbb B_{dr,S}(\mathbb D_S^v(K,W))) 
\xrightarrow{T(D,\mathbb B_{dr})(K,W)}\mathbb D_S\psi_D\mathbb D_S\mathbb B_{dr,S}(K,W)
\xrightarrow{=}\phi_D\mathbb B_{dr,S}(K,W),
\end{eqnarray*}
using definition \ref{TDBdr}.
\item[(ii)]Let $l:S^o\hookrightarrow S$ an open embedding with $S\in\Var(k)$ such that $D=S\backslash S^o\subset S$ is a Cartier divisor.
Let $S=\cup_i S_i$ an open affine cover and 
$i_i:S_i\hookrightarrow\tilde S_i$ closed embedding with $\tilde S_i\in\SmVar(k)$. 
Let $l_I:\tilde S^o_I\hookrightarrow\tilde S_I$ open embeddings such that $\tilde S^o_I\cap S=S^o\cap S_I$
and $\tilde D_I\subset\tilde S_I$ a Cartier divisor such that $D\cap S_I\subset\tilde D_I\cap S$.
We will consider, for $(K,W)\in P_{\mathbb Z_pfil,k}(S^{o,et})$, the canonical isomorphism 
in $D_{\mathbb B_{dr}fil}(S_K^{an,pet}/(\tilde S_{I,K})^{an,pet})$
\begin{eqnarray*}
T(\psi_D,\mathbb B_{dr})(K,W):\mathbb B_{dr,(\tilde S_I)}(\psi_D(K,W))\xrightarrow{:=} 
((\cdots\to\xi_{\tilde E_{d,I}}\cdots\xi_{\tilde E_{1,I}}\psi^u_{\tilde D_I}(i_{I*}j_I^*(K,W)\otimes O_{\tilde S_I}) \\
\otimes_{O_{\tilde S_I}}\mathbb B_{dr,\psi_{\tilde D_I}}\otimes_{\mathbb B_{dr,\tilde S_I}}
\mathbb B_{dr,\tilde S_I\backslash\tilde E_{i_1,I}/S}
\otimes_{\mathbb B_{dr,\tilde S_I}}\cdots\otimes\mathbb B_{dr,\tilde S_I\backslash\tilde E_{i_r,I}/S} \\
\otimes_{\mathbb B_{dr,\tilde S_I}}\mathbb B_{dr,\phi_{\tilde E_{i_{r+1},I}}}
\otimes\cdots\otimes\mathbb B_{dr,\phi_{\tilde E_{i_s,I}}}
\otimes_{\mathbb B_{dr,\tilde S_I}}\mathbb B_{dr,\psi_{\tilde E_{i_{s+1},I}}}
\otimes\cdots\otimes\mathbb B_{dr,\psi_{\tilde E_{i_d,I}}}\to\cdots),\mathbb B_{dr}(t_{IJ})) \\
\xrightarrow{=} 
((\cdots\to\psi^u_{\tilde D_I}\xi_{\tilde E_{d,I}}\cdots\xi_{\tilde E_{1,I}}(i_{I*}j_I^*(K,W)\otimes O_{\tilde S_I})
\otimes_{O_{\tilde S_I}}\mathbb B_{dr,\tilde S_I\backslash\tilde E_{i_1,I}/S}
\otimes_{\mathbb B_{dr,\tilde S_I}}\cdots\otimes\mathbb B_{dr,\tilde S_I\backslash\tilde E_{i_r,I}/S} \\
\otimes_{\mathbb B_{dr,\tilde S_I}}\mathbb B_{dr,\phi_{\tilde E_{i_{r+1},I}}}
\otimes\cdots\otimes\mathbb B_{dr,\phi_{\tilde E_{i_s,I}}}
\otimes_{\mathbb B_{dr,\tilde S_I}}\mathbb B_{dr,\psi_{\tilde E_{i_{s+1},I}}}
\otimes\cdots\otimes\mathbb B_{dr,\psi_{\tilde E_{i_d,I}}}))\to\cdots),\mathbb B_{dr}(t_{IJ})) \\
\xrightarrow{=:}\psi_D\mathbb B_{dr,(\tilde S_I)}(K,W)
\end{eqnarray*}
using definition \ref{VfilKMmap}(vi),
where $(E_1,\ldots,E_d)\in\mathcal S(K)$ is a stratification by Cartier divisor $E_i\subset S^o$, $1\leq i\leq d$, such that 
\begin{equation*}
K_{|E(r)\backslash E(r+1)}:=l_r^*K\in D_{\mathbb Z_p,c}((E(r)\backslash E(r+1))^{et}) 
\end{equation*}
are local systems for all $1\leq r\leq d$, $l_r:E(r)\hookrightarrow S^o$ being the locally closed embeddings,
and $\tilde E_{s,I}\subset\tilde S_I$, $\tilde D_I\subset\tilde S_I$ are (Cartier) divisor 
such that $\bar E_s\cap S_I\subset\tilde E_{s,I}\cap S_I$ and $D\cap S_I\subset\tilde D_I\cap S_I$.
We will also consider, for $(K,W)\in D_{\mathbb Z_pfil,c,k}(S^{et})$, the canonical isomorphism 
in $D_{\mathbb B_{dr}fil}(S_K^{an,pet}/(\tilde S_{I,K})^{an,pet})$
\begin{eqnarray*}
T(\phi_D,\mathbb B_{dr})(K,W):\mathbb B_{dr,(\tilde S_I)}(\phi_D(K,W))=
\mathbb B_{dr,(\tilde S_I)}(\mathbb D^v_S\psi_D\mathbb D^v_S(K,W)) \\
\xrightarrow{T(D,\mathbb B_{dr})(\psi_D\mathbb D^v_S(K,W))}
\mathbb D_S\mathbb B_{dr,(\tilde S_I)}(\psi_D\mathbb D_S^v(K,W)) 
\xrightarrow{\mathbb D_ST(\psi_D,\mathbb B_{dr})(\mathbb D_S^v(K,W))} \\
\mathbb D_S\psi_D\mathbb B_{dr,(\tilde S_I)}(\mathbb D_S^v(K,W))
\xrightarrow{T(D,\mathbb B_{dr})(K,W)}\mathbb D_S\psi_D\mathbb D_S\mathbb B_{dr,(\tilde S_I)}(K,W)
\xrightarrow{=}\phi_D\mathbb B_{dr,(\tilde S_I)}(K,W),
\end{eqnarray*}
using definition \ref{TDBdr}. It gives for $(K,W)\in D_{\mathbb Z_pfil,c,k}(S^{o,et})$, the canonical isomorphisms 
\begin{eqnarray*}
T(\psi_D,\mathbb B_{dr})(K,W):\mathbb B_{dr,(\tilde S_I)}(\psi_D(K,W))\xrightarrow{\sim}
\psi_D\mathbb B_{dr,(\tilde S_I)}(K,W).
\end{eqnarray*}
and
\begin{eqnarray*}
T(\phi_D,\mathbb B_{dr})(K,W):\mathbb B_{dr,(\tilde S_I)}(\phi_D(K,W))\xrightarrow{\sim}
\phi_D\mathbb B_{dr,(\tilde S_I)}(K,W).
\end{eqnarray*}
in $D_{\mathbb B_{dr}fil}(S_K^{an,pet}/(\tilde S_{I,K})^{an,pet})$.
\end{itemize}
\end{defi}

\subsubsection{The geometric $p$-adic Mixed Hodge Modules}

Let $p$ a prime integer. Let $k\subset K\subset\mathbb C_p$ a subfield of a $p$-adic field. 
Denote by $\bar k\subset\mathbb C_p$ its algebraic closure.
Recall $G=\Gal(\bar K,K)\subset\Gal(\bar k,k)$ denotes the Galois group of $K$.

For $S\in\Var(k)$, we denote for short $O_S:=O_{S^{an}_{\mathbb C_p}}$,
$\mathbb B_{dr,S}:=\mathbb B_{dr,S_{\mathbb C_p}}:=\mathbb B_{dr,R_{\mathbb C_p}(S^{an}_{\mathbb C_p})}$ and 
$O\mathbb B_{dr,S}:=O\mathbb B_{dr,S_{\mathbb C_p}}:=O\mathbb B_{dr,R_{\mathbb C_p}(S^{an}_{\mathbb C_p})}$.
where $R_{\mathbb C_p}:\AnSp(\mathbb C_p)\to\AdSp/(\mathbb C_p,O_{\mathbb C_p})$ is the canonical functor (see section 2).

Let $S\in\Var(k)$.
Recall that $S^{et}\subset\Var(k)^{sm}/S$ denote the small etale site.
We then have the morphism of site $\an_S:S^{an,pet}:=S_{\mathbb C_p}^{an,pet}\to S^{et}$ given by the analytical functor
where $S_{\mathbb C_p}^{an,pet}\subset(\AnSp(\mathbb C_p)^{sm}/S)^{pro}$ is the small pro-etale site.
Then, $\PSh_{\mathbb B_{dr},G,fil}(S_{\mathbb C_p}^{an,pet})$ is the category 
whose objects are $\pi_{K/\mathbb C_p}^{*mod}(N,F)$ with $N\in\PSh_{\mathbb B_{dr}fil}(S^{an,pet})$
together with a continuous action of $G$ compatible with the $\mathbb B_{dr,S}$ module structure.

\begin{itemize}
\item Let $S\in\SmVar(k)$. 
The category $C_{\mathcal D(1,0)fil,rh}(S)\times_I D_{\mathbb Z_pfil,c,k}(S^{et})$ is the category 
\begin{itemize}
\item whose set of objects is the set of triples $\left\{((M,F,W),(K,W),\alpha)\right\}$ with 
\begin{eqnarray*}
(M,F,W)\in C_{\mathcal D(1,0)fil,rh}(S), \, (K,W)\in D_{\mathbb Z_pfil,c,k}(S^{et}), \\ 
\alpha:\mathbb B_{dr,S}(K,W)\to F^0DR(S)^{[-]}((M,F,W)^{an}\otimes_{O_S}(O\mathbb B_{dr,S},F))
\end{eqnarray*}
where 
\begin{itemize}
\item we recall that 
\begin{equation*}
DR(S)^{[-]}=DR(S_{\mathbb C_p}^{an})^{[-]}:
C_{\mathcal D(1,0)fil,rh}(S_{\mathbb C_p}^{an})\to C_{\mathbb B_{dr}2fil}(S_{\mathbb C_p}^{an,pet}) 
\end{equation*}
is the De Rahm functor (for $S'\subset S$ a connected component of dimension $d$, $DR(S)^{[-]}_{|S'}=DR(S)_{|S'}[d]$),
\item the functor
\begin{equation*}
\mathbb B_{dr,S}:D_{\mathbb Z_pfil,c,k}(S^{et})\to D_{\mathbb B_{dr}fil}(S_{\mathbb C_p}^{an,pet}) 
\end{equation*}
is the functor from complexes of presheaves with constructible etale cohomology to 
complexes of $\mathbb B_{dr,S}$ modules given in definition \ref{Bdr}
(recall that for $L$ a local system, 
it is given by $\mathbb B_{dr,S}(L):=\an_S^*L\otimes_{\mathbb Q_p}\mathbb B_{dr,S_{\mathbb C_p}}$),
\item $\alpha$ is a morphism in $D_{\mathbb B_{dr},G,fil}(S_{\mathbb C_p}^{an,pet})$, that is
a morphism in $D_{\mathbb B_{dr},fil}(S_{\mathbb C_p}^{an,pet})$ compatible with the action
of the galois group $G=\Gal(\bar K,K)\subset\Gal(\bar k,k)$,
\end{itemize}
\item and whose set of morphisms are 
\begin{equation*}
\phi=(\phi_D,\phi_C,[\theta]):((M_1,F,W),(K_1,W),\alpha_1)\to((M_2,F,W),(K_2,W),\alpha_2)
\end{equation*}
where $\phi_D:(M_1,F,W)\to(M_2,F,W)$ and $\phi_C:(K_1,W)\to (K_2,W)$ are morphisms and
\begin{eqnarray*}
\theta=(\theta^{\bullet},I(F^0DR(S)(\phi^{an}_D\otimes I))\circ I(\alpha_1),I(\alpha_2)\circ I(\mathbb B_{dr,S}(\phi_C))): \\
I(\mathbb B_{dr,S}(K_1,W))[1]\to I(F^0DR(S)((M,F,W)^{an}\otimes_{O_S}(O\mathbb B_{dr,S},F))  
\end{eqnarray*}
is an homotopy, 
$I:C_{\mathbb B_{dr,S},G,fil}(S_{\mathbb C_p}^{an,pet})\to K_{\mathbb B_{dr,S},G,fil}(S_{\mathbb C_p}^{an,pet})$
being the injective resolution functor : for $(N,W)\in C_{\mathbb B_{dr,S},G,fil}(S_{\mathbb C_p}^{an,pet})$, 
$k:(N,W)\to I(N,W)$ with $I(N,W)\in C_{\mathbb B_{dr,S},G,fil}(S_{\mathbb C_p}^{an,pet})$ 
is an injective resolution, and the class $[\theta]$ of $\theta$ does NOT depend of the injective resolution ;
in particular 
\begin{equation*}
F^0DR(S)^{[-]}(\phi^{an}_D\otimes I)\circ\alpha_1=\alpha_2\circ\mathbb B_{dr,S}(\phi_C) 
\end{equation*}
in $D_{\mathbb B_{dr,S},G,fil}(S_{\mathbb C_p}^{an,pet})$, and for
\begin{itemize}
\item $\phi=(\phi_D,\phi_C,[\theta]):((M_1,F,W),(K_1,W),\alpha_1)\to((M_2,F,W),(K_2,W),\alpha_2)$
\item $\phi'=(\phi'_D,\phi'_C,[\theta']):((M_2,F,W),(K_2,W),\alpha_2)\to((M_2,F,W),(K_3,W),\alpha_3)$
\end{itemize}
the composition law is given by 
\begin{eqnarray*}
\phi'\circ\phi:=(\phi'_D\circ\phi_D,\phi'_C\circ\phi_C,
I(F^0DR(S)(\phi^{'an}_C\otimes I))\circ[\theta]+[\theta']\circ I(\mathbb B_{dr,S}(\phi_C))[1]): \\
((M_1,F,W),(K_1,W),\alpha_1)\to((M_3,F,W),(K_3,W),\alpha_3),
\end{eqnarray*}
in particular for $((M,F,W),(K,W),\alpha)\in C_{\mathcal D(1,0)fil,rh}(S)\times_I D_{\mathbb Z_pfil,c,k}(S^{et})$,
\begin{equation*}
I_{((M,F,W),(K,W),\alpha)}=(I_M,I_K,0).
\end{equation*}
\end{itemize}
We have then the full embedding
\begin{eqnarray*}
\PSh_{\mathcal D(1,0)fil,rh}(S)\times_I P_{\mathbb Z_pfil,k}(S^{et})
\hookrightarrow C_{\mathcal D(1,0)fil,rh}(S)\times_I D_{\mathbb Z_pfil,c,k}(S^{et})
\end{eqnarray*}
where the category $\PSh_{\mathcal D(1,0)fil,rh}(S)\times_I P_{\mathbb Z_p,fil}(S^{et})$ is the category 
\begin{itemize}
\item whose set of objects is the set of triples $\left\{((M,F,W),(K,W),\alpha)\right\}$ with 
\begin{eqnarray*}
(M,F,W)\in\PSh_{\mathcal D(1,0)fil,rh}(S), \, (K,W)\in P_{\mathbb Z_pfil,k}(S^{et}), \\ 
\alpha:\mathbb B_{dr,S}(K,W)\to F^0DR(S)^{[-]}((M,F,W)^{an}\otimes_{O_S}(O\mathbb B_{dr,S},F))
\end{eqnarray*}
where $\alpha$ is an isomorphism,
\item and whose set of morphisms are 
\begin{equation*}
\phi=(\phi_D,\phi_C)=(\phi_D,\phi_C,0):((M_1,F,W),(K_1,W),\alpha_1)\to((M_2,F,W),(K_2,W),\alpha_2)
\end{equation*}
where $\phi_D:(M_1,F,W)\to(M_2,F,W)$ and $\phi_C:(K_1,W)\to (K_2,W)$ are morphisms such that
\begin{equation*}
F^0DR(S)^{[-]}(\phi^{an}_D\otimes I)\circ\alpha_1=\alpha_2\circ\mathbb B_{dr,S}(\phi_C) 
\end{equation*}
in $P_{\mathbb B_{dr,S},G,fil}(S_{\mathbb C_p}^{an,pet})$.
\end{itemize}
\item Let $S\in\Var(k)$. Let $S=\cup_{i\in I}S_i$ an open cover such that there
exists closed embeddings $i_i:S_i\hookrightarrow\tilde S_i$ with $\tilde S_I\in\SmVar(k)$.
The category $C_{\mathcal D(1,0)fil,rh}(S/(\tilde S_I))\times_I D_{\mathbb Z_pfil,c,k}(S^{et})$ 
is the category 
\begin{itemize}
\item whose set of objects is the set of triples $\left\{(((M_I,F,W),u_{IJ}),(K,W),\alpha)\right\}$ with
\begin{eqnarray*} 
((M_I,F,W),u_{IJ})\in C_{\mathcal D(1,0)fil,rh}(S/(\tilde S_I)), \, (K,W)\in D_{\mathbb Z_pfil,c,k}(S^{et}), \\ 
\alpha:\mathbb B_{dr,(\tilde S_I)}(K,W)\to
F^0DR(S)^{[-]}(((M_I,F,W),u_{IJ})^{an}\otimes_{O_S}((O\mathbb B_{dr,\tilde S_I},F),t_{IJ}))
\end{eqnarray*}
where 
\begin{itemize}
\item the functor
\begin{equation*}
DR(S)^{[-]}=DR(S_{\mathbb C_p}^{an})^{[-]}:
C_{\mathcal D(1,0)fil,rh}(S_{\mathbb C_p}^{an}/(\tilde S_{I,\mathbb C_p}^{an,pet}))
\to C_{\mathbb B_{dr}2fil}(S_{\mathbb C_p}^{an,pet}/(\tilde S_{I,\mathbb C_p}^{an,pet})) 
\end{equation*}
is the De Rahm functor, 
\item the functor 
\begin{equation*}
\mathbb B_{dr,(\tilde S_I)}:D_{\mathbb Z_pfil,c,k}(S^{et})\to 
D_{\mathbb B_{dr}fil}(S_{\mathbb C_p}^{an,pet}/(\tilde S_{I,\mathbb C_p}^{an,pet})) 
\end{equation*}
is the functor from complexes of presheaves with constructible etale cohomology 
to complexes of $\mathbb B_{dr}$ modules given in definition \ref{Bdr},
\item $\alpha$ is a morphism in $D_{\mathbb B_{dr},G,fil}(S_{\mathbb C_p}^{an,pet}/(\tilde S_{I,\mathbb C_p}^{an,pet}))$,
that is a morphism in $D_{\mathbb B_{dr},fil}(S_{\mathbb C_p}^{an,pet}/(\tilde S_{I,\mathbb C_p}^{an,pet}))$
compatible with the action of the galois group $G=\Gal(\bar{\hat k},\hat k)\subset\Gal(\bar k,k)$
\end{itemize}
\item and whose set of morphisms are 
\begin{equation*}
\phi=(\phi_D,\phi_C,[\theta]):(((M_{1I},F,W),u_{IJ}),(K_1,W),\alpha_1)\to(((M_{2I},F,W),u_{IJ}),(K_2,W),\alpha_2)
\end{equation*}
where $\phi_D:((M_1,F,W),u_{IJ})\to((M_2,F,W),u_{IJ})$ and 
$\phi_C:(K_1,W)\to (K_2,W)$ are morphisms (of filtered complexes) and 
\begin{eqnarray*}
\theta=(\theta^{\bullet},I(F^0DR(S)(\phi_D\otimes I))\circ I(\alpha_1),I(\alpha_2)\circ I(\mathbb B_{dr,S}(\phi_C))): \\
I(\mathbb B_{dr,(\tilde S_I)}(K_1,W))[1]\to 
I(DR(S)(((M_{2I},F,W),u_{IJ})^{an}\otimes_{O_S}((O\mathbb B_{dr,(\tilde S_I)},F),t_{IJ})))  
\end{eqnarray*}
is an homotopy,  
$I:C_{\mathbb B_{dr},G,fil}(S_{\mathbb C_p}^{an}/(\tilde S^{an}_{I,\mathbb C_p}))
\to K_{\mathbb B_{dr},G,fil}(S_{\mathbb C_p}^{an}/(\tilde S^{an}_{I,\mathbb C_p}))$
being the injective resolution functor : for
$((N_I,W),t_{IJ})\in C_{\mathbb B_{dr},G,fil}(S_{\mathbb C_p}^{an}/(\tilde S^{an}_{I,\mathbb C_p}))$, 
\begin{equation*}
k:((N_I,W),t_{IJ})\to I((N_I,W),t_{IJ}) 
\end{equation*}
with $I((N_I,W),t_{IJ})\in C_{\mathbb B_{dr},G,fil}(S_{\mathbb C_p}^{an}/(\tilde S^{an}_{I,\mathbb C_p}))$ 
is an injective resolution, and the class $[\theta]$ of $\theta$ does NOT depend of the injective resolution ;
in particular we have
\begin{equation*}
F^0DR(S)^{[-]}(\phi^{an}_D\otimes I)\circ\alpha_1=\alpha_2\circ\mathbb B_{dr,(\tilde S_I)}(\phi_C) 
\end{equation*}
in $D_{\mathbb B_{dr},G,fil}(S_{\mathbb C_p}/\tilde S_{I,\mathbb C_p}^{an,pet})$, and for
\begin{itemize}
\item $\phi=(\phi_D,\phi_C,[\theta]):(((M_{1I},F,W),u_{IJ}),(K_1,W),\alpha_1)\to(((M_{2I},F,W),u_{IJ}),(K_2,W),\alpha_2)$
\item $\phi'=(\phi'_D,\phi'_C,[\theta']):(((M_{2I},F,W),u_{IJ}),(K_2,W),\alpha_2)\to(((M_{3I},F,W),u_{IJ}),(K_3,W),\alpha_3)$
\end{itemize}
the composition law is given by 
\begin{eqnarray*}
\phi'\circ\phi:=(\phi'_D\circ\phi_D,\phi'_C\circ\phi_C,
I(F^0DR(S)(\phi^{'an}_D\otimes I))\circ[\theta]+[\theta']\circ I(\mathbb B_{dr,(\tilde S_I)}(\phi_C))[1]): \\
(((M_{1I},F,W),u_{IJ}),(K_1,W),\alpha_1)\to(((M_{3I},F,W),u_{IJ}),(K_3,W),\alpha_3).
\end{eqnarray*}
in particular for 
$(((M_I,F,W),u_{IJ}),(K,W),\alpha)\in C_{\mathcal D(1,0)fil,rh}(S/(\tilde S_I))\times_I D_{\mathbb Z_pfil,c,k}(S^{et})$,
\begin{equation*}
I_{(((M_I,F,W),u_{IJ}),(K,W),\alpha)}=((I_{M_I}),I_K,0).
\end{equation*}
\end{itemize}
We have then full embeddings
\begin{eqnarray*}
\PSh^0_{\mathcal D(1,0)fil,rh}(S/(\tilde S_I))\times_I P_{\mathbb Z_pfil,k}(S^{et})
\hookrightarrow C^0_{\mathcal D(1,0)fil,rh}(S/(\tilde S_I))\times_I D_{\mathbb Z_pfil,c,k}(S^{et}) \\
\xrightarrow{\iota^0_{S/\tilde S_I}} C_{\mathcal D(1,0)fil,rh}(S/(\tilde S_I))^0\times_I D_{\mathbb Z_pfil,c,k}(S^{et})
\hookrightarrow C_{\mathcal D(1,0)fil,rh}(S/(\tilde S_I))\times_I D_{\mathbb Z_pfil,c,k}(S^{et})
\end{eqnarray*}
where the category $\PSh^0_{\mathcal D(1,0)fil,rh}(S/(\tilde S_I))\times_I P_{\mathbb Z_pfil,k}(S^{et})$ is the category 
\begin{itemize}
\item whose set of objects is the set of triples $\left\{(((M_I,F,W),u_{IJ}),(K,W),\alpha)\right\}$ with 
\begin{eqnarray*}
((M_I,F,W),u_{IJ})\in\PSh_{\mathcal D(1,0)fil,rh}(S/(\tilde S_I)), \, (K,W)\in P_{\mathbb Z_pfil,k}(S^{et}), \\ 
\alpha:\mathbb B_{dr,(\tilde S_I)}(K,W)\to 
F^0DR(S)^{[-]}(((M_I,F,W),u_{IJ})^{an}\otimes_{O_S}(O\mathbb B_{dr,(\tilde S_I)},F))
\end{eqnarray*}
where $\alpha$ is an isomorphism,
\item and whose set of morphisms are 
\begin{equation*}
\phi=(\phi_D,\phi_C)=(\phi_D,\phi_C,0):(((M_{1I},F,W),u_{IJ}),(K_1,W),\alpha_1)\to
(((M_{2I},F,W),u_{IJ}),(K_2,W),\alpha_2)
\end{equation*}
where $\phi_D:((M_{1I},F,W),u_{IJ})\to((M_{2I},F,W),u_{IJ})$ and $\phi_C:(K_1,W)\to (K_2,W)$ are morphisms such that
\begin{equation*}
F^0DR(S)^{[-]}(\phi^{an}_D\otimes I)\circ\alpha_1=\alpha_2\circ\mathbb B_{dr,(\tilde S_I)}(\phi_C)  
\end{equation*}
in $P_{\mathbb B_{dr,S},G,fil}(S_{\mathbb C_p}^{an,pet}/(\tilde S_{I,\mathbb C_p}^{an,pet}))$.
\end{itemize}
\end{itemize}
Moreover,
\begin{itemize}
\item For 
$(((M_{I},F,W),u_{IJ}),(K,W),\alpha)\in C_{\mathcal D(1,0)fil,rh}(S/(\tilde S_I))\times_I D_{\mathbb Z_pfil,c,k}(S^{et})$, 
we set
\begin{equation*}
(((M_{I},F,W),u_{IJ}),(K,W),\alpha)[1]:=(((M_{I},F,W),u_{IJ})[1],(K,W)[1],\alpha[1]).
\end{equation*}
\item For 
\begin{equation*}
\phi=(\phi_D,\phi_C,[\theta]):(((M_{1I},F,W),u_{IJ}),(K_1,W),\alpha_1)\to(((M_{2I},F,W),u_{IJ}),(K_2,W),\alpha_2)
\end{equation*}
a morphism in $C_{\mathcal D(1,0)fil,rh}(S/(\tilde S_I))\times_I D_{\mathbb Z_pfil,c,k}(S^{et})$, 
we set (see \cite{CG} definition 3.12)
\begin{eqnarray*}
\Cone(\phi):=(\Cone(\phi_D),\Cone(\phi_C),((\alpha_1,\theta),(\alpha_2,0)))
\in D_{\mathcal D(1,0)fil,rh}(S/(\tilde S_I))\times_I D_{\mathbb Z_pfil,c,k}(S^{et}),
\end{eqnarray*}
$((\alpha_1,\theta),(\alpha_2,0))$ being the matrix given by the composition law, together with the canonical maps
\begin{itemize}
\item $c_1(-)=(c_1(\phi_D),c_1(\phi_C),0):(((M_{2I},F,W),u_{IJ}),(K_2,W),\alpha_2)\to\Cone(\phi)$
\item $c_2(-)=(c_2(\phi_D),c_2(\phi_C),0):\Cone(\phi)\to (((M_{1I},F,W),u_{IJ}),(K_1,W),\alpha_1)[1]$.
\end{itemize}
\end{itemize}

\begin{rem}\label{CGremp}
By \cite{CG} theorem 3.25, if 
\begin{equation*}
\phi=(\phi_D,\phi_C,[\theta]):(((M_1,F,W),u_{IJ}),(K_1,W),\alpha_1)\to(((M_2,F,W),u_{IJ}),(K_2,W),\alpha_2)
\end{equation*}
is a morphism in $C_{\mathcal D(1,0)fil,rh}(S/(\tilde S_I))\times_I D_{\mathbb Z_pfil,c,k}(S^{et})$
such that $\phi_D$ is a Zariski local equivalence and $\phi_C$ is an isomorphism then $\phi$ is an isomorphism.
\end{rem}

\begin{defi}\label{falphap}
Let $k\subset{\mathbb C_p}$ a subfield.
\begin{itemize}
\item[(i1)] Let $f:X\to S$ a proper morphism with $S,X\in\SmVar(k)$. Let 
\begin{equation*}
\alpha:\mathbb B_{dr,X}(K,W)\to F^0DR(X)((M,F,W)^{an}\otimes_{O_X}(O\mathbb B_{dr,X},F)) 
\end{equation*}
a morphism in $D_{\mathbb B_{dr},G,fil}(X_{\mathbb C_p}^{an,pet})$, with
\begin{equation*}
(M,F,W)\in C(DRM(X)), \; (K,W)\in D_{\mathbb Z_pfil,c,k,gm}(X^{et}). 
\end{equation*}
We then consider, using definition \ref{TfjBdr1} and definition \ref{TfDRpadic}, 
the map in $D_{\mathbb B_{dr},G,fil}(S_{\mathbb C_p}^{an,pet})$
\begin{eqnarray*}
f_*\alpha=f_*(\alpha):\mathbb B_{dr,S}(Rf_*(K,W))
\xrightarrow{T(f,\mathbb B_{dr})(K,W)}Rf_*\mathbb B_{dr,X}(K,W) \\
\xrightarrow{Rf_*\alpha}Rf_*F^0DR(X)((M,F,W)^{an}\otimes_{O_X}(O\mathbb B_{dr,X},F)) \\
\xrightarrow{Rf_*\iota_{F^0}(-)}F^0Rf_*DR(X)((M,F,W)^{an}\otimes_{O_X}(O\mathbb B_{dr,X},F)) \\
\xrightarrow{F^0T^{B_{dr}}(f,DR)(M,F,W)^{-1}}
F^0DR(S)(\int_f((M,F,W)^{an})\otimes_{O_S}(O\mathbb B_{dr,S},F)) \\
\xrightarrow{F^0DR(S)(T(an,\int_f)(M,F,W)^{-1})}
F^0DR(S)((\int_f(M,F,W))^{an}\otimes_{O_S}(O\mathbb B_{dr,S},F)) \\
\xrightarrow{=}F^0DR(S)((Rf_{*Hdg}(M,F,W))^{an}\otimes_{O_S}(O\mathbb B_{dr,S},F))
\end{eqnarray*}
where $\iota_{F^0}(A)=D_{fil}(\iota_{F^0}(A))$ is the image of the embedding $\iota_{F^0}(A):F^0A\hookrightarrow A$
by the localization functor.
\item[(i2)] Let $j:S^o\hookrightarrow S$ an open embedding with $S\in\SmVar(k)$ and $D=S\backslash S^o$ a (Cartier) divisor. Let 
\begin{equation*}
\alpha:\mathbb B_{dr,S^o}(K,W)\to F^0DR(S^o)((M,F,W)^{an}\otimes_{O_{S^o}}(O\mathbb B_{dr,S^o},F)) 
\end{equation*}
a morphism in $D_{\mathbb B_{dr},G,fil}(S_{\mathbb C_p}^{o,an,pet})$, with
\begin{equation*}
(M,F,W)\in C(DRM(S^o)), \; (K,W)\in D_{\mathbb Z_pfil,c,k}(S^{o,et})^{ad,D}. 
\end{equation*}
We then consider, using definition \ref{TfjBdr2} and the strictness of the $V$-filtration, 
the maps in $D_{\mathbb B_{dr},G,fil}(S_{\mathbb C_p}^{an,pet})$
\begin{eqnarray*}
j_*\alpha=j_*(\alpha):\mathbb B_{dr,S}(j_{*w}(K,W))
\xrightarrow{T(j,\mathbb B_{dr})(K,W)}V_{D0}j_{*w}\mathbb B_{dr,S^o}(K,W)\otimes_{\mathbb B_{dr,S}}\mathbb B_{dr,S^o/S} \\
\xrightarrow{V_{D0}j_*\alpha\otimes I}
V_{D0}j_{*w}(F^0DR(S^o)((M,F,W)^{an}\otimes_{O_{S^o}}(O\mathbb B_{dr,S^o},F)))\otimes_{\mathbb B_{dr,S}}\mathbb B_{dr,S^o/S} \\
\xrightarrow{:=}
(V_{D0}j_{*w}F^0DR(S^o)((M,F,W)^{an}\otimes_{O_S}(O\mathbb B_{dr,S},F)))\otimes_{\mathbb B_{dr,S}} \\
(F^0DR(S)(j_{*Hdg}(O_{S^o},F_b)^{an}\otimes_{O_S}(O\mathbb B_{dr,S},F))) \\
\xrightarrow{w_S\otimes m(M)^{an}}
F^0DR(S)(j_{*Hdg}(M,F,W)^{an}\otimes_{O_S}(O\mathbb B_{dr,S},F))
\end{eqnarray*}
where $m(M):O_{S^o}\otimes_{O_{S^o}} M\xrightarrow{\sim} M, \; m\otimes f\mapsto fm$ 
is the multiplication map structure of the module $M$ and $w_S$ is the wedge product, and
\begin{eqnarray*}
j_!\alpha=j_!(\alpha):\mathbb B_{dr,S}(j_{!w}(K,W))\xrightarrow{T_!(j,\mathbb B_{dr})(K,W)}
V_{D0}\mathbb Dj_*\mathbb D\mathbb B_{dr,S^o}(K,W)\otimes_{\mathbb B_{dr,S}}\mathbb D_S\mathbb B_{dr,S^o/S}) \\
\xrightarrow{(V_{D0}\mathbb D j_*\mathbb D\alpha)\otimes I}
V_{D0}\mathbb Dj_{*w}\mathbb D(F^0DR(S^o)(\mathbb D(M,F,W)^{an}\otimes_{O_{S^o}}(O\mathbb B_{dr,S^o},F)))
\otimes_{\mathbb B_{dr,S}}\mathbb D_S\mathbb B_{dr,S^o/S} \\
\xrightarrow{:=} 
V_{D0}\mathbb Dj_{*w}\mathbb D(F^0DR(S^o)(\mathbb D(M,F,W)^{an}\otimes_{O_{S^o}}(O\mathbb B_{dr,S^o},F))) \\
\otimes_{\mathbb B_{dr,S}}(F^0DR(S)(j_{!Hdg}(O_{S^o},F_b)^{an}\otimes_{O_S}(O\mathbb B_{dr,S},F))) \\
\xrightarrow{w_S\otimes m(M)^{an}}
F^0DR(S)(j_{!Hdg}(M,F,W)^{an}\otimes_{O_S}(O\mathbb B_{dr,S},F))
\end{eqnarray*}
where $m(M):O_{S^o}\otimes_{O_{S^o}} M\xrightarrow{\sim} M, \; h\otimes m\mapsto hm$ is the multiplication map
and $w_S$ is the wedge product.
\item[(i2)']Let $l:S^o\hookrightarrow S$ an open embedding with $S\in\Var(k)$ and $D=S\backslash S^o$ a Cartier divisor.
Let $S=\cup_i S_i$ an open affine cover and 
$i_i:S_i\hookrightarrow\tilde S_i$ closed embedding with $\tilde S_i\in\SmVar(k)$. 
Let $l_I:\tilde S^o_I\hookrightarrow\tilde S_I$ closed embeddings such that $\tilde S^o_I\cap S=S^o\cap S_I$. Let 
\begin{equation*}
\alpha:\mathbb B_{dr,(\tilde S^o_I)}(K,W)\to
F^0DR(S^o)(((M_I,F,W),u_{IJ})^{an}\otimes_{O_{S^o}}((O\mathbb B_{dr,(\tilde S^o_I)},F),t_{IJ})) 
\end{equation*}
a morphism in $D_{\mathbb B_{dr},G,fil}(S_{\mathbb C_p}^{o,an,pet}/(\tilde S_{I,\mathbb C_p}^{o,an,pet}))$, with
\begin{equation*}
((M_I,F,W),u_{IJ})\in C(DRM(S^o))\subset C_{\mathcal D(1,0)fil,rh}(S^o/(\tilde S^o_I)), \; 
(K,W)\in D_{\mathbb Z_pfil,c,k}(S^{o,et})^{ad,D}. 
\end{equation*}
We then consider as in (i2) the maps in $D_{\mathbb B_{dr},G,fil}(S_{\mathbb C_p}^{an,pet}/(\tilde S_{I\mathbb C_p}^{an,pet}))$
\begin{eqnarray*}
l_*\alpha=l_*(\alpha):\mathbb B_{dr,(\tilde S_I)}(l_{*w}(K,W)) \\
\xrightarrow{T(l,\mathbb B_{dr})(K,W)}V_{D0}l_{*w}\mathbb B_{dr,(\tilde S^o_I)}(K,W)
\otimes_{\mathbb B_{dr,(\tilde S_I)}}(\mathbb B_{dr,\tilde S_I^o/\tilde S_I},t_{IJ}) \\
\xrightarrow{V_{D0}l_*\alpha\otimes I}
V_{D0}l_{*w}F^0DR(S^o)(((M_I,F,W),u_{IJ})^{an}\otimes_{O_S}((O\mathbb B_{dr,\tilde S_I^o},F),t_{IJ}))
\otimes_{\mathbb B_{dr,(\tilde S_I)}}(\mathbb B_{dr,\tilde S_I^o/\tilde S_I},t_{IJ}) \\
\xrightarrow{:=} \\
V_{D0}l_{*w}F^0DR(S^o)(((M_I,F,W),u_{IJ})^{an}\otimes_{O_S}((O\mathbb B_{dr,\tilde S_I^o},F),t_{IJ}))
\otimes_{\mathbb B_{dr,(\tilde S_I)}} \\
(F^0DR((\tilde S_I))((l_{I*Hdg}(O_{\tilde S_I^o},F_b),x_{IJ})^{an}\otimes_{O_S}((O\mathbb B_{dr,\tilde S_I},F),t_{IJ}))) \\
\xrightarrow{w_{\tilde S_I}\otimes m(M_I)^{an}}
F^0DR(S)((l_{*Hdg}((M_I,F,W),u_{IJ})^{an})\otimes_{O_S}((O\mathbb B_{dr,\tilde S_I},F),t_{IJ}))
\end{eqnarray*}
and
\begin{eqnarray*}
l_!\alpha=l_!(\alpha):\mathbb B_{dr,(\tilde S_I)}(l_{!w}(K,W)) \\
\xrightarrow{T_!(l,\mathbb B_{dr})(K,W)}
V_{D0}\mathbb Dl_*\mathbb D(\mathbb B_{dr,(\tilde S_I^o)}(K,W))
\otimes_{\mathbb B_{dr,(\tilde S_I})}\mathbb D(\mathbb B_{dr,\tilde S_I^o/\tilde S_I},t_{IJ}) 
\xrightarrow{V_{D0}\mathbb D l_*\mathbb D(\alpha\otimes I)} \\ 
V_{D0}\mathbb Dl_{*w}\mathbb D(F^0DR(S^o)(((M_I,F,W),u_{IJ})^{an}\otimes_{O_{S^o}}((O\mathbb B_{dr,\tilde S_I^o},F),t_{IJ}))) 
\otimes_{\mathbb B_{dr,(\tilde S_I})}\mathbb D(\mathbb B_{dr,\tilde S_I^o/\tilde S_I},t_{IJ}) \\
\xrightarrow{:=}
V_{D0}\mathbb Dl_{*w}\mathbb D(F^0DR(S^o)(((M_I,F,W),u_{IJ})^{an}\otimes_{O_{S^o}}((O\mathbb B_{dr,\tilde S_I^o},F),t_{IJ}))) 
\otimes_{\mathbb B_{dr,(\tilde S_I)}} \\
(F^0DR((\tilde S_I))((l_{I!Hdg}(O_{\tilde S_I^o},F_b),x_{IJ})^{an}\otimes_{O_S}((O\mathbb B_{dr,\tilde S_I},F),t_{IJ}))) \\
\xrightarrow{w_{\tilde S_I}\otimes m(M_I)^{an}}
F^0DR(S)(l_{!Hdg}((M_I,F,W),u_{IJ})^{an}\otimes_{O_S}((O\mathbb B_{dr,\tilde S_I},F),t_{IJ}))
\end{eqnarray*}
\item[(ii0)] Let $f:X\to S$ a morphism with $S,X\in\SmVar(k)$. 
Take a compactification $f:X\xrightarrow{j}\bar X\xrightarrow{\bar f} S$ of $f$ with $\bar X\in\SmVar(k)$,
$j$ an open embedding and $D=\bar X\backslash X$ a divisor 
(see section 2, we can take $D$ a normal crossing divisor but it is unnecessary). Let 
\begin{equation*}
\alpha:\mathbb B_{dr,X}(K,W)\to F^0DR(X)((M,F,W)^{an}\otimes_{O_X}O\mathbb B_{dr,X}) 
\end{equation*}
a morphism in $D_{\mathbb B_{dr},G,fil}(X_{\mathbb C_p}^{an,pet})$, with
\begin{equation*}
(M,F,W)\in C(DRM(X)), \; (K,W)\in D_{\mathbb Z_pfil,c,k,gm}(X^{et})^{ad,D}. 
\end{equation*}
We then consider, using (i1) and (i2) the maps in $D_{\mathbb B_{dr},G,fil}(S_{\mathbb C_p}^{an,pet})$
\begin{eqnarray*}
f_*\alpha=f_*(\alpha):\mathbb B_{dr,S}(Rf_{*w}(K,W))=\mathbb B_{dr,S}(R\bar f_*j_{*w}(K,W)) \\
\xrightarrow{T(\bar f,B_{dr})(-)}R\bar f_*\mathbb B_{dr,\bar X}(j_{*w}(K,W)) 
\xrightarrow{j_*(\alpha)}R\bar f_*F^0DR(X)(j_{*Hdg}(M,F,W)^{an}\otimes_{O_X}(O\mathbb B_{dr,X},F)) \\
\xrightarrow{F^0DR(S)(T(an,\int_f)(-)^{-1})\circ F^0T^{B_{dr}}(\bar f,DR)(-)^{-1}\circ R\bar f_*\iota_{F^0}} \\
F^0DR(S)(R\bar f_{*Hdg}j_{*Hdg}(M,F,W)^{an}\otimes_{O_S}(O\mathbb B_{dr,S},F)) \\
\xrightarrow{=:}F^0DR(S)(Rf_{*Hdg}(M,F,W)^{an}\otimes_{O_S}(O\mathbb B_{dr,S},F))
\end{eqnarray*}
and
\begin{eqnarray*}
f_!\alpha=f_!(\alpha):\mathbb B_{dr,S}(Rf_{!w}(K,W))=\mathbb B_{dr,S}(R\bar f_*j_{!w}(K,W)) \\
\xrightarrow{T(\bar f,B_{dr})(-)}R\bar f_*\mathbb B_{dr,\bar X}(j_{!w}(K,W)) 
\xrightarrow{j_!(\alpha)}R\bar f_*F^0DR(X)(j_{!Hdg}(M,F,W)^{an}\otimes_{O_X}(O\mathbb B_{dr,X},F)) \\
\xrightarrow{F^0DR(S)(T(an,\int_f)(-)^{-1})\circ F^0T^{B_{dr}}(\bar f,DR)(-)^{-1}\circ R\bar f_*\iota_{F^0}} \\
F^0DR(S)(R\bar f_{*Hdg}j_{!Hdg}(M,F,W)^{an}\otimes_{O_S}(O\mathbb B_{dr,S},F)) \\
\xrightarrow{=:}F^0DR(S)(Rf_{!Hdg}(M,F,W)^{an}\otimes_{O_S}(O\mathbb B_{dr,S},F)).
\end{eqnarray*}
\item[(ii)] Let $f:X\to S$ a morphism with $S,X\in\Var(k)$. 
Consider a factorization $f:X\hookrightarrow Y\times S\xrightarrow{p}S$ with $Y\in\SmVar(k)$,
and let $f:X\xrightarrow{j}\bar{X}\hookrightarrow\bar Y\times S\xrightarrow{\bar p}S$ be a compactification of $f$,
with $\bar Y\in\PSmVar(k)$ and $D=\bar Y\backslash Y$ a (Cartier) divisor (e.g. a normal crossing divisor). 
Let $S=\cup_i S_i$ an open affine cover and 
$i_i:S_i\hookrightarrow\tilde S_i$ closed embedding with $\tilde S_i\in\SmVar(k)$. Let 
\begin{equation*}
\alpha:\mathbb B_{dr,(Y\times\tilde S_I)}(K,W)\to 
F^0DR(X)(((M_I,F,W),u_{IJ})^{an}\otimes_{O_{Y\times\tilde S_I}}(O\mathbb B_{dr,Y\times\tilde S_I},F)) 
\end{equation*}
a morphism in $D_{\mathbb B_{dr},G,fil}(X_{\mathbb C_p}^{an,pet}/({Y\times\tilde S_I}_{\mathbb C_p}^{an,pet}))$, with
\begin{equation*}
((M_I,F,W),u_{IJ})\in C(DRM(X))\subset C_{\mathcal D(1,0)fil,rh}(X/(Y\times\tilde S_I)), \; 
(K,W)\in D_{\mathbb Z_pfil,c,k,gm}(X^{et})^{ad,D}. 
\end{equation*}
We then consider, using definitions \ref{TfjBdr1} and (i2)',
the maps in $D_{\mathbb B_{dr},G,fil}(S^{an,pet}_{\mathbb C_p}/(\tilde S_{I,\mathbb C_p}^{an,pet}))$
\begin{eqnarray*}
f_*\alpha=f_*(\alpha):\mathbb B_{dr,(\tilde S_I)}(Rf_{*w}(K,W))=\mathbb B_{dr,(\tilde S_I)}((R\bar f_*j_{*w}(K,W))) \\  
\xrightarrow{T(\bar f,\mathbb B_{dr})(-)}
R\bar{p}_*\mathbb B_{dr,(\bar Y\times\tilde S_I)}(j_{*w}(K,W)) \\
\xrightarrow{j_*\alpha}
R\bar{p}_*F^0DR(X)((j_{*Hdg}((M_I,F,W),u_{IJ})^{an})\otimes_{O_X}((O\mathbb B_{dr,\bar Y\times\tilde S_I},F),t_{IJ})) \\
\xrightarrow{Rp_*\iota_{F^0}}
F^0R\bar{p}_*DR(X)((j_{*Hdg}((M_I,F,W),u_{IJ})^{an})\otimes_{O_X}((O\mathbb B_{dr,\bar Y\times\tilde S_I},F),t_{IJ})) \\
\xrightarrow{F^0T^{B_{dr}}(f,DR)(j_{*Hdg}((M_I,F,W),u_{IJ}))}
F^0DR(S)(\int_{\bar f}(j_{*Hdg}((M_I,F,W),u_{IJ}))^{an}\otimes_{O_S}((O\mathbb B_{dr,\tilde S_I},F),t_{IJ})) \\
\xrightarrow{F^0DR(S)(T(an,\int_f)(-)\otimes I)}
F^0DR(S)((\int_{\bar f}j_{*Hdg}((M_I,F,W),u_{IJ}))^{an}\otimes_{O_S}((O\mathbb B_{dr,\tilde S_I},F),t_{IJ})) \\
\xrightarrow{=}
F^0DR(S)((Rf_{*Hdg}((M_I,F,W),u_{IJ}))^{an}\otimes_{O_S}(O\mathbb B_{dr,\tilde S_I},F)),
\end{eqnarray*}
and
\begin{eqnarray*}
f_!\alpha=f_!(\alpha):\mathbb B_{dr,(\tilde S_I)}(R\bar f_*j_{!w}(K,W)))
\xrightarrow{=}\mathbb B_{dr,(\tilde S_I)}(R\bar f_*j_{!w}(K,W)) \\
\xrightarrow{\mathbb D(T(\bar f,\mathbb B_{dr})(-))^{-1}}
R\bar{p}_*\mathbb B_{dr,(\bar Y\times\tilde S_I)}(j_{!w}(K,W)) \\
\xrightarrow{j_!\alpha}
R\bar{p}_*F^0DR(X)((j_{!Hdg}((M_I,F,W),u_{IJ})^{an})\otimes_{O_X}((O\mathbb B_{dr,\bar Y\times\tilde S_I},F),t_{IJ})) \\
\xrightarrow{Rp_*\iota_{F^0}}
F^0R\bar{p}_*DR(X)((j_{!Hdg}((M_I,F,W),u_{IJ})^{an})\otimes_{O_X}((O\mathbb B_{dr,\bar Y\times\tilde S_I},F),t_{IJ})) \\
\xrightarrow{F^0T^{B_{dr}}(f,DR)(j_{!Hdg}((M_I,F,W),u_{IJ}))}
F^0DR(S)(\int_{\bar f}(j_{!Hdg}((M_I,F,W),u_{IJ}))^{an}\otimes_{O_S}((O\mathbb B_{dr,\tilde S_I},F),t_{IJ})) \\
\xrightarrow{F^0DR(S)(T(an,\int_f)(-)\otimes I)}
F^0DR(S)((\int_{\bar f}j_{!Hdg}((M_I,F,W),u_{IJ}))^{an}\otimes_{O_S}((O\mathbb B_{dr,\tilde S_I},F),t_{IJ})) \\
\xrightarrow{=}
F^0DR(S)((Rf_{!Hdg}((M_I,F,W),u_{IJ}))^{an}\otimes_{O_S}((O\mathbb B_{dr,(\tilde S_I)},F),t_{IJ}).
\end{eqnarray*}
\item[(iii)] Let $l:S^o\hookrightarrow S$ an open embedding with $S\in\Var(k)$ and denote $Z=S\backslash S^o$.
Let $D_1,\cdots,D_d\subset S$ Cartier divisor such that $Z=\cap_{s=1}^dD_s$.
Denote $l_s:D_s\hookrightarrow S$ the closed embeddings.
Let $S=\cup_i S_i$ an open affine cover and 
$i_i:S_i\hookrightarrow\tilde S_i$ closed embedding with $\tilde S_i\in\SmVar(k)$. 
Let $l_{I,s}:\tilde D_{I,s}\hookrightarrow\tilde S_I$ closed embeddings 
such that $\tilde D_{I,s}\cap S=S\cap D_{I,s}$. Let 
\begin{equation*}
\alpha:\mathbb B_{dr,(\tilde S_I)}(K,W)\to
F^0DR(S)(((M_I,F,W),u_{IJ})^{an}\otimes_{O_{S}}((O\mathbb B_{dr,\tilde S_I},F),t_{IJ})) 
\end{equation*}
a morphism in $D_{\mathbb B_{dr},G,fil}(S_{\mathbb C_p}^{an,pet}/(\tilde S_{I,\mathbb C_p}^{an,pet}))$, with
\begin{equation*}
((M_I,F,W),u_{IJ})\in C(DRM(S))\subset C_{\mathcal D(1,0)fil,rh}(S/(\tilde S_I)), \; 
(K,W)\in D_{\mathbb Z_pfil,c,k}(S^{et})^{ad,(D_i)}. 
\end{equation*}
We then have by (i2), the maps in $D_{\mathbb B_{dr},G,fil}(S_{\mathbb C_p}^{an,pet}/(\tilde S_{I,\mathbb C_p}^{an,pet}))$
\begin{eqnarray*}
\Gamma_Z(\alpha):\mathbb B_{dr,(\tilde S_I)}(\Gamma^w_Z(K,W))
\xrightarrow{:=}\mathbb B_{dr,(\tilde S_I)}(\Gamma^w_{D_1}\cdots\Gamma^w_{D_s}(K,W)) \\
\xrightarrow{(I,(l_{1*}\cdots l_{s*}(\alpha)))}
F^0DR(S)((\Gamma^{Hdg}_Z((M_I,F,W),u_{IJ}))^{an}\otimes_{O_S}((O\mathbb B_{dr,\tilde S_I},F),t_{IJ}))
\end{eqnarray*}
and
\begin{eqnarray*}
\Gamma^{\vee}_Z(\alpha):\mathbb B_{dr,(\tilde S_I)}(\Gamma^{\vee,w}_Z(K,W))
\xrightarrow{:=}\mathbb B_{dr,(\tilde S_I)}(\Gamma^{\vee,w}_{D_1}\cdots\Gamma^{\vee,w}_{D_s}(K,W)) \\
\xrightarrow{(I,(l_{1!}\cdots l_{s!}(\alpha)))}
F^0DR(S)((\Gamma^{\vee,Hdg}_Z((M_I,F,W),u_{IJ}))^{an}\otimes_{O_S}((O\mathbb B_{dr,\tilde S_I},F),t_{IJ}))
\end{eqnarray*}
\item[(iv)] Let $f:X\to S$ a morphism with $S,X\in\Var(k)$. 
Consider a factorization $f:X\hookrightarrow Y\times S\xrightarrow{p}S$ with $Y\in\SmVar(k)$. 
Let $S=\cup_i S_i$ an open affine cover and 
$i_i:S_i\hookrightarrow\tilde S_i$ closed embedding with $\tilde S_i\in\SmVar(k)$. Let 
\begin{equation*}
\alpha:\mathbb B_{dr,(\tilde S_I)}(K,W)\to 
F^0DR(S)(((M_I,F,W),u_{IJ})^{an}\otimes_{O_{\tilde S_I}}((O\mathbb B_{dr,\tilde S_I},F),t_{IJ})) 
\end{equation*}
a morphism in $D_{\mathbb B_{dr},G,fil}(S_{\mathbb C_p}^{an,pet}/(\tilde S_{I\mathbb C_p}^{an,pet}))$, with
\begin{equation*}
((M_I,F,W),u_{IJ})\in C(DRM(S))\subset C_{\mathcal D(1,0)fil,rh}(S/(\tilde S_I)), \; 
(K,W)\in D_{\mathbb Z_pfil,c,k}(S^{et})^{ad,(\Gamma_{f,i})}. 
\end{equation*}
We then have by (iii), 
the maps in $D_{\mathbb B_{dr},G,fil}(X_{\mathbb C_p}^{an,pet}/(Y\times\tilde S_{I,\mathbb C_p}^{an,pet}))$
\begin{eqnarray*}
f^!\alpha=f^!(\alpha):\mathbb B_{dr,(Y\times\tilde S_I)}(f^{!w}(K,W))
\xrightarrow{=}\mathbb B_{dr,(Y\times\tilde S_I)}(\Gamma^w_Xp^*(K,W)) \\
\xrightarrow{\Gamma_X(p^*\alpha)}
F^0DR(X)((\Gamma^{Hdg}_Xp^{*mod}((M_I,F,W),u_{IJ}))^{an}\otimes_{O_X}((O\mathbb B_{dr,Y\times\tilde S_I},F),t_{IJ})) \\
\xrightarrow{=}DR(X)(f_{Hdg}^{*mod}((M_I,F,W),u_{IJ})^{an}\otimes_{O_X}((O\mathbb B_{dr,Y\times\tilde S_I},F),t_{IJ}))
\end{eqnarray*}
and
\begin{eqnarray*}
f^*\alpha=f^*(\alpha):\mathbb B_{dr,(Y\times\tilde S_I)}(f^{*w}(K,W))
\xrightarrow{=}\mathbb B_{dr,(Y\times\tilde S_I)}(\Gamma^{\vee,w}_Xp^*(K,W)) \\
\xrightarrow{\Gamma^{\vee}_X(p^*\alpha)}
F^0DR(Y\times S)((\Gamma^{\vee,Hdg}_Xp^{\hat*mod}((M_I,F,W),u_{IJ}))^{an}
\otimes_{O_X}((O\mathbb B_{dr,Y\times\tilde S_I},F),t_{IJ})) \\
\xrightarrow{=}F^0DR(X)(f_{Hdg}^{\hat*mod}((M_I,F,W),u_{IJ})^{an}
\otimes_{O_X}((O\mathbb B_{dr,Y\times\tilde S_I},F),t_{IJ}))
\end{eqnarray*}
with
\begin{eqnarray*}
p^*\alpha:\mathbb B_{dr,(Y\times\tilde S_I)}(p^*(K,W))
\xrightarrow{=}p^{*mod}\mathbb B_{dr,(\tilde S_I)}(K,W) \\
\xrightarrow{p^{*mod}\alpha} 
p^{*mod}F^0DR(S)(((M_I,F,W),u_{IJ})^{an}\otimes_{O_S}((O\mathbb B_{dr,\tilde S_I},F),t_{IJ})) \\
\xrightarrow{=}
F^0DR(Y\times S)(p^{*mod}((M_I,F,W),u_{IJ})^{an}\otimes_{O_{Y\times S}}((O\mathbb B_{dr,Y\times\tilde S_I},F),t_{IJ})).
\end{eqnarray*}
\item[(v)]Let $S\in\Var(k)$. Denote by $\Delta_S:S\hookrightarrow S\times S$ the diagonal closed embedding
and $p_1:S\times S\to S$, $p_2:S\times S\to S$ the projections.
Let $S=\cup_i S_i$ an open affine cover and 
$i_i:S_i\hookrightarrow\tilde S_i$ closed embedding with $\tilde S_i\in\SmVar(k)$. Let 
\begin{eqnarray*}
\alpha:\mathbb B_{dr,(\tilde S_I)}(K,W)\to 
F^0DR(S)(((M_I,F,W),u_{IJ})^{an}\otimes_{O_S}((O\mathbb B_{dr,\tilde S_I},F),t_{IJ})), \\ 
\alpha':\mathbb B_{dr,(\tilde S_I)}(K',W)\to 
F^0DR(S)(((M'_I,F,W),v_{IJ})^{an}\otimes_{O_S}((O\mathbb B_{dr,\tilde S_I},F),t_{IJ})) 
\end{eqnarray*}
two morphisms in $D_{\mathbb B_{dr},G,fil}(S_{\mathbb C_p}^{an,pet}/(\tilde S_{I\mathbb C_p}^{an,pet}))$, with
\begin{eqnarray*}
((M_I,F,W),u_{IJ}),((M'_I,F,W),v_{IJ})\in C(DRM(S))\subset C_{\mathcal D(1,0)fil,rh}(S/(\tilde S_I)), \\ 
(K,W),(K',W)\in D_{\mathbb Z_pfil,c,k}(S^{et}). 
\end{eqnarray*}
We have then, as in (iv), the following map 
in $D_{\mathbb B_{dr},G,fil}(S_{\mathbb C_p}^{an,pet}/(\tilde S_{I\mathbb C_p}^{an,pet}))$
\begin{eqnarray*}
\alpha\otimes\alpha':\mathbb B_{dr,(\tilde S_I)}((K,W)\otimes^{L,w}(K',W)) \\ 
\xrightarrow{T(\otimes,\mathbb B_{dr})((K,W),(K',W))^{-1}}
\mathbb B_{dr,(\tilde S_I)}((K,W)\otimes_{\mathbb B_{dr,S}}\mathbb B_{dr,(\tilde S_I)}(K',W)) \\
\xrightarrow{\alpha\otimes\alpha'} 
F^0DR(S)(((M_I,F,W),u_{IJ})^{an}\otimes_{O_S}((O\mathbb B_{dr,\tilde S_I},F),t_{IJ})) 
\otimes_{\mathbb B_{dr,S}} \\
F^0DR(S)(((M'_I,F,W),v_{IJ})^{an}\otimes_{O_S}((O\mathbb B_{dr,\tilde S_I},F),t_{IJ})) \\
\xrightarrow{w_S}
V_{S_10}\cdots V_{S_r0}\Gamma_S^wF^0DR(S\times S)(p_1^{*mod}((M_I,F,W),u_{IJ})^{an}
\otimes_{O_{S\times S}}p_2^{*mod}((M'_I,F,W),v_{IJ})^{an} \\
\otimes_{O_{S\times S}}((O\mathbb B_{dr,\tilde S_I\times\tilde S_I},F),t_{IJ})) \\
\xrightarrow{(I,m(-)\otimes w_S}
F^0DR(S)(((M_I,F,W),u_{IJ})^{an}\otimes_{O_S}^{Hdg}((M'_I,F,W),v_{IJ})^{an}
\otimes_{O_S}((O\mathbb B_{dr,\tilde S_I},F),t_{IJ})) 
\end{eqnarray*}
where $S=\cap_i S_i$ with $S_i\subset S$ Cartier divisor, and (see (ii)) for $j_i:S\backslash S_i\hookrightarrow S$
the open embedding 
$m(M):V_{S_i0}j_{iw}(M,F,W)\otimes_{O_S}j_{i*Hdg}(O_{S\backslash S_i},F_b)\to j_{i*Hdg}(M,F,W)$
is the multiplication map.
\end{itemize}
\end{defi}

\begin{lem}\label{falphaplem}
\begin{itemize}
\item[(i)]Let $j:S^o\hookrightarrow S$ an open embedding with $S\in\SmVar(k)$ such that
$D=S\backslash S^o=V(s)\subset S$ is a (Cartier) divisor. For $(M,F,W)\in C(DRM(S^o))$,
\begin{equation*}
m(M):V_{D0}j_{*w}(M,F,W)\otimes_{O_S}j_{*Hdg}(O_{S^o},F_b)\to j_{*Hdg}(M,F,W), m\otimes h\mapsto m(M)(m\otimes h):=hm
\end{equation*}
is an isomorphism in $C(DRM(S))$, whose inverse is given by
\begin{equation*}
n(M):j_{*Hdg}(M,F,W)\to V_{D0}j_{*w}(M,F,W)\otimes_{O_S}j_{*Hdg}(O_{S^o},F_b), m\mapsto n(M)(m):=s^rm\otimes 1/s^r.
\end{equation*}
where $r\in\mathbb N$ is such that $s^rm\in\Gamma(W,V_{D0}j_*M)$ for $m\in\Gamma(W,M)$ and $W\subset S$ an open subset.
\item[(i)']Let $l:S^o\hookrightarrow S$ an open embedding with $S\in\Var(k)$ such that
$D=S\backslash S^o=V(s)\subset S$ is a Cartier divisor.
Let $S=\cup_i S_i$ an open affine cover and 
$i_i:S_i\hookrightarrow\tilde S_i$ closed embedding with $\tilde S_i\in\SmVar(k)$. 
Let $l_I:\tilde S^o_I\hookrightarrow\tilde S_I$ closed embeddings such that $\tilde S^o_I\cap S=S^o\cap S_I$.
For $((M_I,F,W),u_{IJ})\in C(DRM(S^o))$
\begin{eqnarray*}
(m(M_I)):V_{D0}l_{*w}((M_I,F,W),u_{IJ})\otimes_{O_S}(l_{I*Hdg}(O_{\tilde S_I^o},F_b),x_{IJ})\to (l_{*Hdg}((M_I,F,W),u_{IJ})
\end{eqnarray*}
is an isomorphism in $C(DRM(S))$ whose inverse is given by
\begin{eqnarray*}
(n(M_I)):l_{*Hdg}((M_I,F,W),u_{IJ})\to V_{D0}l_{*w}((M_I,F,W),u_{IJ})\otimes_{O_S}(l_{I*Hdg}(O_{\tilde S_I^o},F_b),x_{IJ}).
\end{eqnarray*}
\item[(ii)]Let $j:S^o\hookrightarrow S$ an open embedding with $S\in\SmVar(k)$ 
such that $D=S\backslash S^o=V(s)\subset S$ is a (Cartier) divisor. For $(M,F,W)\in C(DRM(S^o))$,
\begin{eqnarray*}
w_S\otimes m(M)^{an}
V_{D0}j_{*w}F^0DR(S^o)(((M,F,W))^{an}\otimes_{O_{S^o}}(O\mathbb B_{dr,S^o},F))\otimes_{\mathbb B_{dr,S}} \\
F^0DR(S)(j_{*Hdg}(O_{S^o},F_b)^{an}\otimes_{O_S}(O\mathbb B_{dr,S},F)) \\
\to F^0DR(S)(j_{*Hdg}(M,F,W)^{an}\otimes_{O_S}(O\mathbb B_{dr,S},F)), \; 
(w_1\otimes m)\otimes(w_2\otimes h)\mapsto (w_1\wedge w_2)\otimes(hm)
\end{eqnarray*}
is an isomorphism in $C(DRM(S))$ whose inverse is 
\begin{eqnarray*}
w_S^{-1}\otimes n(M)^{an}:F^0DR(S)(j_{*Hdg}(M,F,W)^{an}\otimes_{O_S}(O\mathbb B_{dr,S},F)) \\
\to(V_{D0}j_{*w}F^0DR(S^o)(((M,F,W))^{an}\otimes_{O_S}(O\mathbb B_{dr,S},F)))\otimes_{\mathbb B_{dr,S}} \\
(F^0DR(S)(j_{*Hdg}(O_{S^o},F_b)^{an}\otimes_{O_S}(O\mathbb B_{dr,S},F))).
\end{eqnarray*}
\item[(ii)']Let $l:S^o\hookrightarrow S$ an open embedding with $S\in\Var(k)$
such that $D=S\backslash S^o=V(s)\subset S$ is a Cartier divisor. 
Let $S=\cup_i S_i$ an open affine cover and 
$i_i:S_i\hookrightarrow\tilde S_i$ closed embedding with $\tilde S_i\in\SmVar(k)$. 
Let $l_I:\tilde S^o_I\hookrightarrow\tilde S_I$ closed embeddings such that $\tilde S^o_I\cap S=S^o\cap S_I$.
For $((M_I,F,W),u_{IJ})\in C(DRM(S^o))$
\begin{eqnarray*}
(w_{\tilde S_I}\otimes m(M_I)^{an}):
V_{D0}l_{*w}F^0DR(S^o)(((M_I,F,W),u_{IJ})^{an}\otimes_{O_{S^o}}((O\mathbb B_{dr,\tilde S^o_I},F),t_{IJ}))
\otimes_{\mathbb B_{dr,(\tilde S_I)}} \\
(F^0DR((\tilde S_I))((l_{I*Hdg}(O_{\tilde S_I^o},F_b),x_{IJ})^{an}\otimes_{O_S}((O\mathbb B_{dr,\tilde S_I},F),t_{IJ}))) \\
\to F^0DR(S)((l_{*Hdg}((M_I,F,W),u_{IJ})^{an})\otimes_{O_S}((O\mathbb B_{dr,\tilde S_I},F),t_{IJ}))
\end{eqnarray*}
is an isomorphism in $C(DRM(S))$ whose inverse is 
\begin{eqnarray*}
(w_{\tilde S_I}^{-1}\otimes n(M_I)^{an}):
F^0DR(S)((l_{*Hdg}((M_I,F,W),u_{IJ})^{an})\otimes_{O_S}((O\mathbb B_{dr,\tilde S_I},F),t_{IJ})) \\
\to V_{D0}l_{*w}F^0DR(S^o)(((M_I,F,W),u_{IJ})^{an}\otimes_{O_{S^o}}((O\mathbb B_{dr,\tilde S^o_I},F),t_{IJ}))
\otimes_{\mathbb B_{dr,(\tilde S_I)}} \\
F^0DR((\tilde S_I))((l_{I*Hdg}(O_{\tilde S_I^o},F_b),x_{IJ})^{an}\otimes_{O_S}((O\mathbb B_{dr,\tilde S_I},F),t_{IJ})).
\end{eqnarray*}
\end{itemize}
\end{lem}

\begin{proof}
\noindent(i): Follows from the definition of the $V$-filtration and the $F$-filtration.

\noindent(i)':Follows from (i).

\noindent(ii): Follows from (i).

\noindent(ii)':Follows from (ii).
\end{proof}

\begin{prop}\label{falphapprop}
Let $k\subset{\mathbb C_p}$ a subfield.
\begin{itemize}
\item[(i0)]Let $j:S^o\hookrightarrow S$ an open embedding with $S\in\SmVar(k)$ 
and $D=S\backslash S^o$ a (Cartier) divisor. If 
\begin{equation*}
\alpha:\mathbb B_{dr,S^o}(K,W)\to F^0DR(S^o)((M,F,W)^{an}\otimes_{O_{S^o}}(O\mathbb B_{dr,S^o},F)) 
\end{equation*}
is an isomorphism in $D_{\mathbb B_{dr},G,fil}(S_{\mathbb C_p}^{0,an,pet})$, with
\begin{equation*}
(M,F,W)\in C(DRM(S^o)), \; (K,W)\in D_{\mathbb Z_pfil,c,k}(S^{o,et}), 
\end{equation*}
then the maps in $D_{\mathbb B_{dr},G,fil}(S_{\mathbb C_p}^{an,pet})$
\begin{eqnarray*}
j_*\alpha=j_*(\alpha):\mathbb B_{dr,S}(j_{*w}(K,W))\to F^0DR(S)((j_{*Hdg}(M,F,W))^{an}\otimes_{O_S}(O\mathbb B_{dr,S},F))
\end{eqnarray*}
and
\begin{eqnarray*}
j_!\alpha=j_!(\alpha):\mathbb B_{dr,S}(j_{!w}(K,W))\to F^0DR(S)((j_{!Hdg}(M,F,W))^{an}\otimes_{O_S}(O\mathbb B_{dr,S},F))
\end{eqnarray*}
given in definition \ref{falpha} are isomorphism.
\item[(i0)']Let $l:S^o\hookrightarrow S$ an open embedding with $S\in\Var(k)$ and $D=S\backslash S^o$ a Cartier divisor.
Let $S=\cup_i S_i$ an open affine cover and 
$i_i:S_i\hookrightarrow\tilde S_i$ closed embedding with $\tilde S_i\in\SmVar(k)$. 
Let $l_I:\tilde S^o_I\hookrightarrow\tilde S_I$ closed embeddings such that $\tilde S^o_I\cap S=S^o\cap S_I$. If 
\begin{equation*}
\alpha:\mathbb B_{dr,(\tilde S^o_I)}(K,W)\to
F^0DR(S^o)(((M_I,F,W),u_{IJ})^{an}\otimes_{O_{\tilde S^o_I}}(O\mathbb B_{dr,\tilde S^o_I})) 
\end{equation*}
is an isomorphism in $D_{\mathbb B_{dr},G,fil}(S_{\mathbb C_p}^{an,pet}/(\tilde S_{I,\mathbb C_p}^{an,pet}))$, with
\begin{equation*}
((M_I,F,W),u_{IJ})\in C(DRM(S^o)), \; (K,W)\in D_{\mathbb Z_pfil,c,k}(S^{et}), 
\end{equation*}
then the maps in $D_{\mathbb B_{dr},G,fil}(S_{\mathbb C_p}^{an,pet}/(\tilde S_{I,\mathbb C_p}^{an,pet}))$
\begin{eqnarray*}
l_*\alpha=l_*(\alpha):\mathbb B_{dr,(\tilde S_I)}(l_{*w}(K,W))\to
F^0DR(S)(l_{!Hdg}((M_I,F,W),u_{IJ})^{an}\otimes_{O_S}(O\mathbb B_{dr,(\tilde S_I)},F))
\end{eqnarray*}
and
\begin{eqnarray*}
l_!\alpha=l_!(\alpha):\mathbb B_{dr,(\tilde S_I)}(l_{!w}(K,W))\to
F^0DR(S)(l_{!Hdg}((M_I,F,W),u_{IJ})^{an}\otimes_{O_S}(O\mathbb B_{dr,(\tilde S_I)},F))
\end{eqnarray*}
are isomorphisms.
\item[(i)]Let $f:X\to S$ a morphism with $S,X\in\SmVar(k)$. If 
\begin{equation*}
\alpha:\mathbb B_{dr,X}(K,W)\to F^0DR(X)((M,F,W)^{an}\otimes_{O_X}(O\mathbb B_{dr,X},F)) 
\end{equation*}
an isomorphism in $D_{\mathbb B_{dr},G,fil}(X_{\mathbb C_p}^{an,pet})$, with
\begin{equation*}
(M,F,W)\in C(DRM(X)), \; (K,W)\in D_{\mathbb Z_pfil,c,k}(X^{et}), 
\end{equation*}
then the morphisms given in definition \ref{falphap}
\begin{eqnarray*}
f_*\alpha=f_*(\alpha):\mathbb B_{dr,S}(Rf_{*w}(K,W))\to F^0DR(S)((Rf_{*Hdg}(M,F,W))^{an}\otimes_{O_S}(O\mathbb B_{dr,S},F))
\end{eqnarray*}
and
\begin{eqnarray*}
f_!\alpha=f_!(\alpha):\mathbb B_{dr,S}(Rf_{!w}(K,W))\to F^0DR(S)((Rf_{!Hdg}(M,F,W))^{an}\otimes_{O_S}(O\mathbb B_{dr,S},F))
\end{eqnarray*}
are isomorphisms.
\item[(i)']Let $f:X\to S$ a morphism with $S,X\in\QPVar(k)$. 
Consider a factorization $f:X\hookrightarrow Y\times S\xrightarrow{p}S$ with $Y\in\SmVar(k)$. 
Let $S=\cup_i S_i$ an open affine cover and 
$i_i:S_i\hookrightarrow\tilde S_i$ closed embedding with $\tilde S_i\in\SmVar(k)$. If 
\begin{equation*}
\alpha:\mathbb B_{dr,(Y\times\tilde S_I)}(K,W)\to 
F^0DR(X)(((M_I,F,W),u_{IJ})^{an}\otimes_{O_{Y\times\tilde S_I}}(O\mathbb B_{dr,Y\times\tilde S_I},F)) 
\end{equation*}
is an isomorphism in $D_{\mathbb B_{dr},G,fil}(X_{\mathbb C_p}^{an,pet}/({Y\times\tilde S_I}_{\mathbb C_p}^{an,pet}))$, with
\begin{equation*}
((M_I,F,W),u_{IJ})\in C(DRM(X)), \; (K,W)\in D_{\mathbb Z_pfil,c,k}(X^{et}), 
\end{equation*}
then, the maps in $D_{\mathbb B_{dr},G,fil}(S_{\mathbb C_p}^{an,pet}/(\tilde S_{I,\mathbb C_p}^{an,pet}))$
\begin{eqnarray*}
f_*\alpha=f_*(\alpha):\mathbb B_{dr,(\tilde S_I)}(Rf_{*w}(K,W))\to
F^0DR(S)((Rf_{*Hdg}((M_I,F,W),u_{IJ}))^{an}\otimes_{O_S}(O\mathbb B_{dr,(\tilde S_I)},F)),
\end{eqnarray*}
and
\begin{eqnarray*}
f_!\alpha=f_!(\alpha):\mathbb B_{dr,(\tilde S_I)}((Rf_{!w}(K,W)))\to
F^0DR(S)((Rf_{!Hdg}((M_I,F,W),u_{IJ}))^{an}\otimes_{O_S}(O\mathbb B_{dr,(\tilde S_I)},F)),
\end{eqnarray*}
are isomorphisms.
\item[(ii)]Let $f:X\to S$ a morphism with $S,X\in\QPVar(k)$. 
Consider a factorization $f:X\hookrightarrow Y\times S\xrightarrow{p}S$ with $Y\in\SmVar(k)$. 
Let $S=\cup_i S_i$ an open affine cover and 
$i_i:S_i\hookrightarrow\tilde S_i$ closed embedding with $\tilde S_i\in\SmVar(k)$. If 
\begin{equation*}
\alpha:\mathbb B_{dr,(\tilde S_I)}(K,W)\to 
F^0DR(S)(((M_I,F,W),u_{IJ})^{an}\otimes_{O_{\tilde S_I}}(O\mathbb B_{dr,(\tilde S_I)},F)) 
\end{equation*}
is an isomorphism in $D_{\mathbb B_{dr},G,fil}(S_{\mathbb C_p}^{an,pet}/(\tilde S_{I,\mathbb C_p}^{an,pet}))$, with
\begin{equation*}
((M_I,F,W),u_{IJ})\in C(DRM(S)), \; (K,W)\in D_{\mathbb Z_pfil,c,k}(S^{et}), 
\end{equation*}
the maps in $D_{\mathbb B_{dr},G,fil}(X_{\mathbb C_p}^{an,pet}/(Y\times\tilde S_I)_{\mathbb C_p}^{an,pet})$
\begin{eqnarray*}
f^!\alpha=f^!(\alpha):\mathbb B_{dr,(Y\times\tilde S_I)}(f^{!w}(K,W))\to
F^0DR(X)(f_{Hdg}^{*mod}((M_I,F,W),u_{IJ})^{an}\otimes_{O_X}(O\mathbb B_{dr,(Y\times\tilde S_I)},F))
\end{eqnarray*}
and
\begin{eqnarray*}
f^*\alpha=f^*(\alpha):\mathbb B_{dr,(Y\times\tilde S_I)}(f^!(K,W))\to
F^0DR(X)(f_{Hdg}^{\hat*mod}((M_I,F,W),u_{IJ})^{an}\otimes_{O_X}(O\mathbb B_{dr,(Y\times\tilde S_I)},F))
\end{eqnarray*}
given in definition \ref{falphap} are isomorphisms.
\end{itemize}
\end{prop}

\begin{proof}
\noindent(i0): Follows from lemma \ref{falphaplem}(ii).

\noindent(i0)':Follows from lemma \ref{falphaplem}(ii)'.

\noindent(i): Follows from (i0) and on the other hand 
theorem \ref{TfBdrthm} and GAGA for proper morphism of algebraic varieties over a $p$-adic field.

\noindent(i): Follows from (i0)' and on the other hand 
theorem \ref{TfBdrthm} and GAGA for proper morphism of algebraic varieties over a $p$-adic field.

\noindent(ii): Follows from (i0).
\end{proof}

\begin{defi}\label{DHdgpsialphap}
Let $S\in\SmVar(k)$.
Let $D=V(s)\subset S$ a divisor with $s\in\Gamma(S,L)$ and $L$ a line bundle ($S$ being smooth, $D$ is Cartier).
For $\mathcal M=((M,F,W),(K,W),\alpha)\in\PSh_{\mathcal D(1,0)fil,rh}(S)\times_IP_{\mathbb Z_pfil,k}(S^{et})$, 
we then define, using definition \ref{DHdgpsi}, theorem \ref{phipsithmp} and definition {TphipsiBdr},
\begin{itemize}
\item the nearby cycle functor
\begin{equation*}
\psi_D((M,F,W),(K,W),\alpha):=(\psi_D(M,F,W),\psi_D(K,W)[-1],\psi_D\alpha)
\in\PSh_{\mathcal D(1,0)fil,rh}(S)\times_IP_{\mathbb Z_pfil,k}(S^{et}),
\end{equation*}
with
\begin{eqnarray*}
\psi_D\alpha:\mathbb B_{dr,S}(\psi_D(K,W))\xrightarrow{T(\psi_D,\mathbb B_{dr})(K,W)}\psi_D\mathbb B_{dr,S}(K,W)  
\xrightarrow{\psi_D\alpha} \\
\psi_DDR(S)((M,F,W)^{an}\otimes_{O_S}(O\mathbb B_{dr,S},F)) 
\xrightarrow{T^{B_{dr}}(\psi_D,DR)(M,F,W)} \\
DR(S)(\psi_D(M,F,W)^{an}\otimes_{O_S}(O\mathbb B_{dr,S},F)),
\end{eqnarray*}
\item the vanishing cycle functor
\begin{equation*}
\phi_D((M,F,W),(K,W),\alpha):=(\phi_D(M,F,W),\phi_D(K,W)[-1],\phi_D\alpha)
\in\PSh_{\mathcal D(1,0)fil,rh}(S)\times_IP_{\mathbb Z_pfil,k}(S^{et}),
\end{equation*}
with
\begin{eqnarray*}
\phi_D\alpha:\mathbb B_{dr,S}(\phi_D(K,W))\xrightarrow{T(\phi_D,\mathbb B_{dr})(K,W)}\phi_D\mathbb B_{dr,S}(K,W) 
\xrightarrow{\phi_D\alpha} \\
\phi_DDR(S)((M,F,W)^{an}\otimes_{O_S}(O\mathbb B_{dr,S},F)) 
\xrightarrow{T^{B_{dr}}(\phi_D,DR)(M,F,W)} \\
DR(S)(\phi_D(M,F,W)^{an}\otimes_{O_S}(O\mathbb B_{dr,S},F)),
\end{eqnarray*}
\item the canonical maps in $\PSh_{\mathcal D(1,0)fil,rh}(S)\times_IP_{\mathbb Z_pfil,k}(S^{et})$
\begin{eqnarray*}
can(\mathcal M):=(can(M,F,W),can(K,W)):
\psi_D((M,F,W),(K,W),\alpha)\to\phi_D((M,F,W),(K,W),\alpha)(-1),  \\
var(\mathcal M):=(var(M,F,W),var(K,W)): 
\phi_D((M,F,W),(K,W),\alpha)\to\psi_D((M,F,W),(K,W),\alpha).
\end{eqnarray*}
\end{itemize}
\end{defi}

\begin{prop}\label{phipsigmHdgpropalphap}
Let $S\in\SmVar(k)$. Let $D=V(s)\subset S$ a (Cartier divisor). Consider a composition of proper morphisms 
\begin{equation*}
(f:X=X_r\xrightarrow{f_r}X_{r-1}\xrightarrow{f_1}X_0=S)\in\SmVar(k),\; \mbox{proper}, 1\leq i\leq r,
\end{equation*}
and 
\begin{eqnarray*}
(M,F)=H^{n_0}\int_{f_1}\cdots H^{n_r}\int_{f_r}((O_X,F_b), 
H^{n_0}Rf_{1*}\cdots H^{n_r}Rf_{r*}\mathbb Z_{X_{\bar k}}, \\
H^{n_0}f_{1*}\circ\cdots\circ H^{n_r}f_{r*}\alpha(X_{\mathbb C_p})) 
\in\PSh_{\mathcal Dfil,rh}(S)\times_IP_{\mathbb Z_p,k}(S^{et}).
\end{eqnarray*} 
Then,
\begin{eqnarray*}
\psi_D(M,F)=H^{n_0}\int_{f_1}\cdots H^{n_r}\int_{f_r}(\psi_{f^{-1}(D)}(O_X,F_b)), 
H^{n_0}Rf_{1*}\cdots H^{n_r}Rf_{r*}\psi_{f^{-1}(D)}\mathbb Z_{X_{\bar k}}, \\
H^{n_0}f_{1*}\circ\cdots\circ H^{n_r}f_{r*}\psi_{f^{-1}(D)}\alpha(X_{\mathbb C_p})) 
\in\PSh_{\mathcal Dfil,rh}(S)\times_IP_{\mathbb Z_p,k}(S^{et}), 
\end{eqnarray*}
and
\begin{eqnarray*}
\phi_D(M,F)=H^{n_0}\int_{f_1}\cdots H^{n_r}\int_{f_r}(\psi_{f^{-1}(D)}(O_X,F_b)), 
H^{n_0}Rf_{1*}\cdots H^{n_r}Rf_{r*}\phi_{f^{-1}(D)}\mathbb Z_{X_{\bar k}}, \\
H^{n_0}f_{1*}\circ\cdots\circ H^{n_r}f_{r*}\phi_{f^{-1}(D)}\alpha(X_{\mathbb C_p})) 
\in\PSh_{\mathcal Dfil,rh}(S)\times_IP_{\mathbb Z_p,k}(S^{et}), 
\end{eqnarray*}
\end{prop}

\begin{proof}
Immediate from definition.
\end{proof}

Let $S\in\Var(k)$. Let $S=\cup_{i\in I}S_i$ an open cover such that there
exists closed embeddings $i_i:S_i\hookrightarrow\tilde S_i$ with $\tilde S_I\in\SmVar(k)$. We consider
$\mathbb Z_{p,S^{et}}^w\in C_{\mathbb Z_pfil}(S^{et})$ such that 
$j_I^*\mathbb Z_{p,S^{et}}^w=i_I^*\Gamma^{\vee,w}_{S_I}\mathbb Z_{p,\tilde S_I}^{et}$ and set
\begin{eqnarray*}
\alpha(S):\mathbb B_{dr,(\tilde S_I)}(\mathbb Z_{p,S^{et}}^w)\xrightarrow{:=} \\
(\Gamma^{\vee,w}_{S_I}\mathbb Z_{p,(\tilde S_I)}\otimes\mathbb B_{dr,\phi_{\tilde D_{1,I}}}
\otimes_{\mathbb B_{dr,\tilde S_I}}\cdots\otimes_{\mathbb B_{dr,\tilde S_I}}\mathbb B_{dr,\phi_{\tilde D_{d,I}}}, 
\mathbb B_{dr}(t_{IJ})) \\
\xrightarrow{=}
(\Gamma^{\vee,w}_{S_I}\mathbb Z_{p,\tilde S_I}\otimes\mathbb B_{dr,\tilde S_I}
\otimes_{\mathbb B_{dr,\tilde S_I}}\mathbb B_{dr,\phi_{\tilde D_{1,I}}}
\otimes_{\mathbb B_{dr,\tilde S_I}}\cdots\otimes_{\mathbb B_{dr,\tilde S_I}}\mathbb B_{dr,\phi_{\tilde D_{d,I}}}, 
\mathbb B_{dr}(x_{IJ})) \\
\xrightarrow{\alpha((\tilde S_{I,\mathbb C_p})\otimes I)} \\
(V_{\tilde D_{1,I}0}\cdots V_{\tilde D_{d,I}0}\Gamma^{\vee}_{S_I}(F^0DR(\tilde S_I)((O\mathbb B_{dr,\tilde S_I}),F)) 
\otimes_{\mathbb B_{dr,\tilde S_I}}\mathbb B_{dr,\phi_{\tilde D_{1,I}}}
\otimes_{\mathbb B_{dr,\tilde S_I}}\cdots\otimes_{\mathbb B_{dr,\tilde S_I}}\mathbb B_{dr,\phi_{\tilde D_{d,I}}}, 
\mathbb B_{dr}(t_{IJ})) \\
\xrightarrow{(\mathbb DT(\gamma_{S_I},\otimes)(-))
\otimes(DR(S)(\mathbb D\rho_{DR,\tilde D_{1,I}}(O_{\tilde S_I},F_b))\otimes\cdots
\otimes\mathbb D(\rho_{DR,\tilde D_{d,I}}(O_{\tilde S_I},F_b)))} \\
F^0DR(S)(V_{\tilde D_{1,I}0}\cdots V_{\tilde D_{d,I}0}\Gamma^{\vee,h}_{S_I}(O_{\tilde S_I},F_b)
\otimes_{\mathbb B_{dr,\tilde S_I}}O\mathbb B_{dr,\tilde S_I},x_{IJ})\otimes_{\mathbb B_{dr,(\tilde S_I)}} \\
F^0DR(S)((\Gamma^{\vee,Hdg}_{S_I}(O_{\tilde S_I},F_b))^{an}\otimes_{O_{\tilde S_I}}(O\mathbb B_{dr,\tilde S_I},F),x_{IJ}) \\
\xrightarrow{T(DR,\otimes)(-,-)} 
F^0DR(S)(\Gamma^{\vee,h}_{S_I}(O_{\tilde S_I},F_b)\otimes_{O_{\tilde S_I}} 
(\Gamma^{\vee,Hdg}_{S_I}(O_{\tilde S_I},F_b))^{an}\otimes_{O_{\tilde S_I}}(O\mathbb B_{dr,\tilde S_I},F),x_{IJ}) \\
\xrightarrow{F^0DR(S)(m(O_{\tilde S_I}))} 
F^0DR(S)((\Gamma^{\vee,Hdg}_{S_I}(O_{\tilde S_I},F_b))^{an}\otimes_{O_{\tilde S_I}}(O\mathbb B_{dr,\tilde S_I},F),x_{IJ}) \\
\xrightarrow{=:}
F^0DR(S)((\Gamma^{\vee,Hdg}_{S_I}(O_{\tilde S_I},F_b),x_{IJ})^{an}\otimes_{O_S}((O\mathbb B_{dr,(\tilde S_I)},F),t_{IJ})) 
\end{eqnarray*}
is an isomorphism in $D_{\mathbb B_{dr},G}(S_{\mathbb C_p}^{an,pet}/(\tilde S_I)_{\mathbb C_p}^{an,pet})$, 
where $D_1,\ldots,D_d\subset S$ are Cartier divisors such that $S=\cap_{s=1}^dD_s$, 
$\tilde D_{s,I}\subset\tilde S_I$ are Cartier divisor such that $D_s\cap S_I\subset\tilde D_{s,I}\cap S_I$,
$i_I:S_I\hookrightarrow\tilde S_I$ are the closed embeddings, 
$m(O):O\times_O O\xrightarrow{\sim} O, h\otimes f\mapsto hf$ is the multiplication map,
and we use definition \ref{falphap} and proposition \ref{falphapprop}. 

We now give the definition of $p$ adic mixed Hodge modules which is the main definition of this section :

\begin{defi}\label{MHMgmkpdef}
Let $k\subset\mathbb C_p$ a subfield.
\begin{itemize}
\item[(i)]Let $S\in\SmVar(k)$. We denote by 
\begin{eqnarray*}
HM_{gm,k,\mathbb C_p}(S):= \\
<(H^{n_1}\int_{f_1}\cdots H^{n_r}\int_{f_r}(O_X,F_b)(d),
R^{n_1}f_{1*}\cdots R^{n_r}f_{r*}\mathbb Z_{p,X^{et}},H^{n_1}f_{1*}\cdots H^{n_r}f_{r*}\alpha(X)), \\
(f:X=X_r\xrightarrow{f_r}X_{r-1}\to\cdots\xrightarrow{f_1}X_0=S)\in\SmVar(k), \; \mbox{proper}, \; 
n_1,\ldots,n_r,d\in\mathbb Z> \\
\subset PDRM(S)\times_I P_{\mathbb Z_p,k}(S^{et})
\subset\PSh_{\mathcal Dfil,rh}(S)\times_I P_{\mathbb Z_p,k}(S^{et})
\end{eqnarray*}
the full abelian subcategory, where $<,>$ means generated by and $(-)$ the shift of the filtration, 
\begin{equation*}
\alpha(X):\mathbb B_{dr,X}(\mathbb Z_{p,X^{et}}):=\mathbb B_{dr,X_{\mathbb C_p}}
\hookrightarrow DR(X)(O\mathbb B_{dr,X_{\mathbb C_p}}) 
\end{equation*}
is the inclusion quasi-isomorphism in $C_{\mathbb B_{dr},G}(X_{\mathbb C_p}^{pet})$, and we use definition \ref{falphap}. 
We have by proposition \ref{phipsigmHdgpropalphap} for $((M,F),K,\alpha)\in HM_{gm,k,\mathbb C_p}(S)$,
\begin{eqnarray*}
\Gr_k^W\psi_D((M,F),K,\alpha):=\Gr_k^W\psi_D(M,F),\Gr_k^W\psi_DK,\Gr_k^W\psi_D\alpha)\in HM_{gm,k,\mathbb C_p}(S) .
\end{eqnarray*}
and
\begin{eqnarray*}
\Gr_k^W\psi_D((M,F),K,\alpha):=\Gr_k^W\psi_D(M,F),\Gr_k^W\psi_DK,\Gr_k^W\psi_D\alpha)\in HM_{gm,k,\mathbb C_p}(S).
\end{eqnarray*}
for all $k\in\mathbb Z$. We set 
\begin{equation*}
\mathbb Z_{p,S}^{Hdg}:=((O_S,F_b),\mathbb Z_{p,S^{et}},\alpha(S))\in HM_{gm,k,\mathbb C_p}(S)
\end{equation*}
\item[(i)']Let $S\in\Var(k)$. Let $S=\cup_{i\in I}S_i$ an open cover such that there
exists closed embeddings $i_i:S_i\hookrightarrow\tilde S_i$ with $\tilde S_I\in\SmVar(k)$. 
\begin{eqnarray*}
HM_{gm,k,\mathbb C_p}(S):= 
<(R^{n_1}p_{1*Hdg}\cdots R^{n_r}p_{r*Hdg}(\Gamma^{Hdg}_{X_I}(O_{Y\times\tilde X_{r-1,I}},F_b),x_{IJ})(d), \\
R^{n_1}p_{1*}\cdots R^{n_r}p_{r*}T(X/(Y_r\times\tilde X_{r-1,I}))(\mathbb Z_{p,X^{et}}), 
H^{n_1}p_{1*}\cdots H^{n_r}p_{r*}\alpha(X)), \\
(f:X=X_r\xrightarrow{f_r}X_{r-1}\to\cdots\xrightarrow{f_1}X_0=S)\in\Var(k), \; n_1,\ldots,n_r,d\in\mathbb Z> \\
\subset PDRM(S)\times_I P_{\mathbb Z_p,k}(S^{et})\subset
\PSh_{\mathcal Dfil,rh}(S/(\tilde S_I))\times_I P_{\mathbb Z_p,k}(S^{et})
\end{eqnarray*}
the full abelian subcategory, where $<,>$ means generated by and $(-)$ the shift of the filtration, 
$f_i:X_i\hookrightarrow Y_i\times X_{i-1}\xrightarrow{p_i}X_{i-1}$ proper, $Y_i\in\PSmVar(k)$, $X_i$ smooth,
and $\alpha(X)$ is given above.
Note that if $S$ is smooth then this definition of $HM_{gm,k,\mathbb C_p}(S)$ agree with the one given in (i).
\item[(ii)]Let $S\in\Var(k)$.
Take an open cover $S=\cup_iS_i$ such that there are closed embedding $S_I\hookrightarrow\tilde S_I$ with $S_I\in\SmVar(k)$.
We define using the pure case (i) and (i)' the full subcategory of geometric mixed Hodge modules defined over $k$
\begin{eqnarray*}
MHM_{gm,k,\mathbb C_p}(S):= \\
\left\{(((M_I,F,W),u_{IJ}),(K,W),\alpha), \; \mbox{s.t.} \; 
\Gr^W_k(((M_I,F,W),u_{IJ}),(K,W),\alpha)\in HM_{gm,k,\mathbb C_p}(S)\right\} \\
\subset DRM(S)\times_IP_{\mathbb Z_pfil,k}(S^{et})
\subset\PSh_{\mathcal D(1,0)fil,rh}(S/(\tilde S_I))\times_IP_{\mathbb Z_pfil,k}(S^{et})
\end{eqnarray*}
whose object consists of $(((M_I,F,W),u_{IJ}),(K,W),\alpha)\in DRM(S)\times_IP_{\mathbb Z_pfil,k}(S^{et})$ 
such that 
\begin{equation*}
\Gr^W_k(((M_I,F,W),u_{IJ})(K,W),\alpha):=(\Gr^W_k((M_I,F),u_{IJ}),\Gr^W_kK,\Gr^W_k\alpha)
\in HM_{gm,k,\mathbb C_p}(S).
\end{equation*}
where $DRM(S)$ is the category of de Rham modules introduced in section 5 definition \ref{DRMdef}.
The fact that $\alpha$ is an isomorphism implies that the Galois representation of $G$ induced on each $k$-point of $S$
is a de Rham representation. We set 
\begin{equation*}
\mathbb Z_{p,S}^{Hdg}:=((\Gamma_{S_I}^{\vee,Hdg}(O_{\tilde S_I},F_b),x_{IJ}),\mathbb Z_{p,S^{et}}^w,\alpha(S))
\in C(MHM_{gm,k,\mathbb C_p}(S))
\end{equation*}
where $\mathbb Z_{p,S^{et}}^w\in C(P_{fil,k}(S^{et}))$ 
is such that $j_I^*\mathbb Z_{p,S^{et}}^w=i_I^*\Gamma^{\vee,w}_{S_I}\mathbb Z_{p,\tilde S_I^{et}}$
and $\alpha(S)$ given above.
For $S\in\SmVar(k)$ and $D=V(s)\subset S$ a (Cartier) divisor, 
we have for $((M,F,W),(K,W),\alpha)\in MHM_{gm,k,\mathbb C_p}(S)$, using theorem \ref{HSk}, 
\begin{equation*}
\psi_D((M,F,W),(K,W),\alpha),\phi_D((M,F,W),(K,W),\alpha)\in MHM_{gm,k,\mathbb C_p}(S),
\end{equation*}
by the pure case (c.f. (i) and proposition \ref{phipsigmHdgpropalphap}) and the strictness of the $V$-filtration.
\end{itemize}
For $S\in\Var(k)$ we get $D(MHM_{gm,k,\mathbb C_p}(S)):=\Ho_{(zar,et)}(C(MHM_{gm,k,\mathbb C_p}(S)))$ 
after localization with Zariski local equivalence and etale local equivalence.
\end{defi}

We now look at functorialities :

\begin{defi}\label{DHdgjalphap}
Let $k\subset\mathbb C_p$ a subfield. Let $S\in\SmVar(k)$. Let $j:S^o\hookrightarrow S$ an open embedding.
Let $Z:=S\backslash S^o=V(\mathcal I)\subset S$ an the closed complementary subset, 
$\mathcal I\subset O_S$ being an ideal subsheaf. 
Taking generators $\mathcal I=(s_1,\ldots,s_r)$, we get $Z=V(s_1,\ldots,s_r)=\cap^r_{i=1}Z_i\subset S$ with 
$Z_i=V(s_i)\subset S$, $s_i\in\Gamma(S,\mathcal L_i)$ and $L_i$ a line bundle. 
Note that $Z$ is an arbitrary closed subset, $d_Z\geq d_X-r$ needing not be a complete intersection. 
Denote by $j_I:S^{o,I}:=\cap_{i\in I}(S\backslash Z_i)=S\backslash(\cup_{i\in I}Z_i)\xrightarrow{j_I^o}S^o\xrightarrow{j} S$ 
the open embeddings.
Let $(M,F,W)\in MHM_{gm,k,\mathbb C_p}(S^o))$. We then define, using definition \ref{DHdgj} and definition \ref{jwet}
\begin{itemize}
\item the canonical extension 
\begin{eqnarray*}
j_{*Hdg}((M,F,W),(K,W),\alpha):=(j_{*Hdg}(M,F,W),j_{*w}(K,W),j_*\alpha) \\
:=\in MHM_{gm,k,\mathbb C_p}(S), 
\end{eqnarray*}
so that $j^*(j_{*Hdg}((M,F,W),(K,W),\alpha))=((M,F,W),(K,W),\alpha)$,
\item the canonical extension 
\begin{eqnarray*}
j_{!Hdg}((M,F,W),(K,W),\alpha):=(j_{!Hdg}(M,F,W),j_{!w}(K,W),j_!\alpha)
:=\in MHM_{gm,k,\mathbb C_p}(S),  
\end{eqnarray*}
so that $j^*(j_{!Hdg}((M,F,W),(K,W),\alpha))=((M,F,W),(K,W),\alpha)$.
\end{itemize}
Moreover for $((M',F,W),(K',W),\alpha')\in MHM_{gm,k,\mathbb C_p}(S)$,
\begin{itemize}
\item there is a canonical map in $MHM_{gm,k,\mathbb C_p}(S)$
\begin{equation*}
\ad(j^*,j_{*Hdg})((M',F,W),(K',W),\alpha'):((M',F,W),(K',W),\alpha')\to j_{*Hdg}j^*((M',F,W),(K',W),\alpha'), 
\end{equation*}
\item there is a canonical map in $MHM_{gm,k,\mathbb C_p}(S)$
\begin{equation*}
\ad(j_{!Hdg},j^*)((M',F,W),(K',W),\alpha'):j_{!Hdg}j^*((M',F,W),(K',W),\alpha')\to((M',F,W),(K',W),\alpha').
\end{equation*}
\end{itemize}
\end{defi}

For $(M,F,W)\in C(MHM_{gm,k,\mathbb C_p}(S^o))$, 
\begin{itemize}
\item we have the canonical map in $C_{\mathcal D(1,0)fil}(S)\times_IC_{fil}(S^{et})$
\begin{eqnarray*}
T(j_{*Hdg},j_*)((M,F,W),(K,W),\alpha):=(k\circ\ad(j^*,j_*)(-),k\circ\ad(j^*,j_*),0): \\
j_{*Hdg}((M,F,W),(K,W),\alpha)\to (j_*E(M,F,W),j_*E(K,W),\alpha)
\end{eqnarray*}
\item we have the canonical map in $C_{\mathcal D(1,0)fil}(S)\times_IC_{fil}(S^{et})$
\begin{eqnarray*}
T(j_!,j_{!Hdg})((M,F,W),(K,W),\alpha):=(k\circ\ad(j_!,j^*)(-),k\circ\ad(j_!,j^*)(-),0): \\
(j_!(M,F,W),j_!(K,W),j_!\alpha)\to j_{!Hdg}((M,F,W),(K,W),\alpha).
\end{eqnarray*}
\end{itemize}

\begin{prop}\label{jHdgpropadCp}
\begin{itemize}
\item[(i)] Let $S\in\SmVar(k)$.
Let $D=V(s)\subset S$ a divisor with $s\in\Gamma(S,L)$ and $L$ a line bundle ($S$ being smooth, $D$ is Cartier).
Denote by $j:S^o:=S\backslash D\hookrightarrow S$ the open complementary embedding. Then, 
\begin{itemize}
\item $(j^*,j_{*Hdg}):MHM_{gm,k,\mathbb C_p}(S)\leftrightarrows MHM_{gm,k,\mathbb C_p}(S^o)$ is a pair of adjoint functors
\item $(j_{!Hdg},j^*):MHM_{gm,k,\mathbb C_p}(S^o)\leftrightarrows MHM_{gm,k,\mathbb C_p}(S)$ is a pair of adjoint functors.
\end{itemize}
\item[(ii)] Let $S\in\SmVar(k)$.
Let $Z=V(\mathcal I)\subset S$ an arbitrary closed subset, $\mathcal I\subset O_S$ being an ideal subsheaf. 
Denote by $j:S^o:=S\backslash Z\hookrightarrow S$. Then,
\begin{itemize}
\item $(j^*,j_{*Hdg}):D(MHM_{gm,k,\mathbb C_p}(S))\leftrightarrows D(MHM_{gm,k,\mathbb C_p}(S^o))$ is a pair of adjoint functors
\item $(j_{!Hdg},j^*):D(MHM_{gm,k,\mathbb C_p}(S^o))\leftrightarrows D(MHM_{gm,k,\mathbb C_p}(S))$ is a pair of adjoint functors.
\end{itemize}
\end{itemize}
\end{prop}

\begin{proof}
\noindent(i): Follows from proposition \ref{jHdgpropad}.

\noindent(ii):Follows from (i) and the exactness of $j^*$, $j_{*Hdg}$ and $j_{!Hdg}$.
\end{proof}

\begin{defi}\label{gammaHdgalphap}
Let $S\in\SmVar(k)$. Let $Z\subset S$ a closed subset.
Denote by $j:S\backslash Z\hookrightarrow S$ the complementary open embedding. 
\begin{itemize}
\item[(i)] We define using definition \ref{gammaHdg}, definition \ref{gammawet} and definition \ref{falphap}(iii), 
the filtered Hodge support section functor
\begin{eqnarray*}
\Gamma^{Hdg}_Z:C(MHM_{gm,k,\mathbb C_p}(S))\to C(MHM_{gm,k,\mathbb C_p}(S)), \; \; ((M,F,W),(K,W),\alpha)\mapsto \\
\Gamma^{Hdg}_Z((M,F,W),(K,W),\alpha):=(\Gamma_Z^{Hdg}(M,F,W),\Gamma_Z^w(K,W),\Gamma(\alpha)) \\
=\Cone(\ad(j^*,j_{*Hdg})(-):j_{*Hdg},j^*((M,F,W),(K,W),\alpha)\to((M,F,W),(K,W),\alpha)[-1]
\end{eqnarray*}
see definition \ref{DHdgjalphap} for the last equality, together we the canonical map 
\begin{eqnarray*}
\gamma^{Hdg}_Z((M,F,W),(K,W),\alpha):\Gamma^{Hdg}_Z((M,F,W),(K,W),\alpha)\to((M,F,W),(K,W),\alpha).
\end{eqnarray*}
\item[(i)'] Since $j_{*Hdg}:C(MHM_{gm,k,\mathbb C_p}(S^o))\to C(MHM_{gm,k,\mathbb C_p}(S))$ is an exact functor, 
$\Gamma^{Hdg}_Z$ induces the functor
\begin{eqnarray*}
\Gamma^{Hdg}_Z:D(MHM_{gm,k,\mathbb C_p}(S))\to D(MHM_{gm,k,\mathbb C_p}(S)), \\ 
((M,F,W),(K,W),\alpha)\mapsto\Gamma^{Hdg}_Z((M,F,W),(K,W),\alpha)
\end{eqnarray*}
\item[(ii)] We define using definition \ref{gammaHdg}, definition \ref{gammawet} and definition \ref{falphap}(iii) 
the dual filtered Hodge support section functor
\begin{eqnarray*}
\Gamma^{\vee,Hdg}_Z:C(MHM_{gm,k,\mathbb C_p}(S))\to C(MHM_{gm,k,\mathbb C_p}(S)), \; \; ((M,F,W),(K,W),\alpha)\mapsto \\
\Gamma^{\vee,Hdg}_Z((M,F,W),(K,W),\alpha):=(\Gamma_Z^{\vee,Hdg}(M,F,W),\Gamma_Z^{\vee,w}(K,W),\Gamma^{\vee}(\alpha)) \\
=\Cone(\ad(j_{!Hdg},j^*)(-):j_{!Hdg},j^*((M,F,W),(K,W),\alpha) \to ((M,F,W),(K,W),\alpha))
\end{eqnarray*}
see definition \ref{DHdgjalphap} for the last equality, together we the canonical map 
\begin{eqnarray*}
\gamma^{\vee,Hdg}_Z((M,F,W),(K,W),\alpha):((M,F,W),(K,W),\alpha)\to\Gamma_Z^{\vee,Hdg}((M,F,W),(K,W),\alpha).
\end{eqnarray*}
\item[(ii)'] Since $j_{!Hdg}:C(MHM_{gm,k,\mathbb C_p}(S^o))\to C(MHM_{gm,k,\mathbb C_p}(S))$ is an exact functor, 
$\Gamma^{Hdg,\vee}_Z$ induces the functor
\begin{eqnarray*}
\Gamma^{\vee,Hdg}_Z:D(MHM_{gm,k,\mathbb C_p}(S))\to D(MHM_{gm,k,\mathbb C_p}(S)), \\ 
((M,F,W),(K,W),\alpha)\mapsto\Gamma^{\vee,Hdg}_Z((M,F,W),(K,W),\alpha)
\end{eqnarray*}
\end{itemize}
\end{defi}

In the singular case it gives :

\begin{defi}\label{gammaHdgsingalphap}
Let $S\in\Var(k)$. Let $Z\subset S$ a closed subset.
Let $S=\cup_{i=1}^sS_i$ an open cover such that there exist closed embeddings
$i_i:S_i\hookrightarrow\tilde S_i$ with $\tilde S_i\in\SmVar(k)$. Denote $Z_I:=Z\cap S_I$. 
Denote by $n:S\backslash Z\hookrightarrow S$ and $\tilde n_I:\tilde S_I\backslash Z_I\hookrightarrow\tilde S_I$ 
the complementary open embeddings. 
\begin{itemize}
\item[(i)] We define using definition \ref{gammaHdgsing}, definition \ref{gammawet} and definition \ref{falphap}(iii)
the filtered Hodge support section functor
\begin{eqnarray*}
\Gamma^{Hdg}_Z:C(MHM_{gm,k,\mathbb C_p}(S))\to C(MHM_{gm,k,\mathbb C_p}(S)), \\ 
(((M_I,F,W),u_{IJ}),(K,W),\alpha)\mapsto\Gamma^{Hdg}_Z(((M_I,F,W),u_{IJ}),(K,W),\alpha):= \\
:=(\Gamma_Z^{Hdg}((M_I,F,W),u_{IJ}),\Gamma_Z^w(K,W),\Gamma(\alpha))
\end{eqnarray*}
together with the canonical map
\begin{eqnarray*}
\gamma^{Hdg}_Z(((M_I,F,W),u_{IJ}),(K,W),\alpha): \\
\Gamma^{Hdg}_Z(((M_I,F,W),u_{IJ}),(K,W),\alpha)\to(((M_I,F,W),u_{IJ}),(K,W),\alpha).
\end{eqnarray*}
\item[(i)'] By exactness of $\Gamma_Z^{Hdg}$ and $\Gamma_Z^w$ it induces the functor
\begin{eqnarray*}
\Gamma^{Hdg}_Z: D(MHM_{gm,k,\mathbb C_p}(S))\to D(MHM_{gm,k,\mathbb C_p}(S)), \\  
(((M_I,F,W),u_{IJ}),(K,W),\alpha)\mapsto\Gamma^{Hdg}_Z(((M_I,F,W),u_{IJ}),(K,W),\alpha)
\end{eqnarray*}
\item[(ii)] We define using definition \ref{gammaHdgsing}, definition \ref{gammawet} and definition \ref{falphap}(iii)
the dual filtered Hodge support section functor
\begin{eqnarray*}
\Gamma^{\vee,Hdg}_Z:C(MHM_{gm,k,\mathbb C_p}(S))\to C(MHM_{gm,k,\mathbb C_p}(S)), \; \; 
(((M_I,F,W),u_{IJ}),(K,W),\alpha)\mapsto \\
\Gamma^{\vee,Hdg}_Z(((M_I,F,W),u_{IJ}),(K,W),\alpha):=
(\Gamma_Z^{\vee,Hdg}((M_I,F,W),u_{IJ}),\Gamma_Z^{\vee,w}(K,W),\Gamma(\alpha)),
\end{eqnarray*}
together we the canonical map 
\begin{eqnarray*}
\gamma^{\vee,Hdg}_Z(((M_I,F,W),u_{IJ}),(K,W),\alpha): \\
(((M_I,F,W),u_{IJ}),(K,W),\alpha)\to\Gamma_Z^{\vee,Hdg}(((M_I,F,W),u_{IJ}),(K,W),\alpha).
\end{eqnarray*}
\item[(ii)'] By exactness of $\Gamma_Z^{\vee,Hdg}$ and $\Gamma_Z^{\vee,w}$, it induces the functor
\begin{eqnarray*}
\Gamma^{\vee,Hdg}_Z:D(MHM_{gm,k,\mathbb C_p}(S))\to D(MHM_{gm,k,\mathbb C_p}(S)), \\
(((M_I,F,W),u_{IJ}),(K,W),\alpha)\mapsto\Gamma^{\vee,Hdg}_Z(((M_I,F,W),u_{IJ}),(K,W),\alpha) \\
:=(\Gamma_Z^{\vee,Hdg}((M_I,F,W),u_{IJ}),\Gamma_Z^{\vee,w}(K,W),\Gamma(\alpha))
\end{eqnarray*}
\end{itemize}
\end{defi}

This gives the inverse image functor :

\begin{defi}\label{inverseHdgsingalphap}
Let $f:X\to S$ a morphism with $X,S\in\Var(k)$.
Assume there exist a factorization $f:X\xrightarrow{l}Y\times S\xrightarrow{p_S}S$ 
with $Y\in\SmVar(k)$, $l$ a closed embedding and $p_S$ the projection.
Let $S=\cup_{i\in I}$ an open cover such that there exist closed embeddings
$i:S_i\hookrightarrow\tilde S_i$ with $\tilde S_i\in\SmVar(k)$. 
Denote $X_I:=f^{-1}(S_I)$. We have then $X=\cup_{i\in I}X_i$ and the commutative diagrams
\begin{equation*}
\xymatrix{f:X_I\ar[r]^{l_I}\ar[rd] & Y\times S_I\ar[r]^{p_{S_I}}\ar[d]^{i_I':=(I\times i_I)} & S_I\ar[d]^{i_I} \\ 
\, & Y\times\tilde S_I\ar[r]^{p_{\tilde S_I}=:\tilde f_I} & \tilde S_I} 
\end{equation*}
\begin{itemize}
\item[(i)] For $(((M_I,F,W),u_{IJ}),(K,W),\alpha)\in C(MHM_{gm,k,\mathbb C_p}(S))$ 
we set (see definition \ref{gammaHdgsingalphap} for $l$)
\begin{eqnarray*}
f_{Hdg}^{*mod}(((M_I,F,W),u_{IJ}),(K,W),\alpha):= \\
\Gamma_X^{Hdg}((p_{\tilde S_I}^{*mod[-]}(M_I,F,W),p_{\tilde S_I}^{*mod[-]}u_{IJ}),p_S^*(K,W),p_S^*\alpha)(d_Y)[2d_Y]
\in C(MHM_{gm,k,\mathbb C_p}(X)), 
\end{eqnarray*}
\item[(ii)] For $(((M_I,F,W),u_{IJ}),(K,W),\alpha)\in C(MHM_{gm,k,\mathbb C_p}(S)))$ 
we set (see definition \ref{gammaHdgsingalphap} for $l$)
\begin{eqnarray*}
f_{Hdg}^{\hat*mod}(((M_I,F,W),u_{IJ}),(K,W),\alpha):= \\
\Gamma_X^{\vee,Hdg}((p_{\tilde S_I}^{*mod[-]}(M_I,F,W),p_{\tilde S_I}^{*mod[-]}u_{IJ}),p_S^*(K,W),p_S^*\alpha)
\in C(MHM_{gm,k,\mathbb C_p}(X)), 
\end{eqnarray*}
\end{itemize}
\end{defi}

\begin{defi}\label{otimesHdgalphap}
Let $S\in\Var(k)$. 
Let $S=\cup_{i\in I}$ an open cover such that there exist closed embeddings
$i:S_i\hookrightarrow\tilde S_i$ with $\tilde S_i\in\SmVar(k)$. 
We have the following bi-functor
\begin{eqnarray*}
(-)\otimes_{O_S}^{Hdg}(-):D(MHM_{gm,k,\mathbb C_p}(S))^2\to D(MHM_{gm,k,\mathbb C_p}(S)), \\
(((M_I,F,W),u_{IJ}),(K,W),\alpha),(((M'_I,F,W),v_{IJ}),(K',W),\alpha')\mapsto \\
(((M_I,F,W),u_{IJ}),(K,W),\alpha)\otimes_{O_S}^{Hdg}(((M'_I,F,W),v_{IJ}),(K',W),\alpha'):= \\
((M_I,F,W),u_{IJ})\otimes_{O_S}^{Hdg}((M'_I,F,W),v_{IJ}),(K,W)\otimes^{L,w}(K',W),\alpha\otimes\alpha')
\end{eqnarray*}
where the map $\alpha\otimes\alpha'$ is given in definition \ref{falphap}.
\end{defi}

\begin{prop}\label{compDmodDRHdgalphap}
Let $f_1:X\to Y$ and $f_2:Y\to S$ two morphism with $X,Y,S\in\QPVar(k)$. 
\begin{itemize}
\item[(i)]Let $\mathcal M\in C(MHM_{gm,k,\mathbb C_p}(S)))$. Then, 
\begin{equation*}
(f_2\circ f_1)^{!Hdg}(\mathcal M)=f_1^{!Hdg}f_2^{!Hdg}(\mathcal M)\in D(MHM_{gm,k,\mathbb C_p}(X)).
\end{equation*}
\item[(ii)]Let $(M,F,W)\in C(MHM_{gm,k,\mathbb C_p}(S)))$. Then,
\begin{equation*}
(f_2\circ f_1)^{*Hdg}(\mathcal M)=f_1^{*Hdg}f_2^{*Hdg}(\mathcal M)\in D(MHM_{gm,k,\mathbb C_p}(X))
\end{equation*}
\end{itemize}
\end{prop}

\begin{proof}
Immediate from definition.
\end{proof}

\begin{prop}\label{PSkMHMp}
Let $S\in\SmVar(k)$. Let $D=V(s)\subset S$ a (Cartier) divisor, where $s\in\Gamma(S,L)$ 
is a section of the line bundle $L=L_D$ associated to $D$.
\begin{itemize}
\item[(i)] Let $((M,F,W),(K,W),\alpha)\in MHM_{gm,k,\mathbb C_p}(S)$.
We have, using proposition \ref{PSkDRM}, the canonical quasi-isomorphism in $C(MHM_{gm,k,\mathbb C_p}(S))$ :
\begin{eqnarray*}
Is(M):=(Is(M),Is(K),0):((M,F,W),(K,W),\alpha)\to \\
(\psi^u_D((M,F,W),\psi^u(K,W),\psi_D\alpha)\xrightarrow{((c(x_{S^o/S}(M)),can(M)),(c(x(K)),can(K)),0)} \\
(x_{S^o/S}(M,F,W),x_{S^o/S}(K,W),x_{S^o/S}(\alpha))\oplus(\phi^u_D(M,F,W),\phi^u_D(K,W),\phi_D\alpha) \\
\xrightarrow{((exp(s\partial_s+1),var(M)),(0,T-I),var(K)),0)}
(\psi^u_D(M,F,W),\psi^u_D(K,W),\psi_D\alpha).
\end{eqnarray*}
\item[(ii)]We denote by $MHM_{gm,k,\mathbb C_p}(S\backslash D)\times_J MHM_{gm,k,\mathbb C_p}(D)$ 
the category whose set of objects consists of
\begin{equation*}
\left\{(\mathcal M,\mathcal N,a,b),\mathcal M\in MHM_{gm,k,\mathbb C_p}(S\backslash D),
\mathcal N\in MHM_{gm,k,\mathbb C_p}(D), a:\psi_{D1}\mathcal M\to N,b:N\to\psi_{D1}M \right\}
\end{equation*}
The functor (see definition \ref{DHdgpsialphap})
\begin{eqnarray*}
(j^*,\phi_{D},c,v):MHM_{gm,k,\mathbb C_p}(S)\to MHM_{gm,k,\mathbb C_p}(S\backslash D)\times_J MHM_{gm,k,\mathbb C_p}(D), \\
((M,F,W),(K,W),\alpha)\mapsto((j^*(M,F,W),j^*(K,W),j^*\alpha),\phi_{D}((M,F,W),(K,W),\alpha), can(-),var(-))
\end{eqnarray*}
is an equivalence of category.
\end{itemize}
\end{prop}

\begin{proof}
Follows from proposition \ref{PSkDRM}.
\end{proof}

Let $S\in\Var(k)$. Let $S=\cup_{i\in I}S_i$ an open cover such that there
exists closed embeddings $i_i:S_i\hookrightarrow\tilde S_i$ with $\tilde S_I\in\SmVar(k)$.
We have the category $D_{\mathcal D(1,0)fil,rh}(S/(\tilde S_I))\times_I D_{\mathbb Z_pfil,c,k}(S^{et})$  
\begin{itemize}
\item whose set of objects is the set of triples $\left\{(((M_I,F,W),u_{IJ}),(K,W),\alpha)\right\}$ with
\begin{eqnarray*} 
((M_I,F,W),u_{IJ})\in D_{\mathcal D(1,0)fil,rh}(S/(\tilde S_I)), \, (K,W)\in D_{\mathbb Z_pfil,c,k}(S^{et}), \\ 
\alpha:\mathbb B_{dr,(\tilde S_I)}(K,W)\to 
F^0DR(S)^{[-]}(((M_I,F,W),u_{IJ})^{an}\otimes_{O_S}((O\mathbb B_{dr,\tilde S_I},F),t_{IJ}))
\end{eqnarray*}
where $\alpha$ is a morphism in $D_{\mathbb B_{dr},G,fil}(S_{\mathbb C_p}^{an,pet}/(\tilde S_{I,\mathbb C_p}^{an,pet}))$,
\item and whose set of morphisms consists of 
\begin{equation*}
\phi=(\phi_D,\phi_C,[\theta]):(((M_{1I},F,W),u_{IJ}),(K_1,W),\alpha_1)\to(((M_{2I},F,W),u_{IJ}),(K_2,W),\alpha_2)
\end{equation*}
where $\phi_D:((M_1,F,W),u_{IJ})\to((M_2,F,W),u_{IJ})$ and $\phi_C:(K_1,W)\to (K_2,W)$ 
are morphisms and
\begin{eqnarray*}
\theta=(\theta^{\bullet},I(F^0DR(S)(\phi^{an}_D)\otimes I)\circ I(\alpha_1),
I(\alpha_2)\circ I(\mathbb B_{dr,(\tilde S_I)}(\phi_C\otimes I))): \\
I(\mathbb B_{dr,(\tilde S_I)}(K_1,W))[1]\to 
I(F^0DR(S)(((M_{2I},F,W),u_{IJ})^{an}\otimes_{O_S}((O\mathbb B_{dr,\tilde S_I},F),t_{IJ})))  
\end{eqnarray*}
is an homotopy,  
$I:D_{\mathbb B_{dr},G,fil}(S_{\mathbb C_p}^{an,pet}/(\tilde S^{an,pet}_{I,\mathbb C_p}))\to 
K_{\mathbb B_{dr},G,fil}(S_{\mathbb C_p}^{an,pet}/(\tilde S^{an,pet}_{I,\mathbb C_p}))$
being the injective resolution functor, and for
\begin{itemize}
\item $\phi=(\phi_D,\phi_C,[\theta]):(((M_{1I},F,W),u_{IJ}),(K_1,W),\alpha_1)\to(((M_{2I},F,W),u_{IJ}),(K_2,W),\alpha_2)$
\item $\phi'=(\phi'_D,\phi'_C,[\theta']):(((M_{2I},F,W),u_{IJ}),(K_2,W),\alpha_2)\to(((M_{3I},F,W),u_{IJ}),(K_3,W),\alpha_3)$
\end{itemize}
the composition law is given by 
\begin{eqnarray*}
\phi'\circ\phi:=(\phi'_D\circ\phi_D,\phi'_C\circ\phi_C,
I(DR(S)(\phi^{'an}_D\otimes I))\circ[\theta]+[\theta']\circ I(\mathbb B_{dr,(\tilde S_I)}(\phi_C))[1]): \\
(((M_{1I},F,W),u_{IJ}),(K_1,W),\alpha_1)\to(((M_{3I},F,W),u_{IJ}),(K_3,W),\alpha_3),
\end{eqnarray*}
in particular for 
$(((M_I,F,W),u_{IJ}),(K,W),\alpha)\in C_{\mathcal D(1,0)fil,rh}(S/(\tilde S_I))\times_I D_{\mathbb Z_pfil,c,k}(S^{et})$,
\begin{equation*}
I_{(((M_I,F,W),u_{IJ}),(K,W),\alpha)}=((I_{M_I}),I_K,0),
\end{equation*}
\end{itemize}
and also the category 
$D_{\mathcal D(1,0)fil,rh,\infty}(S/(\tilde S_I))\times_I D_{\mathbb Z_pfil,c,k}(S^{et})$ defined in the same way,
together with the localization functor
\begin{eqnarray*}
(D(zar),I):C_{\mathcal D(1,0)fil,rh}(S/(\tilde S_I))\times_I D_{fil,c,k}(S^{et})
\to D_{\mathcal D(1,0)fil,rh}(S/(\tilde S_I))\times_I D_{fil,c,k}(S^{et}) \\
\to D_{\mathcal D(1,0)fil,rh,\infty}(S/(\tilde S_I))\times_I D_{fil,c,k}(S^{et}).
\end{eqnarray*}
Note that if $\phi=(\phi_D,\phi_C,[\theta]):(((M_1,F,W),u_{IJ}),(K_1,W),\alpha_1)\to(((M_2,F,W),u_{IJ}),(K_2,W),\alpha_2)$
is a morphism in $D_{\mathcal D(1,0)fil,rh}(S/(\tilde S_I))\times_I D_{\mathbb Z_pfil,c,k}(S^{et})$
such that $\phi_D$ and $\phi_C$ are isomorphism then $\phi$ is an isomorphism (see remark \ref{CGremp}).
Moreover,
\begin{itemize}
\item For 
$(((M_{I},F,W),u_{IJ}),(K,W),\alpha)\in D_{\mathcal D(1,0)fil,rh}(S/(\tilde S_I))\times_I D_{\mathbb Z_pfil,c,k}(S^{et})$, 
we set
\begin{equation*}
(((M_{I},F,W),u_{IJ}),(K,W),\alpha)[1]:=(((M_{I},F,W),u_{IJ})[1],(K,W)[1],\alpha[1]).
\end{equation*}
\item For 
\begin{equation*}
\phi=(\phi_D,\phi_C,[\theta]):(((M_{1I},F,W),u_{IJ}),(K_1,W),\alpha_1)\to(((M_{2I},F,W),u_{IJ}),(K_2,W),\alpha_2)
\end{equation*}
a morphism in $D_{\mathcal D(1,0)fil,rh}(S/(\tilde S_I))\times_I D_{\mathbb Z_p,fil,c,k}(S^{et})$, 
we set (see \cite{CG} definition 3.12)
\begin{eqnarray*}
\Cone(\phi):=(\Cone(\phi_D),\Cone(\phi_C),((\alpha_1,\theta),(\alpha_2,0)))
\in D_{\mathcal D(1,0)fil,rh}(S/(\tilde S_I))\times_I D_{\mathbb Z_pfil,c,k}(S^{et}),
\end{eqnarray*}
$((\alpha_1,\theta),(\alpha_2,0))$ being the matrix given by the composition law, together with the canonical maps
\begin{itemize}
\item $c_1(-)=(c_1(\phi_D),c_1(\phi_C),0):(((M_{2I},F,W),u_{IJ}),(K_2,W),\alpha_2)\to\Cone(\phi)$
\item $c_2(-)=(c_2(\phi_D),c_2(\phi_C),0):\Cone(\phi)\to (((M_{1I},F,W),u_{IJ}),(K_1,W),\alpha_1)[1]$.
\end{itemize}
\end{itemize}

We have then the following :

\begin{thm}\label{Bekp}
\begin{itemize}
\item[(i)]Let $S\in\Var(k)$. Let $S=\cup_{i\in I}S_i$ an open cover such that there exists
closed embedding $i_i:S_i\hookrightarrow\tilde S_i$ with $\tilde S_i\in\SmVar(k)$. Then the full embedding
\begin{eqnarray*}
\iota_S:MHM_{gm,k,\mathbb C_p}(S)\hookrightarrow
\PSh^0_{\mathcal D(1,0)fil,rh}(S/(\tilde S_I))\times_I P_{\mathbb Z_pfil,k}(S^{et}) 
\hookrightarrow C_{\mathcal D(1,0)fil,rh}(S/(\tilde S_I))\times_I D_{\mathbb Z_pfil,c,k}(S^{et}) 
\end{eqnarray*}
induces a full embedding
\begin{equation*}
\iota_S:D(MHM_{gm,k,\mathbb C_p}(S))\hookrightarrow 
D_{\mathcal D(1,0)fil,rh}(S/(\tilde S_I))\times_I D_{\mathbb Z_pfil,c,k}(S^{et}) 
\end{equation*}
whose image consists of 
$(((M_I,F,W),u_{IJ}),(K,W),\alpha)\in D_{\mathcal D(1,0)fil,rh}(S/(\tilde S_I))\times_I D_{\mathbb Z_pfil,c,k}(S^{et})$ such that 
\begin{equation*}
((H^n(M_I,F,W),H^n(u_{IJ})),H^n(K,W),H^n\alpha)\in MHM_{gm,k,\mathbb C_p}(S) 
\end{equation*}
for all $n\in\mathbb Z$ and such that for all $p\in\mathbb Z$,
the differentials of $\Gr_W^p(M_I,F)$ are strict for the filtrations $F$.
\item[(i)']Let $S\in\Var(k)$. Let $S=\cup_{i\in I}S_i$ an open cover such that there exists
closed embedding $i_i:S_i\hookrightarrow\tilde S_i$ with $\tilde S_i\in\SmVar(k)$. Then,
\begin{eqnarray*}
D(MHM_{gm,k,\mathbb C_p}(S))=
<(\int^{FDR}_f(n\times I)_{!Hdg}(\Gamma_X^{\vee,Hdg}(O_{\mathbb P^{N,o}\times\tilde S_I},F_b),x_{IJ})(d),
Rf_*\mathbb Z_{p,X^{et}}^w,f_*\alpha(X)), \\
(f:X\xrightarrow{l}\mathbb P^{N,o}\times S\xrightarrow{p}S)\in\QPVar(k),\; d\in\mathbb Z> \\
=<(\int^{FDR}_f((\Gamma_X^{\vee,Hdg}(O_{\mathbb P^{N,o}\times\tilde S_I},F_b),x_{IJ})(d),
Rf_*\mathbb Z_{p,X^{et}},f_*\alpha(X)), \\
(f:X\xrightarrow{l}\mathbb P^{N,o}\times S\xrightarrow{p}S)\in\QPVar(k), \; \mbox{proper}, \; X \mbox{smooth}> \\
\subset D_{\mathcal D(1,0)fil,rh}(S/(\tilde S_I))\times_I D_{\mathbb Z_pfil,c,k}(S^{et}) 
\end{eqnarray*}
where $n:\mathbb P^{N,o}\hookrightarrow\mathbb P^N$ are open embeddings, $l$ are closed embedding
and $<,>$ means the full triangulated category generated by and $(-)$ is the shift of the F-filtration.
\item[(ii)]Let $S\in\Var(k)$. Let $S=\cup_{i\in I}S_i$ an open cover such that there exists
closed embedding $i_i:S_i\hookrightarrow\tilde S_i$ with $\tilde S_i\in\SmVar(k)$. Then the full embedding
\begin{eqnarray*}
\iota_S:MHM_{gm,k,\mathbb C_p}(S)\hookrightarrow
\PSh^0_{\mathcal D(1,0)fil,rh}(S/(\tilde S_I))\times_I P_{\mathbb Z_pfil,k}(S^{et})
\hookrightarrow C_{\mathcal D(1,0)fil,rh}(S/(\tilde S_I))\times_I D_{\mathbb Z_pfil,c,k}(S^{et}) 
\end{eqnarray*}
induces a full embedding
\begin{equation*}
\iota_S:D(MHM_{gm,k,\mathbb C_p}(S))\hookrightarrow 
D_{\mathcal D(1,0)fil,\infty,rh}(S/(\tilde S_I))\times_I D_{\mathbb Z_pfil,c,k}(S^{et}) 
\end{equation*}
whose image consists of 
$(((M_I,F,W),u_{IJ}),(K,W),\alpha)\in D_{\mathcal D(1,0)fil,\infty,rh}(S/(\tilde S_I))\times_I D_{\mathbb Z_pfil,c,k}(S^{et})$
such that 
\begin{equation*}
((H^n(M_I,F,W),H^n(u_{IJ})),H^n(K,W),H^n\alpha)\in MHM_{gm,k,\mathbb C_p}(S) 
\end{equation*}
for all $n\in\mathbb Z$ and such that there exist $r\in\mathbb Z$ and an $r$-filtered homotopy equivalence
$((M_I,F,W),u_{IJ})\to ((M'_I,F,W),u_{IJ})$ such that for all $p\in\mathbb Z$
the differentials of $\Gr_W^p(M'_I,F)$ are strict for the filtrations $F$.
\end{itemize}
\end{thm}

\begin{proof}
\noindent(i): We first show that $\iota_S$ is fully faithfull, that is for all
$\mathcal M=(((M_I,F,W),u_{IJ}),(K,W),\alpha),\mathcal M'=(((M'_I,F,W),u_{IJ}),(K',W),\alpha')\in MHM_{gm,k,\mathbb C_p}(S)$ 
and all $n\in\mathbb Z$,
\begin{eqnarray*}
\iota_S:\Ext_{D(MHM_{gm,k,\mathbb C_p}(S))}^n(\mathcal M,\mathcal M'):=
\Hom_{D(MHM_{gm,k,\mathbb C_p}(S))}(\mathcal M,\mathcal M'[n]) \\
\to\Ext_{\mathcal D(S)}^n(\mathcal M,\mathcal M')
:=\Hom_{\mathcal D(S):=D_{\mathcal D(1,0)fil,rh}(S/(\tilde S_I))\times_I D_{fil}(S^{et})}(\mathcal M,\mathcal M'[n])
\end{eqnarray*}
For this it is enough to assume $S$ smooth. We then proceed by induction on $max(\dim\supp(M),\dim\supp(M'))$. 
\begin{itemize}
\item For $\supp(M)=\supp(M')=\left\{s\right\}$, 
it is the theorem for mixed hodge complexes or absolute Hodge complexes, see \cite{CG}. 
If $\supp(M)=\left\{s\right\}$ and $\supp(M')=\left\{s'\right\}$ and $s'\neq s$, 
then by the localization exact sequence
\begin{equation*}
\Ext_{D(MHM_{gm,k,\mathbb C_p}(S))}^n(\mathcal M,\mathcal M')=0=\Ext_{\mathcal D(S)}^n(\mathcal M,\mathcal M')
\end{equation*}
\item Denote $\supp(M)=Z\subset S$ and $\supp(M')=Z'\subset S$.
There exist an open subset $S^o\subset S$ such that $Z^o:=Z\cap S^o$ and $Z^{'o}:=Z'\cap S^o$ are smooth,
and $\mathcal M_{|Z^o}:=((i^*\Gr_{V_{Z^o},0}M_{|S^o},F,W),i^*j^*(K,W),\alpha^*(i))\in MHM_{gm,k}(Z^o)$ and 
$\mathcal M'_{|Z^{'o}}:=((i^{'*}\Gr_{V_{Z^{'o}},0}M'_{|S^o},F,W),i^{'*}j^*K,\alpha^*(i'))\in MHM_{gm,k}(Z^{'o})$ 
are variation of geometric mixed Hodge structure over $k\subset\mathbb C$, 
where $j:S^o\hookrightarrow S$ is the open embedding, and
$i:Z^o\hookrightarrow S^o$, $i':Z^{'o}\hookrightarrow S^o$ the closed embeddings.
Considering the connected components of $Z^o$ and $Z^{'o}$, we way assume that $Z^o$ and $Z^{'o}$ are connected.
Shrinking $S^o$ if necessary, we may assume that either $Z^o=Z^{'o}$ or $Z^o\cap Z^{'o}=\emptyset$,
We denote $D=S\backslash S^o$. Shrinking $S^o$ if necessary, 
we may assume that $D$ is a divisor and denote by $l:S\hookrightarrow L_D$ the zero section embedding.
\begin{itemize}
\item If $Z^o=Z^{'o}$, denote $i:Z^o\hookrightarrow S^o$ the closed embedding.
We have then the following commutative diagram
\begin{equation*}
\xymatrix{\Ext_{D(MHM_{gm,k,\mathbb C_p}(S^o))}^n(\mathcal M_{|S^o},\mathcal M'_{|S^o})
\ar[rr]^{\iota_{S^o}}\ar[d]_{(i^*\Gr_{V_{Z^o},0},i^*,\alpha^*(i))} & \, & 
\Ext_{\mathcal D(S^o)}^n(\mathcal M_{|S^o},\mathcal M'_{|S^o})\ar[d]^{(i^*\Gr_{V_{Z^o},0},i^*,\alpha^*(i))} \\
\Ext_{D(MHM_{gm,k,\mathbb C_p}(Z^o))}^n(\mathcal M_{|Z^o},\mathcal M'_{|Z^o})
\ar[rr]^{\iota_{Z^o}}\ar[u]_{(i_{*mod},i_*,\alpha_*(i))} & \, &
\Ext_{\mathcal D(Z^o)}^n(\mathcal M_{|Z^o},\mathcal M'_{|Z^o})\ar[u]^{(i_{*mod},i_*,\alpha_*(i))}}
\end{equation*}
Now we prove that $\iota_{Z^o}$ is an isomorphism similarly to the proof the the generic case of \cite{Be}.
On the other hand the left and right colummn are isomorphisms.
Hence $\iota_{S^o}$ is an isomorphism by the diagram.
\item If $Z^o\cap Z^{'o}=\emptyset$, we consider the following commutative diagram
\begin{equation*}
\xymatrix{\Ext_{D(MHM_{gm,k,\mathbb C_p}(S^o))}^n(\mathcal M_{|S^o},\mathcal M'_{|S^o})
\ar[rr]^{\iota_{S^o}}\ar[d]_{(i^*\Gr_{V_{Z^o},0},i^*,\alpha^*(i))} & \, & 
\Ext_{\mathcal D(S^o)}^n(\mathcal M_{|S^o},\mathcal M'_{|S^o})\ar[d]^{(i^*\Gr_{V_{Z^o},0},i^*,\alpha^*(i))} \\
\Ext_{D(MHM_{gm,k,\mathbb C_p}(Z^o))}^n(\mathcal M_{|Z^o},0)=0\ar[rr]^{\iota_{Z^o}}\ar[u]_{(i_{*mod},i_*,\alpha_*(i))} & \, &
\Ext_{\mathcal D(Z^o)}^n(\mathcal M_{|Z^o},0)=0\ar[u]^{(i_{*mod},i_*,\alpha_*(i))}}
\end{equation*}
where the left and right column are isomorphism by strictness of the $V_{Z^o}$ filtration
(use a bi-filtered injective resolution with respect to $F$ and $V_{Z^o}$ for the right column).
\end{itemize}
\item We consider now the following commutative diagram in $C(\mathbb Z)$ 
where we denote for short $H:=D(MHM_{gm,k,\mathbb C_p}(S))$
\begin{equation*}
\xymatrix{0\ar[r] & \Hom_{H}^{\bullet}(\Gamma^{\vee,Hdg}_D\mathcal M,\Gamma^{Hdg}_D\mathcal M')
\ar[r]^{\Hom(-,\gamma^{Hdg}_D(\mathcal M'))}\ar[d]^{\iota_S} &
\Hom_{H}^{\bullet}(\Gamma^{\vee,Hdg}_D\mathcal M,\mathcal M')
\ar[r]^{\Hom(-,\ad(j^*,j_{*Hdg})(\mathcal M'))}\ar[d]^{\iota_S} &
\Hom_{H}^{\bullet}(\Gamma^{\vee,Hdg}_D\mathcal M,j_{*Hdg}j^*\mathcal M')\ar[r]\ar[d]^{\iota_S} & 0 \\
0\ar[r] & \Hom_{\mathcal D(S)}^{\bullet}(\Gamma^{\vee,Hdg}_D\mathcal M,\Gamma^{Hdg}_D\mathcal M')
\ar[r]^{\Hom(-,\gamma^{Hdg}_D(\mathcal M'))} &
\Hom_{\mathcal D(S)}^{\bullet}(\Gamma^{\vee,Hdg}_D\mathcal M,\mathcal M')\ar[r]^{\Hom(-,\ad(j^*,j_{*Hdg})(\mathcal M'))} &
\Hom_{\mathcal D(S)}^{\bullet}(\Gamma^{\vee,Hdg}_D\mathcal M,j_{*Hdg}j^*\mathcal M')\ar[r] & 0}
\end{equation*}
whose lines are exact sequence. We have on the one hand,
\begin{equation*}
\Hom_{D(MHM_{gm,k,\mathbb C_p}(S))}^{\bullet}(\Gamma^{\vee,Hdg}_D\mathcal M,j_{*Hdg}j^*\mathcal M')=0=
\Hom_{\mathcal D(S)}^{\bullet}(\Gamma^{\vee,Hdg}_D\mathcal M,j_{*Hdg}j^*\mathcal M')
\end{equation*}
On the other hand by induction hypothesis
\begin{equation*}
\iota_S:\Hom_{D(MHM_{gm,k,\mathbb C_p}(S))}^{\bullet}(\Gamma^{\vee,Hdg}_D\mathcal M,\Gamma^{Hdg}_D\mathcal M')\to
\Hom_{\mathcal D(S)}^{\bullet}(\Gamma^{\vee,Hdg}_D\mathcal M,\Gamma^{Hdg}_D\mathcal M')
\end{equation*}
is a quasi-isomorphism. Hence, by the diagram
\begin{equation*}
\iota_S:\Hom_{D(MHM_{gm,k,\mathbb C_p}(S))}^{\bullet}(\Gamma^{\vee,Hdg}_D\mathcal M,\mathcal M')\to
\Hom_{\mathcal D(S)}^{\bullet}(\Gamma^{\vee,Hdg}_D\mathcal M,\mathcal M')
\end{equation*}
is a quasi-isomorphism.
\item We consider now the following commutative diagram in $C(\mathbb Z)$ 
where we denote for short $H:=D(MHM_{gm,k,\mathbb C_p}(S))$
\begin{equation*}
\xymatrix{0\ar[r] & \Hom_{H}^{\bullet}(\Gamma^{\vee,Hdg}_D\mathcal M,\mathcal M')
\ar[r]^{\Hom(\gamma^{\vee,Hdg}_D(\mathcal M),-)}\ar[d]^{\iota_S} &
\Hom_{H}^{\bullet}(\mathcal M,\mathcal M')
\ar[r]^{\Hom(\ad(j_{!Hdg},j^*)(\mathcal M'),-)}\ar[d]^{\iota_S} &
\Hom_{H}^{\bullet}(j_{!Hdg}j^*\mathcal M,\mathcal M')\ar[r]\ar[d]^{\iota_S} & 0 \\
0\ar[r] & \Hom_{\mathcal D(S)}^{\bullet}(\Gamma^{\vee,Hdg}_D\mathcal M,\mathcal M')
\ar[r]^{\Hom(\gamma^{\vee,Hdg}_D(\mathcal M),-)} &
\Hom_{\mathcal D(S)}^{\bullet}(\mathcal M,\mathcal M')\ar[r]^{\Hom(\ad(j_{!Hdg},j^*)(\mathcal M),-)} &
\Hom_{\mathcal D(S)}^{\bullet}(j_{!Hdg}j^*\mathcal M,\mathcal M')\ar[r] & 0}
\end{equation*}
whose lines are exact sequence. On the one hand, the commutative diagram
\begin{equation*}
\xymatrix{\Hom_{D(MHM_{gm,k,\mathbb C_p}(S))}^{\bullet}(j_{!Hdg}j^*\mathcal M,\mathcal M')\ar[r]^{j^*}\ar[d]^{\iota_{S}} &
\Hom_{D(MHM_{gm,k,\mathbb C_p}(S^o))}^{\bullet}(j^*\mathcal M,j^*\mathcal M')\ar[d]^{\iota_{S^o}} \\
\Hom_{\mathcal D(S)}^{\bullet}(j_{!Hdg}j^*\mathcal M,\mathcal M')\ar[r]^{j^*} &
\Hom_{\mathcal D(S^o)}^{\bullet}(j^*\mathcal M,j^*\mathcal M')}
\end{equation*}
together with the fact that the horizontal arrows $j^*$ are quasi-isomorphism 
by the functoriality given the uniqueness of the $V_S$ filtration for the embedding $l:S\hookrightarrow L_D$, 
(use a bi-filtered injective resolution with respect to $F$ and $V_S$ for the lower arrow)
and the fact that $\iota_{S^o}$ is a quasi-isomorphism by the first two point, show that
\begin{equation*}
\iota_S:\Hom_{D(MHM_{gm,k,\mathbb C_p}(S))}^{\bullet}(j_{!Hdg}j^*\mathcal M,\mathcal M')\to
\Hom_{\mathcal D(S)}^{\bullet}(j_{!Hdg}j^*\mathcal M,\mathcal M')
\end{equation*}
is a quasi-isomorphism. On the other hand, by the third point
\begin{equation*}
\iota_S:\Hom_{D(MHM_{gm,k,\mathbb C_p}(S))}^{\bullet}(\Gamma^{\vee,Hdg}_D\mathcal M,\mathcal M')\to
\Hom_{\mathcal D(S)}^{\bullet}(\Gamma^{\vee,Hdg}_D\mathcal M,\mathcal M')
\end{equation*}
is a quasi-isomorphism. Hence, by the diagram
\begin{equation*}
\iota_S:\Hom_{D(MHM_{gm,k,\mathbb C_p}(S))}^{\bullet}(\Gamma^{\vee,Hdg}_D\mathcal M,\mathcal M')\to
\Hom_{\mathcal D(S)}^{\bullet}(\Gamma^{\vee,Hdg}_D\mathcal M,\mathcal M')
\end{equation*}
is a quasi-isomorphism.
\end{itemize}
This shows the fully faithfulness. We now prove the essential surjectivity : let
\begin{equation*}
(((M_I,F,W),u_{IJ}),(K,W),\alpha)\in C_{\mathcal D(1,0)fil,rh}(S/(\tilde S_I))\times_I C_{fil}(S^{et}) 
\end{equation*}
such that the cohomology are mixed hodge modules and such that the differential are strict.
We proceed by induction on $card\left\{n\in\mathbb Z\right\}, \, \mbox{s.t.} H^n(M_I,F,W)\neq 0$ by taking for 
the cohomological troncation 
\begin{equation*}
\tau^{\leq n}(((M_I,F,W),u_{IJ}),(K,W),\alpha):=
((\tau^{\leq n}(M_I,F,W),\tau^{\leq n}u_{IJ}),\tau^{\leq n}(K,W),\tau^{\leq n}\alpha)
\end{equation*}
and using the fact that the differential are strict for the filtration $F$ and the fully faithfullness.

\noindent(i)':Follows from (i).

\noindent(ii):Follows from (i).Indeed, in the composition of functor
\begin{eqnarray*}
\iota_S:D(MHM_{gm,k,\mathbb C_p}(S))\xrightarrow{\iota_S}
D_{\mathcal D(1,0)fil,rh}(S/(\tilde S_I))\times_I D_{\mathbb Z_pfil,c,k}(S^{et}) \\
\to D_{\mathcal D(1,0)fil,\infty,rh}(S/(\tilde S_I))\times_I D_{\mathbb Z_pfil,c,k}(S^{et}) 
\end{eqnarray*}
the second functor which is the localization functor is an isomorphism on the full subcategory
\begin{equation*}
D_{\mathcal D(1,0)fil,rh}(S/(\tilde S_I))^{st}\times_I D_{\mathbb Z_pfil,c,k}(S^{et})
\subset D_{\mathcal D(1,0)fil,rh}(S/(\tilde S_I))\times_I D_{\mathbb Z_pfil,c,k}(S^{et}) 
\end{equation*}
constisting of complex such that the differentials are strict for $F$, 
and the first functor $\iota_S$ is a full embedding by (i) and 
$\iota_S(D(MHM_{gm,k,\mathbb C_p}(S)))\subset D_{\mathcal D(1,0)fil,rh}(S/(\tilde S_I))^{st}\times_I D_{\mathbb Z_pfil,c,k}(S^{et})$.
\end{proof}

\begin{defi}\label{DHdgalphap}
Let $f:X\to S$ a morphism with $X,S\in\Var(k)$. 
Assume there exist a factorization $f:X\xrightarrow{l}Y\times S\xrightarrow{p_S}S$
with $Y\in\SmVar(k)$, $l$ a closed embedding and $p_S$ the projection.
Let $\bar Y\in\PSmVar(k)$ a smooth compactification of $Y$ with $n:Y\hookrightarrow\bar Y$ the open embedding.
Then $\bar f:\bar X\xrightarrow{\bar l}\bar Y\times_S\xrightarrow{\bar p_S}S$ is a compactification of $f$,
with $\bar X\subset\bar Y\times S$ the closure of $X$ and $\bar l$ the closed embedding,
we denote by $n':X\hookrightarrow\bar X$ the closed embedding so that $f=\bar f\circ n'$.
\begin{itemize}
\item[(i)]For $(((M_I,F,W),u_{IJ}),(K,W),\alpha)\in C(MHM_{gm,k,\mathbb C_p}(X))$, 
we define, using definition \ref{DHdgsing} and theorem \ref{Bekp}, 
\begin{eqnarray*}
Rf_{*Hdg}(((M_I,F,W),u_{IJ}),(K,W),\alpha):=\iota_S^{-1}(\int_f^{Hdg}((M_I,F,W),u_{IJ}),Rf_{*w}(K,W),f_*(\alpha)) \\
\in D(MHM_{gm,k,\mathbb C_p}(S))
\end{eqnarray*}
where $f_*(\alpha)$ is given in definition \ref{falphap}, and since 
\begin{itemize} 
\item by definition 
\begin{equation*}
H^i(\int_{\bar f}^{FDR}\Gr_W^k(n\times I)_{*Hdg}((M_I,F,W),u_{IJ}),R\bar f_*\Gr_W^kn'_{*w}(K,W),\Gr_W^kf_*\alpha)
\in HM_{gm,k,\mathbb C_p}(S) 
\end{equation*}
for all $i,k\in\mathbb Z$, 
hence by the spectral sequence for the filtered complexes 
$\int_{\bar f}^{FDR}(n\times I)_{*Hdg}((M_I,W),u_{IJ})$ and $R\bar f_*n'_{*w}(K,W)$
\begin{eqnarray*}
\Gr_W^kH^i(\int_{f}^{Hdg}((M_I,F,W),u_{IJ}),Rf_{*w}(K,W),f_*\alpha):= \\
(\Gr_W^kH^i\int_{\bar f}^{FDR}((M_I,F,W),u_{IJ}),\Gr_W^kH^iR\bar f_*n'_{*w}(K,W),\Gr_W^kH^if_*\alpha)\in HM_{gm,k,\mathbb C_p}(S) 
\end{eqnarray*}
this gives by definition 
$H^i(\int_f^{Hdg}((M_I,F,W),u_{IJ}),Rf_{*w}(K,W),f_*(\alpha))\in MHM_{gm,k,\mathbb C_p}(S)$ for all $i\in\mathbb Z$. 
\item $\int_{f}^{Hdg}((M_I,F,W),u_{IJ})$ is the class of a complex such that the differential are strict for $F$
by theorem \ref{Sa12} in the complex case.
\end{itemize}
\item[(ii)]For $(((M_I,F,W),u_{IJ}),(K,W),\alpha)\in C(MHM_{gm,k,\mathbb C_p}(X))$, 
we define, using definition \ref{DHdg} and theorem \ref{Bekp}, 
\begin{eqnarray*}
Rf_{!Hdg}(((M_I,F,W),u_{IJ}),(K,W),\alpha):=\iota_S^{-1}(\int_{f!}^{Hdg}((M_I,F,W),u_{IJ}),Rf_{!w}(K,W),f_!(\alpha)) \\
\in D(MHM_{gm,k,\mathbb C_p}(S))
\end{eqnarray*}
where $f_!(\alpha)$ is given in definition \ref{falphap}, and since
\begin{itemize} 
\item by definition 
\begin{equation*}
H^i(\int_{\bar f}^{FDR}\Gr_W^k(n\times I)_{!Hdg}((M_I,F,W),u_{IJ}),
R\bar f_*\Gr_W^kn'_{!w}(K,W),\Gr_W^kf_!\alpha)\in HM_{gm,k,\mathbb C_p}(S) 
\end{equation*}
for all $i,k\in\mathbb Z$, 
hence by the spectral sequence for the filtered complexes 
$\int_{\bar f}^{FDR}(n\times I)_{!Hdg}((M_I,W),u_{IJ})$ and $R\bar f_*n'_{!w}(K,W)$
\begin{eqnarray*}
\Gr_W^kH^i(\int_{f!}^{Hdg}((M_I,F,W),u_{IJ}),Rf_{!w}K,f_!\alpha):= \\
(\Gr_W^kH^i\int_{\bar f}^{FDR}(n\times I)_{!Hdg}((M_I,F,W),u_{IJ}),
\Gr_W^kH^iR\bar f_*n'_{!w}(K,W),\Gr_W^kH^if_!\alpha)\in HM_{gm,k,\mathbb C_p}(S) 
\end{eqnarray*}
this gives by definition 
$H^i(\int_{f!}^{Hdg}((M_I,F,W),u_{IJ}),Rf_{!w}(K,W),f_!(\alpha))\in MHM_{gm,k,\mathbb C_p}(S)$ for all $i\in\mathbb Z$. 
\item $\int_{f!}^{Hdg}((M_I,F,W),u_{IJ})$ is the class of a complex such that the differential are strict for $F$
by theorem \ref{Sa12} in the complex case.
\end{itemize}
\end{itemize}
\end{defi}

\begin{prop}\label{compDmodDRHdgDalphap}
Let $f_1:X\to Y$ and $f_2:Y\to S$ two morphism with $X,Y,S\in\QPVar(k)$. 
\begin{itemize}
\item[(i)]Let $\mathcal M\in C(MHM_{gm,k,\mathbb C_p}(X))$. Then, 
\begin{equation*}
R(f_2\circ f_1)^{Hdg}_*(\mathcal M)=Rf^{Hdg}_{2*}Rf^{Hdg}_{1*}(\mathcal M)\in D(MHM_{gm,k,\mathbb C_p}(S)).
\end{equation*}
\item[(ii)]Let $\mathcal M\in C(MHM_{gm,k,\mathbb C_p}(X))$. Then,
\begin{equation*}
R(f_2\circ f_1)^{Hdg}_!(\mathcal M)=Rf^{Hdg}_{2!}Rf^{Hdg}_{1!}(\mathcal M)\in D(MHM_{gm,k,\mathbb C_p}(S))
\end{equation*}
\end{itemize}
\end{prop}

\begin{proof}
Immediate from definition.
\end{proof}

Let $k\subset K\subset\mathbb C_p$ a subfield of a $p$-adic field $K$.
Definition \ref{inverseHdgsingalphap}, definition \ref{DHdgalphap} and gives by proposition \ref{compDmodDRHdgalphap} 
and proposition \ref{compDmodDRHdgDalphap} respectively, the following 2 functors :
\begin{itemize}
\item We have the following 2 functor on the category of algebraic varieties over $k\subset\mathbb C_p$
\begin{eqnarray*}
D(MHM_{gm,k,\mathbb C_p}(\cdot)):\QPVar(k)\to\TriCat, \; S\mapsto D(MHM_{gm,k,\mathbb C_p}(S)), \\
(f:T\to S)\longmapsto (f^{*Hdg}:(((M_I,F,W),u_{IJ}),(K,W),\alpha)\mapsto \\
f^{!Hdg}(((M_I,F,W),u_{IJ}),(K,W),\alpha):=(f^{*mod}_{Hdg}(((M_I,F,W),u_{IJ})),f^{!w}(K,W),f^!\alpha)).
\end{eqnarray*}
see definition \ref{inverseHdgsing} and definition \ref{falphap} for the equality.
\item We have the following 2 functor on the category of quasi-projective algebraic varieties over $k\subset\mathbb C_p$
\begin{eqnarray*}
D(MHM_{gm,k,\mathbb C_p}(\cdot)):\QPVar(k)\to\TriCat, \; S\mapsto D(MHM_{gm,k,\mathbb C_p}(S)), \\
(f:T\to S)\longmapsto (f_{*Hdg}:(((M_I,F,W),u_{IJ}),(K,W),\alpha)\mapsto Rf_{*Hdg}(((M_I,F,W),u_{IJ}),(K,W),\alpha)).
\end{eqnarray*}
\item We have the following 2 functor on the category of quasi-projective algebraic varieties over $k\subset\mathbb C_p$
\begin{eqnarray*}
D(MHM_{gm,k,\mathbb C_p}(\cdot)):\QPVar(k)\to\TriCat, \; S\mapsto D(MHM_{gm,k,\mathbb C_p}(S)), \\
(f:T\to S)\longmapsto (f_{!Hdg}:(((M_I,F,W)),(K,W),\alpha)\mapsto f_{!Hdg}(((M_I,F,W),u_{IJ}),(K,W),\alpha)).
\end{eqnarray*}
\item We have the following 2 functor on the category of algebraic varieties over $k\subset\mathbb C_p$
\begin{eqnarray*}
D(MHM_{gm,k,\mathbb C_p}(\cdot)):\QPVar(k)\to\TriCat, \; S\mapsto D(MHM_{gm,k,\mathbb C_p}(S)), \\
(f:T\to S)\longmapsto (f^{!Hdg}:(((M_I,F,W),u_{IJ}),(K,W),\alpha)\mapsto \\
f^{*Hdg}(((M_I,F,W),u_{IJ}),(K,W),\alpha):=(f^{\hat*mod}_{Hdg}(((M_I,F,W),u_{IJ})),f^{*w}(K,W),f^*\alpha)).
\end{eqnarray*}
see definition \ref{inverseHdgsing} and definition \ref{falphap} for the equality.
\end{itemize}

\begin{prop}\label{HdgpropadsingCp}
Let $f:X\to S$ with $S,X\in\QPVar(k)$. Then
\begin{itemize}
\item[(i)] $(f^{*Hdg},Rf^{Hdg}_*):D(MHM_{gm,k,\mathbb C_p}(S))\to D(MHM_{gm,k,\mathbb C_p}(X))$ is a pair of adjoint functors.
\begin{itemize}
\item For $(((M_I,F,W),u_{IJ}),(K,W),\alpha)\in C(MHM_{gm,k,\mathbb C_p}(S))$, 
\begin{eqnarray*}
\ad(f^{*Hdg},Rf^{Hdg}_*)(((M_I,F,W),u_{IJ}),(K,W),\alpha):= \\
(\ad(f_{Hdg}^{\hat*mod},Rf^{Hdg}_*)((M_I,F,W),u_{IJ}),\ad(f^{*w},Rf_{*w})(K,W)): \\
(((M_I,F,W),u_{IJ}),(K,W),\alpha)\to Rf^{Hdg}_*f^{*Hdg}(((M_I,F,W),u_{IJ}),(K,W),\alpha) 
\end{eqnarray*}
is the adjonction map in $D(MHM_{gm,k,\mathbb C_p}(S))$. 
\item For $(((N_I,F,W),u_{IJ}),(P,W),\beta)\in C(MHM_{gm,k,\mathbb C_p}(X))$,  
\begin{eqnarray*}
\ad(f^{*Hdg},Rf^{Hdg}_*)(((N_I,F,W),u_{IJ}),(P,W),\beta):= \\
(\ad(f_{Hdg}^{\hat*mod},Rf^{Hdg}_*)((N_I,F,W),u_{IJ}),\ad(f^{*w},Rf_{*w})(P,W)): \\
f^{*Hdg}Rf^{Hdg}_*(((N_I,F,W),u_{IJ}),(P,W),\beta)\to(((N_I,F,W),u_{IJ}),(P,W),\beta) 
\end{eqnarray*}
is the adjonction map in $D(MHM_{gm,k,\mathbb C_p}(X))$ 
\end{itemize}
\item[(ii)]$(Rf^{Hdg}_!,f^{!Hdg}):D(MHM_{gm,k,\mathbb C_p}(X))\to D(MHM_{gm,k,\mathbb C_p}(S))$ is a pair of adjoint functors.
\begin{itemize}
\item For $(((M_I,F,W),u_{IJ}),(K,W),\alpha)\in C(MHM_{gm,k,\mathbb C_p}(S))$, 
\begin{eqnarray*}
\ad(Rf^{Hdg}_!,f^{!Hdg})(((M_I,F,W),u_{IJ}),(K,W),\alpha):= \\
(\ad(f_{Hdg}^{*mod},Rf^{Hdg}_!)((M_I,F,W),u_{IJ}),\ad(f^{!w},Rf_{!w})(K,W)): \\
Rf^{Hdg}_!f^{!Hdg}(((M_I,F,W),u_{IJ}),(K,W),\alpha)\to(((M_I,F,W),u_{IJ}),(K,W),\alpha) 
\end{eqnarray*}
is the adjonction map in $D(MHM_{gm,k,\mathbb C_p}(S))$. 
\item For $(((N_I,F,W),u_{IJ}),(P,W),\beta)\in C(MHM_{gm,k,\mathbb C_p}(X))$, 
\begin{eqnarray*}
\ad(Rf^{Hdg}_!,f^{!Hdg})(((N_I,F,W),u_{IJ}),(P,W),\beta):= \\
(\ad(f_{Hdg}^{*mod},Rf^{Hdg}_!)((N_I,F,W),u_{IJ}),\ad(f^{!w},Rf_{!w})(P,W)): \\
(((N_I,F,W),u_{IJ}),(P,W),\beta)\to f^{!Hdg}Rf^{Hdg}_!(((N_I,F,W),u_{IJ}),(P,W),\beta) 
\end{eqnarray*}
is the adjonction map in $D(MHM_{gm,k,\mathbb C_p}(X))$. 
\end{itemize}
\end{itemize}
\end{prop}

\begin{proof}
Follows from proposition \ref{jHdgpropadCp} after considering a factorization 
$f:X\hookrightarrow\bar Y\times S\xrightarrow{p_S} S$ with $\bar Y\in\PSmVar(k)$.
\end{proof}

\begin{thm}\label{sixMHMkCp}
Let $k\subset\mathbb C_p$ a subfield. 
\begin{itemize}
\item[(i)]We have the six functor formalism on $D(MHM_{gm,k,\mathbb C_p}(-)):\SmVar(k)\to\TriCat$.
\item[(ii)]We have the six functor formalism on $D(MHM_{gm,k,\mathbb C_p}(-)):\QPVar(k)\to\TriCat$.
\end{itemize}
\end{thm}

\begin{proof}
Follows from proposition \ref{HdgpropadsingCp}.
\end{proof}

We give the following version (where the De Rham cohomology is twisted by the p-adic periods)
of the syntomic complex of a p adic analytic space and 
the syntomic cohomology class of an algebraic cycle of a p adic algebraic variety.
\begin{defi}\label{Delkpdef}
\begin{itemize}
\item[(i)] Let $K$ a p adic field. Let $X\in\AnSm(K)$. We have for $d\in\mathbb Z$ the syntomic complex
\begin{equation*}
\mathbb Z_{syn,X}(d):=(\mathbb Z_{p,X}(d)\hookrightarrow 
DR(X)(O\mathbb B_{dr,X})/F_b^d:=(\Omega_X^{\bullet,\leq d}\otimes_{O_X}O\mathbb B_{dr,X}))
\in C(X^{et})
\end{equation*}
Let $D\subset X$ a normal crossing divisor. We have for $d\in\mathbb Z$ the Deligne complexes
\begin{equation*}
\mathbb Z_{syn,(X,D)}(d):=(\mathbb Z_X(d)\hookrightarrow DR(X)(O\mathbb B_{dr,X}(\log D))/F_b^d)
:=(\Omega_X^{\bullet,\leq d}\otimes_{O_X}O\mathbb B_{dr,X}(\log D))\in C(X^{et})
\end{equation*}
and
\begin{equation*}
\mathbb Z_{syn,(X,D)}(d)^{\vee}:=(\mathbb Z_X(d)\hookrightarrow DR(X)(O\mathbb B_{dr,X}(\nul D))/F_b^d)
:=(\Omega_X^{\bullet,\leq d}\otimes_{O_X}O\mathbb B_{dr,X}(\nul D))\in C(X^{et}).
\end{equation*}
Moreover we have as for Deligne complexes canonical products 
\begin{itemize}
\item $(-)\cdot(-):\mathbb Z_{syn,(X,D)}(d)\otimes\mathbb Z_{syn,(X,D)}(d')\to\mathbb Z_{\mathcal D,(X,D)}(d+d')$
\item $(-)\cdot(-):\mathbb Z_{syn,(X,D)}(d)^{\vee}\otimes\mathbb Z_{syn,(X,D)}(d')^{\vee}\to\mathbb Z_{syn,(X,D)}(d+d')^{\vee}$
\end{itemize}
\item[(ii)]Let $K$ a p adic field. Let $X\in\AnSm(K)$. We have for $d\in\mathbb Z$ the syntomic (cohomology) complex
\begin{eqnarray*}
C^{\bullet}_{syn}(X,\mathbb Z(d)):=
\Cone(\Gamma(X,E_{et}(\mathbb Z_{p,X}))\oplus\Gamma(X,F^dE_{et}(DR(X)(O\mathbb B_{dr,X}),F_b)) \\
\hookrightarrow\Gamma(X,E_{et}(DR(X)(O\mathbb B_{dr,X}))))\in C(\mathbb Z_p)
\end{eqnarray*}
Let $D\subset X$ a normal crossing divisor. Denote $U:=X\backslash D$. 
We have for $d\in\mathbb Z$ the syntomic (cohomology) complexes
\begin{eqnarray*}
C^{\bullet}_{syn}((X,D),\mathbb Z(d)):=
\Cone(\Gamma(X,E_{et}(\mathbb Z_{p,X}))\oplus\Gamma(X,F^dE_{et}(DR(X)(O\mathbb B_{dr,X}(\log D)),F_b)) \\
\hookrightarrow\Gamma(X,E_{et}(DR(X)(O\mathbb B_{dr,X}(\log D)))))\in C(\mathbb Z_p)
\end{eqnarray*}
and
\begin{eqnarray*}
C^{\bullet}_{syn}(X,D,\mathbb Z(d)):=
\Cone(\Gamma(X,E_{et}(\mathbb Z_{p,X}))\oplus\Gamma(X,F^dE_{et}(DR(X)(O\mathbb B_{dr,X}(\nul D)),F_b)) \\
\hookrightarrow\Gamma(X,E_{et}(DR(X)(O\mathbb B_{dr,X}(\nul D))))\in C(\mathbb Z_p).
\end{eqnarray*}
\item[(iii)] Let $k\subset K$ an embedding of a field of characteristic zero into a p adic field.
Let $X\in\PSmVar(k)$. We have, for $k\in\mathbb Z$ and $d\in\mathbb Z$, the syntomic cohomology
\begin{eqnarray*}
H_{syn}^k(X_K^{an},\mathbb Z(d)):=\mathbb H^k(X_K^{an},\mathbb Z_{X,syn}(d))=H^kC^{\bullet}_{syn}(X_K^{an},D,\mathbb Z(d))
\end{eqnarray*}
Let $U\in\SmVar(k)$. Let $X\in\PSmVar(k)$ a compactification of $U$ with $D:=X\backslash U$ a normal crossing divisor.
We have, for $k\in\mathbb Z$ and $d\in\mathbb Z$, the syntomic cohomology
\begin{eqnarray*}
H_{syn}^k(U_K^{an},\mathbb Z(d)):=\mathbb H^k(X,\mathbb Z_{(X_K^{an},D_K^{an}),syn}(d))
=H^kC^{\bullet}_{syn}((X_K^{an},D_K^{an}),\mathbb Z(d))
\end{eqnarray*}
and
\begin{eqnarray*}
H_{syn}^k(X,D,\mathbb Z(d)):=\mathbb H^k(X_K^{an},\mathbb Z_{(X_K^{an},D_K^{an}),syn}(d)^{\vee})
=H^kC^{\bullet}_{syn}(X_K^{an},D_K^{an},\mathbb Z(d)).
\end{eqnarray*}
\item[(iv)] Let $k\subset K\subset\mathbb C_p$ an embedding of a field of characteristic zero into a p adic field.
Let $U\in\SmVar(k)$. Let $X\in\PSmVar(k)$ a compactification of $U$ with $D:=X\backslash U$ a normal crossing divisor.
We define the Deligne cohomology of a (higher) cycle $Z\in\mathcal Z^d(U,n)^{\partial=0}$ by
\begin{eqnarray*}
[Z]_{syn}:=\Im(H^{2d-n}(\gamma_{\supp(Z)})([Z])), \\ 
H^k(\gamma_{\supp(Z)}):
\mathbb H^{2d-n}_{syn,\supp(Z)}(X_{\mathbb C_p}^{an},\mathbb Z_{X_{\mathbb C_p}^{an},D_{\mathbb C_p}^{an}}(d))
\to\mathbb H^{2d-n}_{syn}(X_{\mathbb C_p}^{an},\mathbb Z_{X_{\mathbb C_p}^{an},D_{\mathbb C_p}^{an}}(d)) 
\end{eqnarray*}
with $\supp(Z):=p_X(\supp(Z))\subset X$, where $\supp(Z)\subset X\times\square^n$ is the support of $Z$.
\item[(v)]Let $k\subset K$ an embedding of a field of characteristic zero into a p adic field.
Let $U\in\SmVar(k)$. Let $X\in\PSmVar(k)$ a compactification of $U$ with $D:=X\backslash U$ a normal crossing divisor.
We have for $d\in\mathbb Z$ the morphism of complexes
\begin{eqnarray*}
\mathcal R^d_U:\mathcal Z^d(U,\bullet)\to C^{\bullet}_{syn}(X_{\mathbb C_p}^{an},D_{\mathbb C_p}^{an},\mathbb Z(d)), \; 
Z\mapsto\mathcal R^d_U(Z):=(T_{\bar Z},\Omega_{\bar Z},R_{\bar Z})
\end{eqnarray*}
which gives for $Z\in\mathcal Z^d(U,n)^{\partial=0}$, 
\begin{equation*}
[\mathcal R^d_U(Z)]=[Z]_{syn}\in H_{syn}^{2d-n}(U_{\mathbb C_p}^{an},\mathbb Z(d))
\end{equation*}
\end{itemize}
\end{defi}

Let $K$ a p adic field. Let $f:X\to S$ a morphism with $S,X\in\AnSm(K)$.
We have for $d\in\mathbb Z$ the canonical morphism of Deligne complexes
\begin{equation*}
(\ad(f^*,f_*)(\mathbb Z_{p,S}),\Omega_{X/S}^{\leq d}):\mathbb Z_{syn,S}(d)\to f_*\mathbb Z_{syn,X}(d)
\end{equation*}
which induces after taking the canonical flasque resolution of the syntomic complexes 
the morphism in $C(\mathbb Z_p)$
\begin{eqnarray*}
f^*:=(f^*,f^*,\theta(f)^t):C^{\bullet}_{syn}(S,\mathbb Z(d)):= \\
\Cone(\Gamma(S,E_{et}(\mathbb Z_{p,S}))\oplus\Gamma(S,F^dE_{et}(DR(S)(O\mathbb B_{dr,S})))
\hookrightarrow\Gamma(S,E_{et}(DR(S)(O\mathbb B_{dr,S})))) \\
\to C^{\bullet}_{syn}(X,\mathbb Z(d)):= \\
\Cone(\Gamma(X,E_{et}(\mathbb Z_{p,X}))\oplus\Gamma(X,F^dE_{et}(DR(X)(O\mathbb B_{dr,X})))
\hookrightarrow\Gamma(X,E_{et}(DR(X)(O\mathbb B_{dr,X})))) 
\end{eqnarray*}
where $\theta(f)^t$ is the homotopy in the morphism in $D_{fil}(k)\otimes_ID(\mathbb Z_p)$
(where here the comparaison morphisms $\alpha$ are in $D_{\mathbb B_{dr},G,fil}(K)$
instead of $D_{\mathbb B_{dr},G,fil}(\mathbb C_p)$)
\begin{eqnarray*}
(f^*,f^*,\theta(f)^t):
(\Gamma(S,E_{et}(DR(S)(O\mathbb B_{dr,S}),F_b)),\Gamma(S,E_{et}(\mathbb Z_{p,S})),\alpha(S)) \\
\to(\Gamma(X,E_{et}(DR(X)(O\mathbb B_{dr,X}),F_b)),\Gamma(X,E_{et}(\mathbb Z_{p,X})),\alpha(X)),
\end{eqnarray*}
which induces in cohomology for $n\in\mathbb Z$, the morphisms of abelian groups
\begin{equation*}
f^*:H^n_{syn}(S,\mathbb Z(d))\to H^n_{syn}(X,\mathbb Z(d)) \; ;
\end{equation*}
we get dually,
\begin{eqnarray*}
f_*:=(f_*,f_*,\theta(f)): \\
\Cone(\Gamma(X,E_{et}(\mathbb Z_{p,X}))^{\vee}\oplus F^d\Gamma(X,E_{et}(DR(X)(O\mathbb B_{dr,X}),F_b))^{\vee}
\hookrightarrow\Gamma(X,E_{et}(DR(X)(O\mathbb B_{dr,X})))^{\vee}) \\
\to\Cone(\Gamma(S,E_{et}(\mathbb Z_{p,S}))^{\vee}\oplus F^d\Gamma(S,E_{et}(DR(S)(O\mathbb B_{dr,S}),F_b))^{\vee}
\hookrightarrow\Gamma(S,E_{et}(DR(S)(O\mathbb B_{dr,S})))^{\vee}) 
\end{eqnarray*}
where $\theta(f)$ is the homotopy in the morphism in $D_{fil}(k)\otimes_ID(\mathbb Z_p)$
(where here the comparaison morphisms $\alpha$ are in $D_{\mathbb B_{dr},G,fil}(K)$
instead of $D_{\mathbb B_{dr},G,fil}(\mathbb C_p)$)
\begin{eqnarray*}
(f_*,f_*,\theta(f)):
(\Gamma(X,E_{et}(DR(X)(O\mathbb B_{dr,X}),F_b))^{\vee},\Gamma(X,E_{et}(\mathbb Z_{p,X}))^{\vee},\alpha(X)) \\
\to(\Gamma(S,E_{et}(DR(S)(O\mathbb B_{dr,S}),F_b))^{\vee},\Gamma(S,E_{et}(\mathbb Z_{p,S}))^{\vee},\alpha(S)),
\end{eqnarray*}
which induces in homology for $n\in\mathbb Z$, the morphisms of abelian groups
\begin{equation*}
f_*:H_{n,syn}(X,\mathbb Z(d))\to H_{n,syn}(S,\mathbb Z(d)).
\end{equation*}

\begin{thm}\label{Delkp}
Let $k\subset\mathbb C_p$ a subfield.
\begin{itemize}
\item[(i)] Let $U\in\SmVar(k)$. Denote by $a_U:U\to\pt$ the terminal map. 
Let $X\in\PSmVar(k)$ a compactification of $U$ with $D:=X\backslash U$ a normal crossing divisor.
The embedding (see theorem \ref{Bekp})
\begin{equation*}
\iota:D(MHM_{gm,k,\mathbb C_p}(\left\{\pt\right\}))\to D_{fil}(k)\times_ID(\mathbb Z_p) 
\end{equation*}
induces for $k\in\mathbb Z$ and $d\in\mathbb Z$, canonical isomorphisms
\begin{eqnarray*}
\iota(a_{U!Hdg}\mathbb Z^{Hdg}_U):H^k(a_{U!Hdg}\mathbb Z^{Hdg}_U)\xrightarrow{\sim}
H^k_{syn}(X_{\mathbb C_p}^{an},D_{\mathbb C_p}^{an},\mathbb Z(d)), 
\; \mbox{and} \\
\iota(a_{U*Hdg}\mathbb Z^{Hdg}_U):H^k(a_{U*Hdg}\mathbb Z^{Hdg}_U)\xrightarrow{\sim}
H^k_{syn}(U_{\mathbb C_p}^{an},\mathbb Z(d)).
\end{eqnarray*}
\item[(ii)] Let $h:U\to S$ and $h':U'\to S$ two morphism with $S,U,U'\in\SmVar(k)$.
Let $X\in\PSmVar(k)$ a compactification of $U$ with $D:=X\backslash U$ a normal crossing divisor
such that $h:U\to S$ extend to $f:X\to\bar S$.
Let $X'\in\PSmVar(k)$ a compactification of $U'$ with $D':=X'\backslash U'$ a normal crossing divisor
such that $h':U'\to S$ extend to $f':X'\to\bar S$.
The embedding $\iota:D(MHM_{gm,k,\mathbb C_p}(\pt))\to D_{fil}(k)\times_ID(\mathbb Z_p)$ (see theorem \ref{Bekp}) 
induces for $k\in\mathbb Z$ and $d\in\mathbb Z$ a canonical isomorphism
\begin{eqnarray*}
\iota(a_{U'\times_SU!Hdg}\mathbb Z^{Hdg}_{U'\times_SU}):
\Hom_{D(MHM_{gm,k,\mathbb C_p}(S))}(h_{U'!Hdg}\mathbb Z^{Hdg}_{U'},h_{U!Hdg}\mathbb Z^{Hdg}_U(d)[k]) \\
\xrightarrow{RI(-,-)}
\Hom_{D(MHM_{gm,k,\mathbb C_p}(\pt))}(\mathbb Z^{Hdg}_{\pt},a_{U'\times_SU!Hdg}\mathbb Z^{Hdg}_{U'\times_SU}(d)[k])
=H^k(a_{U'\times_SU!Hdg}\mathbb Z^{Hdg}_{U'\times_SU}(d)) \\
\xrightarrow{\sim}
H^k_{\mathcal D}((X'\times_SX)_{\mathbb C_p}^{an},((X'\times_SU)\cup(U'\times_SX))_{\mathbb C_p}^{an},\mathbb Z(d)).
\end{eqnarray*}
\item[(iii)]Let $U\in\SmVar(k)$.  
Let $X\in\PSmVar(k)$ a compactification of $U$ with $D:=X\backslash U$ a normal crossing divisor.
For $[Z]\in\CH^d(U,n)$ and $[Z']\in\CH^{d'}(U,n')$, we have
\begin{equation*}
([Z]\cdot[Z'])_{syn}=[Z]_{syn}\cdot[Z']_{syn}\in H^{2d+2d'-n-n'}(U_{\mathbb C_p}^{an},\mathbb Z(d+d'))
\end{equation*}
where the product on the left is the intersection of higher Chow cycle which is well defined modulo boundary 
(they intersect properly modulo boundary) while the right product of Deligne cohomology classes is induced by 
the product of Deligne complexes 
$(-)\cdot(-):\mathbb Z_{syn,(X,D)}(d)\otimes\mathbb Z_{syn,(X,D)}(d')\to\mathbb Z_{syn,(X,D)}(d+d')$.
\item[(iv)]Let $h:U\to S$,$h':U'\to S$, $h'':U''\to S$ three morphism with $S,U,U',U''\in\SmVar(k)$.
Let $X\in\PSmVar(k)$ a compactification of $U$ with $D:=X\backslash U$ a normal crossing divisor
such that $h:U\to S$ extend to $f:X\to\bar S$.
Let $X'\in\PSmVar(k)$ a compactification of $U'$ with $D':=X'\backslash U'$ a normal crossing divisor
such that $h':U'\to S$ extend to $f':X'\to\bar S$.
Let $X'\in\PSmVar(k)$ a compactification of $U'$ with $D':=X'\backslash U'$ a normal crossing divisor
such that $h':U'\to S$ extend to $f':X'\to\bar S$.
For $[Z]\in\CH^d(U\times_SU',n)$ and $[Z']\in\CH^{d'}(U'\times_SU'',n')$, we have
\begin{equation*}
([Z]\circ[Z'])_{syn}=[Z]_{syn}\circ[Z']_{syn}
\in H^{d''-n''}((U\times_SU'')_{\mathbb C_p}^{an},\mathbb Z(d''-n''))
\end{equation*}
where the composition on the left is the composition of higher correspondence modulo boundary
while the composition on the right is given by (ii).
\end{itemize}
\end{thm}

\begin{proof}
\noindent(i):Standard.

\noindent(ii):Follows on the one hand from (i) and 
on the other hand the six functor formalism on the 2-functor 
$D(MHM_{gm,k,\mathbb C_p}(-)):\SmVar(k)\to\TriCat$ (theorem \ref{sixMHMkC}) gives the isomorphism $RI(-,-)$.

\noindent(iii):Standard.

\noindent(iv):Follows from (iii).

\end{proof}

\section{The algebraic filtered De Rham realizations for Voevodsky relative motives over a field $k$ of characteristic $0$} 

\subsection{The algebraic Gauss-Manin filtered De Rham realization functor}

Let $k$ a field of characteristic zero.
Consider, for $S\in\Var(k)$, the following composition of morphism in $\RCat$ (see section 2)
\begin{eqnarray*}
\tilde e(S):(\Var(k)/S,O_{\Var(k)/S})\xrightarrow{\rho_S}(\Var(k)^{sm}/S,O_{\Var(k)^{sm}/S})
\xrightarrow{e(S)}(S,O_S)
\end{eqnarray*}
with, for $X/S=(X,h)\in\Var(k)/S$,
\begin{itemize}
\item $O_{\Var(k)/S}(X/S):=O_X(X)$, 
\item $(\tilde e(S)^*O_S(X/S)\to O_{\Var(k)/S}(X/S)):=(h^*O_S\to O_X)$.
\end{itemize}
and $O_{\Var(k)^{sm}/S}:=\rho_{S*}O_{\Var(k)/S}$, that is, 
for $U/S=(U,h)\in\Var(k)^{sm}/S$, $O_{\Var(k)^{sm}/S}(U/S):=O_{\Var(k)/S}(U/S):=O_U(U)$

\begin{defi}\label{wtildewGM}
\begin{itemize}
\item[(i)]For $S\in\Var(k)$, we consider the complexes of presheaves 
\begin{equation*} 
\Omega^{\bullet}_{/S}:=
\coker(\Omega_{O_{\Var(k)/S}/\tilde e(S)^*O_S}:
\Omega^{\bullet}_{\tilde e(S)^*O_S}\to\Omega^{\bullet}_{O_{\Var(k)/S}})
\in C_{O_S}(\Var(k)/S) 
\end{equation*}
which is by definition given by 
\begin{itemize}
\item for $X/S$ a morphism $\Omega^{\bullet}_{/S}(X/S)=\Omega^{\bullet}_{X/S}(X)$
\item for $g:X'/S\to X/S$ a morphism, 
\begin{eqnarray*}
\Omega^{\bullet}_{/S}(g):=\Omega_{(X'/X)/(S/S)}(X'):
\Omega^{\bullet}_{X/S}(X)\to g^*\Omega_{X/S}(X')\to\Omega^{\bullet}_{X'/S}(X') \\
\omega\mapsto\Omega_{(X'/X)/(S/S)}(X')(\omega):=g^*(\omega):(\alpha\in\wedge^{k}T_{X'}(X')\mapsto\omega(dg(\alpha)))
\end{eqnarray*}
\end{itemize}
\item[(ii)] For $S\in\Var(k)$, we consider the complexes of presheaves 
\begin{equation*} 
\Omega^{\bullet}_{/S}:=\rho_{S*}\tilde\Omega^{\bullet}_{/S}=
\coker(\Omega_{O_{\Var(k)^{sm}/S}/e(S)^*O_S}:\Omega^{\bullet}_{e(S)^*O_S}\to\Omega^{\bullet}_{O_{\Var(k)^{sm}/S}})
\in C_{O_S}(\Var(k)^{sm}/S) 
\end{equation*}
which is by definition given by 
\begin{itemize}
\item for $U/S$ a smooth morphism $\Omega^{\bullet}_{/S}(U/S)=\Omega^{\bullet}_{U/S}(U)$
\item for $g:U'/S\to U/S$ a morphism, 
\begin{eqnarray*}
\Omega^{\bullet}_{/S}(g):=\Omega_{(U'/U)/(S/S)}(U'):
\Omega^{\bullet}_{U/S}(U)\to g^*\Omega_{U/S}(U')\to\Omega^{\bullet}_{U'/S}(U') \\
\omega\mapsto\Omega_{(U'/U)/(S/S)}(U')(\omega):=g^*(\omega):(\alpha\in\wedge^{k}T_{U'}(U')\mapsto\omega(dg(\alpha)))
\end{eqnarray*}
\end{itemize}
\end{itemize}
\end{defi}

\begin{rem}
For $S\in\Var(k)$, $\Omega^{\bullet}_{/S}\in C(\Var(k)/S)$ 
is by definition a natural extension of $\Omega^{\bullet}_{/S}\in C(\Var(k)^{sm}/S)$. 
However $\Omega^{\bullet}_{/S}\in C(\Var(k)/S)$ does NOT satisfy cdh descent.
\end{rem}

For a smooth morphism $h:U\to S$ with $S,U\in\SmVar(\mathbb C)$, 
the cohomology presheaves  $H^n\Omega^{\bullet}_{U/S}$ of the relative De Rham complex
\begin{equation*}
DR(U/S):=\Omega^{\bullet}_{U/S}:=\coker(h^*\Omega_S\to\Omega_U)\in C_{h^*O_S}(U)
\end{equation*}
for all $n\in\mathbb Z$, have a canonical structure of a complex of $h^*D_S$ modules given by the Gauss Manin connexion : 
for $S^o\subset S$ an open subset, $U^o=h^{-1}(S^o)$, 
$\gamma\in\Gamma(S^o,T_S)$ a vector field and $\hat\omega\in\Omega^{p}_{U/S}(U^o)^c$ a closed form, the action is given by
\begin{equation*}
\gamma\cdot[\hat\omega]=[\widehat{\iota(\tilde\gamma)\partial\omega}], 
\end{equation*}
$\omega\in\Omega^p_U(U^o)$ being a representative of $\hat\omega$ and 
$\tilde\gamma\in\Gamma(U^o,T_U)$ a relevement of $\gamma$ ($h$ is a smooth morphism), so that 
\begin{equation*}
DR(U/S):=\Omega^{\bullet}_{U/S}:=\coker(h^*\Omega_S\to\Omega_U)\in C_{h^*O_S,h^*\mathcal D}(U)
\end{equation*}
with this $h^*D_S$ structure. Hence we get $h_*\Omega^{\bullet}_{U/S}\in C_{O_S,\mathcal D}(S)$ considering this structure.
Since $h$ is a smooth morphism, $\Omega^p_{U/S}$ are locally free $O_U$ modules.

The point (ii) of the definition \ref{wtildew} above gives 
the object in $\DA(S)$ which will, for $S$ smooth, represent the algebraic Gauss-Manin De Rham realisation.
It is the class of an explicit complex of presheaves on $\Var(k)^{sm}/S$. 

\begin{prop}\label{aetfibGM}
Let $S\in\Var(k)$.
\begin{itemize}
\item[(i)] For $U/S=(U,h)\in\Var(k)^{sm}/S$, we have $e(U)_*h^*\Omega^{\bullet}_{/S}=\Omega^{\bullet}_{U/S}$.
\item[(ii)] The complex of presheaves $\Omega^{\bullet}_{/S}\in C_{O_S}(\Var(k)^{sm}/S)$ 
is $\mathbb A^1$ homotopic, in particular $\mathbb A^1$ invariant.
Note that however, for $p>0$, the complexes of presheaves $\Omega^{\bullet\geq p}$ are NOT $\mathbb A^1$ local.
On the other hand, $(\Omega^{\bullet}_{/S},F_b)$ admits transferts 
(recall that means $\Tr(S)_*\Tr(S)^*\Omega^p_{/S}=\Omega^p_{/S}$).
\item[(iii)] If $S$ is smooth, we get $(\Omega^{\bullet}_{/S},F_b)\in C_{O_Sfil,D_S}(\Var(k)^{sm}/S)$
with the structure given by the Gauss Manin connexion. Note that however the $D_S$ structure on the cohomology groups
given by Gauss Main connexion does NOT comes from a structure of $D_S$ module structure on the filtered complex of $O_S$ module.
The $D_S$ structure on the cohomology groups satisfy a non trivial Griffitz transversality (in the non projection cases),
whereas the filtration on the complex is the trivial one.
\end{itemize}
\end{prop}

\begin{proof} 
Similar to the proof of \cite{B4} proposition.
\end{proof}

We have the following canonical transformation map given by the pullback of (relative) differential forms:

Let $g:T\to S$ a morphism with $T,S\in\Var(k)$.
Consider the following commutative diagram in $\RCat$ :
\begin{equation*}
D(g,e):\xymatrix{
(\Var(k)^{sm}/T,O_{\Var(k)^{sm}/T})\ar[rr]^{P(g)}\ar[d]^{e(T)} & \, & 
(\Var(k)^{sm}/S,O_{\Var(k)^{sm}/S})\ar[d]^{e(S)} \\
(T,O_T)\ar[rr]^{P(g)} & \, & (S,O_S)}
\end{equation*}
It gives (see section 2) the canonical morphism in $C_{g^*O_Sfil}(\Var(k)^{sm}/T)$ 
\begin{eqnarray*}
\Omega_{/(T/S)}:=\Omega_{(O_{\Var(k)^{sm}/T}/g^*O_{\Var(k)^{sm}/S})/(O_T/g^*O_S}): \\
g^*(\Omega^{\bullet}_{/S},F_b)=\Omega^{\bullet}_{g^*O_{\Var(k)^{sm}/S}/g^*e(S)^*O_S}\to
(\Omega^{\bullet}_{/T},F_b)=\Omega^{\bullet}_{O_{\Var(k)^{sm}/T}/e(T)^*O_T}
\end{eqnarray*}
which is by definition given by the pullback on differential forms : for $(V/T)=(V,h)\in\Var(k)^{sm}/T$,
\begin{eqnarray*}
\Omega_{/(T/S)}(V/T): 
g^*(\Omega^{\bullet}_{/S})(V/T):=\lim_{(h':U\to S \mbox{sm},g':V\to U,h,g)}\Omega^{\bullet}_{U/S}(U)
\xrightarrow{\Omega_{(V/U)/(T/S)}(V/T)}\Omega^{\bullet}_{V/T}(V)=:\Omega^{\bullet}_{/T}(V/T) \\
\hat\omega\mapsto\Omega_{(V/U)/(T/S)}(V/T)(\omega):=\hat{g^{'*}\omega}.
\end{eqnarray*}
If $S$ and $T$ are smooth, $\Omega_{/(T/S)}:g^*(\Omega^{\bullet}_{/S},F_b)\to(\Omega^{\bullet}_{/T},F_b)$
is a map in $C_{g^*O_Sfil,g^*D_S}(\Var(k)^{sm}/T)$
It induces the canonical morphisms in $C_{g^*O_Sfil,g^*D_S}(\Var(k)^{sm}/T)$:
\begin{eqnarray*}
E\Omega_{/(T/S)}:g^*E_{et}(\Omega^{\bullet}_{/S},F_b)\xrightarrow{T(g,E_{et})(\Omega^{\bullet}_{/S},F_b)}
E_{et}(g^*(\Omega^{\bullet}_{/S},F_b))\xrightarrow{E_{et}(\Omega_{/(T/S)})}E_{et}(\Omega^{\bullet}_{/T},F_b). 
\end{eqnarray*}
and
\begin{eqnarray*}
E\Omega_{/(T/S)}:g^*E_{zar}(\Omega^{\bullet}_{/S},F_b)\xrightarrow{T(g,E_{zar})(\Omega^{\bullet}_{/S},F_b)}
E_{zar}(g^*(\Omega^{\bullet}_{/S},F_b))\xrightarrow{E_{zar}(\Omega_{/(T/S)})}E_{zar}(\Omega^{\bullet}_{/T},F_b). 
\end{eqnarray*}

\begin{defi}\label{TgDRGM}
\begin{itemize}
\item[(i)]Let $g:T\to S$ a morphism with $T,S\in\Var(k)$.
We have, for $F\in C(\Var(k)^{sm}/S)$, the canonical transformation in $C_{O_Tfil}(T)$ :
\begin{eqnarray*}
T^O(g,\Omega_{/\cdot})(F): 
g^{*mod}L_Oe(S)_*\mathcal Hom^{\bullet}(F,E_{et}(\Omega^{\bullet}_{/S},F_b)) \\
\xrightarrow{:=}
(g^*L_Oe(S)_*\mathcal Hom^{\bullet}(F,E_{et}(\Omega^{\bullet}_{/S},F_b)))\otimes_{g^*O_S}O_T \\
\xrightarrow{T(e,g)(-)\circ T(g,L_O)(-)}  
L_O(e(T)_*g^*\mathcal Hom^{\bullet}(F,E_{et}(\Omega^{\bullet}_{/S},F))\otimes_{g^*O_S}O_T) \\ 
\xrightarrow{T(g,hom)(F,E_{et}(\Omega^{\bullet}_{/S}))\otimes I} 
L_O(e(T)_*\mathcal Hom^{\bullet}(g^*F,g^*E_{et}(\Omega^{\bullet}_{/S},F_b))\otimes_{g^*O_S}O_T) \\
\xrightarrow{ev(hom,\otimes)(-,-,-)} 
L_Oe(T)_*\mathcal Hom^{\bullet}(g^*F,g^*E_{et}(\Omega^{\bullet}_{/S},F_b)\otimes_{g^*e(S)^*O_S}e(T)^*O_T) \\
\xrightarrow{\mathcal Hom^{\bullet}(g^*F,E\Omega_{/(T/S)}\otimes I)} 
L_Oe(T)_*\mathcal Hom^{\bullet}(g^*F,E_{et}(\Omega^{\bullet}_{/T},F_b)\otimes_{g^*e(S)^*O_S}e(T)^*O_T) \\
\xrightarrow{m}
L_Oe(T)_*\mathcal Hom^{\bullet}(g^*F,E_{et}(\Omega^{\bullet}_{/T},F_b)
\end{eqnarray*}
where $m(\alpha\otimes h):=h.\alpha$ is the multiplication map.
\item[(ii)] Let $g:T\to S$ a morphism with $T,S\in\Var(k)$, $S$ smooth.
Assume there is a factorization $g:T\xrightarrow{l}Y\times S\xrightarrow{p_S}S$
with $Y\in\SmVar(\mathbb C)$, $l$ a closed embedding and $p_S$ the projection.
We have, for $F\in C(\Var(k)^{sm}/S)$, the canonical transformation in $C_{O_Tfil}(Y\times S)$ :
\begin{eqnarray*}
T(g,\Omega_{/\cdot})(F):
g^{*mod,\Gamma}e(S)_*\mathcal Hom^{\bullet}(F,E_{et}(\Omega^{\bullet}_{/S},F_b)) \\
\xrightarrow{:=} 
\Gamma_TE_{zar}(p_S^{*mod}e(S)_*\mathcal Hom^{\bullet}(F,E_{et}(\Omega^{\bullet}_{/S},F_b))) \\
\xrightarrow{T^O(p_S,\Omega_{/\cdot})(F)} 
\Gamma_TE_{zar}(e(T\times S)_*\mathcal Hom^{\bullet}(p_S^*F,E_{et}(\Omega^{\bullet}_{/Y\times S},F_b))) \\
\xrightarrow{=}
e(T\times S)_*\Gamma_T(\mathcal Hom^{\bullet}(p_S^*F,E_{et}(\Omega^{\bullet}_{/Y\times S},F_b))) \\
\xrightarrow{I(\gamma,\hom)(-,-)}
e(T\times S)_*\mathcal Hom^{\bullet}(\Gamma^{\vee}_Tp_S^*F,E_{et}(\Omega^{\bullet}_{/Y\times S},F_b)).
\end{eqnarray*}
For $Q\in Proj\PSh(\Var(k)^{sm}/S)$,
\begin{eqnarray*}
T(g,\Omega_{/\cdot})(Q):
g^{*mod,\Gamma}e(S)_*\mathcal Hom^{\bullet}(Q,E_{et}(\Omega^{\bullet}_{/S},F_b))\to
e(T\times S)_*\mathcal Hom^{\bullet}(\Gamma^{\vee}_Tp_S^*Q,E_{et}(\Omega^{\bullet}_{/Y\times S},F_b))
\end{eqnarray*}
is a map in $C_{O_Tfil,\mathcal D}(Y\times S)$.
\end{itemize}
\end{defi}

Let $S\in\Var(k)$. We have the canonical map in $C_{O_Sfil}(\Var(k)^{sm}/S)$
\begin{eqnarray*}
w_S:(\Omega^{\bullet}_{/S},F_b)\otimes_{O_S}(\Omega^{\bullet}_{/S},F_b)\to(\Omega^{\bullet}_{/S},F_b)
\end{eqnarray*}
given by for $h:U\to S\in\Var(k)^{sm}/S$ by the wedge product
\begin{eqnarray*}
w_S(U/S):(\Omega^{\bullet}_{U/S},F_b)\otimes_{h^*O_S}(\Omega^{\bullet}_{U/S},F_b)(U)\xrightarrow{w_{U/S}(U)}
(\Omega^{\bullet}_{U/S},F_b)(U)
\end{eqnarray*}
It gives the map
\begin{eqnarray*}
Ew_S:E_{et}(\Omega^{\bullet}_{/S},F_b)\otimes_{O_S}E_{et}(\Omega^{\bullet}_{/S},F_b)\xrightarrow{=}
E_{et}((\Omega^{\bullet}_{/S},F_b)\otimes_{O_S}(\Omega^{\bullet}_{/S},F_b))\xrightarrow{E_{et}(w_S)}
E_{et}(\Omega^{\bullet}_{/S},F_b)
\end{eqnarray*}
If $S\in\SmVar(\mathbb C)$, 
\begin{eqnarray*}
w_S:(\Omega^{\bullet}_{/S},F_b)\otimes_{O_S}(\Omega^{\bullet}_{/S},F_b)\to(\Omega^{\bullet}_{/S},F_b)
\end{eqnarray*}
is a map in $C_{O_Sfil,D_S}(\Var(k)^{sm}/S)$.

\begin{defi}\label{TotimesDRGM}
Let $S\in\Var(k)$.
We have, for $F,G\in C(\Var(k)^{sm}/S)$, the canonical transformation in $C_{O_Sfil}(S)$ :
\begin{eqnarray*}
T(\otimes,\Omega)(F,G):
e(S)_*\mathcal Hom(F,E_{et}(\Omega^{\bullet}_{/S},F_b))\otimes_{O_S}e(S)_*\mathcal Hom(G,E_{et}(\Omega^{\bullet}_{/S},F_b)) \\ 
\xrightarrow{=}
e(S)_*(\mathcal Hom(F,E_{et}(\Omega^{\bullet}_{/S},F_b))\otimes_{O_S}\mathcal Hom(G,E_{et}(\Omega^{\bullet}_{/S},F_b))) \\
\xrightarrow{e(S)_*T(\mathcal Hom,\otimes)(-)}
e(S)_*\mathcal Hom(F\otimes G,E_{et}(\Omega^{\bullet}_{/S},F_b)\otimes_{O_S}E_{et}(\Omega^{\bullet}_{/S},F_b)) \\
\xrightarrow{\mathcal Hom(F\otimes G,Ew_S)}
e(S)_*\mathcal Hom(F\otimes G,E_{et}(\Omega^{\bullet}_{/S},F_b))
\end{eqnarray*}
If $S\in\SmVar(\mathbb C)$, $T(\otimes,\Omega)(F,G)$ is a map in $C_{O_Sfil,\mathcal D}(S)$. 
\end{defi}

\begin{defi}\label{DRalgdefFunctGM}
\begin{itemize}
\item[(i)]Let $S\in\SmVar(\mathbb C)$. We have the functor
\begin{eqnarray*}
C(\Var(k)^{sm}/S)^{op}\to C_{Ofil,\mathcal D}(S), \; \;
F\mapsto e(S)_*\mathcal Hom^{\bullet}(L(i_{I*}j_I^*F),E_{et}(\Omega^{\bullet}_{/S},F_b))[-d_S].
\end{eqnarray*}
\item[(ii)]Let $S\in\Var(k)$ and $S=\cup_{i=1}^l S_i$ an open cover such that there exist closed embeddings
$i_i:S_i\hookrightarrow\tilde S_i$  with $\tilde S_i\in\SmVar(\mathbb C)$. 
For $I\subset\left[1,\cdots l\right]$, denote by $S_I:=\cap_{i\in I} S_i$ and $j_I:S_I\hookrightarrow S$ the open embedding.
We then have closed embeddings $i_I:S_I\hookrightarrow\tilde S_I:=\Pi_{i\in I}\tilde S_i$.
We have the functor
\begin{eqnarray*}
C(\Var(k)^{sm}/S)^{op}\to C_{Ofil,\mathcal D}(S/(\tilde S_I)), \; \;
F\mapsto(e(\tilde S_I)_*\mathcal Hom^{\bullet}(L(i_{I*}j_I^*F),
E_{et}(\Omega^{\bullet}_{/\tilde S_I},F_b))[-d_{\tilde S_I}],u^q_{IJ}(F))
\end{eqnarray*}
where
\begin{eqnarray*}
u^q_{IJ}(F)[d_{\tilde S_J}]:
e(\tilde S_I)_*\mathcal Hom^{\bullet}(L(i_{I*}j_I^*F),E_{et}(\Omega^{\bullet}_{/\tilde S_I},F_b)) \\
\xrightarrow{\ad(p_{IJ}^{*mod},p_{IJ*})(-)}
p_{IJ*}p_{IJ}^{*mod}e(\tilde S_I)_*\mathcal Hom^{\bullet}(L(i_{I*}j_I^*F),E_{et}(\Omega^{\bullet}_{/\tilde S_I},F_b)) \\
\xrightarrow{p_{IJ*}T(p_{IJ},\Omega_{\cdot})(L(i_{I*}j_I^*F))}
p_{IJ*}e(\tilde S_J)_*\mathcal Hom^{\bullet}(p_{IJ}^*L(i_{I*}j_I^*F),E_{et}(\Omega^{\bullet}_{/\tilde S_J},F_b)) \\
\xrightarrow{p_{IJ*}e(\tilde S_J)_*\mathcal Hom(S^q(D_{IJ})(F),E_{et}(\Omega_{/\tilde S_J}^{\bullet,\Gamma},F_b))}  
p_{IJ*}e(\tilde S_J)_*\mathcal Hom^{\bullet}(L(i_{J*}j_J^*F),E_{et}(\Omega^{\bullet}_{/\tilde S_J},F_b)). 
\end{eqnarray*}
For $I\subset J\subset K$, we have obviously $p_{IJ*}u_{JK}(F)\circ u_{IJ}(F)=u_{IK}(F)$.
\end{itemize}
\end{defi}

We then have the following key proposition

\begin{prop}\label{projwachGM}
\begin{itemize}
\item[(i)]Let $S\in\Var(k)$. 
Let $m:Q_1\to Q_2$ be an equivalence $(\mathbb A^1,et)$ local in $C(\Var(k)^{sm}/S)$
with $Q_1,Q_2$ complexes of projective presheaves. Then,
\begin{equation*}
e(S)_*\mathcal Hom(m,E_{et}(\Omega^{\bullet}_{/S},F_b)):
e(S)_*\mathcal Hom^{\bullet}(Q_2,E_{et}(\Omega^{\bullet}_{/S},F_b))\to 
e(S)_*\mathcal Hom^{\bullet}(Q_1,E_{et}(\Omega^{\bullet}_{/S},F_b))
\end{equation*}
is an $2$-filtered quasi-isomorphism. 
It is thus an isomorphism in $D_{O_Sfil,\mathcal D,\infty}(S)$ if $S$ is smooth.
\item[(ii)]Let $S\in\Var(k)$. Let $S=\cup_{i=1}^l S_i$ an open cover such that there exist closed embeddings
$i_i:S_i\hookrightarrow\tilde S_i$  with $\tilde S_i\in\SmVar(\mathbb C)$.
Let $m=(m_I):(Q_{1I},s^1_{IJ})\to (Q_{2I},s^2_{IJ})$ be an equivalence $(\mathbb A^1,et)$ local 
in $C(\Var(k)^{sm}/(\tilde S_I)^{op})$ with $Q_{1I},Q_{2I}$ complexes of projective presheaves. Then,
\begin{eqnarray*}
(e(\tilde S_I)_*\mathcal Hom(m_I,E_{et}(\Omega^{\bullet}_{/\tilde S_I},F_b))): \\
(e(\tilde S_I)_*\mathcal Hom^{\bullet}(Q_{2I},E_{et}(\Omega^{\bullet}_{/\tilde S_I},F_b)),u_{IJ}(Q_{2I},s^2_{IJ}))\to 
(e(\tilde S_I)_*\mathcal Hom^{\bullet}(Q_{1I},E_{et}(\Omega^{\bullet}_{/\tilde S_I},F_b)),u_{IJ}(Q_{1I},s^1_{IJ}))
\end{eqnarray*}
is an $2$-filtered quasi-isomorphism. It is thus an isomorphism in $D_{O_Sfil,\mathcal D,\infty}((\tilde S_I))$.
\end{itemize}
\end{prop}

\begin{proof}
Similar to the proof of \cite{B4} proposition.
\end{proof}

\begin{defi}\label{DRalgdefsingGM}
\begin{itemize}
\item[(i)] We define, using definition \ref{DRalgdefFunctGM}, by proposition \ref{projwachGM}, 
the filtered algebraic Gauss-Manin realization functor defined as
\begin{eqnarray*}
\mathcal F_S^{GM}:\DA_c(S)^{op}\to D_{O_Sfil,\mathcal D,\infty}(S), \; \;
M\mapsto\mathcal F_S^{GM}(M):=e(S)_*\mathcal Hom^{\bullet}(L(F),E_{et}(\Omega^{\bullet}_{/S},F_b))[-d_S]
\end{eqnarray*}
where $F\in C(\Var(k)^{sm}/S)$ is such that $M=D(\mathbb A^1,et)(F)$,
\item[(ii)]Let $S\in\Var(k)$ and $S=\cup_{i=1}^l S_i$ an open cover such that there exist closed embeddings
$i_i:S_i\hookrightarrow\tilde S_i$  with $\tilde S_i\in\SmVar(\mathbb C)$. 
For $I\subset\left[1,\cdots l\right]$, denote by $S_I=\cap_{i\in I} S_i$ and $j_I:S_I\hookrightarrow S$ the open embedding.
We then have closed embeddings $i_I:S_I\hookrightarrow\tilde S_I:=\Pi_{i\in I}\tilde S_i$.
We define, using definition \ref{DRalgdefFunctGM} and proposition \ref{projwachGM}
the filtered algebraic Gauss-Manin realization functor defined as
\begin{eqnarray*}
\mathcal F_S^{GM}:\DA_c(S)^{op}\to D_{Ofil,\mathcal D,\infty}(S/(\tilde S_I)), \; M\mapsto \\
\mathcal F_S^{GM}(M):=((e(\tilde S_I)_*\mathcal Hom^{\bullet}(L(i_{I*}j_I^*F),
E_{et}(\Omega^{\bullet}_{/\tilde S_I}),F_b))[-d_{\tilde S_I}],u^q_{IJ}(F))
\end{eqnarray*}
where $F\in C(\Var(k)^{sm}/S)$ is such that $M=D(\mathbb A^1,et)(F)$.
\end{itemize}
\end{defi}

\begin{prop}\label{keyalgsing1GM}
Let $f:X\to S$ a morphism with $S,X\in\Var(k)$.
Let $S=\cup_{i=1}^l S_i$ an open cover such that there exist closed embeddings
$i_i:S_i\hookrightarrow\tilde S_i$ with $\tilde S_i\in\SmVar(\mathbb C)$. 
Then $X=\cup_{i=1}^lX_i$ with $X_i:=f^{-1}(S_i)$.
Denote, for $I\subset\left[1,\cdots l\right]$, $S_I=\cap_{i\in I} S_i$ and $X_I=\cap_{i\in I}X_i$.
Assume there exist a factorization 
\begin{equation*}
f:X\xrightarrow{l}Y\times S\xrightarrow{p_S} S
\end{equation*}
of $f$ with $Y\in\SmVar(\mathbb C)$, $l$ a closed embedding and $p_S$ the projection.
We then have, for $I\subset\left[1,\cdots l\right]$, the following commutative diagrams which are cartesian 
\begin{equation*}
\xymatrix{
f_I=f_{|X_I}:X_I\ar[r]^{l_I}\ar[rd] & Y\times S_I\ar[r]^{p_{S_I}}\ar[d]^{i'_I} & S_I\ar[d]^{i_I} \\
\, & Y\times\tilde S_I\ar[r]^{p_{\tilde S_I}} & \tilde S_I} \;, \;
\xymatrix{Y\times\tilde S_J\ar[r]^{p_{\tilde S_J}}\ar[d]_{p'_{IJ}} & \tilde S_J\ar[d]^{p_{IJ}} \\
Y\times\tilde S_I\ar[r]^{p_{\tilde S_I}} & \tilde S_I}
\end{equation*}
Let $F(X/S):=p_{S,\sharp}\Gamma_X^{\vee}\mathbb Z(Y\times S/Y\times S)$. 
The transformations maps $(N_I(X/S):Q(X_I/\tilde S_I)\to i_{I*}j_I^*F(X/S))$ and $(k\circ I(\gamma,\hom)(-,-))$, 
for $I\subset\left[1,\cdots,l\right]$, induce an isomorphism in $D_{Ofil,\mathcal D,\infty}(S/(\tilde S_I))$ 
\begin{eqnarray*} 
I^{GM}(X/S): \\
\mathcal F_S^{GM}(M(X/S)):=
(e(\tilde S_I)_*\mathcal Hom(L(i_{I*}j_I^*F(X/S)),E_{et}(\Omega^{\bullet}_{/\tilde S_I},F_b))[-d_{\tilde S_I}],
u^q_{IJ}(F(X/S))) \\
\xrightarrow{(e(\tilde S_I)_*\mathcal Hom(LN_I(X/S),E_{et}(\Omega^{\bullet}_{/\tilde S_I},F_b)))} 
(e(\tilde S_I)_*\mathcal Hom(Q(X_I/\tilde S_I),E_{et}(\Omega^{\bullet}_{/\tilde S_I},F_b))[-d_{\tilde S_I}],
v^q_{IJ}(F(X/S))) \\
\xrightarrow{(k\circ I(\gamma,\hom)(-,-))^{-1}}
(p_{\tilde S_I*}\Gamma_{X_I}E_{zar}(\Omega^{\bullet}_{Y\times\tilde S_I/\tilde S_I},F_b)[-d_{\tilde S_I}],
w_{IJ}(X/S)).
\end{eqnarray*}
\end{prop}

\begin{proof}
Similar to the proof of \cite{B4} proposition.
\end{proof}

\begin{defi}\label{TgDRdefGM}
Let $g:T\to S$ a morphism with $T,S\in\SmVar(\mathbb C)$.
Consider the factorization $g:T\xrightarrow{l}T\times S\xrightarrow{p_S}S$
where $l$ is the graph embedding and $p_S$ the projection.
Let $M\in\DA_c(S)$ and $F\in C(\Var(k)^{sm}/S)$ such that $M=D(\mathbb A^1_S,et)(F)$. 
Then, $D(\mathbb A^1_T,et)(g^*F)=g^*M$. 
\begin{itemize}
\item[(i)]We have then the canonical transformation in $D_{Ofil,\mathcal D,\infty}(T\times S)$
(see definition \ref{TgDRGM}) :
\begin{eqnarray*}
T(g,\mathcal F^{GM})(M):Rg^{*mod[-],\Gamma}\mathcal F_S^{GM}(M):=
g^{*mod,\Gamma}e(S)_*\mathcal Hom^{\bullet}(LF,E_{et}(\Omega^{\bullet}_{/S},F_b)))[-d_T] \\
\xrightarrow{T(g,\Omega_{/\cdot})(LF)} \\ 
e(T\times S)_*\mathcal Hom^{\bullet}(\Gamma_T^{\vee}p_S^*LF,E_{et}(\Omega^{\bullet}_{/T\times S},F_b))[-d_T]
=:\mathcal F_{T\times S}^{GM}(l_*g^*(M,W)). 
\end{eqnarray*}
\item[(ii)]We have then the canonical transformation in $D_{Ofil,\infty}(T)$(see definition \ref{TgDRGM}) :
\begin{eqnarray*}
T^O(g,\mathcal F^{GM})(M,W):Lg^{*mod[-]}\mathcal F_S^{GM}(M):=
g^{*mod}e(S)_*\mathcal Hom^{\bullet}(LF,E_{et}(\Omega^{\bullet}_{/S},F_b)))[-d_T] \\
\xrightarrow{T^O(g,\Omega_{/\cdot})(LF)} \\ 
e(T)_*\mathcal Hom^{\bullet}(g^*LF,E_{et}(\Omega^{\bullet}_{/T},F_b))[-d_T]=:\mathcal F_T^{GM}(g^*M). 
\end{eqnarray*}
\end{itemize}
\end{defi}

\begin{defi}\label{TgDRdefsingGM}
Let $g:T\to S$ a morphism with $T,S\in\Var(k)$.
Assume we have a factorization $g:T\xrightarrow{l}Y\times S\xrightarrow{p_S}S$
with $Y\in\SmVar(\mathbb C)$, $l$ a closed embedding and $p_S$ the projection.
Let $S=\cup_{i=1}^lS_i$ be an open cover such that 
there exists closed embeddings $i_i:S_i\hookrightarrow\tilde S_i$ with $\tilde S_i\in\SmVar(\mathbb C)$
Then, $T=\cup^l_{i=1} T_i$ with $T_i:=g^{-1}(S_i)$
and we have closed embeddings $i'_i:=i_i\circ l:T_i\hookrightarrow Y\times\tilde S_i$,
Moreover $\tilde g_I:=p_{\tilde S_I}:Y\times\tilde S_I\to\tilde S_I$ is a lift of $g_I:=g_{|T_I}:T_I\to S_I$.
Denote for short $d_{YI}:=d_Y+d_{\tilde S_I}$.
Let $M\in\DA_c(S)$ and $F\in C(\Var(k)^{sm}/S)$ such that $M=D(\mathbb A^1_S,et)(F)$.
Then, $D(\mathbb A^1_T,et)(g^*F)=g^*M$. 
We have the canonical transformation in $D_{Ofil,\mathcal D,\infty}(T/(Y\times\tilde S_I))$
\begin{eqnarray*}
T(g,\mathcal F^{GM})(M):Rg^{*mod[-],\Gamma}\mathcal F_S^{GM}(M):= \\
(\Gamma_{T_I}E_{zar}(\tilde g_I^{*mod}e(\tilde S_I)_*\mathcal Hom^{\bullet}(L(i_{I*}j_I^*F),
E_{et}(\Omega^{\bullet}_{/\tilde S_I},F_b)))[-d_Y-d_{\tilde S_I}],\tilde g_J^{*mod}u^q_{IJ}(F)) \\
\xrightarrow{(\Gamma_{T_I}E(T(\tilde g_I,\Omega_{/\cdot})(L(i_{I*}j_I^*(F,W)))))} \\
(\Gamma_{T_I}e(Y\times\tilde S_I)_*\mathcal Hom^{\bullet}(\tilde g_I^*L(i_{I*}j_I^*F),
E_{et}(\Omega^{\bullet}_{/Y\times\tilde S_I},F_b))[-d_Y-d_{\tilde S_I}],\tilde g_J^*u^q_{IJ}(F)_1) \\ 
\xrightarrow{(I(\gamma,\hom(-,-)))} \\
(e(Y\times\tilde S_I)_*\mathcal Hom^{\bullet}(\Gamma_{T_I}^{\vee}\tilde g_I^*L(i_{I*}j_I^*F),
E_{et}(\Omega^{\bullet}_{/Y\times\tilde S_I},F_b))[-d_Y-d_{\tilde S_I}],\tilde g_J^*u^q_{IJ}(F)_2) \\ 
\xrightarrow{(e(Y\times\tilde S_I)_*\mathcal Hom(T^{q,\gamma}(D_{gI})(j_I^*F),
E_{et}(\Omega^{\bullet}_{/Y\times\tilde S_I},F_b)))^{-1}} \\
(e(Y\times\tilde S_I)_*\mathcal Hom^{\bullet}(L(i'_{I*}j_I^{'*}g^*F),
E_{et}(\Omega^{\bullet}_{/Y\times\tilde S_I},F_b))[-d_Y-d_{\tilde S_I}],u^q_{IJ}(g^*F))=:\mathcal F_T^{GM}(g^*M).
\end{eqnarray*} 
\end{defi}

\begin{prop}\label{TgGMprop}
\begin{itemize}
\item[(i)]Let $g:T\to S$ a morphism with $T,S\in\Var(k)$.
Assume we have a factorization $g:T\xrightarrow{l}Y_2\times S\xrightarrow{p_S}S$
with $Y_2\in\SmVar(\mathbb C)$, $l$ a closed embedding and $p_S$ the projection.
Let $S=\cup_{i=1}^lS_i$ be an open cover such that 
there exists closed embeddings $i_i:S_i\hookrightarrow\tilde S_i$ with $\tilde S_i\in\SmVar(\mathbb C)$
Then, $T=\cup^l_{i=1} T_i$ with $T_i:=g^{-1}(S_i)$
and we have closed embeddings $i'_i:=i_i\circ l:T_i\hookrightarrow Y_2\times\tilde S_i$,
Moreover $\tilde g_I:=p_{\tilde S_I}:Y\times\tilde S_I\to\tilde S_I$ is a lift of $g_I:=g_{|T_I}:T_I\to S_I$.
Let $f:X\to S$ a  morphism with $X\in\Var(k)$. Assume that there is a factorization 
$f:X\xrightarrow{l}Y_1\times S\xrightarrow{p_S} S$, with $Y_1\in\SmVar(\mathbb C)$, 
$l$ a closed embedding and $p_S$ the projection. We have then the following commutative diagram
whose squares are cartesians
\begin{equation*}
\xymatrix{f':X_T\ar[r]\ar[d] & Y_1\times T\ar[d]\ar[r] & T\ar[d] \\
f''=f\times I:Y_2\times X\ar[r]\ar[d] & Y_1\times Y_2\times S\ar[r]\ar[d] & Y_2\times S\ar[d] \\
f:X\ar[r] & Y_1\times S\ar[r] & S}
\end{equation*} 
Consider $F(X/S):=p_{S,\sharp}\Gamma_X^{\vee}\mathbb Z(Y_1\times S/Y_1\times S)$ and
the isomorphism in $C(\Var(k)^{sm}/S)$
\begin{eqnarray*}
T(f,g,F(X/S)):g^*F(X/S):=g^*p_{S,\sharp}\Gamma_X^{\vee}\mathbb Z(Y_1\times S/Y_1\times S)\xrightarrow{\sim} \\
p_{T,\sharp}\Gamma_{X_T}^{\vee}\mathbb Z(Y_1\times T/Y_1\times T)=:F(X_T/T).
\end{eqnarray*}
which gives in $\DA(S)$ the isomorphism $T(f,g,F(X/S)):g^*M(X/S)\xrightarrow{\sim}M(X_T/T)$. 
Then, the following diagram in $D_{Ofil,\mathcal D,\infty}(T/(Y_2\times\tilde S_I))$ commutes
\begin{equation*}
\begin{tikzcd}
Rg^{*mod,\Gamma}\mathcal F_S^{GM}(M(X/S))
\ar[rrr,"T(g{,}\mathcal F^{GM})(M(X/S))"]\ar[d,"I^{GM}(X/S)"] & \, & \, & 
\mathcal F_T^{GM}(M(X_T/T))\ar[d,"I^{GM}(X_T/T)"] \\
\shortstack{$g^{*mod[-],\Gamma}(p_{\tilde S_I*}\Gamma_{X_I}
E_{zar}(\Omega^{\bullet}_{Y_1\times\tilde S_I/\tilde S_I},F_b)[-d_{\tilde S_I}],$ \\ $w_{IJ}(X/S))$}
\ar[rrr,"(T(\tilde g_I\times I{,}\gamma)(-)\circ T^O_w(\tilde g_I{,}p_{\tilde S_I}))"] & \, & \, &
\shortstack{$(p_{Y_2\times\tilde S_I*}\Gamma_{X_{T_I}}
E_{zar}(\Omega^{\bullet}_{Y_2\times Y_1\times\tilde S_I/Y_2\times\tilde S_I},F_b)[-d_{Y_2}-d_{\tilde S_I}],$ \\ $w_{IJ}(X_T/T))$}.
\end{tikzcd}
\end{equation*} 
\item[(ii)] Let $g:T\to S$ a morphism with $T,S\in\SmVar(\mathbb C)$.
Let $f:X\to S$ a  morphism with $X\in\Var(k)$. Assume that there is a factorization 
$f:X\xrightarrow{l}Y\times S\xrightarrow{p_S} S$, with $Y\in\SmVar(\mathbb C)$, $l$ a closed embedding and $p_S$ the projection.
Consider $F(X/S):=p_{S,\sharp}\Gamma_X^{\vee}\mathbb Z(Y\times S/Y\times S)$ and
the isomorphism in $C(\Var(k)^{sm}/S)$
\begin{eqnarray*}
T(f,g,F(X/S)):g^*F(X/S):=g^*p_{S,\sharp}\Gamma_X^{\vee}\mathbb Z(Y\times S/Y\times S)\xrightarrow{\sim} \\
p_{T,\sharp}\Gamma_{X_T}^{\vee}\mathbb Z(Y\times T/Y\times T)=:F(X_T/T).
\end{eqnarray*}
which gives in $\DA(S)$ the isomorphism $T(f,g,F(X/S)):g^*M(X/S)\xrightarrow{\sim}M(X_T/T)$. 
Then, the following diagram in $D_{Ofil,\infty}(T)$ commutes
\begin{equation*}
\begin{tikzcd}
Lg^{*mod[-]}\mathcal F_S^{GM}(M(X/S))
\ar[rrr,"T^O(g{,}\mathcal F^{GM})(M(X/S))"]\ar[d,"I^{GM}(X/S)"] & \, & \, & 
\mathcal F_T^{GM}(M(X_T/T))\ar[d,"I^{GM}(X_T/T)"] \\
g^{*mod}L_O(p_{S*}\Gamma_{X}E_{zar}(\Omega^{\bullet}_{Y\times S/S},F_b)[-d_T]
\ar[d,"T_w(\otimes{,}\gamma)(O_{Y\times S})"]\ar[rrr,"(T(g\times I{,}\gamma)(-)\circ T^O_w(g{,}p_S))"] & \, & \, &
p_{Y\times T*}\Gamma_{X_T}E_{zar}(\Omega^{\bullet}_{Y\times T/T},F_b)[-d_T]
\ar[d,"T_w(\otimes{,}\gamma)(O_{Y\times T})"] \\
Lg^{*mod}\int^{FDR}_{p_S}\Gamma_XE(O_{Y\times S},F_b)[-d_Y-d_T]
\ar[rrr,"T^{\mathcal Dmod}(g{,}f)(\Gamma_XE(O_{Y\times S}{,}F_b))"] & \, & \, & 
\int^{FDR}_{p_T}\Gamma_{X_T}E(O_{Y\times T},F_b)[-d_Y-d_T].
\end{tikzcd}
\end{equation*}
\end{itemize}
\end{prop}

\begin{proof}
Follows immediately from definition.
\end{proof}

We have the following theorem:

\begin{thm}\label{mainthmGM}
\begin{itemize}
\item[(i)]Let $g:T\to S$ is a morphism with $T,S\in\Var(k)$. 
Assume there exist a factorization $g:T\xrightarrow{l}Y\times S\xrightarrow{p_S}$
with $Y\in\SmVar(k)$, $l$ a closed embedding and $p_S$ the projection. 
Let $S=\cup_{i=1}^lS_i$ be an open cover such that 
there exists closed embeddings $i_i:S_i\hookrightarrow\tilde S_i$ with $\tilde S_i\in\SmVar(\mathbb C)$.
Then, for $M\in\DA_c(S)$
\begin{eqnarray*}
T(g,\mathcal F^{GM})(M):Rg^{*mod[-],\Gamma}\mathcal F_S^{GM}(M)\xrightarrow{\sim}\mathcal F_T^{GM}(g^*M)
\end{eqnarray*}
is an isomorphism in $D_{O_Tfil,\mathcal D,\infty}(T/(Y\times\tilde S_I))$.
\item[(ii)]Let $g:T\to S$ is a morphism with $T,S\in\SmVar(\mathbb C)$. Then, for $M\in\DA_c(S)$
\begin{eqnarray*}
T^O(g,\mathcal F^{GM})(M):Lg^{*mod[-]}\mathcal F_S^{GM}(M)\xrightarrow{\sim}\mathcal F_T^{GM}(g^*M)
\end{eqnarray*}
is an isomorphism in $D_{O_T}(T)$.
\item[(iii)] A base change theorem for algebraic De Rham cohomology :
Let $g:T\to S$ is a morphism with $T,S\in\SmVar(k)$. 
Let $h:U\to S$ a smooth morphism with $U\in\Var(k)$. Then the map (see definition \cite{B4} section 2)
\begin{equation*}
T^O_w(g,h):Lg^{*mod}Rh_*(\Omega^{\bullet}_{U/S},F_b)\xrightarrow{\sim} Rh'_*(\Omega^{\bullet}_{U_T/T},F_b)
\end{equation*}
is an isomorphism in $D_{O_T}(T)$.
\end{itemize}
\end{thm}

\begin{proof}
Similar to the proof of \cite{B4} theorem.
\end{proof}

\begin{defi}
Let $S\in\Var(k)$ and $S=\cup_{i=1}^l S_i$ an open affine covering and denote, 
for $I\subset\left[1,\cdots l\right]$, $S_I=\cap_{i\in I} S_i$ and $j_I:S_I\hookrightarrow S$ the open embedding.
Let $i_i:S_i\hookrightarrow\tilde S_i$ closed embeddings, with $\tilde S_i\in\SmVar(\mathbb C)$. 
We have, for $M,N\in\DA(S)$ and $F,G\in C(\Var(k)^{sm}/S)$ such that 
$M=D(\mathbb A^1,et)(F)$ and $N=D(\mathbb A^1,et)(G)$, 
the following transformation map in $D_{Ofil,\mathcal D}(S/(\tilde S_I))$
\begin{eqnarray*}
T(\mathcal F_S^{GM},\otimes)(M,N):
\mathcal F_S^{GM}(M)\otimes^{L[-]}_{O_S}\mathcal F_S^{GM}(N):= \\
(e(\tilde S_I)_*\mathcal Hom(L(i_{I*}j_I^*F),E_{et}(\Omega^{\bullet}_{/\tilde S_I},F_b))[-d_{\tilde S_I}],u_{IJ}(F))
\otimes^{[-]}_{O_S} \\
(e(\tilde S_I)_*\mathcal Hom(L(i_{I*}j_I^*G),E_{et}(\Omega^{\bullet}_{/\tilde S_I},F_b))[-d_{\tilde S_I}],u_{IJ}(G)) \\
\xrightarrow{=}
((e(\tilde S_I)_*\mathcal Hom(L(i_{I*}j_I^*F),E_{et}(\Omega^{\bullet}_{/\tilde S_I},F_b))\otimes_{O_{\tilde S_I}} \\
e(\tilde S_I)_*\mathcal Hom(L(i_{I*}j_I^*G),E_{et}(\Omega^{\bullet}_{/\tilde S_I},F_b)))[-d_{\tilde S_I}],
u_{IJ}(F)\otimes u_{IJ}(G)) \\
\xrightarrow{(T(\otimes,\Omega_{/\tilde S_I})(L(i_{I*}j_I^*F),L(i_{I*}j_I^*G)))} \\
(e(\tilde S_I)_*\mathcal Hom(L(i_{I*}j_I^*F)\otimes L(i_{I*}j_I^*G),
E_{et}(\Omega^{\bullet}_{/\tilde S_I},F_{DR}))[-d_{\tilde S_I}],v_{IJ}(F\otimes G)) \\
\xrightarrow{=} 
(e(\tilde S_I)_*\mathcal Hom(L(i_{I*}j_I^*(F\otimes G),E_{et}(\Omega^{\bullet}_{/\tilde S_I},F_b)))[-d_{\tilde S_I}],
u_{IJ}(F\otimes G))=:\mathcal F_S^{GM}(M\otimes N)
\end{eqnarray*}
\end{defi}

\begin{prop}
Let $f_1:X_1\to S$, $f_2:X_2\to S$ two morphism with $X_1,X_2,S\in\Var(k)$. 
Assume that there exist factorizations 
$f_1:X_1\xrightarrow{l_1}Y_1\times S\xrightarrow{p_S} S$, $f_2:X_2\xrightarrow{l_2}Y_2\times S\xrightarrow{p_S} S$
with $Y_1,Y_2\in\SmVar(\mathbb C)$, $l_1,l_2$ closed embeddings and $p_S$ the projections.
We have then the factorization
\begin{equation*}
f_1\times f_2:X_{12}:=X_1\times_S X_2\xrightarrow{l_1\times l_2}Y_1\times Y_2\times S\xrightarrow{p_S} S
\end{equation*}
Let $S=\cup_{i=1}^l S_i$ an open affine covering and denote, 
for $I\subset\left[1,\cdots l\right]$, $S_I=\cap_{i\in I} S_i$ and $j_I:S_I\hookrightarrow S$ the open embedding.
Let $i_i:S_i\hookrightarrow\tilde S_i$ closed embeddings, with $\tilde S_i\in\SmVar(\mathbb C)$. 
We have, for $M,N\in\DA(S)$ and $F,G\in C(\Var(k)^{sm}/S)$ such that 
$M=D(\mathbb A^1,et)(F)$ and $N=D(\mathbb A^1,et)(G)$, 
the following commutative diagram in $D_{Ofil,\mathcal D}(S/(\tilde S_I))$
\begin{equation*}
\begin{tikzcd}
\mathcal F_S^{GM}(M(X_1/S))\otimes^L_{O_S}\mathcal F_S^{GM}(M(X_2/S))
\ar[d,"I^{GM}(X_1/S)\otimes I^{GM}(X_2/S)"]
\ar[rr,"T(\mathcal F_S^{GM}{,}\otimes)(M(X_1/S){,}M(X_2/S))"]  & \, &
\shortstack{$\mathcal F_S^{GM}(M(X_1/S)\otimes M(X_2/S))$ \\ $=\mathcal F_S^{GM}(M(X_1\times_S X_2/S))$}
\ar[d,"I^{GM}(X_{12}/S)"] \\
\shortstack{$(p_{\tilde S_I*}\Gamma_{X_{1I}}
E_{zar}(\Omega^{\bullet}_{Y_1\times\tilde S_I/\tilde S_I},F_b)[-d_{\tilde S_I}],w_{IJ}(X_1/S))\otimes_{O_S}$ \\
$(p_{\tilde S_I*}\Gamma_{X_{2I}}
E_{zar}(\Omega^{\bullet}_{Y_2\times\tilde S_I/\tilde S_I},F_b)[-d_{\tilde S_I}],w_{IJ}(X_2/S))$}
\ar[rr,"(Ew_{(Y_1\times\tilde S_I{,}Y_2\times\tilde S_I)/\tilde S_I})"] & \, &
\shortstack{$(p_{\tilde S_I*}\Gamma_{X_{12I}}
E_{zar}(\Omega^{\bullet}_{Y_1\times Y_2\times\tilde S_I/\tilde S_I},F_b)[-d_{\tilde S_I}],$ \\ $w_{IJ}(X_{12}/S))$}.
\end{tikzcd}
\end{equation*}
\end{prop}

\begin{proof}
Immediate from definition.
\end{proof}

\subsection{The algebraic filtered De Rham realization functor}

Let $k$ a field of caracteristic zero.
We recall (see section 2), for $f:T\to S$ a morphism with $T,S\in\Var(k)$, 
the commutative diagrams of sites (\ref{muf}) and (\ref{Grf})
\begin{equation*}
\xymatrix{
\Var(k)^2/T\ar[rd]^{\rho_T}\ar[rr]^{\mu_T}\ar[dd]_{P(f)} & \, & 
\Var(k)^{2,pr}/T\ar[dd]^{P(f)}\ar[rd]^{\rho_T} & \, \\
\, & \Var(k)^{2,sm}/T\ar[rr]^{\mu_T}\ar[dd]_{P(f)} & \, & \Var(k)^{2,smpr}/T\ar[dd]^{P(f)} \\
\Var(k)^2/S\ar[rd]^{\rho_S}\ar[rr]^{\mu_S} & \, & \Var(k)^{2,pr}/S\ar[rd]^{\rho_S} \\
\, & \Var(k)^{2,sm}/S\ar[rr]^{\mu_S} & \, & \Var(k)^{2,smpr}/S}
\end{equation*}
and
\begin{equation*}
\xymatrix{\Var(k)^{2,pr}/T\ar[rr]^{\Gr_T^{12}}\ar[dd]_{P(f)}\ar[rd]^{\rho_T} & \, & 
\Var(k)/T\ar[dd]^{P(f)}\ar[rd]^{\rho_T} & \, \\
\, & \Var(k)^{2,smpr}/T\ar[rr]^{\Gr_T^{12}}\ar[dd]_{P(f)} & \, & \Var(k)^{sm}/T\ar[dd]^{P(f)} \\
\Var(k)^{2,pr}/S\ar[rr]^{\Gr_S^{12}}\ar[rd]^{\rho_S} & \, & \Var(k)/S\ar[rd]^{\rho_S} & \, \\
\, & \Var(k)^{2,sm}/S\ar[rr]^{\Gr_S^{12}} & \, & \Var(k)^{sm}/S}.
\end{equation*}

For $s:\mathcal I\to\mathcal J$ a functor, with $\mathcal I,\mathcal J\in\Cat$, and
$f_{\bullet}:T_{\bullet}\to S_{s(\bullet)}$ a morphism with 
$T_{\bullet}\in\Fun(\mathcal J,\Var(k))$ and $S_{\bullet}\in\Fun(\mathcal I,\Var(k))$, 
we have then the commutative diagrams of sites (\ref{mufIJ}) and (\ref{GrfIJ})
\begin{equation*}
\xymatrix{\Var(k)^2/T_{\bullet}\ar[rr]^{\mu_{T_{\bullet}}}\ar[dd]_{P(f_{\bullet})}\ar[rd]^{\rho_{T_{\bullet}}} & \, & 
\Var(k)^{2,pr}/T_{\bullet}\ar[dd]^{P(f_{\bullet})}\ar[rd]^{\rho_{T_{\bullet}}} & \, \\
\, & \Var(k)^{2,sm}/T_{\bullet}\ar[rr]^{\mu_{T_{\bullet}}}\ar[dd]_{P(f_{\bullet})} & \, & 
\Var(k)^{2,smpr}/T_{\bullet}\ar[dd]^{P(f_{\bullet})} \\
\Var(k)^2/S_{\bullet}\ar[rr]^{\mu_{S_{\bullet}}}\ar[rd]^{\rho_{S_{\bullet}}} & \, & 
\Var(k)^{2,pr}/S_{\bullet}\ar[rd]^{\rho_{S_{\bullet}}} & \, \\
\, & \Var(k)^{2,sm}/S_{\bullet}\ar[rr]^{\mu_{S_{\bullet}}} & \, & \Var(k)^{2,smpr}/S_{\bullet}}.
\end{equation*}
and
\begin{equation*}
\xymatrix{\Var(k)^{2,pr}/T_{\bullet}
\ar[rr]^{\Gr_{T_{\bullet}}^{12}}\ar[dd]_{P(f_{\bullet})}\ar[rd]^{\rho_{T_{\bullet}}} & \, & 
\Var(k)/T\ar[dd]^{P(f_{\bullet})}\ar[rd]^{\rho_{T_{\bullet}}} & \, \\
\, & \Var(k)^{2,smpr}/T_{\bullet}\ar[rr]^{\Gr_{T}^{12}}\ar[dd]_{P(f_{\bullet})} & \, & 
\Var(k)^{sm}/T_{\bullet}\ar[dd]^{P(f_{\bullet})} \\
\Var(k)^{2,pr}/S_{\bullet}\ar[rr]^{\Gr_{S_{\bullet}}^{12}}\ar[rd]^{\rho_{S_{\bullet}}} & \, & 
\Var(k)/S_{\bullet}\ar[rd]^{\rho_{S_{\bullet}}} & \, \\
\, & \Var(k)^{2,sm}/S_{\bullet}\ar[rr]^{\Gr_{S_{\bullet}}^{12}} & \, & 
\Var(k)^{sm}/S_{\bullet}}.
\end{equation*}

We will use the following map from the property of De Rham modules (see section 5)
together with the specialization map of a filtered D module for a closed embedding (see \cite{B4} section 4.1) :

\begin{defiprop}\label{VfilOlem}
\begin{itemize}
\item[(i)]Let $l:Z\hookrightarrow S$ a closed embedding with $S,Z\in\SmVar(k)$.
Consider an open embedding $j:S^o\hookrightarrow S$. We then have the cartesian square
\begin{equation*}
\xymatrix{S^o\ar[r]^j & S \\ 
Z^o:=Z\times_SS^o\ar[r]^{j'}\ar[u]^{l'} & Z\ar[u]^l}
\end{equation*}
where $j'$ is the open embedding given by base change. Using proposition \ref{VQjO}(ii) or theorem \ref{HSk}, 
the morphisms $Q^{p,0}_{V_Z,V_D}(O_S,F_b)$ for $D\subset S$ a closed subset
of definition-proposition of \cite{B4} induces a canonical morphism in $C_{l^*O_Sfil}(Z)$
\begin{equation*}
Q(Z,j_!)(O_S,F_b):l^*Q_{V_Z,0}j_{!Hdg}(O_{S^o},F_b)\to j'_{!Hdg}(O_{Z^o},F_b),
\end{equation*}
where $V_Z$ is the Kashiwara-Malgrange $V_Z$-filtration and $V_D$ is the Kashiwara-Malgrange $V_D$-filtration,
which commutes with the action of $T_Z$.
\item[(ii)]Let $l:Z\hookrightarrow S$ and $k:Z'\hookrightarrow Z$ be closed embeddings with $S,Z,Z'\in\SmVar(k)$.
Consider an open embedding $j:S^o\hookrightarrow S$. We then have the commutative diagram whose squares are cartesian.
\begin{equation*}
\xymatrix{S^o\ar[r]^j & S \\ 
Z^o:=Z\times_SS^o\ar[r]^{j'}\ar[u]^{l'} & Z\ar[u]^l \\
Z^{'o}:=Z'\times_SS^o\ar[r]^{j''}\ar[u]^{k'} & Z'\ar[u]^{k}}
\end{equation*}
where $j'$ is the open embedding given by base change. Then, 
\begin{eqnarray*}
Q(Z',j_!)(O_S,F_b)=Q(Z',j'_!)(O_Z,F_b)\circ (k^*Q_{V_{Z'},0}Q(Z,j_!)(O_S,F_b)): \\
k^*Q_{V_{Z'},0}l^*Q_{V_Z,0}j_{!Hdg}(O_{S^o},F_b)
\xrightarrow{k^*Q_{V_{Z'},0}Q(Z,j_!)(O_S,F_b)}k^*Q_{V_{Z'},0}j'_{!Hdg}(O_{Z^o},F_b) \\
\xrightarrow{Q(Z',j'_!)(O_Z,F_b)}j''_{!Hdg}(O_{Z^{'o}},F_b)
\end{eqnarray*}
in $C_{k^*l^*O_Sfil}(Z')$ which commutes with the action of $T_{Z'}$.
\item[(iii)] Consider a commutative diagram whose squares are cartesian
\begin{equation*}
\xymatrix{S^{oo}\ar[r]^{j_2} & S^o\ar[r]^{j_1} & S \\ 
Z^{oo}:=Z\times_SS^{oo}\ar[r]^{j'_2}\ar[u]^{l''} & Z^o:=Z\times_SS^o\ar[r]^{j'_1}\ar[u]^{l'} & Z\ar[u]^l}
\end{equation*}
where $j_1$, $j_2$, and hence $j'_1$,$j'_2$ are open embeddings.
We have then the following commutative diagram
\begin{equation*}
\xymatrix{l^*Q_{V_Z,0}j_{1!Hdg}(O_{S^o},F_b)
\ar[rr]^{\ad(j_{2!Hdg},j_2^*)(O_{S^o},F_b)}\ar[d]^{Q(Z,j_!)(O_S,F_b)} & \, &
l^*Q_{V_Z,0}(j_1\circ j_2)^{Hdg}_!(O_{S^{oo}},F_b)\ar[d]^{Q(Z,(j_1\circ j_2)_!)(O_S,F_b)} \\
j'_{1!Hdg}(O_{Z^o},F_b)\ar[rr]^{\ad(j'_{2!Hdg},j_2^{'*})(O_{Z^o},F_b)} & \, &
(j'_1\circ j'_2)_{!Hdg}(O_{Z^{oo}},F_b)}
\end{equation*}
in $C_{l^*O_Sfil}(Z)$ which commutes with the action of $T_{Z}$.
\end{itemize}
\end{defiprop}

\begin{proof}
\noindent(i): By definition of $j_{!Hdg}:\pi_{S^o}(MHM(S^o))\to C(DRM(S)))$,
we have to construct the isomorphism for each complement of a (Cartier) divisor
$j=j_D:S^o=S\backslash D\hookrightarrow S$.
In this case, we have the closed embedding $i:S\hookrightarrow L$ given by the zero section of the line bundle
$L=L_D$ associated to $D$. 
We have then, using definition-proposition of \cite{B4} section 4.1,
the canonical morphism in $PSh_{l^*O_Sfil}(Z)$ which commutes with the action of $T_Z$
\begin{eqnarray*}
Q(Z,j_!)(O_S,F_b):l^*Q_{V_Z,0}j_{!Hdg}(O_{S^o},F_b)\xrightarrow{T_!(l,j)(-)^{-1}}
l^*j_{!Hdg}Q_{V_{Z^o},0}(O_{S^o},F_b)=j'_{!Hdg}(O_{Z^o},F_b).
\end{eqnarray*}
and $V_Z^pT_!(l,j)(-)^{-1}=Q^{p,0}_{V_Z,V_S}(i_{*mod}(O_S,F_b))$.
Now for $j:S^o=S\backslash R\hookrightarrow S$ an arbitrary open embedding, we set
\begin{eqnarray*}
Q(Z,j_!)(O_S,F_b):=\varprojlim_{(D_i),R\subset D_i\subset S} (Q(Z,j_{D_J!})(j_{D_I}^*(O_S,F_b)): 
l^*Q_{V_Z,0}j_{!Hdg}(O_{S^o},F_b)\xrightarrow{\sim}j'_{!Hdg}(O_{Z^o},F_b)
\end{eqnarray*}

\noindent(ii): Follows from \cite{B4} section 4.1.

\noindent(iii): Follows from \cite{B4} section 4.1.
\end{proof}

Using definition-proposition \ref{TgammaHdg} in the projection case, and
the specialization map given in \cite{B4} section 4 and the isomorphism of definition-proposition \ref{VfilOlem},
in the closed embedding case, we have the following canonical map :

\begin{defi}\label{TgHdgO}
Consider a commutative diagram in $\SmVar(k)$ whose square are cartesian
\begin{equation*}
\xymatrix{Z_T\ar[r]^{i'}\ar[dd]^{g'}\ar[rd]^{l'} & T\ar[dd]^g\ar[rd]^l  & 
T\backslash Z_T\ar[l]^{j'}\ar[dd]^{g}\ar[rd]^{l} \\ 
 \, & T\times Z\ar[r]^{I\times i}\ar[ld]^{p_Z} & T\times S\ar[ld]^{p_S} & 
T\times S\backslash (T\times Z)\ar[l]^{I\times j}\ar[ld]^{p_S}\\
Z\ar[r]^i & S & S\backslash Z\ar[l]^j}
\end{equation*}
where $i$ and hence $I\times i$ and $i'$, are closed embeddings,
$j$, $I\times j$, $j'$ are the complementary open embeddings
and $g:T\xrightarrow{l}T\times S\xrightarrow{p_S}S$ 
is the graph factorization, where $l$ is the graph embedding and $p_S$ the projection. 
Then, the map in $C_{l^*O_{T\times S}fil}(T)$
\begin{eqnarray*}
sp_{V_T}(\Gamma_{T\times Z}^{\vee,Hdg}(O_{T\times S},F_b)):
l^*\Gamma_{T\times Z}^{\vee,Hdg}(O_{T\times S},F_b)\xrightarrow{q_{V_T,0}}
l^*Q_{V_T,0}(\Gamma_{T\times Z}^{\vee,Hdg}(O_{T\times S},F_b)) \\
\xrightarrow{Q(T,(I\times j)_!)(O_{T\times S},F_b):=T_!(l,(I\times j))(-)}
\Gamma_{Z_T}^{\vee,Hdg}(O_T,F_b)
\end{eqnarray*}
which commutes with the action of $T_T$, where the first map is given in \cite{B4} section 4.1 
and the last map is studied definition-proposition \ref{VfilOlem}, factors through
\begin{eqnarray*}
sp_{V_T}(\Gamma_{T\times Z}^{\vee,Hdg}(O_{T\times S},F_b)):
l^*\Gamma_{T\times Z}^{\vee,Hdg}(O_{T\times S},F_b)\xrightarrow{n}
l^{*mod}\Gamma_{T\times Z}^{\vee,Hdg}(O_{T\times S},F_b) \\
\xrightarrow{\bar{sp}_{V_T}(\Gamma_{T\times Z}^{\vee,Hdg}(O_{T\times S},F_b))}
\Gamma_{Z_T}^{\vee,Hdg}(O_T,F_b),
\end{eqnarray*}
with for $U\subset T\times S$ an open subset, $m\in\Gamma(U,O_{T\times S})$ and $h\in\Gamma(U_T,O_T)$,
$n(m):=n\otimes 1$ and $\bar{sp}_{V_T}(-)(m\otimes h)=h\cdot sp_{V_T}(m)$ ;
see definition-proposition \ref{TgammaHdg} and theorem \ref{HSk}. Then,
\begin{eqnarray*}
\bar{sp}_{V_T}(\Gamma_{T\times Z}^{\vee,Hdg}(O_{T\times S},F_b)):
l^{*mod}\Gamma_{T\times Z}^{\vee,Hdg}(O_{T\times S},F_b)\to\Gamma_{Z_T}^{\vee,Hdg}(O_T,F_b),
\end{eqnarray*}
is a map in $C_{\mathcal D(1,0)fil}(T)$, i.e. is $D_T$ linear.
We then consider the canonical map in $C_{\mathcal D(1,0)fil}(T)$ 
\begin{eqnarray*}
a(g,Z)(O_S,F_b):g^{*mod}\Gamma_Z^{\vee,Hdg}(O_S,F_b)=l^{*mod}p_S^{*mod}\Gamma_Z^{\vee,Hdg}(O_S,F_b) 
\xrightarrow{l^{*mod}T^{Hdg}(p,\gamma^{\vee})(O_S,F_b)^{-1}} \\
l^{*mod}\Gamma_{T\times Z}^{\vee,Hdg}(O_{T\times S},F_b) 
\xrightarrow{\bar{sp}_{V_T}(\Gamma_{T\times Z}^{\vee,Hdg}(O_{T\times S},F_b))}
\Gamma_{Z_T}^{\vee,Hdg}(O_T,F_b).
\end{eqnarray*}
\end{defi}

\begin{lem}\label{TgHdgOlem}
\begin{itemize}
\item[(i)]For $g:T\to S$ and $g:T'\to T$ two morphism with $S,T,T'\in\SmVar(k)$,
considering the commutative diagram whose squares are cartesian
\begin{equation*}
\xymatrix{Z_{T'}\ar[r]^{i''}\ar[d]^{g'} & T'\ar[d]^{g'}  & T'\backslash Z_{T'}\ar[l]^{j''}\ar[d]^{g'} \\ 
Z_T\ar[r]^{i'}\ar[d]^g & T\ar[d]^{g} & T\backslash Z_T\ar[l]^{j'}\ar[d]^{g} \\
Z\ar[r]^i & S & S\backslash Z\ar[l]^j}
\end{equation*}
we have then
\begin{eqnarray*}
a(g\circ g',Z)(O_S,F_b)=a(g',Z_T)(O_T,F_b)\circ (g^{'*mod}a(g,Z)(O_S,F_b)): \\
(g\circ g')^{*mod}\Gamma_Z^{\vee,Hdg}(O_S,F_b)=g^{'*mod}g^{*mod}\Gamma_Z^{\vee,Hdg}(O_S,F_b)
\xrightarrow{g^{'*mod}a(g,Z)(O_S,F_b)}g^{'*mod}\Gamma_{Z_T}^{\vee,Hdg}(O_T,F_b) \\
\xrightarrow{a(g',Z_T)(O_T,F_b)}\Gamma_{Z_{T'}}^{\vee,Hdg}(O_{T'},F_b).
\end{eqnarray*}
\item[(ii)]For $g:T\to S$ a morphism with $S,T\in\SmVar(k)$,
considering the commutative diagram whose squares are cartesian
\begin{equation*}
\xymatrix{Z'_T\ar[r]^{k'}\ar[d]^g & Z_T\ar[r]^{i'}\ar[d]^{g} & T\ar[d]^{g} \\
Z'\ar[r]^{k} & Z\ar[r]^i & S }
\end{equation*}
we have then the following commutative diagram
\begin{equation*}
\xymatrix{g^{*mod}\Gamma_Z^{\vee,Hdg}(O_S,F_b)
\ar[d]^{a(g,Z)(O_S,F_b)}\ar[rr]^{g^{*mod}T(Z'/Z,\gamma^{\vee,Hdg})(O_S,F_b)} & \, &
g^{*mod}\Gamma_{Z'}^{\vee,Hdg}(O_S,F_b)\ar[d]^{a(g,Z')(O_S,F_b)} \\
\Gamma_{Z_T}^{\vee,Hdg}(O_T,F_b)\ar[rr]^{T(Z'_T/Z_T,\gamma^{\vee,Hdg})(O_T,F_b)} & \, &
\Gamma_{Z'_T}^{\vee,Hdg}(O_T,F_b)}
\end{equation*}
\end{itemize}
\end{lem}

\begin{proof}
\noindent(i):Follows from definition-proposition \ref{VfilOlem} (ii)

\noindent(ii):Follows from definition-proposition \ref{VfilOlem} (iii)
\end{proof}

We can now define the main object :

\begin{defi}\label{wtildew}
\begin{itemize}
\item[(i)]For $S\in\SmVar(k)$, we consider the filtered complexes of presheaves 
\begin{eqnarray*} 
(\Omega^{\bullet,\Gamma,pr}_{/S},F_{DR})\in C_{D_Sfil}(\Var(k)^{2,smpr}/S) 
\end{eqnarray*}
given by, 
\begin{itemize}
\item for $(Y\times S,Z)/S=((Y\times S,Z),p)\in\Var(k)^{2,smpr}/S$, 
\begin{eqnarray*}
(\Omega^{\bullet,\Gamma,pr}_{/S}((Y\times S,Z)/S),F_{DR}):=
((\Omega^{\bullet}_{Y\times S/S},F_b)\otimes_{O_{Y\times S}}\Gamma^{\vee,Hdg}_Z(O_{Y\times S},F_b))(Y\times S) 
\end{eqnarray*}
with the structure of $p^*D_S$ module given by proposition \ref{DRhUS},
\item for $g:(Y_1\times S,Z_1)/S=((Y_1\times S,Z_1),p_1)\to (Y\times S,Z)/S=((Y\times S,Z),p)$ 
a morphism in $\Var(k)^{2,smpr}/S$, denoting for short $\hat{Z}:=Z\times_{Y\times S}(Y_1\times S)$,
\begin{eqnarray*}
\Omega^{\bullet,\Gamma,pr}_{/S}(g):
((\Omega^{\bullet}_{Y\times S/S},F_b)\otimes_{O_{Y\times S}}\Gamma^{\vee,Hdg}_Z(O_{Y\times S},F_b))(Y\times S) \\ 
\xrightarrow{i_{-}} 
g^*((\Omega^{\bullet}_{Y\times S/S},F_b)\otimes_{O_{Y\times S}}\Gamma^{\vee,Hdg}_Z(O_{Y\times S},F_b))(Y_1\times S) \\  
\xrightarrow{\Omega_{(Y_1\times S/Y\times S)/(S/S)}(\Gamma^{\vee,Hdg}_Z(O_{Y\times S},F_b))(Y_1\times S)} 
(\Omega^{\bullet}_{Y_1\times S/S},F_b)\otimes_{O_{Y_1\times S}}g^{*mod}\Gamma^{\vee,Hdg}_Z(O_{Y\times S},F_b))(Y_1\times S) \\ 
\xrightarrow{DR(Y_1\times S/S)(a(g,Z)(O_{Y\times S},F_b))(Y_1\times S)} 
(\Omega^{\bullet}_{Y_1\times S/S},F_b)\otimes_{O_{Y_1\times S}}\Gamma^{\vee,Hdg}_{\hat{Z}}(O_{Y_1\times S},F_b))(Y_1\times S) \\
\xrightarrow{DR(Y_1\times S/S)(T(Z_1/\hat{Z},\gamma^{\vee,Hdg})(O_{Y_1\times S},F_b))(Y_1\times S)} 
(\Omega^{\bullet}_{Y_1\times S/S},F_b)\otimes_{O_{Y_1\times S}}\Gamma^{\vee,Hdg}_{Z_1}(O_{Y_1\times S},F_b))(Y_1\times S), 
\end{eqnarray*}
where 
\begin{itemize}
\item $i_{-}$ is the arrow of the inductive limit,
\item we recall that
\begin{eqnarray*}
\Omega_{(Y_1\times S/Y\times S)/(S/S)}(\Gamma^{\vee,Hdg}_Z(O_{Y\times S},F_b)):
g^*((\Omega^{\bullet}_{Y\times S/S},F_b)\otimes_{O_{Y\times S}}\Gamma^{\vee,Hdg}_Z(O_{Y\times S},F_b)) \\  
\to(\Omega^{\bullet}_{Y_1\times S/S},F_b)\otimes_{O_{Y_1\times S}}g^{*mod}\Gamma^{\vee,Hdg}_Z(O_{Y\times S},F_b))
\end{eqnarray*}
is the map given in \cite{B4} section 4.1, which is $p_1^*D_S$ linear by proposition \ref{TDhwM},
\item the map
\begin{equation*}
a(g,Z)(O_{Y\times S},F_b):g^{*mod}\Gamma^{\vee,Hdg}_{Z}(O_{Y\times S},F_b)\to\Gamma^{\vee,Hdg}_{\hat{Z}}(O_{Y_1\times S},F_b) 
\end{equation*}
is the map given in definition \ref{TgHdgO} 
\item the map
\begin{eqnarray*}
T(Z_1/\hat{Z},\gamma^{\vee,Hdg})(O_{Y_1\times S},F_b):
\Gamma^{\vee,Hdg}_{\hat{Z}}(O_{Y_1\times S},F_b)\to\Gamma^{\vee,Hdg}_{Z_1}(O_{Y_1\times S},F_b) 
\end{eqnarray*}
is given in definition-proposition \ref{TgammaHdg}.
\end{itemize}
\end{itemize}
For $g:((Y_1\times S,Z_1),p_1)\to ((Y\times S,Z),p)$ and $g':((Y'_1\times S,Z'_1),p_1)\to ((Y_1\times S,Z_1),p)$
two morphisms in $\Var(k)^{2,smpr}/S$, we have
\begin{eqnarray*}
\Omega^{\bullet,\Gamma,pr}_{/S}(g\circ g')=
\Omega^{\bullet,\Gamma,pr}_{/S}(g')\circ\Omega^{\bullet,\Gamma,pr}_{/S}(g):
((\Omega^{\bullet}_{Y\times S/S},F_b)\otimes_{O_{Y\times S}}\Gamma^{\vee,Hdg}_Z(O_{Y\times S},F_b))(Y\times S) \\ 
\xrightarrow{\Omega^{\bullet,\Gamma,pr}_{/S}(g)} 
(\Omega^{\bullet}_{Y_1\times S/S},F_b)\otimes_{O_{Y_1\times S}}\Gamma^{\vee,Hdg}_{Z_1}(O_{Y_1\times S},F_b))(Y_1\times S) \\ 
\xrightarrow{\Omega^{\bullet,\Gamma,pr}_{/S}(g')} 
(\Omega^{\bullet}_{Y'_1\times S/S},F_b)\otimes_{O_{Y'_1\times S}}\Gamma^{\vee,Hdg}_{Z'_1}(O_{Y'_1\times S},F_b))(Y'_1\times S),
\end{eqnarray*}
since, denoting for short $\hat{Z}:=Z\times_{Y\times S}(Y_1\times S)$ and $\hat{Z}':=Z\times_{Y\times S}(Y'_1\times S)$
\begin{itemize}
\item we have by lemma \ref{TgHdgOlem}(i)
\begin{equation*}
a(g\circ g',\hat Z')(O_{Y\times S},F_b)=a(g',\hat Z)(O_{Y_1\times S},F_b)\circ g^{'*mod}a(g,Z)(O_{Y\times S},F_b),
\end{equation*}
\item we have by lemma \ref{TgHdgOlem}(ii)
\begin{eqnarray*}
T(Z'_1/\hat{Z}',\gamma^{\vee,Hdg})(O_{Y'_1\times S},F_b)\circ a(g',\hat Z)(O_{Y_1\times S},F_b) \\
=a(g',Z_1)(O_{Y_1\times S},F_b)\circ g^{'*mod}T(Z_1/\hat{Z},\gamma^{\vee,Hdg})(O_{Y_1\times S},F_b).
\end{eqnarray*}
\end{itemize}
\item[(ii)]For $S\in\SmVar(k)$, we have the canonical map $C_{O_Sfil,D_S}(\Var(k)^{sm}/S)$ 
\begin{eqnarray*} 
\Gr(\Omega_{/S}):\Gr_{S*}^{12}(\Omega^{\bullet,\Gamma,pr}_{/S},F_{DR})\to(\Omega^{\bullet}_{/S},F_b) 
\end{eqnarray*}
given by, for $U/S=(U,h)\in\Var(k)^{sm}/S$
\begin{eqnarray*}
\Gr(\Omega_{/S})(U/S):\Gr_{S*}^{12}(\Omega^{\bullet,\Gamma,pr}_{/S},F_{DR})(U/S):=
((\Omega^{\bullet}_{U\times S/S},F_b)\otimes_{O_{U\times S}}\Gamma^{\vee,Hdg}_{U}(O_{U\times S},F_b))(U\times S) \\
\xrightarrow{\ad(i_U^*,i_{U*})(-)(U\times S)}
i^*((\Omega^{\bullet}_{U\times S/S},F_b)\otimes_{O_{U\times S}}\Gamma^{\vee,Hdg}_{U}(O_{U\times S},F_b))(U) \\
\xrightarrow{\Omega_{(U/U\times S)/(S/S)}(-)(U)}
((\Omega^{\bullet}_{U/S},F_b)\otimes_{O_{U}}i_U^{*mod}\Gamma^{\vee,Hdg}_{U}(O_{U\times S},F_b))(U) \\
\xrightarrow{DR(U/S)(a(i_U,U))(U)}
(\Omega^{\bullet}_{U/S},F_b)(U)=:(\Omega^{\bullet}_{/S},F_b)(U/S)
\end{eqnarray*}
where $h:U\xrightarrow{i_U}U\times S\xrightarrow{p_S}S$ is the graph factorization 
with $i_U$ the graph embedding and $p_S$ the projection, note that $a(i_U,U)$ is an isomorphism since for
$j_U:U\times S\backslash U\hookrightarrow U\times S$ the open complementary $i_U^{*mod}j^{Hdg}_{U!}(M,F,W)=0$.
\end{itemize}
\end{defi}

Let $g:T\to S$ a morphism with $T,S\in\SmVar(k)$.
We have the canonical morphism in $C_{g^*D_Sfil}(\Var(k)^{2,smpr}/T)$ 
\begin{eqnarray*}
\Omega^{\Gamma,pr}_{/(T/S)}:g^*(\Omega^{\bullet,\Gamma,pr}_{/S},F_{DR})\to (\Omega^{\bullet,\Gamma,pr}_{/T},F_{DR})
\end{eqnarray*}
induced by the pullback of differential forms : 
for $((Y_1\times T,Z_1)/T)=((Y_1\times T,Z_1),p)\in\Var(k)^{2,smpr}/T$,
\begin{eqnarray*}
\Omega^{\Gamma,pr}_{/(T/S)}((Y_1\times T,Z_1)/T): \\
g^*\Omega^{\bullet,\Gamma,pr}_{/S}((Y_1\times T,Z_1)/T):= 
\lim_{(h:(Y\times S,Z)\to S, \, g_1:(Y_1\times T,Z_1)\to (Y\times T,Z_T),h,g)}
\Omega^{\bullet,\Gamma,pr}_{/S}((Y\times T,Z)/S) \\
\xrightarrow{\Omega^{\bullet,\Gamma,pr}_{/S}(g'\circ g_1)}\Omega^{\bullet,\Gamma,pr}_{/S}((Y_1\times T,Z_1)/S)
\xrightarrow{q(-)(Y_1\times T)}\Omega^{\bullet,\Gamma,pr}_{/T}((Y_1\times T,Z_1)/T), 
\end{eqnarray*}
where $g'=(I_Y\times g):Y\times T\to Y\times S$ is the base change map and 
$q(M):\Omega_{Y_1\times T/S}\otimes_{O_{Y_1\times T}}(M,F)\to\Omega_{Y_1\times T/T}\otimes_{O_{Y_1\times T}}(M,F)$
is the quotient map. 
It induces the canonical morphisms in $C_{g^*D_Sfil}(\Var(k)^{2,smpr}/T)$ :
\begin{eqnarray*}
E\Omega^{\Gamma,pr}_{/(T/S)}:g^*E_{et}(\Omega^{\bullet,\Gamma,pr}_{/S},F_{DR})
\xrightarrow{T(g,E)(-)} E_{et}(g^*(\Omega^{\bullet,\Gamma,pr}_{/S},F_{DR}))
\xrightarrow{E_{et}(\Omega^{\Gamma,pr}_{/(T/S)})}E_{et}(\Omega^{\bullet,\Gamma,pr}_{/T},F_{DR})
\end{eqnarray*}
and
\begin{eqnarray*}
E\Omega^{\Gamma,pr}_{/(T/S)}:g^*E_{zar}(\Omega^{\bullet,\Gamma,pr}_{/S},F_{DR})
\xrightarrow{T(g,E)(-)} E_{zar}(g^*(\Omega^{\bullet,\Gamma,pr}_{/S},F_{DR}))
\xrightarrow{E_{zar}(\Omega^{\Gamma,pr}_{/(T/S)})}E_{zar}(\Omega^{\bullet,\Gamma,pr}_{/T},F_{DR}).
\end{eqnarray*}

\begin{defi}\label{TgDR}
Let $g:T\to S$ a morphism with $T,S\in\SmVar(k)$.
We have, for $F\in C(\Var(k)^{2,smpr}/S)$, the canonical transformation in $C_{\mathcal Dfil}(T)$ :
\begin{eqnarray*}
T(g,\Omega^{\Gamma,pr}_{/\cdot})(F): 
g^{*mod}L_De(S)_*\Gr^{12}_{S*}\mathcal Hom^{\bullet}(F,E_{zar}(\Omega^{\bullet,\Gamma,pr}_{/S},F_{DR})) \\
\xrightarrow{=} 
(g^*L_De(S)_*\mathcal Hom^{\bullet}(F,E_{zar}(\Omega^{\bullet,\Gamma,pr}_{/S},F_{DR})))\otimes_{g^*O_S}O_T \\
\xrightarrow{T(g,\Gr^{12})(-)\circ T(e,g)(-)\circ q}  
e(T)_*\Gr^{12}_{T*}g^*\mathcal Hom^{\bullet}(F,E_{zar}(\Omega^{\bullet,\Gamma,pr}_{/S},F_{DR}))\otimes_{g^*O_S}O_T \\ 
\xrightarrow{(T(g,hom)(-,-)\otimes I)}
e(T)_*\Gr^{12}_{T*}\mathcal Hom^{\bullet}(g^*F,g^*E_{zar}(\Omega^{\bullet,\Gamma,pr}_{/S},F_{DR}))\otimes_{g^*O_S}O_T \\
\xrightarrow{ev(hom,\otimes)(-,-,-)} 
e(T)_*\Gr^{12}_{T*}\mathcal Hom^{\bullet}(g^*F,g^*E_{zar}(\Omega^{\bullet,\Gamma,pr}_{/S},F_{DR}))
\otimes_{g^*e(S)^*O_S}e(T)^*O_T \\
\xrightarrow{\mathcal Hom^{\bullet}(g^*F,(E\Omega^{\Gamma,pr}_{/(T/S)}\otimes m))}
e(T)_*\Gr^{12}_{T*}\mathcal Hom^{\bullet}(g^*F,E_{zar}(\Omega^{\bullet,\Gamma,pr}_{/T},F_{DR}))
\end{eqnarray*}
\end{defi}

The complex of presheaves 
$(\Omega^{\bullet,\Gamma,pr}_{/S},F_{DR})\in C_{D_Sfil}(\Var(k)^{2,smpr}/S)$ 
have a monoidal structure given by the wedge product of differential forms: 
for $p:(Y\times S,Z)\to S\in\Var(k)^{2,smpr}/S$, the map 
\begin{eqnarray*}
DR(-)(\gamma^{\vee,Hdg}_Z(-))\circ w_{Y\times S/S}:
(\Omega^{\bullet}_{Y\times S/S}\otimes_{O_{Y\times S}}(O_{Y\times S},F_b))\otimes_{p^*O_S}
(\Omega^{\bullet}_{Y\times S/S}\otimes_{O_{Y\times S}}(O_{Y\times S},F_b)) \\
\to\Omega^{\bullet}_{Y\times S/S}\otimes_{O_{Y\times S}}\Gamma_Z^{\vee,Hdg}(O_{Y\times S},F_b)
\end{eqnarray*}
factors trough
\begin{eqnarray*}
DR(-)(\gamma^{\vee,Hdg}_Z(-))\circ w_{Y\times S/S}: \\
(\Omega^{\bullet}_{Y\times S/S}\otimes_{O_{Y\times S}}(O_{Y\times S},F_b))\otimes_{p^*O_S}
(\Omega^{\bullet}_{Y\times S/S}\otimes_{O_{Y\times S}}(O_{Y\times S},F_b)) \\
\xrightarrow{DR(-)(\gamma_Z^{\vee,Hdg}(-))\otimes DR(-)(\gamma_Z^{\vee,Hdg}(-))} \\
(\Omega^{\bullet}_{Y\times S/S}\otimes_{O_{Y\times S}}\Gamma_Z^{\vee,Hdg})(O_{Y\times S},F_b)\otimes_{p^*O_S}
\Omega^{\bullet}_{Y\times S/S}\otimes_{O_{Y\times S}}\Gamma_Z^{\vee,Hdg}(O_{Y\times S},F_b) \\
\xrightarrow{(DR(-)(\gamma^{\vee,Hdg}_Z(-))\circ w_{Y\times S/S})^{\gamma}}
\Omega^{\bullet}_{Y\times S/S}\otimes_{O_{Y\times S}}\Gamma_Z^{\vee,Hdg}(O_{Y\times S},F_b)
\end{eqnarray*}
unique up to homotopy, giving the map in $C_{D_Sfil}(\Var(k)^{2,smpr}/S)$:
\begin{eqnarray*}
w_S:(\Omega^{\bullet,\Gamma,pr}_{/S},F_{DR})\otimes_{O_S}(\Omega^{\bullet,\Gamma,pr}_{/S},F_{DR})
\to(\Omega^{\bullet,\Gamma,pr}_{/S},F_{DR})
\end{eqnarray*}
given by for $p:(Y\times S,Z)\to S\in\Var(k)^{2,smpr}/S$,
\begin{eqnarray*}
w_S((Y\times S,Z)/S): \\
(((\Omega^{\bullet}_{Y\times S/S}\otimes_{O_{Y\times S}}\Gamma_Z^{\vee,Hdg})(O_{Y\times S},F_b))\otimes_{p^*O_S}
(\Omega^{\bullet}_{Y\times S/S}\otimes_{O_{Y\times S}}\Gamma_Z^{\vee,Hdg}(O_{Y\times S},F_b)))(Y\times S) \\
\xrightarrow{(DR(-)(\gamma^{\vee,Hdg}_Z(-))\circ w_{Y\times S/S})^{\gamma}(Y\times S)}
(\Omega^{\bullet}_{Y\times S/S}\otimes_{O_{Y\times S}}\Gamma_Z^{\vee,Hdg}(O_{Y\times S},F_b))(Y\times S)
\end{eqnarray*}
which induces the map in $C_{D_Sfil}(\Var(k)^{2,smpr}/S)$
\begin{eqnarray*}
Ew_S:E_{zar}(\Omega^{\bullet,\Gamma,pr}_{/S},F_{DR})\otimes_{O_S} E_{zar}(\Omega^{\bullet,\Gamma,pr}_{/S},F_{DR})
\xrightarrow{=} \\
E_{zar}((\Omega^{\bullet,\Gamma,pr}_{/S},F_{DR})\otimes_{O_S}(\Omega^{\bullet,\Gamma,pr}_{/S},F_{DR}))
\xrightarrow{E_{zar}(w_S)} E_{et}(\Omega^{\bullet,\Gamma,pr}_{/S},F_{DR})
\end{eqnarray*}
by the functoriality of the Godement resolution (see section 2).

\begin{defi}\label{TotimesDR}
Let $S\in\SmVar(k)$.
We have, for $F,G\in C(\Var(k)^{2,smpr}/S)$, the canonical transformation in $C_{\mathcal Dfil}(S)$ :
\begin{eqnarray*}
T(\otimes,\Omega)(F,G): \\
e(S)_*\Gr^{12}_{S*}\mathcal Hom(F,E_{zar}(\Omega^{\bullet,\Gamma,pr}_{/S},F_{DR}))\otimes_{O_S}
e(S)_*\Gr^{12}_{S*}\mathcal Hom(G,E_{zar}(\Omega^{\bullet,\Gamma,pr}_{/S},F_{DR})) \\
\xrightarrow{=} 
e(S)_*\Gr^{12}_{S*}(\mathcal Hom(F,E_{zar}(\Omega^{\bullet,\Gamma,pr}_{/S},F_{DR}))\otimes_{O_S}
\mathcal Hom(G,E_{zar}(\Omega^{\bullet,\Gamma,pr}_{/S},F_{DR}))) \\ 
\xrightarrow{T(\mathcal Hom,\otimes)(-)} 
e(S)_*\Gr^{12}_{S*}\mathcal Hom(F\otimes G,E_{zar}(\Omega^{\bullet,\Gamma,pr}_{/S},F_{DR})\otimes_{O_S}
E_{zar}(\Omega^{\bullet,\Gamma,pr}_{/S},F_{DR}))) \\
\xrightarrow{=} 
e(S)_*\Gr^{12}_{S*}\mathcal Hom(F\otimes G,(E_{zar}(\Omega^{\bullet,\Gamma,pr}_{/S},F_{DR})
\otimes_{O_S}E_{zar}(\Omega^{\bullet,\Gamma,pr}_{/S},F_{DR}))) \\
\xrightarrow{\mathcal Hom(F\otimes G,Ew_S)} 
e(S)_*\Gr^{12}_{S*}\mathcal Hom(F\otimes G,E_{zar}(\Omega^{\bullet,\Gamma,pr}_{/S},F_{DR})).
\end{eqnarray*}
\end{defi}

Let $S\in\Var(k)$. Let $S=\cup_{i=1}^l S_i$ an open affine cover and denote by $S_I=\cap_{i\in I} S_i$.
Let $i_i:S_i\hookrightarrow\tilde S_i$ closed embeddings, with $\tilde S_i\in\Var(k)$. 
For $I\subset\left[1,\cdots l\right]$, denote by $\tilde S_I=\Pi_{i\in I}\tilde S_i$.
We then have closed embeddings $i_I:S_I\hookrightarrow\tilde S_I$ and for $J\subset I$ the following commutative diagram
\begin{equation*}
D_{IJ}=\xymatrix{ S_I\ar[r]^{i_I} & \tilde S_I \\
S_J\ar[u]^{j_{IJ}}\ar[r]^{i_J} & \tilde S_J\ar[u]^{p_{IJ}}}  
\end{equation*}
where $p_{IJ}:\tilde S_J\to\tilde S_I$ is the projection
and $j_{IJ}:S_J\hookrightarrow S_I$ is the open embedding so that $j_I\circ j_{IJ}=j_J$.
This gives the diagram of algebraic varieties $(\tilde S_I)\in\Fun(\mathcal P(\mathbb N),\Var(k))$ which
the diagram of sites $\Var(k)^{2,smpr}/(\tilde S_I)\in\Fun(\mathcal P(\mathbb N),\Cat)$. 
This gives also the diagram of algebraic varieties $(\tilde S_I)^{op}\in\Fun(\mathcal P(\mathbb N)^{op},\Var(k))$ which
the diagram of sites $\Var(k)^{2,smpr}/(\tilde S_I)^{op}\in\Fun(\mathcal P(\mathbb N)^{op},\Cat)$. 
We then get
\begin{eqnarray*}
((\Omega^{\bullet,\Gamma,pr}_{/(\tilde S_I)},F_{DR})[-d_{\tilde S_I}],T_{IJ})
\in C_{D_{(\tilde S_I)}fil}(\Var(k)^{2,smpr}/(\tilde S_I))
\end{eqnarray*}
with
\begin{eqnarray*}
T_{IJ}:(\Omega^{\bullet,\Gamma,pr}_{/\tilde S_I},F_{DR})[-d_{\tilde S_I}]
\xrightarrow{\ad(p_{IJ}^{*mod[-]},p_{IJ*}(-)}
p_{IJ*}p_{IJ}^*(\Omega^{\bullet,\Gamma,pr}_{/\tilde S_I},F_{DR})
\otimes_{p_{IJ}^*O_{\tilde S_I}}O_{\tilde S_J}[-d_{\tilde S_J}] \\
\xrightarrow{m\circ p_{IJ*}\Omega^{\Gamma,pr}_{/(\tilde S_J/\tilde S_I)}[-d_{\tilde S_J}]}
p_{IJ*}(\Omega^{\bullet,\Gamma,pr}_{/\tilde S_J},F_{DR})[-d_{\tilde S_J}].
\end{eqnarray*}
For $(G_I,K_{IJ})\in C(\Var(k)^{2,smpr}/(\tilde S_I)^{op})$, we denote (see section 2)
\begin{eqnarray*}
e'((\tilde S_I))_*\mathcal Hom((G_I,K_{IJ}),(E_{zar}(\Omega^{\bullet,\Gamma,pr}_{/(\tilde S_I)},F_{DR})[-d_{\tilde S_I}],T_{IJ})):= \\
(e'(\tilde S_I)_*\mathcal Hom(G_I,E_{zar}(\Omega^{\bullet,\Gamma,pr}_{/\tilde S_I},F_{DR}))[-d_{\tilde S_I}],u_{IJ}((G_I,K_{IJ})))
\in C_{\mathcal Dfil}((\tilde S_I))
\end{eqnarray*}
with
\begin{eqnarray*}
u_{IJ}((G_I,K_{IJ})):e'(\tilde S_I)_*\mathcal Hom(G_I,E_{zar}(\Omega^{\bullet,\Gamma,pr}_{/\tilde S_I},F_{DR}))[-d_{\tilde S_I}] \\
\xrightarrow{\ad(p_{IJ}^{*mod[-]},p_{IJ*})(-)\circ T(p_{IJ},e)(-)}
p_{IJ*}e'(\tilde S_J)_*p_{IJ}^*\mathcal Hom(G_I,E_{zar}(\Omega^{\bullet,\Gamma,pr}_{/\tilde S_I},F_{DR}))
\otimes_{p_{IJ}^*O_{\tilde S_I}}O_{\tilde S_J}[-d_{\tilde S_J}] \\
\xrightarrow{T(p_{IJ},hom)(-,-)}
p_{IJ*}e'(\tilde S_J)_*\mathcal Hom(p_{IJ}^*G_I,p_{IJ}^*E_{zar}(\Omega^{\bullet,\Gamma,pr}_{/\tilde S_I},F_{DR})) 
\otimes_{p_{IJ}^*O_{\tilde S_I}}O_{\tilde S_J}[-d_{\tilde S_J}] \\
\xrightarrow{m\circ\mathcal Hom(p_{IJ}^*G_I,T_{IJ})}
p_{IJ*}e'(\tilde S_J)_*\mathcal Hom(p_{IJ}^*G_I,E_{zar}(\Omega^{\bullet,\Gamma,pr}_{/\tilde S_J},F_{DR}))[-d_{\tilde S_J}] \\
\xrightarrow{\mathcal Hom(K_{IJ},E_{zar}(\Omega^{\bullet,\Gamma,pr}_{/\tilde S_J},F_{DR}))}
p_{IJ*}e'(\tilde S_J)_*\mathcal Hom(G_J,E_{zar}(\Omega^{\bullet,\Gamma,pr}_{/\tilde S_J},F_{DR}))[-d_{\tilde S_J}].
\end{eqnarray*}
This gives in particular
\begin{eqnarray*}
(\Omega^{\bullet,\Gamma,pr}_{/(\tilde S_I)},F_{DR})[-d_{\tilde S_I}],T_{IJ})\in 
C_{D_{(\tilde S_I)}fil}(\Var(k)^{2,(sm)pr}/(\tilde S_I)^{op}).
\end{eqnarray*}

We now define the filtered De Rahm realization functor.

\begin{defi}\label{DRalgdefFunct}
\begin{itemize}
\item[(i)]Let $S\in\SmVar(k)$.
We have, using definition \ref{wtildew} and definition \ref{RCHhatdef}, the functor
\begin{eqnarray*}
\mathcal F_S^{FDR}:C(\Var(k)^{sm}/S)\to C_{\mathcal Dfil}(S), \; F\mapsto \\ 
\mathcal F_S^{FDR}(F):=e(S)_*\Gr^{12}_{S*}\mathcal Hom^{\bullet}(\hat R^{CH}(\rho_S^*L(F)),
E_{zar}(\Omega^{\bullet,\Gamma,pr}_{/S},F_{DR}))[-d_S]
\end{eqnarray*}
Moreover, the differentials of $\mathcal F_S^{FDR}$ are strict for the filtration by theorem \ref{Sa12}.
\item[(ii)]Let $S\in\Var(k)$ and $S=\cup_{i=1}^l S_i$ an open cover such that there exist closed embeddings
$i_i:S_i\hookrightarrow\tilde S_i$  with $\tilde S_i\in\SmVar(k)$. 
For $I\subset\left[1,\cdots l\right]$, denote by $S_I:=\cap_{i\in I} S_i$ and $j_I:S_I\hookrightarrow S$ the open embedding.
We then have closed embeddings $i_I:S_I\hookrightarrow\tilde S_I:=\Pi_{i\in I}\tilde S_i$.
Consider, for $I\subset J$, the following commutative diagram
\begin{equation*}
D_{IJ}=\xymatrix{ S_I\ar[r]^{i_I} & \tilde S_I \\
S_J\ar[u]^{j_{IJ}}\ar[r]^{i_J} & \tilde S_J\ar[u]^{p_{IJ}}}  
\end{equation*}
and $j_{IJ}:S_J\hookrightarrow S_I$ is the open embedding so that $j_I\circ j_{IJ}=j_J$.
We have, using definition \ref{wtildew} and definition \ref{RCHhatdef}, the functor
\begin{eqnarray*}
\mathcal F_S^{FDR}:C(\Var(k)^{sm}/S)\to C_{\mathcal Dfil}(S/(\tilde S_I)), \; F\mapsto \\
\mathcal F_S^{FDR}(F):=
e'((\tilde S_I))_*\mathcal Hom^{\bullet}(
(\hat R^{CH}(\rho_{\tilde S_I}^*L(i_{I*}j_I^*F)),\hat R^{CH}_{\tilde S_J}(T^q(D_{IJ})(j_I^*F))), \\
(E_{zar}(\Omega^{\bullet,\Gamma,pr}_{/(\tilde S_I)},F_{DR})[-d_{\tilde S_I}],T_{IJ})) \\
:=(e'(\tilde S_I)_*\mathcal Hom^{\bullet}(\hat R^{CH}(\rho_{\tilde S_I}^*L(i_{I*}j_I^*F)), 
E_{zar}(\Omega^{\bullet,\Gamma,pr}_{/\tilde S_I},F_{DR}))[-d_{\tilde S_I}],u^q_{IJ}(F))
\end{eqnarray*}
where we have denoted for short $e'(\tilde S_I)=e(\tilde S_I)\circ\Gr^{12}_{\tilde S_I}$, and
\begin{eqnarray*}
u^q_{IJ}(F)[d_{\tilde S_J}]:
e'(\tilde S_I)_*\mathcal Hom^{\bullet}(\hat R^{CH}(\rho_{\tilde S_I}^*L(i_{I*}j_I^*F)),
E_{zar}(\Omega^{\bullet,\Gamma,pr}_{/\tilde S_I},F_{DR})) \\
\xrightarrow{\ad(p_{IJ}^{*mod},p_{IJ})(-)} 
p_{IJ*}p_{IJ}^{*mod}e'(\tilde S_I)_*\mathcal Hom^{\bullet}(\hat R^{CH}(\rho_{\tilde S_I}^*L(i_{I*}j_I^*F)),
E_{zar}(\Omega^{\bullet,\Gamma,pr}_{/\tilde S_I},F_{DR})) \\
\xrightarrow{p_{IJ*}T(p_{IJ},\Omega^{\gamma,pr}_{\cdot})(-)}  
p_{IJ*}e'(\tilde S_J)_*\mathcal Hom^{\bullet}(p_{IJ}^*\hat R^{CH}(\rho_{\tilde S_I}^*L(i_{I*}j_I^*F)),
E_{zar}(\Omega^{\bullet,\Gamma,pr}_{/\tilde S_J},F_{DR})) \\
\xrightarrow{\mathcal Hom(T(p_{IJ},\hat R^{CH})(Li_{I*}j_I^*F)^{-1},E_{et}(\Omega_{/\tilde S_J}^{\bullet,\Gamma,pr},F_{DR}))} \\
p_{IJ*}e'(\tilde S_J)_*\mathcal Hom^{\bullet}(\hat R^{CH}(\rho_{\tilde S_J}^*p_{IJ}^*L(i_{I*}j_I^*F)),
E_{zar}(\Omega^{\bullet,\Gamma,pr}_{/\tilde S_J},F_{DR})) \\
\xrightarrow{\mathcal Hom(\hat R^{CH}_{\tilde S_J}(T^q(D_{IJ})(j_I^*F)),E_{et}(\Omega_{/\tilde S_J}^{\bullet,\Gamma,pr},F_{DR}))} \\
p_{IJ*}e'(\tilde S_J)_*\mathcal Hom^{\bullet}(\hat R^{CH}(\rho_{\tilde S_J}^*L(i_{J*}j_J^*F)),
E_{zar}(\Omega^{\bullet,\Gamma,pr}_{/\tilde S_J},F_{DR})).
\end{eqnarray*}
For $I\subset J\subset K$, we have obviously $p_{IJ*}u_{JK}(F)\circ u_{IJ}(F)=u_{IK}(F)$.
Moreover, the differentials of $\mathcal F_S^{FDR}$ are strict for the filtration by theorem \ref{Sa12}.
\end{itemize}
\end{defi}

Recall, see section 2, that we have the projection morphisms of sites 
$p_a:\Var(k)^{2,smpr}/(\tilde S_I)^{op}\to\Var(k)^{2,smpr}/(\tilde S_I)^{op}$
given by the functor 
\begin{eqnarray*}
p_a:\Var(k)^{2,smpr}/(\tilde S_I)^{op}\to\Var(k)^{2,smpr}/(\tilde S_I)^{op}, \\
p_a((Y_I\times\tilde S_I,Z_I)/\tilde S_I,s_{IJ}):=
((Y_I\times\mathbb A^1\times\tilde S_I,Z_I\times\mathbb A^1)/\tilde S_I,s_{IJ}\times I), \\ 
p_a((g_I):((Y'_I\times\tilde S_I,Z'_I)/\tilde S_I,s'_{IJ})\to((Y_I\times\tilde S_I,Z_I)/\tilde S_I,s_{IJ}))= \\
(g_I\times I):((Y'_I\times\mathbb A^1\times\tilde S_I,Z'_I\times\mathbb A^1)/\tilde S_I,s'_{IJ}\times I)
\to((Y_I\times\mathbb A^1\times\tilde S_I,Z_I\times\mathbb A^1)/\tilde S_I),s_{IJ}\times I)).
\end{eqnarray*}

We have the following key proposition :

\begin{prop}\label{aetfib}
\begin{itemize}
\item[(i1)]Let $S\in\Var(k)$. Let $S=\cup_{i=1}^l S_i$ an open cover such that there exist closed embeddings
$i_i:S_i\hookrightarrow\tilde S_i$ with $\tilde S_i\in\SmVar(k)$.
The complex of presheaves
$(\Omega^{\bullet,\Gamma,pr}_{/(\tilde S_I)},F_{DR})\in C_{D_{(\tilde S_I)}fil}(\Var(k)^{2,smpr}/(\tilde S_I)^{op})$ 
is $2$-filtered $\mathbb A^1$ homotopic, that is
\begin{equation*}
\ad(p_a^*,p_{a*})(\Omega^{\bullet,\Gamma,pr}_{/(\tilde S_I)},F_{DR}):
(\Omega^{\bullet,\Gamma,pr}_{/S},F_{DR})\to p_{a*}p_a^*(\Omega^{\bullet,\Gamma,pr}_{/(\tilde S_I)},F_{DR})
\end{equation*}
is a $2$-filtered homotopy.
\item[(i2)]Let $S\in\SmVar(k)$. The complex of presheaves
$(\Omega^{\bullet,\Gamma,pr}_{/S},F_{DR})\in C_{D_Sfil}(\Var(k)^{2,smpr}/S)$ 
admits transferts, i.e. 
\begin{equation*}
\Tr(S)_*\Tr(S)^*(\Omega^{\bullet,\Gamma,pr}_{/S},F_{DR})=(\Omega^{\bullet,\Gamma,pr}_{/S},F_{DR}).
\end{equation*}
\item[(iii)]Let $S\in\Var(k)$. Let $S=\cup_{i=1}^l S_i$ an open cover such that there exist closed embeddings
$i_i:S_i\hookrightarrow\tilde S_i$ with $\tilde S_i\in\SmVar(k)$.
Let $m=(m_I):(Q_{1I},K^1_{IJ})\to(Q_{2I},K^2_{IJ})$ be an equivalence $(\mathbb A^1,et)$ local with 
$(Q_{1I},K_{IJ})\to(Q_{2I},K_{IJ})\in C(\Var(k)^{smpr}/(\tilde S_I)^{op})$ complexes of representable presheaves. 
Then, the map in $C_{\mathcal Dfil}((\tilde S_I))$
\begin{eqnarray*} 
M:=(e(\tilde S_I)_*\mathcal Hom^{\bullet}(m_I,E_{zar}(\Omega^{\bullet,\Gamma,pr}_{/\tilde S_I},F_{DR})[-d_{\tilde S_I}])): \\  
e'((\tilde S_I))_*\mathcal Hom^{\bullet}((Q_{2I},K^1_{IJ}),
(E_{zar}(\Omega^{\bullet,\Gamma,pr}_{/(\tilde S_I)},F_{DR})[-d_{\tilde S_I}],T_{IJ})) \\
\to e'((\tilde S_I))_*\mathcal Hom^{\bullet}((Q_{1I},K^1_{IJ}),
(E_{zar}(\Omega^{\bullet,\Gamma,pr}_{/(\tilde S_I)},F_{DR})[-d_{\tilde S_I}],T_{IJ}))
\end{eqnarray*}
is a $2$-filtered quasi-isomorphism. It is thus an isomorphism in $D_{\mathcal Dfil,\infty}((\tilde S_I))$.
\end{itemize}
\end{prop}

\begin{proof}
\noindent(i1): Similar to the proof of \cite{B4}, proposition (ii1)

\noindent(i2): Similar to the proof of \cite{B4}, proposition (ii2) :
Let $\alpha\in\Cor(\Var(\mathbb C)^{2,smpr}/S)((Y_1\times S,Z_1)/S,(Y_2\times S,Z_2)/S)$ irreducible.
Denote by $i:\alpha\hookrightarrow Y_1\times Y_2\times S$ the closed embedding, 
and $p_1:Y_1\times Y_2\times S\to Y_1\times S$, $p_2:Y_1\times Y_2\times S\to Y_2\times S$ the projections.
The morphism $p_1\circ i:\alpha\to Y_1\times S$ is then finite surjective and 
$(Z_1\times Y_2)\cap\alpha\subset Y_1\times Z_2$ (i.e. $p_2(p_1^{-1}(Z_2)\cap\alpha)\subset Z_2$).
Then, the transfert map is given by
\begin{eqnarray*}
\Omega^{\bullet,\Gamma,pr}_{/S}(\alpha):
((\Omega^{\bullet}_{Y_2\times S/S},F_b)\otimes_{O_{Y_2\times S}}
\Gamma_{Z_2}^{\vee,Hdg}(O_{Y_2\times S},F_b))(Y_2\times S) \\
\xrightarrow{i_{-}}
p_2^*((\Omega^{\bullet}_{Y_2\times S/S},F_b)\otimes_{O_{Y_2\times S}}
\Gamma_{Z_2}^{\vee,Hdg}(O_{Y_2\times S},F_b))(Y_1\times Y_2\times S) \\
\xrightarrow{\Omega_{(Y_1\times Y_2\times S/Y_2\times S)/(S/S)}(-)(-)}
((\Omega^{\bullet}_{Y_1\times Y_2\times S/S},F_b)\otimes_{O_{Y_1\times Y_2\times S}}
\Gamma_{Y_1\times Z_2}^{\vee,Hdg}(O_{Y_1\times Y_2\times S},F_b))(Y_1\times Y_2\times S) \\
\xrightarrow{DR(-)(T((Z_1\times Y_2)\cap\alpha/Y_1\times Z_2,\gamma^{\vee,Hdg})(-)(-)} \\
((\Omega^{\bullet}_{Y_1\times Y_2\times S/S},F_b)\otimes_{O_{Y_1\times Y_2\times S}}
\Gamma_{(Z_1\times Y_2)\cap\alpha}^{\vee,Hdg}(O_{Y_1\times Y_2\times S},F_b))(Y_1\times Y_2\times S) \\
\xrightarrow{i_{-}}
i^*((\Omega^{\bullet}_{Y_1\times Y_2\times S/S},F_b)\otimes_{O_{Y_1\times Y_2\times S}}
\Gamma_{(Z_1\times Y_2)\cap\alpha}^{\vee,Hdg}(O_{Y_1\times Y_2\times S},F_b))(\alpha) \\
\xrightarrow{\Omega_{(\alpha/Y_1\times Y_2\times S)/(S/S)}(-)(-)}
((\Omega^{\bullet}_{\alpha/S},F_b)\otimes_{O_{\alpha}}
i^{*mod}\Gamma_{(Z_1\times Y_2)\cap\alpha}^{\vee,Hdg}(O_{Y_1\times Y_2\times S},F_b))(\alpha) \\
\xrightarrow{\Omega_{(\alpha/Y_1\times S)(S/S)}(-)(-)^{tr}}
((\Omega^{\bullet}_{Y_1\times S/S},F_b)\otimes_{O_{Y_1\times S}}
\Gamma_{Z_1}^{\vee,Hdg}(O_{Y_1\times S},F_b))(Y_1\times S).
\end{eqnarray*}

\noindent(ii):Follows from (i) and theorem \ref{DDADM12fil}.
\end{proof}

\begin{prop}\label{projwach}
Let $S\in\Var(k)$.Let $S=\cup_{i=1}^l S_i$ an open cover such that there exist closed embeddings
$i_i:S_i\hookrightarrow\tilde S_i$ with $\tilde S_i\in\SmVar(k)$. 
\begin{itemize}
\item[(i)]Let $m=(m_I):(Q_{1I},K^1_{IJ})\to (Q_{2I},K^2_{IJ})$ be an etale local equivalence local 
with $(Q_{1I},K^1_{IJ}),(Q_{2I},K^2_{IJ})\in C(\Var(k)^{sm}/(\tilde S_I))$
complexes of projective presheaves. Then,
\begin{eqnarray*} 
(e'(\tilde S_I)_*\mathcal Hom^{\bullet}(\hat R_S^{CH}(m_I),
E_{zar}(\Omega^{\bullet,\Gamma,pr}_{/\tilde S_I},F_{DR}))[-d_{\tilde S_I}]): \\ 
e'(\tilde S_I)_*\mathcal Hom^{\bullet}((\hat R^{CH}(\rho_S^*Q_{1I}),\hat R^{CH}(K^1_{IJ})),
(E_{zar}(\Omega^{\bullet,\Gamma,pr}_{/\tilde S_I},F_{DR})[-d_{\tilde S_I}],T_{IJ})) \\
\to e'(\tilde S_I)_*\mathcal Hom^{\bullet}((\hat R^{CH}(\rho_S^*Q_{2I}),\hat R^{CH}(K^2_{IJ})),
(E_{zar}(\Omega^{\bullet,\Gamma,pr}_{/\tilde S_I},F_{DR})[-d_{\tilde S_I}],T_{IJ}))
\end{eqnarray*}
is a filtered quasi-isomorphism. It is thus an isomorphism in $D_{\mathcal Dfil}((\tilde S_I))$.
\item[(ii)]Let $m=(m_I):(Q_{1I},K^1_{IJ})\to (Q_{2I},K^2_{IJ})$ be an equivalence $(\mathbb A^1,et)$ local equivalence local 
with $(Q_{1I},K^1_{IJ}),(Q_{2I},K^2_{IJ})\in C(\Var(k)^{sm}/(\tilde S_I))$
complexes of projective presheaves. Then,
\begin{eqnarray*} 
(e'(\tilde S_I)_*\mathcal Hom^{\bullet}(\hat R_S^{CH}(m_I),
E_{zar}(\Omega^{\bullet,\Gamma,pr}_{/\tilde S_I},F_{DR}))[-d_{\tilde S_I}]): \\ 
e'(\tilde S_I)_*\mathcal Hom^{\bullet}((\hat R^{CH}(\rho_S^*Q_{1I}),\hat R^{CH}(K^1_{IJ})),
(E_{zar}(\Omega^{\bullet,\Gamma,pr}_{/\tilde S_I},F_{DR})[-d_{\tilde S_I}],T_{IJ})) \\
\to e'(\tilde S_I)_*\mathcal Hom^{\bullet}((\hat R^{CH}(\rho_S^*Q_{2I}),\hat R^{CH}(K^2_{IJ})),
(E_{zar}(\Omega^{\bullet,\Gamma,pr}_{/\tilde S_I},F_{DR})[-d_{\tilde S_I}],T_{IJ}))
\end{eqnarray*}
is a filtered quasi-isomorphism. It is thus an isomorphism in $D_{\mathcal Dfil}((\tilde S_I))$.
\end{itemize}
\end{prop}

\begin{proof}
Follows from proposition \ref{aetfib} (see the proof the complex case in \cite{B4} section 6 )
and the fact that the differential of the complexes involved are strict for the F-fitration.
\end{proof}

\begin{defi}\label{DRalgdefsing}
\begin{itemize}
\item[(i)] Let $S\in\SmVar(k)$.
We define using definition \ref{DRalgdefFunct}(i) and proposition \ref{projwach}(ii)
the filtered algebraic De Rahm realization functor defined as
\begin{eqnarray*}
\mathcal F_S^{FDR}:\DA_c(S)\to D_{\mathcal Dfil}(S), \; M\mapsto \\ 
\mathcal F_S^{FDR}(M):=e(S)_*\Gr^{12}_{S*}\mathcal Hom^{\bullet}(\hat R^{CH}(\rho_S^*L(F)),
E_{zar}(\Omega^{\bullet,\Gamma,pr}_{/S},F_{DR}))[-d_S] 
\end{eqnarray*}
where $F\in C(\Var(k)^{sm}/S)$ is such that $M=D(\mathbb A^1,et)(F)$.
\item[(i)'] For the Corti-Hanamura weight structure $W$ on $\DA_c(S)^-$, 
we define using definition \ref{DRalgdefFunct}(i) and proposition \ref{projwach}(ii)
\begin{eqnarray*}
\mathcal F_S^{FDR}:\DA_c^-(S)\to D_{\mathcal D(1,0)fil}^-(S), \; M\mapsto \\
\mathcal F_S^{FDR}((M,W)):= e(S)_*\Gr^{12}_{S*}\mathcal Hom^{\bullet}(\hat R^{CH}(\rho_S^*L(F,W)),
E_{zar}(\Omega^{\bullet,\Gamma,pr}_{/S},F_{DR}))[-d_S] 
\end{eqnarray*}
where $(F,W)\in C_{fil}(\Var(k)^{sm}/S)$ 
is such that $M=D(\mathbb A^1,et)((F,W))$ using corollary \ref{weightst2Cor}. 
Note that the filtration induced by $W$ is a filtration by sub $D_S$ module,
which is a stronger property then Griffitz transversality.
Of course, the filtration induced by $F$ satisfy only Griffitz transversality in general.
\item[(ii)]Let $S\in\Var(k)$ and $S=\cup_{i=1}^l S_i$ an open cover such that there exist closed embeddings
$i_i:S_i\hookrightarrow\tilde S_i$  with $\tilde S_i\in\SmVar(k)$. 
For $I\subset\left[1,\cdots l\right]$, denote by $S_I=\cap_{i\in I} S_i$ and $j_I:S_I\hookrightarrow S$ the open embedding.
We then have closed embeddings $i_I:S_I\hookrightarrow\tilde S_I:=\Pi_{i\in I}\tilde S_i$.
We define, using definition \ref{DRalgdefFunct}(ii) and proposition \ref{projwach}(ii),
the filtered algebraic De Rahm realization functor defined as
\begin{eqnarray*}
\mathcal F_S^{FDR}:\DA_c(S)\to D_{\mathcal Dfil}(S/(\tilde S_I)), \; M\mapsto \\
\mathcal F_S^{FDR}(M):= 
(e'(\tilde S_I)_*\mathcal Hom^{\bullet}(\hat R^{CH}(\rho_{\tilde S_I}^*L(i_{I*}j_I^*F)),
E_{zar}(\Omega^{\bullet,\Gamma,pr}_{/\tilde S_I},F_{DR}))[-d_{\tilde S_I}],u^q_{IJ}(F))
\end{eqnarray*}
where $F\in C(\Var(k)^{sm}/S)$ is such that $M=D(\mathbb A^1,et)(F)$, 
see definition \ref{DRalgdefFunct}.
\item[(ii)'] For the Corti-Hanamura weight structure $W$ on $\DA_c^-(S)$, 
using definition \ref{DRalgdefFunct}(ii) and proposition \ref{projwach}(ii), 
\begin{eqnarray*}
\mathcal F_S^{FDR}:\DA_c^-(S)\to D_{\mathcal D(1,0)fil}^-(S/(\tilde S_I)), \; 
M\mapsto\mathcal F_S^{FDR}((M,W)):= \\ 
(e'(\tilde S_I)_*\mathcal Hom^{\bullet}(\hat R^{CH}(\rho_{\tilde S_I}^*L(i_{I*}j_I^*(F,W))),
E_{zar}(\Omega^{\bullet,\Gamma,pr}_{/\tilde S_I},F_{DR}))[-d_{\tilde S_I}],u^q_{IJ}(F,W))
\end{eqnarray*}
where $(F,W)\in C_{fil}(\Var(k)^{sm}/S)$ 
is such that $(M,W)=D(\mathbb A^1,et)(F,W)$ using corollary \ref{weightst2Cor}.
Note that the filtration induced by $W$ is a filtration by sub $D_{\tilde S_I}$-modules,
which is a stronger property then Griffitz transversality.
Of course, the filtration induced by $F$ satisfy only Griffitz transversality in general.
\end{itemize}
\end{defi}

\begin{prop}\label{FDRwelldef}
For $S\in\Var(k)$ and $S=\cup_{i=1}^l S_i$ an open cover such that there exist closed embeddings
$i_i:S_i\hookrightarrow\tilde S_i$ with $\tilde S_i\in\SmVar(k)$, the functor $\mathcal F_S^{FDR}$ is well defined. 
\end{prop}

\begin{proof}
Similar to the proof of \cite{B4} proposition : follows from proposition \ref{projwach}.
\end{proof}

\begin{rem}
\begin{itemize}
\item[(i)] Let $S\in\SmVar(k)$.
We have, by proposition \ref{aetfib}, for $M\in\DA_c(S)$ the isomorphism in $D_{\mathcal D(1,0)fil,\infty}^-(S)$
\begin{eqnarray*}
\mathcal Hom(-,k)\circ\mathcal Hom(T(\hat R^{CH},R^{CH})(\rho_S^*L(F,W)),-)^{-1}: \\
\mathcal F_S^{FDR}((M,W)):= 
e'(S)_*\mathcal Hom^{\bullet}(\hat R^{CH}(\rho_S^*L(F,W)),E_{zar}(\Omega^{\bullet,\Gamma,pr}_{/S},F_{DR}))[-d_S] \\
\xrightarrow{\sim}
e'(S)_*\mathcal Hom^{\bullet}(L\mu_{S*}\rho_{S*}R^{CH}(\rho_S^*L(F,W)),E_{et}(\Omega^{\bullet,\Gamma,pr}_{/S},F_{DR}))[-d_S] 
\end{eqnarray*}
as it was defined in \cite{B4}.
\item[(ii)]Let $S\in\Var(k)$ and $S=\cup_{i=1}^l S_i$ an open cover such that there exist closed embeddings
$i_i:S_i\hookrightarrow\tilde S_i$  with $\tilde S_i\in\SmVar(k)$. 
We have, by proposition \ref{aetfib}, for $M\in\DA_c(S)$ the isomorphism in $D_{\mathcal D(1,0)fil,\infty}^-(S/(\tilde S_I))$
\begin{eqnarray*}
(\mathcal Hom(-,k))\circ(\mathcal Hom(T(\hat R^{CH},R^{CH})(\rho_{\tilde S_I}^*L(i_{I*}j_I^*(F,W))),-))^{-1}: \\
\mathcal F_S^{FDR}((M,W)):= 
(e'(\tilde S_I)_*\mathcal Hom^{\bullet}(\hat R^{CH}(\rho_{\tilde S_I}^*L(i_{I*}j_I^*(F,W))),
E_{zar}(\Omega^{\bullet,\Gamma,pr}_{/\tilde S_I},F_{DR}))[-d_{\tilde S_I}],u^q_{IJ}(F,W)) \\
\xrightarrow{\sim}
(e'(\tilde S_I)_*\mathcal Hom^{\bullet}(L\mu_{\tilde S_I*}\rho_{\tilde S_I*}R^{CH}(\rho_{\tilde S_I}^*L(i_{I*}j_I^*(F,W))),
E_{et}(\Omega^{\bullet,\Gamma,pr}_{/\tilde S_I},F_{DR}))[-d_{\tilde S_I}],u^q_{IJ}(F,W)) 
\end{eqnarray*}
as it was defined in \cite{B4}.
\end{itemize}
\end{rem}

\begin{prop}\label{keyalgsing1}
Let $f:X\to S$ a morphism with $S,X\in\Var(k)$. Assume there exist a factorization 
\begin{equation*}
f:X\xrightarrow{l}Y\times S\xrightarrow{p_S} S
\end{equation*}
of $f$ with $Y\in\SmVar(k)$, $l$ a closed embedding and $p_S$ the projection.
Let $\bar Y\in\PSmVar(k)$ a compactification of $Y$ with $\bar Y\backslash Y=D$ a normal crossing divisor,
denote $k:D\hookrightarrow \bar Y$ the closed embedding and $n:Y\hookrightarrow\bar Y$ the open embedding.
Denote $\bar X\subset\bar Y\times S$ the closure of $X\subset\bar Y\times S$.
We have then the following commutative diagram in $\Var(k)$
\begin{equation*}
\xymatrix{X\ar[r]^l\ar[d] & Y\times S\ar[rd]^{p_S}\ar[d]^{(n\times I)} & \, \\
\bar X\ar[r]^l & \bar Y\times S\ar[r]^{\bar p_S} & S \\
Z:=\bar X\backslash X\ar[ru]^{l_Z}\ar[u]\ar[r] & D\times S\ar[ru]\ar[u]_{(k\times I)} & \, }.
\end{equation*}
Let $S=\cup_{i=1}^l S_i$ an open cover such that there exist closed embeddings
$i_i:S_i\hookrightarrow\tilde S_i$ with $\tilde S_i\in\SmVar(k)$. 
Then $X=\cup_{i=1}^lX_i$ with $X_i:=f^{-1}(S_i)$.
Denote, for $I\subset\left[1,\cdots l\right]$, $S_I=\cap_{i\in I} S_i$ and $X_I=\cap_{i\in I}X_i$.
Denote $\bar X_I:=\bar X\cap(\bar Y\times S_I)\subset\bar Y\times\tilde S_I$ 
the closure of $X_I\subset\bar Y\times\tilde S_I$, 
and $Z_I:=Z\cap(\bar Y\times S_I)=\bar X_I\backslash X_I\subset\bar Y\times\tilde S_I$.
We have then for $I\subset\left[1,\cdots l\right]$, the following commutative diagram in $\Var(k)$
\begin{equation*}
\xymatrix{X_I\ar[r]^{l_I}\ar[d] & Y\times \tilde S_I\ar[rd]^{p_{\tilde S_I}}\ar[d]^{(n\times I)} & \, \\
\bar X_I\ar[r]^{l_I} & \bar Y\times\tilde S_I\ar[r]^{\bar p_{\tilde S_I}} & \tilde S_I \\
Z_I=\bar X_I\backslash X_I\ar[ru]^{l_{Z_I}}\ar[u]\ar[r] & D\times\tilde S_I\ar[ru]\ar[u]_{(k\times I)} & \, }.
\end{equation*}
Let $F(X/S):=p_{S,\sharp}\Gamma_X^{\vee}\mathbb Z(Y\times S/Y\times S)\in C(\Var(k)^{sm}/S)$. 
We have then the following isomorphism in $D_{\mathcal Dfil}(S/(\tilde S_I))$ 
\begin{eqnarray*} 
I(X/S):\mathcal F_S^{FDR}(M(X/S))\xrightarrow{:=} \\
(e'(\tilde S_I)_*\mathcal Hom(\hat R^{CH}(\rho_{\tilde S_I}^*L(i_{I*}j_I^*F(X/S))), 
E_{zar}(\Omega^{\bullet,\Gamma,pr}_{/\tilde S_I},F_{DR}))[-d_{\tilde S_I}],u^q_{IJ}(F(X/S))) \\ 
\xrightarrow{(\mathcal Hom(\hat R^{CH}_{\tilde S_I}(N_I(X/S)), 
E_{zar}(\Omega^{\bullet,\Gamma,pr}_{/\tilde S_I},F_{DR})))} \\
(e'(\tilde S_I)_*\mathcal Hom(\hat R^{CH}(\rho_{\tilde S_I}^*Q(X_I/\tilde S_I)), 
E_{zar}(\Omega^{\bullet,\Gamma,pr}_{/\tilde S_I},F_{DR}))[-d_{\tilde S_I}],v^q_{IJ}(F(X/S))) \\
\xrightarrow{(\mathcal Hom(\rho_{\tilde S_I*}I_{\delta}((\bar X_I,Z_I)/\tilde S_I),-)[-d_{\tilde S_I}])^{-1}} \\
(p_{\tilde S_I*}E_{zar}((\Omega^{\bullet}_{\bar Y\times\tilde S_I/\tilde S_I},F_b)\otimes_{O_{\bar Y\times\tilde S_I}}
(n\times I)_!^{Hdg}\Gamma^{\vee,Hdg}_{X_I}(O_{Y\times\tilde S_I},F_b))(d_Y+d_{\tilde S_I})[2d_Y+d_{\tilde S_I}],w_{IJ}(X/S)) \\ 
\xrightarrow{=:}\iota_SRf^{Hdg}_!(\Gamma^{\vee,Hdg}_{X_I}(O_{Y\times\tilde S_I},F_b)(d_Y)[2d_Y],x_{IJ}(X/S)).
\xrightarrow{=:}\iota_SRf^{Hdg}_!f^{*mod}_{Hdg}\mathbb Z^{Hdg}_S.
\end{eqnarray*}
\end{prop}

\begin{proof}
Similar to the proof of \cite{B4}, proposition
\end{proof}

\begin{cor}\label{FDRMHM}
Let $S\in\Var(k)$ and $S=\cup_{i=1}^l S_i$ an open cover such that there exist closed embeddings
$i_i:S_i\hookrightarrow\tilde S_i$  with $\tilde S_i\in\SmVar(k)$.
Then for $M\in\DA_c(S)$, $\mathcal F_S^{FDR}\in\iota_S(D(DRM(S)))$, 
where $\iota_S:D(DRM(S))\hookrightarrow D_{\mathcal Dfil}(S/(\tilde S_I))$ is a full embedding by theorem \ref{Bek0}.
\end{cor}

\begin{proof}
There exist by definition of constructible motives an isomorphism $\DA(S)$
\begin{equation*}
w(M):M\xrightarrow{\sim}\Cone(M(X_1/S)\to\cdots\to M(X_r/S)).
\end{equation*}
Hence we have the isomorphism in $D_{\mathcal Dfil}(S/(\tilde S_I))$
\begin{equation*}
\mathcal F_S^{FDR}(w(M)):\mathcal F_S^{FDR}(M)\xrightarrow{\sim}
\Cone(\mathcal F_S^{FDR}(M(X_1/S))\to\cdots\to\mathcal F_S^{FDR}(M(X_r/S))).
\end{equation*}
The result then follows from proposition \ref{keyalgsing1}.
\end{proof}

\begin{prop}\label{FDRHdgwelldef}
For $S\in\Var(k)$ not smooth, the functor (see corollary \ref{FDRMHM}) 
\begin{equation*}
\iota_S^{-1}\mathcal F_S^{FDR}:\DA_c^-(S)^{op}\to D(DRM(S))
\end{equation*}
does not depend on the choice of the open cover $S=\cup_iS_i$
and the closed embeddings $i_i:S_i\hookrightarrow\tilde S_i$ with $\tilde S_i\in\SmVar(k)$.
\end{prop}

\begin{proof}
Similar to the proof of \cite{B4} proposition
\end{proof}

We have the canonical transformation map between the filtered De Rham realization functor and the Gauss-Manin realization functor :

\begin{defi}\label{GMFDRdef}
Let $S\in\Var(k)$ and $S=\cup_{i=1}^l S_i$ an open cover such that there exist closed embeddings
$i_i:S_i\hookrightarrow\tilde S_i$ with $\tilde S_i\in\SmVar(k)$. 
Let $M\in\DA_c(S)$ and $F\in C(\Var(k)^{sm}/S)$ such that $M=D(\mathbb A^1,et)(F)$.
We have, using definition \ref{wtildew}(ii), definition \ref{tus}, proposition \ref{grRCH} and proposition \ref{aetfib},
the canonical map in $D_{O_Sfil,\mathcal D,\infty}(S/(\tilde S_I))$
\begin{eqnarray*}
T(\mathcal F^{GM}_S,\mathcal F^{FDR}_S)(M): \\
\mathcal F_S^{GM}(L\mathbb D_SM)=
(e(\tilde S_I)_*\mathcal Hom^{\bullet}(L(i_{I*}j_I^*\mathbb D_SLF),
E_{zar}(\Omega^{\bullet}_{/\tilde S_I},F_b))[-d_{\tilde S_I}],u^q_{IJ}(F)) \\
\xrightarrow{\sim} 
(e(\tilde S_I)_*\mathcal Hom^{\bullet}(L\mathbb D^0_{\tilde S_I}L(i_{I*}j_I^*F),
E_{zar}(\Omega^{\bullet}_{/\tilde S_I},F_b))[-d_{\tilde S_I}],u^q_{IJ}(F)) \\
\xrightarrow{\mathcal Hom(-,\Gr(\Omega_{\tilde S_I}))^{-1}}
(e(\tilde S_I)_*\mathcal Hom^{\bullet}(L\mathbb D^0_{\tilde S_I}L(i_{I*}j_I^*F),
\Gr^{12}_{\tilde S_I*}E_{zar}(\Omega^{\bullet,\Gamma,pr}_{/\tilde S_I},F_{DR}))[-d_{\tilde S_I}],u^q_{IJ}(F)) \\
\xrightarrow{I(\Gr_{\tilde S_I}^{12*},\Gr^{12}_{\tilde S_I*})(-,-)}
(e(\tilde S_I)_*\mathcal Hom^{\bullet}(\Gr_{\tilde S_I}^{12*}L\mathbb D^0_{\tilde S_I}L(i_{I*}j_I^*F),
E_{zar}(\Omega^{\bullet,\Gamma,pr}_{/\tilde S_I},F_{DR}))[-d_{\tilde S_I}],u^q_{IJ}(F)) \\
\xrightarrow{(\mathcal Hom^{\bullet}(r^{CH}(L(i_{I*}j_I^*F)\circ T(\hat R^{CH},R^{CH})(L(i_{I*}j_I^*F)),
E_{zar}(\Omega^{\bullet,\Gamma,pr}_{/\tilde S_I},F_{DR}))[-d_{\tilde S_I}])} \\
(e'(\tilde S_I)_*\mathcal Hom^{\bullet}(\hat R^{CH}(\rho_{\tilde S_I}^*L(i_{I*}j_I^*F)),
E_{zar}(\Omega^{\bullet,\Gamma,pr}_{/\tilde S_I},F_{DR}))[-d_{\tilde S_I}],u^q_{IJ}(F))=:\mathcal F^{FDR}_S(M)
\end{eqnarray*}
\end{defi}

\begin{prop}\label{TGMFDRprop}
Let $S\in\Var(k)$ and $S=\cup_{i=1}^l S_i$ an open cover such that there exist closed embeddings
$i_i:S_i\hookrightarrow\tilde S_i$ with $\tilde S_i\in\SmVar(k)$. 
\begin{itemize}
\item[(i)] For $M\in\DA_c(S)$ the map in $D_{O_S,\mathcal D}(S/(\tilde S_I))=D_{O_S,\mathcal D}(S)$
\begin{equation*}
o_{fil}T(\mathcal F^{GM}_S,\mathcal F^{FDR}_S)(M):
o_{fil}\mathcal F_S^{GM}(L\mathbb D_SM)\xrightarrow{\sim} o_{fil}\mathcal F^{FDR}_S(M)
\end{equation*}
given in definition \ref{GMFDRdef} is an isomorphism if we forgot the Hodge filtration $F$.
\item[(ii)]For $M\in\DA_c(S)$ and all $n,p\in\mathbb Z$, the map in $\PSh_{O_S,\mathcal D}(S/(\tilde S_I))$
\begin{equation*}
F^pH^nT(\mathcal F^{GM}_S,\mathcal F^{FDR}_S)(M):F^pH^n\mathcal F_S^{GM}(L\mathbb D_SM)\hookrightarrow F^pH^n\mathcal F^{FDR}_S(M)
\end{equation*}
given in definition \ref{GMFDRdef} is a monomorphism.
Note that $F^pH^nT(\mathcal F^{GM}_S,\mathcal F^{FDR}_S)(M)$ is NOT an isomorphism in general :
take for example $M(S^o/S)^{\vee}=D(\mathbb A^1,et)(j_*E_{et}(\mathbb Z(S^o/S)))$ 
for an open embedding $j:S^o\hookrightarrow S$, then 
\begin{equation*}
H^n\mathcal F_S^{GM}(L\mathbb D_SM(S^o/S)^{\vee})=\mathcal F_S^{GM}(\mathbb Z(S^o/S))=j_*E(O_{S^o},F_b)\notin\pi_S(MHM(S)) 
\end{equation*}
and hence is NOT isomorphic to $H^n\mathcal F_S^{FDR}(L\mathbb D_SM(S^o/S)^{\vee})\in\pi_S(MHM(S))$ 
as filtered $D_S$-modules (see remark \ref{remHdgkey}). 
It is an isomorphism in the very particular cases where $M=D(\mathbb A^1,et)(\mathbb Z(X/S))$
or $M=D(\mathbb A^1,et)(\mathbb Z(X^o/S))$ for $f:X\to S$ is a smooth proper morphism and $n:X^o\hookrightarrow X$ is an open subset
such that $X\backslash X^o=\cup D_i$ is a normal crossing divisor 
and such that $f_{|D_i}=f\circ i_i:D_i\to X$ are SMOOTH morphism with $i_i:D_i\hookrightarrow X$ the closed embedding and
considering $f_{|X^o}=f\circ n:X^o\to S$ (see \cite{B4} section 6.1 in the complex case).
\end{itemize}
\end{prop}

\begin{proof}
\noindent(i):Follows from the computation for a Borel-Moore motive.

\noindent(ii):Follows from (i).
\end{proof}

We now define the functorialities of $\mathcal F_S^{FDR}$ with respect to $S$ 
which makes $\mathcal F^{-}_{FDR}$ a morphism of 2 functor.

\begin{defi}\label{TGammaFDR}
Let $S\in\Var(k)$. Let $Z\subset S$ a closed subset.
Let $S=\cup_{i=1}^l S_i$ an open cover such that there exist closed embeddings
$i_i:S_i\hookrightarrow\tilde S_i$ with $\tilde S_i\in\SmVar(k)$.
Denote $Z_I:=Z\cap S_I$. We then have closed embeddings $Z_I\hookrightarrow S_I\hookrightarrow\tilde S_I$.
\begin{itemize}
\item[(i)]For $F\in C(\Var(k)^{sm}/S)$, we will consider the following canonical map 
in $D(DRM(S))\subset D_{\mathcal D(1,0)fil}(S/(\tilde S_I))$
\begin{eqnarray*}
T(\Gamma_Z^{\vee,Hdg},\Omega^{\Gamma,pr}_{/S})(F,W):\\
\Gamma_Z^{\vee,Hdg}\iota_S^{-1}(e'_*\mathcal Hom^{\bullet}(\hat R^{CH}(\rho_{\tilde S_I}^*L(i_{I*}j_I^*(F,W))),  
E_{zar}(\Omega^{\bullet,\Gamma,pr}_{/\tilde S_I},F_{DR}))[-d_{\tilde S_I}],u^q_{IJ}(F,W)) \\ 
\xrightarrow{\mathcal Hom^{\bullet}(\hat R_{\tilde S_I}^{CH}(\gamma^{\vee,Z_I}(L(i_{I*}j_I^*(F,W)))),
E_{zar}(\Omega^{\bullet,\Gamma,pr}_{/\tilde S_I},F_{DR}))} \\ 
\Gamma_Z^{\vee,Hdg}\iota_S^{-1}(e'_*\mathcal Hom^{\bullet}(
\hat R^{CH}(\rho_{\tilde S_I}^*\Gamma^{\vee}_{Z_I}L(i_{I*}j_I^*(F,W))),  
E_{zar}(\Omega^{\bullet,\Gamma,pr}_{/\tilde S_I},F_{DR}))[-d_{\tilde S_I}],u^{q,Z}_{IJ}(F,W)) \\
\xrightarrow{=}  
\iota_S^{-1}(e'_*\mathcal Hom^{\bullet}(\hat R^{CH}(\rho_{\tilde S_I}^*\Gamma^{\vee}_{Z_I}L(i_{I*}j_I^*(F,W))),  
E_{zar}(\Omega^{\bullet,\Gamma,pr}_{/\tilde S_I},F_{DR}))[-d_{\tilde S_I}],u^{q,Z}_{IJ}(F,W)).
\end{eqnarray*}
with $u^{q,Z}_{IJ}(F)$ as in \cite{B4}.
\item[(ii)]For $F\in C(\Var(k)^{sm}/S)$, we have also the following canonical map 
in $D(DRM(S))\subset D_{\mathcal D(1,0)fil}(S/(\tilde S_I))$
\begin{eqnarray*}
T(\Gamma_Z^{Hdg},\Omega^{\Gamma,pr}_{/S})(F,W):\\
\iota_S^{-1}(e'_*\mathcal Hom^{\bullet}(
\hat R^{CH}(\rho_{\tilde S_I}^*L\Gamma_{Z_I}E(i_{I*}j_I^*\mathbb D_S(F,W))),  
E_{zar}(\Omega^{\bullet,\Gamma,pr}_{/\tilde S_I},F_{DR}))[-d_{\tilde S_I}],u^{q,Z,d}_{IJ}(F,W)) 
\xrightarrow{=} \\ 
\Gamma_Z^{Hdg}\iota_S^{-1}(e'_*\mathcal Hom^{\bullet}(
\hat R^{CH}(\rho_{\tilde S_I}^*L\Gamma_{Z_I}E(i_{I*}j_I^*\mathbb D_S(F,W))),  
E_{zar}(\Omega^{\bullet,\Gamma,pr}_{/\tilde S_I},F_{DR}))[-d_{\tilde S_I}],u^{q,Z,d}_{IJ}(F,W)) \\
\xrightarrow{\mathcal Hom^{\bullet}(
\hat R_{\tilde S_I}^{CH}(\gamma^{Z_I}(-)),E_{zar}(\Omega^{\bullet,\Gamma,pr}_{/\tilde S_I},F_{DR}))} \\
\Gamma_Z^{Hdg}\iota_S^{-1}(e'_*\mathcal Hom^{\bullet}(
\hat R^{CH}(\rho_{\tilde S_I}^*L(i_{I*}j_I^*\mathbb D_S(F,W))), 
E_{zar}(\Omega^{\bullet,\Gamma,pr}_{/\tilde S_I},F_{DR}))[-d_{\tilde S_I}],u^q_{IJ}(F,W)) 
\end{eqnarray*}
with $u^{q,Z}_{IJ}(F)$ as in \cite{B4}.
\end{itemize}
\end{defi}

\begin{defi}\label{TgDRdefsing}
Let $g:T\to S$ a morphism with $T,S\in\Var(k)$.
Assume we have a factorization $g:T\xrightarrow{l}Y\times S\xrightarrow{p_S}S$
with $Y\in\SmVar(k)$, $l$ a closed embedding and $p_S$ the projection.
Let $S=\cup_{i=1}^lS_i$ be an open cover such that 
there exists closed embeddings $i_i:S_i\hookrightarrow\tilde S_i$ with $\tilde S_i\in\SmVar(k)$
Then, $T=\cup^l_{i=1} T_i$ with $T_i:=g^{-1}(S_i)$
and we have closed embeddings $i'_i:=i_i\circ l:T_i\hookrightarrow Y\times\tilde S_i$,
Moreover $\tilde g_I:=p_{\tilde S_I}:Y\times\tilde S_I\to\tilde S_I$ is a lift of $g_I:=g_{|T_I}:T_I\to S_I$.
Let $M\in\DA_c(S)^-$ and $(F,W)\in C_{fil}(\Var(k)^{sm}/S)$ such that $(M,W)=D(\mathbb A^1_S,et)(F,W)$.
Then, $D(\mathbb A^1_T,et)(g^*F)=g^*M$ and there exist 
$(F',W)\in C_{fil}(\Var(k)^{sm}/S)$ and an equivalence $(\mathbb A^1,et)$ local $e:g^*(F,W)\to(F',W)$ 
such that $D(\mathbb A^1_T,et)(F',W)=(g^*M,W)$.Denote for short $d_{YI}:=-d_Y-d_{\tilde S_I}$.  
We have, using definition \ref{TgDR} and definition \ref{TGammaFDR}(i), 
the canonical map in $D(DRM(T))\subset D_{\mathcal D(1,0)fil}(T/(Y\times\tilde S_I))$
\begin{eqnarray*}
T(g,\mathcal F^{FDR})(M):g^{\hat*mod}_{Hdg}\iota_S^{-1}\mathcal F_S^{FDR}(M):= \\
\Gamma^{\vee,Hdg}_T\iota_T^{-1}(\tilde g_I^{*mod}
(e'_*\mathcal Hom^{\bullet}(\hat R^{CH}(\rho_{\tilde S_I}^*(L(i_{I*}j_I^*(F,W))),   
E_{zar}(\Omega^{\bullet,\Gamma,pr}_{/\tilde S_I},F_{DR})))[d_{YI}],\tilde g_J^{*mod}u^q_{IJ}(F,W))) \\
\xrightarrow{(T(\tilde g_I,\Omega^{\Gamma,pr}_{/\cdot})(-))} \\
\Gamma^{\vee,Hdg}_T\iota_T^{-1}(e'_*\mathcal Hom^{\bullet}(\tilde g_I^*\hat R^{CH}(\rho_{\tilde S_I}^*L(i_{I*}j_I^*(F,W))), 
E_{zar}(\Omega^{\bullet,\Gamma,pr}_{/Y\times\tilde S_I},F_{DR}))[d_{YI}],\tilde g_J^*u^q_{IJ}(F,W)) \\
\xrightarrow{\mathcal Hom(T(\tilde g_I,R^{CH})(-)^{-1},-)} \\ 
\Gamma^{\vee,Hdg}_T\iota_T^{-1}(e'_*\mathcal Hom^{\bullet}(\hat R^{CH}(\rho_{Y\times\tilde S_I}^*\tilde g_I^*L(i_{I*}j_I^*(F,W))), 
E_{zar}(\Omega^{\bullet,\Gamma,pr}_{/Y\times\tilde S_I},F_{DR}))[d_{YI}],\tilde g_J^*u^q_{IJ}(F,W)) \\  
\xrightarrow{T(\Gamma_T^{\vee,Hdg},\Omega^{\Gamma,pr}_{/S})(F,W)}\\
\iota_T^{-1}(e'_*\mathcal Hom^{\bullet}(\hat R^{CH}(\rho_{Y\times\tilde S_I}^*\Gamma^{\vee}_{T_I}\tilde g_I^*L(i_{I*}j_I^*(F,W))), 
E_{zar}(\Omega^{\bullet,\Gamma,pr}_{/Y\times\tilde S_I},F_{DR}))[d_{YI}],\tilde g_J^{*,\gamma}u^q_{IJ}(F,W)) \\ 
\xrightarrow{(\mathcal Hom(\hat R^{CH}_{Y\times\tilde S_I}(T^{q,\gamma}(D_{gI})(j_I^*(F,W))), 
E_{zar}(\Omega^{\bullet,\Gamma,pr}_{/Y\times\tilde S_I},F_{DR}))[d_{YI}])} \\
\iota_T^{-1}(e'_*\mathcal Hom^{\bullet}(\hat R^{CH}(\rho_{Y\times\tilde S_I}^*L(i'_{I*}j^{'*}_Ig^*(F,W))), 
E_{zar}(\Omega^{\bullet,\Gamma,pr}_{/Y\times\tilde S_I},F_{DR}))[d_{YI}],u^q_{IJ}(g^*(F,W))) \\ 
\xrightarrow{\mathcal Hom(\hat R^{CH}_{Y\times\tilde S_I}(Li'_{I*}j^{'*}_I(e)),-)} \\
\iota_T^{-1}(e'_*\mathcal Hom^{\bullet}(\hat R^{CH}(\rho_{Y\times\tilde S_I}^*L(i'_{I*}j^{'*}_I(F',W))), 
E_{zar}(\Omega^{\bullet,\Gamma,pr}_{/Y\times\tilde S_I},F_{DR}))[d_{YI}],u^q_{IJ}(F',W)) \\
\xrightarrow{=:}\mathcal F_T^{FDR}(g^*M)
\end{eqnarray*}
\end{defi}

\begin{defi}\label{SixTalg}
\begin{itemize}
\item Let $f:X\to S$ a morphism with $X,S\in\Var(k)$. Assume there exist a factorization
$f:X\xrightarrow{l}Y\times S\xrightarrow{p_S}S$ with $Y\in\SmVar(k)$, $l$ a closed embedding and $p_S$ the projection.
We have, for $M\in\DA_c(X)$, the following transformation map in $D(DRM(S))$
\begin{eqnarray*}
T_*(f,\mathcal F^{FDR})(M):\mathcal F_S^{FDR}(Rf_*M)
\xrightarrow{\ad(f_{Hdg}^{\hat*mod},Rf^{Hdg}_*)(-)}Rf^{Hdg}_*f^{\hat*mod}_{Hdg}\mathcal F_S^{FDR}(Rf_*M) \\
\xrightarrow{T(f,\mathcal F^{FDR})(Rf_*M)}Rf^{Hdg}_*\mathcal F_X^{FDR}(f^*Rf_*M)
\xrightarrow{\mathcal F_X^{FDR}(\ad(f^*,Rf_*)(M))}Rf^{Hdg}_*\mathcal F_X^{FDR}(M)
\end{eqnarray*}
Clearly, for $p:Y\times S\to S$ a projection with $Y\in\PSmVar(\mathbb C)$, we have, for $M\in\DA_c(Y\times S)$,
$T_*(p,\mathcal F^{FDR})(M)=T_!(p,\mathcal F^{FDR})(M)[d_Y]$

\item Let $S\in\Var(k)$. Let $Y\in\SmVar(k)$ and $p:Y\times S\to S$ the projection.
We have then, for $M\in\DA(Y\times S)$ the following transformation map in $D(DRM(S))$
\begin{eqnarray*}
T_!(p,\mathcal F^{FDR})(M):p^{Hdg}_!\mathcal F_{Y\times S}^{FDR}(M)
\xrightarrow{\mathcal F_{Y\times S}^{FDR}(\ad(Lp_{\sharp},p^*)(M))}
Rp^{Hdg}_!\mathcal F_{Y\times S}^{FDR}(p^*Lp_{\sharp}M) \\
\xrightarrow{T(p,\mathcal F^{FDR})(Lp_{\sharp}(M,W))}Rp^{Hdg}_!p^{\hat*mod[-]}\mathcal F_S^{FDR}(Lp_{\sharp}M)
\xrightarrow{T(p^{*mod},p^{\hat*mod})(-)}p^{Hdg}_!p^{*mod[-]} \\
\mathcal F_S^{FDR}(Lp_{\sharp}M)\xrightarrow{\ad(Rp^{Hdg}_!,p^{*mod[-]})(\mathcal F_S^{FDR}(Lp_{\sharp}M))}
\mathcal F_S^{FDR}(Lp_{\sharp}M)
\end{eqnarray*}

\item Let $f:X\to S$ a morphism with $X,S\in\Var(k)$. Assume there exist a factorization
$f:X\xrightarrow{l}Y\times S\xrightarrow{p_S}S$ with $Y\in\SmVar(k)$, $l$ a closed embedding and $p_S$ the projection.
We have then, using the second point, for $M\in\DA(X)$ the following transformation map in $D(DRM(S))$
\begin{eqnarray*}
T_!(f,\mathcal F^{FDR})(M):
Rp^{Hdg}_!\mathcal F_X^{FDR}(M,W):=Rp^{Hdg}_!\mathcal F_{Y\times S}^{FDR}(l_*M) \\
\xrightarrow{T_!(p,\mathcal F^{FDR})(l_*M)}\mathcal F_S^{FDR}(Lp_{\sharp}l_*M)
\xrightarrow{=}\mathcal F_S^{FDR}(Rf_!M)
\end{eqnarray*}

\item Let $f:X\to S$ a morphism with $X,S\in\Var(k)$. Assume there exist a factorization
$f:X\xrightarrow{l}Y\times S\xrightarrow{p_S}S$ with $Y\in\SmVar(k)$, $l$ a closed embedding and $p_S$ the projection.
We have, using the third point, for $M\in\DA(S)$, the following transformation map in in $D(DRM(X))$
\begin{eqnarray*}
T^!(f,\mathcal F^{FDR})(M):\mathcal F_X^{FDR}(f^!M)
\xrightarrow{\ad(Rf^{Hdg}_!,f^{*mod}_{Hdg})(\mathcal F_X^{FDR}(f^!M))}
f^{*mod}_{Hdg}Rf^{Hdg}_!\mathcal F_X^{FDR}(f^!M) \\
\xrightarrow{T_!(p_S,\mathcal F^{FDR})(\mathcal F^{FDR}(f^!M))}f^{*mod}_{Hdg}\mathcal F_S^{FDR}(Rf_!f^!M)
\xrightarrow{\mathcal F_S^{FDR}(\ad(Rf_!,f^!)(M))}f^{*mod}_{Hdg}\mathcal F_S^{FDR}(M)
\end{eqnarray*}

\item Let $S\in\Var(k)$. Let $S=\cup_{i=1}^l S_i$ an open cover such that there exist closed embeddings 
$i_i:S_i\hookrightarrow\tilde S_i$ with $\tilde S_i\in\SmVar(k)$. 
We have, using definition \ref{TotimesDR} and the preceding point, 
denoting $\Delta_S:S\hookrightarrow S$ the diagonal closed embedding 
and $p_1:S\times S\to S$, $p_2:S\times S\to S$ the projections, 
for $M,N\in\DA(S)$ and $(F,W),(G,W)\in C_{fil}(\Var(k)^{sm}/S))$ 
such that $(M,W)=D(\mathbb A^1,et)(F,W)$ and $(N,W)=D(\mathbb A^1,et)(G,W)$, 
the following transformation map in $D(DRM(S))$
\begin{eqnarray*}
T(\mathcal F_S^{FDR},\otimes)(M,N):\mathcal F_S^{FDR}(M)\otimes^{Hdg}_{O_S}\mathcal F_S^{FDR}(N) 
:=\Delta_S^{!Hdg}(p_1^{*mod}\mathcal F_S^{FDR}(M)\otimes_{O_{S\times S}}p_2^{*mod}\mathcal F_S^{FDR}(N) \\
\xrightarrow{T^!(p_1,\mathcal F_S^{FDR})(M)\otimes T^!(p_1,\mathcal F_S^{FDR})(M)}
\Delta_S^{!Hdg}(\mathcal F_{S\times S}^{FDR}(p_1^!M)\otimes_{O_{S\times S}}\mathcal F_{S\times S}^{FDR}(p_2^!N) \\
\xrightarrow{(T(\otimes,\Omega)(\hat R^{CH}(\rho_{\tilde S_I\times\tilde S_J}^*L(i_I\times i_J)_*(j_I\times j_J)^*p_1^*F[2d_S]),
\hat R^{CH}(\rho_{\tilde S_I\times\tilde S_J}^*L(i_I\times i_J)_*(j_I\times j_J)^*p_2^*F[2d_S])))} \\
\Delta_S^{!Hdg}(\mathcal F_{S\times S}^{FDR}(p_1^!M\otimes p_2^!N)
\xrightarrow{T^!(\Delta_S,\mathcal F^{FDR})(p_1^!M\otimes p_2^!N)}
\mathcal F_S^{FDR}(\Delta_S^!(p_1^!M\otimes p_2^!N))=\mathcal F_S^{FDR}(M\otimes N)
\end{eqnarray*}
where the last equality follows from the equality in $\DA(S)$
\begin{equation*}
\Delta_S^!(p_1^!M\otimes p_2^!N)=\Delta_S^!p_1^!M\otimes\Delta_S^!p_2^!N=M\otimes N
\end{equation*}
\end{itemize}
\end{defi}

\begin{prop}\label{Tgprop}
Let $g:T\to S$ a morphism with $T,S\in\Var(k)$.
Assume we have a factorization $g:T\xrightarrow{l}Y_2\times S\xrightarrow{p_S}S$
with $Y_2\in\SmVar(k)$, $l$ a closed embedding and $p_S$ the projection.
Let $S=\cup_{i=1}^lS_i$ be an open cover such that 
there exists closed embeddings $i_i:S_i\hookrightarrow\tilde S_i$ with $\tilde S_i\in\SmVar(k)$
Then, $T=\cup^l_{i=1} T_i$ with $T_i:=g^{-1}(S_i)$
and we have closed embeddings $i'_i:=i_i\circ l:T_i\hookrightarrow Y_2\times\tilde S_i$,
Moreover $\tilde g_I:=p_{\tilde S_I}:Y\times\tilde S_I\to\tilde S_I$ is a lift of $g_I:=g_{|T_I}:T_I\to S_I$.
Let $f:X\to S$ a  morphism with $X\in\Var(k)$ such that there exists a factorization 
$f:X\xrightarrow{l}Y_1\times S\xrightarrow{p_S} S$, with $Y_1\in\SmVar(k)$, 
$l$ a closed embedding and $p_S$ the projection. We have then the following commutative diagram
whose squares are cartesians
\begin{equation*}
\xymatrix{f':X_T\ar[r]\ar[rd]\ar[dd]_{g'} & Y_1\times T\ar[rd]\ar[r] & T\ar[rd]\ar[dd]^{g} & \, \\
\, & Y_1\times X\ar[r]\ar[ld] & Y_1\times Y_2\times S\ar[r]\ar[ld] & Y_2\times S\ar[ld] \\
f:X\ar[r] & Y_1\times S\ar[r] & S & \,}
\end{equation*} 
Take a smooth compactification $\bar Y_1\in\PSmVar(\mathbb C)$ of $Y_1$, denote
$\bar X_I\subset\bar Y_1\times\tilde S_I$ the closure of $X_I$, and $Z_I:=\bar X_I\backslash X_I$.
Consider $F(X/S):=p_{S,\sharp}\Gamma_X^{\vee}\mathbb Z(Y_1\times S/Y_1\times S)\in C(\Var(k)^{sm}/S)$ and
the isomorphism in $C(\Var(k)^{sm}/T)$
\begin{eqnarray*}
T(f,g,F(X/S)):g^*F(X/S):=g^*p_{S,\sharp}\Gamma_X^{\vee}\mathbb Z(Y_1\times S/Y_1\times S)\xrightarrow{\sim} \\
p_{T,\sharp}\Gamma_{X_T}^{\vee}\mathbb Z(Y_1\times T/Y_1\times T)=:F(X_T/T).
\end{eqnarray*}
which gives in $\DA(T)$ the isomorphism $T(f,g,F(X/S)):g^*M(X/S)\xrightarrow{\sim}(X_T/T)$.
Then the following diagram in $D(DRM(T))\subset D_{\mathcal D(1,0)fil}(T/(Y_2\times\tilde S_I))$, 
where the horizontal maps are given by proposition \ref{keyalgsing1}, commutes
\begin{equation*}
\begin{tikzcd}
g^{\hat*mod}_{Hdg}\iota_S^{-1}\mathcal F_S^{FDR}(M(X/S))
\ar[dd,"'T(g{,}\mathcal F^{FDR})(M(X/S))"']\ar[rr,"g^{\hat*mod}_{Hdg}I(X/S)"] & \, &
g^{\hat*mod}_{Hdg}Rf^{Hdg}_!(\Gamma^{\vee,Hdg}_{X_I}(O_{Y_1\times\tilde S_I},F_b)(d_{Y_1})[2d_{Y_1}],x_{IJ}(X/S))
\ar[d,"T(p_{\tilde S_I}{,}\gamma^{\vee{,}Hdg})(-)"] \\
\, & \, & Rf^{'Hdg}_!g^{'\hat*mod}_{Hdg}(\Gamma^{\vee,Hdg}_{X_I}(O_{Y_1\times\tilde S_I},F_b)(d_{Y_1})[2d_{Y_1}],x_{IJ}(X/S))
\ar[d,"T(p_{Y_1\times Y_2\times\tilde S_I{,}Hdg}^{\hat*mod}{,}p_{Y_1\times Y_2\times\tilde S_I{,}Hdg}^{*mod})(-)"] \\
\iota_T^{-1}\mathcal F_T^{FDR}(M(X_T/T))\ar[rr,"I(X_T/T)"] & \, & 
Rf^{'Hdg}_!(\Gamma^{\vee,Hdg}_{X_{T_I}}(O_{Y_2\times Y_1\times\tilde S_I},F_b)(d_{Y_{12}})[2d_{Y_{12}}],x_{IJ}(X_T/T)).
\end{tikzcd}
\end{equation*} 
with $d_{Y_{12}}=d_{Y_1}+d_{Y_2}$.
\end{prop}

\begin{proof}
Follows immediately from definition.
\end{proof}

\begin{prop}\label{mainthmprop2}
Let $S\in\Var(k)$. Let $Y\in\SmVar(k)$ and $p:Y\times S\to S$ the projection.
Let $S=\cup_{i=1}^l S_i$ an open cover such that there exist closed embeddings 
$i_i:S_i\hookrightarrow\tilde S_i$ with $\tilde S_i\in\SmVar(k)$.
For $I\subset\left[1,\cdots l\right]$, we denote by $S_I=\cap_{i\in I} S_i$, $j^o_I:S_I\hookrightarrow S$ and
$j_I:Y\times S_I\hookrightarrow Y\times S$ the open embeddings. 
We then have closed embeddings $i_I:Y\times S_I\hookrightarrow Y\times\tilde S_I$.
and we denote by $p_{\tilde S_I}:Y\times\tilde S_I\to\tilde S_I$ the projections.
Let $f':X'\to Y\times S$ a morphism, with $X'\in\Var(k)$ such that there exists a factorization
$f':X'\xrightarrow{l'}Y'\times Y\times S\xrightarrow{p'} Y\times S$ 
with $Y'\in\SmVar(k)$, $l'$ a closed embedding and $p'$ the projection. 
Denoting $X'_I:=f^{'-1}(Y\times S_I)$, we have closed embeddings $i'_I:X'_I\hookrightarrow Y'\times Y\times\tilde S_I$
Consider 
\begin{eqnarray*}
F(X'/Y\times S):=p_{Y\times S,\sharp}\Gamma_{X'}^{\vee}\mathbb Z(Y'\times Y\times S/Y'\times Y\times S)
\in C(\Var(k)^{sm}/Y\times S)
\end{eqnarray*}
and $F(X'/S):=p_{\sharp}F(X'/Y\times S)\in C(\Var(k)^{sm}/S)$, 
so that $Lp_{\sharp}M(X'/Y\times S)[-2d_Y]=:M(X'/S)$. 
Then, the following diagram in $D(DRM(S))\subset D_{\mathcal D(1,0)fil}(S/(Y\times\tilde S_I))$,
where the vertical maps are given by proposition \ref{keyalgsing1}, commutes
\begin{equation*}
\xymatrix{
Rp^{Hdg}!\mathcal F_{Y\times S}^{FDR}(M(X'/Y\times S))
\ar[rr]^{T_!(p,\mathcal F^{FDR})(M(X'/Y\times S))} & \, & 
\mathcal F_S^{FDR}(M(X'/S)) \\
Rp^{Hdg}!Rf^{'Hdg}_!f^{'*mod}_{Hdg}\mathbb Z_{Y\times S}^{Hdg} 
\ar[rr]^{=}\ar[u]^{T(p_{Hdg}^{\hat*mod}{,}p_{Hdg}^{*mod})(-)\circ Rp^{Hdg}!(I(X'/Y\times S))} & \, &
Rf^{Hdg}_!f^{*mod}_{Hdg}\mathbb Z_S^{Hdg}\ar[u]_{I(X'/S)}}. 
\end{equation*}
\end{prop}

\begin{proof}
Immediate from definition.
\end{proof}

\begin{prop}\label{TotimesDRprop}
Let $f_1:X_1\to S$, $f_2:X_2\to S$ two morphism with $X_1,X_2,S\in\Var(k)$. 
Assume that there exist factorizations 
$f_1:X_1\xrightarrow{l_1}Y_1\times S\xrightarrow{p_S} S$, $f_2:X_2\xrightarrow{l_2}Y_2\times S\xrightarrow{p_S} S$
with $Y_1,Y_2\in\SmVar(k)$, $l_1,l_2$ closed embeddings and $p_S$ the projections.
We have then the factorization
\begin{equation*}
f_{12}:=f_1\times f_2:X_{12}:=X_1\times_S X_2\xrightarrow{l_1\times l_2}Y_1\times Y_2\times S\xrightarrow{p_S} S
\end{equation*}
Let $S=\cup_{i=1}^l S_i$ an open affine covering and denote, 
for $I\subset\left[1,\cdots l\right]$, $S_I=\cap_{i\in I} S_i$ and $j_I:S_I\hookrightarrow S$ the open embedding.
Let $i_i:S_i\hookrightarrow\tilde S_i$ closed embeddings, with $\tilde S_i\in\SmVar(k)$. 
We have then the following commutative diagram in $D(DRM(S))\subset D_{\mathcal D(1,0)fil}(S/(\tilde S_I))$
where the vertical maps are given by proposition \ref{keyalgsing1}
\begin{equation*}
\begin{tikzcd}
\mathcal F_S^{FDR}(M(X_1/S))\otimes^{Hdg}_{O_S}\mathcal F_S^{FDR}(M(X_2/S))
\ar[rr,"I(X_1/S)\otimes I(X_2/S)"] \ar[d,"T(\mathcal F_S^{FDR}{,}\otimes)(M(X_1/S){,}M(X_2/S))"]  & \, &
\shortstack{$Rf_{1!}^{Hdg}(\Gamma^{\vee,Hdg}_{X_{1I}}(O_{Y_1\times\tilde S_I},F_b)(d_2)[2d_1],x_{IJ}(X_1/S))\otimes_{O_S}$ \\ 
$Rf_{2!}^{Hdg}(\Gamma^{\vee,Hdg}_{X_{2I}}(O_{Y_2\times\tilde S_I},F_b)(d_1)[2d_2],x_{IJ}(X_2/S))$}
\ar[d,"(Ew_{(Y_1\times\tilde S_I{,}Y_2\times\tilde S_I)/\tilde S_I})"] \\
\mathcal F_S^{FDR}(M(X_1/S)\otimes M(X_2/S)=M(X_1\times_S X_2/S))\ar[rr,"I(X_{12}/S)"] & \, &
Rf_{12!}^{Hdg}(\Gamma^{\vee,Hdg}_{X_{1I}\times_S X_{2I}}(O_{Y_1\times Y_2\times\tilde S_I},F_b)(d_{12})[2d_{12}],x_{IJ}(X_1/S)).
\end{tikzcd}
\end{equation*}
with $d_1=d_{Y_1}$, $d_2=d_{Y_2}$ and $d_{12}=d_{Y_1}+d_{Y_2}$.
\end{prop}

\begin{proof}
Immediate from definition.
\end{proof}

\begin{thm}\label{mainthm}
\begin{itemize}
\item[(i)]Let $g:T\to S$ a morphism, with $S,T\in\Var(k)$. 
Assume we have a factorization $g:T\xrightarrow{l}Y\times S\xrightarrow{p_S}S$
with $Y\in\SmVar(k)$, $l$ a closed embedding and $p_S$ the projection. Let $M\in\DA_c(S)$. 
Then map in $D(DRM(T))$ 
\begin{eqnarray*}
T(g,\mathcal F^{FDR})(M):g^{\hat*mod}_{Hdg}\mathcal F_S^{FDR}(M)\xrightarrow{\sim}\mathcal F_T^{FDR}(g^*M)
\end{eqnarray*} 
given in definition \ref{TgDRdefsing} is an isomorphism.
\item[(ii)] Let $f:X\to S$ a morphism with $X,S\in\Var(k)$. Assume there exist a factorization
$f:X\xrightarrow{l}Y\times S\xrightarrow{p_S}S$ with $Y\in\SmVar(k)$, $l$ a closed embedding and $p_S$ the projection.
Then, for $M\in\DA_c(X)$, the map given in definition \ref{SixTalg}
\begin{equation*}
T_!(f,\mathcal F^{FDR})(M):Rf^{Hdg}_!\mathcal F_X^{FDR}(M)\xrightarrow{\sim}\mathcal F_S^{FDR}(Rf_!M)
\end{equation*}
is an isomorphism in $D(DRM(S))$.
\item[(iii)] Let $f:X\to S$ a morphism with $X,S\in\Var(k)$, $S$ quasi-projective. Assume there exist a factorization
$f:X\xrightarrow{l}Y\times S\xrightarrow{p_S}S$ with $Y\in\SmVar(k)$, $l$ a closed embedding and $p_S$ the projection.
We have, for $M\in\DA_c(X)$, the map given in definition \ref{SixTalg}
\begin{equation*}
T_*(f,\mathcal F^{FDR})(M):\mathcal F_S^{FDR}(Rf_*M)\xrightarrow{\sim}Rf^{Hdg}_*\mathcal F_X^{FDR}(M)
\end{equation*}
is an isomorphism in $D(DRM(S))$.
\item[(iv)] Let $f:X\to S$ a morphism with $X,S\in\Var(k)$, $S$ quasi-projective. Assume there exist a factorization
$f:X\xrightarrow{l}Y\times S\xrightarrow{p_S}S$ with $Y\in\SmVar(k)$, $l$ a closed embedding and $p_S$ the projection.
Then, for $M\in\DA_c(S)$, the map given in definition \ref{SixTalg}
\begin{equation*}
T^!(f,\mathcal F^{FDR})(M):\mathcal F_X^{FDR}(f^!M)\xrightarrow{\sim}f^{*mod}_{Hdg}\mathcal F_S^{FDR}(M)
\end{equation*}
is an isomorphism in $D(DRM(X))$. 
\item[(v)]Let $S\in\Var(k)$. Then, for $M,N\in\DA_c(S)$, the map in $D(DRM(S))$ 
\begin{eqnarray*}
T(\mathcal F_S^{FDR},\otimes)(M,N): 
\mathcal F_S^{FDR}(M)\otimes^{Hdg}_{O_S}\mathcal F_S^{FDR}(N)\xrightarrow{\sim}\mathcal F_S^{FDR}(M\otimes N)
\end{eqnarray*}
given in definition \ref{SixTalg} is an isomorphism.
\end{itemize}
\end{thm}

\begin{proof}
The proof is similar to the complex case : 
follows from \cite{AyoubT} by proposition \ref{Tgprop} and proposition \ref{mainthmprop2}, more precisely :

\noindent(i):follows from proposition \ref{Tgprop} and proposition \ref{keyalgsing1}.

\noindent(ii):follows from proposition \ref{mainthmprop2}

\noindent(iii),(iv): see \cite{B4}.

\noindent(v):follows from proposition \ref{TotimesDRprop}.
\end{proof}

We have the following easy proposition

\begin{prop}
Let $S\in\Var(k)$ and $S=\cup_{i=1}^l S_i$ an open affine covering and denote, 
for $I\subset\left[1,\cdots l\right]$, $S_I=\cap_{i\in I} S_i$ and $j_I:S_I\hookrightarrow S$ the open embedding.
Let $i_i:S_i\hookrightarrow\tilde S_i$ closed embeddings, with $\tilde S_i\in\SmVar(k)$. 
We have, for $M,N\in\DA(S)$ and $F,G\in C(\Var(k)^{sm}/S)$ such that 
$M=D(\mathbb A^1,et)(F)$ and $N=D(\mathbb A^1,et)(G)$, 
the following commutative diagram in $D_{O_Sfil,\mathcal D,\infty}(S/(\tilde S_I))$
\begin{equation*}
\xymatrix{\mathcal F_S^{GM}(L\mathbb D_SM)\otimes^L_{O_S}\mathcal F_S^{GM}(L\mathbb D_SN)
\ar[d]^{T(\mathcal F_S^{GM},\otimes)(L\mathbb D_SM,L\mathbb D_SN)}
\ar[rrrr]^{T(\mathcal F_S^{GM},\mathcal F_S^{FDR})(M)\otimes T(\mathcal F_S^{GM},\mathcal F_S^{FDR})(N)} & \, & \, & \, &
\mathcal F_S^{FDR}(M)\otimes^{Hdg}_{O_S}\mathcal F_S^{FDR}(N)\ar[d]^{T(\mathcal F_S^{FDR},\otimes)(M,N)} \\
\mathcal F_S^{GM}(L\mathbb D_S(M\otimes N))\ar[rrrr]^{T(\mathcal F_S^{GM},\mathcal F_S^{FDR})(M\otimes N)} & \, & \, & \, & 
\mathcal F_S^{FDR}(M\otimes N)}
\end{equation*}
\end{prop}

\begin{proof}
Immediate from definition.
\end{proof}

\section{The Hodge realization functors for relative motives over a field $k$ of characteristic $0$}

\subsection{The Hodge realization functor for relative motives over a subfield $k\subset\mathbb C$}

Let $k\subset\mathbb C$ a subfield.
We have for $f:T\to S$ a morphism with $T,S\in\Var(k)$ we have the commutative diagram of site 
\begin{equation*}
\xymatrix{
\AnSp(\mathbb C)/T_{\mathbb C}^{an}\ar[rr]^{\An_T:=\An_T\circ(\pi_{k/\mathbb C}(-))}\ar[dd]^{P(f)}\ar[rd]^{\rho_T} & \, & 
\Var(k)/T\ar[dd]^{P(f)}\ar[rd]^{\rho_T} & \, \\  
 \, & \AnSp(\mathbb C)^{sm}/T_{\mathbb C}^{an}\ar[rr]^{\An_T}\ar[dd]^{P(f)} & \, & \Var(k)^{sm}/T\ar[dd]^{P(f)} \\  
\AnSp(\mathbb C)/S_{\mathbb C}^{an}\ar[rr]^{\An_S:=\An_S\circ(\pi_{k/\mathbb C}(-))}\ar[rd]^{\rho_S} & \, & 
\Var(k)/S\ar[rd]^{\rho_S} & \, \\  
\, &  \AnSp(\mathbb C)^{sm}/S_{\mathbb C}^{an}\ar[rr]^{\An_S} & \, & \Var(k)^{sm}/S}.  
\end{equation*}
This gives for $s:\mathcal I\to\mathcal J$ a functor with $\mathcal I,\mathcal J\in\Cat$ and
$f:T_{\bullet}\to S_{s(\bullet)}$ a morphism of diagram of algebraic varieties with 
$T_{\bullet}\in\Fun(\mathcal I,\Var(k))$, $S_{\bullet}\in\Fun(\mathcal J,\Var(k))$ 
the commutative diagram of sites 
\begin{equation*}
Dia^{12}(S):=\xymatrix{
\AnSp(\mathbb C)/T_{\bullet,\mathbb C}^{an}
\ar[rr]^{\An_{T_{\bullet}}:=\An_{T_{\bullet}}\circ(\pi_{k/\mathbb C}(-))}
\ar[dd]^{P(f_{\bullet})}\ar[rd]^{\rho_{T_{\bullet}}} & \, & 
\Var(k)/T_{\bullet}\ar[dd]^{P(f_{\bullet})}\ar[rd]^{\rho_{T_{\bullet}}} & \, \\  
 \, & \AnSp(\mathbb C)^{sm}/T_{\bullet,\mathbb C}^{an}\ar[rr]^{\An_{T_{\bullet}}}\ar[dd]^{P(f_{\bullet})} & \, & 
\Var(k)^{sm}/T_{\bullet}\ar[dd]^{P(f_{\bullet})} \\  
\AnSp(\mathbb C)/S_{\bullet,\mathbb C}^{an}
\ar[rr]^{\An_{S_{\bullet}}:=\An_{S_{\bullet}}\circ(\pi_{k/\mathbb C}(-))}\ar[rd]^{\rho_{S_{\bullet}}} & \, & 
\Var(k)/S_{\bullet}\ar[rd]^{\rho_{S_{\bullet}}} & \, \\  
\, &  \AnSp(\mathbb C)^{sm}/S_{\bullet,\mathbb C}^{an}\ar[rr]^{\An_{S_{\bullet}}} & \, & 
\Var(k)^{sm}/S_{\bullet}}.  
\end{equation*}

\subsubsection{The Betti realization functor}

Let $k\subset\mathbb C$ a subfield.

\begin{defi}\label{bettidef}
Let $S\in\Var(k)$.
\begin{itemize}
\item[(i)] The Ayoub's Betti realization functor is 
\begin{equation*}
\Bti_S^*:\DA(S)\to D(S_{\mathbb C}^{an}) \; , \; 
M\in\DA(S)\mapsto\Bti_S^*M=Re(S_{\mathbb C}^{an})_*\An_S^*M=e(S_{\mathbb C}^{an})_*\underline{\sing}_{\mathbb D^*}\An_S^*F
\end{equation*}
where $F\in C(\Var(k)^{sm}/S)$ is such that $M=D(\mathbb A^1,et)(F)$.
\item[(ii)] In \cite{B3}, we define the Betti realization functor as 
\begin{equation*}
\widetilde\Bti_S^*:\DA(S)\to D(S^{an})=D(S^{cw}) \; , \; 
M\mapsto\widetilde\Bti_S^*M=Re(S^{cw})_*\widetilde\Cw_S^*M=e(S^{cw})_*\underline{\sing}_{\mathbb I^*}\widetilde\Cw_S^*F
\end{equation*}
where $F\in C(\Var(k)^{sm}/S)$ is such that $M=D(\mathbb A^1,et)(F)$.
\item[(iii)] For the Corti-Hanamura weight structure on $\DA^-(S)$, we have by functoriality of (i) the functor
\begin{equation*}
\Bti_S^*:\DA^-(S)\to D_{fil}(S_{\mathbb C}^{an}) \; , \; 
M\mapsto(\Bti_S^*M,W):=\Bti_S^*(M,W):=e(S^{an})_*\underline{\sing}_{\mathbb D^*}\An_S^*(F,W)
\end{equation*}
where $(F,W)\in C_{fil}(\Var(k)^{sm}/S)$ is such that $(M,W)=D(\mathbb A^1,et)(F,W)$.
\end{itemize}
Note that by \cite{B3}, $\An_S^*$ and $\widetilde\Cw_S^*$ derive trivially.
\end{defi}

Note that, by considering the explicit $\mathbb D_S^1$ local model for presheaves on $\AnSp(\mathbb C)^{sm}/S_{\mathbb C}^{an}$,
$\Bti_S^*(\DA^-(S))\subset D^-(S^{an})$ ; 
by considering the explicit $\mathbb I_S^1$ local model for presheaves on $\CW^{sm}/S^{cw}$,
$\widetilde\Bti_S^*(\DA^-(S))\subset D^-(S_{\mathbb C}^{an})$.

Let $f:T\to S$ a morphism, with $T,S\in\Var(k)$. We have, 
for $M\in\DA(S)$, $(F,W)\in C_{fil}(\Var(k)^{sm}/S)$ such that $(M,W)=D(\mathbb A^1,et)(F,W)$,
and an equivalence $(\mathbb A^1,et)$ local $e:f^*(F,W)\to (F',W)$ 
with $(F',W)\in C_{fil}(\Var(k)^{sm}/S)$ such that $(f^*M,W)=D(\mathbb A^1,et)(F',W)$ 
the following canonical transformation map in $D_{fil}(T^{an}_{\mathbb C})$:
\begin{eqnarray*}
T^0(f,\Bti)(M,W):f^*\Bti_S^*(M,W):=f^*e(S_{\mathbb C}^{an})_*\underline{\sing}_{\mathbb D^*}\An_S^*(F,W) \\
\xrightarrow{T(f,e)(-)} e(T_{\mathbb C}^{an})_*f^*\underline{\sing}_{\mathbb D^*}\An_S^*(F,W) \\
\xrightarrow{e(T_{\mathbb C}^{an})_*T(f,c)(F,W)}e(T^{an})_*\underline{\sing}_{\mathbb D^*}f^*\An_T^*(F,W)
\xrightarrow{=}e(T_{\mathbb C}^{an})_*\underline{\sing}_{\mathbb D^*}\An_T^*f^*(F,W) \\
\xrightarrow{e(T_{\mathbb C}^{an})_*\underline{\sing}_{\mathbb D^*}\An_T^*e}
e(T_{\mathbb C}^{an})_*\underline{\sing}_{\mathbb D^*}\An_T^*(F',W)=:\Bti_T^*f^*(M,W).
\end{eqnarray*}

\begin{defi}\label{TgBti}
Let $f:T\to S$ a morphism, with $T,S\in\Var(k)$. 
Consider the graph factorization $f:T\xrightarrow{l}T\times S\xrightarrow{p}S$ of $f$
with $l$ the graph closed embedding and $p$ the projection.
We have, for $M\in\DA_c(S)$,  the following canonical transformation map in $D_{fil,c}(T^{an}_{\mathbb C})$:
\begin{eqnarray*}
T(f,\Bti)(M,W):f^{*w}\Bti_S^*(M,W):=l^*\Gamma_T^{\vee,w}p^*\Bti_S^*(F,W) \\
\xrightarrow{T^0(p,\Bti)(-)}l^*\Gamma_T^{\vee,w}\Bti_{T\times S}^*p^*(F,W)
\xrightarrow{\gamma_T^{\vee}(p^*(F,W))}l^*\Gamma_T^{\vee,w}\Bti_{T\times S}^*\Gamma_T^{\vee}p^*(F,W) \\
\xrightarrow{=}l^*\Bti_{T\times S}^*\Gamma_T^{\vee}p^*(F,W)
\xrightarrow{T^0(l,\Bti)(-)}\Bti_T^*l^*\Gamma_T^{\vee}p^*(F,W)=\Bti_T^*f^*(M,W).
\end{eqnarray*}
where we use definition \ref{fw}.
\end{defi}

\begin{defi}\label{TBtiSix}
\begin{itemize}
\item Let $f:X\to S$ a morphism, with $X,S\in\Var(k)$. 
We have, for $M\in\DA_c(X)$, the following transformation map in $D_{fil,c}(S^{an}_{\mathbb C})$
\begin{eqnarray*}
T_*(f,\Bti)(M,W):\Bti_S^*(Rf_*(M,W))\xrightarrow{\ad(f^*,Rf_{*w})(\Bti_S^*(Rf_*(M,W)))}Rf_{*w}f^{*w}\Bti_S^*(Rf_*(M,W)) \\
\xrightarrow{T(f,\Bti)(Rf_*(M,W))}
Rf_{*w}\Bti_X^*(f^*Rf_*(M,W))\xrightarrow{\Bti_X^*(\ad(f^*,Rf_*)(M,W))}Rf_{*w}\Bti_X^*(M,W)
\end{eqnarray*}
Clearly if $l:Z\hookrightarrow S$ is a closed embedding, then $T_*(l,\Bti)(M,W)$ is an isomorphism
since $\ad(l^*,l_*)(-):l^*l_*(M,W)\to (M,W)$ is an isomorphism (see section 3).

\item Let $f:X\to S$ a morphism with $X,S\in\Var(k)$. Assume there exist a factorization
$f:X\xrightarrow{l}Y\times S\xrightarrow{p_S}S$ with $Y\in\SmVar(k)$, $l$ a closed embedding and $p_S$ the projection. 
We have then, for $M\in\DA_c(X)$, using theorem \ref{mainBti} for closed embeddings, 
the following transformation map in $D_{fil}((Y\times S)^{an}_{\mathbb C})$
\begin{eqnarray*}
T_!(f,\Bti)(M):Rf_{!w}\Bti_X^*(M,W)=Rp_{S!w}l_*\Bti_X^*(M,W) \\
\xrightarrow{T_*(l,\Bti)(M,W)}^{-1} Rp_{S!w}\Bti(Y\times S)^*(l_*(M,W)) \\
\xrightarrow{\Bti(Y\times S)^*\ad(Lp_{S\sharp},p_S^*)(l_*(M,W))}
Rp_{S!w}\Bti(Y\times S)^*(p_S^*Lp_{S\sharp}l_*(M,W))\xrightarrow{T(p_S,\Bti)(p_{S\sharp}l_*(M,W))} \\
Rp_{S!w}p_S^*\Bti(Y\times S)^*(Lp_{S\sharp}l_*(M,W))=Rp_{S!w}p_S^{!w}\Bti(Y\times S)^*(Rf_!(M,W)) \\
\xrightarrow{\ad(Rp_{S!w},p_S^{!w})(-)}\Bti(Y\times S)^*(Rf_!(M,W))
\end{eqnarray*}
Clearly, for $f:X\to S$ a proper morphism, with $X,S\in\Var(k)$ we have, for $M\in\DA_c(Y\times S)$,
$T_!(f,\Bti)(M,W)=T_*(f,\Bti)(M,W)$.

\item Let $f:X\to S$ a morphism with $X,S\in\Var(k)$. 
We have, using the second point, for $M\in\DA(S)$, the following transformation map in $D_{fil}(X^{an}_{\mathbb C})$
\begin{eqnarray*}
T^!(f,\Bti)(M,W):\Bti_X^*(f^!(M,W))\xrightarrow{\ad(f_!,Rf^!)(\Bti_X^*(f^!(M,W)))}f^{!w}Rf_{!w}\Bti_X^*(f^!(M,W)) \\
\xrightarrow{T_!(f,\Bti)((f^!(M,W)))}
f^{!w}\Bti_S^*(f_!f^!(M,W))\xrightarrow{\Bti_S^*(\ad(f_!,f^!)(M,W))}f^{!w}\Bti_S^*(M,W)
\end{eqnarray*}

\item Let $S\in\Var(k)$. We have, for $M,N\in\DA(S)$ and $F,G\in C_(\Var(k)^{sm}/S))$  
such that $M=D(\mathbb A^1,et)(F)$ and $N=D(\mathbb A^1,et)(G)$, 
the following transformation map in $D_{fil}(S^{an}_{\mathbb C})$
\begin{eqnarray*}
\Bti_S^*(M,W)\otimes\Bti_S^*(N,W):=
(e(S)_*\underline{\sing}_{\mathbb D^*}\An_S^*(F,W))\otimes(e(S)_*\underline{\sing}_{\mathbb D^*}\An_S^*(G,F)) \\
\xrightarrow{T(\sing_{D^*},\otimes)(\An_S^*(F,W),\An_S^*(G,F))}
e(S)_*\underline{\sing}_{\mathbb D^*}\An_S^*((F,W)\otimes(G,W))=:\Bti_S^*((M,W)\otimes(N,W))
\end{eqnarray*}
\end{itemize}
\end{defi}

\begin{thm}\label{mainBti}
\begin{itemize}
\item[(i)]Let $f:X\to S$ a morphism, with $X,S\in\Var(k)$. For $M\in\DA_c(S)$,
\begin{equation*}
T(f,\Bti)(M,W):f^{*w}\Bti_S^*(M,W)\xrightarrow{\sim}\Bti_X^*f^*(M,W)
\end{equation*}
is an isomorphism in $D_{fil}(X^{an}_{\mathbb C})$.
\item[(ii)] Let $f:X\to S$ a morphism, with $X,S\in\Var(k)$. For $M\in\DA_c(X)$,
\begin{equation*}
T_!(f,\Bti)(M,W):Rf_{!w}\Bti_X^*(M,W)\xrightarrow{\sim}\Bti_S^*Rf_!(M,W)
\end{equation*}
is an isomorphism.
\item[(iii)] Let $f:X\to S$ a morphism, with $X,S\in\Var(k)$. For $M\in\DA_c(X)$,
\begin{equation*}
T_*(f,\Bti)(M,W):Rf_{*w}\Bti_X^*(M,W)\xrightarrow{\sim}\Bti_S^*Rf_*(M,W)
\end{equation*}
is an isomorphism.
\item[(iv)] Let $f:X\to S$ a morphism, with $X,S\in\Var(k)$. For $M\in\DA_c(S)$,
\begin{equation*}
T^!(f,\Bti)(M,W):f^{!w}\Bti_S^*(M,W)\xrightarrow{\sim}\Bti_X^*f^!(M,W)
\end{equation*}
is an isomorphism.
\item[(v)] Let $S\in\Var(k)$. For $M,N\in\DA_c(S)$,
\begin{equation*}
T(\otimes,\Bti)(M,W):\Bti_S^*(M,W)\otimes\Bti_S^*N\xrightarrow{\sim}\Bti_X^*((M,W)\otimes(N,W))
\end{equation*}
is an isomorphism.
\end{itemize}
\end{thm}

\begin{proof}
By functoriality it reduced to the case of Corti-Hanamura motives which is then obvious.
\end{proof}

The main result on the Betti realization functor is the following
\begin{thm}\label{mainBtiCorCor}
\begin{itemize}
\item[(i)] We have $\Bti_S^*=\widetilde\Bti_S^*$ on $\DA^-(S)$
\item[(ii)] The canonical transformations $T(f,\Bti)$, for $f:T\to S$ a morphism in $\Var(k)$,
define a morphism of 2 functor
\begin{equation*}
\Bti_{\cdot}^*:\DA(\cdot)\to D((\cdot)_{\mathbb C}^{an}), \; S\in\Var(k) \mapsto\Bti_S^*:\DA(S)\to D(S_{\mathbb C}^{an})
\end{equation*}
which is a morphism of homotopic 2 functor.
\item[(ii)'] The canonical transformations $T(f,\Bti)$, for $f:T\to S$ a morphism in $\Var(k)$,
define a morphism of 2 functor
\begin{equation*}
\Bti_{\cdot}^*:\DA_c(\cdot)\to D_{fil}((\cdot)_{\mathbb C}^{an}), \; S\in\Var(k) \mapsto\Bti_S^*:\DA(S)\to D_{fil}(S_{\mathbb C}^{an})
\end{equation*}
which is a morphism of homotopic 2 functor.
\end{itemize}
\end{thm}

\begin{proof}
\noindent(i): See \cite{B3}

\noindent(ii) and (ii)':Follows from theorem \ref{mainBti}.
\end{proof}

\begin{rem}
For $X\in\Var(k)$, the quasi-isomorphisms
\begin{equation*}
\mathbb Z\Hom(\bar{\mathbb D_k}_{et}^{\bullet},X)\xrightarrow{\An^*}
\mathbb Z\Hom(\bar{\mathbb D_{\mathbb C}}^n(0,1),X_{\mathbb C}^{an})
\xrightarrow{\Hom(i,X_{\mathbb C}^{cw})}\mathbb Z\Hom([0,1]^n,X_{\mathbb C}^{cw}),
\end{equation*}
where, 
\begin{equation*}
\bar{\mathbb D}^n_{et}:=(e:U\to\mathbb A_k^n,\bar{\mathbb D}^n(0,1)\subset e(U))
\in\Fun(\mathcal V^{et}_{\mathbb A_k^n}(\bar{\mathbb D}^n(0,1)),\Var(k)) 
\end{equation*}
is the system of etale neighborhood of the closed ball $\bar{\mathbb D}_k^n(0,1)\subset\mathbb A_k^n$,
and $i:[0,1]^n\hookrightarrow\bar{\mathbb D}_{\mathbb C}^n(0,1)$ is the closed embedding,
shows that a closed singular chain $\alpha\in\mathbb Z\Hom^n([0,1]^n,X_{\mathbb C}^{cw})$, is homologue to 
a closed singular chain 
\begin{equation*}
\beta=\alpha+\partial\gamma=\tilde{\beta}_{|[0,1]^n}\in\mathbb Z\Hom^n(\Delta^n,X_{\mathbb C}^{cw}) 
\end{equation*}
which is the restriction by the closed embedding 
$[0,1]^n\hookrightarrow U_{\mathbb C}^{cw}\xrightarrow{e}\mathbb A_{\mathbb C}^n$,
where $e:U\to\mathbb A_k^n$ an etale morphism with $U\in\Var(k)$,
of a a complex algebraic morphism $\tilde\beta:U_{\mathbb C}\to X_{\mathbb C}$ defined over $k$. 
Hence $\beta([0,1]^n)=\tilde{\beta}([0,1]^n)\subset X$ is the restriction
of a real algebraic subset of dimension $n$ in $Res_{\mathbb R}(X)$ 
(after restriction a scalar that is under the identification $\mathbb C\simeq\mathbb R^2$)
defined over $k$. 
\end{rem}

\subsubsection{The complex Hodge realization functor for relative motives over a subfield $k\subset\mathbb C$}

Let $k\subset\mathbb C$ a subfield.

Let $S\in\Var(k)$. Let $S=\cup_{i=1}^sS_i$ an open cover such that there exists
closed embedding $i_i:S\hookrightarrow\tilde S_i$ with $\tilde S_i\in\SmVar(k)$.
Recall (see section 5.2) that $D_{\mathcal D(1,0)fil,rh}(S/(\tilde S_I))\times_I D_{fil,c,k}(S_{\mathbb C}^{an})$ is the category  
\begin{itemize}
\item whose set of objects is the set of triples $\left\{(((M_I,F,W),u_{IJ}),(K,W),\alpha)\right\}$ with
\begin{eqnarray*} 
((M_I,F,W),u_{IJ})\in D_{\mathcal D(1,0)fil,rh}(S/(\tilde S_I)), \, (K,W)\in D_{fil,c,k}(S_{\mathbb C}^{an}), \\ 
\alpha:T(S/(\tilde S_I))((K,W)\otimes\mathbb C_{S_{\mathbb C}^{an}})\to DR(S)^{[-]}(((M_I,W),u_{IJ})^{an})
\end{eqnarray*}
where $\alpha$ is an morphism in $D_{fil}(S_{\mathbb C}^{an}/(\tilde S_{I,\mathbb C}^{an}))$,
\item and whose set of morphisms consists of 
\begin{equation*}
\phi=(\phi_D,\phi_C,[\theta]):(((M_{1I},F,W),u_{IJ}),(K_1,W),\alpha_1)\to(((M_{2I},F,W),u_{IJ}),(K_2,W),\alpha_2)
\end{equation*}
where $\phi_D:((M_1,F,W),u_{IJ})\to((M_2,F,W),u_{IJ})$ and $\phi_C:(K_1,W)\to (K_2,W)$ 
are morphisms and
\begin{eqnarray*}
\theta=(\theta^{\bullet},I(DR(S)(\phi^{an}_D))\circ I(\alpha_1),I(\alpha_2)\circ I(\phi_C\otimes I)): \\
I(T(S/(\tilde S_I))((K_1,W)\otimes\mathbb C_{S_{\mathbb C}^{an}}))[1]\to I(DR(S)(((M_{2I},W),u_{IJ})^{an}))  
\end{eqnarray*}
is an homotopy,  
$I:D_{fil}(S_{\mathbb C}^{an}/(\tilde S^{an}_{I,\mathbb C}))\to K_{fil}(S_{\mathbb C}^{an}/(\tilde S^{an}_{I,\mathbb C}))$
being the injective resolution functor, and for
\begin{itemize}
\item $\phi=(\phi_D,\phi_C,[\theta]):(((M_{1I},F,W),u_{IJ}),(K_1,W),\alpha_1)\to(((M_{2I},F,W),u_{IJ}),(K_2,W),\alpha_2)$
\item $\phi'=(\phi'_D,\phi'_C,[\theta']):(((M_{2I},F,W),u_{IJ}),(K_2,W),\alpha_2)\to(((M_{3I},F,W),u_{IJ}),(K_3,W),\alpha_3)$
\end{itemize}
the composition law is given by 
\begin{eqnarray*}
\phi'\circ\phi:=(\phi'_D\circ\phi_D,\phi'_C\circ\phi_C,
I(DR(S)(\phi^{'an}_D))\circ[\theta]+[\theta']\circ I(\phi_C\otimes I)[1]): \\
(((M_{1I},F,W),u_{IJ}),(K_1,W),\alpha_1)\to(((M_{3I},F,W),u_{IJ}),(K_3,W),\alpha_3),
\end{eqnarray*}
in particular for 
$(((M_I,F,W),u_{IJ}),(K,W),\alpha)\in D_{\mathcal D(1,0)fil,rh}(S/(\tilde S_I))\times_I D_{fil,c,k}(S_{\mathbb C}^{an})$,
\begin{equation*}
I_{(((M_I,F,W),u_{IJ}),(K,W),\alpha)}=((I_{M_I}),I_K,0),
\end{equation*}
\end{itemize}
together with the localization functor
\begin{eqnarray*}
(D(zar),I):C_{\mathcal D(1,0)fil,rh}(S/(\tilde S_I))\times_I D_{fil,c,k}(S_{\mathbb C}^{an})
\to D_{\mathcal D(1,0)fil,rh}(S/(\tilde S_I))\times_I D_{fil,c,k}(S_{\mathbb C}^{an}) \\
\to D_{\mathcal D(1,0)fil,rh,\infty}(S/(\tilde S_I))\times_I D_{fil,c,k}(S_{\mathbb C}^{an}).
\end{eqnarray*}
Note that if $\phi=(\phi_D,\phi_C,[\theta]):(((M_1,F,W),u_{IJ}),(K_1,W),\alpha_1)\to(((M_2,F,W),u_{IJ}),(K_2,W),\alpha_2)$
is a morphism in $D_{\mathcal D(1,0)fil,rh}(S/(\tilde S_I))\times_I D_{fil,c,k}(S_{\mathbb C}^{an})$
such that $\phi_D$ and $\phi_C$ are isomorphisms then $\phi$ is an isomorphism (see remark \ref{CGrem}).
Moreover,
\begin{itemize}
\item For 
$(((M_{I},F,W),u_{IJ}),(K,W),\alpha)\in D_{\mathcal D(1,0)fil,rh}(S/(\tilde S_I))\times_I D_{fil,c,k}(S_{\mathbb C}^{an})$, 
we set
\begin{equation*}
(((M_{I},F,W),u_{IJ}),(K,W),\alpha)[1]:=(((M_{I},F,W),u_{IJ})[1],(K,W)[1],\alpha[1]).
\end{equation*}
\item For 
\begin{equation*}
\phi=(\phi_D,\phi_C,[\theta]):(((M_{1I},F,W),u_{IJ}),(K_1,W),\alpha_1)\to(((M_{2I},F,W),u_{IJ}),(K_2,W),\alpha_2)
\end{equation*}
a morphism in $D_{\mathcal D(1,0)fil,rh}(S/(\tilde S_I))\times_I D_{fil,c,k}(S_{\mathbb C}^{an})$, 
we set (see \cite{CG} definition 3.12)
\begin{eqnarray*}
\Cone(\phi):=(\Cone(\phi_D),\Cone(\phi_C),((\alpha_1,\theta),(\alpha_2,0)))
\in D_{\mathcal D(1,0)fil,rh}(S/(\tilde S_I))\times_I D_{fil,c,k}(S_{\mathbb C}^{an}),
\end{eqnarray*}
$((\alpha_1,\theta),(\alpha_2,0))$ being the matrix given by the composition law, together with the canonical maps
\begin{itemize}
\item $c_1(-)=(c_1(\phi_D),c_1(\phi_C),0):(((M_{2I},F,W),u_{IJ}),(K_2,W),\alpha_2)\to\Cone(\phi)$
\item $c_2(-)=(c_2(\phi_D),c_2(\phi_C),0):\Cone(\phi)\to (((M_{1I},F,W),u_{IJ}),(K_1,W),\alpha_1)[1]$.
\end{itemize}
\end{itemize}

Let $S\in\Var(k)$. Let $S=\cup_{i=1}^sS_i$ an open cover such that there exists
closed embedding $i_i:S\hookrightarrow\tilde S_i$ with $\tilde S_i\in\SmVar(k)$.
Consider the category 
\begin{equation*}
(D_{\mathcal D(1,0)fil}(\tilde S_I)\times_ID_{fil,c,k}(\tilde S_{I,\mathbb C}^{an}))\in\Fun(\Gamma(\tilde S_I),\TriCat)
\end{equation*}
such that 
\begin{equation*}
(D_{\mathcal D(1,0)fil}(\tilde S_I)\times_ID_{fil,c,k}(\tilde S_{I,\mathbb C}^{an}))(\tilde S_I)=
D_{\mathcal D(1,0)fil}(\tilde S_I)\times_ID_{fil,c,k}(\tilde S_{I,\mathbb C}^{an})
\end{equation*}
\begin{itemize}
\item whose objects are 
$(((M_I,F,W),(K_I,W),\alpha_I),u_{IJ})\in (D_{\mathcal D(1,0)fil}(\tilde S_I)\times_ID_{fil,c,k}(\tilde S_{I,\mathbb C}^{an}))$ 
such that 
\begin{equation*}
((M_I,F,W),(K_I,W),\alpha_I)\in D_{\mathcal D(1,0)fil}(\tilde S_I)\times_ID_{fil,c,k}(\tilde S_{I,\mathbb C}^{an})
=:\mathcal D(\tilde S_I) 
\end{equation*}
and for $I\subset J$, 
\begin{eqnarray*}
u_{IJ}:((M_I,F,W),(K_I,W),\alpha_I)\to \\
p_{IJ*}((M_J,F,W),(K_J,W),\alpha_J):=(p_{IJ*}(M_J,F,W),p_{IJ*}(K_J,W),p_{IJ*}\alpha_J)
\end{eqnarray*}
are morphisms in $\mathcal D(\tilde S_I)$,
\item whose morphisms $m=(m_I):(((M_I,F,W),(K_I,W),\alpha_I),u_{IJ})\to (((M'_I,F,W),(K'_I,W),\alpha'_I),v_{IJ})$
is a family of morphism such that $v_{IJ}\circ m_I=p_{IJ*}m_J\circ u_{IJ}$ in $\mathcal D(\tilde S_I)$
\end{itemize}
We have then the identity functor
\begin{eqnarray*}
I_S:D_{\mathcal D(1,0)fil}(S/(\tilde S_I))\times_ID_{fil,c,k}(S_{\mathbb C}^{an})\to
(D_{\mathcal D(1,0)fil}(\tilde S_I)\times_ID_{fil,c,k}(\tilde S_{I,\mathbb C}^{an})), \\
(((M_I,F,W),u_{IJ}),(K,W),\alpha)\mapsto (((M_I,F,W),i_{I*}j_I^*(K,W),j_I^*\alpha),(u_{IJ},I,0)), \\ 
m=(m_I,n)\mapsto m=(m_I,i_*j_I^*n)
\end{eqnarray*}
which is a full embedding since by definition for $((M_I,F,W),u_{IJ})\in D_{\mathcal D(1,0)fil}(S/(\tilde S_I))$,
\begin{equation*}
u_{IJ}:(M_I,F,W)\to p_{IJ*}(M_J,F,W) 
\end{equation*}
are filtered Zariski local equivalences, i.e. isomorphisms in $D_{\mathcal D(1,0)fil}(\tilde S_I)$, and hence for 
$((M_I,F,W),u_{IJ}),(K,W),\alpha)\in D_{\mathcal D(1,0)fil}(S/(\tilde S_I))\times_ID_{fil,c,k}(S_{\mathbb C}^{an})$,
\begin{eqnarray*}
(u_{IJ},I,0):((M_I,F,W),i_{I*}j_I^*(K,W),j_I^*\alpha)\to \\
p_{IJ*}((M_J,F,W),i_{J*}j_J^*(K,W),j_J^*\alpha)=(p_{IJ*}(M_J,F,W),i_{I*}j_I^*(K,W),j_I^*\alpha)
\end{eqnarray*}
are isomorphisms in $\mathcal D(\tilde S_I)$.

\begin{defi}\label{IUSm}
For $h:U\to S$ a smooth morphism with $S,U\in\SmVar(k)$ 
and $h:U\xrightarrow{n}X\xrightarrow{f}S$ a compactification of $h$ with $n$ an open embedding, $X\in\SmVar(k)$
such that $D:=X\backslash U=\cup_{i=1}^sD_i\subset X$ is a normal crossing divisor, 
we denote by, using definition \ref{DHdgalpha} and definition \ref{wtildew}
\begin{eqnarray*}
I(U/S):h_{!Hdg}h^{!Hdg}\mathbb Z_S^{Hdg}\xrightarrow{:=} \\
(p_{S*}E_{zar}(\Omega^{\bullet}_{X\times S/S}\otimes_{O_{X\times S}}(n\times I)_{!Hdg}\Gamma_U^{\vee,Hdg}(O_{U\times S},F_b)),
\mathbb D_Sh_*E_{usu}\mathbb Q_{U_{\mathbb C}^{an}},h_!\alpha(U,\delta)) \\
\xrightarrow{((DR(X\times S/S)(\ad((n\times I)_{!Hdg},(n\times I)^*)(-)),0),I,0)} \\
(\Cone((\Omega_{/S}^{\Gamma,pr}(i_{D_i}\times I))_{i\in[1,\ldots,s]}:
p_{S*}E_{zar}(\Omega^{\bullet}_{X\times S/S}\otimes_{O_{X\times S}}\Gamma_X^{\vee,Hdg}(O_{X\times S},F_b))\to \\
(\cdots\to(p_{S*}E_{zar}(\Omega^{\bullet}_{D_I\times S/S}\otimes_{O_{D_I\times S}}\Gamma_{D_I}^{\vee,Hdg}(O_{D_I\times S},F_b)))
\to\cdots)),\mathbb D_Sh_*E_{usu}\mathbb Z_{U_{\mathbb C}^{an}},h_!\alpha(U,\delta)) \\
\xrightarrow{=:}
(\mathcal F_S^{FDR}(\mathbb Z(U/S)),\Bti_S^*\mathbb Z(U/S),\alpha(\mathbb Z(U/S)))
\end{eqnarray*}
the canonical isomorphism in $D_{\mathcal Dfil}(S)\times_ID_{c,k}(S_{\mathbb C}^{an})$, where 
\begin{itemize}
\item we recall that (see section 6.1)
\begin{equation*}
h^{!Hdg}\mathbb Z_S^{Hdg}=
(\Gamma_U^{\vee,Hdg}(O_{U\times S},F_b),\mathbb Z_{U^{an}_{\mathbb C}},\alpha(U))\in HM_{gm,k,\mathbb C}(U), 
\end{equation*}
\item $i_{D_i}:D_i\hookrightarrow X$ are the closed embeddings,
\item $\alpha(\mathbb Z(U/S)):=h_!\alpha(U,\delta):=T^w(h,\otimes)(-)\circ h_!\alpha(U,\delta)$ 
(see definition \ref{falpha}), with 
\begin{eqnarray*}
\alpha(U,\delta):=(DR(U)(\Omega_{(U\times U/U)/(U/pt)}(\Gamma_U^{\vee,Hdg}(O_{U\times U}))))^{-1}\circ\alpha(U):\\
\mathbb C_{U_{\mathbb C}^{an}}\to
DR(U)((p_{U*}E_{zar}(\Omega^{\bullet}_{U\times U/U}\otimes_{O_{U\times U}}\Gamma_U^{\vee,Hdg}(O_{U\times U})))^{an}),
\end{eqnarray*}
by the way we note that the following diagram in $C(U_{\mathbb C}^{an})$ commutes
\begin{equation*}
\xymatrix{\mathbb C_{U_{\mathbb C}^{an}}\ar[r]^{\alpha(U)} & \Omega^{\bullet}_{U_{\mathbb C}^{an}}=:DR(U)(O^{an}_U) \\
p_{U*}E_{usu}\Gamma_U^{\vee}\mathbb C_{{U\times U}_{\mathbb C}^{an}}
\ar[u]^{\ad(\delta^*_U,\delta_{U*})(-)}\ar[r]^{\alpha(U\times U)} &
DR(U)((p_{U*}E_{zar}(\Omega^{\bullet}_{U\times U/U}\otimes_{O_{U\times U}}\Gamma_U^{\vee,Hdg}(O_{U\times U})))^{an})
\ar[u]^{DR(U)(\Omega_{(U\times U/U)/(U/pt)}(\Gamma_U^{\vee,Hdg}(O_{U\times U})))}}
\end{equation*}
\end{itemize}
\end{defi}

\begin{lem}\label{thetam}
Let $S\in\SmVar(k)$. Let $g:U'/S\to U/S$ a morphism with $U/S:=(U,h),U'/S:=(U',h)\in\Var(k)^{sm}/S$.
Let $h:U\xrightarrow{n}X\xrightarrow{f}S$ a compactification of $h$ with $n$ an open embedding, $X\in\SmVar(k)$
such that $D:=X\backslash U=\cup_{i=1}^sD_i\subset X$ is a normal crossing divisor, 
Let $h':U\xrightarrow{n'}X'\xrightarrow{f'}S$ a compactification of $h'$ with $n'$ an open embedding, $X'\in\SmVar(k)$
such that $D':=X\backslash U=\cup_{i=1}^sD_i\subset X$ is a normal crossing divisor and such that
$g:U'\to U$ extend to $\bar g:X'\to X$, see definition-proposition \ref{RCHdef0}.
Then, using definition \ref{IUSm}, 
the following diagram in $D_{\mathcal Dfil}(S)\times_ID_{c,k}(S_{\mathbb C}^{an})$ commutes
\begin{eqnarray*}
\xymatrix{
h'_{!Hdg}h^{'!Hdg}\mathbb Z_S^{Hdg}\ar[rrr]^{I(U'/S)}\ar[d]_{\ad(g_{!Hdg},g^{!Hdg})(h^{!Hdg}\mathbb Z_S^{Hdg})} & \, & \, &
(\mathcal F_S^{FDR}(\mathbb Z(U'/S)),\Bti_S^*\mathbb Z(U'/S),\alpha(\mathbb Z(U'/S)))
\ar[d]^{(\Omega_{/S}^{\Gamma,pr}(R_S^{CH}(g)),\Bti_S^*(g),\theta(g))} \\ 
h_{!Hdg}h^{!Hdg}\mathbb Z_S^{Hdg}\ar[rrr]^{I(U/S)} & \, & \, &
(\mathcal F_S^{FDR}(\mathbb Z(U/S)),\Bti_S^*\mathbb Z(U/S),\alpha(\mathbb Z(U/S)))} 
\end{eqnarray*}
where 
\begin{eqnarray*}
\theta(g):=R_{\mathcal D}([\Gamma_g]):I(\Bti_S^*\mathbb Z(U'/S)\otimes\mathbb C)[1]
\to I(DR(S)(o_F\mathcal F_S^{FDR}(\mathbb Z(U/S))^{an}))
\end{eqnarray*}
is the homotopy given by the third term of the Deligne homology class of the graph $\Gamma_g\subset U'\times_S U$
(see definition \ref{Delkdef}) and 
$o_F:C_{\mathcal Dfil}(S)\to C_{\mathcal D}(S)$ is the forgetful functor 
and we recall (see section 6.1) that
$I:C(S_{\mathcal C}^{an}/\tilde S_{I,\mathbb C}^{an})\to K(S_{\mathcal C}^{an}/\tilde S_{I,\mathbb C}^{an})$
is the injective resolution functor.
\end{lem}

\begin{proof}
Immediate from definition.
\end{proof}

We now define the Hodge realization functor.

\begin{defi}\label{HodgeRealDAsing}
Let $k\subset\mathbb C$ a subfield. Let $S\in\Var(k)$. Let $S=\cup_{i=1}^sS_i$ an open cover such that there exists
closed embedding $i_i:S\hookrightarrow\tilde S_i$ with $\tilde S_i\in\SmVar(k)$.
We define the Hodge realization functor, using definition \ref{DRalgdefFunct}, 
definition \ref{bettidef}, and lemma \ref{thetam}
\begin{eqnarray*}
\mathcal F_S^{Hdg}:=(\mathcal F_S^{FDR},\Bti_S^*\otimes\mathbb Q):
C(\Var(k)^{sm}/S)\to D_{\mathcal D(1,0)fil}(S/(\tilde S_I))^0\times_I D_{fil,c,k}(S_{\mathbb C}^{an}) 
\end{eqnarray*}
first on objects and then on morphisms :
\begin{itemize}
\item for $F\in C(\Var(k)^{sm}/S)$, taking $(F,W)\in C_{fil}(\Var(k)^{sm}/S)$
such that $D(\mathbb A^1,et)(F,W)$ gives the weight structure on $D(\mathbb A^1,et)(F)$,
\begin{eqnarray*}
\mathcal F_S^{Hdg}(F):=(\mathcal F_S^{FDR}(F,W),\Bti_S^*(F,W)\otimes\mathbb Q,\alpha(F)):= \\
(e(S)_*\mathcal Hom((\hat R_{\tilde S_I}^{CH}(\rho_{\tilde S_I}^*Li_{I*}j_I^*(F,W)),
\hat R^{CH}(T^q(D_{IJ})(-))),(E_{zar}(\Omega^{\bullet,\Gamma,pr}_{/\tilde S_I},F_{DR}),T_{IJ})), \\
e(S)_*\underline{\sing}_{\mathbb D^*}\An_S^*L(F,W),\alpha(F)) 
\in D_{\mathcal D(1,0)fil}(S/(\tilde S_I))\times_I D_{fil,c,k}(S_{\mathbb C}^{an})
\end{eqnarray*}
where $\alpha(F)$ is the map in $D_{fil}(S_{\mathbb C}^{an}/(\tilde S_{I,\mathbb C}^{an}))$,
writing for short $DR(S):=DR(S)^{[-]}:=(DR(\tilde S_I)[-d_{\tilde S_I}])$
\begin{eqnarray*}
\alpha(F):T(S/(\tilde S_I))((\Bti_S^*(M,W))\otimes\mathbb C_{S_{\mathbb C}^{an}}):=
(i_{I*}j_I^*((e(S)_*\underline{\sing}_{\mathbb D^*}\An_S^*L(F,W))\otimes\mathbb C_S),I) \\
\xrightarrow{=}
T_S^{-1}(e(\tilde S_I)_*\underline{\sing}_{\mathbb D^*}\An_{\tilde S_I}^*Li_{I*}j_I^*(F,W)\otimes\mathbb C_{S^{an}_{\mathbb C}},
T(p_{IJ},\An)(Li_{I*}j_I^*(F,W))) 
\xrightarrow{=} \\
T_S^{-1}((((\cdot\to\oplus_{(U_{I\alpha},h_{I\alpha})\in V_I}h_{I\alpha !}h^!_{I\alpha}\mathbb C_{\tilde S_{I,\mathbb C}^{an}}
\xrightarrow{\ad(g^{\bullet !}_{I,\alpha,\beta},g^{\bullet}_{I,\alpha,\beta !})(-)} 
\oplus_{(U_{I\alpha},h_{I\alpha})\in V_I}h_{I\alpha !}h^!_{I\alpha}\mathbb C_{\tilde S_{I,\mathbb C}^{an}}\to\cdot),u_{IJ})),W) \\
\xrightarrow{T_S^{-1}(\alpha(\mathbb Z(U_{I\alpha}/\tilde S_I)),\theta(g^{\bullet}_{I,\alpha,\beta}))} \\
T_S^{-1}DR(S)(o_F(e(S)_*\mathcal Hom((\hat R_{\tilde S_I}^{CH}(\rho_{\tilde S_I}^*Li_{I*}j_I^*(F,W)),
\hat R^{CH}(T^q(D_{IJ})(-))),(E_{zar}(\Omega^{\bullet,\Gamma,pr}_{/\tilde S_I},F_{DR}),T_{IJ})))^{an}) \\
\xrightarrow{=:}
DR(S)((o_F\mathcal F^{FDR}_S(M,W))^{an})
\end{eqnarray*}
where
\begin{itemize}
\item $o_F:D_{\mathcal D(1,0)fil}(S/(\tilde S_I))\to D_{\mathcal D0fil}(S/(\tilde S_I))$ is the forgetful functor,
\item $T_S$ is the canonical functor
\begin{equation*}
T_S:D_{fil}(S_{\mathbb C}^{an}/(\tilde S_{I,\mathbb C}^{an}))\to D_{fil}((\tilde S_{I,\mathbb C}^{an})), \;
T_S((K_I,W),u_{IJ})=((K_I,W),u_{IJ}) 
\end{equation*}
which is a full embedding on the full subcategory 
\begin{equation*}
D_{fil}(S_{\mathbb C}^{an}/(\tilde S_{I,\mathbb C}^{an}))^{\sim}:=
D_{usu}(C_{fil}(S_{\mathbb C}^{an}/(\tilde S_{I,\mathbb C}^{an}))^{\sim})
\subset D_{fil}(S_{\mathbb C}^{an}/(\tilde S_{I,\mathbb C}^{an}))
\end{equation*}
where 
$C_{fil}(S_{\mathbb C}^{an}/(\tilde S_{I,\mathbb C}^{an}))^{\sim}\subset C_{fil}(S_{\mathbb C}^{an}/(\tilde S_{I,\mathbb C}^{an}))$
is the full subcategory consisting of $((K_I,W),u_{IJ})$ such that $u_{IJ}:K_I\to p_{IJ*}K_J$ are isomorphisms,
\item we denote by $g^n_{I,\alpha,\beta}:U_{I\alpha}\to U_{I\beta}$ which satisfy 
$h_{I\beta}\circ g^n_{I,\alpha,\beta}=h_{I\alpha}$ the morphisms in the canonical projective resolution
\begin{eqnarray*}
q:Li_{I*}j_I^*(F,W):=((\cdots\to\oplus_{(U_{I\alpha},h_{I\alpha})\in\Var(k)^{sm}/\tilde S_I}\mathbb Z(U_{I\alpha}/\tilde S_I)
\xrightarrow{(\mathbb Z(g^{\bullet}_{I,\alpha,\beta}))} \\
\oplus_{(U_{I\alpha},h_{I\alpha})\in\Var(k)^{sm}/\tilde S_I}\mathbb Z(U_{I\alpha}/\tilde S_I)\to\cdots),W)\to i_{I*}j_I^*(F,W),
\end{eqnarray*}
\end{itemize}
using lemma \ref{thetam}, 
\begin{equation*}
(\alpha(\mathbb Z(U_{I\alpha}/S)),\theta(g^{\bullet}_{I,\alpha,\beta}))
\end{equation*}
being the matrix given inductively by the composition law in 
$D_{\mathcal D(1,0)fil}(\tilde S_I)\times_I D_{fil,c,k}(\tilde S_{I,\mathbb C}^{an})$,
and the fact 
\begin{equation*}
DR(S)((o_F\mathcal F^{FDR}_S(M,W))^{an})\in D_{fil}(S_{\mathbb C}^{an}/(\tilde S_{I,\mathbb C}^{an}))^{\sim}
\end{equation*}
since $\mathcal F^{FDR}_S(M,W)\in D(DRM(S))$ by corollary \ref{FDRMHM}, that is we have the following isomorphism in 
$(D_{\mathcal D(1,0)fil}(\tilde S_I)\times_I D_{fil,c,k}(\tilde S_{I,\mathbb C}^{an}))$,
denoting for short $V_I:=\Var(k)^{sm}/\tilde S_I$
\begin{eqnarray*}
(I^{\bullet}(U_{I\alpha}/\tilde S_I)): 
(((\cdot\to\oplus_{(U_{I\alpha},h_{I\alpha})\in V_I}h_{I\alpha !Hdg}h^{!Hdg}_{I\alpha}\mathbb Q^{Hdg}_{\tilde S_I}
\xrightarrow{\ad(g^{\bullet,!Hdg}_{I,\alpha,\beta},g^{\bullet}_{I,\alpha,\beta !Hdg})(-)} \\
\oplus_{(U_{I\alpha},h_{I\alpha})\in V_I}h_{I\alpha !Hdg}h^{!Hdg}_{I\alpha}\mathbb Q^{Hdg}_{\tilde S_I}\to\cdot),u_{IJ}),W) \\
\xrightarrow{\sim}
I_S(\mathcal F_S^{Hdg}(F):=(\mathcal F_S^{FDR}(F,W),\Bti_{\tilde S_I}^*Li_{I*}j_I^*(F,W)\otimes\mathbb Q,\alpha(F)))
\end{eqnarray*}
\item for $m:F_1\to F_2$ a morphism in $C(\Var(k)^{sm}/S)$, taking $(F_1,W),(F_2,W)\in C_{fil}(\Var(k)^{sm}/S)$
such that $D(\mathbb A^1,et)(F_2,W)$ gives the weight structure on $D(\mathbb A^1,et)(F_2)$ 
$D(\mathbb A^1,et)(F_1,W)$ gives the weight structure on $D(\mathbb A^1,et)(F_1)$
and such that $m:(F_1,W)\to (F_2,W)$ is a filtered morphism, the morphism 
$\mathcal F_S^{Hdg}(m)$ in $D_{\mathcal D(1,0)fil}(S/(\tilde S_I))\times_I D_{fil,c,k}(S_{\mathbb C}^{an})$ is given by
\begin{eqnarray*}
\mathcal F_S^{Hdg}(m):&=&{I_S^{-,-}}^{-1}((I^{\bullet}(U_{I\alpha}/(\tilde S_I)))\circ
(\ad(l_{I\alpha,\beta}^{\bullet !Hdg},l^{\bullet}_{I\alpha,\beta !Hdg})(\mathbb Q_{U_{I\alpha}}^{Hdg})) 
\circ (I^{\bullet}(U_{I\alpha}/(\tilde S_I)))^{-1}) \\
&=&(\mathcal F_S^{FDR}(m),\Bti_S^*(m)\otimes\mathbb Q,\theta(m):=(\theta(l_{I\alpha,\beta}))):
\mathcal F_S^{Hdg}(F_1)\to\mathcal F_S^{Hdg}(F_2)
\end{eqnarray*}
using lemma \ref{thetam}, that is we have the following commutative diagram in 
$(D_{\mathcal D(1,0)fil}(\tilde S_I)\times_I D_{fil,c,k}(\tilde S_{I,\mathbb C}^{an}))$,
denoting for short $V_I:=\Var(k)^{sm}/\tilde S_I$,
\begin{eqnarray*}
\begin{tikzcd}
(((\cdot\to\oplus_{(U_{I\alpha},h_{I\alpha})\in V_I}h_{I\alpha !Hdg}h^{!Hdg}_{I\alpha}\mathbb Q^{Hdg}_{\tilde S_I}
\xrightarrow{A^{Hdg}_{g_{1I,\alpha,\beta}^{\bullet}}} 
\oplus_{(U_{I\alpha},h_{I\alpha})\in V_I}h_{I\alpha !Hdg}h^{!Hdg}_{I\alpha}\mathbb Q^{Hdg}_{\tilde S_I}\to\cdot),u_{IJ}),W)
\ar[r,"(I^{\bullet}(U_{I\alpha}/\tilde S_I))"]
\ar[d,"\ad(l_{I{,}\alpha{,}\beta}^{\bullet !Hdg}{,}l^{I{,}\bullet}_{\alpha{,}\beta !Hdg})(-)"'] & 
\mathcal F_S^{Hdg}(F_1)\ar[d,"\mathcal F_S^{Hdg}(m)=(\mathcal F_S^{FDR}(m){,}\Bti_S^*(m){,}(\theta(l_{I\alpha{,}\beta})))"'] \\
(((\cdot\to\oplus_{(U_{I\alpha},h_{I\alpha})\in V_I}h_{I\alpha !Hdg}h^{!Hdg}_{I\alpha}\mathbb Q^{Hdg}_{\tilde S_I}
\xrightarrow{A^{Hdg}_{g_{2I,\alpha,\beta}^{\bullet}}}
\oplus_{(U_{I\alpha},h_{I\alpha})\in V_I}h_{I\alpha !Hdg}h^{!Hdg}_{I\alpha}\mathbb Q^{Hdg}_{\tilde S_I}\to\cdot),u_{IJ}),W)
\ar[r,"(I^{\bullet}(U_{\alpha}/\tilde S_I))"] & \mathcal F_S^{Hdg}(F_2)
\end{tikzcd}
\end{eqnarray*}
where
\begin{itemize}
\item we denoted for short
$A^{Hdg}_{g_{1I,\alpha,\beta}^{\bullet}}:=
\ad(g^{\bullet,!Hdg}_{1I,\alpha,\beta},g^{\bullet}_{1I,\alpha,\beta !Hdg})(h^{!Hdg}_{I\alpha}\mathbb Z^{Hdg}_{\tilde S_I})$
\item we denoted for short
$A^{Hdg}_{g_{2I,\alpha,\beta}^{\bullet}}:=
\ad(g^{\bullet,!Hdg}_{2I,\alpha,\beta},g^{\bullet}_{2I,\alpha,\beta !Hdg})(h^{!Hdg}_{I\alpha}\mathbb Z^{Hdg}_{\tilde S_I})$
\item we denote by $g^n_{1I,\alpha,\beta}:U_{I\alpha}\to U_{I\beta}$, which satisfy 
$h_{I\beta}\circ g^n_{1I,\alpha,\beta}=h_{I\alpha}$, the morphisms in the canonical projective resolution
\begin{eqnarray*}
q:Li_{I*}j_I^*(F_1,W):=((\cdots\to\oplus_{(U_{I\alpha},h_{I\alpha})\in\Var(k)^{sm}/\tilde S_I}\mathbb Z(U_{I\alpha}/\tilde S_I)
\xrightarrow{(\mathbb Z(g^{\bullet}_{1I,\alpha,\beta}))} \\
\oplus_{(U_{I\alpha},h_{I\alpha})\in\Var(k)^{sm}/\tilde S_I}\mathbb Z(U_{I\alpha}/\tilde S_I)\to\cdots),W)\to i_{I*}j_I^*(F_1,W)
\end{eqnarray*}
\item we denote by $g^n_{2I,\alpha,\beta}:U_{I\alpha}\to U_{I\beta}$, which satisfy 
$h_{I\beta}\circ g^n_{2I,\alpha,\beta}=h_{\alpha}$, the morphisms in the canonical projective resolution
\begin{eqnarray*}
q:Li_{I*}j_I^*(F_2,W):=((\cdots\to\oplus_{(U_{I\alpha},h_{I\alpha})\in\Var(k)^{sm}/\tilde S_I}\mathbb Z(U_{I\alpha}/\tilde S_I)
\xrightarrow{(\mathbb Z(g^{\bullet}_{2I,\alpha,\beta}))} \\
\oplus_{(U_{I\alpha},h_{I\alpha})\in\Var(k)^{sm}/\tilde S_I}\mathbb Z(U_{I\alpha}/\tilde S_I)\to\cdots),W)\to i_{I*}j_I^*(F_2,W)
\end{eqnarray*}
\item we denote by $l_{I\alpha,\beta}^n:U_{I\alpha}\to U_{I\beta}$ which satisfy 
$h_{I\beta}\circ l^n_{I\alpha,\beta}=h_{I\alpha}$ and 
$l^{n+1}_{I\alpha,\beta}\circ g^n_{1I\alpha,\beta}=g^n_{2I\alpha,\beta}\circ l^n_{I\alpha,\beta}$ 
the morphisms in the morphism of canonical projective resolutions
\begin{eqnarray*}
Li_{I*}j_I^*(m):Li_{I*}j_I^*(F_1,W):=
((\cdots\to\oplus_{(U_{I\alpha},h_{I\alpha})\in\Var(k)^{sm}/\tilde S_I}\mathbb Z(U_{I\alpha}/\tilde S_I)\to\cdots),W) 
\xrightarrow{(\mathbb Z(l^{\bullet}_{I\alpha,\beta}))} \\
((\cdots\to\oplus_{(U_{I\alpha},h_{I\alpha})\in\Var(k)^{sm}/\tilde S_I}\mathbb Z(U_{\alpha}/\tilde S_I)\to\cdots),W)
=:Li_{I*}j_I^*(F_2,W),
\end{eqnarray*}
\item the maps $I^{\bullet}(U_{I\alpha})$ are given by definition \ref{IUSm} and lemma \ref{thetam}.
\end{itemize}
\end{itemize}
Obviously $\mathcal F_S^{Hdg}(F[1])=\mathcal F_S^{Hdg}(F)[1]$ and 
$\mathcal F_S^{Hdg}(\Cone(m))=\Cone(\mathcal F_S^{Hdg}(m))$. 
This functor induces by proposition \ref{projwach} and remark \ref{CGrem} the functor
\begin{eqnarray*}
\mathcal F_S^{Hdg}:=(\mathcal F_S^{FDR},\Bti_S^*\otimes\mathbb Q):
\DA(S)\to D_{\mathcal D(1,0)fil}(S/(\tilde S_I))\times_I D_{fil,c,k}(S_{\mathbb C}^{an}), \\ 
M=D(\mathbb A^1,et)(F)\mapsto
\mathcal F_S^{Hdg}(M):=\mathcal F_S^{Hdg}(F)=(\mathcal F_S^{FDR}(M),\Bti_S^*M\otimes\mathbb Q,\alpha(M)),
\end{eqnarray*}
with $\alpha(M)=\alpha(F)$. 
\end{defi}

We now give the functoriality with respect to the five operation using the De Rahm realization case and the Betti realization case :

\begin{prop}\label{TbtiTFDR}
\begin{itemize}
\item[(i)]Let $g:T\to S$ a morphism with $T,S\in\Var(k)$. 
Assume there exists a factorization $g:T\xrightarrow{l}Y\times S\xrightarrow{p}S$, with $Y\in\SmVar(k)$,
$l$ a closed embedding and $p$ the projection.
Let $S=\cup_{i\in I}S_i$ an open cover and $i_i:S_i\hookrightarrow\tilde S_i$ closed embeddings with $\tilde S_i\in\SmVar(k)$.
Then, $\tilde g_I:Y\times\tilde S_I\to\tilde S_I$ is a lift of $g_I=g_{|T_I}:T_I\to S_I$
and we have closed embeddings $i'_I:=i_I\circ l\circ j'_I:T_I\hookrightarrow Y\times\tilde S_I$.
Then, for $M\in DA_c(S)$, the following diagram commutes :
\begin{equation*}
\xymatrix{g^{*w}\Bti_S^*M\otimes\mathbb C\ar[rr]^{g^*(\alpha(M))}\ar[d]_{T(g,bti)(M)} & \, & 
DR(T)^{[-]}(o_F(g^{\hat{*}mod}_{Hdg}\mathcal F_S^{FDR}(M))^{an})
\ar[d]^{DR(T)^{[-]}((T(g,\mathcal F^{FDR})(M))^{an})} \\
\Bti_T^*g^*M\otimes\mathbb C\ar[rr]^{\alpha(g^*M)} &  \, & DR(T)^{[-]}(o_F(\mathcal F_T^{FDR}(g^*M))^{an})},
\end{equation*}
see section 5, definition \ref{TgDRdefsing} and definition \ref{TgBti}
\item[(ii)]Let $f:T\to S$ a morphism with $T,S\in\QPVar(k)$. Then, for $M\in DA_c(T)$,the following diagram commutes :
\begin{equation*}
\xymatrix{Rf_{*w}\Bti_T^*M\otimes\mathbb C\ar[rr]^{f_*(\alpha(M)} & \, & 
DR(S)^{[-]}(o_F(Rf^{Hdg}_*\mathcal F_T^{FDR}(M))^{an}) \\
\Bti_S^*Rf_*M\otimes\mathbb C\ar[u]^{T_*(f,bti)(M)}\ar[rr]^{\alpha(Rf_*M)} & \, &
DR(S)^{[-]}(o_F(\mathcal F_S^{FDR}(Rf_*M))^{an})\ar[u]_{DR(S)^{[-]}((T_*(f,\mathcal F^{FDR})(M))^{an})}}
\end{equation*}
see section 5, definition \ref{SixTalg} and definition \ref{TBtiSix}
\item[(iii)]Let $f:T\to S$ a morphism with $T,S\in\QPVar(k)$. Then, for $M\in DA_c(T)$,the following diagram commutes :
\begin{equation*}
\xymatrix{Rf_{!w}\Bti_T^*M\otimes\mathbb C\ar[d]_{T_!(f,bti)(M)}\ar[rr]^{f_!(\alpha(M))} & \, &  
DR(S)^{[-]}(o_F(Rf^{Hdg}_!\mathcal F^T_{DR}(M))^{an})\ar[d]^{DR(S)^{[-]}((T_!(f,\mathcal F_{FDR})(M))^{an})} \\
\Bti_S^*Rf_!M\otimes\mathbb C\ar[rr]^{\alpha(Rf_!M)} & \, & DR(S)^{[-]}(o_F(\mathcal F^S_{FDR}(Rf_!M))^{an})} 
\end{equation*}
see section 5, definition \ref{SixTalg} and definition \ref{TBtiSix}.
\item[(iv)]Let $f:T\to S$ a morphism with $T,S\in\QPVar(k)$. Then, for $M\in DA_c(S)$,the following diagram commutes :
\begin{equation*}
\xymatrix{f^{!w}\Bti_S^*M\otimes\mathbb C\ar[rr]^{f^!(\alpha(M))} & \, & 
DR(T)^{[-]}(o_F(f^{*mod}_{Hdg}\mathcal F_S^{FDR}(M))^{an}) \\
\Bti_T^*f^!M\otimes\mathbb C\ar[rr]^{\alpha(f^!M)}\ar[u]^{T^!(f,bti)(M)} &  \, &
DR(T)^{[-]}(o_F(\mathcal F_T^{FDR}(f^!M))^{an})\ar[u]_{DR^{[-]}(T)((T^!(g,\mathcal F^{FDR})(M))^{an})}}
\end{equation*}
see section 5, definition \ref{SixTalg} and definition \ref{TBtiSix}.
\item[(v)] Let $S\in\Var(k)$. Then, for $M,N\in DA_c(S)$,the following diagram commutes :
\begin{equation*}
\xymatrix{\Bti_S^*M\otimes\Bti_S^*N\otimes\mathbb C
\ar[rrr]^{\alpha(M)\otimes\alpha(N)}\ar[d]_{T(\otimes,bti)(M,N)} & \, & \, &  
DR(S)(o_F(\mathcal F_S^{FDR}(M)\otimes_{O_S}\mathcal F_S^{FDR}(N))^{an})
\ar[d]^{DR(S)((T(\otimes,\mathcal F^{DR})(M,N))^{an})} \\
\Bti_S^*(M\otimes N)\otimes\mathbb C\ar[rrr]^{(\alpha(M\otimes N))} & \, & \, &  DR(S)((\mathcal F^S_{DR}(M\otimes N))^{an})}
\end{equation*}
see definition \ref{SixTalg} and definition \ref{TBtiSix}.
\end{itemize}
\end{prop}

\begin{proof}

\noindent(i): Follows from the following commutative diagram in 
$(D_{\mathcal D(1,0)fil}(Y\times\tilde S_I)\times_I D_{fil,c,k}(Y\times\tilde S_{I,\mathbb C}^{an}))$,
\begin{eqnarray*}
\begin{tikzcd}
\shortstack{
$(((\to\oplus_{(U_{I\alpha},h_{I\alpha})\in V_I}\tilde g_I^{*Hdg}h_{I\alpha !Hdg}h^{!Hdg}_{I\alpha}\mathbb Z^{Hdg}_{\tilde S_I}
\xrightarrow{A^{Hdg}_{g_{I,\alpha,\beta}^{\bullet}}}$ \\ 
$\oplus_{(U_{I\alpha},h_{I\alpha})\in V_I}h_{I\alpha !Hdg}h^{!Hdg}_{I\alpha}\mathbb Z^{Hdg}_{\tilde S_I}\to),u_{IJ}),W)$}
\ar[rr,"(\tilde g_I^{*Hdg}I^{\bullet}(U_{I\alpha}/\tilde S_I))"]\ar[d,"T^{Hdg}(\tilde g_I{,}h_I)(-)"'] & \, & 
\shortstack{$(g^{\hat*mod}_{Hdg}\mathcal F_T^{FDR}(F)$, \\ $g^{*w}\Bti_S^*(F,W),g^*(\alpha(F)))$}
\ar[d,"(T(g{,}\mathcal F^{FDR})(M){,}T(g{,}\Bti)(M){,}0)"'] \\
\shortstack{
$(((\to\oplus_{(U'_{I\alpha},h_{I\alpha})\in W_I}h'_{I\alpha !Hdg}h^{'!Hdg}_{I\alpha}\mathbb Z^{Hdg}_{Y\times\tilde S_I}
\xrightarrow{A^{Hdg}_{g_{I,\alpha,\beta}^{'\bullet}}}$ \\
$\oplus_{(U'_{I\alpha},h'_{I\alpha})\in W_I}h'_{I\alpha !Hdg}h^{'!Hdg}_{I\alpha}\mathbb Z^{Hdg}_{Y\times\tilde S_I}\to),u_{IJ}),W)$}
\ar[rr,"(I^{\bullet}(U'_{\alpha}/Y\times\tilde S_I))"] & \, &
\shortstack{$(\mathcal F_T^{FDR}(g^*F)$, \\ $\Bti_T^*(g^*F,W),\alpha(g^*F))$}
\end{tikzcd}
\end{eqnarray*}
where, we have denoted for short $V_I:=\Var(k)^{sm}/\tilde S_I$ and $W_I:=\Var(k)^{sm}/Y\times\tilde S_I$,
\begin{itemize}
\item we denoted for short
$A^{Hdg}_{g_{I,\alpha,\beta}^{\bullet}}:=
\ad(g^{\bullet,!Hdg}_{I,\alpha,\beta},g^{\bullet}_{I,\alpha,\beta !Hdg})(h^{!Hdg}_{I\alpha}\mathbb Z^{Hdg}_{\tilde S_I})$
\item we denoted for short
$A^{Hdg}_{g_{I,\alpha,\beta}^{'\bullet}}:=
\ad(g^{'\bullet,!Hdg}_{I,\alpha,\beta},g^{'\bullet}_{I,\alpha,\beta !Hdg})(h^{'!Hdg}_{I\alpha}\mathbb Z^{Hdg}_{Y\times\tilde S_I})$
\item we denote by $g^n_{I,\alpha,\beta}:U_{I\alpha}\to U_{I\beta}$, which satisfy 
$h_{I\beta}\circ g^n_{I,\alpha,\beta}=h_{I\alpha}$, the morphisms in the canonical projective resolution
\begin{eqnarray*}
q:Li_{I*}j_I^*(F,W):=(\cdots\to\oplus_{(U_{I\alpha},h_{I\alpha})\in\Var(k)^{sm}/\tilde S_I}\mathbb Z(U_{I\alpha}/\tilde S_I)
\xrightarrow{(\mathbb Z(g^{\bullet}_{I,\alpha,\beta}))} \\
\oplus_{(U_{I\alpha},h_{I\alpha})\in\Var(k)^{sm}/\tilde S_I}\mathbb Z(U_{I\alpha}/\tilde S_I)\to\cdots)\to i_{I*}j_I^*(F,W)
\end{eqnarray*}
\item we denote by $g^{'n}_{I,\alpha,\beta}:U'_{I\alpha}\to U'_{I\beta}$, which satisfy 
$h'_{I\beta}\circ g^{'n}_{I,\alpha,\beta}=h'_{\alpha}$, the morphisms in the canonical projective resolution
\begin{eqnarray*}
q:Li'_{I*}j_I^{'*}(g^*F,W):=(\cdots\to\oplus_{(U'_{I\alpha},h'_{I\alpha})\in\Var(k)^{sm}/Y\times\tilde S_I}
\mathbb Z(U'_{I\alpha}/Y\times\tilde S_I)
\xrightarrow{(\mathbb Z(g^{'\bullet}_{I,\alpha,\beta}))} \\
\oplus_{(U'_{I\alpha},h'_{I\alpha})\in\Var(k)^{sm}/Y\times\tilde S_I}
\mathbb Z(U'_{I\alpha}/Y\times\tilde S_I)\to\cdots)\to i'_{I*}j_I^{'*}(g^*F,W)
\end{eqnarray*}
\end{itemize}

\noindent(ii): Follows from (i) by adjonction.

\noindent(iii): The closed embedding case is given by (ii) and the smooth projection case follows from (i) by adjonction.

\noindent(iv): Follows from (iii) by adjonction.

\noindent(v):Obvious

\end{proof}

We can now state the following key proposition and the main theorem:

\begin{prop}\label{keyHdg}
Let $k\subset\mathbb C$ a subfield.
\begin{itemize}
\item[(i)] Let $S\in\Var(k)$. Let $S=\cup_i S_i$ an open cover such that there exist closed embeddings
$i_i:S_i\hookrightarrow\tilde S_i$ with $\tilde S_i\in\SmVar(k)$. Then
we have the isomorphism in $D_{\mathcal D(1,0)fil}(S/(\tilde S_I))\times_ID_{fil,c,k}(S_{\mathbb C}^{an})$
\begin{eqnarray*}
\mathcal F_S^{Hdg}(\mathbb Z_S)\xrightarrow{:=}
(\mathcal F_S^{FDR}(\mathbb Z_S,W),\Bti_S^*(\mathbb Z_S,W)\otimes\mathbb Q,\alpha(\mathbb Z_S)) \\
\xrightarrow{((\Omega^{\Gamma,pr}_{/\tilde S_I}(\hat R^{CH}(\ad(i_I^*,i_{I*})(\Gamma^{\vee,w}_{S_I}\mathbb Z_{\tilde S_I})))),I,0)} \\
I_S^{-1}((e(S)_*\mathcal Hom((\hat R^{CH}(\Gamma^{\vee,w}_{S_I}\mathbb Z_{\tilde S_I}),\hat R^{CH}(x_{IJ})),
(E_{zar}(\Omega^{\bullet,\Gamma,pr}_{/\tilde S_I},F_{DR}),T_{IJ})),
T(S/(\tilde S_I))(\mathbb Q_{S_{\mathbb C}^{an}}^w), \\ \Gamma_{S_I}^{\vee,w}\alpha(\tilde S_I,\delta))) 
\xrightarrow{=}
\iota_S((\Gamma^{\vee,Hdg}_{S_I}(O_{\tilde S_I},F_b),x_{IJ}),\mathbb Q_{S_{\mathbb C}^{an}}^w,\alpha(S))
=:\iota_S(\mathbb Q_S^{Hdg})
\end{eqnarray*}
with (see section 6.1) $j_I^*\mathbb Q_{S_{\mathbb C}^{an}}^w=i_I^*\Gamma^{\vee,w}_{S_I}\mathbb Q_{\tilde S_I}$ and
\begin{eqnarray*}
\alpha(S):T(S/(\tilde S_I))((\mathbb Q_{S_{\mathbb C}^{an}}^w)\otimes\mathbb C_{S^{an}_{\mathbb C}})
\xrightarrow{=}(\Gamma_{S_I}^{\vee,w}\mathbb C_{\tilde S^{an}_{I,\mathbb C}},t_{IJ})
\xrightarrow{(\Gamma_{S_I}^{\vee}\alpha(\tilde S_I))}
DR(S)(o_F(\Gamma^{\vee,Hdg}_{S_I}(O_{\tilde S_I},F_b),x_{IJ}))
\end{eqnarray*}
\item[(ii)]Let $f:X\to S$ a morphism with $X,S\in\Var(k)$, $X$ quasi-projective.
Consider a factorization $f:X\xrightarrow{l}Y\times S\xrightarrow{p_S}S$ with $Y=\mathbb P^{N,o}\subset\mathbb P^N$ an open subset,
$l$ a closed embedding and $p_S$ the projection. Let $S=\cup_i S_i$ an open cover such that there exist closed embeddings
$i_i:S_i\hookrightarrow\tilde S_i$ with $\tilde S_i\in\SmVar(\mathbb C)$. 
Recall that $S_I:=\cap_{i\in I} S_i$, $X_I=f^{-1}(S_I)$, and $\tilde S_I:=\Pi_{i\in I}\tilde S_i$. Then,
using proposition \ref{TbtiTFDR}(iii), the maps of definition \ref{SixTalg} and definition \ref{TBtiSix} gives
an isomorphism in $D_{\mathcal D(1,0)fil}(S/(\tilde S_I))\times_ID_{fil,c,k}(S_{\mathbb C}^{an})$
\begin{eqnarray*}
(T_!(f,\mathcal F^{FDR})(\mathbb Z_X,W),T_!(f,\Bti)(\mathbb Z_X,W),0): \\
\mathcal F_S^{Hdg}(M^{BM}(X/S)):=
(\mathcal F_S^{FDR}(Rf_!(\mathbb Z_X,W)),\Bti_S^*Rf_!(\mathbb Z_X,W)\otimes\mathbb Q,\alpha(Rf_!\mathbb Z_X)) \\
\xrightarrow{\sim} 
(Rf_{Hdg!}(\Gamma^{\vee,Hdg}_{X_I}(O_{Y\times\tilde S_I},F_b),x_{IJ}(X/S)),Rf_{!w}\mathbb Q_{X^{an}}^w,
f_!(\alpha(X)))=:\iota_S(Rf_{!Hdg}\mathbb Q_X^{Hdg}).
\end{eqnarray*}
with
\begin{equation*}
\mathbb Q_X^{Hdg}:=((\Gamma^{\vee,Hdg}_{X_I}(O_{Y\times\tilde S_I},F_b),x_{IJ}(X/Y\times S)),
\mathbb Q_{X_{\mathbb C}^{an}}^w,\alpha(X))\in C(MHM_{gm,k,\mathbb C}(X))
\end{equation*}
\end{itemize}
\end{prop}

\begin{proof}
\noindent(i):Follows from proposition \ref{projwach}.

\noindent(ii): Follows from (i) by proposition \ref{TbtiTFDR}(iii),theorem \ref{mainthm}(i) and theorem \ref{mainBti}(i).
\end{proof}

The main theorem of this section is the following :

\begin{thm}\label{main}
Let $k\subset\mathbb C$ a subfield.
\begin{itemize}
\item[(i)] For $S\in\Var(k)$, we have $\mathcal F_S^{Hdg}(\DA_c(S))\subset D(MHM_{gm,k,\mathbb C}(S))$,  
\begin{equation*}
\iota_S:D(MHM_{gm,k,\mathbb C}(S))\hookrightarrow D_{\mathcal D(1,0)fil}(S/(\tilde S_I))\times_I D_{fil,c,k}(S_{\mathbb C}^{an}) 
\end{equation*}
being a full embedding by theorem \ref{Bek}.
\item[(ii)] The Hodge realization functor $\mathcal F_{Hdg}(-)$ define a morphism of 2-functor on $\Var(k)$
\begin{equation*}
\mathcal F^{Hdg}_{-}:\Var(k)\to(\DA_c(-)\to D(MHM_{gm,k,\mathbb C}(-)))
\end{equation*}
whose restriction to $\QPVar(k)$ is an homotopic 2-functor in sense of Ayoub. More precisely,
\begin{itemize}
\item[(ii0)] for $g:T\to S$ a morphism, with $T,S\in\QPVar(k)$, and $M\in\DA_c(S)$, the
the maps of definition \ref{TgDRdefsing} and of definition \ref{TgBti} induce an isomorphism in $D(MHM_{gm,k,\mathbb C}(T))$
\begin{eqnarray*}
T(g,\mathcal F^{Hdg})(M):=(T(g,\mathcal F^{FDR})(M),T(g,bti)(M),0): \\
g^{\hat*Hdg}\mathcal F_S^{Hdg}(M):=
\iota_T^{-1}(g^{\hat{*}mod}_{Hdg}\mathcal F_S^{FDR}(M),g^*\Bti_S(M)\otimes\mathbb Q,g^*(\alpha(M))) \\
\xrightarrow{\sim}
\iota_T^{-1}(\mathcal F_T^{FDR}(g^*M),\Bti_T^*(g^*M)\otimes\mathbb Q,\alpha(g^*M))=:\mathcal F_T^{Hdg}(g^*M),
\end{eqnarray*} 
\item[(ii1)] for $f:T\to S$ a morphism, with $T,S\in\QPVar(k)$, and $M\in\DA_c(T)$,  
the maps of definition \ref{SixTalg} and of definition \ref{TBtiSix} induce an isomorphism in $D(MHM_{gm,k,\mathbb C}(S))$
\begin{eqnarray*}
T_*(f,\mathcal F^{Hdg})(M):=(T_*(f,\mathcal F^{FDR})(M),T_*(f,bti)(M),0): \\
Rf_{Hdg*}\mathcal F_T^{Hdg}(M):=
\iota_S^{-1}(Rf^{Hdg}_*\mathcal F_T^{FDR}(M),Rf_*\Bti_T(M)\otimes\mathbb Q,f_*(\alpha(M))) \\ 
\xrightarrow{\sim}
\iota_S^{-1}(\mathcal F_S^{FDR}(Rf_*M),\Bti_S^*(Rf_*M)\otimes\mathbb Q,\alpha(Rf_*M))=:\mathcal F_S^{Hdg}(Rf_*M),
\end{eqnarray*}  
\item[(ii2)] for $f:T\to S$ a morphism, with $T,S\in\QPVar(k)$, and $M\in\DA_c(T)$, 
the maps of definition \ref{SixTalg} and of definition \ref{TBtiSix} induce an isomorphism in $D(MHM_{gm,k,\mathbb C}(S))$
\begin{eqnarray*}
T_!(f,\mathcal F^{Hdg})(M):=(T_!(f,\mathcal F^{FDR})(M),T_!(f,bti)(M),0): \\
Rf_{!Hdg}\mathcal F_T^{Hdg}(M):=
\iota_S^{-1}(Rf^{Hdg}_!\mathcal F_T^{FDR}(M),Rf_!\Bti_T^*(M)\otimes\mathbb Q,f_!(\alpha(M))) \\ 
\xrightarrow{\sim}
\iota_S^{-1}(\mathcal F_S^{FDR}(Rf_!M),\Bti_S^*(Rf_!M)\otimes\mathbb Q,\alpha(f_!M))=:\mathcal F_S^{Hdg}(Rf_!M),
\end{eqnarray*} 
\item[(ii3)] for $f:T\to S$ a morphism, with $T,S\in\QPVar(k)$, and $M\in\DA_c(S)$,
the maps of definition \ref{SixTalg} and of definition \ref{TBtiSix} induce an isomorphism in $D(MHM_{gm,k,\mathbb C}(T))$
\begin{eqnarray*}
T^!(f,\mathcal F^{Hdg})(M):=(T^!(f,\mathcal F^{FDR})(M),T^!(f,bti)(M),0): \\
f^{*Hdg}\mathcal F_S^{Hdg}(M):=
\iota_T^{-1}(f^{*mod}_{Hdg}\mathcal F_S^{FDR}(M),f^!\Bti_S(M)\otimes\mathbb Q,f^!(\alpha(M))) \\ 
\xrightarrow{\sim}
\iota_T^{-1}(\mathcal F_T^{FDR}(f^!M),\Bti_T^*(f^!M)\otimes\mathbb Q,\alpha(f^!M))=:\mathcal F_T^{Hdg}(f^!M),
\end{eqnarray*}
\item[(ii4)] for $S\in\Var(k)$, and $M,N\in\DA_c(S)$,
the maps of definition \ref{SixTalg} and of definition \ref{TBtiSix} 
induce an isomorphism in $D(MHM_{gm,k,\mathbb C}(S))$
\begin{eqnarray*}
T(\otimes,\mathcal F^{Hdg})(M,N):=(T(\otimes,\mathcal F_S^{FDR})(M,N),T(\otimes,bti)(M,N),0): \\
\iota_S^{-1}(\mathcal F_S^{FDR}(M)\otimes^{Hdg}_{O_S}\mathcal F_S^{FDR}(N),
\Bti_S(M)\otimes\Bti_S(N)\otimes\mathbb Q,\alpha(M)\otimes\alpha(N)) \\
\xrightarrow{\sim}\mathcal F_S^{Hdg}(M\otimes N):=
\iota_S^{-1}(\mathcal F_S^{FDR}(M\otimes N),\Bti_S(M\otimes N)\otimes\mathbb Q,\alpha(M\otimes N)).  
\end{eqnarray*}
\end{itemize}
\item[(iii)] For $S\in\Var(k)$, the following diagram commutes :
\begin{equation*}  
\xymatrix{\Var(k)/S\ar[rrr]^{MH(/S)}\ar[d]_{M(/S)} & \, & \, & D(MHM_{gm,k,\mathbb C}(S))\ar[d]^{\iota^S} \\
\DA(S)\ar[rrr]^{\mathcal F_S^{Hdg}} & \, & \, & D_{\mathcal D(1,0)fil}(S/(\tilde S_I))\times_I D_{fil,c,k}(S_{\mathbb C}^{an})}
\end{equation*}
\end{itemize}
\end{thm}

\begin{proof} 

\noindent(i): 
Let $M\in\DA_c(S)$. There exist by definition of constructible motives an isomorphism in $\DA(S)$ 
\begin{equation*}
w(M):M\xrightarrow{\sim}\Cone(M(X_0/S)[d_0]\xrightarrow{m_1}\cdots\xrightarrow{m_m} M(X_m/S)[d_m]),
\end{equation*}
with $f_n:X_n\to S$ morphisms and $X_n\in\QPVar(k)$.
This gives the isomorphism in $D_{\mathcal D(1,0)fil}(S/(\tilde S_I))\times_I D_{fil,c,k}(S_{\mathbb C}^{an})$
\begin{equation*}
\mathcal F_S^{Hdg}(w(M)):\mathcal F_S^{Hdg}(M)\xrightarrow{\sim}
\Cone(\mathcal F_S^{Hdg}(M(X_0/S))[d_0]\xrightarrow{\mathcal F_S^{Hdg}(m_1)}\cdots
\xrightarrow{\mathcal F_S^{Hdg}(m_1)}\mathcal F_S^{Hdg}(M(X_m/S))[d_m]),
\end{equation*}
On the other hand, by proposition \ref{keyHdg}(i), we have 
\begin{equation*}
\mathcal F_S^{Hdg}(M(X_n/S))\xrightarrow{\sim}Rf_{!Hdg}\mathbb Q_X^{Hdg}\in D(MHM_{gm,k,\mathbb C}(S)).
\end{equation*}
This prove (i).

\noindent(ii0): Follows from theorem \ref{mainthm}(i), proposition \ref{TbtiTFDR}(i) and theorem \ref{mainBti}.

\noindent(ii1): Follows from theorem \ref{mainthm}(iii), proposition \ref{TbtiTFDR}(ii), 
and theorem \ref{mainBti}(iii).

\noindent(ii2):Follows from theorem \ref{mainthm}(ii), proposition \ref{TbtiTFDR}(iii), 
and theorem \ref{mainBti}(ii).

\noindent(ii3): Follows from theorem \ref{mainthm}(iv), proposition \ref{TbtiTFDR}(iv), 
and theorem \ref{mainBti}(iv).

\noindent(ii4):Follows from theorem \ref{mainthm}(v), proposition \ref{TbtiTFDR}(v) and
theorem \ref{mainBti}(v).

\noindent(iii): By (ii), for $g:X'/S\to X/S$ a morphism, with $X',X,S\in\Var(k)$ 
and $X/S=(X,f)$, $X'/S=(X',f')$, we have by adjonction the following commutative diagram
\begin{equation*}
\xymatrix{\mathcal F_S^{Hdg}(M(X'/S)=f'_!f^{'!}\mathbb Z_S=f_!g_!g^!f^!\mathbb Z_S)
\ar[d]_{T_!(f',\mathcal F^{Hdg})(f^{'!}M(X'/S))\circ T^!(f',\mathcal F^{Hdg})(M(X'/S))}
\ar[rr]^{\mathcal F_S^{Hdg}(M(/S)(g)=f_!\ad(g_!,g^!)(f^!\mathbb Z_S))} & \, & 
\mathcal F_S^{Hdg}(M(X/S)=f_!f^!\mathbb Z_S)
\ar[d]_{T_!(f,\mathcal F^{Hdg})(f^!M(X/S))\circ T^!(f,\mathcal F^{Hdg})(M(X/S))} \\
MH(X'/S):=Rf'_{!Hdg}f^{'!Hdg}\mathbb Z^{Hdg}_S=f_{!Hdg}g_{!Hdg}g^{!Hdg}f^{!Hdg}\mathbb Z^{Hdg}_S
\ar[rr]^{f_{!Hdg}\ad(g_{!Hdg},g^{!Hdg})(f^{!Hdg}\mathbb Z^{Hdg}_S)} & \, &  
MH(X/S):=f_{!Hdg}f^{!Hdg}\mathbb Z^{Hdg}_S}.
\end{equation*}
where the left and right columns are isomorphisms by (ii). This proves (iii).
\end{proof}

The theorem \ref{main} gives immediately the following :

\begin{cor}
Let $k\subset\mathbb C$ a subfield.
Let $f:U\to S$, $f':U'\to S$ morphisms, with $U,U',S\in\Var(k)$ irreducible, $U'$ smooth.
Let $\bar{S}\in\PVar(k)$ a compactification of $S$.
Let $\bar{X},\bar{X'}\in\PVar(k)$ compactification of $U$ and $U'$ respectively,
such that $f$ (resp. $f'$) extend to a morphism $\bar f:\bar X\to\bar S$, resp. $\bar{f'}:\bar{X'}\to\bar S$.
Denote $\bar D=\bar X\backslash U$ and $\bar D'=\bar{X'}\backslash U'$ and 
$\bar E=(\bar D\times_{\bar S}\bar X')\cup(\bar X\times_{\bar S}\bar D')$.
Denote $i:\bar D\hookrightarrow\bar X$, $i':\bar D\hookrightarrow\bar X$ denote the closed embeddings
and $j:U\hookrightarrow\bar X$, $j':U'\hookrightarrow\bar X'$ the open embeddings.
Denote $d=\dim(U)$ and $d'=\dim(U')$.
We have the following commutative diagram in $D(\mathbb Z)$
\begin{equation*}
\xymatrix{RHom_{\DA(\bar S)}^{\bullet}(M(U'/\bar S),M((\bar X,\bar D)/\bar S))
\ar[d]^{RI(-,-)}\ar[rr]^{{\mathcal F_S^{Hdg}}^{(-,-)}} & \, & 
RHom_{DMHM(\bar S)}^{\bullet}(f'_{!Hdg}\mathbb Z_{U'}^{Hdg},f_{*Hdg}\mathbb Z_U^{Hdg})
\ar[d]^{RI(-,-)} \\
RHom^{\bullet}(M(\pt),M(\bar X'\times_{\bar S}\bar X,\bar E)(d')[2d'])
\ar[d]^{l}\ar[rr]^{{\mathcal F^{\pt}_{Hdg}}^{(-,-)}} & \, & 
RHom^{\bullet}(\mathbb Z_{\pt}^{Hdg},a_{U'\times_SU!}\mathbb Z_{U\times_SU'}^{Hdg}(d')[2d'])
\ar[d]^{l} \\
\mathcal Z_d(\bar X'\times_{\bar S}\bar X,E,\bullet)
\ar[rr]^{\mathcal R^d_{\bar X'\times_{\bar S}\bar X}} & \, & 
C^{\mathcal D}_{2d+\bullet}(\bar X'\times_{\bar S}\bar X,E,Z(d))}
\end{equation*}
where 
\begin{equation*}
M((\bar X,\bar D)/\bar S):=\Cone(\ad(i_*,i^!):M(\bar D/\bar S)\to M(\bar X/\bar S))
=\bar f_*j_*E_{et}(\mathbb Z(U/U))\in\DA(\bar S)
\end{equation*}
and $l$ the isomorphisms given by canonical embedding of complexes.
\end{cor}

\begin{proof}
The upper square of this diagram follows from theorem \ref{main}(ii).
On the other side, the lower square follows from the absolute case.
\end{proof}

\subsection{The p adic Hodge realization functor for relative motives over a subfield $k\subset\mathbb C_p$}

Let $p$ a prime number. Let $k\subset\mathbb C_p$ a subfield. 

For $S\in\Var(k)$, we have the analytical functor 
\begin{equation*}
\an_S^{*mod}:C_{\mathcal D}(S)\to C_{D(O_{S_{\mathbb C_p}^{an}})}(S_{\mathbb C_p}^{an,pet}), \; 
M\mapsto M^{an}:=\an_S^{*mod}M
\end{equation*} 
given by the morphism of ringed topological spaces 
$\an_S:S_{\mathbb C_p}^{an}\xrightarrow{\an_S}S_{\mathbb C_p}\xrightarrow{\pi_{k/\mathbb C_p}(S)} S$.
For $S\in\Var(k)$, we denote by for short 
$O\mathbb B_{dr,S}:=O\mathbb B_{dr,S_{\mathbb C_p}^{an}}:=O\mathbb B_{dr,R_{\mathbb C_p}(S_{\mathbb C_p}^{an})}$.
where $R_{\mathbb C_p}:\AnSp(\mathbb C_p)\to\AdSp/(\mathbb C_p,O_{\mathbb C_p})$
is the canonical functor (see section 2).
For $\mathcal S\in\Cat$ a site, $p$ a prime number, we recall (see section 2) the functor 
\begin{equation*}
(-)\otimes\mathbb Z_p:C(\mathcal S)\to 
C_{\mathbb Z_p}(\mathcal S)\subset C(\mathbb N\times\mathcal S)=\PSh(\mathcal S,\Fun(\mathbb N,C(\mathbb Z))), \;
K\mapsto K\otimes\mathbb Z_p:=(K\otimes\mathbb Z/p^n\mathbb Z)_{n\in\mathbb N}.
\end{equation*}
Let $S\in\Var(k)$. Let $S=\cup_{i\in I}S_i$ an open cover such that there
exists closed embeddings $i_i:S_i\hookrightarrow\tilde S_i$ with $\tilde S_I\in\SmVar(k)$.
We have the category $D_{\mathcal D(1,0)fil,rh}(S/(\tilde S_I))\times_I D_{\mathbb Z_pfil,c,k}(S^{et})$  
\begin{itemize}
\item whose set of objects is the set of triples $\left\{(((M_I,F,W),u_{IJ}),(K,W),\alpha)\right\}$ with
\begin{eqnarray*} 
((M_I,F,W),u_{IJ})\in D_{\mathcal D(1,0)fil,rh}(S/(\tilde S_I)), \, (K,W)\in D_{\mathbb Z_pfil,c,k}(S^{et}), \\ 
\alpha:\mathbb B_{dr,(\tilde S_I)}(K,W)\to 
F^0DR(S)^{[-]}(((M_I,F,W),u_{IJ})^{an}\otimes_{O_S}((O\mathbb B_{dr,\tilde S_I},F),t_{IJ}))
\end{eqnarray*}
where $\alpha$ is a morphism in $D_{\mathbb B_{dr},G,fil}(S_{\mathbb C_p}^{an,pet}/(\tilde S_{I,\mathbb C_p}^{an,pet}))$,
\item and whose set of morphisms consists of 
\begin{equation*}
\phi=(\phi_D,\phi_C,[\theta]):(((M_{1I},F,W),u_{IJ}),(K_1,W),\alpha_1)\to(((M_{2I},F,W),u_{IJ}),(K_2,W),\alpha_2)
\end{equation*}
where $\phi_D:((M_1,F,W),u_{IJ})\to((M_2,F,W),u_{IJ})$ and $\phi_C:(K_1,W)\to (K_2,W)$ 
are morphisms and
\begin{eqnarray*}
\theta=(\theta^{\bullet},I(F^0DR(S)(\phi^{an}_D)\times I)\circ I(\alpha_1),
I(\alpha_2)\circ I(\mathbb B_{dr,(\tilde S_I)}(\phi_C))): \\
I(\mathbb B_{dr,(\tilde S_I)}(K_1,W))[1]\to 
I(F^0DR(S)(((M_{2I},F,W),u_{IJ})^{an}\otimes_{O_S}((O\mathbb B_{dr,\tilde S_I},F),t_{IJ}))  
\end{eqnarray*}
is an homotopy,  
$I:D_{\mathbb B_{dr},G,fil}(S_{\mathbb C_p}^{an,pet}/(\tilde S^{an,pet}_{I,\mathbb C_p}))\to 
K_{\mathbb B_{dr},G,fil}(S_{\mathbb C_p}^{an,pet}/(\tilde S^{an,pet}_{I,\mathbb C_p}))$
being the injective resolution functor, and for
\begin{itemize}
\item $\phi=(\phi_D,\phi_C,[\theta]):(((M_{1I},F,W),u_{IJ}),(K_1,W),\alpha_1)\to(((M_{2I},F,W),u_{IJ}),(K_2,W),\alpha_2)$
\item $\phi'=(\phi'_D,\phi'_C,[\theta']):(((M_{2I},F,W),u_{IJ}),(K_2,W),\alpha_2)\to(((M_{3I},F,W),u_{IJ}),(K_3,W),\alpha_3)$
\end{itemize}
the composition law is given by 
\begin{eqnarray*}
\phi'\circ\phi:=(\phi'_D\circ\phi_D,\phi'_C\circ\phi_C,
I(DR(S)(\phi^{'an}_D\otimes I))\circ[\theta]+[\theta']\circ I(\mathbb B_{dr,(\tilde S_I)}(\phi_C))[1]): \\
(((M_{1I},F,W),u_{IJ}),(K_1,W),\alpha_1)\to(((M_{3I},F,W),u_{IJ}),(K_3,W),\alpha_3),
\end{eqnarray*}
in particular for 
$(((M_I,F,W),u_{IJ}),(K,W),\alpha)\in C_{\mathcal D(1,0)fil,rh}(S/(\tilde S_I))\times_I D_{\mathbb Z_pfil,c,k}(S^{et})$,
\begin{equation*}
I_{(((M_I,F,W),u_{IJ}),(K,W),\alpha)}=((I_{M_I}),I_K,0),
\end{equation*}
\end{itemize}
and also the category 
$D_{\mathcal D(1,0)fil,rh,\infty}(S/(\tilde S_I))\times_I D_{\mathbb Z_pfil,c,k}(S^{et})$ defined in the same way,
together with the localization functor
\begin{eqnarray*}
(D(zar),I):C_{\mathcal D(1,0)fil,rh}(S/(\tilde S_I))\times_I D_{\mathbb Z_pfil,c,k}(S^{et})
\to D_{\mathcal D(1,0)fil,rh}(S/(\tilde S_I))\times_I D_{\mathbb Z_pfil,c,k}(S^{et}) \\
\to D_{\mathcal D(1,0)fil,rh,\infty}(S/(\tilde S_I))\times_I D_{\mathbb Z_pfil,c,k}(S^{et}).
\end{eqnarray*}
Note that if $\phi=(\phi_D,\phi_C,[\theta]):(((M_1,F,W),u_{IJ}),(K_1,W),\alpha_1)\to(((M_2,F,W),u_{IJ}),(K_2,W),\alpha_2)$
is a morphism in $D_{\mathcal D(1,0)fil,rh}(S/(\tilde S_I))\times_I D_{\mathbb Z_pfil,c,k}(S^{et})$
such that $\phi_D$ and $\phi_C$ are isomorphism then $\phi$ is an isomorphism (see remark \ref{CGremp}).
Moreover,
\begin{itemize}
\item For 
$(((M_{I},F,W),u_{IJ}),(K,W),\alpha)\in D_{\mathcal D(1,0)fil,rh}(S/(\tilde S_I))\times_I D_{\mathbb Z_pfil,c,k}(S^{et})$, 
we set
\begin{equation*}
(((M_{I},F,W),u_{IJ}),(K,W),\alpha)[1]:=(((M_{I},F,W),u_{IJ})[1],(K,W)[1],\alpha[1]).
\end{equation*}
\item For 
\begin{equation*}
\phi=(\phi_D,\phi_C,[\theta]):(((M_{1I},F,W),u_{IJ}),(K_1,W),\alpha_1)\to(((M_{2I},F,W),u_{IJ}),(K_2,W),\alpha_2)
\end{equation*}
a morphism in $D_{\mathcal D(1,0)fil,rh}(S/(\tilde S_I))\times_I D_{\mathbb Z_pfil,c,k}(S^{et})$, 
we set (see \cite{CG} definition 3.12)
\begin{eqnarray*}
\Cone(\phi):=(\Cone(\phi_D),\Cone(\phi_C),((\alpha_1,\theta),(\alpha_2,0)))
\in D_{\mathcal D(1,0)fil,rh}(S/(\tilde S_I))\times_I D_{\mathbb Z_pfil,c,k}(S^{et}),
\end{eqnarray*}
$((\alpha_1,\theta),(\alpha_2,0))$ being the matrix given by the composition law, together with the canonical maps
\begin{itemize}
\item $c_1(-)=(c_1(\phi_D),c_1(\phi_C),0):(((M_{2I},F,W),u_{IJ}),(K_2,W),\alpha_2)\to\Cone(\phi)$
\item $c_2(-)=(c_2(\phi_D),c_2(\phi_C),0):\Cone(\phi)\to (((M_{1I},F,W),u_{IJ}),(K_1,W),\alpha_1)[1]$.
\end{itemize}
\end{itemize}

Let $S\in\Var(k)$. Let $S=\cup_{i=1}^sS_i$ an open cover such that there exists
closed embedding $i_i:S\hookrightarrow\tilde S_i$ with $\tilde S_i\in\SmVar(k)$.
Consider the category 
\begin{equation*}
(D_{\mathcal D(1,0)fil}(\tilde S_I)\times_ID_{\mathbb Z_pfil,c,k}(\tilde S_I^{et}))\in\Fun(\Gamma(\tilde S_I),\TriCat)
\end{equation*}
such that
\begin{equation*}
(D_{\mathcal D(1,0)fil}(\tilde S_I)\times_ID_{\mathbb Z_pfil,c,k}(\tilde S_I^{et}))(\tilde S_I)
=D_{\mathcal D(1,0)fil}(\tilde S_I)\times_ID_{\mathbb Z_pfil,c,k}(\tilde S_I^{et})
\end{equation*}
\begin{itemize}
\item whose objects are 
$(((M_I,F,W),(K_I,W),\alpha_I),u_{IJ})\in (D_{\mathcal D(1,0)fil}(\tilde S_I)\times_ID_{\mathbb Z_pfil,c,k}(\tilde S_I^{et}))$ 
such that 
\begin{equation*}
((M_I,F,W),(K_I,W),\alpha_I)\in D_{\mathcal D(1,0)fil}(\tilde S_I)\times_ID_{\mathbb Z_pfil,c,k}(\tilde S_I^{et})
=:\mathcal D_p(\tilde S_I) 
\end{equation*}
and for $I\subset J$, 
\begin{eqnarray*}
u_{IJ}:((M_I,F,W),(K_I,W),\alpha_I)\to \\
p_{IJ*}((M_J,F,W),(K_J,W),\alpha_J):=(p_{IJ*}(M_J,F,W),p_{IJ*}(K_J,W),p_{IJ*}\alpha_J)
\end{eqnarray*}
are morphisms in $\mathcal D_p(\tilde S_I)$,
\item whose morphisms $m=(m_I):(((M_I,F,W),(K_I,W),\alpha_I),u_{IJ})\to (((M'_I,F,W),(K'_I,W),\alpha'_I),v_{IJ})$
is a family of morphism such that $v_{IJ}\circ m_I=p_{IJ*}m_J\circ u_{IJ}$ in $\mathcal D_p(\tilde S_I)$
\end{itemize}
We have then the identity functor
\begin{eqnarray*}
I_S:D_{\mathcal D(1,0)fil}(S/(\tilde S_I))\times_ID_{\mathbb Z_pfil,c,k}(S^{et})\to
(D_{\mathcal D(1,0)fil}(\tilde S_I)\times_ID_{\mathbb Z_pfil,c,k}(\tilde S_I^{et})), \\
(((M_I,F,W),u_{IJ}),(K,W),\alpha)\mapsto (((M_I,F,W),i_{I*}j_I^*(K,W),j_I^*\alpha),(u_{IJ},I,0)), \\ 
m=(m_I,n)\mapsto m=(m_I,i_*j_I^*n)
\end{eqnarray*}
which is a full embedding since by definition, for $((M_I,F,W),u_{IJ})\in D_{\mathcal D(1,0)fil}(S/(\tilde S_I))$,
\begin{equation*}
u_{IJ}:(M_I,F,W)\to p_{IJ*}(M_J,F,W) 
\end{equation*}
are filtered Zariski local equivalence, i.e. isomorphisms in $D_{\mathcal D(1,0)fil}(\tilde S_I)$, and hence for 
$((M_I,F,W),u_{IJ}),(K,W),\alpha)\in D_{\mathcal D(1,0)fil}(S/(\tilde S_I))\times_ID_{\mathbb Z_pfil,c,k}(S^{et})$,
\begin{eqnarray*}
(u_{IJ},I,0):((M_I,F,W),i_{I*}j_I^*(K,W),j_I^*\alpha)\to \\
p_{IJ*}((M_J,F,W),i_{J*}j_J^*(K,W),j_J^*\alpha)=(p_{IJ*}(M_J,F,W),i_{I*}j_I^*(K,W),j_I^*\alpha)
\end{eqnarray*}
are isomorphisms in $\mathcal D_p(\tilde S_I)$.

\begin{defi}\label{IUSmp}
For $h:U\to S$ a smooth morphism with $S,U\in\SmVar(k)$ 
and $h:U\xrightarrow{n}X\xrightarrow{f}S$ a compactification of $h$ with $n$ an open embedding, $X\in\SmVar(k)$
such that $D:=X\backslash U=\cup_{i=1}^sD_i\subset X$ is a normal crossing divisor, 
we denote by, using definition \ref{DHdgalpha} and definition \ref{wtildew}
\begin{eqnarray*}
I_p(U/S):h_{!Hdg}h^{!Hdg}\mathbb Z_{p,S}^{Hdg}\xrightarrow{:=} \\
(p_{S*}E_{zar}(\Omega^{\bullet}_{X\times S/S}\otimes_{O_{X\times S}}(n\times I)_{!Hdg}\Gamma_U^{\vee,Hdg}(O_{U\times S},F_b)),
\mathbb D_Sh_*E_{et}\mathbb Z_{p,U^{et}},h_!\alpha(U,\delta)) \\
\xrightarrow{((DR(X\times S/S)(\ad((n\times I)_{!Hdg},(n\times I)^*)(-)),0),I,0)} \\
(\Cone((\Omega_{/S}^{\Gamma,pr}(i_{D_i}\times I))_{i\in[1,\ldots,s]}:
p_{S*}E_{zar}(\Omega^{\bullet}_{X\times S/S}\otimes_{O_{X\times S}}\Gamma_X^{\vee,Hdg}(O_{X\times S},F_b))\to \\
(\cdots\to(p_{S*}E_{zar}(\Omega^{\bullet}_{D_I\times S/S}\otimes_{O_{D_I\times S}}\Gamma_{D_I}^{\vee,Hdg}(O_{D_I\times S},F_b)))
\to\cdots)),\mathbb D_Sh_*E_{et}\mathbb Z_{p,U^{et}},h_!\alpha(U,\delta)) \\
\xrightarrow{=:}
(\mathcal F_S^{FDR}(\mathbb Z(U/S)),Rh_!\mathbb Z_{p,U^{et}},\alpha(\mathbb Z(U/S)))
\end{eqnarray*}
the canonical isomorphism in $D_{\mathcal Dfil}(S)\times_ID_{\mathbb Z_p,c,k}(S^{et})$, where 
\begin{itemize}
\item we recall that (see section 6.2)
\begin{equation*}
h^{!Hdg}\mathbb Z_{p,S}^{Hdg}=
(\Gamma_U^{\vee,Hdg}(O_{U\times S},F_b),\mathbb Z_{p,U^{et}},\alpha(U))\in HM_{gm,k,\mathbb C_p}(U), 
\end{equation*}
\item $i_{D_i}:D_i\hookrightarrow X$ are the closed embeddings,
\item $\alpha(\mathbb Z(U/S)):=h_!\alpha(U,\delta):=T^w(h,\otimes)(-)\circ h_!\alpha(U,\delta)$ 
(see definition \ref{falphap}), with 
\begin{eqnarray*}
\alpha(U,\delta):=(DR(U)(\Omega_{(U\times U/U)/(U/pt)}(\Gamma_U^{\vee,Hdg}(O_{U\times U},F_b)))\otimes I)^{-1}\circ\alpha(U):\\
\mathbb B_{dr,U}\to
DR(U)((p_{U*}E_{zar}(\Omega^{\bullet}_{U\times U/U}\otimes_{O_{U\times U}}\Gamma_U^{\vee,Hdg}(O_{U\times U})))^{an}
\otimes_{O_U}(O\mathbb B_{dr,U},F)),
\end{eqnarray*}
by the way we note that the following diagram in $C(U_{\mathbb C_p}^{an,pet})$ commutes
\begin{equation*}
\begin{tikzcd}
\mathbb B_{dr,U}\ar[r,"\alpha(U)"] & 
F^0((\Omega^{\bullet}_{U_{\mathbb C_p}^{an}},F_b)\otimes_{O_U}(O\mathbb B_{dr,U},F))=:F^0DR(U)(O\mathbb B_{dr,U},F) \\
\, & \shortstack{$
E_{et}F^0DR(U)((p_{U*}E_{zar}(\Omega^{\bullet}_{U\times U/U}\otimes_{O_{U\times U}}\Gamma_U^{\vee,Hdg}(O_{U\times U},F_b)))^{an}
\otimes_{O_U}(O\mathbb B_{dr,U},F))$}
\ar[u,"DR(U)(\Omega_{(U\times U/U)/(U/pt)}(\Gamma_U^{\vee{,}Hdg}(O_{U\times U}{,}F_b))\otimes I)"] \\
p_{U*}E_{et}\Gamma_U^{\vee}\mathbb B_{dr,U\times U}
\ar[uu]^{\ad(\delta^{*mod}_U,\delta_{U*})(-)}\ar[r,"\alpha(U\times U)"] &
\shortstack{$
p_{U*}E_{et}F^0DR(U\times U)((\Omega^{\bullet}_{U\times U/U}\otimes_{O_{U\times U}}\Gamma_U^{\vee,Hdg}(O_{U\times U},F_b))^{an}
\otimes_{O_{U\times U}}(O\mathbb B_{dr,U\times U},F))$}
\ar[u,"T^w(p_U{,}\otimes)(-)^{-1}"]
\end{tikzcd}
\end{equation*}
\end{itemize}
\end{defi}

\begin{lem}\label{thetamp}
Let $S\in\SmVar(k)$. Let $g:U'/S\to U/S$ o morphism with $U/S:=(U,h),U'/S:=(U',h)\in\Var(k)^{sm}/S$.
Let $h:U\xrightarrow{n}X\xrightarrow{f}S$ a compactification of $h$ with $n$ an open embedding, $X\in\SmVar(k)$
such that $D:=X\backslash U=\cup_{i=1}^sD_i\subset X$ is a normal crossing divisor, 
Let $h':U\xrightarrow{n'}X'\xrightarrow{f'}S$ a compactification of $h'$ with $n'$ an open embedding, $X'\in\SmVar(k)$
such that $D':=X\backslash U=\cup_{i=1}^sD_i\subset X$ is a normal crossing divisor and such that
$g:U'\to U$ extend to $\bar g:X'\to X$, see definition-proposition \ref{RCHdef0}.
Then, using definition \ref{IUSmp}, 
the following diagram in $D_{\mathcal Dfil}(S)\times_ID_{c,k}(S^{et})$ commutes
\begin{eqnarray*}
\xymatrix{
h'_{!Hdg}h^{'!Hdg}\mathbb Z_{p,S}^{Hdg}\ar[rrr]^{I_p(U'/S)}\ar[d]_{\ad(g_{!Hdg},g^{!Hdg})(h^{!Hdg}\mathbb Z_{p,S}^{Hdg})} 
& \, & \, & (\mathcal F_S^{FDR}(\mathbb Z(U'/S)),e(S^{et})_*C_*(\mathbb Z(U'/S)\otimes\mathbb Z_p),\alpha(\mathbb Z(U'/S)))
\ar[d]^{(\Omega_{/S}^{\Gamma,pr}(R_S^{CH}(g)),Re(S^{et})_*\mathbb Z(g),\theta(g))} \\ 
h_{!Hdg}h^{!Hdg}\mathbb Z_{p,S}^{Hdg}\ar[rrr]^{I_p(U/S)} & \, & \, &
(\mathcal F_S^{FDR}(\mathbb Z(U/S)),e(S^{et})_*C_*(\mathbb Z(U/S)\otimes\mathbb Z_p),\alpha(\mathbb Z(U/S)))} 
\end{eqnarray*}
where 
\begin{equation*}
\theta(g):=R_{\mathcal D}([\Gamma_g]):I(\mathbb B_{dr,S}(h'_!\mathbb Z_{p,{U'}^{et}}))[1]
\to I(F^0DR(S)(\mathcal F_S^{FDR}(\mathbb Z(U/S))^{an}\otimes_{O_S}(O\mathbb B_{dr,S},F)))
\end{equation*}
is the homotopy given by the third term of the syntomic homology class of the graph $\Gamma_g\subset U'\times_S U$,
(see definition \ref{Delkpdef} and we recall (see section 6.2) that 
$I:C_{B_{dr},G}(S_{\mathbb C_p}^{an}/(\tilde S_{I,\mathbb C_p}^{an}))\to 
K_{B_{dr},G}(S_{\mathbb C_p}^{an}/(\tilde S_{I,\mathbb C_p}^{an}))$ 
is the injective resolution functor.
\end{lem}

\begin{proof}
Immediate from definition.
\end{proof}

We can now define the p adic Hodge realization functor for motives :

\begin{defi}\label{HodgeRealDAsingp}
Let $k\subset\mathbb C_p$ a subfield. Let $S\in\Var(k)$. Let $S=\cup_{i=1}^sS_i$ an open cover such that there exists
closed embedding $i_i:S\hookrightarrow\tilde S_i$ with $\tilde S_i\in\SmVar(k)$.
We define the Hodge realization functor as, using definition \ref{DRalgdefFunct},
\begin{eqnarray*}
\mathcal F_S^{Hdg}:=(\mathcal F_S^{FDR},e(S^{et})_*C_*L\otimes\mathbb Z_p):
C(\Var(k)^{sm}/S)\to D_{\mathcal D(1,0)fil}(S/(\tilde S_I))\times_I D_{\mathbb Z_pfil,c,k}(S^{et}), \\ 
F\mapsto\mathcal F_S^{Hdg}(F):=(\mathcal F_S^{FDR}(F,W),e(S^{et})_*C_*(L(F,W)\otimes\mathbb Z_p),\alpha(F)),
\end{eqnarray*}
first on objects and then on morphisms :
\begin{itemize}
\item for $F\in C(\Var(k)^{sm}/S)$, taking $(F,W)\in C_{fil}(\Var(k)^{sm}/S)$
such that $D(\mathbb A^1,et)(F,W)$ gives the weight structure on $D(\mathbb A^1,et)(F)$,
\begin{eqnarray*}
\mathcal F_S^{Hdg}(F):=(\mathcal F_S^{FDR}(F,W),e(S^{et})_*C_*(L(F,W)\otimes\mathbb Z_p),\alpha(F)):= \\
(e(S)_*\mathcal Hom((\hat R_{\tilde S_I}^{CH}(\rho_{\tilde S_I}^*Li_{I*}j_I^*(F,W)),
\hat R^{CH}(T^q(D_{IJ})(-))),(E_{zar}(\Omega^{\bullet,\Gamma,pr}_{/\tilde S_I},F_{DR}),T_{IJ})), \\
e(S^{et})_*C_*(L(F,W)\otimes\mathbb Z_p),\alpha(F)) 
\in D_{\mathcal D(1,0)fil}(S/(\tilde S_I))\times_I D_{\mathbb Z_pfil,c,k}(S^{et})
\end{eqnarray*}
where $\alpha(F)$ is the map in $D_{\mathbb B_{dr},G,fil}(S_{\mathbb C_p}^{an,pet}/(\tilde S_{I,\mathbb C_p}^{an,pet}))$
writing for short $DR(S):=DR(S)^{[-]}:=(DR(\tilde S_I)[-d_{\tilde S_I}])$
\begin{eqnarray*}
\alpha(F):\mathbb B_{dr,(\tilde S_I)}(Re(S^{et})_*((M,W)\otimes^L\mathbb Z_p)):=
\mathbb B_{dr,(\tilde S_I)}((i_{I*}j_I^*e(S^{et})_*C_*(L(F,W)\otimes\mathbb Z_p)),I) \\
\xrightarrow{=}
T_S^{-1}\mathbb B_{dr,(\tilde S_I)}((e(\tilde S_{I,\bar k}^{et})_*C_*(Li_{I*}j_I^*(F,W)\otimes\mathbb Z_p)),I) \\
\xrightarrow{=}
T_S^{-1}((((\cdot\to\oplus_{(U_{I\alpha},h_{I\alpha})\in V_I}
\mathbb B_{dr,\tilde S_I}(h_{I\alpha !}h^!_{I\alpha}\mathbb Z_{p,\tilde S_{I,\bar k}^{et}}) 
\xrightarrow{\mathbb B_{dr,\tilde S_I}(\ad(g^{\bullet !}_{I,\alpha,\beta},g^{\bullet}_{I,\alpha,\beta !})(-))} \\ 
\oplus_{(U_{I\alpha},h_{I\alpha})\in V_I}
\mathbb B_{dr,\tilde S_I}(h_{I\alpha !}h^!_{I\alpha}\mathbb Z_{p,\tilde S_{I,\bar k}^{et}})\to\cdot),u_{IJ})),W) \\
\xrightarrow{T_S^{-1}(\alpha(\mathbb Z(U_{I\alpha}/\tilde S_I)),\theta(g^{\bullet}_{I,\alpha,\beta}))} \\
T_S^{-1}F^0DR(S)((e(S)_*\mathcal Hom((\hat R_{\tilde S_I}^{CH}(\rho_{\tilde S_I}^*Li_{I*}j_I^*(F,W)),
\hat R^{CH}(T^q(D_{IJ})(-))), \\ (E_{zar}(\Omega^{\bullet,\Gamma,pr}_{/\tilde S_I},F_{DR}),T_{IJ})))^{an}
\otimes_{O_S}((O\mathbb B_{dr,\tilde S_I},F),t_{IJ})) \\
\xrightarrow{=:}
F^0DR(S)((\mathcal F^{FDR}_S(M,W))^{an}\otimes_{O_S}((O\mathbb B_{dr,\tilde S_I},F),t_{IJ}))
\end{eqnarray*}
where 
\begin{itemize}
\item $T_S$ is the canonical functor
\begin{equation*}
T_S:D_{\mathbb B_{dr},G,fil}(S_{\mathbb C_p}^{an,pet}/(\tilde S_{I,\mathbb C_p}^{an,pet}))
\to D_{\mathbb B_{dr},G,fil}((\tilde S_{I,\mathbb C_p}^{an,pet})), \; T_S((K_I,W),u_{IJ})=((K_I,W),u_{IJ}) 
\end{equation*}
which is a full embedding on the full subcategory 
\begin{equation*}
D_{\mathbb B_{dr},G,fil}(S_{\mathbb C_p}^{an,pet}/(\tilde S_{I,\mathbb C_p}^{an,pet}))^{\sim}:=
D_{pet}(C_{\mathbb B_{dr},G,fil}(S_{\mathbb C_p}^{an,pet}/(\tilde S_{I,\mathbb C_p}^{an,pet}))^{sim}) \\
\subset D_{\mathbb B_{dr},G,fil}(S_{\mathbb C_p}^{an,pet}/(\tilde S_{I,\mathbb C_p}^{an,pet}))
\end{equation*}
where $C_{\mathbb B_{dr},G,fil}(S_{\mathbb C_p}^{an,pet}/(\tilde S_{I,\mathbb C_p}^{an,pet}))^{\sim}
\subset C_{\mathbb B_{dr},G,fil}(S_{\mathbb C_p}^{an,pet}/(\tilde S_{I,\mathbb C_p}^{an,pet}))$
is the full subcategory consisting of $((K_I,W),u_{IJ})$ such that $u_{IJ}:K_I\to p_{IJ*}K_J$ are isomorphisms,
\item we denote by $g^n_{I,\alpha,\beta}:U_{I\alpha}\to U_{I\beta}$ which satisfy 
$h_{I\beta}\circ g^n_{I,\alpha,\beta}=h_{I\alpha}$ the morphisms in the canonical projective resolution
\begin{eqnarray*}
q:Li_{I*}j_I^*(F,W):=((\cdots\to\oplus_{(U_{I\alpha},h_{I\alpha})\in\Var(k)^{sm}/\tilde S_I}\mathbb Z(U_{I\alpha}/\tilde S_I)
\xrightarrow{(\mathbb Z(g^{\bullet}_{I,\alpha,\beta}))} \\
\oplus_{(U_{I\alpha},h_{I\alpha})\in\Var(k)^{sm}/\tilde S_I}\mathbb Z(U_{I\alpha}/\tilde S_I)\to\cdots),W)\to i_{I*}j_I^*(F,W),
\end{eqnarray*}
\end{itemize}
using lemma \ref{thetamp}, 
\begin{equation*}
(\alpha(\mathbb Z(U_{I\alpha}/S)),\theta(g^{\bullet}_{I,\alpha,\beta}))
\end{equation*}
being the matrix given inductively by the composition law in 
$D_{\mathcal D(1,0)fil}(\tilde S_I)\times_I D_{fil,c,k}(\tilde S_I^{et})$,
and the fact 
\begin{equation*}
F^0DR(S)(\mathcal F^{FDR}_S(M,W)^{an}\otimes_{O_S}((O\mathbb B_{dr,\tilde S_I},F),t_{IJ}))
\in D_{\mathbb B_{dr},G,fil}(S_{\mathbb C_p}^{an,pet}/(\tilde S_{I,\mathbb C_p}^{an,pet}))^{\sim}
\end{equation*}
since $\mathcal F^{FDR}_S(M,W)\in D(DRM(S))$ by corollary \ref{FDRMHM},
that is we have the following isomorphism in 
$(D_{\mathcal D(1,0)fil}(\tilde S_I)\times_I D_{fil,c,k}(\tilde S_I^{et}))$,
denoting for short $V_I:=\Var(k)^{sm}/\tilde S_I$
\begin{eqnarray*}
(I^{\bullet}(U_{I\alpha}/\tilde S_I)): 
(((\cdot\to\oplus_{(U_{I\alpha},h_{I\alpha})\in V_I}h_{I\alpha !Hdg}h^{!Hdg}_{I\alpha}\mathbb Z^{Hdg}_{p,\tilde S_I}
\xrightarrow{\ad(g^{\bullet,!Hdg}_{I,\alpha,\beta},g^{\bullet}_{I,\alpha,\beta !Hdg})(-)} \\
\oplus_{(U_{I\alpha},h_{I\alpha})\in V_I}h_{I\alpha !Hdg}h^{!Hdg}_{I\alpha}\mathbb Z^{Hdg}_{p,\tilde S_I}\to\cdot),u_{IJ}),W) \\
\xrightarrow{\sim}
I_S(\mathcal F_S^{Hdg}(F):=(\mathcal F_S^{FDR}(F,W),e(\tilde S_I^{et})_*C_*(Li_{I*}j_I^*(F,W)\otimes\mathbb Z_p),
\alpha(F)))
\end{eqnarray*}
\item for $m:F_1\to F_2$ a morphism in $C(\Var(k)^{sm}/S)$, taking $(F_1,W),(F_2,W)\in C_{fil}(\Var(k)^{sm}/S)$
such that $D(\mathbb A^1,et)(F_2,W)$ gives the weight structure on $D(\mathbb A^1,et)(F_2)$ 
$D(\mathbb A^1,et)(F_1,W)$ gives the weight structure on $D(\mathbb A^1,et)(F_1)$
and such that $m:(F_1,W)\to (F_2,W)$ is a filtered morphism, the morphism 
$\mathcal F_S^{Hdg}(m)$ in $D_{\mathcal D(1,0)fil}(S/(\tilde S_I))\times_I D_{fil,c,k}(S^{et})$ is given by
\begin{eqnarray*}
\mathcal F_S^{Hdg}(m):&=&{I_S^{-,-}}^{-1}((I^{\bullet}(U_{I\alpha}/(\tilde S_I)))\circ
(\ad(l_{I\alpha,\beta}^{\bullet !Hdg},l^{\bullet}_{I\alpha,\beta !Hdg})(\mathbb Z_{p,U_{I\alpha}}^{Hdg})) 
\circ (I^{\bullet}(U_{I\alpha}/(\tilde S_I)))^{-1}) \\
&=&(\mathcal F_S^{FDR}(m),Re(S^{et})_*\mathbb Z(m),\theta(m):=(\theta(l_{I\alpha,\beta}))):
\mathcal F_S^{Hdg}(F_1)\to\mathcal F_S^{Hdg}(F_2)
\end{eqnarray*}
using lemma \ref{thetam}, that is we have the following commutative diagram in 
$(D_{\mathcal D(1,0)fil}(\tilde S_I)\times_I D_{fil,c,k}(\tilde S_I^{et}))$,
denoting for short $V_I:=\Var(k)^{sm}/\tilde S_I$,
\begin{eqnarray*}
\begin{tikzcd}
(((\cdot\to\oplus_{(U_{I\alpha},h_{I\alpha})\in V_I}h_{I\alpha !Hdg}h^{!Hdg}_{I\alpha}\mathbb Z^{Hdg}_{p,\tilde S_I}
\xrightarrow{A^{Hdg}_{g_{1I,\alpha,\beta}^{\bullet}}} 
\oplus_{(U_{I\alpha},h_{I\alpha})\in V_I}h_{I\alpha !Hdg}h^{!Hdg}_{I\alpha}\mathbb Z^{Hdg}_{p,\tilde S_I}\to\cdot),u_{IJ}),W)
\ar[r,"(I^{\bullet}(U_{I\alpha}/\tilde S_I))"]
\ar[d,"\ad(l_{I{,}\alpha{,}\beta}^{\bullet !Hdg}{,}l^{I{,}\bullet}_{\alpha{,}\beta !Hdg})(-)"'] & 
\mathcal F_S^{Hdg}(F_1)
\ar[d,"\mathcal F_S^{Hdg}(m)=(\mathcal F_S^{FDR}(m){,}Re(S^{et})_*\mathbb Z(m){,}(\theta(l_{I\alpha{,}\beta})))"'] \\
(((\cdot\to\oplus_{(U_{I\alpha},h_{I\alpha})\in V_I}h_{I\alpha !Hdg}h^{!Hdg}_{I\alpha}\mathbb Z^{Hdg}_{p,\tilde S_I}
\xrightarrow{A^{Hdg}_{g_{2I,\alpha,\beta}^{\bullet}}}
\oplus_{(U_{I\alpha},h_{I\alpha})\in V_I}h_{I\alpha !Hdg}h^{!Hdg}_{I\alpha}\mathbb Z^{Hdg}_{p,\tilde S_I}\to\cdot),u_{IJ}),W)
\ar[r,"(I^{\bullet}(U_{\alpha}/\tilde S_I))"] & \mathcal F_S^{Hdg}(F_2)
\end{tikzcd}
\end{eqnarray*}
where
\begin{itemize}
\item we denoted for short
$A^{Hdg}_{g_{1I,\alpha,\beta}^{\bullet}}:=
\ad(g^{\bullet,!Hdg}_{1I,\alpha,\beta},g^{\bullet}_{1I,\alpha,\beta !Hdg})(h^{!Hdg}_{I\alpha}\mathbb Z^{Hdg}_{\tilde S_I})$
\item we denoted for short
$A^{Hdg}_{g_{2I,\alpha,\beta}^{\bullet}}:=
\ad(g^{\bullet,!Hdg}_{2I,\alpha,\beta},g^{\bullet}_{2I,\alpha,\beta !Hdg})(h^{!Hdg}_{I\alpha}\mathbb Z^{Hdg}_{\tilde S_I})$
\item we denote by $g^n_{1I,\alpha,\beta}:U_{I\alpha}\to U_{I\beta}$, which satisfy 
$h_{I\beta}\circ g^n_{1I,\alpha,\beta}=h_{I\alpha}$, the morphisms in the canonical projective resolution
\begin{eqnarray*}
q:Li_{I*}j_I^*(F_1,W):=((\cdots\to\oplus_{(U_{I\alpha},h_{I\alpha})\in\Var(k)^{sm}/\tilde S_I}\mathbb Z(U_{I\alpha}/\tilde S_I)
\xrightarrow{(\mathbb Z(g^{\bullet}_{1I,\alpha,\beta}))} \\
\oplus_{(U_{I\alpha},h_{I\alpha})\in\Var(k)^{sm}/\tilde S_I}\mathbb Z(U_{I\alpha}/\tilde S_I)\to\cdots),W)\to i_{I*}j_I^*(F_1,W)
\end{eqnarray*}
\item we denote by $g^n_{2I,\alpha,\beta}:U_{I\alpha}\to U_{I\beta}$, which satisfy 
$h_{I\beta}\circ g^n_{2I,\alpha,\beta}=h_{\alpha}$, the morphisms in the canonical projective resolution
\begin{eqnarray*}
q:Li_{I*}j_I^*(F_2,W):=((\cdots\to\oplus_{(U_{I\alpha},h_{I\alpha})\in\Var(k)^{sm}/\tilde S_I}\mathbb Z(U_{I\alpha}/\tilde S_I)
\xrightarrow{(\mathbb Z(g^{\bullet}_{2I,\alpha,\beta}))} \\
\oplus_{(U_{I\alpha},h_{I\alpha})\in\Var(k)^{sm}/\tilde S_I}\mathbb Z(U_{I\alpha}/\tilde S_I)\to\cdots),W)\to i_{I*}j_I^*(F_2,W)
\end{eqnarray*}
\item we denote by $l_{I\alpha,\beta}^n:U_{I\alpha}\to U_{I\beta}$ which satisfy 
$h_{I\beta}\circ l^n_{I\alpha,\beta}=h_{I\alpha}$ and 
$l^{n+1}_{I\alpha,\beta}\circ g^n_{1I\alpha,\beta}=g^n_{2I\alpha,\beta}\circ l^n_{I\alpha,\beta}$ 
the morphisms in the morphism of canonical projective resolutions
\begin{eqnarray*}
Li_{I*}j_I^*(m):Li_{I*}j_I^*(F_1,W):=
((\cdots\to\oplus_{(U_{I\alpha},h_{I\alpha})\in\Var(k)^{sm}/\tilde S_I}\mathbb Z(U_{I\alpha}/\tilde S_I)\to\cdots),W) 
\xrightarrow{(\mathbb Z(l^{\bullet}_{I\alpha,\beta}))} \\
((\cdots\to\oplus_{(U_{I\alpha},h_{I\alpha})\in\Var(k)^{sm}/\tilde S_I}\mathbb Z(U_{\alpha}/\tilde S_I)\to\cdots),W)
=:Li_{I*}j_I^*(F_2,W),
\end{eqnarray*}
\item the maps $I^{\bullet}(U_{I\alpha})$ are given by definition \ref{IUSm} and lemma \ref{thetam}.
\end{itemize}
\end{itemize}
Obviously $\mathcal F_S^{Hdg}(F[1])=\mathcal F_S^{Hdg}(F)[1]$ and 
$\mathcal F_S^{Hdg}(\Cone(m))=\Cone(\mathcal F_S^{Hdg}(m))$. 
This functor induces by proposition \ref{projwach} and remark \ref{CGremp} the functor
\begin{eqnarray*}
\mathcal F_S^{Hdg}:=(\mathcal F_S^{FDR},e(S^{et})_*C_*(L\otimes\mathbb Z_p):
\DA(S)\to D_{\mathcal D(1,0)fil}(S/(\tilde S_I))\times_I D_{\mathbb Z_pfil,c,k}(S^{et}), \\ 
M=D(\mathbb A^1,et)(F)\mapsto
\mathcal F_S^{Hdg}(M):=\mathcal F_S^{Hdg}(F)=(\mathcal F_S^{FDR}(M),Re(S^{et})_*(M\otimes^L\mathbb Z_p),\alpha(M)),
\end{eqnarray*}
with $\alpha(M)=\alpha(F)$. 
\end{defi}

We now give the functoriality with respect to the five operation using the De Rahm realization case and the etale realization case :

\begin{prop}\label{TetTFDR}
Let $p$ a prime number. Consider an embedding $k\subset\mathbb C_p$. 
\begin{itemize}
\item[(i)]Let $g:T\to S$ a morphism with $T,S\in\Var(k)$. 
Assume there exists a factorization $g:T\xrightarrow{l}Y\times S\xrightarrow{p}S$, with $Y\in\SmVar(k)$,
$l$ a closed embedding and $p$ the projection.
Let $S=\cup_{i\in I}S_i$ an open cover and 
$i_i:S_i\hookrightarrow\tilde S_i$ closed embeddings with $\tilde S_i\in\SmVar(k)$.
Then, $\tilde g_I:Y\times\tilde S_I\to\tilde S_I$ is a lift of $g_I=g_{|T_I}:T_I\to S_I$
and we have closed embeddings $i'_I:=i_I\circ l\circ j'_I:T_I\hookrightarrow Y\times\tilde S_I$.
Then, for $M=D(\mathbb A^1,et)(F)\in DA_c(S)$, the following diagram commutes :
\begin{equation*}
\begin{tikzcd}
\mathbb B_{dr,(Y\times\tilde S_I)}(g^{*w}Re(S^{et})_*(M\otimes^L\mathbb Z_p))
\ar[r,"g^*(\alpha(M))"]\ar[d,"\mathbb B_{dr{,}(Y\times\tilde S_I)}(T^*(g{,}e)(M\otimes\mathbb Z_p))"] & 
F^0DR(T)^{[-]}((g^{\hat{*}mod}_{Hdg}\mathcal F_S^{FDR}(M))^{an}\otimes_{O_T}((O\mathbb B_{dr,Y\times\tilde S_I},F),t_{IJ}))
\ar[d,"DR(T)^{[-]}((T(g{,}\mathcal F^{FDR})(M))^{an}\otimes I)"] \\
\mathbb B_{dr,(Y\times\tilde S_I)}(Re(T^{et})_*g^*(M\otimes^L\mathbb Z_p))\ar[r,"\alpha(g^*M)"] &     
F^0DR(T)^{[-]}((\mathcal F_T^{FDR}(g^*M))^{an}\otimes_{O_T}((O\mathbb B_{dr,Y\times\tilde S_I},F),t_{IJ})),
\end{tikzcd}
\end{equation*}
see definition \ref{TgDRdefsing} and definition \ref{falphap}.
\item[(ii)]Let $f:T\to S$ a morphism with $T,S\in\Var(k)$. 
Assume there exists a factorization $f:T\xrightarrow{l}Y\times S\xrightarrow{p}S$, with $Y\in\SmVar(k)$,
$l$ a closed embedding and $p$ the projection.
Let $S=\cup_{i\in I}S_i$ an open cover and 
$i_i:S_i\hookrightarrow\tilde S_i$ closed embeddings with $\tilde S_i\in\SmVar(k)$.
Then, for $M=D(\mathbb A^1,et)(F)\in DA_c(T)$,the following diagram commutes :
\begin{equation*}
\begin{tikzcd}
\mathbb B_{dr,(\tilde S_I)}(Rf_{*w}Re(T^{et})_*(M\otimes^L\mathbb Z_p))\ar[r,"f_*(\alpha(M))"] &
F^0DR(S)^{[-]}((Rf^{Hdg}_*\mathcal F_T^{FDR}(M))^{an}\otimes_{O_S}((O\mathbb B_{dr,\tilde S_I},F),t_{IJ})) \\
\mathbb B_{dr,(\tilde S_I)}(Re(S^{et})_*Rf_*(M\otimes^L\mathbb Z_p))
\ar[u,"\mathbb B_{dr{,}(\tilde S_I)}(T_*(f{,}e)(M\otimes\mathbb Z_p))"]\ar[r,"\alpha(Rf_*M)"] &
F^0DR(S)^{[-]}((\mathcal F_S^{FDR}(Rf_*M))^{an}\otimes_{O_S}((O\mathbb B_{dr,\tilde S_I},F),t_{IJ}))
\ar[u,"DR(S)^{[-]}((T_*(f{,}\mathcal F^{FDR})(M))^{an}\otimes I)"]
\end{tikzcd}
\end{equation*}
see definition \ref{SixTalg} and definition \ref{falphap}.
\item[(iii)]Let $f:T\to S$ a morphism with $T,S\in\Var(k)$. 
Assume there exists a factorization $f:T\xrightarrow{l}Y\times S\xrightarrow{p}S$, with $Y\in\SmVar(k)$,
$l$ a closed embedding and $p$ the projection.
Let $S=\cup_{i\in I}S_i$ an open cover and 
$i_i:S_i\hookrightarrow\tilde S_i$ closed embeddings with $\tilde S_i\in\SmVar(k)$.
Then, for $M=D(\mathbb A^1,et)(F)\in DA_c(T)$,the following diagram commutes :
\begin{equation*}
\begin{tikzcd}
\mathbb B_{dr,(\tilde S_I)}(Rf_{!w}Re(T^{et})_*(M\otimes^L\mathbb Z_p))
\ar[d,"\mathbb B_{dr{,}(\tilde S_I)}(T_!(f{,}e)(M\otimes\mathbb Z_p))"]\ar[r,"f_!(\alpha(M))"] &   
F^0DR(S)^{[-]}((Rf^{Hdg}_!\mathcal F^T_{FDR}(M))^{an}\otimes_{O_S}((O\mathbb B_{dr,\tilde S_I},F),t_{IJ}))
\ar[d,"DR(S)^{[-]}((T_!(f{,}\mathcal F_{FDR})(M))^{an}\otimes I)"] \\
\mathbb B_{dr,(\tilde S_I)}(Re(S^{et})_*Rf_!(M\otimes^L\mathbb Z_p))\ar[r,"\alpha(Rf_!M)"] &  
F^0DR(S)^{[-]}((\mathcal F^S_{FDR}(Rf_!M))^{an}\otimes_{O_S}((O\mathbb B_{dr,\tilde S_I},F),t_{IJ}))
\end{tikzcd}
\end{equation*}
see definition \ref{SixTalg} and definition \ref{falphap}.
\item[(iv)]Let $f:T\to S$ a morphism with $T,S\in\Var(k)$. 
Assume there exists a factorization $f:T\xrightarrow{l}Y\times S\xrightarrow{p}S$, with $Y\in\SmVar(k)$,
$l$ a closed embedding and $p$ the projection.
Let $S=\cup_{i\in I}S_i$ an open cover and 
$i_i:S_i\hookrightarrow\tilde S_i$ closed embeddings with $\tilde S_i\in\SmVar(k)$.
Then, for $M=D(\mathbb A^1,et)(F)\in DA_c(S)$,the following diagram commutes :
\begin{equation*}
\begin{tikzcd}
\mathbb B_{dr,(Y\times\tilde S_I)}(f^{!w}Re(S^{et})_*(M\otimes^L\mathbb Z_p))\ar[r,"f^!(\alpha(M))"] &     
F^0DR(T)^{[-]}((f^{*mod}_{Hdg}\mathcal F_S^{FDR}(M))^{an}\otimes_{O_T}((O\mathbb B_{dr,Y\times\tilde S_I},F),t_{IJ})) \\
\mathbb B_{dr,(Y\times\tilde S_I)}(Re(T^{et})_*f^!(M\otimes^L\mathbb Z_p))
\ar[r,"\alpha(f^!M)"]\ar[u,"\mathbb B_{dr{,}(Y\times\tilde S_I)}(T^!(f{,}e)(M\otimes\mathbb Z_p))"] &   
F^0DR(T)^{[-]}((\mathcal F_T^{FDR}(f^!M))^{an}\otimes_{O_T}((O\mathbb B_{dr,Y\times\tilde S_I},F),t_{IJ}))
\ar[u,"DR^{[-]}(T)((T^!(g{,}\mathcal F^{FDR})(M))^{an}\otimes I)"]
\end{tikzcd}
\end{equation*}
see definition \ref{SixTalg} and definition \ref{falphap}.
\item[(v)] Let $S\in\Var(k)$. 
Let $S=\cup_{i\in I}S_i$ an open cover and 
$i_i:S_i\hookrightarrow\tilde S_i$ closed embeddings with $\tilde S_i\in\SmVar(k)$.
Then, for $M,N\in DA_c(S)$,the following diagram commutes :
\begin{equation*}
\begin{tikzcd}
\shortstack{$\mathbb B_{dr,(\tilde S_I)}(Re(S^{et})_*(M\otimes^L\mathbb Z_p))\otimes_{\mathbb B_{dr,S}}$ \\
$\mathbb B_{dr,(\tilde S_I)}(Re(S^{et})_*(N\otimes^L\mathbb Z_p))$}
\ar[rr,"\alpha(M)\otimes\alpha(N)"]
\ar[d,"T(\otimes{,}\mathbb B_{dr})(Re(S^{et})_*M\otimes\mathbb Z_p{,}Re(S^{et})_*N\otimes\mathbb Z_p)"] & \, &   
F^0DR(S)((\mathcal F_S^{FDR}(M)\otimes^{Hdg}_{O_S}\mathcal F_S^{FDR}(N))^{an}\otimes_{O_S}(((O\mathbb B_{dr,\tilde S_I},F),t_{IJ}))
\ar[d,"DR(S)((T(\otimes{,}\mathcal F^{FDR})(M{,}N))^{an}\otimes I)"] \\
\mathbb B_{dr,(\tilde S_I)}(Re(S^{et})_*((M\otimes N)\otimes^L\mathbb Z_p))\ar[rr,"(\alpha(M\otimes N))"] & \, & 
F^0DR(S)((\mathcal F^S_{FDR}(M\otimes N))^{an}\otimes_{O_S}((O\mathbb B_{dr,\tilde S_I},F),t_{IJ}))
\end{tikzcd}
\end{equation*}
see  definition \ref{SixTalg} and definition \ref{falphap}.
\end{itemize}
\end{prop}

\begin{proof}

\noindent(i): Follows from the following commutative diagram in 
$(D_{\mathcal D(1,0)fil}(Y\times\tilde S_I)\times_I D_{fil,c,k}(Y\times\tilde S_I^{et}))$,
\begin{eqnarray*}
\begin{tikzcd}
\shortstack{
$(((\to\oplus_{(U_{I\alpha},h_{I\alpha})\in V_I}\tilde g_I^{*Hdg}h_{I\alpha !Hdg}h^{!Hdg}_{I\alpha}\mathbb Z^{Hdg}_{p,\tilde S_I}
\xrightarrow{A^{Hdg}_{g_{I,\alpha,\beta}^{\bullet}}}$ \\ 
$\oplus_{(U_{I\alpha},h_{I\alpha})\in V_I}h_{I\alpha !Hdg}h^{!Hdg}_{I\alpha}\mathbb Z^{Hdg}_{p,\tilde S_I}\to),u_{IJ}),W)$}
\ar[rr,"(\tilde g_I^{*Hdg}I^{\bullet}(U_{I\alpha}/\tilde S_I))"]\ar[d,"T^{Hdg}(\tilde g_I{,}h_I)(-)"'] & \, &
\shortstack{$(g^{\hat*mod}_{Hdg}\mathcal F_T^{FDR}(F)$, \\ $g^{*w}e(S^{et})_*C_*L(F,W),g^*(\alpha(F)))$}
\ar[d,"(T(g{,}\mathcal F^{FDR})(M){,}T(g{,}e)(M){,}0)"'] \\
\shortstack{
$(((\to\oplus_{(U'_{I\alpha},h_{I\alpha})\in W_I}h'_{I\alpha !Hdg}h^{'!Hdg}_{I\alpha}\mathbb Z^{Hdg}_{p,Y\times\tilde S_I}
\xrightarrow{A^{Hdg}_{g_{I,\alpha,\beta}^{'\bullet}}}$ \\
$\oplus_{(U'_{I\alpha},h'_{I\alpha})\in W_I}h'_{I\alpha !Hdg}h^{'!Hdg}_{I\alpha}\mathbb Z^{Hdg}_{p,Y\times\tilde S_I}\to),u_{IJ}),W)$}
\ar[rr,"(I^{\bullet}(U'_{\alpha}/Y\times\tilde S_I))"] & \, &
\shortstack{$(\mathcal F_T^{FDR}(g^*F)$, \\ $e(T^{et})_*C_*L(g^*F,W),\alpha(g^*F))$}
\end{tikzcd}
\end{eqnarray*}
where, we have denoted for short $V_I:=\Var(k)^{sm}/\tilde S_I$ and $W_I:=\Var(k)^{sm}/Y\times\tilde S_I$,
\begin{itemize}
\item we denoted for short
$A^{Hdg}_{g_{I,\alpha,\beta}^{\bullet}}:=
\ad(g^{\bullet,!Hdg}_{I,\alpha,\beta},g^{\bullet}_{I,\alpha,\beta !Hdg})(h^{!Hdg}_{I\alpha}\mathbb Z^{Hdg}_{p,\tilde S_I})$
\item we denoted for short
$A^{Hdg}_{g_{I,\alpha,\beta}^{'\bullet}}:=
\ad(g^{'\bullet,!Hdg}_{I,\alpha,\beta},g^{'\bullet}_{I,\alpha,\beta !Hdg})(h^{'!Hdg}_{I\alpha}\mathbb Z^{Hdg}_{p,Y\times\tilde S_I})$
\item we denote by $g^n_{I,\alpha,\beta}:U_{I\alpha}\to U_{I\beta}$, which satisfy 
$h_{I\beta}\circ g^n_{I,\alpha,\beta}=h_{I\alpha}$, the morphisms in the canonical projective resolution
\begin{eqnarray*}
q:Li_{I*}j_I^*(F,W):=(\cdots\to\oplus_{(U_{I\alpha},h_{I\alpha})\in\Var(k)^{sm}/\tilde S_I}\mathbb Z(U_{I\alpha}/\tilde S_I)
\xrightarrow{(\mathbb Z(g^{\bullet}_{I,\alpha,\beta}))} \\
\oplus_{(U_{I\alpha},h_{I\alpha})\in\Var(k)^{sm}/\tilde S_I}\mathbb Z(U_{I\alpha}/\tilde S_I)\to\cdots)\to i_{I*}j_I^*(F,W)
\end{eqnarray*}
\item we denote by $g^{'n}_{I,\alpha,\beta}:U'_{I\alpha}\to U'_{I\beta}$, which satisfy 
$h'_{I\beta}\circ g^{'n}_{I,\alpha,\beta}=h'_{\alpha}$, the morphisms in the canonical projective resolution
\begin{eqnarray*}
q:Li'_{I*}j_I^{'*}(g^*F,W):=(\cdots\to\oplus_{(U'_{I\alpha},h'_{I\alpha})\in\Var(k)^{sm}/Y\times\tilde S_I}
\mathbb Z(U'_{I\alpha}/Y\times\tilde S_I)
\xrightarrow{(\mathbb Z(g^{'\bullet}_{I,\alpha,\beta}))} \\
\oplus_{(U'_{I\alpha},h'_{I\alpha})\in\Var(k)^{sm}/Y\times\tilde S_I}
\mathbb Z(U'_{I\alpha}/Y\times\tilde S_I)\to\cdots)\to i'_{I*}j_I^{'*}(g^*F,W)
\end{eqnarray*}
\end{itemize}

\noindent(ii): Follows from (i) by adjonction.

\noindent(iii): The closed embedding case is given by (ii) and the smooth projection case follows from (i) by adjonction.

\noindent(iv): Follows from (iii) by adjonction.

\noindent(v):Obvious
\end{proof}

\begin{prop}\label{keyHdgp}
Let $p$ a prime number. Consider an embedding $k\subset\mathbb C_p$. 
\begin{itemize}
\item[(i)] Let $S\in\Var(k)$. Let $S=\cup_i S_i$ an open cover such that there exist closed embeddings
$i_i:S_i\hookrightarrow\tilde S_i$ with $\tilde S_i\in\SmVar(k)$. Then
we have the isomorphism in $D_{\mathcal D(1,0)fil}(S/(\tilde S_I))\times_ID_{fil,c,k}(S^{et})$
\begin{eqnarray*}
\mathcal F_S^{Hdg}(\mathbb Z_S)\xrightarrow{:=}
(\mathcal F_S^{FDR}(\mathbb Z_S),e(S^{et})_*(\mathbb Z_S\otimes\mathbb Z_p),\alpha(\mathbb Z_S)) \\
\xrightarrow{((\Omega^{\Gamma,pr}_{/\tilde S_I}(\hat R^{CH}(\ad(i_I^*,i_{I*})(\Gamma_{S_I}^{\vee,w}\mathbb Z_{\tilde S_I})))),I,0)} \\
I_S^{-1}((e(S)_*\mathcal Hom((\hat R^{CH}(\Gamma^{\vee,w}_{S_I}\mathbb Z_{\tilde S_I}),\hat R^{CH}(x_{IJ})),
(E_{zar}(\Omega^{\bullet,\Gamma,pr}_{/\tilde S_I},F_{DR}),T_{IJ})),T(S/(\tilde S_I))(\mathbb Z_{p,S^{et}}^w), \\
\alpha(\tilde S_I,\delta))) 
\xrightarrow{=}
\iota_S((\Gamma^{\vee,Hdg}_{S_I}(O_{\tilde S_I},F_b),x_{IJ}),\mathbb Z_{p,S^{et}}^w,\alpha(S))=:\iota_S(\mathbb Z_S^{Hdg})
\end{eqnarray*}
with
\begin{eqnarray*}
\alpha(S):\mathbb B_{dr,(\tilde S_I)}(\mathbb Z_{S^{et}}\otimes\mathbb Z_p)
:=\mathbb B_{dr,(\tilde S_I)}(\Gamma_{S_I}^{\vee,w}\mathbb Z_{p,\tilde S_I^{et}},x_{IJ}) \\
\xrightarrow{(\Gamma_{S_I}^{\vee}\alpha(\tilde S_I))}
F^0DR(S)((\Gamma^{\vee,Hdg}_{S_I}(O_{\tilde S_I},F_b),x_{IJ})\otimes_{O_S}((O\mathbb B_{dr,\tilde S_I},F),t_{IJ}))
\end{eqnarray*}
\item[(ii)]Let $f:X\to S$ a morphism with $X,S\in\Var(k)$, $X$ quasi-projective.
Consider a factorization $f:X\xrightarrow{l}Y\times S\xrightarrow{p_S}S$ with $Y=\mathbb P^{N,o}\subset\mathbb P^N$ an open subset,
$l$ a closed embedding and $p_S$ the projection. Let $S=\cup_i S_i$ an open cover such that there exist closed embeddings
$i_i:S_i\hookrightarrow\tilde S_i$ with $\tilde S_i\in\SmVar(k)$. 
Recall that $S_I:=\cap_{i\in I} S_i$, $X_I=f^{-1}(S_I)$, and $\tilde S_I:=\Pi_{i\in I}\tilde S_i$. Then,
using proposition \ref{TetTFDR}(iii), the map of definition \ref{SixTalg} gives
an isomorphism in $D_{\mathcal D(1,0)fil}(S/(\tilde S_I))\times_ID_{\mathbb Z_pfil,c,k}(S^{et})$
\begin{eqnarray*}
(T_!(f,\mathcal F^{FDR})(\mathbb Z_X),T_!(e,f)(\mathbb Z_{p,X^{et}}))): \\
\mathcal F_S^{Hdg}(M^{BM}(X/S)):=
(\mathcal F_S^{FDR}(Rf_!\mathbb Z_X),e(S^{et})_*Rf_!(\mathbb Z_X\otimes\mathbb Z_p),\alpha(Rf_!\mathbb Z_X)) \\
\xrightarrow{\sim} 
(Rf_{Hdg!}(\Gamma^{\vee,Hdg}_{X_I}(O_{Y\times\tilde S_I},F_b),x_{IJ}(X/S)),Rf_{!w}\mathbb Z_{p,X^{et}}^w,
f_!(\alpha(X)))=:\iota_S(Rf_{!Hdg}(\mathbb Z_{p,X}^{Hdg})).
\end{eqnarray*}
with
\begin{equation*}
\mathbb Z_{p,X}^{Hdg}:=((\Gamma^{\vee,Hdg}_{X_I}(O_{Y\times\tilde S_I},F_b),x_{IJ}(X/S)),
\mathbb Z_{p,X^{et}}^w,\alpha(X))\in C(MHM_{gm,k,\mathbb C_p}(X))
\end{equation*}
\end{itemize}
\end{prop}

\begin{proof}
Follows from proposition \ref{TetTFDR}(iii) and theorem \ref{mainthm}.
\end{proof}

The main theorem of this section is the following :

\begin{thm}\label{mainpadic}
Let $p$ a prime number. Let $k\subset\mathbb C_p$ a subfield. 
\begin{itemize}
\item[(i)] For $S\in\Var(k)$, we have $\mathcal F_S^{Hdg}(\DA_c(S))\subset D(MHM_{gm,k,\mathbb C_p}(S))$,
\begin{equation*}
\iota_S:D(MHM_{gm,k,\mathbb C_p}(S))\hookrightarrow 
D_{\mathcal D(1,0)fil}(S/(\tilde S_I))\times_ID_{\mathbb Z_pfil,c,k}(S^{et})
\end{equation*}
being a full embedding by theorem \ref{Bekp}.
\item[(ii)] The Hodge realization functor $\mathcal F_{Hdg}(-)$ define a morphism of 2-functor on $\Var(k)$
\begin{equation*}
\mathcal F^{Hdg}_{-}:\Var(k)\to(\DA_c(-)\to D(MHM_{gm,k,\mathbb C_p}(-)))
\end{equation*}
whose restriction to $\QPVar(\mathbb C)$ is an homotopic 2-functor in sense of Ayoub. More precisely,
\begin{itemize}
\item[(ii0)] for $g:T\to S$ a morphism, with $T,S\in\QPVar(k)$, and $M\in\DA_c(S)$, the
the map of definition \ref{TgDRdefsing} induces an isomorphism in $D(MHM_{gm,k,\mathbb C_p}(T))$
\begin{eqnarray*}
T(g,\mathcal F^{Hdg})(M):=(T(g,\mathcal F^{FDR})(M),T(g,e)(M\otimes^L\mathbb Z_p),0): \\
g^{\hat*Hdg}\mathcal F_S^{Hdg}(M):=
\iota_T^{-1}(g^{\hat{*}mod}_{Hdg}\mathcal F_S^{FDR}(M),g^*Re(S^{et})_*(M\otimes^L\mathbb Z_p),g^*(\alpha(M))) \\
\xrightarrow{\sim}\iota_T^{-1}(\mathcal F_T^{FDR}(g^*M),Re(T^{et})_*g^*(M\otimes^L\mathbb Z_p),\alpha(g^*M))
=:\mathcal F_T^{Hdg}(g^*M),
\end{eqnarray*} 
\item[(ii1)] for $f:T\to S$ a morphism, with $T,S\in\QPVar(k)$, and $M\in\DA_c(T)$,  
the map of definition \ref{SixTalg} induces an isomorphism in $D(MHM_{gm,k,\mathbb C_p}(S))$
\begin{eqnarray*}
T_*(f,\mathcal F^{Hdg})(M):=(T_*(f,\mathcal F^{FDR})(M),I,0): \\
Rf_{Hdg*}\mathcal F_T^{Hdg}(M):=
\iota_S^{-1}(Rf^{Hdg}_*\mathcal F_T^{FDR}(M),Rf_*Re(T^{et})_*(M\otimes^L\mathbb Z_p),f_*(\alpha(M))) \\ 
\xrightarrow{\sim}\iota_S^{-1}(\mathcal F_S^{FDR}(Rf_*M),Re(S^{et})_*Rf_*(M\otimes^L\mathbb Z_p),\alpha(Rf_*M))
=:\mathcal F_S^{Hdg}(Rf_*M),
\end{eqnarray*}  
\item[(ii2)] for $f:T\to S$ a morphism, with $T,S\in\QPVar(k)$, and $M\in\DA_c(T)$, 
the map of definition \ref{SixTalg} induces an isomorphism in $D(MHM_{gm,k,\mathbb C_p}(S))$
\begin{eqnarray*}
T_!(f,\mathcal F^{Hdg})(M):=(T_!(f,\mathcal F^{FDR})(M),T_!(f,e)(M\otimes^L\mathbb Z_p),0): \\
Rf_{!Hdg}\mathcal F_T^{Hdg}(M):=
\iota_S^{-1}(Rf^{Hdg}_!\mathcal F_T^{FDR}(M),Rf_!Re(T^{et})_*(M\otimes^L\mathbb Z_p),f_!(\alpha(M))) \\ 
\xrightarrow{\sim}\iota_S^{-1}(\mathcal F_S^{FDR}(Rf_!M),Re(S^{et})_*Rf_!M\otimes\mathbb Z_p,\alpha(f_!M))
=:\mathcal F_T^{Hdg}(f_!M),
\end{eqnarray*} 
\item[(ii3)] for $f:T\to S$ a morphism, with $T,S\in\QPVar(k)$, and $M\in\DA_c(S)$,
the map of definition \ref{SixTalg} induces an isomorphism in $D(MHM_{gm,k,\mathbb C_p}(T))$
\begin{eqnarray*}
T^!(f,\mathcal F^{Hdg})(M):=(T^!(f,\mathcal F^{FDR})(M),T^!(f,e)(M\otimes^L\mathbb Z_p),0): \\
f^{*Hdg}\mathcal F_S^{Hdg}(M):=
\iota_T^{-1}(f^{*mod}_{Hdg}\mathcal F_S^{FDR}(M),f^!Re(S^{et})_*(M\otimes^L\mathbb Z_p),f^!(\alpha(M))) \\ 
\xrightarrow{\sim}\iota_T^{-1}(\mathcal F_T^{FDR}(f^!M),Re(T^{et})_*f^!(M\otimes^L\mathbb Z_p),\alpha(f^!M))
=:\mathcal F_T^{Hdg}(f^!M),
\end{eqnarray*}
\item[(ii4)] for $S\in\Var(k)$, and $M,N\in\DA_c(S)$,
the map of definition \ref{SixTalg} induces an isomorphism in $D(MHM_{gm,k,\mathbb C_p}(S))$
\begin{eqnarray*}
T(\otimes,\mathcal F^{Hdg})(M,N):=(T(\otimes,\mathcal F_S^{Hdg})(M,N),I,0): 
\mathcal F_S^{Hdg}(M)\otimes^{Hdg}_{O_S}\mathcal F_S^{Hdg}(N):= \\
\iota_S^{-1}(\mathcal F_S^{FDR}(M)\otimes^{Hdg}_{O_S}\mathcal F_S^{FDR}(N),
Re(S^{et})_*(M\otimes^L\mathbb Z_p)\otimes Re(S^{et})_*(N\otimes^L\mathbb Z_p),\alpha(M)\otimes\alpha(N)) \\
\xrightarrow{\sim}\iota_S^{-1}\mathcal F_S^{Hdg}(M\otimes N):=
\iota_S^{-1}(\mathcal F_S^{FDR}(M\otimes N),Re(S^{et})_*((M\otimes N)\otimes^L\mathbb Z_p),\alpha(M\otimes N)).  
\end{eqnarray*}
\end{itemize}
\item[(iii)] For $S\in\Var(k)$, the following diagram commutes :
\begin{equation*}  
\xymatrix{\Var(k)/S\ar[rrr]^{MH(/S)}\ar[d]_{M(/S)} & \, & \, & D(MHM_{gm,k,\mathbb C_p}(S))\ar[d]^{\iota^S} \\
\DA(S)\ar[rrr]^{\mathcal F_S^{Hdg}} & \, & \, & 
D_{\mathcal D(1,0)fil}(S/(\tilde S_I))\times_I D_{\mathbb Z_pfil,c,k}(S^{et})}
\end{equation*}
\end{itemize}
\end{thm}

\begin{proof} 

\noindent(i):Let $M\in\DA_c(S)$. There exist by definition of constructible motives an isomorphism in $\DA(S)$ 
\begin{equation*}
w(M):M\xrightarrow{\sim}\Cone(M(X_0/S)[d_0]\xrightarrow{m_1}\cdots\xrightarrow{m_m} M(X_m/S)[d_m]),
\end{equation*}
with $f_n:X_n\to S$ morphisms and $X_n\in\QPVar(k)$.
This gives the isomorphism in $D_{\mathcal D(1,0)fil}(S/(\tilde S_I))\times_I D_{\mathbb Z_pfil,c,k}(S^{et})$
\begin{equation*}
\mathcal F_S^{Hdg}(w(M)):\mathcal F_S^{Hdg}(M)\xrightarrow{\sim}
\Cone(\mathcal F_S^{Hdg}(M(X_0/S))[d_0]\xrightarrow{\mathcal F_S^{Hdg}(m_1)}\cdots
\xrightarrow{\mathcal F_S^{Hdg}(m_1)}\mathcal F_S^{Hdg}(M(X_m/S))[d_m]),
\end{equation*}
On the other hand, by proposition \ref{keyHdgp}(i), we have 
\begin{equation*}
\mathcal F_S^{Hdg}(M(X_n/S))\xrightarrow{\sim}Rf_{!Hdg}\mathbb Z_{p,X}^{Hdg}\in D(MHM_{gm,k,\mathbb C_p}(S)).
\end{equation*}
This prove (i).

\noindent(ii0): Follows from theorem \ref{mainthm}(i) and proposition \ref{TetTFDR}(i).

\noindent(ii1): Follows from theorem \ref{mainthm}(iii) and proposition \ref{TetTFDR}(ii).

\noindent(ii2):Follows from theorem \ref{mainthm}(ii) and proposition \ref{TetTFDR}(iii).

\noindent(ii3): Follows from theorem \ref{mainthm}(iv) and proposition \ref{TetTFDR}(iv).

\noindent(ii4):Follows from theorem \ref{mainthm}(v) and proposition \ref{TetTFDR}(v).

\noindent(iii): By (ii), for $g:X'/S\to X/S$ a morphism, with $X',X,S\in\Var(k)$ 
and $X/S=(X,f)$, $X'/S=(X',f')$, we have by adjonction the following commutative diagram
\begin{equation*}
\xymatrix{\mathcal F_S^{Hdg}(M(X'/S)=f'_!f^{'!}\mathbb Z_S=f_!g_!g^!f^!\mathbb Z_S)
\ar[d]_{T_!(f',\mathcal F^{Hdg})(f^{'!}M(X'/S))\circ T^!(f',\mathcal F^{Hdg})(M(X'/S))}
\ar[rr]^{\mathcal F_S^{Hdg}(M(/S)(g)=f_!\ad(g_!,g^!)(f^!\mathbb Z_S))} & \, & 
\mathcal F_S^{Hdg}(M(X/S)=f_!f^!\mathbb Z_S)
\ar[d]_{T_!(f,\mathcal F^{Hdg})(f^!M(X/S))\circ T^!(f,\mathcal F^{Hdg})(M(X/S))} \\
MH(X'/S):=Rf'_{!Hdg}f^{'!Hdg}\mathbb Z^{Hdg}_{p,S}=f_{!Hdg}g_{!Hdg}g^{!Hdg}f^{!Hdg}\mathbb Z^{Hdg}_{p,S}
\ar[rr]^{f_{!Hdg}\ad(g_{!Hdg},g^{!Hdg})(f^{!Hdg}\mathbb Z^{Hdg}_{p,S})} & \, &  
MH(X/S):=f_{!Hdg}f^{!Hdg}\mathbb Z^{Hdg}_{p,S}}.
\end{equation*}
where the left and right columns are isomorphisms by (ii). This proves (iii).
\end{proof}

The theorem \ref{mainpadic} gives immediately the following :

\begin{cor}
Let $p$ a prime number. Let $k\subset\mathbb C_p$ a subfield.
Let $f:U\to S$, $f':U'\to S$ morphisms, with $U,U',S\in\Var(k)$ irreducible, $U'$ smooth.
Let $\bar{S}\in\PVar(k)$ a compactification of $S$.
Let $\bar{X},\bar{X'}\in\PVar(k)$ compactification of $U$ and $U'$ respectively,
such that $f$ (resp. $f'$) extend to a morphism $\bar f:\bar X\to\bar S$, resp. $\bar{f'}:\bar{X'}\to\bar S$.
Denote $\bar D=\bar X\backslash U$ and $\bar D'=\bar{X'}\backslash U'$ and 
$\bar E=(\bar D\times_{\bar S}\bar X')\cup(\bar X\times_{\bar S}\bar D')$.
Denote $i:\bar D\hookrightarrow\bar X$, $i':\bar D\hookrightarrow\bar X$ denote the closed embeddings
and $j:U\hookrightarrow\bar X$, $j':U'\hookrightarrow\bar X'$ the open embeddings.
Denote $d=\dim(U)$ and $d'=\dim(U')$.
We have the following commutative diagram in $D(\mathbb Z)$
\begin{equation*}
\xymatrix{RHom_{\DA(\bar S)}^{\bullet}(M(U'/\bar S),M((\bar X,\bar D)/\bar S))
\ar[d]^{RI(-,-)}\ar[rr]^{{\mathcal F_S^{Hdg}}^{(-,-)}} & \, & 
RHom_{D(MHM_{gm,k,\mathbb C_p}(\bar S))}^{\bullet}(f'_{!Hdg}\mathbb Z_{p,U'}^{Hdg},f_{*Hdg}\mathbb Z_{p,U}^{Hdg})
\ar[d]^{RI(-,-)} \\
RHom^{\bullet}(M(\pt),M(\bar X'\times_{\bar S}\bar X,\bar E)(d')[2d'])
\ar[d]^{l}\ar[rr]^{{\mathcal F^{\pt}_{Hdg}}^{(-,-)}} & \, & 
RHom^{\bullet}(\mathbb Z_{\pt}^{Hdg},a_{U'\times_SU!}\mathbb Z_{p,U\times_SU'}^{Hdg}(d')[2d'])
\ar[d]^{l} \\
\mathcal Z_d(\bar X'\times_{\bar S}\bar X,E,\bullet)
\ar[rr]^{\mathcal R^d_{\bar X'\times_{\bar S}\bar X}} & \, & 
C^{syn}_{2d+\bullet}(\bar X'\times_{\bar S}\bar X,E,Z(d))}
\end{equation*}
where 
\begin{equation*}
M((\bar X,\bar D)/\bar S):=\Cone(\ad(i_*,i^!):M(\bar D/\bar S)\to M(\bar X/\bar S))
=\bar f_*j_*E_{et}(\mathbb Z(U/U))\in\DA(\bar S)
\end{equation*}
and $l$ the isomorphisms given by canonical embedding of complexes.
\end{cor}

\begin{proof}
The upper square of this diagram follows from theorem \ref{mainpadic}(ii).
On the other side, the lower square follows from the absolute case.
\end{proof}


LAGA UMR CNRS 7539 \\
Universit\'e Paris 13, Sorbonne Paris Cit\'e, 99 av Jean-Baptiste Clement, \\
93430 Villetaneuse, France, \\
bouali@math.univ-paris13.fr


\begin{thebibliography}{1}

\bibitem{AyoubB} J.Ayoub, \emph{Note sur les op\'erations de Grothendieck et la realisation de Betti},
Journal of the Institute of Mathematics of Jussieu, Volume 9, Issue 02, April 2010, pp.225-263.

\bibitem{AyoubMHA1} J.Ayoub, \emph{L'algebre de Hopf et le groupe de Galois motiviques d'un corps de caract\'eristique nulle I},
Journal die reine (Crelles Journal), Volume 2014, Issue 693, pp.1-149.

\bibitem{AyoubMHA2} J.Ayoub, \emph{L'algebre de Hopf et le groupe de Galois motiviques d'un corps de caract\'eristique nulle II},
Journal die reine (Crelles Journal), Volume 2014, Issue 693, pp.151-126.

\bibitem{AyoubT} J.Ayoub, \emph{Les six op\'erations de Grothendieck et le formalisme des cycles \'evanescents 
dans le monde motivique I et II}, Soci\'et\'e de Math\'ematiques de France, Ast\'erisque, Volume 314 and 315, 2006.

\bibitem{Be2} A.Beilinson, \emph{How to glue perverse sheaves, in K-theory,arithmetic and geometry}, 
Lecture Notes in Mathematics, Vol.1289, Springer-Verlag, Berlin, 1987.

\bibitem{Be} A.Beilinson, \emph{On the derived category of perverse sheaves, in K-theory, arithmetic and geometry}, 
Lecture Notes in Mathematics, Vol.1289, Springer-Verlag, Berlin, 1987, 27-41.

\bibitem{BR} Bhatwadekar, Rao, \emph{On a question of Quillen} AMS 1983

\bibitem{Bondarko} M.V.Bondarko,\emph{Weight structure vs. t-structure ; weight filtrations, 
spectral sequences, and complexes (for motives and in general)}, 
Journal of K-theory, Volume 6, Issue 3,pp.387-504 , 2010

\bibitem{B3} J.Bouali, \emph{On the realization functor of the derived category of mixed motives}, 
arxiv preprint 1706.04545, 2017-arxiv.org.

\bibitem{B4} J.Bouali, \emph{The Hodge realization functor for relative motives of complex algebraic varieties}, 
arxiv preprint 

\bibitem{CG} J.Cirici, F.Guillen, \emph{Homotopy theory of mixed Hodge Complexes}, 
Preprint (2013) arxiv:1304.6236.

\bibitem{C.D} D.C.Cisinski, F.Deglise, \emph{Triangulated categories of mixed motived},
arxiv preprint 0912.2110, 2009-arxiv.org.

\bibitem{C.H} A.Corti, M.Hanamura, \emph{Motivic decomposition and intersection Chow groups I}, Duke math,
J.103 (2000) 459-522.

\bibitem{Coutino} Coutinho, \emph{A primer of algebraic D-modules}, London Mathematical Society,
Cambrige university press, 1995

\bibitem{DP}H.Esnault, E.Viehweg, \emph{Deligne-Beilinson cohomology} Perspect.Math.vol.4,Academic Press, 43-91, 1988

\bibitem{Fangzhou} F.Jin, \emph{Borel-Moore Motivic homology and weight structure on mixed motives},
Math.Z.283 (2016), no.3, 1149-1183.

\bibitem{Hironaka} H.Hironaka, \emph{Resolution of singularities of an algebraic variety over a field of caracteristic zero},
Ann. of Math. 79(1964), 205-326

\bibitem{LvDmod} R.Hotta, K.Takeuchi, T.Tanisaki, \emph{D-Modules, Perverse Sheaves, and Representation Theory},
Birkhauser Verlag, 2008.

\bibitem{huber} R.Huber, \emph{Etale cohomology of rigid analytic varieties and adic spaces}, Springer, 1996.

\bibitem{Kashiwara} M.Kashiwara, \emph{Vanishing cycle sheaves and holonomic systems of differential equations},
Springer Lecture Note, 1016,(1983),134-142.

\bibitem{Chinois} S.Li, X.Pan, \emph{Logaritmic De Rham comparaison theorem for open rigid spaces},
Forum of Mathematics, Sigma 7 E32, 2019.

\bibitem{milnes} J.Milne, \emph{Etale cohomology}, PMS, Volume 33, 1980.

\bibitem{PS} C.Peters, J.Steenbrink, \emph{Mixed Hodge Structures}, Volume 52, Springer, 2008.

\bibitem{Popa} M.Popa, \emph{The V-filtration on D-modules}

\bibitem{Saito} M.Saito, \emph{Mixed Hodge Modules}, 
Proc. Japan Acad. Ser, A.Math.Sci.,Volume 62,Number 9,360-363, 1986.

\bibitem{Sa2} M.Saito, \emph{Module de Hodge polarizable}, Publ.Res.Inst.Math.Sci 24(1988), no.6,849-995(1989).

\bibitem{Scholze} P.Scholze, \emph{p-adic Hodge theory for rigid analytic varieties}, 
Forum of Mathematics, Pi,1,el,(2013).

\end{thebibliography}
\end{document}